\documentclass[openright]{cambridge7A}
\pdfoutput=1
\usepackage{amsmath,amsthm,amsfonts}
\usepackage{newtxtext}
\usepackage[varvw]{newtxmath}
\let\Top\undefined
\usepackage[utf8]{inputenc}
\usepackage[T1]{fontenc}
\usepackage[english]{babel}
\usepackage[hidelinks,pdfusetitle]{hyperref}
\usepackage[only,llbracket,rrbracket,boxslash]{stmaryrd}
\usepackage[all,cmtip]{xy}
\usepackage{xy4}
\usepackage{graphicx}
\usepackage{makeidx}\makeindex
\usepackage[intoc,refpage]{nomencl}\makenomenclature
\usepackage[capitalize,noabbrev]{cleveref}
\usepackage{theorems}
\usepackage{algo}
\usepackage{macros}
\usepackage{longtable}
\usepackage{satex}
\usepackage{cells}
\usepackage{comment}
\usepackage{algorithm2e}
\usepackage{mathrsfs}
\usepackage{mathpartir}
\usepackage{mathtools}
\usepackage[shortlabels]{enumitem}
\usepackage{url}
\usepackage{fixes}

\usetikzlibrary{arrows}
\usetikzlibrary{calc}

\title{Polygraphs:\\from Rewriting to Higher Categories}
\author{Dimitri Ara, Albert Burroni, Yves Guiraud,\\ Philippe Malbos, François Métayer, Samuel Mimram}
\hypersetup{
  pdftitle={Polygraphs: from Rewriting to Higher Categories},
  pdfauthor={Dimitri Ara, Albert Burroni, Yves Guiraud, Philippe Malbos, François Métayer, Samuel Mimram}
}

\begin{document}
\maketitle

\cleardoublepage

\newcommand{\auteur}[3]{
  \noindent
  \begin{minipage}[t]{.45\textwidth}
    \begin{flushright}
      \textsc{#1} \\
      {\footnotesize\texttt{\href{mailto:#2}{#2}}}
    \end{flushright}
  \end{minipage}
  \qquad
  \begin{minipage}[t]{.50\textwidth}
    #3
  \end{minipage}
}

\thispagestyle{empty}

\begin{center}
  \textbf{Abstract}
\end{center}
Polygraphs are a higher-dimensional generalization of the notion of directed graph. Based on those as unifying concept, this monograph on polygraphs revisits the theory of rewriting in the context of strict higher categories, adopting the abstract point of view offered by homotopical algebra. The first half explores the theory of polygraphs in low dimensions and its applications to the computation of the coherence of algebraic structures. It is meant to be progressive, with little requirements on the background of the reader, apart from basic category theory, and is illustrated with algorithmic computations on algebraic structures. The second half introduces and studies the general notion of $n$-polygraph, dealing with the homotopy theory of those. It constructs the folk model structure on the category of strict higher categories and exhibits polygraphs as cofibrant objects. This allows extending to higher dimensional structures the coherence results developed in the first half.

\vfill

\noindent
This is preprint version of a book published at Cambridge University Press available at \url{http://dx.doi.org/10.1017/9781009498968}.

\vfill

\begin{footnotesize}

\bigskip

\auteur{Dimitri Ara}{dimitri.ara@univ-amu.fr}{Aix Marseille Université, CNRS, I2M,\\Marseille, France}

\bigskip

\auteur{Albert Burroni}{burroni@irif.fr}{Université Paris Cité, CNRS, IRIF,\\F-75013 Paris, France}

\bigskip

\auteur{Yves Guiraud}{yves.guiraud@imj-prg.fr}{Université Paris Cité, CNRS, INRIA, IMJ-PRG,\\F-75013 Paris, France}

\bigskip

\auteur{Philippe Malbos}{malbos@math.univ-lyon1.fr}{Université Claude Bernard Lyon 1, CNRS,\\Institut Camille Jordan,\\F-69622 Villeurbanne, France}

\bigskip

\auteur{Francois Métayer}{metayer@irif.fr}{Université Paris Cité, CNRS, IRIF,\\F-75013 Paris, France}

\bigskip

\auteur{Samuel Mimram}{samuel.mimram@polytechnique.edu}{LIX, CNRS, École polytechnique,\\Institut Polytechnique de Paris,\\F-91120 Palaiseau, France}

\end{footnotesize} 

\cleardoublepage

\setcounter{tocdepth}{1}
\tableofcontents

\cleardoublepage

\setcounter{chapter}{-1}
\chapter{Preface}
\subsection*{Presentations of higher categories}
A group is generally defined by characteristic properties of its
elements, such as the group $\Sigma_n$ of all permutations of
the set $\set{1,...,n}$, or the group of all isometries of a cube.
However, to perform actual computations in a group, we
usually pick a subset of {\em generators} among its elements and write
down certain {\em relations} satisfied by these generators, in such a way
that each element of the group is a product of some generators or their
inverses, and each equality between two elements is derivable from the
relations. This leads to a purely syntactic way of defining a group,
called a {\em group presentation} by generators and relations, which
describes the group as a free group over a set of
generators quotiented by some relations.
Presentations are not limited to groups, and they have been adapted to many
other algebraic structures: monoids, categories, Lawvere theories (this
is the subject of universal algebra), commutative (or not) algebras
(where the relations are usually specified by an ideal), or higher
algebras such as operads, product categories,
and linear monoidal categories, to name a few.

The subject of this monograph is the concept of \emph{polygraph},
which is the notion of presentation adapted to higher
categories, and encompasses the previously mentioned settings as
particular cases. Polygraphs were first introduced by
Street~\cite{street1976limits} under the name of \emph{computads} in their
2-dimensional version, in order to study 2-categorical limits, and then
generalized in arbitrary dimension: this was first published
in~\cite{power1991n}, but the generalization was already known and implicitly
used in~\cite{street1987algebra}. The terminology we adopt here comes
from Burroni~\cite{burroni1991higher, burroni1993higher}, who
independently developed the concept in order to study
generalizations of the word problem and provide an equational presentation of
cartesian categories. The name \emph{polygraph} is meant to suggest a
higher-dimensional analogue of oriented graphs. It should be mentioned
that these ideas were developed on both sides, and informally
circulated long before publication. 

\subsection*{The word problem}
Given a presentation of a group, the completeness of a set of
relations means that for any two sequences $a_1,\ldots,a_n$ and $b_1,\ldots,b_m$ of
generators, also called \emph{words}, their respective products
$a_1...a_n$ and $b_1...b_m$ are
equal in the presented group if and only if one can transform the first into the
second by using the relations. This very naturally led, early on, to
the following algorithmic question, known as the \emph{word problem}:
given a finite (or recursive) presentation and two words as above,
can we decide whether they are equal or not in the presented group? 
From a computational theory point of view: can we implement a computer program that automatically determines the equality of words in a presented group?
However, this question predates by many years the invention of computers. It
was first raised as an important one by Dehn in
1911~\cite{dehn1911unendliche}, who subsequently managed to provide an algorithm
for a certain class of groups (the fundamental groups of closed orientable
surfaces of genus greater than or equal to 2)~\cite{dehn1912transformation}. For
some time it was hoped that the problem could be solved for all groups, but it
was in fact shown to be undecidable around 1955 by
Novikov~\cite{novikov1955algorithmic} and Boone~\cite{boone1958word}.

The word problem has also been considered in other settings where
presentations exist. For monoids, it was first studied by Thue in
1914~\cite{thue1914probleme}, leading to the emergence of the notion of \emph{string
rewriting system}, also called semi-Thue systems. It was much later, in 1947, that the word problem was shown to be undecidable by Post~\cite{post1947recursive} and Markov~\cite{markov1947}. In the case of universal algebra, the undecidability of the word problem follows from the undecidability of conversion in combinatory logic~\cite{curry1972combinatory},
which is closely related to the undecidability of $\beta$-conversion in
$\lambda$-calculus~\cite{church1936unsolvable}.

\subsection*{Rewriting theory}

In most cases where the general word problem is undecidable, we can
nevertheless find many specific presentations with algorithmically
decidable equality. For instance, if we ask a small child whether $2+3+4$ is equal to
$2+2+5$, he will progressively compute both sums and observe that the results
are the same:
\[
  \xymatrix@C=3ex@R=3ex{
    2+3+4\ar[d]\ar@{:}[rr]^?&&\ar[d]2+2+5\\
    5+4\ar[dr]&&\ar[dl]4+5\\
    &9
  }
\]
The starting point of rewriting theory is to provide each relation of
a given presentation
with an orientation taking one side of the equality to a simpler
expression on the other side. The resulting structure is called
a \emph{rewriting system} and the oriented rules are called \emph{rewriting
  rules}. Fundamental examples of rewriting systems are  \emph{string} rewriting systems for
presentations of monoids~\cite{book1993string} and \emph{term} rewriting systems
for presentations of Lawvere theories~\cite{bezem2003term,BaaderNipkow98}.

Each rewriting system now comes with a notion of \emph{normal form},
that is, an expression that cannot be further simplified by applying
rewriting rules. This immediately suggests the \emph{normal
  form algorithm}, consisting, given two expressions, in
applying the two following steps:
\begin{enumerate}
\item simplify the two expressions in order to obtain normal forms, and
\item compare the normal forms.
\end{enumerate}

The normal form algorithm only decides the word problem if the
rewriting system satisfies the two following properties.
\begin{enumerate}
\item Each sequence of rewriting rules eventually reaches a normal
  form after a finite number of steps, in which case we call the
  system \emph{terminating}.
\item If an expression can be reduced in two different ways, we can further
  reduce two expressions into a common one, in which case we call the system
  \emph{confluent}.
\end{enumerate}
A rewriting system is \emph{convergent} when it is both terminating and confluent.

Proofs of termination often rely on embedding the reduction order into a partial
well-founded order over the set of expressions. As for confluence, in
cases where the system is already known to be terminating, it suffices to
check a simpler condition called \emph{local confluence}: this is the content of
\emph{Newman's lemma}~\cite{newman1942theories}.
This lemma holds in \emph{abstract rewriting systems} and therefore does not
depend on the particular formalism under consideration (rewriting on
strings, terms, etc.).  We should also mention that convergent rewriting
systems have independently been discovered in the setting of presentations of
(commutative) algebras in the 1960s by Shirshov~\cite{shirshov1962algorithmic}
and Buchberger~\cite{buchberger1965algorithmus}, where they are known as
\emph{Gröbner bases}. In this context, Newman's lemma is known as the
\emph{diamond lemma}~\cite{bergman1978diamond}.

\subsection*{Tietze transformations and completion}
A presentation of a group---or any other type of algebraic
structure---may be thought of as a particular implementation of it,
thus allowing for computations. As  two
 presentations of the same object may have very different computational
 properties, it is worth considering the family of {\em all}
 presentations of a given object. This has first been done in the case of groups by Tietze in
1908~\cite{tietze1908topologischen} (and later on generalized to other
settings): he introduced a family of operations on group presentations,
now called \emph{Tietze transformations}, such that any two
presentations of  the same group can be transformed into each other by
means of those operations. Precisely, Tietze transformations are
sequences built from two elementary steps and their converse, namely:
\begin{enumerate}
\item to add a new generator which is equal to a product of preexisting
  generators, and
\item to add a relation which is derivable from the already present relations.
\end{enumerate}

Tietze transformations will be used to turn a given rewriting system
into a new one presenting the same structure, possibly with better computational properties.
In particular, the idea behind \emph{completion algorithms} is to
add rules (or generators) to a rewriting system in order to turn it into a
confluent one. These ideas were already present at the beginning of string
rewriting systems~\cite{thue1914probleme} and gained much popularity when they
were developed by Buchberger for Gröbner bases~\cite{buchberger1965algorithmus}
(which eventually lead to very efficient
algorithms~\cite{faugere1999new}), by
Knuth and Bendix for term rewriting systems~\cite{knuth1970simple},
and by Nivat
for string rewriting systems~\cite{Nivat73}.
The new rules to be added to the rewriting system are determined by computing
\emph{critical branchings} (also called \emph{S-polynomials} in the linear
settings), which are minimal obstructions to confluence.

\subsection*{The universality problem}
We have seen that, even though the word problem is undecidable in
general, there are many cases where the presentations are convergent and
the word problem can be decided with the normal form algorithm. This leads to the following question, first formulated by Jantzen in the context of
string rewriting systems~\cite{jantzen1982semi, jantzen1985note}, and sometimes
referred to as the \emph{universality problem} for convergent rewriting: given a
finitely presented monoid with a decidable word problem, does it always admit a
presentation by a finite convergent string rewriting system?

A first answer to this question was brought by Kapur and
Narendran~\cite{kapur1985finite} by considering the Artin presentation
of the positive braid group $B_3^+$ (with two generators $a$ and $b$ and one
relation $aba=bab$), which has a decidable word problem: they showed that one
cannot obtain a convergent presentation from it, by adding or removing
relations. However, this does not settle the general question, because we are only
using Tietze transformations of type~(ii) here, and in fact, the same authors also observed
that, by using a transformation of type~(i), one can indeed obtain a convergent
presentation.

\subsection*{Homotopy and homology}
A complete solution to the universality problem was eventually found
by using ideas coming from algebraic topology, which, in a nutshell,
consists in assigning discrete invariants to continuous shapes. Such
ideas already appear in the works of Euler, but their systematic
development starts with Poincaré~\cite{poincare1895}. Among the first
invariants to be considered are the {\em fundamental group} of a
space, which consists of classes of loops up to continuous
deformations (or homotopies) and the sequence of {\em homology
  groups}, providing information about the ``holes'' of various
dimensions in a given space.

The fact that homological invariants are not just numbers but rather bear a
structure of abelian group was first recognized by Emmy
Noether (see~\cite[p.~478]{Brouwer} or~\cite{hilton1988brief}), and
independently by Vietoris, whose paper~\cite{vietoris1927zusa}
contains the first definition of homology groups ever published.
Homology groups are amenable to effective computations by
using classical tools from linear algebra. On the other hand, one may
consider a group or a monoid as a geometric object, by defining the corresponding
{\em classifying space}. Thus, invariants of spaces may be applied to
groups and monoids. These notions have been vastly generalized over
the years and now apply to algebraic structures whose geometrical content is much less
obvious, such as rings or algebras~\cite{hochschild1945cohomology,
  mac1956homologie}.

\subsection*{Squier's homological and homotopical conditions}
The universality problem was answered negatively by Squier in
1987~\cite{squier1987word, squierotto1987word}.
Squier's argument is based on the homology of monoids and relies on
two fundamental observations. First, any presentation of a monoid
determines a sequence of homology groups, but these groups are
independent of the particular presentation we use: they are invariants
of the monoid itself. Second, a finite convergent presentation of a
monoid~$M$ always yields a third homology group $H_3(M)$ of finite
rank. Now, Squier was able to produce an explicit example of a finitely
presented monoid~$M$ with a decidable word problem, whose third homology
group is not of finite rank and therefore does not admit a finite
convergent presentation.

In more precise terms, the homology of a monoid $M$ is computed by building
a \emph{resolution} of the
trivial $\Z M$-module~$\Z$, that is, an exact sequence
\[
  \xymatrix@C=3ex@R=3ex{
    \cdots\ar[r]^-{d_4}&C_3\ar[r]^{d_3}&C_2\ar[r]^{d_2}&C_1\ar[r]^{d_1}&C_0\ar[r]^\varepsilon&\Z\ar[r]&0
  }
\]
of projective $\Z M$-modules. Tensoring the above sequence by $\Z$ over
$\Z M$ gives a chain complex of abelian groups, which is not exact anymore in the general case;
its homology groups only depend on $M$ and not on the
particular resolution we chose.
Now, if we start with a convergent presentation of $M$, we
obtain a partial resolution of the form
\[
  \xymatrix@C=3ex@R=3ex{
    \cdots\ar[r]^-{d_4}&\Z M[P_3]\ar[r]^{d_3}&\Z M[P_2]\ar[r]^{d_2}&\Z M[P_1]\ar[r]^{d_1}&\Z M[P_0]\ar[r]^-\varepsilon&\Z\ar[r]&0,
  }
\]
where
$P_0$ is the set with one element, $P_1$ is the set of generators of the
presentation, $P_2$ the set of relations, and $P_3$ the set of critical
branchings. For a finite rewriting system, the set $P_3$ of critical branchings
is always finite. Hence the abelian group $\Z M[P_3]\otimes_{\Z M}\Z$
obtained by tensoring by $\Z$ over $\Z M$ has finite
rank and so has the third homology group $H_3(M)$.
Note that the finiteness condition for $H_3$ comes here from the
strictly stronger statement that $M$ is of type \emph{left-$\FP_3$},
which  means that it admits a partial resolution of
length $3$ by finitely generated projective left $\Z M$-modules.

The above algebraic constructions can be interpreted in geometric
terms~\cite{lafont1991church}, at least in the case where the monoid
$M$ is a group. The free resolution $\pair{C_i}{d_i}$ comes from a
cellular decomposition of a {\em contractible} space $X$ on which $M$
acts freely and transitively. Tensoring by $\Z$ over $\Z M$
amounts to quotient $X$ by the action of $M$, thus obtaining a
cellular decomposition of the classifying space of $M$
itself. Starting with a presentation of $M$, the cellular complex we
get is built dimensionwise: in dimension $0$ there is a unique point,
corresponding to the element of $\star\in P_0$. In dimension $1$, each
generator of $P_1$ gives a loop on $\star$. In dimension $2$, each
relation in $P_2$ gives a disk attached to the $1$-dimensional paths
determined by the products of the generators involved.
For instance, with the Artin presentation
of~$B_3^+$ recalled above, we would obtain the space on the left, together with
a disk attached between the paths corresponding to the words $aba$ and $bab$, as
pictured on the right:
\[
  \xymatrix@C=3ex@R=3ex{
    \ar@{-}@(ul,dl)_a\scriptstyle\star\ar@{-}@(ur,dr)^b
  }
  \qquad\qquad\qquad\qquad\qquad
  \begin{tikzpicture}[scale=.5,baseline=(b.base)]
    \coordinate (b) at (0,0);
    \filldraw[color=lightgray,fill=lightgray] (0,0) circle (1);
    \draw (0:1) node {$\scriptstyle\star$};
    \draw (60:1) node {$\scriptstyle\star$};
    \draw (120:1) node {$\scriptstyle\star$};
    \draw (180:1) node {$\scriptstyle\star$};
    \draw (240:1) node {$\scriptstyle\star$};
    \draw (300:1) node {$\scriptstyle\star$};
    \draw (30:1.3) node {$\scriptstyle a$};
    \draw (90:1.3) node {$\scriptstyle b$};
    \draw (150:1.3) node {$\scriptstyle a$};
    \draw (210:1.3) node {$\scriptstyle b$};
    \draw (270:1.3) node {$\scriptstyle a$};
    \draw (330:1.3) node {$\scriptstyle b$};
  \end{tikzpicture}
\]
There may be different ways to fill the gap between two expressions
representing the same element of $M$ by using disks as
above. Geometrically, this corresponds to $3$-dimensional holes in the
corresponding $2$-dimensional complex. These holes have to be filled
by appropriate $3$-dimensional cells, and it turns out that, in the
case of a convergent presentation, a set of $3$-cells coming from
the critical branchings is sufficient for that. Of course, the construction
has to be pursued in higher dimensions in order to obtain the correct topology.

About the same time when Squier studied monoid resolutions,
other authors developed similar ideas:
Anick constructed a resolution of algebras in order
to study Koszulness properties~\cite{anick1986homology}, Brown managed to use
discrete Morse theory in order to ``reduce'' the standard resolution to a small
one~\cite{brown1987finiteness, brown1992geometry}, and Kobayashi extended
Squier's partial resolution into a full one~\cite{kobayashi1990complete}.
Finally, we should mention that Squier subsequently provided a ``homotopical''
variant of his condition, which was published after his
death~\cite{squier1994finiteness}, see also~\cite{lafont1995new}. A finite
presentation is said to be of \emph{finite derivation type} when the \emph{full
  congruence}, which identifies any two witnesses of equality between two given
words, is finitely generated. It can be shown that this property is an invariant
of the monoid and is implied by having a finite convergent
presentation. Moreover, one can construct an explicit example of a monoid that
has a finite presentation with decidable word problem but does not have finite
derivation type.

\subsection*{Polygraphs as higher-dimensional presentations}
A monoid is a category with only one object, and it is thus natural to extend
the notion of presentation, and associated theorems and techniques, from monoids
to small categories. Precisely, a small category $C$ will be presented
by a set $P_0=C_0$ of objects, a set $P_1$ of $1$-generators, that is,
a subset $P_1\subset C_1$ of morphisms generating all morphisms in $C$
by composition, together with a set $P_2$ of relations between certain
pairs of composites of $1$-generators. This pattern generalizes to
higher categories, resulting in the notion of
$n$\nbd-\emph{poly\-graph}~\cite{burroni1993higher,street1976limits},
which consists of the following data:
\begin{itemize}
\item for each $i\in\set{0,...,n}$,  a set $P_i$, freely generating a
  set $P_i^*$ of $i$-dimensional cells, and
\item for each  $i\in\set{1,...,n}$, a pair of maps
  associating to each $i$-generator its source and target in $P_{i-1}^*$. 
\end{itemize}
Given an $(n{+}1)$-polygraph $P$, the source and target maps defined on
$P_{n+1}$ generate an equivalence relation on $P_n^*$, whose quotient
set $C_n$ is the set of $n$-morphisms of the $n$-category
$C$ {\em presented} by $P$.   

Right from the beginning, a polygraph was thought of as a
``higher-di\-men\-sio\-nal rewriting system''~\cite{burroni1993higher,
  eilenberg1986rewrite, street1995higher, street1996categorical} and the
associated theory of rewriting was subsequently developed, in particular by some
of the authors of this book~\cite{guiraud2006termination,guiraud2009higher,GuiraudMalbos10smf,GuiraudMalbos18,lafont2003towards,mimram2014towards}. It
turns out that many of the classical theorems go through but, starting at
dimension $n=3$, there is a major difference with the classical setting: a
finite rewriting system might give rise to an infinite number of critical
pairs, thus preventing easy generalizations of Squier-type theorems.

The fact that polygraphs are higher-dimensional rewriting systems should be
taken in a very strong sense here: they are about rewriting rewriting paths.
Namely, an $(n+1)$-polygraph consists of an $n$-polygraph together with
$(n+1)$\nbd-dimensional rewriting rules, which specify how to rewrite rewriting
paths in the underlying $n$-polygraph.
For instance, consider the abstract rewriting system on the left, which is a
$1$-polygraph:
\[
  \vcenter{
    \xymatrix@C=3ex@R=3ex{
      \ar@(ul,dl)_a\scriptstyle\star\ar@(ur,dr)^b
    }
  }
  \qquad\qquad\qquad\qquad
  \vcenter{
    \xymatrix@C=3ex@R=1ex{
      &\scriptstyle\star\ar[r]^b&\scriptstyle\star\ar@/^/[dr]^a\\
      \scriptstyle\star\ar@/^/[ur]^a\ar@/_/[dr]_b\ar@{}[rrr]|{\Longdownarrow}&&&\scriptstyle\star\\
      &\scriptstyle\star\ar[r]_a&\scriptstyle\star\ar@/_/[ur]_b
    }
  }
\]
Here, we have one object $\star$ and two rewriting rules $a$ and $b$ rewriting
$\star$ to itself. We can obtain a $2$-polygraph by adjoining a $2$-dimensional
rewriting rule which rewrites the path $aba$ into the path $bab$, as pictured on
the right, thus providing a polygraphic presentation of the monoid $B_3^+$: in
this way, we can see that string rewriting is secretly rewriting rewriting paths
in abstract rewriting systems! Similarly, in the next dimension, we can see that
term rewriting, and more generally rewriting of diagrams, is an instance of
rewriting rewriting paths for string rewriting systems.

Several particular generalizations of the notion of polygraph have also been
investigated (for $(n,p)$-categories, linear
higher-categories~\cite{GuiraudHoffbeckMalbos19}, cartesian higher-categories,
etc.). Most of them are particular instances of the notion of
\emph{$T$\nbd-poly\-graph} introduced by Batanin~\cite{batanin1998computads}
in order to define polygraphs adapted to weak higher categories: he defines
a notion of polygraph parametrized by a globular monad~$T$, whose various
instantiations allow recovering the previously mentioned variants of
polygraphs.

\subsection*{Coherence}
Given a presentation of an $n$-category $C$ by an $(n{+}1)$-poly\-graph~$P$,
the elements of $P_{n+1}^*$ witness the equalities between
$n$-morphisms of  $C$. Now, different $(n{+}1)$-cells may witness the
same equality: this defines a congruence on~$P_{n+1}^*$. If we extend
$P$ by a set $P_{n+2}$ generating this congruence, we get a {\em
  coherent presentation} of $M$ by an $(n{+}2)$-polygraph. For instance, from a convergent presentation of a
monoid, we build a coherent presentation by a 
$3$\nbd-polygraph $P$ in which the generators in $P_3$ correspond to the critical branchings.

Coherent presentations provide the ``right'' generalization (in a homotopical sense detailed below) of a structure in higher
dimensions. For instance, the theory of monoids can be described by a
$3$-polygraph. If we extend this polygraph into a coherent one, we obtain a theory
corresponding to pseudo-monoids, which is the expected notion of monoid in a
monoidal $2$-category.

\subsection*{Resolutions with polygraphs}
The construction of a coherent presentation of an $n$-category may be infinitely
pursued in higher dimensions by introducing, for each $m\geq n$, a
set of  $(m{+}1)$-cells generating all desired congruences between
$m$\nbd-cells. More generally, starting with any $\omega$-category
$C$, we may build a polygraph $P$ together with an $\omega$-functor
$p:\free P\to C$ satisfying the following properties.
\begin{enumerate}
\item For each dimension $n\geq 0$, $P_n$ generates $C_n$, in the sense that
  $p_n:\free P_n\to C_n$ is surjective.
\item For any two parallel cells $x$, $y$ in $\free P_n$ such that
  $p_nx=p_ny=u\in C_n$, there is an $(n{+}1)$-cell $z\in \free P_{n+1}$ with
  source $x$ and target $y$  such that $p_{n+1}z=\unit{u}$.
\end{enumerate}
We call such a map $p:\free P\to C$ a {\em polygraphic resolution} of $C$
by the polygraph~$P$~\cite{metayer2003resolutions}.  It turns out that
two polygraphic resolutions of the same $\omega$\nbd-category are
equivalent up to a suitable notion of homotopy. Polygraphic
resolutions are cofibrant replacements in the
``canonical''  model structure on the category of small
$\omega$-categories, as introduced in~\cite{LMW}. As a consequence,
there is a well-defined notion of homology for $\omega$-categories: to
each polygraphic resolution of $C$ by $P$ corresponds a chain complex $\pair{\Z
  P_n}{\partial_n}_{n\geq 0}$ of abelian groups, whose homology only
depends on $C$. This illustrates the general principle according to
which many constructions are better behaved when performed on {\em
  free} objects. In the particular case where $C$ is a monoid $M$  seen
as an $\omega$-category, this
homology coincides with the one computed via free resolutions of $\Z$
by $\Z
M$\nbd-mo\-dules~\cite{lafont2009polygraphic,GuettaHomologyCat}. These
constructions transfer to strict $\omega$-groupoids, or more generally
$\pair{\omega}{n}$-categories~\cite{AraMetWGrp}, where all cells of
dimension strictly above~$n$ are invertible, yielding appropriate notions of polygraphic
resolutions.

\subsection*{Structure of the book(s)}
This book informally splits into two books, which, however strongly
connected, can be read separately, according to one's taste and objectives.

\bigskip\noindent\emph{Low-dimensional book.}
The first book explores the theory of polygraphs in low dimensions and its
applications. It is meant to be very progressive, with little requirements on
the background of the reader, apart from basic category theory,
and is illustrated with algorithmic computations on algebraic structures.
We namely study polygraphs in dimension~1 (\cref{chap:1pol}), in dimension~2
(\cref{chap:2pol,chap:2op,chap:2rewr,chap:2tietze,chap:2linear,chap:2coh,chap:2fdt,chap:2homology}),
and in dimension~3 (\cref{chap:3pol,chap:3term,chap:3coh,chap:trs}).

In all the cases, we introduce the notion of polygraph as well as the associated
notions of generated and presented categories
(\cref{sec:1-pol-cat,sec:pres-set,chap:2pol,sec:3pol}), develop the theory of
rewriting (\cref{sec:ars,chap:2rewr,sec:3pol-rewr}), Tietze transformations and
completion procedures (\cref{sec:pres-set,chap:2tietze}), termination techniques
(\cref{sec:ars,sec:2-red-order,chap:3term}) and coherent presentations
(\cref{chap:2coh,chap:3coh}), with the last notion requiring the introduction of
higher-dimensional notions of polygraphs. We also present the homotopical and
homological invariants (\cref{chap:2fdt,chap:2homology}) they allow to compute,
and introduce variants of the notion of polygraph, namely linear
(\cref{chap:2linear}) and cartesian polygraphs (\cref{chap:trs}).

\bigskip\noindent\emph{Higher-dimensional book.}
The second book goes at a faster pace and supposes that the reader is
familiar with category theory. Moreover, acquaintance with strict higher
categories, as well as the notion of model category, can be helpful, even
though these notions are recalled.
The beginning of this book introduces and studies the general notion of
$n$-polygraph (\cref{chap:n-cat,chap:polygraphs,chap:pol-ex,chap:pol-prop,chap:gen-pol}).
The remainder of the book deals with the homotopy theory of these polygraphs.
We construct the ``folk'' model structure on the category of
$\omega$-categories (\cref{chap:resolutions,chap:w-eq,chap:folk}), in which
polygraphs are precisely the cofibrant objects. This model structure is used
to define a homology theory for $\omega$-categories as a derived functor
(\cref{chap:homology}). Finally, we study the variant of $(\omega,
1)$-polygraphs (\cref{chap:anick}), which allows to formulate a
higher-dimensional generalization of the coherence results developed in the
``low-dimensional book''.

\bigskip\noindent\emph{Appendix.}
The book is followed by a number of chapters containing additional
material. Some of them perform a review of classical---or not---examples
of polygraphs, in order to illustrate their diversity and applications:
2-polygraphs (\cref{chap:2ex}), coherent 2-polygraphs
(\cref{chap:SomeCoherentPresentations}), and 3-polygraphs (\cref{chap:3ex}). Some
other chapters recall elements of classical topics used throughout the
second part of the book: free $n$-categories
(\cref{chap:syntactic_descr}), homology (\cref{Chapter:ComplexesAndHomology}),
locally presentable categories (\cref{chap:loc_pres}), and model categories
(\cref{chap:model-cat}).


\subsection*{Acknowledgments}
The writing of this book started in 2015, under the code name of \emph{polybook}, and the present exposition has
benefited much from interactions with other researchers and friends. We would
like to warmly thank the following people who, voluntarily or not, contributed
to its elaboration: Mathieu Anel, Thibaut Benjamin, Pierre Cagne, Cameron Calk, Cyrille
Chenavier, Pierre-Louis Curien, Benjamin Dupont, Eric Finster, Simon Forest,
Jonas Frey, Andrea Gagna, Léonard Guetta, Amar Hadzihasanovic, Simon Henry, Hugo Herbelin,
Eric Hoffbeck, Cédric Ho Thanh, Jan Willem Klop, Yves Lafont, Chaitanya Leena-Subramaniam, Maxime Lucas, Georges Maltsiniotis,
Paul-André Melliès, Jovana Obradovic, Vincent van Oostrom, Viktorya Ozornova, Jacques
Penon, and all participants of the Higher Category Seminar at IRIF.

String diagrams were first typeset using catex~\cite{catex}, and then
satex~\cite{satex}.

The early stages of the book also benefited from fundings provided by
the projects Focal IDEX and CATHRE ANR.

Our working sessions took place over the years in various locations, and we especially
acknowledge the hospitality of the IHP in Paris, the CIRM in
Marseille, and the Bois-Maréchal in the middle of the Mâconnais vineyard.

\part{Fundamentals of Rewriting}

\chapter{Abstract Rewriting and One-Dimensional Polygraphs}
\chaptermark{One-dimensional polygraphs}
\label{chap:lowdim}
\label{chap:1pol}
We begin by discussing 1-polygraphs, which are simply directed graphs, thought
of here as abstract rewriting systems: they consist of vertices, which represent
the objects of interest, and arrows, which indicate that we can rewrite one
object into another. After formally introducing those in \secr{1-pol-cat}, we
will see in \secr{pres-set} that they provide a notion of \emph{presentation}
for sets, by generators and relations. Of course,
presentations of sets are of little interest in themselves, but they are merely
used here as
a gentle introduction to some of the main concepts discussed in this work:
in particular, we introduce the notion of Tietze
transformations which generate the equivalence between
two presentations of the same set. In
this context, an important question consists in deciding when two objects are
equivalent, \ie represent the same element of the presented set. In order to
address it, we develop the theory of abstract rewriting systems in
\secr{ars}. Most notably, we show that when the rewriting system
satisfies the two
properties of \emph{termination} and \emph{confluence}, equivalence classes of
objects admit a unique canonical representative, the \emph{normal form}, and
equivalence of objects can thus be decided by comparing the associated normal
forms. Finally, in \secr{decreasing-diag}, we detail the more advanced method of
decreasing diagrams, which can be used  to show confluence in the absence of
termination.

\section{The Category of 1-Polygraphs}
\label{sec:1-pol-cat}
\index{polygraph!0-}
\index{0-polygraph}
A \emph{0-polygraph} is simply another name for a set. Since there is not much
to do with those, we move on to 1-polygraphs.

\begin{definition}
  \index{1-polygraph}
  \index{polygraph!1-}
  \index{generator}
  A \emph{1-polygraph}~$P$ consists of a 0-polygraph $P_0$, whose elements
  are called \emph{0-generators}, together with a set~$P_1$ of
  \emph{1-generators} and two functions
  $
    \sce0^P,\tge0^P
    :
    P_1
    \to
    P_0
  $
  respectively associating to each 1-generator its \emph{source} and
  \emph{target} 0-cell. We often write
  $
    \pres{P_0}{P_1}
  $
  for such a polygraph and
  $
    a
    :
    x\rto y
  $
  for a 1-generator $a$ in $P_1$ such that $\sce0^P(a)=x$ and $\tge0^P(a)=y$.
  A 1-polygraph~$P$ is \emph{finite}\index{finite!1-polygraph} when both $P_0$
  and $P_1$ are.
\end{definition}

\index{graph}
\index{1-graph}
The notion of 1-polygraph is simply another name for the notion of
\emph{graph}, by which we always mean a directed multigraph, which we sometimes also call a \emph{1-graph}. Indeed, a polygraph~$P$ as
above is a graph with $P_0$ as set of vertices $P_1$ as set of edges, an edge
$a\in P_1$ having~$\src0^P(a)$ as source and $\tgt0^P(a)$ as
target. Thus, any terminology pertaining to oriented graphs, such as
the notion of {\em path}, immediately applies to 1-polygraphs.

\begin{example}
  \label{ex:graph}
  The directed graph
  \begin{equation}
    \label{eq:ex-graph}
    \vxym{
      x\ar@<.5ex>[r]^a\ar@<-.5ex>[r]_b&y\ar@(ul,ur)^c&z
    }
  \end{equation}
  can be encoded as the 1-polygraph~$P$ with $P_0=\set{x,y,z}$,
  $P_1=\set{a,b,c}$ and
  \begin{align*}
    \sce0(a)&=\sce0(b)=x\pbox,
    &
    \tge0(a)&=\tge0(b)=y\pbox,
    &
    \sce0(c)&=\tge0(c)=y\pbox,
  \end{align*}
  which can be more concisely denoted as
  \[
  P=
  \pres{x,y,z}{a:x\rto y,b:x\rto y,c:y\rto y}
  \pbox.
  \]
\end{example}

\subsection{The category of 1-polygraphs}
A \emph{morphism}
$
f: P\to Q
$
between $1$\nbd-polygraphs~$P$ and~$Q$ consists of a pair of functions
$f_0:P_0\to Q_0$
and
$f_1:P_1\to Q_1$
respectively sending the 0- and 1-cells of~$P$ to those of~$Q$ and
preserving sources and targets:
\begin{align*}
  \sce0^Q\circ f_1&=f_0\circ\sce0^P,
  &
  \tge0^Q\circ f_1&=f_0\circ\tge0^P
  \pbox.
\end{align*}
\nomenclature[Pol1]{$\nPol 1$}{category of 1-polygraphs}%
We write $\nPol 1$ for the category of 1-polygraphs and their
morphisms. Again, this is simply another name for the usual category of directed
graphs and their morphisms.

\section{Presenting Sets}
\label{sec:pres-set}
A 1-polygraph~$P$ can be seen as a {\em presentation} of a set~$X$, in the
following sense. Each element $x$ of $P_0$ denotes an element
$\prescat x$ of $X$, in
such a way that each element of $X$ has at least one ``name'' in
$P_0$, and each element $a:x\to y$ in $P_1$ represents the renaming of
$\prescat x$ by $\prescat y$. The elements of~$P_0$ and~$P_1$ are often respectively called \emph{generators}\index{generator} and
\emph{relations}\index{relation}.

\subsection{$P$-congruence}
\label{sec:P-cong}
\index{congruence!on a set}
The \emph{$P$-congruence} $\approx^P$ associated with a 1-polygraph~$P$ is the
smallest equivalence relation on~$P_0$ such that $x\approx^P y$ for every
1-generator $a:x\to y$ in~$P_1$.

\subsection{The presented set}
\index{presentation!of a set}
The \emph{set $\prescat{P}$ presented} by a 1-polygraph $P$ is the set
$P_0/{\approx^P}$
obtained by quotienting $P_0$ by the $P$-congruence $\approx^P$, what we usually
simply write $P_0/P_1$.
More generally, a set~$X$ is \emph{presented} by a 1-polygraph~$P$ when $X$ is
isomorphic to~$\prescat{P}$, and in this case $P$ is called a
\emph{presentation} of~$X$. Geometrically speaking, $X$ amounts to the
set of connected components of the graph $P$.

\begin{example}
  In \exr{graph}, the relation $\approx^P$ identifies $x$ and $y$, and the
  presented set is the set with two elements, corresponding to the equivalence
  classes $\set{x,y}$ and $\set{z}$.
\end{example}

\noindent
More abstractly, the set presented by a 1-polygraph~$P$ can be characterized by
the following universal property:

\begin{lemma}
  For any set~$X$ and function $f:P_0\to X$ such that $f(x)=f(y)$ for every
  1-generator $a:x\to y$ in~$P_1$, there exists a unique function
  $\prescat{f}:\prescat{P}\to X$ such that $\prescat f\circ q=f$
  \[
  \vxym{
    P_0\ar[d]_-q\ar[r]^f&X\\
    \prescat{P}\ar@{.>}[ur]_{\prescat f}
  }
  \]
  where $q:P_0\to\prescat{P}$ is the function sending an element to its
  equivalence class.
\end{lemma}

\subsection{Tietze transformations}
\label{sec:1pol-tietze}
At this point, a natural question to ask is: when do two polygraphs present the
same set?
For instance, the set with two elements can also be presented by the polygraph
\begin{equation}
  \label{eq:ex-graph'}
  \vxym{
    x\ar[r]^d&x'\ar[r]^e&y&z
  }
\end{equation}
which looks quite different from~\eqref{eq:ex-graph}, and it is not obvious how the two are related.
This question was first studied by Tietze for presentations of
groups~\cite{tietze1908topologischen}, as we shall see in
\chapr{2-tietze}, but similar results already hold for plain sets as we now explain.

\index{Tietze!transformation!of 1-polygraph}
\index{transformation!Tietze}
We call \emph{elementary Tietze transformations} 
the following operations transforming a $1$\nbd-po\-ly\-graph~$P$ into a
1-polygraph~$Q$:
\begin{description}
\item[\tgen] \emph{adding a definable generator}: given $x\in P_0$,
  $y\not\in P_0$, and $a\not\in P_1$, we define
  \[
  Q=\pres{P_0,y}{P_1,a:x\rto y}\text,
  \]
\item[\trel] \emph{adding a derivable relation}: given $x,y\in P_0$ and
  $a\not\in P_1$ such that $x\approx^P y$, we define
  \[
  Q=\pres{P_0}{P_1,a:x\rto y}\text.
  \]
\end{description}
A \emph{Tietze transformation} from~$P$ to~$Q$ is a zigzag of elementary Tietze
transformations, \ie a finite sequence of polygraphs~$(P_i)_{0\leq i\leq n}$
with $P_0=P$ and~$P_n=Q$, together with, for each index~$0\leq i<n$, an
elementary Tietze transformation either from~$P_i$ to~$P_{i+1}$ or from
$P_{i+1}$ to~$P_i$.
The \emph{Tietze equivalence}\index{Tietze!equivalence!of 1-polygraphs}\index{equivalence!Tietze} is the smallest equivalence
relation on 1-polygraphs, identifying any two polygraphs related by an
elementary Tietze transformation and closed by isomorphism; otherwise said, two polygraphs are Tietze
equivalent when there exists a Tietze transformation between them, up
to isomorphism.

\begin{lemma}
  \label{lem:tietze-equiv-1}
  Two Tietze equivalent 1-polygraphs present isomorphic sets.
\end{lemma}
\begin{proof}
  By induction on the length of Tietze transformations, it is enough to show that two polygraphs~$P$ and~$Q$ related
  by an elementary Tietze transformation present the same set. Using the same
  notations as above, in the case of the transformation \tgen, we have
  \[
    \pcat{Q}=(P_0\uplus\set{y})/{\approx^Q}=((P_0\uplus\set{y})/(x\approx y))/{\approx^P}=P_0/{\approx^P}=\pcat{P},
  \]
  where $x\approx y$ denotes the smallest equivalence relation identifying~$x$
  and~$y$. In the case of the transformation \trel, the relations generated
  by~$P_1$ and~$Q_1$ are the same and we have
  \[
    \pcat{Q}=Q_0/{\approx^Q}=P_0/{\approx^P}=\pcat{P}\pbox.\qedhere
  \]
\end{proof}

\noindent
We will see in \thmr{tietze-equiv-1} that the converse also holds: these
operations exactly axiomatize when two finite 1-polygraphs are presenting the
same set.

\begin{example}
  Using the above lemma, one can deduce that the two
  polygraphs~\eqref{eq:ex-graph} and~\eqref{eq:ex-graph'} present the same set,
  by building a series of Tietze transformations relating them:
  \[
  \begin{array}{r@{\qquad}c@{\qquad}r@{\qquad}c@{\qquad}r@{\qquad}c@{\qquad}r}
    \xymatrix@C=2ex@R=2ex{
    x\ar@<.5ex>[r]^a\ar@<-.5ex>[r]_b&y\ar@(ul,ur)^c&z
    }
    &\overset{\trel}\leftsquigarrow&
    \xymatrix@C=2ex@R=2ex{
      x\ar@<.5ex>[r]^a\ar@<-.5ex>[r]_b&y&z
    }
    &\overset{\trel}\leftsquigarrow&
    \xymatrix@C=2ex@R=2ex{
      x\ar[r]^a&y&z
    }\\
    &\overset{\tgen}\rightsquigarrow&
    \xymatrix@C=2ex@R=2ex{
      x'&\ar[l]_-dx\ar[r]^a&y&z
    }
    &\overset{\trel}\rightsquigarrow&
    \xymatrix@C=2ex@R=2ex{
      x'\ar@/_/[rr]_e&\ar[l]_-dx\ar[r]^a&y&z
    }\\
    &\overset{\trel}\leftsquigarrow&
    \xymatrix@C=2ex@R=2ex{
      x'\ar@/_/[rr]_e&\ar[l]_-dx&y&z\pbox{.}
    }
  \end{array}
  \]
  In the first step, $y\approx y$ can be shown without resorting to the relation
  $c:y\to y$ (this is because, by definition, $\approx$ is an equivalence
  relation), and therefore the relation $h$ can be removed using the Tietze
  transformation $\trel$ backward. Other steps can be justified similarly.
  Of course, in this case, it is very easy to compute the sets presented by the
  two polygraphs~\eqref{eq:ex-graph} and~\eqref{eq:ex-graph'} and to see that they
  are isomorphic (both have two elements), but it will no longer be the
  case, when generalizing to higher dimensions.
\end{example}

\subsection{Backward Tietze transformations}
\index{Tietze!transformation!backward}
A Tietze transformation is a zigzag of elementary Tietze transformations. It can
alternatively be seen as a sequence of elementary Tietze transformations or the
following transformations, that we call \emph{backward elementary Tietze
  transformations}, corresponding to using an elementary Tietze transformation
in the ``backward direction'':
\begin{description}
\item[\trgen] \emph{removing a definable generator}: given a polygraph~$P$ of
  the form
  \[
  P=\pres{P_0',x}{P_1',a:x\rto y},
  \]
  where $x$ does not occur in any relation of $P_1'$, we define
  \[
  Q=\pres{P_0'}{P_1'},
  \]
\item[\trrel] \emph{removing a derivable relation}: given a polygraph~$P$ of the
  form
  \[
  P=\pres{P_0}{P_1',a:x\rto y},
  \]
  we define
  \[
  Q=\pres{P_0}{P_1'}
  \]
  whenever $x\approx^Q y$.
\end{description}

\begin{remark}
  Given an elementary Tietze transformation from~$P$ to~$Q$, there is an obvious
  inclusion of~$P$ into~$Q$ that induces a morphism of
  1\nbd-poly\-graphs~$P\to Q$. However, for a backward elementary Tietze
  transformation from~$P$ to~$Q$ there is no canonical morphism~$P\to Q$. For
  instance, consider the transformation
  \[
    \xymatrix@C=2ex@R=2ex{
      x\ar@<.5ex>[r]^a\ar@<-.5ex>[r]_b&y\ar@(ul,ur)^c&z
    }
    \qquad\overset{\trrel}\rightsquigarrow\qquad
    \xymatrix@C=2ex@R=2ex{
      x\ar@<.5ex>[r]^a\ar@<-.5ex>[r]_b&y&z\pbox.
    }
  \]
  The only reasonable choice would be to send the 1-generator~$c:y\to y$ to an
  identity on~$y$, which is not possible with a morphism of 1-polygraph (those
  send 1-generators to 1-generators). This is one of the reasons why we take the
  elementary Tietze transformations (as opposed to the backward ones) as more
  primitive.
\end{remark}

\subsection{Minimal presentations}
\label{sec:1-minimal-pres}
\index{presentation!minimal}
It can be noted that Tietze transformations consisting only of elementary transformations $\tgen$ and $\trel$ make the
presentations larger (in terms of number of generators and relations), whereas those consisting only of
$\trgen$ and $\trrel$ make them smaller. We thus sometimes respectively call \emph{Tietze
  expansions}\index{Tietze!expansion} and \emph{Tietze reductions}\index{Tietze!reduction} these two families of Tietze
transformations and  say that a polygraph~$P$ \emph{Tietze expands} (\resp
\emph{Tietze reduces}) to a polygraph~$Q$ if $Q$ can be obtained from~$P$ by
applying a series of Tietze expansions (\resp Tietze reductions).
One may wonder if, by applying only the second kind of transformations, we
eventually always reach a {\em minimal} presentation \wrt both
generators and relations, and whether two such minimal presentations are necessarily
isomorphic. We will see that it is indeed the case for finite polygraphs. First, note that a
1-polygraph~$P$ without relations (\ie $P_1=\emptyset$)
is always minimal.

\begin{lemma}
  \label{lem:1-disc}
  Any finite 1-polygraph~$P$ Tietze reduces to a
  polygraph isomorphic to~$\smallpres{\prescat{P}}{}$.
\end{lemma}
\begin{proof}
  By induction on the cardinal of~$P_1$, we show that we can remove a
  $1$\nbd-generator using Tietze transformations, unless $P_1$ is empty. Suppose
  that~$P$ contains a non-directed cycle, \ie a non-empty non-directed path
  from a $0$\nbd-gene\-rator~$x$ to itself. We can assume that this path does not use
  the same edge twice; otherwise, we can choose a smaller cycle. Given a
  1-generator $a:x\to y$ occurring in this cycle, there exists a non-directed
  path from~$x$ to~$y$ that is not using $a$. Therefore, we can apply a Tietze
  transformation \trrel{} to remove~$a$. Otherwise, there is no cycle, and
  consider a maximal non-directed path in~$P$. Since~$P$ is finite and acyclic,
  this path will end by a 1-generator $a:x\to y$
  such that either~$x$ or~$y$ is incident to no other edge. Therefore, we can
  use a Tietze transformation \trgen{} to remove~$x$ or~$y$, along with~$a$.
\end{proof}

\noindent
In the case of {\em finite} 1-polygraphs, the above lemma implies
the converse of \lemr{tietze-equiv-1}:

\begin{theorem}
  \label{thm:tietze-equiv-1}
  Two finite 1-polygraphs present isomorphic sets if and only if they are Tietze
  equivalent.
\end{theorem}
\begin{proof}
  Suppose given two polygraphs~$P$ and~$Q$ such that
  $\prescat{P}\isoto\prescat{Q}$. By the previous lemma, $P$ is Tietze equivalent to
  $\smallpres{\prescat{P}}{}$, and similarly~$Q$ is Tietze equivalent
  to~$\smallpres{\prescat{P}}{}$. Finally, the presentations
  $\smallpres{\prescat{P}}{}$ and $\smallpres{\prescat{Q}}{}$ are easily seen to
  be Tietze equivalent because $\prescat{P}$ and $\prescat{Q}$ are isomorphic.
\end{proof}

\begin{remark}
  \label{rem:trans-tietze}
  Note that, given the above definition of Tietze transformations, the
  previous theorem does not generalize to infinite presentations. For instance, the 1-polygraphs
  $\pres{x}{}$ and $\pres{x_i}{a_i:x_i\rto x_0}_{i\in\N}$
  both present the set with one element but are not Tietze equivalent since we
  can only add or remove a finite number of relations using Tietze equivalences
  (the notation on the right means that $i$ ranges over~$\N$ both in
  generators~$x_i$ and relations~$a_i$). In order to overcome this
  counter-example, one might be naively tempted to allow infinite sequences of
  Tietze transformations between 1\nbd-polygraphs, but this does not preserve
  presented sets. For instance, consider the 1-polygraph
  \[
    \pres{x_i,y}{a_i:x_{i+1}\rto x_i,b_i:x_i\rto y}_{i\in\N},
  \]
  \ie the graph
  \[
    \xymatrix@C=4ex@R=4ex{
      x_0\ar[d]|{b_0}&\ar[l]_{a_0}x_1\ar[dl]|{b_1}&\ar[l]_{a_1}x_2\ar[dll]|{b_2}&\ar[l]_{a_2}\cdots\ar[dlll]|{b_3}\\
      y
    }
  \]
  presenting the set with one element. Using Tietze transformations, any finite
  number of relations $b_i$ can be removed from the polygraph, since they are
  derivable. However, if we remove all of them the resulting polygraph presents
  the set with two elements.

  In order to account for infinite presentations, the notion of Tietze
  equivalence has to be generalized as follows. Firstly, we say that a
  1-polygraph~$P$ Tietze expands to~$Q$ if there is a transfinite sequence of
  elementary Tietze expansions from~$P$ to~$Q$; secondly, we define Tietze
  equivalence as the smallest equivalence relation containing Tietze
  expansions. Two (not necessarily finite) 1-polygraphs are Tietze equivalent
  in this sense if and only if they present isomorphic sets. We do not
  dwell further on infinite polygraphs, because we are mostly interested in finite polygraphs in this book; details can be found in~\cite{henry2022tietze}.
\end{remark}

We will see that \lemr{1-disc} does not generalize in dimensions higher
than~$1$, where arbitrary finite sequences of Tietze transformations, interleaving
Tietze reductions and expansions, might be required in order to show that two
polygraphs present the same object. However, an analogous of
\thmr{tietze-equiv-1} will still hold, but its proof has to be carried over
differently, as explained  in \chapr{2-tietze}.

\section{Abstract Rewriting Systems}
\label{sec:ars}
The orientations of the relations do not really matter in a $1$\nbd-polygraph, with
respect to the presented set: if we reverse an edge, the presented set is the
same. This is easily shown using the following series of Tietze transformations:
\[
\pres{P_0}{P_1',x\to y}
\enspace\overset{\trel}\rightsquigarrow\enspace
\pres{P_0}{P_1',x\to y,y\to x}
\enspace\overset{\trrel}\rightsquigarrow\enspace
\pres{P_0}{P_1',y\to x}
\]
which are based on the fact that $\approx$ is an equivalence relation and thus
symmetric.

However, the orientations can still be useful to \emph{decide equality} between
generators, \ie answer the following question:

\begin{quote}
  Given two generators, do they represent the same element of the presented set?
  Or, equivalently, are they related by~$\approx$?
\end{quote}

\noindent
We will see that in good cases, one can come up with canonical representatives
of equivalence classes under $\approx$, in such a way that the representative of
an arbitrary generator can easily be computed. In those situations, the
equivalence of two generators can be tested by checking whether their
representatives are equal or not.
In order to come up with representatives, we use the orientation of the
1-generators. Given two 0-generators $x$ and $y$ such that there is a
1-generator $a:x\to y$, we have $x\approx y$, and the orientation of the
1-generator will be interpreted as indicating that $y$ is a ``more canonical''
representative than~$x$ in the equivalence class under~$\approx$. With respect
to this, the ``most canonical'' elements, which are called normal forms, are
good candidates for being representatives of equivalence classes with good
properties: under reasonable assumptions, it can be shown that every class
admits exactly one such representative.
This point of view is the starting point of \emph{rewriting
  theory}~\cite{BaaderNipkow98, bezem2003term}.

\subsection{Terminology and notations}
\index{rewriting!rule}
\index{rewriting!step}
\index{path!rewriting}
\index{abstract rewriting system}
\index{rewriting!system!abstract}
We have seen that a 1-polygraph~$P$ is simply another name for a graph. Since
people in rewriting theory like to think about it from a different point of
view, they give it yet another name and call it an \emph{abstract rewriting
  system}. In this context, the elements
of~$P_0$ are called \emph{objects} and those of~$P_1$ are called \emph{rewriting
  rules} (or \emph{rewriting steps}). A \emph{rewriting
  path} is simply a path, \ie a sequence
\[
  \xymatrix{
    x_0\ar[r]^{a_1}&x_1\ar[r]^{a_2}&x_2\ar[r]^{a_3}&\ldots\ar[r]^{a_n}&x_n
  }
\] 
of composable rewriting steps. The 0-cells $x_0$ and $x_n$ are respectively
called the \emph{source} and \emph{target} of the path, and we write
$f:x\overset *\to y$ for a path~$f$ from~$x$ to~$y$.
One also writes $x\to y$ (\resp $x\overset *\to y$) when there exists a
rewriting step (\resp a rewriting path) from~$x$ to~$y$, and the notation
$x\overset*\leftrightarrow y$ is often used instead of~$x\approx y$.

\subsection{Normal forms}
\index{normal form}
A 0-cell~$x\in P_0$ is a \emph{normal form} when there is no rule $a:x\to y$
in~$P_1$ with $x$ as source.

We can distinguish the following situations concerning normal forms in
equivalence classes under~$\approx$ of 0-cells in a polygraph~$P$: we say
that~$P$ has
\begin{itemize}
\item the \emph{existing normal form property} when every equivalence class
  contains at least one normal form, \ie for every~$x\in P_0$ there exists a
  normal form $y\in P_0$ such that $x\overset*\leftrightarrow y$,
\item the \emph{unique normal form property} when every equivalence class
  contains at most one normal form, \ie for every normal forms~$x,y\in P_0$,
  $x\overset*\leftrightarrow y$ implies~$x=y$,
\item the \emph{canonical form property} when every equivalence class contains
  exactly one normal form, called the \emph{canonical representative} of the
  class, \ie it satisfies both the existing and the unique normal form property.
\end{itemize}

\begin{example}
  \label{ex:1-fp}
  Consider the following 1-polygraphs:
  \[
    \begin{array}{c@{\qquad\qquad}c@{\qquad\qquad}c}
      \vxym{
        x\ar[r]&y\ar@(ur,dr)
      }
      &
      \vxym{
        x&\ar[l]y\ar[r]&z
      }
      &
      \vxym{
        x&\ar[l]y\ar[r]&z\ar@/^2ex/[ll]
      }
      \\[2ex]
      (1)&(2)&(3)
    \end{array}
  \]
  $(1)$ and $(3)$ have the unique normal form property, $(2)$ and $(3)$ have the
  existing normal form property, and $(3)$ has the canonical form property.
\end{example}

We are interested here in providing practical conditions on~$P$ that ensure
that the canonical form property holds, as well as that we are able to efficiently
compute the canonical form associated to the class of a 0-cell. We will see
that termination of a 1-polygraph implies the existing normal form, that confluence
implies the unique normal form property, and moreover that confluence can be
checked locally for terminating 1-polygraphs.

\subsection{Normalizability}
\label{sec:normalizability}
\index{polygraph!normalizing}
\index{normalizing!1-polygraph}
A polygraph is \emph{normalizing} when every 0-cell~$x$ re\-writes to a normal
form. We sometimes write $\nf x$ for an arbitrary choice of such a normal
form.
From the definition, we deduce the following result.

\begin{lemma}
  \label{lem:1pol-norm-enf}
  A normalizing 1-polygraph has the existing normal form property.
\end{lemma}

\noindent
The converse does not hold, as illustrated in \cref{ex:huet}.

\subsection{Termination}
\label{sec:1-termination}
\index{termination}
\index{well-founded!1-polygraph}
\index{noetherian!1-polygraph}
\index{terminating!1-polygraph}
\index{polygraph!terminating}
In practice, in order to show that a 1-polygraph is normalizing, one often
uses the following property. A 1-polygraph~$P$ is \emph{terminating} (or
\emph{well-founded} or \emph{noetherian} or \emph{strongly normalizing}) when
there is no infinite sequence of rewriting steps
\[
  \vxym{
    x_0\ar[r]^{a_1}&
    x_1\ar[r]^{a_2}&
    x_2\ar[r]^{a_3}&
    \cdots\pbox.
  }
\]
For instance, in \exr{1-fp}, (2) and (3) are terminating but not (1).

Starting from a 0-cell~$x$ in a terminating 1-polygraph, we can define a
sequence of 0-cells by induction by $x_0=x$, and $x_{i+1}$ is the target of an
arbitrary rewriting rule $x_i\to x_{i+1}$ with $x_i$ as source; we stop if there
is no such rewriting rule. Termination ensures that this process will end after
a finite number of steps, and the last 0-cell $x_n$ is necessarily a normal
form. We have just shown the following.

\begin{lemma}
  \label{lem:1pol-term-norm}
  A terminating 1-polygraph is normalizing.
\end{lemma}

\noindent
The converse does not hold, as illustrated in \exr{huet}.

In practice, the termination of a 1-polygraph~$P$ can be shown using the
following lemma. We recall that a poset~$(N,\preceq)$ is \emph{well-founded} when
every decreasing sequence $n_1\succeq n_2\succeq\ldots$ is eventually stationary:
there exists $k\in\N$ such that for every $i,j\in\N$ with $i\geq j\geq k$
one has $n_i=n_j$. Equivalently, the poset is well-founded when there
exists no infinite strictly decreasing sequence $n_1\succ n_2\succ \ldots$ of elements
of~$N$. The typical example of such an order is $(\N,\leq)$, or any ordinal.

\begin{lemma}
  Given a rewriting system~$P$ the following statements are equivalent.
  \begin{enumerate}
  \item The rewriting system $P$ is terminating.
  \item There exists a well-founded order on~$P_0$ such that $x\succ y$ for
    every 1\nbd-gene\-rator~$a:x\to y$ in~$P_1$.
  \item There exists a function~$f:P_0\to N$, where $N$ is a well-founded poset,
    such that $f(x)>f(y)$ for every 1-generator $a:x\to y$ in~$P_1$.
  \end{enumerate}
\end{lemma}
\begin{proof}
  Suppose that~$P$ is terminating. Then the preorder relation on~$P_0$ defined
  by $x\succeq y$ whenever $x\overset *\to y$ is a well-founded partial order
  that shows that~1 implies 2, and taking $f:P_0\to P_0$ to be the identity
  shows that 2 implies~3. Finally, 3 implies 1, for if there was an infinite
  reduction sequence in~$P$, the image of the objects under~$f$ would be an
  infinite strictly decreasing sequence of elements of~$N$.
\end{proof}

\subsection{Well-founded induction}
\label{sec:wfind}
\index{well-founded!induction}
Suppose given a predicate~$\mathcal{P}$ on the $0$\nbd-cells of a
terminating polygraph~$P$. In order to show that $\mathcal{P}$ holds for all the
elements of~$P_0$, it is often useful to use the following \emph{well-founded
  induction principle}: if
\begin{equation}
  \label{eq:wfi-hyp}
  \forall x\in P_0,\qquad
  \big(\pa{\forall y\in P_0,\ x\to y\ \text{implies}\ \mathcal{P}\pa{y}}\quad\text{implies}\quad\mathcal{P}(x)\big)
\end{equation}
then $\forall x\in P_0,\ \mathcal{P}(x)$ holds.

\begin{proposition}
  \label{prop:1-pol-wfi}
  If~$P$ is a terminating 1-polygraph then the well-founded induction
  principle holds.
\end{proposition}
\begin{proof}
  By contradiction, suppose that the well-founded induction principle does not
  hold: there is a predicate~$\mathcal{P}$, such that the
  hypothesis~\eqref{eq:wfi-hyp} holds but not the
  conclusion, \ie $\mathcal{P}(x_0)$ does not hold for
  some~$x_0\in P_0$. By repeated use of~\eqref{eq:wfi-hyp}, we can construct a
  family $(x_i)_{i\in\N}$ of elements of~$P_0$ such that $\mathcal{P}(x_i)$
  does not hold for any $i\in\N$, and $x_0\to x_1\to\cdots$.
  This contradicts the fact that~$P$ is terminating.
\end{proof}

\subsection{Quasi-termination}
\label{sec:quasi-termination}
\index{quasi-!termination}
Following~\cite{dershowitz1987termination}, we introduce the following variant
of the termination condition. We say that a 1-polygraph $P$ is
\emph{quasi-termi\-nating} if every sequence $(x_i)_{i\in\N}$ of 0-cells, with
$x_i\to x_{i+1}$ for every index~$i\in\N$, contains an infinite number of
occurrences of the same 0-cell: there exists a 0-cell $x$ such that for
every $i\in\N$, there exists $j>i$ such that $x_j=x$.

\index{quasi-!normal form}
Let $P$ be a 1-polygraph. A 0-cell $x$ is called a \emph{quasi-normal form}
if for any rewriting step $x\to y$, there exists a rewriting path from $y$ to
$x$.
\index{quasi-!convergence}
If $P$ is quasi-terminating, any 0-cell $x$ rewrites to a quasi-normal
form. Note that this quasi-normal form is neither irreducible nor unique in
general.  We say that $P$ is \emph{quasi-convergent} if it is confluent and it
quasi-terminates.

\begin{example}
  The following 1-polygraph
  \[
    \vxym{
      x\ar[r]&y\ar@/^/[r]&\ar@/^/[l]z
    }
  \]
  is quasi-terminating and quasi-convergent. Both~$y$ and~$z$ are quasi-normal
  forms.
\end{example}

The above termination and normalizability conditions ensure the existing normal
form property. We now investigate conditions implying the unique normal form
property.

\subsection{Joinability}
\index{joinability}
Two 0-cells $x,y\in P_0$ of a polygraph~$P$ are \emph{joinable} when there
exists 0-cell~$z$ such that there are rewriting paths $f:x\overset *\to z$ and
$g:y\overset *\to z$:
\[
  \svxym{
    x\ar@{.>}[dr]_\ast&&y\ar@{.>}[dl]^\ast\pbox.\\
    &z&
  }
\]

\subsection{The Church-Rosser property}
\index{Church-Rosser property}
\index{polygraph!Church-Rosser}
A 1-polygraph~$P$ has the \emph{Church-Rosser property} when any two 0-cells
$x,y\in P_0$ which are equivalent are joinable:
\[
  \svxym{
    x\ar@{.>}[dr]_\ast\ar@{<->}[rr]^\ast&&\ar@{.>}[dl]^\ast y\pbox.\\
    &z&
  }
\]

\begin{proposition}
  \label{prop:1pol-church-rosser}
  A 1-polygraph with the Church-Rosser property has the unique normal form
  property.
\end{proposition}
\begin{proof}
  Suppose given two normal forms $x$ and $y$ such that $x\approx y$. By the
  Church-Rosser property, there exists a 0-cell $z$ and rewriting paths
  $x\overset*\to z$ and $y\overset*\to z$. Since $x$ and $y$ are normal forms,
  these two paths are necessarily empty, and thus $x=y$.
\end{proof}

\noindent
The converse property is not true, as illustrated by the 1-polygraph
\[
  \vxym{
    x&\ar[l]y\ar[r]&z\ar@(ur,dr)
  }
\]
where $x$ and $z$ are equivalent, but cannot be rewritten to a common 0-cell,
even though there is a unique normal form $x$.

In the following, we present more ``local'' properties which imply the
Church-Rosser property, and thus the unique normal form property.

\subsection{Branchings}
\index{branching}
In a 1-polygraph~$P$, a pair $(a,a')$ of coinitial 1-generators
$
a:x\to y
$ and $
a':x\to y'
$
in~$P$ is called a \emph{local branching}; a pair $(f,f')$ of coinitial
rewriting paths
$
f:x\overset*\to y
$ and $
f':x\overset*\to y'
$
is called a \emph{branching}. The 0-cell~$x$ is called the \emph{source} of
the branching.

\subsection{Confluence}
\label{sec:1-confl}
\index{confluence}
\index{polygraph!confluent}
A branching $(f,f')$ as above is \emph{confluent} when $y$ and $y'$ are
joinable:
\[
  \svxym{
    &\ar[dl]_\ast x\ar[dr]^\ast&\\
    y\ar@{.>}[dr]_\ast&&\ar@{.>}[dl]^\ast y'\pbox.\\
    &z&
  }
\]
In this situation, we say that the branching is \emph{confluent}. A
1-polygraph is \emph{confluent} (\resp \emph{locally confluent}) when every
branching (\resp local branching) is confluent. Note that a confluent
1-polygraph is necessarily locally confluent.

The above confluence conditions can be summarized graphically as follows:
\[
  \begin{array}{c@{\qquad\qquad}c@{\qquad\qquad}c}
    \svxym{
    &{}\phantom{x}&\\
    y\ar@{.>}[dr]_\ast\ar@{<->}[rr]^\ast&&\ar@{.>}[dl]^\ast y'\\
    &z&\\
    }
    &
    \svxym{
      &\ar[dl]_\ast x\ar[dr]^\ast&\\
      y\ar@{.>}[dr]_\ast&&\ar@{.>}[dl]^\ast y'\\
      &z&\\
    }
    &
    \svxym{
      &\ar[dl]x\ar[dr]&\\
      y\ar@{.>}[dr]_\ast&&\ar@{.>}[dl]^\ast y'\pbox.\\
      &z&\\
    }
    \\
    \text{Church-Rosser}
    &
    \text{confluence}
    &
    \text{local confluence}
  \end{array}
\]

\begin{proposition}
  \label{prop:1pol-confl-cr}
  A 1-polygraph has the Church-Rosser property if and only if it is confluent.
\end{proposition}
\begin{proof}
  The left-to-right direction is immediate. For the right-to-left direction,
  suppose that $x$ and $y$ are two equivalent 0-cells: this means that there
  exists rewriting paths $f_i:y_i\overset*\to x_i$ and
  $g_i:y_i\overset*\to x_{i+1}$ in~$P$, with $0\leq i<n$, where $x_0=x$ and
  $x_n=y$, forming a diagram as below (ignoring the dotted arrows, $z$ and~$z'$):
  \[
    \xymatrix@!C=1ex@!R=1ex{
      &\ar[dl]_{f_0}^\ast y_0\ar[dr]^{g_0}_\ast&&\ar[dl]_{f_1}^\ast y_1\ar[dr]^{g_1}_\ast&&\ar[dl]_{f_2}^\ast y_2\ar@{}[r]|{\cdots}&y_{n-2}\ar[dr]^{\hspace{-5pt}g_{n-2}}_\ast&&\ar[dl]_{f_{n-1}\hspace{-5pt}}^\ast y_{n-1}\ar[dr]^{g_{n-1}}_\ast\\
      x_0\ar@{.>}@/_/[drrrrr]_\ast&&x_1&&x_2&\ar@{}[r]|{\cdots}&&x_{n-1}\ar@{.>}[dr]_\ast&&\ar@{.>}[dl]^\ast x_n\pbox.\\
      &&&&&z'\ar@{}[u]|{\displaystyle\textsc{ih}}&&&\ar@{.>}[lll]^\ast z\ar@{}[uu]|{\displaystyle\textsc{c}}
    }
  \]
  By induction on~$n\in\N$, we show that~$x_0$ and $x_n$ can be joined. The
  result is immediate when $n=0$, and otherwise the diagram can be completed as
  above using the confluence hypothesis for \textsc{c} and the induction
  hypothesis for \textsc{ih}.
\end{proof}

\noindent
As a direct corollary, we deduce:

\begin{lemma}
  \label{lem:1pol-confl-unf}
  A confluent 1-polygraph has the unique normal form property.
\end{lemma}

\noindent
Confluence is difficult to show in practice, whereas local confluence is much
more tractable. Clearly confluence of a rewriting system implies its local
confluence, and one could hope that both properties are equivalent. This is
however not the case: local confluence does not imply confluence in general, as
illustrated by the following example attributed by Hindley to Kleene,
see~\cite[Figure~6b]{huet1980confluent} and~\cite[Section~1.2]{bezem2003term}.

\begin{example}
  \label{ex:huet}
  Consider the following 1-polygraph:
  \[
  \vxym{
    x'&\ar[l]x\ar@/^/[r]&\ar@/^/[l]y\ar[r]&y'\pbox.
  }
  \]
  It is locally confluent (it is easy to check all the possible cases), but not
  confluent: we have $x\overset *\to x'$ and $x\overset *\to y'$, but there is
  no 0-cell to which both $x'$ and $y'$ rewrite.
\end{example}

\noindent
In the previous example, it can be noted that the rewriting system is not
terminating since there is a directed cycle between the vertices~$x$ and~$y$. It
was shown in a famous lemma by Newman~\cite{newman1942theories}, also known as
the \emph{diamond lemma}, that local confluence and confluence are equivalent
when restricting to terminating rewriting systems, thus providing us with simple
ways of checking for their confluence.

\begin{lemma}
  \label{lem:term-lc-confl}
  \label{lem:newman}
  \index{lemma!Newman}
  \index{Newman's lemma}
  A terminating 1-polygraph is confluent if and only if it is locally
  confluent.
\end{lemma}
\begin{proof}
  We show the right-to-left direction, the other one being immediate.
  We say that a 1-polygraph is \emph{confluent} (\resp \emph{locally
    confluent}) \emph{at a 0-cell~$x$} when every branching (\resp local
  branching) with~$x$ as source is joinable.
  By well-founded induction, whose use is justified by \propr{1-pol-wfi} based
  on the hypothesis that the 1-polygraph is terminating, we show that the
  local confluence property at a vertex~$x$ implies the confluence property
  at~$x$. The base cases are immediate. Otherwise, we have a diagram of the form
  \[
    \xymatrix@!C=1ex@!R=1ex{
      &&\ar[dl]x\ar[dr]&&\\
      &\ar[dl]_\ast y_1\ar@{.>}[dr]_\ast&\txt{\textsc{lc}}&\ar@{.>}[dl]^\ast\ar[dr]^\ast y_1'\\
      y\ar@{.>}[dr]_\ast&\txt{\textsc{ih}}&y''\ar@{.>}[dl]^\ast&\txt{\textsc{ih}}&\ar@{.>}[ddll]^\ast y'\\
      &z\ar@{.>}[dr]_\ast&&&\\
      &&z'
    }
  \]
  which can be closed using the local confluence hypothesis for \textsc{lc} and
  the induction hypothesis for \textsc{ih} (which provides confluence at $y_1$
  and $y_1'$ respectively).
\end{proof}

\begin{remark}
  Showing termination and local confluence is the most usual way of proving that
  an abstract rewriting system is confluent, but it is not the only one. We
  refer to standard rewriting textbooks for other properties which imply
  confluence~\cite{BaaderNipkow98, bezem2003term}. For instance, an abstract
  rewriting system has the \emph{diamond property} when for every pair of
  coinitial rewriting steps $a:x\to y$ and $b:x\to y'$ there exists a pair of
  cofinal rewriting steps (\ie rewriting paths of length one) $a':y\to x$ and
  $b':y'\to x$. Graphically,
  \[
    \svxym{
      &\ar[dl]x\ar[dr]&\\
      y\ar@{.>}[dr]&&\ar@{.>}[dl]y'\pbox.\\
      &z&
    }
  \]
  In this case, the abstract rewriting system is always confluent (this can be
  shown using a variant of the proof of \lemr{term-lc-confl}) even if it is not
  terminating.
\end{remark}

\subsection{Convergence}
\index{convergence}
A 1-polygraph is \emph{convergent} when it is both terminating and
confluent.

\begin{proposition}
\label{prop:1pol-conv-cf}
A convergent 1-polygraph has the canonical form property.
\end{proposition}
\begin{proof}
  Suppose given a convergent 1-polygraph. Since it is terminating, it is
  normalizing by \lemr{1pol-term-norm} and thus has the existing normal form
  property by \lemr{1pol-norm-enf}. Since it is confluent, \lemr{1pol-confl-unf}
  ensures that it also has the unique normal form property.
\end{proof}

\begin{remark}
  A polygraph can have the canonical form property without being convergent:
  \[
    \vxym{
      x\ar@/^/[r]&\ar@/^/[l]y\ar[r]&z\pbox.
    }
  \]
  Here, all the 0-cells are equivalent and $z$ is the only normal form, which
  shows the canonical form property. The polygraph is not terminating (there is
  a cycle between $x$ and $y$) and thus not convergent.
\end{remark}

\subsection{Deciding equality}
\index{word problem!for 1-polygraphs}
\index{decidability!of equality}
\label{sec:1-eq-dec}
Given a finite 1-polygraph~$P$, the \emph{equality decision problem}, or the
\emph{word problem}, for~$P$ consists in answering the following question:
\begin{center}
  Given two 0-cells $x,y\in P_0$, do we have $x\approx y$?
\end{center}
Since we only consider only finite 1-polygraphs, this problem is
\emph{decidable}, meaning that there is a program which takes $P$, $x$ and $y$
as input and outputs whether~$x\approx y$ holds or not. Namely, we can implement
a program which will construct all acyclic paths starting from~$x$, which are in
finite number, and check whether one of those paths ends at~$y$. We will see
that if we assume additional properties on~$P$, this can be performed much more
efficiently.

When the 1-polygraph~$P$ has the canonical form property, the equivalence
class of~$x$ (\resp $y$) contains a unique normal form denoted $\nf{x}$ (\resp
$\nf{y}$), and we have $x\approx y$ if and only if we have $\nf{x}=\nf{y}$. In
this case, the equality decision problem can be decided by comparing normal
forms.
In particular, in the case where the 1-polygraph is convergent, we have seen
in \propr{1pol-conv-cf} that it has the canonical form property, and moreover
the normal form~$\nf{x}$ associated to a $0$\nbd-cell~$x$ can be computed
easily. A maximal path starting from~$x$
\[
  \xymatrix{
    x=x_0\ar[r]^-{a_0}&x_1\ar[r]^{a_1}&x_2\ar[r]^{a_2}&\cdots\ar[r]^{a_{n-1}}&x_n
  }
\]
exists because~$P$ is terminating, and the fact that it is maximal means that
its target is a normal form, \ie $x_n=\nf{x}$. In order to decide whether~$x$
and~$y$ are equivalent, we can thus use the \emph{normal form algorithm}\index{normal form!algorithm}\index{algorithm!normal form} which
consists in
\begin{enumerate}
\item rewriting $x$ as much as possible in order to obtain a normal form~$\nf{x}$,
  and similarly compute a normal form~$\nf{y}$ for~$y$,
\item checking whether $\nf{x}=\nf{y}$ holds or not.
\end{enumerate}
Formally, this is justified as follows:

\begin{proposition}
  \label{prop:conv-eq-dec}
  In a convergent 1-polygraph, two 0-cells $x$ and $y$ are equivalent if and
  only if they have the same normal form: $x\approx y$ if and only
  if~$\nf{x}=\nf{y}$.
\end{proposition}
\begin{proof}
  Since the polygraph is terminating, it is normalizing by
  \lemr{1pol-term-norm}: $x$ rewrites to a normal form~$\nf{x}$, and similarly
  $y$ rewrites to a normal form~$\nf{y}$. If $\nf{x}=\nf{y}$, then clearly $x$
  and $y$ are equivalent:
  \[
    \xymatrix{
      x\ar[r]^-\ast&\nf{x}\ar@{=}[r]&\nf{y}&\ar[l]_-\ast y\pbox.
    }
  \]
  Conversely, suppose that $x$ and $y$ are equivalent, and thus that $\nf{x}$
  and $\nf{y}$ are also equivalent:
  \[
    \xymatrix{
      \nf{x}&\ar[l]_-\ast x\ar@{<->}[r]^-\ast&y\ar[r]^-\ast &\nf{y}\pbox.
    }
  \]
  The confluence of the polygraph implies that it has the Church-Rosser property
  by \propr{1pol-confl-cr}, and thus the unique normal form property by
  \propr{1pol-church-rosser}. Since $\nf{x}$ and $\nf{y}$ are equivalent normal
  forms, we deduce that they are equal.
\end{proof}

\subsection{Deciding confluence}
As a direct corollary of the above proposition, we also have a practical method
for checking whether a terminating 1-polygraph is confluent (and thus
convergent):

\begin{proposition}
  A terminating 1-polygraph is confluent if and only if for every local
  branching $x\to y$ and $x\to z$, we have $\hat y=\hat z$.
\end{proposition}


\section{Decreasing Diagrams}
\label{sec:decreasing-diag}
The main method we have seen so far in order to show the confluence of a
$1$-polygraph is provided by Newman's lemma (\lemr{newman}), which requires
supposing termination of the polygraph.
As a more advanced topic, we explain here the method of \emph{decreasing
  diagrams}, introduced by van Oostrom~\cite{van1994confluence}, see
also~\cite[Section~14.2]{bezem2003term}, which can be used in order to show the
confluence of a $1$-polygraph which is non-terminating. Stronger versions of this
method have been introduced more recently~\cite{felgenhauer2013proof,van2023residuation}.

\subsection{Multisets}
\label{sec:multisets}
Given a set~$A$, a \emph{multiset}\index{multiset} on~$A$ is a function
$\mu:A\to\N$
which is null almost everywhere, \ie the set $\setof{a\in A}{\mu(a)\neq 0}$ is
finite. The set $A$ is called the \emph{domain} of the multiset. Given an
element~$a\in A$, the natural number $\mu(a)$ is called its \emph{multiplicity}
in the multiset: $\mu$ should be thought of as a collection of elements of~$A$
where each element~$a$ occurs~$\mu(a)$ times. We denote by~$\multisets{A}$ the
set of all multisets on~$A$.

We write~$\emptyset$ for the \emph{empty multiset} on~$A$, \ie the constant
function $\emptyset:A\to\N$ equal to~$0$. Given two multisets $\mu$ and~$\nu$
on~$A$, their \emph{union} or \emph{sum} is the multiset~$\mu\mcup\nu$ on~$A$
such that
$
  (\mu\mcup\nu)(a)=\mu(a)+\nu(a)
$
for every element~$a\in A$. The operation $\mcup$ equips $\multisets{A}$ with a
structure of commutative monoid, with $\emptyset$ as neutral element, which
characterizes multisets over~$A$. Given an element $a\in A$, we often write
$\set{a}$ for the multiset with $a$ as only element.
Given two multisets $\mu$ and $\nu$, we say that $\mu$ is \emph{included} in
$\nu$, what we write $\mu\sqsubseteq\nu$ when $\mu(a)\leq\nu(a)$ for every
$a\in A$. This is the case precisely when there is a multiset $\mu'$ such that
$\mu\sqcup\mu'=\nu$. This relation makes $\multisets A$ into a poset which is
well-founded.

A partial order~$\leq$ on a set~$A$ induces an order~$\multisets\leq$ on
$\multisets{A}$, called its \emph{multiset extension}, defined by
$\mu\multisets\leq\nu$ if and only if
\[
  \forall b\in A,
  \qquad
  \mu(b)>\nu(b)
  \quad\text{implies}\quad
  \exists a\in A,
  \quad
  a>b
  \text{ and }
  \mu(a)<\nu(a)
  \pbox.
\]
Let us spell it out: for $\mu$ to be smaller than $\nu$, it is fine to have more
$b$'s as long as $\nu$ has more of something greater than~$b$.
The following result is due to Dershowitz and
Manna~\cite{dershowitz1979proving}:

\begin{proposition}
  \label{prop:multiset-ext-wf}
  Given a well-founded poset $(A,\leq)$, its multiset extension
  $(\multisets{A},\multisets\leq)$ is also well-founded.
\end{proposition}

\subsection{Labeled 1-polygraphs}
\index{labeling}
\index{polygraph!labeled}
A \emph{labeled 1\nbd-poly\-graph} $(P,\labels,\leq,\ell)$ consists of
\begin{itemize}
\item a $1$-polygraph~$P$,
\item a set $\labels$ of \emph{labels} equipped with a well-founded
  ordering~$\leq$, and
\item a function $\ell:P_1\to\labels$ associating a label to each rewriting
  step.
\end{itemize}

\subsection{Lexicographic maximum measure}
\index{lexicographic!maximum measure}
Let $(P,\labels,\leq,\ell)$ be a fixed labeled $1$\nbd-poly\-graph. We write
$\labels^*$ for the sets of words over $\labels$, \ie finite sequences of
elements of $\labels$. The empty word is noted $\emptyword$, and the
concatenation of two words~$w$ and $w$ is noted $ww'$: these operations equip
the sets of words with a structure of monoid. Following~\cite[Definition
3.1]{van1994confluence}, we define the \emph{lexicographic maximum
  measure}~$\lmm{w}$ of a word $w\in\labels^*$ as the multiset defined
inductively by
\begin{align*}
  \lmm{\emptyword}&=\emptyset,
  &
  \lmm{lw}&=\set{l}\mcup\lmm{w^{\not<l}}
  \pbox.
\end{align*}
Above, $w^{\not<l}$ is the subword of~$w$ whose letters are not strictly
below~$l$, which is formally defined by induction by
\begin{align*}
  \emptyword^{\not<l}&=\emptyword,
  &
  (aw)^{\not<l}&=
  \begin{cases}
    w^{\not<l}&\text{if $a<l$,}\\
    aw^{\not<l}&\text{otherwise.}
  \end{cases}  
\end{align*}
Informally, the multiset $\lmm{w}$ thus consists of the letters of~$w$ which are
not dominated by some letter on their left.

The measure $\lmm{\cdot}$ is extended to the set of finite rewriting paths of
$P$ by setting, for every rewriting path $a_1\ldots a_n$,
\[
  \lmm{a_1\ldots a_n}=\lmm{\ell(a_1)\ldots \ell(a_n)}\pbox,
\]
where $\ell(a_1)\ldots \ell(a_n)$ is the product in the monoid
$\labels^\ast$. Finally, the measure $\lmm{\cdot}$ is extended to the set of
finite branchings $(a,b)$ of $P$, by setting
\[
  \lmm{(a,b)}=\lmm{a} \mcup \lmm{b}\pbox.
\]

\subsection{Decreasing diagrams}
\index{decreasing!diagram}
A diagram of rewriting paths of the form
\[
  \vxym{
    \ar[d]_f\ar[r]^g&\ar[d]^{g'}\\
    \ar[r]_{f'}&
  }
\]
is \emph{decreasing} if
\[
  \lmm{ff'}\multisets\leq\lmm{f}\mcup\lmm{g}
  \qqtand
  \lmm{gg'}\multisets\leq\lmm{f}\mcup\lmm{g}
  \pbox.
\]
In the case where $f=a$ and~$g=b$ are both $1$-generators, it can be shown that
the diagram is decreasing if and only if it is of the form
\begin{equation}
  \label{eq:loc-dec-diag}
  \vxym{
    \ar[ddd]_a\ar[rrr]^b&&&\ar[d]^{g'}\\
    &&&\ar[d]^{a'}\\
    &&&\ar[d]^{h_1}\\
    \ar[r]_{f'}&\ar[r]_{b'}&\ar[r]_{h_2}&\\
  }
\end{equation}
where
\begin{itemize}
\item $l<\ell(a)$ for every label~$l$ of a rewriting step in~$f'$,
\item $l<\ell(b)$ for every label~$l$ of a rewriting step in~$g'$,
\item $a'$ is either an identity or a rewriting step labeled by~$\ell(a)$,
\item $b'$ is either an identity or a rewriting step labeled by~$\ell(b)$, and
\item $l<\ell(a)$ or $l<\ell(b)$ for every label~$l$ of a transition in~$h_1$
  (\resp in~$h_2$).
\end{itemize}
\index{polygraph!locally decreasing}
\index{locally decreasing!1-polygraph}
A labeled $1$-polygraph is \emph{locally decreasing} when every local branching
$(a,b)$ can be completed as a locally decreasing
diagram~\eqref{eq:loc-dec-diag}. We can now recall van Oostrom's
theorem~\cite[Theorem~3.7]{van1994confluence}, whose proof follows the one of
Newman's \lemr{newman}:

\begin{theorem}
  \label{thm:ld-confl}
  A locally decreasing $1$-polygraph is confluent.
\end{theorem}

\noindent
This method is complete, in the sense that given a $1$-polygraph with countably
many $0$-cells which is confluent, there is always a way to choose a well-founded
poset~$\labels$ of labels so that the polygraph is locally
decreasing~\cite[Theorem~14.2.32]{bezem2003term}. Moreover, we can always choose
the set~$\labels=\set{0,1}$ with $0<1$ as a set of labels,
see~\cite{endrullis2016decreasing}.


\chapter{Two-Dimensional Polygraphs}
\label{chap:2pol}
This chapter is dedicated to the definition of $2$-polygraphs, which
are a $2$\nbd-dimen\-sional generalization of the $1$-polygraphs presented in
the previous chapter.
Before introducing this notion of $2$-polygraph, we first give in
\secr{free-cat-gpd} a refined viewpoint over $1$-polygraphs. Instead
of merely focusing on the set presented by a $1$-polygraph $P$ as a
set of equivalence classes of $P_0$ modulo the relations in $P_1$, we
now consider the free category generated by $P$, whose set of objects
is $P_0$ and whose morphisms are all the rewriting paths obtained by
composing the elements of $P_1$. A variant of this construction is the
notion of free groupoid generated by a $1$-polygraph $P$, where all
$1$-generators are  supposed to be invertible.

The notion of $2$-polygraph, introduced in~\secr{2pol-cat}, naturally appears as soon as arbitrary,
non necessarily free, small categories are considered. In order to
present such a category~$C$, one starts as above with a polygraph~$P$
such that the elements of $P_1$ generate the morphisms of $C$, but now
we must take account of the relations induced by $C$ among the morphisms of the free
category generated by $\pres{P_0}{P_1}$. These relations will be
generated by a set $P_2$ of $2$-generators, consisting in certain
pairs of morphisms we want to equalize in~$C$, as explained in~\secr{pres-cat}.

Following the same pattern, we finally explain
in~\secr{generating-2-cat} that a
$2$\nbd-poly\-graph can also be seen as a system of generators for a free
$2$-category, thus preparing the study of $3$-polygraphs. We also examine, in~\secr{free-21-category}, the variant where we
freely generate a $(2,1)$-category, that is, a $2$-category in which every $2$-cell is
invertible. This framework allows to invoke explicit witnesses for
the confluence of rewriting paths, which are used in the proofs of
coherence results on abstract rewriting systems.

\section{Generating Categories and Groupoids}
\label{sec:free-cat-gpd}
We now define the free category and the free groupoid generated by a
$1$\nbd-poly\-graph.

\subsection{Underlying polygraph of a category}
\index{underlying!$1$-polygraph}
\label{sec:1-upol}
Any small category~$C$ has an \emph{underlying 1\nbd-poly\-graph}~$P$ with $P_0$
being the set of objects of~$C$, $P_1$ being the set of morphisms of~$C$, and the
source and target of a morphism $f:x\to y$ of~$C$ being, respectively, $x$ and~$y$.
This construction extends in the expected way into a functor
$
\fgf:\Cat\to\nPol1
$
which to every category~$C$ associates its underlying polygraph~$\fgf C$.

\subsection{Freely generated category}
\label{sec:free-1-cat}
\index{free!1-category}
A 1-polygraph~$P$ induces a category~$\freecat{P}$, called the category
\emph{freely generated} by~$P$, or the \emph{free category on $P$}, defined as follows:
\begin{itemize}
\item its objects are the $0$-cells of~$P$,
\item its morphisms from $x$ to $y$ are composable sequences of $1$-generators
  \[
    \vxym{
      x=x_0\ar[r]^-{a_1}&x_1\ar[r]^-{a_2}&x_2\ar[r]^-{a_3}&\ldots\ar[r]^-{a_n}&x_n=y
    }
  \]
  which, more precisely, consist of a sequence $(x_i)_{0\leq i\leq n}$ of
  elements of $P_0$, with $x_0=x$ and $x_n=y$, together with a sequence
  $(a_i)_{0<i\leq n}$ of elements of~$P_1$, for some $n\geq 0$, such that
  $\sce0(a_{i+1})=x_i$ and $\tge0(a_{i+1})=x_{i+1}$ for $0\leq i<n$,
\item identities are morphisms as above with $n=0$, and
\item the composition of two morphisms
  \[
  \vxym{
    x=x_0\ar[r]^-{a_1}&x_1\ar[r]^-{a_2}&x_2\ar[r]^-{a_3}&\ldots\ar[r]^-{a_n}&x_n=y
  }
  \]
  and
  \[
  \vxym{
    y=y_0\ar[r]^-{b_1}&y_1\ar[r]^-{b_2}&y_2\ar[r]^-{b_3}&\ldots\ar[r]^-{b_m}&y_m=z
  }
  \]
  is
  \[
  \vxym{
    x=x_0\ar[r]^-{a_1}&x_1\ar[r]^-{a_2}&\ldots\ar[r]^-{a_n}&x_n=y_0\ar[r]^-{b_1}&y_1\ldots\ar[r]^-{b_m}&y_m=z\pbox{.}
  }
  \]
\end{itemize}
In the terminology of directed graphs, the morphisms of the
category~$\freecat{P}$ are the \emph{directed paths}\index{directed path} in~$P$, with the natural
number~$n$ for a path as above being its \emph{length}\index{length}, identities
are empty paths, and composition is given by
\emph{concatenation}\index{concatenation} of paths. Following the terminology of
rewriting systems, we will also call a morphism in~$\freecat{P}$ a
\emph{rewriting path}.

The category~$\freecat{P}$ can be characterized, up to isomorphism, as the
category satisfying the following universal property:

\begin{lemma}
  \label{lem:free-1-cat}
  For any category~$C$ and morphism of 1-polygraphs $f:P\to\fgf C$, there exists a
  unique functor $\freecat{f}:\freecat{P}\to C$ such that
  the following diagram commutes
  \[
  \vxym{
    P\ar[d]_-i\ar[r]^f&\fgf C\\
    \fgf\freecat{P}\ar@{.>}[ur]_{\fgf\freecat f}
  }
  \]
  where $i:P\to \fgf\freecat{P}$ is the morphism of 1-polygraphs sending a
  0-generator to itself and a 1-generator $a:x\to y$ to the corresponding path
  of length~1.
\end{lemma}

\noindent
By classical theorems~\cite[Section~IV.1]{MacLane98}, this is equivalent to the fact
that the operation which to every polygraph~$P$ associates the freely generated
category~$\freecat{P}$ extends to a functor
$
\freecat-
:
\nPol1
\to
\Cat
$
which is left adjoint to the functor~$\fgf$.

\subsection{Freely generated groupoid}
\label{sec:free-1-gpd}
\index{free!groupoid}
The previous construction can be modified in order to describe the groupoid
freely generated by a polygraph. We recall that a
\emph{groupoid}\index{groupoid}~$C$ is a category in which every morphism is
invertible, \ie for every morphism $f:x\to y$ of~$C$ there exists a morphism
$\finv f:y\to x$ such that $\finv f\circ f=\id_x$ and $f\circ\finv f=\id_y$. We
write $\Gpd$ for the category of groupoids and functors between them.

There is an obvious forgetful functor $\Gpd\to\Cat$, which admits a left
adjoint~\cite[Proposition~5.2.2]{borceux:hcg1}. Given a $1$-polygraph~$P$, we
write $\freegpd{P}$ for the \emph{free groupoid} on the category~$\freecat{P}$. It can
be described as the category whose objects are the elements of $P_0$, and morphisms
from~$x$ to~$y$ are composable sequences
\[
  \vxym{
    x=x_0\ar@{<->}[r]^-{f_1}&x_1\ar@{<->}[r]^-{f_2}&x_2\ar@{<->}[r]^-{f_3}&\ldots\ar@{<->}[r]^-{f_n}&x_n=y,
  }
\]
where $f_{i+1}$ is a 1-generator which is either of the form
$a_{i+1}:x_i\to x_{i+1}$ or $a_{i+1}:x_{i+1}\to x_i$, quotiented by the
equivalence relation identifying
\[
  \vxym{
    \cdots\ar@{<->}[r]^{f_{i-1}}&x\ar[r]^a&y&\ar[l]_ax&\ar@{<->}[l]_{f_{i+1}}\cdots
  }
  \quad\text{with}\quad
  \vxym{
    \cdots\ar@{<->}[r]^{f_{i-1}}&x&\ar@{<->}[l]_{f_{i+1}}\cdots
  }
\]
and
\[
  \vxym{
    \cdots\ar@{<->}[r]^{f_{i-1}}&y\ar@{<-}[r]^a&x&\ar@{<-}[l]_ay&\ar@{<->}[l]_{f_{i+1}}\cdots
  }
  \quad\text{with}\quad
  \vxym{
    \cdots\ar@{<->}[r]^{f_{i-1}}&y&\ar@{<->}[l]_{f_{i+1}}\cdots
  }
\]
which means that we can remove two adjacent occurrences of a $1$-generator~$a$
in different directions. Identities are empty sequences and composition behaves
as in \secr{free-1-cat}. This explicit construction will be presented in more
details in \secr{1-pres-free-gpd}. Note that, in the terminology of graphs, a
morphism as above is called a \emph{non-directed path}\index{path!non-directed}
in~$P$.

The interest of this construction lies in the fact that the morphisms of
$\freegpd{P}$ are ``witnesses'' for the $P$-congruence of $0$-cells, as defined
in \secr{P-cong}:

\begin{lemma}
  Two 0-cells $x,y\in P_0$ are $P$-congruent if and only if there exists a
  morphism $f:x\to y$ in~$\freegpd{P}$.
\end{lemma}

\subsection{Three functors}
To sum up, we have defined the following sequence of functors:
\[
\nPol1
\to
\Cat
\to
\Gpd
\to
\Set,
\]
where
\begin{itemize}
\item $\nPol1\to\Cat$ associates to a 1-polygraph the category it freely
  generates,
\item $\Cat\to\Gpd$ associates to a category the groupoid it freely generates, and
\item $\Gpd\to\Set$ associates to a groupoid~$C$ the corresponding quotient set
  (obtained from the set of objects of~$C$ by identifying any two objects
  between which there is a morphism).
\end{itemize}
The composite functor $\nPol1\to\Set$ was explicitly described in
\secr{pres-set}.

\section{The Category of 2-Polygraphs}
\label{sec:2pol-cat}
We are now ready to introduce $2$-polygraphs, which consist in a $1$-polygraph
together with a set of globular $2$\nbd-gene\-rators between
\ndef{parallel} $1$-cells\index{parallel!1-cells@$1$-cells}, that is, $1$-cells
having same source and target. Those were introduced by Street~\cite{street1976limits}. They are
sometimes also called \emph{linear
  sketches}~\cite[Section~4.6]{barr1990category} when seen as a particular class
of sketches, see \secr{sketches}.

\subsection{Notations on 1-cells}
Given a $1$-polygraph~$P$, we write~$\freecat{P_1}$ for the set of $1$-cells of
the category~$\freecat{P}$, \ie the set of paths in~$P$. The composite of two $1$-cells $u:x\to y$ and $v:y\to z$ is
written $uv:x\to z$, or $u\comp0 v:x\to z$, and the empty path on a $0$-generator~$x$ is written
$\unit{x}$. We also write
$
  \ins1: P_1\to\freecat{P_1}
$
for the canonical inclusion sending a $1$-generator to the corresponding path of
length $1$. Given a $1$-generator $a\in P_1$, we generally simply
write~$a$ instead of $\ins1(a)$.
The source and target functions $\sce0,\tge0:P_1\to P_0$
canonically extend to functions
$
  \freecat{\sce0},\freecat{\tge0}
  :
  \freecat{P_1}\to P_0
$
such that
\begin{equation}
  \label{eq:1-st-ins}
  \freecat{\sce0}\circ\ins1=\sce0
  \qqtand
  \freecat{\tge0}\circ\ins1=\tge0
\end{equation}
respectively sending a path~$f:x\to y$ to its source~$x$ and its target~$y$,
what we often picture as a ``commuting'' diagram of form
\[
  \vxym{
    &\ar@<-.5ex>[dl]_-{\sce0}\ar@<.5ex>[dl]^-{\tge0}P_1\ar[d]^{\ins1}\\
    P_0&\ar@<-.5ex>[l]_-{\freecat{\sce0}}\ar@<.5ex>[l]^-{\freecat{\tge0}}\freecat{P_1}\pbox.
  }
\]

\begin{definition}
\index{polygraph!2-@$2$-}
  A \emph{2-polygraph} consists of
  \begin{itemize}
  \item a 1-polygraph~$P$, \ie
    \[
      \vxym{
        P_0&\ar@<-.5ex>[l]_-{\sce0}\ar@<.5ex>[l]^-{\tge0}P_1,
      }
    \]
  \item a set $P_2$ of \emph{2-generators} together with two functions
    \[
      \sce1,\tge1
      :
      P_2\to\freecat{P_1},
    \]
    associating to each relation its \emph{source} and \emph{target}, which is
    such that
    \begin{equation}
      \label{eq:1-st}
      \freecat{\sce0}\circ \sce1=\freecat{\sce0}\circ\tge1
      \qqtand
      \freecat{\tge0}\circ \sce1=\freecat{\tge0}\circ\tge1
      \pbox.
    \end{equation}
  \end{itemize}
\end{definition}

\noindent
A 2-polygraph thus consists of a diagram of sets and functions
\[
  \vxym{
    &\ar@<-.5ex>[dl]_-{\sce0}\ar@<.5ex>[dl]^-{\tge0}P_1\ar[d]^{\ins1}&\ar@<-.5ex>[dl]_-{\sce1}\ar@<.5ex>[dl]^-{\tge1}P_2\\
    P_0&\ar@<-.5ex>[l]_-{\freecat{\sce0}}\ar@<.5ex>[l]^-{\freecat{\tge0}}\freecat{P_1}
  }
\]
together with compositions and identities on~$\freecat{P_1}$, which ``commutes''
in the sense that the relations \eqref{eq:1-st-ins} and \eqref{eq:1-st} hold.
We often write
$
  \alpha
  :
  u\To v
$
for a $2$-generator $\alpha$ in $P_2$ such that $\sce1(\alpha)=u$ and
$\tge1(\alpha)=v$, and picture it as a $2$-cell
\[
  \xymatrix@C=6ex{
    x\ar@/^3ex/[r]^u\ar@/_3ex/[r]_v\ar@{}[r]|{\phantom\alpha\Longdownarrow\alpha}&y\pbox.
  }
\]
Moreover, we sometimes write
$
\Pres{P_0}{P_1}{P_2}
$
to indicate the generators of a $2$-polygraph~$P$.
A $2$-polygraph is \emph{finite} when the sets $P_0$, $P_1$ and $P_2$ are.
The \emph{underlying 1-polygraph} of a $2$-polygraph~$P$ is
denoted~$\tpol1{P}$.\nomenclature[P1]{$\truncpol{P}{1}$}{underlying
$1$-polygraph of a $2$-polygraph}

\subsection{The category of 2-polygraphs}
\label{sec:Pol2}
A \emph{morphism}
$
  f: P\to Q
$
between two $2$-polygraphs~$P$ and~$Q$ consists of a morphism
$
  f:\tpol1 P\to\tpol1 Q
$  
between the underlying $1$-polygraphs together with a function
$
  f_2: P_2\to Q_2
$  
such that
$
\sce1^Q\circ f_2=f_1\circ\sce1^P
$
and
$
\tge1^Q\circ f_2=f_1\circ\tge1^P
$.
These compose in the expected way, and we write
$\nPol 2$\nomenclature[Pol2]{$\nPol{2}$}{category of $2$-polygraphs} for the category of
$2$-polygraphs and their morphisms.

\section{Presenting Categories}
\label{sec:pres-cat}

\subsection{Quotient categories}
\index{congruence!on a category}
Given a category~$C$, a \emph{congruence}~$\approx$ on~$C$ is an equivalence
relation on the morphisms of~$C$ such that
\begin{itemize}
\item given $u:x\to y$ and $u':x'\to y'$, $u\approx u'$ implies $x=x'$ and
  $y=y'$, and
\item given morphisms $u:x'\to x$, $v,v':x\to y$ and $w:y\to y'$, if
  $v\approx v'$ then $uvw\approx uv'w$:
  \[
    \vxym{
      x\ar@/^/[r]^v\ar@/_/[r]_{v'}\ar@{}[r]|{\vapprox}&y
    }
    \qquad\qquad\text{implies}\qquad\qquad
    \xymatrix{
      x'\ar[r]^u&x\ar@/^/[r]^v\ar@/_/[r]_{v'}\ar@{}[r]|{\vapprox}&y\ar[r]^w&y'\pbox.
    }
  \]
\end{itemize}
In such a situation, one defines the \emph{quotient category}\index{quotient!category}
$
  C/{\approx}
$
as the category whose objects are those of~$C$, whose morphisms are the equivalence classes of morphisms of~$C$, and where composition and identities are induced by those of~$C$.

\subsection{$P$-congruence}
\label{sec:P-congruence}
Given a $2$-polygraph~$P$, the \emph{$P$-congruence}~$\approx^P$, also sometimes
noted $\overset*\Leftrightarrow$, is the smallest congruence on~$\freecat{P}$
such that $u\approx^P v$ for every $2$-generator $\alpha:u\To v$ in~$P_2$.

\subsection{Presented category}
\index{presentation!of a category}
The category~$\prescat{P}$ \emph{presented} by a $2$-polygraph~$P$ is the
category
$
  \prescat{P}
  =
  \freecat{\tpol1 P}/{\approx^P}
$
obtained by quotienting the category freely generated by the underlying
polygraph by the $P$-congruence, what we usually write
$
\freecat{\tpol1{P}}/P_2
$. It can be characterized by the following universal property:

\begin{lemma}
  \label{lem:2-pres}
  Given any category~$C$ and functor $f:\freecat{\tpol1 P}\to C$ such that
  $f(u)=f(v)$ for any $2$-generator $\alpha:u\To v$ in~$P_2$, there exists a
  unique functor $\prescat{f}:\prescat{P}\to C$ such that $\prescat{f}\circ q=f$
  \[
    \vxym{
      \freecat{\tpol1 P}\ar[d]_-q\ar[r]^f&C\\
      \prescat{P},\ar@{.>}[ur]_{\prescat f}
    }
  \]
  where~$q$ is the \emph{quotient functor} which is the identity on objects and
  sends a morphism to its equivalence class under~$\approx^P$.
\end{lemma}

\noindent
Note that, with the notations of the above lemma, two morphisms $u$ and $v$
in~$\freecat{P_1}$ are such that $u\approx^Pv$ if and only if $q(u)=q(v)$.

We say that a category~$C$ is \emph{presented} by a $2$-polygraph~$P$, or
that~$P$ is a \emph{presentation} of~$C$, when~$C$ is isomorphic to~$\prescat{P}$.

\subsection{Presenting monoids}
\index{presentation!of a monoid}
\index{monoid}
\label{sec:pres-mon}
Any monoid~$M$\index{monoid} can canonically be seen as a category with only one
object~$\star$, the morphisms of the category being the elements of~$M$, and
compositions and identities being given by multiplication and unit of the
monoid. This construction extends as a functor $\Mon\to\Cat$ from the category
of monoids to the category of small categories, which is full and faithful.
In the following, by a \emph{presentation of a monoid}, we will always implicitly mean a presentation of the associated
category. Those provide a most abundant source of examples of presentations, see
\apxr{2ex}.

\begin{example}
\label{Example:MonoidWithTwoElements1}
There are exactly two monoids with two elements. They are presented by the following two $2$-polygraphs:
\begin{align*}
P &\: = \:  \Pres{\star}{a:\star\to\star}{\alpha:aa\To 1},\\
Q &\: = \:  \Pres{\star}{a:\star\to\star}{\alpha:aa\To a}.
\end{align*}
The monoid presented by $P$ is $\N/2\N$, see \exr{N2-rewr} for details.
\end{example}

\begin{example}
  \index{permutation}
  \index{symmetric!group}
  The \emph{symmetric group}~$S_n$ is the monoid whose elements are bijections
  on a set with $n$ elements, sometimes also called \emph{permutations}, with composition as multiplication and identities
  as neutral elements. It admits a presentation by a $2$-polygraph~$P$ with
  $P_0=\set\star$, $P_1=\set{a_0,\ldots,a_{n-1}}$, with each $1$-generator having
  $\star$ as source and target, and the relations in $P_2$ are
  \begin{itemize}
  \item $a_ia_i\To 1$, for $0\leq i<n$,
  \item $a_ia_{i+1}a_i\To a_{i+1}a_ia_{i+1}$, for $0\leq i<n-1$, and
  \item $a_ia_j\To a_ja_i$, for $0\leq i<j<n$ with $i+1<j$.
  \end{itemize}
  A generator $a_i$ corresponds here to the transposition exchanging the $i$-th
  and the $(i{+}1)$-th element of the set with $n$ elements, and the reader is
  encouraged to check for himself that these relations make sense, see
  \cref{sec:sym-group,sec:pres-sym} for details.
\end{example}

\subsection{A characterization of presentations}
The following lemma characterizes when a $2$-polygraph is a presentation of a
category~$C$. In practice, it is quite cumbersome to use and more practical
tools will be introduced in \chapr{2rewr}, based on rewriting. We recall that
the underlying 1\nbd-polygraph~$\fgf C$ of a category~$C$ was defined in
\secr{1-upol}.

\begin{lemma}
  \label{lem:2-pres-cond}
  A 2-polygraph~$P$ is a presentation of a category~$C$ if and only if there is
  a morphism of 1-polygraphs
  $
  f:\tpol1 P\to\fgf C
  $
  such that the following three conditions hold.
  \begin{enumerate}
  \item The map $f_0:P_0\to C_0$ is a bijection between $0$-generators and objects
    of~$C$.
  \item For any generator $\alpha:u\To v$ in~$P_2$, $f(u)=f(v)$.
  \item The function $\freecat{f_1}:\freecat{P_1}\to C_1$ induces a bijection
    between $\freecat{P_1}/P_2$ and $C_1$.
  \end{enumerate}
\end{lemma}
\begin{proof}
  By post-composition with the counit $\freecat{(\fgf C)}\to C$ of the
  adjunction described in \secr{free-1-cat}, the functor
  $\freecat{f}:\freecat{\tpol1 P}\to\freecat{(\fgf C)}$ induces a functor
  $\freecat{\tpol1 P}\to C$ that we still write $\freecat{f}$.
  The second condition amounts to require that it induces, by \lemr{2-pres}, a
  well-defined quotient functor $\prescat{\freecat{f}}:\ol P\to C$. The first
  condition amounts to require that this functor is bijective on objects, and
  the third that it is full and faithful.
\end{proof}

\subsection{Models}
\index{model!of a 1-polygraph}
\label{sec:2-models}
Given categories~$C$ and~$S$, the \emph{category of models} of~$C$ in~$S$ is the
category
\[
  \Mod_S(C)
  =
  \Cat(C,S)
\]
of functors from~$C$ to~$S$ and natural transformations between those. We simply
write $\Mod(C)$ in the case $S=\Set$.

\begin{lemma}
  \label{lem:1-model-pres}
  If $P$ is a $2$-polygraph and $C$ category, the category~$\Mod_C(\pcat{P})$
  is isomorphic to the category whose objects consist of
  \begin{itemize}
  \item a family $(f_x)_{x\in P_0}$ of objects of~$C$ indexed by $0$-generators
    in~$P_0$, and
  \item a family $(f_a:f_x\to f_y)_{a:x\to y\in P_1}$ of morphisms in~$C$
    indexed by $1$-generators in~$P_1$,
  \end{itemize}
  such that for every $2$-generator
  $
    \alpha
    :
    a_1\ldots a_m
    \To
    b_1\ldots b_n
  $
  in~$P_2$, we have
  \[
    f_{a_1}\ldots f_{a_m}
    =
    f_{b_1}\ldots f_{b_n}
  \]
  and a morphism~$\phi$ between two objects~$f$ and~$g$ consists of a family
  \[
    (\phi_x:f_x\to g_x)_{x\in P_0}
  \]
  of morphisms of~$C$ indexed by the $0$-generators of~$P$ such that, for every
  $1$\nbd-generator $a:x\to y\in P_1$, the following diagram commutes
  \[
    \svxym{
      f_x\ar[d]_{f_a}\ar[r]^{\phi_x}&g_x\ar[d]^{g_a}\\
      f_y\ar[r]_{\phi_y}&g_y\pbox.
    }
  \]
\end{lemma}

\begin{example}
  \index{isomorphism!walking}
  \label{ex:walking-iso}
  Consider the category~$C$ with two objects~$X$ and~$Y$, and a single morphism
  in each hom-set:
  \[
    \vxym{
      X\ar@(dl,ul)^{\id_X}\ar@/^/[r]^{F}&\ar@/^/[l]^{G}Y\ar@(dr,ur)_{\id_G {\displaystyle.}}
    }
  \]
  This category admits a presentation by the $2$-polygraph
  \[
    P=
    \Pres{x,y}{a:x\to y,b:y\to x}{\alpha:ab\To\id_x,\beta:ba\To\id_y},
  \]
  which can be shown by using \lemr{2-pres-cond}. We define
  a morphism of~$1$\nbd-poly\-graphs $f:\tpol1 P\to\fgf C$ by
  \[
    f(x)=X,
    \qquad\qquad
    f(y)=Y,
    \qquad\qquad
    f(a)=F,
    \qquad\qquad
    f(b)=G
  \]
  and check the conditions of \lemr{2-pres-cond}.
  \begin{enumerate}
  \item The map $f_0:\set{x,y}\to\set{X,Y}$ is a bijection.
  \item The map $\freecat{f}$ preserves the $2$-generator $\alpha:ab\To\id_x$:
    \[
      \freecat{f}(ab)
      =
      \freecat{f}(a)\freecat{f}(b)
      =
      FG
      =
      \id_X
      =
      \id_{\freecat{f}(x)}
      =
      \freecat{f}(\id_x)
      \text,
    \]
    and similarly for~$\beta$.
  \item The morphisms of~$\freecat{P_1}$ are of the form
    \begin{align*}
      (ab)^n&:x\to x,
      &
      (ab)^na&:x\to y,
      &
      (ba)^nb&:y\to x,
      &
      (ba)^n&:y\to y,\\
      \intertext{for some~$n\in\N$, and those are respectively equivalent to}
      \id_x&:x\to x,
      &
      a&:x\to y,
      &
      b&:y\to x,
      &
      \id_y&:y\to y,
    \end{align*}
    by induction on~$n\in\N$, because the presence of $\alpha$ and $\beta$
    respectively imply that we have $ab\approx\id_x$ and $ba\approx\id_y$. These
    are distinct (they have different types) and are in bijection with the
    morphisms of~$C$ (there is one in each hom-set), in a way compatible with
    source and target.
  \end{enumerate}
  By \lemr{1-model-pres}, a model of~$C$ in a category~$S$ consists of
  \begin{itemize}
  \item two objects $f_x$ and $f_y$ of~$S$, and
  \item two morphisms $f_a:f_x\to f_y$ and $f_b:f_y\to f_x$ of~$S$,
  \end{itemize}
  such that $f_b\circ f_a=\id_{f_x}$ and $f_a\circ f_b=\id_{f_y}$.
  The algebras of~$C$ in~$S$ are thus precisely the isomorphisms in~$S$. Otherwise
  said, the category~$C$ represents the functor~$\Cat\to\Set$ sending a category
  to its set of isomorphisms. For this reason, $C$ is sometimes called the
  \emph{walking isomorphism}.
\end{example}

\subsection{Free categories}
\index{free!category}
\label{sec:Z-not-free}
A category~$C$ is \emph{free} when it admits a presentation by a 2-polygraph~$P$
which has no relations (\ie $P_2=\emptyset$): in this case, the category $C$ can
be obtained as the category freely generated by the underlying 1\nbd-poly\-graph
of~$P$. Contrary to the case of sets, see \lemr{1-disc}, not every category
is free: the relations are really needed in order to have a presentation for
every category.

Consider~$\Z$ as the category with one object $\star$, with the morphisms being the
integers with composition given by addition and identity by zero. Suppose given
a presentation without relations of this category. This presentation necessarily
has exactly one 0-generator $\star$, and at least two generators $a,b$: if there
was zero (\resp one) generator, the presented category would be the terminal one
(\resp the one corresponding to the additive monoid $\N$). Writing $P$ for the
underlying polygraph of the presentation, there is an isomorphism
$f:\freecat{P}\to C$. Since~$\Z$ is abelian, we have $f(a)f(b)=f(b)f(a)$ and
therefore $ab=ba$ in~$\freecat{P}$, which does not actually hold
in~$\freecat{P}$. By contradiction, every presentation of~$\Z$ has at least one
relation, \ie $\Z$ is not free. An actual presentation of~$\Z$ is given in
\secr{int-monoid}. For similar reasons, the category corresponding to the monoid
$\N\times\N$ (or in fact any abelian monoid, excepting $\N$ and the
trivial monoid) is not free.

\subsection{Canonical and standard presentations}
\label{sec:2-std-pres}
Any category admits a presentation, in the following way.
Suppose given a category~$C$ and write~$P=\fgf C$ for its underlying
$1$-polygraph, whose $0$-generators are the objects of~$C$ and $1$-generators
are the morphisms of~$C$, see \secr{1-upol}. The identity morphism
$\id_P:P\to\fgf C$ extends, by \lemr{free-1-cat}, as a functor
$\freecat{\id_P}:\freecat{P}\to C$, which sends a $1$-cell in $\freecat{P}$, \ie
a formal composite of morphisms in~$C$, to the result of its composition. The
$2$\nbd-poly\-graph with~$P$ as underlying $1$-polygraph and whose set of
$2$-generators is
\[
  P_2=\setof{(u,v)\in\freecat{P_1}\times\freecat{P_1}}{\freecat{\id_P}(u)=\freecat{\id_P}(v)},
\]
with $\sce1(u,v)=u$ and $\tge1(u,v)=v$, is a presentation of~$C$ called its
\emph{canonical presentation}. 
\index{canonical presentation}
\index{presentation!canonical}

If, in the previous presentation, we restrict the set~$P_2$ to the
$2$\nbd-gene\-rators which are either of the form $(ab,c)$ or of the form
$(\id_x,a)$, with $a,b,c\in P_1$, we obtain another presentation of~$C$, which
is smaller, called its \emph{standard presentation}, detailed in
\secr{cat-std-pres}.
\index{standard!presentation}
\index{presentation!standard}

\section{Generating 2-Categories}
\label{sec:generating-2-cat}
In the same way a $1$-polygraph generates a category which is a ``graph with
compositions'', a $2$-polygraph also generates a $2$-category which is a
``$2$-graph with compositions''. Here, a $2$-graph consists of a graph together
with ``$2$-cells'' which have edges as source and target.

\subsection{2-graphs}
\index{2-graph}
\index{graph!2-}
\index{globular set!2-}
\label{sec:2-graph}
A \emph{2-graph}~$C$, or \emph{2-globular set},
\[
  \xymatrix{
    C_0&\ar@<-.5ex>[l]_{\sce0}\ar@<.5ex>[l]^{\tge0}C_1&\ar@<-.5ex>[l]_{\sce1}\ar@<.5ex>[l]^{\tge1}C_2
  }
\]
consists in sets
\begin{itemize}
\item $C_0$ of \emph{0-cells},
\item $C_1$ of \emph{1-cells} together with functions $\sce0,\tge0:C_1\to C_0$
  respectively associating to each 1-cell its source and target 0-cell, and
\item $C_2$ of \emph{2-cells} together with functions $\sce1,\tge1:C_2\to C_1$
  respectively associating to each 2-cell its source and target 1-cell,
\end{itemize}
such that
\begin{align*}
  \sce0\circ\sce1&=\sce0\circ\tge1,
  &
  \tge0\circ\sce1&=\tge0\circ\tge1\pbox.
\end{align*}
We often write
$
  a
  :
  x
  \to
  y
$
for a $1$-cell~$a$ with $\sce0(a)=x$ and $\tge0(a)=y$, and
\[
  \alpha
  :
  a\To b
  :
  x\to y,
\]
for a $2$-cell~$\alpha$ with
\[ \sce1(\alpha)=a, \quad \tge1(\alpha)=b, \quad
\sce0(a)=\sce0(b)=x, \qtand \tge0(a)=\tge0(b)=y. \]
Any $2$-graph has an \emph{underlying 1-graph} with $C_0$ as vertices and
$C_1$ as edges.

\begin{example}
  The 2-graph~$C$ with
  \begin{align*}
    C_0&=\set{x,y,z},\\
    C_1&=\set{a:x\to y,b_1:y\to z,b_2:y\to z,b_3:y\to z},\\
    C_2&=\set{\alpha:b_1\To b_2,\beta:b_2\To b_3},
  \end{align*}
  can be depicted as
  \[
    \xymatrix@R=1ex{
      x\ar[r]^a&y\ar@/^3ex/[r]^{b_1}\ar@/^1.5ex/@{{}{ }{}}[r]|{\phantom\alpha\Downarrow\alpha}\ar[r]|{b_2}\ar@/_1.5ex/@{{}{ }{}}[r]|{\phantom\beta\Downarrow\beta}\ar@/_3ex/[r]_{b_3}&z\pbox.
    }
  \]
\end{example}

\subsection{2-categories}
\index{2-category}
\label{sec:2-cat}
A \emph{2-category} consists in a 2-graph~$C$ together with
\begin{itemize}
\item for each 0-cell~$x$ an \emph{identity} 1-cell
  \[
  \id_x
  :
  x\to x,
  \]
\item for each 1-cells
  $
  f:x\to y
  $ and $
  g:y\to z
  $,
  a \emph{horizontal composite} 1-cell
  \[
  f\comp0 g
  :
  x\to z,
  \]
\item for each 2-cells
  $
  \alpha:f\To f':x\to y
  $ and $
  \beta:g\To g':y\to z
  $,
  a \emph{horizontal composite} 2-cell
  \[
  \alpha\comp0\beta
  :
  (f\comp0 g)\To(f'\comp0 g')
  :
  x\to z,
  \]
\item for each 1-cell~$f:x\to y$, an \emph{identity} 2-cell
  \[
  \id_f
  :
  f\To f,
  \]
\item for each 2-cells
  $
  \alpha:f\To g:x\to y
  $ and $
  \beta:g\To h:x\to y
  $,
  a \emph{vertical composite} 2-cell
  \[
  \alpha\comp1\beta
  :
  f\To h
  :
  x\to y,
  \]
\end{itemize}
such that
\begin{itemize}
\item the compositions $\comp0$ and $\comp1$ are associative and admit
  identities as neutral elements,
\item the \emph{exchange law} holds: given 1-cells
  $
  f:x\to y
  $ and $
  g:y\to z
  $,
  one has
  \begin{equation}
    \label{eq:2-cat-xch-id}
    \id_f\comp0\id_g
    =
    \id_{f\comp0 g},
  \end{equation}
  and given 2-cells
  \begin{align*}
    \alpha&:f\To f':x\to y,
    &
    \alpha'&:f'\To f'':x\to y,
    \\
    \beta&:g\To g':y\to z,
    &
    \beta'&:g'\To g'':y\to z,
  \end{align*}
  we have
  \begin{equation}
    \label{eq:2-cat-xch}
    (\alpha\comp1\alpha')\comp0(\beta\comp1\beta')
    =
    (\alpha\comp0\beta)\comp1(\alpha'\comp0\beta')
    \pbox.
  \end{equation}
  Graphically, the equations~\eqref{eq:2-cat-xch-id} and~\eqref{eq:2-cat-xch}
  can be pictured as the ``commutation'' of the diagrams
  \[
    \fig{xch-id}
  \]
  and
  \[
    \fig{xch}
  \]
\end{itemize}

Note that every 2-category~$C$ has an underlying category, sometimes denoted
by~$\tcat1C$,
with $C_0$ as set of objects and $C_1$ as set of morphisms.

\subsection{Notation}
In the following, when considering morphisms in $2$-categories, we often omit
writing identities and horizontal compositions~$\comp0$, and write~$\comp{}$ for
vertical composition $\comp1$: for instance, we write
\[
  u\alpha w\comp{}\beta
  \qquad\text{instead of}\qquad
  (\id_u\comp0\alpha\comp0\id_w)\comp1\beta
  \pbox.
\]

\subsection{2-functors}
\index{2-functor}
A \emph{2-functor} $f:C\to D$ between 2-categories~$C$ and~$D$ consists of
functions
\[
  f_0:C_0\to D_0,
  \qquad\qquad
  f_1:C_1\to D_1,
  \qquad\qquad
  f_2:C_2\to D_2,
\]
which are
\begin{itemize}
\item compatible with sources and targets:
  \begin{align*}
    s_0\circ f_1&=f_0\circ s_0,
    &
    t_0\circ f_1&=f_0\circ t_0,
    \\
    s_1\circ f_2&=f_1\circ s_1,
    &
    t_1\circ f_2&=f_1\circ t_1,
  \end{align*}
\item compatible with compositions: for every composable pair of 1-cells $(u,v)$
  (\resp 2-cells $(\alpha,\beta)$),
  \begin{align*}
    f_1(u\comp0 v)&=f_1(u)\comp0f_1(v),
    &
    f_2(\alpha\comp0\beta)&=f_2(\alpha)\comp0f_2(\beta),
    \\
    &&
    f_2(\alpha\comp1\beta)&=f_2(\alpha)\comp1f_2(\beta),
  \end{align*}
\item compatible with identities: for every 0-cell $x$ and 1-cell $u$,
  \begin{align*}
    f_1(\unit{x})&=\unit{f_0(x)},
    &
    f_2(\unit{u})&=\unit{f_1(u)}
    \pbox.
  \end{align*}
\end{itemize}
In the following, we generally omit writing the subscript~$i$ from the
components~$f_i$ of a 2\nbd-func\-tor~$f$.

\subsection{Freely generated 2-categories}
\label{sec:free-2-cat}
\index{free!2-category}
The \emph{free 2-category}~$\freecat{P}$ generated by a 2-polygraph~$P$ is the
2-category~$C$ with $\freecat{\tpol1 P}$ as underlying category, \ie
\begin{itemize}
\item $C_0=P_0$,
\item $C_1=\freecat{P_1}$ with source and target given by
  $\sce0^C=\freecat{\sce0}$ and $\tge0^C=\freecat{\tge0}$,
\end{itemize}
whose $2$-cells are freely generated by~$P_2$, \ie $C_2$ is the smallest set
containing~$P_2$, together with source and target given by the functions $\sce1$
and $\tge1$ of~$P$, such that
\begin{itemize}
\item for every $\alpha,\beta\in C_2$ such that
  $\tge0\circ\tge1(\alpha)=\sce0\circ\sce1(\beta)$, there is an element
  $
    \alpha\comp0\beta
    \in
    C_2
  $,
\item for every $\alpha,\beta\in C_2$ such that $\tge1(\alpha)=\sce1(\beta)$,
  there is an element
  $
    \alpha\comp1\beta
    \in
    C_2
  $,
\end{itemize}
quotiented by the axioms required to form a 2\nbd-category, see \secr{2-cat}. We
write~$\freecat{P_2}$ for the set~$C_2$ of $2$-cells of~$\freecat{P}$.

Another construction for this $2$-category will be presented in
\secr{free-sesqui-2-cat}. It can be characterized by the following universal
property:

\begin{lemma}
  Suppose given a 2-polygraph~$P$, a 2-category~$C$, and functions
  \begin{align*}
    f_0:P_0&\to C_0,
    &
    f_1:P_1&\to C_1,
    &
    f_2:P_2&\to C_2,
  \end{align*}
  such that
  \begin{itemize}
  \item for every $1$-generator $a:x\to y$, we have $f_1(a):f_0(x)\to f_0(u)$,
  \item for every $2$-generator $\alpha:u\To v$, we have
    $f_2(\alpha):\freecat f_1(u)\To\freecat f_1(v)$,
  \end{itemize}
  where $\freecat f_1(a_1\ldots a_n)=f_1(a_1)\ldots f_1(a_n)$.
  Then there exists a 2-functor
  \[
    \freecat{f}:\freecat{P}\to C,
  \]
  which is unique up to isomorphism, such that
  \begin{itemize}
  \item for every $0$-generator $x\in P_0$, $\freecat f(x)=f_0(x)$,
  \item for every $1$-generator $a\in P_1$, $\freecat f(a)=f_1(a)$,
  \item for every $2$-generator $\alpha\in P_2$,
    $\freecat f(\alpha)=f_2(\alpha)$.
  \end{itemize}
\end{lemma}

\subsection{String diagrams}
\index{string diagram}
\index{diagram}
\label{sec:string-diag}
A convenient and intuitive notation for morphisms in free 2-categories is
provided by \emph{string diagrams}. Those were originally introduced by
Feynman~\cite{feynman1949space} and Penrose~\cite{penrose1971applications} in
physics, and formally studied by Joyal and Street~\cite{joyal1991geometry},
see~\cite{baez2009prehistory} for a detailed historical account
and~\cite{selinger2010survey} for a panorama of the possible variations.

Suppose fixed a 2-polygraph~$P$. A 2\nbd-generator
\[
  \alpha
  :
  a_1\ldots a_m
  \To
  b_1\ldots b_n,
\]
where the $a_i:x_{i-1}\to x_i$ and $b_i:y_{i-1}\to y_i$ are 1-generators, can be
thought of as some kind of device with $m$ inputs, with $a_i$ as types, and $n$
outputs, with $b_i$ as types. This suggests that, instead of using the usual
depiction
\[
  \vcenter{
    \xymatrix@C=3ex@R=3ex{
      &x_1\ar[r]^{a_2}&x_2\ar@{.}[r]&x_{m-1}\ar@/^/[dr]^{a_m}\\
      x_0\ar@/^/[ur]^{a_1}\ar@/_/[dr]_{b_1}&&{\phantom\alpha}\Downarrow\alpha&&y_n\\
      &y_1\ar[r]_{b_2}&y_2\ar@{.}[r]&y_{n-1}\ar@/_/[ur]_{b_n}
    }
  }
\]
we can use the alternative graphical notation
\[
  \begin{tikzpicture}[baseline=(current bounding box.center),yscale=-1,scale=0.5,every path/.style={join=round,cap=round},scale=2.5]
    \draw (0.000000,0.700000) node {$\scriptstyle a_1$};
    \draw (1.000000,0.700000) node {$\scriptstyle a_2$};
    \draw (2.000000,0.700000) node {$\scriptstyle a_m$};
    \draw (1.500000,1.000000) node {$\scriptstyle \ldots$};
    \draw (0.000000,0.900000) -- (0.000000,1.200000);
    \draw (1.000000,0.900000) -- (1.000000,1.200000);
    \draw (2.000000,0.900000) -- (2.000000,1.200000);
    \draw (0.000000,1.800000) -- (0.000000,2.100000);
    \draw (1.000000,1.800000) -- (1.000000,2.100000);
    \draw (2.000000,1.800000) -- (2.000000,2.100000);
    \filldraw[fill=white,rounded corners=1pt] (-0.300000,1.200000) -- (2.300000,1.200000) -- (2.300000,1.800000) -- (-0.300000,1.800000) -- cycle;
    \draw (1.000000,1.500000) node {$\alpha $};
    \draw (1.500000,2.000000) node {$\scriptstyle \ldots$};
    \draw (0.000000,2.300000) node {$\scriptstyle b_1$};
    \draw (1.000000,2.300000) node {$\scriptstyle b_2$};
    \draw (2.000000,2.300000) node {$\scriptstyle b_n$};
    \draw (-.5,1.5) node {$\scriptstyle x_0$};
    \draw (0.5,1) node {$\scriptstyle x_1$};
    \draw (0.5,2) node {$\scriptstyle y_1$};
    \draw (2.5,1.5) node {$\scriptstyle y_n$\pbox.};
  \end{tikzpicture}
\]
This notation, called a \emph{string diagram} because the 1-cells become some
kind of strings, is inspired of electric circuits where electronic components
are linked with conductive wires. Namely, a 2-cell will more generally be
depicted as various gates linked together, for instance:
\[
  \satex{sd-ex}\pbox.
\]
This notation is ``dual'' (in the sense of Poincaré duality) to the traditional
one: 0-cells are pictured as 2-dimensional region of the plane, 1-cells are
pictured as wires and 2-cells are pictured as points (the above rectangles, which
are here in order to be able to put labels inside them, should be pictured as
points if we were drawing string diagrams by the book).

String diagrams can be composed in the various ways expected in a
2\nbd-cate\-gory, as we now describe.
Horizontal composition of 2-cells of $\freecat{P}$
\[
  \phi:a_1\ldots a_m\To a'_1\ldots a'_{m'}
  \qquad\text{and}\qquad
  \psi:b_1\ldots b_n\To b'_1\ldots b'_{n'}
\]
amounts to horizontal juxtaposition of the corresponding diagrams, as shown on the left
\begin{align*}
  \satex{sd-hcomp}
  &&
  \satex{sd-vcomp}
\end{align*}
and vertical composition of
\[
  \phi:a_1\ldots a_m\To b_1\ldots b_k
  \qquad\text{and}\qquad
  \psi:b_1\ldots b_k\To c_1\ldots c_n  
\]
is obtained by vertically juxtaposing the corresponding diagrams and linking the
wires, as shown on the right above, identities simply being wires. These
diagrams are to be considered up to ``planar isotopy'', by which we mean
continuous deformations, fixing boundaries, preserving direction, and forbidding
wires to cross. For instance, the three diagrams below are considered to be equal:
\[
  \satex{sd-xch-l}
  =
  \satex{sd-xch-c}
  =
  \satex{sd-xch-r}
\]
This is detailed in~\cite{joyal1991geometry}, where it is proved that the axioms
of 2-categories are satisfied and the following universal property is satisfied:

\begin{theorem}
  String diagrams (up to planar isotopy) are the 2-cells of a 2-category which
  is the free 2-category $\freecat{P}$ on the 2-polygraph~$P$.
\end{theorem}

\index{string diagram!empty}
Given a $0$-cell~$x$, the string diagram corresponding to the
$2$-cell~$\id_{\id_x}$ (the identity on the identity on~$x$) is the \emph{empty}
one. Since this can be confusing, we generally use the diagram
\[
  \satex{empty}
\]
in order to represent it.

\subsection{Strict monoidal categories}
\index{monoidal category!strict}
\index{category!monoidal!strict}
\label{sec:strict-moncat}
\label{sec:PRO}
We introduce here useful categorical structures, which can be seen as particular
cases of $2$-categories.

A \emph{strict monoidal category} $(C,\otimes,\monunit)$ consists of a
category~$C$ together with a \emph{tensor product} functor
$\otimes:C\times C\to C$ and a \emph{unit} object $\monunit\in C$, such that the
product is associative and admits $\monunit$ as unit, which means that for any
objects $x,y,z\in C$, we should have
\[
  (x\otimes y)\otimes z=x\otimes(y\otimes z)
  \qquad\text{and}\qquad
  i\otimes x=x=x\otimes i.
\]
Any $2$-category with only one $0$-cell induces a strict monoidal category
with~$C_1$ as set of objects, $C_2$ as set of morphisms, with composition being given
by $\comp1$ and tensor product by~$\comp0$, and this extends to an isomorphism
between $2$-categories with one fixed $0$-cell and strict monoidal categories:
in this way strict monoidal categories can be seen as particular cases of
$2$-categories (and we will often implicitly make use of this identification in
the following).

Among strict monoidal categories, the following class is a source of
particularly useful examples (see \cref{chap:3pol,chap:3ex}). A
\emph{PRO}\index{PRO} is a strict monoidal category whose set of objects is $\N$
and such that the tensor product of two objects $m$ and $n$ is given by addition
$m\otimes n=m+n$. The terminology is an abbreviation of ``PROduct
category''~\cite{maclane1965categorical}.

\subsection{Bicategories}
\label{sec:bicat}
The axioms for 2-categories are sometimes too strong for situations encountered
in practice: it happens in many situations that composition of 1-cells is not
strictly associative, but only associative up to a coherent invertible 2-cell,
and similarly identity 1-cells are not strictly neutral elements for
composition. This motivates the introduction of the following generalization of
the notion of 2-category.

\index{bicategory}
A \emph{bicategory}~$C$ consists in a 2-graph~$C$ together
with identity 1-cells and 2\nbd-cells, horizontal composition of 1- and 2-cells,
vertical composition of 2\nbd-cells, as in \secr{2-cat}, and moreover natural
families of invertible 2-cells
\begin{align*}
  \alpha_{f,g,h}:(f\comp0 g)\comp0 h&\To f\comp0(g\comp0 h),
  &
  \lambda_f:\unit{x}\comp0 f&\To f,
  \\
  &&
  \rho_f:f\comp0\unit{y}&\To f,  
\end{align*}
indexed by 1-cells $f:x\to y$, $g:y\to z$, and $h:z\to w$, respectively called
\emph{associator} and \emph{left} and \emph{right unitor}, such that
\begin{itemize}
\item vertical composition $\comp1$ is associative on 2-cells and admits
  identities as neutral elements,
\item the exchange law between horizontal and vertical composition holds, and
\item the two following coherence axioms hold for every composable 1-cells $f$,
  $g$, $h$ and $i$:
  \[
    \fig{bicat-pentagon}
  \]
  \[
    \fig{bicat-triangle}
  \]
\end{itemize}

Obviously, any 2-category can be seen as a bicategory where $\alpha_{f,g,h}$,
$\lambda_f$ and~$\rho_f$ are all identity 2-cells. A typical non-trivial example
of bicategory, is the bicategory~$\Span(C)$ of spans in a category~$C$ with
pullbacks, see \secr{span}.

\index{monoidal category}
\index{category!monoidal}
An important particular case of the previous construction is the following
one. We have seen in \cref{sec:strict-moncat} that a $2$-category with one
$0$-cell corresponds to a strict monoidal category. Similarly, a bicategory with
one $2$-cell corresponds to the notion of \emph{monoidal category}, which
generalizes the one of strict monoidal category. The definition is recalled in
\cref{sec:mon-cat}, and more details can be found
in~\cite[Chapter~VII]{MacLane98}.

A fundamental theorem for bicategories is the coherence theorem which ensures that it
can always be ``replaced'' by a 2-category:

\begin{theorem}
  Every bicategory is biequivalent to a 2-category.
\end{theorem}

\noindent
This result is shown for instance in~\cite{power1989general} and a proof of
coherence for monoidal categories, which are equivalent to bicategories with
only one 0-cell, is given in~\cref{sec:mon-coh}.

\section{Coherent Confluence of 1-Polygraphs}
\label{sec:free-21-category}
In \secr{free-1-gpd}, we have seen that the morphisms in the free groupoid
generated by a $1$-polygraph could be seen as representatives of the congruence
associated to the polygraph. A similar construction can be performed for
$2$-polygraphs as follows, from which we will be able to define a notion of
coherent presentation. Here, rather than presenting a category (as in
\secr{pres-cat}) or generating a 2\nbd-cate\-gory (as in \secr{generating-2-cat}), we
think of a 2-polygraph as a presentation of a set by the underlying 1-polygraph
(see \secr{pres-set}) together with additional 2-dimensional coherence data
provided by the 2-generators.

\subsection{Freely generated (2,1)-category}
\label{sec:free-21-cat}
\index{21-category@$(2,1)$-category}
\index{category!21-@$(2,1)$-}
\index{free!21-@$(2,1)$-category}
A \emph{(2,1)-category}~$C$ is a $2$-category in which every $2$-cell is
invertible, that is, for every $2$-cell $\alpha:u\To v$ in $C$ there exists a $2$-cell
$\alpha^-:v\To u$ in $C$ such that
$
\alpha^-\comp1\alpha=\id_u,
$ and $
\alpha\comp1\alpha^-=\id_v  
$.
  
Given a 2-polygraph~$P$, with underlying 1\nbd-polygraph~$\tpol1 P$, one can
generate the \emph{free (2,1)-category}~$\freegpd{P}$ with $\freecat{\tpol1P}$
as underlying category and containing the $2$\nbd-generators in~$P_2$ as
$2$-cells. This construction can be performed in a similar way as in
\secr{free-2-cat}. 
\begin{itemize}
\item The set of 0-cells is the set~$P_0$.
\item The set of 1-cells is~$\freecat{P_1}$. 
\item The set of 2-cells~$\freegpd P_2$ is the set of formal
horizontal and vertical composites of elements of~$P_2$ and identities of
elements of~$\freecat{P_1}$, and their inverses, quotiented by the axioms of
(2,1)-categories.
\end{itemize} 

The following lemma motivates why one can think of 2-cells in
this 2\nbd-cate\-gory as witnesses of equivalence of rewriting paths:

\begin{lemma}
  \label{lem:2pol-2cell-cong}
  Given a 2-polygraph~$P$, two 1-cells $f,g:x\to y$ are $P$\nbd-congruent if and
  only if there exists a 2-cell $\phi:f\To g$ in $\freegpd{P}$.
\end{lemma}

\noindent
For any $2$-polygraph~$P$, there is thus a canonical $2$-functor
$q_P:\freegpd P\to\pcat P$ from the generated $(2,1)$-category to the presented
category, which is the identity on $0$-cells, sends every $1$-cell to its
equivalence class under the $P$-congruence and sends every $2$-cell to the
identity.

Any two presentations of a given category generate equivalent $(2,1)$-cate\-gories
in the following sense:

\begin{lemma}
  \label{lem:21TietzeEquivalence}
  Suppose given two 2-polygraphs~$P$ and~$Q$ both presenting a same
  category~$C$. There exist two $2$-functors
  \begin{align*}
    f&:\freegpd{P}\to\freegpd{Q},
    &
    g&:\freegpd{Q}\to\freegpd{P},    
  \end{align*}
  between the free $(2,1)$-categories generated by $P$ and $Q$ and, for every
  $1$-cells $u$ of $\freegpd{P}$ and $v$ of $\freegpd{Q}$, there exist $2$-cells
  \begin{align*}
    \phi_u&: gf(u) \To u,
    &
    \psi_v&: fg(v) \To v,
  \end{align*}
  in $\freegpd{P}$ and $\freegpd{Q}$, such that the two following conditions are
  satisfied.
  \begin{enumerate}
  \item The $2$-functors $f$ and $g$ induce the identity through the canonical
    projections~$q_P$ and $q_Q$ onto $C$:
    \begin{align*}
      \svxym{
        {\freegpd{P}}
	\ar [r]^{f}
	\ar [d] _-{q_P}
        & {\freegpd{Q}}
	\ar [d] ^-{q_Q} 
        \\
        C
	\ar [r]_{1_C}
        & C,
      }
      &&
      \svxym{
        {\freegpd{Q}}
	\ar [r]^{g}
	\ar [d]_{q_Q}
        & {\freegpd{P}}
        \ar [d]^{q_P}
        \\
        C
	\ar [r]_{1_C}
        & C\pbox.
      }
    \end{align*}
  \item\label{cond:icons} The $2$-cells $\phi_u$ and $\psi_v$ are functorial in $u$ and $v$, that
    is
    \begin{align*}
      \phi_{uu'}&= \phi_u\comp0\phi_{u'},
      &
      \phi_{\unit{x}}&= \unit{\unit{x}},     
    \end{align*}
    for any $0$-composable $1$-cells $u$ and $u'$ and $0$-cell $x$ and
    \begin{align*}
      \psi_{vv'} &= \psi_v\comp0\psi_{v'},
      &
      \psi_{\unit{y}} &= \unit{\unit{y}}, 
    \end{align*}
    for any $0$-composable $1$-cells $v$ and $v'$ and $0$-cell $y$.
  \end{enumerate}
\end{lemma}
\begin{proof}
  We construct the $2$-functor $f$, the case of the $2$-functor $g$ being similar. For a
  $0$-cell~$x$, we set $f(x)=q_{Q}^{-1}q_{P}(x)$. If $a:x\to y$ is a $1$-generator of $P$, we
  choose, in an arbitrary way, a $1$-cell $f(a):f(x)\to f(y)$ in $\freegpd{Q}$ such
  that $q_{Q}f(a)=q_{P}(a)$.
  Then, we extend $f$ to every $1$-cell of
  $\freegpd{P}$ by functoriality. Let $\alpha:u\To u'$ be a $2$\nbd-gene\-rator
  of $P$. Since $P$ is a presentation of $C$, we have $q_{P}(u)=q_{P}(u')$, so
  that $q_{Q}f(u)=q_{Q}f(u')$ holds. Using the fact that $Q$ is a presentation
  of~$C$, we arbitrarily choose a $2$-cell $f(\alpha):f(u)\To f(u')$ in
  $\freegpd{Q}$. Then, we extend $f$ to every $2$-cell of $\freegpd{P}$ by
  functoriality.

  Now, let us define the $2$-cells $\phi_u$, the case of $2$-cells $\psi_v$
  being similar. Let $a$ be a $1$-generator of~$P$. By construction of the $2$-functors $f$
  and~$g$, we have:
  \[
    q_{P} gf(a) = q_{Q} f(a) = q_{P}(a)
    \pbox.
  \]
  Since $P$ is a presentation of $C$, there exists a $2$-cell
  $\phi_a:gf(a)\To a$ in $\freegpd{P}$. We extend $\phi$ to every $1$-cell~$u$
  of~$\freegpd{P}$ by functoriality.
\end{proof}

\begin{remark}
  Condition~\ref{cond:icons} of the lemma exactly expresses that $\phi$ and
  $\psi$ are \emph{icons}\index{icon} (which is the acronym for ``identity
  component oplax natural transformations'', \ie oplax natural transformations
  with whose 1-cell components are identities) in the sense
  of~\cite{lack2010icons}.
\end{remark}

\subsection{Freely generated 2-groupoid}
\index{2-groupoid@$2$-groupoid}
\index{groupoid!2-@$2$-}
\index{20-category@$(2,0)$-category}
\index{category!20-@$(2,0)$-}
\index{free!20-category@$(2,0)$-category}
A \emph{2-groupoid}, also called \emph{(2,0)-cate\-gory}, 
is a $2$-category where both 1-cells and 2-cells are
invertible. As a variant of the previous construction, a $2$-polygraph~$P$ freely
generates a $2$-groupoid, whose set of 0-cells is $P_0$, whose set of 1-cells is
$\freegpd{P_1}$ (the morphisms of the free 1-groupoid generated by the
underlying 1-polygraph, as described in \secr{free-1-gpd}), and whose set of
2-cells is the set of formal horizontal and vertical composites of elements
of~$P_2$ and identities of elements of~$\freegpd{P_1}$, and their inverses,
quotiented by the axioms of 2-groupoids. 

\subsection{Coherent 2-polygraphs}
\index{coherent!2-polygraph@$2$-polygraph}
\index{polygraph!coherent}
Given a $2$-polygraph~$P$, the 1-cells $\freegpd{P}_1$ of the freely generated
2-groupoid may be thought of as witnesses for equivalence of 0\nbd-cells (see
\secr{free-1-gpd}) and the 2-cells $\freegpd{P}_2$ as witnesses of
``equivalences between equivalences'' (as a variant of
\lemr{2pol-2cell-cong}). A 2-polygraph is said to be coherent when equivalences
do not bring essential information, in the sense that there is at most one
equivalence between two given 0-cells, up to equivalence between equivalences.
Formally, a 2-polygraph~$P$ is \emph{coherent} when for every pair of $0$-cells
$x,y\in P_0$ and pair of 1-cells $f,g:x\to y$ in $\freegpd{P}_1$ there is a
$2$-cell $\alpha:f\To g$ in $\freegpd{P}_2$ in the free $2$-groupoid generated
by~$P$.

The notion of coherence can be generalized in the expected way to variants of
$2$-polygraphs where the source and target of $2$-generators are cells in
$\freegpd{P}_2$ (as opposed to $\freecat{P}_2$): those are called
$(2,0)$-polygraphs and formally defined in \cref{sec:np-polygraph}.
\index{cellular extension}
In this context, following the terminology introduced in
\secr{acyclic-extension}, we can also say that $P_2$ is a \emph{cellular
  extension} of the category~$\freegpd P_1$: the 2-polygraph~$P$ is then
coherent precisely when this extension is \emph{acyclic}.

\subsection{Coherent confluence}
\index{coherent!confluence}
\index{confluence!coherent}
Our aim here is to provide techniques to show the coherence of a 2-polygraph by
adapting the rewriting concepts presented in \chapr{lowdim} in order to take
2-dimensional information into account, as an explicit witness of the
commutation of diagrams. This approach first appeared under the terminology of
\emph{commuting diagrams}~\cite[Section 4.3]{Huet87}.

A 2-polygraph~$P$ is \emph{coherently confluent} when for every pair of
coinitial 1\nbd-cells $f:x\to y$ and $f':x\to y'$ in $\freecat{P_1}$ there is a pair
of cofinal 1-cells $g:y\to z$ and $g':y'\to z$ in $\freecat{P}_1$ and a $2$-cell
$\alpha:fg\To f'g'$ in $\freegpd{P}_2$ (in the free $(2,1)$-category it
generates):
\begin{equation}
  \label{eq:coh-confl}
  \svxym{
    &\ar[dl]_fx\ar[dr]^{f'}&\\
    y\ar@{.>}[dr]_g&\overset\alpha\To&\ar@{.>}[dl]^{g'}y'\pbox.\\
    &z&
  }
\end{equation}
In such a situation, we also say that the branching $(f,f')$ is \emph{coherently
joinable}.
Similarly,~$P$ is \emph{locally coherently confluent} when for every pair of
coinitial 1\nbd-gene\-ra\-tors $a:x\to y$ and $b:x\to y'$ in $P_1$ there is a pair of
cofinal 1-cells $g:y\to z$ and $g':y'\to z$ in $\freecat{P}_1$ and a $2$-cell
$\alpha:ag\To bg'$ in $\freegpd{P}_2$:
\[
  \svxym{
    &\ar[dl]_ax\ar[dr]^b&\\
    y\ar@{.>}[dr]_g&\overset\alpha\To&\ar@{.>}[dl]^{g'}y'\pbox.\\
    &z&
  }
\]

We say that a 2-polygraph is \emph{terminating} when the underlying 1-polygraph
is, in the sense of \secr{1-termination}. We can now adapt Newman's lemma (see
\lemr{newman}) to this extended notion of confluence as follows:

\begin{lemma}
  \index{Newman's lemma!coherent}
  \label{L:HNewman-1}
  A terminating 2-polygraph is coherently confluent if and only if it is locally
  coherently confluent.
\end{lemma}
\begin{proof}
  The left-to-right implication is immediate; let us show the right-to-left
  implication. Supposing that the terminating polygraph~$P$ is locally coherently
  confluent, we show that it is coherently confluent at~$x$ (\ie that every
  branching at~$x$ is coherently joinable) by well-founded induction
  on~$x$. Suppose given two coinitial 1-cells $f:x\to y$ and $f':x\to y'$. If one
  of them is the identity, say $f'$, we can close the diagram as in
  \eqref{eq:coh-confl} with $g=\id_y$, $g'=f$, and $\alpha=\id_f$. Otherwise, we
  have $f=af_1$ and $f'=a'f_1'$ for some 1-generators $a$ and $a'$ and 1-cells
  $f$ and~$f'$:
  \[
    \xymatrix@!C=1ex@!R=1ex{
      &&\ar[dl]_ax\ar[dr]^{a'}&&\\
      &\ar[dl]_{f_1}y_1\ar@{.>}[dr]|{f_2}&\overset\alpha\To&\ar@{.>}[dl]|{f_2'}\ar[dr]^{f_1'}y_1'\\
      y\ar@{.>}[dr]_{g_1}&\overset\beta\To&y''\ar@{.>}[dl]|{f_2''}&\overset\gamma\To&\ar@{.>}[ddll]^{g'}y'\pbox.\\
      &z'\ar@{.>}[dr]_{g_2}&&&\\
      &&z
    }
  \]
  By local confluence we deduce the existence of 1-cells $f_2$ and $f_2$ and a
  2-cell $\alpha:af_2\To af_2'$ as above. By induction hypothesis, we obtain the
  existence of 1\nbd-cells $g_1$ and $f_2''$ together with a 2-cell
  $\beta:f_1g_1\To f_2f_2''$. By induction hypothesis again, we deduce the
  existence of 1-cells $g_2$ and $g'$ and a 2-cell
  $\gamma:f_2'f_2''g_2\To f_1'g'$. We have shown the existence of cofinal
  $1$-cells $g_1g_2:y\to z$ and $g':y'\to z$, and of a $2$-cell
  \[
    a\beta g_2\comp{}\alpha f_2''g_2\comp{}a'\gamma
    :
    af_1g_1g_2\To a'f_1'g'
  \]
  and the branching $(af_1,a'f_1')$ is thus coherently joinable.
\end{proof}

\begin{proposition}
  \label{prop:PreSquierOnePolygraphs}
  A terminating and coherently confluent 2-polygraph~$P$ is coherent.
\end{proposition}
\begin{proof}
  Since the polygraph is terminating, for every 0-cell~$x$ there is a normal
  form $\nf{x}$ and a rewriting path
  $
  n_x:x\to\nf x
  $
  in $\freecat{P_1}$. Given a 1-cell $f:x\to y$ in~$\freecat P_1$, since the
  polygraph is convergent we have $\nf{x}=\nf{y}$ and, since the branching
 $(n_x,fn_y)$ is coherently joinable there is a 2-cell~$\alpha_f:fn_y\To n_x$,
  as on the left:
  \begin{align*}
    \svxym{
      x\ar[dr]_{n_x}\ar[rr]^f&\ar@{}[d]|{\overset{\alpha_f}\Leftarrow}&\ar[dl]^{n_y}y\\
      &\nf{x}
    }
    &&
    \svxym{
      y\ar[dr]_{n_y}\ar[rr]^{f^-}&\ar@{}[d]|{\overset{\alpha_{f^-}}\Leftarrow}&\ar[dl]^{n_x}x\pbox.\\
      &\nf{x}
    }    
  \end{align*}
  In the free 2-groupoid generated by~$P$, we also define the 2-cell
  $\alpha_{f^-}=f^-\comp0\alpha_f^-$, which can be pictured as on the right
  above. Now consider a 1-cell $f:x\to y$ in~$\freegpd{P_1}$. It factors as
  \[
    f=f_1^-g_1f_2^-g_2\ldots f_k^-g_k,
  \]
  for some suitably composable morphisms $f_i$ and $g_i$ in $\freecat{P_1}$. We
  then define $\alpha_f:fn_y\To n_x$ as the composite
  \[
    \xymatrix{
      x\ar[d]_{n_x}\ar[r]^{f_1^-}\ar@{}[dr]|{\overset{\alpha_{f_1^-}}\Leftarrow}&
      y_1\ar@{}[dr]|{\overset{\alpha_{g_1}}\Leftarrow}\ar[d]|{n_{y_1}}\ar[r]^{g_1}&
      x_2\ar[d]|{n_{x_2}}\ar[r]&\cdots\ar[r]&
      x_k\ar@{}[dr]|{\overset{\alpha_{f_k^-}}\Leftarrow}\ar[d]|{n_{x_k}}\ar[r]^{f_k^-}&
      y_k\ar@{}[dr]|{\overset{\alpha_{g_k}}\Leftarrow}\ar[d]|{n_{y_k}}\ar[r]^{g_k}&y\ar[d]^{n_y}\\
      \nf x\ar@{=}[r]&\nf x\ar@{=}[r]&\nf x\ar@{=}[r]&\cdots\ar@{=}[r]&\nf x\ar@{=}[r]&\nf x\ar@{=}[r]&\nf x\pbox.
    }
  \]
  Finally, for any pair of parallel 1-cells $f,g:x\to y$ in $\freegpd{P_1}$, the
  composite 2-cell
  \[
    \xymatrix@!C=1ex@!R=1ex{
      &\ar[dl]_fx\ar[dd]|{n_x}\ar[dr]^g&\\
      y\ar@/_3ex/[ddr]_{\id_f}\ar[dr]_{n_y}&\ar@{}[l]|{\overset{\alpha_f^-}\Leftarrow}\ar@{}[r]|{\overset{\alpha_g}\Leftarrow}&\ar[dl]^{n_y}y\ar@/^3ex/[ddl]^{\id_y}\\
      &\nf x\ar[d]|{n_y^-}\\
      &y
    }
  \]
  has type $f\To g$. We thus conclude that the 2-polygraph~$P$ is coherent.
\end{proof}

\noindent
Composing \cref{L:HNewman-1} and
\cref{prop:PreSquierOnePolygraphs}, we obtain the \emph{coherent Squier
theorem} for $1$-polygraphs~\cite[Theorem~5.2]{squier1994finiteness}:

\begin{theorem}
  \index{Squier!theorem!for 1-polygraphs@for $1$-polygraphs}
  \label{thm:SquierHomotopicalOnePolygraphs}
  Let~$P$ be a terminating $2$-poly\-graph. If, for every pair of coinitial
  1-generators $a:x\to y_1$ and $b:x\to y_2$ in $P_1$, there is a pair of
  cofinal $2$-cells $f:y_1\to z$ and $g:y_2\to z$ in $\freecat P_1$ and a
  2-cell $\phi:af\To bg$ in~$\freegpd{P_2}$,
  \[
    \svxym{
      &\ar[dl]_ax\ar[dr]^b&\\
      y_1\ar[dr]_f&\overset\alpha\To&\ar[dl]^gy_2\\
      &z&
    }
  \]
  then $P$ is coherent.
\end{theorem}

\noindent
This theorem is extended to $2$-polygraphs in \chapr{2-Coherent} and to
higher-dimensional polygraphs in \cref{Chapter:ConstructingResolutions}:
we will see that it provides us with a canonical way of extending a convergent
presentation into a coherent one using the homotopical completion procedure, see
\cref{Section:HomotopicalCompletion}.


\chapter{Operations on Presentations}
\label{chap:2-op}
\label{chap:2op}
The usefulness and richness of 2-polygraphs is confirmed by the large number and
variety of categories they present. Some examples of presentations were given in
\cref{sec:pres-cat} and many more are described in \cref{chap:2ex}. In order to
show that a given polygraph is a presentation of a given category, one can
either tackle the issue directly, by using the rewriting tools of
\cref{chap:2rewr}, or take a modular approach, by combining already known
presentations: this is the route taken in the present chapter.

Three significant applications are given. We first address, in
\cref{sec:2-limits}, the presentation of limits and colimits by means
of given presentations of the base categories, and precisely show how to
systematically build presentations of  products, coproducts, and
pushouts.
Next, in \cref{sec:2-pres-loc}, we show how to add
formal inverses to some morphisms of a category at the level of presentations.
Finally, \cref{sec:2-dlaws} is about distributive
laws, in relation to factorization systems on categories. We
introduce a notion of composition along a distributive law between two small categories
sharing the same set of objects and show how to derive a presentation
of this composite from presentations of the components.

\section{Limits and Colimits of Presented Categories}
\label{sec:2-limits}

\subsection{Initial category}
\index{initial category}
The initial category (with no object nor morphism) admits the presentation $\Pres{}{}{}$.

\subsection{Coproducts of categories}
\index{coproduct}
Given two categories~$C$ and~$D$ respectively presented by 2-polygraphs~$P$
and~$Q$, their coproduct~$C\sqcup D$ is presented by the 2-polygraph with
$P_0\sqcup Q_0$ as 0\nbd-gene\-rators, $P_1\sqcup Q_1$ as 1-generators, and
$P_2\sqcup Q_2$ as 2-generators,
with expected source and target maps (this is the coproduct of~$P$ and~$Q$ in
the category of 2-polygraphs).

\subsection{Coproducts with fixed objects}
Given a fixed set~$O$, consider the subcategory~$\Cat_O$ of~$\Cat$ whose objects
are categories with~$O$ as objects, and whose morphisms are functors which are
identities on objects. Given two categories~$C$ and~$D$ in~$\Cat_O$,
respectively presented by 2-polygraphs~$P$ and~$Q$, their coproduct in~$\Cat_O$
is presented by the 2-polygraph with $P_0=Q_0=O$ as 0\nbd-gene\-rators,
$P_1\sqcup Q_1$ as 1-generators and $P_2\sqcup Q_2$ as $2$-generators, with
expected source and target maps. In particular, when $O=\set{\star}$, the
category~$\Cat_O$ is isomorphic to the category of monoids and their coproduct
is the \emph{free product}.

\subsection{Coequalizers}
\index{coequalizer}
Suppose given two categories~$C$ and~$D$ respectively presented by
2-polygraphs~$P$ and~$Q$ and two functors
\[
  f,g : C \to D \pbox.
\]
Reformulating the results of~\cite{bednarczyk1999generalized}, the category
obtained as their coequalizer is presented by the following $2$-polygraph~$R$. We
write $\sim$ for the smallest equivalence relation on $Q_0$ such that
$f(x)\sim g(x)$ for every 0-generator $x\in P_0$, and~$[y]$ for the equivalence
class of a 1-cell $y\in Q_1$.
The coequalizer of~$f$ and~$g$ is the category presented by the polygraph~$R$
such that
\begin{itemize}
\item the 0-generators are the equivalence classes $[x]$ of 0-generators
  in~$Q_0$,
\item the 1-generators are of the form $a:[x]\to[y]$ for 1-generators $a:x\to y$
  in~$Q_1$, and
\item the 2-generators are either of the form
  \begin{itemize}
  \item $\alpha:u\To v$ for 2-generators $\alpha:u\To v$ in~$Q_2$, or
  \item $\alpha_a:u\To v$ for 1-generators $a:x\to y$ in~$P_1$ such that
    $f(a)=\cl u$ and $g(a)=\cl v$.
  \end{itemize}
\end{itemize}

This construction allows to quotient a presented category both on objects and
morphisms. Note that the relations in~$P_2$ are not used in this construction:
we only need a generating graph for~$C$.

\subsection{Pushouts}
\index{pushout}
All finite colimits of presented categories can be constructed from coproducts
and coequalizers~\cite[Section~V.2]{MacLane98}. For instance, the
pushout of the diagram of categories on the left
\begin{align*}
  \xymatrix{
    C&\ar[l]_fB\ar[r]^g&D,
  }
  &&
  \xymatrix{
    B\ar@<+.5ex>[r]^-{i_C\circ f}\ar@<-.5ex>[r]_-{i_D\circ g}&C+D,
  }  
\end{align*}
can be computed as the coequalizer of the diagram on the right, where the
morphisms $i_C:C\to C+D$ and $i_D:D\to C+D$ are the canonical inclusions. In
particular, when both~$f$ and~$g$ are inclusions of polygraphs, the pushout is
simply given by (non-disjoint) union on sets of $n$-generators for $n=0,1,2$.

\subsection{Terminal category}
\index{terminal category}
We now turn to limits. The most simple example is the terminal category which
admits the presentation $\Pres\star{}{}$.

\subsection{Products}
\index{product}
\label{sec:2-pres-prod}
Suppose given two categories~$C$ and~$D$ along with respective presentations by
2-polygraphs~$P$ and~$Q$. A presentation~$R$ of the product category $C\times D$
is given by the 2-polygraph~$R$ with
\begin{itemize}
\item $R_0=P_0\times Q_0$ as set of 0-generators,
\item $R_1=P_1\times Q_0\uplus P_0\times Q_1$ as set of 1-generators with
  \begin{align*}
  (a,y):(x,y)&\to(x',y),
  &
  (x,b):(x,y)&\to(x,y'),
  \end{align*}
  with $a:x\to x'$ in $P_1$ and $y$ in $Q_0$ (\resp $x$ in $P_0$ and
  $b:y\to y'$ in~$Q_1$),
\item $R_2=P_2\times Q_0+P_1\times Q_1+P_0\times Q_2$ as set of 2-generators: a
  2-generator is either
  \[
  (\alpha,y)
  :
  (u,y)\To(u',y)
  :
  (x,y)\to(x',y),
  \]
  with $\alpha:u\To u':x\to x'$ in~$P_2$ and $y\in Q_0$, or
  \[
  (a,b)
  :
  (x,b)(a,y')
  \To
  (a,y)(x',b)
  :
  (x,y)
  \to
  (x',y'),
  \]
  with $a:x\to x'$ in $P_1$ and $b:y\to y'$ in $Q_1$, or
  \[
  (x,\beta)
  :
  (x,v)\To(x,v')
  :
  (x,y)\to(x,y'),
  \]
  with $x\in P_0$ and $\beta:v\To v':y\to y'$ in~$Q_2$.
\end{itemize}
Above, given $u=a_1a_2\ldots a_n$ in $\freecat{P_1}$ (with the $a_i$ being
generators in $P_1$) and $y\in Q_0$, the $1$-cell $(u,y)$ is a notation for
$(u,y)=(a_1,y)(a_2,y)\ldots(a_n,y)$ and the notation $(x,v)$, for $x\in P_0$ and
$v\in\freecat{P_1}$, is similar.

In particular, when $C$ and $D$ are both monoids (or groups), their product in
the above sense is often called their \emph{direct product}.

\section{Localizations of Presented Categories}
\label{sec:2-pres-loc}

\subsection{Free groupoid}
\index{free!groupoid}
\label{sec:1-pres-free-gpd}
The forgetful functor $\Gpd\to\Cat$ witnessing for the fact that a groupoid is a
particular category (with invertible morphisms) admits a left adjoint,
constructing the \emph{free groupoid}~$\freegpd{C}$ (also called the
\emph{enveloping groupoid}) over a category~$C$.
Given a 2-polygraph~$P$ presenting~$C$, the groupoid~$\freegpd{C}$ admits a
presentation (as a category) by the 2-polygraph~$Q$ with
\begin{itemize}
\item $Q_0=P_0$ as set of 0-generators,
\item $Q_1=P_1\sqcup \finv{P_1}$ as set of 1-generators, with
  \[
    \finv{P_1}=\setof{\finv{a}:y\to x}{a:x\to y\in W},
  \]
\item $Q_2=P_2\uplus I_2$ as set of 2-generators with
  \[
  I_2=\setof{a\finv a\To\id_x,\finv aa\To\id_y}{a:x\to y\in W},
  \]
\end{itemize}
where $W=P_1$.

A morphism in~$\freecat{Q_1}$ is \emph{reduced} when it is not of the form
$ua\finv av$ or $u\finv aav$ for some $u,v\in\freecat{Q_1}$ and $a\in P_1$. The
equivalence classes of elements of~$\freecat{Q_1}$ modulo the congruence
generated by $I_2$ contain exactly one reduced morphism, which is often
convenient to choose as a canonical representative; this is detailed in
\exr{finv-confl}.

\subsection{Localization}
\index{localization}
\label{sec:pres-localization}
As a generalization of the previous construction,
given a category~$C$ and a class~$W$ of morphisms of~$C$, we can consider the
category~$\loc CW$, called the \emph{localization} of~$C$ by~$W$, obtained by
formally inverting the morphisms of~$W$, see \cref{sec:localization} for a proper
definition.
Given a category~$C$ presented by a 2-polygraph~$P$ and a set~$W\subseteq P_1$
of $1$-generators, the localization~$\loc C{\pcat W}$ of $C$ by equivalence
classes of elements of~$W$ is presented by the polygraph~$Q$ defined exactly as
in the previous section. In particular, we recover the free groupoid on~$C$ as
$\freegpd{C}=\loc C{P_1}$.


\section{Distributive Laws}
\label{sec:2-dlaws}
In this section, we present a very useful tool in order to build presentations
in a modular fashion. The typical situation we want to address here is when the
category~$E$ we want to present is ``built'' from two subcategories~$C$ and~$D$,
in the sense that every morphism of~$E$ factors a composite of morphisms in the
two subcategories: in this case, we can expect to be able to construct a
presentation of~$E$ from presentations of~$C$ and~$D$.
The way the category~$E$ can be obtained as a composite of~$C$ and~$D$ can be
encoded in a distributive law. This notion was introduced by
Beck~\cite{beck1969distributive}, related to categories and strict factorization
systems by Rosebrugh and Wood~\cite{rosebrugh2002distributive}, and applied to
presentations of categories by Lack~\cite{lack2004composing}. We begin by
recalling this setting and then presenting the generalizations necessary to
handle situations arising in practice.

\subsection{Strict factorization system}
\index{factorization system!strict}
\label{sec:strict-factorization}
A \emph{strict factorization system} on a category~$E$ consists of two
subcategories~$C$ and~$D$ of~$E$, with the same objects as~$E$, such that every
morphism $h$ of $E$ factorizes uniquely as $h=g\circ f$ with $f$ in~$C$ and $g$
in~$D$:
\[
  \xymatrix@C=6ex@R=3ex{
    \ar[rr]^h\ar@{.>}[dr]_{C\owns f}&&\\
    &\ar@{.>}[ur]_{g\in D}
  }\pbox.
\]

This structure can equivalently be encoded through operations which help
expressing every morphism~$E$ as one of $C$ composed with one of~$D$, as we now
explain.

\newcommand{\dll}[2]{{#2}^{#1}} 
\newcommand{\dlm}[3]{\null^{#1}\!#3^{#2}} 
\newcommand{\dlr}[2]{\null^{#2}\!{#1}} 

\subsection{Distributive law}
\label{sec:cat-dlaw}
\index{distributive law}
A \emph{distributive law}~$\dlaw$ between two categories~$C$ and~$D$ having the
same objects is a function, often noted
\[
  \dlaw
  :
  D\otimes C
  \to
  C\otimes D
\]
(the notation as a tensor will be formally justified in \secr{cat-as-monad}),
which to every ``composable'' pair of morphisms
\[
  g:x\to y\in D,
  \qquad\qquad\qquad
  f:y\to z\in C,
\]
associates an object $\dlm gfy$ of $C$ and $D$, and morphisms
\[
  \dll gf:x\to\dlm gfy\in C,
  \qquad\qquad\qquad
  \dlr gf:\dlm gfy\to z\in D,
\]
which can be pictured as
\[
  \xymatrix@C=3ex@R=3ex{
    &y\ar[dr]^{f\in C}&\\
    x\ar[ur]^{D\owns g}\ar@{.>}[dr]_{C\owns\dll gf}&&z\\
    &\dlm gfy\ar@{.>}[ur]_{\dlr gf\in D}&
  }
\]
in a way compatible with compositions
\begin{align*}
  \vcenter{\xymatrix@C=3ex@R=3ex{
    &y\ar[ddrr]^{f_2\circ f_1}\\
    x\ar@{.>}[ddrr]_{\dll g{(f_2\circ f_1)}}\ar[ur]^g&&\\
    &&&z\\
    &&{}\phantom{x}\ar@{.>}[ur]_{\dlr g{f_2\circ f_1}}
  }}
  &=
  \vcenter{\xymatrix@C=3ex@R=3ex{
    &y\ar[dr]^{f_1}\\
    x\ar@{.>}[dr]_{\dll g{f_1}}\ar[ur]^g&&\ar[dr]^{f_2}\\
    &\ar@{.>}[dr]_{\dll{\dlr g{f_1}}{f_2}}\ar@{.>}[ur]|{\dlr g{f_1}}&&z\\
    &&{}\phantom{x}\ar@{.>}[ur]_{\dlr{(\dlr g{f_1})}{f_2}}
  }}
  \\
  \dll g{(f_2\circ f_1)}&=\dll{\dlr g{f_1}}{f_2}\circ\dll g{f_1}
  \\
  \dlr g{f_2\circ f_1}&=\dlr{(\dlr g{f_1})}{f_2}
\end{align*}
\begin{align*}
  \vcenter{\xymatrix@C=3ex@R=3ex{
    &&y\ar[dr]^f\\
    &&&z\\
    x\ar[uurr]^{g_2\circ g_1}\ar@{.>}[dr]_{\dll{g_2\circ g_1}f}&&\\
    &{}\phantom{x}\ar@{.>}[uurr]_{\dlr{g_2\circ g_1}f}
  }}
  &=
  \vcenter{\xymatrix@C=3ex@R=3ex{
    &&y\ar[dr]^f\\
    &\ar[ur]^{g_2}\ar@{.>}[dr]|{\dll{g_2}f}&&z\\
    x\ar[ur]^{g_1}\ar@{.>}[dr]_{\dll{g_1}{(\dll{g_2}f)}}&&\ar@{.>}[ur]_{\dlr{g_2}f}\\
    &{}\phantom{x}\ar@{.>}[ur]_{\dlr{g_1}{\dll{g_2}f}}
  }}
  \\
  \dll{g_2\circ g_1}f&=\dll{g_1}{(\dll{g_2}f)}
  \\
  \dlr{(g_2\circ g_1)}f&=\dlr{g_2}f\circ\dlr{g_1}{\dll{g_2}f}
\end{align*}
and identities
\begin{align*}
  \vcenter{\xymatrix@C=3ex@R=3ex{
    &y\ar[dr]^{\unit{}}\\
    x\ar[ur]^g\ar@{.>}[dr]_{\dll g{\unit{}}}&&y\\
    &x\ar@{.>}[ur]_{\dlr g{\unit{}}}&
  }}
  &=
  \vcenter{\xymatrix@C=3ex@R=3ex{
    &y\ar[dr]^{\unit{}}\\
    x\ar[ur]^g\ar@{.>}[dr]_{\unit{}}&&y\\
    &x\ar@{.>}[ur]_g&
  }}
  &
  \vcenter{\xymatrix@C=3ex@R=3ex{
    &x\ar[dr]^f\\
    x\ar[ur]^{\unit{}}\ar@{.>}[dr]_{\dll{\unit{}}f}&&y\\
    &y\ar@{.>}[ur]_{\dlr{\unit{}}f}
  }}
  &=
  \vcenter{\xymatrix@C=3ex@R=3ex{
    &x\ar[dr]^f\\
    x\ar[ur]^{\unit{}}\ar@{.>}[dr]_f&&y\pbox.\\
    &y\ar@{.>}[ur]_{\unit{}}
  }}
  \\
  \dll g{\unit{}}&=\unit{}
  &
  \dll{\unit{}}f&=f
  \\
  \dlr g{\unit{}}&=g
  &
  \dlr{\unit{}}f&=\unit{}
\end{align*}

\subsection{Composite category}
\index{category!composite}
Given a distributive law~$\dlaw:D\otimes C\to C\otimes D$, we can \emph{compose}
the categories~$C$ and~$D$ along~$\dlaw$ and obtain a new category, noted
$C\otimes_\dlaw D$: it has the same objects as~$C$ and~$D$, a morphism from
$x\to z$ is a pair of morphisms $(f,g)$ with $f:x\to y$ in $C$ and $g:y\to z$
in~$D$ for some object~$y$, identities are pairs of identities, and compositions
are induced in the expected way by the distributive law:
\begin{align*}
  (f',g')\circ(f,g)
  =
  (\dll g{f'}\circ f,g'\circ\dlr g{f'})
  &&
  \vcenter{
    \xymatrix@C=5ex@R=3ex{
      \ar[dr]^f&&\ar[dr]^{f'}&&\\
      &\ar[ur]^g\ar@{.>}[dr]_{\dll g{f'}}&&\ar[ur]^{g'}\\
      &&\ar@{.>}[ur]_{\dlr g{f'}}
    }
  }
\end{align*}
The fact that the axioms of categories are satisfied follows from the axioms of
distributive laws.

\begin{proposition}
  \label{prop:sfs-dlaw}
  Given categories $C,D,E$ with the same objects, the following
  statements are equivalent.
  \begin{enumerate}
  \item The categories $C$ and $D$ form a strict factorization system on~$E$.
  \item There is a distributive law $\dlaw:D\otimes C\to C\otimes D$ such that
    $C\otimes_\dlaw D=E$.
  \end{enumerate}
\end{proposition}
\begin{proof}
  In the case where $C$ and $D$ form a strict factorization system on~$E$, we
  define the distributive law $\dlaw$ which maps a composable pair of morphisms
  $(g,f)\in D\times C$ to the pair of morphisms obtained by factorizing
  $f\circ g\in E$: the axioms of distributive laws follows from the unique
  factorization of morphisms in~$E$, and the functor $C\otimes_\dlaw D\to E$
  which is the identity on objects and sends a composable pair $(f,g)$ to
  $g\circ f$ is easily seen to be an isomorphism.
  Conversely, given the distributive law $\dlaw:D\otimes C\to C\otimes D$, the
  functor $C\to C\otimes_\dlaw D$ which is the identity on objects and sends a
  morphism $f:x\to y$ to $(f,\unit{y})$ is faithful: the category $C$ can be
  seen as a subcategory of~$C\otimes_\dlaw D$, and similarly for $D$. Moreover,
  $C$ and $D$ form a strict factorization system for~$C\otimes_\dlaw D$: for
  every morphism $(f,g)$ of~$C\otimes_\dlaw D$, we have
  $(f,g)=(\unit{},g)\circ(f,\unit{})$, and this is the unique such
  factorization.
\end{proof}

\subsection{Presenting composite categories}
Because of compatibility with composition, we expect that a distributive law is
uniquely determined by the image of pairs of generators for morphisms of the two
subcategories. In the case of presented categories, the composite category can
thus be presented as follows.

\begin{theorem}
  \label{thm:fc-dlaw}
  Suppose given two 2-polygraphs $P,Q$ and a distributive law
  $\dlaw:\pcat{Q}\otimes\pcat{P}\to\pcat{P}\otimes\pcat{Q}$ between the
  presented categories. Then the category~$\pcat{P}\otimes_\dlaw\pcat{Q}$ is
  presented by the polygraph~$R$ with
  \[
    R_0=P_0=Q_0,
    \qquad\qquad
    R_1=P_1\sqcup Q_1,
    \qquad\qquad
    R_2=P_2\sqcup Q_2\sqcup R_2^\dlaw,
  \]
  where $R_2^\dlaw$ contains a 2-generator
  \begin{equation}
    \label{eq:dlaw-alpha-uv}
    \alpha_{u',v'}
    :
    v'u'
    \To
    uv,
  \end{equation}
  for every pair of composable 1-cells $v'\in\freecat{Q_1}$ and
  $u'\in\freecat{P_1}$ such that we have $\dlaw(\cl v',\cl u')=(\cl u,\cl v)$, for some
  $u\in\freecat{P_1}$ and $v\in\freecat{Q_1}$.

  Moreover, if the rewriting relation on 1-cells induced by~$R_2^\dlaw$
  is terminating then one can restrict~$R_2^\dlaw$ to 2-generators of the form
  \begin{equation}
    \label{eq:dlaw-alpha-ab}
    \alpha_{a,b}
    :
    ba
    \To
    uv,
  \end{equation}
  indexed by pairs of 1-generators $a\in P_1$ and $b\in Q_1$.
\end{theorem}
\begin{proof}
  We have a functor $f:\pcat{R}\to\pcat{P}\otimes_\dlaw\pcat{Q}$ which is the
  identity on objects and sends the class of a 1-gene\-rator~$a\in P_1$ (\resp
  $b\in Q_1$) to the morphism $(\cl{a},\unit{})$ (\resp $(\unit{},\cl{b})$).
  This functor is full since every morphism of~$\pcat{P}\otimes_\dlaw\pcat{Q}$
  is of the form $(\cl u,\cl v)$, with~$u\in\freecat{P_1}$ and
  $v\in\freecat{Q_1}$, which is the image of~$\cl{uv}$. Moreover, by the
  rules~$\alpha_{u',v'}$ every morphism $w\in\freecat{R_1}$ is equivalent to one
  of the form $uv$ with~$u\in\freecat{P_1}$ and~$v\in\freecat{Q_1}$, from which
  the faithfulness of the functor follows easily.

  When the rewriting relation generated by~$R'_2$ is terminating, a normal form
  of a morphism~$v'u'$, with~$v'\in\freecat{Q_1}$ and~$u'\in\freecat{P_1}$ is
  necessarily of the form $uv$ with $u\in\freecat{P_1}$ and~$v\in\freecat{Q_1}$
  and therefore a relation of the form \eqref{eq:dlaw-alpha-uv} is derivable for
  every $u'\in\freecat{P_1}$ and $v'\in\freecat{Q_1}$, and we conclude as above.
\end{proof}

\noindent
Note that the first part of the theorem usually gives rise to infinite
presentations (because $R_2^\dlaw$ is infinite), whereas the reduction provided
by the second part produces finite presentations from finite presentations
(provided that the termination condition is satisfied).

\begin{example}
  \label{ex:dlaw-N26}
  The additive monoids $C=\N/2\N$ and $D=\N/3\N$ respectively admit the
  following presentations, see also \cref{sec:N-mod}:
  \begin{align*}
    \Pres\star{a}{aa=1},
    &&
    \Pres\star{b}{bbb=1}
    \pbox.
  \end{align*}
  The product monoid $C\times D$ contains~$C$ and~$D$ as submonoids: an element
  $n\in C$ can be seen as $(m,0)\in C\times D$, and similarly for~$D$. Moreover,
  every element $(m,n)\in C\times D$ can be seen, in a unique way, as a product
  of an element of~$C$ and one of~$D$, namely $(m,n)=(m,0)+(0,n)$, and therefore
  $C$ and $D$ form a factorization system for~$C\times D$. We deduce that
  $C\times D$ admits the presentation
  \[
    \Pres\star{a,b}{aa=1,bbb=1,ba=ab}
    \pbox.
  \]
  The presented monoid is $\N/6\N$ (the generators $a$ and $b$ respectively get
  interpreted as~$3$ and~$2$) and we have embeddings
  \begin{align*}
    \N/2\N&\to\N/6\N&\N/3\N&\to\N/6\N\\
    p&\mapsto 3p&q&\mapsto 2q
  \end{align*}
  which induce the strict factorization system corresponding to the distributive
  law: one readily verifies that every element $n\in\N/6\N$ can be written in a
  unique way as $n=3p+2q$ with $p\in\N/2\N$ and $q\in\N/3\N$.

  We can more generally recover in this way the presentation for products of
  monoids given in \cref{sec:mon-coprod-prod}. Note that the distributive law
  induced between~$C$ and~$D$ is not the only possible one. For instance, the
  presentation
  \[
    \Pres\star{a,b}{aa=1,bbb=1,ba=abb}
  \]
  induces another one, which is not isomorphic (an argument for this is that it
  is not commutative since $ab\neq ba$ can easily be shown, based on the fact
  that the presentation is convergent, see \cref{sec:word-problem}).
\end{example}

\begin{example}
  Starting from two 2-polygraphs~$P$ and~$Q$, we can take their union and add
  relations of the form~\eqref{eq:dlaw-alpha-ab} and hope that the resulting
  2-poly\-graph~$R$ will present a composite category. This is not the case in
  general. For instance, consider the situation with
  \begin{align*}
    P&=\Pres\star{a}{},
    &
    Q&=\Pres\star{b}{bb=1},
    &
    R&=\Pres\star{a,b}{bb=1,ba=1}
    \pbox.
  \end{align*}
  The 2-polygraphs~$P$ and~$Q$ respectively present the monoids $\N$ and
  $\N/2\N$. In the 2-polygraph~$R$, the relation $a=bba=b$ is derivable and
  thus~$R$ presents~$\N/2\N$: the functor~$\pcat{P}\to\pcat{R}$ is not faithful
  and $\pcat{P}$ and $\pcat{Q}$ thus do not form a strict factorization system
  for $\pcat{R}$.
\end{example}

\begin{example}
  The following (counter-)example illustrates the need for the termination
  hypothesis in the second part of \thmr{fc-dlaw}.
  Consider the polygraph~$R$ whose underlying graph is shown on the left,
  together with the relations on the right:
  \[
    \xymatrix@C=3ex@R=3ex{
      x\ar[r]^a\ar@/_5ex/[rrr]_d&x'\ar@/^/[r]^b\ar@/_/[r]_{b'}&y'\ar[r]^c&y
    }
    \qquad\qquad\qquad
    \begin{array}{r@{\ }l@{\qquad}r@{\ }l}
      ab&\To ab',&abc&\To d,\\
      b'c&\To bc,&ab'c&\To d.
    \end{array}
  \]
  We write~$P$ and~$Q$ for the polygraphs, with no relations, whose respective
  underlying graphs are
  \[
    \xymatrix@C=3ex@R=3ex{
      x&x'\ar@/^/[r]^b&y'\ar[r]^c&y
    }
    \qquad\qquad\text{and}\qquad\qquad
    \xymatrix@C=3ex@R=3ex{
      x\ar[r]^a\ar@/_5ex/[rrr]_d&x'\ar@/_/[r]_{b'}&y'&y\pbox.
    }
  \]
  One easily checks that the canonical inclusions $\pcat{P}\to\pcat{R}$ and
  $\pcat{Q}\to\pcat{R}$ are faithful and that every morphism of~$\pcat{R}$
  factorizes uniquely as one from~$\pcat{P}$ followed by one
  from~$\pcat{Q}$. However, if one restricts to relations of the
  form~\eqref{eq:dlaw-alpha-ab}, the only relations left are $ab\To ab'$ and
  $b'c\To bc$, from which the two relations $abc\To d$ and $ab'c\To d$ are not
  derivable. Here, the termination hypothesis of \cref{thm:fc-dlaw} is clearly
  not satisfied since we have the infinite sequence of reductions
  \[
    abc\To ab'c\To abc\To ab'c\To\ldots
  \]
\end{example}

\begin{example}
  \label{ex:simpl-fs}
  \index{simplicial!category}
  \index{category!simplicial}
  \nomenclature[.D]{$\Simplaug$}{augmented simplicial category}
  There is a strict factorization on the augmented simplicial
  category~$\Simplaug$, presented in detail in \cref{sec:simpl-cat}: every
  morphism factorizes as an epimorphism followed by a
  monomorphism. Writing~$\Simplsurj$ (\resp $\Simplinj$) for the subcategory of
  $\Simplaug$, with the same objects, whose morphisms are surjective (\resp
  injective) functions, we thus have
  $\Simplaug=\Simplsurj\otimes_\dlaw\Simplinj$ for some distributive
  law~$\dlaw$. The categories~$\Simplsurj$ and $\Simplinj$ respectively admit
  the presentations
  \begin{align*}
    &\Pres\star{s^n_i:n+1\to n}{s^{n+1}_is^n_j=s^{n+1}_{j+1}s^n_i}_{n\in\N, 0\leq i\leq j<n}\\
    &\Pres\star{d^n_i:n\to n+1}{d^n_jd^{n+1}_i=d^n_id^{n+1}_{j+1}}_{n\in\N, 0\leq i\leq j\leq n}
  \end{align*}
  and by applying \thmr{fc-dlaw}, we recover the presentation of~$\Simplaug$ given
  in \cref{sec:simpl-pres}, see also \cref{sec:2pres-mon} for a 2-dimensional analysis
  of the situation.
\end{example}

\subsection{Spans}
\label{sec:span}
We now briefly recall the construction of span bicategories, which will turn out
to be useful in order to explain the axioms for distributive laws, as well as
provide a rich source of examples for distributive laws between a category and
its opposite. We refer the reader
to~\cite{rosebrugh2002distributive,lack2004composing} for details.

Suppose given a category~$C$ with pullbacks. A \emph{span}\index{span} $(f,g)$
from~$x$ to~$y$ is a pair of coinitial morphisms
\[
  \xymatrix@C=3ex@R=3ex{
    &\ar[dl]_fv\ar[dr]^g&\\
    x&&y
  }
\]
in~$C$. Given a span $(f,g)$ from~$x$ to~$y$ and $(h,i)$ from $y$ to $z$, one
can define a composite span from~$x$ to~$z$ by taking the pullback of the two
arrows in the middle
\[
  \xymatrix@C=3ex@R=3ex{
    &&\ar@{.>}[dl]u\ar@{.>}[dr]\\
    &\ar[dl]_fv\ar[dr]^g&&\ar[dl]_hw\ar[dr]^i\\
    x&&y&&z
  }
\]
and given an object~$x$ one defines the identity span on~$x$ as
\[
  \xymatrix@C=3ex@R=3ex{
    &\ar[dl]_{\unit{x}}x\ar[dr]^{\unit{x}}&\\
    x&&x\pbox.
  }
\]
A morphism~$h$ between two spans $(f,g)$ and $(f',g')$ from~$x$ to~$y$ is a
morphism of~$C$ making the following diagram commute:
\[
  \xymatrix@C=3ex@R=3ex{
    &\ar[dl]_fv\ar[dd]^h\ar[dr]^g&\\
    x&&y\pbox.\\
    &\ar[ul]^{f'}v'\ar[ur]_{g'}&
  }
\]
Because of the way composition was defined, it is generally not strictly
associative, but rather associative up to isomorphism: we can form a
bicategory~$\Span(C)$ whose 0-cells are the objects of~$C$, 1-cells are spans,
and 2-cells are morphisms of spans.
\nomenclature[Span]{$\Span(C)$}{bicategory of spans in a category~$C$}

Of course, when the category~$C$ has pushouts, one can dually define a
bicate\-gory $\CoSpan(C)$ of \emph{cospans}\index{cospan} in~$C$, \ie diagrams
of the form
\[
  \xymatrix@C=3ex@R=3ex{
    &\ar@{<-}[dl]_fv\ar@{<-}[dr]^g&\\
    x&&y
  }
\]
in the category~$C$.

\subsection{Distributive laws between monads}
\label{sec:bicat-dlaw}
Given a bicategory~$\mcB$, one can consider a
\emph{monad}~$(t,\mu,\eta)$\index{monad} in~$\mcB$ (also called a
\emph{monoid}\index{monoid} in~$\mcB$), which consists of an endomorphism
$t:x\to x$ together two 2-cells $\mu:tt\To t$ and $\eta:\unit{x}\To t$,
respectively called multiplication and unit, which are associative and unital:
\begin{align*}
  \xymatrix{
    (tt)t\ar@{=>}[d]_{\mu t}\ar@{=>}[r]^\sim&t(tt)\ar@{=>}[r]^{t\mu}&tt\ar@{=>}[d]^\mu\\
    tt\ar@{=>}[rr]_\mu&&t
  }
  &&
  \xymatrix{
    \unit{}t\ar@{=>}[dr]_\sim\ar@{=>}[r]^{\eta t}&tt\ar@{=>}[d]_\mu&\ar@{=>}[dl]^\sim\ar@{=>}[l]_{t\eta} t1\pbox.\\
    &t&
  }  
\end{align*}
In particular, a monad in the 2-category~$\Cat$ is a monad in the usual sense. We
write $\Mon(\mcB)$ for the category of monads in~$\mcB$, with the expected
notion of morphism.
\nomenclature[Mon(B)]{$\Mon(\mcB)$}{category of monads in a bicategory~$\mcB$}

\index{distributive law!of monads}
Given two monads~$s:x\to x$ and~$t:x\to x$ the composite~$tu$ is not in general
a monad.
The missing piece of data in order to properly compose those was introduced by
Beck~\cite{beck1969distributive}: a \emph{distributive law} between two
monads~$t$ and~$u$ is a 2-cell $\lambda:ut\To tu$ making the diagrams
\[
  \xymatrix@C=3ex@R=3ex{
    u(tt)\ar@{=>}[d]_{u\mu}\ar@{=>}[r]^\sim&(ut)t\ar@{=>}[r]^{\lambda t}&(tu)t\ar@{=>}[r]^\sim&t(ut)\ar@{=>}[r]^{t\lambda}&t(tu)\ar@{=>}[r]^\sim&(tt)u\ar@{=>}[d]^{\mu u}\\
    ut\ar@{=>}[rrrrr]_\lambda&&&&&tu
  }
\]
\[
  \xymatrix@C=3ex@R=3ex{
    (uu)t\ar@{=>}[d]_{\mu t}\ar@{=>}[r]^\sim&u(ut)\ar@{=>}[r]^{u\lambda}&u(tu)\ar@{=>}[r]^\sim&(ut)u\ar@{=>}[r]^{\lambda u}&(tu)u\ar@{=>}[r]^\sim&t(uu)\ar@{=>}[d]^{t\mu}\\
    ut\ar@{=>}[rrrrr]_\lambda&&&&&tu
  }
\]
\begin{align*}
  \xymatrix@C=3ex@R=3ex{
    u\ar@{=>}[d]_\sim\ar@{=>}[r]^\sim&u\unit{}\ar@{=>}[r]^{u\eta}&ut\ar@{=>}[d]^\lambda\\
    \unit{}u\ar@{=>}[rr]_{\eta u}&&tu
  }
  &&
  \xymatrix@C=3ex@R=3ex{
    t\ar@{=>}[d]_\sim\ar@{=>}[r]&\unit{}t\ar@{=>}[r]^{\eta t}&ut\ar@{=>}[d]^\lambda\\
    t\unit{}\ar@{=>}[rr]_{t\eta}&&tu
  }  
\end{align*}
commute. When equipped with a distributive law, one can define a monad structure
on $tu$, called the \emph{composite monad} of~$t$ and~$u$, with multiplication
and unit being respectively
\[
  \xymatrix{
    tutu\ar@{=>}[r]^{t\lambda u}&ttuu\ar@{=>}[r]^{\mu\mu}&tu
  }
  \qqqqtand
  \xymatrix{
    \unit{}\ar@{=>}[r]^{\eta\eta}&tu
  }
\]
(omitting coherence isomorphisms). A more detailed description is given in \cref{sec:dlaw}.

\subsection{Monads in spans}
\label{sec:cat-as-monad}
Interestingly, a monad in $\Span(\Set)$ precisely corresponds to a small
category: $x$ is the set of objects of the category, the endomorphism~$t:x\to x$
is a span of the form
\[
  \xymatrix{
    x&\ar[l]_fu\ar[r]^g&x
  }
\]
providing the underlying graph of the category (where $x$ and $u$ are
respectively the sets of objects and morphisms of
the category and $f$ and $g$ are respectively the source and target functions),
and $\mu$ and $\eta$ respectively describe compositions and identities of the
category.

A distributive law between two categories seen as spans in this way corresponds
precisely to the notion of distributive law defined in \secr{cat-dlaw}, and the
notation~$C\otimes D$ corresponds to the composite of the underlying 1-cells
(\ie graphs) of the monads corresponding to categories~$C$ and~$D$: concretely,
$C\otimes D$ is the graph with the objects of~$C$ (or equivalently~$D$) as
vertices, and pairs $(f,g)$ with $f:x\to y$ in~$C$ and~$g:y\to z$ in~$D$ as
edges~$x\to z$.

\subsection{Categories of spans}
\label{sec:qspan}
\index{span!category of}
\nomenclature[Span(C)]{$\qSpan(C)$}{category of spans in a category~$C$}
\nomenclature[Cospan(C)]{$\qCospan(C)$}{category of cospans in a category~$C$}
Given a category~$C$ with pullbacks, a category~$\qSpan(C)$ can be defined from
the bicategory~$\Span(C)$ by quotienting 1\nbd-cells under isomorphisms and
discarding 2-cells. This provides a rich source of examples of distributive laws between a
category and its opposite, as we now illustrate. A category~$\qCoSpan(C)$ of
cospans can be defined similarly and, of course, satisfies dual results.

There are canonical functors
\[
  C^\op\to\qSpan(C)
  \qquad\qquad\text{and}\qquad\qquad
  C\to\qSpan(C)
\]
respectively sending a morphism $f:x\to y$ to the class of the
span~$(f,\unit{x})$ and~$(\unit{x},f)$, which are both faithful: the categories
$C^\op$ and $C$ can be considered as subcategories of~$\qSpan(C)$. Moreover, in
the category of spans, we have $(f,g)=(\unit{},g)\circ(f,\unit{})$ so that every
morphism is the composite of a morphism in~$C^\op$ followed by one in~$C$. When
this factorization is unique, we have a strict factorization system; in
particular, this is always the case when the category~$C$ has no non-trivial
isomorphism, because in this case the quotient constructing $\qSpan(C)$
from~$\Span(C)$ will be trivial. In such a situation, when we have a
presentation for~$C$, we clearly also have one for $C^\op$, and thus also for
$\qSpan(C)$ by \thmr{fc-dlaw}. The general case corresponds to a generalized
notion of distributive law, presented in \cref{sec:gen-dlaw}.

\begin{example}
  As a simple example, consider the monoid~$\N$. The pullback of two morphisms
  $m$ and $n$
  \[
    \xymatrix@C=3ex@R=3ex{
      &\ar@{.>}[dl]_{m'}\star\ar@{.>}[dr]^{n'}&\\
      \star\ar[dr]_m&&\ar[dl]^n\star\\
      &\star&
    }
  \]
  is given by $m'=\max(m,n)-m$ and $n'=\max(m,n)-n$, and the only isomorphism is
  the identity $0$. From the presentation of~$\N$ as the free monoid on one
  generator, see \cref{sec:pres-nat}, we deduce that a presentation of $\qSpan(\N)$
  is
  $
    \pres{a,b}{ba=1}
  $,
  \index{bicyclic monoid}%
  \index{monoid!bicyclic}%
  \ie this is the bicyclic monoid, see \secr{bicyclic}, where the relation is
  deduced from the fact that the pullback of $1$ with $1$ is given by~$0$
  and~$0$.
\end{example}

\noindent
Other examples are presented in \secr{2-residuation}, in the slightly different
language of residuals (which provide techniques in order to show on the
presentation that the presented category actually has pushouts) and in
\secr{3-dlaws} using 3-polygraphs.

\subsection{Iterated distributive laws}
\label{sec:iterated-dlaw}
\index{distributive law!iterated}
Composing more than two monads in a bicategory can be achieved if one assumes
that there are distributive laws between any pair of monads and every triple of
distributive laws is compatible in the following sense,
see~\cite{cheng2011iterated} for details.

Suppose given three monads $t,u,v:x\to x$ in a bicategory and distributive laws
\[
  \dlaw_{tu}:ut\To tu,
  \qquad\qquad
  \dlaw_{tv}:vt\To tv,
  \qquad\qquad
  \dlaw_{uv}:vu\To uv,
\]
which are \emph{compatible} in the sense that the following diagram commutes:
\[
  \xymatrix@C=6ex@R=3ex{
    &uvt\ar@{=>}[r]^{u\lambda_{tv}}&utv\ar@{=>}[dr]^{\lambda_{tu}v}\\
    vut\ar@{=>}[ur]^{\lambda_{uv}t}\ar@{=>}[dr]_{v\lambda_{tu}}&&&tuv\pbox.\\
    &vtu\ar@{=>}[r]_{\lambda_{tv}u}&tvu\ar@{=>}[ur]_{t\lambda_{uv}}
  }
\]
In this situation, there are distributive laws
\[
  \xymatrix@C=6ex@R=3ex{
    vtu\ar@{=>}[r]^{\lambda_{tv}u}&tvu\ar@{=>}[r]^{t\lambda_{uv}}&tuv
  }
  \qqtand
  \xymatrix@C=6ex@R=3ex{
    uvt\ar@{=>}[r]^{u\lambda_{tv}}&utv\ar@{=>}[r]^{\lambda_{tu}v}&tuv
  }
\]
respectively between $tu$ and $v$, and $t$ and $uv$, which both induce the same
structure of monad on $(tu)v=t(uv)$.

\subsection{More general compositions}
\label{sec:gen-dlaw}
The notion of strict factorization system (or equivalently of distributive law)
is sometimes too restrictive: in many situations, the desirable factorization is
not strictly unique, but only unique up to an isomorphism (or even up to some
subclass of morphisms). For instance, consider the category~$E$ which is the
full subcategory of sets, with finite sets $\set{0,\ldots,n-1}$ as objects for
$n\in\N$ (this category will be denoted~$\Fun$ in
\cref{sec:2pres-mon}). Generalizing the situation of \exr{simpl-fs}, consider
the categories~$C$ and~$D$ which are the subcategories of~$E$ whose morphisms
are respectively surjective and injective functions: the categories~$C$ and~$D$
``almost'' form a factorization system for~$E$. Namely, every function~$h$
factorizes as $h=g\circ f$ where~$f$ is surjective and~$g$ is injective, and
this factorization is ``almost'' unique in the sense that for every other
factorization $h=g'\circ f'$ there exists an isomorphism~$w$ making the
following diagram commute:
\begin{equation*}
  \xymatrix@C=3ex@R=3ex{
    &\ar[dr]^g\ar@{.>}[dd]^w\\
    \ar[ur]^f\ar[dr]_{f'}&&\pbox.\\
    &\ar[ur]_{g'}
  }
\end{equation*}
Writing~$W$ for subcategory~$W$ of isomorphisms of~$E$, we notice that both the
categories~$C$ and~$D$ contain~$W$ as subcategory. Thus, if we know
presentations for both~$C$ and~$D$, we can expect to deduce a presentation
for~$E$ by taking the union of the presentations of~$C$ and~$D$ and adding
distributivity relations as before (\thmr{fc-dlaw}), but we should moreover
identify the presentation of bijections (\ie the subcategories~$W$) in~$C$ and
in~$D$.
In the above situation, note that the category~$W$ ``acts'' on the left (\resp
on the right) on~$C$: for any morphism $w:x'\to x$ (\resp $w:y\to y'$) of~$W$
and $f:x\to y$ of~$C$ one can obtain a new morphism $wf$ (\resp $fw$)
of~$C$ and the situation is the same for~$D$. The distributive law corresponding
to the above situation should now have the form
\[
  \dlaw
  :
  D\otimes_W C
  \to
  C\otimes_W D
\]
where $C\otimes_W D$ is defined as the quotient of $C\otimes D$ above where the
right action of~$W$ on~$C$ is identified with the left action of~$W$ on~$D$.

Other typical situations where we would like to identify subcategories~$W$
of~$C$ and~$D$ is when those are symmetric monoidal categories (in which~$W$
is the actions of symmetric groups) or Lawvere theories (in which
case~$W=\Fun^\op$, see \chapr{trs}).
Proper generalizations of distributive laws techniques in order to encompass
such situations were given by Lack~\cite{lack2004composing} and detailed and
generalized by Cheng~\cite{cheng2011distributive}; we briefly present those
below.

\subsection{}
\label{sec:gen-fs}
\index{factorization system}
A (non-strict) \emph{factorization system} on a category~$E$ consists of
subcategories~$W$, $C$, and~$D$ with the same objects as~$E$ such that
\begin{itemize}
\item $W$ is a subcategory of both~$C$ and~$D$,
\item every morphism of~$E$ factorizes as $g\circ f$ with $f\in C$ and $g\in D$, and
\item any two factorizations $g\circ f$ and $g'\circ f'$ of a given morphism
  in~$E$ are $W$-equivalent: there is a morphism $w\in W$ making the diagram
  \[
    \svxym{
      &\ar[dr]^g\ar@{.>}[dd]|w&\\
      \ar[ur]^f\ar[dr]_{f'}&&\\
      &\ar[ur]_{g'}&
    }
  \]
  commute.
\end{itemize}
Above, the \emph{$W$-equivalence} is the smallest equivalence relation on
composable pairs of morphisms in~$C\otimes D$ such that for every morphisms
$f:x\to y$ in~$C$, $w:y\to y'$ in~$W$ and $g:y'\to z$ in~$D$ the pairs $(fw,g)$
and $(f,wg)$ are $W$\nbd-equi\-va\-lent.

\subsection{}
Given a bicategory~$\mcB$ with coequalizers, we write $\Mod(\mcB)$ for the
bicategory of \emph{bimodules}\index{bimodule} in~$\mcB$ where
\begin{itemize}
\item a 0-cell is a monad in~$\mcB$,
\item given monads $t:x\to x$ and $u:y\to y$, a 1-cell $f:t\to u$ is a
  \emph{bimodule} in~$\mcB$, \ie a 1-cell $f:x\to y$ of~$\mcB$ together with two
  2-cells of~$\mcB$
  \[
    \lambda:tf\To f
    \qqqqtand
    \rho:fu\To u
  \]
  respectively called \emph{left} and \emph{right action} making the following
  diagrams commute:
  \begin{align*}
    \xymatrix{
      ttf\ar@{=>}[d]_{t\lambda}\ar@{=>}[r]^-{\mu f}&tf\ar@{=>}[d]_\lambda&\ar@{=>}[l]_-{\eta f}\ar@{=>}[dl]^{\unit{}}f\\
      tf\ar@{=>}[r]_\lambda&f
    }
    &&
    \xymatrix{
      tfu\ar@{=>}[d]_{t\rho}\ar@{=>}[r]^{\lambda u}&fu\ar@{=>}[d]^\rho\\
      tf\ar@{=>}[r]_\lambda&f
    }
    &&
    \xymatrix{
      f\ar@{=>}[dr]_{\unit{}}\ar@{=>}[r]^{f\eta}&fu\ar@{=>}[d]^\rho&\ar@{=>}[l]_{f\mu}\ar@{=>}[d]^{\rho u}fuu\\
      &f&\ar@{=>}[l]^\rho fu
    }    
  \end{align*}
\item a 2-cell $\phi:f\To g:t\to u$ is a 2-cell $\phi:f\To g$ in~$\mcB$ making
  the following diagram commute:
  \[
    \xymatrix{
      tf\ar@{=>}[d]_{t\phi}\ar@{=>}[r]^\lambda&f\ar@{=>}[d]^\phi&\ar@{=>}[l]_\rho fu\ar@{=>}[d]^{\phi u}\\
      tg\ar@{=>}[r]_\lambda&g&\ar@{=>}[l]^\rho gu
    }
  \]
\item the composite $f\otimes_u g$ of 1-cells $f:t\to u$ and $g:u\to v$ is given
  by the following coequalizer in~$\mcB$:
  \[
    \xymatrix{
      fug\ar@<.8ex>@{=>}[r]^{\rho g}\ar@<-.8ex>@{=>}[r]_{f\lambda}&fg\ar@{:>}[r]&f\otimes_ug
    }
  \]
  and other compositions and identities are the expected ones.
\end{itemize}

In particular, $\Mod(\Set)$ is biequivalent to the usual bicategory of
profunctors, where a 0\nbd-cell is a category and a 1-cell $f:C\to D$ is a
functor $f:C^\op\times D\to\Set$ (called a \emph{profunctor}\index{profunctor}
from $C$ and~$D$).
More interestingly for our matters, consider the bicategory $\Mod(\Span(\Set))$.
By the correspondence given in \secr{cat-as-monad}, the 0-cells are categories
($V$, $W$, \ldots), and a 1-cell $C:V\to W$ is a span
\[
  \xymatrix{
    V_0&\ar[l]_-sC\ar[r]^-t&W_0
  }
\]
that can be seen as a set~$C$ of ``arrows'' with source (\resp target) being an
object of~$V$ (\resp $W$) on which $V$ (\resp~$W$) acts by precomposition (\resp
postcomposition). The horizontal composite $C\otimes_W D$ thus consists of the
set of pairs $(f,g)$ of composable arrows in~$C\times D$ quotiented by the
equivalence relation identifying $(fw,g)$ and $(f,wg)$ for composable $f\in C$,
$w\in W$ and~$g\in D$.

\subsection{}
\label{sec:gen-comp}
Fix a category~$W$. A \emph{category under~$W$} consists of a category~$C$, with
the same objects as~$W$, together with a functor~$f:W\to C$ which is the
identity on objects. We can think of~$C$ as having~$W$ as distinguished
subcategory, at least when the functor~$f$ is faithful.

Given a monad $t:x\to x$ in a bicategory $\mcB$ with colimits, there is always
an equivalence of categories
\[
  \Mon(\Mod(\mcB)(t,t))
  \quad\equivto\quad
  t/\Mon(\mcB(x,x))
  \pbox.
\]
Instantiated to the case where~$\mcB=\Span(\Set)$ and $t$ is a category~$W$,
this says that we have a correspondence between monads in $W$-bimodules and
categories under~$W$.

Given two categories~$C$ and~$D$ under~$W$, we can now define a
\emph{distributive law}\index{distributive law!of categories}
\begin{equation}
  \label{eq:gen-dlaw}
  \dlaw
  :
  D\otimes_WC
  \to
  C\otimes_WD
\end{equation}
as being a distributive law between~$C$ and~$D$, seen as monads in
$W$-bimodules.
Explicitly, it consists of a function which maps a $W$-equivalence class of a
composable pair $(g,f)$ of morphisms $g\in D$ and $f\in C$ to a $W$-equivalence
class of a composable pair $(\dll gf,\dlr gf)$ with $\dll gf\in C$ and
$\dlr gf\in D$ in a way compatible with compositions and identities, in a
similar fashion as for distributive laws, see \cref{sec:bicat-dlaw,sec:dlaw}. As
before, we write~$C\otimes_\dlaw D$ for the resulting composite category.

Every factorization system (in the sense of \secr{gen-fs}) induces a
distributive law in the above sense, and conversely a distributive law induces a
factorization system when the functors $W\to C$ and $W\to D$ are faithful.

\subsection{}
\label{sec:gen-pres}
A generalization of \thmr{fc-dlaw} can also be given as follows. Suppose
given two 2-polygraph $P$ and $Q$ with the same $0$-cells, and write $W=P\cap Q$
for the 2-polygraph such that
\[
  W_0=P_0=Q_0,
  \qquad\qquad
  W_1=P_1\cap Q_1,
  \qquad\qquad
  W_2=P_2\cap Q_2.
\]
Above, we suppose that source and target maps agree in~$P$ and~$Q$ for elements
of~$W_1$ and of $W_2$, and that they induce those of~$W$. The inclusion $W\to P$
induces a functor $\pcat{W}\to\pcat{P}$ making $\pcat{P}$ a category
under~$\pcat{W}$, and similarly for $\pcat{Q}$. Suppose given a distributive law
\[
  \dlaw
  :
  \pcat{Q}\otimes_{\pcat{W}}\pcat{P}
  \to
  \pcat{P}\otimes_{\pcat{W}}\pcat{Q}
\]
between the presented categories. Then the category
$\pcat{P}\otimes_\dlaw\pcat{Q}$ is presented by the polygraph~$R$ with
\[
  R_0=P_0=Q_0,
  \qquad\qquad
  R_1=P_1\cup Q_1,
  \qquad\qquad
  R_2=(P_2\cup Q_2)\sqcup R_2^\dlaw,
\]
(note that some unions are not disjoint) where $R_2^\dlaw$ contains a
2-generator
\[
  \alpha_{u',v'}
  :
  v'u'
  \To
  uv,
\]
for every composable 1-cells $v'\in\freecat{Q_1}$ and $u'\in\freecat{P_1}$ such
that $\dlaw(\cl b,\cl a)=(\cl u,\cl v)$, for some $u\in\freecat{P_1}$ and
$v\in\freecat{Q_1}$. When the rewriting relation induced by~$R_2^\dlaw$ is
moreover terminating, this set can be further reduced as in \thmr{fc-dlaw}.

\begin{example}
  Suppose given a category~$C$ and write $W$ for the subcategory of~$C$, with
  the same objects and the isomorphisms of~$C$ as morphisms. When~$C$ has
  pullbacks, there is a distributive law
  \[
    \dlaw
    :
    C\otimes_WC^\op
    \to
    C^\op\otimes_W C,
  \]
  which to a pair of morphism $(g',f'^\op)$ associates the pullback $(f^\op,g)$:
  \[
    \xymatrix@C=3ex@R=3ex{
      &\ar@{.>}[dl]_fy\ar@{.>}[dr]^g&\\
      \ar[dr]_{g'}x&&z\ar[dl]^{f'}\\
      &y'&
    }
  \]
  and the composite category is $C^\op\otimes_\dlaw C=\qSpan(C)$, with the category
  of isomorphism classes spans described in \secr{qspan},
  see~\cite{lack2004composing,rosebrugh2002distributive,zanasi2015interacting}. Dually,
  when~$C$ has pushouts, the category~$\qCospan(C)$ can be obtained as
  $C\otimes_\dlaw C^\op$ where $\dlaw$ is given by pushout.
\end{example}

\noindent
Other examples and applications are given in \secr{3-dlaws}.


\chapter{String Rewriting and 2-Polygraphs}
\label{chap:2rewr}
\label{Chapter:StringRewritingSystems}

We recast the notion of \emph{string rewriting system} into the language of
polygraphs. This notion, which consists of a set of pairs of words
called \emph{relations} or \emph{rewriting rules} over a fixed
alphabet, can be traced back to Thue. In his 1914 paper~\cite{thue1914probleme}, he introduces the notion of
\emph{word problem}: this is the question of deciding whenever two words are
equivalent with respect to the congruence generated by the relations. He also
shows that the word problem is decidable when the associated rewriting system is
terminating and confluent, and even introduces a completion algorithm in order
to make a system confluent (an accessible presentation of the paper, along with
an English translation can be found in~\cite{power2013thue}). For this reason,
string rewriting systems are also sometimes called \emph{semi-Thue systems}
(the ``semi'' here is to distinguish with \emph{Thue systems} which are defined
in the same way, but where the relations are not oriented). Unexpectedly at the
time, the word problem was shown to be undecidable for those systems in 1947 by
Post~\cite{post1947recursive} and Markov~\cite{markov1947}. Of course, this does
not preclude us from deciding the word problem for subclasses of monoids, and
this is precisely what rewriting is about.
The notion of string rewriting system is a variant of the notion
of presentation for groups, which is adapted to monoids and where the relations are oriented.
Group presentations have been introduced by Dehn~\cite{dehn1911unendliche} in
1911 along with the corresponding word problem for
finitely presented groups and Dehn's algorithm for solving the word problem in
favorable cases. However, the general word problem for groups has been shown
undecidable by Novikov~\cite{novikov1955algorithmic} and Boone~\cite{boone1958word}.
We do not intend to give a complete presentation of those early works, nor
of the recent developments, and we refer the reader to the standard textbooks~\cite{BaaderNipkow98,book1993string,jantzen1988confluent,bezem2003term}
for an in-depth treatment. We rather explain here how string rewriting
systems can be seen as a particular case of 2\nbd-poly\-graphs, and how the
polygraphic rewriting techniques generalize traditional ones.

The notion of string rewriting system---and the more general variant adapted to
categories---is introduced in \cref{sec:2rewr}, where we show that the rewriting
paths form the morphisms of a sesquicategory, in which we can instantiate the
concepts for abstract rewriting systems developed in \cref{chap:lowdim}. In
\cref{sec:word-problem}, we introduce the word problem and show that it can be
efficiently solved for convergent, \ie confluent and
terminating rewriting systems. 
In practice, confluence can be checked by inspecting
the \emph{critical branchings} of the rewriting system, presented in
\cref{sec:critical-branchings}, and termination by introducing a suitable
\emph{reduction order}, as defined in \cref{sec:2-red-order}. The convergence of
a rewriting system is also useful to show that it forms a presentation of a
given category, as  illustrated in \cref{sec:2-presenting}. Finally, in
\cref{sec:2-residuation}, we introduce \emph{residuation} techniques which allow
proving useful properties of categories (such as the existence of pushouts) by
performing computations on their presentations.

\section{String Rewriting Systems}
\label{sec:2rewr}
We have seen in \cref{sec:pres-cat} that a $2$-polygraph~$P$ can be considered
as a notion of presentation for the category~$\pcat{P}$, obtained from the
category freely generated by the underlying 1\nbd-poly\-graph $\tpol 1P$, by
quotienting the $1$-cells under the congruence~$\approx$ generated by~$P_2$. By
\cref{lem:2pol-2cell-cong}, this congruence is the smallest equivalence relation
identifying two $1$-cells~$u$ and~$v$ whenever there is a $2$-cell $\phi:u\To v$
in $\freecat{P_2}$. In such a situation, $u$ and $v$ are thus two representatives of
the same $1$-cell in~$\pcat P$, and if we adopt the point of view developed in
\cref{chap:lowdim}, we can think of $\phi:u\To v$ as indicating that $v$ is a
``more canonical representative'' of the equivalence class than~$u$.
All this suggests that a $2$-polygraph can be considered as a form of rewriting
system, where the objects of interest are the $1$-cells in $\freecat P_1$, and
where the generators $\alpha:u\To v$ in $P_2$ are rewriting rules indicating that
$u$ can be rewritten to~$v$.

If we consider the particular case of a $2$-polygraph~$P$ with only one
$0$\nbd-gene\-rator, say $P_0=\set\star$, the presented category $\pcat P$ has only
$0$-cell and can thus be considered as a monoid, as explained in
\cref{sec:pres-mon}. We will see that, if we restrict to such polygraphs, the
associated notion of rewriting system corresponds precisely to string rewriting
systems, thus establishing $2$-polygraphs as a mild generalization of those, in
which letters are ``typed'' and only well-typed words are considered: in
practice, this extra generality does not bring any major complication and we
develop here the traditional theory of rewriting in full generality.

\subsection{Terminology}
\index{string rewriting system}
\index{rewriting!system!string}
\index{rewriting!rule}
A $2$-polygraph~$P$, when considered as a rewriting system, is sometimes called
a \emph{categorical string rewriting system}, or a \emph{1\nbd-dimensional rewriting system}.
The terminology \emph{string rewriting system} is reserved to the particular case where $P_0=\set{\star}$.
The underlying $1$-polygraph~$\tpol1 P$ is called the \emph{signature} and is
composed of \emph{sorts} (the elements of~$P_0$) and \emph{letters} (the
elements of~$P_1$). The $1$-cells in $\freecat{P_1}$ freely generated by this
signature are called \emph{words} or \emph{strings}, and an identity is
sometimes referred to as an \emph{empty word}. The $2$-generators are the
\emph{rewriting rules} of the rewriting system.

\subsection{Rewriting step}
\index{rewriting!step}
Suppose fixed a $2$-polygraph~$P$. A \emph{rewriting step of $P$}
\[
  \vxym{
    x\ar[r]^u&x'\ar@/^3ex/[r]^v\ar@/_3ex/[r]_{v'}\ar@{}[r]|{\phantom\alpha\Longdownarrow\alpha}&y'\ar[r]^w&y
  }
\]
consists in a $2$-generator
$
  \alpha
  :
  v\To v'
  :
  x'\to y'
$
in $P_2$, together with two $1$-cells
$u:x\to x'$
and
$w:y'\to y$
in $\freecat{P_1}$. Such a rewriting step will be denoted
\begin{equation}
  \label{eq:rewriting-step}
  u\alpha w
  :
  uvw
  \To
  uv'w
  :
  x\to y
\end{equation}
and pictured as
\begin{equation}
  \label{ex:rewriting-step}
  \vxym{
    x\ar@/^3ex/[rrr]^{uvw}\ar@/_3ex/[rrr]_{uv'w}\ar@{}[rrr]|{\textstyle\phantom{u\alpha w}\Longdownarrow u\alpha w}&&&y\pbox.
  }
\end{equation}
The $1$-cell $uvw$ (\resp $uv'w$) in $\freecat{P_1}$ is called its
\emph{source} (\resp \emph{target}). In this situation, we say that $uvw$ is
\emph{reducible}\index{reducible} by~$\alpha$. The pair $(u,w)$ of $1$-cells
in~$\freecat{P_1}$ is sometimes called the \emph{context}\index{context} or
\emph{whisker}\index{whisker} in which the rule $\alpha$ applies to the $1$-cell
$uvw$. We sometimes write
$
u\To v
$
to indicate that there exists a rewriting step of $P$ from~$u$ to~$v$.

\subsection{Rewriting path}
\index{rewriting!path}
\label{subsec:rewr-path}
A \emph{rewriting path of $P$} is a sequence~$\phi$
\begin{equation}
  \label{eq:rp-seq}
  u_1\alpha_1w_1
  ,
  u_2\alpha_2w_2
  ,
  \ldots
  ,
  u_n\alpha_nw_n
\end{equation}
of rewriting steps of $P$
\[
  u_i\alpha_iw_i
  :
  u_iv_iw_i
  \To
  u_iv_i'w_i
  :
  x\to y
\]
which is composable, in the sense that $u_iv_i'w_i=u_{i+1}v_{i+1}w_i$ for
$1\leq i<n$.
\index{rewriting!path!length}
\index{length}
The natural number~$n$ is called the \emph{length} of the rewriting path~$\phi$
and is denoted by~$\sizeof{\phi}$.
The $1$-cells $u_1v_1w_1$ (\resp $u_nv_n'w_n$) are called the
\emph{source} (\resp \emph{target}) of the rewriting path, what we write
\[
  \phi
  :
  u_1v_1w_1
  \To
  u_nv_n'w_n
  :
  x\to y
  \pbox.
\]
By convention, an empty path has a determined source (which is the same as its
target).
We sometimes write
$
u\overset*\To v
$
when there exists a rewriting path from~$u$ to~$v$, in which case we say that
\emph{$u$ rewrites to $v$}. Given two rewriting paths
$\phi:u\To v$
and
$\psi:v\To w$,
we write
$
\phi\comp{}\psi
$
for their concatenation. The rewriting path \eqref{eq:rp-seq} can therefore be
written as a composition of rewriting steps:
\begin{equation}
  \label{eq:rp}
  \phi
  =
  (u_1\alpha_1w_1)
  \comp{}
  (u_2\alpha_2w_2)
  \comp{}
  \ldots
  \comp{}
  (u_n\alpha_nw_n)
  \pbox.
\end{equation}
Given two $1$-cells
$
u:x'\to x
$ and $
w:y\to y'  
$
in $\freecat{P_1}$, we extend the notation~\eqref{eq:rewriting-step} and write
$u\phi w$ for the rewriting path
\begin{equation}
  \label{eq:rp-lr}
  u\phi w
  =
  ((uu_1)\alpha_1(w_1w))\comp{}((uu_2)\alpha_2(w_2w))\comp{}\ldots\comp{}((uu_n)\alpha_n(w_nw))
  \pbox.
\end{equation}
These operations equip the $0$-cells, $1$-cells and rewriting paths in a
polygraph with the structure of a sesquicategory, see \cref{sec:sesquicat}.

\subsection{Support}
\index{support!of a 2-cell}
Any $2$-cell $\phi$ in $\freecat P_2$ can be written as a $1$-composite of
finitely many rewriting steps, of the form \eqref{eq:rp}. We define the
\emph{support} of the $2$-cell $\phi$ as the multiset, denoted by $\msupp[2]f$,
consisting of the $2$-cells $\alpha_i$ occurring in this decomposition. The
support is well-defined because any two decompositions of~$\phi$ in~$\freecat P_2$
into a $1$-composite of rewriting steps involve the same rewriting steps. We
have seen in \cref{sec:multisets} that multiset inclusion is a
well-founded order on supports, allowing us to prove some properties by
induction on the support of $2$-cells.

\subsection{Sesquicategory}
\label{sec:sesquicat}
\index{sesquicategory}
A \emph{sesquicategory}~$C$ consists of
\begin{itemize}
\item a $2$-graph~$C$ (see \cref{sec:2-graph}),
\item a structure of category~$C'$ on the underlying $1$-graph of~$C$,
\item a functor $C(-,-):\opp{C'}\times C'\to\Cat$,
\end{itemize}
such that the composite of the functor $C(-,-)$ with the forgetful functor
$\Cat\to\Set$, which to a category associates its set of objects, coincides with
the hom functor $C'(-,-):\opp{C'}\times C'\to\Set$.

The notion of sesquicategory was introduced by
Street~\cite{street1996categorical}. Let us detail the operations available in
such a structure.
Given \emph{0-cells} $x$, $y$ (\ie objects of~$C'$), we have a category $C(x,y)$
whose objects are the morphisms $f:x\to y$ of~$C'$, called \emph{1-cells},
morphisms
$
\alpha
:
f\To g
:
x\to y
$
are called the \emph{2-cells}, and composition is denoted $\comp{}$ and called
(vertical) composition. Given a 2-cell $\alpha:g\To g':x\to y$ and 1-cells
$f:x'\to x$ and $h:y\to y'$, we have a 2-cell $C(f,h)(\alpha)$ that will be
denoted
\[
f\alpha h
:
fgh\To fg'h
:
x'\to y'
\]
and pictured as
\[
\vxym{
  x'\ar[r]^f&x\ar@/^/[r]^g\ar@{}[r]|{\phantom\alpha\Downarrow\alpha}\ar@/_/[r]_{g'}&y\ar[r]^h&y'\pbox.
}
\]
The functoriality of $C(-,-)$ ensures that this is a proper left and right
action of 1-cells on 2-cells: in a situation such as
\[
  \svxym{
    x''\ar[r]^{f'}&x'\ar[r]^f&x\ar@/^/[r]^g\ar@{}[r]|{\phantom\alpha\Downarrow\alpha}\ar@/_/[r]_{g'}&y\ar[r]^h&y'\ar[r]^{h'}&y''
  }
  \qqqtor
  \svxym{
    x\ar[r]^{\id_x}&x\ar@/^/[r]^g\ar@{}[r]|{\phantom\alpha\Downarrow\alpha}\ar@/_/[r]_{g'}&y\ar[r]^{\id_y}&y\pbox,
  }
\]
we have
\[
f'(f\alpha h)h'=(f'f)\alpha(hh')
\qqqtand
\id_x\alpha\id_y=\alpha
\pbox.
\]
The following observation is a reformulation in the language of polygraphs of
observations originating in~\cite{stell1994modelling,street1996categorical}:

\begin{lemma}
  \label{lem:pol-sesqui}
  Any 2-polygraph~$P$, induces a sesquicategory with $\freecat{\tpol1 P}$ as
  underlying category, rewriting paths of $P$ as 2-cells and left and right actions
  defined as in~\eqref{eq:rp-lr}.
\end{lemma}

Any 2-category~$C$ induces a sesquicategory in the expected way, where $C'$ is
the category underlying~$C$ (with $C_0$ as objects and $C_1$ as morphisms) and
for every $x,y\in C_0$, $C(x,y)$ is the hom-category whose objects are 1-cells
$f:x\to y$ in~$C_1$ and morphisms are 2-cells $\alpha:f\To g:x\to y$ in~$C_2$,
vertical composition~$\comp{}$ is $\comp1$, and the action of 1-cells on 2-cells
is given by $f\alpha g=\unit f\comp0\alpha\comp0 g$. Moreover, the horizontal
composition~$\comp0$ of the original 2-category can be recovered from the
vertical composition and the action since, given $\alpha:f\To f':x\to y$ and
$\beta:g\To g':y\to z$, we have
\[
  \begin{array}{rcccl}
    (\unit x\alpha g)
    \comp{}
    (f'\beta\unit z)
    &=&
    \alpha\comp0\beta
    &=&
    (f\beta\unit z)
    \comp{}
    (\unit x\alpha g')
    \\
    \vcenter{
      \xymatrix{
        x\ar@/^3ex/[r]^f\ar@{{}{ }{}}@/^1.5ex/[r]|{\Downarrow\alpha}\ar[r]_{f'}&y\ar[r]^g\ar@{{}{ }{}}@/_1.5ex/[r]|{\Downarrow\beta}\ar@/_3ex/[r]_{g'}&z
      }
    }
    &=&
    \vcenter{
      \xymatrix{
        x\ar@/^/[r]^f\ar@{{}{ }{}}[r]|{\Downarrow\alpha}\ar@/_/[r]_{f'}&y\ar@/^/[r]^g\ar@{{}{ }{}}[r]|{\Downarrow\beta}\ar@/_/[r]_{g'}&z
      }
    }
    &=&
    \vcenter{
      \xymatrix{
        x\ar[r]^f\ar@{{}{ }{}}@/_1.5ex/[r]|{\Downarrow\alpha}\ar@/_3ex/[r]_{f'}&y\ar@/^3ex/[r]^g\ar@{{}{ }{}}@/^1.5ex/[r]|{\Downarrow\beta}\ar[r]_{g'}&z\pbox.
      }
    }
  \end{array}
\]
In a general sesquicategory, the left and right members of the above equality
are not necessarily equal, and sesquicategories in which this is always the case
are precisely $2$-categories:

\begin{proposition}
  \label{prop:2cat-sesqui}
  A \emph{$2$-category}\index{2-category} is a sesquicategory such that for
  every 2-cells
  $\alpha:f\To f':x\to y$
  and
  $\beta:g\To g':y\to y'$
  we have
  \begin{equation}
    \label{eq:sesqui-exchange}
    (\alpha g)\comp{}(f'\beta)
    =
    (f\beta)\comp{}(\alpha g')
    \pbox.
  \end{equation}
\end{proposition}

\noindent
This explains the name \emph{sesqui}category, meaning a 1\textonehalf-category:
a sesquicategory is almost a 2-category excepting that the exchange law is not
required to hold.

\subsection{Freely generated 2-category}
\index{free!2-category}
\label{sec:free-sesqui-2-cat}
From the alternative description of a $2$\nbd-cate\-gory provided by
\cref{prop:2cat-sesqui}, one can come up with an alternative construction of the
$2$\nbd-cate\-gory~$\freecat{P}$ generated by a $2$\nbd-polygraph~$P$ (see
\secr{free-2-cat}): $\freecat{P}$ is the $2$\nbd-category,
with~$\freecat{\tpol1 P}$ as underlying category, whose $2$-cells are rewriting
paths of $P$ considered up to the congruence generated by~\eqref{eq:sesqui-exchange}.
This means that we do not take in account the order
of rewriting steps operating at disjoint positions and consider them up to the
congruence identifying two rewriting paths of length two of the form
\begin{equation}
  \label{eq:rewriting-exchange}
  \vcenter{
    \xymatrix@C=2.5ex{
      x\ar[r]^{u_1}&x'\ar@/^3ex/[r]^v\ar[r]_{v'}\ar@{{}{ }{}}@/^/[r]|{\phantom\alpha\Downarrow\alpha}&y'\ar[r]^{u_2}&y\ar[r]^w\ar@/_3ex/[r]_{w'}\ar@{{}{ }{}}@/_/[r]|{\phantom\beta\Downarrow\beta}&z\ar[r]^{u_3}&z'
    }
  }
  =
  \vcenter{
    \xymatrix@C=2.5ex{
      x\ar[r]^{u_1}&x'\ar[r]^v\ar@/_3ex/[r]_{v'}\ar@{{}{ }{}}@/_/[r]|{\phantom\alpha\Downarrow\alpha}&y'\ar[r]^{u_2}&y\ar@/^3ex/[r]^w\ar[r]_{w'}\ar@{{}{ }{}}@/^/[r]|{\phantom\beta\Downarrow\beta}&z\ar[r]^{u_3}&z'\pbox.
    }
  }
\end{equation}

It can be shown that the sesquicategory constructed in \cref{lem:pol-sesqui} is
free on the polygraph in the expected sense, akin to \cref{sec:free-2-cat}. One
of the main advantage of considering sesquicategories instead of 2-categories
here is that the $2$-cells are much easier to represent by data structures, thus
making those amenable to mechanized computations: the presence of the quotient
\eqref{eq:rewriting-exchange} makes everything more difficult.

\subsection{Rewriting properties of 2-polygraphs}
\index{abstract rewriting system!of a 2-polygraph}
\label{Subsection:RewritingProperties2Polygraphs}
\label{sec:2-rprop}
Given a $2$-polygraph~$P$, we write here~$P^{\mathrm{rs}}$ for its set of
rewriting steps and $s_1,t_1:P^{\mathrm{rs}}\to\freecat{P_1}$ for the functions
which to a rewriting step respectively associates its source and target.
Any 2-polygraph~$P$ thus induces an abstract rewriting system
\[
  \xymatrix{
    \freecat{P_1}&\ar@<-.5ex>[l]_-{\sce1}\ar@<.5ex>[l]^-{\tge1}P^{\mathrm{rs}},
  }
\]
with $1$-cells as vertices and rewriting steps
\[
u\alpha w
:
uvw
\To
uv'w,
\]
as in~\eqref{eq:rewriting-step}, as edges from $uvw$ to~$uv'w$. We always use
double arrows to denote the edges of this rewriting system. Note that, with this
point of view, the two rewriting paths shown in~\eqref{eq:rewriting-exchange}
are not considered to be equivalent.

This construction allows us to extend the properties of \secr{ars} to
$2$\nbd-poly\-graphs. In particular, a $2$\nbd-polygraph is

\begin{center}
  \emph{terminating} /
  \emph{quasi-terminating} /
  \emph{Church-Rosser} /
 
  \emph{confluent} / 
  \emph{locally confluent} /
  \emph{decreasing}
 
  \emph{convergent} /
  \emph{quasi-convergent}
\end{center}
 when the associated abstract rewriting system
is. Moreover, the properties of \secr{ars} immediately extend to our case. We
list below such constructions and properties, reformulated in the framework of
$2$-polygraphs.

\subsection{Branching}
\index{branching}
\index{local branching}
\index{branching!local}
\index{branching!confluent}
A \emph{branching} in a 2-polygraph~$P$ is a pair $(\phi_1,\phi_2)$ of coinitial
rewriting paths
$\phi_1:u\To v_1$
and
$\phi_2:u \To v_2$
in $\freecat{P_2}$, which we sometimes write
$(\phi_1,\phi_2):u\To(v_1,v_2)$.

Such a branching is \emph{local} when both $\phi_1$ and $\phi_2$ are rewriting steps.
It is \emph{confluent} when there exist cofinal rewriting paths
$\psi_1:v_1\To w$
and
$\psi_2:v_2\To w$
which ``close'' the diagram:
\begin{equation}
  \label{eq:2-branching}
  \svxym{
    &\ar@{=>}[dl]_{\phi_1}u\ar@{=>}[dr]^{\phi_2}&\\
    v_1\ar@{:>}[dr]_{\psi_1}&&\ar@{:>}[dl]^{\psi_2}v_2\pbox.\\
    &w&
  }
\end{equation}
We sometimes write $(\psi_1,\psi_2):(v_1,v_2)\To w$ for such a pair.

The goal of this chapter is to provide conditions which are sufficient to ensure
that a $2$\nbd-polygraph is locally confluent (\secr{2-cb}) and terminating
(\secr{2-red-order}). When both properties are satisfied, we can apply Newman's
\cref{lem:newman} to conclude that it is confluent.

\section{Deciding Equality}
\index{decidability!of equality}
\index{word problem!for 2-polygraphs}
\label{sec:word-problem}
One of the main applications of showing that a $2$-polygraph~$P$ is convergent
is to show that the \emph{equality decision problem}, or \emph{word problem},
for~$P$ is decidable for those. We have already handled this situation in the
case of $1$-polygraphs in~\secr{1-eq-dec}.

\subsection{The word problem}
In the context of $2$-polygraphs, the equality decision problem for a
$2$-polygraph~$P$ is often called the \emph{word problem}
for~$P$, and consists in answering the following question:

\begin{center}
  Given two $1$-cells $u,v\in\freecat{P_1}$, do we have $u\approx v$?
\end{center}

\noindent
Above, we recall that $\approx$ denotes the congruence generated by~$P_2$, see
\cref{sec:P-congruence}. The problem was originally introduced by
Thue~\cite{thue1914probleme}, as well as Dehn~\cite{dehn1910topologie} in the
closely related context of group presentations.

A $2$-polygraph has \emph{decidable word problem} if there is an algorithm
answering the above question: this algorithm consists in a procedure, taking the
$1$-cells $u$ and $v$ as arguments, and answering \code{true} or \code{false}
depending on whether $u\approx v$ holds or not. We should insist on the fact
that we require that the procedure terminates on every input, \ie provides an
answer after some finite amount of time. When there is no such procedure, the
problem is said to be \emph{undecidable}.

\subsection{Undecidability}
We always restrict to finite polygraphs~$P$ (since the algorithm needs to use
this polygraph, the latter must be encoded in a finite way, although we could
consider the more general case of recursively enumerable presentations).
Contrarily to the case of dimension~$1$, the set $\freecat{P_1}$ of $1$-cells
is generally infinite, even though the polygraph~$P$ is supposed
to be finite. Therefore, the argument used in \secr{1-eq-dec} for showing that
the problem is decidable cannot be used anymore: the naive procedure, consisting
in computing the equivalence class of~$u$ and checking whether~$v$ belongs to it or not,
is not guaranteed to terminate since the class might not be finite.
In fact, the word problem was shown by Post~\cite{post1947recursive} and
Markov~\cite{markov1947} to be undecidable in general: there exists a finite
$2$-polygraph~$P$ for which there is no algorithm deciding the word problem. A
concrete example of such a polygraph is given in \secr{tseitin}.

Having a decidable word problem is however a property of the monoid, not of a
particular presentation:
\begin{proposition}
  \label{Proposition:InvarianceDecidabilityWordProblem}
  Let $P$ and $Q$ be two finite Tietze equivalent 2-polygraphs. Then the word
  problem for $P$ is decidable if and only if the word problem for $Q$ is
  decidable.
\end{proposition}
\begin{proof}
  Suppose given two finite Tietze equivalent $2$-polygraphs~$P$ and~$Q$ such
  that the word problem is decidable for~$Q$. Given two parallel $1$-cells $u$
  and $v$ in~$\freecat{P_1}$, the Tietze equivalence allows the effective
  construction for every $1$\nbd-gene\-rator $a\in P_1$ of a $1$-cell
  $[a]\in\freecat{Q_1}$ such that $\pcat a=\pcat{[a]}$. Extending the operation
  $[-]$ as a functor $\freecat P_1\to\freecat Q_1$, we have that $u\approx v$
  in~$P$ if and only if $[u]\approx[v]$ in~$Q$, thus allowing us to conclude.
\end{proof}

\subsection{The normal form algorithm}
\index{algorithm!normal form}
\index{normal form!algorithm}
\label{SS:NormalFormProcedure}
When the $2$-polygraph~$P$ of interest is finite and convergent, the word problem can
be decided as in the case of dimension~$1$ presented in
\secr{1-eq-dec}. Namely, the normal form~$\nf{u}$ of a $1$-cell~$u$ can be
computed by maximally rewriting~$u$, and two $1$-cells $u$ and $v$ are
equivalent if and only if their normal forms $\nf{u}$ and $\nf{v}$ are equal.

We shall now present this algorithm in practice. A $1$-cell~$u\in\freecat{P_1}$
can be encoded as being either
\begin{itemize}
\item a non-empty sequence of elements of~$P_1$, or
\item an element~$x$ of~$P_0$, which we write \code{id($x$)}, representing the
  identity over~$x$.
\end{itemize}
In the following, we will not insist on the handling of identities and
assimilate those to empty lists in order to simplify the writing of algorithms.
We denote by $\wlen{u}$ the length of $u$, with \code{len(id($x$))} being $0$
by convention.
A natural number $i\in\N$ is a \emph{position} in~$u$ when $0\leq i<\wlen{u}$,
and in this case we write $\wlet ui$ for the $i$-th letter of~$u$. Given
$i,k\in\N$ such that $i$ and $i+k-1$ are positions in~$u$, we write $\wsub uik$
for the subword of~$u$ of length~$k$ starting at position~$i$, \ie the word
\[
  \wsub uik
  =
  \wlet ui\wlet u{i+1}\ldots\wlet u{i+k-1}.
\]
We say that a word $v$ \emph{matches}~$u$ at position~$i$ whenever~$v$ is a
subword of~$u$ starting at position~$i$. This can be tested with the following
first procedure:
\begin{lstlisting}
def matches($u$,$i$,$v$) =
  return (sub($u$,$i$,len($v$)) = $v$)
\end{lstlisting}
Given $1$-cells $u,v\in\freecat{P_1}$ such that \code{tgt($u$)} and
\code{src($v$)} are the same, we write \code{concat($u$,$v$)} for their
composition, which is simply the concatenation of the two sequences. More
generally, we allow ourselves to consider the composition
$\code{concat}(u_1,\ldots,u_k)$ of $k$ composable $1$-cells $u_1,\ldots,u_k$.

Given a rule $\alpha:u\To v$ in~$P_2$, we write $\rsrc(\alpha)$ for its source
$u$ and $\rtgt(\alpha)$ for its target~$v$. The normal form of a word~$u$ can be
computed with the following recursive procedure, expressed in a language which
should look familiar to anyone accustomed to imperative programming languages:
\\
\begin{minipage}{\textwidth}
\begin{lstlisting}
def rec normalize($P$,$u$) =
  for $\alpha\in P_2$ do
    $v$ = src($\alpha$)
    for $i$ = $0$ to len($u$)-len($v$) do
      if matches($u$,$i$,$v$) then
        $w_1$ = sub($u$,$0$,$i$)
        $w_2$ = tgt($\alpha$)
        $w_3$ = sub($u$,$i$+len($v$),len($u$)-len($v$)-$i$)
        return normalize($P$,concat($w_1$,$w_2$,$w_3$))
  return $u$
\end{lstlisting}
\end{minipage}
Finally, equality can be decided by the \emph{normal form algorithm} which can
be implemented as
\begin{lstlisting}
def equal($P$,$u$,$v$) =
  return (normalize($P$,$u$) = normalize($P$,$v$))
\end{lstlisting}

\subsection{Complexity of the word problem}
Let us mention the following result obtained by Avenhaus and
Madlener~\cite{avenhaus1977subrekursive,avenhaus1978algorithmische} for
presentations of groups, but the proof can be applied to presentation of
monoids.

\begin{theorem}
  Let $P$ and $Q$ be two Tietze equivalent finite 2-polygraphs. If the word
  problem can be decided for $P$ in time $O(f(n))$, then the word problem for
  $Q$ can be solved in time $O(f(cn))$ for some constant natural number $c>0$.
\end{theorem}

For a finite convergent $2$-poly\-graph~$P$, consider a function $f:\N\to\N$
such that for any $1$-cell $u$ of length~$n=\sizeof{u}$ in~$\freecat{P_1}$, the
leftmost reduction sequence from $u$ to its normal form contains at most~$f(n)$
many steps (here, the leftmost reduction means that we always reduce words with
a reduction which rewrites a subword as much on the left as
possible). In~\cite{book1982confluent}, Book proves that for a finite convergent
and reduced $2$\nbd-polygraph $P$, the normal form for $u$ in $P_1^\ast$ can be
computed in time~$O(n+f(n))$.
We say that a polygraph~$P$ is \emph{length-reducing} when for every rewriting
rule $\alpha:u\To v$ in $P_2$ we have $\sizeof{u}>\sizeof{v}$. As a consequence
of previous result, if a $2$-polygraph $P$ is length-reducing and confluent,
then its word problem is decidable in linear time.

\subsection{Other undecidable problems}
The word problem is far from being the only difficult one for
2-polygraphs~\cite{book1993string,miller1971group}. It is undecidable, given a
finite 2-polygraph~$P$, to determine whether 
\begin{itemize}
\item it is terminating, locally confluent or confluent,
\item there is a finite convergent polygraph presenting the same category,
\item it is presenting the terminal category, a finite category, a free
  category, a cancellative category, or a commutative monoid.
\end{itemize}

\section{Critical Branchings}
\label{sec:2-cb}
In order to show that a $2$-polygraph is confluent using the standard techniques
developed in \cref{sec:1-confl} (by using Newman's lemma, as stated in
\cref{lem:newman}), one has to check that all its local branchings are
confluent. Contrarily to the case of 1-polygraphs, a finite 2-polygraph usually
has an infinite number of local branchings. We however show here that it is
enough to check for the confluence of a finite subset of those, the critical
branchings.

\begin{example}
  \label{ex:aa-a}
  Consider the $2$-polygraph
  \[
  P=
  \Pres{\star}{a}{\alpha:aa\To a}
  \pbox.
  \]
  In order to verify that it is confluent, one can check that all the local
  branchings
  \[
    \xymatrix@!C=1ex@!R=1ex{
      &\ar@{=>}[dl]_-{a^i\alpha a^{n+1-i}}a^{n+3}\ar@{=>}[dr]^-{a^j\alpha a^{n+1-j}}&\\
      a^{n+2}\ar@{:>}[dr]_-{a^j\alpha a^{n-j}}&&a^{n+2}\ar@{:>}[dl]^-{a^i\alpha a^{n-i}}\\
      &a^{n+1}&
    }
  \]
  can be closed, for $n\in\N$ and $0\leq i<j\leq n+1$, where $a^n$ denotes the
  composite $aa\ldots a$ of $n$ instances of~$a$.
\end{example}

\noindent
In order to ease those checks, people are often interested in \emph{critical
  branchings}, which are minimal possible obstructions of confluence, which will
index a finite subset of the above confluence diagrams, enough to ensure the
confluence of the 2-polygraph. In the next
sections, we introduce a classification of local branchings in order to define
the critical ones. Such a classification in view of the critical branching lemma first appeared in \cite{knuth1970simple} for terms rewriting systems and~\cite{Nivat73} for string rewriting systems.

\subsection{Trivial branchings}
\label{sec:2-trivial-branching}
\index{trivial branching}
\index{branching!trivial}
A branching $(\phi_1,\phi_2):u\To(v_1,v_2)$ is \emph{trivial} when
$\phi_1=\phi_2$. Such a branching is always confluent since we can take
$w=v_1=v_2$ and $\phi'_1=\phi'_2=\id_w$ to close the diagram:
\[
  \xymatrix@!C=1ex@!R=1ex{
    &\ar@{=>}[dl]_{\phi_1}u\ar@{=>}[dr]^{\phi_2}&\\
    v_1\ar@{=}[dr]&&\ar@{=}[dl]v_2\pbox.\\
    &w&
  }
\]

\subsection{Orthogonal branchings}
\label{sec:2-orthogonal-branching}
\index{orthogonal branching}
\index{Peiffer branching}
\index{branching!orthogonal}
\index{branching!Peiffer}
A local branching $(\phi_1,\phi_2):u\To(v_1,v_2)$ is \emph{orthogonal} when it is of
the form
\begin{align*}
  u&=u_1vu_2wu_3,
  &
  \phi_1&=u_1\alpha u_2wu_3,
  &
  \phi_2&=u_1vu_2\beta u_3,
\end{align*}
for some words $u_1,v,u_2,w,u_3$ and rules $\alpha:v\To v'$ and $\beta:w\To w'$
(or of the symmetric form, obtained by exchanging the roles of~$\phi_1$
and~$\phi_2$). Such a branching is always confluent:
\[
  \xymatrix@!C=1ex@!R=1ex{
    &\ar@{=>}[dl]_{u_1\alpha u_2wu_3\quad}u_1vu_2wu_3\ar@{=>}[dr]^{\quad u_1vu_2\beta u_3}&\\
    u_1v'u_2wu_3\ar@{=>}[dr]_{u_1v'u_2\beta u_3\quad}&&\ar@{=>}[dl]^{\quad u_1\alpha u_2w'u_3}u_1vu_2w'u_3\pbox.\\
    &u_1v'u_2w'u_3&
  }
\]
Informally, it corresponds to rewriting two independent parts of the word~$u$.
Note that the above diagram corresponds precisely to the
equality~\eqref{eq:rewriting-exchange}. 
Orthogonal branchings are also
sometimes called \emph{Peiffer} branchings in reference to the corresponding notions for spherical diagrams in Cayley complexes associated to presentations of groups~\cite{lyndon2015combinatorial}.

\subsection{Overlapping branchings}
\label{sec:2-overlapping-branching}
\index{overlapping branching}
\index{branching!overlapping}
A local branching is \emph{overlapping} when it is not trivial nor orthogonal.

\subsection{Minimal branchings}
\label{sec:2-minimal-branching}
\index{branching!minimal}
We define a partial order on branchings by setting
$(\phi_1,\phi_2)\ctxleq(\phi_1',\phi_2')$ whenever the second branching can
be obtained by putting the first one in context. Formally, writing $v:x\to y$
for the source of the branching $(\phi_1,\phi_2)$, we have
$(\phi_1,\phi_2)\ctxleq(\phi_1',\phi_2')$ whenever there are words
$u:x'\to x$ and $w:y\to y'$ such that
$\phi_1'=u\phi_1w$ and $\phi_2'=u\phi_2w$.
In such a situation, the confluence of the first branching
$(\phi_1,\phi_2):v\To(v_1,v_2)$, say by $(\psi_1,\psi_2):(v_1,v_2)\To v'$,
implies the confluence of $(\phi_1',\phi_2')$, since we have
\[
  \xymatrix@!C=1ex@!R=1ex{
    &\ar@{=>}[dl]_{u\phi_1w}uvw\ar@{=>}[dr]^{u\phi_2w}&\\
    uv_1w\ar@{=>}[dr]_{u\psi_1w}&&\ar@{=>}[dl]^{u\psi_2w}uv_2w\pbox.\\
    &uv'w&
  }
\]
A branching is \emph{minimal} when it is minimal \wrt this order.

\subsection{Critical branchings}
\label{sec:critical-branchings}
\index{critical!branching}
\index{branching!critical}
A local branching is \emph{critical} when it is overlapping and minimal. The
following lemma is sometimes called the \emph{critical branching lemma}:

\begin{lemma}
  \label{CriticalBranchingLemma}
  \label{lem:2-cb}
  \label{L:CriticalBranchingLemma2Pol}
  A 2-polygraph is locally confluent if and only if all its critical branchings
  are confluent.
\end{lemma}
\begin{proof}
  Suppose given a local branching $(\phi_1,\phi_2)$. If this branching is
  critical, then it is confluent by hypothesis. Otherwise, it is either trivial,
  or orthogonal, or non-minimal. In the first two cases, we have seen that the
  branching is always confluent
  (\cref{sec:2-trivial-branching,sec:2-orthogonal-branching}). When the
  branching is non-minimal, it is greater (\wrt $\ctxleq$) than a critical
  branching (since those are the minimal ones), which is confluent by
  hypothesis, and we have seen that this implies the confluence of the branching
  (\cref{sec:2-minimal-branching}).
\end{proof}

\begin{remark}
  Note that the confluence of a branching $(\phi_1,\phi_2)$ also implies the
  confluence of the branching $(\phi_2,\phi_1)$. We could thus further reduce
  the number of critical branchings by considering them up to symmetry. We will
  refrain from doing so in the following, because it obfuscates the formulation
  of the algorithms without bringing significant improvements (it ``only''
  divides by two the number of critical branchings).
\end{remark}

\subsection{Classification of critical branchings}
\label{sec:cb-class}
Suppose given two rewriting steps
\begin{align*}
  u_1\alpha u_2
  :
  u_1uu_2
  &\To
  u_1u'u_2,
  &
  v_1\beta v_2
  :
  v_1vv_2
  &\To
  v_1v'v_2,
\end{align*}
with $\alpha:u\To u'$ and $\beta:v\To v'$, forming a local branching, \ie
such that $u_1uu_2=v_1vv_2$. We now study when such a branching is critical.

If the $1$-cells $u_1$ and $v_1$ are both non-identities, they are necessarily
of the form $u_1=wu_1'$ and $v_2=wu_2'$ for some non-identity $1$-cell~$w$ and
the branching is thus not minimal. We deduce that either $u_1$ or $v_1$ must be an
identity, and similarly either $u_2$ or $v_2$ must be an identity. Since moreover,
the branching should be overlapping, the situation must be of one of the four
forms given in \figr{cb-classification}, for some $1$-cells $w_1$, $w_2$, and
$w_3$.
\begin{figure}[t]
  \centering
\[
\begin{array}{|c|c|c|c|c|c|c|}
  \hline
  u_1&u&u_2&v_1&v&v_2&\text{Diagram}\\
  \hline
  \emptyword&w_1w_2&w_3&w_1&w_2w_3&\emptyword&
  \vxym{
    \\
    \ar@/^5ex/[rr]^{u'}\ar@{{}{ }{}}@/^2.5ex/[rr]|{\Uparrow\alpha}\ar[r]|{w_1}&\ar@/_5ex/[rr]_{v'}\ar@{{}{ }{}}@/_2.5ex/[rr]|{\Downarrow\beta}\ar[r]|{w_2}&\ar[r]|{w_3}&
    \\
  }
  \\
  \hline
  w_1&w_2w_3&\emptyword&\emptyword&w_1w_2&w_3&
  \vxym{
    \\
    \ar@/_5ex/[rr]_{v'}\ar@{{}{ }{}}@/_2.5ex/[rr]|{\Downarrow\beta}\ar[r]|{w_1}&\ar@/^5ex/[rr]^{u'}\ar@{{}{ }{}}@/^2.5ex/[rr]|{\Uparrow\alpha}\ar[r]|{w_2}&\ar[r]|{w_3}&
    \\
  }
  \\
  \hline
  \emptyword&w_1w_2w_3&\emptyword&w_1&w_2&w_3&
  \vxym{
    \\
    \ar@/^5ex/[rrr]^{u'}\ar@{{}{ }{}}@/^2.5ex/[rrr]|{\Uparrow\alpha}\ar[r]|{w_1}&\ar@/_5ex/[r]_{v'}\ar@{{}{ }{}}@/_2.5ex/[r]|{\Downarrow\beta}\ar[r]|{w_2}&\ar[r]|{w_3}&
    \\
  }
  \\
  \hline
  w_1&w_2&w_3&\emptyword&w_1w_2w_3&\emptyword&
  \vxym{
    \\
    \ar@/_5ex/[rrr]_{v'}\ar@{{}{ }{}}@/_2.5ex/[rrr]|{\Downarrow\beta}\ar[r]|{w_1}&\ar@/^5ex/[r]^{u'}\ar@{{}{ }{}}@/^2.5ex/[r]|{\Uparrow\alpha}\ar[r]|{w_2}&\ar[r]|{w_3}&
    \\
  }
  \\
  \hline
\end{array}
\]  
  \caption{Classification of critical branchings.}
  \label{fig:cb-classification}
\end{figure}
In the two first cases, we suppose that $w_2$ is not an identity (otherwise the
branching is orthogonal). We also suppose that the branching is not trivial,
\ie that we are not in a situation where $w_1$ and $w_3$ are identities and
$\alpha=\beta$. The last two cases are called \emph{inclusion} branchings
because of the relative positions of the rewriting rules as shown in the above
figures.

From the above classification, it should be clear that there is a simple
algorithm, which is detailed below, to compute critical branchings: for any two
rules $\alpha:u\To u'$ and $\beta:v\To v'$, we try to overlap $u$ and $v$ at
various offsets which are small enough to deduce the possible $w_1$, $w_2$, and
$w_3$ for the decompositions of the above form, and remove those which are
trivial. In particular, we have as a consequence:

\begin{lemma}
  \label{lem:2-cb-finite}
  Given a 2-polygraph $P$ with a finite set~$P_2$ of rewriting rules, the number
  of critical branchings is finite.
\end{lemma}

\begin{example}
  Consider the 2-polygraph
  \[
  P
  =
  \Pres{\star}{a,b,c}{\alpha:abc\To a, \beta:ca\To a}
  \pbox.
  \]
  In order to compute the critical branchings, we consider pairs of rules and
  examine how they can overlap. Suppose that we choose $\alpha$ and $\beta$. The
  relative positions of the source $abc$ of $\alpha$ and the source $ca$ of
  $\beta$ can be
  \[
  \begin{array}{c@{}c@{}c@{}c}
    &a&b&c\\
    c&a
  \end{array}
  \qquad
  \begin{array}{c@{}c@{}c@{}c}
    a&b&c\\
    c&a
  \end{array}
  \qquad
  \begin{array}{c@{}c@{}c@{}c}
    a&b&c\\
    c&a
  \end{array}
  \qquad
  \begin{array}{c@{}c@{}c@{}c}
    a&b&c\\
    &c&a
  \end{array}
  \qquad
  \begin{array}{c@{}c@{}c@{}c@{}c}
    a&b&c\\
    &&c&a
  \end{array}
  \]
  Among those, only the first and the last one are valid overlappings, which
  means that the vertically aligned letters are the same, giving rise to the two
  following critical branchings:
  \begin{align*}
    (c\alpha,\beta bc),
    &&
    (\alpha a,ab\beta),
  \end{align*}
  which can also be pictured as
  \begin{align*}
  \vxym{
    \\
    \ar@/_5ex/[rr]_a\ar@{{}{ }{}}@/_2.5ex/[rr]|{\Downarrow\beta}\ar[r]^c&\ar@/^5ex/[rrr]^a\ar@{{}{ }{}}@/^2.5ex/[rrr]|{\Uparrow\alpha}\ar[r]^a&\ar[r]^b&\ar[r]^c&\;,\\
    \\
  }
  &&
  \vxym{
    \\
    \ar@/^5ex/[rrr]^a\ar@{{}{ }{}}@/^2.5ex/[rrr]|{\Uparrow\alpha}\ar[r]^a&\ar[r]^b&\ar@/_5ex/[rr]_a\ar@{{}{ }{}}@/_2.5ex/[rr]|{\Downarrow\beta}\ar[r]^c&\ar[r]^a&\pbox.\\
    \\
  }    
  \end{align*}
  It can be checked that these are the only critical branchings (there is no
  critical branching involving $\alpha$ with $\alpha$ or $\beta$ with
  $\beta$). The first branching is confluent, but not the second one:
  \begin{align*}
    \xymatrix@!C=1ex@!R=1ex{
      &\ar@{=>}[dl]_{\beta bc}cabc\ar@{=>}[dr]^{c\alpha}&\\
      abc\ar@{=>}[dr]_\alpha&&\ar@{=>}[dl]^\beta ca\, ,\\
      &a&
    }
    &&
    \xymatrix@!C=1ex@!R=1ex{
      &\ar@{=>}[dl]_{\alpha a}abca\ar@{=>}[dr]^{ab\beta}&\\
      aa&&aba\pbox.\\
      &&
    }    
  \end{align*}
\end{example}

\begin{example}
  \label{ex:aa-a-confl}
  Consider the 2-polygraph of~\exr{aa-a} again:
  \[
    P=\Pres{\star}{a}{\alpha:aa\To a}
    \pbox.
  \]
  The only critical branching
  \[
    \xymatrix@!C=1ex@!R=1ex{
      &\ar@{=>}[dl]_{\alpha a}aaa\ar@{=>}[dr]^{a\alpha}&\\
      aa\ar@{=>}[dr]_\alpha&&\ar@{=>}[dl]^\alpha aa\\
      &a&
    }
  \]
  is confluent, therefore the 2-polygraph is locally confluent (by \cref{lem:2-cb}). Each rewriting
  step $u\alpha v:uaav\To uav$ has a source whose length is one less than the
  length of the source and therefore the system is terminating. Finally, we
  deduce that the 2-polygraph is convergent (by \cref{lem:newman} and \cref{sec:2-rprop}).
\end{example}

\begin{example}
  \label{ex:finv-confl}
  Given a set~$X$, the free group on this set can be presented, as a
  monoid, by the
  polygraph
  \[
    P=\Pres{\star}{a,\ol a}{\lambda:\ol aa\To 1,\rho:a\ol a\To 1}_{a\in X}
    \pbox.
  \]
  This polygraph is always locally confluent since the two critical branchings
  \begin{align*}
    \svxym{
      &\ar@2[dl]_{\rho a}a\ol aa\ar@2[dr]^{a\lambda}&\\
      a\ar@{=}[rr]&&a,
    }
    &&
    \svxym{
      &\ar@2[dl]_{\lambda a}\ol aa\ol a\ar@2[dr]^{a\rho}&\\
      \ol a\ar@{=}[rr]&&\ol a,
    }
  \end{align*}
  indexed by $a\in X$ are confluent. The rules decrease length and the polygraph
  is also terminating. Normal forms, also sometimes called \emph{reduced words},
  are words which do not contain factors of the form $\ol aa$ or $a\ol a$
  for some $a\in X$.
  This construction easily extends to present the free groupoid on a graph.
\end{example}

\begin{algo}
  \index{critical branching!algorithm}
  \index{algorithm!critical branching}
  \label{algo:2cp}
  The critical branchings of a 2-polygraph~$P$ can be computed thanks to the
  following algorithm, which tries to unify the sources of all pairs of rules
  $(\alpha,\beta)$ in~$P_2$. For simplicity, we suppose here that~$P$ is a
  presentation of a monoid, \ie there is only one possible identity 1-cell
  denoted \code{empty}. The output is a set of pairs
  $(u_1,\alpha,u_2),(v_1,\beta,v_2)$, with $\alpha:u'\To u''$ and $\beta:v'\To v''$
  forming a critical branching
  \[
    \xymatrix@C=10ex{
      u_1u''u_2&\ar@{=>}[l]_-{u_1\alpha u_2}u_1u'u_2=v_1v'v_2\ar@{=>}[r]^-{v_1\beta v_2}&v_1v''v_2\pbox.
    }
  \]
  The procedure in peudo-code is
  \program{cp}
\end{algo}

\section{Reduction Orders}
\label{sec:2-red-order}
In order to show that a 2-polygraph~$P$ is terminating, one has to consider a
well-founded order on 1-cells which is compatible with composition.

\begin{definition}
  \index{reduction!order}
  \index{order!reduction}
  A \emph{reduction order} on a category~$C$ is a partial
  order~$\preceq$ relating pairs of parallel morphisms\index{parallel!morphisms} in~$C$ which is
  \begin{itemize}
  \item well-founded: every weakly decreasing sequence of morphisms is
    eventually stationary,
  \item compatible with composition: for every morphisms $u:x'\to x$,
    $v,v':x\to y$ and $w:y\to y'$, we have that $v\succ v'$ implies
    $uvw\succ uv'w$.
  \end{itemize}
  Given a 2-polygraph~$P$, a reduction order $\preceq$ on~$\freecat{P_1}$ is
  said to be \emph{compatible} with the rules of~$P$ when $u\succ v$ for every
  rule $\alpha:u\To v$ in $P_2$. In this case, the order $\preceq$ is called a
  \emph{termination order}\index{termination!order}\index{order!termination} for~$P$.
\end{definition}

\begin{proposition}
  \label{prop:term-red-ord}
  A 2-polygraph~$P$ is terminating if and only if it admits a termination order.
\end{proposition}
\begin{proof}
  Suppose that $P$ is terminating. Then the following relation $\preceq$ is a
  reduction order compatible with~$P$, where given $1$-cells
  $u,v\in\freecat{P_1}$ we have $u\succeq v$ if and only if $u$ rewrites to~$v$.
  Conversely, in a $2$-polygraph equipped with a reduction order compatible
  with~$P$, every rewriting step is of the form $uvw\To uv'w$ for some rule
  $\alpha:v\To v'$. In such a situation, we have $v\succ v'$ because the order is
  compatible with the rules, and thus $uvw\succ uv'w$ because the order is
  compatible with composition. Therefore, the existence of an infinite sequence
  of reductions in~$P$ would imply the existence of an infinite decreasing
  sequence in the order, contradicting its well-foundedness.
\end{proof}

\begin{remark}
  \label{rem:2-id-nf}
  In a terminating 2-polygraph, an identity is necessarily a normal
  form. Namely, suppose that we have $\id_x\To u$ for some $x\in P_0$ and
  $u\neq\id_x$ in~$P_1$. Then we would have the infinite sequence of rewriting
  steps
  \[
  \id_x
  \To
  u
  =
  \id_xu
  \To
  uu
  =
  \id_xuu
  \To
  \ldots
  \]
\end{remark}

\subsection{Constructing reduction orders}
There is no general rule to construct a reduction order witnessing that a
2\nbd-po\-ly\-graph is terminating: in fact, deciding termination is an
undecidable problem~\cite[Section~2.5]{book1993string}. Fortunately, there is
however a ``standard toolbox'', which we now introduce, from which one is able
to construct orders in many useful cases.

\subsection{Reduction function}
\index{reduction!function}
\label{sec:red-fun}
The most usual method to show the termination of a 2-polygraph~$P$ is to provide
a function $f:P_1^*\to N$, called a \emph{reduction function}, where $(N,\leq)$
is a well-founded poset, which is compatible with composition: for every 1-cells
$u:x'\to x$, $v,v':x\to y$ and $w:y\to y'$, we have that $f(v)>f(v')$ implies
$f(uvw)>f(uv'w)$.
Such a reduction function induces a reduction order~$\preceq$ on~$\freecat P_1$
defined by $u\preceq v$ if and only if $f(u)\leq f(v)$.

A reduction function is \emph{monotone} when $f(u)>f(v)$ for every rule
$\alpha:u\To v$. With such a reduction function, the associated reduction order
is a termination order and thus, by \cref{prop:term-red-ord}:

\begin{lemma}
  \label{lem:red-fun}
  A polygraph equipped with a monotone reduction function is
  terminating.
\end{lemma}

\begin{example}
  \label{ex:length-order}
  The function which to every word $u$ associates its \emph{length}\index{length}~$\sizeof{u}$
  in~$\N$, also called its \emph{degree} in this context, is a reduction
  function. If a polygraph is \emph{length-decreasing}, in the sense that for
  every rule $\alpha:u\To v$ we have $\sizeof u>\sizeof v$, we can thus conclude
  that it is terminating by \cref{lem:red-fun}. This is for instance the
  argument we have been using to show termination in \exr{aa-a-confl}.
\end{example}

\begin{example}
  Functions other than length can also be useful. For instance, consider the
  presentation
  \[
    P=\Pres{\star}{a,b}{\alpha:ba\To aaa}
    \pbox.
  \]
  The 2-polygraph intuitively terminates because each application of a rules
  decreases the number of occurrences of $b$ in a word. In order to formalize
  this, consider the function~$f:P_1^*\to\N$ defined by $f(a)=0$, $f(b)=1$, and
  extended as a morphism of monoids, \ie $f(\emptyword)=0$ and
  $f(uv)=f(u)+f(v)$. We have $f(ba)=1>0=f(a)$ and the function is compatible
  with composition since it is a morphism of monoids. Therefore the 2-polygraph is
  terminating. Note that there are no critical branchings, therefore all of them
  are trivially confluent and the 2-polygraph is convergent.
\end{example}

\begin{remark}
  Instead of reduction functions, one could more generally consider the notion
  of a \emph{reduction 2-functor} which is a 2-functor $f:\freecat{P}\to N$
  where
  \begin{itemize}
  \item $N$ is a 2\nbd-category such that for each pair of objects $x,y\in N$
    the category $N(x,y)$ is a well-founded poset: there is at most one
    morphisms between two objects and every decreasing sequence is eventually
    stationary,
  \item $f$ is injective on parallel 1-cells: for every 1-cells $u,v:x\to y$
    in~$\freecat{P_1}$, $f(u)=f(v)$ implies $u=v$.
  \end{itemize}
\end{remark}

\subsection{Lexicographic product}
\index{lexicographic!product}
\index{product!lexicographic}
Given two posets $(M,\leq_M)$ and $(N,\leq_N)$, one can equip their product
$M\times N$ with an order $\leq_{M\times N}$, called the \emph{lexicographic
  product} of the two orders, such that $(m,n)\leq_{M\times N}(m',n')$ whenever
\begin{itemize}
\item $m<_Mm'$, or
\item $m=m'$ and $n\leq_N n'$.
\end{itemize}
When the two original orders are well-founded, their lexicographic product is
always well-founded. Namely, from every decreasing sequence
\[
  (m_0,n_0)\geq_{M\times N}(m_1,n_1)\geq_{M\times N}(m_2,n_2)\geq_{M\times N}\ldots
\]
the sequence of $(m_i)_{i\in\N}$ is decreasing with respect to $\leq_M$ and thus
eventually stationary, and similarly for the sequence $(n_i)_{i\in\N}$, and thus
the sequence of $(m_i,n_i)_{i\in\N}$ is also eventually stationary. As a
consequence, the lexicographic product of two reduction orders is a reduction
order.

\index{lexicographic!order}
\index{order!lexicographic}
Given a well-founded poset $(N,\leq)$ and $n\in\N$, a partial order~$\leq_n$ can
thus be defined on $N^n$ (the product of $n$ copies of~$N$) by induction on~$n$:
\begin{itemize}
\item $N^0$ is equipped with the trivial order,
\item $N^1=N$ is equipped with the order~$\leq$, and
\item $N^{n+1}=N\times N^n$ is equipped with the lexicographic product of $\leq$
  and $\leq_n$.
\end{itemize}
Finally, one can define a well-founded partial order~$\lexleq$ on
$\coprod_{n\in\N}N^n$, called the \emph{lexicographic order} induced by~$\leq$,
by $u\lexleq v$ whenever both $u$ and $v$ belong to~$N^n$ for some~$n\in\N$ and
$u\leq_n v$. For instance, given $N=\set{a,b}$ with $a\leq b$, one has
$abb\lexleq baa$ and $bbba\lexleq bbaa$.
This is easily adapted to the setting of polygraphs: given a 2\nbd-polygraph~$P$
and a well-founded partial order~$\leq$ on~$P_1$, one can define its
lexicographic extension as above, where $P_1^n$ is now the set of composable
sequences of length~$n$ of elements of~$P_1$. The resulting order $\lexleq$
on~$P_1^*$ is always a reduction order.

The variant of the lexicographic order where letters are compared from right to
left (instead of from left to right) is also useful and called the
\emph{colexicographic order}.

Note that, in the above definition of the lexicographic order, two words of
different lengths are always incomparable. One can define a variant of the
lexicographic order, sometimes called the \emph{dictionary
  order}\index{dictionary order}, which is such that $u\leq v$ whenever $u$ is a
prefix of $v$, or $u$ and $v$ admit respective prefixes $u'$ and $v'$, of the
same length, such that $u'\lexlt v'$. This order is not in general
well-founded, even if the order on the letters is. For instance, with $a<b$, one
has
\[
  b>ab>aab>aaab>\ldots
\]
Another variant of the lexicographic order, not suffering from this problem, is
presented in next section.

\subsection{Deglex order}
\label{sec:deglex}
\index{deglex order}
\index{order!deglex}
Suppose given a polygraph~$P$ equipped with a well-founded partial order
$\preceq_1$ on~$P_1$. We have seen in \cref{ex:length-order} that the length on
words is a reduction function and thus induces a reduction order
on~$\freecat P_1$, as explained in \cref{sec:red-fun}. Moreover, the
lexicographic order induced by~$\preceq_1$ is also a reduction order
on~$\freecat P_1$. By taking the lexicographic product of these two reduction
orders, we obtain a new reduction order $\preceq$ on $\freecat P_1$ called the
\emph{deglex order} associated to $\preceq_1$.
Explicitly, given two words $u=a_1\ldots a_m$ and $v=b_1\ldots b_n$, we have
$u\prec v$ whenever
\begin{itemize}
\item $m<n$, or
\item $m=n$, and there exists~$i$ with $1\leq i\leq n$ such that
  $a_i\preceq_1b_i$ and $a_j=b_j$ for every~$j<i$.
\end{itemize}

\begin{example}
  Consider the 2-polygraph
  \[
    P=\Pres{\star}{a,b}{\alpha:ab\To ba}
    \pbox.
  \]
  If we order the letters by $a>b$, the induced deglex order is a reduction
  order such that $ab>ba$. We can therefore apply \propr{term-red-ord} and
  deduce that the 2-polygraph is terminating. Since there is no critical
  branching, the 2-polygraph is convergent.
\end{example}

\subsection{Derivation}
\label{sec:derivation}
In order to construct a reduction function or a reduction order, one sometimes
needs to propagate information from the left or the right of the word. This idea
is nicely captured by the classical notion of derivation. For instance, consider
the 2-polygraph
\begin{equation}
  \label{eq:ex-mon-der}
  P=\Pres{\star}{a,b,c}{\alpha:cba\To aabc}
  \pbox.
\end{equation}
Most of the simple techniques above (considering the length, the number of
letters, or a deglex order) do not apply here to show the termination of the
rewriting system: for instance, with any deglex order we have $cba<aabc$ because
the second word is longer than the first one. However, one can justify the
termination of the rewriting system by noticing that a rewriting step always
decreases the number of ``occurrences of $c$ on the left of an occurrence
of~$b$''. The purpose of derivation is precisely to formulate such definitions
obtained by propagating information (here, the number of occurrences of~$c$) and
summing over each letter a quantity obtained from the propagated information
(here, the number of occurrence of $c$ on the left for each $b$ and $0$ for
each~$a$ or~$c$). For simplicity, we consider only the case of presentations of
monoids here, but it extends seamlessly to presentations of categories.

\index{bimodule}
\index{action!of a monoid!on a commutative monoid}
Given a monoid~$(M,\cdot,1)$, an \emph{$M$-bimodule~$N$} consists of a
commutative monoid~$(N,+,0)$ together with a function $M\times N\times M\to N$,
called an \emph{action} of $M$ on $N$, the image of a triple $(u,n,v)$ being
written~$u\cdot n\cdot v$, which is
\begin{itemize}
\item linear: for every $u,v\in M$ and $n,n'\in N$,
  \begin{align*}
    u\cdot(n+n')\cdot v
    &=
    u\cdot n\cdot v
    +
    u\cdot n'\cdot v,
    &
    u\cdot 0\cdot v
    &=
    0
    \pbox,
  \end{align*}
\item associative: for every $u',u,v,v'\in M$ and $n\in N$,
  \begin{align*}
    u'\cdot(u\cdot n\cdot v)\cdot v'
    &=
    (u'\cdot u)\cdot n\cdot(v\cdot v'),
    &
    1\cdot n\cdot 1&=n
    \pbox.
  \end{align*}
\end{itemize}
A \emph{derivation}\index{derivation} of $M$ with values in~$N$ is a function
$d:M\to N$ such that, for $u,v\in M$, one has
\[
  d(u\cdot v)
  =
  1\cdot d(u)\cdot v
  +
  u\cdot d(v)\cdot 1
  \qquad\text{and}\qquad
  d(1)=1.
\]
When $N$ is equipped with a partial order, the derivation is \emph{monotone}
when
\begin{itemize}
\item the addition is monotone: $n>n'$ implies $m_1+n+m_2>m_1+n'+m_2$ for every
  $m_1,n,n',m_2\in N$,
\item the action is monotone: $n>n'$ implies $u\cdot n\cdot v>u\cdot n'\cdot v$
  for every $u,v\in M$ and $n,n'\in N$.
\end{itemize}
In the following, we will be mostly interested in derivations in the case
where~$M$ is the monoid $\freecat P_1$ for some fixed $2$-polygraph~$P$ with one
$0$-generator. In such a situation, we say that a derivation $d$ is
\emph{adapted} to~$P$ when $d(u)>d(v)$ for every $2$-generator $\alpha:u\To v$
in $P_2$.

\begin{lemma}
  Suppose given a monoid $N$ equipped with a well-founded partial order,
  together with a structure of $P_1^*$-bimodule. A derivation $d:P_1^*\to N$ of
  $P_1^*$ with values in~$N$ which is monotone and adapted to~$P$, is a
  reduction function.
\end{lemma}
\begin{proof}
  Suppose given words $u,v,v',w$ such that $d(v)>d(v')$, we have
  \begin{align*}
    d(uvw)
    &=
    1\cdot d(u)\cdot vw
    +
    u\cdot d(v)\cdot w
    +
    uv\cdot d(w)\cdot 1
    \\
    &>
    1\cdot d(u)\cdot vw
    +
    u\cdot d(v')\cdot w
    +
    uv\cdot d(w)\cdot 1
    =
    d(uv'w)\pbox.\qedhere
  \end{align*}
\end{proof}

\noindent
Since $P_1^*$ is free, an action onto a monoid~$N$ is specified by its effect on
generators. Namely, any $P_1^*$-bimodule~$N$ induces, by restriction of the
action, two functions
\begin{align*}
  l:P_1\times N&\to N&r:N\times P_1&\to N\\
  (a,n)&\mapsto a\cdot n\cdot 1&(n,b)&\mapsto 1\cdot n\cdot b
\end{align*}
which satisfy for $a,b\in P_1$ and $n\in N$,
\begin{equation}
  \label{eq:lr-compat}
  r(l(a,n),b)
  =
  l(a,r(n,b)),
\end{equation}
both members of the equality being equal to $a\cdot n\cdot b$. Conversely, any
such pair of functions $l:P_1\times N\to N$ and $r:N\times P_1\to N$ satisfying
the above equality extend uniquely as an action. Similarly, a
derivation $d:P_1^*\to N$ is uniquely determined by the function $d:P_1\to N$
obtained as its restriction, and any such function extends uniquely as a
derivation.

The purpose of the action is intuitively to specify which information is
propagated sideways and the derivation determines how the propagated information
is used.

\begin{example}
  The termination of the rewriting system~\eqref{eq:ex-mon-der} can be shown as
  follows. Consider the monoid $(\N\times\N,+,(0,0))$ equipped with the
  componentwise addition, \ie $(m,n)+(m',n')=(m+m',n+n')$, and the partial order
  is the lexicographic of the standard order on~$\N$ by itself. The first
  component of~$d(u)$ will count the number of~$b$ in a word~$u$ and the second
  component the number of~$c$ before a~$b$. The action is given, for
  $(m,n)\in\N\times\N$, by
  \begin{align*}
    a\cdot(m,n)\cdot\emptyword&=(m,n),&\emptyword\cdot(m,n)\cdot a&=(m,n),\\
    b\cdot(m,n)\cdot\emptyword&=(m,n),&\emptyword\cdot(m,n)\cdot b&=(m,n),\\
    c\cdot(m,n)\cdot\emptyword&=(m,n+m),&\emptyword\cdot(m,n)\cdot c&=(m,n)\pbox.
  \end{align*}
  The first column specifies the function $l$ and the second one specifies~$r$,
  and those two functions are easily checked to be compatible in the sense that
  \eqref{eq:lr-compat} holds. The left equation on the last line can be read as:
  given a word~$u$ with $m$ letters $b$ and $n$ occurrences of $c$ before a~$b$,
  the word $au$ has $m$ letters $b$ and $n+m$ occurrences of $c$ before a~$b$;
  other equations are similar. The derivation $d:P_1^*\to\N\times\N$ is defined
  on generators by
  \begin{align*}
    d(a)&=(0,0),
    &
    d(b)&=(1,0),
    &
    d(c)&=(0,0)
    \pbox.
  \end{align*}
  The rule $\alpha$ is decreasing \wrt the derivation as above: we have
  \begin{align*}
    d(cba)
    &=
    1\cdot d(c)\cdot ba+c\cdot d(b)\cdot a+cb\cdot d(a)\cdot 1
    \\
    &=
    (0,0)+(1,1)+(0,0)=(1,1)
    \\
    \shortintertext{and}
    d(aabc)
    &=1\cdot d(a)\cdot abc+a\cdot d(a)\cdot bc+aa\cdot d(b)\cdot c+aab\cdot d(c)\cdot 1\\
    &=(0,0)+(0,0)+(0,0)+(0,0)
    =(0,0)
    \pbox.
  \end{align*}
  Therefore, the rewriting system is terminating (and convergent since it has no
  critical branching).
\end{example}

\section{Constructing Presentations of Categories}
\label{sec:2-presenting}
\label{sec:pres-method}
In order to show that a given category~$C$ is presented by a given
$2$-polygraph~$P$, one must show that the $1$-cells of~$C$ are in bijection with
equivalence classes of $1$\nbd-cells in~$\freecat P_1$ under the congruence
generated by relations in~$P_2$, see \lemr{2-pres-cond} for a formal
statement. Without further hypothesis on the polygraph this is usually difficult
because one has little control over the equivalence classes. However, in the
case where the polygraph~$P$ is convergent, equivalence classes of 1-cells have
normal forms as canonical representatives, which greatly simplifies the
situation. We explore this here by providing a method in order to show that a
given convergent polygraph is a presentation of a given category, based on the
observation that, in this case, \lemr{2-pres-cond} can be reformulated as
follows. We recall that $\fgf C$ denotes the underlying 1-polygraph of a
category~$C$, as defined in \cref{sec:1-upol}.

\begin{lemma}
  \label{lem:2-conv-pres-cond}
  A convergent 2-polygraph~$P$ is a presentation of a category~$C$ if and only
  if there is a morphism of 1-polygraphs $f:\tpol1 P\to \fgf C$ such that
  \begin{enumerate}
  \item $f_0:P_0\to C_0$ is a bijection between the $0$-generators and the
    objects of~$C$,
  \item for any 2-generator $\alpha:u\To v$ in~$P_2$, $f(u)=f(v)$,
  \item the function $\freecat{f_1}:\freecat{P_1}\to C_1$ restricts to a
    bijection between normal forms in~$\freecat{P_1}$ and $C_1$.
  \end{enumerate}
\end{lemma}

\noindent
In the following, we sometimes write $\intp{x}$ instead of $f(x)$ for the image
of a generator in~$P_0$ or~$P_1$ and call it the \emph{interpretation} of~$x$
in~$C$.

\begin{example}
  \label{ex:N2-rewr}
  \label{Example:MonoidWithTwoElements2}
  Let us show that the $2$-polygraph
  \[
    P=\Pres{\star}{a}{\alpha:aa\To 1}
  \]
  of~\cref{Example:MonoidWithTwoElements1} is presenting the monoid $\N/2\N$.
  The polygraph is terminating since the only rule decreases the length of the
  words and confluent since the only critical branching is confluent:
  \[
    \vxym{
      \ar@/_/@{=>}[d]_{\alpha a}aaa\ar@/^/@{=>}[d]^{a\alpha}\\
      a
    }
  \]
  We can thus apply \lemr{2-conv-pres-cond}. There is an obvious bijection
  between $P_0=\set{\star}$ and the only object of the monoid (recall that a
  monoid is considered as a category with only one object), and we define a morphism
  $f:\tpol 1P\to \N/2\N$ by interpreting the generator $a$ as $f(a)=1$. This
  morphism is compatible with the relation, since $f(aa)=1+1=0=f(1)$. Finally,
  the words of $\freecat P_1$ in normal form are $1$ and $a$, and those are in
  bijection with the elements of~$\N/2\N$, allowing us to conclude.
\end{example}

\begin{example}
  \index{isomorphism!walking}
  \label{ex:walking-iso-nf}
  Let us use \cref{lem:2-conv-pres-cond} in order to give a simpler construction
  of the presentation of the category of \exr{walking-iso}. We want to show
  that the posetal category
  \[
    C
    =
    \vxym{
      X\ar@(dl,ul)^{\id_X}\ar@/^/[r]^{F}&\ar@/^/[l]^{G}Y\ar@(dr,ur)_{\id_G}
    }
  \]
  is presented by the 2-polygraph
  \[
    P
    =
    \Pres{x,y}{a:x\to y,b:y\to x}{\alpha:ab\To\id_x,\beta:ba\To\id_y}\pbox.
  \]
  The polygraph is convergent since the rules decrease the length and the two
  critical branchings are confluent:
  \begin{align*}
    \xymatrix@!C=1ex@!R=1ex{
      &\ar@{=>}[dl]_{\alpha a}aba\ar@{=>}[dr]^{a\beta}&\\
      a\ar@{=}[dr]&&\ar@{=}[dl]a,\\
      &a&
    }
    &&
    \xymatrix@!C=1ex@!R=1ex{
      &\ar@{=>}[dl]_{\beta b}bab\ar@{=>}[dr]^{b\alpha}&\\
      b\ar@{=}[dr]&&\ar@{=}[dl]b\pbox.\\
      &b&
    }    
  \end{align*}
  We define a morphism of 1-polygraphs $P'\to \fgf C$ by
  \begin{align*}
    \intp{x}&=X,
    &
    \intp{y}&=Y,
    &
    \intp{a}&=F,
    &
    \intp{b}&=G,
  \end{align*}
  which obviously induces a bijection between $P_0=\set{x,y}$ and
  $C_0=\set{X,Y}$. Finally, the $1$-cells in~$\freecat{P_1}$ which are in normal
  form are the words over the alphabet~$\set{a,b}$ which do not contain $aa$
  nor $bb$ (because~$a$ cannot be composed with~$a$, and similarly for~$b$) nor
  $ab$ nor $ba$ (because the word would not be in normal form since the rule
  $\alpha$ or $\beta$ would apply) as a factor. Thus, there are four normal forms
  in~$\freecat{P_1}$: $\id_x$, $\id_y$, $a$, and $b$. They are respectively sent
  by~$f$ to~$\id_X$, $\id_Y$, $F$, and $G$, and thus we have a bijection between
  normal forms and $1$-cells of~$C$. We conclude that~$P$ is a presentation
  of~$C$, \ie $C\isoto\prescat{P}$.
\end{example}

\begin{remark}
  \index{canonical form}
  \label{rem:2-pres-cf}
  Note that the method given by \lemr{2-conv-pres-cond} would work with any
  notion of ``canonical form'' for the elements of~$\freecat{P_1}$ modulo $\approx$,
  not necessarily corresponding to the normal forms for a rewriting system. Namely, for a
  2-polygraph~$P$ and a category~$C$, suppose given a morphism of 1-polygraphs
  $f:\tpol1 P\to \fgf C$ satisfying the two first conditions of
  \lemr{2-conv-pres-cond} and a set~$N\subseteq\freecat{P_1}$, whose elements
  are called \emph{canonical forms}, such that
  \begin{itemize}
  \item every element of~$\freecat{P_1}$ is equivalent to an element of~$N$,
  \item $f$ induces a bijection between~$N$ and $C_1$,
  \end{itemize}
  then $P$ is a presentation of~$C$. Note that the second condition ensures that
  every element of~$\freecat{P_1}$ is equivalent to a unique canonical form.
\end{remark}

\subsection{The standard presentation}
\index{standard!presentation}
\index{presentation!standard}
\label{sec:cat-std-pres}
Any category~$C$ admits a convergent presentation~$P$, called the \emph{standard
  presentation}, introduced in \secr{2-std-pres}, which is defined by 
\begin{itemize}
\item $P_0$ is the set of $0$-cells of~$C$,
\item $P_1$ contains a $1$-generator $\rep f:x\to y$ for every $1$-cell $f:x\to y$ in~$C$,
\item $P_2$ contains 2-cells of the form
\[
  \eta:\rep{\id_x}\To{\id_x}:x\to x
  \qquad\text{and}\qquad
  \mu_{f,g}:\rep f\rep g\To\rep{fg} : x\to z,
\]
  which can be represented as
\[
  \vxym{
    \\
    x\ar@/_3ex/[rr]_{{\id_x}}\ar@/^3ex/[rr]^{\rep{\id_x}}\ar@{}[rr]|{\phantom{\eta_x}\Longdownarrow\eta_x}&&x
  }
  \qquad\text{and}\qquad
  \vxym{
    &y\ar@/^/[dr]^{\rep g}&\\
    x\ar@/^/[ur]^{\rep f}\ar@/_3ex/[rr]_{\rep{fg}}\ar@{{}{ }{}}@/^2.5ex/[rr]|{\phantom{\mu_{f,g}}\Longdownarrow\mu_{f,g}}&&z,
  }
\]
  for every $0$-cell~$x$ of~$C$ and pair of composable $1$-cells $f:x\to y$ and
  $g:y\to z$ in~$C$.
\end{itemize}
Above, note the subtle distinction between $\rep f\rep g$ and $\rep{fg}$: the source of
the 2\nbd-cell~$\mu_{f,g}$ is a path of length 2 (consisting of the edges $\rep f$
and $\rep g$), whereas its target is a path of length 1 (consisting of the edge
$\rep{fg}$, the generator associated to the composite $1$-cell $fg$).
Similarly, a 2-cell $\eta_x$ has a path of length 1 (\resp 0) as source (\resp target).

The proof that~$P$ presents~$C$ can be performed using the above method.
The 2-polygraph $P$ is terminating because the rules decrease the length of the
1-cells. It is also convergent: its critical branchings are of the form
\begin{align*}
  \xymatrix@!C=1ex@!R=1ex{
    &\ar@{=>}[dl]_{\mu_{f,g}\rep h}\rep f\rep g\rep h\ar@{=>}[dr]^{\rep f\mu_{g,h}}&\\
    \rep{fg}\rep h\ar@{:>}[dr]_{\mu_{fg,h}}&&\ar@{:>}[dl]^{\mu_{f,gh}}\rep f\rep{gh}\,,\\
    &\rep{fgh}&
  }
  &&
  \xymatrix@!C=1ex@!R=1ex{
    &\ar@{=>}[dl]_{\mu_{\id_x,a}}\rep{\id_x}\rep f\ar@{=>}[dr]^{\eta_x\rep f}&\\
    \rep f\ar@{:}[dr]&&\ar@{:}[dl]\rep f\,,\\
    &\rep f&
  }
  &&
  \xymatrix@!C=1ex@!R=1ex{
    &\ar@{=>}[dl]_{\mu_{f,\id_y}}\rep f\rep{\id_y}\ar@{=>}[dr]^{\rep f\eta_y}&\\
    \rep f\ar@{:}[dr]&&\ar@{:}[dl]\rep f\,,\\
    &\rep f&
  }
\end{align*}
for some composable $1$-cells $f:x\to y$, $g:y\to z$ and $h:z\to t$ of~$C$, and
are thus confluent. The normal forms are either empty paths (identities) or
paths of length 1 consisting of a $1$-cell $a$ which is not an identity.
Finally, we define a morphism of 1-polygraphs $\tpol1 P\to\fgf C$ such that the
function on objects $P_0\to C_0$ is the identity and the function on morphisms
$P_1\to C_1$ is the identity, which is obviously compatible with the relations
in~$P_2$. This functor clearly induces a bijection between normal forms and
$1$-cells of~$C$.

A variant where the orientation of the $2$-generators $\eta_x$ is reversed, \ie
$\eta_x:{\id_x}\To\rep{\id_x}$, is more commonly found in the literature. It is
Tietze equivalent to the above one, and thus also a presentation of the
category~$C$, although not a convergent one since identities are not normal
forms (see \cref{rem:2-id-nf}).

\subsection{The simplicial category}
\label{sec:simpl-cat}
\label{sec:simpl-pres}
As a concrete, non-trivial, and useful example, we recall here the well-known
presentation of the \emph{augmented simplicial category}~$\Simplaug$,
\index{simplicial!category}
\index{category!simplicial}
as given in~\cite[Proposition~VII.5.2]{MacLane98}. Its objects are
natural numbers $n\in\N$ and a morphism $f:m\to n$ is a weakly increasing
function $f:\intset{m}\to\intset{n}$, where $\intset{n}$ denotes the finite
ordinal $\set{0,\ldots,n-1}$. We claim that this category admits a presentation
by the 2-polygraph with
\begin{itemize}
\item 0-generators: natural numbers $n\in\N$,
\item 1-generators: for $n\in\N$,
  \[
  s^n_i:n+1\to n,
  \]
  with $0\leq i<n$, and
  \[
  d^n_i:n\to n+1,
  \]
  with $0\leq i\leq n$,
\item 2-generators:
  \begin{align*}
    &\sigma:&s_i^{n+1}s_j^n&\To s_{j+1}^{n+1}s_i^n&&\text{for $0\leq i\leq j<n$,}\\
    &\delta:&d_j^nd_i^{n+1}&\To d_i^nd_{j+1}^{n+1}&&\text{for $0\leq i\leq j\leq n$,}\\
    &\gamma:&d_i^{n+1}s_j^{n+1}&\To s_{j-1}^nd_i^n&&\text{for $0\leq i<j\leq n$,}\\
    &&&\To \id_n&&\text{for $i=j$ or $i=j+1$,}\\
    &&&\To s_j^nd_{i-1}^n&&\text{for $0\leq j+1<i\leq n+1$.}
  \end{align*}
\end{itemize}
We consider the order on generators such that, for $i,j,m,n\in\N$, we have
\begin{itemize}
\item $s_i^n\geq s_j^n$ for $i<j$,
\item $d_j^n\geq d_i^n$ for $j\geq i$,
\item $d_i^n\geq s_j^m$.
\end{itemize}
This order is easily shown to be well-founded, and all the rules are strictly
decreasing according to the associated deglex order. By
\cref{prop:term-red-ord}, the polygraph is thus terminating. For simplicity,
from now on, we omit the superscripts from generators.

The critical branchings of the rewriting system are%
\begin{longtable}{cc}
  $
  \vcenter{\xymatrix@C=0ex@!R=2ex{
    &\ar@{=>}[dl]s_is_js_k\ar@{=>}[dr]&\\
    s_{j+1}s_is_k\ar@{:>}[d]&&\ar@{:>}[d]s_is_{k+1}s_j\\
    s_{j+1}s_{k+1}s_i\ar@{:>}[dr]&&\ar@{:>}[dl]s_{k+2}s_is_j\\
    &s_{k+2}s_{j+1}s_i
  }}
  $
  &
  $
  \vcenter{\xymatrix@C=0ex@!R=2ex{
    &\ar@{=>}[dl]d_kd_jd_i\ar@{=>}[dr]&\\
    d_jd_{k+1}d_i\ar@{:>}[d]&&\ar@{:>}[d]d_kd_id_{j+1}\\
    d_jd_id_{k+2}\ar@{:>}[dr]&&\ar@{:>}[dl]d_id_{k+1}d_{j+1}\\
    &d_id_{j+1}d_{k+2}&\\
  }}
  $
  \\
  for $i\leq j\leq k$
  &
  for $i\leq j\leq k$
  \\
  $
  \vcenter{\xymatrix@C=1ex@R=3ex{
    &\ar@{=>}[dl]d_jd_is_k\ar@{=>}[dr]&\\
    d_id_{j+1}s_k\ar@{:>}[d]&&\ar@{:>}[d]d_js_{k-1}d_i\\
    d_is_{k-1}d_{j+1}\ar@{:>}[dr]&&\ar@{:>}[dl]s_{k-2}d_jd_i\\
    &s_{k-2}d_id_{j+1}
  }}
  $
  &
  $
  \vcenter{\xymatrix@C=1ex@R=3ex{
    &\ar@{=>}[dl]d_jd_is_k\ar@{=>}[dr]&\\
    d_id_{j+1}s_k\ar@{:>}[dr]&&\ar@{:>}[dl]d_js_{k-1}d_i\\
    &d_i
  }}
  $
  \\
  for $i\leq j<k-1$
  &
  for $i<k$, and $j=k-1$ or $j=k$
  \\
  $
  \vcenter{\xymatrix@C=1ex@R=3ex{
    &\ar@{=>}[dl]d_jd_is_k\ar@{=>}[dr]&\\
    d_id_{j+1}s_k\ar@{:>}[d]&&\ar@{:>}[d]d_js_{k-1}d_i\\
    d_is_kd_j\ar@{:>}[dr]&&\ar@{:>}[dl]s_{k-1}d_{j-1}d_i\\
    &s_{k-1}d_id_j\\
  }}
  $
  &
  $
  \vcenter{\xymatrix@C=1ex@R=3ex{
    &\ar@{=>}[dl]d_jd_is_k\ar@{=>}[dd]\\
    d_id_{j+1}s_k\ar@{:>}[dr]&\\
    &d_k
  }}
  $
  \\
  for $i<k<j$
  &
  for $i=j=k$
  \\
  $
  \vcenter{\xymatrix@C=1ex@R=3ex{
    &\ar@{=>}[dl]d_jd_is_k\ar@{=>}[ddd]&\\
    d_id_{j+1}s_k\ar@{:>}[d]\\
    d_is_kd_j\ar@{:>}[dr]&\\
    &d_k
  }}
  $
  &
  $
  \vcenter{\xymatrix@C=1ex@R=3ex{
    &\ar@{=>}[dl]d_jd_is_k\ar@{=>}[dr]&\\
    d_id_{j+1}s_k\ar@{:>}[d]&&\ar@{:>}[d]d_js_kd_{i-1}\\
    d_is_kd_j\ar@{:>}[dr]&&\ar@{:>}[dl]s_kd_{j-1}d_{i-1}\\
    &s_kd_{i-1}d_j\\
  }}
  $
  \\
  for $i=k$ or $i=k+1$, and $k<j$
  &
  for $k+1<i\leq j$
\end{longtable}
The 1-generators
\[
s_i^n:n+1\to n
\quad\text{and}\quad
d_i^n:n\to n+1
\]
of~$P_1$ are respectively interpreted as the morphisms
\[
\intp{s_i^n}:n+1\to n
\quad\text{and}\quad
  \intp{d_i^n}:n\to n+1
\]
of~$\Simplaug$, which are the functions defined by
\[
  \intp{s_i^n}(k)
  =
  \begin{cases}
    k&\text{if $0\leq k\leq i$,}\\
    k-1&\text{if $i<k\leq n$,}
  \end{cases}
\quad\text{and}\quad
  \intp{d_i^n}(k)
  =
  \begin{cases}
    k&\text{if $0\leq k<i$,}\\
    k+1&\text{if $i\leq k<n$.}
  \end{cases}  
\]
For instance, the graphs of $\intp{s^4_2}$ and $\intp{d^3_2}$ are respectively
\[
\xymatrix@R=1ex{
  \ar@{}[r]|4&\cdot\ar@{-}[dr]&\\
  \ar@{}[r]|3&\cdot\ar@{-}[dr]&\cdot&\ar@{}[l]|3\\
  \ar@{}[r]|2&\cdot\ar@{-}[r]&\cdot&\ar@{}[l]|2\\
  \ar@{}[r]|1&\cdot\ar@{-}[r]&\cdot&\ar@{}[l]|1\\
  \ar@{}[r]|0&\cdot\ar@{-}[r]&\cdot&\ar@{}[l]|0
}
\qquad\qquad
\xymatrix@R=1ex{
  \\
  &&\cdot&\ar@{}[l]|3\\
  \ar@{}[r]|2&\cdot\ar@{-}[ur]&\cdot&\ar@{}[l]|2\\
  \ar@{}[r]|1&\cdot\ar@{-}[r]&\cdot&\ar@{}[l]|1\\
  \ar@{}[r]|0&\cdot\ar@{-}[r]&\cdot&\ar@{}[l]|0
}
\]
It can be checked that the interpretation is compatible with the relations
in~$P_2$. For instance, for the rule $d_jd_i\To d_id_{j+1}$, with
$0\leq i\leq j$, we have
\begin{align*}
  \intp{d_jd_i}(k)
  &
  =
  \intp{d_i}\circ\intp{d_j}(k)
  =
  \begin{cases}
    \intp{d_i}(k)&\text{if $k<j$,}\\
    \intp{d_i}(k+1)&\text{if $j\leq k$,}
  \end{cases}
  \\
  &
  =
  \begin{cases}
    k&\text{if $k<i\leq j$,}\\
    k+1&\text{if $i\leq k<j$,}\\
    k+2&\text{if $i\leq j\leq k$,}
  \end{cases}
  =
  \begin{cases}
    \intp{d_{j+1}}(k)&\text{if $k<i$,}\\
    \intp{d_{j+1}}(k+1)&\text{if $i\leq k$,}
  \end{cases}
  \\
  &=
  \intp{d_{j+1}}\circ\intp{d_i}(k)
  =
  \intp{d_id_{j+1}}(k)
  \pbox.
\end{align*}
The normal forms are of the form
\begin{equation}
  \label{eq:delta-nf}
  s_{i_0}^{n+p}s_{i_1}^{n+p-1}\ldots s_{i_p}^nd_{j_0}^nd_{j_2}^{n+1}\ldots d_{j_q}^{n+q},
\end{equation}
for $n,p,q\in\N$, and
\begin{align*}
  n+p>i_0>i_1>\ldots>i_p\geq 0,
  &&
  0\leq j_0<j_1<\ldots<j_q\leq n+q
  \pbox.
\end{align*}
Namely, the rule $\gamma$ imposes that there is no $d_j$ before a $s_i$,
$\sigma$ (\resp $\delta$) imposes that the indices of successive $s_i$ (\resp
$d_i$) are increasing (\resp decreasing).

Every morphism $f:m\to m'$ of $\Simplaug$ is the interpretation of exactly one
such a normal form. We can namely observe that, $f$ being a weakly increasing
function, it is uniquely determined by
\begin{itemize}
\item the set of ``merged'' elements, \ie the set
  \[
  \set{i_0,i_1,\ldots,i_p}
  \quad\subseteq\quad
  \intset{m}
  \]
  of elements such that $f(i_k)=f(i_k+1)$, and
\item its image, or equivalently its complement, \ie the set
  \[
  \set{j_0,j_1,\ldots,j_q}
  \quad\subseteq\quad
  \intset{m'}
  \]
  of elements $j_k$ which are not in the image of~$f$.
\end{itemize}
Finally, with the above notations, and writing $n=m-p=m'-q$, it is easily
checked that $f$ is precisely the interpretation of the normal form
\eqref{eq:delta-nf}.


\section{Residuation}
\label{sec:2-residuation}
The notion of residual, which intuitively specifies what ``remains'' of a
morphism after another one, provides a powerful tool in order to derive
properties of a presented category, from combinatorial properties of its
presentation. Namely, by studying the properties of residuals, through rewriting
systems, one is often able to show interesting properties of the presented
category such as the existence of pushouts, the fact that morphisms are mono, or
that it embeds into its enveloping groupoid. The exposition provided here is
adapted from classical techniques in rewriting theory originating in Lévy's
thesis~\cite{levy1978reductions, huet1991computations},
see~\cite[Section~8.7]{bezem2003term} in the context of term rewriting systems,
\cite{mellies2002axiomatic} for a modern presentation,
\cite[Section~II.4]{dehornoy2015foundations} in the context of presentations of
groups, and~\cite{clerc2015presenting} of which the current presentation is
inspired.

\subsection{Residuation structure}
\index{residuation structure}
A \emph{residuation structure} on a category is a function which to every pair
of coinitial morphisms
$f:x\to y_1$ and $g:x\to y_2$
associates a morphism
\[
  f/g: y_2\to z,
\]
called the \emph{residual} of~$f$ after $g$, satisfying the three
following conditions.
\begin{enumerate}
\item The morphisms $f/g$ and $g/f$ are cofinal and satisfy
  \[
    f(g/f)= g(f/g),
  \]
  \ie
  \[
    \xymatrix@C=3ex@R=3ex{
      &\ar[dl]_fx\ar[dr]^g&\\
      y_1\ar@{.>}[dr]_{g/f}&&\ar@{.>}[dl]^{f/g}y_2\pbox.\\
      &z&
    }
  \]
\item Residuation is compatible with composition: given a morphism $f:x\to y$
  and morphisms $g:x\to z$ and $h:z\to z'$,
  \begin{align*}
    f/\unit{x}&=f,
    &
    f/(gh)&=(f/g)/h,
    \\
    \unit{x}/f&=\unit{y},
    &
    (gh)/f&=(g/f)(h/(f/g))
    \pbox,
  \end{align*}
  \ie
  \[
    \xymatrix{
      x\ar[d]_f\ar[r]^{\unit x}&x\ar@{.>}[d]^f\\
      y\ar@{.>}[r]_{\unit y}&y
    }
    \qquad\qquad
    \xymatrix{
      x\ar[d]_f\ar[r]^g&z\ar@{.>}[d]|{f/g}\ar[r]^h&z'\ar@{.>}[d]^{(f/g)/h}\\
      y\ar@{.>}[r]_{g/f}&y'\ar@{.>}[r]_{h/(f/g)}&y''
    }
  \]
\item Self-residuation is trivial: for every morphism $f:x\to y$,
  \[
    f/f
    =
    \unit{y}
  \]
  \ie
  \[
    \xymatrix@C=3ex@R=3ex{
      &\ar[dl]_fx\ar[dr]^f&\\
      y\ar@{.>}[dr]_{\unit{y}}&&\ar@{.>}[dl]^{\unit{y}}y\pbox.\\
      &y&
    }
  \]
\end{enumerate}
A residuation structure thus provides a witness of confluence for branchings,
which is compatible with the categorical structure.

\begin{remark}
  A residuation structure is precisely a distributive law
  \[
    \ell
    :
    C^\op\otimes C
    \to
    C\otimes C^\op
  \]
  as developed \secr{cat-dlaw}, such that for every morphism $f:x\to y$ in~$C$ we have
  $\ell(f^\op,f)=(\unit{y},\unit{y})$, see also~\secr{qspan}.
\end{remark}

\begin{proposition}
  \label{prop:residual-epi}
  In a category equipped with a residuation structure, every morphism is epi.
\end{proposition}
\begin{proof}
  We show that a morphism~$f:x\to y$ is necessarily epi. Suppose given morphisms
  $g,h:y\to z$ such that $fg=fh$. We have
  \[
    (fg)/f
    =
    (f/f)(g/(f/f))
    =
    \unit{y}(g/\unit{y})
    =
    g
    \pbox.
  \]
  Thus,
  \[
    g=(fg)/f=(fh)/f=h
  \]
  and the morphism~$f$ is epi.
\end{proof}

\begin{proposition}
  \label{prop:residual-pushout}
  In a category equipped with a residuation structure, every pair of coinitial
  morphisms admits a pushout.
\end{proposition}
\begin{proof}
  Given coinitial morphisms $f:x\to y_1$ and $g:x\to y_2$, we claim that the
  morphisms $g/f:y_1\to z$ and $f/g:y_2\to z$ form a pushout cocone. Suppose
  given morphisms $f':y_1\to z'$ and $g':y_2\to z'$ such that $ff'=gg'$:
  \[
    \xymatrix@C=3ex@R=3ex{
      &\ar[dl]_fx\ar[dr]^g&\\
      y_1\ar@/_/[ddr]_{f'}\ar[dr]|{g/f}&&\ar[dl]|{f/g}\ar@/^/[ddl]^{g'}y_2\\
      &z\ar@{.>}[d]|h&\\
      &z'&
    }
  \]
  The morphism $h=f'/(g/f)$ makes the two triangles commute. Namely, we have
  \[
    (g/f)/f'=g/(ff')=g/(gg')=\unit{y_2}/g'=\unit{z'}
  \]
  from which follows the commutation of the left triangle:
  \[
    (g/f)h=(g/f)(f'/(g/f))=f'((g/f)/f')=f'\unit{z'}=f'
    \pbox.
  \]
  Moreover, we have
  \[
    h=f'/(g/f)=(ff')/(f(g/f))=(gg')/(g(f/g))=g'/(f/g)
  \]
  from which we deduce that the right triangle commutes as above, by exchanging
  the roles of $f$ and $g$:
  \[
    (f/g)h=(f/g)(g'/(f/g))=g'((f/g)/g')=g'
    \pbox.
  \]
  Conversely, given a morphism $h:z\to z'$ such that $(g/f)h=f'$ and
  $(f/g)h=g'$, we necessarily have
  \[
    h=((g/f)h)/(g/f)=f'/(g/f)
    \pbox.\qedhere
  \]
\end{proof}

\begin{proposition}
  Suppose given a category~$C$ such that both~$C$ and~$C^\op$ are equipped with
  a residuation structure. Then the canonical functor~$C\to\freegpd{C}$,
  from~$C$ to its enveloping groupoid, is faithful.
\end{proposition}
\begin{proof}
  \index{calculus of fractions}
  Because~$C$ admits a residuation structure, the collection~$W$ of all
  morphisms of~$C$ forms a \emph{calculus of left fractions} in the sense
  of~\cite{gabriel1967calculus}:
  \begin{itemize}
  \item this collection contains identities and is closed under composition,
  \item for any pair of coinitial morphisms~$f$ and $g$ there exists morphisms
    $f'$ and~$g'$ such that $ff'=gg'$ (namely, we can take $f'=g/f$ and
    $g'=f/g$),
  \item for any morphisms $h:x\to y$ and $f,g:y\to z$ such that $hf=hg$, there
    exists a morphism $h':z\to z'$ such that $fh'=gh'$:
    \[
      \xymatrix@C=3ex@R=3ex{
        x\ar[r]^{h}&y\ar@<.5ex>[r]^f\ar@<-.5ex>[r]_g&z\ar@{.>}[r]^{h'}&z'
      }
    \]
    (namely, by \propr{residual-epi} every morphism of~$C$ is epi and we can
    take~$h'=\unit{z}$).
  \end{itemize}
  The enveloping groupoid~$\freegpd{C}$ can thus be described as the category of
  left fractions $\loc CW$. Since~$C^\op$ admits a residuation structure, by
  \propr{residual-epi}, every morphism of~$C$ (and thus of~$W$) is mono, and in
  this case, the canonical functor~$C\to\loc CW$ is easily shown to be faithful.
\end{proof}

\subsection{Residuated presentation}
In practice, it is difficult to directly exhibit a residuation structure on a
category~$C$ and show that it satisfies the required axioms. We provide here a
general methodology in order to show that~$C$ admits a residuation structure in
the case where it is equipped with a presentation by a 2-polygraph. Namely, in
this case, we can specify the residuation structure on generators and extend it
to other morphisms by functoriality.

\index{residuated presentation}
\index{presentation!residuated}
A \emph{residuated presentation}~$P$ is a 2-polygraph together with, for every
pair of coinitial generators
$a:x\to y_1$ and $b:x\to y_2$
in~$P_1$, a morphism
\[
  a/b
  :
  y_1\to z
\]
in~$\freecat P_1$, in such a way that
\begin{itemize}
\item $a/b$ and $b/a$ have the same target,
\item the morphisms $a(b/a)$ and $b(a/b)$ are $P$-congruent:
  \begin{equation}
    \label{eq:res-fill-ax}
    \svxym{
      &\ar[dl]_ax\ar[dr]^b&\\
      y_1\ar[dr]_{b/a}&\overset*\Leftrightarrow&\ar[dl]^{a/b}y_2\\
      &z&
    }
  \end{equation}
\item for every 1-generator $a\in P_1$, we have $a/a=\unit{}$,
\item for every generators
  $
  a:x\to x'
  $ and $
  \alpha:u\To v:x\to y
  $,
  respectively in~$P_1$ and $P_2$, we have
  \[
    a/u= a/v
  \]
  and there is a 2-generator
  \[
    \alpha/a
    :
    u/a\To v/a
    :
    x'\to y'
  \]
  in~$P_2$:
  \begin{equation}
    \label{eq:res-cyl-ax}
    \vcenter{
      \xymatrix@C=8ex@R=4ex{
        x\ar[d]_a\ar@/^/[r]^u\ar@/_/[r]_v\ar@{}[r]|{\Downarrow\alpha}&y\ar@{.>}[d]^{a/u=a/v}\\
        x'\ar@/^/@{.>}[r]^{u/a}\ar@/_/@{.>}[r]_{v/a}\ar@{}[r]|{\Downarrow\alpha/a}&y'\pbox.
      }
    }
  \end{equation}
\end{itemize}

\subsection{Residuation of morphisms}
\label{sec:mor-res}
Suppose fixed a residuated presentation~$P$. We can extend the residuation
operation in order to define the residual $u/v\in\freecat{P_1}$ of a
morphism~$u\in\freecat{P_1}$ after another morphism~$v\in\freecat{P_1}$. The
definition is performed by induction on~$u$ and~$v$ by
\begin{align}
  \label{eq:mor-res}
  u/1&=u,
  &
  u/(vv')&=(u/v)/v',
  \\
  1/u&=1,
  &
  (uu')/v&=(u/v)(u'/(v/u))
  \pbox.\nonumber
\end{align}
We will eventually see in \cref{thm:pcat-res} that, under suitable hypothesis,
this induces a residuation structure on the presented category $\pcat P$.

\begin{example}
  Consider the presentation
  \[
    \Pres{\star}{a,b}{ab=baa}
    \pbox.
  \]
  The relation can be pictured as
  \[
  \vxym{
    \ar[d]_b\ar[rr]^a&&\ar[d]^b\\
    \ar[r]_a&\ar[r]_a&
  }
  \]
  and the only possible residuation structure is defined by
\[
    a/b=aa
\qquad\text{and}\qquad
    b/a=b.
\]
  For instance, we have
  \begin{align*}
    ab/bb&=(ab/b)/b=((a/b)(b/(b/a)))/b=(aa(b/b))/b=aa1/b\\
    &=aa/b=(a/b)(a/(b/a))=aa(a/b)\\
    &=aaaa
    \pbox.
  \end{align*}
  Graphically,
  \[
    \vxym{
      \ar[d]_b\ar[rrrr]^a&&&&\ar[d]|b\ar[r]^b&\ar[d]^1\\
      \ar[d]_b\ar[rr]|a&&\ar[d]|b\ar[rr]|a&&\ar[d]|b\ar[r]|1&\ar[d]^b\\
      \ar[r]_a&\ar[r]_a&\ar[r]_a&\ar[r]_a&\ar[r]_1&
    }
  \]
  and similarly, we have $bb/ab=b$.
\end{example}

It remains to check that the above definition is sound in the sense that we can
always compute a value for the residual using the relations of~\secr{mor-res},
and that residual is uniquely defined, \ie the computed value for $u/v$ does not
depend on the way we bracket~$u$ and~$v$ or the order in which we use the
equalities~\eqref{eq:mor-res}.
Compatibility with bracketing is easily handled:

\begin{lemma}
  Residuation of morphisms is compatible with the axioms of categories.
\end{lemma}
\begin{proof}
  Residuation is compatible with associativity since
  \begin{align*}
    ((uu')u'')/v
    &=
    ((uu')/v)(u''/(v/(uu')))\\
    &=
    (u/v)(u'/(v/u))(u''/((v/u)/u'))\\
    &=
    (u/v)((u'u'')/(v/u))\\
    &=
    (u(u'u''))/v
  \end{align*}
  and
  \[
    u/((vv')v'')
    =
    (u/(vv'))/v''
    =
    ((u/v)/v')/v''
    =
    (u/v)/(v'v'')
    =
    u/(v(v'v''))
    \pbox.
  \]
  Similarly, it is compatible with left and right-unitality since
  \[
    (1u)/v
    =
    (1/v)(u/(v/1))
    =
    u/v
    =
    (u/v)(1/(v/u))
    =
    (u1)/v
  \]
  and
  \[
    u/(1v)
    =
    (u/1)/v
    =
    u/v
    =
    (u/v)/1
    =
    u/(v1)
    \pbox.
    \qedhere
  \]
\end{proof}

\noindent
However, the definition given in \secr{mor-res} is not sound in general, because
it can be vacuous, as illustrated by the following example.

\begin{example}
  Consider the presentation
  \[
    \Pres{\star}{a,b,c,d}{ba=ab,ca=ac,da=abd,cb=bac,db=bd,dc=cd},
  \]
  whose relations can be pictured as
  \[
    \begin{array}{c@{\qquad}c@{\qquad}c}
      \vxym{
        \ar[d]_a\ar[r]^b&\ar[d]^a\\
        \ar[r]_b&
      }
      &
      \vxym{
        \ar[d]_a\ar[r]^c&\ar[d]^a\\
        \ar[r]_c&
      }
      &
      \vxym{
        \ar[d]_a\ar[rr]^d&&\ar[d]^a\\
        \ar[r]_b&\ar[r]_d&
      }
      \\
      \vxym{
        \ar[d]_b\ar[r]^d&\ar[d]^b\\
        \ar[r]_d&
      }
      &
      \vxym{
        \ar[d]_c\ar[r]^d&\ar[d]^c\\
        \ar[r]_d&
      }
      &
      \vxym{
        \ar[d]_b\ar[rr]^c&&\ar[d]^b\\
        \ar[r]_a&\ar[r]_c&
      }
    \end{array}
  \]
  and consider the residuation structure defined by
  \begin{align*}
    a/b&=a,
    &
    a/c&=a,
    &
    a/d&=a,
    &
    b/a&=b,
    &
    b/c&=b,
    &
    b/d&=b,
    \\
    c/a&=c,
    &
    c/b&=ac,
    &
    c/d&=c,
    &
    d/a&=bd,
    &
    d/b&=d,
    &
    d/c&=d
    \pbox.
  \end{align*}
  The process of computing the residual $ac/bd$ by applying, from left to right,
  the relations of \secr{mor-res} defining residuation of morphisms does not
  terminate. Namely, the first two steps of this computation are
  \[
    ac/bd
    =
    aac/d
    =
    ac/bd,
  \]
  which clearly leads to a loop since the left and the right member are the
  same. This can be illustrated as follows:
  \[
    \svxym{
      \ar[d]_a\ar[r]^b&\ar[d]|a\ar[rr]^d&&\ar[d]^a\\
      \ar[dd]_c\ar[r]|b&\ar[d]^a\ar[r]_b&\ar[r]_d&\\
      &\ar[d]^c\\
      \ar[r]_b&
    }
  \]
\end{example}

\subsection{Termination of residuation}
Let $P$~be a fixed residuated presentation. In order to ensure that the
process of computing the residual is well-defined, we follow the technique of
considering ``reversed words'' introduced by
Dehornoy~\cite{dehornoy2000completeness}, and consider the following
polygraph~$Q$ defined from~$P$ by
\begin{align*}
  Q_0&= P_0,
  \\
  Q_1&=\setof{a:x\to y,\finv a:y\to x}{a:x\to y\in P_1},
  \\
  Q_2&=\setof{\finv ab\To v\finv u}{a,b\in P_1,u=a/b,v=b/a},
\end{align*}
where $\finv{(a_1\ldots a_n)}$ is a notation for $\finv a_n\ldots\finv a_1$, and
$\finv a$ should be thought of as a formal inverse for the generator~$a$. A
morphism in~$\freecat{Q_1}$ is a composite of generators in~$P_1$, some of which
might be formally inverted, and rewriting step corresponds to taking residuals
in~$P$:
\[
  \xymatrix@C=3ex@R=3ex{
    &\ar[dl]_ax\ar[dr]^b&\\
    y\ar@{.>}[dr]_{b/a=v}&\Downarrow&\ar@{.>}[dl]^{u=a/b}z\pbox.\\
    &w&
  }
\]
We say that a residuated presentation~$P$ is \emph{terminating} when the
associated $2$-polygraph~$Q$ is terminating in the usual sense.

\begin{lemma}
  \label{lem:res-confl}
  Given a terminating residuated presentation~$P$,
  the associated 2-polygraph~$Q$ is convergent and the residuation operation is
  well-defined on morphisms of~$\freecat{P_1}$. Moreover, for morphisms
  $u:x\to y$ and $v:x\to y'$, the morphisms $u(v/u)$ and $v(u/v)$ are
  $P$-congruent:
  \[
    \xymatrix@C=3ex@R=3ex{
      &\ar[dl]_ux\ar[dr]^v&\\
      y\ar[dr]_{v/u}&\overset*\Leftrightarrow&\ar[dl]^{u/v}y'\pbox.\\
      &z&
    }
  \]
\end{lemma}
\begin{proof}
  The polygraph~$Q$ has no critical pair, it is thus locally confluent by
  \cref{lem:2-cb} and confluent by \cref{lem:newman} since it is assumed to be
  terminating. By well-founded induction, we can show that the normal form of a
  word $\finv uv$ is a word of the form $v'\finv{u'}$ with $v'=v/u$ and
  $u'=u/v$, and that any word of this form is a normal form. The last part of
  the lemma follows by induction from the assumption~\eqref{eq:res-fill-ax}.
\end{proof}

\noindent
In practice, various practical conditions are sufficient to ensure the
termination of the polygraph~$Q$, see~\cite{dehornoy2015foundations,
  clerc2015presenting}. For instance,

\begin{lemma}
  \label{lem:res-1-term}
  Suppose given a function $\omega:P_1\to\N$, which we extend as a function
  $\omega:\freecat{P_1}\to\N$ by $\omega(\unit{})=0$ and
  $\omega(uv)=\omega(u)+\omega(v)$. Suppose moreover that we have
  $\omega(a/b)<\omega(a)$ for every pair of generators $a,b\in P_1$. Then the
  2-polygraph~$Q$ is terminating.
\end{lemma}
\begin{proof}
  We define a function $\omega':\freecat{Q_1}\to\N$ by $\omega'(a)=\omega(a)$
  and $\omega(\finv{a})=0$ for $a\in P_1$, $\omega'(uv)=\omega'(u)+\omega'(v)$,
  $\omega(\unit{})=0$. This function is a reduction order on the 2-polygraph~$Q$
  and we conclude by \lemr{red-fun}.
\end{proof}

\noindent
Finally, we can check that residuation is well defined on morphisms of the
presented category~$\pcat{P}$, \ie that it is compatible with $P$-congruence on
morphisms.

\begin{lemma}
  \label{lem:res-2-compat}
  Given a terminating residuated presentation~$P$, for every morphisms
  $u,v:x\to y$ and $w:x\to x'$ in~$\freecat{P_1}$, we have that
  \[
    u\overset*\Leftrightarrow v
    \qqqqtimpl
    u/w\overset*\Leftrightarrow v/w
    \qtand
    w/u=w/v
    \pbox.
  \]
  Graphically,
  \[
    \xymatrix@C=8ex@R=4ex{
      x\ar[d]_w\ar@/^/[r]^u\ar@/_/[r]_v\ar@{}[r]|{\Updownarrow*}&y\ar@{.>}[d]^{w/u=w/v}\\
      x'\ar@/^/@{.>}[r]^{u/w}\ar@/_/@{.>}[r]_{v/w}\ar@{}[r]|{\Updownarrow*}&y'\pbox.
    }
  \]
\end{lemma}
\begin{proof}
  The assumption that $u$ and $v$ are $P$-congruent means that there exists a
  sequence of $2$-cells of the form
  \[
    \xymatrix@C=5ex@R=3ex{
      x\ar[r]^{u_i'}&x_i\ar@/^/[r]^{u_i}\ar@/_/[r]_{v_i}\ar@{}[r]|{\Updownarrow\alpha_i}&y_i\ar[r]^{u_i''}&y
    }
  \]
  with $1\leq i\leq n$, $u_i,u_i',u_i'',v_i\in\freecat{P_1}$ and
  $\alpha_i\in P_2$, with either $\alpha_i:u_i\To v_i$ or $\alpha_i:v_i\To u_i$,
  such that $u_1'u_1u_1''=u$, $u_{i+1}'u_{i+1}u_{i+1}''=u_i'v_iu_i''$ and
  $u_n'v_nu_n''=v$. We have
  \[
    w/(u_i'u_iu_i'')
    =
    ((w/u_i')/u_i)/u_i''
    =
    ((w/u_i')/v_i)/u_i''
    =
    w/(u_i'v_iu_i'')
  \]
  where the equality $(w/u_i')/u_i=(w/u_i')/v_i$ can be shown by recurrence on
  the length of $w/u_i'$ using axioms~\eqref{eq:res-cyl-ax}. Also, by recurrence
  on $w/u_i'$ and using axioms~\eqref{eq:res-cyl-ax}, we have the existence of a
  2-generator between $u_i/(w/u_i')$ and $v_i/(w/u_i')$:
  \[
    \xymatrix@C=17ex@R=10ex{
      x\ar[d]_w\ar[r]^{u_i'}&x_i\ar@{.>}[d]|{w/u_i'}\ar@/^/[r]^{u_i}\ar@/_/[r]_{v_i}\ar@{}[r]|{\Updownarrow\alpha_i}&y_i\ar@{.>}[d]|{\begin{array}{c}\scriptstyle(w/u_i')/u_i\\[-1ex]\scriptstyle=\\[-1ex]\scriptstyle(w/u_i')/v_i\end{array}}\ar[r]^{u_i''}&y\ar@{.>}[d]|{\begin{array}{c}\scriptstyle((w/u_i')/u_i)/u_i''\\[-1ex]\scriptstyle=\\[-1ex]\scriptstyle((w/u_i')/v_i)/u_i''\end{array}}\\
      x'\ar@{.>}[r]_{u_i'/w}&x_i'\ar@/^/@{.>}[r]^{u_i/(w/u_i')}\ar@/_/@{.>}[r]_{v_i/(w/u_i')}\ar@{}[r]|{\Updownarrow}&y_i'\ar@{.>}[r]_{u_i''/((w/u_i')/u_i)}&y'\pbox.
    }
  \]
  We conclude, by performing a recurrence on~$n$.
\end{proof}

\begin{theorem}
  \label{thm:pcat-res}
  Given a terminating residuated presentation~$P$, the presented
  category~$\pcat{P}$ admits a residuation structure.
\end{theorem}
\begin{proof}
  The residuation operation is well-defined on morphisms in~$\freecat{P_1}$
  modulo $P$-congruence by previous lemmas and immediately satisfies the axioms
  of a residuation structure.
\end{proof}

The axiomatization presented in this section has the advantage of being
relatively simple to state and prove, but more advanced generalizations are
often required in practice.
For instance, in many situations, not every pair of coinitial morphisms $f$ and
$g$ admit a residual, but only those which are \emph{bounded}, \ie for which
there exists $f'$ and $g'$ with $ff'\overset*\Leftrightarrow gg'$.
Also, it is useful to weaken axiom~\eqref{eq:res-cyl-ax} and require that for
every 1-generator $a:x\to x'$ and 2-generator $\alpha:u\To v:x\to y$, we have a
2-cell
\[
  a/\alpha
  :
  a/u
  \To
  a/v
  :
  y\to y'
\]
and a 2-cell
\[
  \alpha/a
  :
  u/a\To v/a
  :
  x'\to y'
\]
in~$\freecat{P_2}$ (or even in $\freegpd{P_2}$):
\begin{equation}
  \label{eq:res-cyl-ax-gen}
  \xymatrix@C=12ex@R=8ex{
    x\ar[d]_a\ar@/^/[r]^u\ar@/_/[r]_v\ar@{}[r]|{\Downarrow\alpha}&y\ar@{.>}@/_/[d]_{a/u}\ar@{}[d]|{\overset{a/\alpha}\To}\ar@{.>}@/^/[d]^{a/v}\\
    x'\ar@/^/@{.>}[r]^{u/a}\ar@/_/@{.>}[r]_{v/a}\ar@{}[r]|{\Downarrow\alpha/a}&y'
  }
\end{equation}
axiom \eqref{eq:res-cyl-ax} being the particular case where we further impose
that $a/\alpha$ is an identity and $\alpha/a$ is a whiskered 2-generator. 
\index{noetherian!category}
In this case, in order for \lemr{res-2-compat} to hold, one has to impose
further termination conditions. Given two cofinal morphisms $f$ and $g$, we
write $g|f$ whenever there exists $h$ with $hg=f$, and in this case we say that
$g$ \emph{divides} $f$ on the right. A category is \emph{right noetherian} when
every infinite sequence $(f_i)$ of cofinal morphisms $f_{i+1}|f_i$ is eventually
stationary. In particular, a category presented by a 2-polygraph whose relations
are homogeneous (\ie preserve the length of words) necessarily has this
property. The following theorem is due to Dehornoy:
see~\cite[Section~II.4]{dehornoy2015foundations} for detailed statement and
proof.

\begin{theorem}
  Given a residuated presentation with generalized
  axiom~\eqref{eq:res-cyl-ax-gen}, whose presented category is right Noetherian,
  every bounded pair of morphisms $u$ and $v$ in~$\freecat{P_1}$ admits a
  residual.
\end{theorem}

\begin{example}
  \index{braid!monoid}
  Consider the positive braid monoid $B_4^+$, see \secr{braid-mon},
  which admits a presentation by a $2$-polygraph with three generators
  $a_0, a_1, a_2$ and three relations
  \begin{align*}
    \alpha_{01}&:a_0a_1a_0\To a_1a_0a_1,
    &
    \alpha_{12}&:a_1a_2a_1\To a_2a_1a_2,
    &
    \alpha_{02}&:a_0a_2\To a_2a_0,
  \end{align*}
  which can respectively be pictured as
  \[
    \xymatrix@C=3ex@R=3ex{
      \ar[dd]_{a_0}\ar[rr]^{a_1}\ar@{}[ddrr]|{\alpha_{01}\Tour}&&\ar[d]^{a_0}\\
      &&\ar[d]^{a_1}\\
      \ar[r]_{a_1}&\ar[r]_{a_0}&
    }
    \qquad\qquad
    \xymatrix@C=3ex@R=3ex{
      \ar[dd]_{a_1}\ar[rr]^{a_2}\ar@{}[ddrr]|{\alpha_{12}\Tour}&&\ar[d]^{a_1}\\
      &&\ar[d]^{a_2}\\
      \ar[r]_{a_2}&\ar[r]_{a_1}&
    }
    \qquad\qquad
    \xymatrix@C=3ex@R=3ex{
      \ar[dd]_{a_0}\ar[rr]^{a_2}\ar@{}[ddrr]|{\alpha_{02}\Tour}&&\ar[dd]^{a_0}\\
      &&\\
      \ar[rr]_{a_2}&&
    }
  \]
  We define residuation on generators by
  \begin{align*}
    a_0/a_0&=\unit{},
    &
    a_1/a_0&=a_1a_0,
    &
    a_2/a_0&=a_2,
    \\
    a_0/a_1&=a_0a_1,
    &
    a_1/a_1&=\unit{},
    &
    a_2/a_1&=a_2a_1,
    \\
    a_0/a_2&=a_0,
    &
    a_1/a_2&=a_1a_2,
    &
    a_2/a_2&=\unit{}
    \pbox.
  \end{align*}
  We can check axiom~\eqref{eq:res-cyl-ax}, \ie that residuation of 1-generators
  is compatible with 2-generators. For instance, for the relation~$\alpha_{01}$,
  the residuals of~$a_0$ ($a_1$ is similar) and $a_2$ after the source and the
  target are respectively
  \[
    \xymatrix@C=3ex@R=3ex{
      \\
      \\
      \ar[rr]^{\unit{}}&&\ar[r]^{a_1}&\ar[r]^{a_0}&\\
      \ar[u]^{a_0}\ar[rr]_{a_0}&&\ar[u]|{\unit{}}\ar[r]_{a_1}\ar@{}[dd]|{\alpha_{01}\Downarrow}&\ar[r]_{a_0}\ar[u]|{\unit{}}&\ar[u]|{\unit{}}\\
      \\
      \ar[dd]_{a_0}\ar[rr]^{a_1}&&\ar[d]|{a_0}\ar[r]^{a_0}&\ar[d]|{\unit{}}\ar[r]^{a_1}&\ar[d]^{\unit{}}\\
      &&\ar[d]|{a_1}\ar[r]|{\unit{}}&\ar[d]|{a_1}\ar[r]|{a_1}&\ar[d]^{\unit{}}\\
      \ar[r]_{a_1}&\ar[r]_{a_0}&\ar[r]_{\unit{}}&\ar[r]_{\unit{}}&\\
      \\
    }
    \qquad\qquad
    \xymatrix@C=3ex@R=3ex{
      \ar[rr]^{a_0}&&\ar[r]^{a_1}&\ar[r]^{a_2}&\ar[r]^{a_0}&\ar[r]^{a_1}&\\
      &&&&&&\ar[u]_{a_0}\\
      &&&&\ar[uu]|{a_1}\ar[rr]|{a_0}&&\ar[u]_{a_1}\\
      \ar[uuu]^{a_2}\ar[rr]_{a_0}&&\ar[uuu]|{a_2}\ar[rr]_{a_1}&\ar@{}[dd]|{\alpha_{01}\Downarrow}&\ar[u]|{a_2}\ar[rr]_{a_0}&&\ar[u]_{a_2}\\
      \\
      \ar[dddd]_{a_2}\ar[rr]^{a_1}&&\ar[dd]|{a_2}\ar[rr]^{a_0}&&\ar[dd]|{a_2}\ar[rr]^{a_1}&&\ar[d]^{a_2}\\
      &&&&&&\ar[d]^{a_1}\\
      &&\ar[dd]|{a_1}\ar[rr]|{a_0}&&\ar[d]|{a_1}\ar[r]|{a_1}&\ar[d]|{\unit{}}\ar[r]|{a_2}&\ar[d]^{\unit{}}\\
      &&&&\ar[d]|{a_0}\ar[r]|{\unit{}}&\ar[d]|{a_0}\ar[r]|{a_2}&\ar[d]^{a_0}\\
      \ar[r]_{a_1}&\ar[r]_{a_2}&\ar[r]_{a_0}&\ar[r]_{a_1}&\ar[r]_{\unit{}}&\ar[r]_{a_2}&\\
    }
  \]
  and we can thus take $a_0/\alpha_{01}=\unit{\unit{}}$,
  $a_2/\alpha_{01}=\unit{a_2a_1a_0}$, $\alpha_{01}/a_0=\unit{a_1a_0}$ and
  $\alpha_{01}/a_2$ to be the 2-cell
  \[
    \xymatrix@C=3ex@R=3ex{
      &\ar[r]^{a_1}&\ar[dr]_{a_0}\ar[r]^{a_2}&\ar[dr]^{a_0}\\
      \ar[ur]^{a_0}\ar[dr]_{a_1}&&&\ar[r]_{a_2}&\ar[dr]^{a_1}\\
      &\ar[dr]_{a_2}\ar[r]^{a_0}&\ar[dr]^{a_2}\ar[ur]|{a_1}&&&\pbox.\\
      &&\ar[r]_{a_0}&\ar[r]_{a_1}&\ar[ur]_{a_2}&
    }
  \]
  Other cases are left to the reader. Being homogeneous, this presentation is
  right Noetherian and residuation is always
  terminating~\cite{dehornoy1997groups}. The category~$B_4^+$ is thus
  residuated. The argument generalize to all positive braid monoids $B_n^+$.
\end{example}

\subsection{Deciding equality}
\index{decidability!of equality}
As a last remark, note that for residuated presentations~$P$ the word problem can be
solved in the following way. Given two morphisms $u,v:x\to y$ in
$\freecat{P_1}$, we have $u\overset*\Leftrightarrow v$ if and only if
\[
  u/v\overset*\Leftrightarrow\unit{y}
  \qqtand
  v/u\overset*\Leftrightarrow\unit{y}
  \pbox.
\]
This follows easily from the fact that residuals corresponds to pushouts cocones
by \propr{residual-pushout}. In particular, when $P$ has no 2-generator with an
identity as source or as target, we have $u\overset*\Leftrightarrow v$ if and
only if
\[
  u/v=\unit{y}
  \qqtand
  v/u=\unit{y}
  \pbox.
\]


\chapter{Tietze Transformations and Completion}
\label{chap:2-tietze}
\label{chap:2tietze}
In this chapter, we introduce a notion of \emph{Tietze transformation} for
2-polygraphs, generalizing the one introduced in \cref{sec:pres-set} for
1-polygraphs. The Tietze transformations are elementary operations on
2-polygraphs, which preserve the presented category, and such that any two
finite $2$-polygraphs presenting the same category can be transformed  into
one  another by applying a series of such transformations. Our notion, introduced
in \cref{sec:2-tietze}, is very close to the one first introduced by Tietze for
presentations of groups~\cite{tietze1908topologischen}. 
We refer to~\cite{lyndon2015combinatorial,magnus04} for more details on the notion of Tietze transformation in combinatorial group theory, see also \cite{magnus82} for a historical account. The notion of Tietze transformation was developed in the polygraphic language in~\cite{GaussentGuiraudMalbos15}.

By using Tietze transformations, one seeks to turn a given
presentation of a category into another one,
possessing better computational properties. In particular, the
\emph{Knuth-Bendix completion} procedure described in~\cref{sec:2kb}
applies those transformations to turn a presentation into a confluent one.

We have seen in~\cref{sec:1-eq-dec}  how convergent presentations
lead to a solution of the word problem: for those, the
equivalence between two words is immediately decided by comparing
their normal forms. In order to tackle the word problem for an
arbitrary presentation, a good strategy thus consists in trying to
transform it into a convergent one by using Tietze transformations.
From this point of view, we  naturally ask ourselves whether a finite presentation of a category
with decidable word problem can always be turned  into a convergent
one by applying Tietze transformations. This problem, called
\emph{universality of convergent presentations}, is introduced in
\cref{sec:universality}. We will see in \cref{chap:2-fdt,chap:2-homology} that
the answer to this question is negative.

\section{Tietze Transformations}
\label{sec:2-tietze}

\subsection{Definition}
\label{sec:2-tietze-def}
\index{Tietze!transformation!of 2-polygraphs}
\index{transformation!Tietze}
The \emph{elementary Tietze transformations} are the following transformations of a
2-polygraph~$P$ into a 2-polygraph~$Q$:
\begin{description}
\item[\tgen] \emph{adding a definable 1-generator}: given $a\not\in P_1$,
  $u:x\to y\in\freecat{P_1}$, and $\alpha\not\in P_2$ we define
  \[
    Q=\Pres{P_0}{P_1,a:x\to y}{P_2,\alpha:a\To u}
    \pbox,
  \]
\item[\trel] \emph{adding a derivable relation}: given $u,v\in\freecat{P_1}$ such
  that $u\approx v$ and $\alpha\not\in P_2$, we define
  \[
    Q=\Pres{P_0}{P_1}{P_2,\alpha:u\To v}
    \pbox.
  \]
\end{description}
\index{Tietze!equivalence!of 2-polygraphs}
\index{equivalence!Tietze}
The \emph{Tietze equivalence} is the smallest equivalence relation on
2-polygraphs which is stable under isomorphisms and Tietze
transformations. We respectively write \trgen{} and \trrel{} for operations
\tgen{} and \trel{} performed backward.

These local transformations completely axiomatize the property of presenting the
same categories. This was first shown by Tietze~\cite{tietze1908topologischen}
for presentations of groups, and the proof extends to the case of 2-polygraphs. 

\begin{theorem}
  \label{thm:2-tietze-equiv}
  Two finite 2-polygraphs~$P$ and~$Q$ present the same category, \ie
  $\prescat{P}\isoto\prescat{Q}$, if and only if they are Tietze equivalent.
\end{theorem}
\begin{proof}
Let $P$ be a 2-polygraph. If the 2-polygraph $Q$ is either isomorphic
to $P$ or obtained by performing transformations  \tgen{} or \trel{} on
  $P$, then $\prescat{Q}$ is isomorphic to $\prescat{P}$. Therefore, any
  two Tietze equivalent 2-polygraphs present the same category.

  Conversely,
  suppose that~$P$ and~$Q$ present the same category~$C$. Up to
  isomorphism, that is, renaming of generators, we may
  suppose that $P_0=Q_0$, $P_1\cap Q_1=\emptyset$ and
  $P_2\cap Q_2=\emptyset$.
  We write
  $q^P:\freecat{P_1}\to(\freecat{P_1}/{\approx_P})=C$ for the quotient functor
  (see \secr{pres-cat}): this functor is full and such that $u\approx^P v$
  precisely when $q^P(u)=q^P(v)$. Similarly, we also consider the quotient
  functor $q^Q:\freecat{Q_1}\to C$.
  Starting from the presentation~$P$, we apply the following series of Tietze
  equivalences.
  \begin{enumerate}
  \item Given a 1-generator~$a\in Q_1$, its image~$q^Q(a)$ is a morphism of~$C$
    and therefore has a representative in~$\freecat P_1$: since $q^P$ is full,
    there exists $u_a\in\freecat P_1$ satisfying $q^P(u_a)=q^Q(a)$. By a
    transformation \tgen{}, we add to~$P$ the 1\nbd-gene\-rator~$a$ and the
    relation $u_a\To a$. Performing this for every generator $a\in Q_1$, we
    obtain the 2-polygraph
    \[
    P'
    =
    \Pres{C_0}{P_1\cup Q_1}{P_2\cup\setof{u_a\To a}{a\in Q_1}}
    \pbox.
    \]
  \item Given a 1-cell $u=a_1\ldots a_n\in\freecat{Q_1}$, by construction
    of~$P'$, we have that $q^{P'}(u_a)=q^Q(a)$, which implies
    \begin{align*}
      q^{P'}(u)&=q^{P'}(a_1)\ldots q^{P'}(a_n)\\
      &=q^{P'}(u_{a_1})\ldots q^{P'}(u_{a_n})\\
      &=q^Q(a_1)\ldots q^Q(a_n)\\
      &=q^{Q}(u)
      \pbox.
    \end{align*}
    For each relation $\alpha:u\To v$ in~$Q_2$, we have $q^Q(u)=q^Q(v)$, which
    implies $q^{P'}(u)=q^{P'}(v)$ by the above, and therefore
    $u\approx^{P'}v$. By a transformation \trel{}, we can thus add to the
    previous 2-polygraph the derivable relation~$\alpha:u\To v$. Performing this
    for every relation $\alpha\in Q_2$, we obtain the 2-polygraph
    \[
    P''
    =
    \Pres{C_0}{P_1\cup Q_1}{P_2\cup Q_2\cup\setof{u_a\To a}{a\in Q_1}}
    \pbox.
    \]
  \item Suppose given a 1-generator~$a\in P_1$. For similar reasons as in first
    step, there exists $v_a\in\freecat{Q_1}$ such that $q^Q(v_a)=q^P(a)$, which
    implies $q^{P''}(v_a)=q^{P''}(a)$, \ie $v_a\approx^{P''} a$. By a
    transformation \trel{}, we can therefore add the derivable relation
    $v_a\To a$. Performing this for every generator $a\in P_1$, we obtain the
    2-polygraph $P'''$ which is
    \[
    \Pres{C_0}{P_1\cup Q_1}{P_2\cup Q_2\cup\setof{v_a\To a}{a\in P_1}\cup\setof{u_a\To a}{a\in Q_1}}
    \pbox.
    \]
  \end{enumerate}
  By exchanging the roles of~$P$ and~$Q$, one shows that~$Q$ is also Tietze
  equivalent to the same polygraph~$P'''$. Therefore, the 2-polygraphs~$P$
  and~$Q$ are Tietze equivalent.
\end{proof}

\begin{remark}
  Similarly to the case of 1-polygraphs (\cref{rem:trans-tietze}), Tietze
  transformations can be extended to account for infinite $2$-polygraphs. The notion
  of Tietze transformation has to be refined in the following way: we say that a
  polygraph~$P$ \emph{Tietze expands} to a polygraph~$Q$ when there is a
  transfinite sequence of Tietze transformations from~$P$ to~$Q$, and we define
  \emph{Tietze equivalence} as the smallest equivalence relation containing
  Tietze expansion. The above proof can be adapted in order to show that two
  polygraphs (of arbitrary cardinality) present the isomorphic categories if and
  only if they are Tietze equivalent.
\end{remark}

\begin{example}
  \label{ex:S3-two-presentations}
  Consider the symmetric group~$S_3$ on 3 elements. We consider it here as a
  monoid (in which elements happen to have inverses). It can be described as the
  category with only one object, whose morphisms are bijections
  $f:\intset{3}\to\intset{3}$, where $\intset{3}$ denotes the set $\set{0,1,2}$
  with three elements, equipped with usual composition and identity. This group
  is generated by the two transpositions $s$ and $t$ whose graphs are
  respectively
  \[
  \svxym{
    \ar@{}[d]|{0}&\ar@{}[d]|{1}&\ar@{}[d]|{2}\\
    \cdot\ar@{-}[dr]&\cdot\ar@{-}[dl]&\cdot\ar@{-}[d]\\
    \cdot&\cdot&\cdot\\
    \ar@{}[u]|{0}&\ar@{}[u]|{1}&\ar@{}[u]|{2}
  }
  \qquad\qquad\qquad\qquad
  \svxym{
    \ar@{}[d]|{0}&\ar@{}[d]|{1}&\ar@{}[d]|{2}\\
    \cdot\ar@{-}[d]&\cdot\ar@{-}[dr]&\cdot\ar@{-}[dl]\\
    \cdot&\cdot&\cdot\\
    \ar@{}[u]|{0}&\ar@{}[u]|{1}&\ar@{}[u]|{2}
  }
  \]
  and working out the relations which are satisfied by those generators, one can
  come up with the following presentation of~$S_3$:
  \[
    P=
    \Pres{\star}{s,t}{ss=1,tt=1,sts=tst}
  \]
  see \cref{sec:pres-sym-group,sec:sym-group-pres} for details.

  The group $S_3$ can also be considered as the group of symmetries of an
  equilateral triangle
  \[
    \fig{S3-sr}
  \]
  Namely, any bijection between the set of vertices determines a unique
  symmetry. As such, it can be generated by a symmetry~$s$ about a vertical axis
  and a rotation~$r$ of angle $2\pi/3$: those respectively correspond to the
  bijections between vertices whose graphs are
  \[
  \svxym{
    \ar@{}[d]|{0}&\ar@{}[d]|{1}&\ar@{}[d]|{2}\\
    \cdot\ar@{-}[dr]&\cdot\ar@{-}[dl]&\cdot\ar@{-}[d]\\
    \cdot&\cdot&\cdot\\
    \ar@{}[u]|{0}&\ar@{}[u]|{1}&\ar@{}[u]|{2}
  }
  \qquad\qquad\qquad\qquad
  \svxym{
    \ar@{}[d]|{0}&\ar@{}[d]|{1}&\ar@{}[d]|{2}\\
    \cdot\ar@{-}[dr]&\cdot\ar@{-}[dr]&\cdot\ar@{-}[dll]\\
    \cdot&\cdot&\cdot\\
    \ar@{}[u]|{0}&\ar@{}[u]|{1}&\ar@{}[u]|{2}
  }
  \]
  Note that the interpretation of~$s$ is the same as previously, and that~$r$
  can be expressed in terms of the previous generators as $r=ts$. Working out the
  relations satisfied by those generators, one obtains the following
  presentation:
  \[
    Q=
    \Pres{\star}{s,r}{rrr=1,ss=1,rsrs=1},
  \]
  which is the usual presentation of the dihedral group~$D_3$, see
  \cref{sec:dihedral-group}.

  Since the two above 2-polygraphs~$P$ and~$Q$ present the same group,
  \cref{thm:2-tietze-equiv} asserts that they are Tietze equivalent. For
  instance, the presentation~$P$ can be transformed into~$Q$ by the following
  series of Tietze transformations. Starting from the 2-polygraph~$P$,
  \begin{description}[labelwidth=\widthof{\trrel}]
  \item[\tgen] add the definable generator $r=ts$:
    \[
    \Pres{\star}{r,s,t}{ss=1,tt=1,sts=tst,r=ts},
    \]
  \item[\trel] add the relation $rrr=1$ (derivable since $rrr=tststs=ttstts=ss=1$):
    \[
    \Pres{\star}{r,s,t}{rrr=1,ss=1,tt=1,sts=tst,r=ts},
    \]
  \item[\trel] add the relation $rsrs=1$ (derivable since $rsrs=tsstss=tt=1$):
    \[
    \Pres{\star}{r,s,t}{rrr=1,ss=1,tt=1,rsrs=1,sts=tst,r=ts},
    \]
  \item[\trel] add the relation $t=rs$ (derivable since $t=tss=rs$):
    \[
    \Pres{\star}{r,s,t}{rrr=1,ss=1,tt=1,rsrs=1,sts=tst,r=ts,t=rs},
    \]
  \item[\trrel] remove the relation $r=ts$ (derivable since $r=rss=ts$):
    \[
    \Pres{\star}{r,s,t}{rrr=1,ss=1,tt=1,rsrs=1,sts=tst,t=rs},
    \]
  \item[\trrel] remove the relation $tt=1$ (derivable since $tt=rsrs=1$):
    \[
    \Pres{\star}{r,s,t}{rrr=1,ss=1,rsrs=1,sts=tst,t=rs},
    \]
  \item[\trrel] remove the relation $sts=tst$, which is derivable since
    \[
    sts=srss=sr=rrrsrss=rrs=rssrs=tst
    \]
    to obtain
    \[
    \Pres{\star}{r,s,t}{rrr=1,ss=1,rsrs=1,t=rs}.
    \]
  \item[\trgen] finally, remove the definable generator~$t$ (which does not
    occur in any relation other than $t=rs$) to obtain~$Q$.
  \end{description}
\end{example}

\begin{example}
  The monoid $(\N/3\N)\times(\N/2\N)$ admits the presentation
  \[
  \Pres{\star}{s,t}{s^3=1,t^2=1,ts=st}
  \]
  but it also admits the presentation
  \[
  \Pres{\star}{r}{r^6=1}
  \]
  (hint: define $r$ by $r=ts$). This shows that we have an isomorphism
  \[
    (\N/3\N)\times(\N/2\N)
    \isoto
    \N/6\N
  \]
  as already noted in \cref{ex:dlaw-N26}. More generally, one can show that the
  monoids
  \[
    (\N/p\N)\times(\N/q\N)
    \isoto
    \N/pq\N
  \]
  are isomorphic when $p$ and $q$ are relatively prime natural numbers.
\end{example}

\subsection{Reduced 2-polygraphs}
\label{sec:2-reduced}
\label{DefinitionReducedPolygraph}
\index{reduced!2-polygraph}
\index{polygraph!reduced}
Tietze equivalences allow one to simplify presentations without changing the
presented category. In particular, one can, without loss of generality, restrict
to the following class of 2-polygraphs, which are often easier to handle than
general 2-polygraphs. 
Those were studied by
Metivier~\cite{metivier1983rewriting} for term rewriting systems, and
Squier~\cite[Theorem~2.4]{squier1987word} for string rewriting systems.

A 2-polygraph~$P$ is
\begin{itemize}
\item \emph{left reduced} when, for every rule $\alpha:u\To v$ in~$P_2$, $u$ is
  not reducible by any rule other than~$\alpha$,
\item \emph{right reduced} when, for every rule $\alpha:u\To v$ in~$P_2$, $v$ is
  not reducible by any rule,
\item \emph{reduced} when it is both left and right reduced.
\end{itemize}
Note that a left reduced 2-polygraph never has inclusion critical branchings, as
defined in \secr{cb-class}, which often simplifies the study of branchings.

\begin{theorem}[{\cite[Theorem~2.4]{squier1987word}}]
  \label{thm:reduced-pres}
  Every convergent 2-polygraph~$P$ is Tietze equivalent to a reduced convergent
  2-polygraph.
\end{theorem}
\begin{proof}
  Starting from the 2-polygraph~$P$, we successively apply the following Tietze
  transformations.
  \begin{enumerate}
  \item Replace every $2$-cell $\alpha:u\To v$ by $\alpha:u\To\nf{u}$, where
    $\nf{u}$ is the normal form of~$u$:
    \[
    \svxym{
      {u}
      \ar@2[rr]^-{\alpha} 
      && {v}
      \ar@2[d]^\ast
      \\
      && {\nf{u}}
    }
    \qquad\overset\trel\rightsquigarrow\qquad
    \svxym{
      {u}
      \ar@2 [rr]
      \ar@2[drr]
      && {v}
      \ar@2 [d]^\ast
      \\
      && {\nf{u}}
    }
    \qquad\overset\trrel\rightsquigarrow\qquad
    \svxym{
      {u} 
      \ar@2 [drr]_-{\alpha} 
      && {v}
      \ar@2 [d]^\ast
      \\
      && {\nf{u}}\pbox.
    }
    \]
  \item If the resulting $2$-polygraph contains parallel $2$-cells\index{parallel!2-cells@$2$-cells}, remove all
    but one:
    \[
    \svxym{
      {u}
      \ar@2 @/^3ex/ [rr] ^-{\alpha_1} ^{}="src"
      \ar@2 @/_3ex/ [rr] _-{\alpha_n} _{}="tgt"
      && {\nf{u}}
      \ar@{.} "src"!<0pt,-10pt>;"tgt"!<0pt,10pt>
    }
    \qquad\qquad\overset\trrel\rightsquigarrow\qquad\qquad
    \svxym{
      {u}
      \ar@2 [rr] ^-{\alpha}
      && {\nf{u}}\pbox.
    }
    \]
  \item Finally, remove, in the resulting $2$-polygraph, every $2$-cell whose
    source is reducible by another $2$-cell:
    \[
    \svxym{
      {uvw}
      \ar@2 [rr] ^-{\alpha}
      \ar@2 [drr] _-{u\beta w} 
      && {\nf{uvw}}
      \\
      && {u\nf{v}w}
      \ar@2 [u] 
    }
    \qquad\qquad\overset\trrel\rightsquigarrow\qquad\qquad
    \svxym{
      {uvw}
      \ar@2 [drr] _-{u\beta w}
      && {\nf{uvw}}
      \\
      && \ar@2 [u]
      {u\nf{v}w}\pbox.
    }
    \]
  \end{enumerate}
  These steps all correspond to Tietze transformations of type \trel{} and the
  resulting polygraph is clearly reduced.
\end{proof}

\begin{example}
  Consider the following presentation of the symmetric group~$S_3$:
  \[
    \Pres{\star}{r,s,t}{
      \begin{array}{r@{\ :\ }r@{\ \To\ }l@{\ }r@{\ :\ }r@{\ \To\ }l@{\ }l}
        \sigma&ss&1&\gamma&sts&tst,&\rho:ts\To r\\
        \tau&tt&1,&\gamma'&sts&sstst
      \end{array}  
    }\pbox.
  \]
  It is not reduced because the target of~$\gamma'$ is not reduced (its normal
  form is~$tst$) and the sources of~$\gamma$ and $\gamma'$ are reducible
  by~$\rho$. Applying the procedure described in the proof of
  \thmr{reduced-pres}, we obtain the following reduced, Tietze equivalent,
  presentation:
  \[
  \Pres{\star}{r,s,t}{\sigma:ss\To 1,\tau:tt\To 1,\gamma:sr\To rt,\rho:ts\To r}\pbox.
  \]
\end{example}

\subsection{The reduced standard presentation}
\label{Example:ReducedStandardPresentation}
\index{standard!presentation!reduced}
\index{presentation!standard!reduced}
In \secr{cat-std-pres}, we have seen that every category~$C$ admits a canonical
presentation, the standard presentation. One can actually achieve a smaller
presentation by not adding identities as 1-generators. The \emph{reduced
  standard polygraphic presentation} of~$C$ is the 2\nbd-polygraph~$R$ where
\begin{itemize}
\item $R_0$ is the set of objects of~$C$,
\item $R_1$ is the set of morphisms of~$C$ which are not identities,
\item $R_2$ contains 2-cells of the form
  \[
    \mu_{a,b}
    :
    ab\To(b\circ a)
    :
    x\to z,
  \]
  for every object~$x$ of~$C$ and pair of composable morphisms $a:x\to y$ and
  $b:y\to z$ in~$C$ such that $b\circ a$ is not an identity, and 2-cells of the
  form
  \[
    \mu'_{a,b}
    :
    ab\To x
    :
    x\to x,
  \]
  for every object~$x$ of~$C$ and pair of composable morphisms $a:x\to y$ and
  $b:y\to x$ in~$C$ such that $b\circ a=\unit{x}$.
\end{itemize}

The proof that this is indeed a presentation of~$C$ can be performed by adapting
the rewriting argument provided in \secr{cat-std-pres}. Another way to show
this, since we know that the standard presentation~$P$ of~$C$ is a presentation
of~$C$, is to show that~$R$ is Tietze equivalent to~$C$. Starting from~$P$, this
can be done by using the following series of Tietze transformations.
\begin{itemize}
\item For each 2-cell $\mu_{a,b}:ab\To(b\circ a):x\to x$ such that
  $b\circ a=\id_x$ is an identity one can add the derivable 2-cell
  $\mu'_{a,b}:ab\To x$ and remove the 2-cell $\mu_{a,b}$ by using
  transformations \trel{}:
  \[
  \svxym{
    &y\ar@/^/[dr]^b&\\
    x\ar@/^/[ur]^a\ar@/_3ex/[rr]_{\unit{x}}\ar@{{}{ }{}}@/^1ex/[rr]|{\phantom{\mu_{a,b}}\Longdownarrow\mu_{a,b}}&&x
  }
  \qquad\qquad\rightsquigarrow\qquad\qquad
  \svxym{
    &y\ar@/^/[dr]^b&\\
    x\ar@/^/[ur]^a\ar@{=}@/_3ex/[rr]_{\phantom{b\circ a}}\ar@{{}{ }{}}@/^1ex/[rr]|{\phantom{\mu'_{a,b}}\Longdownarrow\mu'_{a,b}}&&x\pbox.
  }
  \]
  Note the subtle difference between $\mu_{f,g}$ and $\mu'_{f,g}$: in the first
  case the target is the path of length one consisting of the 1-generator
  $\id_x$, whereas in the second case it is the path of length zero at~$x$.
\item For each $x\in P_0$, remove the 1-generator $\id_x$ along with the 2-cell
  $\eta_x$ by using transformation \tgen{}: this can be done because this
  1-generator does not occur in the source or target of any 2-cell other
  than~$\eta_x$.
\end{itemize}
The resulting 2-polygraph is the reduced standard presentation. In fact, this
presentation is precisely the one that one would obtain by applying the
procedure described in the proof of \thmr{reduced-pres}.

\subsection{Tietze reductions}
\index{Tietze!reduction}
\index{reduction!Tietze}
\label{sec:tietze-reduction}
There are two kinds of Tietze transformations: \tgen{} adding a definable generator
and \trel{} adding a derivable relation. During a Tietze equivalence, those can
also be performed backward: \trgen{} removing a definable generator and \trrel{}
removing a derivable relation. A Tietze equivalence using only the two backward
transformations is called a \emph{Tietze reduction} and consists in making the
presentation smaller by suitable removing generators and relations.
It would be nice if two 2-polygraphs~$P$ and~$Q$ were Tietze equivalent if and
only if they reduce to a common 2-polygraph: this would mean that we do not have
to come up with new generators or relations in order to study Tietze
equivalence. We have seen in \cref{sec:1-minimal-pres} that this holds in the
case of $1$\nbd-poly\-graphs, but we show here that this is not the case for
$2$-polygraphs. This explains why the proof of \cref{thm:2-tietze-equiv}
proceeds by transforming two polygraphs into a bigger one, which contains both,
and not a smaller one.

Consider the presentation
\[
  P=
  \Pres{\star}{a,b}{\alpha:aa\To a,\beta:bb\To b,\gamma:aa\To bb}
  \pbox.
\]
One can apply to it the following Tietze transformations:
\begin{description}
\item[\trel] add the derivable relation $a=b$:
  \[
    \Pres{\star}{a,b}{aa=a,bb=b,aa=bb,a=b}
  \]
  (the relation is derivable by $a=aa=bb=b$),
\item[\trrel] remove the derivable relation $bb=b$:
  \[
    \Pres{\star}{a,b}{aa=a,aa=bb,a=b}
  \]
  (the relation is derivable by $bb=aa=a=b$),
\item[\trrel] remove the derivable relation $aa=bb$:
  \[
    \Pres{\star}{a,b}{aa=a,a=b}
  \]
  (respectively derivable since $a=b$),
\item[\trgen] remove the definable generator $b$:
  \[
    P'=
    \Pres{\star}{a}{aa=a}
    \pbox.
  \]
\end{description}
The polygraphs~$P$ and~$P'$ are thus Tietze equivalent and~$P$ presents the free
monoid with an idempotent element: this monoid has two elements $1$ and $a$,
with multiplication given by $1a=a1=aa=a$.

The polygraph~$P$ is \emph{Tietze minimal}, in the sense that no
non-trivial Tietze reduction can be applied to it; otherwise said, in order to
prove a non-trivial Tietze equivalence, one has to begin by adding definable
generators or derivable relations. Namely, the 1-generator $a$ cannot be removed
along the relation~$\alpha$ because $a$ occurs in the source of~$\gamma$, and
similarly for~$\beta$. Finally, we can show that no relation is derivable by
contradiction as follows.
\begin{itemize}
\item Suppose that $\alpha$ is derivable. This means that the 2-polygraph~$P$ is
  Tietze equivalent to the 2-polygraph
  \[
  Q=
  \Pres{\star}{a,b}{\beta:bb\To b,\gamma:aa\To bb}
  \]
  and we have $\prescat{P}\isoto\prescat{Q}$. The 2-polygraph~$Q$ is not
  convergent, but it is Tietze equivalent to the convergent 2-polygraph
  \[
  Q'=
  \Pres{\star}{a,b}{\beta:bb\To b,\gamma':aa\To b,\delta:ba\To ab}
  \]
  (the relation $\delta$ is derivable by $ba=bba=aaa=ab$, see also
  \cref{ex:KB-compl}). The termination of~$Q'$ can be shown using the deglex
  order generated by $b>a$, and the critical branchings are confluent:
  \begin{align*}
    \svxym{
      &\ar@{=>}[dl]_{\beta b}bbb\ar@{=>}[dr]^{b\beta}&\\
      bb\ar@{=}[rr]&&bb
    }
    &&
    \svxym{
      &\ar@{=>}[dl]_{\gamma' a}aaa\ar@{=>}[dr]^{a\gamma'}&\\
      ba\ar@{=>}[rr]_{\delta}&&ab
    }    
  \end{align*}
  \begin{align*}
    \svxym{
      &\ar@{=>}[ddd]_{\beta a}bba\ar@{=>}[dr]^{b\delta}&\\
      &&bab\ar@{=>}[d]^{\delta b}\\
      &&abb\ar@{=>}[dl]^{a\beta}\\
      &ba
    }
    &&
    \svxym{
      &\ar@{=>}[dl]_{\delta a}baa\ar@{=>}[ddd]^{b\gamma'}&\\
      aba\ar@{=>}[d]_{a\delta}&&\\
      aab\ar@{=>}[dr]_{\gamma'b}&&\\
      &bb\pbox.\\
    }    
  \end{align*}
  The normal forms are $1$, $a$, $b$ and $ab$, \ie there are four morphisms
  in~$\prescat{Q}$ whereas there are only two in~$\prescat{P}$, contradicting
  the isomorphism $\prescat{P}\isoto\prescat{Q}$.
\item By exchanging the role of~$a$ and~$b$ in previous case, and reversing the
  orientation of~$\gamma$ (which does not change the presented category), the
  relation~$\beta$ is not derivable either.
\item Suppose that $\gamma$ is derivable. This means that~$P$ is Tietze
  equivalent to the 2-polygraph
  \[
  Q=
  \Pres{\star}{a,b}{\alpha:aa\To a,\beta:bb\To b}
  \pbox.
  \]
  Again, $Q$ is convergent (rules decrease the length of morphisms and there is
  no critical pair) and every word of the form $ababab\ldots$ is a normal form,
  whereas $\pcat P$ has only two elements.
\end{itemize}
Since the presentation~$P'$ is also (obviously) minimal, we see that there is no
way to show that~$P$ and~$P'$ are Tietze equivalent by Tietze reducing both to a
common 2-polygraph.

\subsection{Tietze transformations up to equivalence}
We have seen in \cref{thm:2-tietze-equiv} that Tietze transformations generate
the following equivalence relation on polygraphs: two polygraphs are equivalent
when they present isomorphic categories.
We consider here the following variant of the notion of equivalence: two
polygraphs are \emph{equivalent} when they present equivalent categories. A
corresponding notion of Tietze transformation can be obtained as a variant of
those presented in \secr{2-tietze-def}, by adding the following kind of
transformation:
\begin{description}
\item[\tobj] \emph{adding an isomorphic 0-generator}: given $x\in P_0$,
  $y\not\in P_0$, $a,b\not\in P_1$, $\alpha,\beta\not\in P_2$, we define
  \[
    Q=\Pres{P_1,y}{P_1,a:x\to y,b:y\to x}{P_2,\alpha:ab\To\id_x,\beta:ba\To\id_y}
    \pbox.
  \]
\end{description}

\section{The Knuth-Bendix Completion Procedure}
\label{sec:2kb}
We have seen that convergent $2$-polygraphs are very convenient to work
with. When given a polygraph which does not have this property, one can in many
cases use Tietze transformations to turn it into one which does,
preserving the presented category. We present here a procedure due to Knuth and
Bendix~\cite{knuth1970simple} (in the setting of term rewriting
systems) whose purpose is to perform this transformation in an automated way:
starting from a 2-polygraph with a reduction order, it adds definable 2-generators until
possibly reaching a convergent 2-polygraph, which is Tietze equivalent to the original
one. We use the terminology of ``procedure'' and not an ``algorithm'', because
there is no guarantee that it will eventually stop, although it very often does
in practice.

This procedure is based on two observations. The first one is that,
in a terminating $2$-polygraph, the completion of any confluent
critical branching can always be chosen to be convergent toward a normal
form. In fact, suppose that a critical pair 
$(\phi_1,\phi_2)$ is closed by $(\phi_1',\phi_2')$ as shown on the
left diagram below:
\[
\svxym{
  &\ar@{=>}[dl]_{\phi_1}u\ar@{=>}[dr]^{\phi_2}&\\
  v_1\ar@{:>}[dr]_{\phi_1'}&&\ar@{:>}[dl]^{\phi_2'}v_2\\
  &w&\\
  &{}\phantom{\nf{w}}&
}
\qquad\qquad\qquad
\svxym{
  &\ar@{=>}[dl]_{\phi_1}u\ar@{=>}[dr]^{\phi_2}&\\
  v_1\ar@{:>}[dr]_{\phi_1'}&&\ar@{:>}[dl]^{\phi_2'}v_2\pbox.\\
  &w\ar@{:>}[d]^\psi&\\
  &\nf{w}&
}
\]
The termination property yields a normalization path $\psi$ from $w$ to a
normal form~$\nf{w}$ of~$w$, so we may close the diagram by the new
pair~$(\psi\circ\phi_1',\psi\circ\phi_2')$ as shown above on the right.

The second observation is that, given a non-confluent critical pair as on
the left below,
\[
\svxym{
  &\ar@{=>}[dl]_{\phi_1}u\ar@{=>}[dr]^{\phi_2}&\\
  v_1&&v_2\\
}
\qquad\quad\rightsquigarrow\qquad\quad
\svxym{
  &\ar@{=>}[dl]_{\phi_1}u\ar@{=>}[dr]^{\phi_2}&\\
  v_1\ar@{:>}[rr]&&v_2\\
}
\quad\text{or}\quad
\svxym{
  &\ar@{=>}[dl]_{\phi_1}u\ar@{=>}[dr]^{\phi_2}&\\
  v_1&&\ar@{:>}[ll]v_2\\
}
\]
we have $v_1\approx v_2$ and it is therefore possible to add the definable
relation $v_1\To v_2$ or $v_2\To v_1$ to the polygraph without changing the
presented category since this is a Tietze transformation of type
\trel{}.  The new 
presentation is ``more confluent'' in the sense that the above critical pair is
now confluent. We are thus tempted to add new rules in this way for
every critical branching. However, newly added rules can create new
non-confluent branchings and we want therefore to add as few of them as
possible. For instance, in the above situation, suppose that we have added a
rule $v_1\To v_2$ and that there was already another reduction $v_1\To v_1'$,
making a non-confluent branching, as shown on the left below:
\[
\svxym{
  &\ar@{=>}[dl]_{\phi_1}u\ar@{=>}[dr]^{\phi_2}&\\
  v_1\ar@{=>}[d]\ar@{=>}[rr]&&v_2\\
  v_1'
}
\quad\rightsquigarrow\quad
\svxym{
  &\ar@{=>}[dl]_{\phi_1}u\ar@{=>}[dr]^{\phi_2}&\\
  v_1\ar@{=>}[d]\ar@{=>}[rr]&&v_2\\
  v_1'\ar@{:>}[urr]
}
\quad\rightsquigarrow\quad
\svxym{
  &\ar@{=>}[dl]_{\phi_1}u\ar@{=>}[dr]^{\phi_2}&\\
  v_1\ar@{=>}[d]&&v_2\pbox.\\
  v_1'\ar@{=>}[urr]
}
\]
In order to make the polygraph confluent, we are now forced to add a new rule
between $v_1'$ and $v_2$, say $v_1'\To v_2$, making the former rule
$v_1\To v_2$ useless: it would have been preferable to directly add 
the rule $v_1'\To v_2$ instead of $v_1\To v_2$. For this reason, given
a critical branching as above, we only add new rules $\nf{v_1}\To\nf{v_2}$ (or
$\nf{v_2}\To\nf{v_1}$), between normal forms $\nf{v_1}$ (\resp $\nf{v_2}$) of
$v_1$ (\resp $v_2$). Finally the newly added
rules must be oriented without breaking the termination
of the original polygraph. This is usually done by orienting rules according to
a reduction order.

\subsection{The completion procedure}
\index{Knuth-Bendix completion}
\index{completion}
\index{procedure!completion}
Suppose given a finite 2-polygraph~$P$, equipped with a total reduction
order~$\preceq$ which is compatible with~$P$, \ie $u\succ v$ for every
2-generator $\alpha:u\To v$ in~$\freecat{P_2}$. By \propr{term-red-ord}, the
polygraph $P$ is necessarily terminating.

The \emph{Knuth-Bendix completion procedure} starts with the 2-polygraph~$P$ and
iteratively transforms it by adding definable relations, as follows.
\begin{enumerate}
\item For every critical branching
  \[
    \vxym{
      v&\ar@{=>}[l]_-{\phi}u\ar@{=>}[r]^-{\psi}&w
    }
  \]
  we compute reduction paths
  $
  \phi':v\overset*\To\nf{v}
  $ and $
  \psi:w\overset*\To\nf{w}
  $
  to some normal forms $\nf{v}$ and $\nf{w}$ of $v$ and $w$ respectively, until
  finding one with $\nf{v}\neq\nf{w}$. If there is none the procedure halts and
  returns the computed polygraph.
\item With the normal forms computed in previous step, we either have
  $\nf{v}\succeq\nf{w}$, in which case we add a 2-generator
  $\alpha:\nf{v}\To\nf{w}$ to~$P$, or $\nf{v}\preceq\nf{w}$, in which case we
  add a 2-generator $\alpha:\nf{w}\To\nf{v}$ to~$P$:
  \[
  \svxym{
    &\ar@{=>}[dl]_\phi u\ar@{=>}[dr]^\psi&\\
    v\ar@{=>}[d]_{\phi'}&&\ar@{=>}[d]^{\psi'}w\\
    \nf{v}\ar@{:>}[rr]_\alpha&&\nf{w}
  }
  \qquad\qquad\qquad
  \svxym{
    &\ar@{=>}[dl]_\phi u\ar@{=>}[dr]^\psi&\\
    v\ar@{=>}[d]_{\phi'}&&\ar@{=>}[d]^{\psi'}w\\
    \nf{v}&&\ar@{:>}[ll]^\alpha\nf{w}\pbox.
  }
  \]
\item Go back to step 1.
\end{enumerate}

If the procedure stops, it returns a $2$-polygraph, which we denote as $\KB(P)$
and call a \emph{Knuth-Bendix completion} of~$P$.
In the case where the procedure does not terminate, it constructs an infinite
sequence $P=P^0,P^1,P^2,\ldots$ of $2$-polygraphs, where $P^{i+1}$ is obtained
from $P^i$ by adding a derivable relation. This sequence is thus increasing, in
the sense that we have $P^i\subseteq P^j$ for $i\leq j$, and thus admits an
inductive limit $\bigcup_iP^i$, which we still denote as $\KB(P)$.

\begin{theorem}[\cite{knuth1970simple, huet1981complete}]
  The Knuth-Bendix completion $\KB(P)$ of a 2\nbd-poly\-graph~$P$ is a
  convergent presentation of the category~$\pcat P$.
\end{theorem}
\begin{proof}
  Since all the rules respect the termination order by construction, the
  reduction order $\preceq$ is a termination order, and the polygraph~$\KB(P)$
  is thus terminating by \cref{prop:term-red-ord}. Moreover, step 1 ensures
  that all the critical branchings are confluent, and the polygraph~$\KB(P)$ is
  thus locally confluent by \cref{lem:2-cb} and confluent by
  \cref{lem:newman}. Finally, the procedure proceeds by adding derivable
  transformations at step 2, \ie by performing Tietze transformations of type
  \trel{}. By \cref{thm:2-tietze-equiv}, the polygraph~$\KB(P)$ thus presents
  the same category as~$P$.
\end{proof}

\noindent
Note that the above theorem applies in both the cases where the procedure
terminates and where it does not. It can moreover be noted that the
$2$-polygraph $\KB(P)$ is finite if and only if the $2$-polygraph $P$ is finite
and the Knuth-Bendix completion procedure halts. For implementation purposes, we
are thus mostly interested in the cases where the procedure computes a result
after a finite amount of time, but for theoretical purposes it is still useful
when it runs indefinitely.
It is also interesting to remark that if the starting $2$-polygraph $P$ is
already convergent, we immediately have $\KB(P)=P$.

\begin{example}
  \label{ex:KB-compl}
  Consider the 2-polygraph
  \[
  P=\Pres{\star}{a,b}{bb\To b,aa\To bb}
  \]
  already encountered in \cref{sec:tietze-reduction}, equipped with the deglex
  order generated by $b>a$, which is compatible with~$P$. The two critical
  branchings are
  \[
  \svxym{
    &\ar@{=>}[dl]bbb\ar@{=>}[dr]&\\
    bb\ar@{:>}[d]&&\ar@{:>}[d]bb\\
    b&&b
  }
  \qquad\qquad\qquad
  \svxym{
    &\ar@{=>}[dl]aaa\ar@{=>}[dr]&\\
    bba\ar@{:>}[d]&&\ar@{:>}[d]abb\\
    ba&&ab
  }
  \]
  and the dotted arrows are chosen normalization 1-cells. In the first case, the
  two normal forms are equal, but not in the second one. Since $ba>ab$, the
  Knuth-Bendix procedure adds a rule $ba\To ab$, thus obtaining the 2-polygraph
  \[
    P=\Pres{\star}{a,b}{bb\To b,aa\To bb,ba\To ab}
    \pbox.
  \]
  Once this new rule added, all the critical pairs are confluent (see
  \secr{tietze-reduction}), so that the procedure halts on the above convergent
  2-polygraph.
\end{example}

\begin{example}
  \label{ex:squier-lafont-monoid}
  Consider the following 2-polygraph from~\cite{lafont1991church}
  \[
    P=
    \Pres{\star}{a,b,c,d}{\alpha_0:ab\To a,\beta:da\To ac},
  \]
  equipped with the deglex order associated to the reverse alphabetic order.
  The Knuth-Bendix completion does not terminate and gives rise to the infinite
  convergent presentation
  \[
    P=
    \Pres{\star}{a,b,c,d}{\alpha_n:ac^nb\To ac^n,\beta:da\To ac}_{n\in\N}
    \pbox.
  \]
  Namely, at the $n$-th step of the procedure the rule $\alpha_{n+1}$ is added
  by closing the critical branching
  \[
    \svxym{
      &\ar@{=>}[dl]_{\beta c^nb}dac^nb\ar@{=>}[dr]^{d\alpha_n}&\\
      ac^{n+1}b\ar@{:>}@/_3ex/[drr]_{\alpha_{n+1}}&&dac^n\ar@{=>}[d]^{\,\beta c^n}\\
      &&ac^{n+1}\pbox.
    }
  \]
  It can be remarked that if we take the converse orientation for rule $\beta$
  \[
    P=
    \Pres{\star}{a,b,c,d}{\alpha_0:ab\To a,\beta:ac\To da}
  \]
  and equip the polygraph with the deglex order associated to the reverse
  alphabetic order, the procedure halts immediately since there is no critical
  branching.
\end{example}

As illustrated in the above example, the procedure depends on many parameters,
each of which can have a strong influence on the output of the procedure, \ie
how small the completed polygraph will be, or even the termination of the
procedure: the termination order, the order in which critical pairs are studied
in step 1, the normal forms $\nf{v}$ and $\nf{w}$ chosen for each critical pair
in step 1.

\subsection{Detailed description of the procedure}
The procedure can be improved so that it produces reduced polygraphs, by
combining it with the procedure presented in~\secr{2-reduced}. It can also be
more efficiently implemented by observing that if a pair is confluent at some
stage, then it is still confluent if new rules are added, therefore one can
restrict step 1 to consider only critical branchings formed by newly added
rules. In a more operational way, close to the presentation given by
Huet~\cite{huet1981complete}, the resulting improved procedure can be described
as follows.

We generalize here slightly the situation of the previous section and suppose given
a 2-polygraph~$P$ together with a reduction order which is not necessarily
compatible with~$P$ (the procedure will reorient the rules anyway) and may be
partial.
The procedure will modify the following variables:
\begin{itemize}
\item a set~$E$ of equations, \ie pairs $u=v$ with $u,v\in\freecat{P_1}$,
  whose initial value is
  \[
  E=
  \setof{u=v}{\alpha:u\To v\in P_2}
  \]
\item a polygraph~$Q$, which is initially the polygraph~$P$ where the set of
  rules has been replaced by the empty set:
  \begin{align*}
    Q_0&=P_0
    &
    Q_1&=P_1
    &
    Q_2&=\emptyset
    \pbox.
  \end{align*}
\end{itemize}
The procedure repeats the following steps until we have $E=\emptyset$:
\begin{enumerate}
\item pick an equation $u=v$ in~$E$ and remove it from~$E$,
\item compute normal forms $\nf{u}$ and $\nf{v}$ of $u$ and $v$,
\item if $\nf{u}=\nf{v}$ then go back to step 1,
\item if neither $\nf{u}\prec\nf{v}$ or $\nf{v}\prec\nf{u}$ then fail,
\item if $\nf{u}\prec\nf{v}$ then exchange $u$ and $v$ (and $\nf{u}$ and
  $\nf{v}$) so that $\nf{u}\succ\nf{v}$,
\item for each rule $\alpha_i:u_i\To v_i$ in~$Q_2$ such that $u_i$ rewrites to
  $u_i'$ by the rule~$u\To v$
  \begin{itemize}
  \item remove $\alpha_i$ from~$Q_2$,
  \item add $u_i'=v_i$ to~$E$,
  \end{itemize}
\item add $\alpha:u\To v$ to~$Q_2$,
\item replace each rule $\alpha_i:u_i\To v_i$ of~$Q_2$ by
  $\alpha_i:u_i\To\nf{v_i}$, where $\nf{v_i}$ is a normal form of~$v_i$ \wrt $Q$
  as computed in previous step,
\item in the polygraph~$Q$, for each critical branching
  \[
    \svxym{
      &\ar@{=>}[dl]_\phi w\ar@{=>}[dr]^\psi&\\
      u'&&v' }
  \]
  where~$\phi$ consists of the rule~$\alpha$ in context, add $u'=v'$ to~$E$.
\end{enumerate}
In the end, \ie when $E=\emptyset$ is reached after a finite number of steps,
the procedure returns the polygraph~$Q$. This polygraph is reduced and Tietze equivalent to
the original polygraph~$P$. It is also possible to reasonably define
a notion of outcome of the procedure when it does not terminate,
see~\cite{huet1981complete} for details.

\begin{example}
  \label{ex:kb-S3-r}
  Consider the presentation
  \[
  \Pres{\star}{a,b,c}{aba\To bab, ba\To c}
  \]
  obtained from the usual presentation of~$B_3^+$, see \secr{braid-mon}, by
  adding a generator~$c$ along with its definition~$ba=c$. We consider the
  deglex order induced by $a>b>c$, which is compatible with the rules. The
  procedure will
  \begin{itemize}
  \item replace $aba\To bab$ by $ac\To cb$,
  \item add the rule $bcb\To cc$ coming from the non-confluent critical
    branching
    \[
    \svxym{
      &\ar@{=>}[dl]bac\ar@{=>}[dr]&\\
      cc&&\ar@{:>}[ll]bcb\pbox,
    }
    \]
  \item add the rule $bcc\To cca$ coming from the non-confluent critical
    branching
    \[
    \svxym{
      &\ar@{=>}[dl]bcba\ar@{=>}[dr]&\\
      cca&&\ar@{:>}[ll]bcc\pbox.
    }
    \]
  \end{itemize}
  We finally obtain the convergent Tietze equivalent presentation
  \[
  \Pres{\star}{a,b,c}{ba\To c,ac\To cb,bcb\To cc, bcc\To cca}\pbox.
  \]
\end{example}

\subsection{The symmetric group}
\index{symmetric!group}
\index{group!symmetric}
\label{sec:pres-sym-group}
Let us work out a fundamental and non-trivial example of a presentation of a
monoid. Given $n\in\N$, we consider the symmetric group~$S_{n+1}$ of
bijections on a set with $n+1$ elements. We claim
that it admits a presentation by the $2$-polygraph
\[
  P
  =
  \Pres{\star}{a_0,\ldots,a_{n-1}}{\alpha_i,\beta_i,\gamma_{i,j}}
\]
where the 2-generators are
\[
  \begin{array}{r@{\ :\ }r@{\ \To\ }l@{\qquad}l}
    \alpha_i&a_ia_i&1&\text{for $0\leq i<n$,}
    \\
    \beta_i&a_{i+1}a_ia_{i+1}&a_ia_{i+1}a_i&\text{for $0\leq i<n-1$,}
    \\
    \gamma_{i,j}&a_ja_i&a_ia_j&\text{for $0\leq i<i+1<j<n$,}
  \end{array}
\]
see \secr{sym-group-pres} for details. Our strategy to show this result is based
on the following two steps.
\begin{enumerate}
\item We use the Knuth-Bendix completion procedure to compute a convergent
  $2$-polygraph~$Q$ presenting the same category as~$P$.
\item We use the techniques presented in \secr{2-presenting} to show that~$Q$ is
  a presentation of~$S_{n+1}$, by showing that elements of $\freecat Q_1$ in
  normal form are in bijection with the elements of~$S_{n+1}$.
\end{enumerate}

In order to apply the Knuth-Bendix procedure, we equip the polygraph~$P$ with
the deglex reduction order~$\preceq$ induced by $a_j>a_i$ whenever
$j>i$, which is compatible with the rules.
After a finite amount of steps, the Knuth-Bendix procedure terminates, producing
the convergent 2-polygraph~$Q$ with the same 0- and 1-generators, and with rules
\[
\begin{array}{r@{\ }c@{\ }r@{\ }c@{\ }l@{\quad}l}
  \alpha_i&:&a_ia_i&\To&1&\text{for $0\leq i<n$,}\\
  \beta^k_i&:&a_{i+k+1}\ldots a_ia_{i+k+1}&\To&a_{i+k}a_{i+k+1}\ldots a_i&\text{for $0\leq k<n$,}\\
  &&&&&\text{and $0\leq i<n-k-1$,}\\
  \gamma_{i,j}&:&a_ja_i&\To&a_ia_j&\text{for $0\leq i<i+1<j<n$,}
\end{array}
\]
where $a_{i+k+1}\ldots a_i$ denotes the sequence of $a_j$, with indices~$j$
decreasing one by one between $i+k+1$ and $i$. The reader is advised to compute
this by himself or refer to~\cite{le1986catalogue} for details.

In an element of~$\freecat Q_1$ in normal form, because of the rules $\alpha_i$
and~$\gamma_{i,j}$, if we have a factor $a_ia_j$, then we have $i<j$ or
$i=j+1$. Taking the rules~$\beta_{i,j}$ in account too, we see that the normal
forms are the words of the form
\[
  w_0w_1w_2\ldots w_{n-1}
  \qquad\qquad\text{with}\qquad\qquad
  w_i=a_ia_{i-1}a_{i-2}\ldots a_{i-k_i}
  \pbox.
\]
Now, let us show that those normal forms are in bijective correspondence with
the elements of~$S_{n+1}$, \ie bijections
$f:\intset{n+1}\to\intset{n+1}$. First, we interpret the generator~$a_i$ as the
bijection $\intp{a_i}:\intset{n+1}\to\intset{n+1}$ which exchanges $i$
and~$i+1$, and can be depicted as
\[
  \satex{trans-i}
\]
To any bijection~$f:\intset{n+1}\to\intset{n+1}$, we associate a 1-cell
$u_f\in\freecat{P_1}$ defined by induction on~$n$. We set $u_f=1$ whenever
$n=0$. Otherwise, we write $f':\intset{n}\to\intset{n}$ for the function
obtained from~$f$ by ``removing'' $n$ from the source of~$f$ and $f(n)$ from its
image, \ie
\[
f'(i)=
\begin{cases}
  f(i)&\text{if $f(i)<f(n)$}\\
  f(i)-1&\text{if $f(i)>f(n)$}\\
\end{cases}
\]
and define
\[
  u_f
  =
  u_{f'}a_{n-1}a_{n-2}\ldots a_{f(n)}
  \pbox.
\]
For instance, consider the bijection~$f:\intset{6}\to\intset{6}$ such that the
images of $0$, $1$, $2$, $3$, $4$, and $5$ are respectively $4$, $1$, $0$, $5$,
$2$, and $3$. Its associated word $u_f$ is $a_0a_1a_0a_3a_2a_4a_3$, which can be
pictured as
\[
\satex{sym_pres_ex_0103243}
\]
Finally, using the above description of normal forms, it can be shown that~$u_f$
is a normal form for any bijection~$f$, and that this provides a bijection
between normal forms and elements of~$S_{n+1}$. Other examples of such
completions for finite groups can be found in~\cite{GaussentGuiraudMalbos15,guiraud2013homotopical,le1986catalogue}.

\subsection{Generated subcategories}
\index{subcategory!generated}
\index{generated!subcategory}
\label{sec:1-gen-subcat}
\newcommand{\gencat}[1]{\langle #1\rangle}
As an application of the previously developed techniques, consider the following
situation. We suppose given a category~$C$ presented by a
2-polygraph~$P$ and a
set~$G$ of morphisms of~$C$, whose elements are called \emph{generators}. The
\emph{category generated} by~$G$, denoted $\gencat{G}$, is the smallest
subcategory of~$C$ which contains the elements of~$G$ as morphisms (and is
closed under identities and composition, source and target of morphisms
in~$C$). Our goal here is to compute a presentation of it: we will provide a
method to perform this in the case where~$C$ admits a suitable convergent
presentation.
Before addressing the general case, we look at the following example:
let $C$ be the monoid $\N/6\N$, presented by
\[
\Pres{\star}{a}{a^6\To 1}
\]
and let us compute the category generated by~$G=\set{a^2}$.

First, note that we can always suppose that each morphism~$f\in G$ admits a
1\nbd-gene\-rator~$b\in P_1$ as a representative. Otherwise, given a
representative~$u\in\freecat{P_1}$ of~$f$, we can apply to~$P$ the Tietze
transformation which consists in adding a new generator~$b$ together with the
rule~$u\To b$. In our example, this amounts to considering the presentation
\[
  \Pres{\star}{a,b}{a^6\To 1,a^2\To b}
  \pbox.
\]
For this reason, we will suppose in the following that the set~$G$ of generating
morphisms is a subset of the 1-generators, \ie $G\subseteq P_1$ (in the above
example, we have $G=\set{b}$). Moreover, we can suppose that the presentation is
convergent and reduced: if it is not the case, we can apply the Knuth-Bendix
procedure, and hope that it succeeds. For instance, with the previous
presentation, consider the deglex order with $a<b$. The critical branchings
\[
\xymatrix@R=2ex{
  &\ar@{=>}[dl]aaa\ar@{=>}[dr]&\\
  ba&&ab
}
\qquad\qquad\qquad
\xymatrix@R=2ex{
  &\ar@{=>}[dl]aaaaaa\ar@{=>}[dr]&\\
  baaaa\ar@{=>}[d]&&1\\
  bbaa\ar@{=>}[d]&&\\
  bbb&&
}
\]
are not confluent, and it can be checked that adding the induced relations
$ba\To ab$ and $b^3\To 1$ makes the presentation convergent and the rule
$a^6\To 1$ superfluous. We thus consider the alternative, convergent,
presentation
\[
\Pres{\star}{a,b}{a^2\To b,ba\To ab,b^3\To 1}
\]
of $\N/6\N$. Using \cref{prop:gen-pres} below, we can finally deduce that the
category~$\gencat{a^2}$ admits the presentation
\[
  \Pres{\star}{b}{b^3\To 1}
  \pbox.
\]
It is thus the monoid $\N/3\N$, as expected.

Below, given $G\subseteq P_1$, we write $\freecat G\subseteq\freecat P_1$ for
the set of morphisms in $\freecat P_1$ which can be expressed as composites of
generators in~$G$.

\begin{proposition}
  \label{prop:gen-pres}
  Suppose given a convergent 2-polygraph~$P$ together with a
  set~$G\subseteq P_1$ of generators, such that for every rule $u\To v$ in~$P_2$
  with $u\in\freecat{G}$ we have $v\in\freecat{G}$. Then the category generated
  by~$G$ admits a presentation by the polygraph~$Q$ where
  \begin{itemize}
  \item $Q_0\subseteq P_0$ consists of the sources and targets of elements
    of~$G$,
  \item $Q_1=G\subseteq P_1$,
  \item $Q_2\subseteq P_2$ consists of the rules $u\To v$ in~$P_2$ such that
    $u\in\freecat{G}$ and $v\in\freecat{G}$.
  \end{itemize}
\end{proposition}
\begin{proof}
  Since $\gencat{G}$ has to be closed under taking the source and target of
  morphisms in~$G$, it contains at least $Q_0$ as objects, and conversely, any
  composite of morphisms in~$G$ will have elements of $Q_0$ as source and
  target; $Q_0$ is thus precisely the set of $0$-cells of~$\gencat G$. The morphisms
  in~$\gencat{G}$ contain the equivalence classes $\eqc{a}$ of 1-generators
  $a\in G$, and since it is closed under composition and identities, its
  morphisms are precisely the equivalence classes of morphisms
  in~$\freecat{G}\subseteq\freecat{P_1}$. Finally, given $u,v\in\freecat{G}$
  such that $\eqc{u}=\eqc{v}$, since~$P$ is convergent both $u$ and $v$ rewrite
  to a common element~$w\in\freecat{P_1}$, \ie $u\overset*\To w$ and
  $v\overset*\To w$. By induction, the two rewriting paths contain only rules
  in~$P_1$ and $w\in\freecat{G}$, thus~$Q_2$ is sufficient to generate the
  required equivalence on elements of~$\freecat{G}$.
\end{proof}

\noindent
As a bonus, note that the 2-polygraph~$Q$ in the previous proposition is
necessarily convergent, because~$P$ is supposed to be so.

\begin{remark}
  Suppose that we start with a $2$-polygraph~$P$ together with a set
  $G\subseteq P_1$, such that the following property is satisfied: for every
  rule $u\To v$ in~$P_2$, $u\in\freecat{G}$ implies $v\in\freecat{G}$. In order
  to be able to apply previous proposition, we need to ensure that~$P$ is
  convergent and, if this is not the case, we can apply the Knuth-Bendix
  completion procedure in order to obtain a convergent polygraph. However, in
  general, the completed polygraph will not satisfy the property anymore. In
  order to improve this, the Knuth-Bendix procedure can be modified in order not
  to produce ``bad rules'', \ie rules of the form $u\To v$ with
  $u\in\freecat{G}$ and $v\in\freecat{P_1}\setminus\freecat{G}$, which prevent
  the resulting polygraph from satisfying the required property. Namely, the
  completion procedure adds new rules of the form $u\To v$ where both $u$ and
  $v$ are normal forms. In the case such a rule is ``bad'', it can be useful to
  add instead a rule $u'\To v$ where $u'\overset\ast\To u$
  and~$u'\in\freecat{P_1}\setminus\freecat{G}$.
\end{remark}

\begin{exo}
  A presentation for the symmetric groups~$S_n$ was constructed
  in~\secr{pres-sym-group}. Deduce from it a presentation for the alternating
  groups~$A_3$ and $A_4$, see \secr{alt-group}.
\end{exo}

\section{Universality of Finite Convergent Rewriting}
\index{universality of convergent rewriting}
\label{sec:universality}
We have seen in \secr{word-problem} that a finite convergent rewriting system always
has decidable word problem.
The question of \emph{universality of convergent
  rewriting} is the converse question, first asked by
Jantzen~\cite{jantzen1982semi}, see
also~\cite{Diekert86,jantzen1985note,jantzen1988confluent,bauer1984finite}:

\begin{quote}
  Given a category~$C$ admitting a finite presentation with decidable word
  problem, does it always admit a finite convergent presentation?
\end{quote}

\noindent
The answer to this question is negative, but showing this requires more tools
than we have at our disposal for now and will be handled in
\cref{chap:2-fdt,chap:2-homology}. We however study here restricted forms of the
question.

\subsection{Universality of Knuth-Bendix completion}
\index{braid!monoid}
\index{monoid!braid}
\label{sec:kapur-narendran}
A more restricted variant of the above question consists in wondering whether it
is always possible to add or remove relations to a 2-polygraph so that it
becomes convergent. Kapur and Narendran~\cite{kapur1985finite} have shown that
this is not the case, by considering the usual presentation of the braid
monoid~$B_3^+$, detailed in \secr{braid-mon}:
\begin{equation}
  \label{eq:B3-pres}
  P
  =
  \Pres{\star}{a,b}{aba\To bab}.
\end{equation}
They show that there is no finite convergent presentation of this monoid on the
same generators, see \cref{prop:kapur-narendran} below.
As a consequence, for such a presentation, the Knuth-Bendix procedure will never
end whichever reduction order or strategy for considering rules is adopted.

\begin{lemma}
  The polygraph~$P$ has decidable word problem.
\end{lemma}
\begin{proof}
  Since the only relation preserves the length of 1-cells, equivalence classes
  contain 1-cells of the same length and are therefore finite.
\end{proof}

\begin{proposition}
  \label{prop:kapur-narendran}
  There is no finite convergent rewriting system which is Tietze equivalent to
  the polygraph~$P$ by a sequence of Tietze transformation consisting only in
  adding or removing derivable relations.
\end{proposition}
\begin{proof}
  First notice that $abbab\overset*\Leftrightarrow babba$ is derivable in~$P$
  since we have
  \[
    abbab\Leftarrow ababa\Rightarrow abbab
    \pbox.
  \]
  More generally, by induction, it can be shown that
  \begin{equation}
    \label{eq:kn-ij}
    a^{i+1}b^{j+2}ab
    \overset*\Leftrightarrow
    bab^{i+2}a^{j+1}
  \end{equation}
  for every $i,j\in\N$. Namely, the base case where $i=j=0$ is handled above,
  and if we suppose that \eqref{eq:kn-ij} holds for some $i$ and $j$, we have
  \[
    a^{i+2}b^{j+2}ab
    \overset*\Leftrightarrow
    abab^{i+2}a^{j+1}
    \Rightarrow
    bab^{i+3}a^{j+1}
  \]
  and
  \[
    a^{i+1}b^{j+2}ab
    \Leftarrow
    a^{i+1}b^{j+1}aba
    \overset*\Leftrightarrow
    bab^{i+2}a^{j+2}\pbox,
  \]
  which constitute the induction step on~$i$ and~$j$, respectively. Another easy
  remark is that, for $n\in\N$, any word~$u$ such that the relation
  $u\overset*\Leftrightarrow b^nab$ (\resp $u\overset*\Leftrightarrow bab^n$) is
  derivable is of the form $u=b^{n-i}aba^i$ (\resp $u=a^ibab^{n-i}$) for some
  $i$ with $0\leq i\leq n$.
  Writing $\eqc{u}$ for the equivalence class of a word $u$ under
  $\overset*\Leftrightarrow$, we thus have
  \begin{align*}
    \eqc{b^nab}&=\setof{b^{n-i}aba^i}{0\leq i\leq n}
    &
    \eqc{bab^n}&=\setof{a^ibab^{n-i}}{0\leq i\leq n}
    \pbox.
  \end{align*}
  We now proceed by contradiction. Suppose given a finite 2-polygraph~$Q$
  convergent and Tietze equivalent to~$P$ by a sequence of Tietze transformations
  consisting only in adding or removing derivable relations. By
  \thmr{reduced-pres}, we can suppose that~$Q$ is reduced. Since,
  $\set{aba,bab}$ is an equivalence class, $Q$ should contain either
  $aba\To bab$ or $bab\To aba$. We suppose that we are in the former case, with the
  other one being similar. Writing~$l$ for the length of the longest left-hand
  side of a rule in~$Q$, and
  \begin{align*}
    u&=a^{l+1}b^{l+2}ab
    &
    v&=bab^{l+2}a^{l+1}
  \end{align*}
  we have $u \overset*\Leftrightarrow v$ and therefore both $u$ and $v$ should
  reduce to a common word. The only factors in those words whose equivalence
  class is not a singleton are of the form $b^nab$ or~$bab^n$. Therefore, we
  must have rules of the form~$b^nab\To w$ or~$bab^n\To w$. By the preceding
  remark, the word $w$ has to be of the form~$b^{n-i}aba^i$ (\resp
  $a^ibab^{n-i}$) with $0<i\leq n$ and therefore is reducible by the rule
  $aba\To bab$, which contradicts the assumption that the rewriting system is
  reduced.
\end{proof}

As an alternative example, it is shown in~\cite{jantzen1985note} that the monoid
(in fact, group) presented by
\[
  \Pres{\star}{a,b}{abba=1}
\]
admits no finite convergent presentation on the same generators.

\subsection{Other Tietze transformations}
The result of \cref{prop:kapur-narendran} can be restated as follows: there is a
finite $2$-polygraph~$P$ which cannot be transformed into a finite convergent
one by using Tietze transformations \trel{} only. However, this does not bring a
definitive answer to the original question of universality of rewriting raised
at the beginning of this section, since it does not rule out the possibility of
turning a presentation into a convergent one by using both transformations
\tgen{} and \trel{}. In fact, this is the case for the presentation
\eqref{eq:B3-pres} of $B_3^+$. Namely, if we use a transformation \tgen{} to
introduce a generator~$c$ and a relation $ba=c$, the resulting presentation can
be completed into a convergent one, this was already detailed in \exr{kb-S3-r}: adding
a superfluous generator allows the Knuth-Bendix procedure to produce a
convergent presentation of~$B_3^+$. Finding a counter-example to the problem of
universality in full generality is much more difficult and will be addressed in
\cref{chap:2-fdt,chap:2-homology}.

The situation encountered for $B_3^+$, where the introduction of a definable
generator improves the properties of the presentation is not an ``isolated
case''. For instance for every natural number $n>3$, the plactic monoid $P_n$ of
type~$A$ does not have a finite presentation on the usual
generators~\cite{kubat2014grobner}. However, if we add the \emph{column
  generators}, we get a finite presentation~\cite{bokut2015new, cain2015finite,
  Hage15}, see Section~\ref{S:CoherentPresentationPlacticMonoids} for details on
convergent presentations of plactic monoids.
Modified Knuth-Bendix completion procedures have been proposed in order to
exploit this and allow for adding generators to handle such
situations~\cite{guiraud2013homotopical}.

\subsection{Conditions for convergence}
Since not every monoid admits a presentation by a finite convergent rewriting
system, a natural question is whether there are natural conditions on monoids
which ensure that this is the case. Diekert~\cite{Diekert86} has addressed this
question in the case of abelian groups: he derived a whole class of finite
string rewriting systems presenting abelian groups with decidable word problem,
which are not Tietze equivalent to a finite convergent string rewriting system
on the same alphabet. Moreover, he constructed necessary and sufficient
conditions for the existence of a convergent presentation for finitely generated
abelian groups. However, the question for general monoids was still open at this
time and new methods had to be introduced to solve this problem, which concerns intrinsic
properties of the presented monoid. In this direction, Squier introduced in
\cite{squier1987word,squier1994finiteness} homotopical and homological
approaches to formulate necessary conditions for a finitely presented monoid to
have a finite convergent presentation. The homotopical construction is presented
in \cref{chap:2-fdt} and the homological one in
\cref{chap:HomologieSquierTheorem}.


\chapter{Linear Rewriting}
\label{chap:2linear}
\label{C:OneDimensionalLinearRewriting}

This chapter presents rewriting techniques for associative
algebras. 
We look here for algorithms turning a given presentation by
generators and relations into a rewriting system by orienting the
latter, thereby producing linear bases of the presented algebra.
In particular, this approach applies to various fundamental decision
problems, such as the word problem, ideal membership, or to compute
quadratic bases, e.g., Poincaré-Birkhoff-Witt bases, Hilbert series,
syzygies of presentations, homology groups, and Poincaré series.
However, if we require rewriting rules to be compatible with the
linear structure, we immediately face the following problem: for any
rule $u\to v$, we also have $-u\to -v$ and thus
\[
v=-u+(u+v)\to -v+(u+v)=u\pbox.
\]
Therefore,
$u\to v$ implies $v\to u$ and thus no rewriting system can be terminating. 
In order to fix this problem, one can either restrict rewriting to be
decreasing with respect to a monomial order, as in the non-commutative
Gröbner basis approach~\cite{bergman1978diamond,Bokut76,Mora94}, or consider the structure of linear polygraph introduced in~\cite{GuiraudHoffbeckMalbos19} with an appropriate notion of reduction. It is the latter notion that we present in this chapter.

We first introduce linear polygraphs as a framework for linear rewriting
in \cref{sec:linear-rewriting}. We then study the confluence
properties of linear polygraphs in
\cref{sec:linear-rewritig-properties}. Finally, in
\cref{SS:PBWGrobner}, we express
Gröbner bases and Poincaré-Birkhoff-Witt bases in the setting of
linear polygraphs.
The polygraphic approach presented in this chapter subsumes many linear rewriting models developed throughout the 20th century. We present a brief historical overview of these works in \cref{S:HistoryLinear}.

The way to define rewriting in associative algebras depends on the definition considered for the associative algebra structure, either as an internal monoid in the category of vector spaces, or a linear category with a single object~\cite{mitchell1972rings}. 
In this chapter, we consider the first point of view, as introduced in~\cite{GuiraudHoffbeckMalbos19}.

\section{Linear Rewriting}
\label{sec:linear-rewriting}

In this section, we introduce the notion of rewriting in associative algebras.

\subsection{Associative algebras}
\nomenclature[Vect]{$\Vect\kk$}{category of vector spaces}
\index{algebra}
\index{associative!algebra}
Suppose fixed a ground field~$\kk$. An \emph{(associative) algebra}~$(A,m,e)$ consists of a
$\kk$-vector space $A$ together with an operation $m:A\otimes A\to A$ and an
element $e\in A$ such that the operation~$m$ is associative and admits~$e$ as
neutral element.
Otherwise said, an algebra is a monoid object (see \cref{ex:monoid}) in the category $\Vect{}$
of vector spaces and linear maps. A \emph{morphism of algebras} $\varphi : (A,m,e) \to (B,m',e')$ is a linear map $\varphi : A \to B$ which is compatible with operations $m$ and $m'$ and the neutral elements:
\begin{align*}
\varphi(m(x,y)) &= m'(\varphi(x),\varphi(y)),
&
\varphi(e) &= e',
\end{align*}
for all $x,y$ in $A$.
We denote by $\Alg$ the category of algebras and their morphisms.
\nomenclature[Alg]{$\Alg$}{category of algebras}

\subsection{Free algebras}
\index{free!algebra}
\index{algebra!free}
\nomenclature[supp]{$\supp{x}$}{support of a cell~$x$}
Given a set $P_0$, we will denote by $\lin{P}_0$ the free algebra over $P_0$. A
\emph{monomial of~$\lin{P}_0$} is an element of the free monoid~$P_0^\ast$
over~$P_0$. The monomials of~$\lin{P}_0$ form a linear basis of the algebra~$\lin{P}_0$, thus every $0$-cell~$p$ of~$\lin{P}_0$ can be uniquely
written as a linear combination
\[
  p = \sum_{i=1}^k\lambda_i u_i
\]
of pairwise distinct monomials~$u_1$, \dots, $u_k$ of~$\lin{P}_0$,
with~$\lambda_1$, \dots, $\lambda_p$ non-zero scalars, called \emph{the canonical
  decomposition of~$p$}. We define the \emph{support of~$p$}\index{support} as the set
$\supp{p}=\set{u_1,\dots,u_k}$.

\subsection{Linear 1-polygraphs}
\index{polygraph!linear}
\index{1-polygraph!linear}
\index{linear!1-polygraph}
A \emph{linear 1-polygraph} consists of a set $P_0$, together with a set $P_1$
equipped with two functions $\src0,\tgt0:P_1\to\lin P_0$.
Such a polygraph is thus characterized by a diagram of sets and functions
\begin{equation}
  \label{eq:lin-1pol}
  \vcenter{
  \xymatrix @=1.5em {
    P_0\ar[d]_{\ins0}
    & P_1
    \ar@<-0.5ex> [dl] _-{\src0}
    \ar@<0.5ex> [dl] ^-{\tgt0}
    \\
    \lin{P}_0
  }}
\end{equation}
where $\lin{P}_0$ is the free algebra over a set~$P_0$
and $\ins0:P_0\to\lin P_0$ is the canonical inclusion. We often write
\[
  \pres{x_i}{\alpha_i:u_i\to v_i}
\]
for a $1$-polygraph with the $x_i$ as elements of~$P_0$ and the $\alpha_i$ as
elements of~$P_1$ with $\src0(\alpha_i)=u_i$ and $\tgt0(\alpha_i)=v_i$.

\subsection{One-dimensional algebras}
\label{sec:1alg}
A \emph{1-algebra} is a category internal to $\Alg$. It thus consists of a
diagram
\[
  \xymatrix{
  A_0&\ar@<-.5ex>[l]_{\src0}\ar@<.5ex>[l]^{\tgt0}A_1
}
\]
comprising two algebras $A_0$ and $A_1$, whose elements are respectively called
\emph{0-} and \emph{1-cells}, with two algebra morphisms $\src0,\tgt0:A_1\to A_0$
respectively providing the \emph{source} and \emph{target} of a $1$-cell,
together with an algebra morphism $i:A_0\to A_1$ which to every $0$-cell $p$ associates
the \emph{identity} $i(p)$ on~$p$, and an algebra morphism $m:A_1\times_{A_0}A_1\to A_1$
which to every pair of composable $1$-cells associates their composite, in such a
way that composition is associative and admits identities as neutral
elements. According to our notations for categories, we set
$m\pair{\phi}{\phi'}=\phi\comp0\phi'$ for any pair $\phi$, $\phi'$ of
composable $1$-cells.

\begin{lemma}
Let $A$ be 1-algebra, then
\begin{itemize}
\item for all composable 1-cells $\phi$ and $\phi'$ in $A$,
\begin{equation}
\label{E:FormuleLinearityComposition}
\phi \comp0 \phi' = \phi - t_0 (\phi) + \phi',
\end{equation}
\item every 1-cell $\phi$ in $A$ is invertible with inverse
\[
\phi^{-} = s_0(\phi) - \phi + t_0(\phi),
\]
\item the product of two 1-cells $\phi,\phi'$ in $A$ decomposes into
\begin{equation}
  \label{E:FormuleLinearityProduct}
  \begin{split}
    \phi\phi'
    &= \phi s_0(\phi') + t_0(\phi)\phi' - t_0(\phi)s_0(\phi')\\
    &= s_0(\phi)\phi'+\phi t_0(\phi')-s_0(\phi)t_0(\phi').
  \end{split}
\end{equation}
\end{itemize}
\end{lemma}
\begin{proof}
For any composable $1$-cells $\phi$ and $\phi'$  in $A$, we have
\[
\phi\comp0 \phi' = (\phi -s_0(\phi') + s_0(\phi')) \comp0 (t_0(\phi)-t_0(\phi) + \phi').
\]
By linearity of the $0$-composition, this implies
\[
\phi \comp0 \phi' = \phi\comp0 t_0(\phi) - s_0(\phi')\comp0 t_0(\phi) + s_0(\phi')\comp0 \phi'
\]
and by neutrality of identities we get~\cref{E:FormuleLinearityComposition}.

The second condition is deduced from the first one. Let $\phi$ be a $1$-cell in $A$, we set
$\phi^{-} = s_0(\phi) - \phi + t_0(\phi)$. We have $s_0(\phi^{-})=t_0(\phi)$ and $t_0(\phi^{-})=s_0(\phi)$. Moreover, from~\cref{E:FormuleLinearityComposition}, we have
$\phi \comp0 \phi^{-} = s_0(\phi)$ and $\phi^{-}\comp0 \phi = t_0(\phi)$. We have thus proved that $\phi^{-}$ is $0$-inverse of $\phi$.

Let us prove the third condition. Let $\phi,\phi'$ be $1$-cells in $A$. The product of these two $1$-cells in $A_1$ decomposes into 
\[
\phi\phi' = (\phi\comp0 t_0(\phi))(s_0(\phi')\comp0 \phi').
\]
With the $0$-composition being an algebra morphism, we deduce that
\[
\phi\phi' = \phi s_0(\phi')\comp0 t_0(\phi)\phi'.
\]
From~\cref{E:FormuleLinearityComposition}, we deduce the first equality in~\cref{E:FormuleLinearityProduct}. The second equality is proved symmetrically.  
\end{proof}

\subsection{Free 1-algebras}
\label{sec:free-1-alg}
A linear $1$-polygraph~$P$ generates a \emph{free 1-algebra}, denoted by~$\lin P$, with
$\lin P_0$ as algebra of $0$-cells and an algebra $\lin P_1$ of $1$-cells that we
now describe. A \emph{1-monomial} of $\lin P$ is a triple
\begin{equation}
  \label{eq:lin-1-mon}
  u\alpha v
\end{equation}
with $u,v\in P_0^*$ monomials and $\alpha\in P_1$. We respectively define the
source and target of such a monomial by
\begin{align*}
  \srcl0(u\alpha v)&=u \src0(\alpha)v\pbox,
  &
  \tgtl0(u\alpha v)&=u \tgt0(\alpha)v\pbox.
\end{align*}
We consider $(\lin P_0\otimes\kk P_1\otimes\lin P_0)\oplus\lin P_0$ the free $\lin P_0$-bimodule on $1$-monomials, and we form the $\lin P_0$-bimodule
\[
  \lin P_1
  =
  (\lin P_0\otimes\kk P_1\otimes\lin P_0)\oplus\lin P_0/\sim,
\]
whose elements are linear combinations of the form
\begin{equation}
  \label{eq:lin-phi}
  \phi=\sum_i\lambda_i\phi_i+\unit p,
\end{equation}
where the $\phi_i$ are distinct monomials and $\unit p$ is a formal identity on
a $0$-cell~$p\in\lin P_0$, quotiented by the relation $\sim$ generated by the relations
\begin{equation}
  \label{eq:lin-xch}
  \phi\srcl0(\psi) + \tgtl0(\phi)\psi-\tgtl0(\phi)\srcl0(\psi) = \srcl0(\phi)\psi + \phi \tgtl0(\psi) - \srcl0(\phi) \tgtl0(\psi),
\end{equation}
where $\phi$ and $\psi$ range over $1$-monomials. 
The relation \eqref{eq:lin-xch} encodes a linear version
of the exchange law. The multiplication of the algebra structure in $\lin P_1$
precisely associates to two $1$-cells $\phi$ and $\psi$ the cell defined by
either member of \eqref{eq:lin-xch}. The source and target maps are the above
functions $\srcl0$ and $\tgtl0$ on monomials, extended by linearity.
Given a $1$-cell $\phi$ in $\lin P$, its \emph{size}\index{size} is the minimum number of
$1$-monomials~$\phi_i$ occurring in a decomposition of the form
\eqref{eq:lin-phi} of~$\phi$. In particular, a monomial is of size~$1$.

If we write $\ins1:P_1\to\lin P_1$ for the canonical inclusion, sending
$\alpha\in P_1$ to the monomial \eqref{eq:lin-1-mon} where $u$ and $v$ are the
empty words, we obtain a diagram
\[
  \xymatrix{
    P_0\ar[d]_{\ins0}
    & P_1
    \ar@<-0.5ex> [dl] _-{\src0}
    \ar@<0.5ex> [dl] ^-{\tgt0}
    \ar[d]^{\ins1}
    \\
    \lin{P}_0
    &
    \ar@<-0.5ex> [l] _-{\srcl0}
    \ar@<0.5ex> [l] ^-{\tgtl0}
    \lin P_1
  }
\]
which ``commutes'' in the sense that we have $\srcl0\circ\ins1=\src0$ and
$\tgtl0\circ\ins1=\tgt0$. In the following, we often simply write respectively
$\src{}(\phi)$ and $\tgt{}(\phi)$ instead of $\srcl0(\phi)$ and $\tgtl0(\phi)$ for
the source and target of a $1$-cell~$\phi$.

The free $1$-algebra $\lin P$ is characterized by the following universal property:

\begin{lemma}
  Suppose given a linear 1-polygraph~$P$, a 1-algebra~$A$ with $\lin{P}_0$ as 
  underlying algebra of $0$-cells and a function $f:P_1\to A_1$
  such that for every 1-generator $\alpha:u\to v$ in $P_1$, we have
  $f(\alpha):u\to v$.
  Then there exists a unique morphism of $1$-algebras
  \[
    \freecat{f}:\lin P \to A
  \]
such that
  $\freecat{f}(\alpha)=f(\alpha)$ for every~$\alpha$ in~$P_1$, seen
  as a 1-cell in~$\lin{P}$.
\end{lemma}

In the sequel, we will use the following decomposition result.

\begin{lemma}
  \label{L:NCellDecomposition}
  Let~$P$ be a linear 1-polygraph. Then, every
  non-identity 1-cell~$\phi$
  of~$\lin P$ admits a decomposition
  $
  \phi = \phi_1 \comp{0} \cdots \comp{0} \phi_k
  $, for some $k\in\N$, 
  where the $\phi_i$ are 1-cells of size~1 in $\lin P$.
\end{lemma}
\begin{proof}
The $1$-cell $\phi$ decomposes into $\phi=\lambda_1\psi_1 + \ldots + \lambda_k\psi_k +\unit p$.
When~$k=1$, the $1$-cell~$\phi$ is of size~$1$. Otherwise, for any $i\in\{1,\ldots,k\}$, we set
\begin{align*}
\alpha_i &= \lambda_1 t(\psi_1) + \ldots + \lambda_i t(\psi_i)\pbox,
&
\beta_i &= \lambda_1 s(\psi_1) + \ldots + \lambda_k s(\psi_k)  
\end{align*}
and $\alpha_0 = \beta_{p+1} =0$. For each~$i\in\{1,\ldots ,k\}$, we define the $1$-cell of size~$1$
\[
\phi_i = \lambda_i \psi_i + \unit p + \unit {\alpha_{i-1}} + \unit {\beta_{i+1}}.
\]
We have $s(\phi_i) = p + \alpha_{i-1} + \beta_i$ and $t(\phi_i) = p + \alpha_i + \beta_{i+1}$, so that $\phi_1\comp0\cdots\comp0\phi_k$ is a well-defined $1$-cell of~$\lin P$. Following relation~\eqref{E:FormuleLinearityComposition}, we deduce
\[
\phi_1 \comp0 \cdots \comp0 \phi_p 
	= \sum_{i=1}^k \lambda_i \psi_i 
	+ \sum_{i=1}^k (1_p + 1_{\alpha_{i-1}} + 1_{\beta_{i+1}}) 
	- \sum_{i=1}^{k-1} (\lambda_i 1_{t(\psi_i)} + 1_p + 1_{\alpha_{i-1}} + 1_{\beta_{i+1}}).
\]
We conclude thanks to $\alpha_{k-1} = \lambda_1 t(\psi_1) + \ldots +  \lambda_{k-1} t(\psi_{k-1})$, and $\beta_{k+1}=0$.
\end{proof}

\subsection{Presentations and ideals of linear polygraphs}
\index{presentation!of an algebra}
Let~$P$ be a linear $1$\nbd-poly\-graph. The \emph{algebra presented by~$P$} is the
quotient algebra
$
  \cl{P} = \lin{P}_0 / P_1
$
of the algebra $\lin{P}_0$ by the congruence generated by the $1$-generators in $P_1$.
We will denote by~$\cl{p}$ the image of a $0$-cell~$p$ of~$\lin{P}_0$ through
the canonical projection.  We say that an algebra~$A$ is \emph{presented
  by~$P$}, or that~$P$ is a \emph{presentation of~$A$}, if~$A$ is isomorphic
to~$\cl{P}$. Two linear $1$-polygraphs $P$ and $Q$ are \emph{Tietze equivalent}
when they present isomorphic algebras: $\cl P\isoto\cl Q$.

We define the \emph{boundary} of a $1$-cell~$\phi$
in the free $1$-algebra~$\lin{P}$, as the $0$-cell
\[
  \bnd(\phi) =  \tgt{0}(\phi) - \src{0}(\phi)
  \pbox.
\]
We denote by~$I(P)$ the ideal of the algebra~$\lin{P}_0$ generated by the
boundaries of the $1$-cells in~$P_1$. Since the algebra~$\lin{P}_0$ is free, the
ideal~$I(P)$ is consists of all the linear combinations
\[
  \sum_{i=1}^k \lambda_i u_i \bnd(\alpha_i) v_i
  \pbox,
\]
where the~$u_i\alpha_i v_i$ are pairwise distinct
$1$-monomials of~$\lin{P}$, and the~$\lambda_i$ are non-zero
scalars, so that the algebra~$\cl{P}$ is isomorphic to the
quotient of~$\lin{P}_0$ by~$I(P)$.

\begin{example}
  \index{Weyl algebra}
  \index{algebra!Weyl}
\label{Example:WeylAlgebras1}
The \emph{Weyl algebra} of dimension $n$ over a field~$\kk$ of characteristic zero is the algebra presented by the linear $1$-polygraph whose $0$-cells are
\[
x_1,\ldots,x_n,\partial_1,\ldots,\partial_n
\]
and with the following $1$-cells:
\begin{align*}
  x_ix_j&\to x_jx_i
  &
  \partial_i\partial_j&\to \partial_j\partial_i
  &
  \partial_i x_j &\to x_j\partial_i
  &&
  \text{for any $1\leq i < j \leq n$,}
  \\
  &&&&
  \partial_ix_i &\to x_i\partial_i +1
  &&
  \text{for any $1\leq i \leq n$.}
\end{align*}
\end{example}

\begin{lemma}
  \label{L:Ideal}
  Let~$P$ be a linear 1-polygraph. For all 0-cells~$p$ and~$q$
  of~$\lin{P}_0$, the following two conditions are equivalent:
  \begin{enumerate}
  \item The 0-cell~$q-p$ belongs to the ideal~$I(P)$.
  \item There exists a 1-cell $\phi:p\to q$ in the free 1-algebra~$\lin{P}$.
  \end{enumerate}
  As a consequence, $I(P)$ exactly contains the 0-cells~$p$ of~$\lin{P}$ such
  that~$\cl{p}=0$ holds in~$\cl{P}$.
\end{lemma}
\begin{proof}   
  Suppose that~$q-p\in I(P)$, that is, 
  \[
q-p = \sum_{1\leq i\leq k} \lambda_i u_i\bnd(\alpha_i) v_i.
  \]
  Then the following $1$-cell~$\phi$ of $\lin{P}$ has source~$p$ and target~$q$:
  \[
    \phi = \sum_{i=1}^k \lambda_i u_i \alpha_i v_i  +  \big( p - \sum_{i=1}^k \lambda_i u_i\src{}(\alpha_i) v_i \big).
  \]

  Conversely, let $\phi:p\fl q$ be a $1$-cell of~$\lin{P}$. Using
  \cref{L:NCellDecomposition}, we decompose~$\phi$ into $1$-cells of size~$1$:
  \[
    \phi = \phi_1\comp0\cdots\comp0 \phi_k
    \qquad\text{with}\qquad
    \phi_i = \lambda_i u_i \alpha_i v_i + h_i.
  \]
  Since $\tgt{}(\phi_i)=\src{}(\phi_{i+1})$, we have
  $q-p = \bnd(\phi_1) + \cdots + \bnd(\phi_p)$. Moreover, since
  $\bnd(\phi_i) = \lambda_i u_i \bnd(\alpha_i) v_i$ we have that
  each~$\bnd(\phi_i)$ belongs to~$I(P)$, and thus so does~$q-p$.

  Finally, if one applies the equivalence to the case~$p=0$, since~$\cl{0}=0$
  holds in~$\cl{P}$, we get that~$q$ is in~$I(P)$ if and only if we
  have~$\cl{q}=0$ in~$\cl{P}$.
\end{proof}

\subsection{Left-monomiality}
\index{monomial!1-polygraph}
\index{polygraph!monomial}
\label{sec:1-left-monomial}
A linear $1$-polygraph $P$ is \emph{left-monomial} if, for every
$1$-generator~$\alpha$ of~$P_1$, the source of~$\alpha$ is a monomial of~$\lin{P}_0$
that does not belong to~$\supp{\tgt{}(\alpha)}$. Note that, from any linear
$1$-polygraph~$P$, one obtains a Tietze equivalent left-monomial linear
$1$-polygraph as follows. For every $1$\nbd-gene\-rator~$\alpha$ in~$P_1$, if the boundary
$\bnd(\alpha)$ is~$0$, discard~$\alpha$, otherwise, replace~$\alpha$ with
\[
  \alpha' : u \to u - \frac{1}{\lambda} \bnd(\alpha),
\]
where~$u$ is any chosen monomial in $\supp{\bnd(\alpha)}$ and~$\lambda$ is the
coefficient of~$u$ in~$\bnd(\alpha)$.

\section{Rewriting Properties of Linear Polygraphs}
\label{sec:linear-rewritig-properties}
In the linear setting, the definition of a rewriting step is more difficult than
in the set-theoretic case, which can be explained as follows. In the
set-theoretic case developed in previous chapters, a $1$-polygraph~$P$ generates
two different objects: a free $1$-category~$\freecat P$ and a free
$1$-groupoid~$\freegpd P$. In this situation, we define a rewriting step as a
size-one $1$-cell of~$\freecat P$, and their compositions generate all the
$1$-cells of~$\freecat P$. But, in the case of associative algebras, there is no
difference between the free $1$-category and the free $1$-groupoid (see also
\cref{T:nAlg}), which is the cause of the problem mentioned in the introduction
of the present chapter. For this reason, we need to adopt a different point of view
to define rewriting steps and positive $1$-cells. Here, we identify, among the
$1$-cells of~$\lin{P}$, a set of positive $1$-cells that will play the same role
as the $1$-cells of~$\freecat P$ with respect to~$\freegpd P$ in the case of
set-theoretic rewriting. When defining this set, we need to ensure that two
conditions are satisfied. Firstly, the set of positive $1$-cells should be big
enough for every $1$-cell of~$\lin{P}$ to factor into a composite of positive
$1$-cells and opposites of positive $1$-cells, as given by
\cref{L:NCellDecomposition} and~\cref{L:FactElem1Cell}. Secondly, the set of
positive $1$-cells should be small enough for preventing a non-trivial $1$-cell
and its inverse to be positive at the same time, so that the polygraph has a
chance to be terminating.

In this section, $P$ denotes a left-monomial linear $1$-polygraph. 

\subsection{Rewriting steps and normal forms}
\index{rewriting!step}
\label{SS:RewritingSteps}
A \emph{rewriting step of~$P$} is a $1$-cell $\lambda\phi + 1_p$ of size~$1$ of
the free $1$-algebra~$\lin{P}$ that satisfies the condition
\[
  \supp{\lambda\src{}(\phi) + p} = \set{\src{}(\phi)}\sqcup\supp{a},
\]
that is, such that $\lambda\neq 0$ and $\src{}(\phi)\notin\supp{p}$. A $1$-cell
of the free $1$-algebra~$\lin{P}$ is called \emph{positive} if it is a (possibly
empty) $0$-composite $\phi_1\comp0\cdots\comp0 \phi_k$ of rewriting steps
of~$P$.

\begin{lemma}
  \label{L:FactElem1Cell}
  Let~$P$ be a left-monomial linear 1-polygraph. Every 1-cell~$\phi$ of
  size~1 of~$\lin{P}$ can be decomposed into $\phi=\psi\comp0\chi^-$, where each
  of~$\psi$ and~$\chi$ is either an identity or a rewriting step of~$P$.
\end{lemma}
\begin{proof}
  Write $\phi=\lambda \phi' + 1_q$, where $\phi':u\to p$ is a $1$-monomial
  of~$\lin{P}$. Let~$\mu$ be the coefficient of~$u$ in~$q$, possibly zero, so
  that $q=\mu u + r$ with~$r$ such that $\supp{r}$ does not contain~$u$. Put
  \[
    \psi = (\lambda + \mu) \phi' + 1_r
    \qquad\text{and}\qquad
    \chi = \lambda 1_p + \mu \phi' + 1_r.
  \]
  The linearity of the $0$-composition of~$\lin{P}$ gives
  $\phi=\psi\comp0\chi^-$.  Moreover, by hypothesis, $u$ does not belong to any
  of $\supp{p}$ or $\supp{r}$. As a consequence, each of the $1$-cells~$\psi$
  and~$\chi$ is either an identity (if $\lambda+\mu=0$ for~$\psi$, if $\mu=0$
  for~$\chi$) or a rewriting step.
\end{proof}

\subsection{Reduced cells and normal forms}
A $0$-cell~$p$ of~$\lin{P}_0$ is called \emph{reduced} if there is no rewriting
step of $P$ of source~$p$.
The reduced $0$-cells of~$\lin{P}_0$ form a linear subspace of the free
algebra~$\lin{P}_0$ which we denote by $\Red(P)$. Because~$P$ is left-monomial,
the set of reduced monomials of~$\lin{P}_0$, denoted by $\Red_m(P)$, forms a
basis of $\Red(P)$.

\index{normal form}
If~$p$ is a $0$-cell of~$\lin{P}_0$, a \emph{normal form of~$p$} is a reduced
$0$-cell~$q$ of~$\lin{P}_0$ such that there exists a positive $1$-cell of
source~$p$ and target~$q$ in the free $1$-algebra $\lin{P}$.

\subsection{Binary relations on free algebras}
Assume that~$\vdash$ is a binary relation on the free monoid~$P_0^\ast$
generated by the set $P_0$. We say that~$\vdash$ is \emph{stable by context} if
$u\vdash u'$ implies $vuw\vdash vu'w$ for all~$u$, $u'$, $v$, and~$w$
in~$P_0^\ast$. We say that $\vdash$ is \emph{compatible with~$P_1$} if
$u\vdash v$ holds for every $1$-cell $\alpha:u\to p$ in $P_1$ and every monomial
$v$ in $\supp{p}$.

The relation~$\vdash$ is extended to the $0$-cells of the free
algebra~$\lin{P}_0$ by setting~$p\vdash q$ when the following two conditions
hold:
\begin{enumerate}
\item $\supp{p}\setminus\supp{q}\neq\emptyset$,
\item for every~$v$ in $\supp{q}\setminus\supp{p}$, there exists~$u$ in
  $\supp{p}\setminus\supp{q}$, such that $u\vdash v$.
\end{enumerate}
As a consequence, if~$u$ is a monomial and~$p$ is a $0$-cell of~$\lin{P}_0$,
then $u\vdash p$ holds if and only if $u\vdash v$ holds for every~$v$
in~$\supp{p}$. Hence, we use the same notation for the relation on~$P_0^\ast$
and for its extension to the $0$-cells of~$\lin{P}_0$.

The relation~$\vdash$ on the $0$-cells of~$\lin{P}_0$ corresponds to the
restriction to finite subsets of~$P_0^\ast$ of the so-called \emph{multiset
  relation} generated by~$\vdash$. We refer to~\cite[Section
2.5]{BaaderNipkow98} for the general definition and the main properties of
multiset relations, and, in particular, the fact that~$\vdash$ is well-founded on
the $0$-cells if and only if it is well-founded on the monomials,
see also \cref{sec:multisets}.

\subsection{The termination order}
\index{order!rewrite}
\index{rewriting!order}
\label{SSS:RewriteOrderTermination}
Define~$\succ_P$ as the smallest transitive binary relation on~$\freecat P_0$
that is stable by context and compatible with~$P_1$. We say that the
polygraph~$P$ \emph{terminates} if the relation~$\succ_P$ is well-founded. In
that case, the reflexive closure~$\succcurlyeq_P$ of the relation~$\succ_P$ is a
well-founded order, called the \emph{termination order of~$P$} (this relation is also
sometimes written $\overset*\to$). This notion of termination order on linear polygraphs corresponds to that defined for 2-polygraphs in \cref{sec:2-red-order}.

Assume that the polygraph~$P$ terminates. Then the minimal $0$-cells for the
termination order of~$P$ are the reduced ones. Moreover, for every non-identity
positive $1$-cell~$p$ of~$\lin{P}_1$, we have~$\src{}(p)\succ_P \tgt{}(p)$. This
implies that the $1$-algebra~$\lin{P}$ contains no infinite sequence of
$0$-composable rewriting steps
\[
  \xymatrix{
    p_0\ar[r]^{\phi_1}&p_1\ar[r]^{\phi_2}&\cdots\ar[r]^{\phi_{n-1}}&a_{n-1}\ar[r]^{\phi_n}&a_n\ar[r]^{\phi_{n+1}}&\cdots
  }
\]
As a consequence, every $0$-cell of~$\lin{P}_0$ admits at least one normal form.
If~$P$ terminates, induction on the well-founded order~$\succ_P$ is called
\emph{noetherian induction}.

\subsection{Monomial orders}
\index{order!monomial}
\index{monomial!order}

A well-founded total order~$\leq$ on the free monoid~$\freecat P_0$ such that the
relation $<$ is stable by context is called a \emph{monomial order}. A
classic example of a monomial order is given, for any well-founded total order
relation~$>$ on~$P_0$, by the \emph{deglex order generated by~$>$}, as already
introduced in \cref{sec:deglex}, which is defined by
\begin{enumerate}
\item $u>_{\text{deglex}}v$ for all monomials~$u$ and~$v$ of~$\freecat{P}_0$ such
  that~$u$ has greater length than~$v$, and
\item $uxv>_{\text{deglex}}uyw$ for all~$x>y$ of~$P_0$, and monomials~$u$, $v$,
  and~$w$ of~$\lin{P}_0$ such that~$v$ and~$w$ have the same length.
\end{enumerate}

\index{leading!term}
\index{leading!coefficient}
\index{leading!monomial}
Given a monomial order~$\preccurlyeq$ on~$\lin{P}_0$. If~$p$ is a non-zero
$1$-cell of~$\lin{P}_0$, the \emph{leading monomial of~$p$} is the maximum
element of~$\supp{p}$ with respect to~$\preccurlyeq$ (or~$0$ if~$\supp{p}$ is
empty), it is denoted by~$\lm_\preccurlyeq(p)$.  The \emph{leading coefficient
  of~$p$} is the coefficient~$\lc_\preccurlyeq(p)$ of~$\lm_\preccurlyeq(p)$
in~$p$, and the \emph{leading term of~$p$} is the element
$\lt_\preccurlyeq(p)=\lc_\preccurlyeq(p)\lm_\preccurlyeq(p)$ of~$\lin{P}_0$.
Observe that, for~$p$ and~$q$ in~$\lin{P}_0$, we have~$p \prec q$ if and only
if either $\lm_\preccurlyeq(p) \prec \lm_\preccurlyeq(q)$ or
($\lt_\preccurlyeq(p)=\lt_\preccurlyeq(q)$ and
$p-\lt_\preccurlyeq(p) \prec q-\lt_\preccurlyeq(q)$).

If there exists a monomial order~$\succ$ on~$\lin{P}_0$ that is compatible
with~$P_1$, then the polygraph~$P$ terminates: the order~$\succ$ is well-founded,
and~$p\succ_P q$ implies~$p\succ q$ for all $0$-cells~$p$ and~$q$. However, the
converse implication does not hold, as illustrated by the following example.

\subsection{Example}
\label{X:xyz=x3+y3+z3}
The following linear $1$-polygraph terminates:
\[
  P = \pres{x,y,z}{\gamma:xyz\to x^3 + y^3 + z^3}.
\]
Indeed, for every monomial~$u$ of~$\lin{P}_1$, denote by~$A(u)$ the number of
factors~$xyz$ that occur in~$u$ and by~$B(u)$ the number of~$y$ that~$u$ contains,
and we consider the function $C(u)=3A(u)+B(u)$. It is sufficient to check that
$C(uxyzv)$ is strictly greater than each of $C(ux^3v)$, $C(uy^3v)$, and
$C(uz^3v)$, for all monomials~$u$ and~$v$ of~$\lin{P}_1$, see
{\cite[Example~3.2.4]{GuiraudHoffbeckMalbos19}} for details. However, no
monomial order on~$\lin{P}_0$ is compatible with~$P_1$, because, for such an
order~$\succ$, one of the monomials~$x^3$, $y^3$, $z^3$ is always greater
than~$xyz$.

\begin{lemma}
\label{L:DecompositionAlgebra}
If~$P$ is a terminating left-monomial linear 1-polygraph, then, as a vector
space, $\lin{P}_0$ admits the decomposition
\[
\lin{P}_0 = \Red(P) + I(P).
\]
\end{lemma}
\begin{proof}
  Since the polygraph~$P$ terminates, every $0$-cell~$p$ of~$\lin{P}_0$ admits
  at least a normal form~$q$. Let us write $p = q + (p-q)$, and note that~$q$
  belongs to~$\Red(P)$, by hypothesis, and that~$p-q$ is in~$I(P)$, by
  \cref{L:Ideal}.
\end{proof}

\subsection{Branchings}
\index{branching}
\label{SS:Branchings}
A \emph{branching} of the polygraph~$P$ is a pair~$(\phi,\psi)$ of positive
$1$-cells of the free $1$-algebra~$\lin{P}$ with the same source, called the \emph{source
  of~$(\phi,\psi)$}. We do not distinguish the branchings~$(\phi,\psi)$
and~$(\psi,\phi)$. A branching~$(\phi,\psi)$ of~$P$ is called \emph{local} if
both~$\phi$ and~$\psi$ are rewriting steps of~$\lin{P}$. For a branching
$(\phi,\psi)$ of~$P$ of source~$p$, define the branching
\[
\lambda u(\phi,\psi)v + q = (\lambda u\phi v + q, \lambda u\psi v + q)
\]
of~$P$ of source $\lambda upv + q$, for all scalar~$\lambda$, monomials~$u$
and~$v$ and $0$-cell~$q$ of~$\lin{P}$. Note that if~$(\phi,\psi)$ is local
and~$\lambda\neq 0$, then $\lambda u(\phi,\psi) + q$ is also local.

\subsection{Classification of local branchings}
Consider a local branching
\[
  (\lambda u_1 \alpha u_2 + p, \mu v_1 \beta v_2 + q)
\]
of~$P$. We have two main possibilities, depending on whether
\[
  u_1 \src{}(\alpha) u_2=v_1\src{}(\beta) v_2
\]
holds or not. Moreover, in the case of equality, there are three different
situations, depending on the respective positions of~$\src{}(\alpha)$
and~$\src{}(\beta)$ in this common monomial. This analysis leads to a partition
of the local branchings of~$P$ into the following four families.
\begin{enumerate}
\item \emph{Trivial} branchings: $\lambda (\phi,\phi) + q$, for all $1$-monomial
  $\phi:u\to p$ of~$\lin{P}$, non-zero scalar~$\lambda$, and $0$-cell~$q$
  of~$\lin{P}$, with $u\notin\supp{q}$.
\item \emph{Additive} branchings: $(\lambda\phi+\mu v+r,\,\lambda u+\mu\psi+r)$,
  for all $1$-monomials $\phi:u\to p$ and $\psi:v\to q$ of~$\lin{P}$, non-zero
  scalars~$\lambda$ and~$\mu$, and $0$-cell~$r$ of~$\lin{P}$, with $u\neq v$ and
  $u,v\notin\supp{r}$.
\item \emph{Multiplicative} branchings: $\lambda(\phi v,u\psi)+r$, for all
  $1$-monomials $\phi:u\fl p$ and $\psi:v\fl q$ of~$\lin{P}$, non-zero
  scalar~$\lambda$, and $0$-cell~$r$ of~$\lin{P}$, with
  $u,v\notin\supp{r}$. 
\item \emph{Overlapping} branchings: $\lambda(\phi,\psi)+r$, for all
  $1$-monomials $\phi:u\fl p$ and $\psi:u\fl q$ of~$\lin{P}$ such
  that~$(\phi,\psi)$ is neither trivial nor multiplicative,
  every non-zero scalar~$\lambda$, and every $0$-cell~$r$ of~$\lin{P}$, with
  $u\notin\supp{r}$.
\end{enumerate}
The \emph{critical branchings of~$P$} are the overlapping branchings of~$P$ such
that~$\lambda=1$ and~$r=0$, and that cannot be factored
$(\phi,\psi)=u(\phi',\psi')v$ in a non-trivial way. Note that an overlapping
branching has a unique decomposition $\lambda u(\phi,\psi)v + r$,
with~$(\phi,\psi)$ critical.

\subsection{Confluence}
\index{confluence}
Assume that~$P$ is a left-monomial linear $1$-polygraph. A
branching~$(\phi,\psi)$ of~$P$ is called \emph{confluent} if there exist
positive $1$-cells~$\phi'$ and~$\psi'$ of~$\lin{P}_1$ as in
\[
  \svxym{
    &\ar[dl]_{\phi}p\ar[dr]^{\psi}&\\
    q_1\ar@{.>}[dr]_{\phi'}&&\ar@{.>}[dl]^{\psi'}q_2\\
    &r
  }
\]
If~$p$ is a $0$-cell of~$\lin{P}_0$, we say that~$P$ is \emph{confluent at~$p$}
(\resp \emph{locally confluent at~$p$}, \resp \emph{critically confluent}) if
every branching (\resp local branching, \resp critical branching) of~$P$ of
source~$p$ is confluent. We say that~$P$ is \emph{confluent} (\resp
\emph{locally confluent}, \resp \emph{critically confluent}) if it is so at
every $0$-cell of~$\lin{P}$. We say that $P$ is~\emph{convergent} when it is
both terminating and confluent.

When the polygraph $P$ is confluent, then every $0$-cell of~$\lin{P}_0$ admits
at most one normal form, and when it is convergent then every $0$-cell~$p$
of~$\lin{P}_0$ has a unique normal form, denoted by~$\rep{p}$, such
that~$\cl{p}=\cl{q}$ holds in~$\cl{P}$ if and only if $\rep{p}=\rep{q}$ holds
in~$\lin{P}_0$. As a consequence, if~$P$ is a convergent presentation of an
algebra~$A$, the assignment of each element~$p$ of~$A$ to the normal form of any
representative of~$p$ in~$\lin{P}$, written~$\rep{p}$ by extension, defines a
section~$A\fl\lin{P}$ of the canonical projection, where~$A$ is seen as a
$1$-algebra with identity $1$-cells only. Note that the section is linear, that
is $\rep{\lambda p + \mu q}=\lambda\rep{p}+\mu\rep{q}$, and it preserves the
unit, that is $\rep{1}=1$. However, in general the equality
$\rep{pq}= \rep{p}\rep{q}$ does not hold.

\begin{proposition}
  \label{P:CharacterisationConfluence}
  Let~$P$ be a terminating left-monomial linear 1-polygraph. The following
  assertions are equivalent:
  \begin{enumerate}
  \item\label{P:CharacterisationConfluence1} The polygraph $P$ is confluent.
  \item\label{P:CharacterisationConfluence2} Every 0-cell of~$I(P)$ admits~$0$ as a normal form.
  \item\label{P:CharacterisationConfluence3} As a vector space, $\lin{P}_0$ admits the direct decomposition
    $\lin{P}_0 = \Red(P) \oplus I(P)$.
  \end{enumerate}
\end{proposition}
\begin{proof}
  {\cref{P:CharacterisationConfluence1}}$\:\dfl\:${\cref{P:CharacterisationConfluence2}}.
  By \cref{L:Ideal}, if~$p$ is in~$I(P)$, then there
  exists a $1$-cell $\phi:p\fl 0$ in~$\lin{P}_1$. Since~$P$ is confluent, this
  implies that~$p$ and~$0$ have the same normal form, if any. And, since~$0$ is
  reduced, this implies that~$0$ is a normal form of~$p$.

  {\cref{P:CharacterisationConfluence2}}$\:\dfl\:${\cref{P:CharacterisationConfluence3}}.
  By \cref{L:DecompositionAlgebra}, it is
  sufficient to prove that $\Red(P)\cap I(P)$ is reduced to~$0$. On the one
  hand, if~$p$ is in~$\Red(P)$, then~$p$ is reduced and, thus, admits itself as
  only normal form. On the other hand, if~$p$ is in~$I(P)$, then~$p$ admits~$0$
  as a normal form by hypothesis.

  {\cref{P:CharacterisationConfluence3}}$\:\dfl\:${\cref{P:CharacterisationConfluence1}}.
  Consider a branching~$(\phi,\psi)$ of~$P$, with
  $\phi:p\fl q$ and $\psi:p\fl r$. Since~$P$ terminates, each of~$q$ and~$r$
  admits at least one normal form, say~$q'$ and~$r'$ respectively. Hence, there
  exist positive $1$-cells $\phi':q\fl q'$ and $\psi':r\fl r'$
  in~$\lin{P}$. Note that the difference~$q'-r'$ is also reduced. Moreover, the
  $1$-cell $(\phi\comp0\phi')^- \comp0 (\psi\comp0\psi')$ has~$q'$ as source
  and~$r'$ as target. This implies, by \cref{L:Ideal}, that~$q'-r'$ also belongs
  to~$I(P)$. The hypothesis gives~$q'-r'=0$ so that~$(\phi,\psi)$ is confluent.
\end{proof}

\begin{theorem}
  \label{T:StandardBasis}
  Let~$A$ be an algebra and~$P$ a convergent presentation of~$A$. Then the
  set~$\Red_m(P)$ of reduced monomials of~$\lin{P}$ is a linear basis of~$A$. As
  a consequence, the vector space~$\Red(P)$, equipped with the product defined
  by $p\cdot q=\rep{pq}$, is an algebra that is isomorphic to~$A$.
\end{theorem}
\begin{proof}
  If~$P$ is convergent, \cref{P:CharacterisationConfluence} implies that the
  following sequence of vector spaces is exact:
  \[
    \xymatrix @C=1.5em {
      0
      \ar[r] 
      & I(P)
      \ar@{>->}[r] 
      & \lin{P}_0
      \ar@{->>}[r]
      & \Red(P)
      \ar[r] 
      & 0.
    }
  \]
  Thus, since the algebra~$\lin{P}_0/I(P)$ is isomorphic to~$\cl{P}$,
  convergence implies that the set~$\Red_m(P)$ is a linear basis of~$\cl{P}$. We
  deduce that~$\Red(P)$ and~$\cl{P}$ are isomorphic as vector spaces. There
  remains to transport the product of~$\cl{P}$ to~$\Red(P)$ to get the result.
\end{proof}

\subsection{Proving confluence}
The techniques developed in previous chapters for proving confluence in
practical cases can be adapted to the setting of linear polygraphs.
Namely, by a direct adaptation of the proof in the set-theoretic case, see
\cref{L:HNewman-1}, one can show an analogous of Newman's lemma: a terminating
left-monomial linear $1$-polygraph~$P$ which is locally confluent is confluent.
The critical branching lemma, see \cref{lem:2-cb}, also generalizes to our
setting. However, compared to the set-theoretic case, one has to add an extra
termination assumption in order to accommodate with the linearity of
contexts. The reason is explained in \cref{rem:lin-cb-term} below, and we defer
the proof to next chapter, where it will be proved in the more general setting of
coherent presentations, see \cref{L:HCriticalLinear}.

\begin{lemma}
  \label{lem:LinearCritical}
  Suppose given a terminating left-monomial linear 1-poly\-graph~$P$. If~$P$ is
  critically confluent, then~$P$ is locally confluent.
\end{lemma}

\noindent
A terminating left-monomial linear $1$-polygraph~$P$ in which all critical
branchings are confluent is thus necessarily confluent.

\begin{example}
  \label{X:BasisReduced}
  Let~$A$ be the algebra presented by the linear $1$-polygraph
  \[
    P=\pres{x,y}{\alpha: xy \fl x^2}.
  \]
  This polygraph terminates, because $xy>x^2$ holds for the deglex order
  generated by~$y>x$. This presentation is also confluent, because it has no
  critical branching, see \cref{lem:LinearCritical}. Hence, the set
  \[
    \Red_m(P) = \setof{y^i x^j}{i,j\in\mathbb{N}}
  \]
  is a linear basis of the algebra~$A$. Moreover, the product defined by
  \[
    y^i x^j \cdot y^k x^l
    =
    \begin{cases}
      y^i x^{j + k + l}
      &\text{if $j\geq k$,}
      \\
      y^{i - j + k} x^{2j + l}
      &\text{if $j\leq k$,}
    \end{cases}
  \]
  turns~$\Red(P)$ into an algebra that is isomorphic to~$A$.

  Now, consider the presentation~$Q=\pres{x,y}{\beta : x^2 \to xy}$
  of~$A$. Termination of~$Q$ follows from the deglex order generated by~$x>y$,
  but~$Q$ is not confluent, since it has a non-confluent critical branching:
  \[
    \svxym{
      &\ar[dl]_{\beta x}x^3\ar[dr]^{x\beta}&\\
      xyx&&x^2y\ar[d]^{\beta y}\\
      &&xy^2
    }
  \]
  Thus the $0$-cell $xyx-xy^2$ is both in~$\Red(Q)$ and~$I(Q)$, proving that the
  sum $\Red(Q)+I(Q)$ is not direct. As a consequence, $\Red_m(Q)$ is not a
  linear basis of~$A$.
\end{example}

\begin{example}
  \index{Weyl algebra}
  \index{algebra!Weyl}
\label{Example:WeylAlgebras2}
The polygraph of \cref{Example:WeylAlgebras1} that presents the Weyl algebra of dimension $n$ is convergent with the following six families of confluent critical branchings:
\begin{align*}
\xymatrix @C=0.5em @R=0.4em{
& x_jx_ix_k
    \ar [r]
& x_jx_kx_i
         \ar@/^2ex/ [dr]
\\
x_ix_jx_k
    \ar@/^2ex/ [ur]
    \ar@/_2ex/ [dr]
&&& x_kx_jx_i
\\
& x_ix_kx_j
    \ar [r]
& x_kx_ix_j
     \ar @/_2ex/ [ur]
}
&&
\xymatrix @C=0.5em @R=0.4em{
& \partial_j\partial_i\partial_k
    \ar [r]
& \partial_j\partial_k\partial_i
         \ar@/^2ex/ [dr]
\\
\partial_i\partial_j\partial_k
    \ar@/^2ex/ [ur]
    \ar@/_2ex/ [dr]
&&& \partial_k\partial_j\partial_i
\\
& \partial_i\partial_k\partial_j
    \ar [r]
& \partial_k\partial_i\partial_j
     \ar @/_2ex/ [ur]
}  
\end{align*}
\begin{align*}
\xymatrix @C=0.5em @R=0.4em{
& x_j\partial_ix_k
    \ar [r]
& x_jx_k\partial_i
         \ar@/^2ex/ [dr]
\\
\partial_ix_jx_k
    \ar@/^2ex/ [ur]
    \ar@/_2ex/ [dr]
&&& x_kx_j\partial_i
\\
& \partial_ix_kx_j
    \ar [r]
& x_k\partial_ix_j
     \ar @/_2ex/ [ur]
}
&&
\xymatrix @C=0.5em @R=0.4em{
& \partial_j\partial_ix_k
    \ar [r]
& \partial_jx_k\partial_i
         \ar@/^2ex/ [dr]
\\
\partial_i\partial_jx_k
    \ar@/^2ex/ [ur]
    \ar@/_2ex/ [dr]
&&& x_k\partial_j\partial_i
\\
& \partial_ix_k\partial_j
    \ar [r]
& x_k\partial_i\partial_j
     \ar @/_2ex/ [ur]
}  
\end{align*}
\[
\xymatrix @C=0.8em @R=0.4em{
& x_i\partial_ix_j + x_j
    \ar [r]
& x_ix_j\partial_i + x_j
         \ar@/^2ex/ [dr]
\\
\partial_ix_ix_j
    \ar@/^2ex/ [ur]
    \ar@/_2ex/ [dr]
&&& x_jx_i\partial_i + x_j
\\
& \partial_ix_jx_i
    \ar [r]
& x_j\partial_ix_i
     \ar @/_2ex/ [ur]
}
\]
\[
\xymatrix @C=0.8em @R=0.4em{
& \partial_j\partial_ix_j
    \ar [r]
& \partial_jx_j\partial_i
         \ar@/^2ex/ [dr]
\\
\partial_i\partial_jx_j
    \ar@/^2ex/ [ur]
    \ar@/_2ex/ [dr]
&&& x_j\partial_j\partial_i + \partial_i
\\
& \partial_ix_j\partial_j + \partial_i
    \ar [r]
& x_j\partial_i\partial_j + \partial_i
     \ar @/_2ex/ [ur]
}
\]
where $1\leq i < j  <k\leq n$. 
\end{example}

\begin{remark}
  \label{rem:lin-cb-term}
  The critical branching lemma for linear $1$-polygraphs given in
  \cref{lem:LinearCritical} differs from its set-theoretic counterpart because
  it requires the polygraph to be terminating, as noted
  in~\cite[Section~4.2]{GuiraudHoffbeckMalbos19}. Indeed, in the set-theoretic case,
  the termination hypothesis is not required, and non-overlapping branchings are
  always confluent, independently of critical confluence. The following two
  counterexamples show that the linear case is different.
The termination assumption comes from the fact that the rewriting steps are modulo the vector space structure.  We refer the reader to \cite{chenavier_dupont_malbos_2021} for an explanation of the linear critical pair lemma in terms of modulo rewriting.
\end{remark}

\begin{example}
  On the one hand, some local branchings can be non-confluent without
  termination, even if critical confluence holds. Indeed, the linear
  $1$\nbd-poly\-graph
  \[
    \pres {x, y, z, t}{\alpha:xy\to xz, \beta:zt\to 2yt}
  \]
  has no critical branching, but it has a non-confluent additive branching:
  \[
    \xymatrix @R=1em {
      && 4xyt 
      \ar [r] ^-*+{4\alpha t}
      & 4xzt
      \ar [r] ^-*+{4x\beta}
      & \cdots
      \\
      & 2xzt
      \ar@/^/ [ur] ^-*+{2x\beta}
      \ar@{.>} [dr] ^(0.55)*+{xzt + x\beta}
      \\
      xyt + xzt 
      \ar@/^/ [ur] ^-*+{\alpha t + xzt}
      \ar@/_/ [dr] _-*+{xyt + x\beta}
      \ar@{} [rr] |-{\sm =}
      && xzt + 2xyt
      \\
      & 3xyt
      \ar@{.>} [ur] _(0.55)*+{\alpha t + 2xyt}
      \ar@/_/ [dr] _-*+{3\alpha t}
      \\
      && 3xzt
      \ar [r] _-*+{3x\beta}
      & 6xyt
      \ar [r] _-*+{6\alpha t}
      & \cdots
    }
  \]
  The only positive $1$-cells of source~$2xzt$ are alternating $0$-compositions
  of $2^k x\beta$ and $2^{k+1}\alpha t$, whose targets are all the $0$-cells
  $2^k xzt$ and $2^{k+1}xyt$, for~$k\geq 1$. Similarly, the only positive
  $1$-cells of source~$3xyt$ have the $0$-cells $3.2^k xyt$ and $3.2^k xzt$ as
  targets, for~$k\geq 0$. The other possible $1$-cells of source~$2xzt$
  and~$3xyt$ are not positive, like the dotted ones. Here, it is the termination
  hypothesis that fails, as testified by the infinite sequences of rewriting
  steps in the previous diagram.
  \end{example}

\begin{example}
  On the other hand, the lack of critical confluence may imply that some
  non-overlapping local branchings are not confluent, even under the hypothesis
  of termination. For example, the linear $1$-polygraph
  \[
    \pres{x, y, z}{\alpha:xy\to 2x, \beta:yz\to z}
  \]
  terminates, but it has a non-confluent orthogonal branching:
  \[
    \xymatrix @R=1em {
      & 6xz
      & 3xz
      \\
      & 3xyz
      \ar [u] ^-*+{3\alpha z}
      \ar [ur] ^-*+{3x\beta}
      \ar@{.>} [dr] ^(0.55)*+{2x\beta + xyz}
      \\
      xyyz + xyz
      \ar@/^/ [ur] ^(0.4)*+{\alpha yz + xyz}
      \ar@/_/ [dr] _(0.4)*+{xy\beta + xyz}
      \ar@{} [rr] |-{\sm =}
      && 2xz + xyz
      \\
      & 2xyz
      \ar@{.>} [ur] _(0.55)*+{\alpha z + xyz}
      \ar [dr] _-{2beta z}
      \ar [d] _-*+{2x\alpha}
      \\
      & 4xz
      & 2xz
    }
  \]
  Here, it is the hypothesis on confluence of critical branchings that is not
  satisfied, since the critical branching~$(\alpha z, x\beta)$ of source~$xyz$
  is not confluent. As a consequence, the only $1$-cells that would close the
  confluence diagram of the Peiffer branching are the dotted ones, which are not
  positive.
\end{example}

\subsection{Reduced convergent presentations}
\index{presentation!reduced}
\index{reduced!presentation}
As with 2-polygraphs in~\cref{sec:2-reduced}, without loss of generality, we can restrict to the class of reduced linear polygraphs.
We say that a left-monomial linear polygraph~$P$ is \emph{left-reduced} if, for
every $1$-cell~$\alpha$ of~$P$, the only rewriting step of~$P$ of
source~$\src{}(\alpha)$ is~$\alpha$ itself. We say that~$P$ is
\emph{right-reduced} if, for every $1$-cell~$\alpha$ of~$P$, the
$0$-cell~$\tgt{}(\alpha)$ is reduced. We say that~$P$ is \emph{reduced} if it is
both left-reduced and right-reduced.

Using the same proof as for 2-polygraphs, \cref{thm:reduced-pres}, we prove that every convergent left-monomial linear $1$-polygraph is Tietze equivalent to a reduced convergent one.

\subsection{Completion of presentations}
The completion procedure, developed by Buchberger for commutative
algebras~\cite{buchberger1965algorithmus} and by Knuth and Bendix for term
rewriting systems~\cite{knuth1970simple}, see \cref{sec:2kb}, adapts to
terminating left-monomial linear $1$-polygraphs as follows, to transform them
into convergent ones.

Fix a left-monomial linear $1$-polygraph~$P$, and a well-founded strict order
that is stable by context and compatible with~$P_1$. For each non-confluent
critical branching~$(\phi,\psi)$ of~$P$, consider $p=r-s$, where~$r$ and~$s$ are
arbitrary normal forms of~$\tgt{}(\phi)$ and~$\tgt{}(\psi)$, respectively. If
$\supp{p}$ contains a maximal element~$u$, add the $1$-cell $u\fl q$ to~$P$,
where~$q$ is defined by $p=\lambda u + q$ and $u\notin\supp{q}$; otherwise, the
procedure fails. After the exploration of all the critical branchings of~$P$,
the procedure, if it has not failed, yields a terminating left-monomial linear
$1$-polygraph~$Q$ such that~$\cl{P}\simeq\cl{Q}$. If~$Q$ is not confluent,
restart with~$Q$. The procedure either stops when it reaches a convergent
left-monomial linear $1$-polygraph or runs forever.

\section{Linear Bases Induced by Monomial Orders}
\label{SS:PBWGrobner}

In this section, we consider linear rewriting systems whose rewriting rules are oriented with respect to a fixed monomial order.
Suppose fixed a monomial order on polynomials (a well-founded order suitably
compatible with multiplication) and an ideal~$I$. A polynomial
$p=\sum_i\lambda_iu_i$ in the ideal~$I$ can of course be interpreted as a
relation
\[
  \sum_i\lambda_iu_i=0
\]
However, supposing that $u_0$ is the monomial which is the greatest with respect
to the fixed order, called the \emph{leading monomial}, it can also be
interpreted as a relation
\[
  u_0=\frac 1{\lambda_0}\sum_{i\neq 0}u_i
\]
which, in turn, can be seen as a rewriting rule transforming the left member~$u_0$
into the right member. A Gröbner basis is then a generating set for the ideal
such that the associated rewriting rules form a confluent rewriting system, and
Buchberger's algorithm to compute a basis can be seen as a form of Knuth-Bendix
completion~\cite{Buchberger87}.

\subsection{Gröbner bases}
\index{Gröbner basis}
\index{basis!Gröbner}
Let~$P_0$ be a set, and let~$I$ be an ideal of the free algebra~$\lin{P}_0$. A
\emph{Gröbner basis} for~$I$ with respect to a monomial order $\preccurlyeq$ on
$\freecat P_0$ is a subset~$\mathcal{G}$ of~$I$ such that the ideals
of~$\lin{P}_0$ generated by~$\lm_\preccurlyeq(I)$ and
by~$\lm_\preccurlyeq(\mathcal{G})$ coincide.

\begin{proposition}
  \label{P:ConvergentGrobner}
  If~$P$ is a convergent left-monomial linear 1-polygraph, and~$\preccurlyeq$
  is a monomial order on~$P_0^\ast$ that is compatible with~$P_1$, then the set
  \[
    \bnd(P_1)=\setof{\bnd(\alpha)}{\alpha\in P_1}
  \] 
  of boundaries of $1$-generators of~$P$ is a Gröbner basis
  for~$(I(P),\preccurlyeq)$.

  Conversely, let~$P_0$ be a set, let~$\preccurlyeq$ be a monomial order
  on~$\lin{P}_0$, and let~$I$ be an ideal of~$\lin{P}_0$ and~$\mathcal{G}$ be a
  subset of~$I$. Define $P(\mathcal{G})$ the linear 1-polygraph whose set of
  0-cells is~$P_0$ and with one 1-cell
  \[
    \alpha_p : \lm(p) \to \lm(p) - \frac{1}{\lc(p)} p,
  \]
  for each~$p$ in~$\mathcal{G}$. If~$\mathcal{G}$ is a Gröbner basis
  for~$(I,\preccurlyeq)$, then~$P(\mathcal{G})$ is a convergent left-monomial
  presentation of the algebra~$\lin{P}/I$, such that $I(P(\mathcal{G}))=I$,
  and~$\preccurlyeq$ is compatible with~$P(\mathcal{G})_1$.
\end{proposition}
\begin{proof}
  If the polygraph~$P$ is convergent, then~$\bnd(\alpha)$ is in~$I(P)$ for every
  $1$-cell~$\alpha$ of~$P$. Since~$\preccurlyeq$ is compatible with~$P_1$, we
  have $\lm(\bnd(\alpha))=\src{}(\alpha)$ for every $1$-cell~$\alpha$
  of~$P$. Now, if~$p$ is in~$I(P)$, it is a linear combination
  \[
    p = \sum_{i} \lambda_i u_i \bnd(\alpha_i) v_i
  \]
  of $1$-cells $u_i\bnd(\alpha_i) v_i$, where~$\alpha_i$ is a $1$-cell
  of~$P$, and~$u_i$ and~$v_i$ are monomials of~$\lin{P}_0$. This implies that
  \[
    \lm(p) = u_i\src{}(\alpha_i)v_i = u_i \lm(\bnd(\alpha_i)) v_i
  \]
  hold for some~$i$. Thus~$\bnd(P_1)$ is a Gröbner basis
  for~$(I(P),\preccurlyeq)$.

  Conversely, assume that~$\mathcal{G}$ is a Gröbner basis
  for~$(I,\preccurlyeq)$. By definition, $\preccurlyeq$ is compatible
  with~$P(\mathcal{G})_1$, so~$P(\mathcal{G})$ terminates, and
  $I(P(\mathcal{G}))=I$ holds, so that the algebra presented by~$P(\mathcal{G})$
  is indeed isomorphic to $\lin{P}/I$. Moreover, the reduced monomials
  of~$\lin{P(\mathcal{G})}$ are the monomials of~$\lin{P}$ that cannot be
  decomposed as~$u\lm(a)v$ with~$a$ in~$\mathcal{G}$, and~$u$ and~$v$ monomials
  of~$\lin{P}$. Thus, if a reduced $0$-cell~$p$ of~$\lin{P(\mathcal{G})}$ is
  in~$I$, its leading monomial must be~$0$, because~$\mathcal{G}$ is a Gröbner
  basis of~$(I,\preccurlyeq)$. As a consequence of
  \cref{P:CharacterisationConfluence}, we get that~$P(\mathcal{G})$ is
  confluent.
\end{proof}

\noindent
By previous proposition and the critical branching lemma
(\cref{lem:LinearCritical}), the notion of Gröbner basis can be related to
confluence as follows. This is sometimes called \emph{Buchberger's criterion}
for determining whether a set of polynomials forms a Gröbner basis with respect to a fixed monomial order.

\begin{proposition}
  \label{C:BuchbergerCriterion}
  Let~$P_0$ be a set, $\preccurlyeq$ be a monomial order on~$P_0^\ast$, and~$I$
  be an ideal of the free algebra~$\lin{P}_0$. A subset~$\mathcal{G}$ of~$I$ is
  a Gröbner basis for~$(I,\preccurlyeq)$ if and only if the linear 1-polygraph
  $\lm(\mathcal{G})$ of \cref{P:ConvergentGrobner} is critically confluent.
\end{proposition}

\subsection{Polygraphs for graded associative algebras}
\nomenclature[gVect]{$\gVect{}$}{category of graded vector spaces}
\index{algebra!graded}
\index{graded algebra}
\label{SS:GradedLinearPolygraphs}

Let denote by~$\gVect{}$ the category of (non-negatively) graded vector spaces over~$\kk$ and graded linear maps of degree~$0$. Recall that a graded vector space~$V$ admits a decomposition $V = \bigoplus_{i\in\Nb} V^{(i)}$, and the elements of~$V^{(i)}$ are said to be \emph{homogeneous of degree~$i$}.
A \emph{graded associative algebra} is an internal monoid in the category $\gVect{}$. Following {\cite[Section 2.2]{GuiraudHoffbeckMalbos19}} we can define a notion of polygraph, called \emph{graded linear polygraphs}, for presentation of graded algebras. Let us expand this notion in low dimensions.

A \emph{graded linear 1-polygraph} is a data $(P_0,P_1)$ made of  
\begin{itemize}
\item a graded linear $0$-polygraph $P_0$, that is, a graded set $P_0 = \coprod_{i\in\Nb} P_0^{(i)}$,
\item a graded cellular extension $P_1$ of the free graded algebra $\lin{P}_0$ generated by $P_0$, meaning that $P_1 = \coprod_{i\in\Nb} P_1^{(i)}$ and that the source and target of each $1$-generator in~$P_1^{(i)}$ are homogeneous of degree~$i$.
\end{itemize}

If~$N\geq 2$, a $1$-polygraph~$P$ is called \emph{$N$-homogeneous} if~$P_0$ is concentrated in degree~$1$ and~$P_1$ is concentrated in degree~$N$. We say \emph{quadratic} and \emph{cubical} instead of $2$-homogeneous and $3$-homogeneous, respectively.

An algebra~$A$ is called \emph{$N$-homogeneous} if it admits a presentation by an
$N$-homogeneous graded linear $1$-polygraph.

\subsection{Poincaré-Birkhoff-Witt bases}
\index{homogeneous!algebra}
\index{PBW basis}
\index{basis!PBW}

Let~$A$ be an $N$-homogeneous algebra, for~$N\geq 2$, let~$P_0$ be a generating
set of~$A$, concentrated in degree~$1$, and let~$\preccurlyeq$ be a monomial
order on~$P^\ast_0$. A \emph{Poincaré-Birkhoff-Witt} (\emph{PBW}) \emph{basis
  for~$(A,P_0,\preccurlyeq)$} is a subset~$\mathcal{B}$ of~$P_0^\ast$ 
  satisfying the following conditions:
\begin{enumerate}
\item $\mathcal{B}$ is a linear basis of~$A$, with~$[u]_{\mathcal{B}}$ denoting
  the decomposition of an element~$u$ of~$P_0^\ast$ in the basis~$\mathcal{B}$.
\item For all~$u$ and~$v$ in~$\mathcal{B}$, we have
  $uv\succcurlyeq [uv]_{\mathcal{B}}$.
\item An element~$u$ of~$P_0^\ast$ belongs to~$\mathcal{B}$ if and only if for
  every decomposition $u=vu'w$ of~$u$ in~$P_0^\ast$ such that~$u'$ has
  degree~$N$, then~$u'$ is in~$\mathcal{B}$.
\end{enumerate}

\begin{proposition}
  \label{P:ConvergentPBW}
  If~$P$ is a convergent left-monomial $N$-homogeneous presentation of an
  algebra~$A$, and~$\preccurlyeq$ is a monomial order on~$\lin{P}_0$ that is
  compatible with~$P_1$, then the set~$\Red_m(P)$ of reduced monomials
  of~$\lin{P}_0$ is a PBW basis for~$(A,P_0,\preccurlyeq)$.

  Conversely, let~$A$ be an $N$-homogeneous algebra, let~$P_0$ be a generating
  set of~$A$ that is concentrated in degree 1, let~$\preccurlyeq$ a monomial
  order on~$\lin{P}_0$, and let~$\mathcal{B}$ be a PBW basis
  of~$(A,P_0,\preccurlyeq)$. Define~$P(\mathcal{B})$ as the linear 1-polygraph
  with $0$-cells~$P_0$ and with one $1$-cell
  \[
    \alpha_{u,v} : uv \to [uv]_{\mathcal{B}}
  \]
  for all~$u$ and~$v$ in~$\mathcal{B}$ such that~$uv$ has degree~$N$
  and~$uv\neq [uv]_{\mathcal{B}}$. Then~$P(\mathcal{B})$ is a convergent
  left-monomial $N$-homogeneous presentation of~$A$, such that
  $\Red_m(P(\mathcal{B}))=\mathcal{B}$, and~$\preccurlyeq$ is compatible
  with~$P(\mathcal{B})_1$.
\end{proposition}
\begin{proof}
  If~$P$ is a convergent left-monomial presentation of~$A$,
  Theorem~\ref{T:StandardBasis} implies that the set~$\Red_m(P)$ of reduced
  monomials of~$P_0^\ast$ is a linear basis of~$A$. The fact that~$\preccurlyeq$
  is compatible with~$P_1$ implies Axiom~2 of a PBW basis, and
  Axiom~3 comes from the definition of a reduced monomial for an
  $N$-homogeneous left-monomial linear $1$-polygraph.

  Conversely, assume that~$\mathcal{B}$ is a PBW basis
  for~$(A,P_0,\preccurlyeq)$. By definition, $P(\mathcal{B})$ is $N$-homogeneous
  and left-monomial, and Axiom~2 of a PBW basis implies
  $\Red_m(P(\mathcal{B}))=\mathcal{B}$. Termination of~$P(\mathcal{B})$ is given
  by Axiom~2 of a PBW basis, because~$\preccurlyeq$ is well-founded. By
  \cref{P:CharacterisationConfluence}, it is sufficient to prove that
  $\Red(P(\mathcal{B}))\cap I(P(\mathcal{B}))=0$ to get confluence: on the one
  hand, a reduced $0$-cell~$a$ of~$\Red(P(\mathcal{B}))$ is a linear combination
  of $0$-cells of~$\mathcal{B}$ so that~$a$ is its only normal form; and, on
  the other hand, if~$a$ belongs to~$I(P(\mathcal{B}))$, then~$a$ admits~$0$ as
  a normal form by \cref{L:Ideal}. Finally, the algebra presented
  by~$P(\mathcal{B})$ is isomorphic to~$\Red(P(\mathcal{B}))$, that is,
  to~$\kk\mathcal{B}$, hence to~$A$, by \cref{T:StandardBasis} and
  because~$\mathcal{B}$ is a linear basis of~$A$.
\end{proof}

\section{Historical Account of Linear Rewriting}
\index{Buchberger algorithm}
\index{algorithm!Buchberger}
\label{S:HistoryLinear}

Gröbner basis theory for ideals in commutative polynomial rings was introduced by Buchberger in \cite{buchberger1965algorithmus}. 
He defined the notion of \emph{$S$-polynomial} to describe the obstructions to local confluence and gave an algorithm for computation of Gröbner bases, \cite{buchberger1965algorithmus, Buchberger70, Buchberger06}, see also \cite{Buchberger87} for an historical account.
In the commutative setting, any ideal of a polynomial ring has a finite Gröbner basis. Indeed, the Buchberger algorithm on a finite family of generators of an ideal always terminates and returns a Gröbner basis of the ideal.
More recently, refined efficient
algorithms have been proposed to achieve this task, such as Faugère's $F_4$ and
$F_5$ algorithms~\cite{faugere1999new,faugere2002new}.

Shirshov introduced in~\cite{shirshov1962algorithmic}  an algorithm to compute a linear basis of a Lie algebra defined by generators and relations. He used the notion of \emph{composition} of elements in a free Lie algebra, that corresponds to the notion of $S$-polynomial in the work of Buchberger. He gave an algorithm to compute bases in free algebras having the computational properties of the Gröbner bases.  He proved that irreducible elements for such a basis form a linear basis of the Lie algebra. This result is called now the \emph{Composition Lemma} for Lie algebras~\cite{bokut2014grobner}.

The Gröbner basis theory has been developed for other types of algebras, such as associative algebras by Bokut in \cite{Bokut76} and by Bergman in \cite{bergman1978diamond}. 
They prove Newman's Lemma for rewriting systems in free associative algebras compatible with a monomial order, stating that local confluence and confluence are equivalent properties.
This result was called \emph{Composition Lemma} by Bokut and \emph{Diamond Lemma}\index{diamond lemma}\index{lemma!diamond} for ring theory by Bergman, see also \cite{Mora94, Ufnarovski95}. 
In general, the Buchberger algorithm does not terminate for ideals in a non-commutative multivariate polynomial ring. Indeed, its termination would give a decision procedure of the undecidable word problem. Even if the ideal is finitely generated it may not have a finite Gröbner basis. However, an infinite Gröbner basis can be computed over a ground field, \cite{Mora94, Ufnarovski98}.
The Buchberger algorithm is the analogue of the Knuth-Bendix completion procedure in a linear setting. Several frameworks unify Buchberger and Knuth-Bendix algorithms; in particular, a Gröbner basis corresponds to a confluent and terminating presentation of an algebra, see \cite{Buchberger87}.

Finally, note that ideas in the style of Gröbner's basis approach appear in many independent works throughout the 20th century. Günter has defined a
similar notion in 1913~\cite{renschuch2003contributions}.
Janet~\cite{Janet20,Janet20a, Janet20b} and Thomas~\cite{Thomas37} developed the notion of involutive bases that are particular cases of Gröbner bases in the context of partial differential algebra. We refer to~\cite{IoharaMalbos20,IoharaMalbos21} for an historical account on involutive bases and their applications to algebraic analysis of linear partial differential systems. 
Hironaka in \cite{hironaka1964resolution} and Grauert in \cite{Grauert72} compute bases of ideals in rings of power series having analogous properties to Gröbner bases but without a constructive method for computing such bases. In \cite{Cohn65}, Cohn gave a method  to decide the word problem by a normal form algorithm based on a confluence property. Much more recently, Gröbner basis theory was developed in various non-commutative contexts such as Weyl algebras, see~\cite{SaitoSturmfelsTakayama00}, or operads~\cite{dotsenko2010grobner, MalbosRen23}.

\part{Coherent Presentations}

\chapter{Coherence by Convergence}
\label{chap:2coh}
\label{chap:2-Coherent}
To any presentation of a category $C$ by a $2$\nbd-poly\-graph~$P$
corresponds a free $(2,1)$\nbd-cate\-gory~$\freegpd P_2$, as defined
in~\cref{sec:free-21-cat}. An \emph{extended presentation} then
consists in a choice of a family $P_3$ of $3$-generators between some
pairs of parallel $2$-cells in~$\freegpd P_2$. Since the category $C$
is already entirely determined by the presentation
$\Pres{P_0}{P_1}{P_2}$, we are mostly interested in the case where
the congruence generated by $P_3$ is the {\em full} relation among pairs of
parallel $2$-cells in $\freegpd P_2$, that is,  when each $2$-sphere is
filled with a $3$-cell generated by $P_3$.
An extended presentation satisfying this property is said to be \emph{coherent}.

Any given presentation~$P$ of a category can be extended into a coherent one by
taking all parallel pairs of $2$-cells as  $3$-generators, but
we are mainly
interested in ``small'' coherent presentations, which are amenable to
computations. The key result in building small coherent
presentations is \cref{thm:SquierHomotopical}, a refined version of
Newman's lemma, called here~\emph{Squier's homotopical theorem}. It
states that a \emph{convergent} presentation $P$ can be extended to a
coherent one by taking for~$P_3$ a family of confluence diagrams of
critical branchings. As a consequence, if $P$ is \emph{finite}
convergent, then $P_3$ can be chosen finite.  

We then introduce a notion of Tietze transformation preserving the
coherence property and the presented category. This, combined with
Squier's homotopical theorem, suggests the following general procedure
to build a coherent extension of a given---not necessarily convergent---presentation $P$: 
\begin{enumerate}
\item Use the Knuth-Bendix completion procedure to compute a convergent
  $2$\nbd-poly\-graph~$Q$ presenting the same category as~$P$.
\item Use Squier's homotopical theorem to extend~$Q$ into a coherent
  presentation~$\tilde Q$.
\item Use Tietze transformations to reduce $\tilde Q$ into a smaller coherent
  presentation~$\tilde P$, which extends~$P$.
\end{enumerate}
The first two steps can be performed at once, using what we call a
\emph{coherent completion procedure}, and a \emph{transfer theorem}
(\cref{thm:HomotopyBasesTransfer}), providing an immediate description of the
coherent extension of~$P$ from the one of~$Q$.

Coherent presentations will prove essential for computing, in low
dimensions, homotopical and homological invariants of presented
categories introduced in~\cref{chap:2-fdt,chap:2-homology}. Moreover, we shall
see in~\cref{chap:ConstructingResolutions} how they extend in any
dimension to \emph{polygraphic resolutions} (\cref{chap:resolutions})
of the presented category. In the language of homotopy theory, these
resolutions are cofibrant replacements of a category by a free $(\omega,1)$-category.
Note also that the rewriting method for calculating coherent
presentations can be applied in many algebraic contexts, as illustrated in \cref{chap:SomeCoherentPresentations}.

This chapter is organized as follows.
In \cref{sec:acyclic-extension}, we introduce the notion of \emph{acyclic
  extension} of a $2$-category, which consists of the additional data of
$3$-generators ``filling all the spheres''. This leads in~\cref{sec:CoherentPresentation}  to  the notion of
\emph{coherent presentation} of a category $C$, that is, 
a $2$-polygraph $P$ presenting $C$ together with an acyclic extension of the free
$(2,1)$-category on $P$. Coherent presentations are then constructed
from convergent ones in \cref{sec:CoherentConfluence}.
The appropriate notion of Tietze transformation between coherent presentations
is studied in \cref{Section:TietzeTransformations31Polygraphes}: this allows us
in \cref{Section:HomotopicalCompletion}
to formulate a coherent variant of the Knuth-Bendix completion procedure, but also a reduction procedure, which can be used to
obtain smaller coherent presentations. Finally, in \cref{sec:coh-alg}, we study coherent
presentations of algebras, thereby defining the proper notion of coherent
extension for the linear polygraphs of \cref{chap:2linear}.

\section{Acyclic Extensions}
\label{Section:AcylicExtensions}
\label{sec:acyclic-extension}

\subsection{Cellular extension of a 2-category}
\index{2-sphere@$2$-sphere}
\index{sphere!2-}
\label{sec:2-spheres}
A \emph{2\nbd-sphere} in a 2-category~$C$ is a
pair~$(\alpha,\beta)$ of parallel $2$-cells in~$C$, \ie satisfying
$\sce1(\alpha)=\sce1(\beta)$ and $\tge1(\alpha)=\tge1(\beta)$.
\index{cellular extension}
A \emph{cellular extension of~$C$} is a set~$X$
equipped with two maps
$
  \sce{2},\tge{2}
  :
  X\to C_2
$
such that, for every~$A$ in $X$, the pair $(\sce{2}(A),\tge{2}(A))$ is a
$2$-sphere of~$C$. More generally, we also call any such element $A$ in $X$ a
\emph{2-sphere} of~$C$.

Every $2$-category~$C$ has two canonical cellular extensions:
\begin{itemize}
\item the empty extension,
\item the full one that contains all the $2$-spheres of~$C$, denoted
  by~$\Sph(C)$.\nomenclature[Sph]{$\Sph(C)$}{$2$-spheres in a $2$-category
  $C$}
\end{itemize}

\subsection{Quotient 2-category}
\index{congruence!on a 2-category}
\label{sec:2-cat-quot}
A \emph{congruence} on a 2-category~$C$ is an
equivalence relation~$\approx$ on the $2$-cells of~$C$ such that
\begin{itemize}
\item given $\varphi:u\dfl v$ and $\varphi':u'\dfl v'$ in $C_2$,
  $\varphi\approx \varphi'$ implies $u=u'$ and $v=v'$,
\item given $1$-cells and $2$-cells of $C$ as in the following diagram
  \[
    \xymatrix{
      x\ar[r]^{u}&x'
      \ar@/^8ex/[rr]
      \ar@{{}{ }{}}@/^6ex/[rr]|{\varphi_1\Downarrow\phantom{\varphi_1}}
      \ar@/^4ex/[rr]
      \ar@{}[rr]|{\psi\Longdownarrow\quad\Longdownarrow\psi'}
      \ar@/_4ex/[rr]
      \ar@{{}{ }{}}@/_6ex/[rr]|{\varphi_2\Downarrow\phantom{\varphi_2}}
      \ar@/_8ex/[rr]
      && y'\ar[r]^{v}&
      y
    }
  \]
  if $\psi\approx\psi'$, then 
  \[
    u\comp0(\varphi_1\comp1 \psi \comp1 \varphi_2)\comp0 v
    \approx
    u\comp0(\varphi_1\comp1 \psi' \comp1 \varphi_2)\comp0 v
    \pbox.
  \]
\end{itemize}

We define the \emph{quotient $2$-category}\index{quotient!2-category} of a
$2$-category~$C$ by a congruence~$\approx$ on $C$ as the $2$\nbd-cate\-gory,
denoted by $C/{\approx}$, whose $0$-cells and $1$-cells are those of~$C$, and
whose $2$-cells are the equivalence classes of $2$-cells of~$C$ modulo the
congruence~$\approx$, composition and identities being induced by those of~$C$.

Given a cellular extension~$X$ of~$C$, the \emph{congruence generated by $X$},
denoted by~$\approx^X$, is defined as the smallest congruence on $C$ such that
$\phi\approx^X\psi$, for every $2$-sphere~$(\phi,\psi)$ in~$X$. That is,
$\approx^X$ is the smallest equivalence relation on the parallel $2$-cells
compatible with all the compositions of~$C$ and relating $\phi$ and $\psi$, for
every $(\phi,\psi)$ in~$X$.

\subsection{Acyclic extension}
\index{acyclic!extension}
\index{extension!acyclic}
We say that a cellular extension~$X$ of a $2$-category~$C$ is
\emph{acyclic}, or equivalently that~$C$ is \emph{$X$-acyclic}, if $\phi\approx^X\psi$ holds for every
$2$\nbd-sphere~$(\phi,\psi)$ of~$C$.  This is equivalent to say that the
equality $\cl{\phi}=\cl{\psi}$ holds in the quotient $2$-category
$C/{\approx^X}$, where $\cl{\phi}$ and $\cl{\psi}$ denote the images of the
$2$-cells under the canonical projection $C\to C/{\approx}^X$.  For instance,
any $2$-category $C$ is $\Sph(C)$-acyclic.

\begin{remark}
  \index{homotopy!relation}
  \index{homotopy!basis}
  \label{rem:TerminologyHomotopyBasis}
  A congruence on a free $(2,1)$-category is called a \emph{homotopy
    relation} by Squier in~\cite{squier1994finiteness}. He noticed that these relations are not really
  the same as usual homotopies in the sense of algebraic topology and justified
  the terminology by saying that, for a homotopy relation generated by a set of
  $2$\nbd-spheres, two ``homotopic'' paths can be transformed into one another
  by a finite sequence of elementary transformation
  steps. In~\cite{squier1994finiteness}, the relation $\approx^X$ is called an
  homotopy relation \emph{generated} by $X$.  The terminology
  \emph{homotopy basis} for an acyclic cellular extension
  was introduced in \cite{kobayashi2000finite,gilbert1998monoid} and
  since then has
  been widely used by various authors. Note also that Squier did not formulate
  his results on the properties of homotopy relations in the
  categorical language 
  we use in the present chapter. Instead of $(2,1)$-categories, he considered
  $2$\nbd-dimen\-sional cellular complexes defined by directed graphs with inverses
  and whose $2$-cells correspond to the exchange relation between compositions
  with respect to $0$- and $1$-composition.
  The categorical formulation of Squier's constructions presented here
  was introduced in~\cite{guiraud2009higher,GuiraudMalbos18}. Another formulation, using
  the structure of monoidal category, is given in~\cite{lafont1995new}.
\end{remark}

\subsection{Transfer theorem for acyclic extensions}
\index{transfer theorem}
Given two presentations $P$ and $Q$ of a $1$-category~$C$, by
\cref{lem:21TietzeEquivalence}, there exist two $2$-functors
\[
f:\freegpd{P}\to\freegpd{Q}
\qqqtand
g:\freegpd{Q}\to\freegpd{P}
\]
and, for every $1$-cell $v$ in $\freegpd{Q}$, there exists a $2$-cell
$\psi_v : fg(v) \To v$ in $\freegpd{Q}$ that satisfies the conditions given in
\cref{lem:21TietzeEquivalence}.
Let us define a cellular extension $X_{Q}$ of the $(2,1)$-category $\freegpd{Q}$ that
contains one $3$-generator
\begin{equation}
  \label{eq:tt-Aa}
  \vcenter{
    \xymatrix @!C @C=2.5em @R=1em{
      & fg(v)
      \ar@2 @/^/ [dr] ^-{\psi_v}
      \ar@3 []!<0pt,-15pt>;[dd]!<0pt,15pt> ^-{A_{\alpha}}
      \\
      fg(u) 
      \ar@2 @/^/ [ur] ^-{fg(\alpha)}
      \ar@2 @/_/ [dr] _-{\psi_u}
      && v	
      \\
      & u
      \ar@2 @/_/ [ur] _-{\alpha}
    }
  }
\end{equation}
for every $2$-generator $\alpha:u\To v$ of $Q$. Furthermore, given a cellular
extension~$X$ of the $(2,1)$\nbd-cate\-gory $\freegpd{P}$, we will denote
by~$f(X)$ the cellular extension of~$\freegpd{Q}$ that contains one $3$-generator
\[
  \xymatrix @C=6em{
    f(u)
    \ar@2 @/^4ex/ [r] ^-{f(\varphi)} _-{}="src"
    \ar@2 @/_4ex/ [r] _-{f(\varphi')} ^-{}="tgt"
    \ar@3 "src"!<0pt,-10pt>;"tgt"!<0pt,10pt> ^-{f(A)}
    & f(v)
  }
\]
for every $3$-generator $A:\varphi\TO\varphi'$ of $X$. Using these notations, we can
formulate the following transfer result among acyclic
extensions of presentations of a given category.

\begin{theorem}
  \label{thm:HomotopyBasesTransfer}
  \label{Theorem:HomotopyBasesTransfer}
  Let $P$ and $Q$ be two presentations of the same category.
  If~$X$ is an acyclic cellular extension of the
  $(2,1)$\nbd-cate\-gory~$\freegpd{P}$, then the cellular extension
  $f(X)\sqcup X_{Q}$ is an acyclic cellular extension of the $(2,1)$-category
  $\freegpd{Q}$.
\end{theorem}
The proof consists in extending the notation on $3$-generators
$A_\alpha$ of \eqref{eq:tt-Aa}, where $\alpha$ is a $2$-generator of
$Q$, in a functorial way, to define a $3$-cell of the shape
\[
  \xymatrix @!C @C=2.5em @R=1em{
    & fg(v)
    \ar@2 @/^/ [dr] ^-{\psi_v}
    \ar@3 []!<0pt,-15pt>;[dd]!<0pt,15pt> ^-{A_\varphi}
    \\
    fg(u) 
    \ar@2 @/^/ [ur] ^-{fg(\varphi)}
    \ar@2 @/_/ [dr] _-{\psi_u}
    && v
    \\
    & u
    \ar@2 @/_/ [ur] _-{\varphi}
  }
\]
for any $2$-cell $\varphi$ in $\freegpd{Q}$. Then, given two parallel $2$-cells
$\varphi,\varphi':u\To v$ of $\freegpd{Q}$, one proves that
$\varphi \approx_{f(X)\sqcup X_{Q}} \varphi'$ by constructing a $3$-cell with
source $\varphi$ and target $\varphi'$ obtained by compositions along $0$-cells,
$1$-cells and $2$-cells of the $3$-cells $A_\alpha$. This construction is based
on the notion of free $(3,1)$-category generated by a cellular extension, which
is the aim of the following section. The full proof of
Theorem~\ref{thm:HomotopyBasesTransfer} will be given in
Section~\ref{SS:ProofHomotopyBasesTransfer}.

\section{Coherent Presentations}
\label{sec:CoherentPresentation}

\subsection{(3,1)-polygraphs}
\label{sec:(3,1)Polygraph}
\index{polygraph!31-@$(3,1)$-}
\index{31-polygraph@$(3,1)$-polygraph}
A \emph{$(3,1)$-polygraph} is a pair $(P,P_3)$ consisting of a
$2$\nbd-poly\-graph~$P$ and a cellular extension $P_3$ of the free
$(2,1)$-category $\freegpd{P}$.
It thus consists of a diagram of sets and functions
\[
\vxym{
&
  \ar@<-.5ex>[dl]_-{\sce0}
  \ar@<+.5ex>[dl]^-{\tge0}
P_1\ar[d]^{\ins1}
&
  \ar@<-.5ex>[dl]_-{\sce1}
  \ar@<+.5ex>[dl]^-{\tge1}
P_2\ar[d]^{\ins2}
&
  \ar@<-.5ex>[dl]_-{\sce2}
  \ar@<+.5ex>[dl]^-{\tge2}
P_3
\\
P_0
&
  \ar@<-.5ex>[l]_-{\freecat{\sce0}}
  \ar@<+.5ex>[l]^-{\freecat{\tge0}}
\freecat{P_1}
&
  \ar@<-.5ex>[l]_-{\freecat{\sce1}}
  \ar@<+.5ex>[l]^-{\freecat{\tge1}}
\freegpd{P_2}
}
\]
together with the compositions and identities of the underlying
$(2,1)$-category
\[
  \xymatrix{
    P_0&
    \ar@<-.5ex>[l]_-{\freecat{\sce0}}
    \ar@<+.5ex>[l]^-{\freecat{\tge0}}
    \freecat{P_1}
    &
    \ar@<-.5ex>[l]_-{\freecat{\sce1}}
    \ar@<+.5ex>[l]^-{\freecat{\tge1}}
    \freegpd{P_2}
}
\]
whose \emph{source} and \emph{target maps} $\sce{i}$ and $\tge{i}$ satisfy the
globular relations
\[
  \freecat{\sce{i}}\circ\sce{i+1} = \freecat{\sce{i}}\circ\tge{i+1}
  \qqtand
  \freecat{\tge{i}}\circ\sce{i+1} = \freecat{\tge{i}}\circ\tge{i+1}
\]
for every $i\in\set{0,1}$. The elements of the cellular extension $P_3$ are
called the \emph{$3$-generators} of the $(3,1)$-poly\-graph~$(P,P_3)$. We write
$A:\phi\TO\psi$ for a $3$-generator $A$ in $P_3$ such that $\sce2(A)=\phi$ and
$\tge2(A)=\psi$, often pictured as
\[
\xymatrix @C=6em{
u
  \ar@2 @/^4ex/ [r] ^-{\phi} _-{}="src"
  \ar@2 @/_4ex/ [r] _-{\psi} ^-{}="tgt"
  \ar@3 "src"!<0pt,-10pt>;"tgt"!<0pt,10pt> ^-{A}
&
v 
}
\]
and, more generally, we will call $A$ a \emph{$3$-generator} of the cellular
extension. A $(3,1)$\nbd-poly\-graph~$P$ will be also denoted by
\[
  \PRes{P_0}{P_1}{P_2}{P_3}
\]
and we will write $\tpol{k} P$ for its underlying $k$-polygraph for
$0\leq k \leq 2$.
Morphisms of $(3,1)$-polygraphs are defined as for $3$-polygraphs
(\cref{sec:Pol3}) and we denote by $\npPol31$ the resulting category.

\begin{example}
  \label{Example:CoherenceAss_3}
  As a simple example, consider the $(3,1)$-polygraph with only one generator in each
  dimension:
  \[
    \PRes\star{a}{\alpha:aa\To a}{A:a\alpha\comp{}\alpha\TO\alpha a\comp{}\alpha}\pbox.
  \]
  The $3$-generator $A$ can be represented by the diagram
  \[
    \vcenter{\xymatrix @! @R=0.15em @C=2em {
        & aa
        \ar@2@/^/ [dr]^{\alpha}
        \ar@3 []!<0pt,-12.5pt>;[dd]!<0pt,12.5pt> ^-*+{A}
        \\
        aaa 
        \ar@2@/^/ [ur]^{a\alpha} 
        \ar@2@/_/ [dr]_{\alpha a}
        && a 
        \\
        & aa
        \ar@2@/_/ [ur]_{\alpha}
      }}
    \pbox.
  \]
  In \cref{sec:CoherentConfluence}, \cref{Example:aa->aCoherent}, we use
  a rewriting argument to show that the $3$-generator $A$ forms an acyclic
  extension of the free $(2,1)$-category generated by~$\Pres{\star}{a}{\alpha}$.
\end{example}

\subsection{Free (3,1)-category}
\index{3-category@$3$-category}
\index{category!3-@$3$-}
The definition of \emph{3-category} is adapted from the one
of $2$-category by replacing the hom-categories and the composition functors by
hom-$2$-categories and composition $2$-functors. We refer the reader to
\chapr{n-cat} for the complete definition of strict $n$-categories for all
$n\geq 0$. In a $3$-category, the $3$-cells can be composed in three different
ways:
\begin{itemize}
\item by~$\comp0$, along their $0$-dimensional boundary:
  \[
    \xymatrix @C=2.5em @!C{
      x
      \ar @/^6ex/ [rr] ^{u} _{}="src1"
      \ar @/_6ex/ [rr] _{u'} ^{}="tgt1"
      && y
      \ar@2 "src1"!<-15pt,-15pt>;"tgt1"!<-15pt,15pt> _{f}="srcA"
      \ar@2 "src1"!<15pt,-15pt>;"tgt1"!<15pt,15pt> ^{f'}="tgtA"
      \ar@3 "srcA"!<15pt,0pt>;"tgtA"!<-15pt,0pt> ^*+{A}
      \ar @/^6ex/ [rr] ^{v} _{}="src2"
      \ar @/_6ex/ [rr] _{v'} ^{}="tgt2"
      && z
      \ar@2 "src2"!<-15pt,-15pt>;"tgt2"!<-15pt,15pt> _{g}="srcB"
      \ar@2 "src2"!<15pt,-15pt>;"tgt2"!<15pt,15pt> ^{g'}="tgtB"
      \ar@3 "srcB"!<15pt,0pt>;"tgtB"!<-15pt,0pt> ^*+{B}
    }
    \quad\longmapsto\quad
    \xymatrix @C=3.2em @!C{
      x
      \ar @/^6ex/ [rr] ^{uv} _{}="src1"
      \ar @/_6ex/ [rr] _{u'v'} _{}="tgt1"
      && z\pbox,
      \ar@2 "src1"!<-20pt,-15pt>;"tgt1"!<-20pt,15pt> _{fg}="srcA"
      \ar@2 "src1"!<20pt,-15pt>;"tgt1"!<20pt,15pt> ^{f'g'}="tgtA"
      \ar@3 "srcA"!<20pt,0pt>;"tgtA"!<-20pt,0pt> ^*+{A\comp0 B}
    }
  \]
\item by~$\comp1$, along their $1$-dimensional boundary:
  \[
    \xymatrix @C=4em @!C{
      x
      \ar @/^8ex/ [rr] ^-{u} _{}="src1"
      \ar [rr] ^-{v} ^{}="tgt1" _{}="src2" 
      \ar @/_8ex/ [rr] _-{w} _{}="tgt2"
      && y
      \ar@2 "src1"!<-20pt,-10pt>;"tgt1"!<-20pt,0pt> _{f}="srcA"
      \ar@2 "src1"!<20pt,-10pt>;"tgt1"!<20pt,0pt> ^{f'}="tgtA"
      \ar@2 "src2"!<-20pt,0pt>;"tgt2"!<-20pt,10pt> _{g}="srcB"
      \ar@2 "src2"!<20pt,0pt>;"tgt2"!<20pt,10pt> ^{g'}="tgtB"
      \ar@3 "srcA"!<20pt,5pt>;"tgtA"!<-20pt,5pt> ^*+{A}
      \ar@3 "srcB"!<20pt,0pt>;"tgtB"!<-20pt,0pt> _*+{B}
    }
    \qquad\longmapsto\qquad
    \xymatrix @C=4.5em @!C{
      x
      \ar @/^8ex/ [rr] ^{u} _{}="src1"
      \ar @/_8ex/ [rr] _{w} _{}="tgt1"
      && y\pbox,
      \ar@2 "src1"!<-20pt,-15pt>;"tgt1"!<-20pt,15pt> _{f\comp1 g} ^{}="srcA"
      \ar@2 "src1"!<20pt,-15pt>;"tgt1"!<20pt,15pt> ^{f'\comp1 g'} _{}="tgtA"
      \ar@3 "srcA"!<10pt,0pt>;"tgtA"!<-10pt,0pt> ^*+{A\comp1 B}
    }
  \]
\item by~$\comp2$, along their $2$-dimensional boundary:
  \[
    \xymatrix @C=4em @!C{
      x
      \ar @/^6ex/ [rr] ^{u} _{}="src"
      \ar @/_6ex/ [rr] _{v} _{}="tgt"
      && y
      \ar@2 "src"!<-30pt,-15pt>;"tgt"!<-30pt,15pt> _{f} ^{}="srcA"
      \ar@2 "src"!<0pt,-10pt>;"tgt"!<0pt,10pt> |*+{g} _{} ^{}="srcB" _{}="tgtA"
      \ar@2 "src"!<30pt,-15pt>;"tgt"!<30pt,15pt> ^{h} _{}="tgtB"
      \ar@3 "srcA"!<10pt,0pt>;"tgtA"!<-10pt,0pt> ^*+{A}
      \ar@3 "srcB"!<10pt,0pt>;"tgtB"!<-10pt,0pt> ^*+{B}
    }
    \qquad\longmapsto\qquad
    \xymatrix @C=4em @!C{
      x
      \ar @/^6ex/ [rr] ^{u} _{}="src1"
      \ar @/_6ex/ [rr] _{v} _{}="tgt1"
      && y\pbox.
      \ar@2 "src1"!<-20pt,-15pt>;"tgt1"!<-20pt,15pt> _{f} ^{}="srcA"
      \ar@2 "src1"!<20pt,-15pt>;"tgt1"!<20pt,15pt> ^{h} _{}="tgtA"
      \ar@3 "srcA"!<10pt,0pt>;"tgtA"!<-10pt,0pt> ^*+{A\comp2 B}
    }
  \]
\end{itemize}
\index{31-category@$(3,1)$-category}
\index{category!31-@$(3,1)$-}
A \emph{$(3,1)$-category} is a $3$-category whose
$2$-cells are invertible \wrt the composition $\comp1$ and whose $3$-cells are
invertible \wrt the composition~$\comp2$ (which implies their invertibility \wrt
$\comp1$).

\index{free!31-category@$(3,1)$-category}
The \emph{free $(3,1)$-category} over a
$(3,1)$-polygraph~$P$ is the $(3,1)$-category, denoted by
$\freegpd{P}$, or $\freegpd{P}_{\leq 2}(P_3)$ in some contexts, whose
\begin{itemize}
\item underlying $2$-category is the free $(2,1)$-category~$\freegpd{\tpol2P}$,
\item set $\freegpd{P_3}$ of $3$-cells consists of all formal compositions \wrt
  $\comp0$, $\comp1$ and $\comp2$ of $3$-generators of~$P$, of their inverses, and
  of identities of $2$-cells, considered up to associativity, identity, exchange,
  and inverse relations.
\end{itemize}
This construction will be detailed in arbitrary dimension $n\geq 1$ in
\cref{sec:np-polygraph}.

We denote by $\finv A$ the inverse \wrt $\comp2$ of a $3$-cell $A$: it satisfies
$A\comp2 A^- = 1_{\src2(A)}$ and $A^- \comp2 A = 1_{\tgt2(A)}$.
Note that if a $3$-cell $A$ is invertible with respect to the composition
$\ast_2$, and its $2$-source and $2$-target are invertible, then it is
invertible with respect the composition $\ast_1$, with inverse given by
\[
t_2(A)^- \ast_1 A^- \ast_1 s_2(A)^-.
\]
Every $3$-cell $A$ of the $(3,1)$-category $\freegpd{P}$ of size $k\geq 1$ has a
decomposition
\[
A = C_1[A_1^{\epsilon_1}] \ast_2 \ldots \ast_2 C_k[A_k^{\epsilon_k}],
\]
with $\epsilon_1,\ldots,\epsilon_k\in\set{-,+}$, and $A_1,\ldots,A_k$ are
$3$-generators of $P$, where for every $1\leq i \leq k$,
$C_i[A_i^{\epsilon_i}]$ denotes a composition of the form
\[
f_2\ast_2(f_1\ast_0 A_i^{\epsilon_i} \ast_0 g_1) \ast_2 g_2,
\]
where $f_j,g_j$ are $j$-cells for $j=1,2$, and where $A_i^+$ is equal to $A_i$

\subsection{Coherent presentations}
\index{extended presentation}
\index{presentation!extended}
\index{coherent!presentation}
\index{coherent!31-polygraph@$(3,1)$-polygraph}
\index{presentation!coherent}
A $(3,1)$-polygraph $P$ is \emph{coherent} when $P_3$ is an acyclic extension of
the free $(2,1)$-category $\freegpd{\tpol2P}$. 
This amounts to requiring that for every pair of parallel
$2$-cells $\phi$ and $\psi$ in $\freegpd P_2$, there is a $3$-cell
$F:\phi\TO\psi$ in $\freegpd P_3$.

An \emph{extended presentation} of a $1$-category~$C$ is a $(3,1)$-poly\-graph~$P$
whose underlying $2$\nbd-polygraph~$\tpol2 P$ is a presentation of $C$.  A
\emph{coherent presentation} of~$C$ is an extended presentation~$P$ of $C$ which
is coherent.

\subsection{Proof of \cref{thm:HomotopyBasesTransfer}}
\label{SS:ProofHomotopyBasesTransfer}
Let us denote by $Y$ the cellular extension $f(X)\sqcup X_Q$. We construct, for
every $2$-cell $\varphi: u \To v$ of $\freegpd{Q}$, a $3$-cell $A_\varphi$ of
the free $(3,1)$-category $\freegpd{Q}(Y)$ with the following shape:
\[
  \xymatrix @!C @C=2.5em @R=1em{
    & fg(v)
    \ar@2 @/^/ [dr] ^-{\psi_v}
    \ar@3 []!<0pt,-15pt>;[dd]!<0pt,15pt> ^-{A_\varphi}
    \\
    fg(u) 
    \ar@2 @/^/ [ur] ^-{fg(\varphi)}
    \ar@2 @/_/ [dr] _-{\psi_u}
    && v
    \\
    & u
    \ar@2 @/_/ [ur] _-{\varphi}
  }
\]
by extending the notation $A_{\alpha}$, where $\alpha$ is a $2$-generator of
$Q$, in a functorial way, according to the following formulas:
\begin{align*}
  A_{1_u}&=1_{\psi_u},
  &
  A_{\varphi\comp0 \varphi'}&=A_\varphi\comp0 A_{\varphi'},
  &
  A_{\varphi^-}&=fg(\varphi)^- \comp1 A_{\varphi}^- \comp1 {\varphi}^-,
\end{align*}
\[
  A_{\varphi\comp1 \varphi'} = \pa{fg(\varphi) \comp1 A_{\varphi'}}\comp2\pa{A_\varphi \comp1 \varphi'}
  \pbox.
\]
We prove that the $3$-cells $A_\varphi$ are well-defined, \ie their
definition is compatible with the relations on $2$-cells, such as the exchange
relation.
Indeed, whenever the composition of the $2$-cells $\varphi_1$, $\varphi_2$, $\varphi_1'$
and $\varphi_2'$ are defined in $\freegpd{Q}$, we have
\begin{align*}
    A_{(\varphi_1\comp0\varphi_2)\comp1(\varphi'_1\comp0\varphi'_2)}
    &=\big((fg(\varphi_1)\comp0 fg(\varphi_2))\comp1(A_{\varphi_1'}\comp0A_{\varphi_2'}\big)\comp2
    \\
    &\hspace{4cm}\big((A_{\varphi_1}\comp0A_{\varphi_2})\comp1 (\varphi_1'\comp0\varphi_2')\big)
    \\
    &=\big((fg(\varphi_1)\comp1A_{\varphi_1'})\comp0(fg(\varphi_2)\comp1 A_{\varphi_2'})\big)\comp2
    \\
    &\hspace{4cm}\big((A_{\varphi_1}\comp1\varphi_1')\comp0(A_{\varphi_2}\comp1\varphi_2')\big)
    \\
    &=\big((fg(\varphi_1)\comp1A_{\varphi_1'})\comp2(A_{\varphi_1}\comp1\varphi_1')\big) \comp0 
    \\
    &\hspace{3.3cm}\big((fg(\varphi_2)\comp1 A_{\varphi_2'})\comp2(A_{\varphi_2}\comp1\varphi_2')\big)
    \\
    &= A_{(\varphi_1\comp1 \varphi'_1)\comp0(\varphi_2\comp1 \varphi'_2)}
    \pbox.
  \end{align*}

Now, let us consider two parallel
$2$-cells $\varphi,\varphi':u\To v$ of $\freegpd{Q}$. The $2$-cells $g(\varphi)$
and~$g(\varphi')$ are parallel in~$\freegpd{P}$ so that, by $X$-acyclicity
of $\freegpd{P}$, there exists a $3$-cell
\[
  \xymatrix @C=5em{
    g(u)
    \ar@2 @/^4ex/ [r] ^-{g(\varphi)} _-{}="src"
    \ar@2 @/_4ex/ [r] _-{g(\varphi')} ^-{}="tgt"
    \ar@3 "src"!<0pt,-10pt>;"tgt"!<0pt,10pt> ^-{A}
    &
    g(v)
  }
\]
in the $(3,1)$-category $\freegpd{P}(X)$. By definition of $Y$ and functoriality of $f$, there exists a $3$-cell
\[
  \xymatrix @C=5em{
    fg(u)
    \ar@2 @/^4ex/ [r] ^-{fg(\varphi)} _-{}="src"
    \ar@2 @/_4ex/ [r] _-{fg(\varphi')} ^-{}="tgt"
    \ar@3 "src"!<0pt,-10pt>;"tgt"!<0pt,10pt> ^-{f(A)}
    &
    fg(v)
  }
\]
in the free $(3,1)$-category $\freegpd{P}(Y)$. Using the $3$-cells $f(A)$,
$A_\varphi$ and $A_{\varphi'}$, we get the following $3$-cell from $\varphi$ to
$\varphi'$ in $\freegpd{P}(Y)$:
\[
  \xymatrix @!C @C=4em{
    {\sm u}
    \ar@2 @/^11ex/ [rrr] ^{\varphi} _{}="src1"
    \ar@2 [r] ^-{\psi_u^-} 
    \ar@2 @/_11ex/ [rrr] _{\varphi'} ^{}="tgt3"
    &
    {\sm fg(u)}
    \ar@2 @/^4ex/ [r] ^(0.7){fg(\varphi)} ^{}="tgt1" _{}="src2"
    \ar@2 @/_4ex/ [r] _(0.7){fg(\varphi')} ^{}="tgt2" _{}="src3"
    & {\sm fg(v)}
    \ar@2 [r] ^-{\psi_v}
    & {\sm v}\pbox.
    \ar@3 "src1"!<-15pt,-10pt>;"tgt1"!<-15pt,10pt> ^-{\psi_u^-\comp1A_\varphi^-}
    \ar@3 "src2"!<-5pt,-10pt>;"tgt2"!<-5pt,10pt> ^-{f(A)}
    \ar@3 "src3"!<-15pt,-10pt>;"tgt3"!<-15pt,10pt> ^-{\psi_u^-\comp1A_{\varphi'}}
  }
\]
This concludes the proof that the
$(2,1)$-category~$\freegpd{Q}$ is $Y$-acyclic.

\subsection{Cofibrant replacements and coherent presentation}
The notion of coherent presentation of a category corresponds to the notion of
cofibrant replacement for the model structure for $2$-categories introduced by Lack
in~\cite{LackFolk2,LackFolkBi}. This will be detailed and generalized in
\cref{chap:folk}, and we only give here a brief overview. 
In this model structure a $2$-category is \emph{cofibrant}\index{cofibrant!2-category@$2$-category} if its
underlying $1$-category is free, and a $2$-functor $F:C\fl D$ is a \emph{weak
  equivalence} if it satisfies the following two conditions.
\begin{enumerate}
\item Every $0$-cell $y$ of $D$ is \emph{equivalent} to a $0$-cell $F(x)$ for
  $x$ in $C$, \ie there exist $1$\nbd-cells $u:F(x)\fl y$ and $v:y\fl F(x)$ and
  invertible $2$-cells $f:u\comp1 v \dfl 1_{F(x)}$ and $g:v\comp1 u\dfl 1_{y}$
  in $D$.
\item For every $0$-cells $x$ and $x'$ in $C$, the induced functor
  \[
    F(x,x'):C(x,x')\fl D(F(x),F(x'))
  \]
  is an equivalence of categories.
\end{enumerate}
In particular, an equivalence of $2$-categories is a weak equivalence. A
\emph{cofibrant replacement}\index{cofibrant!replacement} of a $2$-category $C$ is a cofibrant $2$-category
$\widetilde{C}$ that is weakly equivalent to~$C$. The following theorem is
proved in~\cite[Theorem~1.3.1]{GaussentGuiraudMalbos15}:

\begin{theorem}
  \label{Theorem:CoherentPresentationsCofibrantReplacements}
  Let $P$ be an extended presentation of a category $C$. Then the
  $(3,1)$\nbd-poly\-graph $P$ is a coherent presentation of $C$ if and only if
  the $(2,1)$\nbd-cate\-gory $\cl{P}$ is a cofibrant replacement of $C$.
\end{theorem}

Note that a given category~$C$ may admit other cofibrant replacements than the
$2$-categories presented by coherent presentations of~$C$. For instance,
consider the terminal category $\termcat$: it contains one $0$-cell and the
corresponding identity $1$-cell only. This category $\termcat$ is cofibrant and,
as a consequence, is a cofibrant replacement of itself: this cofibrant
replacement corresponds to the coherent presentation of the terminal category
given by the $(3,1)$-polygraph with one $0$-generator and no higher-dimensional
generators. But the terminal category also admits, as a cofibrant replacement, the
$2$-category with two $0$-cells $x,y$, two $1$-cells $u,v$ as follows
\[
  \xymatrix @C=3em{
    x
    \ar @/^2ex/ [r] ^-{u} _-{}="src"
    &
    y
    \ar @/^2ex/ [l] ^-{v} ^-{}="tgt"
  }
\]
and two invertible $2$-cells $f:uv\dfl 1_x$ and $g:vu\dfl 1_y$. However, this
$2$-category is not presented by a coherent presentation of the terminal
category, since it has two $0$-cells.

\section{Coherent Confluence}
\label{sec:CoherentConfluence}

In this section, we extend to $(3,1)$-polygraphs the results on coherent
confluence given in \cref{sec:free-21-category} for $(2,0)$-polygraphs.

\subsection{Coherent confluence}
\label{sec:squier-coherence}
\index{confluence!coherent}
\index{convergence!coherent}
Let~$P$ be a
$(3,1)$-polygraph. A branching~$(\phi,\psi)$ of~$P$ is \emph{coherently
  confluent} if there exist $2$-cells~$\phi'$ and~$\psi'$ in $\freecat P_2$ and
a $3$-cell~$A$ in $\freegpd P_3$ of the form
\[
  \xymatrix@!C=1ex@!R=1ex{
    &\ar@2[dl]_\phi u\ar@2[dr]^\psi&\\
    v\ar@{:>}[dr]_{\phi'}&\underset{\phantom 1}{\overset A\TO}&\ar@{:>}[dl]^{\psi'}w\pbox.\\
    &w&
  }
\]
We say that~$P$ is \emph{coherently confluent} (\resp \emph{locally coherently
  confluent}, \resp \emph{critically coherently confluent}) when every branching
(\resp local branching, \resp critical branching) of~$P$ is coherently
confluent. We say that~$P$ is \emph{coherently convergent} if it terminates and
is coherently confluent. Note that, given a $2$\nbd-poly\-graph~$P$, by
taking~$P_3=\Sph(\freecat{P})$ to be the set of all $2$-spheres of~$\freecat P$ (see \cref{sec:2-spheres}),
the notions of coherent confluence and coherent convergence in the
$(3,1)$-polygraph $(P,P_3)$ boil down to the ones of confluence and convergence
of the $2$-polygraph defined in
\cref{Subsection:RewritingProperties2Polygraphs}.

The following result essentially amounts to a coherent version of Newman's lemma for
$2$\nbd-polygraphs (\cref{lem:newman}). Its proof is essentially the same as in
the case of $1$\nbd-poly\-graphs (\cref{L:HNewman-1}). It is interesting to note that
the geometric intuition behind it was already present in Newman's original article~\cite[Section~6]{newman1942theories}.

\begin{proposition}
  \label{P:HNewman-2}
  Let~$P$ be a terminating $(3,1)$-polygraph. If~$P$ is locally coherently
  confluent, then~$P$ is coherently confluent.
\end{proposition}

\noindent
As above, we recover the Newman's lemma for $2$-polygraphs, by taking $P_3$ to
be the set of all $2$-spheres. Similarly, the following result is a coherent
version of the critical branching lemma (\cref{lem:2-cb}).

\begin{lemma}
  \label{L:HCritical}
  Let~$P$ be a $(3,1)$-polygraph. If~$P$ is critically coherently confluent,
  then~$P$ is locally coherently confluent.
\end{lemma}
\begin{proof}
  We proceed by case analysis on the type of the local branchings of~$P$. First,
  non-overlapping (\ie trivial and orthogonal) branchings are always coherently
  confluent. Indeed, if~$\phi:u\To v$ is a rewriting step of~$P$, then the
  trivial branching~$(\phi,\phi)$ is coherently confluent because of
  \[
    \xymatrix@!C=1ex@!R=1ex{
      &\ar@2[dl]_\phi u\ar@2[dr]^\phi&\\
      v\ar@2[dr]_{\unit v}&=&\ar@2[dl]^{\unit v}v\pbox.\\
      &v&
    }
  \]
  And, if $\phi:u\To u'$ and $\psi:v\To v'$ are rewriting steps of~$P$, then the
  orthogonal branching $(\phi v,u\psi)$ is coherently confluent thanks to the following equality
  \[
    \xymatrix@!C=1ex@!R=1ex{
      &\ar@2[dl]_{\phi v} uv\ar@2[dr]^{u\psi}&\\
      u'v\ar@2[dr]_{u'\psi}&=&\ar@2[dl]^{\phi v'}uv'\pbox.\\
      &u'v'&
    }
  \]
  Now, assume that~$(\phi,\psi)$ is an overlapping branching, where
  $\phi:u\To v$ and $\psi:u\To w$ are rewriting steps of~$P$. Then we have $u=u_1u'u_2$,
  $\phi=u_1\phi'u_2$ and $\psi=u_1\psi' u_2$, where~$u_1$, $u'$ and~$u_2$ are $1$-cells
  of~$\freecat{P_1}$, and~$\phi'$ and~$\psi'$ are rewriting steps of~$P$ such
  that $(\phi',\psi')$ is a critical branching of~$P$. By hypothesis,
  $(\phi',\psi')$ is critically coherently confluent, from which we deduce the
  existence of $2$\nbd-cells~$\phi''$ and~$\psi''$ of~$\freecat{P_2}$, and of a
  $3$-cell~$F$ of~$\freegpd{P_3}$
  \[
    \xymatrix@!C=1ex@!R=1ex{
      &\ar@2[dl]_{u_1\phi'u_2} u_1u'u_2\ar@2[dr]^{u_1\psi'u_2}&\\
      u_1v'u_2\ar@2[dr]_{u_1\phi''u_2}\ar@3[]!<0ex,0pt>;[rr]!<0ex,0pt>^{u_1Fu_2}&&\ar@2[dl]^{u_1\psi''u_2}u_1w'u_2\\
      &u_1u''u_2&
    }
  \]
  proving that~$(\phi,\psi)$ is coherently confluent.
\end{proof}

\noindent
Again, the critical branching lemma (\lemr{2-cb}) can be recovered by
taking~$P_3$ to be the set of $2$-spheres in a $2$-polygraph.

\begin{proposition}
  \label{prop:PreSquier}
  Let~$P$ be a $(3,1)$-polygraph. If~$P$ is coherently convergent then~$P$ is
  coherent.
\end{proposition}

As for $1$-polygraphs we prove the following coherence result, called
\emph{Squier's homotopical
  theorem}~\cite[Theorem~5.2]{squier1994finiteness}  (see
also~\cite{guiraud2009higher,lafont1995new}).

\begin{theorem}
  \index{Squier!theorem}
  \label{thm:SquierHomotopical}
  \label{Theorem:SquierCompletion2polygraphs}
  Let~$P$ be a convergent 2-poly\-graph, and~$P_3$ be a cellular extension
  of the free $(2,1)$-category~$\freegpd{P}$. If~$P_3$ contains, for every critical
  branching~$(\phi,\psi)$ of~$P$, one 3-generator of the form
  \begin{equation}
    \label{eq:squier-coh-cell}
    \xymatrix@!C=1ex@!R=1ex{
      &\ar@2[dl]_\phi u\ar@2[dr]^\psi&\\
      v\ar@2[dr]_{\phi'}\ar@3[]!<5ex,0ex>;[rr]!<-5ex,0ex>^A&&\ar@2[dl]^{\psi'}w\\
      &u'&
    }
  \end{equation}
  where~$\phi'$ and~$\psi'$ are 2-cells in~$\freecat{P}_2$, then the
  $(3,1)$-polygraph $(P,P_3)$ is coherent.
\end{theorem}

A $3$-generator of the form \eqref{eq:squier-coh-cell}, indexed by a critical branching of $P$, is called a \emph{generating confluence} of the polygraph $P$. \cref{Theorem:SquierCompletion2polygraphs} states that the set of generating confluences of a convergent $2$-polygraph $P$, indexed by all its critical branchings, forms an acyclic extension of the $(2,1)$-category $\freegpd{P}$. 

\begin{example}
\label{Example:aa->aCoherent}
Consider the monoid~$M$ with the convergent presentation 
  \[
    \Pres\star{a}{\alpha:aa\To a}
    \pbox.
  \]
  This $2$-polygraph has exactly one critical branching, whose corresponding
  generating confluence has the form:
  \[
    \xymatrix@!C=1ex@!R=1ex{
      &\ar@2[dl]_{\alpha a}aaa\ar@2[dr]^{a\alpha}&\\
      aa\ar@2[dr]_\alpha\ar@3[]!<5ex,0ex>;[rr]!<-5ex,0ex>^A&&\ar@2[dl]^\alpha aa\pbox.\\
      &a&
    }
  \]
  By \cref{thm:SquierHomotopical}, the $(3,1)$-polygraph
  $\PRes\star{a}{\alpha}{A}$ defined in Example~\ref{Example:CoherenceAss_3} is
  thus a coherent presentation of the monoid $M$.
\end{example}

\subsection{The standard coherent presentation}
\index{standard!presentation}
\index{presentation!standard}
\label{X:StandardCoherentPresentation}
Recall from~\cref{sec:2-std-pres,sec:cat-std-pres} that the \emph{standard presentation} of a
category $C$ is the $2$-polygraph $\Std_2(C)$ such that
\begin{itemize}
\item the $0$-generators are the $0$-cells of~$C$,
\item there is a $1$-generator $\rep{f}:x\fl y$ for every $1$-cell
  $f:x\fl y$ of $C$,
\item there is a $2$-generator $\mu_{f,g}:\rep{f}\rep{g}\dfl\rep{fg}$ for all
  composable $1$-cells $f$ and $g$ of~$C$,
\item there is a $2$-generator $\eta_x:1_x\dfl\rep{1}_x$ for every $0$-cell $x$
  of~$C$.
\end{itemize}
\index{standard!presentation!coherent}
\index{presentation!standard!coherent}
The \emph{standard coherent presentation} $\Std_3(C)$ of $C$ is the presentation
$\Std_2(C)$ extended with the following $3$-generators
\begin{align*}
  \vcenter{
    \xymatrix @C=1em @R=1em {
      & {\rep{fg}\rep{h}}
      \ar@2 @/^/ [dr] ^-{\mu_{fg,h}}
      \\
      {\rep{f}\rep{g}\rep{h}}
      \ar@2 @/^/ [ur] ^-{\mu_{f,g}\rep{h}}
      \ar@2 @/_/ [dr] _-{\rep{f}\mu_{g,h}}
      && {\rep{fgh}}
      \\
      & {\rep{f}\rep{gh}}
      \ar@2 @/_/ [ur] _-{\mu_{f,gh}}
      \ar@3 "1,2"!<0pt,-15pt>;"3,2"!<0pt,15pt> ^-{A_{f,g,h}}
    }
  }
  &&
  \vcenter{
    \xymatrix@C=2em{
      & {\rep{1}_x\rep{f}}
      \ar@2@/^/ [dr] ^-{\mu_{\id{x},f}}
      \\
      {\rep{f}}
      \ar@2@/^/ [ur]  ^{\eta_x \rep{f}}
      \ar@2{=}@/_/ [rr] _-{}="tgt"
      && {\rep{f}}
      \ar@3 "1,2"!<0pt,-15pt>;"tgt"!<0pt,10pt> ^-{L_f}
    }
  }
  &&
  \vcenter{
    \xymatrix@C=2em{
      & {\rep{u}\rep{\id_y}}
      \ar@2@/^/ [dr] ^-{\mu_{f,\id{y}}}
      \\
      {\rep{f}}
      \ar@2@/^/ [ur] ^-{\rep{f}\eta_y}
      \ar@2{=}@/_/ [rr] _-{}="tgt"
      && {\rep{f}}
      \ar@3 "1,2"!<0pt,-15pt>;"tgt"!<0pt,10pt> ^-{R_f}
    }
  }  
\end{align*}
for all $1$-cells $f:x\fl y$, $g:y\fl z$ and $h:z\fl t$ of $C$. Those
$3$-generators can be shown to form an acyclic cellular extension of the free
$(2,1)$\nbd-cate\-gory~$\freegpd{\Std_2(\C)}$ (by first reversing the
orientation of the generators $\eta_x$, as explained in \cref{sec:cat-std-pres},
and then applying \cref{thm:SquierHomotopical}).

\section{Tietze Transformations of (3,1)-Polygraphs}
\label{Section:TietzeTransformations31Polygraphes}
We extend here the notion of Tietze transformation presented in
\secr{1pol-tietze} for 1-polygraphs and in \cref{sec:2-tietze} for 2-polygraphs to
$(3,1)$\nbd-poly\-graphs.

\subsection{Tietze transformations}
\index{Tietze!transformation!of 31-polygraphs@of $(3,1)$-polygraphs}
\index{transformation!Tietze}
If $P$ is a $(3,1)$-polygraph, an \emph{elementary Tietze transformation} on~$P$ is one of
the following operations transforming $P$ into a $(3,1)$\nbd-poly\-graph~$Q$:
\begin{description}
\item[\tgen{}] \emph{adding a definable 1-cell}: $Q$ is obtained from~$P$ by
  adding a $1$-generator $a:x\to y$ together with a $2$-generator
  $\alpha:u\To a$ for some $1$-cell $u:x\to y\in\freecat P_1$:
  \[
    \xymatrix@C=8ex{
      x\ar@/^2.5ex/[r]^u&y
    }
    \qquad\rightsquigarrow\qquad
    \xymatrix@C=8ex{
      x\ar@/^2.5ex/[r]^u\ar@/_2.5ex/[r]_a\ar@{}[r]|{\alpha\Downarrow}&y\pbox,
    }
  \]
\item[\trel{}] \emph{adding a derivable 2-cell}: $Q$ is obtained from~$P$ by adding a
  $2$-generator $\alpha:u\To v$ together with a $3$-generator $A:\phi\TO\alpha$
  for some $3$-cell $\phi:u\To v\in\freegpd P_2$:
  \[
    \xymatrix@C=8ex{
      x\ar@/^2.5ex/[r]^u\ar@/_2.5ex/[r]_v\ar@{}[r]|{\phi\Downarrow}&y
    }
    \qquad\rightsquigarrow\qquad
    \xymatrix@C=8ex{
      x\ar@/^2.5ex/[r]^u\ar@/_2.5ex/[r]_v\ar@{}[r]|{\phi\Downarrow\overset A{\underset{\phantom{A}}\TO}\Downarrow\alpha}&y\pbox,
    }
  \]
\item[\tRel{}] \emph{adding a derivable 3-cell}: $Q$ is obtained from~$P$ by adding a
  $3$-generator $A:\phi\TO\psi$ for some $3$-cell
  $F:\phi\TO\psi\in\freegpd P_3$:
  \[
    \xymatrix@C=8ex{
      x\ar@/^2.5ex/[r]^u\ar@/_2.5ex/[r]_v\ar@{}[r]|{\phi\Downarrow\overset F{\underset{\phantom{F}}\TO}\Downarrow\psi}&y
    }
    \qquad\rightsquigarrow\qquad
    \xymatrix@C=8ex{
      x\ar@/^2.5ex/[r]^u\ar@/_2.5ex/[r]_v\ar@{}[r]|{\phi\Downarrow\overset A{\underset{\phantom{A}}\TO}\Downarrow\psi}&y\pbox.
    }
  \]
\end{description}
A \emph{Tietze transformation} between $(3,1)$-polygraphs~$P$ and~$Q$ is a
finite sequence of polygraphs $P=P_1,P_2,\ldots,P_n=Q$ such that, for
$1\leq i<n$, either $P_{i+1}$ is obtained from $P_i$ by an elementary Tietze
transformation, or $P_i$ is obtained from~$P_{i+1}$ by an elementary Tietze
transformation.
\index{Tietze!equivalence!of 31-polygraphs@of $(3,1)$-polygraphs}
\index{equivalence!Tietze}
Two $(3,1)$-polygraphs are \emph{Tietze equivalent} when there is a Tietze
transformation between them. As in~\cref{sec:1pol-tietze}
and~\cref{sec:2-tietze}, the notion of Tietze equivalence is supposed to be closed
by isomorphism. 

\subsection{Functors induced by Tietze transformations}
For any of the above elementary Tietze transformations from a
$(3,1)$-polygraph~$P$ to a $(3,1)$\nbd-poly\-graph~$Q$, there is a canonical
morphism of polygraphs $P\to Q$, witnessing the inclusion of~$P$ into~$Q$, which
induces a $3$-functor $F:\freegpd P\to\freegpd Q$ between the freely generated
$(3,1)$-categories. This functor always admits a retraction, \ie a
$3$\nbd-func\-tor $G:\freegpd Q\to\freegpd P$ such that
$G\circ F=\id_{\freegpd P}$. For instance, in the case of~\tgen{}, with the same
notations as above, the functor $G$ is such that $G a=u$, $G\alpha=\id_u$, and
$G$ leaves the other generators of $Q$ unchanged.

We recall the following result
from~\cite[Theorem~2.1.3]{GaussentGuiraudMalbos15}:

\begin{theorem}
  \label{thm:31-tietze}
  Two finite $(3,1)$-polygraphs~$P$ and~$Q$ are Tietze equivalent if and only if
  there is an equivalence between the presented 2-categories $\pcat{P}$
  and~$\pcat{Q}$ which induces a bijection between the respective sets of
  0-cells.
\end{theorem}

\noindent
As a consequence, if a $(3,1)$-polygraph $P$ is a coherent presentation of a
category $C$ and if there exists a Tietze transformation from $P$ to a
$(3,1)$\nbd-poly\-graph $Q$, then $Q$ is also a coherent presentation of~$C$.

\subsection{Higher Nielsen transformations}
\index{Nielsen transformation}
\index{transformation!Nielsen}
\label{Subsection:NielsenTransformation}
As a particular subset of Tietze transformation, we identify the following family of
transformations, which will prove useful in the following.
The \emph{elementary Nielsen transformations} on a $(3,1)$\nbd-poly\-graph $P$ are the
following transformations:
\begin{description}
\item[(N1)] the replacement of a $2$-cell by a formal inverse (including in the
  source and target of every $3$-cell),
\item[(N2)] the replacement of a $3$-cell by a formal inverse,
\item[(N3)] the replacement of a $3$-cell $F : \psi\TO\psi'$ by a $3$-cell
  \[
    \tilde{F}:\phi\comp1\psi\comp1\chi\TO\phi\comp1\psi'\comp1\chi
  \]
  where $\phi$ and $\chi$ are $2$-cells of $\freegpd{P}$.
\end{description}
The \emph{Nielsen equivalence} on $(3,1)$-polygraphs is the smallest equivalence
relation identifying any two polygraphs between which there is an elementary
Nielsen transformation. The following is shown
in~\cite[Section~2.1.4]{GaussentGuiraudMalbos15}:

\begin{lemma}
  The elementary Nielsen transformations are Tietze transformations.
\end{lemma}

\subsection{Collapsible generators}
\index{collapsible!generator}
\index{generator!collapsible}
\label{sec:3collapsible}
Given a $(3,1)$-polygraph~$P$, we identify the following families of
redundant generators in the polygraph.
Following the terminology introduced by Brown~\cite{brown1992geometry}, we say
that a $2$-generator~$\alpha$ of~$P$ is \emph{collapsible} if
\begin{itemize}
\item the target of~$\alpha$ is a $1$-generator~$a\in P_1$, and
\item the source of~$\alpha$ is a $1$-cell in which $a$ does not occur.
\end{itemize}
Similarly, a $3$-generator~$A$ of $P$ is \emph{collapsible} if
\begin{itemize}
\item the target of~$A$ is a $2$-generator~$\alpha\in P_2$, and
\item the source of~$A$ is a $2$-cell in which $\alpha$ does not occur.
\end{itemize}
A \emph{$3$-sphere} $\Phi$ is a pair $(F,G)$ of $3$-cells in $\freegpd P_3$ with
the same source and with the same target, where $F$ and $G$ are respectively the
source and target of the $3$\nbd-sphere.
Note that a $3$-sphere can be seen as a $4$-generator in a $(4,1)$-polygraph,
which consists of a $(3,1)$\nbd-poly\-graph~$P$ equipped with a cellular
extension~$P_4$ of the freely generated $(3,1)$\nbd-cate\-gory~$\freegpd{P}$,
see \cref{sec:np-polygraph}. We thus denote by $\Phi:F\TOO G$ a
$3$-sphere from~$F$ to~$G$.
A $3$-sphere $\Phi:F\TOO A$ whose
target $A$ is a $3$-generator is said to be \emph{collapsible}.

Given a collapsible $2$-generator~$\alpha:u\To a$, we write $P/\alpha$ for the
$(3,1)$\nbd-poly\-graph with
\begin{itemize}
\item $P_0$ as $0$-generators,
\item $P_1\setminus\set{a}$ as $1$-generators,
\item $P_2\setminus\set{\alpha}$ as $2$-generators, where every occurrence of~$a$ in the source or target of a $2$-generator has been replaced by~$u$,
\item $P_3$ as $3$-generators, where every occurrence of~$\alpha$ in the source or target of a $3$-generator has been replaced by~$\unit{u}$.
\end{itemize}
Similarly, given a collapsible $3$-generator~$A:\phi\TO\alpha$, we write $P/A$
for the $(3,1)$\nbd-poly\-graph with
\begin{itemize}
\item $P_0$ as $0$-generators,
\item $P_1$ as $1$-generators,
\item $P_2\setminus\set{\alpha}$ as $2$-generators,
\item $P_3\setminus\set{A}$ as $3$-generators, where every occurrence of $\alpha$ in the source or target of a $3$-generator has been replaced by~$\phi$.
\end{itemize}
Similarly, given a collapsible $3$-sphere~$\Phi:F\TOO A$, we write $P/\Phi$ for
the $(3,1)$\nbd-poly\-graph $P_0$, $P_1$, $P_2$, and $P_3\setminus\set{A}$ as
sets of $0$-, $1$-, $2$-, and $3$-generators, respectively.

In the above situation, the target generator of the collapsible cell is said to
be \emph{redundant}, and the polygraph~$P/\alpha$ (\resp $P/A$, \resp $P/\Phi$)
is said to be obtained from~$P$ by \emph{collapsing} $\alpha$ (\resp $A$, \resp
$\Phi$). 

\begin{example}
In the polygraph
\[
  P=\PRes{\star}{a,b,c}{\alpha:ab\To c,\beta:ac\To bc,\gamma:ab\To c}{A:\alpha\TO\gamma}
\]
the $2$-generator~$\alpha$ is collapsible and the polygraph resulting from its
collapse is
\[
  P/\alpha=\PRes{\star}{a,b}{\beta:aab\To bab,\gamma:ab\To ab}{A:\unit{ab}\TO\gamma}
\]
\end{example}

The following is shown in~\cite[Section~2.3]{GaussentGuiraudMalbos15}:

\begin{proposition}
  \label{prop:collapse-Tietze}
  Let $P$ be a $(3,1)$-polygraph.
  Given a collapsible 2\nbd-gene\-rator~$\alpha$ (\resp 3-generator $A$, \resp
  3-sphere $\Phi$) of $P$, the $(3,1)$-polygraphs $P$ and $P/\alpha$ (\resp $P/A$,
  \resp $P/\Phi$) are Tietze equivalent.
\end{proposition}

\begin{remark}
  \label{rem:collapse-Nielsen}
  The class of collapsible generators can be made larger, and thus lead to more
  collapses, by working ``up to Nielsen equivalence''. By this, we mean that one
  can consider that a generator is collapsible in a $(3,1)$-polygraph~$P$, when
  $P$ is Nielsen equivalent to a polygraph~$Q$ in which the corresponding
  generator is collapsible. Namely, one can generalize the above notion of
  collapse to those generators.
\end{remark}

\section{Coherent Completion and Reduction}
\label{Section:HomotopicalCompletion}
Given a convergent $2$-polygraph~$P$, Squier's homotopical theorem (\cref{thm:SquierHomotopical}) provides a way to extend it into a coherent
presentation of the category~$\cl{P}$. When the 2-polygraph~$P$ is not
convergent, we can use the Knuth-Bendix completion procedure (\cref{sec:2kb}) in
order to obtain a convergent presentation of the category~$\cl{P}$ and then
apply Squier's theorem on it in order to finally obtain a coherent presentation
of~$\cl{P}$. We present the \emph{coherent completion
  procedure} from~\cite{guiraud2013homotopical,GaussentGuiraudMalbos15} which
combines the two steps at once: it adds both $2$- and $3$-generators to the
polygraph in order to obtain a coherent convergent presentation.

Often, the resulting coherent presentation of $\pcat{P}$ is not minimal, in the
sense that some of its generators are collapsible. In such a situation, it is
desirable to remove those superfluous generators in order to obtain a smaller
presentation. We also present here techniques to perform this, which, when
combined with the procedure described above, give rise to a \emph{coherent
  completion-reduction procedure}.

\subsection{Family of generating confluences}
\label{sec:FamilyGeneratingConfluences}
Given a $2$-polygraph~$P$, a cellular extension~$X$ of the free $(2,1)$-category~$\freegpd{P}$ containing a $3$-generator
\begin{equation}
  \label{eq:fgc-coh}
  \vcenter{
    \xymatrix@!C=1ex@!R=1ex{
      &\ar@2[dl]_\phi u\ar@2[dr]^\psi&\\
      v\ar@2[dr]_{\phi'}\ar@3[]!<5ex,0ex>;[rr]!<-5ex,0ex>^{A_{\phi,\psi}}&&\ar@2[dl]^{\psi'}w\\
      &u'&
    }}
\end{equation}
for every critical branching $(\phi,\psi)$ of~$P$ is called a \emph{family of
  generating confluences} for~$P$. 

\index{Squier!completion}
A \emph{Squier completion} of the polygraph~$P$ is a $(3,1)$-polygraph, denoted by $\Sr(P)$, obtained from~$P$
by adding a $3$-generator of the form~\eqref{eq:fgc-coh} for every critical
branching $(\phi,\psi)$. By \thmr{SquierHomotopical}, such a $(3,1)$-polygraph is a coherent presentation of the category~$\pcat{P}$.
Note that the notation is slightly abusive since a polygraph $\Sr(P)$ is not entirely determined by $P$. In particular, it depends on a choice of confluences for the critical branchings and the orientation of the $3$-generator $A_{\phi,\psi}$. 
We will see in \cref{Chapter:ConstructingResolutions} that this choice can be encoded as a \emph{$\iota$-contraction} and can be extended in all higher dimensions.

\subsection{Coherent completion of terminating 2-polygraphs}
\label{Subsubsection:HomotopicalCompletion}
By extending the Knuth-Bendix completion procedure, see \cref{sec:2kb}, we define a procedure that computes a coherent presentation of a category~$C$
starting with a terminating, but not necessarily confluent, presentation of $C$,
by suitably adding $2$- and $3$-generators obtained from computing critical
branchings. The procedure is defined as follows.

\index{coherent!completion}
\index{completion!coherent}
\index{procedure!coherent!completion}
\index{procedure!completion!coherent}
Given a terminating $2$-polygraph $P$, equipped with a total termination order,
the \emph{coherent completion} of $P$ is the $(3,1)$\nbd-polygraph obtained from
$P$ by successive applications of Knuth-Bendix and Squier completion steps, as
follows. In this procedure, one considers each critical branching~$(\phi,\psi)$
of~$P$ and performs the following operations:
\begin{itemize}
\item if the branching is confluent, the procedure adds a 3-generator
  \[
    A:\phi\comp1\phi'\TO\psi\comp1\psi'
  \]
  to the polygraph (if such a generator is not already present):
  \[
    \xymatrix @C=3em @R=1em {
      & {v}
      \ar@2@/^/ [dr] ^-{\phi'}
      \ar@3{.>} []!<0pt,-15pt>;[dd]!<0pt,15pt> ^-{A}
      \\
      {u}
      \ar@2@/^/ [ur] ^-{\phi}
      \ar@2@/_/ [dr] _-{\psi}
      && {\nf{v}=\nf{w}}\pbox,
      \\
      & {w}
      \ar@2@/_/ [ur] _-{\psi'}
    }
  \]
\item if the branching is not confluent, the procedure coherently adds a
  $2$\nbd-ge\-ne\-ra\-tor
  \[
    \alpha:\nf v\To\nf w\quad\text{if $\rep{v}>\rep{w}$}
    \qqqtor
    \alpha:\nf w\To\nf v\quad\text{if $\rep{w}>\rep{v}$}
  \]
  and a $3$-generator
  \[
    A:\phi\comp1\phi'\TO\psi\comp1\psi'
  \]
  to the polygraph:
  \[
    \xymatrix @C=3em @R=1em {
      & {v}
      \ar@2 [r] ^-{\phi'}
      \ar@3{.>} []!<10pt,-15pt>;[dd]!<10pt,15pt> ^-{A}
      & {\nf{v}}
      \ar@2{<.>} [dd] ^-{\alpha}
      \\
      {u}
      \ar@2@/^/ [ur] ^-{\phi}
      \ar@2@/_/ [dr] _-{\psi}
      \\
      & {w}
      \ar@2 [r] _-{\psi'}
      & {\nf{w}}\pbox.
    }
  \]
\end{itemize}
In the second case, the procedure adds a new 2-generator $\alpha$, which can in
turn create new critical branchings, which have to be inspected by the
procedure. For this reason, like in the usual Knuth-Bendix procedure, the
process is not guaranteed to terminate. In this situation, this defines an
increasing sequence of $(3,1)$-polygraphs, whose inductive limit is a
potentially infinite $(3,1)$-polygraph.

As a consequence of \cref{thm:SquierHomotopical}, the $(3,1)$-polygraph
constructed using this procedure satisfies the following property~\cite[Theorem
2.2.5]{GaussentGuiraudMalbos15}:

\begin{theorem}
  Let $P$ be a terminating 2-polygraph. Any coherent completion of $P$ is a
  coherent convergent presentation of the category~$\cl{P}$.
\end{theorem}

\subsection{Generic homotopical reduction}
\label{Subsection:GenericHomotopicalReduction}
In order to reduce the size of the $(3,1)$-polygraph obtained by a coherent
completion of a terminating $2$\nbd-poly\-graph, one can identify generators which can
be collapsed, and thus be removed without changing the presented category nor
the coherence of the category (see \cref{prop:collapse-Tietze}). We formalize
here the process of collapsing multiple such generators at once.

\index{collapsible!part}
A \emph{collapsible part} of a $(3,1)$-poly\-graph $P$ is a family~$X$ of its generators that we can collapse together, in the sense introduced in
\cref{sec:3collapsible}. Explicitly, it consists in a triple
$X=(X_2,X_3,X_4)$ made of a family $X_2$ of $2$\nbd-gene\-rators of~$P$, a
family $X_3$ of $3$\nbd-gene\-rators of $P$ and a family~$X_4$ of $3$-spheres of the free $(2,1)$\nbd-cate\-gory~$\freegpd{P}$, such that the following conditions are satisfied:
\begin{itemize}
\item the elements of $X_2$, $X_3$, and $X_4$ are collapsible, potentially up to
  a Nielsen transformation (see \cref{rem:collapse-Nielsen}),
\item no $2$-generator of~$X_2$ is the target of a $3$-generator in $X_3$,
\item no $3$\nbd-gene\-rator of~$X_3$ is the target of a $3$-sphere in~$X_4$,
\item the following relations are well-founded:
  \begin{itemize}
  \item the relation~$<_1$ on~$P_1$ such that $b<_1a$ when there exists a
    $2$-generator $\alpha:u\To a$ in $X_2$ such that~$b$ occurs in~$u$,
  \item the relation $<_2$ on~$P_2$ such that $\beta<_2\alpha$ when there exists
    a $3$\nbd-gene\-rator $A:\phi\TO\alpha$ in $X_2$ such that~$\beta$ occurs
    in~$\phi$,
  \item the relation $<_3$ on~$P_3$ such that $B<_3A$ when there exists a
    $3$-sphere $\Phi:F\TOO A$ in $X_3$ such that~$B$ occurs in~$F$.
  \end{itemize}
\end{itemize}

\index{homotopical!reduction}
\index{reduction!homotopical}
Given such a collapsible part~$X$, one can define a $(3,1)$-polygraph~$P/X$,
obtained by successively collapsing all the elements of~$X$, which is called the
\emph{homotopical reduction of $P$ with respect to $X$}. By construction, the
polygraph~$P/X$ is Tietze equivalent to $P$.

\subsection{Generating triple confluences}
\label{Subsection:TripleConfluences}
The coherent elimination of $3$-gener\-ators of a $(3,1)$-polygraph~$P$ by homotopical
reduction requires a collapsible set of $3$\nbd-spheres of $\freegpd{P}$. When
$P$ is convergent and coherent, its triple critical branchings provide a
convenient way to build such a set.

\index{triple branching}
\index{branching!triple}
\index{critical!branching!triple}
A \emph{local triple branching} is a triple $(\phi,\chi,\psi)$ of $2$-cells
which are rewriting steps with a common source. Similarly to local branchings,
local triple branchings are classified into three families:
\begin{itemize}
\item \emph{trivial} triple branchings have two of the $2$-cells equal,
\item \emph{orthogonal} triple branchings have at least one of their $2$-cells
  that form an orthogonal branching with the other two,
\item \emph{overlapping} triple branchings are the remaining local triple
  branchings.
\end{itemize}
Local triple branchings are ordered by inclusion of their sources, similarly to
branchings.
A \emph{critical} triple branching is an overlapping triple branching that is
minimal for this inclusion. For a reduced $2$-polygraph, such a triple branching
can have two different shapes, where~$\phi$,~$\psi$, and~$\chi$ are $2$-generators:
\[
  \xymatrix@C=4ex{
    \strut
    \ar [r] _-{u_1}
    \ar @/^6ex/ [rrr] _-{}="tgt1"
    & \strut
    \ar [r] |-{u_2} ^-{}="src1"
    \ar @/_6ex/ [rrr] ^-{}="tgt2"
    & \strut
    \ar[r] |-{u_3} _-{} ="src2"
    \ar @/^6ex/ [rrr] _-{}="tgt3"
    & \strut
    \ar[r] |-{u_4} ^-{} ="src3"
    & \strut
    \ar[r] _-{v}
    & \strut
    \ar@2 "src1"!<0pt,5pt>;"tgt1"!<0pt,-5pt> ^-{\phi}
    \ar@2 "src2"!<0pt,-5pt>;"tgt2"!<0pt,5pt> ^-{\psi}
    \ar@2 "src3"!<0pt,5pt>;"tgt3"!<0pt,-5pt> ^-{\chi}
  }
  \qquad\text{or}\qquad
  \xymatrix@C=4ex{
    \strut
    \ar [r] _-{u_1}
    \ar @/^6ex/ [rr] _-{}="tgt1"
    & 
    \ar [r] |-{u_2}
    \ar @/_6ex/ [rrr] ^{}="tgt2"
    & \strut
    \ar[r] ^-{u_3} _{} ="src2"
    & \strut
    \ar[r] |-{u_4}
    \ar @/^6ex/ [rr] _-{}="tgt3"
    & 
    \ar[r] _-{v}
    & \strut\pbox.
    \ar@2 "1,2"!<0pt,5pt>;"tgt1"!<0pt,-5pt> ^-{\phi}
    \ar@2 "src2"!<0pt,-5pt>;"tgt2"!<0pt,5pt> ^-{\psi}
    \ar@2 "1,5"!<0pt,5pt>;"tgt3"!<0pt,-5pt> ^-{\chi}
  }
\]
When the polygraph is not reduced, the other possible type of critical
branchings, with an inclusion of one source into the other one, generates
several other possibilities.

If $P$ is a coherent and convergent $(3,1)$-polygraph, a \emph{generating
  triple confluence of $P$} is a $3$-sphere
\[
  \vcenter{
    \xymatrix@C=3.7ex @R=3ex{
      & {\sm v} 
      \ar@2 @/^/ [rr] ^-{\phi'_1}
      \ar@{} [dr] |-{A}
      && {\sm x'}
      \ar@2 @/^/ [dr] ^-{\psi''}
      \\
      {\sm u}
      \ar@2 @/^/ [ur] ^-{\phi}
      \ar@2 [rr] |-{\chi}
      \ar@2 @/_/ [dr] _-{\psi}
      && {\sm w} 
      \ar@2 [ur] |-{\chi'_1}
      \ar@2 [dr] |-{\chi'_2}
      \ar@{} [rr] |-{C'}
      && {\sm \rep{u}}
      \\
      & {\sm x}
      \ar@2 @/_/ [rr] _-{\psi'_2}
      \ar@{} [ur] |-{B}
      && {\sm v'}
      \ar@2 @/_/ [ur] _-{\phi''}
    }
  }
  \quad
  \vcenter{
    \xymatrix@C=2.6ex{
      \ar@4 [r] ^-*+{\Phi}_*+{\phantom\Phi}&
    }
  }
  \quad
  \vcenter{
    \xymatrix@C=3.7ex @R=3ex{
      & {\sm v}
      \ar@2 @/^/ [rr] ^-{\phi'_1}
      \ar@2 [dr] |-{\phi'_2}
      && {\sm x'}
      \ar@2 @/^/ [dr] ^-{\psi''}
      \\
      {\sm u}
      \ar@2 @/^/ [ur] ^-{\phi}
      \ar@{} [rr] |-{C}
      \ar@2 @/_/ [dr] _-{\psi}
      && {\sm w'}
      \ar@2 [rr] |-{\chi''}
      \ar@{} [ur] |-{B'}
      \ar@{} [dr] |-{A'}
      && {\sm \rep{u}}
      \\
      & {\sm x}
      \ar@2 [ur] |-{\psi'_1}
      \ar@2 @/_/ [rr] _-{\psi'_2}
      && {\sm v'}
      \ar@2 @/_/ [ur] _-{\phi''}
    }
  }
\]
where $(\phi,\chi,\psi)$ is a triple critical branching of $P$ and the $3$-cells are
generated by the generating confluence induced by the critical branchings.

\subsection{Coherent completion-reduction}
\label{Section:HomotopicalCompletionReduction}
\label{Subsection:HomotopicalCompletionReduction}
In practice, we apply homotopical reduction to a coherent
completion $Q$ of a terminating $2$\nbd-poly\-graph $P$. 
In such a situation, one can define a collapsible part~$X$ of~$Q$ whose elements are
\begin{itemize}
\item some of the generating triple confluences of $Q$,
\item the $3$-generators coherently adjoined with a $2$-generator by coherent completion
  to reach confluence,
\item some collapsible $2$-generators or $3$-generators already present in the initial
  presentation~$P$.
\end{itemize}
In practice, the collapsible triple confluences are chosen among those in which
some $3$-generator~$A$ occurs in the source or the target without
$1$-dimensional whiskers and occurs exactly once. Similarly, the collapsible
$3$-generators are chosen among those where a $2$-generator $\alpha$ occurs in
the source or the target without whiskers and occurs exactly once. Finally, the
collapsible $2$-generators are chosen among those of the form $\alpha:u\To a$ or
$\alpha:a\To u$ where $a$ is a generator which does not occur in~$u$. Moreover,
one should check that the conditions of \cref{Subsection:TripleConfluences} are
satisfied. In particular, one should be careful not to select too many such
generators in order for the well-foundedness conditions to be satisfied. An
illustration is given in \cref{ex:B3-KBS} below.

\index{coherent!completion-reduction}
\index{procedure!coherent!completion-reduction}
If $P$ is a terminating $2$-polygraph, the \emph{coherent completion-reduction
  of $P$} with respect to a collapsible part $X$ of its completion~$Q$ is the
$(3,1)$-polygraph the homotopical reduction~$Q/X$ of~$Q$ with respect to~$X$.

\begin{theorem}
  \label{Theorem:HomotopicalCompletionReduction}
  Let $P$ be a terminating 2-polygraph. A coherent completion-reduction of $P$
  is a coherent presentation of the category $\cl{P}$.
\end{theorem}

We refer to \cref{chap:SomeCoherentPresentations} for examples of coherent completion-reduction calculations in the algebraic situations of Artin, plactic, and Chinese monoids. We end this section with a simple example to illustrate the method.

\begin{example}
  \label{ex:B3-KBS}
  Consider the following presentation of the braid monoid~$B_3^+$, already
  encountered in~\cref{ex:kb-S3-r,sec:kapur-narendran}, see also
  \cref{sec:braid-mon}:
  \[
    P
    =
    \Pres\star{a,b,c}{aba\To bab,ba\To c}
  \]
  and equipped with the deglex order induced by $a>b>c$. The $2$-polygraph $P$ is
  terminating and its coherent completion is the $(3,1)$-polygraph:
  \[
    Q=\PRes\star{a,b,c}{\alpha,\beta,\gamma,\delta}{A,B,C,D},
  \]
  where $\alpha:ac\To cb$, $\beta:ba\To c$, $\gamma:bcb\To cc$, $\delta:bcc\To cca$ and $A$, $B$, $C$, and $D$ are the following $3$-generators, induced by completion
  of the critical branchings:
  \[
    \begin{tabular}{c@{\qquad\qquad}c}
      \xymatrix@R=1.5em@C=2em{
        & cc
        \\
        bac
        \ar@2@/^/ [ur] ^{\beta c} _(0.66){}="src"
        \ar@2@/_/ [dr] _{b\alpha} ^(0.66){}="tgt"
        \\
        & bcb
        \ar@2@/_/ [uu] _{\gamma}
        \ar@3 "src"!<0pt,-15pt>;"tgt"!<0pt,15pt> ^{A}
      }
      &
      \xymatrix@R=1.5em@C=2em{
        & cca
        \\
        bcba
        \ar@2@/^/ [ur] ^{\gamma a} _(0.66){}="src"
        \ar@2@/_/ [dr] _{bc\beta} ^(0.658){}="tgt"
        \\
        & bcc
        \ar@2@/_/ [uu] _{\delta}
        \ar@3 "src"!<0pt,-15pt>;"tgt"!<0pt,15pt> ^{B}
      }
    \\[2ex]
      \xymatrix@R=1.5em@C=1em{
        & cccb
        \ar@3 []!<0pt,-20pt>;[dd]!<0pt,20pt> ^{C}
        \\
        bcbcb
        \ar@2@/^/ [ur] ^{\gamma cb}
        \ar@2@/_/ [dr] _{bc\gamma}
        && ccac
        \ar@2@/_/ [ul] _{cc\alpha}
        \\
        & bccc
        \ar@2@/_/ [ur] _{\delta c}
      }
    &
      \xymatrix@R=1.5em@C=1em{
        & cccc
        \ar@3 []!<5pt,-20pt>;[dd]!<5pt,20pt> ^{D}
        & cccba
        \ar@2 [l] _{ccc\beta}
        \\
        bcbcc
        \ar@2@/^/ [ur] ^{\gamma cc}
        \ar@2@/_/ [dr] _{bc\delta}
        \\
        & bccca
        \ar@2 [r] _{\delta ca}
        & ccaca \; .
        \ar@2@/_/ [uu] _{cc\alpha a}
      }
    \end{tabular}
  \]
  The coherent presentation~$Q$ of~$B_3^+$ can be reduced using the collapsible
  part consisting of the following two generating triple confluences
  \begin{align*}
    \vcenter{
      \xymatrix@R=2.5em@C=1.5em{
        & {\sm ccac}
        \ar@2 [r] ^{cc\alpha}
        \ar@3 []!<-25pt,-15pt>;[d]!<-25pt,15pt> ^{Bc}
        & {\sm cccb}
        \\
        {\sm bcbac}
        \ar@2@/^3ex/ [ur] ^{\gamma ac} _(0.68){}="src1"
        \ar@2 [r]^{bc\beta c} ^(0.574){}="tgt1" _(0.574){}="src2"
        \ar@2@/_3ex/ [dr] _{bcb\alpha} ^(0.68){}="tgt2"
        & {\sm bccc}
        \ar@2@/_/ [u] _{\delta c}
        \ar@3 []!<-25pt,-15pt>;[d]!<-25pt,15pt> ^{bcA}
        \\
        & {\sm bcbcb}
        \ar@2@/_/ [u] _{bc\gamma}
      }
    }
    &
    \vcenter{\xymatrix@C=1.5em{\strut \ar@4 [r] ^-*+{\Phi}_-*{\phantom\Phi} &\strut}}
    \vcenter{
      \xymatrix@R=1.5em@C=0.5em{
        & {\sm ccac}
        \ar@2@/^/ [dr] ^{cc\alpha}
        \\
        {\sm bcbac}
        \ar@2@/^/ [ur] ^{\gamma ac}
        \ar@2@/_/ [dr] _{bcb\alpha}
        \ar@{} [rr] |{\rotatebox{90}{=}}
        && {\sm cccb}
        \ar@3 []!<0pt,-20pt>;[dd]!<0pt,20pt> ^{C}
        \\
        & {\sm bcbcb}
        \ar@2 [ur] |{\gamma cb}
        \ar@2@/_/ [dr] _{bc\gamma}
        && {\sm ccac}
        \ar@2@/_/ [ul] _{cc\alpha}
        \\
        && {\sm bccc}
        \ar@2@/_/ [ur] _{\delta c}
      }
    }
    \\
    \shortintertext{and}
    \vcenter{
      \xymatrix@R=2.5em@C=1.5em{
        & {\sm cccba}
        \ar@2 [r] ^{ccc\beta}
        \ar@3 []!<-10pt,-15pt>;[d]!<-10pt,15pt> _{Ca}
        & {\sm cccc}
        \\
        {\sm bcbcba}
        \ar@2@/^3ex/ [ur] ^{\gamma cba}
        \ar@2 [r]_{bc\gamma a} _(0.582){}="src"
        \ar@2@/_3ex/ [dr] _{bcbc\beta} ^(0.68){}="tgt"
        & {\sm bccca}
        \ar@2 [r] _{\delta ca}
        \ar@3 []!<-25pt,-15pt>;[d]!<-25pt,15pt> ^{bcB}
        & {\sm ccaca}
        \ar@2 [ul]_{cc\alpha a}
        \\
        & {\sm bcbcc}
        \ar@2@/_/ [u] _{bc\delta}
      }
    }
    &
    \vcenter{\xymatrix@C=1.5em{\strut \ar@4 [r] ^-*+{\Psi}_-*{\phantom\Psi} &\strut}}
    \vcenter{
      \xymatrix@R=1.5em@C=0.5em{
        & {\sm cccba}
        \ar@2@/^/ [dr] ^{ccc\beta}
        \\
        {\sm  bcbcba}
        \ar@2@/^/ [ur] ^{\gamma cba}
        \ar@2@/_/ [dr] _{bcbc\beta}
        \ar@{} [rr] |{\shortparallel}
        && {\sm cccc}
        \ar@3 []!<10pt,-20pt>;[dd]!<10pt,20pt> ^{D}
        & {\sm cccba}
        \ar@2 [l] _{ccc\beta}
        \\
        & {\sm bcbcc}
        \ar@2 [ur] _{\gamma cc}
        \ar@2@/_/ [dr] _{bc\delta}
        \\
        && {\sm bccca}
        \ar@2 [r] _{\delta ca}
        & {\sm ccaca}
        \ar@2 [uu] |{cc\alpha a}
      }
    }
  \end{align*}
  together with the $3$-generators $A$ and $B$ coherently adjoined with the $2$-generators
  $\gamma$ and~$\delta$ during coherent completion and the $2$-generator
  $\beta:ba\dfl c$ that defines the redundant generator $c$. The generators $\beta$,
  $A$, $B$, $\Phi$, and~$\Psi$ are collapsible up to a Nielsen transformation,
  with respective redundant generators $c$, $\gamma$, $\delta$, $C$, and $D$. We
  conclude that $X$ is collapsible since the relations $<_1$, $<_2$, and $<_3$
  are respectively included in the following well-founded total orders, and are
  thus well-founded:
  \begin{align*}
    c > b > a
    &&
    \delta > \gamma > \beta > \alpha
    &&
    D > C > B > A
    \pbox.
  \end{align*}
  It follows that the homotopical reduction of the coherent presentation $Q$
  with respect to this collapsible part is the following coherent
  $(3,1)$-polygraph:
  \[
    R=\PRes\star{b,a}{aba\To bab}{}
    \pbox.
  \]
  By \cref{Theorem:HomotopicalCompletionReduction}, we recover that the monoid
  $B_3^+$ admits a coherent presentation made of Artin's presentation and no
  $3$-generator. This example is generalized in~\cite{GaussentGuiraudMalbos15} where
  coherent presentations of Artin monoids are constructed, see also
  \cref{chap:SomeCoherentPresentations}.
\end{example}

\section{Coherent Presentations of Associative Algebras}
\label{sec:coh-alg}
In this section, we define the notion of coherent presentation of an
associative algebra by extending the notion of presentation of an algebra
introduced in \cref{chap:2linear}.

\subsection{Extended presentations}
\index{extended presentation}
\index{presentation!extended}
\index{cellular extension!of an algebra}
\index{linear!2-polygraph@$2$-polygraph}
\index{2-polygraph@$2$-polygraph!linear}
A \emph{cellular extension} of a $1$-algebra~$A$, with the notations of
\cref{sec:1alg}, is a set~$X$ equipped with functions $\src1,\tgt1:X\to A_1$
such that $\src0\circ\src1=\src0\circ\tgt1$. A \emph{linear $2$-polygraph}
$(P,P_2)$ consists of a linear $1$\nbd-poly\-graph together with a cellular
extension~$P_2$ of the free $1$-algebra~$\lin P$ generated by~$P$. An
\emph{extended presentation} of an algebra~$A$ is a linear $2$-polygraph whose
underlying linear $1$-polygraph presents~$A$. A linear $2$-polygraph is
\emph{left-monomial} when the underlying linear $1$-polygraph is, in the sense
of \cref{sec:1-left-monomial}.

\index{2-algebra}
\index{algebra!2-}
A \emph{$2$-algebra} is an internal $2$-category in the category~$\Alg$ of
algebras. Note that contrarily to the set-theoretic case, we will not bother
about distinguishing whether we take cells to be invertible or not: it can be
shown that the notion of $2$-algebra coincides with the notion of internal
$2$-groupoid in the category of algebras,
see~\cite{GuiraudHoffbeckMalbos19}. Any linear $2$-polygraph~$P$ freely
generates a linear $2$-algebra that we denote as $\lin P$, and whose algebra of
$2$-cells is in particular written $\lin P_2$.

\subsection{Coherent confluence and convergence}
Let~$P$ be a left-monomial linear $2$-polygraph. A branching~$(\phi,\psi)$
of~$P$ is \emph{coherently confluent} if there exist positive $1$-cells~$\phi'$
and~$\psi'$ in~$\lin{P}_1$ and a $2$-cell~$F$ in~$\lin{P}_2$ as in
\[
  \xymatrix@!C=1ex@!R=1ex{
    &\ar[dl]_{\phi}p\ar[dr]^{\psi}&\\
    q_1\ar@{.>}[dr]_{\phi'}&\underset{\phantom A}{\overset A\To}&\ar@{.>}[dl]^{\psi'}q_2\pbox.\\
    &r
  }
\]
If~$p$ is a $0$-cell of~$\lin{P}_0$, say that~$P$ is \emph{coherently confluent}
(\resp \emph{locally coherently confluent}, \resp \emph{critically coherently
  confluent}) at~$p$ if every branching (\resp local branching, \resp critical
branching) of~$P$ of source~$p$ is coherently confluent. Say that~$P$ is
\emph{coherently confluent} (\resp \emph{locally coherently confluent}, \resp
\emph{critically coherently confluent}) if it is so at every $0$-cell
of~$\lin{P}_0$, and that~$P$ is \emph{coherently convergent} if it is
terminating and coherently confluent.

\begin{lemma}
  \label{L:Confluence2Cell}
  Let~$P$ be a left-monomial linear 2-polygraph with a fixed 0\nbd-cell~$p$, and
  suppose that~$P$ is coherently confluent at every 0-cell~$q$ such that
  $p\overset*\to q$. Let~$\phi$ be a 1-cell of~$\lin{P}_1$ which admits a
  decomposition
  \[
    \xymatrix@C=3ex{
      p_0\ar[r]^\phi&p_k
    }
    \quad=\quad
    \xymatrix@C=3ex{
      p_0\ar[r]^{\phi_1}&p_1\ar[r]^{\phi_2}&\cdots\ar[r]^{\phi_k}&p_k
    }
  \]
  into 1-cells~$\phi_i$ of size~1. If~$p\overset*\to p_i$ holds for
  every~$0\leq i<k$, then there exist positive 1-cells~$\phi'$ and~$\psi$ in~$\lin{P}_1$
  and a 2-cell~$F$ in~$\lin{P}_2$ as in
  \[
    \xymatrix @R=1em {
      & p_k
      \ar@/^/ [dr] ^-{\phi'}
      \ar@2 []!<0pt,-12pt>;[d]!<0pt,3pt> ^-*+{F}
      \\	
      p_0 
      \ar@/^/ [ur] ^-{\phi}
      \ar@/_/ [rr] _-{\psi}
      && p'
    }
  \]
\end{lemma}
\begin{proof}
  Proceed by induction on~$k$. If~$k=0$, then~$\phi$ is an identity, so
  taking~$\phi'=\psi=\unit{p_0}$ and~$F=\unit\phi$ proves the result. Otherwise,
  we construct
  \[
    \xymatrix @R=1.5em {
      && p_k
      \ar@/^/ [dr] ^-{\phi'}
      \ar@2 []!<0pt,-10pt>;[d]!<0pt,10pt> ^-*+{F}
      \\
      & p_1 
      \ar@/^/ [ur] ^(0.4){\phi_2\comp0\cdots\comp0 \phi_k}
      \ar [rr] |-*+{\psi_2}
      \ar [dr] |-*+{\phi'_1}
      \ar@{} []!<-5pt,0pt>;[d]!<-5pt,0pt> |(0.66){\sm =}
      && q_2
      \ar@/^/ [dr] ^-{\phi'_2}
      \ar@2 []!<-15pt,-12.5pt>;[d]!<-15pt,7.5pt> ^-*+{G}
      \\	
      p_0
      \ar@/^/ [ur] ^-{\phi_1}
      \ar@/_/ [rr] _-{\psi_1}
      && q_1 
      \ar@/_/ [rr] _-{\psi'_2}
      && p'\pbox.
    }
  \]
  Apply \cref{L:FactElem1Cell} to the $1$-cell~$\phi_1$ of size~$1$ to get the
  positive $1$-cells~$\phi_1'$ and~$\psi_1$ such that
  $\phi_1=\phi'_1\comp0 \psi_1^-$. We have~$p\overset*\to p_i$ for
  every~$1\leq i<k$, so the induction hypothesis applies to
  $\phi_2\comp0\cdots\comp0\phi_k$, providing the positive $1$-cells~$\phi'$
  and~$\psi_2$, and the $2$-cell~$F$. Then, consider the
  branching~$(\phi_1',\psi_2)$, whose source~$p_1$
  satisfies~$p\overset *\to p_1$: by hypothesis, this branching is coherently
  confluent, giving the positive $1$-cells~$\phi_2'$ and~$\psi_2'$, and the
  $2$-cell~$G$.
\end{proof}

The following result is a formulation of coherent Newman's lemma for linear
polygraphs. The proof is the same as in the set-theoretical case given by
\cref{P:HNewman-2}.

\begin{proposition}
  \label{P:HNewmanLinear}
  Let~$P$ be a terminating left-monomial linear 2-polygraph. If~$P$ is locally
  coherently confluent then it is coherently confluent.
\end{proposition}

The following result is a formulation of the coherent critical branchings lemma,
\cref{L:HCritical}, for linear polygraphs. Due to the linearity of contexts, the
termination is necessary and the proof differs from the set-theoretical case, as
already explained in \cref{rem:lin-cb-term}.

\begin{lemma}
  \label{L:HCriticalLinear}
  Suppose given a terminating left-monomial linear 2-poly\-graph~$P$. If~$P$ is
  critically coherently confluent, then~$P$ is locally coherently confluent.
\end{lemma}
\begin{proof}
  We proceed by noetherian induction on the sources of the local branchings to
  prove that~$P$ is locally coherently confluent at every $0$-cell
  of~$\lin{P}_0$.
  We note that a reduced $0$-cell cannot be the source of a local branching,
  so~$P$ is locally coherently confluent at reduced $0$-cells.
  Now, fix a non-reduced $0$-cell~$p$ of~$\lin{P}_0$, and assume that~$P$ is
  locally coherently confluent at every $0$-cell $q$ with $p\overset*\to
  q$. With a termination-based argument similar to that of \cref{P:HNewman-2},
  we deduce that~$P$ is coherently confluent at every~$q$. Then we proceed by
  case analysis on the type of the local branchings, noting that an aspherical
  branching $\lambda(\phi,\phi)+b$ is always coherently confluent.

  For an additive branching, we construct
  \[
    \xymatrix @R=2.5em @C=1em {
      & {\lambda p + \mu v + r}
      \ar@/^/ [rr] ^-{\phi'_1}
      \ar@{.>} [dr] |-*+{\lambda p + \mu\psi + r}
      && p'
      \ar@/^/ [drrr] ^-{\phi'_2}
      \ar@2 []!<0pt,-30pt>;[dd]!<0pt,30pt> ^-*+{F}
      \\
      {\lambda u + \mu v + r}
      \ar@/^2ex/ [ur] ^-*+{\lambda \phi + \mu v + r}
      \ar@/_2ex/ [dr] _-*+{\lambda u + \mu \psi + r}
      \ar@{} [rr] |(0.45){\sm =}
      && {\lambda p + \mu q + r}
      \ar [ur] |-*+{\phi'}
      \ar [dr] |-*+{\psi'}
      \ar@{} [u] |(0.7){\sm =}
      \ar@{} [d] |(0.7){\sm =}
      &&&& r'\pbox.
      \\
      & {\lambda u + \mu q + r}
      \ar@{.>} [ur] |-*+{\lambda \phi + \mu q + r}
      \ar@/_/ [rr] _-{\psi'_1}
      && {q'}
      \ar@/_/ [urrr] _-{\psi'_2}
    }
  \]
  By linearity of the $0$-composition, we have
  \[
    (\lambda\phi + \mu v + r) \comp0 (\lambda p + \mu \psi + r) 
    = \lambda \phi + \mu \psi + r
    = (\lambda u + \mu \psi + r) \comp0 (\lambda \phi + \mu q + r).
  \]
  Note that the dotted $1$-cells $\lambda p + \mu \psi + r$ and
  $\lambda\phi + \mu q + r$ are not positive in general since~$u$ can be
  in~$\supp{q}$ or~$v$ in~$\supp{p}$. However, those $1$-cells are of size~$1$,
  and \cref{L:FactElem1Cell} applies to both of them to give positive
  $1$-cells~$\phi'_1$, $\psi'_1$, $\phi'$, and~$\psi'$ that satisfy
  \begin{align*}
    \phi'_1 &= (\lambda p + \mu \psi + r) \comp0 \phi'
    &
    \psi'_1 &= (\lambda \phi + \mu q + r) \comp0 \psi'
    \pbox.
  \end{align*}
  Now, $u\overset*\to p$, $v\overset*\to q$, $\lambda\neq 0$ and~$\mu\neq 0$
  imply $\lambda u + \mu v + r \overset*\to \lambda p + \mu q + r$. Thus, the
  branching $(\phi',\psi')$ is coherently confluent by hypothesis, yielding the
  positive $1$-cells~$\phi'_2$ and~$\psi'_2$ and the $2$-cell~$F$.

  Next, in the case of an orthogonal branching, we construct
  \[
    \xymatrix @R=2.5em @C=1.5em {
      & \lambda pv + r
      \ar@/^/ [rr] ^-{\phi'_1}
      \ar@{.>} [dr] |-*+{\lambda p\psi + r}
      && p'
      \ar@/^/ [drr] ^-{\phi'_2}
      \ar@2 []!<0pt,-30pt>;[dd]!<0pt,30pt> ^-*+{H}
      \\
      \lambda uv + r
      \ar@/^2ex/ [ur] ^-*+{\lambda \phi v + r}
      \ar@/_2ex/ [dr] _-*+{\lambda u\psi + r}
      \ar@{} [rr] |(0.45){\sm =}
      && \lambda pq + r
      \ar [ur] |-*+{\phi'}
      \ar [dr] |-*+{\psi'}
      \ar@2 [u]!<0pt,-5pt>;[]!<0pt,25pt> ^-*+{F^-}
      \ar@2 []!<0pt,-25pt>;[d]!<0pt,5pt> ^-*+{G}
      &&& d\pbox.
      \\
      & \lambda uq + r
      \ar@{.>} [ur] |-*+{\lambda \phi q + r}
      \ar@/_/ [rr] _-{\psi'_1}
      &&	q'
      \ar@/_/ [urr] _-{\psi'_2}
    }
  \]
  Use the linearity of the $0$-composition to obtain
  \[
    (\lambda \phi v + r) \comp0 (\lambda p\psi + r) = \lambda\phi\psi + r = (\lambda u\psi + r) \comp0 (\lambda \phi q + r)
    \pbox.
  \]
  Again, the dotted $1$-cells $\lambda \phi q + r$ and $\lambda p\psi + r$ are
  not positive in general: this is the case, for example, if either
  $\supp{uq}\cap\supp{r}$ or $\supp{pv}\cap\supp{r}$ is not empty. Let
  $p = \sum_{i=1}^k \mu_i u_i$ be the canonical decomposition of~$p$. By
  linearity of the $0$-composition, the $1$-cell $\lambda p\psi + r$ admits the
  following decomposition in $1$-cells of size~$1$:
  \[
    \lambda p\psi + r = \psi_1\comp0\cdots\comp0\psi_k
  \]
  with
  \[
    \psi_j
    =
    \sum_{1\leq i < j} \lambda\mu_i u_i q
    + \lambda\mu_j u_j \psi
    + \sum_{j<i\leq k} \lambda\mu_i u_i v + r
    \pbox.
  \]
  We have~$u\overset*\to u_i$ for every~$i$, and~$v\overset*\to b$, giving
  $\lambda uv + v \overset*\to \tgt{}(\psi_j)$ for every~$j$.
  Hence~$\lambda p\psi + r$ is eligible for \cref{L:Confluence2Cell},
  yielding~$\phi'_1$, $\phi'$, and~$F$. The cells~$\psi'_1$, $\psi'$, and~$G$ are
  obtained similarly from $\lambda \phi q + r$. Finally,
  $\lambda uv + r \overset*\to \lambda pq + r$ implies, by induction hypothesis,
  that~$(\phi',\psi')$ is coherently confluent, giving~$\phi'_2$, $\psi'_2$
  and~$H$.

  Finally, for an overlapping branching $(\lambda \phi + r, \lambda \psi + r)$, we construct
  \[
    \xymatrix @R=2.5em @C=1.5em {
      & \lambda p + r
      \ar@/^/ [rr] ^-{\phi'_1}
      \ar@{.>} [dr] |-*+{\lambda \phi' + r}
      \ar@2 []!<0pt,-35pt>;[dd]!<0pt,35pt> ^-*+{F}
      && p'
      \ar@/^/ [drr] ^-{\phi'_2}
      \ar@2 []!<0pt,-35pt>;[dd]!<0pt,35pt> ^-*+{I}
      \\
      \lambda u + r
      \ar@/^2ex/ [ur] ^-*+{\lambda \phi + r}
      \ar@/_2ex/ [dr] _-*+{\lambda \psi + r}
      && \lambda s + r
      \ar [ur] |-*+{\phi''}
      \ar [dr] |-*+{\psi''}
      \ar@2 [u]!<0pt,-5pt>;[]!<0pt,25pt> ^-*+{G}
      \ar@2 []!<0pt,-25pt>;[d]!<0pt,5pt> ^-*+{H}
      &&& t\pbox.
      \\
      & \lambda q + r
      \ar@{.>} [ur] |-*+{\lambda \psi' + r}
      \ar@/_/ [rr] _-{\psi'_1}
      && q'
      \ar@/_/ [urr] _-{\psi'_2}
    }
  \]
  Consider the unique decomposition $(\phi,\psi)=v(\phi_0,\psi_0)w$,
  with~$(\phi_0,\psi_0)$ critical. Since~$(\phi_0,\psi_0)$ is coherently
  confluent by hypothesis, one obtains
  \[
    \xymatrix @R=0.5em {
      & p_0
      \ar@/^/ [dr] ^-{\phi_0'}
      \ar@2 []!<0pt,-12.5pt>;[dd]!<0pt,12.5pt> ^-*+{F_0}
      \\
      u_0
      \ar@/^/ [ur] ^-{\phi_0}
      \ar@/_/ [dr] _-{\psi_0}
      && r_0\pbox.
      \\
      & q_0
      \ar@/_1.4ex/ [ur] _-{\psi_0'}
    }
  \]
  Define the positive $1$-cells $\phi'=v\phi'_0 w$ and $\psi'=v\psi'_0 w$, and
  the $2$-cell~$F=v F_0 w$. As previously, the dotted $1$-cells are not positive
  in general, if~$\supp{c}$ intersects~$\supp{p}$ or~$\supp{q}$, for
  example. However, the $1$-cell~$\phi'$ is positive, so that it is a
  $0$-composite $\phi' = \chi_1\comp0\cdots\comp0 \chi_k$ of rewriting steps. As
  a consequence, we have the chain of reductions
  \[
    u \overset*\to p = \src{}(\chi_1) \overset*\to \cdots \overset*\to \src{}(\chi_k) \overset*\to s
    \pbox.
  \]
  Since we have~$\lambda \neq 0$ and $u\notin\supp{r}$ by hypothesis, the
  inequality
  \[ \lambda u + r \overset*\to \lambda \src{}(\chi_i) + r \]
  holds for
  every~$i$, so that the following decomposition of the $1$-cell
  $\lambda\phi'+r$ satisfies the hypotheses of \cref{L:Confluence2Cell}:
  \[
    \lambda \phi' + r
    =
    \big(\lambda \chi_1 + r \big) \comp1 \cdots \comp1 \big(\lambda \chi_k + r\big)
    \pbox.
  \]
  This gives~$\phi'_1$, $\phi''$, and~$G$. Proceed similarly with the $1$-cell
  $\lambda \psi' + r$ to obtain~$\psi'_1$, $\psi''$, and~$H$. Finally, apply the
  induction hypothesis on~$(\phi'',\psi'')$, since $\lambda u + r \overset*\to \lambda s + r$, to
  get~$\phi'_2$, $\psi'_2$, and~$I$.
\end{proof}

\noindent
Given a terminating left-monomial linear $1$-polygraph~$P$, taking $P_2$ to be
the set of all $2$-spheres, critical coherent confluence (\resp local coherent
confluence) in the linear $2$-polygraph~$(P,P_2)$ is the same as critical
confluence (\resp local confluence) in~$P$. We thus deduce the critical
branching lemma for linear $1$\nbd-poly\-graphs, already announced in
\cref{lem:LinearCritical}, as a particular case.

With a proof similar to the one in the set-theoretical case, see
\cref{thm:SquierHomotopical}, we have the coherent Squier theorem for linear
polygraphs:

\begin{theorem}
  \index{Squier!theorem!for linear polygraphs}
  \label{T:SquierLinearCoh}
  Let~$P$ be a convergent left-monomial linear 1-polygraph and~$P_2$ be a
  cellular extension of~$\lin{P}_1$ that contains a 2-cell
  \[
    \xymatrix@!C=1ex@!R=1ex{
      &\ar[dl]_\phi p\ar[dr]^\psi&\\
      q\ar[dr]_{\phi'}&\underset{\phantom 1}{\overset A\To}&\ar[dl]^{\psi'}r\\
      &p'
    }
  \]
  for every critical branching~$(\phi,\psi)$ of~$P$, with~$\phi'$ and~$\psi'$
  positive 1-cells of~$\lin{P}_1$. Then the 2-polygraph~$(P,P_2)$ is
  coherent.
\end{theorem}

\begin{example}
  \label{X:PP05}
  We consider the quadratic algebra~$A$ presented by
  \[
    \pres{x,y,z}{x^2 + yz=0, x^2 + \lambda zy=0}
  \]
  where~$\lambda$ is a fixed scalar different from $0$ and $1$,
  from~\cite[Section~4.3]{PP05}. Put~$\mu=\lambda^{-1}$. The algebra~$A$ admits
  the presentation
  \[
    P = \pres {x, y, z}{\alpha:yz \to -x^2 , \beta:zy\to - \mu x^2}
    \pbox.
  \]
  The deglex order generated by $z>y>x$ satisfies $yz>x^2$ and $zy>x^2$, proving
  that~$P$ terminates. However, $P$ is not confluent. Indeed, it has two
  critical branchings:
  \[
    \vcenter{\xymatrix @R=0.5em {
        & -x^2y
        \\
        yzy
	\ar @/^/ [ur] ^-{\alpha y}
	\ar @/_/ [dr] _-{y\beta}
        \\
        & - \mu yx^2
      }}
    \qquad\qquad\text{and}\qquad\qquad
    \vcenter{\xymatrix @C=2em @R=0.5em {
        & - \mu x^2z
        \\
        zyz
	\ar @/^/ [ur] ^-{\beta z}
	\ar @/_/ [dr] _-{z\alpha}
        \\
        & -zx^2
      }}
  \]
  and neither of them is confluent because the monomials~$x^2 y$, $yx^2$,
  $x^2z$, and~$zx^2$ are reduced. The adjunction of the $1$-cells
  \[
    \gamma:yx^2\to\lambda x^2y
    \qquad\qquad\text{and}\qquad\qquad
    \delta:zx^2\to\mu x^2z
  \]
  gives a left-monomial linear $1$-polygraph
  \[
    Q = \pres {x, y, z}{
      \begin{array}{r@{\ :\ }r@{\ \to\ }l@{\ }r@{\ :\ }r@{\ \to\ }l}
        \alpha&yz&-x^2,&\gamma&yx^2&\lambda x^2y,\\
        \beta&zy&-\mu x^2,&\delta&zx^2&\mu x^2z
      \end{array}
    }
  \]
  that also presents~$A$, since~$\gamma$ and~$\delta$ induce relations that
  already hold in~$\cl{P}$, and that also terminates, because of $yx^2>x^2y$ and
  $zx^2>x^2z$. Moreover, each one of the four critical branchings of~$Q$ is
  confluent:
  \begin{align*}
    \vcenter{
      \xymatrix @R=0.75em @C=1.25em {
        & -x^2y
        \ar@2 []!<-20pt,-15pt>;[dd]!<-20pt,15pt> ^-*+{A}
        \\
        yzy
        \ar @/^/ [ur] ^-{\alpha y}
        \ar @/_/ [dr] _-{y\beta}
        \\
        & - \mu yx^2
        \ar @/_/ [uu] _-{-\mu\gamma}
      }
    }
    &&
    \vcenter{
      \xymatrix @R=0.75em @C=1.25em {
        & -\mu x^2z
        \ar@2 []!<-20pt,-15pt>;[dd]!<-20pt,15pt> ^-*+{B}
        \\
        zyz
        \ar @/^/ [ur] ^-{\beta z}
        \ar @/_/ [dr] _-{z\alpha}
        \\
        & -zx^2
        \ar @/_/ [uu] _-{-\delta}
      }
    }
    \\
    \vcenter{
      \xymatrix @R=0.75em @C=0.25em {
        & -x^4
        \ar@2 []!<-5pt,-15pt>;[dd]!<-5pt,15pt> ^-*+{C}
        \\
        yzx^2
        \ar@/^/ [ur] ^-{\alpha x^2}
        \ar@/_/ [dr] _-{y\delta}
        && x^2yz
        \ar@/_/ [ul] _-{x^2\alpha}
        \\
        & \mu yx^2z
        \ar@/_/ [ur] _-{\mu \gamma z}
      }
    }
    &&
    \vcenter{
      \xymatrix @R=0.75em @C=0.25em {
        & - \mu x^4
        \ar@2 []!<-5pt,-15pt>;[dd]!<-5pt,15pt> ^-*+{D}
        \\
        zyx^2
        \ar@/^/ [ur] ^-{\beta x^2}
        \ar@/_/ [dr] _-{z\gamma}
        && x^2zy\pbox.
        \ar@/_/ [ul] _-{x^2\beta}
        \\
        & \lambda zx^2y
        \ar@/_/ [ur] _-{\lambda \delta y}
      }
    }
  \end{align*}
  \Cref{T:SquierLinearCoh} implies that the $2$-polygraph
  \[
    \Pres{x,y,z}{\alpha,\beta,\gamma,\delta}{A,B,C,D}
  \]
  is a coherent presentation of~$A$.

  This coherent presentation can be reduced to a smaller one by a collapsing
  mechanism, similar to the one developed in \cref{Section:HomotopicalCompletionReduction}
  in the set-theoretic case, and hinted at in this example. First, some
  $2$-cells may be removed without breaking acyclicity because their boundary
  can also be filled by a composite of other $2$-cells. Here, the ``critical
  $3$-branchings'', where three rewriting steps overlap, reveal two relations
  between $2$-cells:
  \[
    \vcenter{
      \xymatrix @!C @C=3em {
        & -x^2yz
	\ar@/^3ex/ [dr] ^-{-x^2\alpha}
	\ar@2 []!<-40pt,-20pt>;[d]!<-40pt,15pt> ^-*+{Az}
	\ar@2 []!<32.5pt,-37.5pt>;[dd]!<32.5pt,37.5pt> ^-*+{-C}
        \\
        yzyz
	\ar@/^3ex/ [ur] ^-{\alpha yz}
	\ar [r] |-*+{y\beta z}
	\ar@/_3ex/ [dr] _-{yz\alpha}
        & -\mu yx^2z
	\ar [u] _-{-\mu\gamma z}
	\ar@2 []!<-40pt,-15pt>;[d]!<-40pt,20pt> ^-*+{yB}
        & x^4
        \\
        & -yzx^2
        \ar [u] _-{-y\delta}
	\ar@/_3ex/ [ur] _-{-\alpha x^2} 
      }
    }
    \quad\tfl\quad
    \vcenter{
      \xymatrix @!C @C=1.5em @R=1.5em {
        & -x^2yz
        \ar @/^2ex/ [dr] ^{-x^2\alpha}
        \\
        yzyz
	\ar @/^2ex/ [ur] ^{\alpha yz}
	\ar @/_2ex/ [dr] _{y z\alpha}
	\ar@{} [rr] |-{\sm =}
        && x^4\pbox,
        \\
        & -yzx^2
	\ar @/_2ex/ [ur] _{-\alpha x^2}
      }
    }
  \]
  \[
    \vcenter{
      \xymatrix @!C @C=3em{
        & -\mu x^2zy
	\ar@/^3ex/ [dr] ^-{-\mu x^2\beta}
	\ar@2 []!<-40pt,-20pt>;[d]!<-40pt,15pt> ^-*+{By}
	\ar@2 []!<22pt,-37.5pt>;[dd]!<22pt,37.5pt> ^-*+{-\mu D}
        \\
        zyzy
	\ar@/^3ex/ [ur] ^-{\beta zy}
	\ar [r] |-*+{z\alpha y}
	\ar@/_3ex/ [dr] _-{zy\beta}
        & -zx^2y
	\ar [u] _-{-\delta y}
	\ar@2 []!<-40pt,-15pt>;[d]!<-40pt,20pt> ^-*+{zA}
        & \mu^2 x^4
        \\
        & -\mu zyx^2
	\ar [u] _-{-\mu z\gamma}
	\ar@/_3ex/ [ur] _-{-\mu\beta x^2}
      }
    }
    \quad\tfl\quad
    \vcenter{
      \xymatrix @!C @C=1.5em @R=1.5em {
        &-\mu x^2zy
	\ar @/^2ex/ [dr] ^{-\mu x^2\beta}
        \\
        zyzy
	\ar @/^2ex/ [ur] ^{\beta zy}
	\ar @/_2ex/ [dr] _{zy \beta}
	\ar@{} [rr] |-{\sm =}
        && \mu^2 x^4\pbox.
        \\
        &-\mu zyx^2
	\ar @/_2ex/ [ur] _{-\mu\beta x^2}
      }
    }
  \]
  Since the boundaries of~$C$ and~$D$ can also be filled using~$A$ and~$B$ only,
  the $2$\nbd-poly\-graph $\Pres{x,y,z}{\alpha,\beta,\gamma,\delta}{A,B}$ is
  also a coherent presentation of~$A$. Next, the $1$-cells~$\gamma$ and~$\delta$
  are redundant, because the corresponding relations can be derived
  from~$\alpha$ and~$\beta$, as testified by the $2$-cells~$A$ and~$B$:
  removing~$\gamma$ with~$A$, and~$\delta$ with~$B$, proves that~$\lin{P}_1$
  admits an empty acyclic cellular extension, so that
  $\Pres{x,y,z}{\alpha,\beta}{}$ is actually a coherent presentation of~$A$.
\end{example}

\begin{example}[The standard coherent presentation]
  \index{presentation!standard!coherent}
  \label{P:StandardCoherentPresentation}
  Assume that $A=\kk\oplus A_+$ is an augmented algebra, and fix a linear
  basis~$\mathcal{B}$ of~$A_+$.  For~$u$ and~$v$ in~$\mathcal{B}$,
  write~$u\otimes v$ for the product of~$u$ and~$v$ in the free algebra
  over~$\mathcal{B}$, and~$uv$ for their product in~$A$.  Consider the linear
  $1$-polygraph~$\Std(\mathcal{B})_1$ whose $0$-cells are the elements
  of~$\mathcal{B}$, and with a $1$-cell
  \[
    u\otimes v \ofl{u\vert v} uv,
  \]
  for all~$u$ and~$v$ in~$\mathcal{B}$. Note that~$uv$ belongs to the free
  algebra over~$\mathcal{B}$ because~$A$ is augmented. By definition,
  $\Std(\mathcal{B})_1$ is a presentation of~$A$. Moreover,
  $\Std(\mathcal{B})_1$ terminates by a length argument: for all~$u$ and~$v$
  in~$\mathcal{B}$, the monomial $u\otimes v$ is a word of length~$2$ in the
  free monoid over~$\mathcal{B}$, while~$uv$ is a word of length~$1$. Finally,
  $\Std(\mathcal{B})_1$ has one critical branching
  $(u\vert v\otimes w, u\otimes v\vert w)$ for each triple~$(u,v,w)$ of elements
  of~$\mathcal{B}$, and this critical branching is confluent. Thus, extending
  $\Std(\mathcal{B})_1$ with a $2$-cell
  \[
    \xymatrix @R=1em @!C @C=0em {
      & uv\otimes w
      \ar@/^2ex/ [dr] ^-{uv\vert w}
      \ar@2 []!<-10pt, -15pt>;[dd]!<-10pt, 15pt> ^-*+{u\vert v\vert w}
      \\
      u\otimes v\otimes w
      \ar@/^2ex/ [ur] ^-{u\vert v \otimes w}
      \ar@/_2ex/ [dr] _-{u\otimes v\vert w}
      && uvw
      \\
      & u\otimes vw
      \ar@/_2ex/ [ur] _-{u\vert vw}
    }
  \]
  for each triple~$(u,v,w)$ of elements of~$\mathcal{B}$ produces, by
  \cref{T:SquierLinearCoh}, a coherent presentation of~$A$, denoted by
  $\Std(\mathcal{B})_2$. Note that the free $2$-algebra over
  $\Std(\mathcal{B})_2$ does not depend (up to isomorphism) on the choice of the
  basis~$\mathcal{B}$.

  This coherent presentation of~$A$ is extended in every dimension
  in~\cref{SS:StandardPolygraphicResolution} to obtain a polygraphic version of
  the standard resolution of an algebra. As in the previous example, the next
  dimension contains the $3$-cells generated by the ``critical $3$-branchings''
  of~$\Std_1(\mathcal{B})$: there is one such $3$-cell~$u\vert v\vert w\vert x$
  for each quadruple~$(u,v,w,x)$ of elements of~$\mathcal{B}$, with source
\[
  \xymatrix @R=3em @C=0ex @!C {
    & {\sm uv \otimes w\otimes x}
    \ar@/^/ [rr] ^-{uv\vert w\otimes x}
    \ar@2 []!<-5pt,-18pt>;[d]!<-5pt,18pt> ^-{\:u\vert v\vert w\otimes x}
    && {\sm uvw\otimes x}
    \ar@/^2ex/ [dr] ^-{uvw\vert x}
    \ar@2 []!<-10pt,-42.5pt>;[dd]!<-10pt,42.5pt> ^-{\:u\vert vw\vert x}
    \\
    {\sm u\otimes v\otimes w\otimes x}
    \ar@/^2ex/ [ur] ^-{u\vert v\otimes w\otimes x}
    \ar [rr] |-{u\otimes v\vert w\otimes x}
    \ar@/_2ex/ [dr] _-{u\otimes v\otimes w\vert x}
    & \strut
    \ar@2 []!<-5pt,-18pt>;[d]!<-5pt,18pt> ^-{\:u\otimes v\vert w\vert x}
    & {\sm u\otimes vw\otimes x}
    \ar [ur] |-{u\vert vw\otimes x}
    \ar [dr] |-{u\otimes vw\vert x}
    && {\sm uvwx}
    \\
    & {\sm u\otimes v\otimes wx}
    \ar@/_/ [rr] _-{u\otimes v\vert wx}
    && {\sm u\otimes vwx}
    \ar@/_2ex/ [ur] _-{u\vert vwx}
  }
\]
and target
\[
  \xymatrix @R=3em @C=0ex @!C {
    & {\sm uv\otimes w\otimes x}
    \ar@/^/ [rr] ^-{uv\vert w\otimes x}
    \ar [dr] |-{uv\otimes w\vert x}
    \ar@2 []!<-15pt,-42.5pt>;[dd]!<-15pt,42.5pt> ^-{\:1_{u\vert v\otimes w\vert x}}
    && {\sm uvw\otimes x}
    \ar@/^2ex/ [dr] ^-{uvw\vert x}
    \ar@2 []!<-35pt,-18pt>;[d]!<-35pt,18pt> ^-{\:uv\vert w\vert x}
    \\
    {\sm u\otimes v\otimes w\otimes x}
    \ar@/^2ex/ [ur] ^-{u\vert v\otimes w\otimes x}
    \ar@/_2ex/ [dr] _-{u\otimes v\otimes w\vert x}
    && {\sm uv\otimes wx}
    \ar [rr] |-{uv\vert wx}
    & \strut
    \ar@2 []!<-35pt,-18pt>;[d]!<-35pt,18pt> ^-{\:u\vert v\vert wx}
    & {\sm uvwx}\pbox.
    \\
    & {\sm u\otimes v\otimes wx}
    \ar [ur] |-{u\vert v\otimes wx}
    \ar@/_/ [rr] _-{u\otimes v\vert wx}
    && {\sm u\otimes vwx}
    \ar@/_2ex/ [ur] _-{u\vert vwx}
  }
\]
\end{example}


\chapter{Categories of Finite Derivation Type}
\label{chap:2-fdt}
\label{chap:2fdt}
In \cref{chap:2-Coherent}, we have seen a canonical and efficient way
to extend a convergent presentation of a category $C$ by a $2$-polygraph
$P$ into a coherent one. Precisely, the 
$3$-cells used in this extension procedure are in one-to-one correspondence with
the confluence diagrams of critical branchings in $P$
(\cref{thm:SquierHomotopical}). Now if $P$ is finite, so is the set of
its critical branchings and therefore the set of $3$-cells
generating coherence can be taken to be finite. In such a situation, we say that
the polygraph~$P$ has \emph{finite derivation type}, or \emph{FDT}. The
relevance of this concept lies in the following invariance property:
if a category $C$ admits a finite presentation $P$ having finite
derivation type, then \emph{all} finite presentations of $C$ also have FDT
(\cref{thm:TDFTietzeInvariant}). This invariance will prove essential
to show that some finitely presented categories do \emph{not} admit
convergent presentations.

This finiteness condition, introduced by
Squier~\cite{squier1994finiteness}, is of homotopical nature, and is, in
some sense, a refinement of the homological condition introduced
earlier in~\cite{squier1987word}. The latter will be discussed in the
next chapter.
Using these conditions, Squier managed to produce an explicit example
of a finitely presented monoid, with decidable word problem, but
having no finite convergent presentation. This provides a negative
answer to the question of \emph{universality of finite convergent
  rewriting} we raised in
\secr{universality}.
Let us finally emphasize the power of the FDT invariant:
by performing
computations on \emph{one} presentation of a monoid, we are able to deduce properties
of \emph{any} finite presentation of it!

The finiteness condition is introduced in \cref{sec:2fdt-def} and studied in the
case of convergent $2$-polygraphs in \cref{sec:conv-fdt}. In
\cref{sec:IdentitiesAmongRelations}, we define the notion of identities among
relations for $2$\nbd-poly\-graphs, generalizing those already known for presentations
of groups: such identities are described by $2$-spheres of the free
$(2,1)$-category on the $2$-polygraph.

\section{Finite Derivation Type}
\index{finite derivation type!2-polygraph}
\index{finite derivation type!1-category}
\index{FDT (finite derivation type)}
\label{sec:2fdt-def}
A $2$-polygraph $P$ has \emph{finite derivation type}, or~\emph{FDT} for short,
if it is finite and if the $(2,1)$-category $\freegpd{P}$ admits a finite
acyclic cellular extension. Otherwise said, there is a finite coherent
$(3,1)$-polygraph~$Q$ of which~$P$ is the underlying $2$-polygraph. A category~$C$ has \emph{finite derivation type} if it admits a
finite coherent presentation.

\subsection{Tietze invariance of the FDT property}
\index{Tietze!invariance}
Recall from \cref{sec:2-tietze} that two $2$-polygraphs are \emph{Tietze
  equivalent} when they present isomorphic categories. We say that a property
$\mathcal{P}$ on $2$-polygraphs is \emph{Tietze invariant} when for every Tietze
equivalent $2$-polygraphs $P$ and~$Q$, the polygraph~$P$ satisfies the
property~$\mathcal{P}$ if and only if the polygraph~$Q$ does.

Given two Tietze equivalent $2$-polygraphs $P$ and $Q$ whose sets $P_2$ and
$Q_2$ of $2$-generators are finite, consider a finite acyclic cellular extension
$X$ of the free $(2,1)$\nbd-cate\-gory~$\freegpd{P}$. By
\thmr{HomotopyBasesTransfer}, the cellular extension $X$ transfers to a finite
cellular extension of the free $(2,1)$-category
$\freegpd{Q}$.  We may therefore state the following invariance
result, first proved by Squier for
monoids~\cite[Theorem~4.3]{squier1994finiteness} and revisited
in~\cite[Theorem~4.2.3]{GuiraudMalbos18} in polygraphic terms.

\begin{theorem}
  \label{thm:TDFTietzeInvariant}
  Let $P$ and $Q$ be two Tietze equivalent 2-polygraphs such that~$P_2$ and
  $Q_2$ are finite. Then~$P$ has finite derivation type if and only if $Q$ has
  finite derivation type.
\end{theorem}

\noindent
This result shows that the property for a category~$C$ of having finite
derivation type does not depend on the presentation, provided that it is
finite.

The following result will help prove that a presentation admits no finite
acyclic cellular extension, \ie that the presented category does not have finite derivation type.

\begin{proposition}
  \label{Proposition:FiniteBasisExtracted}
  Let $P$ be a 2-polygraph and let $X$ be an acyclic extension of the free
  $(2,1)$-category $\freegpd{P}$. If $\freegpd{P}$ admits a finite acyclic
  cellular extension, then there exists a finite subset of $X$ that is an
  acyclic cellular extension of $\freegpd{P}$.
\end{proposition}
\begin{proof}
  Suppose that $\freegpd{P}$ admits a finite acyclic cellular extension~$Y$ and
  let $A$ be a $3$-generator of $Y$. Since $X$ is an acyclic extension
  of~$\freegpd{P}$, there exists a $3$-cell $F_{A}:\sce2(A)\TO\tge2(A)$ in the
  free $(3,1)$-category $\freegpd{P}(X)$. This induces a $3$-functor between
  free $(3,1)$-categories
  $
    f : \freegpd{P}(Y) \to \freegpd{P}(X)
  $,
  which is the identity on $P$ and such that $f(A)=F_{A}$ for every
  $3$-generator $A$ of~$Y$. Let~$X_{Y}$ be the subset of $X$ containing all the
  $3$-generators occurring in some $3$-cell $F_{A}$, for $A$ in $Y$. Since $Y$
  is finite and each $3$-cell~$F_{A}$ can be written as a composition of
  finitely many $3$-generators of $X$, we deduce that $X_{Y}$ is finite.

  Finally, consider a $2$-sphere $(\phi,\psi)$ of~$\freegpd{P}$. By hypothesis,
  there exists a $3$-cell $A:\phi\TO\psi$ in $\freegpd{P}(Y)$.  By application
  of $f$, one obtains a $3$-cell $f(A):\phi\TO\psi$ in
  $\freegpd{P}(X)$. Moreover, the $3$-cell $f(A)$ is a composite of
  cells~$F_{A}$, and the $3$-cell $f(A)$ is thus in $\freegpd{X}_{Y}$. As a
  consequence, one has $\phi\approx^{X_{Y}}\psi$, so that $X_{Y}$ is a finite
  acyclic cellular extension of~$\freegpd{P}$.
\end{proof}

\section{Convergence and Finite Derivation Type}
\label{sec:conv-fdt}
\cref{thm:SquierHomotopical} states that any family of generating confluences of a convergent $2$-polygraph~$P$ forms an acyclic extension of the free
$(2,1)$-category $\freegpd{P}$. The set of critical branchings of a finite
$2$-polygraph being finite, we deduce that a finite convergent $2$-polygraph has
finite derivation type.  Moreover, from \cref{thm:TDFTietzeInvariant}, the
property of having finite derivation type is Tietze invariant for finite
$2$-polygraphs. We thus obtain a finiteness condition for finitely presented
categories to have a presentation by a finite convergent $2$-polygraph.

\begin{theorem}
  \label{thm:FiniteConvergent=>FDT} 
  If a category admits a finite convergent presentation, then it has finite
  derivation type.
\end{theorem}

\noindent
This result was first proved by Squier for finitely presented
monoids~\cite[Theorem~5.3]{squier1994finiteness}. Several other proofs can be
found in the literature: we refer to~\cite{lafont1995new} for a reformulation of
Squier's arguments and to~\cite{GuiraudMalbos18} for a proof in the polygraphic
language presented in this book.

Now suppose we want to show that some category does not admit a finite
convergent presentation: by \thmr{FiniteConvergent=>FDT} it is
sufficient to prove that it has no finite derivation type.
The first example based on this argument, due to
Squier~\cite{squier1994finiteness}, is presented in
\cref{ex:MonoidSk}. Before that, we turn to a simplified version
introduced by Lafont and Prouté in~\cite{lafont1995new,lafont1991church}.

\begin{example}
  \label{ex:lafont-proute-squier}
  Consider the monoid~$M$ presented by the following $2$\nbd-poly\-graph:
\[
  P
  =
  \Pres{\star}{a,b,c,d,d'}{\alpha_0:ab\To a, \beta:da\To ac,\beta' : d'a\To ac}.
\]
This is a variant of the monoid already encountered in
\cref{ex:squier-lafont-monoid}. 
 It admits a finite presentation and has a decidable word problem,
yet it does not have finite derivation type and, as a consequence, it does not
admit a finite convergent presentation.
To prove these facts, the $2$-polygraph~$P$ is completed, by Knuth-Bendix
procedure (see \cref{sec:2kb}), into the following infinite convergent
$2$\nbd-poly\-graph
\[
  \tilde{P}
  =
  \Pres\star{a,b,c,d,d'}{\alpha_n, \beta, \beta'}_{n\in \N},
\]
with
\begin{align*}
  \alpha_n&:ac^nb\To ac^n,
  &
  \beta&:da\To ac,
  &
  \beta'&:d'a\To ac
  \pbox.
\end{align*}
Event though the polygraph $\tilde{P}$ is infinite, we can
implement an algorithm to normalize $1$-cells in $\tilde{P}$ by iteratively
rewriting those, and therefore decide the word problem in $\tilde{P}$, and
thus in~$P$, by comparing normal forms.
The $2$\nbd-poly\-graph~$\tilde{P}$ has two infinite families of critical
branchings from which we deduce two infinite families of $3$-generators:
\begin{align*}
  \xymatrix@R=3ex@C=3ex{
    & ac^{n+1}b
    \ar@2@/^/ [dr] ^{\alpha_{n+1}} 
    \ar@3 []!<0pt,-20pt>;[dd]!<0pt,20pt> ^{A_n}
    \\
    dac^nb
    \ar@2@/^/ [ur] ^{\beta c^nb}
    \ar@2@/_/ [dr] _{d\alpha_n}
    && ac^{n+1}
    \\
    & dac^n 
    \ar@2@/_/ [ur] _{\beta c^n}
  }
  &&
  \xymatrix@R=3ex@C=3ex{
    & ac^{n+1}b
    \ar@2@/^/ [dr] ^{\alpha_{n+1}} 
    \ar@3 []!<0pt,-20pt>;[dd]!<0pt,20pt> ^{A'_n}
    \\
    d'ac^nb
    \ar@2@/^/ [ur] ^{\beta' c^nb}
    \ar@2@/_/ [dr] _{d'\alpha_n}
    && ac^{n+1}\pbox.
    \\
    & d'ac^n 
    \ar@2@/_/ [ur] _{\beta' c^n}
  }  
\end{align*}
The $3$-generators~$A'_n$ induce a projection functor
$f : \freegpd{\tilde{P}} \to \freegpd{P}$
which is the identity on $0$- and $1$-generators, sends the $2$-generators
$\alpha_0$, $\beta$ and~$\beta'$ to themselves and the image of $\alpha_n$, for
$n>0$, is defined by induction by
\[
  f(\alpha_{n+1})=\beta'^-c^nb\comp1 d'f(\alpha_n)\comp1\beta'c^n\pbox.
\]
This functor is a retract of the canonical inclusion functor
$g:\freegpd P\to\freegpd{\tilde P}$. By the transfer theorem
(\cref{thm:HomotopyBasesTransfer}), the family
\[
  X=\setof{f(A_n)}{n\in\N}
\]
is thus an infinite acyclic cellular extension of the
$(2,1)$-category~$\freegpd{P}$: in this case, the generators of the form
\eqref{eq:tt-Aa} are superfluous because $f\circ g$ is the identity on
$\freegpd{P}$.

By \cref{Proposition:FiniteBasisExtracted}, in order to conclude that the polygraph~$P$ does
not have FDT, it is enough to show that no finite subset of~$X$ forms an acyclic
cellular extension of~$\freegpd{P}$. 
A fully explicit direct proof of this fact is rather tedious. A complete proof is given in {\cite[Section 5]{lafont1995new}} using an abelianized form of the category $\freegpd{P}$ in terms of monoidal groupoids. 
Note also that this can be shown indirectly by a homological argument outlined in the next chapter. Indeed, in \cref{ex:lafont-proute-squier-homology} we show that the third integral homology group of the monoid~$M$ is not of finite type. This shows by \cref{Theorem:TDFimpliesFP3} that the monoid~$M$ does not have FDT, and so in particular, we cannot extract a finite acyclic cellular extension of~$\freegpd{P}$ from~$X$.
\end{example}

\subsection{Squier's monoids}
\index{Squier!monoid}
\index{monoid!Squier}
\nomenclature[Sk]{$S_k$}{Squier monoid}
\label{ex:MonoidSk}
We now recall Squier's original example of a finitely presented monoid that does not admit a finite convergent presentation and studied in~\cite{squier1994finiteness} and~\cite{squier1987word} using homotopical and homological arguments respectively.
Consider, for $k\geq 1$, the monoid~$S_k$ defined in~\cite[Example 4.5]{squier1994finiteness} and presented by
\[
  \Pres{\star}
  {a,b,c,d_i,e_i}
  {\alpha_n, \beta_i,\gamma_i,\delta_i,\varepsilon_i}_{n\in\N,1\leq i\leq k},
\]
where the rules are $\alpha_n:ac^nb\To 1$ for $n\in\N$, and for $1\leq i \leq k$,
\[
\beta_i:d_ia\To acd_i,
\;
\gamma_i:d_ic\To cd_i,
\;
\varepsilon_i:d_ie_i\To 1,
\;
\delta_i:d_ib\To bd_i.
\]
Squier proves the following properties for the monoid~$S_1$ in~\cite[Theorem~6.7,
Corollary~6.8]{squier1994finiteness}. The proof is reworked
in~\cite[Section~6]{lafont1995new} and in~\cite[Section~6]{GuiraudMalbos18} using polygraphs. 

\begin{theorem}
  \label{Theorem:MonoidS1}
  \label{thm:MonoidS1}
The monoid~$S_1$ is a finitely presented monoid that has the following properties.
  \begin{enumerate}
  \item It has a decidable word problem.
  \item It does not have finite derivation type.
  \item It does not have a finite convergent presentation.
  \end{enumerate}
\end{theorem}
\begin{proof}
 The monoid~$S_1$ has the following infinite presentation:
  \[
    P=\Pres{\star}{a,b,c,d,e}{\alpha_n,\beta,\gamma,\delta,\varepsilon}_{n\in\N}
  \]
  with
\[
\alpha_n:ac^nb\To 1,
\;\;
\beta:da\To acd,
\;\;
\gamma:dc\To cd,
\;\;
\varepsilon:de\To 1,
\;\;
\delta:db\To bd.
\]

This presentation is infinite, so the normal-form algorithm of \cref{SS:NormalFormProcedure} cannot be applied to decide the word problem in the monoid $S_1$. However, the sources of the $2$-generators $\alpha_n$ are the elements of the regular language $ac^\ast b$. This implies that the sources of the $2$-generators of the polygraph $P$ form a regular language over the finite set $\{a,b,c,d,e\}$. Following {\cite[Proposition 3.6]{otto1998infinite}}, this implies that the word problem for $S_1$ is decidable, which proves Condition 1.

Condition~3 is a consequence of Condition~2 and \cref{thm:FiniteConvergent=>FDT}.

We sketch the main arguments of the proof of Condition 2, and we refer to~\cite{squier1994finiteness} for the original proof and
  to~\cite[Section~6]{GuiraudMalbos18} for the proof presented here.
We denote by $\gamma_n:dc^n\To c^n d$ the $2$-cell of $\freecat P_2$ defined by
  induction on~$n$ as follows:
  \[
    \gamma_0 = 1_x
    \qquad
    \text{and}
    \qquad
    \gamma_{n+1} = \gamma c^n \comp1 c \gamma_n
    \pbox.
  \]
  For every~$n$, we write $\phi_n:dac^nb\To ac^{n+1}bd$ the following composite
  in~$P_1^*$
  \[
    \xymatrix@C=8ex{
      dac^nb
      \ar@2 [r] ^-{\beta c^n b} 
      & acdc^n b
      \ar@2 [r] ^-{ac\gamma_n b}
      & ac^{n+1}db
      \ar@2 [r] ^-{ac^{n+1}\delta}
      & ac^{n+1}bd\pbox.
    }
  \]
  Considering for every natural number~$n\geq 0$, the following $2$-sphere
  of~$\freegpd{P}_2$:
  \begin{equation}
    \label{alphan}
    \vcenter{
      \xymatrix@C=3ex @R=1.5em {
        & ac^{n+1}bde
        \ar@2 [rr] ^-{ac^{n+1}b\varepsilon} 
        && ac^{n+1}b
        \ar@2@/^/ [dr] ^-{\alpha_{n+1}}
        \\
        dac^nbe
        \ar@2@/^/ [ur] ^-{\phi_n e}
        \ar@2@/_/ [drr] _-{d\alpha_n e}
        &&&& 1
        \\
        && de
        \ar@2@/_/ [urr] _-{\varepsilon}
      }
    }
  \end{equation}
we prove that the monoid~$S_1$ admits the following finite presentation
  \[
Q=\Pres{\star}{a, b, c, d, e} {\alpha_0, \beta, \gamma, \delta, \varepsilon}
    \pbox.
  \]
We also prove that the $2$-polygraph~$P$ is convergent and the Squier completion
  of~$P$ contains a $3$-generator~$A_n$ with shape
  \[
    \xymatrix @!C @C=1.5em {
      & ac^{n+1}bd
      \ar@2@/^/ [dr] ^-{\alpha_{n+1}d}
      \\
      dac^nb
      \ar@2@/^/ [ur] ^-{\phi_n}
      \ar@2@/_/ [rr] _-{d\alpha_n} ^{}="tgt"
      \ar@3 "1,2"!<0pt,-15pt>;"tgt"!<0pt,10pt> ^-*+{A_n}
      && d
    }
  \]
  for every natural number~$n$.
  In order to show that the monoid~$S_1$ does not have FDT,
  by \cref{thm:TDFTietzeInvariant}, it is sufficient to check that the polygraph~$Q$ admits no finite acyclic cellular extension. We
  denote by $g : \freegpd{P} \to \freegpd{Q}$
  the projection that sends the $2$-cells~$\beta$, $\gamma$, $\delta$,
  and~$\varepsilon$ to themselves and whose value on~$\alpha_n$ is given by
  induction on~$n$, thanks to~\eqref{alphan}, \ie
  \[
    g(\alpha_0) = \alpha
    \qquad\text{and}\qquad
    g(\alpha_{n+1}) = (\phi_n e \comp1 ac^{n+1}b\varepsilon)^- \comp1 d f(\alpha_n)e \comp1 \varepsilon.
  \]
  By application of \cref{Theorem:HomotopyBasesTransfer} to the canonical
  inclusion~$f : \freegpd{Q} \fl \freegpd{P}$ and~$g$ defined above, we deduce that the monoid~$S_1$ admits the
  coherent presentation
  \[
    \widetilde{Q} = \PRes{\star}{a,b,c,d,e}{\alpha_0,\beta,\gamma,\delta,\varepsilon}{\tilde A_n}_{n\in\N}
  \]
  where~$\tilde{A}_n$ is the $3$-generator
  \[
    \xymatrix @!C @C=1.5em {
      & ac^{n+1}be
      \ar@2@/^/ [dr] ^-{g(\alpha_{n+1})d}
      \\
      dac^nb
      \ar@2@/^/ [ur] ^-{\phi_n}
      \ar@2@/_/ [rr] _-{dg(\alpha_n)} ^{}="tgt"
      \ar@3 "1,2"!<0pt,-15pt>;"tgt"!<0pt,10pt> ^-*+{\tilde{A}_n}
      && d\pbox.
    }
  \]
By studying the relations among $3$-cells in the $3$-category generated by the coherent presentation $\widetilde{Q}$ in terms of generating critical, see  \cref{Subsection:TripleConfluences}, we show that the polygraph $Q$ is not Tietze equivalent to a polygraph having FDT. Following \cref{thm:TDFTietzeInvariant}, this shows that the monoid $S_1$ does not have FDT.  We refer to~\cite{GuiraudMalbos18} for more details on this proof.
\end{proof}
 
\subsection{Higher-dimensional finite derivation type}
With the aim of characterizing the class of finite presented decidable monoids
admitting a finite convergent presentation, some refinements of the
FDT condition were introduced, such as the $2$-dimensional FDT
property~\cite{PrideGlashanPasku05}, and the 
infinite-deimensional FDT property~\cite{GuiraudMalbos12advances}, see also \cref{sec:higher-FDT}. 
Note that the characterization of this class by finiteness conditions is still an open problem.

\section{Identities Among Relations}
\label{sec:IdentitiesAmongRelations}
The $3$-generators in a coherent presentation are closely related to the notion
of \emph{identity among relations}, which
originates in the work of Peiffer and Reidemeister in combinatorial group
theory~\cite{Peiffer49,Reidemeister49}. This notion is based on the one of crossed
module, introduced by Whitehead, in algebraic topology, for the classification
of homotopy $2$-types~\cite{Whitehead49a,Whitehead49b}. There exist several
formulations of identities for presentations of groups: as
homological $2$-syzygies~\cite{BrownHuebschmann82}, as homotopical
$2$-syzygies~\cite{Loday00}, or as Igusa's
pictures~\cite{kapranov1999hidden,Loday00}. One can also interpret identities
as the critical pairs of a presentation of a group by a
convergent string rewriting system~\cite{cremanns1996groups}. The latter approach
yields an algorithm based on Knuth-Bendix's completion procedure that computes a
family of generators of the module of identities among
relations~\cite{heyworth2003logged}.

In this section, we introduce the notion of identity among relations for a
$2$-polygraph. We relate the property for a polygraph of having a finite
generating set of identities among relations to the property of having abelian
finite derivation type. 
First, let us recall the notion of natural system used in this section.

\subsection{Natural system}
\index{category!of factorizations}
\index{natural system}
\label{sec:natural-system}
The \emph{category of factorizations} of a small category $C$ is the category,
denoted by $\fact{C}$, whose $0$-cells are the $1$-cells of $C$ and whose
$1$-cells from $w$ to $w'$ are pairs $(u,v)$ of $1$-cells of $C$ such that the
following diagram commutes in~$C$:
\[
  \xymatrix @R=2em @C=2em{   
    \cdot\ar[d]_w&\ar[l]_u\cdot\ar[d]^{w'}
    \\
    \cdot\ar[r]_v&\cdot
  }
\]
The triple $(u,w,v)$ is called a \emph{factorization of $w'$}.
\nomenclature[Ab]{$\Ab$}{category of abelian groups}
A \emph{natural system} on $C$ is a functor $D:\fact{C}\to\Ab$ with values in the category $\Ab$ of abelian groups. We will denote by $D_w$ the abelian group which is the
image of $w$ by $D$. The category of natural systems is denoted $\NatSys{C,\Ab}$.
\nomenclature[Nat(C,Ab)]{$\NatSys{C,\Ab}$}{category of natural systems}
We refer to \cref{SS:NaturalSystems} for more details on this notion.

\subsection{Identities among relations}
\label{sec:id-among-rel}
\index{identity among relations}
Let $P$ be a $2$-polygraph. We define the natural system $\Pi(P)$ on the
presented category $\pcat{P}$ of \emph{identities among relations} of $P$ as
follows.
\begin{itemize}
\item If $u$ is a $1$-cell of $\cl{P}$, the abelian group $\Pi(P)_u$ is
  generated by one element~$\iar{\varphi}$, for each $2$-cell $\varphi:v\dfl v$
  of the $(2,1)$-category $\freegpd{P}$ such that $\cl{v}=u$, and subject to the relation
  \begin{equation}
    \label{eq:iar1}
    \iar{\varphi\comp1 \psi} = \iar{\varphi}+\iar{\psi},
  \end{equation}
  for every $2$-cells $\varphi:v\dfl v$ and $\psi:v\dfl v$ of $\freegpd{P}$, with $\cl{v}=u$, and
  \begin{equation}
    \label{eq:iar2}
    \iar{\varphi\comp1 \psi} = \iar{\psi\comp1 \varphi},
  \end{equation}
for every $2$-cells $\varphi:v\dfl w$ and $\psi:w\dfl v$ of $\freegpd{P}$, with $\cl{v}=\cl{w}=u$.
\item If $w'=uwv$ is a factorization in $\cl{P}$, then the homomorphism of groups $\Pi(P)_{(u,v)} : \Pi(P)_w \fl \Pi(P)_{w'}$ is defined by
\[
\Pi(P)_{(u,v)} (\iar{\varphi}) = \iar{\rep{u}\varphi\rep{v}},
\]
where $\rep{u}$ and $\rep{v}$ are any representative $1$-cells of $u$ and $v$ in $P_1^*$, respectively. 
\end{itemize}
Note that the value of $\Pi(P)_{(u,v)}$ does not depend on the choice of the representative $1$-cells $\rep{u}$ and $\rep{v}$.
This proves that $\Pi(P)$ is a natural system on $\cl{P}$. We will often write $\iar{u\varphi v}$ instead of $\iar{\rep{u}\varphi\rep{v}}$.

As a consequence of the defining relations of each group $\Pi(P)_u$, the relations
\begin{align*}
  \iar{1_u} &= 0,
  &
  \iar{\varphi^-} &= -\iar{\varphi},
  & \text{and} \qquad
  \iar{ \psi\comp{1} \varphi \comp{1} \psi^- } = \iar{\varphi} 
\end{align*}
hold for every $1$-cell $u$ and every $2$-cells $\varphi:u\dfl u$ and
$\psi:v\dfl u$ of the free $(2,1)$\nbd-cate\-gory $\freegpd{P}$.

\subsection{Loops and cellular extensions}
In some situations, it is helpful to consider cellular
extensions by the means of \emph{$2$-loops}\index{2-loop}\index{loop!2-} in a $2$-category $C$, \ie $2$-cells~$\varphi$ such that $\src1(\varphi)=\tgt1(\varphi)$.  The following result will be useful in the sequel.

\begin{lemma}
\label{LemmaHomotopiesBoucles}
Let $C$ be a $(2,1)$-category and let $Y$ be a family of 2-loops in~$C$. The following assertions are equivalent. 
\begin{enumerate}
\item The cellular extension $\widetilde{Y} := \{ \widetilde{\beta}:\beta\tfl 1_{\src1(\beta)}, \:\beta\in Y \}$ of $C$ is acyclic.
\item Every 2-loop $\varphi$ in $C$ has a decomposition
\begin{equation}
\label{HomotopyBaseTDFLemmaEq}
\varphi = \left( \psi_1 \comp{1} u_1\beta_1^{\epsilon_1}v_1 \comp{1} \psi_1^- \right) 
	\:\comp{1}\: \cdots \:\comp{1}\: 
	\left( \psi_p \comp{1} u_p\beta_p^{\epsilon_p}v_p \comp{1} \psi_p^- \right)
\end{equation}
with, for every $1\leq i\leq p$, $\beta_i$ in $Y$, $\epsilon_i$ in $\{-,+\}$, $u_i,v_i$ 1-cells of $C$ and $\psi_i$ a 2-cell of~$C$.
\end{enumerate}
\end{lemma}
\begin{proof}
Suppose that $C$ is $\widetilde{Y}$-acyclic. Given a closed $2$-cell $\varphi:w\dfl w$ in~$C$, by hypothesis there exists a $3$-cell $A:\varphi\tfl 1_w$ in $C(\widetilde{Y})$. In the $(3,1)$-category $C(\widetilde{Y})$ the $3$-cell $A$ can be decomposed into
\[
A = A_1\comp2\cdots\comp2 A_k,
\]
where each $A_i$ is a $3$-cell of $C(\widetilde{Y})$ that contains exactly one generating $3$-cell of~$Y$. Thus each $3$-cell $A_i$ has the shape
\[
\psi_i\comp{1} u_i\widetilde{\beta}_i^{\epsilon_i}v_i \comp{1} \psi'_i
\]
with $\beta_i \in Y$, $\epsilon_i \in \{-,+\}$, $u_i,v_i$ $1$-cells of $C$ and $\psi_i, \psi'_i$, $2$-cells of  $C$. By hypothesis on $A$, we have $\varphi=\src2(A)$, hence
$\varphi = \psi_1\comp{1} u_1\src2(\beta_1^{\epsilon_1})v_1 \comp{1} \psi'_1$.
For $\epsilon_1=+$, we have:
\begin{align*}
\varphi \:
	&=\: \psi_1\comp{1} u_1\beta_1v_1 \comp{1} \psi'_1 \\
	&=\: \left( \psi_1\comp{1} u_1\beta_1v_1 \comp{1} \psi_1^- \right) \comp1 \left( \psi_1 \comp{1} \psi'_1 \right) \\
	&=\:  \left( \psi_1\comp{1} u_1\beta_1v_1 \comp{1} \psi_1^- \right) \comp1 \src2(A_2).
\end{align*}
And, for $\epsilon_1=-$, we have:
\begin{align*}
\varphi \:
	&=\: \psi_1 \comp{1} \psi'_1 \\
	&=\: \left( \psi_1 \comp{1} u_1\beta_1^-v_1 \comp{1} \psi_1^- \right) 
		\comp{1} \left( \psi_1 \comp{1} u_1\beta_1v_1 \comp{1} \psi'_1 \right) \\
	&=\: \left( \psi_1 \comp{1} u_1\beta_1^-v_1 \comp{1} \psi_1^- \right) \comp{1} \src2(A_2).
\end{align*}
We proceed by induction on~$k$ to prove that $\varphi$ has a decomposition as in~\eqref{HomotopyBaseTDFLemmaEq}.

Conversely, we assume that every closed $2$-cell $\varphi$ in $C$ has a decomposition as in~\eqref{HomotopyBaseTDFLemmaEq}. Then we have $\varphi \approx^{\widetilde{Y}} 1_{\src1(\varphi)}$ for every closed $2$-cell $\varphi$ in $C$. Let us consider two parallel $2$-cells $\varphi$ and $\psi$ in $C$. Then $\varphi\comp{1} \psi^-$ is a closed $2$-cell, yielding $\varphi\comp{1} \psi^- \approx^{\widetilde{Y}} 1_{s(\varphi)}$. We compose both members by $\psi$ on the right hand to get $\varphi\approx^{\widetilde{Y}} \psi$. Thus~$\widetilde{Y}$ is a homotopy basis of $C$.
\end{proof}

\subsection{Abelian finite derivation type}
\label{subsectionAbelianTDF}
A $(2,1)$-category $C$ is called \emph{abelian} if, for every $1$-cell $u$ of
$C$, the group $\Aut^{C}_u$ of $2$-loops of $C$ with source $u$ is abelian. 
For $C$ an $(2,1)$-category, its \emph{abelianization} $\fab C$ is the quotient of $C$ by the cellular extension that contains one $2$-sphere
$\varphi\comp{1} \psi \TO \psi\comp{1} \varphi$
for all $2$-loops $\varphi$ and $\psi$ of $C$ with the same source.

\index{FDTab@$\FDTAB$ (abelian finite derivation type)}
\index{finite derivation type!abelian}
One says that a $2$-polygraph $P$ has \emph{abelian finite derivation type}, or
$\FDTAB$ for short, when the abelian $(2,1)$-category $\fab{\freegpdab{P}}$
admits a finite acyclic extension.

\begin{proposition}
  \label{Proposition:IsomorphismIarAut}
  Given a 2-polygraph $P$, there exists an isomorphism of natural systems on
  the free category $P^*_1$:
\begin{equation}
\label{E:isormorphismNaturalSystemsLoops}
\Pi(P)_{\pi\circ (-)} \overset{\simeq}{\longrightarrow} \Aut^{\fab{\freegpdab{P}}}_{(-)}
\pbox.
\end{equation}
\end{proposition}
\begin{proof}
For a $1$-cell $u$ of $\fab{\freegpdab{P}}$, we define the morphism of groups
\[
\Phi_u:\Pi(P)_{\cl{u}}\fl\Aut^{\fab{\freegpdab{P}}}_u
\]
given on generators by $\Phi_u(\iar{\varphi}) =
\varphi^\psi=\psi^-\comp 1 \varphi\comp 1\psi$,
  where $\varphi$ is a $2$-loop of $\fab{\freegpdab{P}}$ on a $1$-cell $v$ such
  that $\cl{v}=\cl{u}$ and $\psi:v\To u$ is any $2$-cell of
  $\fab{\freegpdab{P}}$. The morphism $\Phi_u$ is well-defined. Indeed, it is
  independent of the choice of~$\psi$, and its definition is compatible with the
  relations~\eqref{eq:iar1} and~\eqref{eq:iar2} defining $\Pi(P)_{\cl{u}}$.

  For the relation~\eqref{eq:iar1}, let $\varphi_1$ and $\varphi_2$ be $2$-loops of
  $\fab{\freegpdab{P}}$ on a $1$-cell $v$ such that $\cl{v}=\cl{u}$ and let
  $\psi:v\To u$ be a $2$-cell of $\fab{\freegpdab{P}}$. Then,
  \begin{align*}
    \Phi_u(\iar{\varphi_1\comp1\varphi_2}) 
    &= (\varphi_1\comp1\varphi_2)^\psi \\
    &= \varphi_1^\psi \comp1 \varphi_2^\psi \\
    &= \Phi_u(\iar{\varphi_1}) \comp1 \Phi_u(\iar{\varphi_2}) \\
    &= \Phi_u(\iar{\varphi_1}+\iar{\varphi_2})
    \pbox.
  \end{align*}
  For the relation~\eqref{eq:iar2}, we fix $2$-cells
  $\varphi_1:v_1\To v_2$, $\varphi_2:v_2\To v_1$ and $\psi:v_1\To u$, with
  $\cl{v_1}=\cl{v_2}=\cl{u}$. Then,
  \begin{align*}
    \Phi_u(\iar{\varphi_1\comp1\varphi_2})
    &= (\varphi_1\comp1\varphi_2)^\psi \\
    &= (\psi^- \comp1 \varphi_1) \comp1 (\varphi_2 \comp1 \varphi_1) \comp1 (\varphi_1^- \comp1 \psi) \\
    &= (\varphi_2\comp1 \varphi_1)^{\psi^-\comp1 \varphi_1} \\
    &= \Phi_u(\iar{\varphi_2\comp1\varphi_1})
    \pbox.
  \end{align*}
  Thus $\Phi_u$ is a morphism of groups from $\Pi(P)_{\cl{u}}$ to
  $\Aut^{\fab{\freegpdab{P}}}_u$. Moreover, it admits
  $\varphi\mapsto\iar{\varphi}$ as inverse and, as a consequence, is an
  isomorphism.

  Let us prove that $\Phi_u$ is natural in $u$. Let $K$ be a context of 
  $\freecat{P}_1$ such that $v=K[u]$, and prove the equality of the two morphisms
  $\Phi_v\circ\Pi(P)_{\cl{K}}$ and $\Aut^{P}_K \circ \Phi_u$. Let
  $\varphi$ be a $2$-loop of $\fab{\freegpdab{P}}$ with source $u'$ such that
  $\cl{u}'=\cl{u}$. We fix a $2$-cell $\psi:u'\fl u$ in $\fab{\freegpdab{P}}$
  and consider the $2$-cell $K[\psi] : K[u'] \fl v$ of $\fab{\freegpdab{P}}$. Then, we have
  \begin{align*}
    \Phi_v\circ\Pi(P)_{\cl{K}} (\iar{\varphi}) 
    &= (K[\varphi])^{K[\psi]} \\
    &= K[\psi^-] \comp1 K[\varphi] \comp1 K[\psi] \\
    &= K\left[ \psi^-\comp1 \varphi \comp1 \psi \right] \\ 
    &= K[\varphi^\psi] \\
    &= \Aut^{\fab{\freegpdab{P}}}_K \:\circ\: \Phi_u (\iar{\varphi})
    \pbox.\qedhere
  \end{align*}
\end{proof}

\Cref{Proposition:IsomorphismIarAut} characterizes the natural system $\Pi(P)$
on the category $\cl{P}$ up to isomorphism. Using this characterization, we
deduce the following result.

\begin{proposition}
  \label{Proposition:IartfIffFdtab}
  A 2-polygraph $P$ has $\FDTAB$ if and only if the
  natural system~$\Pi(P)$ is finitely generated.
  \end{proposition}
\begin{proof}
  Suppose that the $2$-polygraph $P$ has $\FDTAB$. Then
  the abelian $(2,1)$-category $\fab{\freegpdab{P}}$ admits a finite acyclic
  extension~$X$. Given a $3$-generator $A:\phi\TO\psi$ in~$X$, we write
  $\partial A=\phi\comp2\finv\psi$ and $\partial X=\setof{\partial A}{A\in X}$
  for the set of $2$-loops of~$\fab{\freegpdab{P}}$ associated to $3$-generators in $X$.
    
  By Lemma~\ref{LemmaHomotopiesBoucles}, any $2$-loop
  $\varphi$ can be written in $\fab{\freegpdab{P}}$ as
  \[
    \varphi =
    \pa{\psi_1 \comp{1} u_1\partial A_1^{\epsilon_1}v_1 \comp{1}\finv\psi_1}
    \comp1\ldots\comp1
    \pa{\psi_p \comp{1} u_p\partial A_p^{\epsilon_p}v_p \comp{1}\finv\psi_p},
  \]
  with, for every $1\leq i\leq p$, $A_i$ in $X$, $\epsilon_i$ in $\set{-1,+1}$,
  $u_i,v_i$ $1$-cells of $\freecat P_1$ and $\psi_i$ a $2$-cell of the free
  $(2,1)$-category~$\freegpd{P}$.  As a consequence, for any $\iar{\varphi}$ in
  $\Pi(P)$, we have the following decomposition:
  \[
    \iar{\varphi} 
    = \sum_{i=1}^k (-1)^{\epsilon_i} \iar{\psi_i \comp{1} u_i\partial A_iv_i \comp{1} \psi_i^-}
    = \sum_{i=1}^k (-1)^{\epsilon_i} \cl{u}_i\iar{\partial A_i}\cl{v}_i
    \pbox.
  \]
  Thus, the elements of $\iar{\partial X}$ form a finite generating set for the
  natural system of abelian groups~$\Pi(P)$.

  Conversely, suppose that the natural system~$\Pi(P)$ is
  finitely generated. There exists a finite set~$X$ of $2$-loops of the abelian
  $(2,1)$\nbd-cate\-gory $\fab{\freegpdab{P}}$ such that, for every $1$-cell $\cl{u}$ of
  $\cl{P}$ and every $2$-loop $\varphi$ with source $w$ of $\fab{\freegpdab{P}}$ such that
  $\cl{w}=\cl{u}$, one can write
  \[
    \iar{\varphi} = \sum_{i=1}^p \epsilon_i \cl{u}_i\iar{\alpha_i}\cl{v}_i,
  \]
  with, for every $1\leq i\leq p$, $\alpha_i$ in $X$, $\epsilon_i$ an integer
  and $u_i,v_i$ $1$-cells of $\cl{P}$ such that, for every representative
  $\rep{u}_i$ of $\cl{u}_i$ and $\rep{v}_i$ of $\cl{v}_i$ in $\fab{\freegpdab{P}}$,
  $\rep{u}_i\alpha_i\rep{v}_i$ is a $2$-loop of $\fab{\freegpdab{P}}$ whose source $w_i$
  that satisfies $\cl{w}_i=\cl{w}$. We fix, for every~$i$, a $2$-cell
  $\psi_i:w\dfl w_i$ in $\freegpd{P}$. Then, the properties of $\Pi(P)$ imply:
  \begin{align*}
    \iar{\varphi}
    &= \sum_{i=1}^p \iar{\psi_i \comp{1} \rep{u}_i\alpha_i^{\epsilon_i}\rep{v}_i \comp{1} \psi_i^-} \\
    &= \iar{ \big( \psi_1 \comp{1} \rep{u}_1\alpha_1^{\epsilon_1}\rep{v_1} \comp{1} \psi_1^- \big)
      \comp{1} \ldots \comp{1}
      \big( \psi_p \comp{1} \rep{u}_p\alpha_p^{\epsilon_p}\rep{v_p} \comp{1} \psi_p^- \big) }\pbox.
  \end{align*}
  We use the isomorphism~\cref{E:isormorphismNaturalSystemsLoops} and Lemma~\ref{LemmaHomotopiesBoucles} to deduce
  that the cellular extension
  $\setof{A_\alpha:\alpha\To 1_{s(\alpha)}}{\alpha\in X}$ of $\fab{\freegpdab{P}}$
  is acyclic, proving that the $2$-polygraph $P$ has $\FDTAB$.
\end{proof}

The following result states that the property of being finitely generated
for~$\Pi(P)$ is Tietze invariant for polygraphs $P$ having a finite set of
$2$-gener\-ators~\cite[Proposition 2.3.5]{GuiraudMalbos10smf}.

\begin{proposition}
Let $P$ and $Q$ be two Tietze equivalent 2-polygraphs such that $P_2$ and $Q_2$ are finite. Then the natural system~$\Pi(P)$ is finitely generated if and only if the natural system $\Pi(Q)$ is finitely generated.
\end{proposition}

\noindent
From this result and \cref{Proposition:IartfIffFdtab}, we deduce that the property $\FDTAB$ is Tietze invariant for finite polygraphs. As a consequence, we can define a category $\FDTAB$ if it admits a presentation by a finite $2$-polygraph having~$\FDTAB$.

\medskip

We conclude this chapter with a remarkable property of the natural system of identities among relations from~{\cite[Proposition 2.4.2]{GuiraudMalbos10smf}}, which is a consequence of Squier's homotopical theorem. By \cref{Theorem:SquierCompletion2polygraphs}, the set of generating confluences of a convergent $2$-polygraph $P$ forms an acyclic extension of the $(2,1)$-category $\freegpd{P}$. Following the proof of \cref{Proposition:IartfIffFdtab}, we transform this extension into a generating set for the natural system $\Pi(P)$, proving the following result.

\begin{theorem}
Let $P$ be a convergent 2-polygraph.
The natural system $\Pi(P)$ is generated by the generating confluences of $P$.
\end{theorem}


\chapter{Homological Syzygies and Confluence}
\label{chap:2-homology}
\label{chap:2homology}
\label{chap:HomologieSquierTheorem}
The main purpose of algebraic topology is the classification of topological
spaces and continuous maps by means of discrete algebraic invariants
preserving homotopy equivalence. Among those invariants, a
particularly important one is homology, which assigns to each space a sequence of
abelian groups. Starting from very geometric insights, 
homology has developed into a whole body of concepts and methods known
as homological algebra and has been applied to the
study of various algebraic structures, including groups and
monoids~\cite{MacLane95}. For instance, the homology of a monoid is defined by
first building a \emph{resolution} of it, that is, an exact sequence of left-modules
over the ring generated by the monoid, ending at the trivial
module. Of course the soundness of this definition is based on the fact that the
homology does not depend on the  choice of the resolution.

Squier showed in his 1987 article~\cite{squier1987word} that a convergent presentation $P$ of a
monoid $M$ yields a partial resolution generated by the set~$P_1$ of
generators in dimension~$1$, by the set $P_2$ of rules in
dimension~$2$ and by the critical branchings in dimension~$3$. 
If moreover the presentation $P$  is finite, the Squier resolution is
finitely generated up to dimension~$3$. In this case, we say that the
monoid $M$ is
of \emph{homological type left-$\FP_3$}. This property readily implies that the third integral homology group $H_3(M,\Z)$ of the monoid is finitely
generated. Therefore, a monoid whose third homology group is not
finitely generated does not admit a finite convergent presentation.
By explicitly exhibiting an example of this type,
Squier first provided a negative answer to the question of universality of
convergent rewriting.

The homological finiteness condition is of course linked to the homotopical one 
discussed in the previous chapter. Indeed, we will prove that for a monoid, the 
property of having $\FDT$ implies the property of having left-$\FP_3$
(\cref{Theorem:TDFimpliesFP3}). In this sense, the homological finiteness
condition is weaker than its homotopical counterpart. It however has the
advantage of being simpler to compute.

We begin by introducing the finite homological type for monoids in
\cref{sec:homological-finiteness}. We study the case $n=2$ in \cref{sec:FP2} and
show that finitely presented monoids are left-$\FP_2$
(\cref{Proposition:MonoidsFP2}). Then, we study the case~$n=3$ in
\cref{sec:FP3}, and show that monoids with finite convergent
presentations have the left-$\FP_3$ property (\cref{Theorem:SquierAb1}). 
We illustrate these results with several examples.
The constructions of this chapter will be generalized in any homological dimension for $1$-categories in \cref{sec:H-iar}. The case of monoids treated in this chapter corresponds to $1$-categories with a single object.

The homological notions used in this chapter are recalled in \cref{Chapter:ComplexesAndHomology}. In particular, homology of monoids is recalled in \cref{sec:mon-homology}.
In this chapter, we study the homological type left-$\FP_3$ relative to left modules, but the homological type right-$\FP_3$ relative to right modules is treated in the same way. We refer to \cref{S:CategoriesFiniteHomologicalType} and \cref{Section:AnotherHomologicalType} for relationships between homology types according to the module categories considered.

\section{Monoids of Finite Homological Type}
\label{sec:homological-finiteness}

\subsection{Monoid ring}
\index{monoid!ring}
Let~$M$ be a monoid. The \emph{ring generated by~$M$} is the free abelian group
over~$M$, denoted by~$\Z M$. Its elements are formal sums $\sum_{u\in M}n_uu$ of
elements $u$ of~$M$ with coefficients $n_u\in\Z$, finitely many of which are
non-zero, and it is equipped with the canonical extension of the product of~$M$:
\[
\pa{\sum_{u\in M}n_u u}
\pa{\sum_{v\in M}n_v v}
=
\sum_{u,v\in M}n_un_v uv
=
\sum_{w\in M}
\pa{\sum_{uv=w}n_un_v}w
\pbox.
\]
This construction coincides with the one of the free $\Z$-module, thus the
notation.

\subsection{Free modules}
\index{free!module}
Given a monoid~$M$ and a set~$X$, we write $\freemod{\Z M}{X}$ for the free left
$\Z M$-module generated by~$X$: its elements are formal sums of the form
\[
\sum_{u\in M, x\in X}n_{u,x}u[x]
\] 
with $n_{u,x}\in\Z$, finitely many of which
are non-zero, and other operations are defined in the expected way.
Any function $f:X\to C$, where $C$ is a $\Z M$-module extends uniquely as a
morphism of $\Z M$-modules $f:\freemod{\Z M}X\to C$.
Note that any $\Z M$-module is also canonically a $\Z$-module.

\subsection{Resolutions}
\index{resolution}
\index{projective!resolution}
\index{trivial module}
\index{module!trivial}
\label{SS:ResolutionsHomologicalSyzygies}
If~$M$ is a monoid, the \emph{trivial $\Z M$-module} is the abelian group~$\Z$
equipped with the trivial action~$un=n$, for every~$u$ in~$M$ and~$n$ in~$\Z$.
A \emph{partial resolution} of length~$n$ of this trivial module consists of a
chain complex
\begin{equation}
  \label{eq:resolution-n}
  \xymatrix{
    C_{n} \ar[r] ^-{d_n}
    & C_{n-1} \ar[r] ^-{d_{n-1}}
    & \cdots \ar[r] ^-{d_2}
    & C_1 \ar[r] ^-{d_1} 
    & C_0 \ar[r] ^-{d_0} 
    & \Z \ar[r]
    & 0
  }
\end{equation}
where, for~$0\leq k\leq n$, the $C_k$ are left $\Z M$-modules, and
the $d_k$ are $\Z M$\nbd-linear maps making the complex exact. By convention, $C_{-1}=\Z$ and
$d_{-1}:\Z\to 0$ is the terminal map. The main properties on resolutions that we use in
this chapter are recalled in \cref{S:Resolution}.

A contracting homotopy of a chain complex of the form~\eqref{eq:resolution-n} is a sequence
\[
  \xymatrix{
    C_{n}
    & C_{n-1} \ar[l]_-{i_n}
    & \cdots \ar[l]_-{i_{n-1}}
    & C_1 \ar[l]_-{i_2} 
    & C_0 \ar[l]_-{i_1} 
    & \Z \ar[l]_-{i_0}
  }
\]
where the $i_k$ are $\Z$-linear maps for
$0\leq k\leq n$, and such that
\[
  d_k\circ i_k+i_{k-1}\circ d_{k-1}=\id_{C_{k-1}}
\]
holds for $0\leq k\leq n$, see~\cref{SS:ContractingHomotopy} for details.
By convention, $i_{-1}:0\to\Z$ is the initial map. Any
chain complex of the form \eqref{eq:resolution-n} equipped with a contracting
homotopy is necessarily a partial resolution, see \cref{prop:contr-hom-exact}.

\subsection{Homological type left-$\FP_n$}
\nomenclature[FPn]{$\FP_n$}{homological type}
\index{homological!left-FP8@left-$\FP_n$}
\index{left-FP8@left-$\FP_n$}
\index{FP8@$\FP_n$}
\label{Definition:LeftFPn}
A monoid~$M$ has \emph{homological type left-$\FP_n$} (where $\FP_n$ stands for ``finitely $n$-presented''), for a natural number~$n$, if there exists a
partial resolution of length~$n$ of the trivial $\Z M$-module~$\Z$ of the
form~\eqref{eq:resolution-n}, where the $C_i$ are projective modules which are finitely generated.
A monoid~$ M$ has \emph{homological type left-$\FP_{\infty}$} if it has
homological type left\nbd-$\FP_n$ for all $n\geq 0$.

We will use the following characterization given by \cref{lemme_pl_n}: a monoid~$ M$ has homological type left-$\FP_n$ if and only if there exists a free, finitely generated partial resolution of the trivial $\Z M$-module~$\Z$ of length~$n$:
 \begin{equation}
 \label{E:ResolutionFPn}
\xymatrix{
F_n \ar[r] 
& F_{n-1} \ar[r] 
& \cdots \ar[r] 
& F_0 \ar[r] 
& \Z.
}
\end{equation}

\subsection{Finiteness homological type and homology}
For a monoid $M$,
having homological type left-$\FP_n$ implies a finiteness property on its homology modules. First, we recall the definition of homology of a monoid and we refer to \cref{sec:mon-homology} for more details.
Given a free resolution 
\[
\xymatrix{
\cdots \ar[r]
& F_{n+1} \ar[r] ^-{d_{n+1}} 
& F_n \ar[r] 
& \cdots \ar[r]
& F_1 \ar[r] ^-{d_1}
& F_0 \ar[r] ^-{\varepsilon}
& \Z
}
\]
of the trivial $\Z M$-module  $\Z$ by left $\Z M$-modules, the
operation of tensoring by the trivial right $\Z M$-module~$\Z$ gives the following complex of $\Z$-modules:
\[
\xymatrix@C=1.5em{
\cdots \ar[r]
& \Z\otimes_{\Z M}F_{n+1} \ar[r] ^-{\widetilde{d}_{n+1}} 
& \Z\otimes_{\Z M}F_n \ar[r] 
& \cdots \ar[r]
& \Z\otimes_{\Z M}F_1 \ar[r] ^-{\widetilde{d}_1}
& \Z\otimes_{\Z M}F_0
}
\]
where $\widetilde{d}_k$ denotes the map $\id_{\Z}\otimes_{\Z M} d_{k}$, for all $k\geq 1$.
The \emph{$n$-th homology group}\index{homology!group} of $M$ with integral coefficient $\Z$ is defined as the following $\Z$-module: 
\[
\mathrm{H}_n(M,\Z) = \ker \widetilde{\bnd}_n / \im \widetilde{\bnd}_{n+1},
\]
with the convention that $\widetilde{d}_0=0$. By definition, for any monoid $M$, we have $\mathrm{H}_0(M,\Z)\simeq \Z$.

Now, suppose that the monoid $M$ has homological type left-$\FP_n$ and
consider a resolution of $M$ of the form~\cref{E:ResolutionFPn}. Then the $\Z$-modules $\Z\otimes_{\Z M}F_i$ are finitely generated for $0\leq i \leq n$. This proves the following result.

\begin{proposition}
  \label{prop:FPn-H-fg}
  If a monoid $M$ has homological type left-$\FP_n$ for some $n\in\N$, then the
  groups $H_k(M,\Z)$ are finitely generated for $0\leq k\leq n$.
\end{proposition}

\subsection{Homological type left-$\FP_0$}
\label{Subsection:MonoidHomologicalTypeFP0}
Let $M$ be a monoid. We write $P_0=\set{\star}$ for a set with one element. We
have that $\freemod{\Z M}{P_0}\simeq\Z M$ and sometimes implicitly
identify the elements of these
two modules. The \emph{augmentation map} of $\Z M$ is the morphism of
$\Z M$-modules
\[
  \varepsilon
  :
  \freemod{\Z M}{P_0}
  \to
  \Z
\]
defined by $\varepsilon(u)=1$ for any $u$ in $\Z M$. The
augmentation map is clearly surjective and thus the sequence
\[
  \xymatrix{
    \freemod{\Z M}{P_0}
    \ar [r] ^-{\varepsilon}
    & \Z
    \ar [r]
    & 0
  }
\]
is exact. It follows that every monoid has homological type left-$\FP_0$.

\subsection{Homological type left-$\FP_1$}
Let~$P$ be a presentation of a monoid~$M$. We define a free partial resolution
of length~$1$ of the trivial $\Z M$-module~$\Z$ by $\Z M$-modules
\[
  \xymatrix{
    \freemod{\Z M}{P_1}
    \ar [r] ^-{d_1}
    &
    \freemod{\Z M}{P_0}
    \ar [r] ^-{\varepsilon}
    & \Z
    \ar [r]
    & 0
  }
\]
where
the morphism $\varepsilon$ is the augmentation map and the morphism $d_1$ is
defined, on any generator $[a]$, by
\[
d_1([a]) = \cl{a} - 1.
\]

A \emph{section} of the canonical projection $\pi : P_1^\ast \to M$ is a map $M \to P_1^\ast$ sending every~$u$ in~$M$ to a $1$-cell~$\rep{u}$ of~$P_1^\ast$ such that $\pi(\rep{u})=u$. In general, we do not assume that the chosen section is functorial, i.e., that $\rep{uv}=\rep{u}\rep{v}$ holds in~$P_1^\ast$. However, we assume that~$\rep{1}=1$. For a $1$-cell~$u$ of~$P_1^\ast$, we simply write~$\rep{u}$ for~$\rep{\cl{u}}$.

\begin{proposition}
  \label{Proposition:MonoidsFP1}
  If a monoid $M$ is finitely generated, then it has homological
  type left\nbd-$\FP_1$, and thus the group $H_1(M,\Z)$ is finitely generated.
\end{proposition}
\begin{proof}
  We first note that the sequence is a chain complex. Indeed, exactness at~$\Z$
  was already observed in \cref{Subsection:MonoidHomologicalTypeFP0}. Moreover,
  we have
  \[
    \varepsilon d_1[a] = \varepsilon (\cl{a}) - \varepsilon(1) = 1 - 1 = 0
  \]
  for every $1$-generator~$a$ of~$P$. In order to prove exactness at
  $\freemod{\Z M}{P_0}$, as explained in \cref{SS:ResolutionsHomologicalSyzygies}, we construct contracting homotopies
\[
i_0 : \Z \fl \freemod{\Z M}{P_0}
\qquad\text{and}\qquad
i_1 : \freemod{\Z M}{P_0}\fl \freemod{\Z M}{P_1}
\]
  as follows. The morphism $i_0$ is simply defined by $i_0(1)=1$ and
  extended by linearity. As for $i_1$ we first need to extend the
  bracket map~$[-]:P_1\to\freemod{\Z M}{P_1}$ to a map  $[-]:\freecat
  P_1\to\freemod{\Z M}{P_1}$. This is done by induction on the length
  of the words in $\freecat P_1$ by setting
\[
[1]= 0
\qquad\text{and}\qquad
[aw] = [a]+\cl{a}[w],
\]
for $a\in P_1$ and $w\in\freecat P_1$ (technically, we extend the map as a
  derivation, see \cref{sec:derivation}).  It follows that the
  equation
  \begin{equation}
    \label{eq:derivation1}
    d_1([w]) = \cl{w} - 1
  \end{equation}
   holds for all elements $w$ of
  $\freecat P_1$, not just for generators. We reason by induction on the
  length of the words in $\freecat P_1$. One first has
  $d_1([1])=0=\cl{1}-1$. Let now $a\in P_1$ and $w\in\freecat P_1$
  such that the equation \eqref{eq:derivation1} holds for $w$. Then
\[
d_1([aw])=d_1([a])+\cl{a}d_1([w])=\cl{a}-1+\cl{a}(\cl{w}-1)=\cl{aw}-1.
\]
Now, we choose a section and define the morphism $i_1$ by setting 
\[
i_1(u)=[\rep{u}]
\]
and extending it by linearity.
Finally, for any $u\in M$, we have $i_0\varepsilon(u) = 1$ and
\[
    d_1i_1 (u) = d_1[\rep{u}] = \cl{\rep{u}}-1=u-1
    \pbox.
  \]
  Thus, $d_1 i_1+i_0\varepsilon = \id_{\Z M}$ and $i_0$, $i_1$ are
  a contracting homotopies.
\end{proof}

\section{Monoids Having Homological Type Left-\texorpdfstring{$\FP_2$}{FP2}}
\label{sec:FP2}

\subsection{Presentations and partial resolutions of length 2}
\label{Section:ReidemesterFoxJacobian}
Let~$P$ be a presentation of a monoid~$M$. We define a partial resolution of
length~$2$ of the trivial $\Z M$-module~$\Z$ by free $\Z M$-modules
\[
  \xymatrix{
    \freemod{\Z M}{P_2}
    \ar [r] ^-{d_2}
    &
    \freemod{\Z M}{P_1}
    \ar [r] ^-{d_1}
    & \freemod{\Z M}{P_0}
    \ar [r] ^-{\varepsilon}
    & \Z
    \ar [r]
    & 0. 
  }
\]
The morphisms $\varepsilon$ and $d_1$ are those defined in the previous
section. The morphism~$d_2$ is defined, on generators of  $\freemod{\Z M}{P_2}$, by
\[
  d_2([\alpha]) = [\src{}(\alpha)] - [\tgt{}(\alpha)]
  \pbox,
\]
for every $\alpha$ in $P_2$, and called the \emph{Reidemeister-Fox Jacobian}\index{Reidemeister-Fox Jacobian} of the
presentation~$P$.

\subsection{Normalization strategies}
\label{sec:normalization-strat}
Let $P$ be a $2$-polygraph with a given section.
A \emph{normalization strategy}\index{normalization strategy}~$\sigma$ for~$P$ is a map
\[
\sigma : P_1^\ast \to \freegpd P_2
\]
that sends every $1$-cell~$w$ of~$P_1^\ast$ to a $2$-cell
\[
\sigma(w) : w\To\rep{w}
\]
in~$\freegpd P_2$, such that $\sigma(\rep w)=\unit{\rep w}$ holds for every $1$-cell~$w$ of~$P_1^\ast$.
A normalization strategy~$\sigma$ is a \emph{left} (resp. \emph{right}) one if it also satisfies
\[
\sigma(wv) = \sigma(w)v\comp1\sigma(\rep{w}v)
\qquad
\text{\big(resp.} 
\quad
\sigma(wv) = w\sigma(v)\comp1\sigma(w\rep{v})
\;\;
\big)
\]
that is
\[
\sigma(wv) \:=
\vcenter{\xymatrix @C=3.5em @R=2.5em {
& \cdot  
	\ar@/^/ [dr] ^-{v}
	\ar@2 []!<0pt,-15pt>;[d]!<0pt,-15pt> ^-*+{\sigma(\rep{w}v)}
\\
\cdot
	\ar@/^5ex/ [ur] ^-{w} _-{}="src"
	\ar@/_/ [ur] _-{\rep{w}} ^-{}="tgt"
	\ar@2 "src"!<4pt,-6.5pt>;"tgt"!<-4pt,7.5pt> ^-{\sigma(w)}
	\ar@/_/ [rr] _-{\rep{wv}}
&& \cdot
}}
\quad
\big( \text{resp.\ } \sigma(wv) \:=
\vcenter{\xymatrix @C=3.5em @R=2.5em {
& \cdot  
	\ar@/^5ex/ [dr] ^-{v} _-{}="src"
	\ar@/_/ [dr] _-{\rep{v}} ^-{}="tgt"
	\ar@2 "src"!<-8pt,-8.5pt>;"tgt"!<0pt,3.5pt> ^-{\sigma(v)}
	\ar@2 []!<-22pt,-15pt>;[d]!<-22pt,4.5pt> ^-*+{\sigma(w\rep{v})}
\\
\cdot
	\ar@/^/ [ur] ^-{w}
	\ar@/_/ [rr] _-{\rep{wv}}
&& \cdot
}}
\;\big).
\]
A $2$-polygraph~$P$ always admits left and right normalization strategies. 
Let us prove this in the left case, the right case being treated in the same way.
Let us arbitrarily choose a $2$-cell $\sigma(wa):wa\To\rep{wa}$ in~$\freegpd P_2$, for every~$1$-cell~$w$ of~$P_1^\ast$ and every $1$-generator~$a$ of~$P$, such that~$\rep{w}=w$ and $\rep{wa}\neq wa$. Then we extend~$\sigma$ into a left normalization strategy by setting $\sigma(w)=1_w$ if $\rep{w}=w$ (which implies $\sigma(1)=1$), and
\[
\sigma(w) = \sigma(v)a \comp1 \sigma(\rep{v}a)
\]
if $\rep{w}\neq w$ and $w=va$ with~$v$ in~$P_1^\ast$ and~$a$ in~$P_1$.

\begin{proposition}
  \label{Proposition:LowSquierResolution}
  Let~$M$ be a monoid and let~$P$ be a presentation of~$M$. The sequence of
  $\Z M$-modules
  \[
    \xymatrix{
      \freemod{\Z M}{P_2}
      \ar [r] ^-{d_2}
      &
      \freemod{\Z M}{P_1}
      \ar [r] ^-{d_1}
      & \freemod{\Z M}{P_0}
      \ar [r] ^-{\varepsilon}
      & \Z
      \ar [r]
      & 0
    }
  \]
  is a partial free resolution of length~2 of~$\Z$.
\end{proposition}
\begin{proof}
  In \cref{Proposition:MonoidsFP1}, we have proved the exactness at $\Z$ and
  $\freemod{\Z M}{P_0}$, and exactness at $\freemod{\Z M}{P_1}$ remains to be
  shown. The equation $d_1d_2 = 0$ is a consequence of \cref{eq:derivation1}.
  Indeed, we have
  \[
 d_1d_2 [\alpha] = d_1[s(\alpha)] - d_1[t(\alpha)] = \cl{s(\alpha)} - \cl{t(\alpha)} = 0,
  \]
  for every $2$-generator~$\alpha$ of~$P$, where the last equality comes
  from the equality $\cl{\src{}(\alpha)}=\cl{\tgt{}(\alpha)}$, which holds because~$P$ is a
  presentation of~$M$.

  In order to prove the exactness at $\freemod{\Z M}{P_1}$, we construct a
  contracting homotopy of the complex. The morphisms of $\Z$-modules $i_0$ and
  $i_1$ are defined in the proof of \cref{Proposition:MonoidsFP1}, and the
  morphism of $\Z$-modules
  \[
    i_2 : \freemod{\Z M}{P_1} \fl \freemod{\Z M}{P_2}
  \]
  is defined by fixing a left normalization strategy~$\sigma$ for the
  $2$\nbd-poly\-graph~$P$. Namely, we define the morphism of $\Z$-modules~$i_2$ by its
  value on generic elements
  \[
    i_2(u[a]) = [\sigma(\rep{u}a)],
  \]
  where the bracket~$[-]$ is extended to every $2$-cell of the free
  $(2,1)$-category~$\freegpd{P}$ by the following relations
  \begin{align*}
    [\unit{u}] &= 0,
    &
    [u\phi v] &= \cl{u}[\phi],
    &
    [\phi\comp1\psi] &= [\phi]+[\psi],
  \end{align*}
  for all $1$-cells~$u$ and~$v$ and $2$-cells~$\phi$ and~$\psi$ of~$\freegpd{P}$
  such that the composite~$\phi\comp1\psi$ are defined.

  We have, on the one hand,
  \[
    i_1 d_1(u[a]) = i_1(u\cl{a} - u) =  [\rep{ua}]  - [\rep{u}]
  \]
  and, on the other hand, 
  \[
    d_2 i_2 (u[a]) = d_2[\sigma(\rep{u}a)] = [\rep{u}a] - [\rep{ua}] = u[a] + [\rep{u}] - [\rep{ua}]
    \pbox.
  \]
  For the equality in the middle, one proves that
  $d_2[\phi]=[\src{}(\phi)]-[\tgt{}(\phi)]$ holds for every $2$-cell~$\phi$
  of~$\freegpd{P}$ by induction on the size of~$\phi$. Hence we have
  \[
    d_2 i_2 + i_1d_1=\id_{\freemod {\Z M}{P_1}},
  \]
  thus concluding the proof.
\end{proof}

\noindent
The previous proposition allows us to deduce:

\begin{proposition}
  \label{Proposition:MonoidsFP2}
  If a monoid $M$ admits a finite presentation, then it has homological
  type left\nbd-$\FP_2$, and thus the group $H_2(M,\Z)$ is finitely generated.
\end{proposition}

\subsection{Homological 2-syzygies}
The kernel of the morphism $d_2$ defined in \cref{Section:ReidemesterFoxJacobian} is called the $\Z M$-module of \emph{homological 2-syzygies} of the $2$-poly\-graph~$P$.
Using natural systems as modules, we will establish in \cref{sec:H-iar} an isomorphism between the homological $2$-syzygies and the identities among relations for a $1$-category presented by a $2$-polygraph.

\section{Homological Type Left-\texorpdfstring{$\FP_3$}{FP3} and Confluence}
\label{sec:FP3}
\label{SquierHomologicalTheorem}

\subsection{Coherent presentations and partial resolutions of length 3}
Let~$P$ be a coherent presentation of a monoid~$M$. Let us extend the partial
resolution of \cref{Proposition:LowSquierResolution} into the resolution of
length~$3$
\[
  \xymatrix@C=2em{
    \freemod {\Z M}{P_3}
    \ar [r] ^-{d_3}
    &
    \freemod {\Z M}{P_2}
    \ar [r] ^-{d_2}
    &
    \freemod {\Z M}{P_1}
    \ar [r] ^-{d_1}
    &\freemod{\Z M}{P_0}
    \ar [r] ^-{\varepsilon}
    & \Z
    \ar [r]
    & 0\pbox.
  }
\]
The boundary map~$d_3$ is defined, for every $3$-cell~$A$ of~$P$, by
\[
  d_3[A] = [\src2(A)] - [\tgt2(A)]
  \pbox.
\]
The bracket notation~$[-]$ is extended to $3$-cells of~$\freegpd{P}$ by setting
\begin{align*}
  [uFv] &= \cl{u}[F]
  &
  [F\comp1 G] &= [F] + [G]
  &
  [F\comp2 G] &= [F] + [G]
\end{align*}
for all $1$-cells~$u$ and~$v$ and $3$-cells~$F$ and~$G$ of~$\freegpd{P}$ such
that the composites are defined. In particular, the latter relation implies
$[1_\phi]=0$ for every $2$-cell~$\phi$ of~$\freegpd{P}$. We check, by induction
on the size, that $d_3[F] = [\src2(F)] - [\tgt2(F)]$ holds for every
$3$-cell~$F$ of~$\freegpd{P}$.

\begin{proposition}
  \label{Proposition:LowSquierResolution2}
  Let~$P$ be a coherent presentation of a monoid~$M$. The sequence of
  $\Z M$-modules
  \[
    \xymatrix{
      \freemod {\Z M}{P_3}
      \ar [r] ^-{d_3}
      &
      \freemod {\Z M}{P_2}
      \ar [r] ^-{d_2}
      &
      \freemod {\Z M}{P_1}
      \ar [r] ^-{d_1}
      & \freemod{\Z M}{P_0}
      \ar [r] ^-{\varepsilon}
      & \Z
      \ar [r]
      & 0
    }
  \]
  is a partial free resolution of length~3 of~$\Z$.
\end{proposition}
\begin{proof}
  We proceed with the same notations as in the proof of
  \cref{Proposition:LowSquierResolution}, with the extra hypothesis
  that~$\sigma$ is a left normalization strategy for~$P$. This implies that
  $i_2(u[v])=[\sigma(\rep{u}v)]$ holds for all~$u$ in~$M$ and~$v$ in~$\freecat P_1$, by
  induction on the length of~$v$.
  We have $d_2d_3=0$ because $\src1\src2=\src1\tgt2$ and
  $\tgt1\src2=\tgt1\tgt2$. Then, we define the following morphism of
  $\Z$-modules $i_3 : {\freemod {\Z M}{P_2}} \fl {\freemod {\Z M}{P_3}}$ by
  setting, for $u\in M$ and $\alpha\in P_3$,
  \[
    i_3(u[\alpha]) =  [\sigma(\rep{u}\alpha)]
  \]
  where $\sigma(\rep{u}\alpha)$ is a $3$-cell of~$\freegpd{P}$ with the
  following shape, with $v=\src{}(\alpha)$ and $w=\tgt{}(\alpha)$:
  \[
    \xymatrix @!C @C=2em @R=2em {
      & {\rep{u}w}
      \ar@2@/^/ [dr] ^-{\sigma(\rep{u}w)}
      \\
      {\rep{u}v}
      \ar@2@/^/ [ur] ^-{\rep{u}\alpha}
      \ar@2@/_/ [rr] _-{\sigma(\rep{u}v)} ^-{}="tgt"
      && {\rep{uv}}
      \ar@3 "1,2"!<-2pt,-11pt>;"tgt"!<-2pt,6pt> ^-{\sigma(\rep{u}\alpha)}
    }
  \]
  Let us note that such a $3$-cell necessarily exists in~$\freegpd{P}$
  because~$P_3$ is an acyclic cellular extension of~$\freegpd{P}$. Then we have,
  on the one hand,
  \[
    i_2d_2(u[\alpha]) = i_2(u[v]-u[w]) = [\sigma(\rep{u}v)] - [\sigma(\rep{u}w)]
  \]
  and, on the other hand,
  \begin{align*}
    d_3i_3(u[\alpha]) 
    &= [\rep{u}\alpha \comp1 \sigma(\rep{u}w)] - [\sigma(\rep{u}v)]\\
    &= u[\alpha] + [\sigma(\rep{u}w)] - [\sigma(\rep{u}v)]
    \pbox.
  \end{align*}
  Hence $d_3 i_3 + i_2 d_2=\id_{\freemod{\Z M}{P_2}}$, concluding the proof.
\end{proof}

\begin{remark}
  The proof of \cref{Proposition:LowSquierResolution2} uses the fact that~$P_3$
  is an acyclic cellular extension to produce, for every $2$-cell~$\alpha$
  of~$P_2$ and every~$u$ in~$M$, a $3$-cell $\sigma(\rep{u}\alpha)$ with the
  required shape. The hypothesis on~$P_3$ could thus be modified to only require
  the existence of such a $3$-cell in~$\freegpd{P}$: however, it is proved
  in~\cite{GuiraudMalbos12advances} that this implies that~$P_3$ is an acyclic
  cellular extension.
\end{remark}

\noindent
The previous proposition has the following consequence, already noted
in~\cite{pride1995low}, \cite[Theorem~3.2]{cremanns1994finite}, and
\cite[Theorem~3]{lafont1995new}:

\begin{theorem}
\label{Theorem:TDFimpliesFP3}
Let~$M$ be a finitely presented monoid. If~$M$ has finite derivation type, then it has homological type left-$\FP_3$, and thus the group $H_3(M,\Z)$ is finitely generated.
\end{theorem}

\index{homological!finiteness condition}
\index{Squier!theorem!homological}
\noindent
By \cref{thm:FiniteConvergent=>FDT} and \cref{prop:FPn-H-fg}, this implies the following homological finiteness condition for finite convergence~\cite[Theorem~4.1]{squier1987word}:

\begin{theorem}
\label{Theorem:SquierAb1}
If a monoid $M$ admits a finite convergent presentation, then it has homological type left-$\FP_3$, and thus the group $H_3(M,\Z)$ is finitely generated.
\end{theorem}

\noindent
The construction of this chapter will be generalized in
\cref{chap:ConstructingResolutions} to produce a free resolution  of
infinite length,  involving $n$-fold critical branchings for every natural number $n$~(\cref{T:AbelianResolution}). 

\begin{example}
  \label{Example:Monoidaa}
  Consider the monoid~$M$ with the convergent presentation
  \[
    P
    =
    \Pres{\star}{a}{\mu : aa\dfl a}
    \pbox.
  \]
  Writing $\nf w$ for the normal form of a word~$w$, we have $\nf{w}=a$ for every
  non-identity $1$-cell $w\in\freecat P_1$. With the leftmost normalization
  strategy~$\sigma$, we get, writing the $2$-cell~$\mu$ as a string
  diagram~$\satex{mu-vsmall}$:
  \begin{align*}
    \sigma(a) &= 1_a
    &
    \sigma(aa) &= \satex{mu}
    &
    \sigma(aaa) &= \mu a \comp1 \mu = \satex{mon-assoc-l}
    \pbox.
  \end{align*}
  The presentation has exactly one critical branching, whose corresponding
  generating confluence can be written in two equivalent ways
  \[
    \vcenter{
      \xymatrix @!C @R=1.1em @C=1.7em {
        & aa
        \ar@2@/^/ [dr]^{\satex{mu}}
        \ar@3 []!<-15pt,-12.5pt>;[dd]!<-15pt,12.5pt> ^-*+{\satex{asso}}
        \\
        aaa 
        \ar@2@/^/ [ur]^{\satex{1mu}} 
        \ar@2@/_/ [dr]_{\satex{mu1}}
        && a 
        \\
        & aa
        \ar@2@/_/ [ur]_{\satex{mu}}
      }}
    \qquad\qquad\text{or}\qquad\qquad
    \vcenter{
      \xymatrix@C=3em{
        {\satex{mon-assoc-r}}
	\ar@3 [r] ^-*+{\satex{asso}}
        & {\satex{mon-assoc-l}}
      }
    }
    \pbox.
  \]
  The $\Z M$-module~$\ker d_2$ is generated by
  \begin{align*}
    d_3\spa{\satex{asso}}
    &=
    \spa{\satex{mon-assoc-r}}
    - 
    \spa{\satex{mon-assoc-l}}
    \\
    &=
    \spa{\satex{1mu}} + \spa{\satex{mu}}
    - 
    \spa{\satex{mu1}} - \spa{\satex{mu}}
    \\
    &= a \spa{\satex{mu}} - \spa{\satex{mu}}
    \pbox.
  \end{align*}
\end{example}

We will see in \cref{chap:ConstructingResolutions} that the construction of the
partial resolution in \cref{Proposition:LowSquierResolution2} can be extended to arbitrary length. We only provide here a
small generalization~\cite[Theorem~3.2]{squier1987word}, which is enough to
imply a negative answer to the universality of finite convergent rewriting, see
\cref{ex:lafont-proute-squier-homology}.

\subsection{A short exact sequence}
\cref{Theorem:SquierCompletion2polygraphs} states that any set $P_3$ of generating confluences of a convergent $2$-polygraph $P$, indexed by all its critical branchings, forms an acyclic extension of the $(2,1)$-category $\freegpd{P}$.
Following \cref{Proposition:LowSquierResolution2}, this induces a partial free resolution of length~$3$ of~$\Z$ by $\Z M$-modules
  \[
    \xymatrix{
      \freemod {\Z M}{P_3}
      \ar [r] ^-{d_3}
      &
      \freemod {\Z M}{P_2}
      \ar [r] ^-{d_2}
      &
      \freemod {\Z M}{P_1}
      \ar [r] ^-{d_1}
      & \freemod{\Z M}{P_0}
      \ar [r] ^-{\varepsilon}
      & \Z
      \ar [r]
      & 0.
    }
  \]
We have also shown in \cref{Subsection:TripleConfluences} that the triple generating confluences of $P$ generate the relations among the $3$-cells of the free $(3,1)$-category $\freegpd{(P,P_3)}$.
We will show in \cref{S:SquierResolution,S:Abelianization} that this allows us to extend the previous resolution with a boundary map $d_4 : \freemod {\Z M}{P_4} \to \freemod {\Z M}{P_3}$ defined on the free module generated by a set of $4$-chains $P_4$ indexed by generating triple confluences. In particular, when there are no critical triples, we recover the following result shown by Squier in {\cite[Theorem 3.2]{squier1987word}}, see also \cref{C:SyzygiesWithoutTriple}.

\begin{proposition}
  \label{Proposition:LowSquierResolution3}
  Suppose given a convergent 2-polygraph~$P$ without critical 3-branching, and write $P_3$ for a set of~2-spheres containing a confluence
  diagram for every critical branching of~$P$. Then the sequence of $\Z M$-modules
  \[
    \xymatrix@C=1.5em{
      0
      \ar [r]
      &
      \freemod {\Z M}{P_3}
      \ar [r] ^-{d_3}
      &
      \freemod {\Z M}{P_2}
      \ar [r] ^-{d_2}
      &
      \freemod {\Z M}{P_1}
      \ar [r] ^-{d_1}
      & \freemod{\Z M}{P_0}
      \ar [r] ^-{\varepsilon}
      & \Z
      \ar [r]
      & 0
    }
  \]
  is a partial resolution of length~4 of~$\Z$.
\end{proposition}

In the rest of this section, we show that this result turns out to be very useful for constructing examples of finitely presented monoids having an infinite third integral homology group while having a decidable word problem.

\begin{example}
  Consider the monoid $M$ presented by the $2$-polygraph
  \[
    P = \Pres\star{a,b,c}{\alpha_n : ac^nb \dfl 1}_{n\in\N}
    \pbox.
  \]
  The polygraph $P$ is convergent without critical branchings. Hence, by Squier's
  \cref{thm:SquierHomotopical} it can be extended into a coherent presentation
  with an empty set of $3$-generators. Following
  \cref{Proposition:LowSquierResolution2}, we have a partial free resolution of
  length~$3$ of~$\Z$ by free $\Z M$-modules:
  \[
    \xymatrix{
      0	
      \ar [r]
      & \freemod{\Z M}{P_2}
      \ar [r] ^-{d_2}
      & \freemod{\Z M}{P_1}
      \ar [r] ^-{d_1}
      & \freemod{\Z M}{P_0}
      \ar [r] ^-{\varepsilon}
      & \Z
      \ar [r]
      & 0
      \pbox.
    }
  \]
  We have
  \[
    \widetilde{d}_1(a)=\widetilde{d}_1(b)=\widetilde{d}_1(c)=0
  \]
  and 
  \[
    \widetilde{d}_2([\alpha_n]) = [a] + n[c] + [b]
  \]
  for all $n\geq 0$. As a consequence $H_1(M,\Z)=\Z$ and
  $\mathrm{H}_2(M,\Z)=\ker \widetilde{d}_2$ is the free $\Z$-module generated by
  \[
    [\alpha_n] - n[\alpha_1] + (n-1)[\alpha_0]
  \]
  for $n\geq 2$. Since $H_2(M)$ is not finitely generated, this shows that the finitely generated monoid $M$ cannot be finitely presented, by \cref{Proposition:MonoidsFP2,prop:FPn-H-fg}.
\end{example}

\begin{example}
  \label{ex:squier-lafont-monoid-homology}
  Consider the monoid $M$ presented by the following
  $2$-poly\-graph considered in~\cite{lafont1991church}:
  \[
    P=
    \Pres{\star}{a,b,c,d}{\alpha_0:ab\To a,\beta:da\To ac}.
  \]
  We have seen in \cref{ex:squier-lafont-monoid} that using the Knuth-Bendix
  completion procedure, this polygraph can be completed into the following
  convergent polygraph with infinitely many $2$-generators:
  \[
    \widetilde{P} = \pres{a,b,c,d}{\alpha_n : ac^nb \dfl ac^n, \beta : da \dfl ac}_{n\in\N}
    \pbox.
  \]
  There are infinitely many critical branchings, indexed by $n\in\N$:
  \[
    \xymatrix@!C @R=1.6em@C=1.5em {
      & ac^{n+1}b
      \ar@2@/^/ [dr] ^{\alpha_{n+1}} 
      \ar@3 []!<0pt,-20pt>;[dd]!<0pt,20pt> ^{A_n}
      \\
      dac^nb
      \ar@2@/^/ [ur] ^{\beta c^nb}
      \ar@2@/_/ [dr] _{d\alpha_n}
      && ac^{n+1}
      \\
      & dac^n 
      \ar@2@/_/ [ur] _{\beta c^n}
    } 
  \]
  Denoting by $\widetilde P_3$ the set of $3$-generators $\setof{A_n}{n\in\mathbb{N}}$, by
  \cref{thm:SquierHomotopical}, $\widetilde P_3$ extends $\widetilde{P}$ into a coherent
  presentation.  This system has no critical $3$-branching, thus by
  \cref{Proposition:LowSquierResolution3}, we have an exact sequence
  \[
    \xymatrix@C=1.5em{
      0	
      \ar [r]
      & \freemod{\Z M}{\widetilde P_3}
      \ar [r] ^-{d_3}
      & \freemod{\Z M}{\widetilde P_2}
      \ar [r] ^-{d_2}
      & \freemod{\Z M}{\widetilde P_1}
      \ar [r] ^-{d_1}
      & \freemod{\Z M}{\widetilde P_0}
      \ar [r] ^-{\varepsilon}
      & \Z
      \ar [r]
      & 0
      \pbox.
    }
  \]
To calculate homology groups of $M$, we consider the maps $\widetilde d_k := 1_{\Z}\otimes_{\Z M} d_k$ defined on $\Z[\widetilde P_k]$ and with values in $\Z[\widetilde P_{k-1}]$. We have
  \begin{align*}
    \widetilde d_1(a)&=\widetilde d_1(b)=\widetilde d_1(c)=\widetilde d_1(d)=0,
    \\
    \widetilde d_2([\alpha_n])&=[a]+n[c]+[b]-([a]+n[c])=[b],
    \\
    \widetilde d_2([\beta])&=[d]+[a]-([a]+[c])=[d]-[c],
    \\
    \widetilde d_3([A_n])&=[\beta]+[\alpha_{n+1}]-([\alpha_n]+[\beta])=[\alpha_{n+1}]-[\alpha_n].
  \end{align*}
  Thus
  \begin{align*}
    H_0(M,\Z)&=\Z,
    &
    H_1(M,\Z)&=\Z^2,
    &
    H_i(M,\Z)&=0,
    &
    \text{for $i=2,3$.}
  \end{align*}
At this stage, therefore, we cannot use the finiteness condition of \cref{Theorem:SquierAb1} to conclude the existence of a convergent presentation for the monoid $M$. As noted in {\cite[Section 3.5]{lafont1991church}}, we can nevertheless construct a finite convergent presentation of $M$ with another orientation of the rule $\beta$. Indeed, the following polygraph presents the monoid M and has no critical branching:
\[
\Pres{\star}{a,b,c,d}{\alpha_0:ab\To a,\gamma:ac\To da}.
\]
It is therefore trivially convergent.    
\end{example}

\begin{example}
  \label{ex:lafont-proute-squier-homology}
  Consider the monoid~$M$ of \cref{ex:lafont-proute-squier} presented by the
  following $2$-polygraph:
  \[
    P
    =
    \Pres\star{a,b,c,d,d'}{\alpha_0:ab\To a,\beta:da \To ac,\beta':d'a \To ac}
    \pbox.
  \]
  We have seen that, by using the Knuth-Bendix completion procedure it can be
  completed into an infinite convergent polygraph, from which we deduce the
  following coherent presentation
  \[
    \widetilde P=
    \PRes{\star}{a,b,c,d,d'}{\alpha_n, \beta, \beta'}{A_n,A'_n}_{n\in \N},
  \]
  with
  \begin{align*}
    \alpha_n&:ac^nb\To ac^n,
    &
    \beta&:da\To ac,
    &
    \beta'&:d'a\To ac,
  \end{align*}
  and
  \begin{align*}
    A_n&:\beta c^nb\comp1\alpha_{n+1}\TO d\alpha_n\comp1\beta c^n,
    &
    A'_n&:\beta' c^nb\comp1\alpha_{n+1}\TO d'\alpha_n\comp1\beta' c^n.
  \end{align*}
  There are no critical $3$-branching and thus by
  \cref{Proposition:LowSquierResolution3} we have a partial resolution of
  length~$4$
  \[
    \xymatrix@C=1.5em{
      0
      \ar [r]
      & \freemod{\Z M}{\widetilde P_3}
      \ar [r] ^-{d_3}
      & \freemod{\Z M}{\widetilde P_2}
      \ar [r] ^-{d_2}
      & \freemod{\Z M}{\widetilde P_1}
      \ar [r] ^-{d_1}
      & \freemod{\Z M}{\widetilde P_0}
      \ar [r] ^-{\varepsilon}
      & \Z
      \ar [r]
      & 0
      \pbox.
    }
  \]
  The computations are similar to those of
  \cref{ex:squier-lafont-monoid-homology}. The map $\widetilde d_1$ is zero on the $\Z$-module $\Z\widetilde P_1$, and we have
\[
    \widetilde d_2([\alpha_n])=[b],
\qquad
    \widetilde d_2([\beta])=[d]-[c],
\qquad
    \widetilde d_2([\beta'])=[d']-[c],
\]
\[
    \widetilde d_3([A_n])=[\alpha_{n+1}]-[\alpha_n],
 \qquad
    \widetilde d_3([A'_n])=[\alpha_{n+1}]-[\alpha_n].
\]
  We deduce that
  \begin{align*}
    H_0(M,\Z)&=\Z,
    &
    H_1(M,\Z)&=\Z^2,
    &
    H_2(M,\Z)&=0,
  \end{align*}
and the $\Z$-module $H_3(M,\Z)$ is freely generated by the infinite family
\[
([A_n]-[A'_n])_{n\geq 0}.
\]
Following~\cref{Theorem:SquierAb1}, we deduce that the monoid $M$ does not have a finite convergent presentation.
This example thus exhibits a finitely presented monoid, with a decidable word problem, which does not admit a finite convergent presentation.
It therefore illustrates the fact that string rewriting theory is not universal for deciding the word problem in monoids, see \cref{sec:universality}.
\end{example}

\begin{example}
  Consider the coherent presentation of the monoid $B_3^+$ given in
  \cref{ex:B3-KBS}. By \cref{Proposition:LowSquierResolution2}, it induces a
  resolution of the trivial $\Z B_3^+$-module $\Z$, from which we can compute
  the following homology groups:
  \begin{align*}
    H_0(B_3^+,\Z)&=H_1(B_3^+,\Z)=H_2(B_3^+,\Z)=\Z,
    &
    H_3(B_3^+,\Z)&=0.
  \end{align*}
\end{example}

\subsection{Remark}
Note that in combinatorial group theory several examples of finitely presented groups with a decidable word problem that do not have homological type $\FP_3$ were discovered before Squier's work on homology of monoids.
In particular, Stallings constructed in \cite{stallings1963finitely} a finitely
presented group whose $3$-dimensional homology group
with integer coefficients is not finitely generated and thus it does not have homological type left-$\FP_3$. 
  The group is presented by
  \[
    \Pres\star{a,b,c,d,e}{[d,a],[e,a],[d,b],[e,b],[\finv{a}d,c],[\finv{a}e,c],[\finv{b}a,c]}
  \]
  where $[u,v]$ denotes the relation $uv=vu$.
Bieri proved that this group has a decidable word problem \cite{bieri1976homological}.
It was not yet known that its word problem cannot be solved by the normal form algorithm.

\subsection{Remarks on other homological finiteness conditions}
\label{Section:AnotherHomologicalType}
In the definition of homological type left-$\FP_n$ for a
monoid~$M$ (\cref{Definition:LeftFPn}), changing left modules to right modules, bimodules, or natural systems
gives the definitions of the homological types \emph{right-$\FP_n$},
\emph{bi-$\FP_n$} and~$\FP_n$. We refer the
reader to~\cite[Section 5.2]{GuiraudMalbos12advances} for the relations between
these different finiteness conditions, see also \cref{SS:CategoriesFiniteHomologicalType}. 
In particular, for~$n=3$, all of these
homotopical conditions are consequences of the finite derivation type
property defined in \cref{sec:2fdt-def}. 
Moreover, all these homological finiteness properties are necessary conditions for finite convergence.
The proof is similar to the one for the left-$\FP_3$ property given
in \cref{SquierHomologicalTheorem}. 
In particular, for the
right-$\FP_3$ property, we consider right modules and, to get the contracting
homotopy, we construct a right normalization strategy~$\sigma$ by defining a
$3$-cell $\sigma(\alpha\rep{u})$ with shape
\[
  \xymatrix @!C @C=2em @R=2em {
    & {w\rep{u}}
    \ar@2@/^/ [dr] ^-{\sigma(w\rep{u})}
    \\
    {v\rep{u}}
    \ar@2@/^/ [ur] ^-{\alpha\rep{u}}
    \ar@2@/_/ [rr] _-{\sigma(v\rep{u})} ^-{}="tgt"
    && {\rep{vu}}
    \ar@3 "1,2"!<-2pt,-11pt>;"tgt"!<-2pt,6pt> ^-{\sigma(\alpha\rep{u})}
  }
\]
for any $2$-generator $\alpha:v\dfl w$ and element~$u$ in the monoid $M$.

To conclude this chapter, we summarize in the following theorem the properties of the family of monoids $S_k$, for $k\geq 0$, defined in \cref{ex:MonoidSk}, which is Squier's original example~\cite{squier1987word,squier1994finiteness}. This family illustrates the homological and homotopical finiteness conditions for the convergence of string rewriting systems studied in this and the previous chapter, see \cref{Theorem:SquierAb1,thm:FiniteConvergent=>FDT}.

\begin{theorem}
For $k\geq 1$, the monoid~$S_k$ is a finitely presented monoid that has the following properties.
\begin{enumerate}
\item It has a decidable word problem~\cite[Example 4.5]{squier1987word}.
\item For $k=1$, it does not have finite derivation type~\cite[Theorem 6.7]{squier1994finiteness}.
\item For $k=1$, it has homological type left-$\FP_\infty$~\cite[Example 4.5]{squier1987word}.
\item For $k\geq 2$, it does not have homological type left-$\FP_3$~\cite[Example 4.5]{squier1987word}.
\item It does not have a finite convergent presentation.
\end{enumerate}
\end{theorem}

Conditions 1, 2, and 5 are seen in \cref{Theorem:MonoidS1} for $k=1$.
\Cref{Theorem:TDFimpliesFP3} proves that finite derivation type implies homological type left-$\FP_3$. Conditions 2 and 3 on monoid $S_1$ prove that the converse implication is false in general.  Note, however, that in the special case of groups, the property of having finite derivation type is equivalent to the homological finiteness condition left-$\FP_3$~\cite{cremanns1996groups}. The latter result is based on the Brown-Huebschmann isomorphism between identities among relations and homological syzygies~\cite{BrownHuebschmann82}, see also \cref{Theorem:IsomorphismPiH2}.

\part{Diagram Rewriting}


\chapter{Three-Dimensional Polygraphs}
\label{chap:3pol}
We have seen in \cref{chap:lowdim} that $1$-polygraphs provide a notion of
presentation for sets and in \cref{chap:2pol} that $2$-polygraphs provide a
notion of presentation for categories. We go on climbing the dimensional ladder
and establish $3$-polygraphs as a notion of presentation for $2$-categories,
see \cref{sec:3pol}. As expected, those consist in generators for $0$-,
$1$-, and $2$-dimensional cells, together with relations between freely generated
$2$-cells, which are represented by generating $3$-cells. As
particular cases, let us mention the notions of presentation of
monoidal category (when there is only one $0$\nbd-gene\-rator) and of PRO (when
there is only one $0$\nbd-gene\-rator and one $1$-generator). This includes
$2$-categories encoding theories for fundamental algebraic structures such as
monoids, groups, etc. Note that in the point of view on
$(3,1)$-polygraphs we adopt here,
the $3$-cells encode relations, as opposed to~\cref{chap:2-Coherent}
where they encode coherences between relations.

Any $3$-polygraph induces an abstract rewriting system, so that all
general rewriting concepts still make sense in this setting:
confluence, termination, etc.
However, more specific tools have to be adapted to this context:
the notion of \emph{critical branching} is defined for $3$-polygraphs in
\cref{sec:3pol-rewr}, along with the proof that confluence of critical
branchings implies the local confluence of the polygraph (\cref{lem:3-cb-lc}). In the case where the
polygraph is terminating (techniques to show this will be presented in
\cref{chap:3term}), local confluence implies confluence, and we thus have
a systematic method to show the convergence of a $3$-polygraph.
When this is the case, normal forms give canonical representatives for
$2$-cells modulo the congruence generated by $3$-cells, and we explain how to
exploit this to show that a given $3$-polygraph is a presentation of a
given $2$-category in \cref{sec:3pol-pres}. There is however a major difference with the case of
$2$-dimensional polygraphs: a finite
convergent polygraph might give rise to an infinite number of critical
branchings~(\cref{sec:3pol-indexed}). This prevents us from making direct generalizations of homotopical or
homological finiteness conditions (\cref{chap:2-fdt,chap:2-homology}) from $2$-
to $3$-polygraphs. Finally, in \cref{sec:3-dlaws}, we provide some techniques for
combining presentations of $2$-categories and for building presentations
of $2$-categories in a modular way.

\section{Three-Dimensional Polygraphs}
\index{polygraph!3-@$3$-}
\index{3-polygraph@$3$-polygraph}
\label{sec:3pol}

\subsection{Definition}
A \emph{3-polygraph} $(P,P_3)$ consists of a $2$-polygraph~$P$ together with a
cellular extension~$P_3$ of the $2$-category~$\freecat{P}$ freely generated
by~$P$, the elements of~$P_3$ being referred to as \emph{3-generators}.
Explicitly, a $3$-polygraph consists of a diagram
\[
  \xymatrix{
    &\ar@<-.5ex>[dl]_-{\sce0}\ar@<.5ex>[dl]^-{\tge0}P_1\ar[d]^{\ins1}&\ar@<-.5ex>[dl]_-{\sce1}\ar@<.5ex>[dl]^-{\tge1}P_2\ar[d]^{\ins2}&\ar@<-.5ex>[dl]_-{\sce2}\ar@<.5ex>[dl]^-{\tge2}P_3\\
    P_0&\ar@<-.5ex>[l]_-{\freecat{\sce0}}\ar@<.5ex>[l]^-{\freecat{\tge0}}\freecat{P_1}&\ar@<-.5ex>[l]_-{\freecat{\sce1}}\ar@<.5ex>[l]^-{\freecat{\tge1}}\freecat{P_2}
  }
\]
in $\Set$, together with a structure of $2$-category on the $2$-graph
\[
  \xymatrix{
    P_0&\ar@<-.5ex>[l]_-{\freecat{\sce0}}\ar@<.5ex>[l]^-{\freecat{\tge0}}\freecat{P_1}&\ar@<-.5ex>[l]_-{\freecat{\sce1}}\ar@<.5ex>[l]^-{\freecat{\tge1}}\freecat{P_2}&{}\phantom{P_3}
  }
\]
such that, for $n\in\set{0,1,2}$,
\begin{itemize}
\item $\freecat{P_n}$ is the set of $n$-cells of the $n$-category freely
  generated by the underlying $n$-polygraph,
\item $\ins{n}:P_n\to\freecat{P_n}$ is the
  canonical inclusion,
\item $\freecat{\sce{n}}$ and~$\freecat{\tge{n}}$ are the respective canonical 
  extensions of~$\sce{n}$ and~$\tge{n}$, satisfying
  \[
    \freecat{\sce{n}}\circ\ins{n}=\sce{n}
    \qquad\text{and}\qquad
    \freecat{\tge{n}}\circ\ins{n}=\tge{n},
  \]
\item the globular identities are satisfied:
  \[
    \freecat{\sce{n}}\circ\sce{n+1}=\freecat{\sce{n}}\circ\tge{n+1}
    \qquad\text{and}\qquad
    \freecat{\tge{n}}\circ\sce{n+1}=\freecat{\tge{n}}\circ\tge{n+1}.
  \]
\end{itemize}
We write $A : \phi \TO \psi$ for a $3$-generator $A\in P_3$ with $\sce2(A)=\phi$
and $\tge2(A)=\psi$. A $3$-polygraph~$P$ is often concisely denoted
\[
  \PRes{P_0}{P_1}{P_2}{P_3}\pbox.
\]
\index{underlying!2-polygraph}
We write~$\tpol2{P}$ for the \emph{underlying $2$-polygraph} of a
$3$-polygraph~$P$. Contrarily to previous chapters, we always respectively
denote by~$\comp0$ and~$\comp1$ the horizontal and vertical compositions of a
$2$-category.

\begin{example}
  \index{polygraph!of monoids}
  \label{ex:ass}
  The $3$-polygraph~$\MonPoly$ is
    \[
    \MonPoly
    =
    \PRes{\star}{a}{\mu:a\comp0 a\To a,\, \eta:\unit\star\To a}{A,L,R}
  \]
  where the sources and targets of the $3$-generators are given by
  \begin{align*}
    A&:(\mu\comp0 a)\comp1\mu\TO(a\comp0\mu)\comp1\mu
    &
    L&:(\eta\comp0 a)\comp1\mu\TO a
    \\
    &&
    R&:(a\comp0\eta)\comp1\mu\TO a    
  \end{align*}
  Using string diagrams (see \secr{string-diag}), the $2$-generators
  of~$\MonPoly$ are pictured as
  \begin{align*}
    \mu&=\satex{mu}
    &
    \eta&=\satex{eta}
  \end{align*}
  and its $3$-generators~$A$, $L$, and~$R$ respectively as
  \begin{align*}
    \satex{A2}
    &:
    \satex{mon-assoc-l}
    \TO
    \satex{mon-assoc-r}
    &
    \satex{L2}
    &:
    \satex{mon-unit-l}
    \TO
    \satex{mon-unit-c}
    &
    \satex{R2}
    &:
    \satex{mon-unit-r}
    \TO
    \satex{mon-unit-c}
  \end{align*}
\end{example}

\noindent
Many other examples of~$3$-polygraphs are given in \chapr{3ex}.

\subsection{Presented 2-category}
\index{presentation!of a 2-category}
Let~$P$ be a $3$-polygraph. The \emph{2-category presented by~$P$} is the
$2$-category, denoted by~$\pcat{P}$, obtained by quotienting the free
$2$-category over~$\tpol2{P}$ by the congruence $\approx^P$ generated by~$P_3$
on $2$-cells, as described in~\secr{2-cat-quot}:
\[
  \pcat{P}
  =
  \freecat{\tpol2{P}}/P_3
\]
If~$C$ is a $2$-category, we say that~$P$ \emph{presents~$C$} if~$C$ is
isomorphic to~$\pcat{P}$.

In particular, when the set $P_0$ is reduced to one element, the category presented
by~$P$ has one $0$-cell and is thus a strict monoidal category (see \cref{sec:strict-moncat}). Moreover,
when both~$P_0$ and~$P_1$ are reduced to one element, the set $\freecat P_1$ of $1$-cells is the free
monoid on one generator, \ie $\N$, and the presented category is a PRO (see \cref{sec:PRO}). This is
for instance the case in \cref{ex:ass}.

\index{Tietze!equivalence!of 3-polygraphs}
\index{equivalence!Tietze}
Two $3$-polygraphs $P$ and~$Q$ are said \emph{Tietze equivalent} when the
presented $2$-categories are isomorphic: $\pcat P\isoto\pcat Q$.

\subsection{Models}
\index{model!of a 2-category}
Given a $2$-category~$C$, the category of \emph{models} (or \emph{algebras})
of~$C$ in a $2$\nbd-cate\-gory~$S$ is the category $\nCat2(C,S)$ of $2$-functors
$C\to S$ and oplax $2$-natural transformations between those (see
\cref{oplax-transformation} for a general definition).
More explicitly, given a $3$-polygraph~$P$, a \emph{model} of~$\pcat{P}$ in a
$2$-category~$C$ consists of
  \begin{itemize}
  \item a family
    \[
      (f_x)_{x\in P_0}
    \]
    of $0$-cells of~$C$ indexed by the $0$-generators of~$P$,
  \item a family
    \[
      (f_a:f_x\to f_y)_{a:x\to y\in P_1}
    \]
    of~$1$-cells of~$C$ indexed by the $1$-generators of~$P$; if
    $u=a_1\ldots a_n$ is a $1$-cell of $\freecat{P}$, we write $f_u$ for the
    $1$-cell $f_{a_1}\ldots f_{a_n}$,
  \item a family
    \[
      (f_\alpha:f_u\To f_v)_{\alpha:u\To v\in P_2}
    \]
    of~$2$-cells of~$C$ indexed by the $2$-generators of~$P$; the notation
    $f_{\phi}$ is extended to any $2$-cell of~$\freecat{P}$ by
    $f_{\phi\comp0\psi}=f_\phi\comp0 f_\psi$,
    $f_{\phi\comp1\psi}=f_\phi\comp1f_\psi$ and $f_{\unit{u}}=\unit{f_u}$,
  \end{itemize}
  such that, for every $3$-generator
  $
  A
  :
  \phi
  \TO
  \psi
  $
  of $P$, we have
  \[
    f_\phi
    =
    f_\psi
    \pbox.
  \]

\begin{example}
  \index{polygraph!of monoids}
  \label{ex:ass-alg}
  \label{ex:monoid}  
  Let $C$ be a monoidal category. The models of the $3$-polygraph
  $\MonPoly$ of \cref{ex:ass} in~$C$ (considered as a $2$-category with only one
  $0$-cell) are precisely monoids in~$C$ in the following sense.
  A \emph{monoid}\index{monoid} in a monoidal category~$C$ consists of
  an
  object~$x$ of~$C$ and two morphisms
  \begin{align*}
    m&:x\otimes x\to x
    &
    e&:i\to x
  \end{align*}
  such that the diagrams
  \begin{align*}
    \xymatrix{
      x\otimes x\otimes x\ar[d]_{x\otimes m}\ar[r]^{m\otimes x}&x\otimes x\ar[d]^m\\
      x\otimes x\ar[r]_m&x
    }
    &&
    \xymatrix{
      x\otimes x\ar[dr]_m&\ar[l]_-{e\otimes x}x\ar[d]^{\unit{x}}\ar[r]^-{x\otimes e}&x\otimes x\ar[dl]^m\\
      &x
    }
  \end{align*}
  commute. A morphism $f:(x,m,e)\to(x',m',e')$ between monoids is a morphism
  $f:x\to x$ of~$C$ such that the diagrams
  \begin{align*}
    \xymatrix{
      x\otimes x\ar[d]_m\ar[r]^{f\otimes f}&x'\otimes x'\ar[d]^{m'}\\
      x\ar[r]_f&x'
    }
    &&
    \xymatrix{
      i\ar[d]_e\ar[r]^{\unit{i}}&i\ar[d]^{e'}\\
      x\ar[r]_f\ar[r]_f&x'
    }
  \end{align*}
  commute.
  In particular, a monoid in~$\Set$ equipped with cartesian product as
  tensor product and terminal set as unit is precisely a monoid in the usual
  sense (by Mac Lane's coherence theorem, \cref{thm:moncat-coh}, we can always
  consider that it forms a strict monoidal category).

  As another application, the algebras of $\MonPoly$ in the $2$-category~$\Cat$
  (of categories, functors, and natural transformations) are precisely
  monads\index{monad}.
\end{example}

\subsection{2-categories admitting a presentation}
\label{sec:2-cat-no-pres}
We should first note that, contrarily to the case of $1$-categories
(\cref{sec:2-std-pres}), not every $2$-category admits a presentation by a
$3$-polygraph. Namely, any $2$-category~$C$ is presented by a
$3$-polygraph~$P$ is a quotient of the free
$2$-category~$\freecat{\tpol2{P}}$ by the congruence generated by~$P_3$. Since
there is no quotient on $1$-cells, the underlying $1$-category of~$C$ is always
free (on the underlying $1$-polygraph of~$P$). Therefore, only $2$-categories
whose underlying $1$-category is free may have a presentation.

For instance, the category~$\Z$ corresponding to the additive monoid of
integers is not free, by \cref{sec:Z-not-free}.  Therefore
the $2$-category with~$\Z$ as underlying category and only identity $2$-cells
admits no presentation by a $3$-polygraph. Extensions of the notion of
polygraph aimed at addressing this problem have been proposed 
in~\cite{curien2017coherent,dupont2022coherent,mimram:LIPIcs.FSCD.2023.16}.

\subsection{The canonical and standard presentations}
\index{canonical presentation}
\index{standard!presentation}
\index{presentation!canonical}
\index{presentation!standard}
Any $2$\nbd-cate\-gory~$C$ whose underlying category is freely generated by a
$1$-polygraph~$Q$ admits a presentation by a $3$\nbd-polygraph.
The \emph{canonical presentation of~$C$} is the $3$\nbd-poly\-graph~$P$ with
\begin{itemize}
\item $Q$ as underlying $1$-polygraph,
\item the set $P_2=C_2$ of all $2$-cells of~$C$ as $2$-generators,
\item the subset $P_3$ of $\freecat{P_2}\times\freecat{P_2}$ of pairs of
  parallel $2$-cells whose evaluation as $2$-cells of~$C$ are
  equal.
\end{itemize}
An analogous of the standard presentation (\cref{sec:2-std-pres}) can also be
defined for $2$-categories, giving rise to slightly smaller presentations.

\subsection{The category of 3-polygraphs}
\label{sec:Pol3}
A \emph{morphism}~$f:P\to Q$ between $3$-polygraphs~$P$ and~$Q$
consists of a morphism $f:\tpol2 P\to\tpol2 Q$ between the underlying
$2$-polygraphs (see \cref{sec:Pol2}) together with a function $f_3:P_3\to Q_3$
such that $\src2\circ f_3=f_2\circ\src2$ and $\tgt2\circ f_3=f_2\circ\tgt2$.
These compose in the expected way, and we write
$\nPol 3$\nomenclature[Pol3]{$\nPol{3}$}{category of $3$-polygraphs} for the category of
$3$-polygraphs and their morphisms.

\section{Rewriting Properties of 3-Polygraphs}
\label{sec:3pol-rewr}
A $3$-polygraph~$P$ can be seen as a $3$-dimensional rewriting system: its
underlying $2$-polygraph generates a $2$-category, whose $2$\nbd-cells are the
``terms'' which get rewritten by the $3$-generators. For this reason, the
elements of~$P_3$ are sometimes called \emph{rewriting rules}. We now formalize
this point of view.

\subsection{Occurrences of 2-generators}
\index{occurrence!of a $2$-generator}
\label{sec:2-occurrences}
Let $P$ be a $2$-polygraph. Given a $2$-cell $\phi$ in $\freecat P_2$
and a $2$-gene\-rator $\alpha\in P_2$, we write $\sizeof[\alpha]\phi$ for the
number of occurrences of $\alpha$ in~$\phi$. It can be formally defined as
follows.

We write $N$ for the 2-category with one $0$-cell $\star$, one $1$-cell
$\unit\star$, $\N$ as set of $2$-cells, horizontal and vertical compositions
being given by addition and the identity 2-cell by $0$.
Given a $2$-generator $\alpha\in P_2$, there exists a unique $2$\nbd-func\-tor
\[
  \sizeof[\alpha]{{-}}:\freecat{P_2}\to N
\]
such that $\sizeof[\alpha]\alpha=1$ and $\sizeof[\alpha]\beta=0$ for every
$2$-generator $\beta\in P_2$ with $\beta\neq\alpha$. Given a 2-cell $\phi$ in
$\freecat{P_2}$, the natural number $\sizeof[\alpha]\phi$ is called the
\emph{number of occurrences} of the generator $\alpha$ in~$\phi$.
We also write
\[
  \sizeof{{-}}:\freecat{P_2}\to N
\]
for the 2-functor such that $\sizeof\alpha=1$ for every generator
$\alpha\in P_2$. Given a morphism~$\phi$, we have
\[
  \sizeof\phi=\sum_{\alpha\in P_2}\sizeof[\alpha]\phi
\]
and this quantity is called the \emph{number of generators} in~$\phi$ or the
\emph{size}\index{size} of~$\phi$.

\subsection{Contexts}
\index{context}
\label{sec:2-context}
Let $P$ be a 2-polygraph and $u,v:x\to y$ be two parallel 1-cells  in
$\freecat{P_1}$. We write $P[X]$ for the 2-polygraph with the same $0$- and
$1$-cells as $P$ and with $P_2\sqcup\set{X}$ as $2$-cells, with $\src1(X)=u$ and
$\tgt1(X)=v$.
A \emph{context}~$K$ of type $(u,v)$ in~$P$ is a cell $K$ in $\freecat{P[X]_2}$
in which the generator~$X$ occurs exactly once, \ie such that $\sizeof[X]{K}=1$.

\begin{lemma}
  \label{lem:2-context}
  Any context~$K$ of type~$(u,v)$ can be written in the form
  \[
    K=\psi\comp1(w\comp0 X\comp0 w')\comp1\psi'
  \]
  for some 1-cells $w:x'\to x$ and $w':y\to y'$ and 2-cells
  $\phi:u'\To wuw'$ and $\phi':wvw'\To v'$. Graphically, $K$ can be depicted as
  \[
    \vxym{
      \\
      x'\ar[r]^w
      \ar@/^11ex/[rrrr]^{u'}\ar@/_11ex/[rrrr]_{v'}
      \ar@{{}{ }{}}@/^8ex/[rrrr]|{\psi\Longdownarrow\phantom{\psi}}
      \ar@{{}{ }{}}@/_8ex/[rrrr]|{\psi'\Longdownarrow\phantom{\psi'}}
      &
      x\ar@/^5ex/[rr]_{u}\ar@/_5ex/[rr]^{v}
      \ar@{}[rr]|{\displaystyle X\Longdownarrow\phantom{X}}
      &&y\ar[r]^{w'}&y'\pbox.
      \\
      \\
    }
  \]
\end{lemma}

Given a context~$K$ as in previous lemma and a $2$-cell $\phi:u\To v$ in
$\freecat{P_2}$, we write $C[\phi]$ for the following $2$-cell of
$\freecat{P_2}$:
\[
  K[\phi]=\psi\comp1(w\comp0\phi\comp0 w')\comp1\psi'\pbox.
\]

\subsection{Rewriting steps}
\index{rewriting!step}
Let $P$ be  a $3$-polygraph. Given a context~$K$ in~$P_2$ of type $(u,v)$
as in \lemr{2-context} and a $3$-cell
\[
  F:\phi\TO\phi':u\To v:x\to y
\]
in $\freecat{P_3}$, we extend the previous notation and write
\[
  K[F]=\psi\comp1(w\comp0 F\comp0 w')\comp1\psi'
\]
Graphically,
\[
  \vxym{
    \\
    x'\ar[r]^w
    \ar@/^11ex/[rrrr]^{u'}\ar@/_11ex/[rrrr]_{v'}
    \ar@{{}{ }{}}@/^8ex/[rrrr]|{\psi\Longdownarrow\phantom{\psi}}
    \ar@{{}{ }{}}@/_8ex/[rrrr]|{\psi'\Longdownarrow\phantom{\psi'}}
    &
    x\ar@/^5ex/[rr]_u\ar@/_5ex/[rr]^v
    \ar@{}[rr]|{\displaystyle\phi\phantom{'}\Longdownarrow\overset F\TO\Longdownarrow\phi'}
    &&y\ar[r]^{w'}&y'\pbox.
    \\
    \\
  }
\]

A \emph{rewriting step} is a $3$-cell of the form $K[A]$ for some context~$K$
and $3$\nbd-gene\-rator~$A:\phi\TO\phi'$. The $2$-cells $K[\phi]$ and $K[\phi']$
are respectively called the source and target of the rewriting step.
Using the axioms of $3$-categories, one shows that every $3$-cell of~$P^*$ is
a composite of rewriting steps:

\begin{lemma}
  Any 3-cell $F$ of $\freecat{P}$ can be decomposed as
  \[
    F
    =
    F_1\comp2 F_2\comp2\ldots\comp2 F_k
  \]
  where the $F_i$'s are rewriting steps. Moreover, the number $k\in\N$ is the
  same for all such decompositions.
\end{lemma}

\noindent
The number~$k$ in the previous lemma is called the \emph{length} of~$F$.

\subsection{Termination and confluence}
Given a $3$-polygraph~$P$, we write here~$P^{\mathrm{rs}}$ for its set of
rewriting steps and $s_2,t_2:P^{\mathrm{rs}}\to\freecat{P_2}$ for the
functions respectively 
taking a rewriting step to its source and target.
\index{abstract rewriting system!of a 3-polygraph}
Any $3$-polygraph~$P$ induces an abstract rewriting system
\[
  \xymatrix{
    \freecat{P_2}&\ar@<-.5ex>[l]_-{\sce2}\ar@<.5ex>[l]^-{\tge2}P^{\mathrm{rs}}
  }
\]
 with $2$-cells as vertices and rewriting steps as edges.
A $3$\nbd-poly\-graph is said to be
\emph{terminating}, \emph{Church-Rosser}, \emph{confluent}, \emph{locally
  confluent}, or \emph{convergent} when the associated abstract rewriting system
is, see \cref{sec:ars}.

The termination of a $3$-polygraph can be shown in a similar way as for
$2$\nbd-poly\-graphs (\secr{2-red-order}) by considering suitable reduction
orders: this will be detailed in \chapr{3term}.
In the next section, we study local
confluence through critical branchings, generalizing the techniques introduced
for $2$-polygraphs.

\subsection{Branchings}
Let~$P$ be a $3$-polygraph. A \emph{branching}\index{branching} is a
branching of the underlying abstract rewriting system. It consists of a pair
$(F_1,F_2)$ of $3$-cells
$F_1:\phi\TO\psi_1$ and $F_2:\phi\TO\psi_2$
in $\freecat{P}$ with the same source. We say that the $2$-cell~$\phi$ is the
\emph{source} of $(F_1,F_2)$. A branching $(F_1,F_2)$ of~$P$ is
\emph{local}\index{local branching}\index{branching!local} if both~$F_1$
and~$F_2$ are rewriting steps; it is \emph{confluent}\index{branching!confluent}
if there exist cofinal $3$-cells $F_1':\psi_1\TO\chi$ and $F_2':\psi_2\TO\chi$
in~$\freecat{P}$.

We say that a local branching~$(F_1,F_2)$ of~$P$ is \emph{trivial} if
$F_1=F_2$. We say that the branching $(F_1,F_2)$ (\resp $(F_2,F_1)$) is
\emph{orthogonal} if $F_1$ and $F_2$ are of the form
\begin{align*}
  F_1
  &=
  \psi
  \comp1
  (w_1\comp0 A_1\comp0 w_1')
  \comp1
  \psi'
  \comp1
  (w_2\comp0\phi_2\comp0 w_2')
  \comp1
  \psi''
  \\
  F_2
  &=
  \psi
  \comp1
  (w_1\comp0\phi_1\comp0 w_1')
  \comp1
  \psi'
  \comp1
  (w_2\comp0 A_2\comp0 w_2')
  \comp1
  \psi''
\end{align*}
where $A_1:\phi_1\TO\phi_1'$ and $A_2:\phi_2\TO\phi_2$ are
$3$-generators,~$w_1$ , $w_1'$, $w_2$, and~$w_2'$ are $1$-cells of~$\freecat{P}$, and~$\psi$, $\psi'$, and~$\psi''$ are $2$-cells of~$\freecat{P}$. The orthogonal situation can be pictured as 
\[
\xymatrix@C=8ex@R=3.5ex{
  \\
  \ar@{{}{ }{}}@/^6ex/[rrrr]|{\psi\Longdownarrow\phantom{\psi}}
  &
  x_1\ar@/^4ex/[rr]_{u_1}\ar@/_4ex/[rr]^{v_1}
  \ar@{}[rr]|{\displaystyle\phi_1\phantom{'}\Longdownarrow\overset{A_1}\TO\Longdownarrow\phi_1'}
  &&
  y_1\ar[dr]_{w_1'}
  &
  \\
  x
  \ar@{{}{ }{}}[rrrr]|{\psi'\Longdownarrow\phantom{\psi'}}
  \ar@/^14.5ex/[rrrr]^u
  \ar@/_14.5ex/[rrrr]_v
  \ar[ur]_{w_1}
  \ar[dr]^{w_2}
  &&&&y\pbox.
  \\
  \ar@{{}{ }{}}@/_6ex/[rrrr]|{\psi''\Longdownarrow\phantom{\psi''}}
  &
  x_2\ar@/^4ex/[rr]_{u_2}\ar@/_4ex/[rr]^{v_2}
  \ar@{}[rr]|{\displaystyle\phi_2\phantom{'}\Longdownarrow\overset{A_2}\TO\Longdownarrow\phi_2'}
  &&
  y_2\ar[ur]^{w_2'}
  &
}
\]

Local branchings are ordered by the relation~$\sqsubseteq$ generated by
\[
(F_1, F_2) \sqsubseteq (\phi \comp{i} F_1, \, \phi \comp{i} F_2) 
\quad\text{and}\quad
(F_1, F_2) \sqsubseteq (F_1 \comp{i} \phi, \, F_2 \comp{i} \phi),
\]
where~$\phi$ ranges over the $2$-cells of~$\freecat{P}$ and~$i$ over~$\set{0,1}$
such that the involved composites are defined. A branching of~$P$ is called
\emph{critical} if it is local, orthogonal, and minimal for~$\sqsubseteq$. We say
that~$P$ is \emph{critically confluent} if all its critical branchings are
confluent.

As in the case of $2$-polygraphs, \cref{L:CriticalBranchingLemma2Pol},
the critical branching lemma holds for $3$-polygraphs:
\begin{lemma}
  \label{lem:3-cb-lc}
  \label{L:CriticalBranchingLemma3Pol}
  A 3-polygraph is locally confluent if and only if all its critical branchings
  are confluent.
\end{lemma}

\noindent
As a direct corollary of this lemma and Newman's lemma (\cref{lem:newman}), we
may state the following proposition, which is used to show the
convergence of polygraphs in the vast majority of cases.

\begin{proposition}
  \label{prop:3-term-cp}
  A 3-polygraph which is terminating and has all its critical branchings
  confluent is convergent.
\end{proposition}

\begin{example}
  \label{ex:ass-cp}
  The $3$-polygraph $\MonPoly$ of \cref{ex:ass} has five critical
   branchings. All of them are confluent, as shown by the string diagrams below:
  \[
        \vcenter{
          \xymatrix @R=1em @C=1em {
            &{ \satex{mon-cp1-1} }
            \ar@3 [rr] 
            ^-{\satex{mon-cp1-12}} 
            _-{}="1" 
            && { \satex{mon-cp1-2} }
            \ar@3@/^/ [dr] ^-{\satex{mon-cp1-23} }
            \\
            { \satex{mon-cp1-0} }
            \ar@3@/^/ [ur] ^-{\satex{mon-cp1-01}}
            \ar@3@/_/ [drr] _-*+{\satex{mon-cp1-04}}
            &&&& { \satex{mon-cp1-3} }
            \\
            && { \satex{mon-cp1-4} }
            \ar@3@/_/ [urr] _-*+{\satex{mon-cp1-43}}
          }
        }
  \]
  \begin{align*}
      \vcenter{
        \xymatrix @R=1em @C=1em {
          & {\satex{mon-cp2-1}}
          \ar@3@/^/ [dr] ^-{\satex{mon-cp2-12}}
          \\
          {\satex{mon-cp2-0}}
          \ar@3@/^/ [ur] ^-{\satex{mon-cp2-01}}
          \ar@3@/_/ [rr] _-{\satex{mon-cp2-02}} ^-{}="1"
          && {\satex{mon-cp2-2}}
        }
      }
      &&
      \vcenter{
        \xymatrix @R=1em @C=1em {
          & { \satex{mon-cp3-1} }
          \ar@3@/^/ [dr] ^-{\satex{mon-cp3-01}}
          \\
          { \satex{mon-cp3-0} }
          \ar@3@/^/ [ur] ^-{\satex{mon-cp3-12}}
          \ar@3@/_/ [rr] _-{\satex{mon-cp3-02}} ^-{}="1"
          && { \satex{mon-cp3-2} }
        }
      }
    \end{align*}
    \begin{align*}
      \vcenter{
        \xymatrix @R=1em @C=1em {
          & {\satex{mon-cp4-1}}
          \ar@3@/^/ [dr] ^-{\satex{mon-cp4-12}}
          \\
          {\satex{mon-cp4-0}}
          \ar@3@/^/ [ur] ^-{\satex{mon-cp4-01}}
          \ar@3@/_/ [rr] _-{\satex{mon-cp4-02}} ^-{}="1"
          && {\satex{mon-cp4-2}}
        }
      }
      &&
      \vcenter{
        \xymatrix @C=1.6em {
          {\satex{mon-cp5-0}}
          \ar@3@/^4.5ex/ [rr] ^{\satex{mon-cp5-01}} _{}="1"
          \ar@3@/_4.5ex/ [rr] _{\satex{mon-cp5-10}} ^{}="2"
          && {\satex{mon-cp5-1}}
        }
      }    
  \end{align*}
  For instance, the source of the first critical branching is
  \[
    \satex{mon-cp1-0} = (\mu\comp0 a\comp0 a)\comp1(\mu\comp0 a)\comp1\mu
  \]
  which can be rewritten by using the $3$-generator
  \[
    \satex{A2}:\satex{mon-assoc-l}\TO\satex{mon-assoc-r}
  \]
  in two ways, yielding the two rewriting steps:
  \begin{align*}
    \satex{mon-cp1-01}& = (A\comp0 a)\comp1\mu
    &
    \satex{mon-cp1-04}& = (\mu\comp0 a\comp0 a)\comp1 A
    \pbox.
  \end{align*}
  We will see in \cref{ex:ass-term} that this polygraph is terminating, based on
  the observation that rewriting either removes $2$-generators or moves subtrees
  to the right. By \cref{prop:3-term-cp}, it is thus convergent.
\end{example}

\subsection{The Knuth-Bendix completion procedure}
\index{Knuth-Bendix completion}
\index{completion}
\index{procedure!completion}
A Knuth-Bendix completion procedure can be defined for $3$-polygraphs. We do not
detail it much, because it is very similar to the one for $2$-polygraphs given
in \cref{sec:2kb}: starting with a $3$-polygraph~$P$ and a total reduction order adapted to the polygraph
(\cref{def:3-red-order}), we compute the critical branchings and, for each of
those branchings, normalize both members, and add a new rule when the normal
forms differ, oriented according to the reduction order. As for $2$-polygraphs,
this procedure might not terminate because we keep on adding new rules. However,
there is a new potential source of non-termination for $3$-polygraphs: we will
see in \cref{sec:3pol-indexed} that a finite polygraph might give rise to an
infinite number of critical branchings, and the completion procedure will have
to examine each of them.

\section{Constructing Presentations}
\label{sec:3pol-pres}
When a $3$-polygraph~$P$ is convergent, the normal forms provide canonical
representatives of equivalence classes of $2$-cells in $\freecat P_2$ modulo the
congruence generated by~$P_3$. We explain here that this can be exploited to
show that~$P$ presents a given $2$-category~$C$, by showing
that the $2$-cells of~$C$ are in bijection with the normal forms of the
polygraph. The following proposition thus generalizes the method proposed
in~\secr{pres-method} to construct presentations of categories.

\begin{proposition}
  \label{sec:3-pres-method}
  \label{prop:3-pres-method}
  Let $P$ be a convergent 3-polygraph and $C$ a
  2-category whose underlying category is isomorphic to the category
  freely generated by the underlying 1-polygraph of~$P$, \ie we have
  isomorphisms $f_0:P_0\to C_0$ and $f_1:\freecat{P_1}\to C_1$.
  Let moreover
  $
  f_2:P_2\to C_2
  $
  be a function compatible with source and target, \ie $f_2$ sends a 2-generator
  $\alpha:u\To v$ in $P_2$ to a 2-cell $f_2(\alpha):f_1(u)\To f_1(v)$. We
  extend $f_2$ to the 2-cells in $\freecat P_2$ by functoriality by
  \begin{align*}
    f_2(\phi\comp0\psi)&=f_2(\phi)\comp0 f_2(\psi)\pbox,
    \\
    f_2(\phi\comp1\psi)&=f_2(\phi)\comp1 f_2(\psi)\pbox,
     \\                    
    f_2(\unit{u})&=\unit{f_1(u)}\pbox.
  \end{align*}
  Suppose finally that
  \begin{itemize}
  \item for any 3-generator $A:\phi\TO\psi$ in~$P_3$, we have
    $f_2(\phi)=f_2(\psi)$,
  \item the function $\freecat{f_2}:\freecat{P_2}\to C_2$ restricts to a
    bijection between normal forms in $\freecat{P_2}$ and $C_2$.
  \end{itemize}
  Then~$P$ is a presentation of~$C$.
\end{proposition}
\begin{proof}
  Let us write $\fgf C$ for the $2$-polygraph with $P$ as underlying
  $1$-polygraph and whose set of $2$-generators is the set $C_2$ of $2$-cells
  of~$C$. The triple of morphisms $(f_0,f_1,f_2)$ precisely corresponds to a
  morphism of $2$-polygraphs $f:\tpol2 P\to\fgf C$. This morphism induces a
  functor $\freecat f:\freecat{\tpol2 P}\to C$ from the freely generated
  $2$-category. The first condition ensures that it induces a quotient functor
  $\pcat f:\pcat P\to C$, and the second condition ensures that $\pcat f$ is a
  bijection. Namely, $\pcat f$ is a bijection in dimensions $0$ and $1$ by
  hypothesis. Moreover, the $2$-cells of $\pcat P$ are in bijection with
  $2$-cells of $\freecat P_2$ in normal form because the polygraph is convergent
  (\cref{prop:1pol-conv-cf}), and those are in bijection with the $2$-cells
  of~$C$ by hypothesis.
\end{proof}

\subsection{A presentation of $\Simplaug$}
\index{simplicial!category}
\index{category!simplicial}
\label{sec:ass-simpl-pres}
As a detailed example of the above method, we show here that the
$3$-polygraph~$\MonPoly$ of \cref{ex:ass} presents the augmented simplicial
category~$\Simplaug$.
This category was introduced in \cref{sec:simpl-cat}: its objects are natural
numbers and morphisms are non-decreasing functions. It is moreover monoidal,
with tensor product given on objects by addition (such a monoidal category is
called a PRO, see \cref{sec:PRO}). As such, it can be considered as a
$2$-category (see \cref{sec:strict-moncat}), of which we now make an explicit
description.

The $2$-category $\Simplaug$ has one $0$-cell $\star$, the $1$-cells are natural
numbers and the $2$-cells $f:m\to n$ are non-decreasing maps from $\intset{m}$
to $\intset{n}$, where $\intset{n}=\set{0,\dots, n-1}$ for~$n\geq 0$. The
vertical composition of $2$-cells is the usual composition of functions, with
identities as neutral elements. The horizontal composition of $1$-cells is given
by addition, with $0$ as neutral element, and the horizontal composition of
$2$-cells $f:m\to n$ and $f':m'\to n'$ is given by
\begin{equation}
  \label{eq:simpl-tensor}
  (f\comp0 f')(i)
  =
  \begin{cases}
    f(i)&\text{if $0\leq i<m$,}\\
    f(i-m)+n&\text{if $m\leq i<m+m'$.}
  \end{cases}
\end{equation}

We have seen in \cref{ex:ass-cp} that the polygraph~$\MonPoly$ is convergent.
The normal forms in $\freecat{P_2}$ can be characterized as follows. Given
$n\in\N$, we define the \emph{right comb} $\mu_n:n\To 1$ in $\freecat{P_2}$ by
induction:
\begin{align*}
  \mu_0&=\eta=\satex{ass-mu0}
  &
  \mu_1&=a=\satex{ass-mu1}
  &
  \mu_{n+2}&=(a\comp0\mu_{n+1})\comp1\mu=\satex{ass-mun2}\pbox.
\end{align*}
A \emph{right forest} is a horizontal composite of right combs, \ie a morphism
of the form
\[
  \mu_{n_1}\comp0\mu_{n_2}\comp0\ldots\comp0\mu_{n_k}
  =
  \satex{ass-forest}
\]
for some $k\in\N$ called the \emph{width} of the right forest, and
$(n_1,\ldots,n_k)\in\N^k$.

\begin{lemma}
  \label{lem:mu-mu-mu}
  Given $n,n'\in\N$, $(\mu_n\comp0\mu_{n'})\comp1\mu$ rewrites to $\mu_{n+n'}$.
\end{lemma}
\begin{proof}
  By induction on~$n$. For $n=0$, we have
  $(\mu_0\comp0\mu_{n'})\comp1\mu=\mu_{n'+1}$. Otherwise,
  $(\mu_{n+1}\comp0\mu_{n'})\comp1\mu$ rewrites in one step to
  $(\mu_n\comp0\mu_{n'+1})\comp1\mu$ and we conclude using the induction
  hypothesis.
\end{proof}

\begin{lemma}
  \label{lem:ass-nf}
  The 2-cells of $\freecat{P_2}$ in normal form are precisely the right
  forests.
\end{lemma}
\begin{proof}
  The right forests are easily checked to be normal forms.
  We now show that, conversely, every normal form is a right forest, by showing
  that every $2$\nbd-cell~$\phi$ in $\freecat{P_2}$ rewrites to a right
  forest. The proof is done by induction on the size of~$\phi$. If the size of
  $\phi$ is $0$ then it is an identity, which is a right forest. Otherwise it
  can be decomposed as
  \[
    \phi
    =
    \psi\comp1(a^i\comp0\mu\comp0 a^j)
    =
    \satex{ass-dec}
  \]
  ($a^i$ denoting the horizontal composition of $i$ instances of~$a$) where, by
  induction, $\psi$ is a right forest
  \[
    \psi
    =
    \mu_{n_1}\comp0\mu_{n_2}\comp0\ldots\comp0\mu_{n_{i+1+j}}
    \pbox.
  \]
  By \lemr{mu-mu-mu}, $\phi$ rewrites to the right forest
  \[
    \mu_{n_1}\comp0\ldots\comp0\mu_{n_{i-1}}\comp0\mu_{n_i+n_{i+1}}\comp0\mu_{n_{i+2}}\comp0\ldots\comp0\mu_{n_{i+1+j}}
  \]
  which can be graphically depicted as
  \begin{align*}
    \phi
    &=
    \satex{ass-phi-ind-l}
    \\
    &=
    \satex{ass-phi-ind-r}
  \end{align*}
  and concludes the proof.
\end{proof}

The underlying $1$-category of~$\Simplaug$ is the additive monoid $\N$, seen as a
category. It is thus the free category on the $1$-polygraph $\pres{\star}{a}$,
which is the underlying $1$-polygraph of $\MonPoly$.
We define a function $f_2:\MonPoly_2\to{(\Simplaug)}_2$ where
$f_2(\mu):2\To 1$ and $f_2(\eta):0\To 1$
are both  terminal arrows. Consider the $3$-generator
\[
  A
  :
  \phi\TO\psi
  :
  a^3\To a
  :
  \star\to\star
\]
with $\phi=(\mu\comp0 a)\comp1\mu$ and $\psi=(a\comp0\mu)\comp1\mu$. We
necessarily have $f_2(\phi)=f_2(\psi)$ (where $f_2$ is extended by
functoriality, as in \cref{prop:3-pres-method}) because both $2$-cells are of
type $3\To 1$ in~$\Simplaug$, and the object $1$ is terminal in this category. A
similar reasoning can be held for the two other $3$-generators $L$, and $R$.

\begin{lemma}
  \label{lem:ass-nf-bij}
  The function $f_2$ induces a bijection between the normal forms in~$\freecat{P_2}$
  and the 2-cells of~$\Simplaug$.
\end{lemma}
\begin{proof}
  Given a $2$-cell $g:n\to k$ in~$\Simplaug$, \ie a non-decreasing function
  $g:\intset{n}\to\intset{k}$, and $j\in\intset{k}$, we write $n_j$ for the
  cardinal of the set $f^{-1}(j)$. Note that the function $g$ is entirely
  determined by the tuple $(n_1,\ldots,n_k)\in\N^k$ since, for $i\in\intset{n}$,
  $g(i)$ is the unique element of~$\intset{k}$ satisfying
  \[
    \sum_{0\leq j<g(i)}n_j
    \leq
    i
    <
    \sum_{0\leq j\leq g(i)}n_j
  \]
  and conversely any $k$-tuple of natural numbers $(n_1,\ldots,n_k)$ determines
  an increasing function from $\sum_jn_j$ to $k$ in this way. Given a right
  forest $\mu_{n_1}\comp0\ldots\comp0\mu_{n_k}$, one easily checks that its
  image under~$f$ is the non-decreasing function with associated $k$-tuple 
  $(n_1,\ldots,n_k)$, thus establishing a bijection between forests of width~$k$
  and non-decreasing functions with codomain~$k$.
\end{proof}

\noindent
In order to illustrate the above proof, consider the function $g:5\to 3$ whose
graph is depicted on the left below:
\begin{align*}
  \vcenter{
    \xymatrix@C=3ex@R=5ex{
      0\ar@{-}[d]&1\ar@{-}[dl]&2\ar@{-}[dll]&3\ar@{-}[dr]&4\ar@{-}[d]\\
      0&&1&&2
    }
  }
  &&\qquad
  \satex{ass_ex}
\end{align*}
The associated sequence in $\N^3$ is $(3,0,2)$ and the associated normal form is
the right forest $\mu_3\comp0\mu_0\comp0\mu_2$ as pictured on the right. Note
the clear correspondence between the two figures.

Let us sum up the results we have obtained in this section for the augmented
simplicial category $\Simplaug$. We have
\begin{enumerate}
\item constructed two $1$-cells $\mu:2\to 1$ and $\eta:0\to 1$ in the
  $2$-category $\Simplaug$: both are terminal $2$-cells,
\item shown that they generate the $2$-category: every morphism of $\Simplaug$
  can be written as a (horizontal or vertical) composite of those $2$-cells,
\item shown that those two $2$-cells satisfy the axioms $A$, $L$ and $R$ of
  \exr{ass} (reading $3$-cells are equalities), expressing that $\mu$ is
  associative and admits~$\eta$ as left and right unit,
\item shown that this set of axioms is complete: if two composites of $\mu$ and
  $\eta$ give rise to the same $2$-cell in $\Simplaug$ then one can show that they
  are equal using axioms $A$, $L$, $R$, and axioms for $2$-categories.
\end{enumerate}
Since the generators and the relations are precisely those of monoids, we can
deduce that $\Simplaug$ impersonates the notion of monoid, in the sense that a
$2$\nbd-func\-tor $\Simplaug\to C$ to some $2$-category~$C$ determines a monoid,
as explained in \exr{ass-alg}.

\subsection{Presentations from polygraphs with canonical forms}
\label{sec:3-pres-cf}
The observation made for $2$-polygraphs in \cref{rem:2-pres-cf} generalizes to
$3$-polygraphs. Namely, in the proof of \cref{prop:3-pres-method}, we did not
fully use the convergence of the polygraph, only the existence of
canonical representatives of equivalence classes provided by normal
forms. This suggests the following generalization of the above method, which
applies in cases where, even though it is difficult to construct a convergent
presentation, one can still directly come up with a notion of canonical
form. Many applications of this methodology were studied by
Lafont~\cite{lafont2003towards}.

\begin{proposition}
  \label{prop:3-pres-method-cf}
  Suppose given
  \begin{itemize}
  \item a 3-polygraph~$P$,
  \item a 2-category~$C$ whose underlying category is isomorphic to the free
    category on the underlying 1-polygraph of~$P$,
  \item a function $f_2:P_2\to C_2$ which is compatible with source and target,
  \item a set $\nf{\freecat{P_2}}\subseteq\freecat{P_2}$ of 2-cells called
    \emph{canonical forms},
  \end{itemize}
  such that
  \begin{enumerate}
  \item for any 2-generator $A:\phi\TO\psi$ in $P_3$, we have
    $f_2(\phi)=f_2(\psi)$,
  \item every 2-cell in $\freecat{P_2}$ is equivalent, \wrt the congruence
    $\approx^P$ generated by~$P_3$, to a canonical form,
  \item $f_2$ restricts to a bijection between canonical forms $u\To v$ in
    $\freecat{P_2}$ and 2-cells $f(u)\To f(v)$ in~$C$, where $f_2$ is implicitly extended to 2-cells in $\freecat P_2$ by
  functoriality.
  \end{enumerate}
   Then $P$ is a presentation of the 2-category~$C$.
\end{proposition}
\begin{proof}
  As in the proof of \cref{prop:3-pres-method}, by the first condition, $f_2$ induces
  a $2$-functor $f:\pcat P\to C$, which is bijective on $0$- and $1$-cells. By
  surjectivity of $f_2$ in the third condition, for every $2$\nbd-cell~$\phi$
  in~$C_2$, there is a canonical form $\nf\phi$ such that
  $f_2(\nf\phi)=\phi$. We have $f(\pcat{\nf\phi})=f_2(\nf\phi)=\phi$ and $f$ is
  thus surjective on $2$-cells. Moreover, consider two parallel $2$-cells $\phi$
  and $\psi$ of~$\pcat P_2$, such that $f(\phi)=f(\psi)$. The second condition
  ensures that there are canonical forms $\nf\phi$ and $\nf\psi$ such that
  $\phi=\pcat{\nf\phi}$ and $\psi=\pcat{\nf\psi}$. We have
  \[
    f_2(\nf\phi)=f(\pcat{\nf\phi})=f(\phi)=f(\psi)=f(\pcat{\nf\psi})=f_2(\nf\psi)
    \pbox.
  \]
  By injectivity of $f_2$ in the third condition, we have
  $\nf\phi=\nf\psi$. Thus,
  $
  \phi=\pcat{\nf\phi}=\pcat{\nf\psi}=\psi
  $
  and $f$ is injective on $2$-cells.
\end{proof}

\noindent
In particular, when~$P$ is a convergent $3$-polygraph, we can take
$\nf{\freecat{P_2}}$ to be the set of $2$-cells in normal form, thus recovering
\cref{prop:3-pres-method} as a particular case of previous proposition.

\begin{example}
  \label{ex:cf-not-nf}
  In order to show that $\MonPoly$ presents $\Simplaug$, we could have  chosen the
  following alternative definition of right combs:
  \begin{align*}
    \mu_0&=\eta
    &
    \mu_{n+1}&=(a\comp0\mu_n)\comp1\mu    
  \end{align*}
  Right forests are obtained as horizontal composites of such right combs, and
  we consider those as canonical forms.
  The definition is mostly the same as before except that we have
  \[
    \mu_1
    =
    (a\comp0\eta)\comp1\mu
    =
    \satex{ass-mu1-alt}
    \pbox.
  \]
  It can be shown that every $2$-cell is equivalent to a canonical form using a
  variant of the proof of \lemr{ass-nf}, and one can construct a bijection
  between canonical forms and $2$-cells in $\Simplaug$ using a variant of the
  proof of \cref{lem:ass-nf-bij}. We can thus conclude that $\MonPoly$ is a
  presentation of $\Simplaug$ by \cref{prop:3-pres-method-cf}.
  Note that the canonical form associated to $\id_a$ is~$\mu_1$, so there is no
  hope to obtain canonical forms as normal forms for some convergent
  $3$-polygraph, because no terminating polygraph can rewrite identities.
\end{example}

\begin{remark}
  As a variant of the previous example, consider the
  category~$\Simplaug^2=\Simplaug\times\Simplaug$. The monoidal structure on~$\Simplaug$
  induces one on $\Simplaug^2$ given on objects by
  \[
    (m_1,n_1)\otimes(m_2,n_2)
    =
    (m_1\otimes m_2,n_1\otimes n_2)
  \]
  and similarly on morphisms. The underlying monoid of objects of this category
  is $\N\times\N$ which is abelian and thus not free (see \secr{Z-not-free}).
  Therefore, by \secr{2-cat-no-pres}, there
  is no presentation of it (seen as a $2$-category induced by the monoidal
  structure, see \secr{PRO}) by a $2$-polygraph. This situation is detailed
  in~\cite{curien2017coherent}. The same argument applies to most products of
  $2$-categories, but for degenerated cases.
\end{remark}

\subsection{Presenting categories}
\index{presentation!of a $2$-category}
\label{sec:3-pres-2-pres}
Let $P$ be a $3$-polygraph presenting a monoidal category~$C$, seen as a
$2$-category: this presentation is of the form
\[
  P
  =
  \PRes\star{P_1}{P_2}{P_3}
  \pbox.
\]
The monoidal category~$C$ has an underlying category obtained by forgetting the
tensor product and unit object. This category admits, as a category, a
presentation by the following $2$-polygraph~$Q$ constructed from~$P$:
\begin{itemize}
\item $Q_0=\freecat P_1$ is the set of $1$-cells $u:\star\to\star$
  in~$\freecat{P_1}$,
\item $Q_1=\freecat{P_1}P_2\freecat{P_2}$ contains a $1$-generator
  \[
    u\alpha w
    :
    uvw
    \to
    uv'w
  \]
  for every $1$-cells $u,w\in\freecat{P_1}$ and $2$-generator $\alpha:v\To v'$
  in $P_2$,
\item $Q_2$ contains a $2$\nbd-generator
  \[
    uAw
    :
    u\phi w
    \To
    u\psi w
  \]
  for every $1$-cells $u,w\in\freecat{P_1}$ and $3$-generator $A:\phi\TO\psi$,
  where $u\phi w$ and $u\psi w$ are seen as elements of $\freecat{Q_2}$ in the
  expected functorial way, it also contains a $2$-generator
  \[
    X_{u,\alpha,u',\beta,u''}
    :
    u\alpha u'wu''\comp{}uv'u'\beta u''
    \To
    uvu'\beta u''\comp{}u\alpha u'w'u''
  \]
  for every $1$-cells $u,u',u''\in\freecat{P_1}$ and $2$-generators
  $\alpha:v\To v'$ and $\beta:w\To w'$ in $P_2$ (which encodes the exchange
  law).
\end{itemize}

\begin{example}
  \index{simplicial!category}
  \index{category!simplicial}
  \label{ex:simpl-3-2-pres}
  We have seen above that the augmented simplicial category~$\Simplaug$ was
  presented, as a monoidal category, by the $3$-polygraph~$\MonPoly$, defined in
  \exr{ass}. We deduce the presentation of~$\Simplaug$, as a category, by the
  $2$-polygraph~$P$ with
  \begin{itemize}
  \item $0$-generators: for $i\in\N$, a generator $a^i$,
  \item $1$-generators: for $i,j\in\N$,
    \begin{align*}
      a^i\mu a^j&:a^{i+2+j}\to a^{i+1+j}
      &
      a^i\eta a^j&:a^{i+j}\to a^{i+1+j}      
    \end{align*}
  \item $2$-generators: for $i,j\in\N$,
    \[
      \begin{array}{r@{\ :\ }r@{\ \To\ }l}
      a^iAa^j&a^i\mu a^{j+1}\comp{}a^i\mu a^j&a^{i+1}\mu a^j\comp{}a^i\mu a^j\\
      a^iLa^j&a^i\eta a^{j+1}\comp{}a^i\mu a^j&a^{i+j}\\
      a^iRa^j&a^{i+1}\eta a^j\comp{}a^i\mu a^j&a^{i+j}\\
      X_{a^i,\mu,a^j,\mu,a^k}&a^i\mu a^{j+2+k}\comp{}a^{i+1+j}\mu a^k&a^{i+2+j}\mu a^k\comp{}a^i\mu a^{j+1+k}\\
      X_{a^i,\mu,a^j,\eta,a^k}&a^i\mu a^{j+k}\comp{}a^{i+1+j}\eta a^k&a^{i+2+j}\eta a^k\comp{}a^i\mu a^{j+1+k}\\
      X_{a^i,\eta,a^j,\mu,a^k}&a^i\eta a^{j+2+k}\comp{}a^{i+1+j}\mu a^k&a^{i+j}\mu a^k\comp{}a^i\eta a^{j+1+k}\\
      X_{a^i,\eta,a^j,\eta,a^k}&a^i\eta a^{j+k}\comp{}a^{i+1+j}\eta a^k&a^{i+j}\eta a^k\comp{}a^i\eta a^{j+1+k}        
      \end{array}
    \]
  \end{itemize}
  It can be checked that we precisely recover the presentation for the
  simplicial category given in~\secr{simpl-cat}, up to renaming the
  $1$-generators $a^i\mu a^j$ to $s^{i+j+1}_i$ and $a^i\eta a^j$ to~$d^{i+j}_i$.
\end{example}

\section{Indexed Critical Branchings}
\label{sec:3pol-indexed}
There is a major difference between rewriting in $3$-polygraphs compared to the
case of $2$-polygraphs studied in previous chapters: contrarily to presentations
of categories (see \lemr{2-cb-finite}), the number of critical branchings of a
finite $3$-polygraph can be infinite. We begin with an example of this
phenomenon, originally observed by Lafont~\cite{lafont2003towards}.

\subsection{Presenting the theory for symmetries}
\index{category!of permutations}
\index{polygraph!of permutations}
\index{permutation}
\nomenclature[S]{$\Bij$}{category of permutations}
\label{sec:sym-pres}
The category $\Bij$ is the category whose objects are natural numbers and a
morphism $f:m\to n$ is a bijection (also called a \emph{permutation}) from
$\intset{m}$ to $\intset{n}$, the ordinals with $m$ and $n$ elements,
respectively, with usual compositions and identities. Here, all morphisms are
in fact endomorphisms. This category is monoidal with tensor product given by
addition on objects (this is a PRO) and as for $\Simplaug$ on morphisms,
see~\eqref{eq:simpl-tensor} in \cref{sec:ass-simpl-pres}.

Starting from the fact that any bijection can be decomposed as a composite of
transpositions, we expect that this monoidal category, seen as a $2$-category,
admits a presentation by the following $3$-polygraph~$P$:
\[
  \PRes{\star}{a}{\gamma:a\comp0a\to a\comp0a}{I,Y}\pbox.
\]
Here, the $2$-generator $\gamma$ corresponds to the transposition on a set with
two elements and is usually pictured as
\[
  \satex{gamma}
\]
The two $3$-generators express
\begin{itemize}
\item the \emph{involutivity} of the transposition:
  \begin{align*}
    I
    :
    \gamma\comp1\gamma
    &\TO
    \id_a\comp0\id_a
    \\
  \intertext{which can be represented as}
  \satex{sym-l}
  &\TO
  \satex{sym-r}
  \end{align*}
\item the \emph{Yang-Baxter} relation~$Y$ of type
  \begin{align*}
    (\gamma\comp0\id_a)\comp1(\id_a\comp0\gamma)\comp1(\gamma\comp0\id_a)
    &\TO
    (\id_a\comp0\gamma)\comp1(\gamma\comp0\id_a)\comp1(\id_a\comp0\gamma)
    \\
    \intertext{which can be represented as}
    \satex{yb-l}
    &\TO
    \satex{yb-r}    
  \end{align*}
\end{itemize}
In this polygraph, it can be noted that,
for any $2$-cell $\phi:a^{m+1}\To a^{n+1}$,
the $2$-cell
\[
  (\gamma\comp0a^{m+1})\comp1(a\comp0\gamma\comp0a^m)\comp1(\gamma\comp0\phi)\comp1(a\comp0\gamma\comp0a^n)\comp1(\gamma\comp0a^{n+1})
\]
of $\freecat{P_2}$ can be rewritten in two ways using the
rule~$Y$:
\[
  \satex{sym-icp-l}
  \quad\OT\quad
  \satex{sym-icp}
  \quad\TO\quad
  \satex{sym-icp-r}
  \pbox.
\]
This gives rise to a critical branching when $\phi$ is either $\id_a$ or of the
form $\gamma^{n+1}$ (the vertical composite of $n+1$ instances of $\gamma$). The
critical branchings of the rewriting system are thus
\begin{align*}
  \satex{sym-cp1}
  &&
  \satex{sym-cp2}
  &&
  \satex{sym-cp3}
  &&
  \satex{sym-cp4}
  &&
  \satex{sym-cp5}
  \,\pbox.
\end{align*}
In the first case, the $2$-cell can be rewritten by $I$ in two different ways, in
the second and third case, the $2$-cells can be rewritten both by $I$ and $Y$. The
fourth and fifth cases can be rewritten by~$Y$ in two ways, as described above.
A finite $3$-polygraph can thus give
rise to an infinite number of critical branchings. Still, they can be checked to
be confluent. In the first four cases, this can be checked directly:
\begin{align*}
    \xymatrix@R=5ex{
      \ar@3@/_/[d]\satex{sym-cp1u}\ar@3@/^/[d]\\
      \satex{sym-cp1d}
    }
    &&
    \xymatrix@R=2ex{
      \satex{sym-cp2r}\ar@3[d]&\satex{sym-cp2u}\ar@3[l]\ar@3[d]\\
      \satex{sym-cp2rr}\ar@3[r]&\satex{sym-cp2d}
    }
    &&
    \xymatrix@R=2ex{
      \ar@3[d]\satex{sym-cp3u}\ar@3[r]&\satex{sym-cp3r}\ar@3[d]\\
      \satex{sym-cp3d}&\satex{sym-cp3rr}\ar@3[l]
    }
\end{align*}
\[
  \xymatrix@R=2ex{
    \satex{sym-cp4l}\ar@3[d]&\ar@3[l]\satex{sym-cp4u}\ar@3[r]&\ar@3[d]\satex{sym-cp4r}\\
    \satex{sym-cp4ll}\ar@3[r]&\satex{sym-cp4d}&\ar@3[l]\satex{sym-cp4rr}
  }
\]
In the fifth case, we have, depending on whether $n$ is odd or even:
\[
  \gamma^{n+1}\overset*\TO\unit{a^2}
  \qquad\qquad\text{or}\qquad\qquad
  \gamma^{n+1}\overset*\TO\gamma
\]
and in the two cases the local branching is confluent:
\begin{itemize}
\item if $\gamma^{n+1}\overset*\TO\unit{a^2}$:
  \[
    \xymatrix@R=2ex{
      \satex{sym-cp5l}\ar@3[d]_\ast&&\ar@3[ll]\satex{sym-cp5u}\ar@3[rr]&&\ar@3[d]^\ast\satex{sym-cp5r}\\
      \satex{sym-cp5ll}\ar@3[r]&\satex{sym-cp5lll}\ar@3[r]&\satex{sym-cp5d}&\ar@3[l]\satex{sym-cp5rrr}&\ar@3[l]\satex{sym-cp5rr}
    }
  \]
\item if $\gamma^{n+1}\overset*\TO\gamma$:
  \[
    \xymatrix@R=2ex{
      \satex{sym-cp5l2}\ar@3[d]&\satex{sym-cp5l}\ar@3[l]_\ast&\ar@3[l]\satex{sym-cp5u}\ar@3[r]&\ar@3[r]^\ast\satex{sym-cp5r}&\ar@3[d]\satex{sym-cp5r2}\\
      \satex{sym-cp5l3}\ar@3[r]&\satex{sym-cp5l4}\ar@3[r]&\satex{sym-cp5dd}&\ar@3[l]\satex{sym-cp5r4}&\ar@3[l]\satex{sym-cp5r3}
    }
  \]
\end{itemize}
Finally, the rewriting system can be shown to be terminating
(\cref{ex:sym-term,ex:sym-term2}) and thus convergent.

The normal forms can be characterized as follows. We first define, for every
$n\in\N$, a $2$-cell
\[
  \gamma_n
  :
  a^{n+1}
  \to
  a^{n+1}
\]
by
\begin{align*}
  \gamma_0&=\unit{a}
  &
  \gamma_{n+1}&=(\gamma_n\comp0 a)\comp1(a^n\comp0\gamma)
  \pbox.
\end{align*}
This cell can be seen as a generalization of $\gamma$ since $\gamma_1=\gamma$,
and $\gamma_n$ is most naturally pictured as
\[
  \satex{gamman}
\]
\ie for low values of $n$,
\begin{align*}
  \gamma_0&=\satex{gamma0}
  &
  \gamma_1&=\satex{gamma1}
  &
  \gamma_2&=\satex{gamma2}
  &
  \gamma_3&=\satex{gamma3}
  \,\pbox.
\end{align*}
The normal forms of the rewriting system can be characterized as follows, from
which we can conclude that this is indeed a presentation of~$\Bij$.

\begin{proposition}
  A 2-cell~$\phi$ is a normal form if and only if it is either
  \begin{itemize}
  \item $\phi=\unit{a^0}$, or
  \item there exists a normal form $\psi:n\to n$ and $m\in\N$ with $0\leq m\leq n$ such
    that~$\phi$ is
    \[
      (a\comp0\psi)\comp1(\gamma_m\comp0 a^{n-m})
      :
      m+1
      \to
      m+1
    \]
    what we write
    \[
      \phi
      =
      \Gamma_m\psi
    \]
    \ie graphically,
    \[
      \phi=\satex{sym-nf}
      \pbox.
    \]
  \end{itemize}
\end{proposition}
\begin{proof}
  We call \emph{canonical forms} the $2$-cells of the above form. Those are in
  normal form: this can be checked directly. Conversely, we show that any
  $2$\nbd-cell~$\phi$ rewrites to a canonical form by induction on its size. If
  the size of $\phi$ is $0$ then $\phi$ is of the form $\unit{a^n}$ and we have
  $\phi=\Gamma_0\ldots\Gamma_0\Gamma_0\unit{a^0}$ (with $n$ occurrences of
  $\Gamma_0$). Otherwise, $\phi$ is of the form
  \[
    \phi
    =
    \phi'\comp1(a^i\comp0\gamma\comp0a^j)
    =
    \satex{sym-dec}
  \]
  where $\phi'$ rewrites to a canonical form $\Gamma_m\psi$. Then depending on
  the respective values of $m$ and $i$, four generic situations are possible,
  and in all of them $\phi$ rewrites to a canonical form (by using the induction
  hypothesis in the first and in the last case):
  \begin{align*}
    \satex{sym-ind1-l}
    &\TO
    \satex{sym-ind1-r}
    \\
    \satex{sym-ind2-l}
    &\TO
    \satex{sym-ind2-r}
    \\
    \satex{sym-ind3-l}
    &\TO
    \satex{sym-ind3-r}
    \\
    \satex{sym-ind4-l}
    &=
    \satex{sym-ind4-r}
  \end{align*}
  This concludes the proof.
\end{proof}

\index{Lehmer code}
Given a bijection $f:[n]\to[n]$ its \emph{Lehmer code} is a sequence of~$n$
natural numbers $(k_0,k_1,\ldots,k_{n-1})$ such that $0\leq k_i<n-i$ for
every index~$i$, where $k_i$ is the cardinal of the set
$\setof{j>i}{f(j)<f(i)}$. It can be shown that this induces a bijection between
permutations of~$[n]$ and such
sequences, see~\cite{laisant1888numeration,lehmer1960teaching}. For instance,
consider the permutation of~$[4]$ whose images are $(2~0~4~3~1)$, pictured on
the left:
\begin{align}
  \label{eq:lehmer-ex}
  \vcenter{
    \xymatrix@C=1ex{
      0\ar@{-}[drr]&1\ar@{-}[dl]&2\ar@{-}[drr]&3\ar@{-}[d]&4\ar@{-}[dlll]\\
      0&1&2&3&4
    }
  }
  &&
  \satex{lehmer-ex}
\end{align}
The associated Lehmer code is $(2,0,2,1,0)$. We can finally conclude that the
polygraph~$P$ defined in~\ref{sec:sym-pres} is a presentation of the PRO of
symmetries, see~\cite{lafont2003towards} for details.

\begin{theorem}
  The polygraph~$P$ is a presentation of~$\Bij$. 
\end{theorem}
\begin{proof}
  Following the method described in \secr{3-pres-method}, we interpret the
  morphism $\gamma:2\to 2$ as the transposition $\intset2\to\intset2$, and this
  interpretation is compatible with the two rules.
  It is easy to see that the interpretation of a normal form
  $\Gamma_{k_0}\ldots\Gamma_{k_{n-1}}\unit{a^0}$ is the bijection whose Lehmer
  code is $(k_0,k_1,\ldots,k_n)$. For instance, the normal form associated to the
  bijection \cref{eq:lehmer-ex} is
  $\Gamma_2\Gamma_0\Gamma_2\Gamma_1\Gamma_0\id_{a^0}$, which is pictured on the
  right of \cref{eq:lehmer-ex}. This clearly establishes a bijection between
  $2$-cells in normal form in $\freecat{\tpol2 P}$ and $2$-cells of $\Bij$.
\end{proof}

\subsection{Classification of critical branchings}
\label{sec:3cb-classification}
Critical branchings in $3$-poly\-graphs are classified
in~\cite{guiraud2009higher}. This case is more difficult than in dimension~$2$ mainly because,
as initially noted in~\cite{lafont2003towards} and observed in previous section, a finite
$3$-polygraph may have an infinite number of critical branchings.
However, an analysis of the possible shapes of these critical
branchings yields a sufficient condition for confluence that only requires to consider a
finite subset of them.

Assume that~$P$ is a $3$-polygraph. By examination of the different
possibilities, the critical branchings of~$P$ are classified as
follows~\cite[Section~5.1.1]{guiraud2009higher}.
\begin{enumerate}
\item \emph{Inclusion} critical branchings,
  with the following source, if~$\chi$ is the source of a $3$-generator of~$P$, and
  $\phi \comp1 u\chi v \comp1\psi$ is the source of another one:
  \[
    \satex{cp3-inclusion}
  \]
\item \emph{Regular} critical branchings, with the following source,
  if~$\phi\comp1 u\chi$ and~$\chi v\comp1\psi$ (or $\phi\comp1\chi v$
  and~$u\chi\comp1\psi$) are the sources of two $3$-generators of~$P$:
  \[
    \satex{cp3-regular1}
    \qquad\text{or}\qquad
    \satex{cp3-regular2}
\]
\item Instances of \emph{left-indexed} critical branchings, with the following
  source, if~$\phi\comp1 u\chi$ and~$v\chi\comp1\psi$ are the sources of two
  $3$-generators of~$P$, and $\zeta:wu\to xv$ is a $2$-cell of~$P^*$:
  \[
    \satex{cp3-indexed-l}
  \]
\item Instances of \emph{right-indexed} critical branchings, with the following
  source, if~$\phi\comp1\chi u$ and~$\chi v\comp1\psi$ are the sources of
  $3$-generators of~$P$, and $\zeta:uw\dfl vx$ is a $2$-cell of~$P^*$:
  \[
    \satex{cp3-indexed-r}
  \]
\item Instances of \emph{multi-indexed} critical branchings, in all the other
  cases: one has a $3$-generator with a source of the form
  \[
    \phi\comp1(u_0 \comp0 \chi_1\comp0 u_1\comp0 \chi_2\comp0 \cdots \comp0 u_{n-1} \comp0 \chi_n \comp0 u_n) 
  \]
  and another $3$-generator with a source of the form
  \[
    (v_0\comp0 \chi_1\comp0 v_1\comp0 \chi_2\comp0 \cdots \comp0 v_{n-1} \comp0 \chi_n \comp0 v_n)\comp1\psi
  \]
  so that the source of the branching is of the form
  \[
    \satex{cp3-mi}
  \]
\end{enumerate}
For example, in the presentation of \cref{sec:sym-pres}, the four first
branchings are regular whereas the generic family of branchings is
right-indexed (by $\gamma^{n+1}$).

An \emph{instance} of a left- or right-indexed branching is a left- or
right-indexed branching as above with a particular value for the $2$-cell
$\zeta$. It is a \emph{normal instance} when $\zeta$ is in normal form.

\subsection{Indexed polygraphs}
\index{indexed polygraph}
\index{polygraph!indexed}
\label{sec:3pol-non-indexed}
We say that a $3$-polygraph~$P$ is \emph{non-indexed} if it has inclusion or
regular critical branchings only, \emph{left-indexed} (\resp
\emph{right-indexed}) if it has inclusion, regular or left-indexed (\resp
right-indexed) critical branchings only, and \emph{finitely indexed} if each of
its indexed critical branchings has a finite number of reduced instances.
Then we have the following results, which apply to the presentation of
permutations in \cref{sec:sym-pres}.

\begin{proposition}[{\cite[Proposition~5.1.3]{guiraud2009higher}}]
  If~$P$ has a finite number of 3-cells, then it has a finite number of
  inclusion and regular critical branchings.
\end{proposition}

\begin{proposition}[{\cite[Proposition~5.3.1]{guiraud2009higher}}]
  If~$P$ is terminating and left-indexed (\resp right-indexed), then~$P$ is
  confluent if and only if all its inclusion and regular critical branchings,
  and all the reduced instances of its left-indexed (\resp right-indexed)
  critical branchings are confluent.
\end{proposition}

\section{Distributive Laws}
\label{sec:3-dlaws}
The notions introduced in \secr{2-dlaws} for combining presentations using
distributive laws generalize easily to presentations of monoidal categories as
we now briefly explain, following~\cite{lack2004composing}.

\subsection{Monoidal categories as monads}
We have seen in \cref{sec:cat-as-monad} that a category corresponds precisely to
a monad in $\Span(\Set)$, and our aim is to generalize this situation to strict
monoidal categories. The main starting point is that, in a strict monoidal
category, the set of objects forms a monoid, with tensor as product.

The category $\Mon$ of monoids and their morphisms has small limits and it
therefore makes sense to consider the bicategory $\Span(\Mon)$ of spans
internals to~$\Mon$ as explained in \cref{sec:span}: a $0$-cell of this
bicategory is a monoid, a $1$-cell from~$A$ to~$B$ is a diagram of the form
\[
  \xymatrix{
    A&\ar[l]C\ar[r]&B
  }
\]
in~$\Mon$ and composition is given by pullback. It can then be observed that a
monad in this bicategory (see \cref{sec:bicat-dlaw}) is precisely a monoidal
category. In particular, a PRO corresponds to a monad on the
monoid~$(\N,+,0)$. This point of view allows us to compose monoidal categories
through distributive laws between the corresponding monads. We briefly present
it below.

\subsection{Distributive laws}
\index{distributive law!of monoidal categories}
Given two monoidal categories~$C$ and~$D$ with the same monoid of objects, a
\emph{distributive law} between them is a distributive law between the corresponding
monads in $\Span(\Mon)$, in the sense of \cref{sec:bicat-dlaw}.
It consists of a distributive law
\[
  \dlaw
  :
  D\otimes C
  \to
  C\otimes D
\]
in the sense of \cref{sec:cat-dlaw}, between the underlying categories, which is
compatible with tensor product in the sense that, for $f_1,f_2,f_1',f_2'$ and
$g_1,g_2,g_1',g_2'$ morphisms of $C$ and $D$ respectively,
\[
  \begin{array}{l}
    \dlaw(g_1',f_1')=(f_1,g_1)\\
    \dlaw(g_2',f_2')=(f_2,g_2)
  \end{array}
  \qquad\text{implies}\qquad
  \dlaw(g_1'\otimes g_2',f_1'\otimes f_2')=(f_1\otimes f_2,g_1\otimes g_2)
  \pbox.
\]

\subsection{Factorization systems}
\index{factorization system!monoidal}
Given monoidal categories $C$, $D$, and $E$, with the same monoid of objects, we
have $C\otimes_\dlaw D=E$ precisely when~$C$ and~$D$ form a \emph{monoidal
  factorization system} for~$E$, \ie $C$ and~$D$ are monoidal subcategories of~$E$
such that every morphism of~$E$ factorizes uniquely a morphism of~$C$ followed
by a morphism of~$D$.

In the case where~$C$ (\resp $D$) admits a presentation by a 3-polygraph~$P$
(resp.~$Q$), with $P_0=Q_0=\set{\star}$ and $P_1=Q_1$, the category~$E$ admits a
presentation by the 3-polygraph~$R$ with
\begin{align*}
  R_0&=\set{\star}
  &
  R_1&=P_1=Q_1
  &
  R_2&=P_2\sqcup Q_2
  &
  R_3&=P_3\sqcup Q_3\sqcup R_3^\dlaw  
\end{align*}
where $R_3^\dlaw$ presents the distributive law similarly to
\cref{thm:fc-dlaw}.

\subsection{Composing PROPs}
Generalized composition of categories, as described in \secr{gen-dlaw}, also
extends to this setting. A typical situation where this is useful is the case of
PROPs: those are strict symmetric monoidal categories, with $(\N,+,0)$ as monoid
of objects. When composing two PROPs, one would like to identify the symmetric
structures already present in both of them.

The category~$\Bij$ of finite cardinals and bijections is the free PROP
(this follows from the presentation constructed in \cref{sec:sym-pres}, see also \cref{sec:free-sym}),
and a PROP~$C$ can thus be seen as a monad on the additive
monoid~$\N$ in $\Span(\Mon)$ equipped with a functor $\Bij\to C$, \ie an object
on the category on the left of \eqref{eq:slice-mon} below.
As explained in \cref{sec:gen-comp}, we have an isomorphism
\begin{equation}
  \label{eq:slice-mon}
  \Bij/\Mon(\Span(\Mon)(\N,\N))
  \isoto
  \Mon(\Mod(\Span(\Mon))(\Bij,\Bij))
\end{equation}
which thus allows one to consider a PROP as a monad on bimodules of spans of
monoids over~$\Bij$, and two PROPs~$C$ and~$D$ can be composed along a
distributive law
\[
  \dlaw
  :
  D\otimes_\Bij C
  \to
  C\otimes_\Bij D
\]
between the corresponding monads: the composite defined in this way always gives
rise to a PROP and identifies the symmetry structure in the composed PROPs as
expected~\cite{lack2004composing,zanasi2015interacting}.
When the two PROPs~$C$ and~$D$ are presented, one can obtain a presentation of
their composite using a direct generalization of \cref{sec:gen-pres}.

Similarly, Lawvere theories can be seen as monads in bimodules of spans of
monoids over $\Fun$ and can be composed along distributive
laws~\cite{cheng2011distributive}, see also~\secr{trs-dlaw}.

\begin{example}
  We write $\Fun$ (\resp $\Funsurj$, \resp $\Funinj$) for the PROP of finite
  cardinals and functions (\resp surjective functions, \resp injective
  functions), see \secr{2pres-mon}. The PROPs~$\Funsurj$ and~$\Funinj$ are
  subcategories of~$\Fun$, any morphism $h\in\Fun$ factorizes as $h=g\circ f$
  with $f\in\Funsurj$ and $g\in\Funinj$, and for any other factorization
  $h=g'\circ f'$ there exists a permutation~$w\in\Bij$ making the following
  diagram commute:
  \[
    \xymatrix@C=3ex@R=3ex{
      &y\ar@{.>}[dd]^w\ar[dr]^g&\\
      x\ar[ur]^f\ar[dr]_{f'}&&z\\
      &y'\ar[ur]_{g'}&
    }
  \]
  \ie we have $\Fun=\Funsurj\otimes_\dlaw\Funinj$. The categories~$\Funsurj$ and
  $\Funinj$ respectively admit presentations by the polygraphs~$P$ and~$Q$ with
  generators
  \begin{align*}
    P_0=Q_0&=\set{\star}
    &
    P_1=Q_1&=\set{a}
    &
    P_2&=\set{\satex{gamma},\satex{mu}}
    &
    Q_2&=\set{\satex{gamma},\satex{eta}}    
  \end{align*}
  where the relations in~$P_3$ are
  \begin{align*}
    \satex{yb-l}
    &\TO
    \satex{yb-r}
    &
    \satex{mu-gamma-l}
    &\TO
    \satex{mu-gamma-r}
    &
    \satex{gamma-mu-l}
    &\TO
    \satex{gamma-mu-r}
    \\
    \satex{mon-assoc-l}
    &\TO
    \satex{mon-assoc-r}
    &
    \satex{mon-com-l}
    &\TO
    \satex{mon-com-r}    
  \end{align*}
  and the relations in~$Q_3$ are
  \begin{align*}
    \satex{yb-l}
    &\TO
    \satex{yb-r}
    &
    \satex{eta-gamma-l}
    &\TO
    \satex{eta-gamma-r}
    &
    \satex{gamma-eta-l}
    &\TO
    \satex{gamma-eta-r}    
  \end{align*}
  From those, we deduce the presentation of~$\Fun$ by the polygraph~$R$ with
  \begin{align*}
    R_0&=\set{\star}
    &
    R_1&=\set{a}
    &
    R_2&=P_2\sqcup Q_2=\set{\satex{gamma},\satex{mu},\satex{eta}}    
  \end{align*}
  and the relations are
  \begin{itemize}
  \item the common relation for symmetry:
    \[
      \satex{yb-l}
      \TO
      \satex{yb-r}
    \]
  \item the relations of~$P$:
    \begin{align*}
      \satex{mu-gamma-l}
      &\TO
      \satex{mu-gamma-r}
      &
      \satex{gamma-mu-l}
      &\TO
      \satex{gamma-mu-r}
      \\
      \satex{mon-com-l}
      &\TO
      \satex{mon-com-r}
      &
      \satex{mon-assoc-l}
      &\TO
      \satex{mon-assoc-r}
    \end{align*}
  \item the relations of~$Q$:
    \begin{align*}
      \satex{eta-gamma-l}
      &\TO
      \satex{eta-gamma-r}
      &
      \satex{gamma-eta-l}
      &\TO
      \satex{gamma-eta-r}
    \end{align*}
  \item the relations generated by the distributive law:
    \begin{align*}
      \satex{mon-unit-l}
      &\TO
      \satex{mon-unit-c}
      &
      \satex{mon-unit-r}
      &\TO
      \satex{mon-unit-c}
    \end{align*}
  \end{itemize}
\end{example}


\chapter{Termination of 3-Polygraphs}
\label{chap:3term} 
This chapter presents techniques for proving the termination of $3$-polygraphs,
generalizing those already introduced for $2$-polygraphs in
\cref{sec:2-red-order}. A first method, described
in~\cref{sec:3red-order}, is based on a certain type of well-founded
orders called \emph{reduction orders}. We then turn
in~\cref{sec:3-term-func-interp} to \emph{functorial interpretations}:
these amount to construct a functor from the underlying category to
another category which already bears a reduction order.
This covers quite
a few useful examples. To address more complex cases, we present in
\cref{sec:3-term-der} a powerful technique, due to
Guiraud~\cite{guiraud2006termination,guiraud2009higher}, based on the
construction of a derivation from the polygraph. Here, termination is
obtained  by specifying quantities on $2$-cells
which decrease during rewriting, based on information propagated by the
$2$-cells themselves.

\section{Reduction and Termination Orders}
\label{sec:3red-order}
We begin by extending the notion of reduction order introduced in
\secr{2-red-order}, from 2-polygraphs to 3-polygraphs.

\begin{definition}
  \index{reduction!order}
  \index{termination!order}
  \index{order!reduction}
  \index{order!termination}
  \label{def:3-red-order}
  Given a 2-category~$C$, a \emph{reduction order} is a partial order $\succeq$
  on pairs of parallel 2-cells which is
  \begin{itemize}
  \item well-founded: every weakly decreasing sequence of 2-cells is
    eventually stationary,
  \item compatible with 0-composition: for every 2-cells
    \begin{align*}
      \phi&:u\To u':x'\to x
      &
      \psi_1,\psi_2&:v\To v':x\to y
      &
      \phi'&:w\To w':y\to y'
    \end{align*}
    which can be represented as
    \[
      \xymatrix@C=8ex{
        x'\ar@/^3ex/[r]^u\ar@/_3ex/[r]_{u'}\ar@{}[r]|{\phantom\phi\Longdownarrow\phi}&x\ar@/^3ex/[r]^v\ar@/_3ex/[r]_{v'}\ar@{}[r]|{\psi_1\Longdownarrow\ \Longdownarrow\psi_2}&y\ar@/^3ex/[r]^w\ar@/_3ex/[r]_{w'}\ar@{}[r]|{\phantom{\phi'}\Longdownarrow{\phi'}}&y'
      }
    \]
    we have that
    \[
      \psi_1\succ\psi_2
      \qquad\text{implies}\qquad
      \phi\comp0\psi_1\comp0\phi'\succ\phi\comp0\psi_2\comp0\phi'
      \pbox,
    \]
  \item compatible with 1-composition: for every 2-cells
    \begin{align*}
      \phi&:u'\To u:x\to y
      &
      \psi_1,\psi_2&:u\To v:x\to y
      &
      \phi'&:v\To v':x\to y
    \end{align*}
    which can be represented as
    \[
      \xymatrix@C=14ex{
        x
        \ar@/^7ex/[r]^{u'}\ar@{{}{ }{}}@/^5ex/[r]|{\phantom\phi\Longdownarrow\phi}
        \ar@/^2.5ex/[r]|u\ar@/_2.5ex/[r]|v\ar@{}[r]|{\psi_1\Longdownarrow\quad\Longdownarrow\psi_2}
        \ar@/_7ex/[r]_{v'}\ar@{{}{ }{}}@/_5ex/[r]|{\phantom{\phi'}\Longdownarrow{\phi'}}
        &y
      }
    \]
    we have that
    \[
      \psi_1\succ\psi_2
      \qquad\text{implies}\qquad
      \phi\comp1\psi_1\comp1\phi'\succ\phi\comp1\psi_2\comp1\phi'
      \pbox.
    \]
  \end{itemize}
  Given a 3-polygraph~$P$, a reduction order $\succeq$ on the 2-category
  $\freecat{P_2}$ is said to be \emph{compatible} with the rules of~$P$ when
  $\phi\succ\phi'$ for every rule $A:\phi\TO\phi'$ in~$P_3$. In this case, the
  order $\succeq$ is called a \emph{termination order} for~$P$.
\end{definition}

In an arbitrary 3-polygraph~$P$, we write $\freecat\TO$ for the relation on
parallel 2\nbd-cells of $\freecat{P_2}$ such that $\phi\freecat\TO\psi$ whenever
$\phi$ rewrites to $\psi$, or equivalently whenever there is a 3-cell
$F:\phi\TO\psi$ in $\freecat{P_3}$.

\begin{proposition}
  \label{prop:3-term-red-ord}
  \label{P:TerminationOrder}
  Given a 3-polygraph~$P$, the following statements are equivalent.
  \begin{enumerate}
  \item The polygraph~$P$ is terminating.
  \item The relation~$\freecat\TO$ is a termination order.
  \item The polygraph~$P$ admits a termination order.
  \end{enumerate}
\end{proposition}
\begin{proof}
  Similar to the proof of \propr{term-red-ord}.
\end{proof}

\section{Functorial Interpretations}
\index{functorial interpretation}
\label{sec:3-term-func-interp}
The main criterion for showing the termination of a 3-polygraph is given by
constructing a suitable interpretation of its 2-cells in 2-categories for which
a reduction order is known. This generalizes the technique of reduction
functions introduced in \secr{red-fun}.

\newcommand{\tfun}[1]{\left[#1\right]}

\begin{proposition}
  \label{P:TerminationInterpretation}
  \label{prop:term-interp}
  Let~$P$ be a 3-polygraph. The following statements are equivalent.
  \begin{enumerate}
  \item\label{prop:term-interp1} The 3-polygraph~$P$ terminates.
  \item\label{prop:term-interp2} There exists a 2-category~$C$, equipped with a reduction
    order~$\succeq$, and a 2\nbd-func\-tor $\tfun-:\freecat{P_2}\to C$, such
    that $\tfun\phi\succ\tfun\psi$ for every 3-generator $A:\phi\To\psi$ in
    $P_3$.
  \end{enumerate}
\end{proposition}
\begin{proof}
  \cref{prop:term-interp1} $\Rightarrow$ \cref{prop:term-interp2}.
  If the polygraph~$P$ terminates, we can take $C=\freecat{P_2}$ and
  $\tfun-=\id_{P_2^*}$. By \cref{P:TerminationOrder},
  taking~$\freecat{\TO}$ for~$\succeq$ gives a termination order for~$P$, hence
  a reduction order on~$C$ such that $\tfun\phi\succ\tfun\psi$ for every
  3-generator $A:\phi\To\psi$ in~$P_3$.

  \cref{prop:term-interp2} $\Rightarrow$ \cref{prop:term-interp1}.
  Conversely, $\tfun-$ being a functor, and~$\succeq$ being compatible with the
  compositions of~$C$ imply that, for every rewriting step $F:\phi\TO\psi$, we
  have $\tfun f\succ\tfun g$ in~$C$. Now, assume that~$P$ does not
  terminate. Then there exists an infinite sequence of composable rewriting steps
  in~$P$, yielding an infinite decreasing sequence for~$\succeq$, which is
  excluded by well-foundedness of~$\succeq$.
\end{proof}

\subsection{The number of generators}
A first very simple situation in which a 3\nbd-poly\-graph terminates is when each
rewriting rule (and thus each rewriting step) decreases the number of
2-generators in 2-cells (we recall that the number~$\sizeof\phi$ of 2-generators in a 2-cell
$\phi$ was formally defined in \cref{sec:2-occurrences}). Namely, the usual
order $\geq$ on natural numbers is a reduction order because it is well-founded
and addition is strictly increasing. Thus, the following statement is
a direct application of \propr{term-interp}.

\begin{proposition}
  \label{prop:3-term-nb-gen}
  A 3-polygraph~$P$ such that $\sizeof\phi>\sizeof{\phi'}$ for every rewriting rule
  $A:\phi\TO\phi'$ in~$P_3$, is terminating.
\end{proposition}

\begin{example}
  \label{ex:ass-term}
  Consider the 3-polygraph~$\MonPoly$ of monoids described in \cref{ex:ass,ex:ass-cp}, and
  write~$\MonPoly'$ for the polygraph obtained from $\MonPoly$ by removing the
  rewriting rule~$A$ corresponding to associativity.
  For both rewriting rules $L$ and $R$, the number of generators of the source
  is $2$ and the number of generators of the target is $0$. By
  \propr{3-term-nb-gen}, the polygraph~$\MonPoly'$ is thus terminating.

  However, in the rule~$A$, the number of generators in the
  source and in the target are both~$2$. Therefore the above proposition does not
  apply to $\MonPoly$. We see below a more
  general method which is able to handle this case, see \exr{mon-term}.
\end{example}

\subsection{A 2-category of posets}
\nomenclature[Ord]{$\Ord$}{2-category of posets}
We define a 2-category~$\Ord$ with one $0$\nbd-cell~$\star$, $1$-cells are
posets, and $2$-cells are weakly increasing functions. The composition of two
$1$-cells $(X,\leq_X)$ and $(Y,\leq_Y)$ is given by their cartesian product
$(X\times Y,\leq_{X\times Y})$, where the product order is such that
$(x,y)\leq_{X\times Y}(x',y')$ if and only if $x\leq_X x'$ and $y\leq_Y
y'$. Likewise, the horizontal composition of $2$-cells is given by their
cartesian product, and their vertical composition is the usual composition of
functions.

We write~$\succeq$ for the pointwise order on $2$-cells: given $f,g:X\To Y$, we
have $f\succ g$ if and only if $f(x)>g(x)$ for every $x\in X$. This order
is always compatible with horizontal and vertical composition of~$2$-cells.

\begin{lemma}
  \label{lem:pointwise-order-wf}
  Given posets~$X$ and~$Y$ such that~$X$ is non-empty and~$Y$ is well-founded,
  the order $\succeq$ on functions $X\To Y$ is well-founded.
\end{lemma}
\begin{proof}
  Fix an arbitrary element~$x$ in~$X$, which is supposed to be non-empty. An
  infinite strictly decreasing sequence of functions $f_0\succ f_1\succ\ldots$
  would induce an infinite strictly decreasing sequence $f_0(x)>f_1(x)>\ldots$
  of elements of~$Y$. This is excluded by well-foundedness of~$Y$.
\end{proof}

This provides us with the following technique for showing termination of
polygraphs, first considered in~\cite{lafont2003towards}:

\begin{proposition}
  \label{P:TerminationByFunctor}
  Let $P$ be a 3-polygraph and $(X,\leq)$ a non-empty well-founded
  po\-set. Suppose that to each~2-generator $\alpha:u\To v$ in~$P_2$
  whose source has length~$m=\sizeof{u}$ and target has
  length~$n=\sizeof{v}$ is assigned  a strictly increasing function
  $\tfun{\alpha}:X^m\to X^n$, and write $\tfun-:\freecat{P_2}\to\Ord$ for the induced
  2\nbd-func\-tor. If, for every rewriting rule $A:\phi\TO\psi$ in~$P_3$, we
  have $\tfun\phi\succ\tfun\psi$ then the polygraph~$P$ terminates.
\end{proposition}
\begin{proof}
  Write~$C$ for the full sub-$2$-category of~$\Ord$ whose $1$-cells are the
  powers~$X^n$ of the poset~$X$. The order on~$X^n$ is well-founded as a product
  of well-founded orders, and therefore the induced order $\succeq$ on $2$-cells
  $X^m\To X^n$ is also well-founded by \lemr{pointwise-order-wf}. For every
  $2$-cell $\phi$, we have that $\tfun\phi$ is a strictly increasing function
  and thus the order is compatible with composition of~$2$-cells in the sense
  of~\defr{3-red-order}, it is thus a reduction order on~$C$. We conclude with
  \propr{term-interp}.
\end{proof}

\begin{example}
  \index{polygraph!of monoids}
  \label{ex:mon-term}
  Consider the 3-polygraph~$\MonPoly$ of monoids introduced
  in~\exr{monoid}. Graphically, its rules are
  \begin{align*}
    A:\satex{mon-assoc-l}
    &\TO
    \satex{mon-assoc-r}
    &
    L:\satex{mon-unit-l}
    &\TO
    \satex{mon-unit-c}
    &
    R:\satex{mon-unit-r}
    &\TO
    \satex{mon-unit-c}
    \pbox.
  \end{align*}
  Writing~$X=\N\setminus\set{0}$ equipped with the usual order, we define
  \[
    \tfun{\satex{mu}}(i,j) = 2i+j
    \qquad\text{and}\qquad
    \tfun{\satex{eta}}()= 1.
  \]
  By \cref{P:TerminationByFunctor}, the following strict inequalities imply that the $3$-polygraph $\MonPoly$ terminates:
  \begin{align*}
    \tfun{\:\satex{mon-assoc-l}\:}(i,j,k)
    = 4i+2j+k
    &> 2i+2j+k
    = \tfun{\:\satex{mon-assoc-r}\:}(i,j,k) 
    \\
    \tfun{\:\satex{mon-unit-l}\:}(i)
    = 2+i
    &> i
    = \tfun{\:\satex{mon-unit-c}\:}(i)
    \\
    \tfun{\:\satex{mon-unit-r}\:}(i)
    = 2i+1
    &> i
    = \tfun{\:\satex{mon-unit-c}\:}(i)
    \pbox.
\end{align*}
\end{example}

\begin{example}
  \index{polygraph!of permutations}
  \label{ex:sym-term}
  Consider the 2-polygraph of permutations introduced in \secr{sym-pres}. It has
  one $0$-generator $\star$, one $1$-generator $a$, one $2$-generator
  \[
    \satex{gamma}
  \]
  and its rules are
  \begin{align*}
    N:
    \satex{sym-l}
    &\TO
    \satex{sym-r}
    &
    Y:
    \satex{yb-l}
    &\TO
    \satex{yb-r}
    \,\pbox.
  \end{align*}
  We consider the well-founded poset~$X=\N\setminus\set{0}$ equipped with the
  usual order and consider the interpretation
  \[
    \tfun{\satex{gamma}}(i,j)=(i+j,i)\pbox.
  \]
  By \cref{P:TerminationByFunctor}, we ensure that the polygraph
  is terminating by showing that the rewriting rules are strictly decreasing:
  \begin{align*}
    \tfun{\:\satex{sym-l}\:}(i,j)=(2i+j,i+j)&>(i,j)=\tfun{\:\satex{sym-r}\:}(i,j)
    \\
    \tfun{\:\satex{ybs-l}\:}(i,j,k)=(2i+j+k,i+j,i)&>(i+j+k,i+j,i)=\tfun{\:\satex{ybs-r}\:}(i,j,k)\pbox.
  \end{align*}
\end{example}

\begin{example}
  By combining the two previous examples, one can construct a convergent
  rewriting system corresponding to the theory of commutative monoids (see
  \cref{sec:pres-com-mon}) which is terminating~\cite{lafont2003towards}.
\end{example}

\begin{remark}
  \label{R:TerminationByFunctorFails}
  Consider a 3-polygraph~$P$ containing a $3$-generator
  \[
    A:\phi\TO\psi:u\To v:x\to y
  \]
  such that $\sizeof{v}=0$, \ie $x=y$ and $v=\unit{x}$. Then $X^0$ is reduced to
  one element and we cannot have $[\phi]\succ[\psi]$. Therefore, in this case,
  we cannot use \cref{P:TerminationByFunctor} to show the termination of the
  polygraph.
\end{remark}

\begin{remark}
  Proposition~\ref{P:TerminationByFunctor} can be generalized by taking a
  possibly different well-founded poset~$X(a)$ for each 1-generator $a$ of the
  $3$\nbd-poly\-graph~$P$. In that case, the interpretation of each 2-generator
  $\alpha:u\to v$ is replaced by an increasing map $\tfun{\alpha}:X(u)\to X(v)$,
  where~$X$ is extended to a $1$-functor from the free
  $1$-category~$\freecat{P_1}$ to the underlying $1$-category of~$\Ord$.
\end{remark}

\section{Termination by Derivations}
\label{sec:3-term-der}
The above technique is often not applicable for rewriting systems such that some
2-generators have multiple outputs (with the notable exception of
\exr{sym-term}). We present here another technique for showing the termination
of 3-polygraphs, which is due to
Guiraud~\cite{guiraud2006termination,guiraud2009higher} and is based on the
following intuition. As suggested by the string diagrammatic representation, we
can think of a 2-cell in a polygraph as some kind of electric circuit built
from basic components, the 2-generators. For each of those components, we
are going to specify how the current is transmitted from inputs to outputs (in
both directions, from top to bottom and from bottom to top), as well as how much
heat it emits when the current flows through. Finally, if the rewriting
rules are such that rewriting a circuit strictly decreases its heat, and the heats are
taken in a well-founded order, then we will be able to conclude that the
rewriting system is terminating.

\subsection{The category of contexts}
\index{context}
\index{category!of contexts}
\newcommand{\contexts}[1]{\mathcal{K}_{#1}}
Given a $2$-polygraph~$P$, we write $\contexts{P}$ for the \emph{category of
  contexts} of~$P$. The objects of this category are the $2$-cells
of~$\freecat{P_2}$ and a morphism from a $2$-cell $\phi:u\To v:x\to y$ to a
$2$-cell $\phi':u'\To v':x'\to y'$ is a context~$K$ of type $(u,v)$, as defined
in \secr{2-context}, such that $K[\phi]=\phi'$.

\subsection{Natural system}
\index{natural system}
A \emph{natural system} on a $2$-polygraph~$P$ is a functor
$N:\contexts{P}\to\Ab$, associating an abelian group $N_\phi$ to every $2$-cell
$\phi\in\freecat{P_2}$ and a morphism of groups $N_K:N_\phi\to N_{\phi'}$ to
every context $K$ such that $K[\phi]=\phi'$. This notion is a 2-categorical
variant of the one already encountered in \cref{sec:natural-system}.

By abuse of notation, given an element $n\in N_\phi$ and a context $K$ of
suitable type, we sometimes write $K[n]$ instead of $N_K(n)$. Moreover, given
$i$\nbd-com\-posable morphisms $\phi$ and $\phi'$, for $i\in\set{0,1}$, and elements
$n\in N_\phi$ and $n'\in N_{\phi'}$, we often write $n\comp i\phi'$ (\resp
$\phi\comp in'$) instead of $N_{K'}(n)$ (\resp $N_K(n')$) where $K'$ is the
context $X\comp i\phi'$ (\resp $K$ is the context $\phi\comp iX$).

\subsection{Derivation}
\index{derivation}
\label{sec:2derivation}
Given a natural system~$N$ on a $2$-polygraph~$P$, a \emph{derivation}~$d$
of~$P$ into~$N$ is a function which to every $2$-cell $\phi$ in $\freecat{P_2}$
associates an element of the group~$N_\phi$, in such a way that
\begin{align*}
  d(\phi\comp i\psi)
  &=
  d(\phi)\comp i\psi
  +
  \phi\comp i d(\psi)
\end{align*}
for suitably $i$-composable $2$-cells $\phi$ and $\psi$ in $\freecat{P_2}$, with
$i\in\set{0,1}$.
Note that such a derivation~$d$ is uniquely determined by the images
$d(\alpha)\in N_\alpha$ of the $2$-generators $\alpha\in P_2$.

\begin{lemma}
  \label{lem:der-id}
  Given a derivation~$d$ as above and a 1-cell $u\in\freecat{P_1}$, we have
  $d(\unit{u})=0$.
\end{lemma}
\begin{proof}
  We have
  \[
    d(\unit{u})
    =
    d(\unit{u}\comp1\unit{u})
    =
    d(\unit{u})\comp1\unit{u}
    +
    \unit{u}\comp1 d(\unit{u})
    =
    d(\unit{u})+d(\unit{u})
  \]
  from which we conclude.
\end{proof}

\begin{example}
  We write~$Z:\contexts P\to\Z$ for the \emph{trivial natural system} which
  sends every object to $\Z$ and every morphism to the identity on~$\Z$. Fix a
  $2$\nbd-gene\-rator $\alpha$ in~$P_2$. The operation introduced in
  \secr{2-occurrences}, which to a $2$-cell $\phi$ in~$\freecat{P_2}$ associates
  the number~$\sizeof[\alpha]\phi$ of occurrences of~$\alpha$ in~$\phi$, is the
  derivation of~$P$ into the trivial natural system such that
  $\sizeof[\alpha]\alpha=1$ and $\sizeof[\alpha]\beta=0$ for $\beta\in P_2$ such
  that $\beta\neq\alpha$.
\end{example}

\subsection{A natural system of interest}
\newcommand{\currents}[3]{\mathcal{O}(#1,#2,#3)}
Suppose fixed
\begin{itemize}
\item a $2$-functor $X:\freecat{P_2}\to\Ord$,
\item a $2$-functor $Y:(\freecat{P_2})^\co\to\Ord$, where $(\freecat{P_2})^\co$
  is the $2$-category obtained from $\freecat{P_2}$ by formally changing the
  direction of all $2$-cells,
\item a commutative monoid~$(M,+,0)$ whose addition is strictly increasing.
\end{itemize}
We define a natural system $N:\contexts P\to\Ab$ as follows.
\begin{itemize}
\item To every $2$-cell~$\phi:u\To v$ in $\freecat{P_2}$, $N$ associates the
  monoid~$N_\phi$ of functions $X_u\times Y_v\to M$ with addition being induced
  pointwise by the one in~$M$.
\item For every $2$-cell~$\phi:u\To v:x\to y$ in $\freecat{P_2}$, and
  every pair of
  $1$-cells $w:x'\to x$, $w':y\to y'$ as in
  \[
    \xymatrix{
      x'\ar[r]^w&x\ar@/^3ex/[r]^u\ar@/_3ex/[r]_v\ar@{}[r]|{\phi\Longdownarrow\phantom\phi}&y\ar[r]^{w'}&y'
    }
  \]
  the image of the context $K=w\comp0-\comp0 w'$ is the group morphism
  \[
    N_K:N_\phi\to N_{w\comp0\phi\comp0w'}
  \]
  which sends a function
  \[
    f:X_u\times Y_v\to M
  \]
  to the function
  \[
    N_K(f):X_w\times X_u\times X_{w'}\times Y_w\times Y_v\times Y_{w'}\to M
  \]
  obtained by precomposing~$f$ with the canonical projection
  \[
    X_w\times X_u\times X_{w'}\times Y_w\times Y_v\times Y_{w'}
    \to
    X_u\times Y_v
    \pbox.
  \]
\item For every $2$-cell $\phi:u\To v$ and $2$-cells $\psi:u'\To u$ and
  $\psi':v\To v'$ as in
  \[
    \xymatrix@C=14ex{
      x
      \ar@/^7ex/[r]^{u'}
      \ar@{{}{ }{}}@/^5ex/[r]|{\psi\Longdownarrow\phantom\psi}
      \ar@/^2.5ex/[r]|{u}
      \ar@{}[r]|{\phi\Longdownarrow\phantom\phi}
      \ar@/_2.5ex/[r]|{v}
      \ar@{{}{ }{}}@/_5ex/[r]|{\psi'\!\!\Longdownarrow\!\!\phantom{\psi'}}
      \ar@/_7ex/[r]_{v'}
      &y
    }
  \]
  the image of the context $K=\psi\comp1-\comp1\psi'$ is the group morphism
  \[
    N_K:N_\phi\to N_{\psi\comp1\phi\comp1\psi'}
  \]
  which sends a function
  \[
    f:X_u\times X_v\to M
  \]
  to the function
  \[
    N_K(f):X_{u'}\times Y_{v'}\to M
  \]
  defined by
  \[
    N_K(f)=f\circ(X_\psi\times Y_{\psi'})
    \pbox.
  \]
\end{itemize}
The above conditions entirely determine the derivation~$N$, which we often denote by
$\currents XYM$ to make clear the dependency on~$X$, $Y$,
and~$M$.
Note that such a natural system is entirely determined by the data of
\begin{itemize}
\item the posets $X_a$ and $Y_a$ for every $a$ in $P_2$,
\item the functions $X_\alpha:X_u\to X_v$ and $Y_\alpha:Y_v\to Y_u$ for every
  $2$-generator $\alpha:u\To v$ in~$P_2$,
\end{itemize}
where $X_{a_1\ldots a_n}=X_{a_1}\times\ldots\times X_{a_n}$.

\begin{remark}
  We could generalize the definition to the case where~$X$ and~$Y$ are
  objects in an arbitrary cartesian category and $M$ is a commutative monoid
  internal to this category.
\end{remark}

Given a $2$-cell $\phi:u\To v$ in~$\freecat{P_2}$, the monoid $N_\phi$ is
canonically equipped with the order such that, for elements
$f,g:X_u\times Y_v\to M$ of $N_\phi$ we have
\[
  f>g
  \qquad\text{if and only if}\qquad
  f(x,y)>g(x,y)
  \text{ for every $(x,y)\in X_u\times Y_v$.}
\]
The above construction was introduced by Guiraud~\cite{guiraud2006termination,
  guiraud2019rewriting, guiraud2009higher}. It is the basis of the
following useful termination criterion.

\begin{theorem}
  \label{thm:3-der-term}
  Consider a 3-polygraph~$P$. Suppose given
  \begin{itemize}
  \item two 2-functors $X:\freecat{P_2}\to\Ord$ and
    $Y:(\freecat{P_2})^\co\to\Ord$ such that for every $1$-generator $a\in P_1$
    the posets $X_a$ and $Y_a$ are non-empty, and $X_\phi\geq X_\psi$ and
    $Y_\phi\geq Y_\psi$ for every $3$-generator $A:\phi\To\psi$ in~$P_3$,
  \item a well-founded partially ordered commutative monoid~$(M,+,0)$ such that
    addition is strictly increasing,
  \item a derivation~$d$ from the underlying 2-polygraph of~$P$ to
    $\currents XYM$ such that $d(\phi)>d(\psi)$ for every 3-generator
    $A:\phi\To\psi$ in~$P_3$.
  \end{itemize}
  Then the polygraph~$P$ is terminating.
\end{theorem}
\begin{proof}
  Let $K[\alpha]:\chi\TO\chi'$ be a rewriting step, for some $2$-generator
  $\alpha:\phi\TO\phi'$ in $P_2$ and context~$K\in\contexts P$. By
  \lemr{2-context}, the context~$K$ can be written in the form
  \[
    K=\psi\comp1(w\comp0-\comp0w')\comp1\psi'
  \]
  for suitably typed $1$-cells $w$ and $w'$ in $\freecat{P_1}$ and $2$-cells
  $\psi$ and $\psi'$ in $\freecat{P_2}$. Using the definition of derivations,
  see \secr{derivation}, and \lemr{der-id}, we have that $d(K[\chi])$ is equal
  to
  \[
    d(\psi)\comp1(w\comp0\chi\comp0w')\comp1\psi'
    +
    \psi\comp1(w\comp0 d(\chi)\comp0w')\comp1\psi'
    +
    \psi\comp1(w\comp0\chi\comp0w')\comp1 d(\psi')
  \]
  and similarly for $d(K[\chi'])$. By hypothesis, we have $d(\chi)>d(\chi')$,
  and therefore
  \[
    \psi\comp1(w\comp0 d(\chi)\comp0w')\comp1\psi'
    >
    \psi\comp1(w\comp0 d(\chi)\comp0w')\comp1\psi'
    \pbox.
  \]
  Moreover, since $X$ and $Y$ are decreasing on generators, by functoriality of
  $X$ and $Y$ we have $X_\chi\geq X_{\chi'}$, and thus
  \[
    d(\psi)\comp1(w\comp0\chi\comp0w')\comp1\psi'
    \geq
    d(\psi)\comp1(w\comp0\chi'\comp0w')\comp1\psi'
  \]
  and similarly
  \[
    \psi\comp1(w\comp0\chi\comp0w')\comp1 d(\psi')
    \geq
    \psi\comp1(w\comp0\chi'\comp0w')\comp1 d(\psi')
    \pbox.
  \]
  Finally, since addition is strictly increasing, we deduce
  \[
    d(K[\chi])>d(K[\chi'])
    \pbox.
  \]
  An infinite sequence of rewriting steps starting from a $2$-cell
  $\phi:u\To v$, would thus induce a strictly decreasing sequence
  \[
    f_0>f_1>f_2>\ldots
  \]
  of elements of $d(\phi)$, \ie functions $X_u\times Y_v\to M$. Since $X_u$ and
  $Y_v$ are supposed to be non-empty, we can pick an element
  $(x,y)\in X_u\times Y_v$, and we would have a strictly decreasing sequence
  \[
    f_0(x,y)>f_1(x,y)>f_2(x,y)>\ldots
  \]
  of elements of~$M$, which is excluded by hypothesis. The polygraph~$P$ is thus
  terminating.
\end{proof}

In order to give some intuition, let us consider a $3$-polygraph~$P$. A $2$-generator
\[
  \alpha:a_1a_2\ldots a_m\To b_1b_2\ldots b_n
\]
in~$P_2$, where the $a_i$ and $b_i$ are $1$-generators in~$P_1$, can be seen as
an operation with $m$ inputs and $n$ outputs
\[
  \satex{sd-alpha}
\]
which, as the figure suggests, can be thought of as a building piece of some
electrical circuit. The poset $X_{a_i}$ (\resp $Y_{a_i}$, $X_{b_i}$, $Y_{b_i}$)
is the set of possible values for the currents flowing into $a_i$ (\resp out
from $a_i$, out from $b_i$, into~$b_i$). The function
\[
  X_\alpha:X_{a_1}\times\ldots\times X_{a_m}\to X_{b_1}\times\ldots\times X_{b_n}
\]
then indicates, given currents flowing into the inputs $a_i$, what currents we
get from the outputs~$b_i$. Similarly, the function
\[
  Y_\alpha:Y_{b_1}\times\ldots\times Y_{b_n}\to Y_{a_1}\times\ldots\times Y_{a_m}
\]
indicate the current we obtain from the $a_i$ if we use the device ``upside down''
and flow currents into the $b_i$. Finally, the monoid~$M$ can be thought of as
the possible values for ``heat'' emitted by our electrical circuit and the function
\[
  d(\alpha):X_{a_1}\times\ldots\times X_{a_m}\times Y_{b_1}\times\ldots\times Y_{b_n}\to M
\]
indicates, given currents flowing into the $a_i$ and into the $b_i$, the heat
that our circuit produces. The fact that it is a derivation amounts to impose
that the heat produced by a circuit is the sum of the heat emitted by its
components. Finally, the hypotheses of \thmr{3-der-term} ensure that rewriting a
circuit will always transform it into a ``colder'' (\ie less heat-emitting)
circuit.

\begin{example}
  \index{polygraph!of permutations}
  \label{ex:sym-term2}
  Following~\cite[Section~5.4]{guiraud2009higher}, let us apply the above
  technique to show that the rewriting system for permutations, already
  considered in \cref{ex:sym-term}, is terminating. We suppose here that
  \begin{itemize}
  \item the poset~$X_a$ associated to the $1$-generator $a$ is $\N$ equipped
    with the usual order, and $X$ is defined on the $2$-generator by
    \[
      X\pa{\satex{gamma}}(i,j)=(j+1,i)\pbox,
    \]
  \item $Y_a=\set{*}$ is the terminal poset (reduced to one element~$*$), and
    $Y$ is defined on the $2$-generator by
    \[
      Y\pa{\satex{gamma}}(*,*)=(*,*)\pbox,
    \]
  \item the monoid $M$ is the additive monoid $\N$,
  \item the derivation~$d$ is defined on the $2$-generator by
    \[
      d\pa{\satex{gamma}}(i,j)=i
    \]
    (more precisely, the derivation takes four arguments $(i,j,k,l)$, but the
    two last arguments are necessarily equal to $*$, the only element of~$Y_a$,
    and are thus omitted).
  \end{itemize}
  The rewriting rules make $X$ weakly decrease
  \[
    \begin{array}{r@{\ =\ }c@{\ =\ }l}
      X\pa{\:\satex{sym-l}\:}(i,j)
      &
      (i+1,j+1)
      \geq
      (i,j)
      &
      X\pa{\:\satex{sym-r}\:}(i,j)
      \\
      X\pa{\:\satex{yb-l}\:}(i,j,k)
      &
      (k+2,j+1,i)
      &
      X\pa{\:\satex{yb-r}\:}(i,j,k)
    \end{array}
  \]
  as well as obviously $Y$, and make the derivation strictly decrease
  \[
    \begin{array}{r@{\ =\ }c@{\ >\ }c@{\ =\ }l}
      d\pa{\:\satex{sym-l}\:}(i,j)
      &
      i+j+1
      &
      0
      &
      d\pa{\:\satex{sym-r}\:}(i,j)
      \\
      d\pa{\:\satex{yb-l}\:}(i,j,k)
      &
      2i+j+1
      &
      2i+j
      &
      d\pa{\:\satex{yb-r}\:}(i,j,k)
      \pbox.
    \end{array}
  \]
  By \thmr{3-der-term}, the $3$-polygraph is thus terminating.
\end{example}

\noindent
A presentation which cannot be handled with the techniques of
\secr{3-term-func-interp}, see \cref{R:TerminationByFunctorFails}, but can be
handled with derivations, is given in \cref{sec:pearls}.


\chapter{Coherent Presentations of 2-Categories}
\label{chap:3coh}
In this chapter, we generalize the definitions and results of
\cref{chap:2-Coherent} from categories to $2$-categories, following~\cite{guiraud2009higher,GuiraudMalbos12mscs}. In
\cref{sec:2coh-pres}, we introduce the notion of coherent presentation of a
$2$-category by a $(4,2)$-polygraph, where the $4$-generators encode the
relations among relations.
We explain in \cref{Section:CaseOf3Polygraphs} that, in the case of convergent
polygraphs, we can construct the Squier completion, which is a coherent
completion whose $4$-generators come from confluence diagrams for the
critical branchings. We show in \cref{sec:pearls} that, contrarily to the case of categories, a
$2$-category presented by a finite convergent $3$\nbd-poly\-graph is not necessarily of finite
derivation type.
In \cref{sec:32PRO}, we develop a $3$-dimensional generalization of the notion
of PRO, for which coherent presentations can be given by
$(4,2)$-polygraphs. This allows us, in \cref{sec:mon-coh}, to use the
constructions of coherent presentations to obtain coherence results
such as Mac Lane's coherence theorem for monoidal categories, as well as
generalizations to symmetric and braided monoidal categories in
\cref{Section:CoherenceSym}.

\section{Coherent Presentations of 2-Categories}
\label{sec:2coh-pres}

\subsection{(4,2)-polygraphs}
\label{sec:(4,2)Polygraph}
\index{42-polygraph@$(4,2)$-polygraph}
\index{polygraph!42-@$(4,2)$-}
A \emph{$(4,2)$-polygraph} is a pair $(P,P_4)$ consisting of a
$3$-polygraph~$P$, and a cellular extension $P_4$ of the free $3$-category $\freegpd{P}$ over~$P$.
It thus consists of a diagram of sets and functions
\[
\vxym{
&
  \ar@<-.5ex>[dl]_-{\sce0}
  \ar@<+.5ex>[dl]^-{\tge0}
P_1\ar[d]^{\ins1}
&
  \ar@<-.5ex>[dl]_-{\sce1}
  \ar@<+.5ex>[dl]^-{\tge1}
P_2\ar[d]^{\ins2}
&
  \ar@<-.5ex>[dl]_-{\sce2}
  \ar@<+.5ex>[dl]^-{\tge2}
P_3\ar[d]^{\ins3}
&
  \ar@<-.5ex>[dl]_-{\sce3}
  \ar@<+.5ex>[dl]^-{\tge3}
P_4
\\
P_0
&
  \ar@<-.5ex>[l]_-{\freecat{\sce0}}
  \ar@<+.5ex>[l]^-{\freecat{\tge0}}
\freecat{P_1}
&
  \ar@<-.5ex>[l]_-{\freecat{\sce1}}
  \ar@<+.5ex>[l]^-{\freecat{\tge1}}
\freecat{P_2}
&
  \ar@<-.5ex>[l]_-{\freecat{\sce2}}
  \ar@<+.5ex>[l]^-{\freecat{\tge2}}
\freegpd{P_3}
}
\]
together with the compositions and identities of the underlying
$(3,2)$-category
\[
  \xymatrix{
    P_0&
    \ar@<-.5ex>[l]_-{\freecat{\sce0}}
    \ar@<+.5ex>[l]^-{\freecat{\tge0}}
    \freecat{P_1}
    &
    \ar@<-.5ex>[l]_-{\freecat{\sce1}}
    \ar@<+.5ex>[l]^-{\freecat{\tge1}}
    \freecat{P_2}
    &
    \ar@<-.5ex>[l]_-{\freecat{\sce2}}
    \ar@<+.5ex>[l]^-{\freecat{\tge2}}
    \freegpd{P_3},
}
\]
whose \emph{source} and \emph{target maps} $\sce{i}$ and $\tge{i}$ satisfy the
globular relations
\[
  \freecat{\sce{i}}\circ\sce{i+1} = \freecat{\sce{i}}\circ\tge{i+1}
  \qqtand
  \freecat{\tge{i}}\circ\sce{i+1} = \freecat{\tge{i}}\circ\tge{i+1},
\]
for every $0\leq i \leq 2$. The elements of the cellular extension $P_4$ are
called the \emph{4-generators} of the polygraph~$P$. We write
$\Lambda: F \ttfl G$ for a $4$-generator $\Lambda$ in~$P_4$ such that
$\sce2(\Lambda)=F$ and $\tge2(\Lambda)=G$. Given a $4$-polygraph $P$, we write
$\tpol3 P$ for its underlying $3$-polygraph.

\subsection{Coherence}
\index{coherent!42-polygraph@$(4,2)$-polygraph}
A $(4,2)$-polygraph~$P$ is \emph{coherent} when for any parallel $3$\nbd-cells
\index{parallel!3-cells@$3$-cells}
$F,G:\phi\TO\psi$ of $\freegpd{\tpol3 P}$, the free $(3,2)$-category generated
by the underlying $3$-polygraph of~$P$, we have $F\simeq^PG$, \ie $F$
and~$G$ are related by the congruence generated by the cells in~$P_4$.

Given a $(3,2)$-category~$C$ a cellular extension~$X$ of~$C$ is \emph{acyclic}
when $F\simeq^XG$ holds for every pair of parallel $3$-cells $F$ and
$G$ of~$C$. Here
$\simeq^X$ is the congruence generated by~$X$ (which is defined in the expected way,
generalizing the definition of~\cref{sec:2-cat-quot}). Given a $(3,2)$\nbd-poly\-graph~$P$, a
$(4,2)$-polygraph~$(P,P_4)$ is thus coherent precisely when $P_4$ is an acyclic extension of $\freegpd{P}$.

\subsection{Polygraphs of finite derivation type}
\index{finite derivation type!3-polygraph@$3$-polygraph}
\index{finite derivation type!2-category@$2$-category}
\index{FDT (finite derivation type)}
The property of finite derivation type defined in \cref{chap:2-fdt} for
$2$-polygraphs is extended to $3$-polygraphs as follows.
One says that a $3$-polygraph $P$ has \emph{finite derivation type} when it is
finite and when the free $(3,2)$-category $\freegpd{P}$ admits
a finite acyclic cellular extension. As in the case of presentation of $1$-categories, given
two presentations of the same $2$-category by finite $3$-polygraphs, the following
result proves that both have finite derivation type or none at all.

\begin{theorem}
  \label{TDFTietzeInvariant3Pol}
  Let $P$ and $Q$ be two Tietze equivalent 3-polygraphs such that~$P_2$ and~$Q_2$ are finite.
  Then~$P$ has finite derivation type if and only if $Q$ has finite derivation type.
\end{theorem}
\begin{proof}
  The proof is similar to the one given in the case of $1$-categories by
  \cref{thm:TDFTietzeInvariant}.
\end{proof}

As a consequence of \cref{TDFTietzeInvariant3Pol}, one can say that a $2$-category has \emph{finite derivation type} when it admits a presentation by a $3$-polygraph having finite derivation type.
The property of having finite derivation type is invariant by Tietze equivalence for finite $3$-polygraphs. This is not the case for infinite ones as shown by the following example~\cite[Section~4.3.10]{guiraud2009higher}.

\begin{example}
  \label{Example:infinitePolTDF}
  Consider the $3$-polygraph~$P$ with one $0$-generator, one
  $1$\nbd-gene\-rator, three $2$-generators $\satex{itdf-alpha-vsmall}$,
  $\satex{itdf-beta-vsmall}$, $\satex{itdf-gamma-vsmall}$
  and the following two $3$-generators:
  \begin{align*}
    A&:
    \satex{itdf-A-l}
    \TO
    \satex{itdf-A-r}\;,
    &
    B&:
    \satex{itdf-B-l}
    \TO
    \satex{itdf-B-r}
    \pbox.
  \end{align*}
  We prove that the $(3,2)$-category $\freegpd{P}$ admits an empty acyclic extension
  and thus has finite derivation type.

  The $3$-polygraph $P$ is Tietze equivalent to the $3$\nbd-poly\-graph $Q$
  defined the same way as~$P$ except for the orientation of the $3$-cell~$A$:
  \begin{align*}
    A&:
    \satex{itdf-A-r}
    \TO
    \satex{itdf-A-l}\;,
    &
    B&:
    \satex{itdf-B-l}
    \TO
    \satex{itdf-B-r}
    \,\pbox.
  \end{align*}
  In this polygraph, we introduce the notation $\satex{itdf-gammak-small}$ for the
  $2$-cell defined by induction on the natural number~$k$ by
  \begin{align*}
    \satex{itdf-gamma1}&=\satex{itdf-gamma}\;,
    &
    \satex{itdf-gammak1}&=\satex{itdf-gammak1-def}\pbox.
  \end{align*}
  The polygraph $Q$ is not convergent, but we can complete it into the infinite
  $3$-polygraph $Q_\infty = Q \sqcup \setof{B_k}{k\geq 1}$, where~$B_0$ is $B$
  and $B_k$ is the following $3$-cell:
  \[
    B_k:
    \satex{itdf-Bk-l}
    \TO
    \satex{itdf-Bk-r}
    \pbox.
  \]
  It can be shown that the $3$-polygraph $Q_\infty$ does not have finite derivation
  type. In particular the $(3,2)$-category $\freegpd{Q_{\infty}}$ has an
  infinite acyclic extension $\setof{\Lambda_k}{k\in\N}$ with
  \[
    \satex{itdf-Lambda}
  \]
  which cannot be reduced to a finite one.
\end{example}

\subsection{Generating confluences}
\label{Subsection:BranchingHomotopyBases}
\cref{thm:SquierHomotopicalOnePolygraphs,thm:SquierHomotopical} state that 
the set of critical branchings of a convergent $n$-polygraph $P$ generates an acyclic extension of the $(n,n{-}1)$-category $\freegpd{P}$ when $n\leq 2$. The proof of this result can be
extended to $3$-polygraphs as follows.  Given a convergent $3$-polygraph $P$, a
\emph{family of generating confluences of~$P$} is a cellular extension of the
free $(3,2)$\nbd-cate\-gory~$\freegpd{P}$ that contains exactly one $4$-cell $\Lambda$
of the form
\[
  \xymatrix @R=1.25em@C=3em @!C{
    & {\chi}
    \ar@3 @/^/ [dr] ^{F'}
    \ar@4 []!<0pt,-15pt>;[dd]!<0pt,15pt> ^-*+{\Lambda}
    \\
    {\phi}
    \ar@3 @/^/ [ur] ^{F}
    \ar@3 @/_/ [dr] _{G}
    && {\phi'},
    \\
    & {\psi}
    \ar@3 @/_/ [ur] _{G'}
  }
\]
for every critical branching $(F,G)$ of $P$.
\index{Squier!completion}
We define the \emph{Squier
  completion} of the $3$-polygraph $P$ as the $(4,2)$-polygraph denoted by
$\Sr(P)$ and defined by $\Sr(P)=(P,P_4)$, where $P_4$ is a chosen family of
generating confluences of $P$. As in the case were $n=2$, see
\cref{thm:SquierHomotopical}, we have~\cite[Proposition~4.3.4]{guiraud2009higher}:

\begin{theorem}
  \label{Theorem:SquierCompletion3polygraphs}
  Given a convergent presentation of a 2-category $C$ by a
  3\nbd-poly\-graph~$P$, any Squier completion of $P$ is coherent
  presentation of~$C$.
\end{theorem}

\noindent
As a consequence of \cref{Theorem:SquierCompletion3polygraphs}, a finite convergent $3$-polygraph with a finite set of critical branchings has finite derivation type. 
In particular, a terminating polygraph with no critical branching has finite derivation type.
However, this result fails to generalize to $n$-categories when $n\geq 2$, see \cref{S:PolygraphsFDTn} and~\cite{guiraud2009higher}.
A counterexample to show this for $n=3$ is developed in~\cref{Subsection:main_counter_example}.

\section{Squier's Completion of 3-Polygraphs}
\label{Section:CaseOf3Polygraphs}

\subsection{Non-indexed 3-polygraphs}
\label{Subsection:NonIndexed}
Recall from \cref{sec:3pol-non-indexed} that a $3$\nbd-poly\-graph is
non-indexed when each of its critical branchings is an inclusion one or a
regular one.  It can be proved that a $3$-polygraph with a finite set of
$3$-cells has a finite number of inclusion and regular critical
branchings~\cite[Proposition 5.1.3]{guiraud2009higher}. As a consequence, we
have the following finiteness condition in the non-indexed
case~\cite[Theorem~5.1.4]{guiraud2009higher}:

\begin{theorem}
\label{Thm:3PolyNonIndexe}
A finite, convergent, and non-indexed 3-polygraph has finite derivation type. 
\end{theorem}

\subsection{Confluence in indexed 3-polygraphs}
\label{Subsection:Indexed}
Now let us consider the problem of finite-convergence for finitely indexed $3$\nbd-poly\-graphs (those for which each indexed critical branching has a finite number of normal instances). The situation is more complicated than in the non-indexed case. However, we have the following confluence result~\cite[Proposition 5.3.1]{guiraud2009higher}:

\begin{proposition}
  \label{Theorem:RightIndexedConfluence}
  Let $P$ be a terminating right-indexed (\resp left-indexed)
  3\nbd-poly\-graph. Then $P$ is confluent if and only if every inclusion critical
  branching, every regular critical branching, and every instance of every
  right-indexed (\resp left-indexed) critical branching is confluent.
\end{proposition}
\begin{proof}
  Suppose that $P$ is a terminating right-indexed $3$-polygraph (the
  left-indexed case is similar) such that all its inclusion critical branchings, regular critical
  branchings, and all the normal instances of its right-indexed critical
  branchings are confluent. It is sufficient to prove that every non-normal
  instance of its right-indexed critical branchings is confluent.
  Let us consider an instance of right-indexed critical branching. With the
  notations of \cref{sec:3cb-classification}, it is of the form
  \[
    \xymatrix@R=-4ex{
      &\satex{cp3-indexed-r1}\\
      \satex{cp3-indexed-r}\ar@3[ur]\ar@3[dr]\\
      &\satex{cp3-indexed-r2}
    }
  \]
  for some $2$-cell $\zeta$ in $\freecat P_2$.
  If $\zeta$ is not a normal form, it admits a normal form
  $\nf\zeta$, because~$P$ terminates. There is another instance of the above critical branching with
  $\nf\zeta$ in place of $\zeta$. Since $\nf\zeta$ is a normal form, this is a
  normal instance, so that, by hypothesis, it is confluent. This ensures
  the confluence of the original branching as follows:
  \[
    \xymatrix@R=2ex@C=-6.6ex{
      &\satex{cp3-indexed-r1}\ar@3[rr]&&\satex{cp3-indexed-rn1}\ar@3[dr]\\
      \satex{cp3-indexed-r}\ar@3[ur]\ar@3[dr]\ar@3[rr]&&\satex{cp3-indexed-rn}\ar@3[ur]\ar@3[dr]&&\satex{cp3-indexed-rnn}.\\
      &\satex{cp3-indexed-r2}\ar@3[rr]&&\satex{cp3-indexed-rn2}\ar@3[ur]
    }
  \]
  In this way, we prove that the polygraph $P$ is confluent.
\end{proof}

\subsection{Acyclic extensions of indexed 3-polygraphs}
Let $P$ be a locally confluent and right-indexed (\resp left-indexed)
$3$-polygraph. Suppose that a confluence has been chosen for each inclusion and
regular critical branching and each normal instance of each right-indexed (\resp
left-indexed) critical branching. Let $P_4$ be the cellular extension of the
$(3,2)$-category $\freegpd{P}$ corresponding to these confluence diagrams. We
can prove that if $P$ is convergent and right-indexed (\resp left-indexed), then
$P_4$ forms an acyclic extension of the $(3,2)$-category $\freegpd{P}$, \ie the
$(4,2)$-polygraph~$(P,P_4)$ is
coherent~\cite[Proposition~5.3.3]{guiraud2009higher}.
The proof follows the same scheme as the proof given for
\cref{Theorem:SquierCompletion2polygraphs,Theorem:SquierCompletion3polygraphs}.
It is the same for trivial, inclusion, and regular critical branchings. For
right-indexed (\resp left-indexed) critical branchings, we follow a reasoning
similar to the proof of~\cref{Theorem:RightIndexedConfluence}.
We thus have the following result~\cite[Theorem~5.3.4]{guiraud2009higher}:

\begin{theorem}
  \label{Theorem:FinitelyIndexedFDT}
  A finite, convergent, and finitely indexed 3-polygraph has finite derivation
  type.
\end{theorem}

\noindent
In the next section, we present an illustration of this result with a $3$-polygraph which
is finite, convergent, and right-indexed, and thus has an infinite number of
critical branchings. Yet, the polygraph has finite derivation type thanks to finite
indexation.

\subsection{Example: the 3-polygraph of permutations}
\index{polygraph!of permutations}
\label{Example:PolygraphPermutations}
Consider the $3$\nbd-poly\-graph~$P$ presenting the PRO $\Bij$ of whose morphisms are
permutations, which is introduced in~\cref{sec:sym-pres}.
This polygraph has one $0$-cell, one $1$-cell, one $2$-cell
$\satex{gamma-vsmall}$ and the following two $3$-cells:
\begin{align*}
  I:\satex{sym-l}&\TO\satex{sym-r}\;,
  &
  Y:\satex{yb-l}&\TO\satex{yb-r}
  \;\pbox.
\end{align*}
This polygraph was shown to be terminating in \cref{ex:sym-term,ex:sym-term2}.
We have seen in \cref{sec:sym-pres} that it has three regular and one right-indexed critical
branchings, with the following sources:
\begin{align*}
  \satex{sym-cp1}\;,
  &&
  \satex{sym-cp2}\;,
  &&
  \satex{sym-cp3}\;,
  &&
  \satex{sym-cpi}
 \; \,\pbox.
\end{align*}
From \cref{Theorem:RightIndexedConfluence}, we know that, to show the confluence
of the polygraph, it is sufficient to prove that the three regular critical
branchings are confluent and that each normal instance of the right-indexed one
is. This is, in fact, what we have been doing in \cref{sec:sym-pres}: we only
briefly recall here those confluence diagrams.
First, the three regular critical branchings are confluent:
\[
  \xymatrix@C=7ex{
    \satex{sym-cp1u}
    \ar@3 @/^5ex/ [rr] ^-{I} _-{~}="1"
    \ar@3 @/_5ex/ [rr] _-{I} ^-{~}="2"
    \ar@4 "1"!<0pt,-3ex>;"2"!<0pt,3ex> ^-*+{\Gamma}
    &&
    \satex{sym-cp1d}\;,
  }
\]
\begin{align*}
  \xymatrix@R=3ex@C=3ex{
    & \satex{sym-cp3r}
    \ar@3 [r] ^-{Y} _-{~}="3"
    & \satex{sym-cp3rr}
    \ar@3 [dr] ^-{I}
    \\
    \satex{sym-cp3u}
    \ar@3 [ur] ^-{Y}
    \ar@3 [rrr] _-{I} ^-{~}="4"
    \ar@4 "3"!<0ex,-3ex>;"4"!<0ex,3ex> ^-*+{\Delta}
    &&&
    \satex{sym-cp3d}\;,
  }
  &&
  \xymatrix@R=3ex@C=3ex{
  \satex{sym-cp2u}
  \ar@3 [rrr] ^-{I} _-{~}="1"
  \ar@3 [dr] _-{Y}
  &&&
  \satex{sym-cp2d}\;\pbox.
  \\
  &\satex{sym-cp2r}
  \ar@3 [r] _-{Y} ^-{~}="2"
  \ar@4"1"!<0ex,-3ex>;"2"!<0ex,3ex> ^-*+{\Theta}
  & \satex{sym-cp2rr}
  \ar@3 [ur] _-{I}
  }
\end{align*}
From the characterization of the set of normal forms given in \cref{sec:sym-pres},
we deduce that there are two normal instances of the right-indexed critical
branching: for $k=\satex{id-vsmall}\,$ and $k=\satex{gamma-vsmall}\,$. We
check that both are confluent. For $k=\satex{id-vsmall}\,$, we have:
\[
  \xymatrix@R=-4ex@C=9ex{
    & \satex{sym-cp4l}
    \ar@3 [r] ^-{I} _-{~}="1"
    & \satex{sym-cp4ll}
    \ar@3 [dr] ^-{I}
    \\
    \satex{sym-cp4u}
    \ar@3 [ur] ^-{Y}
    \ar@3 [dr] _-{Y}
    &&& \satex{sym-cp4d}\;\pbox.
    \\
    & \satex{sym-cp4r}
    \ar@3 [r] _-{I} ^-{~}="2"
    & \satex{sym-cp4rr}
    \ar@3 [ur] _-{I} 
    \ar@4 "1"!<0ex,-5ex>;"2"!<0ex,5ex> ^-*+{\Lambda}
  }
\]
For $k=\satex{gamma-vsmall}\,$, we have:
\[
  \xymatrix@R=-4ex@C=7ex{
    & \satex{sym-cp5-2}
    \ar@3 [r] ^-{Y} 
    & \satex{sym-cp5-3}
    \ar@3 [r] ^-{Y} 
    & \satex{sym-cp5-4}
    \ar@3 [dr] ^-{Y}
    \\
    \satex{sym-cp5-1}
    \ar@3 [ur] ^-{Y}
    \ar@3 [dr] _-{Y}
    &&&& \satex{sym-cp5-5}\,\pbox.
    \\
    & \satex{sym-cp5-8}
    \ar@3 [r] _-{Y} 
    & \satex{sym-cp5-7}
    \ar@3 [r] _-{Y} 
    & \satex{sym-cp5-6}
    \ar@3 [ur] _-{Y} 
    \ar@4 "1,3";"3,3" ^-*+{\Lambda'}
  }
\]
The $3$-polygraph $P$ is finite, convergent, and finitely indexed, by
\cref{Theorem:FinitelyIndexedFDT}, it follows that it has finite derivation
type. More precisely, the five $4$-cells $\Gamma$, $\Delta$, $\Theta$, $\Lambda$,
and $\Lambda'$ form an acyclic extension of the $(3,2)$-category $\freegpd{P}$.

\subsection{Main counterexample: the polygraph of pearls}
\index{pearl}
\index{polygraph!of pearls}
\label{Subsection:main_counter_example}
\label{sec:pearls}
Let us mention a $3$-polygraph, studied in~\cite{guiraud2009higher}, which illustrates the fact
that, without finite indexation, finiteness, and convergence are not
sufficient to ensure finiteness of derivation type.
We consider the $3$-polygraph $P$ of \emph{pearls} with one
$0$-cell~$\star$, one $1$-cell~$a$, three $2$-cells $\satex{pearl-vsmall}$,
$\satex{cap-vsmall}$, and $\satex{cup-vsmall}$, and the following four
$3$-cells:
\begin{align*}
  A:\satex{pearl-cap-l}&\TO\satex{pearl-cap-r}\;,\;
  &
  B:\satex{pearl-cup-l}&\TO\satex{pearl-cup-r}\;,\;
  &
  C:\satex{pearl-zz-l}&\TO\satex{pearl-zz-c}\;,\;
  &
  D:\satex{pearl-zz-r}&\TO\satex{pearl-zz-c}
  \,\pbox.
\end{align*}
We define by induction on the natural number $k$ the $2$-cell $\satex{pearl}^k$ as
follows:
\begin{align*}
  \satex{pearl0}&=\satex{pearl-id}\;,
  &
  \satex{pearlk1}&=\satex{pearlk1-def}
  \pbox.
\end{align*}

Let us show that the polygraph is terminating. The rules $C$ and
$D$ make the number of $2$-generators strictly decrease, while this number
is invariant by the rules $A$ and $B$, so that we only have to show
that the rules~$A$ and~$B$ are terminating. In order to show this, we
apply \cref{thm:3-der-term} and use the notations of this theorem in the
following. The posets associated to the $1$\nbd-gene\-rator are $X_a=Y_a=\N$
equipped with the usual order. The interpretations of the $2$-generators are
\begin{align*}
  X(\satex{pearl-vsmall})(i)&=i+1,
  &
  X(\satex{cap-vsmall})()&=(0,0),
  &
  X(\satex{cup-vsmall})(i,j)&=(),
  \\
  Y(\satex{pearl-vsmall})(i)&=i+1,
  &
  X(\satex{cup-vsmall})()&=(0,0),
  &
  X(\satex{cap-vsmall})(i,j)&=()
  \pbox.
\end{align*}
We take $M$ to be the monoid $(\N,+,0)$ and define the derivation~$d$ by
\begin{align*}
  d(\satex{pearl-vsmall})(i,j)&=0,
  &
  d(\satex{cap-vsmall})(i,j)&=i,
  &
  d(\satex{cup-vsmall})(i,j)&=i
  \pbox.
\end{align*}
For the rule $A$, we have, using the properties of derivation,
\begin{align*}
  d(\freecat{\src2}(A))&=d(\satex{pearl-cap-l-vsmall})=d(\satex{cap-vsmall})\comp1\satex{pearl1-vsmall}+\satex{cap-vsmall}\comp1 d(\satex{pearl1-vsmall}),
  \\
  d(\freecat{\tgt2}(A))&=d(\satex{pearl-cap-r-vsmall})=d(\satex{cap-vsmall})\comp1\satex{1pearl-vsmall}+\satex{cap-vsmall}\comp1 d(\satex{1pearl-vsmall}),
\end{align*}
so that
\begin{align*}
  d(\freecat{\src2}(A))(i,j)
  &=d(\satex{cap-vsmall})(i+1,j)+d(\satex{pearl-vsmall})(0,i)
  \\
  &=(i+1)+0
  \\
  &>i+0
  \\
  &=d(\satex{cap-vsmall})(i,j+1)+d(\satex{pearl-vsmall})(0,j)
  \\
  &=d(\freecat{\tgt2}(A))(i,j)
\end{align*}
and therefore $d(\freecat{\src2}(A))>d(\freecat{\tgt2}(A))$. Similarly,
$d(\freecat{\src2}(B))>d(\freecat{\tgt2}(B))$.
By \cref{thm:3-der-term}, we thus deduce that the polygraph is terminating.

The $3$-polygraph $P$ has four regular critical branchings whose sources are
\begin{align*}
  \satex{pearl-cp-nun}\;,
  &&
  \satex{pearl-cp-unu}\;,
  &&
  \satex{pearl-cp-onu}\;,
  &&
  \satex{pearl-cp-oun}
  \,\pbox.
\end{align*}
It also has one right-indexed critical branching, generated by the
$3$-cells~$A$ and~$B$, with source
\[
  \satex{pearl-onouf}
  \,\pbox.
\]
Thus $P$ is a terminating and right-indexed $3$-polygraph. By application of
\cref{Theorem:RightIndexedConfluence}, the confluence of $P$ can be shown by proving that
its four regular critical branchings and all normal instances of its
right-indexed critical branchings are confluent. For the regular ones, we have
the following confluence diagrams:
\[
  \begin{array}{c@{\qquad}c}
    \xymatrix{
      {\satex{pearl-nun}}
      \ar@3 @/^6ex/ [rr] ^{C}_{}="1"
      \ar@3 @/_6ex/ [rr] _{D} ^{}="2"
      && {\satex{pearl-cap-id}}\; ,
      \ar@4 "1"!<0ex,-3ex>;"2"!<0ex,3ex> |-*+{\scriptstyle CD}
    }
    &
    \xymatrix{
      {\satex{pearl-unu}}
      \ar@3 @/^6ex/ [rr] ^{D} _{}="1"
      \ar@3 @/_6ex/ [rr] _{C} ^{}="2"
      && {\satex{pearl-cup-id}}\; ,
      \ar@4 "1"!<0ex,-3ex>;"2"!<0ex,3ex> |-*+{\scriptstyle DC}
    }
    \\
    \xymatrix@C=3ex{
      & {\satex{pearl-nou}}
      \ar@3 [r] ^{B} _-{}="source"
      & {\satex{pearl-nuo}}
      \ar@3 [dr] ^{C} \\
      {\satex{pearl-onu}}
      \ar@3 [ur] ^{A}
      \ar@3 [rrr] _{C} ^-{}="target"
      &&& {\satex{pearl2}} \; ,
      \ar@4 "source"!<0ex,-3ex> ; "target"!<0ex,3ex> |-*+{\scriptstyle AC}
    }
    &
    \xymatrix@C=3ex{
      & {\satex{pearl-uon}}
      \ar@3 [r] ^{A} _-{}="source"
      & {\satex{pearl-uno}}
      \ar@3 [dr] ^{D} \\
      {\satex{pearl-oun}}
      \ar@3 [ur] ^{B}
      \ar@3 [rrr] _{D} ^-{}="target"
      &&& {\satex{pearl2}}\,\pbox.
      \ar@4 "source"!<0ex,-3ex>; "target"!<0ex,3ex> |-*+{\scriptstyle BD}
    }
  \end{array}
\]
From the characterization of normal forms of the polygraph given
in~\cite[Section~5.5.2]{guiraud2009higher}, the normal instances of the
right-indexed critical branching $AB\pa{\satex{pearl-f-small}}$ are
the instances corresponding to the following $2$-cells
\begin{align*}
  \satex{pearl-f}&=\satex{pearl-u1n}\;,\;\;
  &
  \satex{pearl-f}&=\satex{pearl-n0i}\;,\;\;
  &
  \satex{pearl-f}&=\satex{pearl-u0i}\;,\;\;
  &
  \satex{pearl-f}&=\satex{pearl-f0n}\!\text,
\end{align*}
where, in the latter, $n\in\N$ and $\satex{pearl-f0-small}$ ranges over the
set $N_0$, the subset of $\freecat P_2$ consisting of normal forms of~$P$ with
degenerate source and target, which are characterized by the following two
construction rules:
\[
  \satex{pearl-f0}
  \qquad
  =
  \qquad
  \satex{empty}
  \qquad\text{or}\qquad
  \satex{pearl-k-f0}
\]
(on the left, this is the empty diagram). Now we check that, for each one of
these $2$-cells, the corresponding critical branching
$AB\pa{\satex{pearl-f-small}}$ is confluent. Let us note that, for
the first three cases, there are several possible confluence diagrams, because
they also contain regular critical branchings of $P$.
\begin{itemize}
\item For $\satex{pearl-f}=\satex{pearl-u1n-small}$, we choose the following one:
\[
  \xymatrix@R=0ex{
    & \satex{pearl-cp5-1}
    \ar@3 [r] ^{D}
    & \satex{pearl-cp5-2}
    \ar@3 [r] ^{B}
    & \satex{pearl-cp5-3}
    \ar@3 [dr] ^{C}
    \\
    \satex{pearl-cp5-0}
    \ar@3 [ur] ^{A} 
    \ar@3 [dr] _{B} 
    &&&& \satex{pearl-cp5-4}\,\pbox.
    \\
    & \satex{pearl-cp5-7}
    \ar@3 [r] _{C}
    & \satex{pearl-cp5-6}
    \ar@3 [r] _{A}
    & \satex{pearl-cp5-5}
    \ar@3 [ur] _{D}
    \ar@4 "1,3"!<0ex,-6ex>;"3,3"!<0ex,6ex> |-*+{\scriptstyle AB\left(\satex{pearl-u1n-small}\right)}
  }
\]
\item For $\satex{pearl-f}=\satex{pearl-n0i}\:$:
\[
  \xymatrix@R=0ex{
    && \satex{pearl-cp6-1}
    \ar@3 [drr] ^{D} 
    \\
    \satex{pearl-cp6-0}
    \ar@3 [urr] ^{A}
    \ar@3 [dr] _{B}
    &&&& \satex{pearl-cp6-2}\,\pbox.
    \\
    & \satex{pearl-cp6-5}
    \ar@3 [r] _{A} 
    & \satex{pearl-cp6-4}
    \ar@3 [r] _{D}
    & \satex{pearl-cp6-3}
    \ar@3 [ur] _{A}
    \ar@4 "1,3"!<0ex,-3ex>;"3,3"!<0ex,3ex> |-*+{\scriptstyle AB\left(\satex{pearl-n0i-small}\right)} 
  }
\]
\item For $\satex{pearl-f}=\satex{pearl-u0i}\:$:
\[
  \xymatrix@R=0ex{
    & \satex{pearl-cp7-1}
    \ar@3 [r] ^{B} 
    & \satex{pearl-cp7-2}
    \ar@3 [r] ^{C}
    & \satex{pearl-cp7-3}
    \ar@3 [dr] ^{B}
    \\
    \satex{pearl-cp7-0}
    \ar@3 [ur] ^{A}
    \ar@3 [drr] _{B}
    &&&& \satex{pearl-cp7-4}\,\pbox.
    \\
    && \satex{pearl-cp7-5}
    \ar@3 [urr] _{C} 
    \ar@4 "1,3";"3,3" |-*+{\scriptstyle AB\left(\satex{pearl-u0i-small}\right)}
  }
\]
\item Finally, for $\satex{pearl-f} = \satex{pearl-f0n}$:
\[
  \xymatrix@R=0ex{
    {\satex{pearl-n-of0o-u}}
    \ar @3 @/^6ex/ [rr] ^-{A} _-{~}="1"
    \ar @3 @/_6ex/ [rr] _-{B} ^-{~}="2"
    && {\satex{pearl-n-1f0oo-u}}\,\pbox.
    \ar@4 "1"!<0ex,-2ex>;"2"!<0ex,2ex> |-{\scriptstyle AB\pa{\satex{pearl-f0n-small}}} 
  }
\]
\end{itemize}
It follows that the $3$-polygraph $P$ is convergent and right-indexed, and the
following $4$-cells form an acyclic extension of $\freegpd{P}$:
\[
CD, \;
DC, \;
AC, \; 
BD, \; 
AB\left(\,\satex{pearl-u1n-small}\,\right), \;
AB\left(\,\satex{pearl-n0i-small}\,\right), \;
AB\left(\,\satex{pearl-u0i-small}\,\right),\;
AB\left(\,\satex{pearl-f0n-small}\,\right),
\]
where $\satex{pearl-f0-small}$ is in $N_0$ and $n$ is in $\Nb$. It can be
observed that the $4$-cells $AB\pa{\,\satex{pearl-u1n-small}\,}$,
$AB\pa{\,\satex{pearl-n0i-small}\,}$ and
$AB\pa{\,\satex{pearl-u0i-small}\,}$ are superfluous. Namely, the
$3$-spheres forming their boundaries are also the boundaries of $4$-cells of
$\freegpd{Q}$ where $Q$ is the $(4,2)$\nbd-poly\-graph obtained from the
$3$-polygraph~$P$ by adding the $4$-generators $AC$ and $BD$.

Let us denote by $X_0$ the family made of the $4$-cells $CD$,
$DC$, $AC$, and $BD$. Then, for every natural number
$n$, one defines:
\[
X_{n+1} = X_n \sqcup \setof{ AB\pa{\satex{pearl-f0n-small}}}{\satex{pearl-f0-small} \in N_0}\pbox.
\]
Thus, the following set of $4$-cells forms an acyclic extension of the $(3,2)$\nbd-cate\-gory $\freegpd{P}$:
\[
  X = \bigcup_{n\in\Nb} X_n
  \pbox.
\]
It can be shown that this infinite number of confluence
diagrams cannot be filled by a finite cellular extension and thus that the
$3$-polygraph $P$ does not have finite derivation
type~\cite[Theorem~5.5.7]{guiraud2009higher}:

\begin{theorem}
  \label{T:Pearl3NotFDT}
  The above 3-polygraph $P$ does not have finite derivation type.
\end{theorem}

\noindent
We will see in \cref{sec:pol-pearls} a generalization of this result to
$n$-polygraphs with $n\geq 3$.

\section{(3,2)-PROs}
\label{sec:32PRO}
We generalize, in dimension $3$, the notion of PRO introduced in
\cref{sec:PRO}, as well as introduce symmetric and coherent variants. This will
be used in subsequent sections to show coherence theorems such as Mac Lane's
coherence theorem for monoidal categories.

\subsection{3-PROs}
\index{3-monoid@$3$-monoid}
\index{monoid!3-@$3$-}
\index{3-PRO@$3$-PRO}
\index{PRO!3-@$3$-}
A \emph{3-monoid} is a $3$-category~$C$ where there is exactly one
$0$\nbd-cell~$\star$. In such a $3$-category, the set $C_1$ is canonically a
monoid when equipped with $0$-composition as multiplication and $\unit{\star}$
as unit.
A \emph{3-PRO} is a $3$-monoid whose monoid of $1$-cells is the additive
monoid~$\N$. Note that the underlying $2$-category of a $3$-PRO is always a PRO,
as defined in \cref{sec:PRO}, thus the name.
A $3$-PRO which is also a $(3,2)$-category is called a \emph{$(3,2)$-PRO}: this
is a $3$-PRO in which every $3$-cell is invertible with respect to
composition~$\comp2$.

\subsection{3-PROPs}
\label{sec:3PROP}
A \emph{symmetry} on a $3$-monoid~$C$ is an invertible transformation
\[
  \gamma_{a,b}:a\comp0 b\to b\comp0 a,
\]
indexed by $1$-cells $a,b\in C_1$, which is natural in both components and makes
the following diagrams commute:
\begin{align*}
  \svxym{
    a\ar@{=}[dr]\comp0\unit\star\ar[rr]^{\gamma_{a,\unit\star}}&&\ar@{=}[dl]\unit\star\comp0a\\
    &a&
  }
  &&
  \svxym{
    &c\comp0a\comp0 b\ar[dr]^{c\comp0\gamma_{a,b}}&\\
    a\comp0 b\comp0 c\ar[ur]^{\gamma_{a\comp0 b,c}}\ar[dr]_{\gamma_{a,b\comp0 c}}&&c\comp0 b\comp0 a\pbox.\\
    &b\comp0 c\comp0 a\ar[ur]_{\gamma_{b,c}\comp0 a}&
  }  
\end{align*}
A \emph{3-PROP}\index{3-PROP@$3$-PROP}\index{PROP!3-@$3$-} (\resp \emph{$(3,2)$-PROP}) is a 3-PRO (\resp $(3,2)$-PRO) equipped with a symmetry.

\newcommand{\PP}{C} 

\subsection{Algebras over 3-PRO(P)s}
\index{algebra!of a $3$-PRO}
\label{algebrasOverPro(p)}
The $2$-category~$\Cat$ of categories, functors, and natural transformations
is monoidal when equipped with the cartesian product as tensor product. By
Mac Lane's coherence theorem, it can be considered as a strict monoidal
category or, equivalently, as a $3$-PRO with
categories as $1$-cells, functors as $2$-cells, natural transformations as
$3$-cells, cartesian product as $0$-composition, composition of functors as
$1$-composition, vertical composition of natural transformations as
$2$-composition.

If $\PP$ is a $3$\nbd-PRO, a \emph{$\PP$-algebra} is a $3$-functor from $\PP$ to
$\Cat$. If $\PP$ is a $3$\nbd-PROP, we moreover require that this $3$-functor
preserves the symmetry. Given a $\PP$\nbd-alge\-bra~$A$, we often write $A_\star$
instead $A(\star)$.
If $A$ and $B$ are $\PP$-algebras, a \emph{morphism of $\PP$-algebras} from $A$ to
$B$ is a natural transformation from $A$ to $B$, \ie a pair $(F,\phi)$ where
$F:A_\star \fl B_\star$ is a functor and $\phi$ is a map sending every $2$-cell
$f:m\dfl n$ in $\PP$ to a natural isomorphism with the following shape:
\[
\xymatrix@R=1em{
  & {B_\star^m}
  \ar@/^/ [dr] ^-{B(f)}
  \ar@2 []!<0pt,-15pt>;[dd]!<0pt,15pt> _-{\phi_f} 
  \\
  {A_\star^m} 
  \ar@/^/ [ur] ^-{F^m}
  \ar@/_/ [dr] _-{A(f)}
  && {B_\star^n}
  \\
  & {A_\star^n}
  \ar@/_/ [ur] _-{F^n}
}
\]
such that the following relations hold:
\begin{itemize}
\item for every $2$-cells $f:m\dfl n$ and $g:p\dfl q$ of $\PP$, we have
  \[
    \phi_{f\comp0 g} = \phi_f\times \phi_g,
  \]
  \ie graphically,
  \[
    \vcenter{
      \xymatrix@R=1em@C=2ex{
        & {B_\star^{m+p}}
        \ar@/^/ [dr] ^-{B(f\comp0 g)}
        \ar@2 []!<0pt,-10pt>;[dd]!<0pt,10pt> |-{\phi_{f\comp0 g}} 
        &
        \\
        {A_\star^{m+p}} 
        \ar@/^/ [ur] ^-{F^{m+p}}
        \ar@/_/ [dr] _-{A(f\comp0 g)}
        && {B_\star^{n+q}}
        \\
        & {A_\star^{n+q}}
        \ar@/_/ [ur] _-{F^{n+q}}
        &
      }
    }
    =
    \vcenter{
      \xymatrix@R=1em@C=2ex{
        &   {B_\star^m\times B_\star^p}
        \ar@/^/ [dr] ^-{B(f)\times B(g)}
        \ar@2 []!<0pt,-10pt>;[dd]!<0pt,10pt> |-{\phi_{f}\times\phi_{g}}
        \\
        {A_\star^m\times A_\star^p}
        \ar@/^/ [ur] ^-{F^{m}\times F^{p}}
        \ar@/_/ [dr] _-{A(f)\times A(g)}
        && {B_\star^{n}\times B_\star^{q}}\,,
        \\
        & {A_\star^{n}\times A_\star^{q}}
        \ar@/_/ [ur] _-{F^{n}\times F^{q}}
      }
    }
  \]
\item for every $2$-cells $f:m\dfl n$ and $g:n\dfl p$ in $\PP$, we have
  \[
    \phi_{f\comp1 g} = (\phi_f\comp1 B(g)) \comp2 (A(f)\comp1 \phi_g),
  \]
  \ie graphically,
  \[
    \vcenter{
      \xymatrix@R=1em{
        & {B_\star^{m}}
        \ar@/^/ [dr] ^-{B(f\comp1 g)}
        \ar@2 []!<0pt,-10pt>;[dd]!<0pt,10pt> |-{\phi_{f\comp1 g}}  
        \\
        {A_\star^{m}} 
        \ar@/^/ [ur] ^-{F^{m}}
        \ar@/_/ [dr] _-{A(f\comp1 g)}
        && {B_\star^{p}}
        \\
        & {A_\star^{p}}
        \ar@/_/ [ur] _-{F^{p}}
      }
    }
    =
    \vcenter{
      \xymatrix@R=1em{
        & {B_\star^m}
	\ar@/^/ [dr] ^-{B(f)}
	\ar@2 []!<0pt,-10pt>;[dd]!<0pt,10pt> _-{\phi_{f}}  
        \\
        {A_\star^m} 
	\ar@/^/ [ur] ^-{F^{m}}
	\ar@/_/ [dr] _-{A(f)}
        && {B_\star^n}
	\ar@/^/ [dr] ^-{B(g)}
	\ar@2 []!<0pt,-10pt>;[dd]!<0pt,10pt> ^-{\phi_{g}}  
        \\ 
        & {A_\star^n}
	\ar@/_/ [dr] _-{A(g)}
	\ar [ur] |-{F^n}
        && {B_\star^p}\,,
        \\
        && {A_\star^p}
	\ar@/_/ [ur] _-{F^p}
      }
    }
  \]
\item for every $3$-cell $\alpha:f\tfl g:m\dfl n$ in $\PP$, we have
  \[
    \phi_f\comp2 (A(\alpha) \comp1 F^n) = (F^m \comp1 B(\alpha)) \comp2 \phi_g,
  \]
  \ie graphically,
  \[
    \vcenter{
      \xymatrix@R=1em{
        & {B_\star^m}
        \ar@/^/ [dr] ^-{B(f)}
        \ar@2 []!<0pt,-10pt>;[dd]!<0pt,10pt> ^-{\phi_f}  
        \\
        {A_\star^m} 
        \ar@/^/ [ur] ^-{F^m}
        \ar [dr] ^(0.4){A(f)} _-{}="src"
        \ar@/_7ex/ [dr] _-{A(g)} ^-{}="tgt"
        \ar@2 "src";"tgt" |-{A(\alpha)}
        && {B_\star^n}
        \\
        & {A_\star^n}
        \ar@/_/ [ur] _-{F^n}
      }
    }
    =
    \vcenter{
      \xymatrix@R=1em{
        & {B_\star^m}
        \ar@/^7ex/ [dr] ^-{B(f)} _-{}="src"
	\ar [dr] _(0.6){B(g)} ^-{}="tgt"
	\ar@2 []!<0pt,-10pt>;[dd]!<0pt,10pt> _-{\phi_g} 
	\ar@2 "src";"tgt" |-{B(\alpha)}
        \\
        {A_\star^m} 
	\ar@/^/ [ur] ^-{F^m}
	\ar@/_/ [dr] _-{A(g)}
        && {B_\star^n}\,.
        \\
        & {A_\star^n}
	\ar@/_/ [ur] _-{F^n}
      }
    }
  \]
\end{itemize}
The $\PP$-algebras and their morphisms form a category denoted by $\Alg(\PP)$.

\subsection{The coherence problem for algebras over a \pdfm{3}-PRO(P)}
\label{coherenceProblemPro(p)}
Let $\PP$ be a $3$-PRO(P) and $A$ be a $\PP$-algebra. A \emph{$\PP$-diagram}
in $A$ is the image~$A(\Delta)$ of a $3$-sphere $\Delta$ in $\PP$, \ie a pair
$(\alpha,\beta)$ of $3$-cells with the same source, and with the same target,
where $\alpha$ (\resp $\beta$) is the \emph{source} (\resp \emph{target}) of the
$3$-sphere and is denoted by $\src{}(\Delta)$ (\resp $\tgt{}(\Delta)$). A
$\PP$-diagram $A(\Delta)$ in $A$ \emph{commutes} if the relation
\[
  A(\src{}(\Delta))= A(\tgt{}(\Delta))
\]
is satisfied in $\Cat$.

\index{coherence!for algebras over a $3$-PRO}
The \emph{coherence problem} for algebras over a $3$-PRO(P) is the following
question:
\begin{center}
\emph{Given a $3$-PRO(P) $\PP$, does every $\PP$-diagram commute in every $\PP$-algebra?}  
\end{center}

\index{aspherical!$3$-PRO}
A $3$-PRO is \emph{aspherical} when there is at most one $3$-cell between two
given $2$\nbd-cells, \ie any two parallel $3$-cells are equal.
As a consequence of this definition, we have the following sufficient condition
for giving a positive answer to the coherence problem:

\begin{proposition}
  \label{Propaspherical=>coherence}
  If $\PP$ is an aspherical 3-PRO(P), then every $\PP$-diagram commutes in
  every $\PP$-algebra.
\end{proposition}

\subsection{Presentations of (3,2)-PROs}
\index{presentation!of a 32-PRO@of a $(3,2)$-PRO}
A \emph{presentation} of a $(3,2)$-PRO $\PP$ is a $(4,2)$-polygraph~$P$
such that $\PP \isoto \freegpd{\tpol 3P}/P_4$,
\ie $\PP$ is isomorphic to the $(3,2)$\nbd-cate\-gory generated by the underlying
$3$\nbd-poly\-graph~$\tpol 3P$, quotiented by the congruence generated by the
$4$-generators.
By definition of a $(3,2)$-PRO, in the case where we have a presentation as
above, the $3$-polygraph $P$ necessarily has exactly one $0$-cell and one
$1$\nbd-cell. A presentation~$P$ of $\PP$ is called \emph{coherently
  convergent} rather than \emph{convergent} when $P$
is a convergent $3$\nbd-poly\-graph and $P_4$ is a cellular extension of
generating confluences of~$P$, see~\cite[Section~2.1.1]{GuiraudMalbos12mscs}.

\newcommand{\PRO}[1]{\mathrm{#1}}

\begin{example}
  \index{polygraph!of monoids}
  \label{ex:as-cat}
  Consider the $(4,2)$-polygraph
  \[
    P=
    \PREs{\star}{a}{\mu:aa\To a}{A:(\mu\comp0a)\comp1\mu\TO(a\comp0\mu)\comp1\mu}{\Gamma}.
  \]
  The $3$-generator $\mu$ is often pictured as $\satex{mu-vsmall}$ and $A$ as
  \[
    \satex{alpha2}:\satex{mon-assoc-l}\TO\satex{mon-assoc-r}\pbox.
  \]
  Similarly, the $4$-generator~$\Gamma$ is pictured as $\satex{aleph-vsmall}$
  and its boundary is given by
  \[
    \vcenter{
      \xymatrix @R=0em @C=1em {
        &{ \satex{mon-cp1-small-1} }
        \ar@3 [rr] 
        ^-{\satex{mon-cp1-small-12}} 
        _-{}="1" 
        && { \satex{mon-cp1-small-2} }
        \ar@3@/^/ [dr] ^-{\satex{mon-cp1-small-23} }
        \\
        { \satex{mon-cp1-small-0} }
        \ar@3@/^/ [ur] ^-{\satex{mon-cp1-small-01}}
        \ar@3@/_/ [drr] _-*+{\satex{mon-cp1-small-04}}
        &&&&\satex{mon-cp1-small-3}\pbox.
        \\
        && { \satex{mon-cp1-small-4} }
        \ar@3@/_/ [urr] _-*+{\satex{mon-cp1-small-43}}
        \ar@4 "1"!<0ex,-4.5ex>;[]!<0ex,5.5ex> ^-*+{\satex{aleph}}
      }
    }
  \]

  Consider the category~$\Simplsurj$ whose objects are natural numbers and
  morphisms from~$m$ to~$n$ are surjective increasing functions $[m]\to[n]$
  where $[n]$ denotes the set $\set{0,\ldots,n-1}$. This is a monoidal
  subcategory of the augmented simplicial category $\Simplaug$
  (see \cref{sec:simpl-cat,sec:ass-simpl-pres}), already encountered in
  \cref{ex:simpl-fs}. As a variant of the presentation of~$\Simplaug$, see
  \cref{sec:ass-simpl-pres}, the underlying $3$\nbd-poly\-graph of~$P$ can be
  shown to present the monoidal category~$\Simplsurj$. Moreover, the rewriting
  system is convergent and the boundary of the $4$-generator~$\Gamma$ shown above is a confluence diagram for the only critical branching of the rewriting system so that,
  by~\cref{Theorem:SquierCompletion3polygraphs}, $P_4$ forms an acyclic extension of the $(3,2)$-category $\freegpd{\tpol 3P}$, \ie $P$ is a coherent presentation of $\Simplsurj$.

  We write $\PRO{AsCat}$ for the $(3,2)$-PRO presented by~$P$, \ie 
  $\PRO{AsCat} = \freegpd{\tpol 3P}/P_4$.
  The category of its algebras $\Alg(\PRO{AsCat})$ is isomorphic to the category
  of associative categories: we recall that an \emph{associative category} is a
  category~$C$ equipped with a bifunctor $\otimes:C\times C\to C$ and a natural
  transformation $\alpha_{a,b,c}:(a\otimes b)\otimes c\to a\otimes(b\otimes c)$
  satisfying the usual coherence law (see \cref{sec:mon-cat}). Namely, the
  correspondence between an associative category $(C,\otimes,\alpha)$ and a
  $3$-functor $A:\PRO{AsCat}\to\Cat$ is given by
  \begin{align*}
    A\pa{\:\satex{1}\:} &= C,
    &
    A\pa{\satex{mu}} &= \otimes,
    &
    A\pa{\satex{alpha}} &= \alpha
    \pbox.
  \end{align*}
  This correspondence is well-defined since the coherence diagram satisfied by
  associative categories corresponds to the
  $4$-cell~$\satex{aleph-vsmall}$.

  Since $P_4=\set{\satex{aleph-vsmall}}$ is an acyclic extension of  $\freegpd{\tpol3\polyAs}$, we have that $\PRO{AsCat}$ is an aspherical
  $(3,2)$-PRO. As a consequence, in every associative category $C$, every
  $\PRO{AsCat}$-diagram is commutative. This fact can be informally restated as:
  every diagram built in $C$ from the functor $\otimes$ and the natural
  transformation~$\alpha$ is commutative.
\end{example}

\subsection{Coherence in algebras over (3,2)-PROs}
By definition, a $3$-PRO $\PP$ is aspherical if, for every
presentation $P$ of $\PP$, the cellular extension $P_4$ of $\freegpd{\tpol3P}$
is acyclic. The latter condition is satisfied by any convergent
presentation of~$\PP$ yielding the following sufficient condition for giving a
positive answer to the coherence problem for
$\PP$-algebras~\cite[Theorem~2.1.2]{GuiraudMalbos12mscs}:

\begin{theorem}
  \label{TheoremProAspherical}
  If a $(3,2)$-PRO $\PP$ admits a convergent presentation then every
  $\PP$-diagram commutes in every $\PP$-algebra.
\end{theorem}

\section{Coherence in Monoidal Categories}
\label{sec:mon-coh}
The coherence problems in monoidal categories can be formulated in terms of
asphericity problems for $(3,2)$-categories. This section briefly reviews this
approach in the case of monoidal categories. Symmetric and
braided monoidal categories are handled in the next section.

\subsection{Monoidal categories}
\index{monoidal category}
\label{sec:mon-cat}
A \emph{monoidal category} is a category $C$, equipped with two functors
\begin{align*}
  \otimes&:C\times C\fl C,
  &
  e&:1\fl C,  
\end{align*}
and three natural isomorphisms
\begin{align*}
  \alpha_{x,y,z} &: (x\otimes y)\otimes z \fl x\otimes(y\otimes z),
  &
  \lambda_x &: e \otimes x \fl x,
  &
  \rho_x &: x\otimes e \fl x,
\end{align*}
such that the following two diagrams commute in $C$:
\[
  \xymatrix@C=-3.5ex@R=3ex{
    &&(x\otimes y)\otimes(z\otimes t)\ar@/^/[drr]^{\alpha_{x,y,z\otimes t}}\\
    ((x\otimes y)\otimes z)\otimes t\ar@/^/[urr]^{\alpha_{x\otimes y,z,t}}\ar[dr]_{\alpha_{x,y,z}\otimes t}&&&&x\otimes(y\otimes(z\otimes t))\\
    &(x\otimes(y\otimes z))\otimes t\ar[rr]_{\alpha_{x,y\otimes z,t}}&&x\otimes((y\otimes z)\otimes t)\ar[ur]_{x\otimes\alpha_{y,z,t}}\\
    &(x\otimes e)\otimes y\ar[dr]_{\alpha_{x,e,y}}\ar[rr]^{\rho_x\otimes y}&&x\otimes(e\otimes y)\\
    &&x\otimes(e\otimes y).\ar[ur]_{x\otimes\lambda_y}
  }
\]

\subsection{The \pdfm{3}-PRO of monoidal categories}
\label{2ProMonCat}
\label{sec:moncat-coh}
Consider the $(4,2)$-polygraph
\[
  P=\PREs{\star}{a}{\mu:aa\To a,\eta:\unit\star\to a}{A,L,R}{\Gamma,\Delta},
\]
whose $2$-generators $\mu$ and $\eta$ are respectively pictured as
$\satex{mu-vsmall}$ and $\satex{eta-vsmall}$, whose $3$\nbd-gene\-rators
$A$, $L$, and $R$ are respectively pictured as
\begin{align*}
  \satex{A2}&:\satex{mon-assoc-l}\TO\satex{mon-assoc-r},
  &
  \satex{L2}&:\satex{mon-unit-l}\TO\satex{mon-unit-c}\;,
  &
  \satex{R2}&:\satex{mon-unit-r}\TO\satex{mon-unit-c}\;,
\end{align*}
whose $4$-generators~$\Gamma$ and $\Delta$ are respectively
\[
  \vcenter{
    \xymatrix @R=1em @C=1em {
      && \satex{mon-cp1s-4}
      \ar@{{}{ }{}}[]_{}="2"
      \ar@3@/^/ [drr] ^-*+{\satex{mon-cp1s-43}}
      \\
      \satex{mon-cp1s-0}
      \ar@3@/_/ [dr] _-{\satex{mon-cp1s-01}}
      \ar@3@/^/ [urr] ^-*+{\satex{mon-cp1s-04}}
      &&&&
      \satex{mon-cp1s-3}
      \\
      &\satex{mon-cp1s-1}
      \ar@3 [rr] 
      _-{\satex{mon-cp1s-12}} 
      ^-{}="1" 
      &
      &
      \satex{mon-cp1s-2}
      \ar@3@/_/ [ur] _-{\satex{mon-cp1s-23}}
      \ar@4 "2"!<0ex,-7ex>;"1"!<0ex,5ex> ^-*+{\satex{aleph}}
    }
  }
\]
and
\[
  \vcenter{
    \xymatrix @R=1em @C=1em {
      &\satex{mon-cp3s-0}
      \ar@3@/_/ [dr] _-{\satex{mon-cp3s-01}}
      \ar@3[rr] ^-{\satex{mon-cp3s-02}} _{}="3"
      &&
      \satex{mon-cp3s-2}
      \,.
      \\
      &&\satex{mon-cp3s-1}
      \ar@3@/_/ [ur] _-{\satex{mon-cp3s-12}}
      \ar@{{}{ }{}}[]^{}="4"
      \ar@4 "3"!<1.4pt,-2ex>;"4"!<0ex,4ex> ^-*+{\satex{beth}}
    }
  }
\]
We denote by $\PRO{MonCat}$ the $(3,2)$-PRO presented by this polygraph.
It is easily seen that the category of small monoidal categories and monoidal
functors is isomorphic to the category
$\Alg(\PRO{MonCat})$~\cite[Lemma~2.3.2]{GuiraudMalbos12mscs}.

Note that the underlying $3$-polygraph~$\tpol3P$ of~$P$ is the polygraph of
monoids defined in \cref{ex:ass}, which was shown to be terminating in
\cref{ex:ass-term}. Its five critical branchings are computed in
\cref{ex:ass-cp} and shown to be confluent.
Consider the cellular extension~$X$ of $\freegpd{\tpol3 P}$ with five $4$-cells:
the $4$\nbd-gene\-rators $\Gamma$ and $\Delta$, as well as
\begin{align*}
  \vcenter{
    \xymatrix @R=1em @C=1em {
      & {\satex{mon-cp2s-1}}
      \ar@3@/^/ [dr] ^-{\satex{mon-cp2s-12}}
      \\
      {\satex{mon-cp2s-0}}
      \ar@3@/^/ [ur] ^-{\satex{mon-cp2s-01}}
      \ar@3@/_/ [rr] _-{\satex{mon-cp2s-02}} ^-{}="1"
      && {\satex{mon-cp2s-2}}\,\pbox\;,
      \ar@4 "1,2"!<-5pt,-20pt>;"1"!<-5pt,10pt> ^-*+{\Lambda_1}
    }
  }
  &&
  \vcenter{
    \xymatrix @R=1em @C=1em {
      & {\satex{mon-cp4s-1}}
      \ar@3@/^/ [dr] ^-{\satex{mon-cp4s-12}}
      \\
      {\satex{mon-cp4s-0}}
      \ar@3@/^/ [ur] ^-{\satex{mon-cp4s-01}}
      \ar@3@/_/ [rr] _-{\satex{mon-cp4s-02}} ^-{}="1"
      && {\satex{mon-cp4s-2}}\,\pbox\;,
      \ar@4 "1,2"!<-5pt,-20pt>;"1"!<-5pt,10pt> ^-*+{\Lambda_2}
    }
  }
  &&
  \vcenter{
    \xymatrix @C=1.6em {
      {\satex{mon-cp5s-0}}
      \ar@3@/^4.5ex/ [rr] ^{\satex{mon-cp5s-01}} _{}="1"
      \ar@3@/_4.5ex/ [rr] _{\satex{mon-cp5s-10}} ^{}="2"
      && {\satex{mon-cp5s-1}}\,\pbox\;.
      \ar@4 "1"!<-5pt,-13.5pt>;"2"!<-5pt,13.5pt> ^-*+{\Lambda_3}
    }
  }  
\end{align*}
By \cref{Theorem:SquierCompletion3polygraphs}, $X$ forms an acyclic extension of
$\freegpd{\tpol3P}$ since its elements are a choice of confluence diagrams for
the five critical branchings. It can be shown that $\Lambda_1$, $\Lambda_2$, and
$\Lambda_3$ are superfluous in this cellular extension, \ie for each $4$\nbd-cell
$\Lambda_i$, we have $\cl{\src{}(\Lambda_i)}= \cl{\tgt{}(\Lambda_i)}$ in
$\PRO{MonCat}$~\cite[Section~2.3.3]{GuiraudMalbos12mscs}. Therefore
$\set{\Gamma,\Delta}$ is still an acyclic extension, \ie the polygraph~$P$ is
coherent.
We have thus proved~\cite[Theorem~5.2.2]{guiraud2009higher},
\cite[Proposition~2.3.3]{GuiraudMalbos12mscs}:

\begin{proposition}
  \label{Mon_4HomotopyBase}
  The above $(4,2)$-polygraph~$P$ is coherent.
\end{proposition}

Mac Lane's coherence theorem~\cite[Theorem~VII.2.1]{MacLane98} states that, in a
monoidal category, every diagram whose arrows are built up from instances of
$\otimes$, $\alpha$, $\lambda$, and $\rho$ commute. From
Proposition~\ref{Mon_4HomotopyBase}, we can deduce this theorem, which can be
reformulated as follows:

\begin{theorem}
  \label{thm:moncat-coh}
  The 3-PRO $\PRO{MonCat}$ is aspherical.
\end{theorem}

\section{Coherence in Symmetric and Braided Monoidal Categories}
\label{Section:CoherenceSym}

\subsection{Symmetric monoidal categories}
\index{symmetric!monoidal category}
\index{monoidal category!symmetric}
A \emph{symmetric monoidal category} is a monoidal category
$(C,\otimes,e,\alpha,\lambda,\rho)$ equipped with a natural isomorphism
\[
  \gamma_{x,y} : x\otimes y \longrightarrow y\otimes x,
\]
called the \emph{symmetry} and such the following two diagrams commute
in~$C$:
\[
  \vcenter{
    \xymatrix@C=1em{
      & y\otimes x
      \ar[dr] ^-{\gamma_{y,x}}
      \\
      x\otimes y 
      \ar[ur] ^-{\gamma_{x,y}}
      \ar@{=}[rr] ^-{}="1"
      && x\otimes y,
    }
  }
\]
\[
  \vcenter{
    \xymatrix{
      & x\otimes (y\otimes z) 
      \ar[rr] ^-{\gamma_{x,y\otimes z}}
      && (y\otimes z)\otimes x
      \ar[dr] ^-{\alpha}
      \\
      (x\otimes y) \otimes z 
      \ar[ur] ^-{\alpha_{x,y,z}}
      \ar[dr] _-{\gamma_{x,y}\otimes z}
      &&&& y \otimes (z\otimes x)\pbox.
      \\
      & (y\otimes x)\otimes z 
      \ar[rr] _-{\alpha_{y,x,z}}
      && y\otimes(x\otimes z)
      \ar[ur] _-{y\otimes\gamma_{x,z}}
    }
  }
\]

\subsection{PROPs}
\index{PROP}
Recall from \cref{sec:PRO} that a PRO is a strict monoidal category whose monoid
of objects is~$(\N,+,0)$. We now introduce the following symmetric variant. A
\emph{PROP} is a strict symmetric monoidal category whose monoid of objects
is~$(\N,+,0)$. In the following, we consider PROs and PROPs as $2$-categories
with one $0$-cell. In particular, the underlying $2$-category of a $3$-PROP, as
defined in \cref{sec:3PROP}, is a PROP.

\subsection{PROPs as PROs}
\index{symmetry}
\index{Yang-Baxter relation}
\label{Subsection:Presentations2Props}
PROPs can be characterized among PROs as follows, see \cref{sec:free-sym} and~\cite[Proposition A.3 and
Corollary A4]{Guiraud04}.
A PRO $\PP$ is a PROP if and only if it contains a $2$-cell $\gamma:2\dfl 2$,
represented by~$\satex{gamma-vsmall}$, such that the following relations
hold:
\begin{itemize}
\item \emph{involutivity} of the symmetry
  \[
    \gamma\comp1\gamma = \id_2,
  \]
  which can be pictured as
\[
\satex{sym-l} = \satex{sym-r}\;,
\]
\item the \emph{Yang-Baxter relation}
  \[
    (\gamma\comp0 1)\comp1(1\comp0\gamma)\comp1(\gamma\comp0 1) = (1\comp0\gamma)\comp1(\gamma\comp0 1)\comp1(1\comp0\gamma),
  \]
  which can be pictured as
\[
\satex{yb-l} = \satex{yb-r}\;,
\]
\item for every $2$-cell $\phi:m\dfl n$ of $\PP$, the \emph{left and right
    naturality relations} for~$\phi$
  \begin{align*}
    (\phi \comp0 1) \comp1 \gamma_{n,1} &= \gamma_{m,1} \comp1 (1\comp0 \phi),
    \\
    (1\comp0 \phi) \comp1 \gamma_{1,n} &= \gamma_{1,m} \comp1 (\phi\comp0 1),
  \end{align*}
  with the inductively defined notations:
  \begin{align*}
    \gamma_{0,1} &= \gamma_{1,0} = \id_1,
    &
    \gamma_{n+1,1} &= (n \comp0 \gamma) \comp1 (\gamma_{n,1}\comp0 1),
    \\
    &&
    \gamma_{1,n+1} &= (\gamma \comp0 n) \comp1 (1\comp0\gamma_{1,n})
    \pbox.
  \end{align*}
\end{itemize}
If we represent $\phi$ by $\satex{coh-phi-small}$\,, any $\gamma_{n,1}$
by~$\satex{coh-gamma1-small}$\,, and any~$\gamma_{1,n}$ by~$\satex{coh-gamma2-small}$\,, the
naturality relations for $\phi$ are
\begin{align*}
  \satex{coh-gamma-nat1-l}&=\satex{coh-gamma-nat1-r}\;,
  &
  \satex{coh-gamma-nat2-l}&=\satex{coh-gamma-nat2-r}\;.
\end{align*}

\subsection{The PROP of permutations}
\index{polygraph!of permutations}
The initial PROP is the PROP of \emph{permutations}, denoted by $\Sym$ and introduced in \cref{sec:sym-pres}, whose
$2$-cells from $n$ to $n$ are the permutations of $\set{0,\dots,n-1}$ and with
no $2$-cell from $m$ to $n$ if $m\neq n$. The $2$-PROP $\Sym$ is presented by
the $3$-polygraph $P$ of permutations defined
in~\cref{sec:sym-pres} whose $3$-cells correspond to the
involutivity and Yang-Baxter relations:
\begin{align}
  \label{eq:sym-yb}
  \satex{sym-l}&\tfl \satex{sym-r}\;,
  &
  \satex{yb-l}&\tfl \satex{yb-r}\;.
\end{align}
There is an isomorphism between the category of small categories and functors
and the category $\Alg(\Perm)$.

\subsection{Presentations of PROPs}
\index{presentation!of a PROP}
\index{free!PROP}
\label{sec:free-PROP}
Let $P$ be a $2$-polygraph with one $0$-cell and one $1$-cell. We denote by~$SP$
the $3$-polygraph obtained from $P$ by adjoining a $2$-cell
$\satex{gamma-vsmall}:2\dfl 2$ and the following $3$-cells:
\begin{itemize}
\item the symmetry $3$-cell and the Yang-Baxter $3$-cell~\eqref{eq:sym-yb},
\item two $3$-cells for every $2$-cell $\phi=\satex{phi-small}$ of $P$,
  corresponding to the naturality relations for $f$:
  \begin{align*}
    \satex{coh-gamma-nat1-l}&\TO\satex{coh-gamma-nat1-r}\;,
    &
    \satex{coh-gamma-nat2-l}&\TO\satex{coh-gamma-nat2-r}\; . 
  \end{align*}
\end{itemize}
The \emph{free PROP generated by $P$} is the $2$-category, denoted by $P^S$,
presented by the $3$-polygraph $SP$, see also \cref{sec:pres-sym}.  We define a \emph{presentation of a PROP
  $\PP$} as a pair $(P,P_3)$ made of a $2$-polygraph $P$ with one $0$-cell and
one $1$-cell and a cellular extension $P_3$ of the free $2$-PROP $P^S$, such
that $\PP \simeq P_2^S/P_3$.

\subsection{Presentations of (3,2)-PROPs}
\index{presentation!of a 32-PROP@of a $(3,2)$-PROP}
Let $(P,P_3)$ be a presentation of a PROP. We denote by~$Q$ the
$(4,2)$-polygraph obtained from the $3$-polygraph $SP$ by adjoining the
$3$-cells of $P_3$ and a cellular extension made of the following two $4$-cells
for each $3$-generator $A:\phi\TO\psi$ in $P_3$, corresponding to the naturality
relations for~$A$:
\begin{align*}
  \xymatrix@R=1em@C=1em{
    & \satex{coh-gamma-3nat1-2}
    \ar@3 [dr]
    \ar@4 []!<0pt,-30pt>;[dd]!<0pt,30pt> 
    \\
    \satex{coh-gamma-3nat1-1}
    \ar@3 [ur] ^-{A}
    \ar@3 [dr]
    && \satex{coh-gamma-3nat1-3}
    \\
    & \satex{coh-gamma-3nat1-4}
    \ar@3 [ur] _-{A}
  }
  &\;\;&
  \xymatrix@R=1em@C=1em{
    & \satex{coh-gamma-3nat2-2}
    \ar@3 [dr]
    \ar@4 []!<0pt,-30pt>;[dd]!<0pt,30pt> 
    \\
    \satex{coh-gamma-3nat2-1}
    \ar@3 [ur] ^-{A}
    \ar@3 [dr]
    && \satex{coh-gamma-3nat2-3}\,.
    \\
    & \satex{coh-gamma-3nat2-4}
    \ar@3 [ur] _-{A}
  }
\end{align*}
The \emph{free $(3,2)$-PROP generated by $P$} is the $(3,2)$-category, denoted
by $P^S$, presented by the $(4,2)$-polygraph~$Q$ defined above: $P^S = \freegpd{\tpol3Q}/Q_4$.

We define a \emph{presentation of a $(3,2)$-PROP $\PP$} as a pair $(P,P_4)$,
where $P$ is a presentation of a PROP and~$P_4$ is a cellular extension of the
free $(3,2)$-PROP $P^S$ generated by $P$, such that
$\PP \isoto P^S/P_4$. A presentation~$P$ of a $(3,2)$-PROP is called \emph{convergent} when the
$3$-polygraph $SP$ is convergent.

\subsection{Application to symmetric monoidal categories}
\index{symmetric!monoidal category}
\index{monoidal category!symmetric}
\label{Subsection:ApplicationCoherenceSym}
\label{2PROPSymCat}

Let $\Sym\Cat$ be the $(3,2)$-PROP presented by the polygraph~$P$ given as
follows.
\begin{itemize}
\item $P_0=\set{\star}$, $P_1=\set{a}$,
\item $P_2$ is the $2$-polygraph, containing two
  $2$-cells~$\satex{mu-vsmall}$ and~$\satex{eta-vsmall}$\;,
\item $P_3$ is the cellular extension of the free $2$-PROP $P_2^S$ generated by
  $P_2$ containing the three $3$-cells
  \begin{align*}
    \xymatrix@C=8ex{
      \satex{mon-assoc-l}
      \ar@3[r] ^-{\satex{A}}
      &
      \satex{mon-assoc-r}\;,
    }
    &&
    \xymatrix@C=5ex{
      \satex{mon-unit-l}
      \ar@3[r] ^-{\satex{L}}
      &
      \satex{mon-unit-c}\;,
    }  
    &&
    \xymatrix@C=5ex{
      \satex{mon-unit-r}
      \ar@3[r] ^-{\satex{R}}
      &
      \satex{mon-unit-c}\;,
    }
  \end{align*}
  plus the following extra $3$-cell:
  \[
    \xymatrix{
      \satex{mon-com-l}
      \ar@3[r] ^-{\satex{C}}
      &
      \satex{mon-com-r}\;,
    }
  \]
\item $P_4$ is the cellular extension of the free $(3,2)$-PROP $P_3^S$ generated
  by $P_3$ containing the two $4$-cells
  \[
    \satex{sym-4cell1}
    \quad
    \satex{sym-4cell2}
  \]
  plus the following two extra $4$-cells:
  \[
    \satex{sym-4cell4}
    \quad
    \satex{sym-4cell3}
  \]
\end{itemize}
The category of small symmetric monoidal categories and symmetric monoidal
functors is isomorphic to the category $\Alg(\Sym\Cat)$.
A convergent presentation of the $(3,2)$-PROP $\Sym\Cat$ is constructed
in~\cite[Section~3.2]{GuiraudMalbos12mscs}. The coherence theorem for symmetric
monoidal categories~\cite{MacLane63} can be deduced from this construction: the
$(3,2)$-PROP $\Sym\Cat$ is
aspherical~\cite[Corollary~3.3.6]{GuiraudMalbos12mscs}.

\subsection{Braided monoidal categories}
\index{braided!monoidal category}
\index{monoidal category!braided}
A \emph{braided monoidal category} is a monoidal category
$(C, \otimes,e,\alpha,\lambda,\rho)$ equipped with a natural isomorphism
\[
\beta_{x,y} : x\otimes y \longrightarrow y\otimes x,
\]
called the \emph{braiding} and such that the following diagrams commute in $C$:
\[
  \xymatrix@C=1em{
    & x\otimes (y\otimes z) 
    \ar[rr] ^-{\beta_{x,y\otimes z}}
    && (y\otimes z)\otimes x
    \ar[dr] ^-{\alpha_{y,z,x}}
    \\
    (x\otimes y) \otimes z 
    \ar[ur] ^-{\alpha_{x,y,z}}
    \ar[dr] _-{\beta_{x,y}\otimes z}
    &&&& y \otimes (z\otimes x),
    \\
    & (y\otimes x)\otimes z 
    \ar[rr] _-{\alpha_{y,x,z}}
    && y\otimes(x\otimes z)
    \ar[ur] _-{y\otimes\beta_{x, z}}
  }
\]
\[
  \xymatrix@C=1em{
    & x\otimes (y\otimes z) 
    \ar[rr] ^-{\beta^-_{y\otimes z,x}}
    && (y\otimes z)\otimes x
    \ar[dr] ^-{\alpha_{y,z,x}}
    \\
    (x\otimes y) \otimes z 
    \ar[ur] ^-{\alpha_{x,y,z}}
    \ar[dr] _-{\beta^-_{y,x}\otimes z}
    &&&& y \otimes (z\otimes x)\pbox.
    \\
    & (y\otimes x)\otimes z 
    \ar[rr] _-{\alpha_{y,x,z}}
    && y\otimes(x\otimes z)
    \ar[ur] _-{y\otimes\beta^-_{z,x}}
  }
\]

\subsection{Generalized coherence theorems}
Contrarily to the case of monoidal and symmetric monoidal categories, we do not
have that every diagram commutes in a braided monoidal category. For instance,
the morphisms $\beta_{x,y}$ and $\beta^{-}_{y,x}$, from $x\otimes y$ to
$y\otimes x$, have no reason to be equal. In fact, they are equal if and only if
$\beta$ is a symmetry, hence if and only if all diagrams commute.  As a
consequence, the coherence problem for braided monoidal categories requires a
generalized version of the coherence problem we have considered so far. The
\emph{generalized coherence problem} is the following one:

\begin{center}
  \emph{Given a $(3,2)$-PROP~$\PP$, decide, for any 3-sphere $\alpha$ of
    $\PP$,\\
    whether or not the diagram $A(\alpha)$ commutes in every $\PP$-algebra $A$.}
\end{center}

A solution for the generalized coherence problem is a decision procedure for the
equality of $3$-cells of $\PP$. For the coherence problems considered so far,
this decision procedure answers yes for every
$3$-sphere. A method to study the generalized coherence theorem of $3$-PROPs is
given in~\cite[Section~4]{GuiraudMalbos12mscs} and illustrated
on the $(3,2)$-PROP of braided monoidal categories. In this way, we recover the
coherence result of Joyal and Street~\cite{joyal1993braided}: a diagram~$\Delta$
commutes if and only if both $\src{}(\Delta)$ and $\tgt{}(\Delta)$ have the same
associated braid.



\chapter{Term Rewriting Systems}
\label{chap:trs}
\label{Chapter:LawvereBurroni}
The study of \emph{universal algebra}, that is, the description of
algebraic structures by means of symbolic expressions subject to
equations,  dates back to the end of the
19th century~\cite{whitehead1898treatise}. It was motivated by the large number of
fundamental mathematical structures fitting into this framework: groups, rings,
lattices, and so on. From the 1970s on, the algorithmic aspect
became prominent and led to the notion of term rewriting system. This
chapter briefly revisits these ideas from a polygraphic viewpoint,
introducing only what is strictly necessary for understanding. 
We refer the reader to standard textbooks such
as~\cite{BaaderNipkow98,bezem2003term} for a proper study of this
vast topic.

In \cref{sec:lt-pres}, we begin by introducing term rewriting systems as
presentations of Lawvere theories, which are particular cartesian
categories. Some classical results on Lawvere theories are recalled in
\cref{sec:models}. The theory of rewriting in this context is explored in
\cref{sec:trs-rewriting} by defining rewriting steps, critical branchings, and
the original Knuth-Bendix completion procedure~\cite{knuth1970simple}. In \cref{sec:trs-3pol}, we show
that a term rewriting system can also be described by a $3$-polygraph in which
variables are handled explicitly, \ie by taking into account their duplication
and erasure. Finally, in \cref{sec:cart-pol}, we give a precise
meaning to the statement that
term rewriting systems are   ``cartesian polygraphs''.

\section{Presentations of Lawvere Theories}
\label{sec:lt-pres}

\subsection{Signatures}
\index{signature}
\label{sec:trs-sig}
In the context of term rewriting systems, a \emph{signature}~$P$ consists of
\begin{itemize}
\item a set~$P_0$ of \emph{sorts},
\item a set~$P_1$ of \emph{operations} together with functions
  $\src0:P_1\to\freecat{P_0}$ and $\tgt0:P_1\to P_0$
  respectively associating to each operation the sorts of its inputs and of its
  output, where $\freecat{P_0}$ denotes the free monoid over~$P_0$.
\end{itemize}
A \emph{morphism} $f:P\to Q$ between signatures~$P$ and~$Q$ consists of two functions
$f_0:P_0\to Q_0$ and $f_1:P_1\to Q_1$ such that
$\src0\circ f_1=\freecat{f_0}\circ \src0$ and $\tgt0\circ f_1=f_0\circ\tgt0$
(where $\freecat{f_0}:\freecat{P_0}\to\freecat{Q_0}$ is the extension of $f_0$ as a morphism of monoids).
We write~$\nCPol1$ for the resulting category (this notation will be justified in \cref{sec:cart-pol} below).

Given an operation $\alpha$, we write
$
  \alpha
  :
  a_1\ldots a_n
  \to
  a
$
to indicate that its source is $\src0(\alpha)=a_1\ldots a_n$ and target is $\tgt0(\alpha)=a$.
The natural number~$n$ is called the \emph{arity} of~$\alpha$.
A signature is \emph{mono-sorted} when $P_0$ is reduced to one element: in this
case, $\freecat{P_0}=\N$ and $\tgt0$ is the terminal function.

\subsection{Terms}
\index{variable}
\index{term}
Given a signature~$P$, a \emph{variable} is a symbol of the form~$x_i^u$ with
$u\in\freecat{P_0}$ and $i\in\N$.
A \emph{term} on a signature~$P$ is a ``well-typed'' tree whose nodes are
decorated in operations and leaves are decorated in variables. Formally, the
family of sets $\freecat{P_1}(u,a)$ of \emph{terms} from $u\in\freecat{P_0}$ to
$a\in P_0$ is the smallest family, indexed by~$u$ and~$a$, such that
\begin{itemize}
\item given $a_1,\ldots,a_n\in P_0$ and $1\leq i\leq n$, we have a variable term
  \[
    x_i^{a_1\ldots a_n}
    \qin
    \freecat{P_1}(a_1\ldots a_n,a_i)
    \pbox,
  \]
\item given an operation $\alpha:a_1\ldots a_n\to a$ in~$P_1$,
  $u\in\freecat{P_0}$, and terms $\phi_i\in\freecat{P_1}(u,a_i)$ for
  $1\leq i\leq n$, we have a composite term
  \[
    \alpha(\phi_1,\ldots,\phi_n)
    \qin
    \freecat{P_1}(u,a)
    \pbox.
  \]
\end{itemize}
We write $\phi:u\to a$ to indicate that $\phi$ is a term in
$\freecat{P_1}(u,a)$. In the following, we sometimes omit writing the
superscripts from variables.

\subsection{Substitutions}
\label{sec:trs-subst}
\index{substitution}
Given sorts $u\in\freecat{P_0}$ and $a_1,\ldots,a_n\in P_0$, a
\emph{substitution} $\sigma:u\to a_1\ldots a_n$ is an $n$-tuple of terms
$\sigma=\uple{\sigma_1,\ldots,\sigma_n}$, where the $\sigma_i:u\to a_i$ are
terms with $1\leq i\leq n$. Given a term $\phi:a_1\ldots a_n\to a$, we write
\[
  \phi\cdot\sigma
  :
  u\to a
\]
for the term obtained from~$\phi$ by replacing each variable $x_i$ by~$\sigma_i$:
this term is defined inductively by
\begin{align*}
  x_i\cdot\sigma&=\sigma_i,
  &
  \alpha(\phi_1,\ldots,\phi_m)\cdot\sigma&=\alpha(\phi_1\cdot\sigma,\ldots,\phi_m\cdot\sigma)
  \pbox.
\end{align*}

\subsection{The generated category}
\label{sec:law-gen-cat}
Given a signature~$P$, we write $\freecat{P}$ for the category with
$\freecat{P_0}$ as objects and substitutions $\sigma:u\to v$ as morphisms.
Given two substitutions $\sigma:u\to v$ and $\tau:v\to w$ with
$\tau=\uple{\tau_1,\ldots,\tau_n}$, their composite is the substitution
\[
  \tau\circ\sigma
  =
  \uple{\tau_1\cdot\sigma,\ldots,\tau_n\cdot\sigma},
\]
and given an object $u=a_1\ldots a_n$ the identity on $u$
is~$\uple{x^u_1,\ldots,x^u_n}$. Note that, given a term~$\phi:w\to a$, we have
\begin{align*} 
  \phi\cdot(\tau\circ\sigma)
  &=
  (\phi\cdot\tau)\cdot\sigma,
  &
  \phi\cdot\uple{x^u_1,\ldots,x^u_n}
  &=
  \phi
  \pbox. 
\end{align*}
We write $\freecat{P_1}$ for the set of all morphisms of~$\freecat{P}$ and
$
  \freecat{\src0},\freecat{\tgt0}
  :
  \freecat{P_1}
  \to
  P_0
$
for the source and target functions.

\subsection{Cartesian categories}
\index{cartesian!category}
\index{category!cartesian}
In a category~$C$, a \emph{cartesian product} of two objects~$u$ and~$v$ is an
object, usually noted $u\times v$, together with morphisms
$\pi_1:u\times v\to u$ and $\pi_2:u\times v\to v$, called \emph{projections},
such that for every object~$w$ and morphisms $\phi:w\to u$ and $\psi:w\to v$,
there exists a unique morphism $\uple{\phi,\psi}:w\to u\times v$ satisfying
$\pi_1\circ\uple{\phi,\psi}=\phi$ and $\pi_2\circ\uple{\phi,\psi}=\psi$:
\[
  \xymatrix@C=4ex@R=4ex{
    &w\ar@/_/[ddl]_\phi\ar@{.>}[d]|-{\uple{\phi,\psi}}\ar@/^/[ddr]^\psi&\\
    &\ar[dl]^{\pi_1}u\times v\ar[dr]_{\pi_2}&\\
    u&&v\pbox.
  }
\]
An object~$1$ is \emph{terminal} in a category when for every object~$u$ there
exists a unique morphism $u\to 1$.
A category is \emph{cartesian} when it has a terminal object and every pair of
objects admits a cartesian product. In the following, for simplicity, we suppose
fixed a choice of a product for any pair of objects in~$C$, which we can suppose to
be strictly associative and unital by Mac Lane's coherence theorem
(\cref{thm:moncat-coh}).

A \emph{morphism} $f:C\to D$ of cartesian categories, also called a
\emph{cartesian functor}, is a functor which preserves cartesian products and the
terminal object. Here, we only consider functors for which this preservation is
strict, by which we mean that $f(u\times v)=f(u)\times f(v)$ and $f(1)=1$. We write
$\Cart$, or $\nCart 1$, for the category of cartesian categories.
\nomenclature[Cart]{$\Cart$}{category of cartesian categories}

\begin{lemma}
  The category~$\freecat{P}$ is cartesian.
\end{lemma}
\begin{proof}
  The product of two objects $u=a_1\ldots a_p$ and $v=b_1\ldots b_q$ is given by
  their concatenation~$uv$, and the canonical projections are
  \begin{align*}
    \uple{x^{uv}_1,\ldots,x^{uv}_p}&:uv\to u,
    &
    \uple{x^{uv}_{p+1},\ldots,x^{uv}_{p+q}}&:uv\to v
    \pbox.    
  \end{align*}
  Finally, given two morphisms
  \begin{align*}
    \sigma&=\uple{\sigma_1,\ldots,\sigma_p}:w\to u,
    &
    \tau&=\uple{\tau_1,\ldots,\tau_q}:w\to v,   
  \end{align*}
  the associated universal morphism is
  \[
    \uple{\sigma,\tau}
    =
    \uple{\sigma_1,\ldots,\sigma_p,\tau_1,\ldots,\tau_q}
    :
    w
    \to
    uv
    \pbox.
    \qedhere
  \]
\end{proof}

\nomenclature[F]{$\Fun$}{category of functions}
Let us describe an important case of the above construction. We write $\Fun$ for
the category whose objects are natural numbers and morphisms $m\to n$ are
functions $[m]\to[n]$ where $[n]=\set{0,\ldots,n-1}$ is a set with $n$ elements,
see also \cref{sec:2pres-mon}. We write $I:\Fun\to\Set$ for the canonical
inclusion functor. Given a set~$P_0$, we (abusively) write $\Fun/P_0$ for the comma
category $\commacat{I}{P_0}$ of~$I$ over the set~$P_0$.

\begin{lemma}
  Given a signature~$P$ such that~$P_1=\emptyset$, we have
  $
    \freecat{P}
    \isoto
    (\Fun/P_0)^\op
  $.
  In particular, if~$P_0=\set{\star}$ then $\freecat{P}\isoto\Fun^\op$.
\end{lemma}

\noindent
We now describe the universal property satisfied by the
construction~$\freecat{P}$.

\subsection{Lawvere theories}
\index{Lawvere theory}
\label{sec:lawvere-theory}
Suppose fixed a set~$P_0$ of sorts. A \emph{$P_0$-sorted} \emph{Lawvere theory}
(or \emph{algebraic theory}) is a cartesian category~$C$ equipped with functor
\[
  (\Fun/P_0)^\op\to C,
\]
which preserves finite products and is the identity on objects. A morphism
between two Lawvere theories~$C$ and~$D$ is a functor $f:C\to D$ making the
following diagram commute:
\[
  \xymatrix@C=3ex@R=3ex{
    C\ar[rr]^f&&D\\
    &\ar[ul](\Fun/P_0)^\op\ar[ur]
  }
\]
We write $\Law_{P_0}$ for the resulting category.
\nomenclature[Law]{$\Law$}{category of Lawvere theories}

\subsection{The free Lawvere theory}
\index{free!Lawvere theory}
Consider the subcategory~$\category{S}_{P_0}$ of the category $\nCPol1$ of
signatures, where objects are the signatures having~$P_0$ as sorts, and
morphisms are those which are identity on sorts. There is a forgetful
functor
\[
  W_0
  :
  \Law_{P_0}
  \to
  \category{S}_{P_0}
\]
sending a Lawvere theory~$C$ to the signature~$P$ with
\[
  P_1=\coprod_{u\freecat{P_0},a\in P_0}C(u,a)
\]
as operations, with source and target being respectively given by the
indices~$u$ and~$a$ of the coproduct.
The cartesian category generated by a signature introduced in
\cref{sec:law-gen-cat} can be shown to be freely generated in the following
sense.

\begin{proposition}
  \label{prop:free-lt}
  The functor~$W_0$ admits a left adjoint
  \[
    L_0
    :
    \category{S}_{P_0}
    \to
    \Law_{P_0},
  \]
  such that the image of a signature~$P$ is the Lawvere theory~$\freecat{P}$.
\end{proposition}

\subsection{Congruence}
\index{congruence!on a Lawvere theory}
A \emph{congruence} on a Lawvere theory~$C$ is a relation~$\approx$ on parallel
morphisms, such that
\begin{itemize}
\item given morphisms $f:u'\to u$, $h:v\to v'$, and $g,g':u\to v$, 
  \[
    g\approx g'
    \qquad\text{implies}\qquad
    f\comp{}g\comp{}h
    \approx
    f\comp{}g'\comp{}h
    \pbox,
  \]
\item given morphisms $f,f':w\to u$ and $g,g':w\to v$,
  \[
    f\approx f'
    \qtand
    g\approx g'
    \qquad\text{implies}\qquad
    \uple{f,g}\approx\uple{f',g'}
    \pbox.
  \]
\end{itemize}
Given such a congruence, we write $C/{\approx}$ for the associated
\emph{quotient} Lawvere theory, obtained from~$C$ by quotienting morphisms
under~$\approx$.

\subsection{Term rewriting systems}
\label{sec:trs}
\index{term!rewriting system}
A \emph{term rewriting system} $P$ consists of a signature
$(P_0,\src0,\tgt0,P_1)$ together with a set~$P_2$ of \emph{rewriting rules} (or
\emph{relations}) equipped with \emph{source} and \emph{target}
functions~$\src1,\tgt1:P_2\to\freecat{P_1}$ such that the associated morphisms
are parallel, \ie $\freecat{\src0}\circ\src1=\freecat{\src0}\circ\tgt1$ and
$\freecat{\tgt0}\circ\src1=\freecat{\tgt0}\circ\tgt1$. A rewriting rule~$A$ with
source (\resp target) $\phi:u\to a$ (\resp $\psi:u\to a$) is often denoted
\[
  A
  :
  \phi\To\psi
  :
  u\to a
  \pbox.
\]
The Lawvere theory \emph{presented} by a term rewriting system~$P$ is
$
  \pcat{P}
  =
  \freecat{P}/P_2
$,
\ie the category obtained from~$\freecat{P}$ by quotienting morphisms under the
congruence~$\approx^P$ generated by~$P_2$, and we say that~$P$ is a
\emph{presentation} of~$\pcat{P}$.

A morphism~$P\to Q$ between term rewriting systems consists of a morphism
between the underlying signatures together with a function~$P_2\to Q_2$ which is
compatible with source and target. We write $\nCPol2$ for the category of term
rewriting systems.

\subsection{Models}
\index{model!of a Lawvere theory}
\label{sec:lt-model}
A \emph{model} of a Lawvere theory~$C$ is a functor $C\to\Set$ which preserves
finite products. In the case where~$C$ admits a presentation~$P$, a model
amounts to the data of
\begin{itemize}
\item a set $\intp{a}$ for every sort $a\in P_0$,
\item a function
  $\intp{\alpha}:\intp{a_1}\times\ldots\times\intp{a_n}\to\intp{a}$ for every
  operation $\alpha:a_1\ldots a_n\to a$ in $P_1$,
\end{itemize}
such that $\intp{\phi}=\intp{\psi}$ for every relation $A:\phi\To\psi$, where
$\intp{-}$ is extended to terms by
\[
  \intp{\alpha(\phi_1,\ldots,\phi_n)}=\intp{\alpha}\circ\uple{\intp{\phi_1},\ldots,\intp{\phi_n}}
\]
and
$\intp{x^{a_1\ldots a_n}_i}:\intp{a_1}\times\ldots\times\intp{a_n}\to\intp{a_i}$
is the canonical projection. We sometimes abusively speak of a model of a
signature (\resp term rewriting system) to mean a model of the generated (\resp
presented) Lawvere theory.

\begin{example}
  \label{ex:trs-group}
  \index{group}
  \index{theory!of groups}
  The theory of groups is presented by~$P$ with
  \begin{align*}
    P_0&=\set{a},
    &
    P_1&=\set{\mu:2\to 1,\eta:0\to 1,\iota:1\to 1},
  \end{align*}
  and rewriting rules
  \begin{align*}
    \mu(\eta,x_1)&\To x_1,
    &
    \mu(x_1,\eta)&\To x_1,
    &
    \mu(\mu(x_1,x_2),x_3)&\To\mu(x_1,\mu(x_2,x_3)),
    \\
    \mu(\iota(x_1),x_1)&\To\eta,
    &
    \mu(x_1,\iota(x_1))&\To\eta,
  \end{align*}
  where, given $n\in\N$, we write $n$ instead of $a^n$ for an element
  of~$\freecat{P_0}$. A model for this theory is a group.
\end{example}

\noindent
Of course, as a variation of the previous example, usual algebraic structures have
an associated Lawvere theory: groups, rings, modules, vector spaces, algebras,
lattices, etc. Most are mono-sorted, apart from the theory for modules (as well
as the one of vector spaces) which has two sorts: one corresponding to the ring
of scalars and one to the abelian group. As a notable exception, there is no
Lawvere theory corresponding to fields: intuitively, this is because the inverse
operation is only partially defined ($0$ is not invertible). Below, we give an
example of a Lawvere theory of more computational nature.

\begin{example}
  \label{ex:ski}
  \index{combinatory logic}
  Combinatory logic was introduced by
  Schönfinkel~\cite{schonfinkel1924bausteine} and
  Curry~\cite{curry1930grundlagen} as an algebraic way of capturing binding and
  substitution. It can be presented by the mono-sorted term rewriting system~$P$
  with operations
  \begin{align*}
    \alpha&:2\to 1,
    &
    \sigma&:0\to 1,
    &
    \kappa&:0\to 1,
    &
    \iota&:0\to 1
  \end{align*}
  ($\alpha$ should be read as an ``application'', and the constants $\sigma$,
  $\kappa$, and $\iota$ are usually respectively denoted $S$, $K$, and $I$) and
  relations
  \begin{align*}
    \alpha(\alpha(\alpha(\sigma,x_1),x_2),x_3)
    &\To
    \alpha(\alpha(x_1,x_2),\alpha(x_1,x_3))\text,
    \\
    \alpha(\alpha(\kappa,x_1),x_2)
    &\To
    x_1\text,
    \\
    \alpha(\iota,x_1)
    &\To
    x_1\text.
  \end{align*}
  A combinatory term~$\phi$ in~$\freecat{P_1}$ can be interpreted as a
  $\lambda$\nbd-term~$\intp{\phi}$ by
  \begin{align*}
    \intp{\alpha(\phi,\psi)}&=\intp{\phi}\intp{\psi}\text,
    &
    \intp{\sigma}&=\lambda xyz.(xz)(yz)\text,
    &
    \intp{\kappa}&=\lambda xy.x\text,
    &
    \intp{\iota}=\lambda x.x\text,
  \end{align*}
  and conversely every $\lambda$-term can be interpreted as a combinatory logic
  term, giving rise to a correspondence between the morphisms in the presented category
  and $\lambda$\nbd-terms modulo $\beta$-reduction, although the details are subtle,
  see~\cite{selinger2002lambda} for a survey on the subject. A model of this
  Lawvere theory is called a \emph{combinatory algebra}.
\end{example}

\subsection{Tietze transformations}
\index{Tietze!transformation!of term rewriting systems}
\index{transformation!Tietze}
The \emph{elementary Tietze transformations} consist, starting from a presentation~$P$, in
\begin{description}
\item[\tgen] \emph{adding a superfluous operation}: given a term $\phi:u\to a$ in
  $\freecat{P_1}$, we construct the presentation~$P'$ such that
  \begin{align*}
    P'_0&=P_0,
    &
    P'_1&=P_1\sqcup\set{\alpha:u\To a},
    &
    P'_2&=P_2\sqcup\set{A:\phi\To\alpha},    
  \end{align*}
\item[\trel] \emph{adding a superfluous relation}: given two terms $\phi$ and $\psi$
  such that $\phi\approx^P\psi$, we construct the presentation~$P'$ such that
  \begin{align*}
    P'_0&=P_0,
    &
    P'_1&=P_1,
    &
    P'_2&=P_2\sqcup\set{A:\phi\To\psi}. 
  \end{align*}
\end{description}
The \emph{Tietze equivalence} is the smallest equivalence relation on
presentations such that~$P$ is Tietze equivalent to~$P'$ whenever there exists a
Tietze transformation from~$P$ to~$P'$. The proof of the following theorem
carries over as in the case of polygraphs (see \thmr{2-tietze-equiv}).

\begin{theorem}
  Two presentations with finite sets of operations and relations present
  isomorphic categories if and only if they are Tietze equivalent.
\end{theorem}

\subsection{Composing presentations}
\label{sec:trs-dlaw}
\index{distributive law!of Lawvere theories}
In \cref{sec:2-dlaws}, we have seen that we could compose presented categories,
when given a distributive law between them, and this was extended to
presentations of monoidal categories in \cref{sec:3-dlaws}. We mention
here that this also generalizes to Lawvere theories: the corresponding notion of
distributive law is studied in~\cite{cheng2011distributive}.

\section{More on Models}
\label{sec:models}
In this section, we briefly recall some of the classical theory of Lawvere
theories, as initiated by Lawvere in his PhD thesis~\cite{lawvere1963functorial},
see~\cite{adamek1994localp, adamek2011algebraic} for an in-depth
presentation. For the sake of simplicity, we only handle here the mono-sorted case.

\subsection{Models}
Given a Lawvere theory~$C$, we have seen in \secr{lt-model} that a model is a
functor $C\to\Set$ which preserves finite limits. A morphism between models is a
natural transformation and we write $\Mod(C)$ for the category of models.

\subsection{Free models}
Given a Lawvere theory~$C$, there is a forgetful functor
\[
  U
  :
  \Mod(C)
  \to
  \Set
\]
which to a model $M:C\to\Set$ associates $M(1)$.

\begin{theorem}
  The functor~$U$ is monadic: it admits a left adjoint and the category
  $\Mod(C)$ is equivalent to the category of $T$-algebras, where $T$ is the
  monad associated to the adjunction.
\end{theorem}

\subsection{Monads}
\index{monad!finitary}
\newcommand{\Mnd}{\category{Mnd}}
By the above theorem, every Lawvere theory~$C$ canonically induces a monad~$T$
on~$\Set$. An explicit description of this monad can be given by the following
coend formula:
\[
  TX
  =
  \int^nC(n,1)\times X^n.
\]
Not every monad arises in this way, and those which do can be characterized as
being \emph{finitary}, \ie preserving filtered colimits. Writing $\Mnd$ for the
category of monads on~$\Set$, we have the following equivalence of categories.

\begin{theorem}
  The category~$\Law$ of Lawvere theories is equivalent to the full subcategory
  of~$\Mnd$ whose objects are finitary monads.
\end{theorem}

\noindent
We have already explained above how to associate a finitary monad to a Lawvere
theory. Conversely, given such a monad~$T$, the opposite category of the Kleisli
category is always a Lawvere theory. These constructions give rise to the
equivalence stated in the theorem.

\subsection{The Birkhoff theorem}
\index{Birkhoff theorem}
\index{HSP theorem}
We now turn to a different approach to models of a Lawvere theory. Fix a
signature~$P$. In the context of model theory, its models are sometimes called
\emph{structures}\index{structure}.
Given a presentation~$Q$ on this signature (\ie $Q_0=P_0$ and $Q_1=P_1$), we
have a quotient functor
\[
  \freecat{P}\to\pcat{Q}=\freecat{P}/Q_2,
\]
which induces, by precomposition, a functor
\[
  \Mod(Q)\to\Mod(P)
\]
between the categories of models. The functor $P\to Q$ being surjective on
objects and full, the induced functor between models is full and faithful, and
we can thus consider $\Mod(Q)$ as a full subcategory of~$\Mod(P)$. Conversely,
given a full subcategory~$\C$ of~$\Mod(P)$, one may wonder whether there is a
set~$P_2$ of relations such that~$\C$ is precisely the category of models of the
Lawvere theory presented by $(P,P_2)$. The following theorem, due to
Birkhoff~\cite{birkhoff1935structure}, see~\cite{adamek2011algebraic},
and sometimes called the \emph{HSP theorem}, gives a
characterization of those situations. 

\begin{theorem}
  Given a signature~$P$ and a full subcategory~$\C$ of $\Mod(P)$, $\C$ is the
  category of models of a term rewriting system~$Q$ on the signature~$P$ if and
  only if it is closed under
  \begin{description}[labelwidth=\widthof{(H)}]
  \item[(H)] \emph{homomorphic images}: given a regular epimorphism $f:M\to N$ in
    $\Mod(P)$ (\ie $f$ is an epi which can be obtained as a coequalizer) with
    $M\in\C$, the object~$N$ also belongs to~$\C$,
  \item[(S)] \emph{subalgebras}: given a monomorphism $f:M\to N$ in~$\Mod(P)$ with
    $N\in\C$, the object~$M$ also belongs to~$\C$, and
  \item[(P)] \emph{products}: given $M,N\in\C$, their product $M\times N$ in $\Mod(P)$
    also belongs to~$\C$.
  \end{description}
\end{theorem}

\section{Term Rewriting}
\label{sec:trs-rewriting}
Up to now, we have been using term rewriting systems as a notion of
presentation, for which the orientation of the rules does not really matter. We
now introduce the rewriting structure, following what we have done for
$1$-polygraphs (\cref{sec:ars}), $2$-polygraphs (\cref{chap:2rewr}), and
$3$-polygraphs (\cref{sec:3pol-rewr}).

\subsection{Occurrences}
\newcommand{\occ}[2]{o_{#1}(#2)}
\index{occurrence}
Given a term $t:a_1\ldots a_n\to a$ and an index $i$, with $1\leq i\leq n$, the
number $\occ it$ of \emph{occurrences} of the $i$-th variable $x_i$
into~$t$ is defined by induction on~$t$ by
\begin{align*}
  \occ i{x_i}&=1,
  &
  \occ i{x_j}&=0,
  &
  \occ i{\alpha(\phi_1,\ldots,\phi_k)}&=\sum_{i=1}^k\occ i{\phi_i},
\end{align*}
for $j\neq i$. A variable~$x_i$ is \emph{linear} in a term~$t$ when it occurs exactly
once, \ie $\occ it=1$.

\subsection{Contexts}
\label{sec:trs-context}
\index{context}
A \emph{context} $\kappa: a_1\ldots a_na\to b$ is a term such that the
variable~$x_{n+1}$ (of type~$a$) is linear in~$C$. Given a
term~$\phi:a_1\ldots a_n\to a$, we write
\[
  \kappa\cdot\phi
  :
  a_1\ldots a_n
  \to
  b
\]
for the term obtained from~$\kappa$ by substituting~$\phi$ for~$x_{n+1}$, \ie
\[
  \kappa\cdot\phi
  =
  \kappa\cdot\uple{x_1,\ldots,x_n,\phi},
\]
with the notations of~\secr{trs-subst}.

Given two contexts $\kappa:a_1\ldots a_na\to b$ and $\rho:a_1\ldots a_nb\to c$,
their \emph{composite} is the context
\[
  \rho\circ\kappa
  :
  a_1\ldots a_na\to c,
\]
obtained from~$\rho$ by replacing the variable $x_{n+1}$ (of type~$b$)
by~$\kappa$, \ie
\[
  \rho\circ\kappa
  =
  \rho\cdot\kappa
  =
  \rho\cdot\uple{x_1,\ldots,x_n,\kappa}
\]
and the \emph{identity context} is
\[
  x_{n+1}^{a_1\ldots a_na}
  :
  a_1\ldots a_na\to a
  \pbox.
\]
Note that given a term~$\phi:a_1\ldots a_n\to a$, we have
\begin{align*}
  (\rho\circ\kappa)\cdot\phi
  &=
  \rho\cdot(\kappa\cdot\phi),
  &
  x_{n+1}^{a_1\ldots a_na}\cdot\phi
  &=
  \phi
  \pbox.  
\end{align*}
Given a fixed $u\in\freecat{P_0}$, we can thus build a category with~$P_0$ as
set of objects, a morphism $\kappa:a\to b$ being a context $\kappa:ua\to b$.

It is also useful to introduce a \emph{base change} operation on contexts. Given
a context $\kappa:a_1\ldots a_na\to b$ and a substitution
$\sigma:u\to a_1\ldots a_n$, we write
\[
  \sigma^*(\kappa)
  :
  ua
  \to
  b
\]
for the context defined by
\[
  \sigma^*(\kappa)
  =
  \kappa\cdot\uple{\sigma_1,\ldots,\sigma_n,x_{n+1}^{ua}}
  \pbox.
\]

\subsection{Contexts and substitutions}
\label{sec:trs-ctxt-sub}
Given a context $\kappa$, a term $\phi$ and a substitution~$\sigma$ of
appropriate type, an expression of the form $\kappa\cdot\phi\cdot\sigma$ is
always implicitly bracketed as $\kappa\cdot(\phi\cdot\sigma)$.
These operations are compatible with the categorical structures in the sense
that, for suitably typed contexts and substitutions, we have
\begin{align*}
  \rho\cdot(\kappa\cdot\phi\cdot\tau)\cdot\sigma
  &=
  (\rho\circ\sigma^*(\kappa))\cdot\phi\cdot(\tau\circ\sigma),
  &
  \unit{}\cdot\phi\cdot\unit{}
  &=
  \phi
  \pbox.  
\end{align*}
Graphically, it is sometimes convenient to depict the
term~$\kappa\cdot\phi\cdot\sigma$ as
\[
  \begin{tikzpicture}
    \draw (-.5,.5) rectangle (.5,1);
    \draw (0,.75) node {$\sigma$};
    \draw (-.5,.5) -- (0,0) -- (.5,.5) node[right] {\,.};
    \draw (0,.25) node{$\phi$};
    \draw (-.5,1) to[bend right=40] (0,-.5) to[bend right=40] (.5,1);
    \draw (0,-.25) node{$\kappa$};
  \end{tikzpicture}
\]

\subsection{$k$-ary contexts}
Generalizing the construction of \secr{trs-context}, a \emph{$k$\nbd-ary
  context}, for $k\in\N$, is a term
\[
  \kappa:a_1\ldots a_na_1'\ldots a_k'\to a,
\]
such that the variables $x_{n+1}$, \ldots, $x_{n+k}$ are linear
in~$\kappa$. Given terms
\[
  \phi_i
  :
  a_1\ldots a_n\to a'_i,
\]
with $1\leq i\leq k$, we write
\[
  \kappa\cdot(\phi_1,\ldots,\phi_k)
  =
  \kappa\uple{x_1,\ldots,x_n,\phi_1,\ldots,\phi_k}
  \pbox.
\]
In the following, we will only use binary contexts, \ie the case where $k=2$.

\subsection{Rewriting steps}
\index{rewriting!step}
A \emph{rewriting step}
\[
  \kappa\cdot A\cdot\sigma
  :
  \kappa\cdot\phi\cdot\sigma
  \To
  \kappa\cdot\psi\cdot\sigma
  :
  u'\to a'
\]
consists of a rewriting rule
\[
  A
  :
  \phi\To\psi
  :
  u\to a,
\]
together with a substitution and a context
\[
  \sigma
  :
  u'\to u
  \qqtand
  \kappa
  :
  u'a\to a'
  \pbox.
\]
In this case, we say that the term $\kappa\cdot\phi\cdot\sigma$ \emph{rewrites}
in one step to the term $\kappa\cdot\psi\cdot\sigma$. A \emph{rewriting path} is
a sequence of composable rewriting steps.

\subsection{Rewriting}
\index{termination}
\index{rewriting!system!terminating}
\index{abstract rewriting system!of a term rewriting system}
A term rewriting system~$P$ induces an abstract rewriting system with the terms
in~$\freecat{P_1}$ as vertices and rewriting steps as edges. The rewriting
system is \emph{terminating}, (\emph{locally}) \emph{confluent}, etc., when the
associated abstract rewriting system is.

\subsection{Critical branchings}
\index{branching}
A \emph{branching} is a pair
\begin{equation}
  \label{eq:trs-branching}
  \xymatrix@C=7ex{
    \kappa_1\cdot\psi_1\cdot\sigma_1
    &
    \ar@{=>}[l]_-{\kappa_1\cdot A_1\cdot\sigma_1}
    \kappa_1\cdot\phi_1\cdot\sigma_1
    =
    \kappa_2\cdot\phi_2\cdot\sigma_2
    \ar@{=>}[r]^-{\kappa_2\cdot A_2\cdot\sigma_2}
    &
    \kappa_2\cdot\psi_2\cdot\sigma_2
  }
\end{equation}
of coinitial rewriting steps. We can identify the following families of
branchings.
\begin{enumerate}
\item A branching is \emph{trivial} when it is of the form
  \[
    \xymatrix@C=8ex{
      \psi
      &
      \ar@{=>}[l]_-{\kappa\cdot A\cdot\sigma}
      \phi
      \ar@{=>}[r]^-{\kappa\cdot A\cdot\sigma}
      &
      \psi.
    }
  \]
  Such a branching is clearly confluent.
\item A branching is \emph{parallel orthogonal} when
  it is of the form
  \[
    \xymatrix@C=14ex{
      \phi_1
      &
      \ar@{=>}[l]_-{\kappa\cdot(A_1\cdot\sigma_1,\psi_2\cdot\sigma_2)}
      \phi
      \ar@{=>}[r]^-{\kappa\cdot(\psi_1\cdot\sigma_1,A_2\cdot\sigma_2)}
      &
      \phi_2,
    }
  \]
  for some suitably-typed binary context~$\kappa$, rewriting rules
  $A_1:\psi_1\To\psi_1'$ and $A_2:\psi_2\To\psi_2'$, and substitutions
  $\sigma_1$ and $\sigma_2$.
\item A branching is \emph{inclusion orthogonal} when it is of the form
  \[
    \xymatrix@C=24ex{
      \phi_1
      &
      \ar@{=>}[l]_-{\kappa\cdot A_1\cdot\uple{\sigma_1,\ldots,\sigma_n}}
      \phi
      \ar@{=>}[r]^-{\kappa\cdot\psi_1\cdot\uple{\sigma_1,\ldots,\kappa'\cdot A_2\cdot\sigma',\ldots,\sigma_n}}
      &
      \phi_2,
    }
  \]
  for some suitably-typed contexts~$\kappa$ and $\kappa'$, rewriting rules
  $A_1:\psi_1\To\psi_1'$ and $A_2:\psi_2\To\psi_2'$, and substitutions
  $\sigma=\uple{\sigma_1,\ldots,\sigma_n}$ and $\sigma'$, such that
  $\sigma_i=\kappa'\cdot\psi_2\cdot\sigma'$ for some index $i$ with
  $1\leq i\leq n$.
\item A branching is \emph{non-minimal} when it is of the form
  \[
    \xymatrix@C=14ex{
      \phi_1
      &
      \ar@{=>}[l]_-{\kappa\cdot(\kappa_1\cdot A_1\cdot\sigma_1)\cdot\sigma}
      \phi
      \ar@{=>}[r]^-{\kappa\cdot(\kappa_1\cdot A_1\cdot\sigma_1)\cdot\sigma}
      &
      \phi_2,
    }
  \]
  for some suitably-typed context~$\kappa$ and substitution~$\sigma$, not both
  identities, and of suitable types.
\item A branching is \emph{critical}\index{branching!critical}\index{critical!branching} when it is not of
  any of the above forms.
\end{enumerate}

The above definition of critical branching makes it easy to show the critical
branchings lemma, stated below. Moreover, those can be efficiently computed,
see~\cite[Section~6.2]{BaaderNipkow98} for a presentation of the classical
algorithms.

\begin{lemma}
  \label{lem:trs-critical-branching}
  A term rewriting system is locally confluent if and only if all its critical
  branchings are.
\end{lemma}

\begin{example}
  \index{monoid}
  \index{theory!of monoids}
  \label{ex:trs-monoids}
  The theory of monoids, with 0-generators $P_0=\set{a}$, 1\nbd-gene\-rators
  $P_1=\set{\mu:2\to 1,\eta:0\to 1}$ and 2-generators
  \begin{align*}
    \mu(\eta,x_1)&\To x_1,
    &
    \mu(x_1,\eta)&\To x_1,
    &
    \mu(\mu(x_1,x_2),x_3)&\To\mu(x_1,\mu(x_2,x_3)),
  \end{align*}
  is locally confluent since its five critical branchings, whose source is shown
  below, are confluent:
  \[
    \begin{array}{c@{\qquad}c@{\qquad}c}
    \mu(\mu(\mu(x_1,x_2),x_3),x_4),
    &
    \mu(\mu(\eta,x_1),x_2),
    &
    \mu(\mu(x_1,\eta),x_2),
    \\
    &
    \mu(\mu(x_1,x_2),\eta),
    &
    \mu(\eta,\eta).       
    \end{array}
  \]
  The rewriting system can be shown to be terminating and is thus convergent.
\end{example}

\subsection{Reduction orders}
\index{reduction!order}
\index{order!reduction}
Given a signature $P$, a \emph{reduction order}~$\succ$ is an order on $P_1^*$
which is
\begin{itemize}
\item well-founded,
\item closed under application: for every $\alpha\in P_1$ of arity~$n$ and
  every terms
  $\phi_1,\ldots,\phi_n\in P_1^*$ and $\phi_i'\in P_1^*$, we have that $\phi_i\succ\phi_i'$ implies
  \[
    \alpha(\phi_1,\ldots,\phi_{i-1},\phi_i,\phi_{i+1},\ldots,\phi_n)
    \succ
    \alpha(\phi_1,\ldots,\phi_{i-1},\phi_i',\phi_{i+1},\ldots,\phi_n)
    \pbox,
  \]
\item closed under substitution: given terms $\phi,\phi'\in\freecat P_1$ and
  substitution $\sigma$, we have $\phi\succ\phi'$ implies
  $\phi\cdot\sigma\succ\phi'\cdot\sigma$.
\end{itemize}
Given a rewriting system $(P,P_2)$ on a signature~$P$, a \emph{termination}
order~$\succ$ is a reduction order on~$P$ such that $\phi\succ\psi$ for every
rewriting rule $A:\phi\To\psi$ in~$P$. As in the case of $2$-polygraphs
(\cref{prop:term-red-ord}), we have~\cite[Theorem~5.2.3]{BaaderNipkow98}:

\begin{proposition}
  A rewriting system is terminating if and only if it admits a termination
  order.
\end{proposition}

\subsection{Completion}
\label{sec:trs-completion}
\index{Knuth-Bendix completion}
\index{completion}
\index{procedure!completion}
Given a rewriting system~$P$ equipped with a termination order~$\succ$, we can
turn it into a convergent rewriting system by the following Knuth-Bendix
completion procedure~\cite{knuth1970simple}. It is very similar to the one
already given in \cref{sec:2kb} and consists in iteratively applying the
following steps.
\begin{itemize}
\item For every critical branching $\phi_1\Leftarrow\phi\To\phi_2$, compute
  normal forms~$\nf\phi_1$ and $\nf\phi_2$ for~$\phi_1$ and~$\phi_2$,
  respectively.
\item If~$\nf\phi_1=\nf\phi_2$ for every possible branching, the procedure halts.
\item Otherwise, there is a critical branching with $\nf\phi_1\neq\nf\phi_2$:
  \begin{itemize}
  \item if $\nf\phi_1\succ\nf\phi_2$, we add the rule $\nf\phi_1\To\nf\phi_2$
    to~$P$,
  \item if $\nf\phi_2\succ\nf\phi_1$, we add the rule $\nf\phi_2\To\nf\phi_1$
    to~$P$.
  \end{itemize}
\end{itemize}

As in the case of $2$-polygraphs, the procedure is not guaranteed to stop. In
the case it does, the resulting rewriting system is a convergent presentation of
the Lawvere theory $\pcat P$.
When the procedure does not terminate, the above steps produce an infinite
sequence of rewriting systems~$P^i$ with $P^0=P$, by iteratively adding rules,
and the inductive limit $\bigcup_iP^i$ is always a convergent presentation
of~$\pcat P$.

\begin{example}
  The presentation of groups given in \cref{ex:trs-group} is not locally
  confluent. By applying the Knuth-Bendix completion procedure, one can arrive
  at the following convergent presentation with the same 0- and 1-generators
  and whose relations are those of \exr{trs-group} together with
  \begin{align*}
    \mu(\iota(x_1,\mu(x_1,x_2)))&\To x_2,
    &
    \iota(\eta)&\To\eta,
    &
    \iota(\mu(x_1,x_2))&\To\mu(\iota(x_2),\iota(x_1)),
    \\
    \mu(x_1,\mu(\iota(x_1,x_2)))&\To x_2,
    &
    \iota(\iota(x_1))&\To x_1,
  \end{align*}
  see~\cite[Example~1]{knuth1970simple}.
\end{example}

\subsection{Non-linearity and confluence}
Because of the presence of variables which can potentially duplicate terms when substituted, one
should be careful when the rewriting system is not terminating. For instance, contrary to the
case of polygraphs studied in previous chapters, it
is not true that a rewriting system without critical branchings is always
confluent. Namely, consider the mono-sorted rewriting system due to
Huet~\cite{huet1980confluent}, with generators
\begin{align*}
  \tau&:0\to 1,
  &
  \phi&:0\to 1,
  &
  \omega&:0\to 1,
  &
  \sigma&:1\to 1,
  &
  \varepsilon&:2\to 1,
\end{align*}
(which should respectively be read as ``true'', ``false'', ``infinity'',
``successor'', and ``equality'') and relations
\begin{align*}
  \omega&\To\sigma(\omega),
  &
  \varepsilon(x_1,x_1)&\To\tau,
  &
  \varepsilon(x_1,\sigma(x_1))&\To\phi.
\end{align*}
The first rule clearly makes the rewriting system non-terminating.
There is no critical branching, yet the system is not confluent:
\[
  \xymatrix@C=3ex{
    \tau
    &
    \ar@{=>}[l]
    \varepsilon(\omega,\omega)
    \ar@{=>}[r]
    &
    \varepsilon(\omega,\sigma(\omega))
    \ar@{=>}[r]
    &
    \phi.
  }
\]
It can however be shown that a rewriting system which is \emph{left-linear} (\ie where
no variable occurs twice in the left member of a rewriting rule) and without
critical branchings is always confluent~\cite[Section~6.4]{BaaderNipkow98}.

\section{Term Rewriting Systems and 3-Polygraphs}
\label{sec:trs-3pol}
We now explain that presentations of Lawvere theories can be seen as particular
3-poly\-graphs~\cite{burroni1993higher,fox1976coalgebras}.
This is based on the idea, familiar to people working on linear logic, that a cartesian category is a
monoidal category in which every object can be duplicated and erased, see for
instance~\cite[Section~6]{mellies2009categorical}. For the sake of simplicity,
we consider only strict monoidal categories here, as justified by Mac Lane's
coherence theorem (\cref{thm:moncat-coh}) although many results extend
seamlessly to the general case.

\subsection{Underlying monoidal category}
Suppose given a cartesian category~$C$. It can be equipped with a structure of
symmetric monoidal category. The unit object is the terminal object. The tensor
product of two objects is their cartesian product, and the tensor product of two
morphisms $\phi:u\to u'$ and $\psi:v\to v'$ is the morphism
\[
  \phi\times\psi:u\times v\to u'\times v',
\]
obtained by the universal property of the product:
\[
  \xymatrix@C=3ex@R=3ex{
    &\ar[dl]_{\pi_1}u\times v\ar@{.>}[dd]|-{\phi\times\psi}\ar[dr]^{\pi_2}&\\
    u\ar[dd]_\phi&&v\ar[dd]^\psi\\
    &\ar[dl]^{\pi_1'}u'\times v'\ar[dr]_{\pi_2'}&\\
    u'&&v',
  }
\]
where the morphisms $\pi_1$, $\pi_2$, $\pi_1'$, $\pi_2'$ are the
projections. Finally, the symmetry
\[
  \gamma_{u,v}:v\times u\to u\times v
\]
is defined by
\[
  \xymatrix@C=3ex@R=3ex{
    &v\times u\ar@/_/[ddl]_{\pi_2}\ar@{.>}[d]|-{\gamma_{u,v}}\ar@/^/[ddr]^{\pi_1}&\\
    &\ar[dl]^{\pi_1}u\times v\ar[dr]_{\pi_2}&\\
    u&&v.
  }
\]
In general, this monoidal structure is not strict, but Mac Lane's coherence
theorem ensures that there is no harm in considering it to be strict up to
monoidal equivalence of categories.

A cartesian monoidal category is a symmetric monoidal category which is induced
by a cartesian category as above. We now show that those can be characterized
among symmetric monoidal categories as being those in which every object is
equipped with a structure of commutative comonoid in a natural way.

\subsection{Comonoids}
Suppose given a strict monoidal category $(C,\otimes,\monunit)$.
A \emph{comonoid}\index{comonoid} $(u,\delta,\varepsilon)$ in~$C$ consists of an
object~$u$ together with morphisms
\begin{align*}
  \delta&:u\to u\otimes u,
  &
  \varepsilon&:\monunit\to u,  
\end{align*}
satisfying the usual associativity and unitality axioms
\begin{align*}
  (\delta\otimes\unit u)\circ\delta
  &=
  (\unit u\otimes\delta)\circ\delta,
  &
  (\varepsilon\otimes\unit u)\circ\delta
  &=
  \unit u,
  &
  (\unit u\otimes\varepsilon)\circ\delta
  &=
  \unit u.
\end{align*}
This structure is dual to the one of monoid (see \cref{ex:monoid}).
In the case where the monoidal category is equipped with a symmetry gamma, the
comonoid is \emph{commutative} when it satisfies
\[
  \gamma_{u,u}\circ\delta=\delta
  \pbox.
\]
The following theorem is detailed in various places, \eg
\cite{mellies2009categorical}:

\begin{theorem}
  \label{thm:mon-cart}
  In a symmetric monoidal category $(C,\otimes,\monunit)$, the tensor product is
  a cartesian product if and only if there are natural transformations of
  components
  \begin{align*}
    \delta_u&:u\to u\otimes u,
    &
    \varepsilon_u&:u\to\monunit,
  \end{align*}
  which are monoidal, \ie for every objects $u,v\in C$ we have
\begin{align*}
\xymatrix@C=3ex@R=3ex{
u\otimes v\ar[dr]_-{\delta_u\otimes\delta_v}\ar[rr]^-{\delta_{u\otimes v}}&&\ar[dl]^-{u\otimes \gamma_{u,v}\otimes v}u\otimes v\otimes u\otimes v\\
&u\otimes u\otimes v\otimes v&
}
&&
\xymatrix@C=3ex@R=3ex{
\monunit \ar@/^2ex/[rr]^{\varepsilon_\monunit}\ar@/_2ex/[rr]_{\unit{\monunit}}&&\monunit
}    
\end{align*}
  and such that $(u,\delta_u,\varepsilon_u)$ is a commutative comonoid for every
  object~$u$.
\end{theorem}
\begin{proof}
  Suppose that~$C$ is a cartesian category. Given an object~$u$, the comonoid
  morphisms $\delta_u$ are induced by the universal property of the product:
  \[
    \xymatrix@C=3ex@R=3ex{
      &\ar@/_/[ddl]_{\unit u}u\ar@{.>}[d]|{\delta_u}\ar@/^/[ddr]^{\unit u}&\\
      &\ar[dl]^{\pi_1}u\times u\ar[dr]_{\pi_2}&\\
      u&&u
    }
  \]
  and $\varepsilon_u:u\to 1$ is the terminal morphism. The verification of
  axioms of commutative comonoids and naturality is left to the reader.

  Conversely, suppose that~$C$ is a symmetric monoidal category equipped with
  natural transformations $\delta$ and $\varepsilon$ as in the statement of the
  theorem. For any pair of objects~$u$ and~$v$, we claim that their cartesian
  product is $u\otimes v$ equipped with projections
  \begin{align*}
    \unit{u}\otimes\varepsilon_{v}&:u\otimes v\to u,
    &
    \varepsilon_{v}\otimes\unit{v}&:u\otimes v\to v
    \pbox.    
  \end{align*}
  Given morphisms $\phi:w\to u$ and $\psi:w\to v$, we claim that the universal
  morphism is $\chi=(\phi\otimes\psi)\circ\delta_w$:
  \[
    \xymatrix@C=8ex@R=3ex{
      &\ar@/_/[dddl]_{\phi}w\ar@{.>}[d]|{\delta_w}\ar@/^/[dddr]^{\psi}&\\
      &w\otimes w\ar@{.>}[d]|{\phi\otimes\psi}&\\
      &\ar[dl]^{\unit{u}\otimes\varepsilon_v}u\otimes v\ar[dr]_{\varepsilon_u\otimes\unit{v}}&\\
      u&&v.
    }
  \]
  Namely, we have
  \begin{align*}
    (\unit{u}\otimes\varepsilon_v)\circ(\phi\otimes\psi)\circ\delta_w
    &=
    (\unit{u}\circ\phi)\otimes(\varepsilon_v\circ\psi)\circ\delta_w,
    &&\text{(interchange law)}
    \\
    &=
    (\phi\otimes\varepsilon_w)\circ\delta_w,
    &&\text{(naturality of~$\varepsilon$)}
    \\
    &=
    \phi\circ(\unit{w}\otimes\varepsilon_w)\circ\delta_w,
    &&\text{(interchange law)}
    \\
    &=
    \phi,
    &&\text{(axiom of comonoids)}
  \end{align*}
  so that the triangle on the left commutes, and similarly for the one on the
  right. Conversely, given a morphism~$\chi:w\to u\otimes v$ such that
  $\pi_1\circ\chi=\phi$ and $\pi_2\circ\chi=\psi$, we have
  \begin{align*}
    \chi
    &=
    (((\unit{u}\otimes\varepsilon_{u})\circ\delta_u)\otimes((\varepsilon_{v}\otimes\unit{v})\circ\delta_v))\circ\chi,
    &&\text{(axiom of comonoids)}
    \\
    &=
    (\unit{u}\otimes\varepsilon_{u}\otimes\varepsilon_{v}\otimes\unit{v})\circ(\delta_u\otimes\delta_v)\circ\chi,
    &&\text{(interchange)}
    \\
    &=
    (\unit{u}\otimes\varepsilon_{u}\otimes\varepsilon_{v}\otimes\unit{v})\circ(\unit{u}\otimes\gamma_{u,v}\otimes\unit{v})\circ\delta_{u\otimes v}\circ\chi,
    &&\text{($\delta$ is monoidal)}\\
    &=
    (\unit{u}\otimes\varepsilon_{v}\otimes\varepsilon_{u}\otimes\unit{v})\circ\delta_{u\otimes v}\circ\chi,
    &&\text{($\gamma$ natural)}.
  \end{align*}
  This concludes the proof.
\end{proof}

\begin{remark}
  \nomenclature[MonCat]{$\MonCat$}{category of monoidal categories}
  In a more general way, it can be shown that the forgetful functor
  $\Cart\to\MonCat$ from cartesian categories to monoidal categories admits a
  right adjoint
  \[
    \operatorname{Comon}
    :
    \Cart
    \to
    \MonCat
  \]
  which associates to a monoidal category~$C$, the category of comonoids
  in~$C$. This has been rediscovered many times and can be traced back
  to Fox~\cite{fox1976coalgebras}.
\end{remark}

As a consequence of this theorem, the free cartesian category on a presented
symmetric monoidal category can be presented as follows~\cite{burroni1993higher}:

\begin{theorem}
  \label{thm:pres-free-cart}
  Let $P$ be a 3-polygraph presenting a symmetric monoidal
  category~$\prescat{P}$ (in particular, $P_0=\set{\star}$ is reduced to one
  element). The free cartesian category on~$\pcat{P}$ is presented by the
  3-polygraph~$Q$ such that
  \begin{align*}
    Q_0&= P_0,\\
    Q_1&= P_1,\\
    Q_2&= P_2\sqcup\setof{\delta_a:a\to aa,\varepsilon_a:a\to\unit{}}{a\in P_1},\\
    Q_3&= P_3\sqcup Q_3',
  \end{align*}
  where $Q_3'$ consists of the generators
  \begin{align*}
    A_a&:\delta_a\comp{}\delta_aa\TO\delta_a\comp{}a\delta_a,
    &
    L_a&:\delta_a\comp{}\varepsilon_aa\TO\unit{a},
    &
    D_\alpha&:\alpha\comp{}\delta_v\TO\delta_u\comp{}\alpha\alpha,
    \\
    &&
    R_a&:\delta_a\comp{}a\varepsilon_a\TO\unit{a},
    &
    E_\alpha&:\alpha\comp{}\varepsilon_v\TO\varepsilon_u,
  \end{align*}
  indexed by 1-generators $a$ in $P_1$ and 2-generators $\alpha:u\To v$
  in~$P_2$. Here, given $u\in\freecat{P_1}$, the 2-cells $\delta_u$ and
  $\varepsilon_u$ are defined by induction on~$u$ by
  \begin{align*}
    \delta_\star&=\unit\star,
    &
    \delta_{au}&=\delta_a\delta_u\comp{}a\gamma_{u,a}\comp{}u,
    &
    \varepsilon_\star&=\unit\star,
    &
    \varepsilon_{au}&=\varepsilon_a\varepsilon_u
    \pbox.
  \end{align*}
\end{theorem}

\noindent
Graphically, the morphisms $\delta_u$ and $\varepsilon_u$ can be respectively
depicted as
\begin{align*}
  \satex{delta-u},
  &&
  \satex{eps-u}.  
\end{align*}
and satisfy
\begin{align*}
  \satex{fc-delta-star-l}&=\satex{empty}\,,
  &
  \satex{fc-delta-au-l}&=\satex{fc-delta-au-r},
  &
  \satex{fc-eps-star-l}&=\satex{empty}\,,
  &
  \satex{fc-eps-au-l}&=\satex{fc-eps-au-r}.
\end{align*}
The relations are
\begin{align*}
  \satex{fc-A-l}&\overset{A_a}\TO\satex{fc-A-r},
  &
  \satex{fc-L-l}&\overset{L_a}\TO\satex{fc-L-r},
  &
  \satex{fc-D-l}&\overset{D_\alpha}\TO\satex{fc-D-r},
  \\
  &&
  \satex{fc-R-l}&\overset{R_a}\TO\satex{fc-R-r},
  &
  \satex{fc-E-l}&\overset{E_\alpha}\TO\satex{fc-E-r}.
\end{align*}

More generally, the free cartesian category on a presented monoidal
category~$C$, can be obtained by first presenting the free symmetric monoidal
category on~$C$, see~\cref{sec:free-PROP}, and then applying the above construction.

\subsection{Presentations of Lawvere theories by polygraphs}
The above construction can be used to translate a term rewriting system
presenting a Lawvere theory~$C$ to a polygraph presenting the underlying
monoidal category of~$C$.
First, we can define a functor
\[
  U:\nCPol1\to\nPol2,
\]
which to every signature~$P$ associates the 2-polygraph~$UP$ defined by
\[
  (UP)_0=\set{\star},
  \qquad\qquad
  (UP)_1=P_0,
  \qquad\qquad
  (UP)_2=P_1
  \pbox.
\]
This functor induces an isomorphism between $\nCPol1$ and the full subcategory
of~$\Pol_2$ whose objects are polygraphs with $\star$ as only
$0$-generator. However, note that, given a signature~$P$, the categories
$\freecat{P}$ and $\freecat{(UP)}$ are generally not isomorphic: the former is
cartesian whereas the latter is generally only monoidal. In order to address
this discrepancy, \thmr{mon-cart} suggests that we consider the 2-polygraph
obtained from~$UP$ by formally adding a symmetry, see
\cref{sec:free-PROP,thm:free-smc}, and a natural structure of commutative
comonoid for every object, see \thmr{pres-free-cart}. We thus define a functor
\[
  L:\nCPol1\to\nPol3,
\]
where $LP$ is the polygraph obtained from~$UP$ (seen as a 3-polygraph by the
canonical inclusion $\Pol_2\to\Pol_3$ adding an empty set of 3-generators) by
performing those constructions.

\begin{proposition}
  Given a signature $P\in\nCPol1$, the monoidal category~$\pcat{LP}$ presented
  3-polygraph~$LP$ is the cartesian category~$\freecat{P}$ generated by~$P$.
\end{proposition}

\noindent
As a variant of the above construction, one can show~\cite{burroni1993higher}:

\begin{theorem}
  For every term rewriting system~$P$, there is a 3-polygraph~$Q$ such that
  $\pcat{P}$ is isomorphic to~$\pcat{Q}$ (as monoidal categories). Moreover,
  when~$P$ is finite, the polygraph~$Q$ can also be chosen to be finite.
\end{theorem}

\begin{example}
  We have described, in \cref{ex:trs-group}, a term rewriting system
  corresponding to the theory of groups. By applying the above construction, we
  obtain the following $3$-polygraph~$P$ which presents the same Lawvere theory,
  considered as a monoidal category. We have $P_0=\set\star$, $P_1=\set{a}$
  (thus $\freecat{P_1}\isoto\N$), the 2-generators are those coming from the
  original term rewriting system
  \begin{align*}
    \mu&:2\to 1,
    &
    \eta&:0\to 1,
    &
    \iota&:1\to 1,
  \end{align*}
  respectively pictured as
  \begin{align*}
    \satex{mu}\;,
    &&
    \satex{eta}\;,
    &&
    \satex{sigma}\;,
  \end{align*}
  as well as those corresponding to the cartesian structure
  \begin{align*}
    \delta&:1\to 2,
    &
    \varepsilon&:1\to 0,
    &
    \gamma&:2\to 2,
  \end{align*}
  respectively pictured as
  \begin{align*}
    \satex{delta}\;,
    &&
    \satex{eps}\;,
    &&
    \satex{gamma}\;,
  \end{align*}  
  and the relations are those coming from the term rewriting system
  \begin{align*}
    \satex{mon-unit-l}&\TO\satex{mon-unit-c}\;,
    &
    \satex{mon-unit-r}&\TO\satex{mon-unit-c}\;,
    &
    \satex{mon-assoc-l}&\TO\satex{mon-assoc-r}\;,
    \\
    \mu a\comp{}\mu&\TO\unit{a},
    &
    a \mu\comp{}\mu&\TO\unit{a},
    &
    \mu a\comp{}\mu&\TO a\mu\comp{}\mu,
    \\
    \satex{hopf-l}&\TO\satex{hopf-c}\;,
    &
    \satex{hopf-r}&\TO\satex{hopf-c}\;,
    \\
    \delta\comp{}\iota a\comp{}\mu&\TO\varepsilon,
    &
    \delta\comp{}a \iota\comp{}\mu&\TO\varepsilon,
  \end{align*}
  in addition to those corresponding to the cartesian structure (omitted here).
  Note that the use of~$\delta$ in the two last relations is due to the fact that
  the variable~$x_1$ is used twice in the corresponding relation on terms. It
  turns out that this is the polygraph for cocommutative Hopf algebras, see
  \secr{hopf}.
\end{example}

\begin{example}
  The 3-polygraph~$P$ corresponding to the Lawvere theory of commutative monoids
  is the polygraph of bicommutative bialgebras, see \secr{bialgebra}.
\end{example}

\section{Cartesian Polygraphs}
\index{cartesian!polygraph}
\index{polygraph!cartesian}
\label{sec:cart-pol}
As the notations used in \cref{sec:lt-pres} are meant to suggest, term rewriting
systems can be seen as particular instances of a notion of \emph{cartesian
  polygraph}, adapted to presenting cartesian categories~\cite{cartpol}.
We briefly review this notion here.

\subsection{Cartesian 0-polygraphs}
The category~$\nCPol0$ of cartesian $0$-polygraphs is the category of sets (as
for regular polygraphs, see~\cref{sec:1-pol-cat}).

\subsection{Cartesian 1-polygraphs}
The category $\nCPol1$ of cartesian $1$-polygraphs is the category of
signatures, see \secr{trs-sig}.

\subsection{Cartesian 2-polygraphs}
The category $\nCPol2$ of cartesian $2$-polygraphs is the category of term rewriting
systems, see \secr{trs}. Given a $2$-polygraph~$P$, we write $\tpol1P$ for the
underlying $1$-polygraph.

\subsection{Cartesian 2-categories}
In order to define cartesian $3$-polygraphs, we first need to introduce the following
notion. A $2$-category~$C$ is \emph{cartesian} if its underlying category is
cartesian and for every pair of $2$-cells
\begin{align*}
  F&:\phi\To\phi':w\to u,
  &
  G&:\psi\To\psi':w\to v,
\end{align*}
with same $0$-source (\resp $0$-target), there exists a unique morphism
\[
  \uple{F,G}
  :
  \uple{\phi,\phi'}
  :
  w\to u\times v,
\]
such that $\uple{F,G}\comp0\pi_1=F$ and $\uple{F,G}\comp0\pi_2=G$.
Graphically,
\[
  \xymatrix@C=25ex@R=5ex{
    &w
    \ar@/_4ex/[ddl]_{\phi}
    \ar@/_2ex/@{{}{ }{}}[ddl]|{\overset{F}{\Longrightarrow}}
    \ar[ddl]^{\phi'}
    \ar@/_2.5ex/[d]_>>>{\uple{\phi,\psi}}
    \ar@{}[d]|{\overset{\uple{F,G}}{\Longrightarrow}}
    \ar@/^2.5ex/[d]^>>>{\uple{\phi',\psi'}}
    \ar[ddr]_\psi
    \ar@/^2ex/@{{}{ }{}}[ddr]|{\overset{G}{\Longrightarrow}}
    \ar@/^4ex/[ddr]^{\psi'}
    &\\
    &\ar[dl]^{\pi_1}u\times v\ar[dr]_{\pi_2}&\\
    u&&v\pbox.
  }
\]
\nomenclature[Cart2]{$\nCart2$}{category of cartesian $2$-categories}
We write $\nCart2$ for the category of cartesian $2$-categories, morphisms being
$2$\nbd-func\-tors whose underlying functor is cartesian.

\subsection{Lawvere 2-theories}
\index{Lawvere theory!$2$-}
Given a set~$P_0$, the cartesian category $(\Fun/P_0)^\op$ can canonically be
seen as cartesian 2-category with only identity 2-cells. A \emph{Lawvere
  2-theory} is a cartesian 2-category equipped with a cartesian 2-functor
\[
  (\Fun/P_0)^\op
  \to
  C,
\]
which preserves finite products and is the identity on
objects~\cite{yanofsky2001coherence}.

\subsection{The generated cartesian (2,1)-category}
\index{free!cartesian $(2,1)$-category}
Given a cartesian $2$-poly\-graph~$P$, we write $\freegpd{P}$ for the cartesian
$(2,1)$-category it generates. It has the category $\freecat{\tpol1P}$ freely
generated by the underlying signature (\ie 1-polygraph) as underlying category
and its $2$\nbd-cells are generated under composition and inverses by the
elements of~$P_2$, with source and target indicated by $\src1$ and~$\tgt1$. We
write $\freegpd{P_2}$ for the set of $2$-cells of $\freegpd{P}$.
The cartesian $(2,1)$-category $\freegpd{P}$ is canonically a Lawvere
$2$-theory, with $P_0$ as sorts.

\subsection{Cartesian (3,1)-polygraphs}
A \emph{cartesian $(3,1)$-polygraph} consists of
\begin{itemize}
\item a $2$-polygraph~$P$,
\item a set $P_3$ of $3$-generators together with functions
  $
    \src2,\tgt2
    :
    P_3
    \to
    \freegpd P_2
  $
  such that
  $
  \freecat{\src1}\circ\src2=\freecat{\tgt1}\circ\src2
  $ and $
  \freecat{\src1}\circ\tgt2=\freecat{\tgt1}\circ\tgt2
  $.
\end{itemize}
A \emph{morphism} $f:P\to Q$ of cartesian $(3,1)$-polygraphs consists of a
morphism $\tpol2P\to\tpol2Q$ between the underlying $2$-polygraphs together with
a function $P_3\to Q_3$ which commutes with source and with target.

\subsection{Congruence}
\index{congruence!on a cartesian $2$-category}
A \emph{congruence}~$\approx$ on a cartesian $2$-category~$C$ is a congruence on
the underlying $2$-category such that, for every $1$-cells
\begin{align*}
  F,F'&:\phi\To\phi':w\to u,
  &
  G,G'&:\psi\To\psi':w\to v,
\end{align*}
we have that
\[
  F\approx F'
  \qtand
  G\approx G'
  \qquad\text{implies}\qquad
  \uple{F,G}
  \approx
  \uple{F',G'}
  \pbox.
\]
Given a $3$-polygraph~$P$, the \emph{$P$-congruence} $\approx^P$ is the smallest
congruence such that $F\approx^P G$, for every $3$-generator $\Lambda:F\TO G$.

\subsection{Coherent presentation}
\index{coherent!presentation!of a cartesian category}
\index{presentation!coherent}
A $(3,1)$-polygraph~$P$ is a \emph{coherent presentation} of a cartesian
category~$C$ when $C$ is the cartesian category presented by the underlying
$2$-polygraph, \ie $\pcat{\tpol2P}=C$, and for every parallel $2$-cells
$F,G:\phi\To\psi$ in $\freegpd P_2$ one has $F\approx^PG$.

\subsection{Cartesian Squier homotopical theorem}
\index{Squier!theorem!for cartesian polygraphs}
An analogous of Squier's homotopical theorem
(\cref{thm:SquierHomotopical,Theorem:SquierCompletion3polygraphs}) can be
formulated in this context: the cartesian category presented by a
convergent term rewriting system~$P$ admits a coherent presentation by
the $(3,1)$\nbd-poly\-graph $(P,P_3)$, where~$P_3$
consists of a confluence diagram for every critical branching of~$P$. This can be
used to recover various coherence results (such as Mac
Lane's coherence theorem for monoidal categories) through term rewriting
systems~\cite{beke2011categorification, cohen2009coherence, mimram:LIPIcs.FSCD.2023.16, yanofsky2001coherence}.

\subsection{Homological invariants}
Finally, we shall briefly mention that the homological tools developed for
2-polygraphs in \cref{chap:2homology} can be adapted to term rewriting systems~\cite{Malbos2004,malbos2016homological}. In particular, the homology of a 
Lawvere theory can be used in order to obtain lower bounds on generators and
relations that any presentation should
have~\cite{ikebuchi2021,ikebuchi2022lower,malbos2016homological}.

\part{Polygraphs}

\chapter{Higher Categories}
\label{chap:n-cat}
The remaining chapters of this book present a general theoretical background
underlying all constructions encountered so far. Thus, from this point on, we
shall assume that the reader is well acquainted with the basics of
category theory, as developed in~\cite{MacLane98}.

\bigskip

Among the many existing notions of higher categories, the notion
of strict globular $n$-category that we shall describe is, in some
sense, the most basic one.  The earliest published reference to the
concept appears to be~\cite{brown1981equivalence}, where it is
motivated by the study of higher homotopies. Precisely, one looks here
for higher-dimensional analogues of the fundamental groupoid of a
space. Grothendieck soon afterward realized the need for a weak version of
infinity-groupoids to fulfill this purpose, see~\cite{GrothendieckPstacks,MaltsiniotisWCat}.
The theory has been further developed in the highly influential
paper~\cite{street1987algebra}, advocating the use of strict
higher categories as coefficients for non-abelian cohomology. The same
paper introduces the notion of {\em freely generated $\omega$-category over
a computad}---here called a polygraph---which is central in the
present work, together with the definition of oriented simplices or
orientals. Orientals yield a nerve functor from strict $\omega$-categories to
simplicial sets, turning the former into models of homotopy types, as
developed in~\cite{AraMaltsiThomstr,GagnaStrhot}.

The present work stresses yet another aspect of strict
$\omega$-categories, especially the free ones, as higher-dimensional
rewriting ``spaces'', in the spirit
of~\cite{GuiraudMalbos12advances}.

In this chapter, we set the essential definitions and notations. Starting
with a description of the basic ``shapes'', that is, the presheaf
category $\oGlob$ of globular
sets, we define a family of operations endowing a globular set with a
structure of $\omega$-category. We then prove that the category
$\oCat$ of strict \oo-categories is exactly the category of
algebras of the monad induced by the forgetful functor from $\oCat$ to
$\oGlob$. We finally define important subcategories of $\oCat$ obtained by
requiring cells to be invertible above a given dimension.

\section{Globular Sets}
\label{sec:gset}
  
\subsection{Globes}
\index{globe!category of}
\nomenclature[O]{$\glob$}{category of globes}
We first define the small category $\glob$ of \emph{globes}: its objects are the
integers $0,1,\ldots$ and its morphisms are generated by a double sequence
\[
  \cosce{n},\cotge{n}: n\to n+1,
\]
where $n\in\N$,
\[
  \xymatrix{
    0\doubr{\cosce{0}}{\cotge{0}}& 1 \doubr{\cosce{1}}{\cotge{1}}&
    \cdots\doubr{\cosce{n-1}}{\cotge{n-1}} & n \doubr{\cosce{n}}{\cotge{n}}&
    n+1 \doubr{\cosce{n+1}}{\cotge{n+1}}& \cdots \pbox,
    }
  \]
quotiented by the equations
\begin{align}
  \cosce{n+1}\circ\cosce{n} & = \cotge{n+1}\circ\cosce{n},\label{eq:coglob1}\\
  \cosce{n+1}\circ\cotge{n} & = \cotge{n+1}\circ\cotge{n}.\label{eq:coglob2}
\end{align}
As a consequence of these equations, whenever $0\leq m<n$, the hom-set
$\hom{\glob}{m}{n}$
contains exactly two morphisms:
\begin{align*}
  \cossce mn & = \cosce{n-1}\circ\cdots\circ\cosce m,\\
  \cottge mn & = \cotge{n-1}\circ\cdots\circ\cotge m.
\end{align*}

\index{globular set}
A \emph{globular set} is then a presheaf on $\glob$, that is, a functor
$X:\opp{\glob}\to\Set$.
Thus, a globular set $X$ amounts to a sequence of sets
$X(n)$ of $n$-dimensional globes, for each $n\geq 0$, together with source and
target maps
\begin{displaymath}
  X(\cosce{n}),X(\cotge{n}): X(n+1)\to X(n)
\end{displaymath}
satisfying the {\em globular relations}, dual to~(\ref{eq:coglob1}) and~(\ref{eq:coglob2}).
Let us denote
$X(\cosce{n})$ by $\sce n$ and $X(\cotge{n})$ by $\tge n$. Whenever $m \le n$, we set
\begin{align*}
  \ssce{m}{n} & = \sce{m}\circ\cdots\circ\sce{n-1},\\
  \ttge{m}{n} & = \tge{m}\circ\cdots\circ\tge{n-1},
\end{align*}
thus $\ssce mn,\ttge mn:X(n)\to X(m)$.

\nomenclature[Glob8]{$\nGlob n$}{category of $n$-globular sets}
\nomenclature[Glob9]{$\oGlob$}{category of globular sets}
Globular sets and natural transformations between them define a category denoted
by $\oGlob$. For any globular set $X$ and integer $n$, we shall denote $X(n)$ by
$X_n$. The elements of $X_n$ are called \emph{$n$-cells}\index{cell}.
For any $n$-cell $x$ and $m \le n$, the notations $\sce{m}(x)$ and $\tge{m}(x)$
will stand for $\ssce{m}{n}(x)$ and $\ttge{m}{n}(x)$, respectively. The
$m$\nbd-cell $\sce{m}(x)$ will be called the
\ndef{$m$-source}\index{source!of an $n$-cell} of $x$ (or simply the
\ndef{source}\index{source} if $m = n - 1$) and the $m$-cell $\tge{m}(x)$
the \ndef{$m$-target}\index{target!of an $n$-cell} of $x$ (or simply the \ndef{target} if
$m = n - 1$). Two $n$-cells $x$, $y$ are called
\ndef{parallel}\index{parallel!8-cells@$n$-cells} if either $n=0$ or 
$\sce{n-1}(x)=\sce{n-1}(y)$ and~$\tge{n-1}(x)=\tge{n-1}(y)$ otherwise.

\nomenclature[On]{$\glob^{(n)}$}{category of globes of dimension $\le n$}
Let us denote $\glob^{(n)}$ the full subcategory of $\glob$ whose objects are
$0,\ldots,n$. Then the category $\nGlob{n}$ of \emph{$n$-globular
sets}\index{globular set!8-@$n$-}\index{8-globular set@$n$-globular set} is by
definition the category of presheaves on $\glob^{(n)}$. 
If $0\leq m<n$, the canonical inclusion
\begin{displaymath}
  \glob^{(m)}\to\glob^{(n)}
\end{displaymath}
gives rise, by precomposition, to a truncation functor
\begin{displaymath}
  \trunc^n_m:\nGlob n\to \nGlob m.
\end{displaymath}
Also, for each $n\geq 0$, the canonical inclusion $\glob^{(n)}\to \glob$ gives
rise to a truncation functor
\begin{displaymath}
  \trunc_n:\oGlob\to \nGlob n
\end{displaymath}
making all the following triangles commute:
\begin{displaymath}
  \begin{xy}
    \xymatrix{
      \oGlob \ar[d]_{\trunc_n}\ar[rd]^{\trunc_m}& \\
      \nGlob n \ar[r]_{\trunc_m^n}& \nGlob m\pbox.
    }
  \end{xy}
\end{displaymath}
For any $n$-globular set $X$, and $m<n$, the notation $\trunc_m(X)$
will stand for $\trunc_m^n(X)$.
Remark that $\oGlob$ is the projective limit of the diagram
\begin{displaymath}
\xymatrix{
\nGlob 0 
&\nGlob 1\ar[l]_-{\trunc_0^1} 
&\cdots\ar[l]
& \nGlob {n-1} \ar[l]
& \nGlob n \ar[l]_-{\trunc_{n-1}^{n}}
&\cdots\ar[l]\pbox.}
\end{displaymath}

\subsection{Globes and spheres}
As for any presheaf category, we get a Yoneda embedding
\begin{displaymath}
  \yoneda:\glob\to\oGlob
\end{displaymath}
defined on objects by $Y(m)(n)=\glob(m,n)$.
\nomenclature[On]{$\globe n$}{$n$-globe (as a globular set)}
\nomenclature[On2]{$\sphere n$}{$n$-sphere (as a globular set)}
For each $n\geq 0$, we call \emph{$n$-globe}\index{globe} and denote by $\globe n$
the representable globular set $\yoneda(n)$. The $n$-globe has exactly two
$i$-cells in dimensions $0\leq i<n$, one $n$-cell, and no $i$-cell for
$i>n$. The sub-globular set of $\globe n$ having the same cells as
$\globe n$ in all dimensions $i\neq n$ and no $n$-cell will be called
the \emph{$n$-sphere}\index{sphere}, and denoted by $\sphere n$. Remark that
$\sphere 0$ is the initial globular set with no cells at all. Globes
and spheres come with a family of canonical
inclusion morphisms
\begin{displaymath}
  \gencof n:\sphere n\to\globe n
\end{displaymath}
which we shall encounter in numerous occasions. For example, the case
$n=2$ may be pictured as
\[
 \{
    \xymatrix@C=2pc{
      \bullet
      \ar@/^2ex/[r]_{}="0"
      \ar@/_2ex/[r]_{}="1"
      &
      \bullet
    }
    \}
    \into
    \{
    \xymatrix@C=2pc{
      \bullet
      \ar@/^2ex/[r]_{}="0"
      \ar@/_2ex/[r]_{}="1"
      \ar@2"0";"1"
      &
      \bullet
    }
    \}.
\]

\section{Strict \pdfm{n}-Categories}
\index{8-category@$n$-category}
\index{category!8-@$n$-}
\index{9-category@$\omega$-category}
\index{category!9-@$\omega$-}
\nomenclature[Cat8]{$\nCat n$}{category of $n$-categories}
\nomenclature[Cat9]{$\ooCat$}{category of $\omega$-categories}

\subsection{Definition}
A \emph{strict $\omega$-category} is given by a globular set $C$ together
with a family of partial binary composition operations $(\comp i)_{i\in\N}$ and identity
operations $(\Unit{i}{})_{i\in\N\setminus\set{0}}$  subject to the
following conditions:
\begin{itemize}
\item if $0\leq i<k$ and $x$, $y$ are $k$-cells such that
  $\tge i(x)=\sce i(y)$ (in which case we say that $x$ and $y$ are {\em
  $i$-composable}\index{composable cells}) there is a $k$-cell $x\comp i y$,
\item if $k>0$ and $x$ is a $(k - 1)$-cell, there is a $k$-cell $\Unit{k}{x}$,
  and more generally, if $i \ge 0$ and $x$ is an $i$-cell, we may define
  recursively on $k>i$ a $k$-cell $\Unit{k}{x}$
  by~$\Unit{k}{x}=\Unit{k}{\Unit{k-1}{x}}$.
\end{itemize}
Compositions and units are subject to:
\begin{enumerate}
\item positional conditions prescribing the source and target of composites and
  units, namely
  \begin{itemize}
  \item if $0\leq i<j$, then
    \[
    \sce j(x\comp i y)=\sce j(x)\comp i \sce j(y)
    \qand
    \tge j(x\comp i y)=\tge j(x)\comp i \tge j(y), 
   \]
  \item if $0\leq j\leq i$, then
    \[
      \sce j (x\comp i y)=\sce j (x) \qand \tge j(x\comp i y)=\tge j (y),
    \]
  \item if $0\leq i<k$ and $x$ is an $i$-cell, then
    \[ \sce i(\Unit{k}{x})= x = \tge i(\Unit{k}{x}), \]
  \end{itemize}
\item computational conditions of
  \begin{itemize}
  \item associativity: if $i<k$ and $x$, $y$, $z$ are $k$-cells such that
    $\tge i(x)=\sce i(y)$ and $\tge i(y)=\sce i(z)$, then 
    \[ 
      (x\comp i y)\comp i z= x\comp i (y\comp i z),
    \]
  \item neutrality of units: if $0\leq i<k$ and $x$ is a $k$-cell, then
    \[ \Unit{k}{\sce i(x)}\comp i x= x\comp i \Unit{k}{\tge i(x)}= x, \]
  \item exchange: if $i<j<k$ and $x$, $y$, $z$, $v$ are $k$-cells such that
    $\tge j(x)=\sce j(y)$, $\tge j(z)=\sce j(v)$ and $\tge i(x)=\sce i(z)$, then
    also $\tge j(y)=\sce j(v)$, and
    \[
      (x\comp j y)\comp{i} (z\comp j v)=(x\comp i z) \comp j (y\comp i v),
    \]
 \item compatibility of units: if $0\leq i<j<k$ and $x$, $y$ are $i$-composable
    $j$-cells, then
    \[
      \Unit{k}{x\comp i y}=\Unit{k}{x}\comp i\Unit{k}{y}.
    \]
  \end{itemize}
\end{enumerate}

\index{9-functor@$\omega$-functor}
\index{functor!9-@$\omega$-}
Let $C$, $D$ be two strict $\omega$-categories. An {\em
$\omega$-functor} $f:C\to D$ is a morphism of the underlying globular sets
which preserves the compositions and units.  Strict $\omega$-categories and
strict $\omega$-functors build a (large) category we denote by $\oCat$. If
we restrict the above construction to cells of dimension at most $n$,
we get the category of strict $n$-categories, denoted by $\nCat{n}$.

From now on we will drop the adjective ``strict'' and we will speak of
``\oo-categories'' and ``$n$-categories'' when we mean ``strict
\oo-categories'' and ``strict $n$-categories''. 

\begin{remark}
 The structure of $n$-category is sometimes presented by the alternative set
 of operations and axioms described
 in~\cref{chap:syntactic_descr}.
\end{remark}

\subsection{\pdfoo-categories as models of a projective sketch}
\label{paragr:ocatdiagrams}
The above axioms for \oo-categories can be presented in diagrammatic
form as follows. 

For any globular set $X$ and $0\leq i<n$, there is a pullback square
 in $\Set$:
 \begin{displaymath}
   \begin{xy}
     \xymatrix{\plb{X_n}{X_n}{X_i}\ar[d]_-{\lproj_i^n}\ar[r]^-{\rproj_i^n}& X_n\ar[d]^{\ssce{i}{n}}\\
                      X_n\ar[r]_{\ttge{i}{n}} & X_i\pbox.}
   \end{xy}
 \end{displaymath}
The operations of compositions and units become maps:
\begin{displaymath}
 \comp i:\plb{X_n}{X_n}{X_i}\to X_n
\end{displaymath}
and
\begin{displaymath}
  \UNIT{n}: X_i\to X_n
\end{displaymath}
for $0\leq i<n$.

The positional conditions for compositions amount to the commutation of the following
diagrams:
\begin{displaymath}
   \begin{xy}
     \xymatrix{\plb{X_n}{X_n}{X_i}\ar[d]_{\plb{\ssce{j}{n}}{\ssce{j}{n}}{X_i}}\ar[r]^-{\comp{i}}& X_n\ar[d]^{\ssce{j}{n}}\\
                     \plb{X_j}{X_j}{X_i} \ar[r]_-{\comp i} & X_i}
   \end{xy}
\qquad\qquad
\begin{xy}
     \xymatrix{\plb{X_n}{X_n}{X_i}\ar[d]_{\plb{\ttge{j}{n}}{\ttge{j}{n}}{X_i}}\ar[r]^-{\comp{i}}& X_n\ar[d]^{\ttge{j}{n}}\\
                     \plb{X_j}{X_j}{X_i} \ar[r]_-{\comp i} & X_i}
   \end{xy}
 \end{displaymath}
for $0\leq i<j<n$ and
 \begin{displaymath}
   \begin{xy}
     \xymatrix{\plb{X_n}{X_n}{X_i}\ar[d]_{\lproj_i^n}\ar[r]^-{\comp{i}}& X_n\ar[d]^{\ssce{j}{n}}\\
                      X_n\ar[r]_{\ssce{j}{n}} & X_j}
   \end{xy}
\qquad\qquad
\begin{xy}
     \xymatrix{\plb{X_n}{X_n}{X_i}\ar[d]_{\rproj_i^n}\ar[r]^-{\comp{i}}& X_n\ar[d]^{\ttge{j}{n}}\\
                      X_n\ar[r]_{\ttge{j}{n}} & X_j}
   \end{xy}
 \end{displaymath}
for $0\leq j\leq i < n$.

As for units, if $0\leq i<n$,  the positional conditions amount to the commutations of
\begin{displaymath}
  \begin{xy}
    \xymatrix{X_n\ar[rd]_{\ssce{i}{n}} & X_i \ar[d]^{\id}\ar[l]_{\UNIT{n}}\ar[r]^{\UNIT{n}}& X_n\ar[ld]^{\ttge{i}{n}}\\
                            & X_i & \pbox.}
  \end{xy}
\end{displaymath}

Now, each axiom is expressed by the commutation of a diagram in $\Set$
involving arrows derived from the source, target, compositions, and
unit arrows by means of universal constructions. 
\begin{itemize}
\item Associativity of compositions amounts to the commutation of
\begin{displaymath}
  \begin{xy}
    \xymatrix{\plb{X_n}{(\plb{X_n}{X_n}{X_i})}{X_i}\ar[dd]_{\alpha}
      \ar[rr]^(.6){\plb{X_n}{(\comp{i})}{X_i}}&&\plb{X_n}{X_n}{X_i}\ar[dr]^{\comp{i}} & \\
   & & & X_n\\
\plb{(\plb{X_n}{X_n}{X_i})}{X_n}{X_i} \ar[rr]_(.6){\plb{(\comp{i})}{X_n}{X_i}}&& \plb{X_n}{X_n}{X_i}\ar[ur]_{\comp{i}} & }
  \end{xy}
\end{displaymath}
where $\alpha$ is the canonical bijection between both pullbacks.
\item Let $0\leq i<j<n$. We build the diagram
  \begin{displaymath}
    \begin{xy}
      \xymatrix{\plb{(\plb{X_n}{X_n}{X_j})}{(\plb{X_n}{X_n}{X_j})}{X_i}\ar[rr]^-{\rproj}\ar[dd]_{\lproj}\ar@{.>}[rd]^-{\lambda}&
      & \plb{X_n}{X_n}{X_j}\ar[d]^{\lproj_j^n}\\
&\plb{X_n}{X_n}{X_i}\ar[r]^-{\rproj_i^n}\ar[d]_{\lproj_i^n}&X_n\ar[d]^{\ssce{n}{i}}\\
\plb{X_n}{X_n}{X_j}\ar[r]_-{\lproj_j^n}&X_n\ar[r]_{\ttge{n}{i}}&X_i}
    \end{xy}
  \end{displaymath}
where both solid squares are pullbacks. There is a unique universal
arrow~$\lambda$ making the whole diagram commute. Similarly, by
replacing the left projection $\lproj_j^n$ with the right projection
$\rproj_j^n$ in the above diagram, we get a universal arrow
\begin{displaymath}
  \rho:\plb{(\plb{X_n}{X_n}{X_j})}{(\plb{X_n}{X_n}{X_j})}{X_i} \to
  \plb{X_n}{X_n}{X_i}.
\end{displaymath}
Consider now the diagram
\begin{displaymath}
    \begin{xy}
      \xymatrix{\plb{(\plb{X_n}{X_n}{X_j})}{(\plb{X_n}{X_n}{X_j})}{X_i}\ar[rr]^-{\rho}\ar[dd]_{\lambda}\ar@{.>}[rd]^-{\theta}&
      & \plb{X_n}{X_n}{X_i}\ar[d]^{\comp{i}}\\
&\plb{X_n}{X_n}{X_j}\ar[r]^-{\rproj_j^n}\ar[d]_{\lproj_j^n}&X_n\ar[d]^{\ssce{n}{j}}\\
\plb{X_n}{X_n}{X_i}\ar[r]_-{\comp{i}}&X_n\ar[r]_{\ttge{n}{j}}&X_j\pbox.}
    \end{xy}
  \end{displaymath}
The positional conditions on compositions and the globular relations
ensure that the outer square commutes, whereas the small solid square
is a pullback by definition. Therefore, we get a unique universal
arrow $\theta$ as shown in the diagram.

Now the exchange rule amounts to the commutation of
  \begin{displaymath}
    \begin{xy}
      \xymatrix{\plb{(\plb{X_n}{X_n}{X_j})}{(\plb{X_n}{X_n}{X_j})}{X_i}\ar[r]^-{\theta}\ar[d]_{\plb{(\comp{j})}{(\comp{j})}{X_i}}  & \plb{X_n}{X_n}{X_j}\ar[d]^{\comp{j}}\\
                        \plb{X_n}{X_n}{X_i} \ar[r]_-{\comp{i}}& X_n\pbox.}
    \end{xy}
  \end{displaymath}
\item Let $0\leq i<n$. The positional conditions on units imply that
  the following diagram of solid arrows commutes:
  \begin{displaymath}
    \begin{xy}
      \xymatrix{X_n\ar[d]_{\ssce{i}{n}}\ar[rrd]^{\id}\ar@{.>}[rd]_{\iota}& & \\
                       X_i \ar[rd]_{\UNIT{n}}& \plb{X_n}{X_n}{X_i}\ar[r]_-{\rproj_i^n}\ar[d]_{\lproj_i^n} & X_n\ar[d]^{\ssce{i}{n}}\\
                             & X_n\ar[r]_{\ttge{i}{n}} & X_i\pbox.}
    \end{xy}
  \end{displaymath}
Therefore, there is a unique universal arrow $\iota$ making the whole
diagram commute. Now the first axiom for left units amounts to the fact
that $\iota$ equalizes the pair $\pair{\rproj_i^n}{\comp{i}}$, that is, the
commutation of
\begin{displaymath}
  \begin{xy}
    \xymatrix{X_n\ar[r]^(.4){\iota} &\plb{X_n}{X_n}{X_i}\ar@<2pt>[r]^(.6){\rproj_i^n}\ar@<-2pt>[r]_(.6){\comp{i}} & X_n\pbox.}
  \end{xy}
\end{displaymath}
The first axiom for right units is treated similarly.
\item Finally, the compatibility of compositions with units amounts to
  the commutation of the following diagram
  \begin{displaymath}
    \begin{xy}
      \xymatrix{\plb{X_m}{X_m}{X_i}\ar[d]_{\plb{\UNIT{n}}{\UNIT{n}}{X_i}}\ar[r]^-{\comp{i}} & X_m\ar[d]^{\UNIT{n}}\\
                          \plb{X_n}{X_n}{X_i}\ar[r]^-{\comp{i}} & X_n}
    \end{xy}
  \end{displaymath}
whenever $0\leq i<m<n$.

\end{itemize}

\begin{proposition}\label{prop:ocatlimsketch}
  The category $\oCat$ is a category of models of a projective sketch.
\end{proposition}

\begin{proof}
The above diagrammatic presentation of the axioms of
$\omega$-categories defines a sketch whose models are actual $\omega$-categories. 
\end{proof}

\begin{corollary}\label{coro:ooCat_bicomp}
The category $\oCat$ is complete and cocomplete.
\end{corollary}

\begin{proof}
  This follows from the previous result by using
  Proposition~\ref{prop:mod_bicompl}.
\end{proof}

\subsection{Truncation functors}
\nomenclature[Unm]{$\trunc^n_m$}{truncation functor for $n$-categories}
\nomenclature[Um]{$\trunc_m$}{truncation functor for $\omega$-categories}
As for globular sets, whenever $0\leq m < n$, we get \ndef{truncation
functors}\index{truncation functor}\index{functor!truncation}
\begin{displaymath}
  \trunc^n_m:\nCat n\to\nCat m
\end{displaymath}
and
\begin{displaymath}
  \trunc_n:\oCat\to\nCat n
\end{displaymath}
making all triangles 
\begin{displaymath}
  \begin{xy}
    \xymatrix{\oCat \ar[d]_{\trunc_n}\ar[rd]^{\trunc_m}& \\
     \nCat n \ar[r]_{\trunc_m^n}& \nCat m}
  \end{xy}
\end{displaymath}
commute.

Here again, $\oCat$ appears as the projective limit of the diagram
\begin{displaymath}
  \xymatrix{\nCat 0 & \nCat 1 \ar[l]_{\trunc_0^1}& \cdots \ar[l]& \nCat n\ar[l] & \nCat{n+1} \ar[l]_{\trunc_n^{n+1}}&\cdots \ar[l]\pbox{.}}
\end{displaymath}
Remark that, by abuse of language, we use the same notation for
truncation functors among $n$-categories and among $n$-globular
sets.

\nomenclature[V]{$\fgf$}{forgetful functor from $\omega$-categories to globular sets}
\nomenclature[Vn]{$\fgf_n$}{forgetful functor from $n$-categories to $n$-globular sets}
By construction, $\omega$-categories are globular sets with structure,
whence a forgetful functor
\begin{displaymath}
  \fgf:\oCat\to\oGlob
\end{displaymath}
which restricts for each $n\in\N$ to
\begin{displaymath}
  \fgf_n:\nCat{n}\to\nGlob{n}.
\end{displaymath}
These forgetful functors commute with the above truncation functors,
that is, the following diagram commutes whenever $0\leq m<n$:
\begin{equation}
  \label{diag:trunc-fgf-commute}
  \vcenter{
  \begin{xy}
    \xymatrix{\nCat{n} \ar[r]^{\trunc_m^n}\ar[d]_{\fgf_n}& \nCat{m}\ar[d]^{\fgf{m}}\\
                      \nGlob{n} \ar[r]_{\trunc_m^n}&\nGlob{m}\pbox.}
  \end{xy}
  }
\end{equation}

\begin{prop}\label{prop:catglob}
  The forgetful functor $\fgf:\oCat\to\oGlob$ admits
  a left adjoint
  $\freecatg:\oGlob\to\oCat$. Likewise, for each $n\in\N$,
  $\fgf_n:\nCat{n}\to\nGlob{n}$ admits a left adjoint $\freecatg_n:\nGlob{n}\to\nCat{n}$.
  \end{prop}

\begin{proof}
The categories $\oGlob$ and $\oCat$ are categories of models of
projective sketches $S$ and~$S'$, respectively,
whereas the functor $\fgf$ is the one induced on models by the inclusion morphism
$S\hookto S'$. Therefore $\fgf$ admits a
left adjoint (see~Theorem~\ref{thm:sketch-morphism}). The same arguments hold for $\nGlob{n}$ and $\nCat{n}$, where~$n\in\N$.
\end{proof}

\begin{proposition}\label{prop:filcolim}
  The forgetful functor $\fgf : \oCat \to \oGlob$ preserves filtered colimits.
\end{proposition}
\begin{proof}
 This is a general property of functors induced by a morphism of
 projective sketches involving only finite cones (see~\cite[Chapter 4,
 Theorem 4.4]{barr1985toposes}).
 Concretely, the colimit $C$ of a filtered
 diagram in $\oCat$ is obtained by taking the colimit $X$ of the
 underlying diagram in $\oGlob$ and defining a structure of
 \oo-category on $X$ in the obvious way, so that $\fgf(C)=X$. Thus,
 the stronger statement that $\fgf$ {\em creates} filtered colimits holds.
\end{proof}

\begin{paragr}[Globes and spheres]
\label{paragraph:GlobesSpheres}
\nomenclature[On]{$\globe n$}{$n$-globe (as an $\omega$-category)}
\nomenclature[On2]{$\sphere n$}{$n$-sphere (as an $\omega$-category)}
By abuse of notation, the free $\omega$-category $\freecatg(\globe n)$
generated by $\globe n$ will be still denoted by $\globe n$ and called the
\emph{$n$-globe}\index{globe}. Likewise
we denote by $\sphere n$ the free $\omega$-category $\freecatg(\sphere
n)$, and call it the \emph{$n$-sphere}\index{sphere}. 
Remark that, in the case of globes and spheres, the free functor only
adds new identity cells to the ones already present in the globular
globes and spheres. 

For any $\omega$-category $C$, the set of \emph{$n$-globes of $C$} is
the hom-set $\hom{\oCat}{\globe n}{C}$, which amounts to the set $C_n$
of $n$-cells in $C$. Likewise, the set of \emph{$n$-spheres of~$C$}
is the hom-set $\hom{\oCat}{\sphere{n+1}}{C}$, which amounts to the set
of pairs of parallel $n$-cells in $C$, that is,  parallel in the
underlying globular set. 
\end{paragr}

\section{Basic Examples}\label{sec:basic_ex}

Let us first mention a few immediate examples of \oo-categories.
\begin{itemize}
\item Sets: as $\nCat 0=\nGlob 0=\Set$ and $\nCat n$ naturally embeds in $\oCat$
(see~\cref{paragr:adjunctions} below), any set $S$ can be viewed as an
\oo-category, precisely the \oo-category~$C$ whose $0$-cells are the
elements of $S$ and whose $n$-cells, $n>0$,  are all of the form $\Unit{n}{x}$ for
$x\in S$.
\item Monoids: any monoid $M$, being a $1$-category with a unique
object, can be seen as an \oo-category whose $n$-cells are
identities for all $n\geq 2$.
\item Commutative monoids: to any {\em commutative} monoid $\pair{A}{+}$ we may associate the
  \oo-category $C$  defined by $C_0=\set{\star}$,
  $C_1=\set{\Unit{1}{\star}}$, $C_2=A$ and having only identity cells in higher
  dimensions. The source and target maps are uniquely determined,
  whereas compositions are defined for any pair $\pair uv$ of
  $2$-cells by $u\comp 0 v=u\comp 1 v=u+v$. The axioms of
  \oo-categories are easily checked. Conversely, for any
  \oo-category $C$ such that $C_0=\set{\star}$,
  $C_1=\set{\Unit{1}{\star}}$ and $C_n$ has only identity cells for
  $n>2$, the $\comp 0$ and $\comp 1$ compositions on $C_2$ coincide
  and are commutative operations, so that $\pair{C_2}{\comp 0}$ (or
  $\pair{C_2}{\comp 1}$) becomes an abelian monoid.
\end{itemize}

\section{More Properties of \pdfm\oCat}

\subsection{Local presentability}
\label{paragr:oCat_loc_pres}
By definition, the category $\oGlob$ of globular sets is a category of
presheaves, so that limits and colimits are computed pointwise. 
On the other hand, we already noticed that $\oCat$ 
is the category of models of a projective sketch and is hence 
complete and cocomplete (see~\cref{coro:ooCat_bicomp}). The forgetful
functor $\fgf:\oCat\to\oGlob$, being a right adjoint, preserves limits. Thus
limits in $\oCat$ are computed as in $\oGlob$. Colimits however are hard to
compute, even in $\nCat 1$. 

\subsection{Enrichment}
For each $n\geq 0$, the category $\nCat n$ has a monoidal structure defined by
its cartesian product and terminal object. Thus, the notion of
$\nCat n$-enriched category makes sense: in fact, $\nCat n$-enriched categories
are just $(n + 1)$-categories. As for $n=\omega$, it turns out that $\oCat$ is
enriched over itself.

 \begin{prop}\label{prop:catglobmonad}
   The forgetful functor $\fgf:\oCat\to \oGlob$ is monadic.
 \end{prop}
 \begin{proof}
First remark that, for any $0\leq i<n$, 
the correspondence
\begin{displaymath}
  X\mapsto \plb{X_n}{X_n}{X_i}
\end{displaymath}
is functorial from $\oGlob$ to $\Set$. Indeed, the Yoneda embedding yields a
pushout
\begin{displaymath}
  \begin{xy}
    \xymatrix{\globe i \ar[r]^{\yoneda(\cottge{i}{n})}\ar[d]_{\yoneda(\cossce{i}{n})}&\globe n\ar[d]\\
\globe n \ar[r]& \psh{\globe n}{\globe n}{\globe i}}
  \end{xy}
\end{displaymath}
in $\oGlob$ and we get a natural bijection
\begin{displaymath}
  \plb{X_n}{X_n}{X_i} \isoto \hom{\oGlob}{\psh{\globe n}{\globe n}{\globe i}}{X}
\end{displaymath}
so that our correspondence is the object part of the (representable) functor
\begin{displaymath}
  \hom{\oGlob}{\psh{\globe n}{\globe n}{\globe i}}{-}:\oGlob\to\Set.
\end{displaymath}

   Now, by Proposition~\ref{prop:catglob}, the functor $\fgf$ admits a
   left adjoint $\freecatg$. By Beck's monadicity
   theorem~\cite{MacLane98}, it is sufficient to prove that
   $\fgf$ creates absolute coequalizers.
Thus, let $C$, $D$ be
   $\omega$-categories, $f,g:C\to D$ a pair of morphisms, $X=\fgf C$,
   $Y=\fgf D$ the underlying globular sets. Suppose
   $w:Y\to Z$ is an absolute coequalizer of the pair $u=\fgf f$,
   $v=\fgf g$. We must
   prove the existence of a unique $\omega$-category $E$, and a unique
   morphism $h:D\to E$ such that
   \begin{itemize}
   \item $\fgf E=Z$,
   \item $\fgf h=w$
  \end{itemize}
and check that $h$ is the coequalizer of the pair $\pair fg$ in
$\oCat$.

Let us first define a structure of $\omega$-category on the globular
set $Z$.
By hypothesis, the diagram
\begin{equation}
    \xymatrix{X\doubr{u}{v}& Y\ar[r]^{w} & Z}
\label{eq:coeq}
\end{equation}
is an absolute coequalizer. Thus, by applying the functor
$\hom{\oGlob}{\psh{\globe n}{\globe n}{\globe i}}{-}$ defined in the
preliminary remark, we
still get a coequalizer, now in $\Set$:
\begin{displaymath}
  \begin{xy}
    \xymatrix{\plb{X_n}{X_n}{X_i} \ar@<2pt>[r]^{u'}\ar@<-2pt>[r]_{v'}& \plb{Y_n}{Y_n}{Y_i}\ar[r]^{w'} & \plb{Z_n}{Z_n}{Z_i}\pbox.}
  \end{xy}
\end{displaymath}
Likewise by applying the functor $\hom{\oGlob}{\globe n}{-}$, we get
another coequalizer diagram
\begin{displaymath}
  \begin{xy}
    \xymatrix{X_n \doubr{u_n}{v_n}& Y_n\ar[r]^{w_n} & Z_n\pbox.}
  \end{xy}
\end{displaymath}
Consider now the following diagram:
\begin{displaymath}
  \begin{xy}
    \xymatrix{\plb{X_n}{X_n}{X_i} \doubr{u'}{v'}\ar[d]_{\comp{i}}&
      \plb{Y_n}{Y_n}{Y_i}\ar[r]^{w'} \ar[d]_{\comp{i}}& \plb{Z_n}{Z_n}{Z_i}\ar@{.>}[d]_{?}\\
X_n \doubr{u_n}{v_n}& Y_n\ar[r]^{w_n} & Z_n\pbox.
}
  \end{xy}
\end{displaymath}
As $u$ and $v$ come from morphisms in $\oCat$, the squares on $u'$,
$u_n$, and $v'$, $v_n$ commute. Therefore, $w_n\circ(\comp{i})$
coequalizes the pair $\pair{u'}{v'}$, and there exists a unique map
from $\plb{Z_n}{Z_n}{Z_i}$ to $Z_n$
making the right-hand square commute. This map defines the
$i$-composition among $n$-cells in $Z$, still denoted by $\comp i$ . By a similar argument we may
define the unit map $\UNIT{n}:Z_{n-1}\to Z_n$. 

It remains to check that the compositions and units just defined on
$Z$ satisfy the axioms of $\omega$-categories. This amounts to check
the commutation of all diagrams expressing these axioms. We shall
treat the axiom of associativity in detail and leave the remaining
axioms as exercises. By applying appropriate functors to the
coequalizer diagram~(\ref{eq:coeq}), we get the following diagram in
$\Set$:
\begin{displaymath}
  \fig[width=\textwidth]{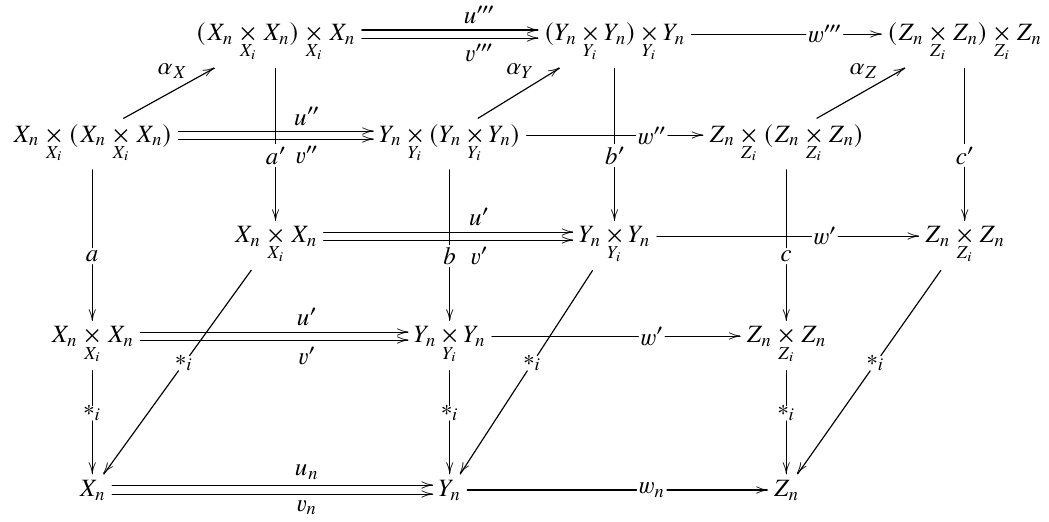}
\end{displaymath}
where $a=\plb{X_n}{(\comp i)}{X_i}$, $a'=\plb{(\comp i)}{X_n}{X_i}$, and $b$, $b'$,
$c$, $c'$ are defined accordingly.

Because the coequalizer~(\ref{eq:coeq}) is absolute, all horizontal
lines are also coequalizer diagrams. Now $w_n\circ (\comp i)\circ b$
coequalizes the pair $\pair{u''}{v''}$, therefore 
\begin{equation}
  (\comp i)\circ c\circ w''= w_n\circ (\comp i)\circ b.
\label{eq:assoc1}
\end{equation}
On the other hand, by naturality of $\alpha$,
\begin{align*}
  (\comp i)\circ c' \circ \alpha_Z \circ w''& = (\comp i)\circ c'
  \circ w'''\circ \alpha_Y \\
 & = (\comp i)\circ w'\circ b'\circ \alpha_Y\\
& = w_n\circ (\comp i)\circ b'\circ \alpha_Y.
\end{align*}
and by the associativity of $\comp i$ on $Y$, the latter expression is
equal to $w_n\circ (\comp i)\circ b$ so that we get 
\begin{equation}
  \label{eq:assoc2}
  (\comp i)\circ c' \circ \alpha_Z \circ w''= w_n\circ (\comp i)\circ b.
\end{equation}
Now as $w''$ is a coequalizer,~\eqref{eq:assoc1} and~\eqref{eq:assoc2}
imply
\begin{displaymath}
  (\comp i)\circ c'\circ \alpha_Z=(\comp i)\circ c.
\end{displaymath}
This is the commutation of the rightmost pentagon and $\comp i$ is
associative on~$Z$.
The other axioms are proved in a similar way. Thus we have defined an
$\omega$\nbd-category $E$ with underlying globular set $Z=\fgf E$.  

Now the commutation of the square involving $w_n$ and $w'$
in the above diagram expresses the preservation of the composition $\comp i$
by $w$. The preservation of units holds by a similar argument. Therefore, we
get a unique morphism $h:D\to E$ in $\oCat$ such that $\fgf h=w$.

Finally, we must show that the diagram
\begin{displaymath}
  \xymatrix{C \doubr{f}{g}& D \ar[r]^{h}& E}
\end{displaymath}
is itself a coequalizer in $\oCat$. Thus, let $K$ be an
$\omega$-category and $k:D\to K$ a morphism such that $kf=kg$. 
We have
to prove the existence of a unique morphism $\ell:E\to K$ such that
$\ell h=k$. Now, if $T=\fgf K$ and $t=\fgf k$, there is a unique
morphism $s:Z\to T$ such that $sw=t$:
\begin{displaymath}
  \xymatrix{X\doubr{u}{v} & Y\ar[r]^{w}\ar[rd]_{t} & Z\ar[d]^{s}\\
                      & & T \pbox.}
\end{displaymath}
 We look for an $\ell$ such that
$\fgf \ell=s$. Uniqueness is obvious, and existence reduces to the
observation that $s$ preserves compositions and units. As for
compositions, consider the following diagram:
\begin{displaymath}
  \xymatrix{\plb{Y_n}{Y_n}{Y_i}\ar[r]^{w'}\ar[d]_{\comp i} &\plb{Z_n}{Z_n}{Z_i}\ar[r]^{s'}\ar[d]_{\comp i} & \plb{T_n}{T_n}{T_i}\ar[d]_{\comp i}\\
                    Y_n \ar[r]_{w_n}& Z_n\ar[r]_{s_n}&T_n \pbox.}
\end{displaymath}
The left-hand square commutes because $w=\fgf h$, and the outer square
commutes because $sw=t$ and $t=\fgf k$. But $w'$ is a coequalizer
map, whence the right-hand square also commutes: this shows that $s$
preserves compositions, as required. A similar argument proves the
preservation of units, and we get the unique morphism $\ell:E\to K$ such that
$s=\fgf \ell$.
\end{proof}

\begin{remark}
 Instead of using Beck's criterion to prove the monadicity of~$\oCat$ over
 $\oGlob$, we could have used a less known criterion in terms of
 sketches due to Lair. Precisely, globular sets are models
 of a projective sketch $S$ with underlying category $\opp{\glob}$
 and no cones, whereas $\omega$-categories are models of a projective sketch
 $S'$, with an obvious sketch inclusion $S\hookto S'$, inducing the
 forgetful functor $\fgf:\oCat\to\oGlob$ between the corresponding
 categories of models. Now (i) the base of each cone of $S'$ already belongs to $S$ and (ii) each object of
 $S'$ not in $S$ is the tip of at least one cone of $S'$. By
 Theorem~\ref{thm:Lair}, these two conditions ensure the monadicity of
 $\fgf$.
 \end{remark}

\begin{remark}
Recall that monadicity is not transitive. For example, the forgetful
functor $\Cat\to \Graph$ is monadic, as well as the forgetful
functor $\Graph\to \Set^2$ taking the graph $\xymatrix{X_0  &
  \doubl{}{} X_1}$ to
the pair $\pair{X_0}{X_1}$. However the composite $\Cat\to \Set^2$ is
\emph{not} monadic. Consider in fact the following categories:
\begin{itemize}
\item the category $C$, freely generated
on the graph having a set of five vertices $V_C=\set{0,1,2,3,4}$ and a
set of two edges $E_C=\set{a:0\to 1;b:3\to 4}$,
\item the subcategory $D$ of $C$ obtained by removing the isolated
  vertex $2$.
\end{itemize}
 Define two morphisms $u,v:C\to D$ such that $u$ and $v$ are
 both retractions of the inclusion $D\to C$, $u(2)=1$ and $v(2)=3$. We
 leave it as an exercise to check that the forgetful functor $\Cat\to
 \Set^2$ takes the pair $\pair uv$ to a pair $\pair fg$ whose
 coequalizer $e$ in $\Set^2$ is split, but the coequalizer $w$ of
 $\pair uv$ in $\Cat$ is not sent to $e$. The generating graphs for
 $C$, $D$ and of the coequalizer $E$ of the pair $\pair uv$ are
 represented in the picture below:
 \begin{center}
   \begin{tikzpicture}
     \clip (-3.5,2) rectangle (7,4.5);
     \draw (-3.12,4.2) node[anchor=north west] {$\bullet^0\to\bullet^1$};
     \draw (-2.23,3.8) node[anchor=north west] {$\bullet^2 $};
     \draw (-2.23,3.4) node[anchor=north west] {$\bullet^3\to\bullet^4$};
     \draw (0.9,4.2) node[anchor=north west] {$\bullet^0\to\bullet^1 $};
     \draw (1.8,3.4) node[anchor=north west] {$\bullet^3\to\bullet^4$};
     \draw (4.52,3.8) node[anchor=north west] {$\bullet^0\to\bullet\to\bullet^4 $};
     \draw [dotted] (-3.36,4.28) rectangle (-0.58,2.62);
     \draw [dotted] (0.55,4.28) rectangle (3.33,2.62);
     \draw [dotted] (4.19,4.28) rectangle (6.97,2.62);
     \draw (-0.36,4) node[anchor=north west] {$\overset{u}{\to}$};
     \draw (-0.36,3.5) node[anchor=north west] {$\underset{v}{\to}$};
     \draw (3.46,3.9) node[anchor=north west] {$\overset{w}{\to}$};
     \draw (-2.18,2.5) node[anchor=north west] {$C$};
     \draw (1.66,2.5) node[anchor=north west] {$D$};
     \draw (5.18,2.5) node[anchor=north west] {$E$};
   \end{tikzpicture}
 \end{center}
\end{remark}

\begin{paragr}[Adjunctions]
  \label{paragr:adjunctions}
  \label{sec:adjunctions}
  If $0\leq m<n$, there is a canonical inclusion functor
  \begin{displaymath}
   \incl_n^m:\nCat{m} \to \nCat{n}
  \end{displaymath}
taking an $m$-category $C$ to the $n$-category $D=\incl_n^m(C)$ such
that $\trunc_m^n(D)=C$ and all $i$-cells of $D$ are units for $m<i\leq
n$. Likewise, we get a canonical inclusion
\begin{displaymath}
   \incl^m:\nCat{m} \to \oCat.
  \end{displaymath}
Therefore, each $m$-category may be naturally identified with an $n$-category for
any $m<n\leq\omega$. 

The functor $\incl_n^m$ (\resp $\incl^m$) has a
right adjoint, namely the truncation functor~$\trunc_m^n$ (\resp $\trunc_m$). 

Now $\incl_n^m$  also admits a left adjoint
\begin{displaymath}
  \truncbis_m^n:\nCat{n}\to\nCat{m}.
\end{displaymath}
Let $C$ be an $n$-category and
$D=\truncbis_m^n(C)$. Up to dimension $m - 1$, $D$ coincides with
$\trunc_{m-1}^n(C)$, whereas $D_m$ is the quotient of $C_m$ modulo the
congruence generated by $C_{m+1}$. Precisely, two parallel $m$-cells
$x$, $y$ in $C_m$ are congruent modulo $C_{m+1}$ if and only if there
is a sequence $x_0=x,x_1,\ldots,x_p=y$ of $m$-cells of $C_m$ and a
sequence $z_1,\ldots,z_p$ of $(m + 1)$-cells in $C_{m+1}$ such that, for
each $i=1,\ldots,p$, either $z_i=x_{i-1}\to x_i$ or $z_i:x_i\to
x_{i-1}$. Note that the source and target maps, as well as
compositions on $C_m$ are compatible with the congruence
relation. Therefore, $D$ is a well-defined $m$-category, as
expected. Also, the action of~$\truncbis_m^n$ on morphisms is
immediate, and clearly functorial. Likewise, $\incl^m$ admits a
left adjoint $\truncbis_m:\oCat\to \nCat{m}$.

Let us finally remark that the truncation functor $\trunc_m^n$
(resp. $\trunc_n$) also
admits a right adjoint $\inclbis_n^m:\nCat{m}\to \nCat{n}$
(resp. $\inclbis^m:\nCat{m}\to\oCat$): let $C$ be an $m$-category,
$D=\inclbis_n^m(C)$ is the $n$-category such that (i) $\trunc_m^n(D)=C$,
(ii)  for each pair $x,y$ of parallel $m$-cells in $D_m$, there is
exactly one $(m + 1)$-cell $z:x\to y$ in $D_{m + 1}$, and (iii) all
$i$-cells of $D$ are units whenever $i>m+1$. The
$\omega$-category~$\inclbis^m(C)$ is defined accordingly.
To sum up, omitting the indices, we get a series of adjunctions
between inclusions and truncation functors:
\begin{displaymath}
  \truncbis \quad\dashv\quad \incl \quad\dashv\quad \trunc \quad\dashv\quad
  \inclbis.
\end{displaymath}

\end{paragr}

\section{\pdfm{(n,p)}-Categories}
\label{sec:npcat}

\subsection{Invertible cells}
\label{subsec:inverses}
\index{inverse}
\index{invertible cell}
\index{cell!invertible}
Let $0\leq i<k\leq n\leq\omega$. Let $C$ be an $n$-category and $u$ a
$k$-cell of $C$. A $k$-cell $v$ of $C$ is a {\em $\comp i$-inverse to
  $u$} if $v$ is left and right $i$-composable with $u$ and $u\comp i
v=\Unit{k}{\sce{i}(u)}$ and $v\comp i u=\Unit{k}{\tge{i}(u)}$. If such
a $k$-cell $v$ exists, it is necessarily unique. In that case, we call
$u$ an {\em $\comp i$-invertible} cell. A $k$-cell $u$ is called
simply {\em invertible} if it is $\comp{k-1}$-invertible. 

\begin{lemma}\label{lemma:invertible_2_cells}
  If a $2$-cell is  $\comp 0$-invertible, then it is also $\comp 1$-invertible.
\end{lemma}
\begin{proof}
  Let $u$ be a $\comp 0$-invertible $2$-cell and $v$ its $\comp
  0$-inverse, so that $u\comp 0 v=\Unit{2}{\sce{0}(u)}$ and $v\comp 0
  u=\Unit{2}{\tge{0}(u)}$. This implies that the $1$-cells
  $\sce{1}(u)$ and $\tge{1}(u)$ are $\comp 0$-invertible, with
  $\sce{1}(v)$ and $\tge{1}(v)$ as respective $\comp 0$-inverses. Let
  \[
    v'=\Unit{2}{\tge{1}(u)}\comp 0 v\comp 0\Unit{2}{\sce{1}(u)}.
  \]
  We claim that $v'$ is a $\comp 1$-inverse to $u$. In fact,
  \[
    \sce{1}(v')=\tge{1}(u)\comp 0 \sce{1}(v)\comp 0\sce{1}(u)=\tge{1}(u)
  \]
  so that $u$ and $v'$ are $\comp 1$-composable. Moreover, by using
  the exchange rule
  \begin{eqnarray*}
    u\comp 1 v' & = & (u\comp 0 \Unit{2}{\sce{1}(v)} \comp 0
                      \Unit{2}{\sce{1}(u)})\comp 1
                      (\Unit{2}{\tge{1}(u)}  \comp 0 v\comp
                      0\Unit{2}{\sce{1}(u)})\\
                & = & (u\comp 1 \Unit{2}{\tge{1}(u)})\comp 0 (
                      \Unit{2}{\sce{1}(v)}\comp 1 v)\comp 0 \Unit{2}{\sce{1}(u)}\\
                & = & u\comp 0 v\comp 0 \Unit{2}{\sce{1}(u)}\\
                & = & \Unit{2}{\sce{1}(u)}.
  \end{eqnarray*}
 Likewise, one checks that $v'\comp 1  u=\Unit{2}{\tge{1}(u)}$.
\end{proof}

\begin{corollary}\label{corollary:invertible_cells}
  Let $i\leq j<k$. Each $\comp i$-invertible $k$-cell
  is also $\comp j$-invertible. 
\end{corollary}
\begin{proof}
Let $i\leq j<k\leq n$, and let $C$ be an $n$-category. We define a
$2$-category~$D$ by $D_0=C_i$, $D_1=C_j$, $D_2=C_k$ with the obvious
source and target maps, units and compositions induced by $C$. The
statement then immediately follows by applying
Lemma~\ref{lemma:invertible_2_cells} to the $2$-cells of $D$.
\end{proof}

\begin{definition} 
  \index{8p-category@$(n,p)$-category}
  \index{9p-category@$(\omega,p)$-category}
  \index{category!8p-@$(n,p)$-}
  \index{category!9p-@$(\omega,p)$-}
  \nomenclature[Cat8p]{$\npCat np$}{category of $(n,p)$-categories}
  \index{8-groupoid@$n$-groupoid}
  \index{9-groupoid@$\omega$-groupoid}
  \index{groupoid!8-@$n$-}
  \index{groupoid!9-@$\omega$-}
  \nomenclature[Gpdn]{$\nGpd n$}{category of $n$-groupoids}
 Let $0\leq p\leq n\leq \omega$. The category $\npCat np$ is the full
 subcategory of $\nCat{n}$ having as objects the $n$-categories whose
 $k$-cells are invertible for all $k>p$. These objects are called {\em
   $\pair np$-categories}. In particular, the objects of
 $\npCat n0$ are the {\em $n$-groupoids}, where $k$-cells are
 invertible for all~$k>0$. We also denote $\npCat n0$ by $\nGpd n$. 
\end{definition}

\begin{prop}\label{prop:invertible_cells_in_catnp}
  Let $p\leq i<k\leq n$ and let $C$ be an $(n,p)$-category. Each $k$-cell
  of $C$ is $\comp i$-invertible.
\end{prop}
\begin{proof}
As in the proof of Corollary~\ref{corollary:invertible_cells}, the
statement reduces to the fact that, in each $2$-category $C$ all whose
   $1$-cells are $\comp 0$-invertible, any $2$-cell $u$ is $\comp
   1$-invertible if and only if it is $\comp 0$-invertible. The ``if''
   direction follows from
   \cref{lemma:invertible_2_cells}. Conversely, suppose that all
   $1$-cells are invertible, and let $u$ be a $\comp 1$-invertible
   $2$-cell, with $\comp 1$-inverse $v$. By hypothesis, $\sce{1}(u)$
   and $\tge{1}(u)$ have $\comp 0$\nbd-inverses $v_1^{-}$ and $v_1^{+}$,
   respectively. One easily checks that the required $\comp 0$-inverse
   to $u$ is given by
   \[
     v'=\Unit{2}{v_1^{+}}\comp 0 v \comp 0 \Unit{2}{v_1^{-}}
     \]
(see also~\cite[§1.3.]{AraMetWGrp}).
\end{proof}

\begin{prop}\label{prop:npcat_properties}
  For each $n$ and $p$ such that $0\leq p\leq n\leq \omega$, the
  category~$\npCat np$ is complete, cocomplete, and monadic over $\nGlob
  n$. Moreover, the inclusion functor $\npCat np\to \nCat n$ admits a left adjoint.
\end{prop}
\begin{proof}
  Like $\nCat n$, the category $\npCat np$ is the category of models of a
  projective sketch $S_{n,p}$, hence it is complete and cocomplete. The
  inclusion functor $\npCat np\to\nCat n$ is induced by the morphism of
  corresponding sketches $S_n\to S_{n,p}$, hence
  admits a left adjoint. Finally the monadicity of~$\npCat np$ over~$\nGlob
  n$ is proved as in \cref{prop:catglobmonad}.
\end{proof}


\chapter{Polygraphs}
\label{chap:polygraphs}

The notion of $2$-polygraph, already introduced in \cref{chap:2pol},
first appears in~\cite{street1976limits} under the name of {\em
  computad}, as an essential tool in proving the existence of limits
in $2$-categories. Although its relevance to rewriting theory was recognized by
Eilenberg and Street from the very beginning~\cite{eilenberg1986rewrite}, this point of view is not
explicitly mentioned in the literature until the early 1990s. The
general notion of $n$\nbd-compu\-tad explicitly appears in~\cite{power1991n}, and
independently
in~\cite{burroni1991higher} and~\cite{burroni1993higher} under the
name of \emph{polygraph}. We adopt here Burroni's presentation and terminology.
The source of Burroni's approach can be traced back in his work on graphical
algebras~\cite{burroni1981algrap}, where he presents a ``concept of
dimension in formal languages''.
Let us mention that~\cite{batanin1998computads} introduces a wide generalization of
the notion of $n$-computad attached to a finitary monad $T$ on globular
sets, presented in more detail in \cref{chap:gen-pol}.
We deal in this chapter with the particular case where the monad $T$ comes
from the adjunction between~$\oCat$ and~$\oGlob$.

\section{Main Definitions}

Throughout this section, we denote by $n$ a natural number.

\subsection{Cellular extensions}
\label{sec:cell-ext}
\index{cellular extension}
Given an $n$-category~$C$, a \emph{cellular extension of $C$} is a family
\[
(X_i:\sphere{n+1}\to C)_{i\in I}
\]
of $n$\nbd-spheres in~$C$ indexed by a set~$I$. This amounts to a
family of pairs of parallel $n$-cells in $C$. Equivalently, it can also be seen as an \oo-functor
\[
  X:\coprod_{i\in I}\sphere{n+1}\to C.
\]
Note that, in order for $X$
to make sense, we identify the $n$-category $C$ with
its image in $\oCat$ by the inclusion functor defined in~\cref{sec:adjunctions}. 
A morphism
\[
f
:
(C,(X_i)_{i\in I})
\to
(D,(Y_j)_{j\in J})
\]
between two cellular extensions of $n$-categories consists of a pair $\pair{g}{h}$,
where $g:C\to D$ is a morphism in $\nCat{n}$ and $h:I\to J$ is a map
such that, for each $i\in I$, $g\circ X_i=Y_{h(i)}$.
\nomenclature[Catnp]{$\nCatp n$}{category of $n$-categories with a cellular extension}
We write $\nCatp n$ for the resulting category. More abstractly, the category
$\nCatp n$ is the pullback of $\nCat n$ and $\nGlob{n+1}$ over $\nGlob n$ in
$\CAT$\nomenclature[CAT]{$\CAT$}{category of possibly large categories} (which denotes the
category of possibly large categories and functors):
\begin{equation}
  \label{eq:ncap-pullback}
  \vxym{
    \nCatp n\ar@{.>}[d]\ar@{.>}[r]&\nGlob{n+1}\ar[d]^{\trunc_n^{n+1}}\\
    \nCat n\ar[r]_{\fgf_n}&\nGlob n\pbox.
  }
\end{equation}
In the above diagram, the forgetful functor $\fgf_n$ and the truncation functor
$\trunc^{n + 1}_{n}$ are those defined in \chapr{n-cat}, and the forgetful functor from $\nCatp n$ to $\nCat n$ takes a cellular
extension $\pair CX$ to the $n$-category $C$. Finally, the horizontal
dotted arrow takes a cellular extension $\pair CX$ to the
$(n+1)$-globular set extending $\fgf_n(C)$ with the set of $(n+1)$-cells
determined by $X$.

\subsection{Freely generated category}
\label{sec:cell-ext-free-cat}
Consider now the forgetful functor
\[
\fgfplus_n:\nCat{n+1}\to\nCatp n
\]
which to an $(n+1)$-category~$C$ associates the pair $(\trunc^{n+1}_n(C),(X_x)_{x\in C_{n+1}})$
where for each cell $x\in C_{n+1}$,
$X_x$ is the $n$-sphere $(\sce n(x),\tge n(x))$. This functor $\fgfplus_n$ is in fact the universal arrow from $\nCat{n+1}$ to
$\nCatp n$ resulting from the commutation of the diagram
\begin{displaymath}
  \begin{xy}
    \xymatrix{\nCat{n+1}\ar[r]^{\fgf_{n+1}}\ar[d]_{\trunc^{n+1}_{n}} & \nGlob{n+1}\ar[d]^{\trunc^{n+1}_{n}}\\
  \nCat n\ar[r]_{\fgf_n}& \nGlob n}
  \end{xy}
\end{displaymath}
and the pullback property of~\eqref{eq:ncap-pullback}.
\begin{proposition}\label{prop:ncatp-leftadjoint}
  The functor $\fgfplus_n$ admits a left adjoint 
$\freeplus_n$:
\[
\vxym{
  \nCat{n+1}\ar@/_5ex/[ddr]\ar@/^5ex/[drr]^-{\fgf_{n+1}}\ar@/_3ex/[dr]_-{\fgfplus_n}\ar@{}[dr]|{\rotatebox[origin=c]{-45}{$\bot$}}&\\
  &\ar@/_3ex/[ul]_-{\freeplus_n}\nCatp n\ar[d]\ar[r]&\nGlob{n+1}\ar[d]^-{\trunc_n^{n+1}}\\
  &\nCat n\ar[r]_-{\fgf_n}&\nGlob n \pbox.
}
\]
\end{proposition}

\begin{proof}
  Let $C$ be an $n$-category and $(X_i:\sphere{n+1}\to C)_{i\in
  I}$    a cellular extension, so that $\pair CX$ is an object of
$\nCatp n$. By considering $C$, $\sphere{n+1}$ and $\globe{n+1}$ as
$(n + 1)$\nbd-categories, and remembering that $\nCat{n+1}$ is cocomplete, we define $\freeplus_n\pair CX$
as the pushout given by  the following diagram in $\nCat{n+1}$:
\begin{displaymath}
  \vxym{\coprod_{i\in I}\sphere{n+1}\ar[r]^-X\ar[d]_{\coprod_{i\in I} \gencof{n+1}} & C\ar[d]\\
\coprod_{i\in I}\globe{n+1}\ar[r] & \freeplus_n\pair CX\pbox.
}
\end{displaymath}
This construction yields a functor $\freeplus_n: \nCatp{n} \to \nCat{n+1}$.
We claim that $\freeplus_n$ is left adjoint to $\fgfplus_n$. Consider in fact $\pair
CX$ an object of $\nCatp{n}$, $D$ an $(n + 1)$\nbd-category, and
$f: \pair CX\to \fgfplus_n(D)$
a morphism in $\nCatp{n}$. Recall that $f$ is a pair $\pair gh$ where
$g:C\to \trunc_n(D)$ is an $n$-functor and $h:X\to D_{n+1}$ is a map
preserving the globular structure. This amounts to a commutative
diagram in~$\nCat{n+1}$ of the form
\begin{displaymath}
  \vxym{\coprod_{i\in I}\sphere{n+1}\ar[r]^-X\ar[d]_{\coprod_{i\in I} \gencof{n+1}}
    & C\ar[d]^{\tilde g}\\
\coprod_{i\in I}\globe{n+1}\ar[r]_-{\tilde h} & D\pbox,
}
\end{displaymath}
where $\tilde g$ and $\tilde h$ are the $(n+1)$-functors built from
$g$ and $h$, respectively. The pushout property then gives a unique
morphism
$f^*: \freeplus_n\pair CX\to D$
making the following diagram commutative:
\begin{displaymath}
  \vxym{\coprod_{i\in I}\sphere{n+1}\ar[r]^X\ar[d]_{\coprod_{i\in I} \gencof{n+1}}
    & C\ar[d]\ar@/^2ex/[rdd]^{\tilde g}& \\
\coprod_{i\in I}\globe{n+1}\ar[r] \ar@/_2ex/[rrd]_{\tilde h}& \freeplus_n\pair CX\ar[rd]^{f^*}&\\
& & D\pbox.
}
\end{displaymath}
The correspondence $f\mapsto f^*$ is then a natural isomorphism
\[
  \xymatrix{\Hom{\nCatp{n}}{\pair CX}{\fgfplus_n (D)} \ar[r]^{\isoto}&\Hom{\nCat{n+1}}{\freeplus_n\pair CX}{D}\pbox{,}}
  \]
which ends the proof.
\end{proof}

\begin{remark}
 There are several approaches to the construction of the above functor
 $\freeplus_n$. On the abstract side, its existence comes from the fact that
 $\fgfplus_n$ is a limit and filtered colimit preserving functor between
 locally presentable categories,
 see~\cite[14.6]{gabrielulmer1971lokprc}, \cite[Theorem~4.1]{barr1985toposes},
 or~\cite[1.66]{adamek1994localp}.
 More concretely, a purely syntactic construction of $\freeplus_n$ based on a type
 system is given in~\cite{metayer2008cofibrant}. An alternative
 construction is given in \cref{chap:syntactic_descr}.
\end{remark}

  Given a cellular extension $\pair{C}{X}$ of an $n$-category
$C$, we call $\freeplus_n\pair{C}{X}$ the \emph{freely
  generated $(n + 1)$-category} on this extension and denote it by
$C[X]$.
In practice, the universal property of $C[X]$ will be used by applying
the following lemma.

\begin{lemma}\label{lemma:ext-universal}
  Let $\pair{C}{(X_i)_{i\in I}}$ be a cellular extension of an
  $n$-category~$C$.
  For every $(n+1)$-category~$D$, every morphism
  $g:C\to\trunc_n^{n+1}(D)$ of $\nCat{n}$ and every map
 $h:I\to D_{n+1}$ making $\pair{g}{h}$ a morphism of
  cellular extensions, there exists a unique morphism $\overline{g}: C[X]\to D$
 in $\nCat{n+1}$ such that $\fgfplus_n(\overline{g})=\pair{g}{h}$.
\end{lemma}
\begin{proof}
 The statement of the lemma is a mere rephrasing of the fact that
 $\freeplus_n$ is left adjoint to $\fgfplus_n$.   
\end{proof}

\subsection{Quotient category}
\label{sec:quotient-category}
\index{quotient!8-category@$n$-category}
For any cellular extension $\pair CX$ of an $n$-cate\-gory~$C$, the \emph{quotient}
$n$-category $C/X$ is the $n$-category obtained from~$C$ by identifying $n$-cells under the
smallest congruence (\wrt compositions and identities) containing all
spheres $\pair xy$ in~$X$.
More precisely,
$C/X$ is nothing but~$\truncbis_n^{n+1}\!(C[X])$ were $\truncbis$ is the
truncation functor already described in \secr{adjunctions}. 

\subsection{\pdfm{n}-polygraphs}
\label{subsec:npolyg}
\index{polygraph!8-@$n$-}
\index{8-polygraph@$n$-polygraph}
\nomenclature[Poln]{$\nPol n$}{category of $n$-polygraphs}
\nomenclature[FP]{$\freecatpol(P)$}{free $\omega$-category on a polygraph $P$}
We already gave an explicit description of $n$-polygraphs for
$0\leq n \leq 3$ in previous chapters. We now turn to the general
definition. Thus, the category $\nPol n$ of $n$-polygraphs is defined by
induction on $n$, together with a functor $\freecatpol_n:\nPol
n\to\nCat n$.
\begin{itemize}
\item The category $\nPol 0$ is $\Set=\nCat 0$ and
  $\freecatpol_0=\id$.
\item Given $\nPol n$ and $\freecatpol_n:\nPol n\to \nCat n$, the
  category $\nPol{n+1}$ is defined by the following pullback in $\CAT$
  \begin{displaymath}
    \begin{xy}
      \xymatrix{\nPol{n+1}\ar[r]^{J_n}\ar[d]_{\trunc^{n+1}_{n}} & \nCatp n \ar[d]\\
                       \nPol n\ar[r]_{\freecatpol_n} & \nCat n\pbox,}
    \end{xy}
  \end{displaymath}
whereas $\freecatpol_{n+1}$ is $\freeplus_nJ_n$:
\begin{displaymath}
  \xymatrix{\nPol{n+1}\ar[r]_{J_n}\ar@/^3ex/[rr]^{\freecatpol_{n+1}} & \nCatp n \ar[r]_{\freeplus_n}& \nCat{n+1}\text{.}}
\end{displaymath}
\end{itemize}
More explicitly, an $(n + 1)$-polygraph $P^{(n+1)}$ is a pair $\pair{P^{(n)}}{X}$
where $P^{(n)}$ is an $n$-polygraph and $X$ is a cellular extension of
the $n$-category $\freecatpol_n(P^{(n)})$.  The first projection of
the pullback yields a truncation functor $\nPol{n+1}\to \nPol n$ we
still denote by $\trunc_n^{n+1}$ as in the case of $n$-globular sets and
$n$-categories.  Note also that the following square commutes:
\begin{displaymath}
  \begin{xy}
    \xymatrix{\nPol n\ar[d]_{\freecatpol_n} &\ar[l]_{\trunc_{n}^{n+1}} \ar[d]^{\freecatpol_{n+1}}\nPol{n+1}\\
                      \nCat n & \ar[l]^{\trunc_{n}^{n+1}}\nCat{n+1}\pbox.}
  \end{xy}
\end{displaymath}
For each $n$-polygraph $P$ and $k<n$, the polygraph
$\trunc_k^{k+1}\circ\cdots\circ\trunc_{n-1}^n (P)$ will be denoted by $\truncpol{P}{k}$.
\nomenclature[P1]{$\truncpol{P}{k}$}{underlying $k$-polygraph of a polygraph}
\index{underlying!8-polygraph@$n$-polygraph}
%
Thus, the data defining an $n$-polygraph $P$ may be displayed in the
following diagram in $\Set$:
\begin{displaymath}
  \xymatrix@C=6ex{
  P_0\ar[d]^{\ins0}&\ar@<-.5ex>[dl]_(.4){\sce
    0}\ar@<.5ex>[dl]^(.4){\tge{0}}P_1\ar[d]^{\ins1}&\ar@<-.5ex>[dl]_(.4){\sce
    1}\ar@<.5ex>[dl]^(.4){\tge{1}}P_2\ar[d]^{\ins2}&\ar@<-.5ex>[dl]_(.4){\sce 2}\ar@<.5ex>[dl]^(.4){\tge{2}}\ar@{}[r]|{\textstyle\ldots}&&\ar@<-.5ex>[dl]_(.4){\sce{n-2}}\ar@<.5ex>[dl]^(.4){\tge{n-2}}P_{n-1}\ar[d]^{\ins{n-1}}&\ar@<-.5ex>[dl]_(.4){\sce{n-1}}\ar@<.5ex>[dl]^(.4){\tge{n-1}}P_{n}\\
  P_0^*&\ar@<-.5ex>[l]_-{\ssce{0}{*}}\ar@<.5ex>[l]^-{\ttge{0}{*}}P_1^*&\ar@<-.5ex>[l]_-{\ssce{1}{*}}\ar@<.5ex>[l]^-{\ttge{1}{*}}P_2^*&\ar@<-.5ex>[l]_-{\ssce{2}{*}}\ar@<.5ex>[l]^-{\ttge{2}{*}}\ar@{}[r]|{\textstyle\ldots}&{\phantom{P_{n-2}^*}}&\ar@<-.5ex>[l]_-{\ssce{n-2}{*}}\ar@<.5ex>[l]^-{\ttge{n-2}{*}}P_{n-1}^*\\
}
\end{displaymath}
where, for each $0\leq k<n-1$, $\ssce{k}{*}\circ \ins{k+1}=\sce{k}$ and
$\ttge{k}{*}\circ \ins{k+1}=\tge{k}$.
Also, at each level $k\leq n$, $k>0$, we add a new set $P_k$ of {\em $k$-generators}\index{generator!8-@$n$-}
together with source and target maps
\begin{displaymath}
  \sce{k-1},\tge{k-1} : P_{k}\to P_{k-1}^*
\end{displaymath}
such that
 \begin{equation}
  \label{eq:pol-globular}
  \ssce{k-2}{*}\circ \sce{k-1}=\ssce{k-2}{*}\circ \tge{k-1}
  \qquad\qquad
  \ttge{k-2}{*}\circ \tge{k-1}=\ttge{k-2}{*}\circ \tge{k-1}
\end{equation}
for $k>2$,
thus defining a $k$-cellular extension of the $(k-1)$-category $\freecat{\truncpol{P}{k-1}}$
\begin{displaymath}
 \xymatrix@C=8ex{
  P_0^*&\ar@<-.5ex>[l]_-{\sce0^*}\ar@<.5ex>[l]^-{\tge0^*}P_1^*&\ar@<-.5ex>[l]_-{\sce1^*}\ar@<.5ex>[l]^-{\tge1^*}P_2^*&\ar@<-.5ex>[l]_-{\sce2^*}\ar@<.5ex>[l]^-{\tge2^*}\ar@{}[r]|{\textstyle\ldots}&{\phantom{P_{-k2}^*}}&\ar@<-.5ex>[l]_-{\sce{k-2}^*}\ar@<.5ex>[l]^-{\tge{k-2}^*}P_{k-1}^*\\
}
\end{displaymath}
already defined at the previous level. The vertical map $\ins k$ is the natural
inclusion of the set $P_k$ of $k$-generators into the set $P_k^*$  of
$k$-cells of the free $k$-category on the cellular
extension of $\freecat{\truncpol{P}{k-1}}$ by $\sce{k-1},\tge{k-1} : P_{k}\to
P_{k-1}^*$. Formally, $\ins k$ is the $k$-dimensional component of the
unit of the monad $\fgfplus_{k-1}\freeplus_{k-1}$ on this cellular extension.
Concretely, the cells in $P_k^*$ are all formal compositions of 
$k$-generators in $P_k$ and units on cells in $P_{k-1}^*$, quotiented
by the axioms of $k$-categories.

\begin{prop}\label{prop:generating_cells}
  Let $P$ be an $n$-polygraph and let $0<k<n$. Suppose that $A\subset \free
  P_k$ is a subset of $\free P_k$ such that
  \begin{itemize}
  \item $\ins k(P_n)\subset A$,
   \item $A$ contains all cells of the form $\Unit{k}{u}$ for $u\in
     \free P_{k-1}$,
    \item for each $i\in\set{0,...,k-1}$ and each pair $u,v\in A$ of $\comp
      i$-composable cells,  $A$ contains their $\comp i$-composition $u\comp i v$.
    \end{itemize}
    Then $A=\free P_k$.
  \end{prop}
  \begin{proof}
    Let $C$ be the $k$-category whose $(k - 1)$-truncation is $P_{\leq
      k-1}^*$ and such that $C_k=A$, with source and target maps $A\to
    \free P_{k-1}$  given by the restriction of
    $\sce{k-1},\tge{k-1}:\free P_k\to \free P_{k-1}$. By hypothesis,
    $C$ is a sub-$k$-category of $\free P_{\leq k}$, and there is an
    inclusion morphism $j:C\to \free P_{\leq k}$. On the other hand,
    the inclusion $\ins k:P_k\to A$ determines a unique morphism
    $f:\free P_{\leq k}\to C$ by
    Lemma~\ref{lemma:ext-universal}. Now $j\circ f:\free P_{\leq
      k}\to \free P_{\leq k}$ is the identity on $\free P_{\leq k}$, by the
    uniqueness property from Lemma~\ref{lemma:ext-universal}. This
    implies that $j_k\circ f_k: \free P_k\to \free P_k$ is the
    identity on $\free P_k$, whence $j_k$ is surjective. Therefore
    $A=\free P_k$. 
  \end{proof}

  \subsection{Structural induction}
\label{subsec:structural_induction}
Proposition~\ref{prop:generating_cells} allows 
reasoning by structural induction on the cells of freely generated
$\omega$-categories. Precisely, in order to prove that a certain
property $A$  holds for all cells in $\free P_k$, it suffices to check
that
\begin{itemize}
\item $A$ holds for all units $\unit v$ where $v\in \free P_{k-1}$,
 \item $A$ holds for all generating cells of the form $\ins k (a)$ for
   $a\in P_k$,
 \item whenever $A$ holds for two $i$-composable cells $u$
   and $v$, then $A$ holds for $w=u\comp i v$.  
 \end{itemize}
 
\subsection{\pdfoo-polygraphs}
\label{subsec:omegapol}
\index{polygraph!9-@$\omega$-}
\index{9-polygraph@$\omega$-polygraph}
\nomenclature[Polw]{$\oPol$}{category of $\omega$-polygraphs}
\nomenclature[Um]{$\trunc_m$}{truncation functor for $\omega$-polygraphs}
The category $\oPol$ of $\omega$-polygraphs (or simply polygraphs) is
the projective limit of the following diagram in $\CAT$: 
\begin{displaymath}
  \xymatrix{\nPol 0 & \nPol 1 \ar[l]_{\trunc_0^1}& \cdots \ar[l]& \nPol n\ar[l] & \nPol{n+1} \ar[l]_{\trunc_n^{n+1}}&\cdots \ar[l]}.
\end{displaymath}
Thus, as in the case of globular sets and $\omega$-categories, we get
a family of truncation functors $\trunc_n: \oPol\to \nPol{n}$. 
Also, keeping the above notations, for each $\omega$-polygraph $P$ and
integer $n$, the $n$-polygraph $\trunc_n(P)$ will be denoted
by~$\truncpol{P}{n}$. Likewise,  for each $n$, any morphism $f:P\to Q$ in
$\oPol$ gives rise by truncation to a morphism $\truncpol{f}{n}: \truncpol
Pn\to\truncpol Qn$.

\section{Three Adjunctions}\label{sec:adj_triangle}

This section investigates some fundamental adjunctions between
the categories $\oPol$, $\oCat$ and $\oGlob$. We have already examined
the monadic adjunction between $\oCat$ and $\oGlob$ (see
Propositions~\ref{prop:catglob} and~\ref{prop:catglobmonad}). We now
turn to two further important pairs of adjoint functors.

\begin{paragr}
\nomenclature[GC]{$\catpol(C)$}{standard resolution of an $\omega$-category
$C$}
 For each $n\in\N$, we first define a functor $\catpol_n: \nCat n \to \nPol n$
together with a natural transformation 
$\varepsilon : \freecatpol_n\catpol_n\to \id$.

Thus, let $C$ be an $n$-category. The $n$-polygraph $P=\catpol_n(C)$, as
well as $\varepsilon_C$, are defined dimensionwise as follows:
\begin{itemize}
\item The set $P_0=(\catpol_n(C))_0$ is just $C_0=(\freecatpol_n(P))_0$, and $(\varepsilon_C)_0$ is the
  identity.
\item Suppose $P$ and $\varepsilon_C$ have been defined up to dimension
  $k<n$. The set~$P_{k+1}$ of $(k + 1)$-generators then consists in triples
  $p=(z,x,y)$ where $x$, $y$ are parallel $k$-cells in $P_k^*=(\freecatpol_n\catpol_n(C))_k$ and $z$
  is a $(k + 1)$-cell in $C_{k+1}$ of the form
  $z:\varepsilon_C(x)\to\varepsilon_C(y)$. The source and target maps 
  \begin{displaymath}
    \sce{k},\tge{k} : P_{k+1}\to P_k^*
  \end{displaymath}
are given by $\sce{k}(p)=x$ and $\tge{k}(p)=y$, and $\varepsilon_C$
extends to $P_{k+1}$ by
$\varepsilon_C(p)=z$. Thus $P$ is now defined up to dimension $k + 1$, and
by applying \cref{lemma:ext-universal}, the map $\varepsilon_C$
extends to $P_{k+1}^*$, yielding a $(k+1)$-functor.
\end{itemize}
Likewise, $\catpol_n$ is defined on morphisms. Functoriality of $\catpol_n$ and naturality of
$\varepsilon$ immediately follow from the construction.
\end{paragr}

\begin{lemma}
  For each $n\in\N$, the functor $\catpol_n$ is right adjoint to $\freecatpol_n$.
\end{lemma}
\begin{proof}
  It is sufficient to check that the natural transformation
  $\varepsilon$, which becomes of course the counit of the adjunction,
  satisfies the following universal property: for any $n$-functor
    $f:\freecatpol_n(P)\to C$, where $P$ is an $n$-polygraph and $C$ an
    $n$\nbd-category, there is a unique morphism $g:P\to \catpol_n(C)$
in $\nPol n$ such that the following triangle commutes:
\begin{displaymath}
  \vxym{ & \freecatpol_n\catpol_n(C)\ar[d]^{\varepsilon_C}\\
               \freecatpol_n(P)\ar[r]_f\ar[ru]^{\freecatpol_n(g)}  & C
             \pbox.}
\end{displaymath}
Here again $g$ is built by induction on all dimensions $0\leq k \leq
n$. For $k=0$, we must have $g_0=f_0$. Suppose now that $g$ has been
defined up to dimension $k<n$, satisfying the commutation
condition. Let $p\in P_{k+1}$ be a $(k + 1)$\nbd-generator of $P$,
$u=\sce{k}(p)$, and $v=\tge{k}(p)$ in $P_k^*$. The induction hypothesis
and the definition of $\varepsilon$ imply that
$g(p)=(f(p),u,v)$.
Now Lemma~\ref{lemma:ext-universal} applies, and we may
extend~$\freecatpol_n(g)$ up to a $(k+1)$-functor still satisfying the
commutation condition. 
\end{proof}

\begin{paragr}
Consider now the diagram
\begin{displaymath}
  \xymatrix{\nPol 0 \ar@<-1ex>[d]_{\freecatpol_0}& \nPol 1 \ar[l]_{\trunc_0^1}\ar@<-1ex>[d]_{\freecatpol_1}& \cdots \ar[l]&
    \nPol n\ar[l] \ar@<-1ex>[d]_{\freecatpol_n}& \nPol{n+1} \ar[l]_{\trunc_n^{n+1}}\ar@<-1ex>[d]_{\freecatpol_{n+1}}&\cdots \ar[l]\\
\nCat 0\ar@<-1ex>[u]_{\catpol_0} & \nCat 1 \ar[l]_{\trunc_0^1}\ar@<-1ex>[u]_{\catpol_1}& \cdots \ar[l]& \nCat n\ar[l]\ar@<-1ex>[u]_{\catpol_n} & \nCat{n+1} \ar[l]_{\trunc_n^{n+1}}\ar@<-1ex>[u]_{\catpol_{n+1}}&\cdots \ar[l]}
\end{displaymath}
where the squares involving $\freecatpol$ and the squares involving $\catpol$ commute. The projective
limit of the top row is the category $\oPol$ and the projective limit of the bottom
row is the category $\oCat$. Therefore we get a pair of functors
$\freecatpol:\oPol\to \oCat$ and $\catpol:\oCat \to \oPol$ such that, for each $n$,
$\freecatpol_n\circ \trunc_n=\trunc_n\circ \freecatpol$ and $\catpol_n\circ
\trunc_n=\trunc_n\circ \catpol$. Moreover $\freecatpol$ is left adjoint to $\catpol$.  For
any polygraph $P$, $\freecatpol(P)$ is the \emph{freely generated
$\omega$-category}\index{free!9-category@$\omega$-category} on $P$ and will be denoted by $\free P$.
\nomenclature[P*]{$P^\ast$}{$\omega$-category generated by a polygraph $P$}
\end{paragr}

  \begin{paragr}
Let us now define a pair of adjoint functors between categories
$\oPol$ and $\oGlob$.  Consider first the functor 
  $\poltoglob : \oPol \to \oGlob$
which takes a polygraph $P$ to a globular set $X$ by keeping only
``hereditary globular'' generators. Precisely, let $P$ be a polygraph,
we define the globular set $X=\poltoglob(P)$ dimensionwise, such that
for each $n\in\N$, $X_n\subset P_n$; recall that $\ins{k}:P_k\to \free
P_k$ denotes the canonical insertion of $k$-generators into $k$-cells:
\begin{itemize}
\item For $n=0$, $X_0=P_0$.
\item Let $n>0$ and suppose we have defined $X_k\subset P_k$ for all
  $k<n$, together with source and target maps building an
  $(n - 1)$-globular set. Let $X_n\subset P_n$ be the set of
  $n$-generators $a$ of $P$ such that $\sce{n-1}(a)$ and
  $\tge{n-1}(a)$ belong to $\ins{n-1}(X_{n-1})$ and define source and
  target maps $\sce{n-1}^X,\tge{n-1}^X:X_n\to X_{n-1}$ as the unique
  maps such that $\ins{n-1}\sce{n-1}^X(a)=\sce{n-1}(a)$ and
  $\ins{n-1}\tge{n-1}^X(a)=\tge{n-1}(a) $ for each~$a\in X_n$. This
  extends $X$ to an $n$-globular set, as shown in the following diagram:
  \begin{displaymath}
    \xymatrix{X_n \ar@<2pt>[d]^{\tge{n-1}^X}\ar@<-2pt>[d]_{\sce{n-1}^X}\ar@{^{(}->}[r]& P_n \ar@<2pt>[rd]^{\tge{n-1}}\ar@<-2pt>[rd]_{\sce{n-1}}& \\
                      X_{n-1} \ar@{^{(}->}[r]& P_{n-1}\ar[r]_{\ins{n-1}} & \free P_{n-1}\pbox.}
  \end{displaymath}
 \end{itemize}
The previous construction is clearly functorial and defines the
required functor~$\poltoglob$. Remark that $\poltoglob$ admits a left adjoint
$\globtopol:\oGlob\to\oPol$ which takes the globular set $X$ to a
polygraph $P$ such that $P_n=X_n$, in other words, $\globtopol$
defines a natural inclusion of $\oGlob$ into $\oPol$. 
\end{paragr}

\begin{lemma}\label{lemma:adj_triangle}
  There is a natural isomorphism  $\phi:\poltoglob\catpol\to \fgf $, that is,
  the diagram 
  \begin{displaymath}
    \xymatrix{\oCat \ar[r]^{\catpol}\ar[d]_{\fgf}& \oPol\ar[ld]^{\poltoglob}\\
                    \oGlob &}
                \end{displaymath}
commutes up to a natural isomorphism.
\end{lemma}
\begin{proof}
  Let $C$ be an $\omega$-category, and $X=\poltoglob\catpol(C)$. For each
  $n\in\N$, let $\phi^C_n:X_n\to C_n$ be the composition of the
  following maps
  \begin{displaymath}
    \xymatrix{X_n \ar@{^{(}->}[r]& \catpol(C)_n \ar[r]^(.4){\ins{n}}& \freecatpol\catpol(C)_n \ar[r]^(.6){(\varepsilon_C)_n}& C_n\pbox.}
  \end{displaymath}
The family $(\phi^C_n)_{n\in\N}$ defines a globular morphism
$\phi^C:\poltoglob\catpol(C)\to \fgf(C)$, natural in $C$. Thus we get a natural
transformation $\phi:\poltoglob\catpol\to\fgf$.  
 
Let us now define $\chi^C_n:C_n\to X_n$ by induction on $n$ such that
$\phi^C_n\circ \chi^C_n=\unit{C_n}$: 
\begin{itemize}
\item For $n=0$, $X_0=C_0$ and $\phi^C_0=\unit{C_0}=\unit{X_0}$, so that
  $\chi^C_0:C_0\to X_0$ is also $\unit{C_0}=\unit{X_0}$.
\item Suppose $n>0$ and $\chi^C_k$ has been defined up to $k=n - 1$,
  and let $z\in C_n$. Let $u=\sce{n-1}(z)$ and $v=\tge{n-1}(z)$ in $C_{n-1}$.
  By induction hypothesis, 
  $\chi^C_{n-1}(u)$ and $\chi^C_{n-1}(v)$ belong to $X_{n-1}$. Let
  $x=\ins{n-1}\chi^C_{n-1}(u)$, $y=\ins{n-1}\chi^C_{n-1}(v)$ in $\freecatpol\catpol(C)_{n-1}$
  and define $a=\chi^C_n(z)=(z,x,y)$. By construction $a\in X_{n}$ and $\phi^C_n(a)=z$.
\end{itemize}
It remains to prove that $\phi^C_n$ is injective. We reason again by
induction on $n$.
\begin{itemize}
\item For $n=0$, $\phi^C_0$ is an identity, hence injective.
\item Suppose $n>0$ and $\phi^C_{n-1}$ injective. Let
  $a_i=(z_i,x_i,y_i)\in X_n$ for $i=0,1$ such that
  $\phi^C_n(a_0)=\phi^C_n(a_1)$. Thus $z_0=z_1$. Also
  \begin{align*}
    \phi^C_{n-1}(\sce{n-1}^X(a_0)) & = \sce{n-1}(\phi^C_n(a_0))\\
                                            & =\sce{n-1}(\phi^C_n(a_1))\\
                                             & = \phi^C_{n-1}(\sce{n-1}^X(a_1))    
  \end{align*}
and because $\phi^C_{n-1}$ is injective, 
\begin{displaymath}
  \sce{n-1}^X(a_0)=\sce{n-1}^X(a_1).
\end{displaymath}
Now
\begin{align*}
  x_0 &= \sce{n-1}(a_0)\\
        & = \ins{n-1}\sce{n-1}^X(a_0)\\
        & = \ins{n-1}\sce{n-1}^X(a_1)\\
        & = \sce{n-1}(a_1)\\
        & = x_1.  
\end{align*}
Likewise $y_0=y_1$, and we get $a_0=a_1$. Hence $\phi^C_n$ is
injective, and we are done.\qedhere
\end{itemize}
\end{proof}

\section{\pdfm{(n,p)}-Polygraphs}
\label{sec:np-polygraph}
\index{8p-polygraph@$(n,p)$-polygraph}
\index{polygraph!8p-@$(n,p)$-}
\nomenclature[Polnp]{$\npPol np$}{category of $(n,p)$-polygraphs}

In the same way as an $n$-polygraph generates an $n$-category, we may
define, for each $n\geq p$, a notion of {\em $\pair np$-polygraph} generating an $\pair
np$-category. The particular case where $n=3$ and $p=1$ has been
already used to define coherent presentations of categories in \cref{chap:2-Coherent}. Remark that the construction of the left adjoint in
Proposition~\ref{prop:ncatp-leftadjoint} applies to the case where
$\nCat n$ is replaced by~$\npCat np$, as in the following diagram:
\[
\vxym{
  \npCat{n+1}p\ar@/_5ex/[ddr]\ar@/^5ex/[drr]\ar@/_3ex/[dr]_-{\fgfplus'_n}\ar@{}[dr]|{\rotatebox[origin=c]{-45}{$\bot$}}&\\
  &\ar@/_3ex/[ul]_-{L'_n}\nCatp{n,p}\ar[d]\ar[r]&\nGlob{n+1}\ar[d]\\
  &\npCat np\ar[r]&\nGlob n\pbox.
}
\]
Here $\nCatp{n,p}$ is defined as above by the pullback square of the
obvious forgetful functors. We may now define the category $\npPol np$
of $\pair np$-polygraphs, together with a functor
  $\freecatpolbis_n:\npPol np\to \npCat np$
taking an $\pair np$-polygraph to the $\pair np$\nbd-category it
generates. The definition is by induction on $n\geq p$.
\begin{itemize}
\item For $n=p$, $\npPol np$ is just $\nPol p$, and
  $\freecatpolbis_n$ is $\freecatpol_p$, as $\npCat np$ is just $\nCat
  p$.
\item Given $\npPol np$ and $\freecatpolbis_n:\npPol
  np\to\npCat np$, the category $\npPol{n+1}p$ is defined by the
  following pullback in $\CAT$:
  \[
\xymatrix{\npPol{n+1}p\ar[r]^{J'_n}\ar[d]_{U'_{n+1,n}}& \nCatp{n,p} \ar[d]\\
                       \npPol np\ar[r]_{\freecatpolbis_n} & \npCat np\pbox,}
  \]
  whereas $\freecatpolbis_{n+1}$ is $L'_nJ'_n$:
  \[
    \xymatrix{\npPol{n+1}p\ar[r]_{J'_n}\ar@/^3ex/[rr]^{\freecatpolbis_{n+1}} & \nCatp{n,p} \ar[r]_{L'_n}& \npCat{n+1}p\pbox.}
  \]
\end{itemize}
\nomenclature[P**]{$\freegpd P$}{$(n, p)$-category generated by an $(n,p)$-polygraph}
Concretely, let $(P_k)_{0\leq k\leq n}$ the sequence of
$k$-dimensional generators of an $\pair np$\nbd-polygraph $P$, and $C$ be
the $\pair np$-category generated by $P$, we shall denote
$C_k$ by $\freecat P_k$ for $0\leq k\leq p$ and by $\freegpd P_k$ for
$p<k\leq n$. Note that in the latter case, $\freegpd P_k$ contains all composites {\em and
  inverses} of the generators in $P_k$. As for $n$-polygraphs in \cref{subsec:npolyg}, the data defining an $(n,p)$-polygraph $P$ may be displayed in the
following diagram in $\Set$:
\begin{displaymath}
  \xymatrix@C=5ex{
  P_0
  \ar[d]^{\ins0}
  &\ar@<-.5ex>[dl]\ar@<.5ex>[dl]{\textstyle\ldots}
  &\ar@<-.5ex>[dl]\ar@<.5ex>[dl]P_p\ar[d]^{\ins{p}}
    &\ar@<-.5ex>[dl]_(.4){\sce p}\ar@<.5ex>[dl]^(.4){\tge{p}}P_{p+1}\ar[d]^{\ins{p+1}}
    &\ar@<-.5ex>[dl]_(.4){\sce {p+1}}\ar@<.5ex>[dl]^(.4){\tge{p+1}}P_{p+2}\ar@{}[r]|{\textstyle\ldots}&
    \ar@<-.5ex>[dl] \ar@<.5ex>[dl] P_{n-1}\ar[d]^{\ins{n-1}}&\ar@<-.5ex>[dl]_(.4){\sce{n-1}}\ar@<.5ex>[dl]^(.4){\tge{n-1}}P_{n}
\\
  P_0^*
  &\ar@<-.5ex>[l]\ar@<.5ex>[l]{\textstyle\ldots}
  &
  \ar@<-.5ex>[l]\ar@<.5ex>[l]P_{p}^*
  &\ar@<-.5ex>[l]_-{\ssce{p}{*}}\ar@<.5ex>[l]^-{\ttge{p}{*}}\freegpd P_{p+1}\ar@{}[r]|{\textstyle\ldots}&
  {\phantom{P_{n-2}^*}}{\textstyle\ldots}
  &\ar@<-.5ex>[l] \ar@<.5ex>[l] \freegpd P_{n-1}
\pbox.}
\end{displaymath}

\index{polygraph!9p-@$(\omega,p)$-}
Finally, as in the case of (plain) polygraphs, we may define a
category $\npPol{\omega}p$ of $\pair{\omega}p$-polygraphs as the
projective limit of the system $(\npPol np,U_{n+1,n})_{n\geq p}$.

\chapter{Properties of the Category of \pdfm{n}-Polygraphs}
\label{chap:pol-prop}

In this chapter, we establish the main properties of the category
$\nPol n$ of $n$\nbd-polygraphs. We first show how to compute limits and colimits and prove
that $\nPol n$ is complete and cocomplete for any $n\geq 0$. The
behavior of the cartesian product deserves a special attention in that it does {\em
  not} correspond to the product of generators. The monomorphisms
(resp.~epimorphisms) in $\nPol n$ are then characterized as injective
(resp.~surjective) maps between generators. The linearization of
polygraphic expressions plays a
central role in proving these facts. Whereas $\nPol n$ is a
presheaf category for $n\in\set{0,1,2}$, it already fails to be cartesian
closed for~$n\geq 3$, as proved in~\cite{makkai3}, the culprit for
this defect being as usual the Eckmann-Hilton phenomenon. The categories
$\nPol n$ are however locally presentable for all $n\in
\N\cup\set{\omega}$. We introduce the technical notion of {\em context}, in
relation with $n$-dimensional rewriting, and use it to prove that if an
\oo-category is freely generated by a polygraph then this polygraph is
unique up to isomorphism. Finally, we show how to define rewriting
properties of $n$-polygraphs and to prove coherence results by rewriting on
$(n-1)$-categories presented by convergent $n$-polygraphs.

\section{Limits and Colimits}

\begin{paragr}[Terminal object]\label{paragr:termpol} 
The category $\oPol$ has a terminal
  object $\termpol$, defined as the image of the terminal
  $\omega$-category $\termcat$ by the right adjoint functor $\catpol:\oCat\to
  \oPol$. Concretely, $(\termpol)_0$ consists in a single $0$-cell,
  whereas for each $n>0$, the set $(\termpol)_n$ of $n$-generators
  consists in all pairs $\pair uv$ of parallel cells in
  $(\termpol)_{n-1}^*$. Thus $\termpol$ has exactly one $0$-generator,
  one $1$-generator, and infinitely many $n$-generators for each
  $n\geq 2$. A detailed description in the case~$n=2$ can be found in~\cite{street1996categorical}.
\end{paragr}

\begin{paragr}[Products]\label{paragr:pol_product}
  Let $P$, $Q$ be a pair of polygraphs, and $\cano{P}:P\to \termpol$,
  $\cano{Q}:Q\to \termpol$ the canonical morphisms from $P$, $Q$ to the
  terminal object. We define a polygraph $R$, together with morphisms
  $p:R\to P$, $q:R\to Q$, by induction on the dimension.
  \begin{itemize}
\item For $n=0$, $R_0=P_0\times Q_0$ and $p_0:R_0\to P_0$, $q_0:R_0\to
  Q_0$ are the usual projection maps.
  \item Let $n>0$ and suppose $R$, $p$, $q$ have been defined up to dimension
    $n-1$. The set~$R_n$ of $n$-generators consists of
    quadruples 
    \[
    c=(a,b,u,v)\in P_n\times Q_n\times R_{n-1}^*\times R_{n-1}^*
  \]
such that $\cano{P}_n(a)=\cano{Q}_n(b)$ and $u$, $v$ are parallel cells
satisfying the equations $p_{n-1}^*u = \sce{n-1}(a)$, $q_{n-1}^*u =
\sce{n-1}(b)$, $p_{n-1}^*v = \tge{n-1}(a)$, and $q_{n-1}^*v = \tge{n-1}(b)$.
 The projection maps are defined by  $p_nc=a$ and $q_nc=b$.
The source and target maps $\sce{n-1},\tge{n-1}:R_n\to
  R_{n-1}^*$ are defined by $\sce{n-1}(r)=u$ and~$\tge{n-1}(r)=v$.
\end{itemize}
Now the following square is a pullback in $\oPol$:
\begin{displaymath}
  \vxym{R\ar[r]^p\ar[d]_q & P\ar[d]^{\cano{P}}\ \\
             Q \ar[r]_{\cano{Q}}& \termpol\pbox.}
\end{displaymath}
In fact, let $S$ be a polygraph and $f:S\to P$, $g:S\to Q$ morphisms
such that the following diagram commutes:
\begin{equation}
  \vxym{S\ar[r]^f\ar[d]_g & P\ar[d]^{\cano{P}}\ \\
             Q \ar[r]_{\cano{Q}}& \termpol\pbox.}
\label{eq:pullbackpol}
\end{equation}
We show that there is a unique morphism $h:S\to R$ making the
following diagram commute, for each $n\geq 0$:
  \begin{equation}
  \vxym{\truncpol Sn\ar[rd]|{\truncpol{h}{n}}\ar@/^/[rrd]^{\truncpol{f}{n}}\ar@/_/[rdd]_{\truncpol{g}{n}}& & \\
&\truncpol Rn\ar[r]_{\truncpol{p}{n}}\ar[d]^{\truncpol{q}{n}} &
    \truncpol Pn\ar[d]^{!^P_{\leq n}}\ \\
             &\truncpol Qn \ar[r]_{!^Q_{\leq n}}&
             \mathbf{1}^{\mathrm{pol}}_{\leq n}\pbox.}
\label{eq:pullbackpoln}
\end{equation}
Let us build $h$ and prove its uniqueness by induction on $n$.
\begin{itemize}
\item If $n=0$, $R_0=P_0\times Q_0$,
  $h_0:S_0\to R_0$ takes $c\in S_0$ to $\pair{f_0c}{g_0c}$ and this
  choice is unique.
\item Let $n\geq 0$ and suppose that $h$ has been defined up to
  dimension $n$ such that~(\ref{eq:pullbackpoln}) commutes. We extend $h$ to
  dimension $n + 1$ as follows: given $c\in S_{n+1}$, we define $h_{n+1}c=(a,b,u,v)$
  where 
  \begin{align}
    \label{eq:pullback1}
    a &= f_{n+1}c,\\
\label{eq:pullback2}
   b &= g_{n+1}c,\\
\label{eq:pullback3}
u& = h_n^*\sce{n}(c),\\
\label{eq:pullback4}
v&=h_n^*\tge n(c).
  \end{align}
The equations~(\ref{eq:pullback1}) to~(\ref{eq:pullback4}) ensure that
$h_{n+1}c\in R_n$ and that the globular relations 
  $\sce{n}(h_{n+1}c)=h_n^*\sce{n}(c)$ and
  $\tge{n}(h_{n+1}c)=h_n^*\tge{n}(c)$ hold, so that $h$ extends to a
  morphism up to dimension $n + 1$. Moreover, this extension now satisfies
  the required commutation conditions $f=ph$ and $g=qh$ up to
  dimension~$n + 1$. Conversely,  the commutation conditions
  imply~(\ref{eq:pullback1}) and~(\ref{eq:pullback2}), and the
  requirement that $h$ be a morphism implies~(\ref{eq:pullback3})
  and~(\ref{eq:pullback4}). Hence, the above choice for $h_{n+1}c$
  is unique and we are done. 
 \end{itemize}
 Therefore $R$ is the cartesian product $P\times Q$ in $\oPol$.
It should be emphasized that generally, from $n=2$ on,  the map
$(a,b,u,v)\mapsto \pair ab$ from $R_n$ to $P_n\times Q_n$ is not
surjective, and from $n=3$ on, not injective either. For
example, if
\begin{eqnarray*}
  P&=&\Pres{\star_1,\star_2}{a,b:\star_1\to\star_2}{f:a\to b},\\
  Q &=& \Pres{\star_1,\star_2,\star_3}
{a:\star_1\to\star_2,b:\star_2\to\star_3,c:\star_1\to\star_3}{g:a\comp 0
  b\to c}
\end{eqnarray*}
and  $R=P\times Q$, then $R_2=\emptyset$ and $P_2\times Q_2=\set{\pair{f}{g}}$, thus showing
that $R_2\to P_2\times Q_2$ is not surjective. Another example of
non-injectivity is given in~\ref{subsec:polnotcc} below.
A similar proof shows that $\oPol$ has all small
products $\prod_{i\in I} P^{i}$.
\end{paragr}

\begin{paragr}[Equalizers]
 Let $P$, $Q$ be polygraphs and $f,g:P\to Q$ be two morphisms in
 $\oPol$. For each $n\in\N$, we define a subset of $P_n$ by
 \begin{displaymath}
   R_n=\setof{a\in P_n}{f(a)=g(a)}.
 \end{displaymath}
As the inclusions $j_n:R_n\to P_n$ commute with the source and target
maps, this defines a polygraph $R$ together with a morphism $j:R\to P$ in
$\oPol$. Now $j:R\to P$ is clearly an equalizer of the pair
$\pair{f}{g}$ in $\oPol$.
\end{paragr}

\begin{paragr}[Coproducts]
  The coproducts in $\oPol$ are built by taking the coproducts of the
  corresponding sets of $n$-generators in each dimension $n$. Thus, if $(P^{i})_{i\in I}$
  is a family of polygraphs, we have for each dimension $n$, $(\coprod_{i\in I}P^{i})_n=\coprod_{i\in I}P^{i}_n$.
\end{paragr}

\begin{paragr}[Coequalizers]
  Let $P$, $Q$ be two polygraphs and $f,g:P\to Q$ be a pair of
  morphisms in $\oPol$. We first build, for each $n\in \N$, a
  coequalizer $k_n:Q_n\to R_n$ of the pair $\pair{f_n}{g_n}$ in
  $\Set$:
  \begin{displaymath}
    \vxym{
  P_n\ar@<.5ex>[r]^-{f_n}\ar@<-.5ex>[r]_-{g_n}&Q_n\ar[r]^{k_n}&R_n\pbox.
}
  \end{displaymath}
Concretely, there is a binary relation on $Q_n$ defined by $b\leadsto_n
b'$ if and only if there is an $a\in P_n$ such that $b=f_n(a)$ and
$b'=g_n(a)$. If $\sim_n$ is the smallest equivalence relation on
$Q_n$ containing $\leadsto_n$, then $R_n$ is the set $Q_n/{\sim_n}$ of
equivalence classes of $\sim_n$ and $k_n$ is the canonical surjection.

We now define, by induction on $n$, a polygraph $R$ with $R_n$ as the
set of $n$\nbd-generators, together with a morphism $k:Q\to R$ in $\oPol$
whose $n$\nbd-dimensional component is the above $k_n$.
\begin{itemize}
\item For $n=0$, we just take $R_0$ and $k_0$ as above, and there is
  nothing to prove.
\item Let  $n>0$ and suppose that the polygraph $R$ has been
  constructed up to dimension $n-1$, together with the morphism
  $k:Q\to R$ and that, for all $0\leq i<n$,
  $k_i^*f_i^*=k_i^*g_i^*$. Taking $R_n$ as above, we have to define
  source and target maps $\sce{n-1},\tge{n-1}:R_n\to
  R_{n-1}^*$. Consider the following diagram
  \begin{displaymath}
    \vxym{P_n \ar@<.5ex>[r]^-{f_n}\ar@<-.5ex>[r]_-{g_n}\ar[d]_{\sce{n-1}}& Q_n\ar[r]^{k_n}\ar[d]_{\sce{n-1}} & R_n\ar@{.>}[d]^{\sce{n-1}}\\
                P_{n-1}^* \ar@<.5ex>[r]^-{f_{n-1}^*}\ar@<-.5ex>[r]_-{g_{n-1}^*}& Q_{n-1}^*\ar[r]_{k_{n-1}^*} & R_{n-1}^*}
  \end{displaymath}
in which all solid arrows are already given, and let $c\in R_n$. Let
$b$ be a representative of $c$ in $Q_n$, that is, $c=k_n(b)$ and set
$\sce{n-1}(c)=k_{n-1}^*\sce{n-1}(b)\in R_{n-1}^*$. We must show that
$\sce{n-1}(c)$ so defined does not depend on the choice of the representative $b$. Thus let $b'\in Q_n$
such that $k_n(b')=c$. By definition, $b\sim_n b'$. In order to show
that $k_{n-1}^*\sce{n-1}(b)=k_{n-1}^*\sce{n-1}(b')$, it is sufficient
to check it in the base case where $b\leadsto_n b'$. In this case,
there is a generator~$a\in P_n$ such that $f_n(a)=b$ and $g_n(a)=b'$,
therefore
\begin{align*}
  k_{n-1}^*\sce{n-1}(b) & = k_{n-1}^*\sce{n-1}f_n(a) \\
                                    &= k_{n-1}^*f_{n-1}^*\sce{n-1}(a)\\
                                   & = k_{n-1}^*g_{n-1}^*\sce{n-1}(a)\\
                                   & = k_{n-1}^*\sce{n-1}g_n(a) \\
                                   & = k_{n-1}^*\sce{n-1}(b')
\end{align*}
so that the source map on $R_n$ is well-defined. Of course the target
map $\tge{n-1}$ is defined accordingly. As for the globular relations,
remark that if $n\geq 2$, by induction hypothesis,
\begin{align*}
  \tge{n-2}\sce{n-1} k_n& = \tge{n-2}k_{n-1}^*\sce{n-1}\\
                                      & =
                                            k_{n-2}^*\tge{n-2}\sce{n-1}\\
                                      & =
                                            k_{n-2}^*\tge{n-2}\tge{n-1}\\
                                      & =
                                            \tge{n-2}k_{n-1}^*\tge{n-1}\\
                                       & = \tge{n-2}\tge{n-1} k_n
\end{align*}
and because $k_n$ is surjective, this implies 
$\tge{n-2}\sce{n-1}=\tge{n-2}\tge{n-1}$. Likewise
$\sce{n-2}\sce{n-1}=\sce{n-2}\tge{n-1}$. Thus, the polygraph $R$ is
now defined up to dimension~$n$. By using
Lemma~\ref{lemma:ext-universal}, we finally extend $k$ to a morphism of
polygraphs up to dimension $n$, so that $k_n^*f_n^*=k_n^*g_n^*$. 
\end{itemize}
Having defined $k:Q\to R$, a similar induction process shows that $k$
is indeed a coequalizer of the pair $\pair fg$ in $\oPol$.
\end{paragr}

\begin{prop}
  The category $\oPol$ is complete and cocomplete.
\end{prop}
\begin{proof}
  The category $\oPol$ has all small products and equalizers of all
  pairs of morphisms, hence has all small limits (see~\cite[Chapter~V,
  Theorem~1]{MacLane98}). Dually, $\oPol$ has coproducts and
  coequalizers, hence has all small colimits.
\end{proof}

\begin{remark}
  \label{rmk:polsize}
  \label{rem:polsize}
  To each polygraph $P$ we may associate a set
  \[
    \sizeof P=\coprod_{k\in\N}P_k
    \]
  consisting of all generators of $P$. This correspondence is the
  object part of a functor $\sizeof{{-}}:\oPol \to \Set$.
  In fact, a morphism $f:P\to Q$ in $\oPol$ takes, for each $k\in \N$,
  the set $P_k$ to the set $Q_k$, whence a map $\sizeof f:\sizeof P\to \sizeof Q$.
Now, the above construction of the colimits in $\oPol$ shows 
that this functor $\sizeof{{-}}$ preserves all small colimits. As a
morphism $f$ of polygraphs is entirely determined by its components
$f_n$, the functor $\sizeof{{-}}$ is faithful, making $\oPol$ a {\em
  concrete} category.
\end{remark}

\section{Morphisms in \oPol}

We need a few technical preliminaries before
characterizing  monomorphisms and epimorphisms in $\oPol$. First, let $C$ be
an \oo-category and $X\subset C_n$ a set of $n$-cells of $C$:  we say
that $X$ is {\em closed under divisors}\index{divisor} if, for any $i$-composable
$n$-cells $u,v\in C_n$ such that $u\comp i v\in X$, then $u\in X$ and
$v\in X$. We may then state the following crucial property of cellular extensions:
\begin{lemma}\label{lemma:cellextinj}
  Let $n\in\N$, and
  $f:\pair CX\to\pair DY$ a morphism in $\nCatp{n}$, where $f=\pair
  gh$, with $g:C\to D$ a morphism in $\nCat{n}$ and $h:X\to Y$ a
  source and target preserving map. Let $\free f:\freeplus_n\pair CX\to
  \freeplus_n\pair DY$ be the induced morphism in $\nCat{n+1}$. 
Suppose that the maps $g_n:C_n\to D_n$ and $h:X\to Y$ are injective
and that the image of $g_n$ is closed under divisors. Then $\free
f_{n+1}$ is injective and its image is closed under divisors.
\end{lemma}
\begin{proof}
 
 We restrict here to the general line of reasoning and refer
 to~\cite[Section~2.3]{lucas2015coherence} for a complete proof. Note
 also that~\cite{makkai2005wordcomp} contains a thorough analysis of the shape
 of freely generated cells, essentially encompassing the present
 material. Thus, let $\pair CX$ and $\pair DY$ as in the above
 statement. Let $C_{n+1}$ (\resp $D_{n+1}$)
the set of $(n+1)$-cells in $\freeplus_n\pair CX$ (\resp $\freeplus_n\pair
DY$). According to~\cite[Section 4.1]{metayer2008cofibrant}, there are
sets of well-typed formal expressions $E^C$ and~$E^D$ endowed with
binary relations~$\sim^C$ and $\sim^D$ generating congruence
relations $\simeq^C$ and $\simeq^D$ such that  $C_{n+1}=E^C/{\simeq^C}$ and
$D_{n+1}=E^D/{\simeq^D}$. Moreover, $f$ induces a map
$\overline{f}:E^C\to E^D$ such that the following diagram commutes,
the vertical maps being the canonical surjections:
\begin{displaymath}
  \vxym{E^C\ar[r]^{\overline{f}}\ar[d] & E^D\ar[d]\\
              C_{n+1} \ar[r]_{f_{n+1}^*}& D_{n+1}\pbox.}
\end{displaymath}
Then structural induction on formal expressions shows that
for each $a, b\in E^C$ such that $\overline{f}(a)\sim^D
\overline{f}(b)$, then  $a\sim^C b$.  Thus, whenever
$\overline{f}(a)\simeq^D\overline{f}(b)$, $a\simeq^C b$. Therefore
$f_{n+1}^*$ is injective. Moreover, the image of $f_{n+1}^*$ is still
closed by divisors.
\end{proof}

\begin{proposition}\label{prop:injectivity}
  Let $P$, $Q$ be polygraphs, $f:P\to Q$ be a morphism in~$\oPol$ and
  $n\in \N$. If for all $k\leq n$, $f_k:P_k\to Q_k$ is injective, then
  $f_n^*:P_n^*\to Q_n^*$ is injective.
\end{proposition}

\begin{proof}
Suppose $f:P\to Q$ is a morphism in $\oPol$ such that, for all $k\leq n$, the map
$f_k:P_k\to Q_k$ is injective.  Then, by applying
Lemma~\ref{lemma:cellextinj} dimensionwise, we
check the two following properties by induction on~$k\in\set{0,...,n-1}$.  
\begin{itemize}
\item The maps $f_k^*:P^*_k\to Q^*_k$ and $f_{k+1}:P_{k+1}\to Q_{k+1}$ are
  injective.
\item The image of $f^*_k$ is closed by divisors.
\end{itemize}
Therefore, by applying once more Lemma~\ref{lemma:cellextinj}, $f^*_n$ is injective.
\end{proof}

Let $\pair CX$ be a cellular extension of an $n$-category $C$ by $X$
and $C[X]$ the freely generated $(n + 1)$-category on this
extension. As all categories $\nCat{m}$ canonically embed into
$\oCat$, we may view $C$ and $C[X]$ as \oo-categories, and consider
the canonical \oo-functor $j:C\to C[X]$.  Let us call an \oo-functor
{\em injective} if its underlying globular map is a monomorphism in
$\oGlob$, that is, injective in each dimension. Then the
following result follows immediately as a particular case of Lemma~\ref{lemma:cellextinj}.

\begin{proposition}\label{prop:canonical_injection}
  For any cellular extension $\pair CX$, the canonical \oo\nbd-functor
  $j:C\to C[X]$ is injective.
\end{proposition}

\begin{remark}
In fact, Proposition~\ref{prop:canonical_injection} holds for a more
general interpretation of the notion of cellular extension, namely one
when $C$ is {\em any }  \oo-category and $X$ is a set of cells of {\em
  any dimensions} freely attached to $C$. We refer
to~\cite[Section 4, p.~36]{makkai2005wordcomp} for a complete proof of
this generalized statement. When translated into the present language, Makkai's theorem states
precisely the following: for any \oo-category $C$ and any set $X$, together with a
family of morphisms $f_x: \sphere{n_x}\to C$, $x\in X$, the morphism
$j:C\to C[X]$ in the pushout square
\[
\vxym{
\coprod_{x\in X}\sphere{n_x}\ar[r]^-{[f_x]_{x\in X}}
  \ar[d]_{\coprod_{x\in X}\gencof{n_x}}& C\ar[d]^j\\
\coprod_{x\in X}\globe{n_x}\ar[r] & C[X]
}
\]
is injective.
The proof of Makkai goes along the same lines as the one of
Lemma~\ref{lemma:cellextinj} and involves a precise analysis of the formal
expressions denoting the cells of~$C[X]$. 
\end{remark}

\begin{paragr}[Linearization]
  \label{paragr:linear}
  \label{sec:linearization}
To each pair $\pair Xn$ such that $X$ is a set and $n\in \N$ we may
associate an \oo-category $C(X,n)$ whose only non-trivial cells are in
dimension $n$, where
\begin{itemize}
\item $C(X,n)_n=X$ if $n=0$,
\item $C(X,n)_n=\N[X]$, the free abelian monoid generated by $X$ if $n>0$.
\end{itemize}
Note that $\N[X]$ consists of linear combinations of the form
$u=\sum_{x\in X} n_x x$ where $n_x\in \N$ and $n_x=0$ for all but a
finite number of indices $x$. Moreover, for all $i<n$, the
$i$-composition of such $n$-cells is given by $u\comp i v=u+v$.
Remark that, {\em except for $n=1$}, the category $C(X,n)$ is freely
generated by a polygraph~$P(X,n)$ whose generators are given by
$P(X,n)_n=X$, $P(X,n)_0=\set{\ast}$ if~$n>0$ and
$P(X,n)_i=\emptyset$ otherwise,
the source and target maps being  uniquely determined by these data.

Let now $P$ be a polygraph and $n\in \N$. To any $n$-generator $a\in
P_n$ and any $n$-cell $u\in \free P_n$ we may unambiguously attach a natural
number $\wht{a}{u}$, the {\em weight} of $a$ in $u$, measuring the
number of occurrences of $a$ in $u$. Consider indeed the $n$-polygraph
$\truncpol Pn$ which coincides with $P$ up to dimension $n$, and
let~$X=P_n$: the above defined \oo-category $C(X,n)$ can be seen as an
$n$-category. Now, by Lemma~\ref{lemma:ext-universal}, the unique morphism
in $\nCat{n-1}$ taking $\free{(\truncpol{P}{n-1})}$ to the terminal object
extends uniquely to a morphism $\lambda:\free{(\truncpol{P}{n})}\to
C(X,n)$ whose restriction to $P_n$ is the identity $P_n\to X$. We may
therefore define the natural numbers $\wht{a}{u}$ by
\[
\lambda_n(u)=\sum_{a\in P_n} \wht{a}{u} a
\]
for any $n$-cell $u\in \free P_n$. Note that $\lambda_n(u)$ may be
seen as a {\em multiset} $\msupp u$ on~$P_n$ mapping
each generator $a\in P_n$ to $\wht{a}{u}$.
\end{paragr}

\begin{paragr}[Support]
  \label{paragr:support}
  \label{sec:support}
  \index{support}
  \nomenclature[supp]{$\supp{x}$}{support of a cell~$x$}
The above notion of weight leads to the technical notion of {\em
  support}. Let $P$ be a polygraph and $u$ an $n$-cell in $\free P_n$.
We first define the set of $n$-generators actually occurring in $u$ by
\begin{displaymath}
  \supptop{n}{u}=\setof{a\in P_n}{\wht{a}{u}\neq 0}.
\end{displaymath}
Now the \emph{total support} of $u$ is defined by induction on the
dimension of $u$:
\begin{itemize}
\item For $n=0$, $u\in P_0^*=P_0$ and $\supp{u}=\set{u}$.
\item For $n>0$,
$\supp{u}=U\cup S\cup T\cup S'\cup T'$
where
\begin{align*}
  U&=\supptop{n}{u},\\
  S&=\supp{\sce{n-1}(u)},\\
  T&=\supp{\tge{n-1}(u)}\\
  S'&=\cup_{a\in\supptop{n}{u}}
  \supp{\sce{n-1}(a)},\\
  T'&=\cup_{a\in\supptop{n}{u}}
    \supp{\tge{n-1}(a)}.
  \end{align*}
\end{itemize}
In other words, the total support of $u$ consists in all generators in
$\cup_{0\leq i\leq n} P_i$ needed to express $u$.
Likewise, for any subset $A\subset \coprod_{k\in\N}P_k$ of generators of
$P$, we define the  {\em support of $A$}  by
  $\supp{A}=\cup_{a\in A}\supp{a^*}$.
\end{paragr}

\begin{paragr}[Subpolygraph]
  \label{paragr:subpol}
  \index{subpolygraph}
  Let $P$ be a polygraph and $a\in
  P_n$ an $n$-generator. Let~$a^*$ denote the corresponding $n$-cell
  in $P_n^*$. The {\em subpolygraph of $P$
    generated by $a$}, denoted $Q=\subpol aP$, is defined as follows: for each $k\leq n$,
  the set of $k$\nbd-generators of $Q$ is $Q_k=P_k\cap \supp{a^*}$, where
  $\supp{a^*}$ is the total support of $a^*$ defined
  in~\cref{paragr:support}, whereas the source and target maps
  $\sce{k-1},\tge{k-1}:Q_k\to Q_{k-1}^*$ are obtained by restriction of
  $\sce{k-1},\tge{k-1}:P_k\to P_{k-1}^*$ for $k>0$. Of course
  $Q_k=\emptyset$ for $k>n$. Thus, we get a canonical morphism
$k_a:\subpol aP\to P$
such that $Q_n=\set{a}$ and $k_a(a)=a$.
More generally, for any subset $A\subset \coprod_{k\in\N}P_k$ of generators of
$P$, the {\em subpolygraph of $P$ generated by $A$} is the polygraph
$Q=\subpol AP$ such that, for each $k\in \N$, $Q_k=P_k\cap\supp{A}$,
and built dimensionwise, together with a canonical morphism
$k_A:\subpol AP\to P$
by using at each level the appropriate restriction to $Q$ of the source and
target maps of $P$.
\end{paragr}

\begin{proposition}\label{prop:monos_in_pol}
 A morphism $f:Q\to P$ in $\oPol$ is a monomorphism if and only if
 $f_n:Q_n\to P_n$ is injective in each dimension $n$.
\end{proposition}

\begin{proof}
In one direction, notice that the functor $\oPol\to \Set^\N$
taking a polygraph $P$ to the sequence $(P_n)_{n\in\N}$ is faithful by
construction, hence reflects monomorphisms. Therefore, if $f_n$ is
injective for all $n$, then $f$ is a monomorphism.

Conversely, suppose that $f:Q\to P$ is a monomorphism. Suppose also
that there is a $k\in \N$ such that $f_k$ is not injective, and let
$n$ be the smallest such integer. By hypothesis, there are two distinct
$n$-generators $a,a'\in Q_n$ such that $f_n(a)=f_n(a')$. Hence
\begin{displaymath}
  f_{n-1}^*(\sce{n-1}(a))=\sce{n-1}(f_n(a))=\sce{n-1}(f_n(a'))=f_{n-1}^*(\sce{n-1}(a')).
\end{displaymath}
But as $f_{n-1}$ is injective, so is $f_{n-1}^*$ by
Lemma~\ref{prop:injectivity}. Hence
$\sce{n-1}(a)=\sce{n-1}(a')$. Likewise
$\tge{n-1}(a)=\tge{n-1}(a')$. As a consequence, there is a unique
morphism $h:\subpol{a}{Q}\to\subpol{a'}{Q}$ taking $a$ to $a'$. Let
$k=k_a$ and $k'=k_b\circ h$. We now have a diagram
\begin{displaymath}
  \vxym{\subpol{a}{Q}\ar@<.5ex>[r]^-{k}\ar@<-.5ex>[r]_-{k'} & Q \ar[r]^f& P}
\end{displaymath}
such that $fk=fk'$ but $k\neq k'$. This contradicts the hypothesis,
and we are done.  
\end{proof}

\begin{corollary}\label{corol:monos_in_pol}
The functor $\freecatpol:\oPol\to \oCat$ preserves monomorphisms.
\end{corollary}
\begin{proof}
  Suppose that $f:Q\to P$ is a monomorphism in $\oPol$. By
  Proposition~\ref{prop:monos_in_pol}, $f_n$ is injective in each
  dimension $n$, and by Lemma~\ref{prop:injectivity}, so is~$f_n^*$. Now the
  functor $C\mapsto (C_n)_{n\in\N}$ from $\oCat$ to
  $\Set^{\N}$ is faithful and reflects monomorphisms. Hence
  $f^*:Q^*\to P^*$ is a monomorphism.
\end{proof}

\begin{remark}
  Proposition~\ref{prop:monos_in_pol} shows that for each subset $A$
  of generators of a polygraph $P$, the canonical morphism
  $k_A:\subpol AP\to P$ is a monomorphism. Conversely, any
  monomorphism $f:Q\to P$ in $\oPol$ factorizes as $f=k_A\circ h$,
  where $A=\setof{f(q)}{q\in Q_k,k\in\N}$ and $h:Q\to\subpol AP$ is an isomorphism.
\end{remark}

\begin{proposition}\label{prop:epis_in_pol}
  A morphism $f:P\to Q$ in $\oPol$ is an
  epimorphism if and only if $f_n$ is surjective in each dimension
  $n$.
\end{proposition}

\begin{proof}
  As above, the functor $P\mapsto (P_n)_{n\in\N}$ is faithful and
reflects epimorphisms. Therefore if $f_n$ is surjective in each
dimension, then $f$ is an epimorphism.
Conversely, suppose that $f:P\to Q$ is an epimorphism in
$\oPol$. Consider the abelianization functor $\ab:\oCat\to \pCh[\Z]$
(see \cref{chap:homology}). By
Proposition~\ref{prop:ab_leftadj}, $\ab$ is a left adjoint. Recall
that $\freecatpol:\oPol\to\oCat$ is also a left adjoint. Therefore the
composition
$\ab\circ\freecatpol:\oPol\to \pCh[\Z]$ is a left adjoint and
preserves epimorphisms. As a consequence, 
the map  $(\ab(f^*))_n:\Z[P_n]\to\Z[Q_n]$ is surjective for each $n$,
but $(\ab(f^*))_n$ is nothing but the linearization of $f_n$, hence
$f_n$ itself is surjective.
\end{proof}

\begin{proposition}\label{prop:isomonoepi}
  A morphism $f:P\to Q$ in $\oPol$ is an isomorphism if and only if, for each
  $n\geq 0$,  it induces a bijection $f_n:P_n\to Q_n$.
\end{proposition}
\begin{proof}
  In one direction, if $f$ is an isomorphism, so is $\sizeof f$ by
  functoriality, whence also all maps $f_n$ for $n\geq 0$.
Conversely, suppose $f_n$ is a bijection in all dimensions $n$. We define $g:Q\to
  P$ inverse to $f$ in $\oPol$ by induction on the dimension.
  \begin{itemize}
  \item For $n=0$,  let $g_0:Q_0\to P_0$ be the inverse map
    $(f_0)^{-1}$ of $f_0$.
    \item Let $n\geq 0$ and suppose $g$ has been defined up to
      dimension $n$ such that $g_k\circ f_k=(\id_{P})_k$ and $f_k\circ
      g_k=(\id_{Q})_k$ for all $0\leq k\leq n$. Define
      \[
        g_{n+1}=(f_{n+1})^{-1}:Q_{n+1}\to P_{n+1}.
        \]
     Let $a\in
      Q_{n+1}$. We have to check that
      $\sce{n}(g_{n+1}(a))=\free g_n(\sce{n}(a))$.
      Now
    \[\free f_n(\sce{n}(g_{n+1}(a)) =
         \sce{n}(f_{n+1}g_{n+1}(a))= \sce{n}(a)= \free f_n(\free g_n(\sce{n}(a)).
    \]
   By induction hypothesis, $\free f_n$ is a bijection, whence the
   desired equality.
  \end{itemize} 
Likewise $\tge{n}(g_{n+1}(a))=\free g_n(\tge{n}(a))$.
 Therefore, $g$ is defined, and is inverse to
 $f$ up to dimension $n+1$, whence the result.
 \end{proof}

\begin{remark}
  \label{rmk:polsize-preservation}
  As an immediate consequence of the above results, the isomorphisms
  in $\oPol$ are exactly the morphisms which are monomorphisms and epimorphisms.
Moreover the
  functor $\sizeof{{-}}$ from \cref{rmk:polsize} preserves and
  reflects monomorphisms, epimorphisms, and isomorphisms.
\end{remark}

\begin{paragr}[Subobject classifier]
The category $\oPol$ has a {\em subobject classifier}, that is, an
object $\subobj$ together with a monomorphism $\true:\termpol\to \subobj$
such that for every monomorphism $f:P\to Q$ in $\oPol$, there is a
unique morphism $\charfunc{f}:Q\to \subobj$ making the following
diagram a pullback square:
\begin{equation}
  \vxym{P\ar[r]^{\cano{P}}\ar[d]_f & \termpol\ar[d]^{\true} \\
              Q \ar[r]_{\charfunc{f}}& \subobj\pbox.}
\label{eq:subobj_pullback}
\end{equation}
As usual, these conditions determine $\subobj$ up to isomorphism.
Let us now  build $\true:\termpol\to \subobj$ by induction on the
 dimension. To avoid overloaded notation, we denote $\termpol$ by $T$
 throughout the construction.
 \begin{itemize}
 \item For $n=0$, $T_0$ is a singleton, whereas
   $\subobj_0=T_0+T_0$ elements and 
   $\true_0$ is the left inclusion $T_0\to T_0+T_0$.
\item Let $n\geq 0$ and suppose $\subobj$ and $\true$ have been
  defined up to dimension $n$. Recall from~\cref{paragr:termpol} that the set
  of $(n + 1)$-generators of $T$ is 
  \begin{displaymath}
    T_{n+1}=\setof{\pair{u}{v}}{u\in T_n^*,v\in
      T_n^*,\sce{n-1}(u)=\sce{n-1}(v), \tge{n-1}(u)=\tge{n-1}(v)}
  \end{displaymath}
with source and target maps defined by $\sce{n}\pair uv=u$ and
$\tge{n}\pair uv=v$. Consider now the set
\begin{displaymath}
  S_{n+1}=\setof{\pair{u}{v}}{u\in \subobj_n^*,v\in
      \subobj_n^*,\sce{n-1}(u)=\sce{n-1}(v), \tge{n-1}(u)=\tge{n-1}(v)}
\end{displaymath}
and the following two subsets of $S_{n+1}$:
\begin{align*}
  S_{n+1}^0 & = \setof{\pair{\true^*_n u}{\true^*_n v}}{\pair uv\in
                  T_{n+1}},\\
  S_{n+1}^1 & = S_{n+1}\setminus S_{n+1}^0.
\end{align*}
The set of $(n + 1)$-generators of $\subobj$ is then
\begin{equation}
  \subobj_{n+1}= S_{n+1}^0+S_{n+1}^0+S_{n+1}^1
\label{eq:subobj}
\end{equation}
and $\true_{n+1}$ sends $T_{n+1}$ to the first copy of $S_{n+1}^0$ by 
\begin{displaymath}
  \true_{n+1}\pair uv=\pair{\true^*_n u}{\true^*_n v}.
\end{displaymath}
The source and target maps $\sce{n},\tge{n}:\subobj_{n+1}\to
\subobj_n^*$ are naturally given by $\sce{n}\pair uv=u$ and
$\tge{n}\pair uv=v$, making $\true$ a morphism of $(n + 1)$-polygraphs.
\end{itemize}
Suppose now that $f:P\to Q$ is a monomorphism. We define a morphism
$\charfunc f:Q\to \subobj$ by induction on the dimension.
\begin{itemize}
\item  If $n=0$, $(\charfunc{f})_0:Q_0\to \subobj_0=T_0+T_0$ sends a
  $0$-cells $u$ of $Q$ to the left component if and only if $a\in \im
  f_0$.
\item Suppose $\charfunc f$ has been defined up to dimension $n$, and
  let $a\in Q_{n+1}$, with $\sce{n}(a)=u$ and $\tge{n}(a)=v$. The pair 
$c=\pair{(\charfunc{f}^*)_n(u)}{(\charfunc{f}^*)_n(v)}$ is by
induction a pair of parallel $n$-cells of $\subobj_n^*$ and three
cases are possible:
\begin{itemize}
\item If $c\in S_{n+1}^0$ and $a\in\im f_{n+1}$, $(\charfunc{f})_{n+1}$
  sends $a$ to $c$ in the first $S_{n+1}^0$ component of~(\ref{eq:subobj}).
\item If $c\in S_{n+1}^0$ and $a\notin\im f_{n+1}$, $(\charfunc{f})_{n+1}$
  sends $a$ to $c$ in the second $S_{n+1}^0$ component of~(\ref{eq:subobj}).
\item If $c\in S_{n+1}^1$, $(\charfunc{f})_{n+1}$
  sends $a$ to $c$ in the $S_{n+1}^0$ component of~(\ref{eq:subobj}).
\end{itemize}
\end{itemize}
This defines $\charfunc f$ as a morphism of $(n + 1)$-polygraphs.
We easily check that $\charfunc f$ so defined is the unique morphism
such that~(\ref{eq:subobj_pullback}) is a pullback square.

\end{paragr}

\subsection{A counterexample}

Let us end this review of morphisms in $\oPol$ by the following small
observation. To each \oo-category $C$ corresponds an
\oo-groupoid $\freegrpd C$ in $\oGpd$ obtained by formally inverting
all $n$-cells, $n>0$,  in~$C$. The correspondence $C\mapsto \freegrpd
C$ is in fact the object part of the left adjoint to the inclusion
$\oGpd\to \oCat$ and the unit of the associated monad yields a
morphism $\eta_C:C\to \freegrpd C$. In case $C=\freecat P$ is freely
generated by a polygraph~$P$, one could expect $\eta_C$ to be
injective. However, this is proved wrong by the following
counterexample. Let
\[
P=\Pres{\star}{x:\star\to\star}{a,b:x\to\unit{\star},
  c:\unit{\star}\to x}
\]
and $C=\freecat P$. Consider the $2$-cells $u,v\in C_2$ given by
$u=a\comp 1 c\comp 1 b$,  $v=b\comp 1 c\comp 1 a$ and define
$u'=\eta_C u$, $v'=\eta_C v$. Whereas $u\neq v$ in $C$, the presence
of a strict inverse $a^{-1}$ for $a$ in $\freegrpd C$ implies that
$u'=v'$, as
\[
  \begin{split}
    u' & = a\comp 1 c\comp 1 b = a\comp 1 (c\comp 1 a) \comp 1 (a^{-1}\comp
    1 b)
    \\
       &  = a \comp 1 (a^{-1}\comp 1 b) \comp 1 (c\comp 1 a) 
       = b\comp 1 c\comp 1 a
       = v'.
  \end{split}
\]

\section{Is \nPol n a Topos?}\label{sec:pol_topos}
As the category of $n$-polygraphs is complete and cocomplete and has a
subobject classifier, it is natural to ask if $\nPol{n}$ is a
topos. In fact, $\nPol 0$ is the category of sets and  $\nPol{1}$ is
the category of graphs, hence both are
presheaf categories and thus are topoi. As for $n=2$, 
Carboni and Johnstone~\cite{carboni1995connected} (corrected
in~\cite{carboni2004corrigenda}) have proved that $\nPol{2}$ is also a presheaf
category and thus a topos. However, from $n=3$ on, $\nPol n$ fails to
be cartesian {\em closed}, as shown by Makkai and
Zawadowski~\cite{makkai3} (see also~\cite[Section~6, p.~57]{makkai2005wordcomp}), and therefore is not even an elementary
topos. Let us finally  mention that Batanin~\cite{batanin2002computads}
defines a notion of $T$-computad for each monad $T$ on
globular sets and gives sufficient conditions on $T$ for the category
of $T$\nbd-com\-pu\-tads to be a presheaf category. Of course, \cite{makkai3}
implies that these conditions do not hold when $T$ is the monad of
strict $\omega$-categories on globular sets.

\subsection{The category \texorpdfstring{$\nPol 2$}{Pol2} is a presheaf category}
The category $\nPol 2$  is a category of presheaves over the small
category whose objects are $p_0$, $p_1$, and $p_{m,n}$ for $m,n\in\N$, morphisms
are generated by $\src{}:p_0\to p_1$, $\tgt{}:p_0\to p_1$ and 
\begin{align*}
  \src{m,n}^i&:p_1\to p_{m,n}&&\text{for $m,n\in\N$ with $0\leq i<m$,}
  \\
  \tgt{m,n}^i&:p_1\to p_{m,n}&&\text{for $m,n\in\N$ with $0\leq i<n$,}
  \\
  \sigma_{m,n}^i&:p_0\to p_{m,n}&&\text{for $m,n\in\N$ with $0\leq i\leq m$,}
  \\
  \tau_{m,n}^i&:p_0\to p_{m,n}&&\text{for $m,n\in\N$ with $0\leq i\leq n$,}
\end{align*}
subject to the relations
\begin{align*}
  \sigma_{m,n}^i&=\src{m,n}^i\circ\src{}
  &
  \sigma_{m,n}^{i+1}&=\src{m,n}^i\circ\tgt{}
  &
  \sigma_{m,n}^0&=\tau_{m,n}^0
  \\
  \tau_{m,n}^i&=\tgt{m,n}^i\circ\src{}
  &
  \tau_{m,n}^{i+1}&=\tgt{m,n}^i\circ\tgt{}
  &
  \sigma_{m,n}^{m}&=\tau_{m,n}^{n}
\end{align*}
for every index such that the morphisms are defined.
A presheaf~$P$ over this category then corresponds to a polygraph with
$P(p_0)$ as $0$-cells, $P(p_1)$ as $1$-cells
\[
  a:x\to y
\]
with $x=P(\src{}(a))$ and $y=P(\tgt{}(a))$, and $P(p_{m,n})$ as $2$-cells
\[
  \alpha:a_1\ldots a_m\to b_1\ldots b_n
\]
with $a_i=P(\src{m,n}^i(\alpha))$ and $b_i=P(\tgt{m,n}^i(\alpha))$.

\subsection{The category \texorpdfstring{$\nPol 3$}{Pol3} is not cartesian closed}\label{subsec:polnotcc}

The original argument by Makkai and Zawadowski being quite intricate,
we give here the simpler proof due to Cheng~\cite{cheng2012direct},
based on an explicit counterexample: we shall describe  $3$-polygraphs
$P$, $Q$, $R$, and $S$ such that the diagram
\[
\vxym{
  P\ar@<.5ex>[r]^-{f_1}\ar@<-.5ex>[r]_-{f_2}&Q\ar[r]^g&R
}
\]
is a coequalizer in $\nPol 3$ not preserved by the functor
$-\times S:\nPol 3\to\nPol 3$. Therefore $-\times S$ does not preserve
colimits, hence admits no right adjoint and $\nPol 3$ is not cartesian
closed. Now $P$, $Q$, $R$, $S$ are as follows:
\begin{itemize}
\item $P_0=\set{\star}$, $P_1=\emptyset$, $P_2=\set{a}$, and
  $P_3=\emptyset$ so that there is no choice for the source and target
  maps,
\item $Q_0=\set{\star}$, $Q_1=\emptyset$, $Q_2=\set{b_1,b_2}$ where
  $f_ia=b_i$ for $i\in\set{1,2}$, and $Q_3=\set{b_3}$ where $b_3:b_1\comp 1 b_2\to \unit{2}(\star)$,
\item $R_0=\set{\star}$,
  $R_1=\emptyset$, $R_2=\set{c}$ where $c=gb_1=gb_2$, and $R_3=\set{d}$
  where $gb_3=d$ and $d:c\comp 1 c \to \unit{2}(\star)$,
\item  $S=Q$.
 \end{itemize}
Let $P'=P\times S$, $Q'=Q\times S$, and $R'=R\times S$. By the
construction of products described in \cref{paragr:pol_product},
$P'$, $Q'$, and $R'$ have a single $0$-generator $\star'=\pair{\star}{\star}$
and no $1$-generator. As for generators in dimensions $2$ and $3$,
\begin{itemize}
\item $P'_2$ has two elements $a'_i=(a,b_i,\unit{1}(\star'),\unit{1}(\star'))$ for
  $i\in\set{1,2}$ and $P'_3=\emptyset$,
\item $Q'_2$ has four elements $b'_{ij}=(b_i,b_j,\unit{1}(\star'),\unit{1}(\star'))$
  for $\pair ij\in\set{1,2}\times\set{1,2}$,
  \item $Q'_3$ has {\em two}
  elements 
  \begin{align*}
    b'_{31}&=(b_3,b_3,b'_{12}\comp 1 b'_{21},\unit{2}(\star')),\\
    b'_{32}&=(b_3,b_3,b'_{11}\comp 1 b'_{22},\unit{2}(\star')),
  \end{align*}
\item $R'_2$ has two elements
  $c'_i=(c,b_i,\unit{1}(\star'),\unit{1}(\star'))$ for $i\in\set{1,2}$,
  \item $R'_3$ has a single element $d'=(d,b_3,c'_1\comp 1 c'_2,\unit{2}(\star'))$. 
\end{itemize}
Consider now the coequalizer diagram
\[
\vxym{
  P'\ar@<.5ex>[r]^-{f'_1}\ar@<-.5ex>[r]_-{f'_2}&Q'\ar[r]^{g'}&R''
}
\]
where $f'_i=f_i\times \id_S$. As $P'_3=\emptyset$ and $Q'_3$ has two
elements, $R''_3$ also has two elements, whereas $R'_3$ only has
one. Therefore $R'$ is not isomorphic to $R''$ and we are done.
Note that the key of the above argument is the Eckmann-Hilton
phenomenon, according to which $b_1\comp 1 b_2=b_2 \comp 1 b_1$: 
this in turn implies that $b'_{31}$ and $b'_{32}$ are defined in such a
way that both projections of $b'_{12}\comp 1 b'_{21}$ and~$b'_{11}\comp 1 
b'_{22}$ match the actual source of $b_3$.

\subsection{Presheaf subcategories of polygraphs}

Although $\oPol$ is not a presheaf category, there are several
interesting full subcategories of $\oPol$ which are presheaf
categories, as shown in~\cite{henry2018regular,henry2017nonunital}. The general
principle is to restrict the shape of generators to sufficiently
``regular''  ones. Important examples are {\em many-to-one}
polygraphs, where the target of each generator is a generator,
or {\em non-unital} polygraphs, where the source and target of each
generator cannot be identities.
An alternative approach is taken
in~\cite{hadzihasanovic2019representable,hadzihasanovic2018weak,hadzihasanovic2020combinatorial}, where a category of {\em
  regular polygraphs} is {\em defined} as the presheaf category on a
small category of certain
globular shapes, then shown to be equivalent to a full subcategory of~$\oPol$.


\section{Local Presentability}
This section closely
follows~\cite[Section~5]{makkai2005wordcomp}. \Cref{chap:loc_pres}
recalls everything we need about locally presentable categories in the
present section.
Let us call a polygraph $P$ {\em finite}\index{polygraph!finite} when
$\sizeof P=\coprod_k P_k$ is a finite
set.  Given any polygraph $Q$, we may then consider the set of finite
subpolygraphs of $Q$, in the sense of~\cref{paragr:subpol}.  Let $I_Q$ be
the subcategory of $\oPol$ whose objects are the finite
subpolygraphs of $Q$, the morphisms are the canonical inclusions
$P\to P'$ for all $P$, $P'$ such that $\sizeof P\subset \sizeof{P'}$ and $D:I_Q\to \oPol$ the corresponding inclusion
functor. Let $P$, $P'$ be objects in $I_Q$, together with their
canonical inclusion morphisms $f:P\to Q$ and $f':P'\to Q$. Consider the
finite subpolygraph $P''$ of $Q$ generated by $\sizeof
P\cup\sizeof{P'}$ and $h:P''\to Q$ the corresponding inclusion
morphism. Clearly $f$ and $f'$ factor through $f''$ as in the diagram
\[
  \xymatrix@R=2ex{
    P \ar[rd]\ar@/^/[rrd]^f& & \\
    & P''\ar[r]|{f''} & Q\\
    P'\ar[ru]\ar@/_/[rru]_{f'}& &
  }
  \]
  so that $I_Q$ is a directed poset. Now, by taking the colimit of
  the sets $\sizeof P$ for~$P\in I_Q$ in $\Set$, we get
  \[
    \colim_{P\in I_Q}\sizeof{D(P)}=\bigcup_{F\in I_Q}\sizeof P =\sizeof Q 
  \]
  but the functor $\sizeof{{-}}$ preserves small colimits, as
  remarked in \cref{rmk:polsize}, so that
  \[
    \sizeof{\colim_{P\in I_Q}D(P)}=\sizeof Q.
  \]
  Therefore the canonical map $j:\colim_{P\in I_Q}D(P)\to Q$ is such
  that $\sizeof j$ is an identity. By
  Remark~\ref{rmk:polsize-preservation}, the functor $\sizeof{{-}}$
  reflects isomorphisms, hence $j$ is an isomorphism in $\oPol$. Thus
  we have  proved the following statement:
  \begin{proposition}~\label{prop:canonical-colim}
    Any polygraph is a canonical colimit of
  its finite subpolygraphs.
\end{proposition}
Let now $P$ be a finite polygraph and $\pair{Q^{(i)}}{f_{ij}}$ a
filtered system in $\oPol$  indexed by a small category $I$, with
colimit $Q=\colim_{i}Q^{(i)}$ and $f_i:Q^{(i)}\to Q$ the canonical
morphisms for $i\in I$. Let $f:P\to Q$ be a morphism in
$\oPol$. Because the functor $\sizeof{{-}}$ preserves colimits and
$\sizeof P$ is finite, there is an $i\in I$ and a map $h:\sizeof P\to
\sizeof{Q^{(i)}}$ such that
\begin{itemize}
\item $\sizeof f$ factors as in the following
diagram
\[
  \xymatrix{
    & \sizeof{Q^{(i)}}\ar[d]^{\sizeof{f_i}}\\
    \sizeof P\ar[r]_{\sizeof f}\ar[ru]^{h} & \sizeof{Q}\pbox,
    }
  \]
  \item for all $a\in \sizeof P$, $b\in \sizeof{Q^{(i)}}$, and $b'\in
    \sizeof{Q^{(i)}}$ such that 
    \[ \sizeof{f_i}(b)=\sizeof{f_i}(b')=\sizeof f(a), \]
    we have $b=b'$.
\end{itemize}
A double induction on the pair $\pair np$, where $n$ is the highest
dimension $k$ such that $P_k\neq \emptyset$ and $p$ is the cardinal of
$P_n$, shows with the help of Lemma~\ref{lemma:ext-universal} that there is a morphism
$g:P\to Q^{(i)}$ is $\oPol$ factoring $f$ as in the diagram
\[
  \xymatrix{
    & Q^{(i)}\ar[d]^{f_i}\\
   P\ar[r]_f\ar[ru]^{g} & Q \pbox.
    }
  \]
  Thus
  \[
    \hom{\oPol}{P}{\colim_iQ^{(i)}}\simeq \colim_i\hom{\oPol}{P}{Q^{(i)}}
  \]
  and we have proved the following result:
  \begin{proposition}~\label{prop:fp-polygraphs}
      A finite polygraph is a finitely presentable object in $\oPol$.
    \end{proposition}
    
  One easily checks that the isomorphism classes of finite polygraphs
  form a set. As a consequence of
  Propositions~\ref{prop:canonical-colim} and~\ref{prop:fp-polygraphs}
  we get the following theorem:
  \begin{theorem}~\label{thm:pol-is-fp}
    The category $\oPol$ is locally finitely presentable.
  \end{theorem}

\section{Contexts}
\label{sec:contexts}
We recall here very briefly the notion of {\em context}, based on the presentation
of~\cite[Section~5, p.~191]{metayer2008cofibrant}, and refer to this
article for detailed proofs. Let~$n\geq 1$, $P$ a polygraph and
$x=\pair uv$ an ordered pair of parallel $(n - 1)$-cells in~$\free
P$. We call here such a pair $x$ an {\em $n$-type} (a convenient
alternative terminology for ``$n$-sphere'' as defined in~\secr{cell-ext}).
The polygraph $P$ may be extended to a polygraph $P[x]$
by adjoining a new generator $x$ in dimension $n$, such that $\sce{n-1}x=u$ and
$\tge{n-1}x=v$. We call the $n$-cell $\indet x=\free x$ (where $\free
x$ is short for $\ins{n} (x)$) of $\free{P[x]}$ an {\em
  $n$-indeterminate\index{indeterminate} of type $x$ over $P$}.

\begin{defi}\label{def:contexts}
  Let $n\geq 1$, $P$ a polygraph and $x$ an $n$-type of $P$. An
  $n$-cell $u\in \free{P[x]}$ is a {\em context}\index{context} of type $x$ if $\wht{x}{u}=1$.
\end{defi}

A context $u$ of type $x$ will be denoted $u=\varctx{c}{x}$, where
$\indet{x}=\free x$. An $n$\nbd-con\-text $u=\varctx cx$ of type $x$ is {\em
  thin}\index{context!thin}\index{thin context} whenever $\wht{a}{u}=0$ for all $n$-generators $a\in
P[x]_n\setminus\set{x}$, that is, all generators in $\supp{u}$ are of
dimension $<n$ but $x$ itself, 
and $\varctx cx$ is {\em trivial} if $\varctx
cx=\indet x$. There is a well-defined operation of {\em substitution}\index{substitution}
in $n$-contexts: let $z$ be an $n$-cell of $\free P$ of type $x=\pair uv$, that
is, such that $\sce{n-1}z=u$ and $\tge{n-1}z=v$, and $\varctx cx$ an
$n$-context of type $x$ over $P$. From Lemma~\ref{lemma:ext-universal}, we get a morphism
\[
  \subst z : \free{P[x]}\to \free P
  \]
in $\oCat$ taking $\indet x$ to $z$ and leaving other generators unchanged. Then the substitution of~$\indet x$ by~$z$ in~$c$, noted $\ctx cz$ may be defined as
$\subst z(\varctx cx)$.

Recall from~\cite{metayer2008cofibrant} that thin $n$-contexts can always be
expressed in the (non-unique) form:
\begin{equation}
  \label{eq:thin_context}
  \varctx cx= u_{n-1}\comp{n-2}(... \comp 1(u_1\comp 0 \indet x\comp 0
 v_1)\comp 1...)\comp{n-2} v_{n-1}
\end{equation}
where $u_i$ and $v_i$ are identities over $i$-dimensional cells. This
remark leads to the following technical result:

\begin{lemma}\label{lemma:deltacontext}
  For each $n>1$ there is a map $\partial$ taking each thin
  $n$-context~$\varctx cx$ to an
  $(n - 1)$-context $\dvarctx c{x'}$ of type $x'=\pair{\sce{n-2}\indet
    x}{\tge{n-2}\indet x}$ such that
  \begin{itemize}
  \item for each $n$-cell $z$ of type $x$, $\sce{n-1}\ctx
    cz=\dctx{c}{\sce{n-1}z}$ and $\tge{n-1}\ctx
    cz=\dctx{c}{\tge{n-1}z}$,
    \item if $\dvarctx c{x'}$ is trivial, then so is $\varctx cx$.
  \end{itemize}
\end{lemma}
\begin{proof}
  See~\cite[p.~193]{metayer2008cofibrant}.
\end{proof}

\begin{remark}\label{rk:context_promotion}
 Conversely, let $\varctx cy$ be an $(n - 1)$-context of type
$y=\pair uv$, where $u$, $v$ are parallel $(n - 2)$-cells, and $\indet
x$ be an $n$-indeterminate such that $\sce{n-2}(\indet x)=u$ and
$\tge{n-2}(\indet x)=v$, there is a unique thin $n$-context $\ivarctx
cx$ such that~$\partial\ivarctx cy=\varctx cy$. 
\end{remark}

\subsection{Composition}
Let now $u=\varctx cx$ be an $n$-context of type $x$ over $P$, and
$\varctx dy$ be an $n$-context of type
$y=\pair{\sce{n-1}u}{\tge{n-1}u}$ over $P[x]$. The previously defined
substitution process yields an $n$-context $\ctx d{\varctx cx}$ of type
$x$ over $P$.

\begin{remark}\label{rk:ctx-prom-comp}
  Let $u=\varctx cx$ be an $(n - 1)$\nbd-context of type $x$,
$\varctx dy$ be an $(n - 1)$\nbd-context of type
$y=\pair{\sce{n-2}u}{\tge{n-2}u}$, and $\varctx ex=\ctx d{\varctx cx}$
the composed context as above. For each $n$-indeterminate $\indet x$
such that
\[
\pair{\sce{n-2}\indet x}{\tge{n-2}\indet
  x}=\pair{\sce{n-2}u}{\tge{n-2}u},
\]
the thin $n$-context defined in \cref{rk:context_promotion}  satisfies
the equation 
\begin{equation}
  \label{eq:comp-prom-ctx}
  \ivarctx ex=\ictx{d}{\ivarctx{c}{x}}.
\end{equation}
\end{remark}

\begin{lemma}\label{lemma:composcontext}
  For any $n$-cell $z$ of type $x$, if $\ctx d{\ctx cz}=z$ then both
  contexts $\varctx cx$ and $\varctx dy$ are trivial.
\end{lemma}
\begin{proof}
 We reason by induction on $n$. If $n=1$, by computing the weights on
 both sides of the equality $\ctx d{\ctx cz}=z$, we see that $\varctx
 cx$ and $\varctx dy$ are thin, but thin $1$-contexts are trivial. Let
 $n>1$ and suppose that the result holds in dimension~$n - 1$. Let
 $\varctx cx$ and $\varctx dy$ be $n$-contexts and $z$ an $n$-cell as
 in the statement, such that $\ctx d{\ctx cz}=z$. By computing the
 weights on both sides of this equality, we see that $\varctx cx$ and
 $\varctx dy$ are thin contexts. Thus,
 Lemma~\ref{lemma:deltacontext} yields $(n - 1)$\nbd-contexts
 $\dvarctx{c}{x'}$ and $\dvarctx{d}{y'}$ such that, for $z'=\sce{n-1}z$,
 $z'=\dctx{d}{\dctx{c}{z'}}$. By induction hypothesis,
 $\dvarctx{c}{x'}$ and $\dvarctx{d}{y'}$ are trivial, and so are
 $\varctx cx$ and~$\varctx dy$ by Lemma~\ref{lemma:deltacontext}.
\end{proof}

\section{Basis Uniqueness} 
\label{sec:basisunique}

Let $P$ be a polygraph and $C=\freecat P$ the free \oo-category it
generates. We shall prove that the generators of $P$ are entirely
determined by $C$ (see also~\cite[Section~4.(8.3)]{makkai2005wordcomp}).

\begin{defi}\label{defi:irreducible}
  Let $P$ be a polygraph and $n>0$. An $n$-cell $u\in \free P_n$ is
  {\em irreducible} if it is not a unit and whenever $u=v\comp i w$,
  then either $u=v$ and~$w=\Unit{n}{\tge{i}(v)}$ or $u=w$ and $v=\Unit{n}{\sce{i}(w)}$.
\end{defi}

\begin{lemma}\label{lemma:irreducible}
  An $n$-cell $u\in \free P_n$ is irreducible if and only if it is a
  generating cell of the form $u=\free a$ for $a\in P_n$.
\end{lemma}

\begin{proof}
  Suppose that $u$ is irreducible. By structural induction on $u$
  (see \cref{subsec:structural_induction}),
  either $u=\Unit nv$, which contradicts the hypothesis, or $u=\free
  a$ with $a\in \free P_n$, in which case we get the result, or
  $u=v\comp i w$. By definition, $u=v$ or $u=w$, so that the induction hypothesis
  applies and $u$ is a generating cell.

  Conversely, suppose that $u=\free a$ and $u=v\comp i w$. By
  computing the weights on both sides, we may suppose without loss of
  generality that $\wht av=1$ and~$\wht aw=0$. Then, there is a context
  $\varctx cx$ of type $x=\pair{\sce{n-1}\free a}{\tge{n-1}\free a}$
  such that $v=\ctx c{\free a}$. Let $y=\pair{\sce{n-1}v}{\tge{n-1}v}$
  and denote by $\varctx dy$ the context $\indet y\comp i w$ of type
  $y$. By substitution, $\ctx d{\ctx{c}{\free a}}=u=\free a$. By
  Lemma~\ref{lemma:composcontext}, both contexts $\varctx cx$ and
  $\varctx dy$ are trivial. In particular, $\indet y\comp i w=\indet
  y$. By repeated applications of Lemma~\ref{lemma:deltacontext}, we
  see that $\sce{k}(w)$ must be a unit cell for all $k>i$,  which
  implies $w=\Unit{n}{\tge{i}(\indet y)}=\Unit{n}{\tge{i}(v)}$. Therefore $u=v\comp i w=v$ and we are done.
\end{proof}

We may now state the main result of this section:

\begin{prop}\label{prop:uniquebasis}
  Let $P$, $Q$ two polygraphs such that there is an isomorphism
  $f:\freecat P\to\freecat Q$ in $\oCat$. Then there is a morphism
  $g:P\to Q$ in $\oPol$ such that, for each $n\in \N$, $g_n:P_n\to Q_n$ is a
  bijection, and $f=\freecat g$.
\end{prop}
\begin{proof}
  Let $a\in P_n$ be an $n$-generator. By
  Lemma~\ref{lemma:irreducible}, $\free a$ is irreducible, and as
  $f:\freecat P\to \freecat Q$ is an isomorphism, so is $f(\free
  a)$. Therefore, by Lemma~\ref{lemma:irreducible} again,  there is a
  $b=g_n(a)\in Q_n$ such that $f(\free a)=\free b$. This defines the
  required map $g_n:P_n\to Q_n$. By construction, $f=\freecat g$.
\end{proof}

\section{Rewriting Properties of \pdfm{n}-Polygraphs}
\label{sec:n-rewriting}
The theory of rewriting developed in
\cref{chap:lowdim,chap:2rewr,chap:3pol} extends to
$n$\nbd-poly\-graphs in a seamless way. Throughout this section, we fix an
$n$-polygraph~$P$, with $n\geq 2$. Recall from \cref{subsec:npolyg}, that,
for any $k\leq n$, we denote by~$\tpol kP$ the $k$-polygraph obtained by
truncating $P$ to dimension
$k$.
By definition, the set $P_n$ of $n$-generators
is a cellular extension of the category~$\freecat{\tpol {n - 1}P}$ freely generated by the $(n - 1)$-polygraph~$\tpol{n - 1}P$ and we
think of~$P$ as a presentation of the $(n - 1)$-category
\[
  \pcat P=\freecat{\tpol {n - 1}P}/P_n
\]
defined as the quotient $(n - 1)$-category, in the sense of
\secr{quotient-category}. In this setting, the elements of $P_n$ are called
\emph{rewriting rules of $P$}\index{rewriting!rule}.

\subsection{Rewriting step}
\label{sec:n-rewriting-step}
Let $x$, $y$ be two parallel $(n - 1)$-cells in $\freecat{\tpol
  {n - 1}P}$. A~{\em rewriting step of $P$} from $x$ to $y$ is an $n$-cell 
$w\in\free P_n$ with source $x$ and target $y$ of the form
  $w=\ctx c{a}$
where $a$ is a rewriting rule of $P$ and $\ctx c{\indet w}$ is a thin $n$-context
of type
$\pair{\sce{n-1}(a)}{\tge{n-1}(a)}$. From~\eqref{eq:thin_context},
each rewriting step can be expressed in a (non-unique) form as
\[
  \ctx ca= u_{n-1}\comp{n-2}(... \comp 1(u_1\comp 0 a\comp 0 v_1)\comp 1...)\comp{n-2} v_{n-1}
\]
where $u_i$ and $v_i$ are identities over $i$-dimensional cells, for every $1\leq i \leq n-1$.
We denote by $\rsteps{P}_n$ the set of rewriting steps of the polygraph $P$.
\nomenclature[P1]{$\rsteps{P}_n$}{set of rewriting steps of a polygraph $P$}

\subsection{Rewriting path}
\label{sec:n-rewriting-path}
A \emph{rewriting path of $P$}\index{rewriting!path} is a sequence
\begin{equation}
  \label{eq:rewr-path}
  \phi=w_1,...,w_k
\end{equation}
of rewriting steps of $P$ such that
$\tge{n-1}(w_i)=\sce{n-1}(w_{i+1})$, for every $1\leq i \leq k-1$. 
Therefore, a rewriting path
$\phi$ of type~\eqref{eq:rewr-path} yields an $n$-cell of $\freecat P_n$
\begin{equation}
  \label{eq:rewr-path-cell}
  \overline{\phi}=w_1\comp{n-1}...\comp{n-1}w_k.
\end{equation}
In the case where $k=0$, the rewriting path $\phi$  is said to
be \emph{empty}, and the corresponding $n$-cell $\overline{\phi}$ is
of the form $\unit{u}$ for some
$u\in\freecat P_{n-1}$. Note that for any rewriting path $ \phi=w_1,...,w_k$, the
$(n - 1)$-cells $x=\sce{n-1}(w_1)$ and $y=\tge{n-1}(w_k)$ are respectively
the source and target of $\overline\phi$, so that there is no ambiguity
in writing $\phi:x\to y$. 
As in~\cref{subsec:rewr-path}, if $\phi:x\to y$ and $\psi:y\to z$ are
rewriting paths, we denote their concatenation by $\phi\comp{}\psi$. In
particular, the sequence~(\ref{eq:rewr-path}) may be denoted
\begin{equation}
  \label{eq:rewr-concat}
 \phi= w_1\comp{}...\comp{}w_k.
\end{equation}

The following result is now an easy consequence of Proposition~\ref{prop:generating_cells}.

\begin{proposition}\label{prop:rewr-path-cell}
  For each $n$-cell $z$ in $\freecat P_n$, there is a rewriting path
  $\phi$ of $P$ such that $z=\overline\phi$.
\end{proposition}

\begin{remark}
  In Proposition~\ref{prop:rewr-path-cell}, the rewriting path $\phi$
  is not uniquely determined by $z$. However, if $\phi=w_1\comp{}...\comp{}w_k$
  and $\phi'=w'_1\comp{}...\comp{}w'_{k'}$ are two rewriting paths such that
  $\overline{\phi}=\overline{\phi'}$, then $k'=k$. Therefore, the {\em
    length} of a rewriting path $\phi$ only depends on
  $\overline{\phi}$. Moreover, given $\phi$ and $\phi'$ as above, there is
  a permutation $\sigma$ of $\set{1,...,k}$ such that for each
  $i\in\set{1,...,k}$, $w_i$ and $w'_{\sigma(i)}$ are thin contexts
  $w_i=\ctx{c}{a}$ and $w'_{\sigma(i)}=\ctx{c'}{a}$
  over the {\em same} generator $a\in P_{n-1}$.
\end{remark}

\subsection{Rewriting properties}
\index{abstract rewriting system!of an 8-polygraph@of an $n$-polygraph}
Let $P$ be an $n$-polygraph~$P$. Its two top
dimensions build a  $1$-polygraph $(\freecat P_{n-1},\rsteps{P}_n)$, \ie an abstract rewriting system, whose $0$-generators
are the $(n - 1)$-cells of~$\freecat P_{n-1}$ and $1$-generators are rewriting steps of $P$. We thus extend
the rewriting notions defined on $1$-polygraphs in \secr{ars} to $n$-polygraphs as follows.
We say that the $n$-polygraph~$P$ is

\begin{center}
  \emph{terminating} /
  \emph{quasi-terminating} /
  \emph{Church-Rosser} /
 
  \emph{confluent} / 
  \emph{locally confluent} /
  \emph{decreasing}
 
  \emph{convergent} /
  \emph{quasi-convergent}
\end{center}
 when the $1$-polygraph $(\freecat P_{n-1},\rsteps{P}_n)$
is. In particular, the following general result still applies here, see \cref{lem:newman}:

\begin{proposition}\label{prop:local-confl}
  A terminating polygraph~$P$ is confluent if and only if it is locally
  confluent.
\end{proposition}

\noindent
The above proposition naturally leads us to investigate local branchings.

\subsection{Classification of local
  branchings}\label{subsec:class-loc-branch}

A \emph{branching of $P$} is a pair $\pair{\phi}{\psi}$ of rewriting
paths $\phi:x\to y$, $\psi:x\to z$ of $P$ with common source $x$. Such a
branching is {\em local} if both $\phi$ and $\psi$ are of length~$1$,
that is, are rewriting {\em steps} of~$P$. A local branching
$\pair{\phi}{\psi}$ is
\begin{itemize}
\item {\em trivial} \index{trivial branching}
  \index{branching!trivial} if $\psi=\phi$,
\item {\em orthogonal}   \index{orthogonal branching}
  \index{branching!orthogonal} if there are rewriting steps $w:u\to v$
  and $w':u'\to v'$ such that $\overline{\phi}=w\comp{n-2}\unit{u'}$
  and $\overline{\psi}=\unit{u}\comp{n-2} w'$ (or
  $\overline{\phi}=\unit{u'}\comp{n-2}w$ and $\overline{\psi}=
  w'\comp{n-2}\unit{u}$),
  \item {\em overlapping} \index{overlapping branching}
    \index{branching!overlapping} if it is neither trivial nor orthogonal.
\end{itemize}

\subsection{Minimal branching}

The notion of {\em minimal branching}  extends in arbitrary dimension $n\geq
2$ as follows. There is a binary relation $\ctxleq$ on the set of local
branchings defined by
$\pair{\phi}{\psi}\ctxleq \pair{\phi'}{\psi'}$
if and only if there is an $(n - 1)$-context~$\varctx cy$ of type
$\pair{\sce{n-2}(\overline\phi)}{\tge{n-2}(\overline\phi)}$
(which is the same as $\pair{\sce{n-2}(\overline\psi)}{\tge{n-2}(\overline\psi)}$)
such that $\overline{\phi'}=\ictx{c}{\overline\phi}$ and $\overline{\psi'}=\ictx{c}{\overline\psi}$,
where $\ivarctx cx$ is the thin $n$-context defined
in~\cref{rk:context_promotion}. The relation
$\ctxleq$ is clearly reflexive. Antisymmetry follows from
\cref{lemma:composcontext}, whereas transitivity is a consequence
of~\eqref{eq:comp-prom-ctx}. Therefore $\ctxleq$ is a partial
order and we may define a minimal branching as a minimal element of
the set of local branchings with
respect to this order.

\subsection{Critical branching}
A branching is \emph{critical} when it is overlapping and minimal. 

\begin{lemma}
  \label{lem:n-cb}
  Given $n\geq 2$, an $n$-polygraph is locally confluent if and only if all its
  critical branchings are confluent.
\end{lemma}
\begin{proof}
  Same proof as for~\cref{lem:2-cb}.
\end{proof}

\section{Polygraphs with Finite Derivation Type}
\index{finite derivation type!8-polygraph@$n$-polygraph}
\index{finite derivation type!8-category@$n$-category}
\label{S:PolygraphsFDTn}
The property of finite derivation type has been studied in \cref{chap:2-fdt} for $2$\nbd-poly\-graphs and in \cref{chap:3coh} for $3$-polygraphs. It can be extended to any $n$-polygraphs as follows. An $n$-polygraph $P$ has \emph{finite derivation type} if it is finite and if the free $(n,n - 1)$\nbd-cate\-gory~$\freegpd{P}$ admits a
finite acyclic cellular extension, that is, a finite cellular extension generating
all $n$-spheres of $\freegpd{P}$.
An $(n - 1)$-category has \emph{finite derivation type} when it admits a presentation by an $n$-polygraph with finite derivation type.
As in the case of $1$- and $2$-categories, given two presentations of the same $(n - 1)$-category by finite $n$-polygraphs, both are of finite derivation type or neither is. The proof of this result given in~\cite[Proposition~3.3.4]{guiraud2009higher} is similar to the proof in the case of $1$- and $2$-categories given in~\cref{thm:TDFTietzeInvariant,TDFTietzeInvariant3Pol}, respectively.

\begin{proposition}
\label{TDFTietzeInvariant}
Let $P$ and $Q$ be Tietze equivalent finite $n$-polygraphs. Then the polygraph~$P$ has finite derivation type if and only if $Q$ has.
\end{proposition}

\subsection{Squier's coherence theorem}
For $n\leq 3$, we have shown that, for a convergent $n$-polygraph $P$, the set of critical branchings generate a homotopy basis of the $(n,n - 1)$-category $\freegpd{P}$, see \cref{thm:SquierHomotopicalOnePolygraphs,thm:SquierHomotopical,Theorem:SquierCompletion3polygraphs}. 
The proofs of these results extend to higher-dimensional polygraphs as follows, see~\cite[Proposition 4.3.4]{guiraud2009higher}.

\begin{theorem}
\index{Squier!theorem}
\label{Theorem:SquierCompletionNpolygraphs}
Let~$P$ be a convergent $n$-polygraph, and~$P_{n+1}$ be a cellular extension
  of the free $(n,n-1)$-category~$\freegpd{P}$. If~$P_{n+1}$ contains, for every critical
  branching~$(\phi,\psi)$ of~$P$, one $(n+1)$-generator of the form
\[
    \xymatrix@!C=1ex@!R=1ex{
      &\ar@2[dl]_\phi u\ar@2[dr]^\psi&\\
      v\ar@2[dr]_{\phi'}\ar@3[]!<5ex,0ex>;[rr]!<-5ex,0ex>^A&&\ar@2[dl]^{\psi'}w\\
      &u'&
    }
\]
  where~$\phi'$ and~$\psi'$ are $n$-cells in~$\freecat{P}_n$, then the
  $(n+1,n-1)$-polygraph $(P,P_{n+1})$ is coherent.
\end{theorem}

Following \cref{Theorem:SquierCompletionNpolygraphs}, for every $n\geq 1$, a finite convergent $n$-polygraph with a finite set of critical branchings has finite derivation type. 
A $1$-category having a finite convergent presentation therefore has finite derivation type, see Theorem~\ref{thm:FiniteConvergent=>FDT}.
Note, however, that this result fails to generalize to $n$-categories for $n\geq 2$ as shown with the following counterexample.

\subsection{A counterexample}
\label{sec:pol-pearls}
We have constructed in \cref{sec:pearls} a finite convergent
$3$-polygraph $\Pearl_3$ which does not have finite derivation type. By
shifting dimensions on the polygraph $\Pearl_3$, we obtain an $n$-polygraph $\Pearl_n$, for any $n\geq 3$. It has exactly the same cells and compositions in dimensions~$n - 3$, $n - 2$, $n - 1$, and $n$ as $\Pearl_3$ has in dimensions $0$, $1$, $2$, and $3$; on top of that, it has one cell in each dimension up to $n - 4$ and no other possible compositions, except with degenerate cells.
By construction, the polygraph $\Pearl_n$ is finite and convergent, yet it still fails to have finite derivation type.
We have thus proved that, for every natural number $n\geq 2$, there exists an $n$-category which does not have finite derivation type and admits a presentation by a finite convergent $(n + 1)$-polygraph. 
We refer the reader to \cite{guiraud2009higher} for more details on this result.

\chapter{A Catalogue of \pdfm{n}-Polygraphs}
\label{chap:pol-ex}

We have already presented a wealth of low-dimensional examples of
$n$\nbd-poly\-graphs in Parts~I, II, and~III, and many more will be
found in~\cref{chap:2ex,chap:3ex}.
The present chapter concentrates on some useful
families of $n$-polygraphs based on familiar shapes: cylinders, cubes,
and simplices, namely Street's {\em orientals} defined in the
seminal paper~\cite{street1987algebra}. These families are crucial in
the development of a homotopy theory of $\omega$-categories. 

We shall explain two methods for generating the above families. The first
one is based on a direct definition of the cylinder polygraph $\globe{1}
\tensor P$ of a polygraph~$P$. The second is based on Steiner's theory of
augmented directed complexes~\cite{steiner2004omega}, which is a very
powerful tool to build polygraphs using chain complexes. In particular, it
allows to define a tensor product for polygraphs (or even
$\omega$-categories) from which we can recover the cylinder polygraph,
as well as a join operation.

\section{First Examples of \pdfm{n}-Polygraphs}
\label{sec:pol_examples}

\subsection{Monoids}
Of particular importance is the $2$-polygraph associated to a presentation of a
monoid by generators and rewriting rules (see~\cref{chap:2rewr} and~\chapr{2ex}). Recall
that given $M$ a monoid presented by a set of generators~$P_1$ and a set $P_2$
of rewriting rules of the form $\alpha:w\To w'$ where $w,w'\in P_1^*$, we get a
$2$-polygraph $P=\Pres{P_0}{P_1}{P_2}$, where $P_0=\set{\star}$:
\begin{displaymath}
  \vxym{P_0 \ar[d]^{\ins0}& \ar@<-.5ex>[dl]_(.4){\sce
    0}\ar@<.5ex>[dl]^(.4){\tge{0}}P_1\ar[d]^{\ins1} & \ar@<-.5ex>[dl]_(.4){\sce
    1}\ar@<.5ex>[dl]^(.4){\tge{1}}P_2.\\
 P_0^* & P_1^* \ar@<-.5ex>[l]_-{\ssce{0}{*}}\ar@<.5ex>[l]^-{\ttge{0}{*}}& }
\end{displaymath}
The source and target maps are defined by $\sce{0}(a)=\tge{0}(a)=\star$ for each $a\in P_1$,
and $\sce{1}(\alpha)=w$, $\tge{1}(\alpha)=w'$ for each rewrite rule
$\alpha:w\To w'$.
Remark that in this case we recover $M$ as the quotient category $C/X$ where $C$
is the free category $\xymatrix{P_0^* & \doubl{}{}P_1^*}$ and $X=P_2$.

\subsection{Terminal polygraph}
\label{paragr:pol_terminal}
The functor $\catpol:\oCat \to\oPol$ from~\ref{subsec:omegapol} also provides
many natural examples of polygraphs. In particular $\termpol=\catpol(\termcat)$
where $\termcat$ is the terminal $\omega$-category has infinitely many
generators in all dimensions $k\geq 2$.  It is in fact a terminal object in the
category $\oPol$, as immediately implied by the adjunction.

\subsection{Globular sets}
Each $n$-globular set $X$ may be seen as an $n$-polygraph. Precisely, for each
$n\in\N$ there is a natural inclusion functor $\nGlob n\to \nPol n$ taking the
globular set $X$ to the polygraph $P$ whose set of $k$-generators is just
$P_k=X_k$ for $0\leq k\leq n$. For example, the $n$-globe $\globe{n}$ is a
polygraph with exactly two generators in dimensions $0\leq k<n$ and a unique
generator in dimension $n$. Likewise, the $n$-spheres $\sphere{n}$ provide
another example of a series of $n$-polygraphs.

\section{Tensoring Polygraphs by \texorpdfstring{$\globe{1}$}{O1}}
\label{sec:tensor}
\begin{paragr}
Before we turn to more complex examples of polygraphs, we need to 
describe a useful general construction, namely the {\em tensor
  product} of a polygraph by the $1$-globe $\globe{1}$. In order to
check the coherence of the following construction, we shall rely on the
existence and the properties of the functor $\CatCyl:\oCat\to\oCat$
taking an \oo-category $C$ to the \oo-category $\CatCyl C$ of ``small
cylinders internal to $C$''. We refer to \cref{chap:w-eq} for a
detailed account of this functor $\CatCyl$. As for the
present argument, we just need to know that $\CatCyl$ comes equipped
with natural transformations
\[ \projcyltop, \projcylbot : \CatCyl \to \id_{\ooCat} \]
and maps
$\ppal{-}:(\CatCyl C)_n\to C_{n+1}$ taking an $n$-cylinder $u\in
(\CatCyl C)_n$ to its {\em principal cell} $\ppal u\in C_{n+1}$, such
that, for any $n$-cylinder $u\in(\CatCyl C)_n$, the source and target
of $\ppal u$ are given by
\[
\ppal{u} : \projcyltop_C(u)\comp 0 \ppal{\tge{0}(u)}\comp 1\cdots
\comp{n-1}\ppal{\tge{n-1}(u)}\to
\ppal{\sce{n-1}(u)}\comp{n-1}\cdots\comp{1}\ppal{\sce{0}(u)}\comp{0} \projcylbot_C(u).
\] 

Let us first
introduce a few notations to be used throughout the construction.
\begin{itemize}
\item The $0$-generators of $\globe{1}$ will be denoted by $\sgl$ and $\tgl$.
\item The $1$-generator of $\globe{1}$ will be denoted by $\pgl$, so that
  $\pgl:\sgl\to\tgl$.
\item For any set $A$ of symbols, we denote by $\sgltens{A}$ (\resp
  $\tgltens{A}$, $\pgltens{A}$) the set of all symbols of the form $\sgltens{a}$
  (resp.\ $\tgltens{a}$, $\pgltens{a}$), where $a\in A$.
\item Given a symbol $a$, we define the following formal expressions:
  \begin{align*}
  \laxsce a0 &=\sgltens{a},\\
  \laxtge a0& =\tgltens{a} ,\\
  \laxsce a1&=(\sgltens{a})\comp 0 (\pgltens{\tge{0}(a)}),\\
  \laxtge a1&= (\pgltens{\sce{0}(a)})\comp 0 (\tgltens{a}) ;
  \end{align*}
and more
generally, for each integer $i>1$: 
  \begin{align*}
\laxsce{a}{i} & = \Laxsce{a}{i}, \\
\laxtge{a}{i} & = \Laxtge{a}{i}.
\end{align*}
These expressions will eventually denote actual cells whenever the indeterminate $a$
is interpreted appropriately.
\end{itemize}
\end{paragr}

\begin{paragr}
Let now $P$ be a polygraph. We shall build  a new
polygraph $Q=\globe{1}\tensor P$ endowed with a morphism
$h:\free P\to \CatCyl{\free Q}$
giving rise to the following families of maps:
\begin{itemize}
\item morphisms $\tinc,\binc:P\to Q$ in $\oPol$ respectively taking an $n$-cell
  $u\in \free P_n$ to 
$\tinc(u)=\sgltens{u}=\projcyltop_{\free Q}(h(u))$
and
$\binc(u)=\tgltens{u}=\projcylbot_{\free Q}(h(u))$,
\item  maps $\pgltens{-}:\free P_{n}\to \free Q_{n+1}$ defined as the
  composite 
\[
\xymatrix{\free P_n \ar[r]^{h_n}& \CatCyl{\free Q}_n\ar[r]^{\ppal{-}} & \free Q_{n+1}}
\]
\end{itemize}
such that 
  \begin{equation}
    \label{eq:laxboundaries}
\pgltens{u}:\laxsce{u}{n}\to \laxtge{u}{n}    
  \end{equation}
for any $n$-cell $u\in \free P_n$. 
\end{paragr}

\begin{paragr}
We now define $Q$, $\binc$, $\tinc$, and $\pgltens{-}$ by
simultaneous induction on the dimension.
\begin{itemize}
\item For $n=0$, the set of $0$-generators of $Q$ is defined by 
$Q_0=\sgltens{P_0} \uplus \tgltens{P_0}$
whereas $\tinc_0$ and $\binc_0$ are given by $\tinc_0(a)=\sgltens{a}$
and $\binc_0(a)=\tgltens{a}$ for~$a\in P_0$.
\item The set of $1$-generators of $Q$ is defined by
  $Q_1=\sgltens{P_1} \uplus \tgltens{P_1} \uplus \pgltens{P_0}$,
whereas $\tinc_1,\binc_1:P_1\to Q_1$ and $\pgltens{-}:P_0\to Q_1$ are given by
the obvious canonical inclusions. The source and target maps 
$\sce{0}^Q,\tge{0}^Q:Q_1\to\free Q_0$ are given, for each $a\in P_1$
by
$\sce{0}^Q(\sgltens{a}) = \sgltens{\sce{0}^P(a)} $,
$\sce{0}^Q(\tgltens{a}) = \tgltens{\sce{0}^P(a)} $,
$\tge{0}^Q(\sgltens{a}) = \sgltens{\tge{0}^P(a)} $ and $\tge{0}^Q(\tgltens{a}) = \tgltens{\tge{0}^P(a)} $,
so that $\tinc$ and $\binc$ are as expected morphisms up to dimension
1. As for $a\in P_0$, we define $\sce{0}^Q(\pgltens{a}) = \sgltens{a}$ and $\tge{0}^Q(\pgltens{a}) = \tgltens{a}$.
Because $\free P_0=P_0$ there is nothing more to verify here,
and~(\ref{eq:laxboundaries}) holds.
\item Suppose now $Q$ has been defined up to dimension $n$, as well as morphisms
$\tinc,\binc:P\to Q$ in $\oPol$ up to dimension $n$, a morphism $h:\free P\to
\CatCyl{\free Q}$ in~$\oCat$ up to
dimension $n-1$ and corresponding maps $\pgltens{-}:\free P_k\to \free Q_{k+1}$
satisfying equations~(\ref{eq:laxboundaries}) for all
$k\in\set{0,...,n-1}$.
The set of $(n+1)$\nbd-gener\-ators of $Q$ is defined by
$Q_{n+1}=\sgltens{P_{n+1}}\uplus\tgltens{P_{n+1}}\uplus\pgltens{P_n}$
whereas $\tinc_{n+1},\binc_{n+1}:P_{n+1}\to Q_{n+1}$ and $\pgltens{-}:P_n\to Q_{n+1}$
are given by the obvious canonical inclusions. The source and target
maps are defined as above for the generators of the form $\sgltens{a}$
and $\tgltens{a}$, namely:
\begin{align*}
  \sce{n}^Q(\sgltens{a}) & = \sgltens{\sce{n}^P(a)}, \\
\sce{n}^Q(\tgltens{a}) & = \tgltens{\sce{n}^P(a)}, \\
\tge{n}^Q(\sgltens{a}) & = \sgltens{\tge{n}^P(a)}, \\
\tge{n}^Q(\tgltens{a}) & = \tgltens{\tge{n}^P(a)}.
\end{align*}
As for generators of the form $\pgltens{a}$, with $a\in P_n$, the
induction hypothesis implies that $\laxsce{a}{n}$ and $\laxtge{a}{n}$
denote well-defined cells in $\free Q_n$ and moreover that these two
cells are parallel. Hence, the $n$-source and $n$-target of
$\pgltens{a}$ may be defined by
$\sce{n}^Q(\pgltens{a}) = \laxsce{a}{n}$ and $\tge{n}^Q(\pgltens{a}) = \laxtge{a}{n}$,
whence
$\pgltens{a}:\laxsce{a}{n}\to\laxtge{a}{n}$.
It follows that the polygraph $Q$ is now defined up to dimension
$n+1$, together with the morphisms $\tinc,\binc:P\to Q$, and
that~(\ref{eq:laxboundaries}) holds for the $(n+1)$-{\em generators} of
$P$. This implies that the morphism $h$ can be extended in dimension
$n$ by a map $h_{n}:P_n\to (\CatCyl{\free Q})_n$ commuting to source
and target maps. By the universal property of polygraphs
(Lemma~\ref{lemma:ext-universal}), $h$ extends as a morphism $\free
P\to \CatCyl{\free Q}$ up to dimension $n$, and by composing with the
principal cell map $\ppal{-}:(\CatCyl{\free Q})_n\to \free Q_{n+1}$,
we get a map
\[
\pgltens{-}:\free P_n\to\free Q_{n+1}
\]
satisfying~(\ref{eq:laxboundaries}) for all cells $u\in \free
P_n$. Thus, the induction is complete. 
\end{itemize}
\end{paragr}

\section{Families of Polygraphs}\label{sec:families}

\subsection{Cylinders}
\nomenclature[Cyln]{$\Cyl{n}$}{polygraphic $n$-cylinder}
\index{cylinders (as polygraphs)}
Let $n\geq 0$. The {\em free-standing $n$-cylinder}  is by
definition the polygraph
\[
\Cyl n= \globe{1}\tensor\globe{n}.
\]
Now recall that the $n$-globe $\globe{n}$ has $2n+1$ generators,
namely the only $n$\nbd-gener\-ator~$\gener{n}$ together with the $2n$ generators of
the form $\gener{i}^{-}=\sce{i}(\gener{n})$ and
$\gener{i}^{+}=\tge{i}(\gener{n})$ for $i\in\set{0,...,n-1}$. Thus,
$\Cyl 0$ is just $\globe 1$, whereas for $n>0$ 
the generators of~$\Cyl n$ are listed below.
\begin{itemize}
\item There are four $0$-generators, namely $\sgltens{\gener{0}^{-}}$,
  $\sgltens{\gener{0}^{+}}$, $\tgltens{\gener{0}^{-}}$, and
  $\tgltens{\gener{0}^{+}}$.
\item For $0<i\leq n-1$, there are six $i$-generators, namely
  $\sgltens{\gener{i}^{-}}$, $\sgltens{\gener{i}^{+}}$,
  $\tgltens{\gener{i}^{-}}$, $\tgltens{\gener{i}^{+}}$,
  $\pgltens{\gener{(i{-}1)}^{-}}$, and
  $\pgltens{\gener{(i{-}1)}^{+}}$.
\item There are four $n$-generators, namely $\sgltens{\gener{n}}$,
  $\tgltens{\gener{n}}$, $\pgltens{\gener{(n{-}1)}^{-}}$, and $\pgltens{\gener{(n{-}1)}^{+}}$.
\item There is only one $(n+1)$-generator $\pgltens{\gener{n}}$.
\end{itemize}
Therefore $\Cyl n$ has exactly $6n+3$ generators. Moreover, the source
and target of these generators are given by the
formulas~(\ref{eq:laxboundaries}).
Here are pictures of $\Cyl{0}$, $\Cyl{1}$, and $\Cyl{2}$:
  \[
   \xymatrix@R=3pc{
   \gener{0}^- \ar[d]_{\gener{1}} \\
   \gener{0}^+
   \pbox,
   }
   \qquad
   \qquad
    \xymatrix@C=3pc@R=3pc{
      \gener{0}^- \otimes \gener{0}^- \ar[r]^{\gener{0}^- \otimes \gener{1}}
      \ar[d]_{\gener{1} \otimes \gener{0}^-} &
      \gener{0}^- \otimes \gener{0}^+ \ar[d]^{\gener{1} \otimes \gener{0}^+} \\
      \gener{0}^+ \otimes \gener{0}^- \ar[r]_{\gener{0}^+ \otimes \gener{1}}
      y & \gener{0}^+ \otimes \gener{0}^+
      \ar@{}[u];[l]_(.30){}="s"
      \ar@{}[u];[l]_(.70){}="t"
      \ar@2"s";"t"_{\gener{1} \otimes \gener{1}}
      \pbox{,}
    }
    \]
    \[
    \xymatrix@C=5pc@R=5pc{
      \gener{0}^- \otimes \gener{0}^-
      \ar@/^3ex/[r]^(0.70){\gener{0}^- \otimes \gener{1}^-}_{}="0"
      \ar@/_3ex/[r]_(0.70){\gener{0}^- \otimes \gener{1}^+}_{}="1"
      \ar[d]_{}="f"_{\gener{1} \otimes \gener{0}^-}
      \ar@2"0";"1"_{\gener{0}^- \otimes \gener{2}\,\,}
      &
      \gener{0}^- \otimes \gener{0}^+
      \ar[d]^{\gener{1} \otimes \gener{0}^+} \\
      \gener{0}^+ \otimes \gener{0}^-
      \ar@{.>}@/^3ex/[r]^(0.30){\gener{0}^+ \otimes \gener{1}^-}_{}="0"
      \ar@/_3ex/[r]_(0.30){\gener{0}^+ \otimes \gener{1}^+}_{}="1"
      \ar@{:>}"0";"1"_{\gener{0}^+ \otimes \gener{2}\,\,}
      &
      \gener{0}^+ \otimes \gener{0}^+
      \ar@{}[u];[l]_(.40){}="x"
      \ar@{}[u];[l]_(.60){}="y"
      \ar@<-1.5ex>@/_1ex/@{:>}"x";"y"_(0.60){\gener{1} \otimes \gener{1}^-\,}_{}="0"
      \ar@<1.5ex>@/^1ex/@2"x";"y"^(0.40){\!\gener{1} \otimes \gener{1}^+}_{}="1"
      \ar@{}"1";"0"_(.05){}="z"
      \ar@{}"1";"0"_(.95){}="t"
      \ar@3{>}"z";"t"_{\gener{1} \otimes \gener{2}}
      \pbox{.}
    }
  \]

\subsection{Cubes}
\nomenclature[Cubn]{$\cub{n}$}{polygraphic $n$-cube}
\index{cubes (as polygraphs)}
The tensor product construction above also leads to the definition of
the {\em polygraphic $n$-cubes}. Precisely, this 
family $(\cub n)_{n\in \N}$ of  polygraphs is defined by
\begin{itemize}
\item $\cub 0=\globe{0}$,
\item $\cub{n+1}=\globe{1}\tensor \cub n$.
\end{itemize}
It follows from the above construction that, for each $i\in\set{0,...,n}$, the set of
$i$-generators of $\cub n$ has  $\binom{n}{i} 2^{n-i}$ elements, and
therefore the total number of generators in $\cub n$ is
\[
\sum_{i=0}^n\binom{n}{i}2^{n-i}= 3^n.
\] 
A convenient way to encode these generators is by labeling them by
words on the alphabet $\set{-,1,+}$, the $i$-generators being those
with exactly $i$ occurrences of the letter $1$. For example, $\cub 2$
looks like
\[
\xymatrix{ {\scriptstyle [--]}\ar[r]^{[-1]}\ar[d]_{[1-]} &
  {\scriptstyle [-+]}\ar[d]^{[1+]}\ar@{=>}[ld]|{[11]}\\
                {\scriptstyle [+-]} \ar[r]_{[+1]}& {\scriptstyle [++]}\pbox{,}}
\]
whereas $\cub 3$ looks like
\[
\xymatrix@R+1.1pc{
{\scriptstyle[---]} \ar[rr]|{[--1]}\ar[rd]|{[-1-]}\ar[dd]_{[1--]}& & {\scriptstyle[--+]}\ar@/_4ex/@{:>}[lldd]\ar@{=>}[ld]|{[-11]}\ar[rd]|{[-1+]}\ar@{.>}[dd]|(.3){[1-+]} & \\
& {\scriptstyle[-+-]} \ar@{=>}[ld]|{[11-]}\ar[rr]|(.3){[-+1]}\ar[dd]|(.7){[1+-]}&& {\scriptstyle[-++]} \ar@{:>}[ld]\ar@/^4ex/@{=>}[lldd]_(.7){[1+1]}\ar[dd]|{[1++]} \\
{\scriptstyle[+--]}\ar[rd]|{[+1-]}\ar@{.>}[rr]|(.3){[+-1]} & & {\scriptstyle[+-+]}\ar@3{.>}[ul]|{[111]}\ar@{:>}[ld]\ar@{.>}[rd]|{[+1+]} & \\
& {\scriptstyle[++-]} \ar[rr]|{[++1]}& & {\scriptstyle[+++]} \pbox{,}
  }
\]
where the $3$-cell $[111]$ goes from the composition of the front faces
\[
  \sce{2}[111]=([-11]\comp 0[1{+}{+}])\comp 1([-1-]\comp 0 [1{+}1])\comp
  1([11-]\comp 0[{+}{+}1])
\]
to the composition of the back faces
\[
  \tge{2}[111]=([{-}{-}1]\comp 0 [11+])\comp 1 ([1{-}1]\comp 0[+1+]) \comp 1 ([1{-}{-}]\comp
  0 [+11]).
  \]

\subsection{Simplices}
\label{paragr:simpl}

The correspondence $P\mapsto \globe{1}\tensor P$ from $\oPol$ to $\oPol$
is easily seen to be functorial. Now, for each polygraph $P$, we can
form the following pushout in $\oCat$:
\[
\xymatrix{\free P\ar[d]_{!}\ar[r]^-{\tinc} & \free{(\globe 1\tensor P)}\ar[d]\\
\termcat \ar[r]& C(P) \pbox.
}
\]
This defines a functor $P\mapsto C(P)$ from $\oPol$ to
$\oCat$. Because $\tinc$ is a cofibration of $\oCat$
(see~\cref{paragr:cof_triv_fib}), $\termcat\to C(P)$ is also a
cofibration and as $\termcat$ is cofibrant, so is $C(P)$. By
Theorem~\ref{thm:cofib}, there is a polygraph $S(P)$ such that
$C(P)=S(P)^*$. It is now possible to define a family
$(\simpl{n})_{n\in\N}$ of polygraphs by
\nomenclature[On]{$\simpl{n}$}{$n$-th oriental}
\index{oriental}
\begin{itemize}
\item $\simpl{0}=\globe{0}$,
\item $\simpl{n+1}=S(\simpl{n})$.
\end{itemize}
 It turns out that $\simpl{n}$ is precisely the polygraphic $n$-th
 simplex, or {\em $n$-th oriental}, first
 introduced in~\cite{street1987algebra}. For each $0\leq i\leq n$ the set of
$i$-generators (``$i$-faces'') of $\simpl n$  has $\binom{n+1}{i+1}$
elements. The $i$-generators of $\simpl n$  may be conveniently
encoded by strictly increasing sequences of integers
$\smp{n_0,...,n_i}$, where $0\leq n_0<n_1<...<n_i\leq n$. For example,
$\simpl 2$ and $\simpl 3$ may be pictured as follows:

\[
  \begin{xy}
\xymatrix @=1.5em{
&\overset{\smp{0}}{\cdot} 
           \ar[ldd]_{\smp{0,1}} 
            \ar[rdd]^{\smp{0,2}}
         &
        \ar@2{->}(13,-11)*{};(5,-18)*{}|{\scriptstyle{\smp{0,1,2}}}\\
&&\\
\underset{\makebox[0pt]{$\scriptstyle{\smp{1}}\ \ \ $}}{\cdot} 
\ar[rr]_{\smp{1,2}}
&&
\underset{\makebox[0pt]{$\ \ \ \scriptstyle{\smp{2}}$}}{\cdot}  }
  \end{xy}
\hspace{4em}
  \begin{xy}
    \xymatrix @=1.2em{
&&\overset{\smp{0}}{\cdot}
\ar[lldddd]_{\smp{0,1}}
\ar[rrrddd]^{\smp{0,3}}
\ar[rddddd]|(.2){\smp{0,2}}
&&& \\
&&&&& \\
&&&
\ar@2{.>}@<-5pt>[ld]_(.7)*-<3pt>{\scriptstyle{\smp{0,1,3}}}
\ar@2{->}@<3pt>[d]|{\smp{0,2,3}}
&& \\
&&&
\ar@3{.>}@<4pt>[l]^(.4){\smp{0,1,2,3}}
&&{\hbox to 2mm{\hfill\vbox to 2mm{}}\cdot\scriptstyle{\smp{3}}}
\\
*{\scriptstyle{\smp{1}}\textstyle{\cdot}\hbox to 3mm{\hfill\vbox to 4mm{}}}
\ar[rrrd]_{\smp{1,2}}
\ar@{.>}[rrrrru]|(.75){\smp{1,3}}&&
\ar@2{.>}[rd]!<-3ex,4ex>^*-<3pt>{\scriptstyle{\smp{1,2,3}}}
\ar@2{->}@<-6pt>[l]_*-<3pt>{\scriptstyle{\smp{0,1,2}}}&&& \hbox to 2mm{\hfill\vbox to 10mm{}}\\
&&&\underset{\smp{2}}{\cdot}
\ar[rruu]_{\smp{2,3}}&&
}
\end{xy}
  \]


\section{Construction of Polygraphs via Steiner's Theory}
\label{sec:steiner}

We shall now present Steiner's theory of augmented directed
complexes~\cite{steiner2004omega} and use it to construct some
polygraphs. Other similar formalisms include Street's \emph{parity
complexes}~\cite{street1991parity,street1994parity} and Johnson's
\emph{pasting schemes}~\cite{johnson1989combinatorics}.

\begin{paragraph}[Augmented directed complexes]
  In this section, by ``chain complex'' we will always mean ``chain complex of
  abelian groups in non-negative degree''
  (see~\cref{paragr:def_chain_complex}). Recall that an
  \ndef{augmented chain complex}\index{augmented chain complex} $(K, d, e)$
  is a chain complex~$(K, d)$ endowed with an
  \ndef{augmentation}~$e$, that
  is, with a map of abelian groups $e : K_0 \to \Z$ such that $ed_1 = 0$. A
  morphism from an augmented chain complex $(K, d_K, e_K)$ to a second
  augmented chain complex $(L, d_L, e_L)$ is a morphism of chain complex $f$
  from $(K, d_K)$ to~$(L, d_L)$ such that $e_Lf_0 = e_K$.

  \index{augmented directed complex}\index{ADC (augmented directed complex)}
  An \emph{augmented directed complex} is an augmented chain complex $(K, d,
  e)$ equipped with, for every $n \ge 0$, a submonoid $\adcmon{K_n}$ of
  \ndef{positive chains}\index{positive chains of an augmented directed
  complex} of the abelian group~$K_n$. A morphism $f:K\to L$
  between two augmented directed complexes consists of a morphism of
  the underlying augmented chain complexes that respects the positive chains
  in the sense that $f(\adcmon{K_n})\subseteq\adcmon{L_n}$ for every $n \ge
  0$. We will denote by $\ADC$ the resulting category.
  \nomenclature[ADC]{$\ADC$}{category of augmented directed complexes}
\end{paragraph}

\subsection{From $\omega$-categories to augmented directed complexes}
  \nomenclature[.lambda]{$\lambda(C)$}{augmented directed complex associated
to an $\omega$-category $C$}
We define a functor $\lambda:\oCat\to\ADC$ by sending
an $\omega$-category $C$ to the following augmented directed complex
$\lambda(C)$:
\begin{itemize}
  \item For every $n \ge 0$, the abelian group $(\lambda(C))_n$ is generated
    by elements $\adcgen{x}$, for $x$ an $n$-cell of $C$, subject to the
    relations
  \[
    \adcgen{x\comp i y}
    =
    \adcgen{x}+\adcgen{y}
  \]
  for pairs of $i$-composable $n$-cells~$x$ and~$y$.
  \item For every $n \ge 1$ and every $n$-cell $x$ of $C$, we set
  \[
    d_n(\adcgen{x})
    =
    \adcgen{\tge{n-1}(x)}-\adcgen{\sce{n-1}(x)}.
  \]
  \item If $x$ is a $0$-cell of $C$, we set
  \[
    e(\adcgen{x})
    =
    1.
  \]
  \item Finally, for every $n \ge 0$, the submonoid $\adcmon{(\lambda(C))_n}$
    is the submonoid of $(\lambda(C))_n$ generated by the elements
    $\adcgen{x}$ for $x\in C_n$.
\end{itemize}
The globular relations easily imply that $d$ is indeed a differential so
that $\lambda(C)$ is indeed an augmented directed complex.

Given an \oo-functor $f:C\to D$, the morphism of augmented directed complexes
$\lambda(f) : \lambda(C) \to \lambda(D)$ is defined on generators by
\[
  (\lambda(f))(\adcgen{x})
  =
  \adcgen{f(x)}.
\]

\subsection{From chain complexes to $\omega$-categories}
\label{paragr:def_mu}
\nomenclature[.muK]{$\mu(K)$}{$\omega$-category associated to a chain
complex $K$}
 
We now define a functor $\mu:\pCh[\Z] \to \oCat$, where $\pCh[\Z]$ denotes
the category of chain complexes. Given a chain complex $K$, the associated
$\omega$-category $\mu(K)$ is defined as follows:
\begin{itemize}
  \item For $n \ge 0$, an $n$-cell of $\mu(K)$ consists of a table
  \[
      x =
      \begin{pmatrix}
        x_0^- & \cdots & x_n^- \\
        \noalign{\vskip 3pt}
        x_0^+ & \cdots & x_n^+ \\
      \end{pmatrix}
      ,
  \]
  where
  \begin{itemize}
    \item $x_i^\e$, for $0 \le i \le n$ and $\e = \pm$, belongs to $K_i$,
    \item $x_n^- = x_n^+$,
    \item $d(x^\e_i) = x_{i-1}^+ - x_{i-1}^-$ for $1 \le i \le n$ and $\e =
      \pm$.
  \end{itemize}
  \item For $n \ge 1$, the source and target of such a cell $x$ are given
    by
  \[
      \sce{}(x) =
      \begin{pmatrix}
        x_0^- & \cdots & x_{n-1}^- \\
        \noalign{\vskip 3pt}
        x_0^+ & \cdots & x_{n-1}^- \\
      \end{pmatrix}
      \qquad\text{and}\qquad
      \tge{}(x) =
      \begin{pmatrix}
        x_0^- & \cdots & x_{n-1}^+ \\
        \noalign{\vskip 3pt}
        x_0^+ & \cdots & x_{n-1}^+ \\
      \end{pmatrix}
      .
  \]
  \item For $n \ge 0$, the unit cell of such a cell $x$ is given by
  \[
    \unit{x} =
      \begin{pmatrix}
        x_0^- & \cdots & x_{n}^- & 0 \\
        \noalign{\vskip 3pt}
        x_0^+ & \cdots & x_{n}^+ & 0 \\
      \end{pmatrix}.
  \]
  \item Finally, if
  \[
    x =
      \begin{pmatrix}
        x_0^- & \cdots & x_n^- \\
        \noalign{\vskip 3pt}
        x_0^+ & \cdots & x_n^+ \\
      \end{pmatrix}
    \qquad\text{and}\qquad
    y =
      \begin{pmatrix}
        y_0^- & \cdots & y_n^- \\
        \noalign{\vskip 3pt}
        y_0^+ & \cdots & y_n^+ \\
      \end{pmatrix}
  \]
  are two $n$-cells of $\mu(K)$ such that the $\tge{i}(x) = \sce{i}(y)$ for
  some $0 \le i < n$, then
  \[
      x \comp{i} y =
      \begin{pmatrix}
        x_0^- & \cdots & x_i^- & x_{i+1}^- + y_{i+1}^- & \cdots & x_n^- + y_n^- \\
        \noalign{\vskip 3pt}
        y_0^+ & \cdots & y_i^+ & x_{i+1}^+ + y_{i+1}^+ & \cdots & x_n^+ + y_n^+ \\
      \end{pmatrix}.
  \]
\end{itemize}
One checks that $\mu(K)$ is indeed an \oo-category.

If $f : K \to L$ is a morphism of chain complexes, then the action on
$n$-cells of the \oo-functor $\mu(f)$ is defined by
  \[
      \begin{pmatrix}
        x_0^- & \cdots & x_n^- \\
        \noalign{\vskip 3pt}
        x_0^+ & \cdots & x_n^+ \\
      \end{pmatrix}
      \mapsto
      \begin{pmatrix}
        f(x_0^-) & \cdots & f(x_n^-) \\
        \noalign{\vskip 3pt}
        f(x_0^+) & \cdots & f(x_n^+) \\
      \end{pmatrix}.
  \]

\begin{remark}
  The construction $\mu$ of the previous paragraph actually lands into the
  category $\ooCat(\Ab)$ of \oo-categories internal to abelian groups, that
  is, of \oo-categories whose set of $n$-cells is endowed with a structure
  of abelian group, and whose operations (sources, target, units, and
  compositions) are compatible with these structures of abelian groups on
  cells. More precisely, the functor $\mu$ naturally lifts to a functor
  $\pCh[\Z] \to \ooCat(\Ab)$. A theorem of Bourn~\cite{bourn1990another}
  states that this functor is an equivalence of categories. This is sometimes
  called the \emph{globular Dold-Kan correspondence}.
\end{remark}

\subsection{From augmented directed complexes to $\omega$-categories}
\label{paragr:def_nu}
\nomenclature[.nuK]{$\nu(K)$}{$\omega$-category associated to an augmented
directed complex}

We now define a functor $\nu:\ADC\to\oCat$ as a subfunctor of the
functor
\[ \ADC \xto{U} \pCh[\Z] \xto{\mu} \oCat, \]
where $U$ denotes the obvious forgetful functor. If $K$ is an augmented
directed complex, an $n$-cell
\[
  \begin{pmatrix}
    x_0^- & \cdots & x_n^- \\
    \noalign{\vskip 3pt}
    x_0^+ & \cdots & x_n^+ \\
  \end{pmatrix}
\]
of $\mu(K)$ belongs to the sub-\oo-category $\nu(K)$ if
\begin{itemize}
  \item $x_i^\e$, for $0 \le i \le n$ and $\e = \pm$, is a positive chain
    (that is, is in $\adcmon{K_i}$),
  \item $e(x^\e_0) = 1$ for $\e = \pm$.
\end{itemize}
One checks that the operations of the \oo-category $\mu(K)$ restrict to
$\nu(K)$ and that if $f : K \to L$ is a morphism of augmented directed
complexes, then the \oo-functor $\mu(f) : \mu(K) \to \mu(L)$ restricts to an
\oo-functor $\nu(f) : \nu(K) \to \nu(L)$.

\begin{proposition}\label{prop:adj_lambda_nu}
  The functors
  \[
    \lambda : \oCat \to \ADC
    \qquad
    \qquad
    \nu : \ADC \to \oCat
  \]
  define a pair of adjoint functors.
\end{proposition}

\begin{proof}
  We will only define the components of the adjunction morphisms and leave
  the verification to the reader. If $C$ is an \oo-category, the unit of the
  adjunction at $C$ is the \oo-functor $\eta_C : C \to
  \nu(\lambda(C))$ defined on $n$-cells by
  \[
  x
  \mapsto
  \begin{pmatrix}
    \adcgen{\sce{0}(x)} & \cdots & \adcgen{\sce{n}(x)} \\
    \noalign{\vskip 3pt}
    \adcgen{\tge{0}(x)} & \cdots & \adcgen{\tge{n}(x)} \\
  \end{pmatrix}.
  \]
  If $K$ is an augmented chain complex, the counit of the adjunction $K$ is
  the morphism $\varepsilon_K : \lambda(\nu(K)) \to K$
  defined on generating $n$-chains by
  \[
    \adcgenbig{
    \begin{pmatrix}
      x_0^- & \cdots & x_n^- \\
      \noalign{\vskip 3pt}
      x_0^+ & \cdots & x_n^+ \\
    \end{pmatrix}
    }
    \mapsto x_n,
  \]
  where $x_n$ denotes the $n$-chain $x_n^- = x_n^+$.
\end{proof}

\begin{remark}\label{rem:adj_mu}
  Similarly, the functors
  \[
    U\lambda : \oCat \to \pCh[\Z]
    \qquad
    \qquad
    \mu : \pCh[\Z] \to \oCat,
  \]
  where $U : \ADC \to \pCh[\Z]$ denotes the forgetful functor,
  define a pair of adjoint functors.
\end{remark}

\begin{paragraph}[Augmented directed complexes with basis]
  \index{basis!augmented directed complex}
  \index{augmented directed complex!basis}
A \ndef{basis}~$B$ of an augmented directed complex is a sequence of
subsets $B_n\subseteq K_n$, indexed by $n \ge 0$, such that
\begin{itemize}
\item $B_n$ is a basis of the $\Z$-module $K_n$,
\item $B_n$ generates $\adcmon{K_n}$ as a submonoid of $K_n$.
\end{itemize}
The data of such a basis gives, for every $n \ge 0$, an isomorphism of
abelian groups between $K_n$ and $\Z^{(B_n)}$ restricting to an isomorphism
of monoids between $\adcmon{K_n}$ and~$\N^{(B_n)}$. The elements of $B_n$
can be recovered as the minimal non-zero elements of $\adcmon{K_n}$ for the
order on $K_n$ defined by $x \le y$ if $y - x \in \adcmon{K_n}$. This shows
that if such a basis exists it is unique.
\end{paragraph}

\begin{prop}\label{prop:lambda_pol}
  Let $P$ be a polygraph. Then, for each $n\geq 0$,  the family
  $([x])_{x\in P_n}$ forms a basis of $\lambda(P^\ast)$.
\end{prop}

\begin{proof}
  As any $n$-cell of $P^\ast$ can be expressed as a composition of units and
  $n$-generators  (see Proposition~\ref{prop:generating_cells}),
the set $P_n$ generates the monoid $\lambda(P^\ast)_n^+$. Let us show
  that it is a $\Z$-basis of $\lambda(P^\ast)_n$.
  We have to prove that the morphism $\beta : \freemod\Z{P_n} \to
  \lambda(P^\ast)_n$ which, given $x$ in $P_n$, sends $[x]$
  in~$\freemod\Z{P_n}$ to $[x]$ in $\lambda(P^\ast)_n$ is an
  isomorphism.
  By \ref{sec:linearization}, we have a map $\gamma :
  \free{P}_n \to \freemod\Z{P_n}$ sending an $n$-generator $x$
  of $P$ to $[x]$ and compositions to sums. We thus get a morphism $\bar{\gamma}
  : \lambda(\free{P})_n \to \freemod\Z{P_n}$. We claim that
  $\bar{\gamma}$ is an inverse of $\beta$. As $\freemod\Z{P_n}$ is generated
  by the $[x]$, it suffices to check equality $\bar{\gamma} \circ \beta =
  \id{}$ on these elements, which is true by definition. As
  $\lambda(P^\ast)_n$ is generated by the $[x]$, we similarly get the
  equality~$\beta \circ \bar{\gamma} = \id{}$.
\end{proof}

\begin{paragraph}
  \index{support}
  \nomenclature[supp]{$\supp{x}$}{support of a cell~$x$}
Let $K$ be an augmented directed complex with a basis $B$. Let $n \ge 0$ and
let $x$ be an $n$-chain of $K$. We can write
\[ x = \sum_{b \in B_n} n_b b, \]
where the $n_b$ are integers, in a unique way. The \ndef{support} $\supp{x}$
of $x$ is the set of $b$ such that $n_b$ is non-null. The \ndef{negative
support} $\suppneg{x}$ and \ndef{positive support} $\supppos{x}$ of $x$ are
the sets
\[
  \suppneg{x} = \setof{b}{n_b < 0}
  \qqtand
  \supppos{x} = \setof{b}{n_b > 0}.
\]
We define
\[
  x^- = \sum_{b \in \suppneg{x}} (-n_b) b
  \qqtand
  x^+ = \sum_{b \in \supppos{x}} n_b b.
\]
We have
\[ x = x^+ - x^-. \]
Actually, $x^-$ and $x^+$ are the only positive $n$-chains $y$ and $z$ with
disjoint support such that $x = z - y$.

If now $n \ge 1$ and $x$ is still an $n$-chain, we define the positive
$(n-1)$-chains $d^-(x)$ and $d^+(x)$ to be
\[
  d^-(x) = (d(x))^-
  \qqtand
  d^+(x) = (d(x))^+.
\]
More generally, if $0 \le i \le n$, we define two positive $i$-chains
$d^-_i(x)$ and $d^+_i(x)$ by
\[
  d^-_i(x) = (d^-)^{n-i}(x)
  \qqtand
  d^+_i(x) = (d^+)^{n-i}(x).
\]
\end{paragraph}

\begin{paragraph}[Unital augmented directed complexes]
Let $K$ be an augmented directed complex with basis $B$. For every element
of the basis $B$, we define a table
\[
  \atom{x} =
      \begin{pmatrix}
        d_0^-(x) & \cdots & d_{n-1}^-(x) & x \\
        \noalign{\vskip 3pt}
        d_0^+(x) & \cdots & d_{n-1}^+(x) & x \\
      \end{pmatrix}.
\]
This table defines an $n$-cell of $\nu(K)$ if and only if $e(d_0^\e(x)) = 1$
for $\e = \pm$.

This motivates the following definition. An augmented directed complex with
basis is said to be \ndef{unital}\index{augmented directed
complex!unital}\index{unital augmented directed complex} if for every $n \ge
0$ and every $n$-chain $x$ of the basis, we
have $e(d_0^-(x)) = 1$ and $e(d_0^+(x)) = 1$.
If $K$ is a unital augmented directed complex, the cells of the form
$\atom{x}$, for $x$ in the basis of $K$, are called \ndef{atoms}.
\index{atom!of an augmented directed complex}
\index{augmented directed complex!atom}
\end{paragraph}

\subsection{Loop-free augmented directed complexes}
Let $K$ be an augmented directed complex~$K$ with basis~$B$. We say that $K$
is \ndef{loop-free}\index{augmented directed
complex!loop-free}\index{loop-free!augmented directed complex} if there exists a partial order on $\coprod_{n \ge 0}
B_n$ such that, for every $n \ge 1$, every $x$ in $B_n$ and every $0 \le i < n$,
any element of the support of $d^-_i(x)$ is smaller than any element of the
support of $d^+_i(x)$.

\index{strongly loop-free!augmented directed complex}
\index{augmented directed complex!strongly loop-free}
Similarly, we say that $K$ is \ndef{strongly loop-free} if there exists a
partial order~$\preceq$ on $\coprod_{n \ge 0} B_n$ such that, for every $n \ge
1$, every $x$ in $B_n$, every $y$ in the support of~$d^-x$ and every $z$
in the support of $d^+x$, one has
  \[ y \preceq x \preceq z. \]
As the terminology suggests, one can show that a strongly loop-free
augmented directed complex is loop-free.

\begin{remark}
  The definition of ``loop-free'' given in the previous paragraph is not the
  one from \cite{steiner2004omega} but the one used in
  \cite{SteinerOpetopes}. The two definitions can be shown to be equivalent.
\end{remark}

\begin{paragraph}[Steiner complexes]
  We will say that an augmented directed complex is a \ndef{Steiner
  complex}\index{Steiner!complex}
  if it is unital and loop-free. Similarly, we will say that it is a \ndef{strong
    Steiner complex}\index{Steiner!complex!strong}\index{strong Steiner!complex}
  if it is unital and strongly loop-free.
\end{paragraph}

\begin{paragraph}[Loop-free polygraphs]
  Let $P$ be a polygraph. We say that $P$ is
  \ndef{loop-free}\index{loop-free!polygraph}\index{polygraph!loop-free} if
  there exists a partial order on the set $\coprod_{n \ge 0} P_n$ such that,
  for every $n \ge 1$, every $x$ in $P_n$ and every $0 \le i < n$, any
  element of the support of $\sce{i}(x)$ is strictly smaller than any
  element of the support of $\tge{i}(x)$.

  \index{strongly loop-free!polygraph}
  \index{polygraph!strongly loop-free}
  Similarly, we will say that $P$ is \ndef{strongly
  loop-free} if there exists a partial order~$\preceq$ on $\coprod_{n \ge 0}
  P_n$ such that, for every $n \ge 1$, every $x$ in $P_n$, every $0 \le i <
  n$, every $y$ in the support of $\sce{i}(x)$ and every $z$ in the support
  of $\tge{i}(x)$, one has~$y \preceq x \preceq z$.

  We will say that an \oo-category is a \ndef{Steiner
  \oo-category}\index{Steiner!9-category@$\omega$-category}\index{9-category@$\omega$-category!Steiner} (resp.{} a
  \ndef{strong Steiner \oo-category}\index{strong Steiner!9-category@$\omega$-category}\index{9-category@$\omega$-category!strong Steiner}) if it is
  generated by loop-free polygraph (resp.~by a strongly loop-free
  polygraph).
\end{paragraph}

\interbreak

We can now state a reformulation of the main theorems
of~\cite{steiner2004omega}:

\begin{theorem}
  The adjunction
  \[
    \xymatrix@C=2pc@R=2pc{
      \oCat\ar@/^/[r]^\lambda\ar@{}[r]|-{\bot}&\ar@/^/[l]^\nu\ADC
    }
  \]
  restricts to an equivalence of categories between
  \begin{itemize}
  \item the full subcategory of $\oCat$ spanned by Steiner's \oo-categories
    (resp.{} by strong Steiner's \oo-categories), and
  \item the full subcategory of~$\ADC$ on Steiner complexes (resp.{} by strong
    Steiner complexes).
  \end{itemize}
  Moreover, if $K$ is a Steiner complex, the \oo-category $\nu(K)$ is
  freely generated in the sense of polygraphs by its atoms.
\end{theorem}

\begin{remark}
  The description of the full subcategories of $\oCat$ appearing in the
  previous theorem is different in \cite{steiner2004omega}. Nevertheless,
  it is shown in~\cite{AGOV} that they are equivalent.
\end{remark}

We will now use Steiner's theory to construct some polygraphs and in
particular recover the ones defined in the previous section.

\begin{paragraph}[Orientals]
  Consider the functor
  \[
    \nc : \SSet \to \pCh[\Z]
  \]
  sending a simplicial set (see~\cref{paragr:def_simpl_sets}) to its
  normalized chain complex (see~\cref{paragr:homology_SSet}). This functor
  naturally lifts to a
  functor
  \[
    \nadc : \SSet \to \ADC.
  \]
  Indeed, if $X$ is a simplicial set, then the chain complex $\nc(X)$ can be
  equipped with the following structure of  augmented directed complex.
  \begin{itemize}
    \item The augmentation $e : \freemod{\Z}{X_0} \to \Z$ sends
      $[x]$, for $x$ a $0$-simplex, to~$1$.
    \item The submonoid $\adcmon{\nc(X)}_n$, for $n \ge 0$, is generated by canonical
      basis of $\nc(X)_n$.
  \end{itemize}
\end{paragraph}
Steiner showed that the composite
\[ \Simpl \hookto \SSet \xto{\nadc} \ADC, \]
where the first functor denotes the Yoneda embedding,
lands into strong Steiner complexes. In particular, by post-composing by
$\nu : \ADC \to \ooCat$, we get a functor
\[ \simpl{} : \Simpl \to \ooCat \]
landing into \oo-categories generated by polygraphs. This functor is the
so-called \ndef{cosimplicial object of orientals}\index{oriental}. In particular, for $n \ge
0$, we recover the $n$-th oriental $\simpl{n}$ as defined
in~\cref{paragr:simpl}.
\nomenclature[On]{$\simpl{n}$}{$n$-th oriental}

More generally, it is shown in \cite{AraMaltsiCondE} that if the simplicial
set $X$ is a simplicial complex (this means, first, that the $(n+1)$-faces of
any $n$-simplex are distinct and, second, that for any set $E$ of $n+1$
$0$-simplices, there is at most one non-degenerate $n$-simplex whose set of
$0$-simplices is $E$), then $\nadc(X)$ is a strong Steiner complex.
Under this condition, we thus get an \oo-category $\nu(\nadc(X))$ generated
by a polygraph, which deserves to be called the \ndef{oriental} associated
to $X$.

\subsection{Globes}

The $n$-category $\globe{n}$ is freely generated by a globular set. In
particular, it is generated by a polygraph $P$ whose generators are the cells of
this globular set. If we denote by $x$ the principal $n$-cell of $\globe{n}$, then
\[ \sce{0}(x) \preceq \cdots \preceq \sce{n-1}(x) \preceq \sce{n}(x) = x =
  \tge{n}(x) \preceq \tge{n-1}(x) \preceq \cdots \preceq \tge{0}(x) \]
is a total order on generators of $P$ showing that $P$ is strongly
loop-free. This shows that $\globe{n}$ is a strong Steiner \oo-category, so
that we have
  \[ \globe{n} \simeq \nu(\lambda(\globe{n})). \]

Let us describe explicitly the augmented directed complex
$\lambda(\globe{n})$.
\begin{itemize}
  \item The chains are defined by
  \[
    \lambda(\globe{n})_i =
    \begin{cases}
      \freemod{\Z}{\{\sce{i}(x), \tge{i}(x)\}} & \text{if $0 \le i < n$,} \\
      \freemod{\Z}{\{x\}} & \text{if $i = n$,} \\
      0 & \text{if $i > n$.}
    \end{cases}
  \]
  \item If $0 < i < n$, then
    \[ 
      d(\sce{i}(x)) = \tge{i-1}(x) - \sce{i-1}(x)
      \qqtand
      d(\tge{i}(x)) = \tge{i-1}(x) - \sce{i-1}(x),
    \]
    and
    \[ d(x) = \tge{n-1}(x) - \sce{n-1}(x). \]
  \item The augmentation is defined by
    \[ e(\sce{0}(x)) = 1 \qqtand e(\tge{0}(x)) = 1. \]
  \item The monoids of positive chains are given by
  \[
    \adcmon{\lambda(\globe{n})}_i =
    \begin{cases}
      \freemod{\N}{\{\sce{i}(x), \tge{i}(x)\}} & \text{if $0 \le i < n$,} \\
      \freemod{\N}{\{x\}} & \text{if $i = n$,} \\
      0 & \text{if $i > n$.}
    \end{cases}
  \]

  One can easily check that $\nu(\lambda(\globe{n}))$ is indeed isomorphic to
  $\globe{n}$ without invoking Steiner's theorem.

\end{itemize}

\begin{paragraph}[Tensor product of augmented directed complexes]
Let $K$ and $L$ be two augmented directed complexes. We define their
\ndef{tensor product} $K \otimes L$ in the following way:
\begin{itemize}
  \item For $n \ge 0$,
    \[
      (K \otimes L)_n
      =
      \bigoplus_{\substack{i+j = n\\i \ge 0,\, j \ge 0}}K_i \otimes L_j.
    \]
  \item If $x$ is in $K_i$ and $y$ is in $K_j$ with $i + j > 0$, then
    \[ d(x \otimes y) = d(x) \otimes y + (-1)^i x \otimes d(y), \]
    where by convention $d(z) = 0$ is $z$ is in $K_0$ or $L_0$.
  \item If $x$ is in $K_0$ and $y$ is in $L_0$, then
    \[ e(x \otimes y) = e(x)e(y). \]
  \item For $n \ge 0$, the submonoid $\adcmon{(K \otimes L)}_n$ is generated
    by the chains of the form $x \otimes y$ with $x$ in $K^+_i$ and $y$ in
    $L^+_j$
    with $i + j = n$.
\end{itemize}

This tensor product defines a (non-symmetric) monoidal structure on
$\ADC$ which is biclosed in the sense that, for $X$ an object, the functors $X
\otimes \var$ and $\var \otimes X$ both admit a right adjoint. The unit
of this tensor product is $\lambda(\globe{0})$.
\end{paragraph}

\begin{proposition}
The tensor product of two strong Steiner complexes is a strong Steiner complex.
\end{proposition}

\begin{proof}
  This is \cite[Example 3.10]{steiner2004omega}.
\end{proof}

\begin{paragraph}[Tensor product of strong Steiner \pdfoo-categories]
  Let $C$ and $D$ be two strong Steiner \oo-categories. We define their tensor
  product by
  \[ C \otimes D = \nu(\lambda(C) \otimes \lambda(D)). \]
  Steiner's theory and the previous proposition imply that $C \otimes D$
  is still a strong Steiner \oo-category. Moreover, if $C$ is generated by a
  polygraph $P$ and $D$ is generated by a polygraph $Q$, the \oo-category $C
  \otimes D$ is generated by a polygraph~$P \otimes Q$ whose $n$-generators
  are given by the formula
  \[
    (P \otimes Q)_n = \coprod_{i+j = n} P_i \times Q_j.
  \]
\end{paragraph}

\begin{theorem}
  \index{Gray tensor product}
  \index{tensor product of $\omega$-categories}
  \index{9-category@$\omega$-category!tensor product}
  There exists a unique (up to unique isomorphism) biclosed monoidal
  structure on $\oCat$ that extends the monoidal structure given by the
  tensor product of Steiner \oo-categories.
\end{theorem}

\begin{proof}
  This is stated in \cite[Section 7]{steiner2004omega}. A detailed proof
  can be found in~\cite[Appendix A]{AraMaltsiJoint}.
\end{proof}

\begin{remark}
  The tensor product given by the previous theorem is called the \ndef{Gray
    tensor product}. It was defined for $2$-categories by Gray in
  \cite{GrayFCT}. It was first extended to \oo-categories by Al-Agl and
  Steiner \cite{AlAglSteiner}. Crans then gave alternate descriptions
  \cite{CransThesis}.
\end{remark}

\begin{proposition}
  The tensor product of two \oo-categories generated by polygraphs is
  generated by a polygraph.
\end{proposition}

\begin{proof}
  See \cite[Theorem 1.35]{AmarThesis} or \cite[Proposition
  5.1.2.7]{LucasThesis}.
\end{proof}

\begin{remark}
  This proposition shows that it makes sense to talk of the tensor product
  of two polygraphs. In particular, identifying $\globe{1}$ with its generating
  polygraph, for any polygraph $P$, we get a polygraph $\globe{1} \otimes P$.
  One can show that this polygraph is the same as the one defined in
  Section~\ref{sec:tensor}.
\end{remark}

\begin{paragraph}[Cylinders and cubes]
\nomenclature[Cyln]{$\Cyl{n}$}{polygraphic $n$-cylinder}
\index{cylinders (as polygraphs)}
\nomenclature[Cubn]{$\cub{n}$}{polygraphic $n$-cube}
\index{cubes (as polygraphs)}
  It follows from the previous remark that one has
  \[ (\Cyl{n})^\ast \simeq \globe{1} \otimes \globe{n} \simeq
    \nu(\lambda(\globe{1}) \otimes \lambda(\globe{n}))
  \]
  and
  \[
    (\cub{n})^\ast \simeq \globe{1} \otimes \cdots \otimes \globe{1} \simeq
    \nu(\lambda(\globe{1}) \otimes \cdots \otimes \lambda(\globe{1})),
  \]
  where $\globe{1}$ and $\lambda(\globe{1})$ both appear $n$ times.
\end{paragraph}

\begin{paragraph}[Join of augmented directed complexes]
Let $K$ and $L$ be two augmented directed complexes. We define their
\ndef{join} $K \join L$ in the following way:
\begin{itemize}
  \item For $n \ge 0$,
    \[
      (K \join L)_n
      =
      \bigoplus_{\substack{i+1+j = n\\i \ge -1,\, j \ge -1}} K_i \otimes L_j,
    \]
    where by convention $K_{-1} = \Z$ and $L_{-1} = \Z$.
  \item If $x$ is in $K_i$ and $y$ is in $K_j$ with $i + 1 + j > 0$, then
    \[ d(x \otimes y) = d(x) \otimes y + (-1)^{i+1} x \otimes d(y), \]
    where by convention $d(z) = e(z)$ if $z$ is in $K_0$ or $L_0$, and $d(n)
    = 0$ for $n$ in~$K_{-1}$ or $L_{-1}$.
  \item If $x$ is in $K_0$ and $y$ is in $L_0$, then
    \[ e(x \otimes 1) = e(x) \qqtand e(1 \otimes y) = e(y). \]
  \item For $n \ge 0$, the submonoid $\adcmon{(K \join L)}_n$ is generated
    by the chains of the form $x \otimes y$ with $x$ in $\adcmon{K}_i$ and
    $y$ in $\adcmon{L}_j$ with $i + 1 + j = n$.
\end{itemize}

The join defines a (non-symmetric) monoidal structure on $\ADC$. The unit
is the initial augmented directed complex, the null complex.
This monoidal structure is not biclosed but only locally biclosed in some
appropriate sense (see \cite[paragraph~5.7]{AraMaltsiJoint}). This locally
biclosedness is equivalent to the fact that the joint commutes with
(non-empty) connected colimits in each variable.
\end{paragraph}

\begin{proposition}
The join of two strong Steiner complexes is a strong Steiner complex.
\end{proposition}

\begin{proof}
  This is \cite[Corollary 6.21]{AraMaltsiJoint}.
\end{proof}

\begin{paragraph}[Join of strong Steiner \pdfoo-categories]
  Let $C$ and $D$ be two strong Steiner \oo-categories. We define their
  \ndef{join} by
  \[ C \join D = \nu(\lambda(C) \join \lambda(D)). \]
  Steiner's theory and the previous proposition imply that $C \join D$
  is still a strong Steiner \oo-category. Moreover, if $C$ is generated by a
  polygraph $P$ and $D$ is generated by a polygraph $Q$, the \oo-category $C
  \join D$ is generated by a polygraph~$P \join Q$ whose $n$-generators
  are given by the formula
  \[
    (P \join Q)_n = \coprod_{\substack{i+1+j = n\\i \ge -1,\, j \ge -1}} P_i \times Q_j,
  \]
  where by convention $P_{-1}$ and $Q_{-1}$ are both singletons.
\end{paragraph}

\begin{theorem}
  \index{join of $\omega$-category}\index{9-category@$\omega$-category!join}
  There exists a unique (up to unique isomorphism) monoidal
  structure on $\oCat$, called the \ndef{join}, that extends the monoidal
  structure given by the join of Steiner \oo-categories and whose monoidal
  product commutes with (non-empty) connected colimits in each variable.
\end{theorem}

\begin{proof}
  This is~\cite[Theorem 6.29]{AraMaltsiJoint}.
\end{proof}

\begin{proposition}
  The join of two \oo-categories generated by polygraphs is generated by a
  polygraph.
\end{proposition}

\begin{proof}
  This is~\cite[Corollary 7.6]{AraLucasFolkMon}.
\end{proof}

\begin{remark}
  This proposition shows that it makes sense to talk of the join
  of two polygraphs. In particular, identifying $\globe{0}$ with its generating
  polygraph, for any polygraph $P$, we get a polygraph $\globe{0} \join P$.
  One can show that it is canonically isomorphic to the polygraph $S(P)$
  defined in~\cref{paragr:simpl}.
\end{remark}

\begin{paragraph}[Orientals]
  It follows from the previous paragraph that one has
  \[ (\simpl{n})^\ast \simeq \globe{0} \join \cdots \join \globe{0}
     \simeq \nu(\lambda(\globe{0}) \join \cdots \join \lambda(\globe{0})),
  \]
  where $\globe{0}$ and $\lambda(\globe{0})$ both appear $n$ times. Note that
  the full cosimplicial object $\simpl{} : \Simpl \to \ooCat$ can be
  recovered from this definition of the orientals. Indeed, as $\globe{0}$ is
  a terminal object in $\oCat$, it is canonically endowed with a monoid
  structure for the monoidal structure given by the join on $\oCat$. By the
  universal property of the augmented simplicial category $\SimplAug$
  (see~\cite[Chapter VII, Section 5, Proposition 1]{MacLane98}),
  this monoid structure induces a functor $\SimplAug \to \oCat$, whose
  restriction to $\Simpl$ gives back the cosimplicial object of orientals.
\end{paragraph}


\chapter{Generalized Polygraphs}
\label{chap:gen-pol}

For each $n\in\N\cup\set{\omega}$, strict $n$-categories are the algebras
of the monad $T=\fgf_n\freecatg_n$ induced by the forgetful functor
$\fgf_n:\nCat n\to\nGlob n$ and its left adjoint $\freecatg_n$, as shown in \cref{chap:n-cat}. However,
the notion of {\em strict} $n$-category is sometimes too restrictive, whence
the need for a notion of {\em weak} $n$-category. One proposal, due to
Penon~\cite{penon1999approche}, defines weak $n$-categories as
algebras of another monad $P$ on $\nGlob n$, which in some sense
``relaxes'' the above monad $T$. In the same vein,
Batanin~\cite{batanin2002penon} describes a general process consisting
in replacing equalities with coherence cells, of which Penon's
construction is a typical instance. 
In~\cite{batanin1998computads}, Batanin generalizes the notion of
polygraph to a notion of \emph{$T$-polygraph} (that he calls
\emph{$T$-computad}), where $T$ is any finitary monad on globular sets.

This chapter starts with our presentation of Batanin's ideas, the main
point being the fairly general adjunction result
of~\cref{subsec:pullback_monadic}. This immediately applies to $\pair
np$-polygraphs, as a straightforward generalization of $n$-polygraphs.

We then turn to two key examples of
this general setting: the monad of weak
$n$-categories~(\cref{sec:weak-polygraphs}), which was the initial
motivation for the general construction, and the monad associated to
linear polygraphs~(\cref{sec:linear-polygraphs}), of special interest
in the present book.

\section{\pdfm{T}-Polygraphs}
\label{sec:T-polygraphs}

The definition of polygraphs presented in \cref{chap:polygraphs}
strongly relies on \cref{prop:ncatp-leftadjoint}, which
asserts the existence of a left adjoint to a certain functor
\[
  \fgfplus_n:\nCat{n+1}\to \nCatp n.
\]
This functor $\fgfplus_n$ is in turn based on the monad of strict
$\omega$-categories on globular sets. Following Batanin's original
idea~\cite{batanin2002computads}, we shall briefly explain how this
construction adapts to an arbitrary finitary monad $T$ on
globular sets, giving rise to the general notion of
``$T$-polygraph''. 

\subsection{Globular algebras}
\label{subsec:globular_algebras}

Let $m\in \N\cup\set{\omega}$ and $T_m$ be a finitary monad on~$\nGlob
m$, that is, a monad whose underlying endofunctor preserves filtered colimits. Let $\nTAlg m$ denote the category of $T_m$-algebras. 
The category $\nTAlg m$ comes with a forgetful functor $\fgf_m:\nTAlg m\to
\nGlob m$, right adjoint to a functor $\fgfad_m:\nGlob m\to \nTAlg m$. Now
for each $n<m$, there is a truncation functor
\[
  \trunc_{n}^{m}:\nGlob m\to \nGlob n
  \]
  right adjoint to the canonical inclusion 
  \[\incl_n^m:\nGlob n\to\nGlob m.\]
Therefore, we get a pair of adjoint
functors $\trunc_{n}^{m}\fgf_m:\nTAlg m\to \nGlob n$ and
$\fgfad_m\incl_n^m:\nGlob n\to \nTAlg m$ whose composition gives a monad $T_n=\trunc_{n}^{m}\fgf_m\fgfad_m\incl_n^m $
on $\nGlob n$. Let again $\nTAlg n$ denote the category of
$T_n$-algebras and $\fgf_n$,  $\fgfad_n$ the corresponding adjoint
functors. There is then a unique comparison functor $K:\nTAlg m\to\nTAlg n$ making the following diagrams commute:
\[
  \xymatrix{\nTAlg m\ar[r]^{\fgf_m}\ar[d]_{K}& \nGlob m\ar[d]^{\trunc_{n}^{m}}\\
    \nTAlg n\ar[r]_{\fgf_n}& \nGlob n \pbox,}
  \qquad\qquad\qquad
   \xymatrix{ \nGlob m\ar[r]^{\fgfad_m}& \nTAlg m\ar[d]^K\\
    \nGlob n\ar[r]_{\fgfad_n}\ar[u]^{\incl_n^m}& \nTAlg n \pbox.}
\]
Thinking of $K$ as a truncation functor, we shall denote it from now
on by the same letter $\trunc_n^m$ as the corresponding truncation
between globular sets. 

\subsection{Freely adjoining cells}
\label{subsec:freecells}

Following the pattern of \cref{chap:polygraphs}, we turn to the
special case of the above situation where $m=n+1$, and consider the
category $\nTAlgp n$ defined by the following pullback square in
$\CAT$:
\[
  \xymatrix{\nTAlgp n \ar[r]\ar[d]& \nGlob{n+1}\ar[d]^{\trunc_{n}^{n+1}}\\
  \nTAlg n\ar[r]_{\fgf_n} & \nGlob n\pbox.}
\]
As the square
\[
  \xymatrix{\nTAlg{n+1}\ar[r]^{\fgf_{n+1}}\ar[d]_{\trunc_{n}^{n+1}}& \nGlob {n+1}\ar[d]^{\trunc_{n}^{n+1}}\\
    \nTAlg n\ar[r]_{\fgf_n}& \nGlob n}
\]
commutes, there is a unique functor $\fgfplus_n:\nTAlg{n+1}\to \nTAlgp n$
such that the following diagram commutes:
\begin{equation}
  \label{diag:functor_w}
  \vcenter{
   \xymatrix{\nTAlg{n+1}\ar[rd]|{\fgfplus_n}\ar@/^/[rrd]^{\fgf_{n+1}}\ar@/_/[rdd]_{\trunc_{n}^{n+1}}&&\\
    &\nTAlgp n \ar[r]\ar[d]& \nGlob{n+1}\ar[d]^{\trunc_{n}^{n+1}}\\
    &\nTAlg n\ar[r]_{\fgf_n} & \nGlob n \pbox.}
  }
\end{equation}
As in the definition of ``ordinary'' polygraphs, the crucial step
will be the construction of a left adjoint $\freeplus_n$ for $\fgfplus_n$. Before
proving this fact, we shall need a few preliminary results on
pullbacks in $\CAT$.

\begin{remark}
  \label{rmk:pullbacks_in_cat}
  Pullbacks in $\CAT$, sometimes called \emph{strict pullbacks}\index{strict
  pullback}, are
  generally badly behaved and do not preserve equivalences of categories.
  Therefore, given categories
  $\ct A$, $\ct B$, $\ct C$ and functors $U:\ct A\to \ct C$, $V:\ct
  B\to \ct C$, one usually defines the \ndef{pseudo-pullback}
  \[
    \xymatrix{\ct P\ar[r]\ar[d] & \ct A\ar[d]^U\\
    \ct B\ar[r]_V & \ct C}
\]
by taking for  objects of $\ct P$ the triples $(a,b,\phi)$ where $a$
is an object of $\ct A$, $b$ an object of $\ct B$ and $\phi$ an
isomorphism from $Ua$ to $Vb$. Note that the above square only commutes up
to a canonical isomorphism.
However, in case $U$ is an {\em
  isofibration}, that is, if for any object $a$ in $\ct A$ and any
isomorphism $g:Ua\to c$ in $\ct C$ there is an isomorphism $f:a\to a'$
such that $Uf=g$, 
it turns out that such a $\ct P$ is equivalent to the strict pullback
(exercise!).
\end{remark}

\subsection{Definition of $\nTAlgp n$}
\label{subsec:ntalgp}
In the present setting, the truncation functor
$\trunc_{n}^{n+1}:\nGlob{n+1}\to \nGlob{n}$ is in fact an
isofibration, and we shall slightly depart here from Batanin's presentation by
defining $\nTAlgp n$ as the {\em strict} pullback of
$\trunc_{n}^{n+1}$ and $\fgf_n$.
Thus, $\nTAlgp n$ has objects pairs $\pair XA$ where $X$ is an
$(n+1)$-globular set and $A$ a $T_n$-algebra such that
$\trunc_{n}^{n+1}X=\fgf_n A$. Morphisms are defined accordingly.

\subsection{Pullbacks of monadic functors}
\label{subsec:pullback_monadic}

Let $\trunc:\ct A\to \ct C$ and $\fgf :\ct B\to \ct C$ be two
functors, and consider their strict pullback in $\CAT$:
\[
    \xymatrix{\ct P\ar[r]^{\fgfpb}\ar[d]_{\truncpb} & \ct A\ar[d]^\trunc\\
    \ct B\ar[r]_\fgf & \ct C \pbox.}
\]
Suppose in addition that
\begin{itemize}
\item $\fgf$ is strictly monadic, with left adjoint $\fgfad$, meaning
  that $\ct B$ is isomorphic to the category of algebras of the
  associated monad $\fgf\fgfad$,
\item $\ct A$ is cocomplete,
\item $\trunc$ is an isofibration, 
\item $\trunc$ admits a left adjoint $\incl$ such that
  $\trunc\incl=\unit{\ct{C}}$ and for each object $c$, the unit $\adju c:c\to
  \trunc\incl c$ is $\unit c$, whence also, for each object $a$ in
  $\ct A$, $\trunc(\adjc a)=\unit{\trunc a}$,
  \item $\trunc$ has a right adjoint.
\end{itemize}
Then the functor $\fgfpb$ is also monadic. We shall denote by
$\adju{}'$ and $\adjc{}'$ the unit and counit of the adjunction
between $\fgf$ and $\fgfad$.

Before proving the
statement, let us point out that the above hypotheses are immediately
satisfied in case $\fgf=\fgf_n$ and $\trunc=\trunc_n^{n+1}$ as defined
in~\ref{subsec:freecells}. 

As for the existence of a left adjoint for $\fgfpb$, let $a$ be an object
of $\ct A$ and define a pair $p=\pair{a^+}{b}$, where $a^+$ is an object
of $\ct A$ and $b$ an object of $\ct B$ as follows: taking first
$b=\fgfad\trunc a$, in order for $p$ to be an object of $\ct P$, we need to
define $a^+$ such that $\trunc a^+=\fgf b$. Consider first the following pushout square in $\ct A$, which
exists because of the cocompleteness assumption:
\begin{equation}
  \label{eq:aplus}
  \vcenter{
  \xymatrix{\incl\trunc a\ar[r]^{\adjc a}\ar[d]_{\incl(\adju{\trunc a}')} & a\ar[d] \\
  \incl\fgf \fgfad\trunc a\ar[r] & a' \pbox.}
  }
\end{equation}
 Now $\trunc$ being left adjoint, it preserves pushouts and moreover
$\trunc(\adjc a)=\unit{\trunc a}$. Therefore, by applying $\trunc$ to the above
square, the bottom arrow becomes an isomorphism $\phi:\fgf b\to \trunc
a'$. Now $\trunc$ being an
isofibration, we may choose an object $a^+$ and an isomorphism
$\psi:a^+\to a'$ such that $\trunc\psi=\phi$. We therefore get a
pushout square
\begin{equation}
  \label{eq:aplusbis}
  \vcenter{
  \xymatrix{\incl\trunc a\ar[r]^{\adjc a}\ar[d]_{\incl(\adju{\trunc a}')} & a\ar[d] \\
  \incl\fgf \fgfad\trunc a\ar[r] & a^+}
  }
\end{equation}
whose top and bottom arrows are taken to identities by $\trunc$.
The construction $\fgfadl:a\mapsto
\pair{a^+}{\fgfad\trunc a}$ is clearly functorial.

It remains to show that
$\fgfadl$ is in fact left adjoint to $\fgfpb$. Let us define two
natural transformations
\[\adjubis{}:\unit{\ct A}\to\fgfpb\fgfadl\]
and
\[\adjcbis{}:\fgfadl\fgfpb\to\unit{\ct P}\]
as follows.

For an object $a$ in $\ct A$, as $\fgfpb\fgfadl a=a^+$ we
may define $\adjubis a$ as the right vertical arrow in the pushout
square~(\ref{eq:aplusbis}). This clearly defines a natural
transformation from $\unit{\ct{A}}$ to $\fgfpb\fgfadl$.

Let now $p=\pair ab$ be an object of $\ct P$, that is, $\trunc a=\fgf
b$. By definition, $\fgfadl\fgfpb p=\pair{a^+}{\fgfad\fgf b}$ and we thus
look for a morphism
\[
  \adjcbis{p}:\pair{a^+}{\fgfad\fgf b}\to \pair{a}{b}.
  \]
The second component is immediately given by the counit
$\adjc{b}':\fgfad\fgf b\to b$. As for the first component, we need to
build a morphism
\[
t_a:a^+\to a
\]
in $\ct A$. Note first that the following triangle
\[
  \xymatrix{\trunc a \ar[d]_{\adju{\trunc a}'}\ar[rd]^{\unit{\trunc a}}& \\
  \fgf\fgfad\trunc a \ar[r]_{\fgf(\adjc{b}')}& \trunc a}
  \]
commutes because $\trunc a=\fgf b$ and the triangular identity between
$\adju{}'$ and $\adjc{}'$. By applying the functor $\incl$ to the
above triangle, we get the commutation of
\[
  \xymatrix{\incl\trunc a \ar[d]_{\incl(\adju{\trunc a}')}\ar[rd]^{\unit{\incl\trunc a}}& \\
  \incl\fgf\fgfad\trunc a \ar[r]_{\incl\fgf(\adjc{b}')}& \incl\trunc a\pbox,}
\]
therefore also the commutation of
\[
  \xymatrix{\incl\trunc a \ar[d]_{\incl(\adju{\trunc
        a}')}\ar[rd]^{\unit{\incl\trunc a}}\ar[rr]^{\adjc a}& &
    a\ar[d]^{\unit a}\\
  \incl\fgf\fgfad\trunc a \ar[r]_{\incl\fgf(\adjc{b}')}& \incl\trunc
  a\ar[r]_{\adjc a}&a \pbox.}
\]
As the outer square above commutes, the universal property of the
pushout yields a unique morphism $t_a:a^+\to a$ such that the
following diagram commutes:
\begin{equation}
  \label{eq:apluster}
  \vcenter{
  \xymatrix{\incl\trunc a\ar[r]^{\adjc a}\ar[d]_{\incl(\adju{\trunc
        a}')} & a\ar[d]_{\adjubis a}\ar@/^/[rdd]^{\unit a} & \\
    \incl\fgf \fgfad\trunc a\ar[r]\ar@/_/[rrd]_{\adjc a\circ\incl\fgf(\adjc{b}')} & a^+\ar[rd]|{t_a}&\\
  && a \pbox.}
  }
\end{equation}
Now, by applying $\trunc$ to the above diagram, one gets
$\trunc(t_a)=\fgf(\adjc{b}')$, whence $\adjcbis p$ is indeed a
morphism in $\ct P$. The naturality of the construction is
immediate. It remains to establish the triangular identities between
$\adjubis{}$ and $\adjcbis{}$. Let first $p=\pair ab$ be an object of
$\ct P$, so that $a=\fgfpb p$. The commutation of the right triangle
in~(\ref{eq:apluster}) reads
\[
  \fgfpb(\adjcbis p)\circ \adjubis{\fgfpb p}=\unit{\fgfpb p}
\]
which gives the first identity. Finally, let $a$ be an object of $\ct
A$ and apply $\fgfadl$ to the right triangle in~(\ref{eq:apluster}),
one obtains
\[
  \fgfadl(t_a)\circ \fgfadl(\adjubis{a})=\unit{\fgfadl a}
\]
and we have to show that
\[
  \fgfadl(t_a)=\adjcbis{\fgfadl a}
\]
which reduces to the equality between both components:
\[
   \truncpb\fgfadl(t_a)=\truncpb(\adjcbis{\fgfadl a})
   \qand
    \fgfpb\fgfadl(t_a)=\fgfpb(\adjcbis{\fgfadl a}).
\]
 The first one comes from the fact that
 $\truncpb\fgfadl=\fgfad\trunc$ and by applying this last functor to
 the bottom triangle in~(\ref{eq:apluster}), and the second from the
 definition of~$\fgfadl$ and~$\adjcbis{}$ according to which
 $\fgfpb\fgfadl(t_a)=t_{a^+}=\fgfpb(\adjcbis{\fgfadl a})$. Therefore
 \[
   \adjcbis{\fgfadl a}\circ \fgfadl(\adjubis{a})=\unit{\fgfadl a}
 \]
 and the second triangular identity is proved. Thus $\fgfadl$ is
 left adjoint to $\fgfpb$.

The monadicity of $\fgfpb$ now follows from Beck's criterion. Let
$p=\pair ab$, $p'=\pair{a'}{b'}$ in $\ct P$ and 
$u,v:p\to p'$ be  a pair of parallel morphisms such that
\begin{equation}
  \label{diag:abscoeq1}
  \xymatrix{\fgfpb p\doubr{\fgfpb u}{\fgfpb v} & \fgfpb p'\ar[r]^e & a''}
\end{equation}
is an absolute coequalizer in $\ct A$. It follows that
\begin{equation}
  \label{diag:abscoeq2}
  \xymatrix{\fgf b\doubr{\trunc\fgfpb u}{\trunc\fgfpb v} & \fgf b'\ar[r]^{\trunc e} & \trunc a''}
\end{equation}
is also an absolute coequalizer in $\ct C$. As $\fgf$ is strictly
monadic, it creates such coequalizers. There is therefore a unique
morphism $f:b'\to b''$ in $\ct B$ such that $\fgf f=\trunc e$ and
\[
  \xymatrix{b\doubr{\truncpb u}{\truncpb v} & b'\ar[r]^f & b''}
 \]
is a coequalizer in $\ct B$. Clearly $\pair{a''}{b''}$ is in $\ct P$
and the morphism
\[ \pair ef: \pair{a'}{b'}\to \pair{a''}{b''} \] is the
coequalizer of $u$ and $v$ we are looking for.


\interbreak
 
The existence of a left adjoint for the functor~$\fgfplus_n$ 
of~(\ref{diag:functor_w}) is then a consequence of the following
result in general category theory:
\begin{prop}\label{prop:adjtriangle}
  Let $X:\ct A\to\ct B$ and $Y:\ct B\to\ct C$ be two functors and
  $Z=YX$, and suppose that
  \begin{itemize}
  \item $\ct A$ is cocomplete,
  \item $Z$ has a left adjoint,
  \item $Y$ is monadic.
   \end{itemize}
  Then the functor $X$ admits a left adjoint.
\end{prop}

\begin{proof}
This is essentially a simpler, less general version of Dubuc's adjoint
triangle theorem~\cite{Dubuc1968triangle}.   
\end{proof}

\subsection{Existence of the functor $\freeplus_n$}
This above statement immediately applies to the triangle
\[
  \xymatrix{\nTAlg{n+1}\ar[rrd]^{\fgf_{n+1}}\ar[rd]_{\fgfplus_n} & &\\
  &\nTAlgp{n}\ar[r]_{\fgfpb_n} & \nGlob{n+1} \pbox.}
  \]
In fact, the cocompleteness of $\nTAlg{n+1}$ comes from the finitary assumption
on the monad $T$, the functor $\fgf_{n+1}$, being the forgetful functor from
$(n+1)$-globular algebras to $(n+1)$-globular sets  admits a
left adjoint, and we have just shown that $\fgfpb_n$ is
monadic. 
Therefore $\fgfplus_n$ admits a left adjoint $\freeplus_n$.

\subsection{Definition of $T$-polygraphs}
We may now define by induction on $n\geq 0$ a category $\nTPol n$,
together with  a functor
\[
  \freecatpol_n:\nTPol n\to \nTAlg n,
  \]
  following the pattern of Section~\ref{subsec:npolyg}.
  \begin{itemize}
  \item For $n=0$, $\nTPol 0=\nGlob 0=\Set=\nTAlg 0$ and $F_0$ is the
    identity functor.
  \item Let $n\geq 0$ and suppose we have defined
    \[
\freecatpol_n:\nTPol n\to \nTAlg n.
    \]
     The category $\nTPol{n+1}$ is then given by the following pullback square
    in $\CAT$:
     \[
      \xymatrix{\nTPol{n+1}\ar[r]^{J_n}\ar[d] & \nTAlgp{n}\ar[d]\\
      \nTPol n\ar[r]_{\freecatpol_n} & \nTAlg n \pbox;}
  \]
  and the functor $\freecatpol_{n+1}$ is the composite
  \[
    \freeplus_nJ_n:\nTPol{n+1}\to \nTAlg{n+1}
  \]
  where $\freeplus_n$ is the above defined left adjoint to $\fgfplus_n$.
\end{itemize}
If we start with a finitary monad $T$ on $\oGlob$, we get 
    by construction a sequence of canonical functors in $\CAT$
    \[
      \xymatrix{\nTPol 0 & \nTPol 1\ar[l] & \nTPol 2\ar[l]& \cdots \ar[l] }
      \]
whose projective limit $\nTPol{\omega}$ is by definition the category of
$\omega$-polygraphs with respect to the monad $T$.

Finally, for each $n\geq 0$, there is a functor
\[
  \catpol_n: \nTAlg n \to \nTPol n
  \]
  and a natural transformation
  \[
    \varepsilon : \freecatpol_n\catpol_n\to \id
  \]
  such that $\freecatpol_n$ is left adjoint to $\catpol_n$ and
  $\varepsilon$ is the counit of this adjunction.

 The construction of $\catpol_n$ in the general case  being
 essentially the same as the one explained in Section~\ref{sec:adj_triangle}
 for the particular case where $\nTAlg{n}=\nCat n$, we shall only very briefly sketch
 the induction step yielding $\catpol_{n+1}$ from $\catpol_n$.

 Thus, suppose we have defined $\catpol_n$ as a right adjoint to
 $\freecatpol_n$, with counit $\varepsilon$, and let $C$ be an object
 of $\nTAlg{n+1}$, with underlying $(n+1)$-globular set $X$. Let
 $C'=\trunc_n^{n+1}(C)$. We define
 $\catpol_{n+1}(C)$ as a pair $\pair{P}{D^+}$ where $P$ is in $\nTPol{n}$
 and $D^+$ in $\nTAlgp{n}$ by taking $P=\catpol_n(C')$
 and $D^+=\pair{D}{Z}$, with $D=\freecatpol_n(P)$ and $Z$ the  $(n+1)$-globular
 set defined as follows. Up to dimension $n$, $Z$ coincides with the
 underlying $n$-globular set of $\freecatpol_n(P)$, whereas $Z_{n+1}$
 consists of triples $(z,x,y)\in X_{n+1}\times Z_n\times Z_n$ such
 that $z:\varepsilon^{C'}_n (x)\to\varepsilon^{C'}_n(y)$.

 Now
 $\varepsilon^{C'}:\freecatpol_n\catpol_n(C')\to C'$ extends to a
 natural transformation 
 \[ \theta:D^+\to \fgfplus_n(C) \]
 by sending the generator $(z,x,y)$ to $z$.

 Thus, by adjunction, we get $\theta^*:\freeplus_n(D^+)\to C$, but
 $\freeplus_n(D^+)$ is precisely $\freecatpol_{n+1}\catpol_{n+1}(C)$, so that
 $\theta^*$ defines the counit
 $\varepsilon^C:\freecatpol_{n+1}\catpol_{n+1}(C)\to C$.

\subsection{Basic examples}
Besides the basic case of strict \oo-categories, an immediate example
of the above construction is given by $\pair np$-polygraphs introduced
in \cref{sec:np-polygraph}.
We just remark here that,
again, these polygraphs are particular instances of $T$-polygraphs,
where the monad $T$ is the one induced by the forgetful functor
\[
  V:\npCat np\to \nGlob n.
\]

 \section{Polygraphs for Weak \pdfm{n}-Categories}
\label{sec:weak-polygraphs}

The general construction of Section~\ref{sec:T-polygraphs} applies to
Penon's monad on $n$\nbd-globular sets,
from~\cite{penon1999approche}. This section basically follows
Penon's approach, but for a correction pointed out by Cheng and
Makkai in~\cite{ChengMakkai2009penon}: whereas Penon works over {\em
  reflexive} globular sets, it turns out that the category of plain globular sets,
that is, our category $\nGlob n$, yield more examples of weak
$n$-categories. In particular, braided monoidal categories fit in the
latter setting, but not in the former.

\nomenclature[Magn]{$\nMag n$}{category of $n$-magmas}
Given $n\in\N\cup\set{\omega}$, we first consider the category $\nMag
n$ of {\em $n$-magmas}, whose objects are $n$-globular sets endowed
with the same family of binary composition operations as (strict)
$n$-categories, satisfying the same ``positional'' conditions with
respect to source and target maps, but without requiring identities,
associativity and exchange. The morphisms of $\nMag n$ are the
globular maps preserving all compositions. Of course any $n$-category
is a particular $n$-magma. Now, let $M$ be an $n$-magma and $C$ be an
$n$-category seen as an $n$-magma, a morphism
\[
  f:M\to C
\]
in $\nMag n$ is called a {\em categorical stretching}\index{categorical
stretching}---``étirement
catégorique'' in Penon's terminology, or ``trivial fibration'' in
Batanin's---if for any $k<n$ and any pair
$\pair ab$ of parallel $k$-cells in $M$ such that $f(a)=f(b)$, there
is a $(k+1)$\nbd-cell $c$ in $M$ such that $f(c)=\Unit{k+1}{f(a)}$, and if
moreover, in case $n\neq\omega$, $f_n$ is injective. A {\em
  trivialization} of a categorical stretching $f$ is a map $[,]_f$
choosing a cell $c=[a,b]_f$ with $f(c)=\Unit{k+1}{f(a)}$ for each pair
$\pair ab$ as above. There is now a category $\ct{Q}$ whose objects
are pairs $\pair{f}{[,]_f}$, where $f$ is a categorical stretching and
$[,]_f$ a trivialization of $f$, and whose morphisms are commutative
squares
\[
  \xymatrix{
    M\ar[r]^u\ar[d]_f & M'\ar[d]^{f'}\\
    C\ar[r]_v & C'
    }
  \]
  in $\nMag n$ such that $u[a,b]_f=[va,vb]_{f'}$ for all parallel
  $k$-cells $a$, $b$, such that $f(a)=f(b)$. As each $n$-magma $M$ has
  an underlying $n$-globular set, the correspondence taking $f:M\to C$
  to $M$ induces a forgetful functor
  \[
    U:\ct{Q}\to\nGlob n.
    \]
 Now, as shown in~\cite{penon1999approche}, this functor $U$ admits a
 left adjoint $F$, thus defining a monad $P=UF$ on $\nGlob n$. The
 algebras of this monad $P$ are precisely the weak $n$-categories we
 were looking for.

 Finally, this monad $P$ satisfies the hypotheses of
 Section~\ref{sec:T-polygraphs}, and therefore produces an appropriate
 notion of $P$-polygraphs.

\section{Linear Polygraphs}

\subsection{Linear polygraphs as $T$-polygraphs}
\nomenclature[Alg]{$\Alg$}{category of algebras}
\nomenclature[Algn]{$\nAlg n$}{category of $n$-algebras}
\index{linear!8-polygraph@$n$-polygraph}
\index{polygraph!linear}
\label{sec:linear-polygraphs}
\label{S:LinearPolygraphs}

In \cref{C:OneDimensionalLinearRewriting} we introduced the
notion of $1$-dimensional linear polygraphs as rewriting systems for
associative unital algebras over a field $\kk$. 
These are again special
cases of the general construction~\ref{sec:T-polygraphs}, for a
certain monad on $n$-globular sets we now briefly describe. Let
$n\in\N\cup\set{\omega}$ and $\kk$ be a field. We denote by~$\Alg$ the
category of unital and associative $\kk$-algebras. The category $\nAlg n$ of $n$-algebras is the category of $n$-categories internal to
$\Alg$. Consider the forgetful functor
\[U:\nAlg n\to \nGlob n\]
obtained by composing both functors $\nAlg n\to\nCat n$ and
$\nCat n\to \nGlob n$. 
The categories $\nGlob n$, $\nCat n$, $\Set$, and $\Alg$ are
models of projective sketches $S$, $S'$, $T$, and $T'$, respectively,
where $T$ is of course the trivial sketch on the terminal
category. Moreover, there are inclusion morphisms of sketches $i:S\to
S'$ and $j:T\to T'$, and both of them satisfy Lair's conditions of
\cref{thm:Lair}. These induce a morphism
\[i\tensor j:S\tensor T\to S'\tensor T'\]
still satisfying Lair's conditions (see~\cite[p.~18]{ageron1992logic} about
tensoring sketches). Now the models of
$S'\tensor T'$ in $\Set$ are also the models of $S'$ in
$\Modset{T'}=\Alg$, in other words 
\[\Modset{S'\tensor T'}=\nAlg n\]
whereas $S\tensor T\simeq S$. Therefore our functor $U$ is precisely
the one induced by $i\tensor j$ on models:
\[U=\Modset{i\tensor j}:\Modset{S'\tensor T'}\to \Modset{S}.\]
As a consequence of \cref{thm:Lair}, the functor $U$ is
monadic. It is also easily seen to preserve filtered colimits, whence
inducing a finitary monad $L$ on $\nGlob n$. Thus the machinery
of~\ref{sec:T-polygraphs} applies, and produces the category of
$L$-polygraphs, that is, the {\em linear
  polygraphs} we wanted.

\subsection{Bimodules}
\index{globular bimodule}
\index{bimodule!globular}
\nomenclature[BMod(A)]{$\BMod(A)$}{category of bimodules over an algebra~$A$}
\label{subsec:bimod}
The description of $n$-algebras and linear $n$\nbd-poly\-graphs can be made more
explicit by introducing the following notion of globular bimodules. For an
algebra~$A$, we denote by $\BMod(A)$ the category of bimodules over $A$ and
their morphisms. We consider the category $\GBMod(A)$ of \emph{globular
$A$-bimodules}, that is, functors
\[
X:\opp{\glob}\to \BMod(A),
\]
and their morphisms.
Let us define the category $\BModg$ whose objects are pairs~$(A,M)$ made of
an algebra~$A$ and a globular $A$-bimodule~$M$, and whose morphisms
from~$(A,M)$ to~$(B,N)$ are pairs~$(F,G)$ made of a morphism $F:A\fl B$ of
algebras and a morphism $G:M\fl N$ of bimodules, that is,
\[
G(ama')=F(a)G(m)F(a') 
\]
holds for all~$a$ and~$a'$ in~$A$ and~$m$ in~$M$.

\interskip

The following result gives a characterization of the category of $\omega$-algebras.

\begin{theorem}[{\cite[Theorem 1.3.3]{GuiraudHoffbeckMalbos19}}]
\label{T:nAlg}
The category~$\nAlg{\omega}$ is isomorphic to the full subcategory of
$\BModg$ whose objects are the pairs~$\pair AM$ such that~$M_0$  is equal to~$A$, with its canonical $A$-bimodule structure, and that satisfy, for all $n$-cells~$a$ and~$b$ of~$M$, the relation
\begin{equation}
\label{E:nAlgExchangeOmegaAlgebras}
a \sce{0}(b) + \tge{0}(a) b - \tge{0}(a) \sce{0}(b) = \sce{0}(a) b + a \tge{0}(b) - \sce{0}(a) \tge{0}(b).
\end{equation}
\end{theorem}

\subsection{Explicit construction}
The key step in the explicit inductive construction of $T$-polygraphs at
level $n$ is the concrete description of the left adjoint~$\freeplus_n$ to the forgetful functor $\fgfplus_n:\nTAlg{n+1}\to\nTAlgp{n}$ of~(\ref{diag:functor_w}). In the present case of
linear polygraphs, this forgetful functor is
\[
\fgfplus_n:\nAlg{n+1}\to \nAlgp{n}
\]
and its left adjoint $\freeplus_n$ takes an extended
$n$-algebra~$(A,X)$ to the $(n+1)$\nbd-algebra $A[X]$ constructed as
follows.
First, we consider the $A_0$-bimodule
\[
M = \big( A_0 \otimes \kk[X] \otimes A_0 \big) \oplus A_{n}
\]
obtained by the direct sum of the free $A_0$-bimodule with basis~$X$ and of a copy of~$A_{n}$, equipped with its canonical $A_0$-bimodule structure. Thus~$M$ contains linear combinations of elements $axb$, for~$a$ and~$b$ in~$A_0$ and~$x$ in~$X$, and of an $n$-cell~$c$ of~$A$. We define the source, target, and identity maps 
\[
  \xymatrix{M\doubr{\src{}}{\tgt{}} & A_{n} \ar@/_3.5ex/[l]_{\reflex}}
\]
by
\begin{align*}
  \src{}(axb)&= a \src{}(x) b,
  &
  \src{}(c) &= c,
  &
  \reflex(c) &= c,
  \\
  \tgt{}(axb)&= a \tgt{}(x) b,
  &
  \tgt{}(c) &= c,
\end{align*}
for all~$x$ in~$X$, $a$ and~$b$ in~$A_0$, and~$c$ in~$A_{n-1}$. 
Then we define the $A_0$\nbd-bimodule $A[X]_{n+1}$ as the quotient of~$M$ by the $A_0$-bimodule ideal generated by all the elements
\[
\big( a \sce{0}(b) + \tge{0}(a) b - \tge{0}(a) \sce{0}(b) \big) 
	\: - \: 
	\big( \sce{0}(a) b + a \tge{0}(b) - \sce{0}(a) \tge{0}(b) \big),
\] 
where~$a$ and~$b$ range over $A_0 \otimes \kk X \otimes A_0$. We prove that the source and target maps are compatible with the quotient, so that, by \cref{T:nAlg}, the $A_0$-bimodule $A[X]_{n+1}$ extends~$A$ into a uniquely defined $(n+1)$-algebra~$A[X]$.

\part{Homotopy Theory of Polygraphs}
\label{part:cof}

\chapter{Polygraphic Resolutions}
\label{chap:resolutions}

The purpose of this chapter is to introduce the notion of a polygraphic
resolution of an \oo-category. This notion was introduced by Métayer
\cite{metayer2003resolutions} to define a homology theory for \oo-categories,
that is now known as the \emph{polygraphic homology}. It was then showed by
himself and Lafont \cite{lafont2009polygraphic} that this homology recovers
the classical homology of monoids for \oo-categories coming from monoids.
It is now known by the work of Lafont, Métayer, and Worytkiewicz~\cite{LMW} that
these polygraphic resolutions are resolutions in the sense of a model
category structure on $\ooCat$, the so-called \emph{folk model structure},
that we will present in the next chapters.

Roughly speaking, a polygraphic resolution is a non-abelian version of a
resolution of a module in the sense of homological algebra. More precisely,
if~$C$ is an \oo-category, a polygraphic resolution of $C$ is a polygraph
$P$ endowed with an \oo-functor $P^\ast \to C$ that is a \emph{trivial
fibration}, meaning in a nutshell that it is surjective with source and target
fixed at all levels. Technically, we will define these trivial fibrations
by a right lifting property with respect to inclusions of spheres into
disks, in a very similar way as trivial fibrations of topological
spaces are defined.

The chapter is organized as follows. In the first section, we introduce the
notion of a weak factorization system in a general category. We explain the
 \emph{small object argument}. In the second section, we define
\emph{cofibrations} and \ndef{trivial cofibrations} of \oo-categories as the
classes appearing in the weak factorization system generated by the set
$\setgencof$ of inclusions of spheres into disks. We also introduce the
class of \emph{relative polygraphs}, which are  cell complexes
generated by~$\setgencof$, and we explain why polygraphs almost
tautologically correspond to relative polygraphs of the form $\varnothing
\to C$. In the third section, we introduce the notion of a polygraphic
resolution. We show that every \oo-category admits such a resolution and we
give the example of the so-called \emph{canonical resolution}. Finally, in the
last section, we study the uniqueness of these resolutions, showing that there
is always a map between two such resolutions and postponing to a later
chapter the fact that two such maps are homotopic in some appropriate sense.

\section{Weak Factorization Systems}
\label{sec:wfs}

The purpose of this section is to introduce some basic results on weak
factorization systems, which will be applied in the next section to
a factorization system giving rise to ``polygraphic resolutions''. This
factorization system has the additional property of being generated by
morphisms between finitely presentable objects. This additional hypothesis
will allow us to avoid the use of ordinals and cardinals.

\medbreak

In the section, we fix a category $\C$. Our case of interest is $\C=\ooCat$.

\begin{paragr}[Lifting properties]
  \index{lifting property}
  \index{left lifting property}
  \index{right lifting property}
  \index{lift}
  Let $f : X \to Y$ and $g : Z \to T$ be two morphisms of $\C$. One says
  that $f$ \ndef{has the left lifting property} with respect to $g$ or that
  $g$ \ndef{has the right lifting property} with respect to $f$ if for every
  commutative square
  \[
  \xymatrix{
    X \ar[d]_f \ar[r] & Z \ar[d]^g \\
    Y \ar[r] & T \\
  }
  \]
  there exists a \ndef{lift}, that is, a morphism $h : Y \to Z$ making
  the two triangles
  \[
  \xymatrix{
    X \ar[d]_f \ar[r] & Z \ar[d]^g \\
    Y \ar[r] \ar@{.>}[ur]^h & T \\
  }
  \]
  commute. More generally, one says that $f$ has the left lifting property
  with respect to a class of maps $\clI$ if it has the left lifting
  property with respect to every morphism in $\clI$, and similarly for the
  right lifting property. We will denote by $\lorth{\clI}$ and
  $\rorth{\clI}$ the class of maps having the left or right lifting property
  with respect to a class $\clI$.
\end{paragr}

\begin{paragr}[Stability properties of $\lorth{\clI}$]\label{paragr:stab_lI}\label{sec:stab_lI}
  Suppose $\C$ is cocomplete and let $\clI$ be a class of morphisms of $\C$.
  One checks that $\lorth{\clI}$ contains the class of isomorphisms and
  is stable under
  \begin{enumerate}
    \item sums: if $(i_k : X_k \to Y_k)_{k \in K}$ is a (small) family of
    elements of $\lorth{\clI}$, then the sum
    \[
      \coprod_{k \in K} i_k : \coprod_{k \in K} X_k \to \coprod_{k
      \in K} Y_k
    \] 
    belongs to $\lorth{\clI}$,
    \item pushouts: if
      \[
        \xymatrix{
          X \ar[d]_i \ar[r] & Z \ar[d]^j \\
          Y \ar[r] & T
        }
      \]
      is a pushout square and $i$ belongs to $\lorth{\clI}$, then so does $j$,
    \item countable compositions: if
    \[
     \xymatrix@C=1.5pc{
     X_0 \ar[r]^{i_1} & X_1 \ar[r]^{i_2} & \cdots \ar[r]^{i_n} & X_n
     \ar[r]^{i_{n+1}} & \cdots
     }
    \]
    is a diagram of elements of $\lorth{\clI}$, then the morphism
    \[ X_0 \to \limind_{n \ge 0} X_n \]
    belongs to $\lorth{\clI}$,
    \item retracts: if
    \[
      \xymatrix{
        X \ar[d]_i \ar[r] \ar@/^3ex/[rr]^{\unit{X}} & Z \ar[d]_j \ar[r] & X \ar[d]^i \\
        Y \ar[r] \ar@/_3ex/[rr]_{\unit{Y}} & T \ar[r] & Y
      }
    \]
    is a commutative diagram and $j$ belongs to $\lorth{\clI}$, then so does
    $i$.
  \end{enumerate}
\end{paragr}

\begin{rem}
  More generally, under the same assumption and with the same notation as in
  the paragraph above, the class $\lorth{\clI}$ is stable under transfinite
  compositions (see~\cite[Definition 2.1.1]{hovey2007model}). Note that
  being stable by transfinite compositions and pushouts implies being stable
  by sums.
\end{rem}

\begin{paragr}[Weak factorization systems]
  \index{weak factorization system}
  \index{factorization system!weak}
  A \ndef{weak factorization system} on $\C$ is a pair $(\clL, \clR)$ of
  classes of morphisms of $\C$ satisfying the following conditions:
  \begin{enumerate}
    \item Every morphism $f$ of $\C$ factors as $f = pi$, where $p$ is in
    $\clR$ and $i$ is in $\clL$.
    \item We have
    \[ \clL = \lorth{\clR} \qqtand \clR = \rorth{\clL}. \]
  \end{enumerate}
\end{paragr}

\index{finitely!presentable object}
Recall that an object~$A$ of a category is \emph{finitely presentable} when the
functor $\Hom\C A-:\C\to\Set$ preserves filtered colimits. In practice, we are
mostly interested in preservation of colimits of diagrams in~$\C$ of the form
\[
  \xymatrix{
    X_0\ar[r]&X_1\ar[r]&X_2\ar[r]&\cdots
  }
\]
In this case, the fact that the canonical arrow $\colim_n\Hom\C
A{X_n}\to\Hom\C A{\colim_nX_n}$ is an isomorphism means in particular that
any morphism $f:A\to\colim_nX_n$ factors through some $X_n$. For instance,
the finitely presentable sets are the finite sets.

\begin{prop}[Small object argument: first version]
  \label{prop:small_obj_arg_1}
  \index{small object argument}
  Suppose $\C$ is cocomplete and let $I$ be a \emph{set} of morphisms of $\C$
  whose sources are finitely presentable objects. Then $(\lrorth{I},\rorth{I})$
  is a weak factorization system on $\C$.
\end{prop}

\begin{proof}
  The second condition of the definition of a weak factorization system
  follows from the general equality $rlr(I) = r(I)$. Let us sketch a
  proof for the first one.  Let $f : X \to Y$ be a morphism of $\C$. Let $S$
  be the set
  of commutative squares of the form
  \[
    \xymatrix{
      A_s \ar[r] \ar[d]_{i_s} & X \ar[d]^f \\
      B_s \ar[r] & Y \pbox{,}
    }
  \]
  where $i_s$ is an element of $I$. Summing all these squares, we get a
  commutative square
  \[
    \xymatrix{
      \coprod_{s \in S} A_s \ar[r] \ar[d]_{\coprod_{s \in S}i_s} & X \ar[d]^f \\
      \coprod_{s \in S} B_s \ar[r] & Y \pbox{.}
    }
  \]
  Taking the pushout of the top left corner, we get a factorization of $f$:
  \[
    \xymatrix{
      \coprod_{s \in S} A_s \ar[r] \ar[d]_{\coprod_{s \in S}i_s} & X
      \ar@/^3ex/[ddr]^f
      \ar[d]^{j_1} \\
      \coprod_{s \in S} B_s \ar[r] \ar@/_3ex/[drr] & X_1 \ar[rd]^{p_1} \\
      & & Y \pbox{.}
    }
  \]
  By \cref{paragr:stab_lI}, the morphism $j_1$ belongs to $\lrorth{I}$, but there
  is no reason for $p_1$ to belong to $\rorth{I}$: we thus apply the same
  procedure to~$p_1$, obtaining a factorization $p_1 = p_2j_2$, where
  $j_2 : X_1 \to X_2$ belongs to $\lrorth{I}$. Going on by induction, we get a
  sequence of morphisms
  \[
    \xymatrix@C=1.5pc{
      X \ar[r]^{j_1} & X_1 \ar[r]^{j_2} & \cdots \ar[r]^{j_n} & X_n
      \ar[r]^{j_{n+1}} & \cdots
    }
  \]
  belonging to $\lrorth{I}$ and hence a morphism
  \[
    j_\infty : X \longrightarrow X_\infty = \limind_n X_n
  \]
  in $\lrorth{I}$ giving rise to a factorization
  \[ f = p_\infty j_\infty, \]
  where $p_\infty$ is the morphisms induced by the $p_n$'s. To conclude
  the proof, it suffices to show that $p_\infty$ belongs to
  $\rorth{I}$. Consider a commutative square
  \[
    \xymatrix{
      A \ar[r]^k \ar[d]_i & X_\infty \ar[d]^{p_\infty} \\
      B \ar[r]_l & Y \pbox{,}
    }
  \]
  where $i$ is in $I$. Since by our additional assumption, the object $A$ is
  finitely presentable, the morphism $k$ factors through some $X_n$ and the
  diagram factors as
  \[
    \xymatrix{
      A \ar[r] \ar[d]_i & X_n \ar[r] &  X_\infty \ar[d]^{p_\infty} \\
      B \ar[rr]_l & & Y \pbox{,}
    }
  \]
  where $X_n \to X_\infty$ is the canonical morphism. The composite from $X_n$
  to $Y$ is~$p_n$, and this diagram defines a square appearing in the definition
  of $X_{n+1}$. This means that there exists a lift
  \[
    \xymatrix{
      A \ar[r] \ar[d]_i & X_n \ar[r]^{j_{n+1}} & X_{n+1} \ar[r] &  X_\infty \ar[d]^{p_\infty} \\
      B \ar[rrr]_l \ar@{.>}[urr] & & & Y \pbox{,}
    }
  \]
  thereby proving the result.
\end{proof}

\begin{rem}\label{rem:hyp_small_obj_arg}
  The same conclusion holds with the weaker hypothesis that there exists a
  regular cardinal $\kappa$ such that the sources of the morphisms in $I$ are
  $\kappa$\nbd-presentable (see Section~\ref{sec:def_kappa-small}), the
  above statement being the case $\kappa = \aleph_0$.  This condition is
  automatic if $\C$ is locally presentable. See also
  \cite[Section~10.5]{HirMC} for weaker assumptions for the small object
  argument. The proof of these more general results is basically the same
  except that one has to proceed by transfinite induction.
\end{rem}

\begin{paragr}[$\clI$-cells]
  \index{cellular extension!for factorization systems}
  Suppose $\C$ is cocomplete and let $\clI$ be a class of morphisms of~$\C$.
  The class $\cellom{\clI}$ of \ndef{countable $\clI$-cellular extensions}
  is the smallest class of morphisms of $\C$ containing $\clI$ and stable
  under sums, pushouts, and countable compositions. Every element of
  $\cellom{\clI}$ can be obtained as a countable composition of pushouts of
  sums of elements of~$\clI$.
\end{paragr}

\begin{rem}
  When $\clI$ is a set and the source of the elements of $\clI$ are finitely
  presentable, by~\cite[Proposition A.6]{MaltsiniotisWCat}, the class
  $\cellom{\clI}$ is equal to the more classical class $\cell{\clI}$ of
  $\clI$-cellular extensions, which is defined as $\cellom{\clI}$ but asking
  also for stability under transfinite compositions.
\end{rem}

\begin{prop}[Small object argument: refined version]
  \label{prop:small_obj_arg_2}
  \index{small object argument}
  Suppose $\C$ is cocomplete and let $I$ be a \emph{set} of morphisms of
  $\C$ whose sources are finitely presentable. Then every morphism of $\C$
  factors as $f = pi$, where $p$ is in $\rorth{I}$ and $i$ is in
  $\cellom{I}$. Moreover, every element of $\lrorth{I}$ is a retract of an
  element of~$\cellom{I}$.
\end{prop}

\begin{proof}
  The first point was actually proven in the proof
  of~\cref{prop:small_obj_arg_1}. The second point follows from the following
  lemma, the so-called ``retract lemma''.
\end{proof}

\begin{lemma}[Retract lemma]\label{lemma:retr}
  \index{retract!lemma}
  \index{lemma!retract}
  Suppose we have a factorization $f = pi$, in a category $\C$, where $f$
  has the left lifting property with respect to $p$. Then $f$ is a retract
  of $i$.
\end{lemma}

\begin{proof}
  Denote $i : X \to Y$ and $p : Y \to Z$. By hypothesis, there exists $h :
  Z \to Y$ making the square
  \[
  \xymatrix{
    X \ar[d]_f \ar[r]^i & Y \ar[d]^p \\
    Z \ar@{=}[r] \ar[ur]^h & Z \\
  }
  \]
  commute. The commutative diagram
  \[
    \xymatrix{
      X \ar@{=}[r] \ar[d]_f & X \ar@{=}[r] \ar[d]_i & X \ar[d]^f \\
      Z \ar[r]_h & Y \ar[r]_p & Z
    }
  \]
  proves the result.
\end{proof}

\begin{rem}
  Proposition~\ref{prop:small_obj_arg_2} holds under the weaker hypothesis
  described in Remark~\ref{rem:hyp_small_obj_arg} if one replaces
  $\cellom{I}$ with $\cell{I}$.
\end{rem}

\section{Cofibrations and Trivial Fibrations}
\label{sec:oocat-cof-tfib}

We will now apply the theory of weak factorization systems described in the
previous section to the category $\ooCat$ and the set of inclusions of
boundaries of globes. This will lead to a notion of cofibrations and
trivial fibrations of \oo-categories.

\begin{paragr}[Cofibrations, trivial fibrations, and relative polygraphs]
  \label{paragr:cof_triv_fib}
  \index{cofibration!of 9-categories@of $\omega$-categories}
  \index{trivial fibration!of 9-categories@of $\omega$-categories}
  \index{fibration!trivial!of 9-categories@of $\omega$-categories}
  \index{9-category@$\omega$-category!cofibration}
  \index{9-category@$\omega$-category!trivial fibration}
  \index{relative!polygraph}
  \index{polygraph!relative}
  Recall that for every $n \ge 0$ we denote by $\gencof{n} : \sphere{n} \to
  \globe{n}$ the inclusion of the $(n-1)$-sphere into the $n$-globe. We set
  \[ \setgencof \qeq \setof{\gencof{n}}{n \ge 0}. \]
  From this set $\setgencof$, we obtain three classes of \oo-functors:
  \begin{itemize}
    \item the class $\lrorth{\setgencof}$ of \ndef{cofibrations},
    \item the class $\rorth{\setgencof}$ of \ndef{trivial fibrations}, and
    \item the class $\cellom{\setgencof}$ of \ndef{relative polygraphs}.
  \end{itemize}

  The category $\ooCat$ being locally presentable
  (see~\cref{paragr:oCat_loc_pres}) and the spheres being finitely
  presentable \oo-categories, we can apply the small object argument
  (Propositions~\ref{prop:small_obj_arg_1} and~\ref{prop:small_obj_arg_2})
  and we obtain:
\end{paragr}

\begin{prop}\label{prop:factor_cof_fib_triv}
  The following assertions hold:
  \begin{enumerate}
    \item The pair $(\{\text{cofibrations}\}, \{\text{trivial fibrations}\})$ is a
    weak factorization system on $\ooCat$.
    \item Every \oo-functor $f$ factors as $f = pi$, where $i$ is a relative
    polygraph and $p$~is a trivial fibration.
    \item The cofibrations are the retracts of relative polygraphs.
  \end{enumerate}
\end{prop}

\medbreak

We will now describe more concretely trivial fibrations and relative polygraphs.

\begin{paragr}[Trivial fibrations]\label{paragr:triv_fib}
  \index{fibration!trivial}
  By definition, an \oo-functor $f : C \to D$ is a trivial fibration if, for
  every $n \ge 0$, it has the right lifting property with respect to
  $\gencof{n} : \sphere{n} \to \globe{n}$. The data of a commutative diagram
  \[
    \xymatrix{
      \sphere{n} \ar[d]_{\gencof{n}} \ar[r] & C \ar[d]^f \\
      \globe{n} \ar[r] & D
    }
  \]
  is equivalent to the data of a pair $(x, y)$ of parallel $(n-1)$-cells of
  $C$ and of an $n$-cell $b : f(x) \to f(y)$ in $D$. (To make sense of this
  when $n = 0$, one has to consider that every \oo-category has a unique cell
  of dimension $-1$.) A lift for such a square is given by an $n$-cell $a :
  x \to y$ such that $f(a) = b$.

  In other words, $f : C \to D$ is a trivial fibration if and only if the
  following conditions hold:
  \begin{enumerate}
    \item $f$ is surjective on objects,
    \item for every $n \ge 0$, every pair $x, y$ of parallel $n$-cells of
    $C$ and every $(n+1)$\nbd-cells $b : f(x) \to f(y)$ in $D$, there exists an
    $(n+1)$-cell $a : x \to y$ in $C$ such that $f(a) = b$.
  \end{enumerate}
\end{paragr}

\begin{paragr}[Relative polygraphs]
  \index{relative polygraph}
  \index{polygraph!relative}
  By definition, an \oo-functor $f : C \to D$ is a relative polygraph if
  $f$ can be obtained as a countable composition
  \[
     \xymatrix@C=1.5pc{
     C = C_0 \ar[r]^-{i_1} & C_1 \ar[r]^{i_2} & \cdots \ar[r]^{i_n} & C_n
     \ar[r]^{i_{n+1}} & \cdots \pbox{,}
     }
  \]
  i.e., if $f$ is the canonical morphism $C_0 \to \limind_n C_n$,
  where $i_n$, for $n \ge 1$, is part of a pushout square
  \[
    \xymatrix{
      \coprod_{k \ge 0} E_k \times \sphere{k} \ar[r] \ar[d]_{\coprod_{k \ge
      0} E_k \times \gencof{k}} &
      C_{n-1} \ar[d]^{i_n} \\
      \coprod_{k \ge 0} E_k \times \globe{k} \ar[r] & C_n \pbox{,}
    }
  \]
  the $E_k$'s being sets (depending on $n$). In other words, $C_n$ is
  obtained from~$C_{n-1}$ by freely adding cells (of any dimension) and
  $i_n : C_{n-1} \to C_n$ is the canonical \oo-functor. This means that $D$
  is obtained from $C$ by adding cells (in any order) in several steps.
  Using the fact that the forgetful functor from $\nCat{m}$ to $\nCat{n}$,
  when $m > n$, respects colimits (see~\cref{paragr:adjunctions}) and in
  particular pushouts, one gets that the new cells can always be attached
  dimension by dimension. This means that we can suppose that~$i_n$, for $n \ge
  1$, is part of a pushout square of the simpler form
  \[
    \xymatrix{
      \coprod_{F_n} \sphere{n} \ar[r]
         \ar[d]_{\coprod_{F_n} \gencof{n}} &
      C_{n-1} \ar[d]^{i_n} \\
      \coprod_{F_n} \globe{n} \ar[r] & C_n \pbox{,}
    }
  \]
  where $F_n$ is any set. In other words, $i_n$ is the canonical \oo-functor
  associated to the cellular extension $(C_{n-1}, F_n)$. In particular, we
  get:
\end{paragr}

\begin{prop}
  An \oo-category $C$ is generated by a polygraph if and only if the unique
  functor from the initial \oo-category to $C$ is a relative polygraph.
\end{prop}

\begin{remark}
  We will see in \chapr{folk} that the previous proposition can be
  strengthened by saying that $C$ is generated by a polygraph if and only if
  the unique \oo-functor from the initial \oo-category to $C$ is a
  cofibration. In other words, the ``cofibrant objects'' are exactly the
  \oo-categories generated by polygraphs.
\end{remark}

\begin{proposition}
  Every cofibration is a monomorphism.
\end{proposition}

\begin{proof}
  Cofibrations are retracts of countable compositions of canonical
  \oo-func\-tors associated to cellular extensions. As retracts and
  countable compositions of monomorphic \oo-functors are monomorphisms, the
  result follows from the fact that the canonical \oo-functor associated to
  a cellular extension is a monomorphism
  (Proposition~\ref{prop:canonical_injection}).
\end{proof}

\section{Polygraphic Resolutions}
\label{S:PolygraphicResolutionsOmegaCategories}

Polygraphic resolutions are the \oo-categorical version of the free resolutions
of homological algebra, see Section~\ref{SS:ResolutionModules}. 
We will see shortly that every \oo-category admits a
polygraphic resolution and that such a resolution is in some sense unique up
to homotopy.

\begin{paragr}[Polygraphic resolutions]\label{paragr:pol_res}
  \index{resolution!polygraphic}
  \index{polygraphic!resolution}
  A \ndef{polygraphic resolution} of an \oo-category~$C$ is a pair $(P, p)$,
  where $P$ is a polygraph and $p : P^\ast \to C$ is a trivial fibration.
\end{paragr}

\begin{prop}\label{prop:pol_cell}
  Every \oo-category admits a polygraphic resolution.
\end{prop}

\begin{proof}
  Let $C$ be an \oo-category. Consider the unique \oo-functor $\varnothing_C
  : \varnothing \to C$ from the initial \oo-category $\varnothing$ to $C$.
  By Proposition~\ref{prop:factor_cof_fib_triv}, this $\omega$-functor
  factors as $p\circ\varnothing_{P^\ast}$, where $P$ is a polygraph,
  $\varnothing_{P^\ast} : \varnothing \to P^\ast$ is the unique such
  morphism and $p$ is a trivial fibration, thereby proving the result.
\end{proof}

The previous proposition shows the existence of a polygraphic resolution for any
\oo-category by abstract non-sense. Here is a canonical choice of such a
resolution:

\begin{paragr}[The canonical resolution]
  \index{canonical resolution}
  \index{resolution!canonical}
  Let $C$ be an \oo-category. The counit of the adjunction between the
  categories of polygraphs and \oo-categories gives an \oo-functor
  $\varepsilon_C : G(C)^\ast \to C$, where $G : \oCat \to \oPol$ denotes the
  functor described in \cref{subsec:omegapol}. This \oo-functor is a
  trivial fibration, essentially by definition. Indeed, it is bijective on
  objects by definition and, if $x$ and~$y$ are two parallel $k$\nbd-cells of
  $G(C)^\ast$, for some $k \ge 0$, and $z$ is a $(k+1)$-cell of~$G(C)^\ast$
  from~$\varepsilon_C(x)$ to~$\varepsilon_C(y)$, then, by definition of $G(C)$,
  there is a generating $(k+1)$\nbd-cell $p = (z, x, y)$ in $G(C)$ from $x$ to
  $y$ such that $\varepsilon_C(p) = z$. This means that $(G(C), \varepsilon_C)$ is
  a polygraphic resolution of $C$. This polygraphic resolution is
  called the \ndef{canonical resolution} of $C$. Note that it is functorial
  in $C$.
\end{paragr}

\section{Uniqueness of Polygraphic Resolutions}
\label{sec:pol-res-unique}

In this section, we prove that polygraphic resolutions are unique up to a
non-canonical homotopy in an appropriate sense. 

\begin{prop}
  Let $f : C \to D$ be an \oo-functor, let $(P, p)$ and $(Q, q)$ two
  polygraphic resolutions of $C$ and $D$, respectively. There exists a
  (non-canonical) \oo-functor $g : P^\ast \to Q^\ast$ making the square
  \[
    \xymatrix{
      P^\ast \ar[d]_{p} \ar[r]^g & Q^\ast \ar[d]^q \\
      C \ar[r]_f & D
    }
  \]
  commute.
\end{prop}

\begin{proof}
  Consider the commutative square
  \[
    \xymatrix{
      \varnothing \ar[d] \ar[rr] & & Q^\ast \ar[d]^q \\
      P^\ast \ar[r]_p & C \ar[r]_f & D \pbox{,}
    }
  \]
  where the unlabeled arrows are the unique such \oo-functors.
  By Proposition~\ref{prop:pol_cell}, the \oo-functor $\varnothing \to
  P^\ast$ is a cofibration. Since by definition, the \oo-functor $q$ is a
  trivial fibration, the square admits a lift $f$, thereby proving the
  result.
\end{proof}

\begin{rem}
  \label{rem:pol-res-func}
  In the particular case where $C = D$ and $f$ is the identity \oo-functor,
  we get the existence of an \oo-functor
  \[
    \xymatrix@C=1pc@R=2pc{
      P^\ast \ar[dr]_{p} \ar[rr]^g & & Q^\ast \ar[dl]^q \\
      & C
    }
  \]
  between any two polygraphic resolutions of $C$.
\end{rem}

\begin{remark}
  We will prove in \cref{chap:folk}, using the ``folk model
  structure'' on $\ooCat$, that the \oo-functor $g$ of the previous
  proposition is unique up to some appropriate notion of homotopy (see
  Proposition~\ref{prop:pol_res_uniq_hom}).
\end{remark}


\chapter[Toward the Folk Model Structure]{Toward the Folk Model Structure: \pdfoo-Equivalences}
\label{chap:w-eq}

In this chapter, we introduce all the notions and tools that will allow us to
define and establish the existence, in the next chapter, of the folk model
category structure on $\ooCat$.  We particularly focus on the concept of
\oo-equivalences, that will be the weak equivalences of this model
structure. Essentially all the material of this chapter is extracted from Lafont,
Métayer, and Worytkiewicz~\cite{LMW}, although cylinders first appeared in
work by Métayer~\cite{metayer2008cofibrant}.

The class of \oo-equivalences is the appropriate generalization to
\oo-categories of the class of equivalences of ordinary categories. In
particular, an \oo-equiv\-alence between $1$\nbd-categories is nothing but an
equivalence of categories. To define this notion, we need to generalize the
concept of an invertible cell (or isomorphism). This leads to the notion of
a reversible cell, which is, in intuitive terms, a cell admitting an inverse up to cells
admitting inverses, up to cells admitting inverses, etc. Another fundamental
tool is the \oo-category $\CatCylRev(C)$ of reversible cylinders in an
\oo-category $C$ that we will lead to a sensible notion of homotopy.

The first section is devoted to \oo-equivalences. We start
by defining the notion of a \emph{reversible cell} in an \oo-category. Following
\cite{LMW}, we define it by coinduction. We then introduce the class of
\emph{\oo-equivalences}. We prove some basic stability properties of this class. We
observe that trivial fibrations in the sense of the previous chapter are
\oo-equivalences. In the second section, we introduce the \oo-category
$\CatCyl(C)$ of \emph{cylinders} in an \oo-category $C$. The purpose of this
\oo-category is to allow the definition of \oo-functors playing the role of
homotopies. This leads to the definition of an \emph{oplax transformations},
generalizing the $1$-categorical natural transformations and the
$2$-categorical oplax transformations. The third section is about the
sub-\oo-category $\CatCylRev(C) \subset \CatCyl(C)$ of
\emph{reversible cylinders}. This is a fundamental tool to establish
properties of \oo-equivalences. We show that this \oo-category behaves as
a homotopical path object and we use its properties to show that
the class of \oo-equivalences satisfies the $2$-out-of-$3$ property, a
result that turns out to be non-trivial. In the fourth section, we introduce a
``coherent'' version of reversible cells leading to the notion of
a \emph{fibration} of \oo-categories. We show that an \oo-functor is a
trivial fibration if and only if it is both an \oo-equivalence and a
fibration. Finally, in the last section, we study the class of \oo-functors
having the left lifting properties with respect to fibrations (that could be
called trivial cofibrations). We show that it consists of \oo-functors
being both an \oo-equivalence and a cofibration. To do so, we introduce the
notion of an \emph{immersion}, which is a kind of strong deformation
retract, with respect to the ``path object'' of reversible cylinders.

\section{\pdfoo-Equivalences}

The notion of an \oo-equivalence is the higher-dimensional generalization of
the notion of an equivalence of categories. In particular, \oo-equivalent
\oo-categories can be considered as being close to be equal.
The definition is a bit involved and requires the introduction of the
auxiliary concept of a reversible cell.

\begin{paragr}[Reversible cells]
  \label{paragr:rev_cells}
  \index{reversible!cell}
  \index{cell!reversible}
  Given an \oo-category $C$, the notion of a reversible cell of $C$ is
  defined by \emph{coinduction} in the following way. An $n$-cell $u : x \to
  y$ of~$C$, for some $n > 0$, is \ndef{reversible} if there exists an
  $n$-cell $\weakinv{u} : y \to x$ in $C$ and $(n+1)$-cells
  \[
    \alpha : \unit{x} \to u \comp{n-1} \weakinv{u}
    \qqqtand
    \beta : \weakinv{u} \comp{n-1} u \to \unit{y}
  \]
  in $C$ that are both reversible.

  Concretely, this means that an $n$-cell $u : x \to y$ of $C$ is reversible
  if and only if it belongs to a set $X$ of cells of $C$ having the
  following property: for every~$n' > 0$ and every $n'$-cell $u' : x'
  \to y'$ in $X$, there exists an $n'$-cell $\weakinv{u'} : y' \to x'$ and $(n'+1)$-cells
  $\alpha' : \unit{x'} \to u' \comp{n'-1} \weakinv{u'}$ and $\beta' :
  \weakinv{u'} \comp{n'-1} u' \to \unit{y'}$ in $X$.
\end{paragr}

\begin{remark}
  Strictly speaking, coinduction allows to define algebraic structures. What we have
  really defined in the previous paragraph is a ``reversibility structure''
  where all the ``there exists'' are replaced by actual choices. This
  structure can somehow be flattened into a kind of tree. Then one can define a
  reversible cell as a cell that can be endowed with a reversibility
  structure.

  We will often say that ``we reason by coinduction''. This basically means that
  we are defining a ``reversibility structure'' according to its actual
  definition, producing choices of $\weakinv{u}$, $\alpha$, and $\beta$
  as in the definition. In particular, if $R$ is a set of cells of an
  \oo-category $C$, ``proving by coinduction'' that cells of $R$ are
  reversible will consist in producing, for every $n$-cell $u$ of $R$,
  a formula giving $\weakinv{u}$, $\alpha$, and $\beta$ as in the definition
  of a reversible cell \emph{assuming} that the $(n+1)$-cells of $R$ are
  reversible.
\end{remark}

\begin{paragr}[Weak inverses]\label{paragr:w_inv}
  Let $n \ge 0$. If an $n$-cell $u : x \to y$ in some \oo-category is
  reversible, then any $\weakinv{u} : y \to x$ as in the definition will be
  called a \ndef{weak inverse}\index{weak inverse}\index{inverse!weak} of
  $u$.

  If $C$ is a $1$-category, a reversible $1$-cell of $C$ is the
  same thing as an isomorphism of $C$ and a weak inverse is an inverse. This
  follows from the fact that units are reversible (see the next lemma).
  Similarly, if $C$ is a $2$-category, a $1$-cell $u : x \to y$ of $C$ is
  reversible if and only if it is an equivalence, that is, if and only if there
  exists a $1$-cell $\weakinv{u} : y \to x$ and $2$-cells $\alpha : \unit{x}
  \to u \comp{0} \weakinv{u}$ and $\beta : \weakinv{u} \comp{0} u \to
  \unit{y}$ that are isomorphisms.
\end{paragr}

\begin{lemma}\label{lemma:w_inv}
  Let $C$ be an \oo-category.
  \begin{enumerate}
    \item If $x$ is a cell of $C$, then $\unit{x}$ is reversible.
    \item If $u : x \to y$ is a reversible cell of $C$, then any weak
      inverse $\weakinv{u} : y \to x$ of $u$ is reversible and has $u$ as a
      weak inverse.
    \item If $u : x \to y$ and $v : y \to z$ are two reversible $n$-cells of
    $C$ for an $n > 0$, then $u \comp{n-1} v$ is reversible.
    \item More generally, if $u$ and $v$ are two reversible $n$-cells of $C$
      such that the composition $u \comp{i} v$ is defined for some $0 \le i
      < n$, then this composition is reversible.
  \end{enumerate}
\end{lemma}

\begin{proof}
  We proceed by coinduction.
  \begin{enumerate}
    \item Let $x$ be an $n$-cell of $C$. Set $u = \weakinv{u} =
      \unit{x}$. Then $u \comp{n} \weakinv{u} = \weakinv{u} \comp{n} u =
      \unit{\unit{x}}$ is reversible by coinduction, and so is $\unit{x}$ by
      definition.
    \item This is immediate by the symmetry in the definition of a
      reversible cell.
    \item By definition, there exist cells $\weakinv{u} : y \to x$,
      $\weakinv{v} : z \to y$ and reversible cells
      $\alpha : \unit{x} \to u \comp{n-1} \weakinv{u}$,
      $\beta : \weakinv{u} \comp{n-1} u \to \unit{y}$,
      $\gamma : \unit{y} \to v \comp{n-1} \weakinv{v}$,
      $\delta : \weakinv{v} \comp{n-1} v \to \unit{z}$. We get cells
      \begin{align*}
        \alpha \comp{n} (u \comp{n-1} \gamma \comp{n-1} \weakinv{u}) & :
        \unit{x} \to
        u \comp{n-1} v \comp{n-1} \weakinv{v} \comp{n-1} \weakinv{u} \\
        (\weakinv{v} \comp{n-1} \beta \comp{n-1} v) \comp{n} \delta & :
        \weakinv{v} \comp{n-1} \weakinv{u} \comp{n-1} u \comp{n-1} v \to \unit{z}.
      \end{align*}
      By coinduction, it suffices to show that the cells
      $\alpha \comp{n} (u \comp{n-1} \gamma \comp{n-1} \weakinv{u})$ and
      $(\weakinv{v} \comp{n-1} \beta \comp{n-1} v) \comp{n} \delta$
      are reversible. It thus suffices to show that the whiskering of a
      reversible cell by any cell is reversible. This can be shown by a new
      coinduction but also follows from the next proposition.
    \item The case where $i = n - 1$ is the previous assertion. If $i < n -
      1$, then, using the exchange law, this follows from the case $i = n -
      1$ and the case of a whiskering (see the proof of the previous
      assertion). \qedhere
  \end{enumerate}
\end{proof}

\begin{proposition}\label{prop:funct_oo-equiv}
  Let $f : C \to D$ be an \oo-functor. If $u$ is a reversible cell of~$C$,
  then $f(u)$ is a reversible cell of $D$.
\end{proposition}

\begin{proof}
  Let $u : x \to y$ be a reversible $n$-cell of $C$. By definition, there
  exists a cell $\weakinv{u} : y \to x$ and two reversible cells $\alpha :
 \unit{x} \to u \comp{n-1} \weakinv{u}$ and $\beta : \weakinv{u} \comp{n-1}
  u \to \unit{x}$. By coinduction, the cells
  $f(\alpha) :
  \unit{f(x)} \to f(u) \comp{n-1} f(\weakinv{u})$ and $f(\beta) :
  f(\weakinv{u}) \comp{n-1} f(u) \to \unit{f(y)}$ are reversible. This means
  that $f(u)$ is indeed reversible.
\end{proof}

\begin{paragr}[\pdfoo-equivalent cells]\label{paragr:oeq-cells}
  Let $n \ge 0$. Two $n$-cells $x$ and $y$ of an \oo-category~$C$ are
  \ndef{\oo-equivalent}\index{9-equivalence@\oo-equivalence!of cells}\index{cell!9-equivalence@\oo-equivalence}
  if there exists a reversible cell $u : x \to y$ in $C$. We denote this
  relation by $x\oeq y$.

  We chose to define the notion of being \oo-equivalent in terms of the
  notion of a reversible cell but we could also have defined this notion
  directly using coinduction: two parallel $n$-cells $x$ and $y$ are
  \oo-equivalent if there exists $(n+1)$\nbd-cells $u : x \to y$ and
  $\weakinv{u} : y \to x$ such that $\unit{x}$ and $u \comp{n} \weakinv{u}$
  are \oo-equivalent, and $\weakinv{u} \comp{n} u$ and $\unit{y}$ are
  \oo-equivalent.

  It follows from~\cref{paragr:w_inv} that \oo-equivalent objects
  in a $1$-category are isomorphic objects and that \oo-equivalent objects
  in a $2$-category are equivalent objects.
\end{paragr}

\begin{prop}\label{prop:omegaeq_cong}
  Let $C$ be an \oo-category. The relation ``being \oo-equivalent'' is a
  congruence relation on the set of cells of $C$ in the sense that:
 \begin{enumerate}
   \item This relation is an equivalence relation.
   \item This equivalence relation is compatible with
     compositions: if $u$ and $u'$, and $v$ and $v'$ are \oo-equivalent cells,
     then $u \comp{i} v$ is \oo-equivalent to $u' \comp{i} v'$ when these
     compositions make sense.
 \end{enumerate}
\end{prop}

\begin{proof}
  This is a direct consequence of Lemma~\ref{lemma:w_inv}.
\end{proof}

\begin{proposition}\label{prop:img_w-equiv}
  Let $f : C \to D$ be an \oo-functor. If two cells $x$ and $y$ of~$C$ are
  \oo-equivalent, then so are $f(x)$ and $f(y)$.
\end{proposition}

\begin{proof}
  This follows from Proposition~\ref{prop:funct_oo-equiv}
\end{proof}

\begin{paragr}
  We now turn to a result of paramount importance in the construction
  of the folk model structure, namely the ``division Lemma''~(\cite[Lemma
  4.6]{LMW}). As the proof is quite intricate, it will be convenient
  to introduce the following terminology: for any property $\mathcal P$ applying to $n$-cells, we say that there
is a {\em weakly unique} $n$-cell satisfying $\mathcal P$ whenever any
two $n$-cells satisfying $\mathcal P$ are $\omega$-equivalent. 
\end{paragr}

\begin{lemma}[Division lemma]\label{lemma:weakdiv}
  \index{division lemma}
  \index{lemma!division}
  Any reversible 1-cell $u:x\to y$ satisfies the following {\em left division property}:
  \begin{itemize}
  \item For any 1-cell $w : x \to z$, there is a weakly unique
    1-cell $v : y \to z$ such that $u \comp 0 v \oeq w$.
  \item For any pair of parallel 1-cells $a, b : y \to z$ and any 2-cell $w : u
    \comp 0 a \to u \comp 0 b$, there is a weakly unique 2-cell $v :
    a\to b$ such that $u \comp 0 v \oeq w$.
 \item More generally, for any $n > 0$, any pair of parallel $n$-cells $a,
   b$ such that $\sce 0(a)=\sce 0(b)=y$, $\tge 0(a)=\tge 0(b)=z$ and any
   $(n+1)$-cell $w : u \comp 0 a \to u \comp 0 b$, there is a weakly unique
   $(n+1)$-cell $v : a \to b$ such that $u \comp 0 v \oeq w$.   
 \end{itemize}
 Likewise, reversible 1-cells satisfy the corresponding right
 division property.  
\end{lemma}

\begin{proof}
  Let $u:x\to y$ be a reversible $1$-cell. By definition there is a
  weak inverse $\wki u:y\to x$ together with a reversible $2$-cell
  $\beta:\wki{u}\comp 0 u\to \unit{y}$, as well as a reversible
  $2$-cell $\wki{\beta}:\unit{y}\to \wki{u}\comp 0 u$ such that
  $\beta\comp 1\wki{\beta}\oeq \unit{\wki{u}\comp 0 u}$.
  \begin{itemize}
  \item In the first case, define $v=\wki{u}\comp 0 w:y\to
    z$. Then $u\comp 0 v=u\comp 0 \wki{u}\comp 0 w\oeq
    w$. Moreover, for any $v':y\to z$ such that $u\comp 0 v'\oeq w$,
    we get
    $\wki{u}\comp 0 u\comp 0 v' \oeq \wki{u}\comp 0 w = v$,
    whence $v'\oeq v$ and weak uniqueness holds. Likewise, the right
    division property holds in the first case.
  \item   Let $a,b:y\to z$ be $1$-cells and $w:u\comp 0 a\to
    u\comp 0 b$. Suppose there is a $2$-cell $v:a\to b$ such that $u\comp 0
    v\oeq w$. Then, by applying exchange and compatibility of
    $\omega$-equivalence with compositions, one gets
    \begin{eqnarray*}
      (\wki{\beta}\comp 0 a)\comp 1 (\wki{\beta}\comp 0 w) \comp1 (\beta\comp 0 b) & \oeq &
                                                               (\wki{\beta}\comp 0 a)\comp 1 (u\comp 0 v)\comp 1 (\beta\comp 0 b)\\
                                                      & \oeq & (\wki{\beta}\comp 1 \beta)\comp 0 v\\
                                                      & \oeq & v,
    \end{eqnarray*}
    which implies weak uniqueness for $v$. As for existence, 
    define
    \[
      v=(\wki{\beta}\comp 0 a)\comp 1 (\wki{\beta}\comp 0 w) \comp 1 (\beta\comp 0 b).
    \]
    By definition, $v:a\to b$. Consider now
    \[v'=(\beta\comp 0 a)\comp 1 v \comp 1 (\wki{\beta}\comp 0 b).\]
    Again, exchange and compatibility yield
    \[v'\oeq \wki{u}\comp 0 w\]
    but also
    \[v'\oeq \wki{u}\comp 0 (u\comp 0 v). \]
    By applying weak uniqueness to left division by $\wki{u}$, one
    gets $w\oeq u\comp 0 v$, whence the result. Right division is
    proved accordingly in this case.
   \item The general case is proved by induction on $n$. The case
     $n=1$ has been just proved above. Let now $n>1$ and suppose that
     the property of left and right division holds for any reversible
     $1$-cell in any $\omega$-category $C$ up to dimension $n-1$. Let
     $a$, $b$ be parallel $n$-cells such that
     $\sce{0}(a)=\sce{0}(b)=y$,  $\tge{0}(a)=\tge{0}(b)=z$ and
     $w:u\comp 0 a\to u\comp 0 b$ an $(n+1)$-cell in $C$. Consider
     \[w'=(\wki{u}\comp 0 w)\comp 1 (\beta\comp 0\tge{1}(a)).\]
     We have
     \[\sce{n}(w') =(\beta\comp 0\sce{1}(a))\comp 1 a\quad\hbox{and}\quad
       \tge{n}(w') =(\beta\comp 0\sce{1}(a))\comp 1 b. \]
   Now $\beta\comp 0\sce{1}(a)$ may be seen as a reversible $1$-cell in
   the $\omega$-category $C(y,z)$, so that the induction hypothesis
   applies and yields a weakly unique $(n+1)$-cell $v:a\to b$ such
   that
   \begin{equation}
     \label{eq:division}
     (\beta\comp 0\sce{1}(a))\comp 1 v\oeq w'.
   \end{equation}
     One then shows that any $(n+1)$-cell $v':a\to b$ such that $u\comp 0 v'\oeq w$
     satisfies the equation~(\ref{eq:division}). Therefore, by
     induction, $v'\oeq v$ and weak uniqueness is proved. It remains
     to check that the above $(n+1)$-cell $v$ satisfies $u\comp 0
     v\oeq w$. Rewriting~(\ref{eq:division}) using the exchange law, we get
     \[(\wki{u}\comp 0 u\comp 0 v)\comp 1 (\beta\comp 0 \tge{1}(a))\oeq
       (\wki{u}\comp 0 w)\comp 1 (\beta\comp 0 \tge{1}(a)).\]
     By induction applied to the reversible $1$-cell $\beta\comp 0
     \tge{1}(a)$ of $C(y,z)$, weak uniqueness for right division
     implies
     \[\wki{u}\comp 0 u\comp 0 v\oeq \wki{u}\comp 0  w,\]
     and finally, by weak uniqueness applied to $\wki{u}$, $u\comp 0
     v\oeq w$, as required.
     \qedhere
  \end{itemize}
\end{proof}


\begin{paragr}[\pdfoo-equivalences]\label{paragr:def_oo-equiv}
  An \oo-functor $f : C \to D$ is an
  \ndef{\oo-equivalence}\index{9-equivalence@\oo-equivalence} if the following
  conditions are satisfied:
  \begin{enumerate}
    \item For every object $y$ of $D$, there exists an object $x$ of $C$
      such that $f(x)$ is \oo-equivalent to $y$.
    \item For every $n \ge 0$, every pair of parallel $n$-cells $x$ and $y$
      of $C$ and every $(n+1)$\nbd-cell $v : f(x) \to f(y)$, there exists an
      $(n+1)$-cell $u : x \to y$ such that $f(u)$ is \oo-equivalent to $v$.
  \end{enumerate}

  If $f : C \to D$ is a $1$-functor between $1$-categories, then $f$ is an
  \oo-equivalence if and only if it is an equivalence of categories.
  Similarly, if $f : C \to D$ is a $2$-functor, then $f$ is an
  \oo-equivalence if and only if $f$ is a biequivalence. Note that an
  \oo-equivalence between two $2$-categories is not the same thing as a
  biequivalence as an \oo-equivalence is required to be a strict
  $2$-functor as opposed to a bifunctor.
\end{paragr}

\begin{rem}
  In general, an \oo-equivalence does \emph{not} admit an inverse, in any
  reasonable sense, that is a \emph{strict} \oo-equivalence. Morally, it
  admits a \emph{weak} \oo-functor as a weak inverse in some sense. For
  instance, the obvious \oo-functor from the ``pseudo-$2$-triangle''
  \[
    \xymatrix{
      & \cdot
      \\
      \cdot \ar[r] \ar[ur]_{}="s" & \cdot \ar[u]
      \ar@{}"s";[]_(0.05){}="ss"_(0.85){}="tt"
      \ar@2"ss";"tt"^(0.30)\alpha
      \pbox{,}
    }
  \]
  where $\alpha$ is a reversible cell, to the commutative triangle, is
  easily seen to be an \oo-equivalence, but there is no \oo-equivalence from
  the commutative triangle to the pseudo-$2$-triangle.
\end{rem}

\begin{proposition}\label{prop:triv_fib_eq}
  Trivial fibrations are \oo-equivalences.
\end{proposition}

\begin{proof}
  This follows immediately from the characterization of trivial fibrations
  given at the end of~\cref{paragr:triv_fib} and the fact that
  identities are reversible cells.
\end{proof}

An \oo-equivalence is injective up to \oo-equivalence of cells in the
following sense:

\begin{proposition}\label{prop:oo-equiv_inj}
  Let $f : C \to D$ be an \oo-equivalence and let $x$ and $y$ be two parallel
  cells of $C$. If $f(x)$ and $f(y)$ are \oo-equivalent, then $x$ and
  $y$ are \oo-equivalent.
\end{proposition}

\begin{proof}
  By definition, there exists a cell $v : f(y) \to f(x)$ such that
  $\unit{f(x)}$ and $f(u) \comp{n-1} v$, and $v \comp{n-1} f(u)$ and
  $\unit{f(y)}$ are \oo-equivalent cells. Since $f$ is an \oo-equivalence,
  there exists a cell $\weakinv{u} : y \to x$ such that $f(\weakinv{u})$ is
  \oo-equivalent to $v$. Using Proposition~\ref{prop:omegaeq_cong}, we get
  that $f(\unit{x})$ and $f(u \comp{n-1} \weakinv{u})$, and $f(\weakinv{u}
  \comp{n-1} u)$ and $f(\unit{y})$ are \oo-equivalent cells. The result thus
  follows by coinduction.
\end{proof}

\begin{prop}\label{prop:w_eq_comp}
  The composition of two \oo-equivalences is an \oo-equiv\-alence.
\end{prop}

\begin{proof}
  Let $f : C \to D$ and $g : D \to E$ be two \oo-equivalences. Let us prove
  that $gf : C \to E$ is an \oo-equivalence.
  \begin{enumerate}
    \item Let $z$ be an object of $E$. Since $g$ is an \oo-equivalence,
      there exists an object~$y$ of~$D$ such that $g(y)$ and $z$ are
      \oo-equivalent. Similarly, since $f$ is an \oo-equivalence, there
      exists an object $x$ of $C$ such that $f(x)$ and $y$ are
      \oo-equivalent. It follows from
      Proposition~\ref{prop:funct_oo-equiv} that $g(f(x))$ and $g(y)$ are
      \oo-equivalent. Hence, by transitivity of the relation of
      \oo-equivalence, $gf(x)$ and $z$ are \oo-equivalent.
    \item Let $n \ge 0$, let $x$ and $y$ be two parallel $n$-cells of
      $C$, and let $w : gf(x) \to gf(y)$ be an $(n+1)$-cell of $E$. Since
      $g$ is an \oo-equivalence, there exists an $(n+1)$\nbd-cell $v : f(x) \to
      f(y)$ of $D$ such that $g(v)$ and $w$ are \oo-equivalent. Similarly,
      since $f$ is an \oo-equivalence, there exists an $(n+1)$-cell $u : x
      \to y$ such that $f(u)$ and $v$ are \oo-equivalent. It follows from
      \cref{prop:img_w-equiv} that $gf(u)$ and $g(v)$ are \oo-equivalent,
      and hence that $gf(u)$ and $w$ are \oo-equivalent,
      thereby proving the result.  \qedhere
  \end{enumerate}
\end{proof}

\begin{paragraph}[2-out-of-3 property]
Recall that a class of maps $\clW$ in a category $\C$ is said to satisfy the
\ndef{2-out-of-3 property}\index{2-out-of-3 property} if for any commutative triangle
\[
\xymatrix{
  X \ar[rr]^f \ar[dr]_h & & Y \ar[dl]^g \\
                    & Z
}
\]
in $\C$, if two morphisms among $f$, $g$, and $h$ are in $\clW$, then so is
the third one.
For instance, isomorphisms in a category satisfy the 2-out-of-3 property.
More generally, any reasonable notion of ``equivalence'' in a category
should satisfy this property.
\end{paragraph}

\medskip

The 2-out-of-3 property is made of three different properties, depending on
which of two morphisms among $f$, $g$, and $h$ are assumed to be in
$\clW$. If $\clW$ is the class of \oo-equivalences, the previous proposition
gives the case where these two morphisms are $f$ and $g$. The following
proposition is the case where they are $g$ and $h$:

\begin{proposition}\label{prop:2-out-of-3-easy}
  Let $f : C \to D$ and $g : D \to E$ be two \oo-functors. If $g$ and $gf$
  are \oo-equivalences, then so is $f$.
\end{proposition}

\begin{proof}
  Let us prove that $f$ is an \oo-equivalence.
  \begin{enumerate}
    \item Let $y$ be an object of $D$. Consider the object $g(y)$ of $E$.
    Since $gf$ is an \oo-equivalence, there exists an object $x$ of $C$ such
    that $gf(x)$ is \oo-equivalent to $g(y)$. By
    Proposition~\ref{prop:oo-equiv_inj}, this implies that $f(x)$ and $y$
    are \oo-equivalent.
    \item Let $n \ge 0$, let $x$ and $y$ be two parallel $n$-cells of $C$
    and let $v : f(x) \to f(y)$ be an $(n+1)$-cell of $D$. Consider the
    $(n+1)$-cell $g(v) : gf(x) \to gf(y)$ of~$E$. Since $gf$ is an
    \oo-equivalence, there exists an $(n+1)$-cell $u : x \to y$ such that
    $gf(u)$ and $g(v)$ are \oo-equivalent. By
    Proposition~\ref{prop:oo-equiv_inj}, this implies that $f(u)$ and $v$
    are \oo-equivalent. \qedhere
  \end{enumerate}
\end{proof}

The remaining case of the 2-out-of-3 property is much harder to prove and
requires the introduction of the \oo-category of cylinders, which is the
topic of the next section.

We end the section with two easy stability conditions of the class of
\oo-equiv\-alences.

\begin{proposition}\label{prop:w-eq_retr}
  A retract of an \oo-equivalence is an \oo-equivalence.
\end{proposition}

\begin{proof}
  Consider a commutative diagram
  \[
      \xymatrix{
        C \ar[d]_{f} \ar[r]_i \ar@/^3ex/[rr]^{\unit{C}} & C' \ar[d]_{f'}
        \ar[r]_p & C \ar[d]^{f} \\
        D \ar[r]^j \ar@/_3ex/[rr]_{\unit{D}} & D' \ar[r]^q & D
      }
  \]
  where $f'$ is an \oo-equivalence. Let us prove that $f$ is an
  \oo-equivalence.
  \begin{enumerate}
    \item Let $y$ be an object of $D$. Since $f'$ is an \oo-equivalence, there
      exists a cell~$x'$ in $C'$ such that $f'(x') = y'$, where $y' = j(y)$.
      We thus get an object $x = p(x')$ in $C$. As $f(x) = fp(x') = qf'(x)$,
      the cell $f(x)$ is \oo-equivalent to $qy' = qj(y) = y$.
    \item Let $x$ and $y$ be parallel cells of $C$ and let $v : f(x) \to
      f(y)$ be a cell of $D$. Consider the cell $v' = j(v) : jf(x) \to
      jf(y)$. We have $v' : f'i(x) \to f'i(y)$ and, since $f'$ is an
      \oo-equivalence, there exists a cell $u' : i(x) \to i(y)$ such that
      $f(u')$ is \oo-equivalent to $v'$. Setting $u = pf(u')$, we get as in
      the previous point that $f(u)$ is \oo-equivalent to $v$.
      \qedhere
  \end{enumerate}
\end{proof}

\begin{proposition}\label{prop:w-eq_filtered}
  A filtered colimit of \oo-equivalences is an \oo-equivalence. In
  particular, a countable composition (or more generally a transfinite
  composition) of \oo-equivalences is an \oo-equivalence.
\end{proposition}

\begin{proof}
  Let $I$ be a filtered category, let $F, G : I \to \ooCat$ be two functors,
  and let $f : F \to G$ a natural transformation such that for every
  $i$ in $I$, $f_i : F(i) \to G(i)$ is an \oo-equivalence. Let us prove
  that $f_\infty = \colim f : \colim F \to \colim G$ is an \oo-equivalence. Recall
  (see Proposition~\ref{prop:filcolim}) that the forgetful functor from
  \oo-categories to globular sets respects filtered colimits. This means
  that for every $n \ge 0$, the set of $n$-cells of $\colim F$ is the
  colimit of the sets of $n$-cells of the~$F(i)$. If $x$ is an $n$-cell of
  $F(i)$, we will denote by $[x]$ the corresponding $n$-cell of $\colim F$.
  \begin{enumerate}
    \item Let $[y]$ be an object of $\colim G$ coming from an object $y$ of $G(i)$
      for some~$i$ in $I$. Since $f_i : F(i) \to G(i)$ is an
      \oo-equivalence, there exists an object~$x$ of $F(i)$ such that
      $\alpha_i(x)$ is \oo-equivalence to $y$. By definition, we have
      $f_\infty([x]) = [y]$. Moreover, since the canonical \oo-functor $G(i)
      \to \colim G$, as any \oo-functor, sends \oo-equivalent cells to
      \oo-equivalent cells, $[x]$ and $[y]$ are \oo-equivalent.
    \item Let $n \ge 0$ and let $[x]$ and $[y]$ be two parallel $n$-cells of
      $\colim F$, where $x$ is an $n$-cell of $F(i)$ and $y$ is an $n$-cell
      of $F(j)$. Let $[u] : f_\infty([x]) \to f_\infty([y])$ be an
      $(n+1)$-cell of $\colim G$, where $u$ is an $(n+1)$-cell of $G(k)$.
      Using the fact that $I$ is filtered, one can suppose that $i = j = k$
      and $u : f_i(x) \to f_i(y)$. Since $f_i$ is an \oo-equivalence, there
      exists $v : x \to y$ such that $f_i(u)$ is \oo-equivalent to $u$. This
      implies that $[u] : [x] \to [y]$ and that $f_\infty([u])$ is
      \oo-equivalent to~$[v]$.
      \qedhere
  \end{enumerate}
\end{proof}

\section{The \pdfoo-Category of Cylinders}

In this section, we fix an \oo-category $C$.

\begin{paragraph}[Cylinders]\label{paragr:def_cyl}
  For $n \ge 0$, an \ndef{$n$-cylinder in $C$}\index{cylinder!in an \oo-category} is
  given by two $n$-cells $x$ and~$y$ of $C$, and a sequence
  \[
    \alpha_0^-, \alpha_0^+, \dots, \alpha_{n-1}^-, \alpha_{n-1}^+,
    \alpha_n^- = \alpha^{}_n = \alpha_n^+
  \]
  where $\alpha_j^\e$, for $0 \le j \le n$ and $\e = \pm$,
  is a $(j+1)$-cell whose source and target are
  \begin{align*}
    \alpha_j^- & :
    \sce{j}(x) \comp{0} \alpha_0^+ \comp{1} \cdots \comp{j-1} \alpha_{j-1}^+
    \to
    \alpha_{j-1}^- \comp{j-1} \cdots \comp{1} \alpha_0^- \comp{0}
    \sce{j}(y),
    \\
    \alpha_j^+ & :
    \tge{j}(x) \comp{0} \alpha_0^+ \comp{1} \cdots \comp{j-1} \alpha_{j-1}^+
    \to
    \alpha_{j-1}^- \comp{j-1} \cdots \comp{1} \alpha_0^- \comp{0}
    \tge{j}(y).
  \end{align*}
  Note that for $n = j$, we get
  \[
    \alpha_n
    :
    x \comp{0} \alpha_0^+ \comp{1} \cdots \comp{n-1} \alpha_{n-1}^+
    \to
    \alpha_{n-1}^- \comp{n-1} \cdots \comp{1} \alpha_0^- \comp{0} y.
  \]
  We say that such a cylinder is a \ndef{cylinder from $x$ to $y$}
  and we write $\alpha : x \cylto y$.

  Here are pictures of a $0$-cylinder, a $1$-cylinder, and a $2$-cylinder
  in $C$:
  \[
   \xymatrix@R=3pc{
   x \ar[d]_{\alpha_0} \\
   y
   }
   \qquad
   \qquad
    \xymatrix@C=3pc@R=3pc{
      \sce{0}(x) \ar[r]^x \ar[d]_{\alpha^-_0} &
      \tge{0}(x) \ar[d]^{\alpha^+_0} \\
      \sce{0}(y) \ar[r]_y & \tge{0}(y)
      \ar@{}[u];[l]_(.30){}="s"
      \ar@{}[u];[l]_(.70){}="t"
      \ar@2"s";"t"_{\alpha_1}
    }
   \qquad
   \qquad
    \xymatrix@C=3pc@R=3pc{
      \sce{0}(x)
      \ar@/^2ex/[r]^(0.70){\sce{1}(x)}_{}="0"
      \ar@/_2ex/[r]_(0.70){\tge{1}(x)}_{}="1"
      \ar[d]_{}="f"_{\alpha^-_0}
      \ar@2"0";"1"_{x\,\,}
      &
      \tge{0}(x)
      \ar[d]^{\alpha^+_0} \\
      \sce{0}(y)
      \ar@{.>}@/^2ex/[r]^(0.30){\sce{1}(y)}_{}="0"
      \ar@/_2ex/[r]_(0.30){\tge{1}(y)}_{}="1"
      \ar@{:>}"0";"1"_{y\,\,}
      &
      \tge{0}(y)\,.
      \ar@{}[u];[l]_(.40){}="x"
      \ar@{}[u];[l]_(.60){}="y"
      \ar@<-1.5ex>@/_1ex/@{:>}"x";"y"_(0.60){\alpha^-_1\,}_{}="0"
      \ar@<1.5ex>@/^1ex/@2"x";"y"^(0.40){\!\!\!\alpha^+_1}_{}="1"
      \ar@{}"1";"0"_(.05){}="z"
      \ar@{}"1";"0"_(.95){}="t"
      \ar@3{>}"z";"t"_{\alpha}
    }
  \]
\end{paragraph}

\begin{paragraph}[Inductive definition of cylinders]
  Alternatively, the notion of an $n$\nbd-cyl\-inder~$\alpha : x \cylto y$, where
  $x$ and $y$ are $n$-cells of $C$, can defined inductively in the following
  way. A $0$-cylinder $\alpha : x \to y$ is a $1$-cell $\alpha_0 : x \to y$.
  For $n > 0$, an $n$-cylinder $\alpha : x \cylto y$ consists of two $1$-cells
  \[
    \alpha_0^- : \sce{0}(x) \to \sce{0}(y)
    \qqqtand
    \alpha_0^+ : \tge{0}(x) \to \tge{0}(y)
  \]
  and an $(n-1)$-cylinder
  \[ [\alpha] : [x \comp{0} \alpha_0^+] \cylto [\alpha_0^- \comp{0} y] \]
  in the \oo-category $C(\sce{0}(x), \tge{0}(y))$, where
  $[x \comp{0} \alpha_0^+]$ and $[\alpha_0^- \comp{0} y]$ denote the
  $n$-cells $x \comp{0} \alpha_0^+$ and $\alpha_0^- \comp{0} y$ of $C$ seen
  as $(n-1)$-cells of $C(\sce{0}(x), \tge{0}(y))$.

  The equivalence between the two definitions easily follows by induction.
  This second definition is the one used \cite{LMW}. Its main advantage
  is to allow inductive arguments on cylinders.
\end{paragraph}

\medbreak

We will see that the cylinders in $C$ organize themselves in an
\oo-category. We now define the associated operations.

\begin{paragraph}[Source and target of a cylinder]
  Let $x$ and $y$ be two $n$-cells for some~$n > 0$, and let $\alpha : x
  \cylto y$ be an $n$-cylinder. We define the \ndef{source} of this cylinder
  to be the $(n-1)$-cylinder
  \[
    \sce{n-1}(\alpha) : \sce{n-1}(x) \cylto \sce{n-1}(y)
  \]
  given by the cells
  \[
    \alpha_0^-, \alpha_0^+, \dots, \alpha_{n-2}^-, \alpha_{n-2}^+,
    \alpha_{n-1}^-.
  \]
  Similarly, the \ndef{target} of such a cylinder is the $(n-1)$-cylinder
  \[
    \tge{n-1}(\alpha) : \tge{n-1}(x) \cylto \tge{n-1}(y)
  \]
  defined by the cells
  \[
    \alpha_0^-, \alpha_0^+, \dots, \alpha_{n-2}^-, \alpha_{n-2}^+,
    \alpha_{n-1}^+.
  \]
  Geometrically, the source of a cylinder is the ``back face'' of this
  cylinder and the target is the ``front face''. For instance, the source
  and target of a $2$-cylinder
  \[
    \xymatrix@C=3pc@R=3pc{
      \sce{0}(x)
      \ar@/^2ex/[r]^(0.70){\sce{1}(x)}_{}="0"
      \ar@/_2ex/[r]_(0.70){\tge{1}(x)}_{}="1"
      \ar[d]_{}="f"_{\alpha^-_0}
      \ar@2"0";"1"_{x\,\,}
      &
      \tge{0}(x)
      \ar[d]^{\alpha^+_0} \\
      \sce{0}(y)
      \ar@{.>}@/^2ex/[r]^(0.30){\sce{1}(y)}_{}="0"
      \ar@/_2ex/[r]_(0.30){\tge{1}(y)}_{}="1"
      \ar@{:>}"0";"1"_{y\,\,}
      &
      \tge{0}(y)
      \ar@{}[u];[l]_(.40){}="x"
      \ar@{}[u];[l]_(.60){}="y"
      \ar@<-1.5ex>@/_1ex/@{:>}"x";"y"_(0.60){\alpha^-_1\,}_{}="0"
      \ar@<1.5ex>@/^1ex/@2"x";"y"^(0.40){\!\!\!\alpha^+_1}_{}="1"
      \ar@{}"1";"0"_(.05){}="z"
      \ar@{}"1";"0"_(.95){}="t"
      \ar@3{>}"z";"t"_{\alpha}
    }
  \]
  are the $1$-cylinders
  \[
    \vcenter{
      \xymatrix{
        \sce{0}(x) \ar[r]^{\sce{1}(x)} \ar[d]_{\alpha^-_0} &
        \tge{0}(x) \ar[d]^{\alpha^+_0} \\
        \sce{0}(y) \ar[r]_{\sce{1}(y)} & \tge{0}(y)
        \ar@{}[u];[l]_(.30){}="s"
        \ar@{}[u];[l]_(.70){}="t"
        \ar@2"s";"t"_{\alpha^-_1}
      }
    }
    \qqtand
    \vcenter{
      \xymatrix{
        \sce{0}(x) \ar[r]^{\tge{1}(x)} \ar[d]_{\alpha^-_0} &
        \tge{0}(x) \ar[d]^{\alpha^+_0} \\
        \sce{0}(y) \ar[r]_{\tge{1}(y)} & \tge{0}(y)
        \ar@{}[u];[l]_(.30){}="s"
        \ar@{}[u];[l]_(.70){}="t"
        \ar@2"s";"t"_{\alpha^+_1}
        \pbox{.}
      }
    }
  \]
\end{paragraph}

\begin{paragraph}[Unit of a cylinder]
  Let $\alpha : x \cylto y$ be an $n$-cylinder in $C$. We define the \ndef{unit}
  of this cylinder to be the $(n+1)$-cylinder
  \[
    \unit{\alpha} : \unit{x} \cylto \unit{y}
  \]
  defined by the cells
  \[
    \alpha_0^-, \alpha_0^+, \dots, \alpha_{n-1}^-, \alpha_{n-1}^+,
    \alpha_n, \alpha_n, \unit{\alpha_n}.
  \]
  Geometrically, this cylinder is obtained by gluing two copies of $\alpha$
  and putting a trivial cell in the middle. For instance, the unit of the
  $1$-cylinder
  \[
    \xymatrix{
      \sce{0}(x) \ar[r]^x \ar[d]_{\alpha^-_0} &
      \tge{0}(x) \ar[d]^{\alpha^+_0} \\
      \sce{0}(y) \ar[r]_y & \tge{0}(y)
      \ar@{}[u];[l]_(.30){}="s"
      \ar@{}[u];[l]_(.70){}="t"
      \ar@2"s";"t"_{\alpha_1}
    }
  \]
  is the $2$-cylinder
  \[
    \xymatrix@C=3pc@R=3pc{
      \sce{0}(x)
      \ar@/^2ex/[r]^(0.70){x}_{}="0"
      \ar@/_2ex/[r]_(0.70){x}_{}="1"
      \ar[d]_{}="f"_{\alpha^-_0}
      \ar@2"0";"1"_{\unit{x}\,\,}
      &
      \tge{0}(x)
      \ar[d]^{\alpha^+_0} \\
      \sce{0}(y)
      \ar@{.>}@/^2ex/[r]^(0.30){y}_{}="0"
      \ar@/_2ex/[r]_(0.30){y}_{}="1"
      \ar@{:>}"0";"1"_{\unit{y}\,\,}
      &
      \tge{0}(y)
      \ar@{}[u];[l]_(.40){}="x"
      \ar@{}[u];[l]_(.60){}="y"
      \ar@<-1.5ex>@/_1ex/@{:>}"x";"y"_(0.60){\alpha_1\,}_{}="0"
      \ar@<1.5ex>@/^1ex/@2"x";"y"^(0.40){\!\alpha_1}_{}="1"
      \ar@{}"1";"0"_(.05){}="z"
      \ar@{}"1";"0"_(.95){}="t"
      \ar@3{>}"z";"t"_{\unit{\alpha_1}}
      \pbox{.}
    }
  \]
\end{paragraph}

\begin{paragraph}[Composition of cylinders]
  Let $\alpha : x \cylto y$ and $\beta : z \cylto t$ be two $n$\nbd-cylinders in
  $C$ for some $n > 0$.
  Let $0 \le i < n$ and suppose that $\alpha$ and $\beta$ are composable in
  dimension $i$, that is, that we have
  \[
    \tge{i}(\alpha) : \tge{i}(x) \cylto \tge{i}(y)
    =
    \sce{i}(\beta) : \sce{i}(x) \cylto \sce{i}(y).
  \]
  We define the \ndef{composition in dimension $i$} of these two $n$-cylinders
  to be the $n$-cylinder
  \[
    \gamma = \alpha \comp{i} \beta
    :
    x \comp{i} z \cylto y \comp{i} t
  \]
  defined by, for $j<i$,
  \begin{align*}
    \gamma^\e_j &= \alpha^\e_j = \beta^\e_j\\
    \gamma^-_i&= \alpha^-_i\\
    \gamma^+_i &= \beta^+_i\\
    \gamma^-_{i+1} & =
    \big(\sce{i+1}(x) \comp{0} \beta_0^+ \comp{1} \cdots \comp{i-1}
    \beta_{i-1}^+ \comp{i} \beta_{i+1}^\e\big)\\
    &\qquad\quad \comp{i+1}
    \big(\alpha^\e_{i+1} \comp{i} \alpha_{i-1}^- \comp{i-1} \cdots \comp{1}\alpha_0^- \comp{0} \sce{i+1}(t)\big)
    \\
    \gamma^+_{i+1} & =
    \big(\tge{i+1}(x) \comp{0} \beta_0^+ \comp{1} \cdots \comp{i-1} \beta_{i-1}^+ \comp{i} \beta_{i+1}^\e\big)
    \\
    &\qquad\quad\comp{i+1}
    \big(\alpha^\e_{i+1} \comp{i} \alpha_{i-1}^- \comp{i-1} \cdots \comp{1} \alpha_0^- \comp{0} \tge{i+1}(t)\big)
  \end{align*}
  if $i + 1 < n$, and
  \[
    \gamma^\e_j = \big(\sce{i+1}(x) \comp{0} \beta_0^+ \comp{1} \cdots \comp{i-1}
      \beta_{i-1}^+ \comp{i} \beta_j^\e\big) \comp{i+1}
      \big(\alpha^\e_j \comp{i} \alpha_{i-1}^- \comp{i-1} \cdots \comp{1}
      \alpha_0^- \comp{0} \tge{i+1}(t)\big)
  \]
  for $i + 1 < j \le n$. In particular, for $j = n$, we have
  \[
    \gamma^\e_n = \big(\sce{i+1}(x) \comp{0} \beta_0^+ \comp{1} \cdots \comp{i-1}
      \beta_{i-1}^+ \comp{i} \beta_n\big) \comp{i+1}
      \big(\alpha_n \comp{i} \alpha_{i-1}^- \comp{i-1} \cdots \comp{1}
      \alpha_0^- \comp{0} \tge{i+1}(t)\big).
  \]
  Here are some pictures of composable cylinders in low dimensions:
  \[
    \xymatrix{
      \cdot \ar[r]^x \ar[d]_{\alpha^-_0} &
      \cdot \ar[d] \ar[r]^z & \cdot \ar[d]^{\beta^+_0} \\
      \cdot \ar[r]_y & \cdot
      \ar[r]_t
      \ar@{}[u];[l]_(.30){}="s"
      \ar@{}[u];[l]_(.70){}="t"
      \ar@2"s";"t"_{\alpha_1}
      & \cdot
      \ar@{}[u];[l]_(.30){}="s'"
      \ar@{}[u];[l]_(.70){}="t'"
      \ar@2"s'";"t'"_{\beta_1}
    }
  \]
  \[
    \xymatrix@C=3pc@R=3pc{
      \cdot
      \ar@/^2ex/[r]_{}="0"
      \ar@/_2ex/[r]_{}="1"
      \ar[d]_{}="f"_{\alpha^-_0}
      \ar@2"0";"1"_{x\,\,}
      &
      \cdot
      \ar@/^2ex/[r]_{}="0"
      \ar@/_2ex/[r]_{}="1"
      \ar@2"0";"1"_{z\,\,}
      \ar[d]
      &
      \cdot
      \ar[d]^{\beta_0^+}
      \\
      \cdot
      \ar@{.>}@/^2ex/[r]_{}="0"
      \ar@/_2ex/[r]_{}="1"
      \ar@{:>}"0";"1"_{y\,\,}
      &
      \cdot
      \ar@{}[u];[l]_(.40){}="x"
      \ar@{}[u];[l]_(.60){}="y"
      \ar@<-1.5ex>@/_1ex/@{:>}"x";"y"_(0.60){\alpha^-_1\!}_{}="0"
      \ar@<1.5ex>@/^1ex/@2"x";"y"^(0.40){\!\alpha^+_1}_{}="1"
      \ar@{}"1";"0"_(.05){}="z"
      \ar@{}"1";"0"_(.95){}="t"
      \ar@3{>}"z";"t"_{\alpha}
      \ar@{.>}@/^2ex/[r]_{}="0"
      \ar@/_2ex/[r]_{}="1"
      \ar@{:>}"0";"1"_{t\,\,}
      &
      \cdot
      \ar@{}[u];[l]_(.40){}="x"
      \ar@{}[u];[l]_(.60){}="y"
      \ar@<-1.5ex>@/_1ex/@{:>}"x";"y"_(0.60){\beta^-_1\!}_{}="0"
      \ar@<1.5ex>@/^1ex/@2"x";"y"^(0.40){\!\beta^+_1}_{}="1"
      \ar@{}"1";"0"_(.05){}="z"
      \ar@{}"1";"0"_(.95){}="t"
      \ar@3{>}"z";"t"_{\beta}
    }
    \qquad
    \qquad
    \xymatrix@C=4pc@R=4pc{
      \cdot
      \ar@/^3ex/[r]_{}="0"
      \ar[r]_{}="01"
      \ar@/_3ex/[r]_{}="1"
      \ar[d]_{}="f"_{\alpha^-_0 = \beta^-_0}
      \ar@2"0";"01"_{x\,\,}
      \ar@2"01";"1"_{z\,\,}
      &
      \cdot
      \ar[d]^{\alpha^+_0 = \beta^+_0} \\
      \cdot
      \ar@{.>}@/^3ex/[r]_{}="0"
      \ar@{.>}[r]_{}="01"
      \ar@/_3ex/[r]_{}="1"
      \ar@{:>}"0";"01"_{y\,\,}
      \ar@{:>}"01";"1"_{t\,\,}
      &
      \cdot
      \ar@{}[u];[l]_(.35){}="x"
      \ar@{}[u];[l]_(.65){}="y"
      \ar@<-3ex>@/_1ex/@{:>}"x";"y"_(0.60){\alpha^-_1\!}_{}="0"
      \ar@<0ex>@{:>}"x";"y"_{}="01"
      \ar@<3ex>@/^1ex/@2"x";"y"^(0.40){\!\beta^+_1}_{}="1"
      \ar@3{>}"01";"0"_{\alpha_2}
      \ar@3{>}"1";"01"_{\beta_2}
    }
  \]
\end{paragraph}

\subsection{The \pdfoo-category of cylinders}
\index{9-category@$\omega$-category!of cylinders}
We define the \oo-category $\CatCyl(C)$ of \ndef{cylinders
in~$C$} to be the
\oo-category whose $n$-cells are $n$-cylinders in $C$ and whose sources,
targets, units, and compositions are given according to the formulas given in
the previous paragraphs.
Tedious calculations show that:

\begin{theorem}\label{thm:catcyl}
  If $C$ is an \oo-category, then $\CatCyl(C)$ is indeed an \oo-category.
\end{theorem}

\begin{proof}
  See \cite[Appendix~A]{metayer2003resolutions}.
\end{proof}

\begin{paragraph}[The free-standing $n$-cylinder]
  \label{sec:free-n-cyl}
Theorem~\ref{thm:catcyl} is a key ingredient in the 
construction $\globe{1}\tensor (-)$ presented
in~\cref{sec:tensor}. In particular, we get for each $n \ge 0$ a polygraph
$\Cyl n=\globe{1}\tensor\globe{n}$. The associated \oo-category
$\free{\Cyl{n}}$ is the free-standing $n$-cylinder, and by construction, for
any \oo-category $C$, the set of $n$-cylinders in $C$ is just
$\Hom{\oCat}{\free{\Cyl{n}}}{C}$. 
\end{paragraph}

\begin{remark}
  The \oo-category $\CatCyl(C)$ can be defined more conceptually using the
  Gray tensor product of \oo-categories. Indeed, the Gray tensor product is
  a biclosed monoidal structure on \oo-category. This means that we have
  functors
  \[
    \HomOpLax : \ooCat^\op \times \ooCat \to \ooCat
    \quad\text{and}\quad
    \HomLax : \ooCat^\op \times \ooCat \to \ooCat
  \]
  and isomorphisms
  \[
      \ooCat(C \otimes D, E)
      \simeq \ooCat(C, \HomOpLax(D, E))
  \]
  and
  \[
      \ooCat(C \otimes D, E)
      \simeq \ooCat(D, \HomLax(C, E)),
  \]
  natural in $C$, $D$, and $E$ in $\ooCat$.
  (See for instance \cite[Appendix A]{AraMaltsiJoint} for definitions.) It
  can be shown that the \oo-category $\Gamma(C)$ is canonically isomorphic
  to~$\HomLax(\globe{1}, C)$ (see \cite[Section B.1]{AraMaltsiJoint}).
\end{remark}

\begin{paragraph}[Functoriality of the cylinder \pdfoo-category]
  Let $f : C \to D$ be an \oo-functor. If $\alpha : x \cylto y$ is an
  $n$-cylinder in $C$, we define an $n$-cylinder $f(\alpha) : f(x)
  \cylto f(y)$ in $D$ by applying $f$ to the components of $\alpha$, that
  is, by the cells
  \[
    f(\alpha_0^-), f(\alpha_0^+), \cdots f(\alpha_n^-), f(\alpha_n^+),
    f(\alpha_{n+1}).
  \]
  One checks that this assignment defines an \oo-functor $\CatCyl(f) :
  \CatCyl(C) \to \CatCyl(D)$. Moreover, this assignment is functorial in
  $f$. In other words, we have defined an endofunctor
  \[
    \CatCyl : \ooCat \to \ooCat.
  \]
\end{paragraph}

\begin{paragraph}[Projections]
  Let $\alpha : x \cylto y$ be a cylinder in $C$. We set
  \[ \projcyltop(\alpha) = x \qqqtand \projcylbot(\alpha) = y. \]
  Geometrically, $\projcyltop(\alpha)$ and $\projcylbot(\alpha)$ are,
  respectively, the top face and the bottom face of the cylinder $\alpha$.
  The formulas defining the structure of the \oo-category $\CatCyl(C)$ make
  transparent the fact that these assignments define \oo-functors
  \[ \projcyltop^{}_C, \projcylbot_C : \CatCyl(C) \to C. \]
  Moreover, these \oo-functors are natural in $C$ and we have defined
  natural transformations
  \[ \projcyltop, \projcylbot : \CatCyl \to \id_{\ooCat}. \]
\end{paragraph}

\begin{paragraph}[Trivial cylinders]
  Let $x$ be an $n$-cell of $C$. The \emph{trivial cylinder} on $x$ is the
  $n$-cylinder $\trivcyl(x) : x \cylto x$ defined by the cells
  \[
    \unit{\sce{0}(x)}, \unit{\tge{0}(x)},
    \dots,
    \unit{\sce{n-1}(x)}, \unit{\tge{n-1}(x)},
    \unit{x}.
  \]
  For instance, the trivial cylinder on a $2$-cell
  \[
    \xymatrix@C=3pc@R=3pc{
      a
      \ar@/^2ex/[r]^{u}_{}="0"
      \ar@/_2ex/[r]_{v}_{}="1"
      \ar@2"0";"1"_{x\,\,}
      &
      b
    }
  \]
  can be pictured as
  \[
    \xymatrix@C=3pc@R=3pc{
      a
      \ar@/^2ex/[r]^(0.70){u}_{}="0"
      \ar@/_2ex/[r]_(0.70){v}_{}="1"
      \ar[d]_{}="f"_{\unit{a}}
      \ar@2"0";"1"_{x\,\,}
      &
      b
      \ar[d]^{\unit{b}} \\
      a
      \ar@{.>}@/^2ex/[r]^(0.30){u}_{}="0"
      \ar@/_2ex/[r]_(0.30){v}_{}="1"
      \ar@{:>}"0";"1"_{x\,\,}
      &
      b
      \ar@{}[u];[l]_(.40){}="x"
      \ar@{}[u];[l]_(.60){}="y"
      \ar@<-1.5ex>@/_1ex/@{:>}"x";"y"_(0.60){\unit{u}\,}_{}="0"
      \ar@<1.5ex>@/^1ex/@2"x";"y"^(0.40){\unit{v}}_{}="1"
      \ar@{}"1";"0"_(.05){}="z"
      \ar@{}"1";"0"_(.95){}="t"
      \ar@3{>}"z";"t"_{\unit{x}}
      \pbox{.}
    }
  \]
  One checks that this assignment defines an \oo-functor $\trivcyl_C :
  C \to \CatCyl(C)$ and this \oo-functor is natural in $C$. In other
  words, we have defined a natural transformation
  \[ \trivcyl : \id_{\ooCat} \to \CatCyl\pbox. \]
\end{paragraph}

\begin{paragraph}[Oplax transformations]
  \label{oplax-transformation}
  \index{oplax!transformation}
  \index{transformation!oplax}
  If $C$ is an \oo-category, then the diagonal \oo-func\-tor $C \to C
  \times C$ factors as
  \[
    \xymatrix@C=3pc{
      C \ar[r]^-{\trivcyl_C} & \CatCyl(C)
      \ar[r]^-{(\projcyltop^{}_C, \projcylbot_C)} & C \times C
      \pbox{.}
    }
  \]
  As every factorization of the diagonal, this factorization induces some
  kind of notion of homotopies. These homotopies will be called
  \ndef{oplax transformations}.
  In other words, if $f, g : C \to D$ are two \oo-functors, an oplax
  transformation
  $\alpha$ from $f$ to $g$, denoted by $\alpha : f \tod g$, is an \oo-functor
  $\alpha : C \to \CatCyl(D)$ making the diagram
  \[
    \xymatrix@R=1pc{
      & & D \\
      C \ar[r]^-{\alpha} \ar@/^3ex/[urr]^f \ar@/_3ex/[drr]_g
      & \CatCyl(D) \ar[ru]_{\projcyltop_D} \ar[rd]^{\projcylbot_D} \\
      & & D
    }
  \]
  commute.
\end{paragraph}

\begin{paragraph}[Algebraic description of oplax transformations]
  \label{paragr:oo-trans_alg}
  Let $\alpha : f \tod g$ be an oplax transformation. If $x$ is an
  $n$-cell of $C$, we will denote by $\alpha_x$ the principal cell of the
  $n$-cylinder~$\alpha(x)$. Thus, $\alpha_x$ is an $(n+1)$-cell
  of $D$. One can show that the cells $\alpha_x$ of $D$, where $x$ varies
  among the cells of~$C$, fully determine~$\alpha$. Better, an
  oplax transformation can be fully defined in terms of these $\alpha_x$.
 
  More precisely, if $f, g : C \to D$ are two \oo-functors, the data of an
  oplax transformation $\alpha : f \tod g$ is equivalent to the data of, for
  every $n$-cell $x$ of $C$, an $(n+1)$-cell
  \[
    \alpha_x :
    f(x) \comp{0} \alpha_{\tge{0}(x)} \comp{1} \cdots \comp{n-1} \alpha_{\tge{n-1}(x)}
    \to
    \alpha_{\sce{n-1}(x)} \comp{n-1} \cdots \comp{1} \alpha_{\sce{0}(x)} \comp{0} g(x)
  \]
  such that the following relations are satisfied:
  \begin{itemize}
    \item for every cell $x$ of $C$, we have
      \[ \alpha_{\unit{x}} = \unit{\alpha_x}, \]
    \item for every $n > i \ge 0$ and every pair of $n$-cells $x$ and $y$
      such that $x \comp{i} y$ is well-defined, we have
      \begin{align*}
        \alpha_{x \comp{i} y} & = \big(f(\sce{i+1}(x)) \comp{0}
        \alpha_{\tge{0}(y)} \comp{1} \cdots \comp{i-1}
        \alpha_{\tge{i-1}(y)} \comp{i} \alpha_y\big) \comp{i+1} \\
        & \phantom{=}\hspace{1ex}
        \big(\alpha_x \comp{i} \alpha_{\sce{i-1}(x)} \comp{i-1} \cdots
        \comp{1} \alpha_{\sce{0}(x)} \comp{0} g(\tge{i+1}(y))\big).        
      \end{align*}
  \end{itemize}
\end{paragraph}

\begin{paragraph}[Some operations on oplax transformations]
  \label{paragr:op_oo-trans}
  Let $f : C \to C$ be an \oo-functor. The \oo-functor $\trivcyl_C : C \to
  \CatCyl(C)$ defines an oplax transformation from~$f$ to~$f$ that we will
  call the \ndef{unit oplax transformation} of $f$ and that we will denote by
  $\unit{f} : f \tod f$. If $x$ is a cell of $C$, we have
  \[ (\unit{f})_x = \unit{x}. \]

  If $\alpha : f \tod f'$ is an oplax transformation between \oo-functors
  from $C$ to $D$
  and $g : D \to E$ is an \oo-functor, we define an oplax
  transformation $g\alpha : gf \tod gf'$ by the composition
  \[
    \xymatrix{C \ar[r]^-\alpha & \CatCyl(D) \ar[r]^-{\CatCyl(g)} & \CatCyl(E)
    \pbox{.}}
  \]
  If $x$ is a cell of $C$, we have
  \[ (g\alpha)_x = g(\alpha_x). \]

  Similarly, if $f : C \to D$ is an \oo-functor and $\alpha : g \tod g'$ is
  an oplax transformation between \oo-functors from $D$ to $E$, we
  define an oplax transformation $\alpha f : gf \tod g'f$ by the composition
  \[
    \xymatrix{C \ar[r]^-f & D \ar[r]^-{\alpha} & \CatCyl(E)}
  \]
  and, for $x$ a cell of $C$, we have
  \[ (\alpha f)_x = \alpha_{f(x)}. \]
\end{paragraph}

\begin{remark}\label{rem:2-Cat_oplax}
  If $f, g, h : C \to D$ are three \oo-functors and $\alpha : f \tod g$ and
  $\beta : g \tod h$ are oplax transformations, one can define in a natural
  way a composite $\beta \circ \alpha : f \tod h$. We will only need the
  existence of this composition and therefore we do not give its precise
  definition. One can show that \oo-categories, \oo-functors, and oplax
  transformations, with the operations defined in the previous paragraph and
  this operation $\circ$, form a sesquicategory (see \cref{sec:sesquicat} for a
  definition). However, this sesquicategory is \emph{not} a $2$-category.
  We refer the reader to \cite[Appendix C]{AraMaltsiJoint} for more details.
\end{remark}

\section{The \pdfoo-Category of Reversible Cylinders}

\begin{paragraph}[Reversible cylinders]
  \index{reversible!cylinder}
  \index{cylinder!reversible}
  Let $C$ be an \oo-category. An $n$-cylinder $\alpha : x \cylto y$ is said to
  be \emph{reversible} if the cells
  \[
    \alpha_0^-, \alpha_0^+, \dots, \alpha_{n-1}^-, \alpha_{n-1}^+,
    \alpha_n
  \]
  are reversible (see~\cref{paragr:rev_cells}).
\end{paragraph}

\begin{paragraph}[The \pdfoo-category of reversible cylinders]
  \label{paragr:def_catcylrev}
  \index{9-category@$\omega$-category!of cylinders!reversible}
  Let $C$ be an \oo-category and let $\CatCylRev(C)_n$, for $n \ge 0$, be the set of
  reversible $n$-cylinders in $C$. By definition, we have an inclusion
  $\CatCylRev(C)_n \subset \CatCyl(C)_n$. One immediately checks that these
  $\CatCylRev(C)_n$ define a subcategory $\CatCylRev(C)$ of $\CatCyl(C)$.
  Moreover, as the image of a reversible cell by an \oo-functor is
  reversible, the functor $\CatCyl : \ooCat \to \ooCat$ induces a functor
  \[
    \CatCylRev : \ooCat \to \ooCat.
  \]
  Since units are reversible cells, the trivial cylinder on a cell is
  reversible and the natural transformation $\trivcyl : \id_{\ooCat} \to
  \CatCyl$ factors by a natural transformation
  \[ \trivcyl : \id_{\ooCat} \to \CatCylRev. \]
  Finally, by restriction, the natural transformations
  $\projcyltop, \projcylbot : \CatCyl \to \id_{\ooCat}$
  define natural transformations
  \[ \projcyltop, \projcylbot : \CatCylRev \to \id_{\ooCat}. \]
\end{paragraph}

\begin{paragraph}[Reversible oplax transformations]
  \index{oplax!transformation!reversible}
  \index{reversible!oplax transformation}
  Let $f, g : C \to D$ be two \oo-func\-tors. We say that an
  oplax transformation $\alpha : f \tod g$ is
  \ndef{reversible} if the
  \oo-functor $\alpha : C \to \CatCyl(D)$ factors through the inclusion
  $\CatCylRev(D) \subset \CatCyl(D)$. In other words, we have a
  factorization
  \[
    \xymatrix@C=3pc{
      D \ar[r]^-{\trivcyl_D} & \CatCylRev(D)
      \ar[r]^-{(\projcyltop^{}_D, \projcylbot_D)} & D \times D
    }
  \]
  of the diagonal of $D$ and a reversible oplax transformation from $f$ to $g$
  is an \oo-functor $\alpha : C \to \CatCylRev(D)$ making the diagram
  \[
    \xymatrix@R=1pc{
      & & D \\
      C \ar[r]^-{\alpha} \ar@/^3ex/[urr]^f \ar@/_3ex/[drr]_g
      & \CatCylRev(D) \ar[ru]_{\projcyltop_D} \ar[rd]^{\projcylbot_D} \\
      & & D
    }
  \]
  commute. 

  In the algebraic definition of an oplax transformation given
  in~\cref{paragr:oo-trans_alg}, an oplax transformation $\alpha$ is
  reversible if and only if $\alpha_x$ is reversible for every cell~$x$.

  One immediately checks that the various operations on oplax transformations
  defined in~\cref{paragr:op_oo-trans} restrict to reversible
  oplax transformations. In other words, the unit oplax transformation of an
  \oo-functor is reversible and if we can compose horizontally a reversible
  oplax transformation and an \oo-functor (in both directions) we get a
  reversible oplax transformation.
\end{paragraph}

\medbreak

We now introduce some terminology to state the ``transport lemma''
which is of crucial importance to prove the existence of the folk model
structure.

\begin{paragraph}[Incomplete cylinders]
  Let $C$ be an \oo-category. For $n \ge 0$, a \ndef{bottom-incomplete
  $n$-cylinder}\index{bottom-incomplete
  $n$-cylinder}\index{cylinder!bottom-incomplete} in $C$ consists of the
  same data as an $n$-cylinder $\alpha : x \cylto y$ except that $y$ and
  $\alpha_n$ are not given.
  More formally, it consists of an object $x$ of~$C$ if $n = 0$ and,
  if $n > 0$, of an $n$-cell $x$, two parallel $(n-1)$-cells $y^-$
  and $y^+$ and a sequence
  \[
    \alpha_0^-, \alpha_0^+, \dots, \alpha_{n-1}^-, \alpha_{n-1}^+
  \]
  where $\alpha_j^\e$ is a $(j+1)$-cell whose source and target are
  \begin{align*}
    \alpha_j^- & :
    \sce{j}(x) \comp{0} \alpha_0^+ \comp{1} \cdots \comp{j-1} \alpha_{j-1}^+
    \to
    \alpha_{j-1}^- \comp{j-1} \cdots \comp{1} \alpha_0^- \comp{0}
    \sce{j}(y^-),
    \\
    \alpha_j^+ & :
    \tge{j}(x) \comp{0} \alpha_0^+ \comp{1} \cdots \comp{j-1} \alpha_{j-1}^+
    \to
    \alpha_{j-1}^- \comp{j-1} \cdots \comp{1} \alpha_0^- \comp{0}
    \tge{j}(y^+).
  \end{align*}
  Here are pictures of bottom-incomplete cylinders in dimensions $0$, $1$, and $2$:
  \[
   \xymatrix@R=3pc{x\ \cdot\ar@{}[d] \\
   &
   }
   \qquad
   \qquad
    \xymatrix@C=3pc@R=3pc{
      \sce{0}(x) \ar[r]^x \ar[d]_{\alpha^-_0} &
      \tge{0}(x) \ar[d]^{\alpha^+_0} \\
      y^-
      & y^+
      \ar@{}[u];[l]_(.30){}="s"
      \ar@{}[u];[l]_(.70){}="t"
    }
   \qquad
   \qquad
    \xymatrix@C=3.5pc@R=3pc{
      \cdot
      \ar@/^2ex/[r]_{}="0"
      \ar@/_2ex/[r]_{}="1"
      \ar[d]_{}="f"_{\alpha^-_0}
      \ar@2"0";"1"_{x\,\,}
      &
      \cdot
      \ar[d]^{\alpha^+_0} \\
      \cdot
      \ar@{.>}@/^2ex/[r]^(0.30){y^-}_{}="0"
      \ar@/_2ex/[r]_(0.30){y^+}_{}="1"
      &
      \cdot
      \ar@{}[u];[l]_(.40){}="x"
      \ar@{}[u];[l]_(.60){}="y"
      \ar@<-1.5ex>@/_1ex/@{:>}"x";"y"_(0.60){\alpha^-_1\,}_{}="0"
      \ar@<1.5ex>@/^1ex/@2"x";"y"^(0.40){\!\alpha^+_1}_{}="1"
      \ar@{}"1";"0"_(.05){}="z"
      \ar@{}"1";"0"_(.95){}="t"
    }
  \]

  We say that such a bottom-incomplete $n$-cylinder is
  \ndef{reversible} if the cells
  \[
    \alpha_0^-, \alpha_0^+, \dots, \alpha_{n-2}^-, \alpha_{n-2}^+,
    \alpha_{n-1}^-
  \]
  are reversible.

  \index{top-incomplete $n$-cylinder}
  \index{cylinder!top-incomplete}
  We define similarly the notion of a \ndef{top-incomplete $n$-cylinder} and of
  \ndef{reversible top-incomplete $n$-cylinder}.
\end{paragraph}

\begin{lemma}[Transport lemma]
  \index{transport lemma}
  \index{lemma!transport}
  \label{lemma:transport}
  Let $C$ be an \oo-category.
  \begin{enumerate}
    \item\label{item:bi_ext}
    Any reversible bottom-incomplete $n$-cylinder extends (non-uniquely) to a
    reversible $n$-cylinder.
  \item Consider a reversible bottom-incomplete $n$-cylinder.
    \begin{enumerate}
      \item If $\alpha :x \cylto y$ and $\alpha' : x \cylto y'$ are two
        extensions as in~\ref{item:bi_ext}, then $y$ and $y'$ are
        $\omega$-equivalent.
      \item If $\alpha : x \cylto y$ is an extension as in~\ref{item:bi_ext}
        and $y'$ is an $n$-cell \oo-equivalent to~$y$, then there exists an
        extension $\alpha : x \cylto y'$ as in~\ref{item:bi_ext}.
    \end{enumerate}
  \end{enumerate}
\end{lemma}

\begin{proof}
  The case $n = 0$ is obvious. For $n=1$, consider the incomplete
  $1$-cylinder defined by  $(x,y^-,y^+,\alpha_0^-,\alpha_0^+)$ and let 
 $\weakinv{\alpha_0^-}$ be a weak inverse of $\alpha_0^-$. The $1$-cell
 \[y=\weakinv{\alpha_0^-}\comp 0 x\comp 0\alpha_0^+\]
satisfies the relation $x\comp 0\alpha_0^+\oeq\alpha_0^-\comp 0 y $ where
$\oeq$ denotes \oo-equivalence as defined in~\cref{paragr:oeq-cells}. This is witnessed by a reversible $2$-cell
$\alpha_1:x\comp 0\alpha_0^+\to\alpha_0^-\comp 0 y$, which gives the desired
extension. Conditions 2(a) and 2(b) immediately follow from the remark
that $y\oeq y'$ if and only if $\alpha_0^-\comp 0 y \sim \alpha_0^-\comp 0 y'$.
  
 Let us now suppose $n > 1$ and let
  \[ (x, y^-, y^+, \alpha_0^-, \alpha_0^+, \dots, \alpha_{n-2}^-, \alpha_{n-2}^+,
  \alpha_{n-1}^-) \]
be a reversible bottom-incomplete $n$-cylinder. We have to prove that there
is a weakly unique $n$-cell $y$ such that
  \begin{equation}
    \label{eq:incomplete_cyl}
    x \comp{0} \alpha_0^+ \comp{1} \cdots \comp{n-1} \alpha_{n-1}^+
    \quad\oeq\quad
    \alpha_{n-1}^- \comp{n-1} \cdots \comp{1} \alpha_0^- \comp{0} y,
  \end{equation}
  where $y : y^- \to y^+$.
  Let $w=x \comp{0} \alpha_0^+ \comp{1} \cdots \comp{n-1}
  \alpha_{n-1}^+$ denote the left member
  of~(\ref{eq:incomplete_cyl}). We prove by induction on
  $k$ such that $1 \le k \le n$ that there is a weakly unique $n$-cell $v^k$ in $C$
  such that
  \begin{equation}
    \label{eq:incomplete_cyl_k}
    w\quad \oeq\quad \alpha_{n-1}^- \comp{n-1} \cdots \comp{1}
    \alpha_{n-k}^- \comp{n-k} v^k.
  \end{equation}
  For $k=1$, this equation reads
  \begin{equation}
    \label{eq:incomplete_cyl_1}
    w\quad\oeq\quad
  \alpha_{n-1}^-\comp{n-1}v^1.
  \end{equation}
  It can be seen as an equation on $1$-cells in the \oo-category $C'$ of
  $(n-2)$\nbd-cells of~$C$. Therefore the division lemma
  (Lemma~\ref{lemma:weakdiv}) applies and yields a
  weakly unique $1$-cell $v^1$ of $C'$, which is also an $n$-cell of $C$
  satisfying~(\ref{eq:incomplete_cyl_1}).
  Let~$k$ such that $1 \le k < n$ and suppose
 that~(\ref{eq:incomplete_cyl_k}) has a weakly unique solution~$v^k$. By applying again
 the division lemma, there is a weakly unique $n$-cell
  $v^{k+1}$ solution of
  \begin{equation}
    \label{eq:inducstep}
     v^k\quad\oeq\quad \alpha_{n-k-1}^-\comp{n-k-1} v^{k+1}.
  \end{equation}
  Therefore, substituting $v^k$ by the second member
  of~(\ref{eq:inducstep}) in~(\ref{eq:incomplete_cyl_k}) one gets
  \[
     w\quad \oeq\quad \alpha_{n-1}^- \comp{n-1} \cdots \comp{1}
     \alpha_{n-k-1}^- \comp{n-k-1} v^{k+1},
  \]
  which completes the induction. Now $v^n=y$ provides a weakly unique
  solution of~(\ref{eq:incomplete_cyl}), as required.
\end{proof}

\begin{remark}
  A similar result holds for reversible top-incomplete cylinders.
\end{remark}

\begin{paragraph}[Alternative description of the transport lemma]
  In practice, the transport lemma will be used in the following way. Let
  $\alpha : x_0 \cylto x_1$, $\beta : y_0 \cylto y_1$ be two parallel
  $n$-cylinders and let $u : x_0 \to y_0$ be an $(n+1)$-cell in an
  \oo-category~$C$. The data of $\alpha$, $\beta$, and $u$ is equivalent to
  the data of a bottom-incomplete $(n+1)$-cylinder (informally, $\alpha$,
  $\beta$, and $u$ corresponds respectively to the back face, the front face,
  and the top face of the bottom-incomplete cylinder). Moreover, this
  bottom-incomplete cylinder is reversible if and only if the cylinders
  $\alpha$ and $\beta$ are reversible. In this case, an extension, as in the
  transport lemma, corresponds to an $(n+1)$-cell $v : x_1 \to y_1$ and a
  reversible $(n+1)$-cylinder $\Lambda : u \cylto v$, whose source and
  target cylinders are $\alpha$ and $\beta$. We will often write $\Lambda :
  u \cylto v : \alpha \to \beta$ in such a situation. The existence and
  uniqueness properties of the second part of the transport lemma apply to
  $v$.
\end{paragraph}

\begin{proposition}\label{prop:projcyl_triv_fib}
  Let $C$ be an \oo-category.
  \begin{enumerate}
    \item The projections $\projcyltop^{}_C, \projcylbot_C : \CatCylRev(C)
      \to C$ are trivial fibrations.
    \item The \oo-functor $\trivcyl_C : C \to \CatCylRev(C)$ is an
      \oo-equivalence.
  \end{enumerate}
\end{proposition}

\begin{proof}
  \begin{enumerate}
    \item Let us prove the result for $\projcyltop$, the other proof
      being analogous. The surjectivity on objects is clear. Let $n \ge 0$ and
      let $\alpha : x \cylto y$ and $\beta : z \cylto t$ be two parallel
      reversible $n$-cylinders. Suppose we have an $n$-cell between their
      images by the \oo-functor $\projcyltop$, i.e., an $n$-cell $u : x \to
      z$. We are precisely in the situation of the previous paragraph and
      the transport lemma gives a reversible $(n+1)$-cylinder $\Lambda :
      \alpha \to \beta : u \cylto v$, showing that $\projcyltop$ is indeed a
      trivial fibration.
    \item Since the composite
      \[
        \xymatrix{C \ar[r]^-{\trivcyl_C} & \CatCylRev(C)
          \ar[r]^-{\projcyltop_C} & C
        }
      \]
      is the identity and hence an \oo-equivalence, and that the trivial
      fibration~$\projcyltop_C$ is an \oo-equivalence
      (see~\cref{prop:triv_fib_eq}), it follows from one
      of the known cases of the 2-out-of-3 property
      (\cref{prop:2-out-of-3-easy}) that $\trivcyl_C$ is an \oo-equivalence.
      \qedhere
  \end{enumerate}
\end{proof}

\begin{paragraph}[Mapping path space factorization]
  \label{paragr:mapcyl}
  Let $f : C \to D$ be an \oo-functor. We define an \oo-category
  $\MapCyl{f}$ by the pullback square
  \[
    \xymatrix{
      \MapCyl{f} \ar[d]_{\leftproj} \ar[r]^-{\rightproj} & \CatCylRev(D)
      \ar[d]^{\projcyltop} \\
      C \ar[r]_f & D \pbox{.}
    }
  \]
  Explicitly, an $n$-cell of $\MapCyl{f}$ is a pair $(x, \gamma : f(x)
  \cylto y)$, where $x$ and $y$ are $n$-cells of $C$ and $D$, and $\gamma$
  is a reversible $n$-cylinder. Note that since trivial fibrations are
  stable under pullback, by the previous proposition, the projection
  $\leftproj : \MapCyl{f} \to C$ is a trivial fibration.

  The \oo-functor $f : C \to D$ factors as
  \[
    \xymatrix{
      C \ar[r]^{\incmapcyl{f}} & \MapCyl{f} \ar[r]^{\projmapcyl{f}} & D\pbox,
    }
  \]
  where $\incmapcyl{f}$ is defined by
  \[ \incmapcyl{f}(x) = (x, \iota(f(x)) : f(x) \cylto f(x)) \]
  and $\projmapcyl{f}$ by
  \[ \projmapcyl{f}(x, \alpha : f(x) \cylto y) = y. \]
  As the triangle
  \[
    \xymatrix{
      C \ar[dr]_{\id_C} \ar[r]^{\incmapcyl{f}} & \MapCyl{f}
      \ar[d]^{\leftproj} \\
                                               & C
    }
  \]
  is commutative and $\leftproj$ is a trivial fibration, by one of the known
  cases of the 2-out-of-3 property (\cref{prop:2-out-of-3-easy}),
  $\incmapcyl{f}$ is an \oo-equivalence.
\end{paragraph}

\begin{proposition}\label{prop:w-eq_mapcyl}
  An \oo-functor $f : C \to D$ is an \oo-equivalence if and only if 
  $\projmapcyl{f} : \MapCyl{f} \to D$ is a trivial fibration.
\end{proposition}

\begin{proof}
  If $\projmapcyl{f}$ is a trivial fibration and hence an \oo-equivalence, then, since
  $\incmapcyl{f}$ is always an \oo-equivalence, so is
  $f = \projmapcyl{f}\incmapcyl{f}$.

  Suppose conversely that $f$ is an \oo-equivalence.
  \begin{enumerate}
    \item Let $y$ be an object of $D$. As $f$ is an \oo-equivalence, there
      exists an object $x$ of $C$ and a reversible $1$-cell $u : f(x) \to y$.
      This $1$-cell defines a reversible $0$-cylinder $u : f(x) \cylto y$ sent
      to $y$ by $\projmapcyl{f}$.
    \item Let $n > 0$ and let $(x, \alpha : f(x) \cylto y)$ and $(x',
      \alpha' : f(x') \cylto y')$ be two parallel $n$-cells of $\MapCyl{f}$.
      Suppose that we have a cell between their image by~$\projmapcyl{f}$,
      i.e., an $n$-cell $v : y \to y'$. By the transport lemma
      (Lemma~\ref{lemma:transport}), there exist an $(n+1)$-cell $v' : f(x)
      \to f(x')$ and an $(n+1)$-cylinder
      $\Lambda : v' \cylto v : \alpha \to \alpha'$. Since $f$ is an
      \oo-equivalence, there exists $u : x \to x'$ such that $f(u)$ is
      \oo-equivalent to $v'$. Using the transport lemma again, we get
      a reversible $(n+1)$-cylinder $\Lambda' : f(u)
      \cylto v : \alpha \to \alpha'$. This means that $(u, \Lambda' : f(u)
      \cylto v)$ is an $(n+1)$-cell of $\MapCyl{f}$ sent to $v$ by
      $\projmapcyl{f}$, thereby proving the result.
      \qedhere
  \end{enumerate}
\end{proof}

\begin{theorem}\label{thm:w-eq-2-3}
  The class of \oo-equivalences satisfies the 2-out-of-3 property.
\end{theorem}

\begin{proof}
  The only remaining case is the following one (see
  Propositions~\ref{prop:w_eq_comp} and~\ref{prop:2-out-of-3-easy} for the
  other ones). Let $f : C \to D$ and $g : D \to E$ be two \oo-functors.
  Suppose $f$ and $gf$ are \oo-equivalences. Let us prove that $g$ is an
  \oo-equivalence.
  \begin{enumerate}
    \item Let $z$ be an object of $E$. Since $gf$ is an \oo-equivalence, we
      get an object $x$ of $C$ such that $gf(x)$ and $z$ are \oo-equivalent
      and $y = f(x)$ shows that $g$ is surjective up to \oo-equivalence.
    \item
  Let $y_0$ and $y_1$ be two parallel $n$-cells of $D$
  and let $w : g(y_0) \to g(y_1)$ be an $(n+1)$-cell. By the previous
  proposition, since $f$ is an \oo-equivalence, $\projmapcyl{f} : \MapCyl{f}
  \to D$ is a trivial fibration. We can thus lift the pair of parallel
  cells~$y_0, y_1$ to a pair of parallel cells
  \[
    (x_0, \alpha : f(x_0) \cylto y_0)
    \qqtand
    (x_1, \beta : f(x_1) \cylto y_1)
  \]
  of $\MapCyl{f}$. By applying $g$, we get a pair of parallel cells
  \[
    (x_0, g(\alpha) : gf(x_0) \cylto g(y_0))
    \qqtand
    (x_1, g(\beta) : gf(x_1) \cylto g(y_1))
  \]
  of $\MapCyl{gf}$. Since $gf$ is an
  \oo-equivalence, by the previous proposition, the projection
  $\projmapcyl{gf} : \MapCyl{gf} \to E$ is a trivial fibration. We can thus
  lift the cell $w : g(y_0) \to g(y_1)$ to a cell 
  \[
    (u : x_0 \to x_1, \Lambda : gf(u) \cylto w : g(\alpha) \to g(\beta))
  \]
  of $\MapCyl{gf}$.
  Consider the bottom-incomplete cylinder defined by $\alpha : f(x_0) \cylto
  y_0$, $\beta : f(x_1) \cylto
  y_1$ and $f(u) : f(x_0) \to f(x_1)$. By the transport lemma
  (Lemma~\ref{lemma:transport}), it can be extended and we get
  a cell $v : y_0 \to y_1$ and a cylinder $\Delta : f(u) \cylto v : \alpha
  \to \beta$. Applying the uniqueness part of the transport lemma to the cylinders
  \[
    g(\Delta) : gf(u) \cylto g(v) : g(\alpha) \to g(\beta)
  \]
  and
  \[
    \Lambda : gf(u) \cylto w : g(\alpha) \to g(\beta),
  \]
  we get that $g(v)$ is \oo-equivalent to $w$, thereby proving the result.
  \qedhere
  \end{enumerate}
\end{proof}

\section{Coherent Reversible Cells and Fibrations}

\begin{paragraph}[The free-standing reversible cell]
  \label{paragr:def_R1}
  We define an \oo-category $\Rev{1}$ generated by a polygraph in the
  following way. The \oo-category $\Rev{1}$ has two objects, $0$ and $1$.
  There are two generating $1$-cells
  \[
    r : 0 \to 1
    \qqtand
    \weakinv{r} : 1 \to 0
  \]
  and four generating $2$-cells
  \[
    r_- : \unit{0} \to r \comp{0} \weakinv{r}, \quad
    \weakinv{r}_- : r \comp{0} \weakinv{r} \to \unit{0}, \quad
    r_+ : \weakinv{r} \comp{0} r \to \unit{1}, \quad
    \weakinv{r}_+ : \unit{1} \to \weakinv{r} \comp{0} r.
  \]
  More generally, for $j \ge 2$, there are $2^j$ generating $j$-cells
  \[ r_{l_1, \dots, l_{j-1}} \qqtand \weakinv{r}_{l_1, \dots, l_{j-1}}, \]
  where for $1 \le k \le j - 1$, $l_k = \pm$. The source and target of the
  generators are given by
  \begin{align*}
    r_{l_1, \dots, l_{j-2}, -}
    &:&
    \unit{\sce{j-2}(r_{l_1, \dots, l_{j-2}})}
    &\to
    r_{l_1, \dots, l_{j-2}} \comp{j-2} \weakinv{r}_{l_1, \dots, l_{j-2}},
    \\
    \weakinv{r}_{l_1, \dots, l_{j-2}, -}
    &:&
    r_{l_1, \dots, l_{j-2}} \comp{j-2} \weakinv{r}_{l_1, \dots, l_{j-2}}
    &\to
    \unit{\sce{j-2}(r_{l_1, \dots, l_{j-2}})},
    \\
    r_{l_1, \dots, l_{j-2}, +}
    &:&
    \weakinv{r}_{l_1, \dots, l_{j-2}} \comp{j-2} r_{l_1, \dots, l_{j-2}}
    &\to
    \unit{\tge{j-2}(r_{l_1, \dots, l_{j-2}})},
    \\
    \weakinv{r}_{l_1, \dots, l_{j-2}, +}
    &:&
    \unit{\tge{j-2}(r_{l_1, \dots, l_{j-2}})}
    &\to
    \weakinv{r}_{l_1, \dots, l_{j-2}} \comp{j-2} r_{l_1, \dots, l_{j-2}}.
  \end{align*}

  This definition was made so that $\Rev{1}$ is, in some sense, freely
  generated by a reversible $1$-cell $r$.
  In particular, $\Rev{1}$ comes with a \ndef{canonical inclusion}
  $\globe{1} \into \Rev{1}$ corresponding to the $1$-cell $r$ and, if $C$ is an
  \oo-category, a $1$-cell $u$ of $C$ is reversible if and only if there
  exists a dotted arrow making the triangle
  \[
    \xymatrix{
      \globe{1} \ar@{^{(}->}[d] \ar[r]^u & C \\
      \Rev{1} \ar@{.>}[ur]
    }
  \]
  commute. In particular, two objects $x$ and $y$ of $C$ are \oo-equivalent
  if and only if there exists a dotted arrow making the triangle
  \[
    \xymatrix{
      \sphere{1} \ar@{^{(}->}[d] \ar[r]^{(x, y)} & C \\
      \Rev{1} \ar@{.>}[ur] & \pbox{,}
    }
  \]
  where the vertical inclusion is the composite $\sphere{1} \into \globe{1}
  \into \Rev{1}$, commute.

  Similarly, for every $n \ge 1$, one can define an \oo-category $\Rev{n}$,
  endowed with a \ndef{canonical inclusion} $\globe{n} \into \Rev{n}$
  corresponding to an $n$-cell $r$, modeled on the definition of a
  reversible $n$-cell. This \oo-category has, by definition, the following
  properties: an $n$-cell $u$ of an \oo-category $C$ is reversible if and
  only if there exists a dotted arrow making the triangle
  \[
    \xymatrix{
      \globe{n} \ar@{^{(}->}[d] \ar[r]^u & C \\
      \Rev{n} \ar@{.>}[ur]
    }
  \]
  commute; in particular, two parallel $(n-1)$-cells $x$ and $y$ of $C$ are
  \oo-equivalent if and only if there exists a dotted arrow making the
  triangle
  \[
    \xymatrix{
      \sphere{n} \ar@{^{(}->}[d] \ar[r]^{(x, y)} & C \\
      \Rev{n} \ar@{.>}[ur]
    }
  \]
  commute.
\end{paragraph}

\begin{paragraph}[Incoherence]
  If $x$ and $y$ are two objects of a $2$-category, it is well known that $x$
  and $y$ are equivalent if and only if there exists an adjoint equivalence
  between $x$ and $y$. In other words, if there exist $1$-cells $u : x \to y$, $v
  : y \to x$ and $2$-cells $\eta : \unit{x} \to uv$ and $\e : vu \to \unit{y}$
  that are isomorphisms, then we can choose these cells such that the
   triangular identities
  \[
    (\eta \comp{} u)(u \comp{} \e) = \unit{u}
    \qtand
    (v \comp{} \eta)(\e \comp{} v) = \unit{v}
  \]
  hold. A $4$-tuple $(u, v, \eta, \e)$ as above is a witness of the fact that
  $x$ and $y$ are \oo-equivalent. We can think of this witness as
  \ndef{incoherent} if the triangular identities do not hold and as
  \ndef{coherent} if they hold.

  Similarly, if $x$ and $y$ are two objects of an \oo-category $C$, the data
  of an \oo-functor $\Rev{1} \to C$ making the triangle
  \[
    \xymatrix{
      \sphere{1} \ar@{^{(}->}[d] \ar[r]^{(x, y)} & C \\
      \Rev{1} \ar[ur]
    }
  \]
  commute is witness of the fact that $x$ and $y$ are \oo-equivalent.
  This witness is \ndef{incoherent} as, for instance, it does not satisfy the
  triangular identities up to \oo-equivalence. Technically, this boils down
  to the fact that $R_1$ is not contractible (see~\cref{paragr:homology_R1}).
\end{paragraph}

\begin{paragraph}[Coherence]\label{paragr:coherence}
  \label{paragr:def_globeinv}
  Let $n \ge 0$. Consider the \oo-functor
  \[ \sphere{n+1} \to \globe{n} \]
  corresponding to the collapsing of the two non-trivial $n$-cells of
  $\sphere{n+1}$. By Proposition~\ref{prop:factor_cof_fib_triv}, this
  \oo-functor factors (non-uniquely) as a cofibration followed by a trivial
  fibration. Let us fix such a factorization
  \[
    \xymatrix{
      \sphere{n+1} \ar[r]^j & \globeinv{n+1} \ar[r]^p & \globe{n} \pbox.
    }
  \]

  Let $x$ and $y$ be two parallel $n$-cells of an \oo-category $C$. The fact
  that $p$ is a trivial fibration makes it reasonable to think of a dotted
  arrow making the triangle
      \[
        \xymatrix{
          \sphere{n+1} \ar[d]_j \ar[r]^{(x, y)} & C \\
          \globeinv{n+1} \ar@{.>}[ur]
        }
      \]
  commute as a \ndef{coherent} witness of the fact that $x$ and $y$ are
  \oo-equivalent.
\end{paragraph}

\begin{remark}\label{rem:choice_Jn}
  As already mentioned, the \oo-category $\globeinv{n+1}$ introduced in the
  previous paragraph is not uniquely defined. Some choices are wiser than
  others. For instance, we can choose $\globeinv{n+1}$ so that its underlying
  $(n+1)$-category is obtained from $\globe{n+1}$ by adding an $(n+1)$-cell
  $v : y \to x$ in the other direction than the principal cell $u : x \to y$
  of $\globe{n+1}$. We can even fix its underlying $(n+2)$\nbd-categ\-ory by
  saying that it is generated by two $(n+2)$-cells $\eta : \unit{x} \to u
  \comp{n} v$ and $\e : v \comp{n} u \to \unit{y}$.
  We do not assume that these additional properties hold in the remaining of
  the text. (We will use the existence of these better choices in
  \cref{sec:folk_trans}.)
\end{remark}

The following proposition shows that, as in the case of $2$-categories,
there exist coherent witnesses if and only if there exist incoherent
witnesses.

\begin{proposition}
  Let $n \ge 0$ and let $x$ and $y$ be two parallel $n$-cells of an
  \oo-category~$C$. The following assertions are equivalent:
  \begin{enumerate}
    \item\label{item:coh_a} The cells $x$ and $y$ are \oo-equivalent.
    \item\label{item:coh_b} There exists a dotted arrow making the triangle
      \[
        \xymatrix{
          \sphere{n+1} \ar@{^{(}->}[d] \ar[r]^{(x, y)} & C \\
          \Rev{n+1} \ar@{.>}[ur]
        }
      \]
      commute.
    \item\label{item:coh_c} There exists a dotted arrow making the triangle
      \[
        \xymatrix{
          \sphere{n+1} \ar[d]_j \ar[r]^{(x, y)} & C \\
          \globeinv{n+1} \ar@{.>}[ur]
        }
      \]
      commute.
    \item\label{item:coh_d} There exists a factorization
      \[
        \xymatrix{
          \sphere{n+1} \ar[r]^k & K \ar[r]^q & \globe{n}
        }
      \]
      of the \oo-functor $\sphere{n+1} \to \globe{n}$
      of~\cref{paragr:coherence} with $q$ a trivial fibration, and a
      dotted arrow making the triangle
      \[
        \xymatrix{
          \sphere{n+1} \ar[d]_k \ar[r]^{(x, y)} & C \\
          K \ar@{.>}[ur]
        }
      \]
      commute.
  \end{enumerate}
\end{proposition}

\begin{proof}
  The equivalence between \ref{item:coh_a} and \ref{item:coh_b} is true by the
  definition of $\Rev{n+1}$.

  By definition, \ref{item:coh_c} implies \ref{item:coh_d}. The converse
  follows by considering the solid commutative square
  \[
    \xymatrix{
      \sphere{n+1} \ar[d]_j \ar[r]^k & K \ar[d]^q \\
      \globeinv{n+1} \ar[r]_p \ar@{.>}[ur] & \globe{n}
    }
  \]
  that admits a lift as $j$ is a cofibration and $q$ is a trivial fibration.

  Let us prove that \ref{item:coh_d} implies \ref{item:coh_a}. The
  \oo-functor $k : \sphere{n+1} \to K$ defines a pair of parallel $n$-cells
  in $K$. These two cells are sent to the principal cell of $\sphere{n}$ by
  $q : K \to \globe{n}$. Since this \oo-functor $q$ is a trivial fibration,
  by Proposition~\ref{prop:oo-equiv_inj}, there exists a reversible
  $(n+1)$-cell $r$ in $K$ between these two parallel $n$-cells. By
  commutativity of the triangle of the hypothesis, we have $f(r) : x \to y$,
  showing that $x$ and $y$ are \oo-equivalent.

  To conclude the proof, let us prove that \ref{item:coh_a} implies
  \ref{item:coh_d}. Consider the \oo-functor $x : \globe{n} \to C$
  corresponding to the $n$-cell $x$ and the associated \oo-category
  $\MapCyl{x}$ (see~\cref{paragr:mapcyl}). We get a factorization

  \[
    \xymatrix@C=2.5pc{
      \sphere{n+1} \ar[r]^-{(k_1, k_2)} & \MapCyl{x} \ar[r]^\leftproj & \globe{n}
    }
  \]
  of the \oo-functor $\sphere{n+1} \to \globe{n}$ of the assertion by choosing two parallel
  $n$\nbd-cells~$k_1$ and $k_2$ of $\MapCyl{x}$ in the following way. Denote by
  $p$ the principal cell of~$\globe{n}$.  We set $k_1 = (p, \trivcyl(x) : x
  \cylto x)$. As for the cell $k_2$, we choose a reversible
  $(n+1)$-cell $u : x \to y$. Using this $u$, one easily defines a
  reversible cylinder $\alpha_u : x \cylto y$ and we set $k_2 = (p, \alpha_u
  : x \cylto y$). This ends the definition of the desired factorization.
  Moreover, by definition, the triangle
  \[
    \xymatrix{
      \sphere{n+1} \ar[d]_{(k_1, k_2)} \ar[r]^{(x, y)} & C \\
      \MapCyl{x} \ar[ur]_{\projmapcyl{x}}
    }
  \]
  commutes, thereby proving the result.
\end{proof}

\begin{proposition}\label{prop:char_w-eq}
  Let $f : C \to D$ be an \oo-functor. The following assertions are
  equivalent:
  \begin{enumerate}
    \item The \oo-functor $f$ is an \oo-equivalence.
    \item For every $n \ge 0$, every commutative square
      \[
        \xymatrix{
          \sphere{n} \ar[d]_{\gencof{n}} \ar[r] & C \ar[d]^f \\
          \globe{n} \ar[r] & D
        }
      \]
      factors as
      \[
        \xymatrix{
          \sphere{n} \ar[d]_{\gencof{n}} \ar[r]^{\gencof{n}} & \globe{n}
          \ar[r] \ar[d] & C \ar[d]^f \\
          \globe{n} \ar[r] & \Rev{n+1} \ar[r] & D \pbox{,}
        }
      \]
      where the two \oo-functors $\globe{n} \to \Rev{n}$ are the two
      components of the inclusion $\sphere{n+1} \into \globe{n+1} \into
      \Rev{n+1}$.
    \item For every $n \ge 0$, every commutative square
      \[
        \xymatrix{
          \sphere{n} \ar[d]_{\gencof{n}} \ar[r] & C \ar[d]^f \\
          \globe{n} \ar[r] & D
        }
      \]
      factors as
      \[
        \xymatrix{
          \sphere{n} \ar[d]_{\gencof{n}} \ar[r]^{\gencof{n}} & \globe{n}
          \ar[r] \ar[d] & C \ar[d]^f \\
          \globe{n} \ar[r] & \globeinv{n+1} \ar[r] & D
        }
      \]
      where the two \oo-functors $\globe{n} \to \globeinv{n}$ are
      the two components of the \oo-functor $j : \sphere{n+1} \to
      \globeinv{n+1}$.
  \end{enumerate}
\end{proposition}

\begin{proof}
  Using the defining property of $\Rev{n+1}$, the second assertion
  translates directly to the definition of an \oo-equivalence. Indeed, if $n
  = 0$, the data of a square as in the assertion boils down to the data of
  an object $y$ of $D$ and the data of a factorization boils down to the
  data of an object $x$ in $C$ and of an incoherent witness of the fact that
  $f(x)$ and $y$ are \oo-equivalent. Similarly, if~$n > 0$, the data of the
  square boils down to the data of two $(n-1)$-cells $x$ and $x'$ of $C$ and
  an $(n+1)$-cell $v : f(x) \to f(x')$ in $D$, and the data of the
  factorization boils down to the data of an $(n+1)$-cell $u : x \to y$ of
  $C$ and an incoherent witness that $f(u)$ and $v$ are \oo-equivalent. The
  third assertion translates similarly replacing incoherent witnesses by
  coherent witnesses and the result follows from the previous proposition.
\end{proof}

\begin{paragraph}\label{paragr:fib}
  We will denote by
  \[
    \gentrivcof{n} : \globe{n} \to \globeinv{n+1}
  \]
  the composite
  \[
    \xymatrix{\globe{n} \ar@{^{(}->}[r]^{\leftincl} & \sphere{n+1} \ar[r]^j & \globeinv{n+1}
    \pbox{,}}
  \]
  where $\leftincl : \globe{n} \into \sphere{n+1}$ denotes the \oo-functor
  corresponding to ``source $n$\nbd-cell'' of $\sphere{n+1}$. Note that
  $\gentrivcof{n}$ is defined as a composite of two cofibrations and is
  hence a cofibration. Moreover, it is an \oo-equivalence. Indeed, we have a
  commutative triangle
  \[
    \xymatrix{
      \globe{n} \ar[r]^{\gentrivcof{n}} \ar[dr]_{\id_{\globe{n}}}
      & \globeinv{n+1} \ar[d]^p \\
      & \globe{n} \pbox{,}
    }
  \]
  where $p$ is trivial fibration and hence an \oo-equivalence, and this
  follows from the 2-out-of-3 property (and more precisely, the easy case of
  Proposition~\ref{prop:2-out-of-3-easy}). We set
  \[
    \setgentrivcof = \setof{\gentrivcof{n}}{n \ge 0}.
  \]

  An \oo-functor $f : C \to D$ is a \ndef{fibration}
  \index{fibration!of 9-categories@of $\omega$-categories}
  \index{9-category@$\omega$-category!fibration}
  if it has the right
  lifting property with respect to $\setgentrivcof$. Concretely, this means
  that $f$ is a fibration if, for every $n \ge 0$, every $n$-cell $u$ of
  $C$, every $n$-cell $v$ of $D$, every coherent witness that $f(u)$ and $v$
  are \oo-equivalent can be lifted to a coherent witness that $u$ and some
  $n$-cell $u'$ such that $f(u') = v$ are \oo-equivalent.
\end{paragraph}

\begin{theorem}\label{thm:triv_fib}
  An \oo-functor $f$ is a trivial fibration if and only if it is both an
  \oo-equivalence and a fibration.
\end{theorem}

\begin{proof}
  As elements of $\setgentrivcof$ are cofibrations, \ie elements of
  $\lrorth{\setgencof}$, we have
  \[
    \rorth{\setgencof} = \rlrorth{\setgencof} \subseteq
    \rorth{\setgentrivcof},
  \]
  meaning that trivial fibrations are fibrations. As
  we already know that trivial fibrations are \oo-equivalence
  (Proposition~\ref{prop:triv_fib_eq}), this establishes one implication.

  Conversely, let $f : C \to D$ be an \oo-functor being both an
  \oo-equivalence and a fibration. Fix $n \ge 0$ and consider a commutative
  square
  \[
    \xymatrix{
      \sphere{n} \ar[d]_{\gencof{n}} \ar[r] & C \ar[d]^f \\
      \globe{n} \ar[r] & D \pbox{.}
    }
  \]
 As $f$ is an \oo-equivalence, by Proposition~\ref{prop:char_w-eq} this
 square factors as
  \[
    \xymatrix{
      \sphere{n} \ar[d]_{\gencof{n}} \ar[r]^{\gencof{n}} & \globe{n}
      \ar[r] \ar[d]_{\gentrivcof{n}} & C \ar[d]^f \\
      \globe{n} \ar[r] & \globeinv{n+1} \ar[r] & D \pbox{.}
    }
  \]
  But as $f$ is a fibration, the right square of the factorization admits a
  lift, showing that the initial square admits a lift and hence that $f$ is
  a trivial fibration.
 \end{proof}

\section{Immersions}

\begin{paragraph}\label{paragr:immersion}
  \index{immersion!of 9-categories@of $\omega$-categories}
  An \oo-functor $i : C \to D$ is an \ndef{immersion} if it admits a
  retraction $r : D \to C$, so that we have $ri = \id_{C}$, and a
  reversible oplax transformation $\alpha : ir \tod \id_{D}$ such that $\alpha
  \comp{} i = \unit{i}$ (see~\cref{paragr:op_oo-trans}).
  Diagrammatically, we have the following commutative diagrams:
  \[
    \vcenter{
      \xymatrix{
        C \ar[r]^i \ar@/_3ex/[rr]_{\unit{C}}\ar@/^3ex/@{{}{ }{}}[rr]^{\phantom{\unit{C}}} & D \ar[r]^r & C\\
      }}
    \qquad
    \vcenter{
      \xymatrix@R=1pc{
        & & D \\
        D \ar[r]^-{\alpha} \ar@/^3ex/[urr]^{ir} \ar@/_3ex/[drr]_{\id_D}
        & \CatCylRev(D) \ar[ru]^{\projcyltop_D} \ar[rd]_{\projcylbot_D} \\
        & & D
      }}
    \qquad
    \vcenter{
      \xymatrix{
        C \ar[r]^i \ar[d]_i & D \ar[d]^\alpha \\
        D \ar[r]_-{\trivcyl_D} & \CatCylRev(D) \pbox.
      }}
  \]
  The last square can also be written
  \[
      \xymatrix@C=3pc{
        C \ar[r]^i \ar[d]_{\trivcyl_C} & D \ar[d]^\alpha \\
        \CatCylRev(C) \ar[r]_-{\CatCylRev(i)} & \CatCylRev(D) \pbox{.}
      }
  \]
  This notion of an immersion is an \oo-categorical version of the notion of
  a ``strong deformation retract'' in topology.
\end{paragraph}

\begin{proposition}
  An \oo-functor $f : C \to D$ is an immersion if and only if the
  commutative square
  \[
    \xymatrix{
      C \ar[r]^{\incmapcyl{f}} \ar[d]_f & \MapCyl{f} \ar[d]^{\projmapcyl{f}} \\
      D \ar[r]_{\id_D} & D
    }
  \]
  admits a lift.
\end{proposition}

\begin{proof}
  By definition, we have a pullback square
  \[
    \xymatrix{
      \MapCyl{f} \ar[d]_{\leftproj} \ar[r]^-{\rightproj} & \CatCylRev(D)
      \ar[d]^{\projcyltop} \\
      C \ar[r]_f & D \pbox{.}
    }
  \]
  This means that an \oo-functor $h : D \to \MapCyl{f}$ corresponds to
  the data of an \oo-functor $g : D \to C$ and an \oo-functor $\alpha : D
  \to \CatCylRev(D)$ such that $\projcyltop\alpha = fg$, i.e., a reversible
  oplax transformation whose source is $fg$. The commutativity of the
  right-lower triangle of
  \[
    \xymatrix{
      C \ar[r]^{\incmapcyl{f}} \ar[d]_f & \MapCyl{f} \ar[d]^{\projmapcyl{f}} \\
      D \ar[r]_{\id_D} \ar[ru]^h & D \pbox{,}
    }
  \]
  using the equality $\projmapcyl{f} = \projcylbot\rightproj$, means that
  $\projcylbot\alpha = \id_D$, i.e., that the target of $\alpha$ is~$\id_D$.
  Finally, the commutativity of the left-upper triangle, using the fact
  that $\incmapcyl{f} = (\id_{C}, \trivcyl_D f) : C \to C \times_D \CatCylRev(D)$,
  means that $(g, \alpha)f = (\id_{C}, \trivcyl_Df)$, i.e., that $gf = \id_C$ and
  that $\alpha \comp{} f = \unit{f}$, thereby proving the result.
\end{proof}

\begin{corollary}\label{coro:w-eq_cof_imm}
  If an \oo-functor is both an \oo-equivalence and a cofibration, then it
  is an immersion.
\end{corollary}

\begin{proof}
  If $f : C \to D$ is an \oo-equivalence, then, by
  Proposition~\ref{prop:w-eq_mapcyl}, $\projmapcyl{f} :
  \MapCyl{f} \to D$ is a trivial fibration. In particular, if $f$ is
  additionally a cofibration, then the square of the previous proposition admits
  a lift, showing that $f$ is indeed an immersion.
\end{proof}

\begin{proposition}\label{prop:imm_w-eq}
  Immersions are \oo-equivalences.
\end{proposition}

\begin{proof}
  Let $i : C \to D$ be an immersion, and let $r$ and $\alpha$ be as in the
  definition of~\cref{paragr:immersion}.
  \begin{enumerate}
    \item Let $y$ be an object of $D$. The reversible cell
      $\alpha_y : ir(y) \to y$ shows that, for $x = r(y)$, $i(x)$ and $y$
      are \oo-equivalent.
    \item Let $x$ and $y$ be two parallel $n$-cells of $C$ and let $v : i(x)
      \to i(y)$ be an $(n+1)$\nbd-cell of $D$. Consider the cell $r(v)$. As
      $ri = \id_{C}$, we have $r(v) : x \to y$. Consider now the reversible
      cell $\alpha_v$. A priori, its source and target are given by
      \[
        \alpha_v
        :
        ir(v) \comp{0} \alpha_{\tge{0}(v)} \comp{1} \cdots \comp{n} \alpha_{\tge{n}(v)}
        \to
        \alpha_{\sce{n}(v)} \comp{n} \cdots \comp{1}
        \alpha_{\sce{0}(v)}\comp{0} v.
      \]
      As the cells $\sce{0}(v), \tge{0}(v), \dots, \sce{n}(v), \tge{n}(v)$
      are in the image of $i$, it follows from the relation $\alpha \comp{} i =
      \unit{i}$ that the cells $\alpha_{\sce{0}(v)}, \alpha_{\tge{0}(v)},
      \dots, \alpha_{\sce{n}(v)}, \alpha_{\tge{n}(v)}$ are units, so that we
      have $\alpha_v : ir(v) \to v$. This shows
      that $i(u)$, where $u = r(v)$, and $v$ are \oo-equivalent.
      \qedhere
  \end{enumerate}
\end{proof}

\begin{proposition}\label{prop:immersion_pushout}
  The class of immersions is closed under pushouts.
\end{proposition}

\begin{proof}
  Let $i : C \to D$ be an immersion, and let $r$ and $\alpha$ be as in the
  definition of~\cref{paragr:immersion}. Consider a pushout square
  \[
    \xymatrix{
      C \ar[d]_i \ar[r]^f & C' \ar[d]^{i'} \\
      D \ar[r]_g & D' \pbox{.}
    }
  \]
  Using the universal property of this square, we get a retraction $r'$ of $i'$:
  \[
    \xymatrix@C=1pc@R=1pc{
      C \ar[dd]_i \ar[rr]^f && C' \ar[dd]^{i'} \ar@/^2ex/[dddr]^{\id_{C'}} \\
      \\
      D \ar[rr]_g \ar@/_2ex/[drrr]_{fr} && D' \ar@{.>}[dr]^{r'} \\
                                      &&& C' \pbox{.}
    }
  \]
  Similarly, we get a reversible oplax transformation $\alpha'$:
  \[
    \xymatrix@C=0.5pc@R=1pc{
      C \ar[dd]_i \ar[rr]^f && C' \ar[dd]^{i'} \ar[dr]^-{\trivcyl_{C'}} \\
      &&& \CatCylRev(C') \ar[dd]^{\CatCylRev(i')} \\
      D \ar[rr]_g \ar[dr]_-{\alpha} && D' \ar@{.>}[dr]^-{\alpha'} &
      \\
      & \CatCylRev(D) \ar[rr]_-{\CatCylRev(g)} && \CatCylRev(D') \pbox{.}
    }
  \]
  Indeed, the outer diagram commutes as
  \[
    \begin{split}
      \CatCylRev(i') \trivcyl_C' f = \CatCylRev(i') \CatCylRev(f) \trivcyl_C
      = \CatCylRev(g) \CatCylRev(i) \trivcyl_C = \CatCylRev(g) \alpha i.
    \end{split}
  \]
  Using the uniqueness part of the same universal property, one checks that
  $\projcyltop \alpha = i'r'$ and $\projcylbot \alpha = \id_{D'}$, i.e.,
  that $\alpha : i'r' \tod \id_{D'}$. The fact that $\alpha' \comp{} i' =
  \unit{i'}$ is expressed in one of the commutative squares defining
  $\alpha'$. This proves that $i' : C' \to D'$ is indeed an immersion.
\end{proof}

\begin{theorem}\label{thm:w-eq_cof_closed}
  The class of \oo-functors that are both \oo-equivalences and cofibrations
  is closed under pushouts.
\end{theorem}

\begin{proof}
  The class of cofibrations being closed under pushouts
  (see~\cref{paragr:stab_lI}), all we have to show is that the pushout of an
  \oo-functor $f$ being both an \oo-equivalence and a cofibration is
  an \oo-equivalence. But by Corollary~\ref{coro:w-eq_cof_imm}, such an $f$
  is an immersion. Hence, by the previous proposition, a pushout of $f$ is
  still an immersion, and hence an \oo-equivalence by
  Proposition~\ref{prop:imm_w-eq}.
\end{proof}

\begin{coro}\label{coro:lrJ}
  Every element of $\lrorth{\setgentrivcof}$ is both an \oo-equivalence and
  a cofibration.
\end{coro}

\begin{proof}
  Denote by $\clW$ the class of \oo-equivalences and by $\clC$
  the class of cofibrations. By Propositions~\ref{prop:w-eq_retr}
  and~\ref{prop:w-eq_filtered}, the class $\clW$ is closed under retracts and
  countable compositions. The same properties hold for the class $\clC =
  \lrorth{\setgencof}$ (see~\cref{paragr:stab_lI}). Using the
  previous proposition, we get that the class $\clW \cap \clC$ is closed
  under retracts, pushouts, and countable compositions. It follows, using the
  small object argument (Proposition~\ref{prop:small_obj_arg_2}), that to
  show the inclusion $\lrorth{\setgentrivcof} \subset \clW \cap \clC$, it
  suffices to show the inclusion $\setgentrivcof \subset \clW \cap \clC$.
  But we already noted in~\cref{paragr:fib} that for every $n$, the
  \oo-functor $\gentrivcof{n} : \globe{n} \to \globeinv{n+1}$ is both an
  \oo-equivalence and a cofibration.
\end{proof}

\begin{remark}\label{rem:lrJ}
  We will see in the proof of Theorem~\ref{thm:folk} that it follows
  formally from what we proved so far that, conversely, an \oo-functor being
  both an \oo-equivalence and a cofibration is in $\lrorth{\setgentrivcof}$.
\end{remark}

\chapter{The Folk Model Structure}
\index{folk model structure}
\index{model structure!folk}
\label{chap:folk}

This chapter is about proving the existence of the so-called folk model
category structure on $\ooCat$, following Lafont, Métayer, and
Worytkiewicz~\cite{LMW}. This model category structure is a generalization
of a model category structure on~$\Cat$ whose weak equivalences are the
equivalences of categories, a folklore result, whence the name. The
analogous result for $2$-categories was proved by Lack~\cite{LackFolk2,
LackFolkBi}.

The \emph{folk model category structure} is a model category structure on $\ooCat$
whose weak equivalences are the $\omega$-equivalences and whose cofibrant
resolutions are the polygraphic resolutions. It is the natural homotopical
framework in which the notion of polygraphic resolutions lives. As a
convincing evidence of this, we will see in the next chapter that Métayer's
polygraphic homology can be expressed as a derived functor with respect to
the folk model category structure.

The chapter is organized as follows. The first section is devoted to the
proof of the existence of the folk model category structure. The hard
work was done in the previous chapter, and we are basically just assembling
various results. Our proof differs in one point from the original one:
we avoid the use of Smith's theorem.  In the second section, we
prove, still according to \cite{LMW}, that the \oo-category $\CatCylRev(C)$
of reversible cylinders in an \oo-category $C$ forms a path object for $C$ in
the sense of the folk model category structure. We deduce from this fact that
polygraphic resolutions are unique in a stronger sense than the one proved
before. In the last section, we transfer the folk model structure
on~$\ooCat$ to the category of $(n, p)$-categories. We give explicit
descriptions of the resulting structures for various special cases:
$n$-categories with particular small values of $n$, \oo-groupoids, and $(\omega,
1)$-categories.

\section{The Folk Model Structure on \pdfm{\ooCat}}

The purpose of this section is to prove the existence of the so-called
``folk model structure'' on $\ooCat$. We start by recalling the definition
of a model category, see \chapr{model-cat} for a more detailed presentation.

\begin{paragraph}
  \index{model category}
  \index{weak equivalence}
  \index{fibration}
  \index{cofibration}
  \index{trivial fibration}
  \index{fibration!trivial}
  \index{trivial cofibration}
  \index{cofibration!trivial}
  A \ndef{model category} is a category $\M$ endowed with three classes of
  maps: the \ndef{weak equivalences}, the \ndef{cofibrations}, and the
  \ndef{fibrations}; these data are required to satisfy the following
  axioms:
  \begin{enumerate}
    \item the category $\M$ is finitely complete and finitely cocomplete,
    \item the class of weak equivalences satisfies the 2-out-of-3 property,
    \item the class of weak equivalences, cofibrations, and fibrations are
      closed under retracts,
    \item cofibrations have the left lifting property with respect to
      \ndef{trivial fibrations} (that is, maps that are both a fibration and
      a weak equivalence); trivial cofibrations (that is, maps that are both
      a cofibration and a weak equivalence) have the left lifting property
      with respect to fibrations, and
    \item every map of $\M$ factors as a cofibration followed by a trivial
      fibration, and as a trivial cofibration followed by a fibration.
  \end{enumerate}
\end{paragraph}

\begin{theorem}\label{thm:folk}
  The classes of \oo-equivalences (see~\cref{paragr:def_oo-equiv}),
  cofibrations (see~\cref{paragr:cof_triv_fib}) and fibrations
  (see~\cref{paragr:fib}) define a model structure on $\ooCat$ known as the
  ``folk model structure''.
\end{theorem}

\begin{proof}
  Let us denote by $\clW$ the class of weak equivalences, by $\clC$ the
  class of cofibrations and by $\clF$ the class of fibrations.
  By definition, we have
  \[ \clC = \lrorth{\setgencof} \qtand \clF = \rorth{\setgentrivcof} \]
  (see \cref{paragr:cof_triv_fib,paragr:fib}).
  We showed that
  \begin{enumerate}
    \item\label{item:folk_a} $\clW$ satisfies the 2-out-of-3 property
      (Theorem~\ref{thm:w-eq-2-3}) and is closed under retracts
      (Proposition~\ref{prop:w-eq_retr}),
      \item\label{item:folk_b}
      $\lrorth{\setgentrivcof} \subset \clW \cap \clC$
      (Corollary~\ref{coro:lrJ}),
    \item\label{item:folk_c} $\clW \cap \clF = \rorth{\setgencof}$
      (Theorem~\ref{thm:triv_fib}).
  \end{enumerate}
  We will now see that these properties formally imply the theorem.

  Let us first prove that the inclusion in \ref{item:folk_b}
  is actually an equality. Let
  $f$ be an \oo-functor in $\clW \cap \clC$. Applying the small object
  argument (Proposition~\ref{prop:small_obj_arg_2}) to $f$ and the set
  $\setgentrivcof$ yields a factorization~$f = pi$, where $p$ is in $\clF =
  \rorth{\setgentrivcof}$ and $i$ in $\lrorth{\setgentrivcof}$.
  By~\ref{item:folk_b}, $i$ is in $\clW$, and so
  is $p$ by~\ref{item:folk_a}. It follows from~\ref{item:folk_c} that $p$ is
  in~$\rorth{\setgencof}$. As $f = pi$ is in $\clC = \lrorth{\clI}$, it has
  the left lifting property with respect to $p$ and, by the retract lemma
  (Lemma~\ref{lemma:retr}), $f$~is a retract of~$i$. As $i$ is in
  $\lrorth{\setgentrivcof}$, so is $f$, thereby proving the desired
  equality.

  Let us now check that we indeed have a model structure. The
  category~$\ooCat$ is complete and cocomplete. We have already proved that
  the class of \oo-equiv\-alences satisfies the 2-out-of-3 property and is
  closed under retracts. The fact that the classes of cofibrations and fibrations
  can be defined by lifting conditions implies that they are closed under
  retracts as well. The equalities
  \[
      \clC = \lrorth{\setgencof}, \quad \clF = \rorth{\setgentrivcof}, \quad
      \clW \cap \clC = \lrorth{\setgentrivcof},
      \qtand
      \clW \cap \clF = \rorth{\setgencof}
  \]
  show that the required lifting properties are fulfilled and, by applying
  the small object argument, that the required factorization properties are
  fulfilled as well.
\end{proof}

\begin{remark}
  \index{combinatorial model category}
  \index{model category!combinatorial}
  \label{rem:comb_model_str}
  The folk model structure is what is called a ``combinatorial model
  structure'': a model category $\C$ is said to be
  \ndef{combinatorial} if, first, the category $\C$ is locally presentable
  (see \chapr{loc_pres}) and, second, there exist \emph{sets} (by
  opposition to classes)
  of morphisms $\setgencof$ and $\setgentrivcof$ such that the class of
  cofibration is $\lrorth{\setgencof}$ and the class of trivial cofibrations
  is $\lrorth{\setgentrivcof}$. The sets $\setgencof$ and~$\setgentrivcof$
  (which are not unique) are then said to \ndef{generate} the model
  category $\C$.
\end{remark}

\begin{paragraph}
  Let $\C$ be a model category. An object $X$ of $\C$ is said to be
  \ndef{cofibrant}\index{cofibrant!object} if the unique morphism from the initial object of $\C$ to
  $X$ is a cofibration; similarly, $X$ is said to be
  \ndef{fibrant}\index{fibrant!object} if the unique morphism from $X$ to
  the terminal object of~$\C$ is a fibration.
\end{paragraph}

\begin{proposition}
  Every \oo-category is fibrant in the folk model structure.
\end{proposition}

\begin{proof}
  Let $C$ be an \oo-category. Consider a diagram
  \[
     \xymatrix{
       D \ar[r]^f \ar[d]_i & C \\
       E & \pbox{,}
     }
  \]
  where $i$ is a trivial cofibration. By Corollary~\ref{coro:w-eq_cof_imm},
  $i$ admits a retraction $r$ and the triangle
  \[
     \xymatrix{
       D \ar[r]^f \ar[d]_i & C \\
       E  \ar[ur]_{fr}
     }
  \]
  is commutative, thereby proving that $C$ is fibrant.
\end{proof}

\begin{theorem}
  \label{thm:cofib}
  \label{thm:pol-cof-obj}
  The cofibrant objects of the folk model structure are exactly the
  \oo-categories generated by polygraphs.
\end{theorem}

\begin{proof}
 Let $C$ be a cofibrant \oo-category. By
  Proposition~\ref{prop:factor_cof_fib_triv}, the \oo-functor $\varnothing
  \to C$ is a retract of a relative polygraph $D \to E$. In particular,
  $\varnothing$ is a retract of~$D$ and $C$ is a retract of $E$. This means
  that $D \simeq \varnothing$ and hence that $D$ is generated by a
  polygraph. This shows that $C$ is a retract of an \oo-category generated
  by a polygraph.

  To conclude the proof, it suffices to show that a retract of an
  \oo-category generated by a polygraph is also generated by a
  polygraph. This immediately follows from~\cite[Theorem
  7.1]{metayer2008cofibrant}, stating that the full 
  subcategory of $\ooCat$ consisting of those \oo-categories that are
  generated by polygraphs is Cauchy-complete. We briefly describe the
  general idea of the argument, and refer
  to~\cite{metayer2008cofibrant} for a complete proof. Thus, let $P$
  be a polygraph and $h:\free P\to \free P$ be an idempotent morphism
  in $\ooCat$.  We have to build a polygraph $Q$, together with
  morphisms $r:\free P\to \free Q$ and $s:\free Q\to \free P$ such
  that $rs=\unit{\free{Q}}$ and $sr=h$. The polygraph $Q$ and the maps
  $r$, $s$ are defined by simultaneous induction on the dimension. In
  dimension~$0$, $Q_0=\setof{h(x)}{x\in \free P_0=P_0}$, the map $s$
  is given by the obvious inclusion of $Q_0$ into $P_0$ and
  $r(x)=h(x)$ for all $x\in P_0$. Let now~$n>0$ and suppose we have
  defined $Q$, $r$, and $s$ up to dimension $n-1$, satisfying the
  above equations. We must now define the set $Q_n$
  of $n$-generators in $Q$, together with source and target maps
  $\sce{n-1},\tge{n-1}:Q_n\to \free Q_{n-1}$.  The crucial observation
  is that $h$ induces a partition of $P_n$ in three subsets:
  \[
    P_n=P_n^0\uplus P_n^1\uplus P_n^2
  \]
  where $P_n^0$ is the set of those $n$-generators whose image by $h$
  is an identity, $P_n^1$~is the set of those generators $a\in P_n$ such that
  $h(\free a)=\ctx c{\free a}$ where $\varctx cx$ is a thin
  $n$-context~(see \cref{sec:contexts}), and $P_n^2$ is the set of remaining
  $n$-generators. We may now define
  \[
    Q_n=\setof{h(\free a)}{a\in P_n^1}
  \]
  which gives an obvious inclusion $i:Q_n\to \free P_n$. By induction
  hypothesis, the morphism $r$ is defined up to dimension $n-1$, so
  that we can define
  \[
    \sce{n-1},\tge{n-1}:Q_n\to \free Q_{n-1}
    \]
  by requiring the commutation of the following square
  \[
    \xymatrix{Q_n \ar[d]_i\ar[rr]^{\sce{n-1},\tge{n-1}}& &\free Q_{n-1}\\
    \free P_n \ar[rr]^{\sce{n-1},\tge{n-1}}& &\free
    P_{n-1}\ar[u]^{r}}
    \]
  for source and target maps. Thus the $\omega$-category $\free Q$ is now
  defined up to dimension $n$ and one easily checks that $i:Q_n\to
  \free P_n$ induces an extension of $s$ up to dimension $n$ from
  $\free Q$ to $\free P$ satisfying $hs=s$.

  Finally, we must extend $r$ up to dimension $n$. The main difficulty
  is to find where the generators belonging to
  $S_n^2$ should be sent. Therefore, we first consider an
  auxiliary subpolygraph $R$ of $P$,  identical to $P$ up to dimension
  $n-1$, and whose set of $n$-generators is just $R_n=P_n^0\uplus
  P_n^1$. The inclusion $j:R\to P$ now induces a morphism $\free j$ of
  $n$-categories from $\free R\to \free P$, and one shows the
  existence of $n$-morphisms $h':\free R\to \free R$ and $k:\free P\to
  \free R$ such that the following square commutes:
  \[
    \xymatrix{\free R\ar[r]^{h'}\ar[d]_{\free{j}} & \free R\ar[d]^{\free{j}}\\
    \free P\ar[r]_h \ar[ur]|k& \free P}
    \]
  The purpose of the previous step was precisely to get rid of
  $S_n^2$. Now one builds $r':\free R\to \free Q$ and $s':\free Q\to
  \free R$, extending $r$, $s$ up to dimension $n$  is such a way that
  $s'r'=h'$ and $r's'=\unit{\free{Q}}$. Note that the proof of this
  last equation relies on the fine properties of contexts, as
  presented in \cref{sec:contexts}. We obtain the desired
  retraction $r$ up to dimension $n$ by taking $r=r'k$.
 \end{proof}

\begin{remark}
  The previous theorem easily implies that the folk model structure is not
  ``cartesian'' (meaning that the cartesian product functor is not a Quillen
  bifunctor in the sense of \cite[Definition 4.2.1]{hovey2007model}).
  Indeed, in a cartesian model structure, the product of two cofibrant
  objects is cofibrant. But it is not true that the product of two
  \oo-categories generated by a polygraph is generated by a polygraph. For
  instance, the product $\globe{1} \times \globe{1}$ is not generated by a
  polygraph.
\end{remark}

\section{The Path Objects of Cylinders}

The goal of this section is to prove that, for every \oo-category $C$, the
\oo-cate\-gory~$\CatCylRev{C}$ (see~\cref{paragr:def_catcylrev}), endowed
with the maps
\[
  \xymatrix@=2.5pc{
    C \ar[r]^-{\trivcyl_C} &
    \CatCylRev(C) \ar[r]^-{(\projcyltop, \projcylbot)} &
    C \times C \pbox{,}
  }
\]
is a path object for $C$ in the folk model structure in the following sense:

\begin{paragraph}[Path objects]
  Let $\M$ be a model category and let $X$ be an object of~$\M$. A
  \ndef{path object}\index{path!object} for $X$ in $\M$ is an object $P$ of $\M$, endowed with
  a factorization
  \[ \xymatrix{ X \ar[r]^r & P \ar[r]^-p & X \times X} \]
  of the diagonal of $X$, where $r$ a weak equivalence of $\M$ and $p$ a
  fibration of~$\M$.
\end{paragraph}

\medskip

As we already know that $\trivcyl_C : C \to \CatCylRev(C)$ is an
\oo-equivalence (Proposition~\ref{prop:projcyl_triv_fib}), all we have to
prove is that $(\projcyltop, \projcylbot) : \CatCylRev(C) \to C \times C$ is a
fibration. To do so, we will use the following variation on the notion of an
immersion:

\begin{paragraph}[Strong immersions]\label{paragr:def_strong_immersion}
  An \oo-functor $i : C \to D$ is a \ndef{strong immersion}\index{immersion!strong} if it admits a
  retraction $r : D \to C$ and reversible oplax transformations $\alpha : ir
  \tod \id_{D}$ and $\alpha' : \id_{D} \tod ir$ such that
  $\alpha \comp{} i = \unit{i}$ and $\alpha' \comp{} i = \unit{i}$.
\end{paragraph}

\begin{proposition}
  Let $C$ be an \oo-category. Then the \oo-functor
  \[ (\projcyltop, \projcylbot) : \CatCylRev(C) \to C \times C \]
  has the right lifting properties with respect to strong immersions.
\end{proposition}

\begin{proof}
  Let $i : D \to E$ be a strong immersion and let $r$, $\alpha$, and
  $\alpha'$ be as in~\cref{paragr:def_strong_immersion}. Consider a
  commutative square
  \[
    \xymatrix{
      D \ar[d]_i \ar[r]^-{\gamma} &
      \CatCylRev(C) \ar[d]^{(\projcyltop, \projcylbot)} \\
      E \ar[r]_-{(f, g)} & C \times C \pbox{.}
    }
  \]
  In other words, we have two \oo-functors $f, g : E \to C$ and a reversible
  oplax transformation $\gamma : fi \tod gi$. A lift of such a square amounts
  to a reversible oplax transformation $\delta : f \tod g$ such that $\delta
  \comp{} i = \gamma$. One defines $\delta$ as the composite
  \[
    \xymatrix{
      f \ar@2[r]^-{f \comp{} \alpha'} & fir \ar@2[r]^-{\gamma \comp{} r} & gir
      \ar@2[r]^-{g \comp{} \alpha} & g
    }
  \]
  (see Remark~\ref{rem:2-Cat_oplax}). We have
  \[
    \begin{split}
      \delta \comp{} i &
      = \big((g \comp{} \alpha)(\gamma \comp{} r)(f \comp {} \alpha')\big) \comp{} i
      \\
      &
      = (g \comp{} \alpha \comp{} i)(\gamma \comp{} ri)(f \comp {} \alpha' \comp{} i)
      \\
      &
      = \unit{gi} (\gamma \comp{} \unit{D}) \unit{fi}
      = \gamma,
    \end{split}
    \]
    thereby proving the result.
\end{proof}

\begin{proposition}
  For every $n \ge 0$, the \oo-functor $\gentrivcof{n} : \globe{n} \to
  \globeinv{n+1}$ of~\cref{paragr:def_globeinv} is a strong immersion.
\end{proposition}

\begin{proof}
  Let $n \ge 0$. By Corollary~\ref{coro:w-eq_cof_imm}, every trivial
  cofibration is an immersion. This means that there exists a section $r :
  \globeinv{n+1} \to \globe{n}$ and a reversible oplax transformation $\alpha
  : ir \tod \unit{D}$. A dual proof shows that there exists a section $r' :
  \globeinv{n+1} \to \globe{n}$ and a reversible oplax transformation $\alpha'
  : \unit{D} \tod ir'$. (More formally, one can apply
  Corollary~\ref{coro:w-eq_cof_imm} to $(\gentrivcof{n})^\mathrm{o}$,
  where ${-}^\mathrm{o} : \ooCat \to \ooCat$ denotes the duality of
  $\ooCat$ consisting in reverting the orientation of $n$-cells for every $n
  \ge 1$.) To conclude the proof it suffices to check that $r = r'$.
  Consider the reversible oplax transformation $r \comp{} \alpha' : r \tod
  rir' = r'$. As in $\globe{n}$, the only reversible cells are the
  identities, this oplax transformation has to be the identity, showing that
  $r = r'$.
\end{proof}

\begin{corollary}
  For every \oo-category $C$, the \oo-functor
  \[ (\projcyltop, \projcylbot) : \CatCylRev(C) \to C \times C \]
  is a fibration.
\end{corollary}

\begin{proof}
  This is an immediate consequence of the two previous propositions.
\end{proof}

\begin{remark}
  The same proof shows that $(\projcyltop, \projcylbot) : \CatCyl(C) \to C
  \times C$ is also a fibration.
\end{remark}

We thus have showed:

\begin{theorem}
  For every \oo-category $C$, the \oo-category $\CatCylRev(C)$ (endowed with
  the \oo-functors described before) is a path object for $C$ in the folk
  model structure.
\end{theorem}

\begin{remark}
  Although this is not related to the goal of the section, let us mention
  that one can prove in a very similar way that for every \oo-functor $f : C \to
  D$, the \oo-functor $\projmapcyl{f} : \MapCyl{f} \to D$ is a fibration.
  Indeed, one checks that it has the right lifting property with respect to
  (not necessarily strong) immersions and the assertion thus follows from
  Corollary~\ref{coro:w-eq_cof_imm}.
\end{remark}

\begin{proposition}\label{prop:pol_res_uniq_hom}
  Let $f : C \to D$ be an \oo-functor and let $(P, p)$ and $(Q, q)$ be polygraphic
  resolutions (see~\cref{paragr:pol_res}) of $C$ and $D$, respectively. If
  $g, g' : P^\ast \to Q^\ast$ are two \oo-functors such that the two squares
  \[
    \xymatrix@C=2pc@R=2pc{
      P^\ast \ar[d]_{p} \ar@/^2ex/[r]^g \ar@/_2ex/[r]_{g'}
      & Q^\ast \ar[d]^q \\
      C \ar[r]_f & D
    }
  \]
  commute, then there exists a reversible oplax transformation $\alpha : g
  \tod g'$.
\end{proposition}

\begin{proof}
  The fact that $\CatCylRev(A)$ is a path object for any \oo-category $A$
  implies that $u : A \to B$ is ``right homotopic'' to $v : A \to B$ in the
  sense of model categories if and only if there exists a reversible
  oplax transformation from $u$ to $v$. The result thus follows from general
  properties of cofibrant resolutions in model categories (see for
  instance~\cite[Proposition 8.1.25]{HirMC}).
\end{proof}

\begin{rem}
  In particular, in the case where $C = D$ and $f$ is the identity
  \oo-functor, we get a diagram
  \[
    \xymatrix@C=1pc@R=2pc{
      P^\ast \ar[dr]_{p} \ar@/^2ex/[rr]^{g}_{}="g" \ar@/_2ex/[rr]_{g'}_{}="g'"
      \ar@2"g";"g'"^\alpha
      & & Q^\ast \ar[dl]^q \\
      & C
    }
  \]
  from which one can deduce that the morphism given in \remr{pol-res-func}
  between any two resolutions is unique up to a (non-canonical) reversible
  oplax transformation.
\end{rem}

\section{The Folk Model Structure on \pdfm{\nCat{n}} and \pdfm{\npCat{n}{p}}}
\label{sec:folk_trans}

The purpose of this section is to transfer the folk model structure on
$\ooCat$ to subcategories of $\ooCat$ such as the $\nCat{n}$ or more generally
$\npCat{n}{p}$.
To do so, we will use the following classical transfer lemma:

\begin{lemma}
  Let $(\M, \clW, \clC, \clF)$ be a combinatorial model category generated
  by $I$ and $J$.  Let $\C$ be a locally presentable category, and let
  \[ F : \M \to \C, \quad G : \C \to \M \]
  be a pair of adjoint functors. Suppose that
  \[ G(\lrorth{F(J)}) \subset \clW, \]
  where $\clW$ denotes the class of weak equivalences of $\M$.
  Then $F(I)$ and $F(J)$ generate a combinatorial model structure on $\C$,
  whose class of weak equivalences is $G^{-1}(\clW)$ and whose class of
  fibrations is $G^{-1}(\clF)$.
\end{lemma}

\begin{proof}
  See for instance \cite[Theorem 3.3]{CransTrans}.
\end{proof}

\begin{proposition}\label{prop:abs_folk}
  Let $\C$ be a reflective subcategory of $\ooCat$ closed under pushouts and
  filtered colimits. Suppose that, for every \oo-category $C$ in $\C$, the
  \oo-category $\CatCyl(C)$ is still in $\C$. Then there exists a model
  structure on $\C$ whose weak equivalences are the \oo-equivalences between
  objects of $\C$ and whose fibrations are the ``folk'' fibrations
  (i.e., the one defined in~\cref{paragr:fib}) between objects of $\C$.
  This model structure is generated by $F(\setgencof)$ and
  $F(\setgentrivcof)$, where $F : \ooCat \to \C$ denotes the left adjoint of
  the inclusion functor.
\end{proposition}

\begin{proof}
  As $\C$ is a reflective subcategory closed under filtered colimits of a
  locally presentable category, the category $\C$ is itself locally
  presentable. Let us check the hypothesis of the previous lemma: we have to
  show that $\lrorth{F(\setgentrivcof)}$ is included in the class of \oo-equivalences.
  By the small object argument applied to $F(\setgentrivcof)$ in $\C$,
  elements of $\lrorth{F(\setgentrivcof)}$ are transfinite compositions of
  pushouts of elements of $F(\setgentrivcof)$, these colimits being taken in
  $\C$. By hypothesis, the inclusion functors to $\ooCat$ preserves these
  colimits and they can be computed in $\ooCat$. By the following lemma, the
  elements of $F(\setgentrivcof)$ are immersions. By
  Proposition~\ref{prop:immersion_pushout}, a pushout of such an immersion
  is an immersion, and hence is an \oo-equivalence by
  Proposition~\ref{prop:imm_w-eq}. As, by
  Proposition~\ref{prop:w-eq_filtered}, \oo-equivalences are stable under
  transfinite compositions, this concludes the proof.
\end{proof}

\begin{lemma}
  Under the hypothesis of the previous proposition, if $i$ is an immersion,
  so is $F(i)$.
\end{lemma}

\begin{proof}
  Let $D$ be an \oo-category. By hypothesis, $\CatCylRev F(D)$ is in $\C$.
  For $E$ an \oo-category, denote by $\e_E : E \to F(E)$ the canonical
  \oo-functor. The universal property of $F$ gives a dotted arrow
  \[
    \xymatrix{
      \CatCylRev(D) \ar[r]^-{\CatCylRev(\e_D)} \ar[d]_{\e^{}_{\CatCylRevS(D)}} &
      \CatCylRev F(D) \\
      F\mkern1mu\CatCylRev(D) \ar@{.>}[ur]_{\tau^{}_D}
    }
  \]
  making the triangle commute. If $f, g : C \to D$ are two \oo-functors and
  $\alpha : f \tod g$ is a reversible oplax transformation, then
  \[
    \xymatrix{
      F(C) \ar[r]^-{F(\alpha)} & F\mkern1mu\CatCylRev(D) \ar[r]^{\tau_D} &
      \CatCylRev F(D)
    }
  \]
  defines a reversible oplax transformation. By abuse of notation, we will denote it
  by~$F(\alpha)$. One checks that we have $F(\alpha) : F(f) \tod F(g)$ as
  expected.

  Let now $i : C \to D$ be an immersion and let $r$ and $\alpha$ be as in the
  definition of~\cref{paragr:immersion}. By functoriality, $F(r)$ is a
  retraction of $F(i)$. The previous paragraph gives a reversible
  oplax transformation $F(\alpha) : F(i)F(r) \tod \unit{F(D)}$. One checks
  that $F(\alpha) \comp{} F(i) = \id_{F(i)}$, thereby proving the result. (We
  refer the reader to~\cite[Lemma 6.2]{LMW} for a more detailed proof.)
\end{proof}

From now on, we fix $n$ to be either an integer or $\omega$, and $p \le n$.

\begin{proposition}
  If $C$ is an $(n, p)$-category, then so is $\CatCylRev(C)$.
\end{proposition}

\begin{proof}
  Let $k > n$ and let $\alpha : x \cylto y$ be a $k$-cell of
  $\CatCylRev(C)$. Such a cell is a unit in $\CatCylRev(C)$ if and only if $x$,
  $y$, and $\alpha_k$ are units in $C$ and $\alpha_{k-1}^- = \alpha_{k-1}^+$.
  As $C$ is an $n$-category, the $k$-cells $x$, $y$, $\alpha^-_{k-1}$,
  $\alpha^+_{k-1}$ and the $(k+1)$-cell $\alpha_k$ are units. As
  \[
    \begin{split}
      \sce{k-1}(\alpha^-_{k-1})
      & =
      \sce{k-1}(x) \comp{0} \alpha_0^+ \comp{1} \cdots \comp{k-2} \alpha_{k-2}^+
      \\
      & =
      \tge{k-1}(x) \comp{0} \alpha_0^+ \comp{1} \cdots \comp{k-2} \alpha_{k-2}^+
      = \sce{k-1}(\alpha^+_{k-1}),
    \end{split}
  \]
  we have $\alpha^-_{k-1} = \alpha^+_{k-1}$, showing that $\alpha : x \cylto
  y$ is indeed a unit.

  Let now $k > p$ and let $\alpha : x \cylto y$ be a $k$-cell of
  $\CatCylRev(C)$. As $C$ is an $(n, p)$\nbd-category, the $k$-cells $x$ and $y$
  are invertible. We now define an inverse $\beta : x^{-1} \cylto
  y^{-1}$ of $\alpha : x \cylto y$ in $\CatCylRev(C)$ in the following way:
  \[
    \beta^\e_l = \alpha^\e_l \quad\text{for $l \le k-2$},
  \]
  \[ \beta^-_{k-1} = \alpha^+_{k-1}, \qquad \beta^+_{k-1} = \alpha^-_{k-1}, \]
  \[
    \begin{split}
      \beta_k & = (x^{-1} \comp{0} \alpha_0^+ \comp{1} \cdots \comp{k-2}
      \alpha^+_{k-2}\big) \\
              & \qquad \comp{k-1} \alpha_k^{-1} \comp{k-1} (\alpha^-_{k-2}
              \comp{k-2} \cdots \comp{1} \alpha_0^- \comp{0} y^{-1}),
    \end{split}%
  \]%
  the inverse of the $(k+1)$-cells $\alpha_k$ existing by hypothesis. One
  checks that this indeed defines a revertible cylinder and that this
  cylinder is an inverse~of $\alpha : x \cylto y$.
\end{proof}

\begin{theorem}
\label{T:FolkModelStructureCatNP}
  There exists a model structure on $\npCat{n}{p}$, known as ``the folk
  model structure'', whose weak equivalences are the \oo-equivalences
  between $(n, p)$\nbd-categories and whose fibrations are the ``folk''
  fibrations (i.e., the one defined in~\cref{paragr:fib}) between
  $(n, p)$-categories. This model structure is generated by $F(\setgencof)$ and
  $F(\setgentrivcof)$, where $F : \ooCat \to \npCat{n}{p}$ denotes the left
  adjoint of the inclusion functor.
\end{theorem}

\begin{proof}
  The result will be a consequence of Proposition~\ref{prop:abs_folk} once
  we check that the  hypotheses hold. The previous proposition gives us one of
  the hypotheses and we are left to check that the subcategory
  $\npCat{n}{p}$ of $\ooCat$ is stable under pushout and filtered colimits.
  But it is actually stable under any colimits as the inclusion functor
  $\npCat{n}{p} \into \ooCat$ admits a right adjoint, namely the functor
  taking an \oo-category $C$ to the $(n,p)$-category obtained from $C$ by throwing out
  $k$-cells for $k > n$, and non-invertible $l$-cells for $l > p$.
\end{proof}

Let us now describe more precisely this model structure in some
specific cases.

\begin{paragraph}[The folk model structure on $\nCat{n}$]
  The \oo-category $\nCat{n}$, which is nothing but $\npCat{n}{n}$, is
  endowed with a folk model structure by the previous theorem. Its weak
  equivalences, fibrations and hence trivial fibrations are inherited from
  the folk model structure on $\ooCat$. Let us describe its cofibrations. By
  the previous theorem, the class of cofibrations is
  $\lrorth{F(\setgencof)}$, where $F : \ooCat \to \nCat{n}$ denotes the left
  adjoint to the inclusion functor. For $k \ge 0$, we have
  \[
    F(\globe{k}) =
    \begin{cases}
      \globe{k} & \text{if $k \le n$,} \\
      \globe{n} & \text{if $k > n$,}
    \end{cases}
  \]
  and, similarly,
  \[
    F(\sphere{k}) =
    \begin{cases}
      \sphere{k} & \text{if $k \le n+1$,} \\
      \globe{n} & \text{if $k > n+1$.}
    \end{cases}
  \]
  This implies that
  \[
    F(\gencof{k})
    \begin{cases}
      \gencof{k} & \text{if $k \le n$,} \\
      \gencofbar{n+1} & \text{if $k \le n+1$,} \\
      \unit{\globe{n}} & \text{if $k > n+1$,}
    \end{cases}
  \]
  where $\gencofbar{n+1} : \sphere{n+1} \to \globe{n}$ corresponds to the
  collapsing of the two non-trivial $n$-cells of $\sphere{n+1}$. This shows
  that the class of cofibrations is generated by
  \[ \gencof{0}, \dots, \gencof{n}, \gencofbar{n+1}. \]
  If $C$ is an $n$-category, then the \oo-functor from the empty
  $n$-category to $C$ is a cellular extension for these generators if and
  only if $C$ is presented by an $(n+1)$-polygraph, that is, if and only if its
  underlying $(n-1)$-category is freely generated by a polygraph. In
  particular, any $n$-category admits such an $n$-category as a fibrant
  replacement in $\nCat{n}$. Actually, one can show that an $n$-category is
  cofibrant if and only if it is presented by an $(n+1)$-polygraph.

  Similarly, if one chooses carefully the $\globeinv{k+1}$ (see
  Remark~\ref{rem:choice_Jn}), one gets that $F(\gentrivcof{k}) =
  \unit{\globe{n}}$ if $k \ge n$ so that trivial cofibrations of $\nCat{n}$
  are generated by
  \[
     \gentrivcofbar{1} = F(\gentrivcof{1}), \dots,
     \gentrivcofbar{n-1} = F(\gentrivcof{n-1}).
  \]
\end{paragraph}

\begin{paragraph}[The folk model structure on $\nCat{0}$, $\nCat{1}$, and $\nCat{2}$]
  The weak equivalence of the folk model structure on $\Set = \nCat{0}$ are
  the \oo-equivalences between sets, that is, the bijections. It follows
  from the discussion in the previous paragraph that its class of
  cofibrations is generated by
  \[ \varnothing \to \{x\} \qand \{x, y\} \to \{z\}. \]
  This implies that any map is cofibration. In particular, every set is
  cofibrant. Similarly, its class of trivial cofibrations is generated by an empty set
  of generators, showing that any map is fibration.

  The weak equivalences of the folk model structure on $\Cat = \nCat{1}$ are
  the \oo-equivalences between categories, that is, the equivalences of
  categories. Its class of cofibrations is generated by the obvious functors
  of the form
  \[ \varnothing \into \cdot, \qquad
    \{\xymatrix@C=0.8pc{\cdot & \cdot}\}
    \into \{\xymatrix@C=0.8pc{\cdot \ar[r] & \cdot}\},
    \quad
    \{
    \xymatrix@C=2pc{
      \cdot
      \ar@/^2ex/[r]_{}="0"
      \ar@/_2ex/[r]_{}="1"
      &
      \cdot
    }
    \}
    \to
    \{
    \xymatrix@C=1.5pc{
      \cdot
      \ar[r]
      &
      \cdot
    }
    \}.
  \]
  This implies that cofibrations are exactly the functors injective on
  objects. In particular, every category is cofibrant. Similarly, the class
  of trivial cofibrations is generated by $\gentrivcofbar{1}$. If one
  chooses $\globeinv{1}$ according to Remark~\ref{rem:choice_Jn}, one gets
  that $\gentrivcofbar{1}$ is the inclusion functor
  \[
    \{x\} \into \{\xymatrix@C=1pc{x \ar[r]^{\sim} & y}\},
  \]
  where the symbol $\sim$ denotes an isomorphism. This means that an
  \oo-functor $f : C \to D$ is a fibration in $\Cat$ if and only if it is
  an \emph{iso-fibration}, that is, if and only if for any object $x$ of $C$
  and any isomorphism $v : f(x) \to y$ of $D$, there exists an isomorphism
  $u : x \to x'$ in $C$ such that $f(u) = v$. We have thus recovered the
  classical folk model structure on $\Cat$, as described for instance in
  \cite{RezkFolk}.

  Let us now move on to $\nCat{2}$. The weak equivalences are the
  $2$-equivalences of $2$-categories. Its class of cofibrations is generated by
  the obvious $2$-functors of the form
  \[ \varnothing \into \cdot, \qquad
    \{\xymatrix@C=0.8pc{\cdot & \cdot}\}
    \into \{\xymatrix@C=0.8pc{\cdot \ar[r] & \cdot}\},
    \quad
    \{
    \xymatrix@C=2pc{
      \cdot
      \ar@/^2ex/[r]_{}="0"
      \ar@/_2ex/[r]_{}="1"
      &
      \cdot
    }
    \}
    \into
    \{
    \xymatrix@C=2pc{
      \cdot
      \ar@/^2ex/[r]_{}="0"
      \ar@/_2ex/[r]_{}="1"
      \ar@2"0";"1"
      &
      \cdot
    }
    \}.
  \]
  and
  \[
    \{
    \xymatrix@C=2pc{
      \cdot
      \ar@/^2ex/[r]_{}="0"
      \ar@/_2ex/[r]_{}="1"
      \ar@<-1ex>@2"0";"1"
      \ar@<1ex>@2"0";"1"
      &
      \cdot
    }
    \}
    \to
    \{
    \xymatrix@C=2pc{
      \cdot
      \ar@/^2ex/[r]_{}="0"
      \ar@/_2ex/[r]_{}="1"
      \ar@2"0";"1"
      &
      \cdot
    }
    \}.
  \]
  Its class of trivial cofibrations is generated by
  \[
    \{x\} \into F(\globeinv{1})
    \qand
    \{
    \xymatrix@C=2pc{
      \cdot
      \ar@/^2ex/[r]^u
      &
      \cdot
    }
    \}
    \into
    \{
    \xymatrix@C=2pc{
      \cdot
      \ar@/^2ex/[r]^u_{}="0"
      \ar@/_2ex/[r]_v_{}="1"
      \ar@2"0";"1"^(0.40)*[@]{\sim\mskip-4mu}
      &
      \cdot
    }
    \}
    \pbox.
  \]
  Choosing $\globeinv{1}$ according to Remark~\ref{rem:choice_Jn}, one gets
  that $F(\globeinv{1})$ is the free-standing adjoint equivalence: it has
  two objects $x$ and $y$, two generating $1$\nbd-cells $u : 0 \to 1$ and $v : 1
  \to 0$, two generating $2$-cells $\eta : \unit{x} \to uv$ and
  $\e : vu \to \unit{y}$ and these $2$-cells satisfy the
  triangular identities: $(\eta \comp{} u)(u \comp{} \e) = \unit{u}$ and $(v
  \comp{} \eta)(\e \comp{} v) = \unit{v}$. We have thus recovered the folk
  model structure introduced by Lack in \cite{LackFolk2}, with a correction
  in \cite{LackFolkBi}.
\end{paragraph}

\begin{paragraph}[The folk model structure on $\oGpd$]
  The category $\oGpd$, which is nothing but $\npCat{\omega}{0}$, is endowed
  with a folk model structure. The weak equivalences are the
  \oo-equivalences between \oo-groupoids. They can be characterized as the
  \oo-functors inducing bijections on connected components and homotopy
  groups but we will not enter into that.
  Denote by $F : \ooCat \to \oGpd$ the left adjoint to the
  inclusion functor. This functor sends an \oo-category $C$ to the
  \oo-groupoid obtained by formally inverting every cell of $C$. Set
  \[ \globegpd{n} = F(\globe{n}). \]
  The class of cofibrations of $\oGpd$ is generated by the inclusion
  \[ F(\gencof{n}) : \spheregpd{n} \into \globegpd{n}, \]
  where $\spheregpd{n}$ denotes the underlying $(n-1)$-groupoid of
  $\globegpd{n}$. Less trivially, the class of trivial cofibrations is
  generated by the set $\setgentrivcofgpd$ of \oo-functors
  \[ \gentrivcofgpd{n} : \globegpd{n} \into \globegpd{n+1}, \]
  where $n \ge 1$, corresponding to the source of the principal cell of
  $\globegpd{n+1}$. (Note that these \oo-functors are \emph{not} the
  $F(\gentrivcof{n})$.) To prove this, one can proceed as follows. First,
  one notes that these \oo-functors are trivial cofibrations, so that
  the class $\setgentrivcofgpd$ is included in the class of trivial
  cofibration. Second, one checks that an \oo-equivalence between \oo-groupoids
  having the right lifting property with respect to the $\gentrivcofgpd{n}$
  is a trivial fibration. Third, one concludes
  using a similar argument as in the proof of \cref{thm:folk}. Let
  $f$ be a trivial cofibration in~$\oGpd$. Using the small object argument,
  one can factor it as $f = pi$, where $i$ is in
  $\lrorth{\setgentrivcofgpd}$ and $p$ is in $\rorth{\setgentrivcofgpd}$. By
  the first point, $i$ is a trivial cofibration and so $p$ is a weak
  equivalence by the 2-out-of-3 property. This implies that $p$ is a trivial
  fibration by the second point. Thus, $f$ has the left lifting property
  with respect to $p$ and, by the retract lemma, $f$ is a retract of
  $i$, showing that $f$ is in~$\lrorth{\setgentrivcofgpd}$.

  We refer the reader to~\cite{AraMetWGrp} for more details on the folk
  model structure on $\oGpd$. In this paper, this model structure is called
  the Brown-Golasiński model structure as it coincides, through the
  equivalence of categories between \oo-groupoids and crossed complexes,
  with the model structure on crossed complexes introduced by Brown and
  Golasiński in \cite{BrownGolas}.
\end{paragraph}

\begin{paragraph}[The folk model structure on $\npCat{\omega}{1}$]
\label{P:DescriptionFolkModelStructureOmegaUnCat}
  The description of the folk model structure on $\npCat{\omega}{1}$ is very
  close to the one of $\oGpd$. Its weak equivalences are the
  \oo-equivalences between $(\omega, 1)$-categories. Setting
  \[ \globegpd{n} = F(\globe{n}), \]
  where $F : \ooCat \to \npCat{\omega}{1}$ denotes the left adjoint to the
  inclusion functor, the class of cofibrations of $\npCat{\omega}{1}$ is
  generated by the inclusions
  \[ F(\gencof{n}) : \spheregpd{n} \into \globegpd{n}, \]
  where $\spheregpd{n}$ denotes the underlying $(n-1, 1)$-category
  of $\globegpd{n}$. The class of trivial cofibrations can be shown to be
  generated by
  \[ \gentrivcofgpd{1} = F(\gentrivcof{1}) : \globegpd{0} \to
  F(\globeinv{1}) \]
  and by the
  \[ \gentrivcofgpd{n} : \globegpd{n} \into \globegpd{n+1}, \]
  for $n \ge 1$, corresponding to the source of the principal cell of
  $\globegpd{n+1}$.
\end{paragraph}


\chapter{Homology of \pdfoo-Categories}
\label{chap:homology}

This chapter is about Métayer's polygraphic homology of \pdfoo-categories
\cite{metayer2003resolutions}. This homology theory was first defined in the
following way: the polygraphic homology of an \oo-category is the homology
of the abelianization of any of its polygraphic replacements. Métayer then
showed with Lafont that for every monoid, considered
as an \oo-category, its polygraphic homology coincides with its classical homology as a monoid
\cite{lafont2009polygraphic}. This result was then generalized to
$1$-categories by Guetta \cite{GuettaHomologyCat}.

In this chapter, we prove that the polygraphic homology is the left derived
functor of a linearization functor from $\ooCat$ to the category $\pCh[\Z]$
of chain complexes, where $\ooCat$ is endowed with \oo-equivalences
and $\pCh[\Z]$ with quasi-isomorphisms.

In the first section of the chapter, we define the abelianization functor,
and we prove that an oplax transformation induces a chain homotopy after
abelianization. In the second section, we show that the abelianization
functor is a left Quillen functor for the folk model category structure
on $\ooCat$ and the projective structure on $\pCh[\Z]$. We define polygraphic
homology using this derived functor. The fact that polygraphic resolutions
are cofibrant resolutions for the folk model category structure proves
that this homology is indeed Métayer's polygraphic homology. In the third
section, we show that polygraphic homology generalizes monoid homology.
The proof we present is based on a conceptual reinterpretation of Lafont and
Métayer's proof by Guetta. Finally, in the last section, we compute some
examples of polygraphic homology.

In the whole chapter, we assume some familiarity with homological algebra,
and in particular homology of monoids. We refer the reader to
\cref{Chapter:ComplexesAndHomology} for a quick introduction.

\section{The Abelianization Functor}
\label{sec:abfunctor}

\begin{paragr}[Abelianization of an \pdfoo-category]
  \index{abelianization}
  We now define the \ndef{abelianization functor}
  \[
  \ab : \ooCat \to \pCh[\Z],
  \]
  where $\pCh[\Z]$ denotes the category of chain complexes of abelian groups
  in non-negative degree, see~\cref{paragr:def_ChR}. Using the notation of
  Section~\ref{sec:steiner}, this functor is just the composite
  \[ \ooCat \xto{\lambda} \ADC \xto{U} \pCh[\Z]. \]
  Nevertheless, we will give a direct definition.

  Let $C$ be an \oo-category. We define the chain complex $\ab(C)$ in the
  following way. For every $n \ge 0$, the group $\ab(C)_n$ is generated by
  elements $[x]$, for every $n$-cell~$x$ of $C$, subject to the
  relations $[x \comp{i} y] = [x] + [y]$, for every pair of $n$\nbd-cells~$x$
  and $y$ such that $x \comp{i} y$ is defined. It follows that if $z$ is an
  $m$-cell for~$m < n$, then $[\unit{z}] = 0$, where $\unit{z}$ denotes the
  iterated identity in dimension~$n$. For $n \ge 1$, we define $d_n$ by
  \[
    d_n([x]) = [\tge{n-1}(x)] - [\sce{n-1}(x)]
    \pbox,
  \]
  where $x$ is an $n$-cell of $C$. The axioms of \oo-categories giving the
  sources and targets of compositions imply that these maps are well-defined
  and the globular relations that they define a chain complex.

  If $f : C \to D$ is an \oo-functor, then we define a chain map
  \[ \ab(f) : \ab(C) \to \ab(D) \]
  by setting 
  \[ \ab(f)_n([x]) = [f(x)] \]
  for $n \ge 0$ and $x$ an $n$-cell of $C$. One checks that this is
  well-defined, that $\ab(f)$ is a chain map and that $\ab$ is indeed a
  functor.
\end{paragr}

\begin{proposition}\label{prop:ab_leftadj}
  The functor $\ab$ is a left adjoint. In particular, it preserves colimits.
\end{proposition}

\begin{proof}
  We saw in Remark~\ref{rem:adj_mu} that $\ab$, which is $U\lambda$ with the
  notation of this remark, admits the functor $\nu : \pCh[\Z] \to \ooCat$ of
  \cref{paragr:def_nu} as a left adjoint.
\end{proof}

\begin{prop}\label{prop:ab_pol}
  Let $P$ be a polygraph. For every $n \ge 0$, the abelian group~$\ab(P^\ast)_n$
  is free with basis the elements of the form $[x]$, where $x$ is in $P_n$.
\end{prop}

\begin{proof}
  This is part of Proposition~\ref{prop:lambda_pol}.
\end{proof}

\begin{proposition}\label{prop:ab_homot}
  Let $f, g : C \to D$ be two \oo-functors and let $\alpha : f \tod g$ be an
  oplax transformation. Then, setting
  \[ \ab(\alpha)_n([x]) = [\alpha_x] \]
  for $n \ge 0$ and $x$ in $C_n$, we obtain a chain homotopy
  (see~\cref{paragr:chain_homot})
  \[
    \ab(\alpha) : \ab(f) \tod \ab(g)
    \pbox.
  \]
\end{proposition}

\begin{proof}
  Let us first check that the map $\ab(\alpha)_n$ is well-defined. Let $x$
  and $y$ be two $n$\nbd-cells such that $x \comp{i} y$ is defined for some $i
  < n$. We have to check the equality $\ab(\alpha)_n([x \comp{i} y]) =
  \ab(\alpha)_n([x]) + \ab(\alpha)_n([y])$, that is,
  $[\alpha_{x \comp{i} y}] = [\alpha_x] + [\alpha_y]$. But this is exactly
  what we get by linearizing the formula
  \[
        \begin{split}
          \alpha_{x \comp{i} y} & = \big(f(\sce{i+1}(x)) \comp{0}
            \alpha_{\tge{0}(y)} \comp{1} \cdots \comp{i-1}
            \alpha_{\tge{i-1}(y)} \comp{i} \alpha_y\big) \comp{i+1} \\
            & \qquad\quad\qquad
            \big(\alpha_x \comp{i} \alpha_{\sce{i-1}(x)} \comp{i-1} \cdots
              \comp{1} \alpha_{\sce{0}(x)} \comp{0}
            g(\tge{i+1}(y))\big)\pbox,
        \end{split}
  \]
  keeping in mind that, as we are linearizing in dimension $n+1$, we have
  $[z] = 0$ for $z$ an $m$-cell with $m < n + 1$.

  Let us now prove that $\ab(\alpha)$ is indeed a chain homotopy. If $x$ is
  $0$-cell of~$C$, we have $\alpha_x : f(x) \to g(x)$ and so $d_1[\alpha_x] =
  [g(x)] - [f(x)]$, showing that $d_1(\ab(\alpha)_0([x])) = \ab(g)_0([x]) -
  \ab(f)_0([x])$. Similarly, if $x$ is an $n$-cell of~$C_n$ for $n > 0$, we have
  \[
    \alpha_x :
    f(x) \comp{0} \alpha_{\tge{0}(x)} \comp{1} \cdots \comp{n-1} \alpha_{\tge{n-1}(x)}
    \to
    \alpha_{\sce{n-1}(x)} \comp{n-1} \cdots \comp{1} \alpha_{\sce{0}(x)} \comp{0} g(x)
  \]
  and, by linearizing,
  \[
    d_{n+1}[\alpha_x] =
    [\alpha_{\sce{n-1}(x)}] + [g(x)]
    -[f(x)] - [\alpha_{\tge{n-1}(x)}]
    .
  \]
  As
  \[
    \begin{split}
    [\alpha_{\tge{n-1}(x)}] -
    [\alpha_{\sce{n-1}(x)}]
    & =
    \ab(\alpha)_n([\tge{n-1}(x)]) -
    \ab(\alpha)_n([\sce{n-1}(x)])
    \\
    & =
    \ab(\alpha)_n([\tge{n-1}(x)] - [\sce{n-1}(x)])
    \\
    & =
    \ab(\alpha)_n(d_n[x]),
    \end{split}
  \]
  we get
  \[
    d_{n+1}(\ab(\alpha)_n([x])) +
    \ab(\alpha)_{n-1}(d_n[x])
    = \ab(g)_n([x]) - \ab(f)_n([x]),
  \]
  thereby showing the result.
\end{proof}

\section{Deriving the Abelianization Functor}

In this section, we will show that the abelianization functor
\[ \ab : \ooCat \to \pCh[\Z] \]
can be derived, thus producing a homology theory for \oo-categories.

\medbreak

We start with some preliminaries on derived functors.

\begin{paragraph}[Homotopy category]
\index{localizer}
\index{relative!category}
\index{weak equivalence}
\index{homotopy!category}
  A \ndef{localizer} (also called \ndef{relative category}) is a category
  $\C$ endowed with a class $\clW$ of morphisms called \ndef{weak
  equivalences}. The \ndef{homotopy category} of a localizer $(\C, \clW)$ is
  the category $\loc{\C}{\clW}$ obtained from $\C$ by formally inverting
  arrows in $\clW$. We will also denote this category by $\Ho(\C)$, making
  implicit the class $\clW$. There is a canonical functor~$p : \C \to
  \Ho(\C)$.

  In particular, any model category $\M$ has an underlying localizer $(\M,
  \clW)$ and thus a homotopy category $\Ho(\M)$.
\end{paragraph}

\begin{paragraph}[Left derived functors]
  Let $(\C, \clW_\C)$ and $(\D, \clW_\D)$ be two localizers, and let
  $F : \C \to \D$ be a functor. The \ndef{(total) left derived
  functor}\index{left derived functor}\index{functor!left derived} of $F$, if
  it exists, is the universal pair consisting of a functor
  \[ \Lder F : \Ho(\C) \to \Ho(\D) \]
  and a natural transformation
  \[
  \xymatrix@C=2pc@R=2pc{
    \C\ar[d]_{p_\C}\ar[r]^F&\D\ar[d]^{p_\D}\\
    \Ho(\C)
    \ar@{.>}[r]_{\Lder{F}}
    \ar@{}[ur]_(0.35){}="a"_(0.65){}="b"
    \ar@2"a";"b"^{\lambda}
    &
    \Ho(\D)
    \pbox.
  }
  \]

  By abuse of language, one often refers to $\Lder F$ as the left derived
  functor of~$F$.
\end{paragraph}

\medbreak

One important use of model categories is to provide tools to prove the
existence of derived functors. In particular, the so-called left Quillen
functors can be left derived.

{
\renewcommand\N{\mathcal{N}}

\begin{paragraph}[Left Quillen functors]
\index{Quillen!functor}
\index{functor!Quillen}
  Let $\M$ and $\N$ be two model categories. A left adjoint functor $F : \M
  \to \N$ is called a \ndef{left Quillen functor}\index{left Quillen
  functor}\index{Quillen!functor!left} if it sends cofibrations of $\M$ to cofibrations of $\N$ and
  trivial cofibrations of $\M$ to trivial
  cofibrations of~$\N$.

  Similarly, a right adjoint functor $G : \N \to \M$ is said to be a
  \ndef{right Quillen functor}\index{right Quillen
  functor}\index{Quillen!functor!right} if it sends fibrations to fibrations
  and trivial fibrations to trivial fibrations.

  If
  \[ F : \M \to \N \qquad\qquad G : \N \to \M \]
  is a pair of adjoint functors, then $F$ is a left Quillen functor if and
  only if $G$ is a right Quillen functor. The pair $(F, G)$ is then called a
  \ndef{Quillen pair}\index{Quillen!pair} or a \ndef{Quillen
  adjunction}\index{Quillen!adjunction}\index{adjunction!Quillen}.
\end{paragraph}

\begin{theorem}[Quillen]\label{thm:Quillen_derived}
  A left Quillen functor $F : \M \to \N$ admits a left derived functor
  $\Lder F : \Ho(\M) \to \Ho(\N)$. Moreover, if $X$ is an object of $\M$,
  then $\Lder F(p_\M(X))$ is canonically isomorphic to $p_\N(F(Q))$, where
  $(Q, Q \to X)$ is a cofibrant replacement, \ie a cofibrant object $Q$
  endowed with a weak equivalence~$Q \to X$.

  Similarly, a right Quillen functor admits a right derived functor that can
  be computed using fibrant replacement.
\end{theorem}

\begin{remark}
  One actually only needs $F$ to send trivial cofibrations between
  cofibrant objects to weak equivalences for the previous theorem to apply.
\end{remark}
}

We will now apply Quillen's result to the abelianization functor
\[ \ab : \ooCat \to \pCh[\Z]. \]
We first introduce a model category structure on $\pCh[\Z]$.

\begin{paragr}[The projective model structure on chain complexes]
  \index{model structure!on chain complexes}
  \index{model structure!projective}
  \index{projective!model structure}
  \index{chain complex!model structure}
The category of chain complexes $\pCh[\Z]$ can be endowed with the so-called
\ndef{projective model structure} (see \cite[Chapter II, Section 4, page
11]{quillen1967homotopical}):

\begin{itemize}
  \item the weak equivalences are the quasi-isomorphisms
  (see~\cref{paragr:def_qis}),
  \item the cofibrations are the monomorphisms $f$ such that, for every $n
  \ge 0$, the cokernel of $f_n$ is projective,
  \item the fibrations are the morphisms $f$ such that, for every $n > 0$,
  $f_n$ is an epimorphism.
\end{itemize}
In particular, every chain complex is fibrant for this model structure and
the cofibrant chain complexes are exactly the ones that are projective in
every degree.
\end{paragr}

\begin{theorem}
  The functor $\ab : \ooCat \to \pCh[\Z]$ is a left Quillen functor, where
  $\ooCat$ is endowed with the folk model structure and $\pCh[\Z]$ is endowed
  with the projective model structure.
\end{theorem}

\begin{proof}
  We have to show that $\ab$ preserves cofibrations and trivial
  cofibrations. Since the cofibrations and trivial cofibrations of $\ooCat$
  are generated by sets~$\setgencof$ and $\setgentrivcof$, it suffices to
  prove that these sets are sent to cofibrations and weak equivalences,
  respectively.

  Let $\gencof{n} : \sphere{n} \to \globe{n}$ be an element of
  $\setgencof$. Proposition~\ref{prop:ab_pol} gives a concrete description
  of the morphism $\ab(\gencof{n}) : \ab(\sphere{n}) \to \ab(\globe{n})$:
  for $k$ such that~$0 \le k < n$, the morphism $\ab(\gencof{n})_k$ can be
  identified with the identity of $\Z^2$, for $k = n$, it can be identified
  with the unique morphism $0 \to \Z$ and, for $k > n$, with the identity of
  $0$. All these morphisms are monomorphisms of cokernel either $0$ or $\Z$.
  This means that $\ab(\gencof{n})$ is a cofibration.

  Let now $\gentrivcof{n} : \globe{n} \to \globeinv{n+1}$ be an element of
  $\setgentrivcof$. By Corollary~\ref{coro:w-eq_cof_imm}, this \oo-functor
  admits an inverse up to reversible oplax transformations.
  (Now that we have the folk model structure, this also
  follows from the so-called ``Whitehead Theorem'', see for
  instance~\cite[Theorem 7.5.10]{HirMC}, as $\gentrivcof{n}$ is a weak
  equivalence between cofibrant and fibrant objects.)
  But since by Proposition~\ref{prop:ab_homot}, oplax transformations are
  sent to chain homotopies, the morphism~$\ab(\gentrivcof{n})$ is a homotopy
  equivalence (see~\cref{paragr:def_homot_equiv}) and thus a quasi-isomor\-phism,
  thereby proving the result.
\end{proof}

\begin{paragraph}[Polygraphic homology of \pdfoo-categories]
  \index{polygraphic!homology}
  \index{homology!polygraphic}
  By the previous theorem, the functor
  \[ \ab : \ooCat \to \pCh[\Z] \]
  can be derived to a functor
  \[ \Lab : \Ho(\ooCat) \to \Ho(\pCh[\Z]). \]
  In particular, for $n \ge 0$, by post-composing by the functor
  \[ \HG_n : \Ho(\pCh[\Z]) \to \Ab \]
  (induced by the $H_n : \pCh[\Z] \to \Ab$ functor
  of~\cref{paragr:def_funct_Hn}),
  we get a functor
  \[ \HPol_n : \Ho(\ooCat) \to \Ab. \]
  If $C$ is an \oo-category, by definition, the \ndef{$n$-th polygraphic
  homology group} of $C$ is the abelian group $\HPol_n(C)$.

  Concretely, by \cref{thm:Quillen_derived}, the group $\HPol_n(C)$ is
  computed in the following way:
  \[
    \HPol_n(C) = \HG_n(\ab(P^\ast)),
  \]
  where $(P, P^\ast \to C)$ is any polygraphic resolution of $C$.
\end{paragraph}

\section{Comparison with Homology of Monoids}
\label{Section:ComparisonWithHomologyMonoids}

In the previous section, we defined the polygraphic homology of an
\oo-categ\-ory~$C$. In this section, we will show, following
\cite{lafont2009polygraphic}, that when $M$ is a monoid, the polygraphic
homology of $M$, seen as an \oo-category, coincides with the homology of the
monoid $M$. A similar result holds for $1$-categories, the homology of
$1$-categories being defined using simplicial sets (see
Section~\ref{sec:simpl-H}), as proved in~\cite{GuettaHomologyCat}.

\begin{thm}
  Let $M$ be a monoid. The polygraphic homology of $M$ seen as an
  \oo-category coincides with the homology of the monoid $M$.
\end{thm}

\begin{proof}
  Let $p : P^\ast \to M$ be a polygraphic resolution of $M$, considered as a
  $1$\nbd-category, such that $P^0$ consists of a unique element that we will
  denote by~$\ast$. For instance, one could take
  the standard resolution. By definition, the polygraphic homology of $M$ is
  the homology of the abelianization of $P^\ast$. We will now associate to
  this polygraphic resolution a resolution by free left $\Z M$\nbd-modules of~$\Z$
  (endowed with the trivial action).

  Denote by $e$ the unit element of $M$ and by $\coslice{M}{e}$ the coslice
  category. Explicitly, an object of $\coslice{M}{e}$ is an element of $M$,
  and if $m$ and $m'$ are two objects, a morphism from $m$ to $m'$ is an
  element $n$ of $M$ such that $nm = m'$. Denote by~$\coslice{P^\ast}{e}$ the
  \oo-category defined by the pullback square
  \[
    \xymatrix{
      \coslice{P^\ast}{e} \ar[r] \ar[d]_q & P^\ast \ar[d]^p \\
      \coslice{M}{e} \ar[r] & M \pbox{,}
    }
  \]
  where the bottom horizontal morphism is the forgetful functor.
  Let us describe the cells of this \oo-category. An object of
  $\coslice{P^\ast}{e}$ consists of an element $m$ of $M$ that we will
  denote by $(m, \ast)$. A $1$-cell from
  an object $m$ to an object $m'$ consists of a $1$-cell $x$ such that
  $p(x)m = m'$. A $1$-cell is thus uniquely determined by a pair $(m, x)$,
  where $m$ is in $M$ and $x$ is a $1$-cell of $P^\ast$. The source of such
  a $1$-cell is $p(x)m$ and its target is $m'$. Similarly, an $i$-cell for
  $i > 1$ can be described as a pair $(m, x)$, where $m$ is in $M$ and $x$
  is an $i$-cell of $P^\ast$. Its source is $(m, \src{i-1}(x))$ and its target
  $(m, \tgt{i-1}(x))$.

  One can show that the \oo-category $\coslice{P^\ast}{e}$ is freely
  generated in the sense of polygraphs by its cells of the form $(m,
  a)$, where $m$ is in $M$ and $a$ is in $P_i$ for~$i \ge 0$. This
  implies that the linearization of this \oo-category in degree $i$ is
  \[ \ab(\coslice{P^\ast}{e})_i \isoto \freemod{\Z}{M \times P_i} \isoto
  \freemod{\Z M}{P_i}\pbox. \]
  In particular, it has a natural structure of left $\Z M$-module. If $x$
  is a cell, setting $[x] = [(e, x)]$, we have $[(m, x)] = m[x]$.
  By definition, if $m$ is in $M$ and $x$ is an $i$-cell for $i > 0$,
  we have
  \[
    d_i(m[x]) =
    \begin{cases}
      m[\ast] - mp(x)[\ast] & \text{if $i = 1$,} \\
    m[\tgt{i-1}(x)] - m[\src{i-1}(x)] & \text{if $i > 1$,}
    \end{cases}
  \]
  and the map $d_i$ is thus $\Z M$-linear. This shows that
  $\ab(\coslice{P^\ast}{e})$ is a complex of free left $\Z M$-modules.

  The unique \oo-functor from $\coslice{P^\ast}{e}$ to the terminal
  \oo-category induces a morphism of complexes $\ab(\coslice{P^\ast}{e})$ to
  $\Z$. In other words, by sending $m[\ast]$ in $\ab(\coslice{P^\ast}{e})_0$
  to $1$ in $\Z$, we get an augmented complex of left $\Z M$-modules
  \[
    \xymatrix{
      \Z & \ar[l] \freemod{\Z M}{P_0} & \ar[l] \freemod{\Z M}{P_1} & \ar[l]
      \freemod{\Z M}{P_2} & \ar[l] \cdots \pbox{.}
    }
  \]
  We will see that this complex is exact. Assuming this, let
  us end the proof. As the complex is exact, we have built, as announced, a
  resolution of $\Z$ by free left $\Z M$-module. The homology of the monoid
  $M$ is thus the homology of the complex
  \[
    \xymatrix{
      \Z \otimes_{\Z M} \freemod{\Z M}{P_0} & \ar[l] \Z \otimes_{\Z M}
      \freemod{\Z M}{P_1} & \ar[l] \Z \otimes_{\Z M} \freemod{\Z M}{P_2} &
      \ar[l] \cdots \pbox{,}
    }
  \]
  which is canonically isomorphic to the complex
  \[
    \xymatrix{
     \freemod{\Z}{P_0} & \ar[l] \freemod{\Z}{P_1} & \ar[l] \freemod{\Z}{P_2}
     & \ar[l] \cdots \pbox{,}
    }
  \]
  which is nothing but $\ab(P^\ast)$. The homology of the monoid $M$ is thus
  the homology of $\ab(P^\ast)$, that is, the polygraphic homology of $M$.

  To end the proof, it thus suffices to show that the augmented complex
  introduced at the beginning of the previous paragraph is exact. We have a
  unique \oo-functor $r : \coslice{P^\ast}{e} \to \termobj{}$, where
  $\termobj{}$ denotes the terminal \oo-category. We have to show that
  $\ab(r)$ is a quasi-isomorphism. Consider the \oo-functor $e : \termobj{}
  \to \coslice{P^\ast}{e}$ corresponding to the object $(e, \ast)$ of
  $\coslice{P^\ast}{e}$. The composition $re$ is the identity of
  $\termobj{}$. We will construct an oplax transformation $\alpha : er \tod
  \id_{\coslice{P^\ast}{e}}$. Using Proposition~\ref{prop:ab_homot}, we will
  get that $\ab(r)$ is a homotopy equivalence and thus a quasi-isomorphism.
  Let us construct this oplax transformation $\alpha$. First, note that there
  exists an oplax transformation $\beta : qer \tod q$, where $q :
  \coslice{P^\ast}{e} \to \coslice{M}{e}$ is the \oo-functor introduced at
  the beginning of the proof. (Note that $qer$ is the constant \oo-functor of
  value $e$.) Indeed, as $\coslice{M}{e}$ admits $e$ as an initial object, we
  have a natural transformation $\gamma : e \tod \id_{\coslice{M}{e}}$,
  where $e$ denotes the constant endofunctor of $\coslice{M}{e}$ of value
  $e$. The oplax transformation $\beta$ is thus $\gamma \comp{} q$. Second,
  note that $q$ is a trivial fibration, as it is defined by pulling back the
  trivial fibration $p$. The existence of the transformation $\alpha$ thus
  follows from the following lemma that concludes the proof.
\end{proof}

\begin{lemma}
  Let $p : C \to D$ and $f, g : B \to C$ be three \oo-functors and let $\beta : pf
  \tod pg$ be an oplax transformation. If $p$ is a trivial fibration and $B$
  is a cofibrant \oo-category, then there exists an oplax transformation
  $\alpha : f \tod g$ such that~$\beta = p \comp{} \alpha$.
\end{lemma}

\begin{proof}
  The oplax transformation $\beta$ corresponds to an \oo-functor $B \to
  \CatCyl(D)$. Consider the pullback square
  \[
    \xymatrix{
      \CatCyl(D) \times_{D \times D} C \times C \ar[r] \ar[d]
      & \CatCyl(D) \ar[d]^{(\projcyltop, \projcylbot)} \\
      C \times C \ar[r]_{p \times p} & D \times D \pbox{.}
    }
  \]
  As the source of $\beta$ is $pf$ and its target is $pg$, we get an
  \oo-functor
  \[  (\beta, f \times g) : B \to \CatCyl(D) \times_{D \times D} C \times C\pbox{.} \]
  The naturality square
  \[
    \xymatrix{
      \CatCyl(C) \ar[r]^{\CatCyl(p)} \ar[d]_{(\projcyltop, \projcylbot)} &
      \CatCyl(D) \ar[d]^{(\projcyltop, \projcylbot)} \\
      C \times C \ar[r]_{p \times p} & D \times D
    }
  \]
  induces an \oo-functor
  \[
    q : \CatCyl(C) \to \CatCyl(D) \times_{D \times D} C \times C \pbox{.}
  \]
  An oplax transformation $\alpha$ as in the statement exactly corresponds to
  a dotted arrow making the triangle
  \[
    \xymatrix@C=3pc{
      & \CatCyl(C) \ar[d]^q \\
    B \ar@{.>}[ur] \ar[r]_-{(\beta, f \times g)} & \CatCyl(D) \times_{D \times D} C \times C
    }
  \]
  commute. As $B$ is cofibrant, to get the result it suffices to show that
  $q$ is a trivial fibration.

  Let us prove this:
  \begin{enumerate}
    \item An object of $\CatCyl(D) \times_{D \times D} C \times C$ is given
    by two objects $x$ and $x'$ of $C$ and a $1$-cell $v : p(x) \to p(x')$.
    As $p$ is a trivial fibration, there exists $u : x \to x'$ such that
    $p(u) = v$, showing that $q$ is surjective on objects.
    \item Let $n \ge 1$ and let $\gamma : x \cylto y$ and $\delta : z \cylto
    t$ be two parallel $n$-cylinders of~$C$. Suppose we have an $(n+1)$-cell
    from $q(\gamma)$ to $q(\delta)$. This means that we have
    two $n$-cells $u : x \to z$ and $v : y \to t$ of $C$ and an
    $(n+1)$\nbd-cylinder $\Gamma : p(\gamma) \to p(\delta) : p(u) \cylto p(v)$.
    As $p$ is a trivial fibration, the $(n+2)$\nbd-cell of the $(n+1)$-cylinder $\Gamma$
    can be lifted to $C$ yielding an $(n+1)$-cylinder $\Lambda : \gamma \to
    \delta : u \cylto v$ of $C$ such that $p(\Lambda) = \Gamma$. This shows
    that $q$ is a trivial fibration, thereby proving the lemma.
    \qedhere
  \end{enumerate}
\end{proof}

\section{Examples}
\subsection{Polygraphic homology of $\Rev{1}$}
\label{paragr:homology_R1}

We will compute the polygraphic homology of the free-standing reversible
cell $\Rev{1}$, introduced in \cref{paragr:def_R1}. We will see that its
homology is non-trivial. This implies that $\Rev{1}$ is not weakly
contractible.

Since $\Rev{1}$ is freely generated by a polygraph, its homology is the homology
of its linearization $C = \Ab(\Rev{1})$. We have
\[ C_0 = \freemod{\Z}{0, 1} \]
and
\[
  C_j
    = \freemod{\Z}{\{ r_{l_1, \dots, l_{j-1}}, \weakinv{r}_{l_1, \dots,
    l_{j-1}} \mid l_k = \pm, 1 \le k < j \}},
\]
with
\[ d_1([r]) = [\tge{0}(r)] - [\sce{0}(r)] = [1] - [0], \]
and similarly,
\[ d_1([\weakinv{r}]) = [0] - [1] = -d([r]), \]
\[
  \begin{split}
   d_j([r_{l_1, \dots, l_{j-2}, -}])
    & = [\tge{j-1}(r_{l_1, \dots, l_{j-2}, -})] - [\sce{j-1}(r_{l_1, \dots,
    l_{j-1}},-)] \\
    & = [r_{l_1, \dots, l_{j-2}} \comp{j-2} \weakinv{r}_{l_1, \dots,
    l_{j-2}}] - [\unit{\sce{j-2}(r_{l_1, \dots, l_{j-2}})}] \\
    & = [r_{l_1, \dots, l_{j-2}}] + [\weakinv{r}_{l_1, \dots, l_{j-2}}],
  \end{split}
\]
and similarly,
\[
  \begin{split}
  d_j([\weakinv{r}_{l_1, \dots, l_{j-2}, -}])
    & = -d_j([r_{l_1, \dots, l_{j-2}, -}]),
    \\
  d_j([r_{l_1, \dots, l_{j-2}, +}])
    & = -d_j([r_{l_1, \dots, l_{j-2}, -}]),
    \\
  d_j([\weakinv{r}_{l_1, \dots, l_{j-2}, +}])
    & = d_j([r_{l_1, \dots, l_{j-2}, -}]).
  \end{split}
\]
Let us compute the homology of this chain complex. One has
\[ \HPol_0(\Rev{1}) = \Z[0,1]/\lspan{[1] - [0]} \simeq \Z. \]
A $1$-chain $a[r] + b[\weakinv{r}]$ is a cycle if and only if
\[ d_1(a[r] + b[\weakinv{r}]) = (a - b)[r] = 0, \]
that is, if and only if $a = b$. But $[r] + [\weakinv{r}]$ is a boundary.
Hence
\[ \HPol_1(\Rev{1}) = 0. \]
Similarly, a $2$-chain $a[r_-] + b[\weakinv{r}_-] + c [r_+] +
d[\weakinv{r}_+]$ is a cycle if and only if one has $a - b - c + d = 0$. The
abelian group of $2$-cycles is thus isomorphic to~$\Z^3$. As for the abelian
group of $2$-boundaries, it is spanned by $[r_-] + [\weakinv{r}_-]$ and
$[r_+] + [\weakinv{r}_+]$. It is thus isomorphic to $\Z^2$. Computing the
quotient, one gets that
\[ \HPol_2(\Rev{1}) \simeq \Z. \]
This already proves that $\Rev{1}$ is not weakly contractible.
More generally, for~$n \ge 2$, one checks that
\[ \HPol_n(\Rev{1}) \simeq \Z^{(2^n - 2^{n-2}) - 2^{n-1}} = \Z^{2^{n-2}}. \]

\subsection{Thomason homology and the homology of $K(\N,2)$}

In this chapter, we defined a homology theory for \oo-categories: the
polygraphic homology. There is a second homology theory for \oo-categories,
that we will call the Thomason homology, which can be morally defined in the
following way. Let $C$ be an \oo-category. Consider the weak \oo-groupoid
obtained by weakly inverting all the cells of $C$. This weak \oo-groupoid
corresponds to a homotopy type and, morally, the Thomason homology of $C$ is
the homology of this homotopy type. One way to give a precise definition of
this homology is to use Street nerve $N : \ooCat \to \SSet$, introduced in
\cite{street1987algebra}, which extends the usual nerve functor
(see~\cref{paragr:nerve}) to \oo-categories. If $C$ is an \oo-category, for
every $n \ge 0$, we set
\[ \HThom_n(C) = \HG_n(NC). \]
This is the \ndef{Thomason homology}\index{Thomason homology} of $C$. When
$C$ is a category, we recover the classical homology of categories
(generalizing the classical homology of monoids) recalled in
Section~\ref{sec:simpl-H}.  The main result of~\cite{GuettaHomologyCat} says
that the polygraphic homology and the Thomason homology of a category agree.
It is tempting to think that these two homologies always agree. This is not
the case!

Here is a counter-example. Recall (see~Section~\ref{sec:basic_ex}) that
to any abelian monoid $M$, one can associate a $2$-category with only one
object, one $1$-cell (the identity of the unique object), and $M$ as the set
of $2$-cells. Let us denote this $2$-category by $K(M, 2)$.
By~\cite[Theorem 4.7]{AraThmB}, if $M$ is an abelian group $A$, then the
Thomason homology of $K(A, 2)$ is the homology of the corresponding
Eilenberg-Mac Lane space, that is, of any CW-complex whose homotopy groups
are trivial except the second one which is $A$. In particular, by classical
results, the Thomason homology of $K(\Z, 2)$ is the homology of
$\mathbb{C}\mathrm{P}^{\infty}$, the infinite-dimensional complex projective
space, and is hence $\Z$ in even degree and null in odd degree. Moreover,
by~\cite[Theorem 4.9]{AraThmB}, the inclusion $2$-functor $K(\N, 2) \hookto
K(\Z, 2)$ induces an isomorphism in homology. We thus have
\[
\HThom_n(K(\N, 2))
=
\begin{cases}
  \Z & \text{if $n$ is even,} \\
  0 & \text{if $n$ is odd.}
\end{cases}
\]
But $K(\N, 2)$ is freely generated by the unique $2$-polygraph $P$ with
\[ P_0 = \{\star\}, \quad P_1 = \emptyset, \quad P_2 = \{ \alpha \}. \]
Its polygraphic homology is hence the homology of its linearization
  \[
    \xymatrix{
      \Z & \ar[l] 0 & \Z \ar[l] & 0 \ar[l] & \ar[l] \cdots
    }
  \]
and we have
\[
\HPol_n(K(\N, 2))
=
\begin{cases}
  \Z & \text{if $n = 0$ or $n = 2$,} \\
  0 & \text{otherwise.}
\end{cases}
\]
In particular, $\HPol_4(K(\N, 2)) \not\simeq \HThom_4(K(\N, 2))$.



\newcommand{\Er}{\mathcal{E}}
\newcommand{\Srr}{\mathcal{S}}

\newcommand{\polygraphe}{\mathrm}
\newcommand{\AsPoly}{\polygraphe{As}}
\newcommand{\Sq}{\polygraphe{Sq}}

\chapter{Resolutions by \texorpdfstring{$(\omega,1)$}{(omega,1)}-Polygraphs}
\label{Chapter:ConstructingResolutions}
\label{chap:ConstructingResolutions}
\label{chap:anick}
Anick and Green constructed the first explicit free resolutions for algebras from a presentation of relations by non-commutative Gröbner bases~\cite{Anick85, anick1986homology, AnickGreen87, Green99}.
Their constructions provide resolutions to compute homological invariants, such as homology groups, Hilbert and Poincaré series of algebras presented by generators and relations given by a Gröbner basis. The chains of these resolutions are defined by iterated overlaps of the leading terms of the Gröbner basis and the differentials are constructed by noetherian induction. Similar methods for calculating free resolutions for monoids and algebras, inspired by string rewriting mechanisms, have been developed in numerous works~\cite{brown1992geometry,Groves90,kobayashi1990complete,kobayashi2005grobner}.
A purely polygraphic approach to the construction of these resolutions by rewriting has been developed in~\cite{GuiraudMalbos12advances} using the notion of $(\omega,1)$\nbd-poly\-graphic resolution, where the mechanism for proving the acyclicity of the resolution relies on the construction of a normalization strategy extended in all dimensions. The construction of polygraphic resolutions by rewriting has also been applied to the case of associative algebras in~\cite{GuiraudHoffbeckMalbos19} and shuffle operads in~\cite{MalbosRen23}, introducing in each case a notion of polygraph adapted to the algebraic structure.

In this chapter, we show how to construct a polygraphic resolution of a category from a convergent presentation of that category and how to deduce an abelian version of such a resolution.
The notion of polygraphic resolution of an $\omega$\nbd-category was introduced in \cref{S:PolygraphicResolutionsOmegaCategories}: it consists of a polygraph which is weakly equivalent to the category.
We consider here a variant of this notion adapted to
$(\omega,1)$-categories, related to the folk model structure on the category~$\npCat{\omega}{1}$ of $(\omega,1)$-categories.
The chapter is organized as follows. 
In \cref{S:ResolutionsContractions}, we introduce the notion of contraction with respect to a unital section, which we often use to show that an $(\omega,1)$\nbd-poly\-graph is acyclic.
In \cref{S:SquierResolution}, we show how to compute a cofibrant replacement of a category in the category $\npCat{\omega}{1}$ from one of its convergent presentations. This construction extends to higher dimensions the one given in low dimensions in \cref{chap:2-Coherent} in terms of coherent presentations. 
In \cref{S:Abelianization}, we explain how to deduce an abelian resolution
from a resolution by an $(\omega,1)$-polygraph, thus again extending the constructions given in low dimensions in \cref{chap:2-homology}.
We deduce from this resolution several homological and homotopical finiteness conditions for finite convergence.
In \cref{S:FinitenessConditions}, we extend the results about finite
derivation type presented in \cref{chap:2-fdt}.

\section{Polygraphic Resolutions and Contractions}
\label{S:ResolutionsContractions}
In this section, we consider the folk model structure on $\npCat{\omega}{1}$ constructed in \cref{T:FolkModelStructureCatNP} and described in \cref{P:DescriptionFolkModelStructureOmegaUnCat}, and the notion of oplax transformation between $\omega$-functors as defined in \cref{paragr:oo-trans_alg}. 

\subsection{Polygraphic resolution of an $(\omega,1)$-category}
\index{polygraphic!resolution!of an $(\omega,1)$-category}
A \emph{polygraphic resolution of a category~$C$ in~$\npCat{\omega}{1}$} is a pair $(P,p)$ made of an $(\omega,1)$-polygraph~$P$ and a trivial fibration $p : \freegpd{P} \to C$, where~$\freegpd{P}$ is the free $(\omega,1)$-category generated by $P$.
Expanding the definition, $P$ is a polygraphic resolution of~$C$ if and only if it presents~$C$ and, for every~$n\geq 2$, the extension~$P_{n+1}$ of~$\freegpd{\tpol nP}$ is acyclic.

\subsection{Unital sections and essential cells}
Let~$P$ be an $(\omega,1)$-polygraph. For~$u$ a $1$-cell of the quotient category~$\cl{P}$, we denote by $\freegpd{P}_u$ the corresponding fiber of the canonical projection $\pi:\freegpd{P}\twoheadrightarrow\cl{P}$. By definition, $\freegpd{P}_u$ is an $\omega$-groupoid, whose $0$-cells are the representatives of~$u$ in~$\freegpd{P}$. 
To avoid confusion, we keep the dimensions of the $(\omega,1)$-category~$\freegpd{P}$ when talking about the cells and compositions of~$\freegpd{P}_u$.

A \emph{unital section of~$P$}\index{unital section} is a family
\[
  \iota=
  \left(\iota_u : \termobj{}\fl\freegpd{P}_u\right)_{u\in\cl{P}}
\]
of $\omega$-functors, satisfying $\iota_{1_x} = \unit{\unit{x}}$ for every $0$\nbd-generator~$x$ of~$P$. 
Such a family of functors assigns to every $1$-cell~$u$ of~$\cl{P}$ a representative $1$-cell~$\iota_u$ in~$\freegpd{P}$, in such a way that identities are mapped to identities.
A unital section of~$P$ is almost a functorial section of the canonical projection $\pi:\freegpd{P}\fl\cl{P}$, except that it is not defined in dimension~$0$ and no specific compatibility with the $0$-composition is required.

Fix a unital section~$\iota$ of~$P$. If~$\phi$ is an $n$-cell of~$\freegpd{P}$, we will write~$\rep{\phi}$ for $\iota\pi(\phi)$ when no confusion occurs. Note that~$\rep{\phi}$ is an identity if~$n\geq 2$. A $1$-cell~$u$ of~$\freegpd{P}$ is \emph{$\iota$-reduced} if $u=\rep{u}$ holds. A non-$\iota$-reduced $1$-cell~$u$ of~$\freegpd{P}$ is \emph{$\iota$-essential} if~$u=av$, with~$a$ a $1$-generator of~$P$ and~$v$ an $\iota$-reduced $1$-cell of~$\freegpd{P}$.

\subsection{Contractions}
\index{contraction}
\label{SS:Contractions}
Let~$P$ be an $(\omega,1)$-polygraph, and~$\iota$ be a unital section of~$P$. An \emph{$\iota$-contraction of~$P$} is a family
\[
  \sigma=
\left(
\vcenter{\xymatrix @R=1em @C=1em {
{\freegpd{P}_u}
	\ar@/^/ [rr] ^-{\id_{\freegpd{P}_u}} _-{}="src"
	\ar@/_/ [dr] _-{\varepsilon}
&& {\freegpd{P}_u}
\\
& \termobj{}
	\ar@/_/ [ur] _-{\iota}
	\ar@2 "src"!<-0.5ex,-2.2ex>;[]!<-0.5ex,2.2ex> ^-*+{\sigma_u}
}}
\right)_{u\in\cl{P}}
\]
of oplax transformations such that $\sigma_{\sigma_\phi}=\id_{\sigma_\phi}$ and $\sigma_{\iota_u}=\id_{\iota_u}$ for every cell~$\phi$ in~$\freegpd{P}$ and $1$-cell $u$ in $\pcat{P}$, where $\sigma_\psi$ is a short notation for $(\sigma_{\rep\psi})_\psi$.
An $\iota$-contraction is thus almost an oplax transformation from~$\id_{\freegpd{P}}$ to~$\iota\varepsilon$, but, like~$\iota\varepsilon$, it is not defined on $0$-cells and no specific compatibility with the $0$-composition is required.

Fix an $\iota$-contraction~$\sigma$ of~$P$. By definition of~$\sigma$, for all~$n\geq 1$, $n$-cell~$\phi$ of~$\freegpd{P}$, and~$1\leq k<n$,
\[
\src k(\sigma_\phi) = \phi\comp1\sigma_{t_1(\phi)}\comp2\cdots\comp k\sigma_{t_k(\phi)}
\quad\text{and}\quad
\tgt k(\sigma_\phi) = 
\begin{cases}
\rep{\phi} 
	&\text{if~$k=1$,} \\
\sigma_{s_k(\phi)}
	&\text{otherwise.}
\end{cases}
\]

An $n$-cell~$\phi$ of~$\freegpd{P}$ is \emph{$\sigma$-reduced} if it is an identity or in the image of~$\sigma$.

\subsection{Sided contractions}
\label{SS:SidedContractions}

We say that an $\iota$-contraction~$\sigma$ is \emph{right} if, for all~$n\geq 1$ and $n$-cells~$\phi$, $\psi$ of~$\freegpd{P}$ of respective $1$-sources~$u$ and~$v$, it satisfies
\begin{equation}
\label{E:RightContraction}
\sigma_{\phi\psi} = u\sigma_\psi \comp1 \sigma_{\phi\rep{v}}.
\end{equation}
Symmetrically, an $\iota$-contraction is \emph{left} if for all~$n\geq 1$ and $n$-cells~$\phi$, $\psi$ of~$\freegpd{P}$ of respective $1$-sources~$u$ and~$v$, it satisfies
\begin{equation}
\label{E:LeftContraction}
\sigma_{\phi\psi} = \sigma_\phi v \comp1 \sigma_{\rep{u}\psi}.
\end{equation}
In the sequel, we will consider right $\iota$-contractions, however, the definitions and results admit a left version.

If~$\sigma$ is a right $\iota$-contraction of~$P$, and~$n\geq 1$, a non-$\sigma$-reduced $n$-cell~$\phi$ of~$\freegpd{P}$ is \emph{$\sigma$-essential} if there exists an $n$-cell~$\alpha$ of~$P$ and an $\iota$-reduced $1$-generator~$u$ of~$\freegpd{P}$ such that $\phi=\alpha u$.

\begin{lemma}[{\cite[Corollary~3.3.5]{GuiraudMalbos12advances}}]
\label{L:FreeRightContraction}
Let~$P$ be an $(\omega,1)$-polygraph, and~$\iota$ be a unital section of~$P$. A right $\iota$-contraction~$\sigma$ of~$P$ is uniquely and entirely determined by its values on the $\iota$-essential $1$-cells of~$\freegpd{P}$ and, for every~$n\geq 1$, on the $\sigma$-essential $n$-cells of~$\freegpd{P}$. 
\end{lemma}
\begin{proof}
If~$\sigma$ is a right $\iota$-contraction, then its values are prescribed on every cell of~$\freegpd{P}$ that is not $\iota$-essential or $\sigma$-essential. Now, the values of~$\sigma$ on $\iota$-essential and $\sigma$-essential cells of~$\freegpd{P}$ can be chosen freely (with correct source and target), provided that these values make~$\sigma$ compatible with all the defining relations of the structure of $(\omega,1)$-category, and in particular with exchange relations between the $0$-composition and the other compositions. It turns out that~\eqref{E:RightContraction} imposes compatibility with these exchange relations.
\end{proof}

\begin{theorem}
  \label{T:ResolutionContraction}
Let~$P$ be an $(\omega,1)$-polygraph, and~$\iota$ be a unital section of~$P$. The canonical projection $\pi : \freegpd{P}\fl\cl{P}$ is a trivial fibration in $\npCat{\omega}{1}$ if and only if $P$ admits a right $\iota$-contraction.
\end{theorem}
\begin{proof}
Assume that $\pi : \freegpd{P}\fl\cl{P}$ is a trivial fibration. Let us define a right $\iota$-contraction~$\sigma$ of~$P$ thanks to \cref{L:FreeRightContraction}. If~$au$ is an $\iota$-essential $1$-cell of the free $(\omega,1)$-category~$\freegpd{P}$, then~$\pi(au)=\pi(\rep{au})$, so that, by definition of~$\cl{P}$, there exists a $2$-generator
\[
\sigma_{au}:au\fl\rep{au} 
\]
in~$\freegpd{P}$. Assume that~$\sigma$ is defined on the $n$-cells of~$\freegpd{P}$, for~$n\geq 1$, and let~$\alpha u$ be a $\sigma$-essential $n$-cell of~$\freegpd{P}$. The $n$-cells $\src{}(\sigma_{\alpha u})$ and $\tgt{}(\sigma_{\alpha u})$ are parallel, so, by hypothesis, there exists an $(n{+}1)$-cell
\[
\sigma_{\alpha u} : \src{}(\sigma_{\alpha u}) \fl \tgt{}(\sigma_{\alpha u})
\]
in~$\freegpd{P}$.

Conversely, let~$\sigma$ be an $\iota$-contraction of~$P$, and~$\phi$, $\psi$ be parallel $n$-cells of~$\freegpd{P}$, for~$n\geq 1$. We have $\tgt{}(\sigma_\phi)=\sigma_{\src{}(\phi)}=\sigma_{\src{}(\psi)}=\tgt{}(\sigma_\psi)$ by hypothesis, so that the $(n{+}1)$-cell $\sigma_\phi\comp n\sigma_\psi^-$ is well defined, with source $\src{}(\sigma_\phi)$ and target $\src{}(\sigma_\psi)$. The fact that $\tgt k(\phi)=\tgt k(\psi)$ holds for every $0\leq k<n$ implies that
\[
(\sigma_\phi \comp n \sigma_\psi)^- \comp{n-1} \sigma_{t_{n-1}(\phi)}^- \comp{n-2} \cdots \comp0 \sigma_{t_0(\phi)}^-
\]
is a well-defined $(n{+}1)$-cell of~$\freegpd{P}$, with source~$\phi$ and target~$\psi$, thus proving that~$P_{n+1}$ is acyclic.
\end{proof}

This theorem shows that to prove that an $(\omega,1)$-polygraph $P$ is a polygraphic resolution of the category $\cl{P}$, it suffices to provide it with a right $\iota$-contraction.
In the next section, we show how to construct such a contraction from a convergent presentation of the category $\cl{P}$.

\section{Polygraphic Resolutions from Convergence}
\label{S:SquierResolution}

\subsection{The cells of the Squier polygraphic resolution}

Assume that~$P$ is a $2$\nbd-poly\-graph. Define~$\Sq(P)$ as the graded set $(\Sq_n(P))_{n\geq 0}$ where
\begin{enumerate}
\item $\Sq_0(P) = P_0$ and $\Sq_1(P)=P_1$, 
\item for~$n\geq 2$, $\Sq_n(P)$ is the set of tuples $(u_1,\dots,u_n)$, written $u_1\vert\cdots\vert u_n$, of non-identity reduced $1$-cells of~$P^*$ such that
\begin{itemize}
\item $u_1$ is a $1$-generator of~$P$,
\item for every $1\leq i<n$, the $1$-cell $u_i u_{i+1}$ is not reduced, 
\item for every $1\leq i<n$, every proper left-factor of $u_i u_{i+1}$ is reduced.
\end{itemize}
\end{enumerate}

\subsection{Interpretation in the reduced case}
\label{SS:InterpretationReducedCase1}

Assume that~$P$ is a reduced $2$\nbd-poly\-graph. Then~$u_1\vert u_2$ is a $2$-generator of~$\Sq(P)$ if and only if $u_1$ is a $1$\nbd-gene\-rator of~$P$ and $u_1u_2$ is the source of a $2$-generator of~$P$. 

From the classification of critical branchings of a $2$-polygraph given in \cref{sec:cb-class}, the critical branchings of the reduced polygraph~$P$ are of the form 
\[
\xymatrix @C=3em {
\strut 
	\ar [r] _-{u_1}
	\ar@/^6ex/ [rrr] _-{}="t1"
& \strut
	\ar [r] _{v_2} ^{}="s1"
& \strut
	\ar [r] |{w_2} 
	\ar@/_6ex/ [rr] ^{}="t2"
& \strut
	\ar[r] ^-{u_3} _(0.25){}="s2"
& \strut
\ar@2 "s1"!<0pt,7.5pt>;"t1"!<0pt,-7.5pt> ^-{\alpha}
\ar@2 "1,4";"t2"!<0pt,7.5pt> ^-{\beta}
}
\]
where~$\alpha$ and~$\beta$ are $2$-generators of~$P$, and~$u_1v_2$, $w_2$, and~$u_3$ are reduced non-identity $1$-cells of~$P^*$, with~$u_1$ a $1$-generator of~$P$. Putting $u_2=v_2w_2$ induces a one-to-one correspondence between the $3$-generators $u_1\vert u_2\vert u_3$ of~$\Sq_3(P)$ and the critical branchings of~$P$ whose source is $u_1u_2u_3$.

For~$n\geq 3$, define the \emph{critical $n$-branchings of~$P$}\index{critical!branching!8-@$n$-} as the non-ordered families $(\alpha_1,\dots,\alpha_n)$ of rewriting steps of~$P$ with the same source, overlapping in a non-trivial and minimal way. Conducting a similar analysis as in \cref{sec:cb-class} shows that the critical $3$-branchings of~$P$ fall in one of the two cases
\begin{equation}
\label{E:w1ResolCT1}
\xymatrix @C=2.25em {
\strut
	\ar [r] _-{u_1}
	\ar@/^6ex/ [rrrr] _-{}="tgt1"
& \strut
	\ar@{} [r] ^{}="src1"
	\ar [r] _-{v_2}
& \strut 
	\ar [r] |-{w_2}
	\ar@/_6ex/ [rrr] ^-{}="tgt2"
& \strut
	\ar [r] |-{x_2}
	\ar@{} [r] _-{} ="src2"
	\ar@/^6ex/ [rrr] _-{}="tgt3"
& \strut
	\ar [r] |-{u_3} ^-{} ="src3"
& \strut
	\ar [r] _-{u_4}
& \strut
\ar@2 "1,3"!<0pt,5pt>;"tgt1"!<0pt,-7.5pt> ^-{\alpha}
\ar@2 "src2"!<0pt,-7.5pt>;"tgt2"!<0pt,7.5pt> ^-{\beta}
\ar@2 "src3"!<0pt,7.5pt>;"tgt3"!<0pt,-7.5pt> _-{\gamma}
}
\end{equation}
or
\begin{equation}
\label{E:w1ResolCT2}
\xymatrix @C=2.25em {
\strut
	\ar [r] _-{u_1}
	\ar @/^6ex/ [rrr] _-{}="tgt1"
& 
	\ar [r] _-{v_2} ^{}="src1"
& \strut
	\ar [r] |-{w_2}
	\ar@/_6ex/ [rrr] ^{}="tgt2"
& \strut
	\ar [r] ^-{v_3} _{}="src2"
& 
	\ar [r] |-{w_3}
	\ar@/^6ex/ [rr] _{}="tgt3"
& \strut
	\ar [r] _-{u_4}
& \strut
\ar@2 "src1"!<0pt,7.5pt>;"tgt1"!<0pt,-7.5pt> ^-{\alpha}
\ar@2 "src2"!<0pt,-7.5pt>;"tgt2"!<0pt,7.5pt> ^-{\beta}
\ar@2 "1,6"!<0pt,5pt>;"tgt3"!<0pt,-7.5pt> _-{\gamma}
}
\end{equation}
where~$\alpha$, $\beta$, and~$\gamma$ are $2$-generators of~$P$, and~$u_1v_2$, $w_2$, $u_3$, and~$u_4$ in \eqref{E:w1ResolCT1} or~$u_1v_2$, $w_2$, $v_3w_3$, and~$u_4$ in \eqref{E:w1ResolCT2} are reduced non-identity $1$-cells of~$P^*$, with~$u_1$ a $1$-generator of~$P$. Putting $u_2=v_2w_2x_2$, in \eqref{E:w1ResolCT1}, or $u_2=v_2w_2$ and $u_3=v_3w_3$, in \eqref{E:w1ResolCT2}, induces a one-to-one correspondence between the $4$-cells $u_1\vert u_2\vert u_3\vert u_4$ of~$\Sq_4(P)$ and the critical $3$-branchings of~$P$. 
This observation generalizes to establish a bijection between the $(n{+}1)$-generators of~$\Sq(P)$ and the critical $n$-branchings of~$P$.

\begin{theorem}
\label{T:SquierResolution}
Let~$P$ be a convergent 2-polygraph. There exists a unique structure of $(\omega,1)$\nbd-poly\-graph on~$\Sq(P)$, and unique unital section~$\iota$ and right $\iota$-contraction~$\sigma$ of~$\Sq(P)$, that satisfy $\iota_u = \rep{u}$, for every 1\nbd-cell~$u$ of~$P^*$, and
\begin{equation}
\label{E:SquierContraction}
\sigma_{(u_1\vert\cdots\vert u_n) u_{n+1}} =
\begin{cases}
u_1\vert\cdots\vert u_{n+1}
	&\text{if $u_1\vert\cdots\vert u_{n+1} \in \Sq_{n+1}(P)$,}
\\
1_{(u_1\vert\cdots\vert u_n) u_{n+1}}
	&\text{if~$u_n u_{n+1}$ is reduced,}
\end{cases}
\end{equation}
for all $n$-generator $u_1\vert\cdots\vert u_n$ of~$\Sq_n(P)$ with~$n\geq 1$, and reduced 1-cell~$u_{n+1}$ of~$P^*$.
Moreover, this structure makes $\Sq(P)$ a polygraphic resolution of the category~$\prescat{P}$.
\end{theorem}

\index{polygraphic!resolution!Squier}
\noindent
The polygraphic resolution $\Sq(P)$ thus constructed is called the \emph{Squier polygraphic resolution} of the category~$\prescat{P}$ with respect to the presentation $P$.

\begin{proof}
If the condition~\eqref{E:SquierContraction} is satisfied, then the source and target maps of~$\Sq(P)$ are imposed by the first case and the definition of an $\iota$-contraction. Indeed, writing $\ul{u}=u_1\vert\cdots\vert u_{n-1}$, we must have
\[
\src{}(u_1\vert\cdots\vert u_n) 
	= \src{}(\sigma_{\ul{u}u_n}) = \ul{u}u_n \comp1\sigma_{\tgt1(\ul{u})u_n}\comp2\cdots\comp{n-1}\sigma_{\tgt{n-1}(\ul{u})u_n},
\]
and
\[
\tgt{}(u_1\vert\cdots\vert u_n) 
	= \tgt{}(\sigma_{\ul{u}u_n}) = 
\begin{cases}
\rep{u_1u_2}
	&\text{if~$n=2$,}
\\
\sigma_{\src{}(\ul{u})u_n}
	&\text{otherwise.}
\end{cases}
\]
Then we prove, using the definition of the source and target of an $\iota$-contraction, that these source and target maps satisfy the globular relations. Next, according to \cref{L:FreeRightContraction}, it is necessary and sufficient to define~$\sigma$ on the $\iota$-essential and $\sigma$-essential cells of~$\freegpd{\Sq(P)}$. 

The $\iota$-essential $1$-cells are the $1$-cells~$u_1u_2$, where~$u_1$ is a $1$-generator of~$P$, $u_2$ is a reduced $1$-cell of~$P^*$, and~$u_1u_2$ is not reduced. If~$u_1\vert u_2$ is a $2$-generator of~$\Sq(P)$, then~\eqref{E:SquierContraction} imposes $\sigma_{u_1u_2}=u_1\vert u_2$. Otherwise, there exists a proper factorization $u_2=v_2w_2$ such that $u_1\vert v_2$ is a $2$-generator of~$\Sq(P)$, and~\eqref{E:SquierContraction} reads $\sigma_{(u_1\vert v_2)w_2}=1_{(u_1\vert v_2)w_2}$. This last equality imposes that the source and target of $\sigma_{(u_1\vert v_2)w_2}$ must be equal, giving the value of~$\sigma$ on~$u_1u_2$:
\[
\sigma_{u_1u_2}
	= \tgt{}(\sigma_{(u_1\vert v_2)w_2})
	= \src{}(\sigma_{(u_1\vert v_2)w_2})
	= (u_1\vert v_2)w_2 \comp1 \sigma_{\rep{u_1v_2}w_2}.
\]

Now, fix~$n\geq 2$. The $\sigma$-essential $n$-cells of~$\freegpd{\Sq(P)}$ are the $\ul{u}u_{n+1}$, where $\ul{u}=u_1\vert\cdots\vert u_n$ is an $n$-generator of~$\Sq(P)$, and~$u_{n+1}$ is a reduced $1$-cell of~$P^*$. We distinguish three cases. First, if $\ul{u}\vert u_{n+1}$ is an $(n+1)$-generator of~$\Sq(P)$, then~\eqref{E:SquierContraction} imposes $\sigma_{\ul{u}u_{n+1}}=\ul{u}\vert u_{n+1}$. Second, if $u_n u_{n+1}$ is reduced, then~\eqref{E:SquierContraction} gives $\sigma_{\ul{u}u_{n+1}}=1_{\ul{u}u_{n+1}}$. Otherwise, there exists a proper factorization $u_{n+1}=v_{n+1}w_{n+1}$ such that $\ul{u}\vert v_{n+1}$ is an $(n+1)$-generator of~$\Sq(P)$. In that case, \eqref{E:SquierContraction} implies that the source and the target of $\sigma_{(\ul{u}\vert v_{n+1})w_{n+1}}$ are equal. On the one hand, we have
\[
\src{}(\sigma_{(\ul{u}\vert v_{n+1})w_{n+1}})
= (\ul{u}\vert v_{n+1}) w_{n+1} \comp1 \sigma_{\tgt1(\ul{u}\vert
v_{n+1})w_{n+1}} \comp2 \cdots \comp n \sigma_{\tgt n(\ul{u}\vert v_{n+1})
w_{n+1}},
\]
and, on the other hand, we obtain
\begin{align*}
\tgt{}(\sigma_{(\ul{u}\vert v_{n+1})w_{n+1}})
	&= \sigma_{\src{}(\ul{u}\vert v_{n+1}) w_{n+1}} 
	= \sigma_{\src{}(\sigma_{\ul{u} v_{n+1}}) w_{n+1}} \\
	&= \sigma_{\ul{u}u_{n+1} \comp1 \sigma_{\tgt1(\ul{u})v_{n+1}} w_{n+1}
  \comp2 \cdots \comp n \sigma_{\tgt n(\ul{u})v_{n+1}} w_{n+1}}.
\end{align*}
Using the compatibility of~$\sigma$ with the compositions~$\comp{i}$ for $1\leq i \leq n$, we develop the latter expression, by induction on~$n$, to obtain a composite $(n{+}1)$-cell containing~$\sigma_{\ul{u}u_{n+1}}$, $\sigma_{\sigma_{t_n(\ul{u})v_{n+1}}w_{n+1}}$, and lower-dimensional invertible cells. Thus, we obtain a relation between two composite $(n{+}1)$-cells that defines~$\sigma_{\ul{u}u_{n+1}}$ in terms of the other involved cells. 

Finally, we apply \cref{T:ResolutionContraction} to conclude that~$\Sq(P)$ is a polygraphic resolution of the category~$C$.
\end{proof}

\subsection{Interpretation in the reduced case}
\label{SS:InterpretationReducedCase2}

Assume that~$P$ is a reduced convergent $2$-polygraph, and let us examine the first dimensions of~$\Sq(P)$.

The $2$-generators~$a\vert u$ of~$\Sq_2(P)$, for~$a$ a $1$-generator of~$P$ and~$u$ a reduced $1$-cell of~$P^*$ such that~$au$ is not reduced, have the shape
\[
  a\vert u:au \To \rep{au}.
\]
The $\iota$-contraction~$\sigma$ is given, on a $1$-cell~$au$ of~$P^*$ with~$a\in P_1$ and~$u$ reduced, by
\[
\sigma_{au} =
\begin{cases}
a\vert u
	&\text{if $a\vert u\in\Sq_2(P)$,} 
\\
1_{au} 
	&\text{if~$au$ is reduced,}
\\
(a\vert v)w \comp1 \sigma_{\rep{av}w}
	&\text{if~$u=vw$ with~$a\vert v\in\Sq_2(P)$.}
\end{cases}
\]
On more general $1$-cells, $\sigma$ is defined by the fact that it is a right $\iota$-contraction, and the relation
\[
  \sigma_{uv} = u\sigma_v \comp1 \sigma_{u\rep{v}}.
\]
By construction, the $3$-generators of~$\Sq_3(P)$ have the shape
\[
\xymatrix @R=1em {
& {\rep{au}v}
	\ar@2 @/^/ [dr] ^-{\sigma_{\rep{au}v}}
\\
{auv} 
	\ar@2 @/^/ [ur] ^-{(a\vert u)v}
	\ar@2 @/_2ex/ [rr] _-{\sigma_{auv}} ^-{}="tgt"
&& {\rep{auv}}
\ar@3 "1,2"!<-5pt,-10pt>;"2,2"!<-5pt,0pt> ^-*+{a\vert u\vert v}
\pbox.
}
\]
The $\iota$-contraction~$\sigma$ is defined on the $2$-generators $(a\vert u)v$ by~\eqref{E:SquierContraction}. The simple cases are $\sigma_{(a\vert u)v} = a\vert u\vert v$, if the latter belongs to~$\Sq_3(P)$, and $\sigma_{(a\vert u)v}=1_{(a\vert u)v}$ if~$uv$ is reduced. The more complicated case is the definition of~$\sigma_{(a\vert u)vw}$ when~$a\vert u\vert v$ belongs to~$\Sq_3(P)$. In this situation, the relation $\sigma_{(a\vert u\vert v)w}=1_{(a\vert u\vert v)w}$ implies $s(\sigma_{(a\vert u\vert v)w})=\tgt{}(\sigma_{(a\vert u\vert v)w})$, which develops into
\[
\vcenter{\xymatrix @C=4.4em @R=4.4em {
{\rep{au}vw}
	\ar@2@/^/ [r] ^-{\sigma_{\rep{au}v}w}
	\ar@3 []!<10pt,-15pt>;[d]!<10pt,25pt> ^(0.1){(a\vert u\vert v)w}
& {\rep{auv}w}
	\ar@2@/^/ [d] ^-{\sigma_{\rep{auv}w}}
	\ar@3 []!<-10pt,-25pt>;[d]!<-10pt,15pt> _(0.9){\sigma_{\sigma_{auv}w}}
\\
auvw
	\ar@2@/^/ [u] ^-{(a\vert u)vw}
	\ar@2 [ur] |-*+{\sigma_{auv}w}
	\ar@2@/_/ [r] _-{\sigma_{auvw}}
& {\rep{auvw}}
}}
\,\, = \,\,
%
\vcenter{\xymatrix @C=4.4em @R=4.4em {
{\rep{au}vw}
	\ar@2@/^/ [r] ^-{\sigma_{\rep{au}v}w}
	\ar@2 [dr] |-*+{\sigma_{\rep{au}vw}}
	\ar@3 []!<10pt,-25pt>;[d]!<10pt,15pt> ^(0.9){\sigma_{(a\vert u)vw}}
& {\rep{auv}w}
	\ar@2@/^/ [d] ^-{\sigma_{\rep{auv}w}}
	\ar@3 []!<-10pt,-15pt>;[d]!<-10pt,25pt> _(0.1){\sigma_{\sigma_{\rep{au}v}w}}
\\
auvw
	\ar@2@/^/ [u] ^-{(a\vert u)vw}
	\ar@2@/_/ [r] _-{\sigma_{auvw}}
& {\rep{auvw}}\pbox.
}}
\]
Finally, the $4$-generators of~$\Sq_4(P)$ have the same shape as this last defining equation, but in the case where~$vw$ is not reduced and with all proper left-factors reduced:
\[
\vcenter{\xymatrix @C=4.5em @R=4.5em {
{\rep{au}vw}
	\ar@2@/^/ [r] ^-{\sigma_{\rep{au}v}w}
	\ar@3 []!<10pt,-15pt>;[d]!<10pt,25pt> ^(0.1){(a\vert u\vert v)w}
& {\rep{auv}w}
	\ar@2@/^/ [d] ^-{\sigma_{\rep{auv}w}}
	\ar@3 []!<-10pt,-25pt>;[d]!<-10pt,15pt> _(0.9){\sigma_{\sigma_{auv}w}}
\\
auvw
	\ar@2@/^/ [u] ^-{(a\vert u)vw}
	\ar@2 [ur] |-*+{\sigma_{auv}w}
	\ar@2@/_/ [r] _-{\sigma_{auvw}}
& {\rep{auvw}}
}}
\!\!
\oqfl{a\vert u\vert v\vert w}
\!\!
\vcenter{\xymatrix @C=4.5em @R=4.5em {
{\rep{au}vw}
	\ar@2@/^/ [r] ^-{\sigma_{\rep{au}v}w}
	\ar@2 [dr] |-*+{\sigma_{\rep{au}vw}}
	\ar@3 []!<10pt,-25pt>;[d]!<10pt,15pt> ^(0.9){\sigma_{(a\vert u)vw}}
& {\rep{auv}w}
	\ar@2@/^/ [d] ^-{\sigma_{\rep{auv}w}}
	\ar@3 []!<-10pt,-15pt>;[d]!<-10pt,25pt> _(0.1){\sigma_{\sigma_{\rep{au}v}w}}
\\
auvw
	\ar@2@/^/ [u] ^-{(a\vert u)vw}
	\ar@2@/_/ [r] _-{\sigma_{auvw}}
& {\rep{auvw}}\pbox.
}}
\]

\section{Abelianization of Polygraphic Resolutions}
\label{S:Abelianization}

In \cref{Section:ComparisonWithHomologyMonoids}, we saw that the polygraphic homology of a monoid, seen as an $\omega$-category, coincides with its integral homology. In the same spirit, in this section we show how to deduce the homology of a small category with coefficients in natural systems from one of its resolutions in $\npCat{\omega}{1}$.
In particular, we show how to associate a resolution of natural systems to a polygraphic resolution of a category in~$\npCat{\omega}{1}$, and we illustrate this construction with examples of polygraphic resolutions calculated from convergent polygraphs.

\subsection{Free natural systems}
\index{free!natural system}
\index{natural system!free}
A natural system on a category $C$ is a functor from the category $\fact{C}$ of factorization of $C$ and with values in the category $\Ab$, see \cref{SS:NaturalSystems} for details.
For a family~$X$ of $1$-cells of~$C$, we denote by $\fnat{C}{X}$ the free natural system on~$C$ generated by~$X$ and given by 
\[
  \fnat{C}{X} = \bigoplus_{x\in X} \fact{C}(x,-).
\]

Fix an $(\omega,1)$-polygraph~$P$ presenting~$C$. We consider the free natural system $\fnat{C}{P_0}$ generated by the identity $1$-cells~$1_x$, for~$x\in P_0$. If~$u$ is a $1$-cell of~$C$, then $\fnat{C}{P_0}_u$ is the free abelian group generated by the pairs $(v,w)$ of $1$-cells of~$C$ such that $\tgt{}(v)=\src{}(w)=x$ and $vw=u$. 

We also consider, for every natural number~$n\geq 1$, the free natural system $\fnat{C}{P_n}$ generated by one copy of the $1$-cell~$\cl{\alpha}$ of~$C$ for each $n$-generator~$\alpha$ of~$P$. If~$u$ is a $1$-cell of~$C$, then $\fnat{C}{P_n}_u$ is the free abelian group generated by the triples $(v,\alpha,w)$, denoted by $v[\alpha]w$, made of an $n$-generator~$\alpha$ of~$P$, and $1$-cells~$v$ and~$w$ of~$C$, such that the composite~$v\cl{\alpha}w$ is well defined in~$C$ and equal to~$u$

The mapping of every $1$-generator~$x$ of~$P$ to the element~$[x]$ of $\fnat{C}{P_1}_{\cl{x}}$ is extended into a derivation of~$P_1^*$ into $\fnat{C}{P_1}$ by putting
\[
[1_u] = 0 
\qquad\text{and}\qquad
[uv] = [u]\cl{v} + \cl{u}[v],
\]
for all composable $1$-cells $u$ and $v$ in $C$.
Here, the natural system $\fnat{C}{P_1}$ on~$C$ is seen as a natural system on~$P_1^*$ by composition with the canonical projection $p : P_1^* \fl \cl{P}$. 

For $n>1$, the mapping of every $n$-generator~$\alpha$ of~$P$ to the element~$[\alpha]$ of $\fnat{C}{P_n}_{\cl{\alpha}}$ is extended to associate to every $n$-cell~$\phi$ of~$\freegpd{P}$ the element~$[\phi]$ of $\fnat{C}{P_n}_{\cl{\phi}}$, defined by induction on the size of~$\phi$ as follows:
\begin{align*}
  [1_\phi] &= 0,
  &
  [\phi^-] &= -[\phi],
  &
  [\phi\comp k \psi] &= 
	\begin{cases}
	[\phi]\cl{\psi} + \cl{\phi}[\psi] &\text{if $k=0$,} \\
	[\phi] + [\psi] &\text{otherwise.}
	\end{cases}  
\end{align*}

\subsection{Abelianization of polygraphic resolutions}
\label{SS:Abelianization}
\index{abelianization!of an $(\omega,1)$-polygraph}
\index{91-polygraph@$(\omega,1)$-polygraph!abelianization}
Let~$P$ be an $(\omega,1)$\nbd-poly\-graph. We denote by~$\fnat{\cl{P}}{P}$ the complex
\[
\xymatrix@C=4ex{
\cdots\, 
	\ar [r] 
& \fnat{\cl{P}}{P_n}
	\ar [r] ^-{\bnd_n}
& \fnat{\cl{P}}{P_{n-1}}
	\ar [r] 
& \,\cdots\,
	\ar [r] ^-{\bnd_1}
& \fnat{\cl{P}}{P_0}
	\ar [r] ^-{\varepsilon}
& \Z
	\ar [r]
&
0
}
\]
of natural systems on~$\cl{P}$, whose boundary maps are defined as follows. The augmentation morphism~$\varepsilon$ is defined, on every pair~$(u,v)$ of composable $1$\nbd-cells of~$\cl{P}$, by
\[ 
  \varepsilon(u,v) = 1.
\]
For $n\geq 1$, the morphism~$\bnd_n$ of natural systems on~$\cl{P}$ is given, on the generator~$[\alpha]$ corresponding to an $n$-generator~$\alpha$ of~$P$, by
\[
\bnd_n[\alpha] =
\begin{cases}
 (1,\cl{\alpha}) - (\cl{\alpha},1) 
	&\text{if~$n=1$,}
\\
[\tgt{}(\alpha)] - [\src{}(\alpha)]
	&\text{otherwise.}
\end{cases}
\]
By induction on the size of cells of~$\freegpd{P}$, we prove, for every $n$-cell~$\phi$ in~$\freegpd{P}$, with~$n\geq 1$, that
\[
\bnd_n[\phi] =
\begin{cases}
(1,\cl{\phi}) - (\cl{\phi},1) 
	&\text{if $n=1$,}
\\
[\tgt{}(\phi)] - [\src{}(\phi)] 
	&\text{otherwise.}
\end{cases}
\]
As a consequence, we have $\varepsilon\bnd_1=0$ and $\bnd_n\bnd_{n+1}=0$, for every $n\geq 1$, proving that $\fnat{\cl{P}}{P}$ is indeed a chain complex.

\begin{theorem}
\label{T:AbelianResolution}
If~$P$ is a polygraphic resolution of a category~$C$, then~$\fnat{C}{P}$ is a resolution by free natural systems on~$C$ of the constant natural system~$\Z$.
\end{theorem}
\begin{proof}
Let~$\iota$ be a unital section of~$P$. 
By \cref{T:ResolutionContraction}, $P$ admits a right $\iota$\nbd-contraction~$\sigma$. 
Let us consider the following families of morphisms of $\Z$\nbd-mod\-ules, indexed by a $1$-cell~$w$ of~$C$:
\begin{align*}
i_{-1} &: \Z \fl \fnat{C}{P_0}_w
&
i_n &: \fnat{C}{P_n}_u \fl \fnat{C}{P_{n+1}}_w 
\end{align*}
defined by
\begin{itemize}
\item $i_{-1}(1)=(w,1)$,
\item $i_0(u,v)=u[\rep{v}]$, for all $1$-cells $u,v$ of $C$ such that $w=uv$,
\item $i_n(u[\alpha]v)=u[\sigma_{\alpha\rep{v}}]$, for all~$n\geq 1$, $n$-generator $\alpha$ in $P$, and $1$-cells $u,v$ of $C$ such that $w=u\cl{\alpha}v$.
\end{itemize}

By induction on the size of the $n$-cells of the free $(\omega,1)$-category~$\freegpd{P}$, using the properties of a right $\iota$-contraction, we prove that
\[
i_n (u[\phi]v) = u[\sigma_{\phi\rep{v}}]
\]
holds for all~$n\geq 1$, $n$-cell~$\phi$ of~$\freegpd{P}$, and $1$-cells~$u,v$ of~$C$ such that the composite~$u\cl{\phi}v$ is well defined. We deduce that $(\sigma_n)_{n\geq 1}$ is a contracting homotopy for the complex $\fnat{C}{P}$.
\end{proof}

\subsection{Homological syzygies}
\label{SS:HomologicalSyzygiesCategory}
For~$n\geq 1$, the kernel of the differential~$\bnd_n$ defined in \cref{SS:Abelianization} is the natural system on $\cl{P}$ defined pointwise by 
\[
(h_n)_w = \ker 
\big( \xymatrix{
 \fnat{\cl{P}}{P_n}_w
	\ar [r] ^-{\bnd_n}
& \fnat{\cl{P}}{P_{n-1}}_w
}
\big),
\]
for every $1$-cell $w$ in $\cl{P}$. It is denoted by~$h_n(P)$, and its elements are called the \emph{homological $n$-syzygies of~$P$}\index{homological!syzygy}.
As a consequence of \cref{T:SquierResolution,T:AbelianResolution}, we obtain the following result.

\begin{theorem}
\label{T:SyzygiesGenerators}
Let~$C$ be a category and~$P$ a convergent presentation of~$C$. Then, for every~$n\geq 2$, the natural system~$h_n(P)$ is generated by the elements
\[
\bnd_n[u_1\vert\cdots\vert u_n] = [u_1\vert\cdots\vert u_{n-1}]\cl{u}_n + [\sigma_{\tgt{}(u_1\vert\cdots\vert u_{n-1})u_n}] - [\sigma_{\src{}(u_1\vert\cdots\vert u_{n-1})u_n}]
\]
where $u_1\vert\cdots\vert u_n$ ranges over the $n$-generators of~$\Sq(P)$, and~$\sigma$ is the right $\iota$\nbd-contrac\-tion associated to~$\Sq(P)$.
\end{theorem}

We now state a consequence of \cref{T:SquierResolution,T:AbelianResolution} for reduced presentations without critical $3$-branchings.
In \cref{SS:InterpretationReducedCase1,SS:InterpretationReducedCase2} we give an interpretation of the $n$-generators of the resolution $\Sq(P)$ when $P$ is a reduced convergent $2$-polygraph. In particular, there is a one-to-one correspondence between the $2$-generators of the resolution $\Sq(P)$ and the critical branchings on the one hand and between the $3$-generators and the critical $3$-branchings on the other. Moreover, when $P$ admits no critical $3$-branchings, it also has no critical $n$-branchings for $n\geq 3$, and so $\Sq_n(P)$ is empty for $n\geq 3$. The result is thus as follows.

\begin{corollary}
\label{C:SyzygiesWithoutTriple}
Let $C$ be a category, and $P$ a reduced convergent presentation of~$C$ without critical 3-branchings. 
Then the sequence
\[
\xymatrix@C=4ex{
0
	\ar [r] 
& \fnat{C}{\Sq_2(P)}
	\ar [r] ^-{\bnd_3}
&\fnat{C}{P_{2}}
	\ar [r] ^-{\bnd_2}	
&\fnat{C}{P_{1}}
	\ar [r] ^-{\bnd_1}
& \fnat{C}{P_0}
	\ar [r] ^-{\varepsilon}
& \Z
	\ar [r]
&
0
}
\]
is a partial resolution of length 4 of the trivial natural system $\Z$.
\end{corollary}

This result was proved by Squier for presentation of monoids by string rewriting systems using a direct method based on the characterization of critical $3$-branchings {\cite[Theorem 3.2]{squier1987word}}.
It is useful for proving examples of categories or monoids of finite type that do not admit a finite convergent presentation, as we did in \cref{chap:2-homology}.

\subsection{Convergence and homology of categories}
The construction given in this section allows us to calculate the homology of a category $C$ from a presentation of this category by a convergent $2$-polygraph $P$. Indeed, by \cref{T:SquierResolution}, $\Sq(P)$ is a polygraphic resolution of the category $C$, and by \cref{T:AbelianResolution},~$\fnat{C}{\Sq(P)}$ is a resolution by free natural systems on~$C$ of the constant natural system~$\Z$.
The homology of $C$ with coefficient in a contravariant natural system $D$ on $C$, as defined in \cref{SS:HomologyWithCoefficients}, is thus given by
\[
H_\ast(C,D) = \mathrm{Tor}^{FC}(D,\mathbb{Z}) = \mathrm{H}_\ast(D\otimes_{F C}{\Sq(P)}).
\]

\subsection{Example: the reduced standard polygraphic resolution}
\label{SS:StandardPolygraphicResolution}
Let~$C$ be a category. To simplify the example, assume that, if a composite morphism~$fg$ of~$C$ is an identity, then so are~$f$ and~$g$. The \emph{reduced standard presentation of~$C$} is the $2$-polygraph $\ol{\Std}_2(C)$ whose $0$-generators are objects of~$C$, with one $1$-generator~$\rep{f}$ for each non-identity morphism~$f$ of~$C$, and one $2$-generator
\[
f\vert g : \rep{f}\rep{g} \dfl \rep{fg}
\] 
for each pair $(f,g)$ of composable non-identity morphisms in~$C$. Without the simplifying hypothesis on~$C$, the target of~$f\vert g$ is replaced by $1_x$ if~$fg=1_x$ in~$C$.

The $2$-polygraph~$\ol{\Std}_2(C)$ is reduced and convergent, and applying \cref{T:SquierResolution} extends it into a polygraphic resolution of~$C$, denoted by~$\ol{\Std}(C)$ and called the \emph{reduced standard polygraphic resolution of~$C$}\index{standard!polygraphic resolution!reduced}\index{polygraphic!resolution!standard!reduced}. For~$n\geq 2$, the $n$-generators of~$\Std(C)$ are the  $f_1\vert\cdots\vert f_n$, such that each~$f_i$ is a non-identity morphism of~$C$ and each $(f_i,f_{i+1})$ is composable.

The source and target of the $3$-generators of~$\ol{\Std}(C)$ are given by
\[
\xymatrix @C=3em @R=0.7em {
& {\rep{fg}\rep{h}}
	\ar@2@/^/ [dr] ^-{fg\vert h}
\\
{\rep{f}\rep{g}\rep{h}}
	\ar@2@/^/ [ur] ^-{(f\vert g)\rep{h}}
	\ar@2@/_/ [dr] _-{\rep{f}(g\vert h)}
&& {\rep{fgh}}
\\
& {\rep{f}\rep{gh}}
	\ar@2@/_/ [ur] _-{f\vert gh}
\ar@3 "1,2"!<-5pt,-17pt>;"3,2"!<-5pt,17pt> ^-*+{f\vert g\vert h}
}
\]
and a $4$-generator $a\vert b\vert c\vert d$ has source
\[
  \vcenter{\xymatrix @C=2em@R=2em {
& {\rep{fg}\rep{h}\rep{k}}
	\ar@2 [rr] ^-{(fg\vert h)\rep{k}}
	\ar@{} [dr] |{(f\vert g\vert h)\rep{k}}
&& {\rep{fgh}\rep{k}}
	\ar@2@/^2ex/ [dr] ^-{fgh\vert k}
\\
{\rep{f}\rep{g}\rep{h}\rep{k}}
	\ar@2@/^2ex/ [ur] ^-{(f\vert g)\rep{h}\rep{k}}
	\ar@2 [rr] |{\rep{f}(g\vert h)\rep{k}}
	\ar@2@/_2ex/ [dr] _-{\rep{f}\rep{g}(h\vert k)}
&& {\rep{f}\rep{gh}\rep{k}}
	\ar@2 [ur] |{(f\vert gh)\rep{k}}
	\ar@2 [dr] |{\rep{f}(gh\vert k)}
	\ar@{} [rr] |-{f\vert gh\vert k}
&& {\rep{fghk}}
\\
& {\rep{f}\rep{g}\rep{hk}}
	\ar@2 [rr] _-{\rep{f}(g\vert hk)}
	\ar@{} [ur] |{\rep{f}(g\vert h\vert k)}
&& {\rep{f}\rep{ghk}}
	\ar@2@/_2ex/ [ur] _-{f\vert ghk}
}}
\]
and target
\[
  \vcenter{\xymatrix @C=2em@R=2em {
& {\rep{fg}\rep{h}\rep{k}}
	\ar@2 [rr] ^-{(fg\vert h)\rep{k}}
	\ar@2 [dr] |{\rep{fg}(h\vert k)}
&& {\rep{fgh}\rep{k}}
	\ar@2@/^2ex/ [dr] ^-{fgh\vert k}
\\
{\rep{f}\rep{g}\rep{h}\rep{k}}
	\ar@2@/^2ex/ [ur] ^-{(f\vert g)\rep{h}\rep{k}}
	\ar@2@/_2ex/ [dr] _-{\rep{f}\rep{g}(h\vert k)}
	\ar@{} [rr] |-{=}
&& {\rep{fg}\rep{hk}}
	\ar@2 [rr] |{fg\vert hk}
	\ar@{} [ur] |-{fg\vert h\vert k}
	\ar@{} [dr] |-{f\vert g\vert hk}
&& {\rep{fghk}}.
\\
& {\rep{f}\rep{g}\rep{hk}}
	\ar@2 [ur] |{(f\vert g)\rep{hk}}
	\ar@2 [rr] _-{\rep{f}(g\vert hk)}
&& {\rep{f}\rep{ghk}}
	\ar@2@/_2ex/ [ur] _-{f\vert ghk}
}}
\]
(with the arrows of $3$-cells removed for clarity).

For $n$-cells, $n\geq 2$, we prove, by induction on~$n$, that the source and target of $n$-generators are composites of the $(n-1)$-cells
\[
d_i(f_1\vert\cdots\vert f_n) = 
\begin{cases}
\rep{f}_1(f_2\vert\cdots\vert f_n) 
	& \text{if $i=0$,}
\\
f_1\vert\cdots\vert f_i f_{i+1}\vert\cdots\vert f_n
	& \text{if $1\leq i\leq n-1$,}
\\
(f_1\vert\cdots\vert f_{n-1})\rep{f}_n
	& \text{if $i=n$,}
\end{cases}
\]
with $k$-cells, for $1<k<n-1$.
More precisely, the source of $f_1\vert\cdots\vert f_n$ contains one copy of each $d_i(f_1\vert\cdots\vert f_n)$ for~$n-i$ even, and its target, one copy of each $d_i(f_1\vert\cdots\vert f_n)$ for~$n-i$ odd.

\Cref{T:AbelianResolution} applied to~$\ol{\Std}(C)$ gives a free resolution
\[
\xymatrix@C=3.1ex{
\cdots
	\ar [r]
& \fnat{C}{\ol{\Std}_n(C)}
	\ar [r] ^-{\bnd_{n}}
& \fnat{C}{\ol{\Std}_{n-1}(C)}
	\ar [r]	
& \cdots
	\ar [r]
& \fnat{C}{C_0}
	\ar [r] ^-{\varepsilon}
& \Z
	\ar [r]
&
0
}
\]
with differential defined by
\begin{align*}
  \bnd_n[f_1\vert\cdots\vert f_n] &\:=\: 
  (-1)^n f_1[f_2\vert\cdots\vert f_n] + \\
  &\hspace{1cm}
  \sum_{i=1}^{n-1} (-1)^{n-i} [d_i(f_1\vert\cdots\vert f_n)]
  + [f_1\vert\cdots\vert f_{n-1}]f_n.
\end{align*}

\subsection{The associative polygraphic resolution}
\label{SS:associativePolygraphicResolution}
Let~$A$ be the monoid with one non-trivial idempotent element, that is presented by the following $2$-polygraph:
\[
  \AsPoly_2 = \Pres{a_0}{a_1}{a_2:a_1 a_1\dfl a_1}.
\]
This polygraph is reduced and convergent, see~\cref{ex:aa-a-confl}, with one critical $n$-branching for every~$n\geq 2$. Thus, the reduced standard polygraphic resolution $\AsPoly_{\omega}=\Sq(\AsPoly_2)$ of~$A$, given by \cref{T:SquierResolution}, has one $n$-generator~$a_n$ for every~$n\geq 0$, corresponding to the product $a_1\vert\cdots\vert a_1$ of~$n$ copies of~$a_1$. 
The $3$-generator~$a_3$ of~$\AsPoly_{\omega}$ is given in classical notation and in string diagrams respectively as follows:
\begin{align*}
  a_2a_1 \comp1 a_2 &\otfl{a_3} a_1a_2 \comp1 a_2
  &
  \satex{mon-assoc-small-l}
  &\otfl{\satex{alpha-small}}
  \satex{mon-assoc-small-r}.
\end{align*}
The $4$-generator~$a_4$ of~$\AsPoly_{\omega}$ is
\[
  \satex{As-a4}
\]
which, contracting by one dimension, can also be pictured as Mac Lane's pentagon or Stasheff's polytope~$K_4$:
\[
  \xymatrix @R=1em @C=1em {
    &{ \satex{mon-cp1-vsmall-1} }
    \ar@3 [rr] 
    ^-{\satex{mon-cp1-vsmall-12}} 
    _-{}="1" 
    && { \satex{mon-cp1-vsmall-2} }
    \ar@3@/^/ [dr] ^-{\satex{mon-cp1-vsmall-23} }
    \\
    { \satex{mon-cp1-vsmall-0} }
    \ar@3@/^/ [ur] ^-{\satex{mon-cp1-vsmall-01}}
    \ar@3@/_/ [drr] _-*+{\satex{mon-cp1-vsmall-04}}
    &&&& { \satex{mon-cp1-vsmall-3}\text.}
    \\
    && { \satex{mon-cp1-vsmall-4} }
    \ar@3@/_/ [urr] _-*+{\satex{mon-cp1-vsmall-43}}
    \ar@4 "1"!<-2.5ex,-5ex>;[]!<-2.5ex,6ex> ^-*+{\satex{aleph-small}}
  }
\]
Finally, the $5$-generator~$a_5$ of~$\AsPoly_{\omega}$ has the shape of Stasheff's polytope~$K_5$, its source being
\[
  \satex{K5-src}
\]
and its target being given by a symmetric composite $4$-cell, see~\cite[Section~6.1]{GuiraudMalbos12advances}. \Cref{T:AbelianResolution}, applied to~$\AsPoly_{\omega}$, yields a resolution 
\[
\xymatrix @C=1em{
\cdots\:
	\ar [r]
& {\fnat{A}{\:\satex{aleph-vsmall}\:}}
	\ar [r] ^-{\bnd_4}
& {\fnat{A}{\:\satex{alpha-vsmall}\:}}
	\ar [r] ^-{\bnd_3}
& {\fnat{A}{\:\satex{mu-vsmall}\:}}
	\ar [r] ^-{\bnd_2}
& {\fnat{A}{\:\satex{1-vsmall}\:}}
	\ar [r] ^-{\bnd_1}
& {\fnat{A}{\ast}}
	\ar [r] ^-{\varepsilon}
& \Z
	\ar [r]
& 0
}
\]
of~$\Z$ by free natural systems on~$A$. 
Computing this differential on each $n$\nbd-cell of~$\AsPoly_{\omega}$ gives generators of the natural systems of homological $n$-syzygies of~$\AsPoly_{\omega}$. For example, $h_2(\AsPoly)$ is generated by
\[
\bnd_3 [\:\satex{alpha-small}\:]
	= \left[ {\satex{mon-assoc-small-l}} \right] - \left[ {\satex{mon-assoc-small-r}} \right]
	= \big[\satex{mu-small}\big] a - a \big[\satex{mu-small}\big]
\]
while $h_3(\AsPoly)$ is generated by 
\begin{align*}
\bnd_4\spa{\satex{aleph-small}} \:
	&=\: \left[\:\satex{mon-cp1-vsmall-01}\:\right]
		+ \left[\:\satex{mon-cp1-vsmall-12}\:\right]
		+ \left[\:\satex{mon-cp1-vsmall-23}\:\right]
		- \left[\:\satex{mon-cp1-vsmall-04}\:\right]
		- \left[\:\satex{mon-cp1-vsmall-43}\:\right] \\
	&=\: a\left[\:\satex{alpha-small}\:\right] - \left[\:\satex{alpha-small}\:\right] + \left[\:\satex{alpha-small}\:\right] a.
\end{align*}
Similarly, $h_4(\AsPoly)$ is generated by $\bnd_5[a_5]$, which is equal, by definition, to 
\[
\left[ \:\satex{mon-K5-vsmall-Y}\: \right]
+ \left[ \:\satex{mon-K5-vsmall-X}\: \right]
+ \left[ \:\satex{mon-K5-vsmall-Z}\: \right]
- \left[ \:\satex{mon-K5-vsmall-Y2}\: \right]
- \left[ \:\satex{mon-K5-vsmall-X2}\: \right]
- \left[ \:\satex{mon-K5-vsmall-Z2}\: \right]
\]
and reduces to 
\[
  \bnd_5[a_5] = \left[ \:\satex{aleph-small}\: \right] a - a \left[
  \:\satex{aleph-small}\: \right].
\]

\subsection{The category of monotone surjections}
\label{sec:simpl-surj-coh}
We denote by~$\Simplsurj$ the subcategory of the simplicial category whose
objects are the natural numbers and whose morphisms from~$m$ to~$n$ are the
monotone surjections from $\set{0,\dots,m}$ to $\set{0,\dots,n}$. This category, studied in~\cite{LivernetRichter11}, see also {\cite[Section 7.5, Exercise~3.(a)]{MacLane98}}, admits a presentation by the $2$-polygraph~$P$ with the natural numbers as $0$-generators, with one $1$-generator $x^n_i:n+1\fl n$ for all natural numbers $1\leq i\leq n$, and one $2$-generator
\[
\xymatrix @!C @R=1em{
& n+1
	\ar@/^/ [dr] ^-{x^n_j} 
	\ar@2 []!<0pt,-15pt>;[dd]!<0pt,15pt> ^-{x^n_{i,j}}
\\
n+2 
	\ar@/^/ [ur] ^-{x^{n+1}_i}
	\ar@/_/ [dr] _-{x^{n+1}_{j+1}}
&& n
\\
& n+1 
	\ar@/_/ [ur] _-{x^n_i}
}
\]
for all natural numbers $0\leq i\leq j\leq n + 1$. The $1$-generator~$x_i^n$ represents the map 
\[
x^n_i (j) \:=\: 
\begin{cases}
j 
	&\text{if $j\leq i$,} 
\\
j-1 
	&\text{if $j>i$.}
\end{cases}
\]
This is a variant of the presentation constructed in \cref{sec:simpl-cat}.
Thereafter, we drop the exponents of the $1$-generators and $2$-generators of~$P$, simply writing~$x_i$ and~$x_{i,j}$.

The $2$-polygraph~$P$ is convergent. Indeed, for termination, given a $1$-cell $u=x_{i_1}\dots x_{i_k}$ of~$P^*$, we define the natural number $\nu(u)$ as the number of pairs $(i_p,i_q)$ such that $i_p\leq i_q$, with $1\leq p<q\leq k$. In particular, we have $\nu(x_i x_j) = 1$ and $\nu(x_{j+1} x_i)=0$ when $i\leq j$, 
giving $\nu(\src{}(x_{i,j}))>\nu(\tgt{}(x_{i,j}))$. Moreover, we have $\nu(w u w') > \nu(w v w')$ when $\nu(u)>\nu(v)$ holds. Thus, for every non-identity $2$-cell $a:u\dfl v$ of~$P^*$, the strict inequality $\nu(u)>\nu(v)$ is satisfied, giving termination. Moreover, the $2$-polygraph~$P$ has one critical branching $(x_{i,j} x_k, x_i x_{j,k})$ for all possible $0\leq i\leq j\leq k\leq n+2$, which is confluent.

\Cref{T:SquierResolution}, applied to~$P$, gives a polygraphic resolution~$\Sq(P)$ of~$\Simplsurj$, 
 whose $3$-cells are given, in classical notation and in string diagrams (with $x_i = \:\satex{1-vsmall}_i\:$ and $x_{i,j}=\satex{tau-vsmall}_{i,j}$) respectively, by
\[
\vcenter{
\xymatrix @!C @C=-1em {
& x_{j+1} x_i x_k 
	\ar@2 [rr] ^-{x_{j+1} x_{i,k}} _{}="src"
&& x_{j+1} x_{k+1} x_i
	\ar@2 @/^2ex/ [dr] ^-{x_{j+1,k+1} x_i}
\\
x_i x_j x_k 
	\ar@2 @/^2ex/ [ur] ^-{x_{i,j} x_k} 
	\ar@2 @/_2ex/ [dr] _-{x_i x_{j,k}}
&&&& x_{k+2} x_{j+1} x_i
\\
& x_i x_{k+1} x_j
	\ar@2 [rr] _-{x_{i,k+1} x_j} ^{}="tgt"
&& x_{k+2} x_i x_j
	\ar@2 @/_2ex/ [ur] _-{x_{k+2} x_{i,j}}
\ar@3 "src"!<0pt,-25pt>;"tgt"!<0pt,25pt> ^-{x_{i,j,k}}
}}
\]
\[
  \vcenter{
    \xymatrix @C=7em {
      \satex{yb_l-ssmall}_{i,j,k}
      \ar@3 [r] ^-*+{\satex{yb-ssmall}_{i,j,k}}
      &
      \satex{yb_r-ssmall}_{i,j,k}
    }
  }
\]
The $(\omega,1)$-polygraph~$\Sq(P)$ has one $4$-generator~$x_{i,j,k,l}$ for every
possible index $0\leq i\leq j\leq k\leq l \leq n+3$, given in string diagrams and omitting the subscripts, by the following diagrams:
\[
  \satex{yb-cp}
\]
Then, \cref{T:AbelianResolution} gives, in particular, generators for the natural systems of homological $n$-syzygies of~$P$. For example, $h_2(P)$ is generated by the elements
\begin{align*}
\bnd_3\left[\satex{yb-ssmall}_{i,j,k}\right] 
	&= \left[\satex{yb_l-ssmall}_{i,j,k} \right] 
        - \left[\satex{yb_r-ssmall}_{i,j,k} \right]
        \\
	&= 
\left\{
\begin{array}{cl}
  & \pa{\spa{\satex{tau-ssmall}_{i,j}} x_k - x_{k+2} \spa{\satex{tau-ssmall}_{i,j}}}\\[1ex]
  + & \pa{x_{j+1} \spa{\satex{tau-ssmall}_{i,k}} - \spa{\satex{tau-ssmall}_{i,k+1}} x_j}\\[1ex]
  + & \pa{\spa{\satex{tau-ssmall}_{j+1,k+1}} x_i - x_i \spa{\satex{tau-ssmall}_{j,k}}}\pbox.
\end{array}
\right.  
\end{align*}

\section{Categories of Finite Homological Type}
\label{sec:finiteness-conditions}
\label{S:FinitenessConditions}

In \cref{chap:2-fdt}, we introduced the notion of a finite derivation type for categories. In this section, we show how to refine this notion in higher dimensions and how to relate it with finite homological type.

\subsection{Higher-dimensional finite derivation type}
\label{sec:higher-FDT}
\index{finite derivation type!8-@$n$-}
\index{finite derivation type!9-@$\infty$-}
\index{FDT8 (finite 8-derivation type)@$\FDT_n$ (finite $n$-derivation type)}
\index{FDT9 (finite 9-derivation type)@$\FDT_\infty$ (finite $\infty$-derivation type)}

Let~$C$ be a category and $n\in\Nbinfty$. We say that~$C$ has \emph{finite $n$-derivation type}, $\FDT_n$ for short, if it admits a polygraphic resolution~$P$ in~$\npCat{\omega}{1}$ such that~$P_k$ is finite for every~$k\leq n$. In particular, $C$ has~$\FDT_1$ if it is finitely generated, $\FDT_2$ if it is finitely presented, and $\FDT_3$ if it has finite derivation type as defined in \cref{sec:2fdt-def}. By definition, $\FDT_{\infty}$ implies $\FDT_n$, and $\FDT_{n+1}$ implies $\FDT_n$, for every $n\geq 0$.

As an immediate consequence of Theorem~\ref{T:SquierResolution}, we deduce the following condition for finite convergence.

\begin{theorem}
A category with a finite convergent presentation has~$\FDT_{\infty}$.
\end{theorem}

\subsection{Categories of finite homological type}
\index{FP8@$\FP_n$}
\index{homological!type $\FP_n$}
A category $C$ is of \emph{homological type~$\FP_n$} if the constant natural system $\Z$ on $C$ is of homological type~$\FP_n$, see \cref{S:CategoriesFiniteHomologicalType} for a summary on the notion of finite homological type for categories.
As a consequence of \cref{T:SquierResolution,T:AbelianResolution}, we obtain the following implications.

\begin{theorem}
\label{T:FDT=>FP}
Let~$C$ be a category and $n\in\Nbinfty$. If~$C$ has~$\FDT_n$, then it is of homological type~$\FP_n$. In particular, if~$C$ has a finite convergent presentation, then it is of homological type~$\FP_{\infty}$.
\end{theorem}

This result generalizes~\cite[Theorem~3.2]{cremanns1994finite} and~\cite[Theorem~3]{lafont1995new}, stating that if a monoid has $\FDT$, then it is~$\FP_3$
(see also~\cite{pride1995low}). It also generalizes Squier's homological theorem~\cite[Theorem~4.1]{squier1987word}, that says that a monoid admitting a finite convergent presentation is~$\FP_3$, and the extensions of Squier's result in \cite{brown1992geometry,Groves90,kobayashi1990complete} that prove that such a monoid is $\FP_\infty$.

\subsection{Example}
The $2$-polygraph $\AsPoly_2$ defined in \cref{SS:associativePolygraphicResolution} extends to a polygraphic resolution $\AsPoly_{\omega}$ having one $n$-generator for every $n\geq 0$. Hence the polygraph $\AsPoly_2$ and the presented monoid have~$\FDT_{\infty}$.

\section{Homological Syzygies and Identities Among Relations}
\label{sec:H-iar}
In this section we establish an isomorphism between the natural systems of homological $2$-syzygies and identities among relations for a category presented by a $2$-polygraph. This is an extension to category presentations of a Brown-Huebschmann theorem in group theory that states an isomorphism between the modules of identities among relations and homological $2$-syzygies for group presentations~\cite{BrownHuebschmann82}.

In this section, $P$ denotes a $2$-polygraph. We aim to build an isomorphism between the natural system $\Pi(P)$ of identities among relations of $P$ defined in~\cref{sec:id-among-rel} and the natural system $h_2(P)$ of its homological $2$-syzygies defined in \cref{SS:HomologicalSyzygiesCategory}.

\begin{lemma}
  Let $P$ be a 2-polygraph. For every 2-loop $\psi$ of $\freegpd{P}$, we have $[\psi]=0$ in~$\fnat{\cl{P}}{P_2}$ if and only if $\iar{\psi}=0$ holds in $\Pi(P)$.
\end{lemma}

\begin{proof}
  To prove that $\iar{\psi}=0$ implies $[\psi]=0$, we check that the relations~\eqref{eq:iar1} and~\eqref{eq:iar2} defining $\Pi(P)$ are also satisfied in $\fnat{\cl{P}}{P_2}$. The first relation is given by the definition of the map $[-]$. The second relation is checked as follows:
  \[
    [\psi\comp1 \phi] = [\psi]+[\phi] = [\phi]+[\psi] = [\phi\comp1 \psi].
  \]
  Conversely, let us consider a $2$-loop $\psi$ in $\freegpd{P}$ with source $w$ such that $[\psi]=0$. We decompose $\psi$ into
  \[
    \psi = u_1 \alpha_1^{\epsilon_1} v_1 \comp1 \cdots \comp1 u_p \alpha_p^{\epsilon_p} v_p,
  \]
  where $\alpha_i$ is a $2$-generator, $u_i,v_i$ are $1$-cells of $\freegpd{P}$, and $\epsilon_i\in\{-,+\}$. Then we get
  \[
 0 = [\psi] = \epsilon_1 \cl{u}_1 [\alpha_1] \cl{v}_1 + \ldots + \epsilon_p \cl{u}_p [\alpha_p] \cl{v}_p.
  \]
  Since the natural system $\fnat{\cl{P}}{P_2}$ is freely generated by the elements $[\alpha]$ of $\fnat{\cl{P}}{P_2}_{\cl{\alpha}}$, for $\alpha$ a $2$-generator, this implies the existence of a self-inverse permutation $\tau$ of $\{1,\dots,p\}$ such that the following relations are satisfied:
  \begin{align*}
    \alpha_i  &= \alpha_{\tau(i)},
    &
    \cl{u}_i &= \cl{u}_{\tau(i)},
    &
    \cl{v}_i &= \cl{v}_{\tau(i)},
    &
    \epsilon_i &= -\epsilon_{\tau(i)}.
  \end{align*}
  Let us denote, for every $1\leq i\leq p$, the source and target of $\alpha_i^{\epsilon_i}$ by $w_i$ and $w'_i$, respectively. They satisfy $\cl{w}_i=\cl{w}'_i$. We also fix a section $\rep{\;\cdot\;}$ and a left strategy~$\sigma$ for the $2$-polygraph $P$, so that $\rep{u}=\rep{v}$ for all $1$-cells $u$ and $v$ such that $\cl{u}=\cl{v}$.
  
  For every $1\leq i\leq p$, we denote by $\psi_i$ the following $2$-cell of $\freegpd{P}$:
  \[
    \psi_i = \sigma^-_{u_i w_i v_i} \comp1 u_i\alpha_i^{\epsilon_i}v_i \comp1 \sigma_{u_i w'_i v_i}.
  \]
  Using the facts that $w=u_1w_1v_1=u_p w'_p v_p$ and that $u_iw'_iv_i=u_{i+1} w_{i+1} v_{i+1}$ for every $1\leq i<p$, we can write the $2$-loop $\psi$ as the following composite:
  \[
    \psi = \sigma_w \comp1 \psi_1 \comp1 \cdots \comp1 \psi_p \comp1 \sigma^-_w.
  \]
  As a consequence, we get
  \[
    \iar{\psi} 
    = \iar{\sigma^-_w\comp1 \psi\comp1 \sigma_w} 
    = \iar{\psi_1} + \ldots + \iar{\psi_p}.
  \]
  In order to conclude, we prove that the equality $\iar{\psi_{\tau(i)}} = -\iar{\psi_i}$ holds, for every $1\leq i\leq p$. Since $\sigma$ is a left strategy, we have
  \[
    \sigma_{u_i w_i v_i} = \sigma_{u_i} w_i v_i \comp1 \sigma_{\rep{u}_i w_i} v_i \comp1 \sigma_{\rep{u_i w}_i v_i}
  \]
  and, using the fact that $\rep{u_i w}_i=\rep{u_i w}'_i$, we have
  \[
    \sigma_{u_i w'_i v_i} = \sigma_{u_i} w'_i v_i \comp1 \sigma_{\rep{u}_i w'_i} v_i \comp1 \sigma_{\rep{u_i w}_i v_i}.
  \]
  This gives
\begin{align*}
\iar{\psi_i} &= \iar{ \sigma^-_{\rep{u_i w_i}v_i} \!\comp1\!\sigma^-_{\rep{u}_i w_i} v_i  \!\comp1\! \sigma^-_{u_i} w_i v_i  \!\comp1 u_i\alpha_i^{\epsilon_i}v_i  \comp1\! \sigma_{u_i} w'_i v_i  \!\comp1 \!\sigma_{\rep{u}_i w'_i} v_i  \!\comp1 \! \sigma_{\rep{u_i w_i}v_i}}\\
&
    = \iar{ 
      \sigma^-_{\rep{u}_i w_i} v_i \comp1 \sigma^-_{u_i} w_i v_i
      \comp1 u_i\alpha_i^{\epsilon_i}v_i 
      \comp1 \sigma_{u_i} w'_i v_i \comp1 \sigma_{\rep{u}_i w'_i} v_i
    }\\
&
    = \iar{ 
      \sigma^-_{\rep{u}_i w_i}
      \comp1 \rep{u}_i\alpha_i^{\epsilon_i} 
      \comp1 \sigma_{\rep{u}_i w'_i}
    } \cl{v}_i.
  \end{align*}
Now, let us compute $\iar{\psi_{\tau(i)}}$. We already know that $\alpha_{\tau(i)}=\alpha_i$ and $\epsilon_{\tau(i)}=-\epsilon_i$. As a consequence, we get $w_{\tau(i)}=w'_i$ and $w'_{\tau(i)}=w_i$. Moreover, we have $\rep{u}_{\tau(i)}=\rep{u}_{i}$, so that we have:
\begin{align*}
\iar{\psi_{\tau(i)}} 
	\:=\: &\left\lfloor 
		\sigma^-_{\rep{u_i w_i}v_{\tau(i)}}
			\star_1 \sigma^-_{\rep{u}_i w'_i} v_{\tau(i)} 
			\star_1 \sigma^-_{u_i} w'_i v_{\tau(i)} \right.\\
		&\left.\star_1 u_{\tau(i)} \alpha_i^{-\epsilon_i} v_{\tau(i)} 
		\star_1 \sigma_{u_i} w_i v_{\tau(i)} 
			\star_1 \sigma_{\rep{u}_i w_i} v_{\tau(i)} 
			\star_1 \sigma_{\rep{u_i w_i} v_{\tau(i)}}\right\rfloor \\
& = \iar{
      \sigma^-_{\rep{u}_i w'_i}  
      \comp1 \rep{u}_i\alpha_i^{-\epsilon_i} 
      \comp1 \sigma_{\rep{u}_i w_i} 
    } \cl{v}_i
    = -\iar{\psi_i}.
\end{align*}
  This implies $\iar{\psi}=0$, thus concluding the proof.
\end{proof}

\begin{lemma}
\label{Lemma:h2Pi}
For every element $a$ in $h_2(P)$, there exists a 2-loop $\psi$ in~$\freegpd{P}$ such that $a=[\psi]$ holds.
\end{lemma}

\begin{proof}
Let $w$ be the $1$-cell of $\cl{P}$ such that $a$ belongs to $\fnat{\cl{P}}{P_2}_w$ and let $P_3$ be an acyclic extension of the $(2,1)$-category~$\freegpd{P}$. Since $d_2(a)=0$, by acyclicity of $P_3$ and \cref{T:AbelianResolution}, there exists $b$ in $\fnat{\cl{P}}{P_3}_w$ such that $a=d_3(b)$. By definition of $\fnat{\cl{P}}{P_3}_w$, we can write
\[
b =\epsilon_1 u_1 [\alpha_1] v_1 + \ldots + \epsilon_p u_p [\alpha_p] v_p,
\]
with, for every $1\leq i\leq p$, $\alpha_i\in P_3$, $u_i,v_i\in\cl{P}$ and $\epsilon_i\in\{-,+\}$ such that $u_i\cl{\alpha}_iv_i=w$ holds. We fix a section $\rep{\;\cdot\;}$ of $P$ and we choose $2$-cells
\[
\phi_i \::\: \rep{w} \:\dfl\: \rep{u}_i s_1(\alpha_i^{\epsilon_i}) \rep{v}_i
\qquad\text{and}\qquad
\psi_i \::\: \rep{u}_i t_1(\alpha_i^{\epsilon_i}) \rep{v}_i \:\dfl\: \rep{w}.
\]
Let $A$ be $3$-cell of $\freegpd{P}_3$ defined by
\[
A 
	= \big( \phi_1\comp1 \rep{u}_1\alpha_1^{\epsilon_1} \rep{v}_1 \comp1 \psi_1 \big)
			\comp1 \cdots
			\comp1 \big( \phi_k\comp1 \rep{u}_k\alpha_k^{\epsilon_k} \rep{v}_k \comp1 \psi_k \big). 		
\]
By definition of $[\cdot]$ on $3$-cells, we have 
\[
[A] 
	= \sum_{i=1}^p \big[ \phi_i\comp1 \rep{u}_i\alpha_i^{\epsilon_i}\rep{v}_i \comp1 \psi_i \big]
	= \sum_{i=1}^p \big( [1_{\phi_i}] + \epsilon_i u_i[\alpha_i] v_i + [1_{\psi_i}] \big)
	= b.
\]
Finally, we get
\[
a = d_3[A] = [s(A)]-[t(A)] = [s(A)\comp1 t(A)^-]. 
\]
Hence $\psi=s(A)\comp1 t(A)^-$ is a $2$-loop of $\freegpd{P}$ that satisfies $a=[\psi]$.
\end{proof}

When~$P$ is a convergent $2$-polygraph, we have seen that the natural systems $h_2(P)$ and $\Pi(P)$ on $\cl{P}$ are generated by a family of generating confluences of~$P$. The following result from~{\cite[Theorem 5.6.5]{GuiraudMalbos12advances}} states that, more generally, the natural systems $h_2(P)$ and $\Pi(P)$ are isomorphic, as proved by Brown-Huebschmann for presentations of groups in \cite{BrownHuebschmann82}.

\begin{theorem}
\label{Theorem:IsomorphismPiH2}
Let $P$ be a 2-polygraph. The natural systems $\Pi(P)$ and $h_2(P)$ are isomorphic.
\end{theorem}

\begin{proof}
We define a morphism of natural systems $\Phi : \Pi(P) \fl h_2(P)$ by setting 
$\Phi\iar{\phi} = [\phi]$ for every identity $\phi$ of $P$. 
This definition is correct, since the defining relations of $\Pi(P)$ also hold in $\fnat{\cl{P}}{P_2}$, and thus in $h_2(P)$. 
Moreover, $\Phi$ is a morphism of natural systems, since we have
\[
\Phi(u\iar{\phi}v) = \Phi(\iar{\rep{u}\phi\rep{v}}) = [\rep{u}\phi\rep{v}] = u[\phi]v = u\Phi(\iar{\phi})v,
\]
for every $2$-loop $\phi$ in $\freegpd{P}$ and $1$-cells $u$, $v$ in $\cl{P}$ such that $\rep{u}\phi\rep{v}$ is defined. 

Now, let us define a morphism of natural systems $\Psi : h_2(P) \fl \Pi(P)$. Let $a$ be an element of $h_2(P)_w$. By \cref{Lemma:h2Pi}, there exists a $2$-loop $\psi:u\dfl u$ in $\freegpd{P}$ such that $a=[\psi]$ and $w=\cl{u}$. We define $\Psi(a) = \iar{\psi}$.
This definition does not depend on the choice of $\psi$. Indeed, let us assume that $\phi:v\dfl v$ is a $2$-loop such that $a=[\phi]$. It follows that $\cl{v}=\cl{u}$, and we can choose a $2$-cell $\xi:u\dfl v$ in $\freegpd{P}$. Then we have
\[
a = [\psi] = [\phi] = [\xi\comp1 \phi\comp1 \xi^-].
\]
As a consequence, we get
\[
[\psi\comp1 \xi^-\comp1 \phi^-\comp1 \xi] = [\psi] - [\xi \comp1 \phi \comp1 \xi^-] = 0.
\]
Thus
\[
0 = \iar{\psi\comp1 \xi^- \comp1 \phi^- \comp1 \xi} = \iar{\psi} - \iar{\xi\comp1 \phi\comp1 \xi^-} = \iar{\psi} -\iar{\phi}.
\]
Finally, the relations $\Psi\Phi=1_{\Pi(P)}$ and $\Phi\Psi=1_{h_2(P)}$ are direct consequences of the definitions of $\Phi$ and $\Psi$.
\end{proof}

The following result relates the low-dimensional finiteness properties seen in this chapter and in \cref{subsectionAbelianTDF} for the property $\FDTAB$.

\begin{theorem}
\label{Theorem:TDFAB<=>FP_3}
Let $P$ be a finite 2-polygraph. The following conditions are equivalent.
\begin{enumerate}
\item \label{Theorem:TDFAB<=>FP_3_1} The category $\prescat{P}$ is of homological type $\FP_3$.
\item \label{Theorem:TDFAB<=>FP_3_2} The natural system $h_2(P)$ on $\prescat{P}$ is finitely generated.
\item \label{Theorem:TDFAB<=>FP_3_3} The natural system $\Pi(P)$ on $\prescat{P}$ is finitely generated.
\item \label{Theorem:TDFAB<=>FP_3_4} The category $\prescat{P}$ has $\FDTAB$.
\end{enumerate}
\end{theorem}
\begin{proof}
The equivalence between \ref{Theorem:TDFAB<=>FP_3_1} and \ref{Theorem:TDFAB<=>FP_3_2} comes from the definition of the property $\FP_3$. The equivalence between \ref{Theorem:TDFAB<=>FP_3_2} and \ref{Theorem:TDFAB<=>FP_3_3} is a consequence of \cref{Theorem:IsomorphismPiH2}. The equivalence between \ref{Theorem:TDFAB<=>FP_3_3} and \ref{Theorem:TDFAB<=>FP_3_4} is given by \cref{Proposition:IartfIffFdtab}.
\end{proof}

Note that, following \cref{T:FDT=>FP}, the property $\FDT_3$ implies $\FP_3$. We expect the reverse implication to be false in general, which amounts to proving that $\FDTAB$ does not imply $\FDT_3$, since $\FP_3$ is equivalent to $\FDTAB$ for finitely presented categories. This question is still open.

\subsection{Identities among relations for higher polygraphs}
\label{SS:IARnPolygraphs}
We conclude this chapter by mentioning some results on identities among relations for $n$\nbd-poly\-graphs for any $n\geq 0$. 
As in \cref{subsectionAbelianTDF} for the case $n=2$, we call an $(n,n-1)$\nbd-cate\-gory~$C$ \emph{abelian} if, for every $(n-1)$-cell $u$ of~$C$, the group $\Aut^{C}_u$ of $n$-loops of~$C$ with source $u$ is abelian. 
For $C$ an $(n,n-1)$-category, its \emph{abelianization}~$\fab C$ is the quotient of~$C$ by the cellular extension that contains one $n$-sphere
$\varphi\comp{n-1} \psi \to \psi\comp{n-1} \varphi$
for every $n$-loops $\varphi$ and $\psi$ of~$C$ with the same source.

Identities among relations for $2$-polygraphs, as defined in~\cref{sec:id-among-rel}, were extended to the structure of an $n$-polygraph $P$ in~\cite{GuiraudMalbos10smf} to form a natural system~$\Pi(P)$ on the $(n - 1)$-category $\cl{P}$. This definition is based on a generalization of a result proved by Baues and Jibladze, see~\cite{BauesJibladze02} for the case $n=2$, stating that an $(n,n-1)$-category is abelian if and only if it is linear, where $(n,n-1)$-categories correspond to a notion of ``globular crossed module'' for $n$-categories.
When the polygraph~$P$ is convergent, the natural system $\Pi(P)$ is generated by the generating confluences of $P$~\cite[Proposition 2.4.2]{GuiraudMalbos10smf}.

\index{FDTab@$\FDTAB$ (abelian finite derivation type)}
\index{finite derivation type!abelian}
A notion of \emph{abelian finite derivation type}, $\FDTAB$ for short, can also be defined for $n$-polygraphs. An $n$-polygraph $P$ has $\FDTAB$ when the abelian $(n,n-1)$\nbd-category 
$\fab{\freegpd{P}}$ admits a finite acyclic cellular extension. 
As for case $n=2$ in \cref{Proposition:IartfIffFdtab}, we prove that an $n$\nbd-polygraph is $\FDTAB$ if and only if the natural system~$\Pi(P)$ of identities among relations of $P$ is finitely generated~\cite[Proposition 5.7.2]{GuiraudMalbos12advances}.

\addcontentsline{toc}{part}{Appendix}
\part*{Appendix}

\appendix
\setcounter{tocdepth}{1}
\addtocontents{toc}{\protect\setcounter{tocdepth}{1}}


\chapter{A Catalogue of 2-Polygraphs}
\label{chap:pres-mon}
\label{chap:2ex}
In this chapter, we list some examples of presentations of monoids and
categories by a 2\nbd-poly\-graph.

\section{Presentations of Monoids}
The richest source of presentations of categories can be found in presentations
of monoids (and groups), which are seen here as a particular instance of a
category.

\subsection{Monoid}
\label{sec:monoid}
A \emph{monoid}\index{monoid} $(M,\times,1)$ is a set equipped with a binary
operation~$\times$ which is associative and admits~$1$ as neutral element. A
morphism between two monoids is a function between the underlying sets which
preserves multiplication and neutral element.

The following lemma shows that one can always see a monoid as a particular case
of a category with only one element, conventionally denoted $\star$.

\begin{lemma}
  \nomenclature[Mon]{$\Mon$}{category of monoids}
  \label{lem:mon-cat}
  The category $\Mon$ of monoids is isomorphic to the full subcategory of~$\Cat$
  whose objects are the categories with~$\star$ as only object.
\end{lemma}
\begin{proof}
  To a monoid $(M,\times,1)$ we can associate a category $BM$, sometimes called
  the \emph{delooping} of~$M$, with $\star$ as only object, the elements
  $a\in M$ as morphisms $a:\star\to\star$, composition being given by
  $b\circ a=a\times b$, with $1$ as identity on~$\star$. This construction is
  easily extended as a functor which is an isomorphism of categories.
\end{proof}

\noindent
A presentation of a monoid consists of a $2$-polygraph~$P$ whose set of
$0$\nbd-gene\-ra\-tors is reduced to one element $P_0=\set{\star}$. A group
being a monoid, presentations of groups provide many such
examples~\cite{coxeter1972generators}; we mostly restrict here to those which
are ``really monoids'', meaning that the presence of inverses does not play a
crucial role in the presentation.

\subsection{Notations}
In a polygraph $P=\Pres{\star}{P_1}{P_2}$,
we generally omit mentioning the source and target of 1-generators since they
are necessarily~$\star$. Moreover, as far as we are concerned with the presented
category, the names of the 2-generators, as well as their orientation, will not
be relevant. The 2-polygraph
\[
  \Pres{\star}{a:\star\to\star}{\alpha:aa\To 1}
\]
will thus often be simply noted
\begin{equation}
  \label{eq:pres-bool}
  \Pres{\star}{a}{aa=1}
  \pbox.
\end{equation}
Traditionally, the set~$P_0$ is even omitted, \ie the above presentation is
noted $\pres{a}{aa=1}$, we will however refrain from doing so in order to avoid
a possible confusion with a 1-polygraph. We sometimes write the indices at the
bottom right of the presentation. For instance
\[
\Pres{\star}{a_i}{a_ja_i=a_ia_j}_{i,j\in\N}
\]
denotes a presentation with
\[
  P_1=\setof{a_i}{i\in\N}
  \qqtand
  P_2=\setof{a_ja_i=a_ia_j}{i,j\in\N}.
\]

\subsection{Natural numbers}
\label{sec:pres-nat}
\index{monoid!of natural numbers}
\index{natural number}
\nomenclature[N]{$\N$}{natural numbers}
The additive monoid~$\N$ of natural numbers admits the presentation
\[
  \Pres{\star}{a}{}
\]
where~$a$ corresponds to the natural number~$1$, and more generally $a^n$ to
$n\in\N$.

\subsection{Cyclic monoids}
\label{sec:N-mod}
\index{cyclic monoid}
\index{monoid!cyclic}
The additive monoid of booleans~$\N/2\N$ admits the
presentation~\eqref{eq:pres-bool} above. More generally, given $n\in\N$, the
additive monoid~$\N/n\N$ admits the presentation
\[
  \Pres{\star}{a}{a^n=1}
  \pbox.
\]

\subsection{Booleans}
\index{boolean}
\index{monoid!boolean}
\nomenclature[B]{$\Bool$}{booleans}
The set~$\Bool$ of booleans respectively consists of the two elements~$\bot$
(standing for false) and~$\top$ (standing for true). We respectively write
$\land$, $\lor$, and $\times$ for conjunction, disjunction, and exclusive
disjunction. The monoid~$(\Bool,\times,\bot)$ is isomorphic to $\N/2\N$ and thus
admits~\eqref{eq:pres-bool} as presentation. The monoids~$(\Bool,\lor,\bot)$ and
$(\Bool,\land,\top)$ are isomorphic and admit the presentation
\[
  \Pres{\star}{a}{aa=a}
  \pbox.
\]

\subsection{Free monoids}
\index{free!monoid}
\index{monoid!free}
Fix a set~$X$. A \emph{word}~$u$ over~$X$ is a finite sequence $u=a_1\ldots a_n$
of elements~$a_i$ of $X$. Given two words $u=a_1\ldots a_m$ and
$v=b_1\ldots b_n$, their \emph{concatenation} is the word
$uv=a_1\ldots a_mb_1\ldots b_n$. The \emph{free monoid}~$X^*$ over~$X$ can be
described as the monoid of words over~$X$, with concatenation as composition and
empty word as unit. It admits the presentation
\[
  \Pres{\star}{X}{}
  \pbox.
\]

\subsection{Free commutative monoids}
\label{sec:multiset}
\index{multiset}
\index{free!commutative monoid}
\index{monoid!free!commutative}
Fix a set~$X$. Recall from \cref{sec:multisets} that a \emph{multiset} over~$X$ is a function $\mu:X\to\N$, assigning a
\emph{multiplicity} to an element $a$ of~$X$, such that every element~$a$ of $X$
has a null multiplicity, excepting for a finite number. Given two multisets $\mu$
and $\nu$, we write $\mu\mcup\nu$ for their pointwise sum. The constant function equal
to~$0$ is called the \emph{empty multiset}. The \emph{free commutative monoid}
over~$X$ can be described as the monoid of multisets with $\mcup$ as
composition and empty multiset as unit. It admits the presentation
\[
  \Pres{\star}{X}{ab=ba}_{a,b\in X}
  \pbox.
\]
For instance, the multiplicative monoid~$\N\setminus\set{0}$ is the free
commutative monoid over the set~$\N$, since an element of this monoid can be
interpreted as the multiset of its prime factors. It admits the presentation
\[
  \Pres{\star}{a_i}{a_ia_j=a_ja_i}_{i,j\in\N}
\]
where the 1-generator $a_i$ corresponds to the $i$-th prime number, \ie $a_0=2$,
$a_1=3$, $a_2=5$, etc.

\subsection{Partially commutative monoids}
\index{partially commutative monoid}
\index{monoid!partially commutative}
Suppose given an alphabet $X$ and $I\subseteq X\times X$ a reflexive and
symmetric relation, called \emph{independence}. The monoid with presentation
\[
  \Pres{\star}{X}{ab=ba}_{(a,b)\in I}
\]
is called a \emph{trace monoid}, a \emph{partially commutative monoid}, or a
\emph{heap monoid}. This family of monoids was introduced by Cartier and
Foata~\cite{cartier1969problemes} in order to study combinatorial problems of
rearrangements and also studied in computer
science~\cite{diekert1995book,fajstrup2016directed,mazurkiewicz1977concurrent} for the following
reason. The elements of~$X$ can be interpreted as actions, and thus $X^*$ as the
set of possible sequences of actions. Sometimes, the order in which two
actions~$a$ and~$b$ are performed does not matter, for instance when $a$ and $b$
consist in reading or writing at disjoint positions in the memory in a computer program: in the end,
executing~$a$ then~$b$, or~$b$ then~$a$, will lead to the same state. This kind
of situation is naturally modeled by having $(a,b)\in I$ and occurs when
considering concurrent processes.

\subsection{Baumslag–Solitar monoids}
\index{Baumslag–Solitar monoid}
\index{monoid!Baumslag–Solitar}
Given natural numbers~$m$ and~$n$, the Baumslag-Solitar monoid
$BS(m,n)$~\cite{baumslag1962some} is the monoid presented by
\[
  \Pres{\star}{a,b}{ab^m=b^na}
\]
which can be seen as a variant of the commutativity relation described above,
$BS(1,1)$ being the free commutative monoid on two generators. The enveloping
groups of those monoids provide examples of non-Hopfian groups, e.g., $BS(2,3)$.

\subsection{Finite subsets}
\index{powerset}
\index{monoid!powerset}
\nomenclature[P]{$\powerset$}{powerset}
Given a set~$I$, the set $\powerset_{\mathrm{fin}}(I)$ of finite subsets of~$I$
is a monoid with union as multiplication and empty set as neutral element. It
admits as presentation
\[
  \Pres{\star}{a_i}{a_ia_j=a_ja_i,a_ia_i=a_i}_{i,j\in I}
  \pbox.
\]
It is thus the free idempotent commutative monoid on the set~$I$.

\subsection{Coproduct and product}
\label{sec:mon-coprod-prod}
\index{coproduct!of monoids}
\index{product!of monoids}
\index{monoid!coproduct}
\index{monoid!product}
\index{free product}
The monoids $\N\sqcup\N$ and $\N\times\N$ respectively admit the presentations
\[
  \Pres{\star}{a,b}{}
  \qqqtand
  \Pres{\star}{a,b}{ba=ab}
  \pbox.
\]
More generally, given 2-polygraphs $P$ and $Q$, respectively presenting
monoids~$M$ and~$N$,
\begin{itemize}
\item their coproduct (or \emph{free product}) $M+N$ is presented by
  \[
    \Pres{\star}{P_1,Q_1}{P_2,Q_2}
    \pbox,
  \]
\item their product $M\times N$ is presented by
  \[
    \Pres{\star}{P_1,Q_1}{P_2,Q_2,R_2}
  \]
  where
  \[
    R_2=\setof{ba=ab}{a\in P_1,b\in Q_1}
    \pbox.
  \]
\end{itemize}
These constructions are detailed and generalized in \chapr{2-op}.

\subsection{Bicyclic monoid}
\label{sec:bicyclic}
\index{bicyclic monoid}
\index{monoid!bicyclic}
\index{Dyck word}
\index{well-bracketed word}
The \emph{bicyclic monoid} is the monoid presented by
\[
  \Pres{\star}{a,b}{ab=1}
  \pbox.
\]
If we read $a$ as an ``opening bracket'' and $b$ as a ``closing bracket'', then
the words in the equivalence class of~$1$ are precisely the \emph{well-bracketed
  words} (also called \emph{Dyck words}). From this property follow applications
in combinatorics and computer science (e.g., this monoid is the syntactic monoid
of the language of Dyck words).

\subsection{Integers}
\label{sec:int-monoid}
\index{integer}
\index{monoid!of integers}
The additive group $\Z$ admits, as a monoid, the presentation
\[
  \Pres{\star}{a,b}{ab=1,ba=1}
\]
where $a$ and $b$ can respectively be interpreted as $1$ and $-1$.

\subsection{Enveloping group}
\label{sec:env-group-pres}
\index{group!enveloping}
\index{enveloping!group}
The forgetful functor $\Grp\to\Mon$, from the category of groups to the category
of monoids, admits a left adjoint, under which the image of a monoid is called
the associated \emph{enveloping group}.
For instance, the additive group~$\Z$ is the enveloping group of the additive
monoid~$\N$. As such, the previous example is an instance of a more general
construction: given a polygraph~$P$ presenting a monoid~$M$, its enveloping
group~$G$ is presented as a monoid by
\[
  \Pres{\star}{P_1,\finv{P_1}}{P_2,R_2}
\]
where
\[
  \finv{P_1}=\setof{\finv{a}}{a\in P_1}
  \qqtand
  R_2=\setof{a\finv{a}=1,\finv{a}a=1}{a\in P_1}
  \pbox.  
\]
This construction is detailed and generalized in \secr{1-pres-free-gpd}.
By abuse of notation, given a word $u=a_1\ldots a_n$ in $\freecat{P_1}$, we
sometimes write $\finv u=\finv a_n\ldots\finv a_1$.

Two distinct monoids can generate isomorphic enveloping groups. For instance,
consider the monoid presented by
\[
  \Pres{\star}{a,b}{ab=bba,ba=aab}\pbox.
\]
This monoid is not the trivial one since for instance the classes of $1$, $a$,
and $b$ contain only one element (there is no derivable relation involving
those, the generating relations being between words of length~$2$ and~$3$). We
can observe that the relation $ab=baab$ is derivable since $ab=bba=baab$, and
similarly the relation $ba=abba$ is also derivable. In the enveloping group, we
thus have $ba=baab\finv{b}\finv{a}=ab\finv{b}\finv{a}=1$, and similarly $ab=1$.
Therefore, $1=ab=bba=b1=b$ and similarly $1=a$. The presented group is therefore
the trivial one, which is also the enveloping group of the trivial monoid.
Another example is given in \cref{sec:wirtinger}.

\subsection{Group presentation}
\index{presentation!of a group}
As a variant of the previous construction, by a \emph{group presentation}, one
usually understands a pair
\[
  \pres{P_1}{P_2}
\]
where~$P_1$ is a set and $P_2\subseteq\freecat{(P_1\cup\finv{P_1})}$, which
implicitly means the group presented by the 2-polygraph
\[
  \Pres{\star}{P_1,\finv{P_1}}{\finv{a}a=1,a\finv a=1,u=1}_{a\in P_1,u\in P_2}
  \pbox.
\]
Here, we only consider relations of the form $u=1$, since having a
relation $u=v$ is equivalent (by Tietze transformations) to having the
relation~$\finv uv=1$.

\subsection{Positive rational numbers}
We have seen a presentation of the multiplicative monoid
$\N\setminus\set0$ in \cref{sec:multiset}. The multiplicative group $\mathbb{Q}^{>0}$ of strictly
positive rational numbers is its enveloping group and, from the construction of
\cref{sec:env-group-pres}, we deduce that it admits the presentation
\[
  \Pres{\star}{a_i,\finv a_i}{a_ia_j=a_ja_i,\finv a_ia_i=1,a_i\finv a_i=1}_{i,j\in\N}
  \pbox.
\]

\subsection{Non-negative rational numbers}
\label{sec:Qplus}
Let us detail an example of a presentation coming
from~\cite[Section~5.7]{johnson1997presentations}. The \emph{additive} monoid
$\mathbb{Q}^+$ of non-negative rational numbers can be presented by the
2-polygraph
\[
  P=\Pres{\star}{a_i}{(a_{i+1})^{i+1}=a_i}_{i\in\N\setminus\set{0}}
  \pbox.
\]
Here, any pair of generators commute, \ie the relations $a_ja_i=a_ia_j$ are
derivable since, supposing $j>i$, we have
\[
  a_ja_i
  =
  a_ja_j^{j(j-1)\ldots(i+1)}
  =
  a_j^{j(j-1)\ldots(i+1)}a_j
  =
  a_ia_j
  \pbox.
\]
The interpretation of a generator $a_i$ in~$\mathbb{Q}^+$ is $f(x_i)=1/i!$,
which extends as a morphism of monoids $f:\freecat{P}\to\mathbb{Q}^+$. It is
compatible with the relations since
\[
  f(a_{i+1}^{i+1})
  =
  \frac{i+1}{(i+1)!}
  =
  \frac1{i!}
  =
  f(a_i)
\]
and thus induces a morphism $f:\prescat{P}\to\mathbb{Q}^+$, which we prove to be an
isomorphism. The function~$f$ is surjective since for
$(p,q)\in\N\times\pa{\N\setminus\set{0}}$, one has
\[
  \frac pq
  =
  \frac{p(q-1)!}{q!}
  =
  f\pa{a_q^{p(q-1)!}}
  \pbox.
\]
Let us show that it is injective. Suppose given two words $u,v\in\freecat{P_1}$
such that $f(u)=f(v)$. Since the generators commute, the words $u$
and $v$ are, up to equivalence, of the form
\[
  u=a_1^{m_1}a_2^{m_2}\ldots a_k^{m_k}
  \qqqtand
  v=a_1^{n_1}a_2^{n_2}\ldots a_{k}^{n_k}
\]
with $k\in\N$ and $m_i,n_i\in\N$ for $0\leq i\leq k$. We can also suppose that
$m_k\neq 0$, up to exchanging the roles of $u$ and $v$ and lowering~$k$. Suppose
moreover that $k>1$ (this hypothesis is not innocuous since it will lead to a
contradiction). If $n_k\geq k$, we can use the relation $a_k^k=a_{k-1}$ to
transform $u$ into a word $u'$ equivalent to $u$ and satisfying $0\leq
m_k<k$. For this reason, we can also suppose, without loss of generality, that
$0<m_k<k$ and $0\leq n_k<k$. Up to exchanging the roles of~$u$ and~$v$, we can
also suppose $m_k\geq n_k$. We have
\[
  0
  =
  f(u)-f(u')
  =
  \frac{m_1-n_1}{1!}
  +
  \frac{m_2-n_2}{2!}
  +
  \ldots
  +
  \frac{m_k-n_k}{k!}
\]
and therefore
\begin{align*}
  0
  &=
    (k-1)!\pa{f(u)-f(u')}
  \\
  &=
    (k-1)!
    \pa{
    \frac{m_1-n_1}{1!}
    +
    \frac{m_2-n_2}{2!}
    +
    \ldots
    +
    \frac{m_{k-1}-n_{k-1}}{(k-1)!}
    }
    +
    \frac{m_k-n_k}{k}
  \\
  &=
    m+\frac{m_k-n_k}{k}
\end{align*}
with $m\in\Z$. Since $0<m_k<k$, $0\leq n_k<k$ and $m_k\geq n_k$, we have
$0<(m_k-n_k)/k<1$, and thus a contradiction. Therefore, we have $k=1$ and
$m_k=f(u)=f(v)=n_k$, \ie $u=v$. The morphism $f:\prescat{P}\to\mathbb{Q}^+$ is
thus injective. From the results of \secr{env-group-pres}, we can deduce the
following presentation of the additive monoid~$\mathbb{Q}$, which is its
enveloping group:
\[
  \Pres{\star}{a_i,b_i}{a_{i+1}^{i+1}=a_i,a_ib_i=1,b_ia_i=1}_{i\in\N\setminus\set{0}}
  \pbox.
\]

\subsection{Symmetric groups}
\label{sec:sym-group}
\label{sec:sym-group-pres}
\index{symmetric!group}
\index{group!symmetric}
\index{permutation}
\nomenclature[S]{$S_n$}{symmetric group}
The \emph{symmetric group}~$S_{n+1}$ is the group of bijections (or
\emph{permutations}) from the set $\set{0,\ldots,n}$ to itself, with
multiplication given by composition and unit by identity. More geometrically, it
can also be defined as the group of symmetries of an $n$-simplex. Considered as
a monoid, it admits a presentation where the set of generators is
$P_1=\set{a_0,\ldots,a_{n-1}}$ and the relations are
\begin{itemize}
\item $a_ia_i=1$, for $0\leq i<n$,
\item $a_ia_{i+1}a_i=a_{i+1}a_ia_{i+1}$, for $0\leq i<n-1$,
\item $a_ia_j=a_ja_i$, for $0\leq i<j<n$ with $j>i+1$,
\end{itemize}
as found out by Moore~\cite{moore1896concerning}.
Here, the generator $a_i$ should be interpreted as the transposition exchanging
$i$ with $i+1$, written $(i(i+1))$, whose graph can be pictured as
\[
  \satex{sym-ai}
  \,\pbox.
\]
The presentation is studied in detail in~\secr{pres-sym-group}.

The above presentation is based on the transpositions $(i(i+1))$ as
generators~$a_i$, but other choices of generators are possible and will give rise
to other presentations~\cite[Section~6.2]{coxeter1972generators}, for instance:
\begin{itemize}
\item the generators $a_i=(in)$, with $0\leq i<n$,
  induce a presentation with relations
  \begin{align*}
    a_i^2&=1
    &
    (a_ia_{i+1})^3&=1
    &
    (a_ia_{i+1}a_ia_j)^2&=1
  \end{align*}
  for $1\leq i,j\leq n-1$ with $j\neq i$ and $j\neq i+1$, where $a_n=a_0$ by
  convention,
\item the generators $a_0=(1n)$ and $a_i=(0in)$ for $1\leq i<n$
  induce a presentation with relations
  \begin{align*}
    a_0^2&=1
    &
    a_j^3&=1
    &
    (a_ia_j)^2&=1
  \end{align*}
  for $0\leq i<j<n$,
\item the two generators $a=(01)$ and $b=(012\ldots n)$ induce a presentation
  with relations
  \begin{align*}
    a^2&=1
    &
    b^{n+1}&=1
    &
    (ba)^n&=1
    &
    (ab^nab)^3&=1
    &
    (ab^{n+1-j}ab^j)&=1
  \end{align*}
  for $1\leq j<n-1$.
\end{itemize}

\subsection{Alternating groups}
\label{sec:alt-group}
\index{alternating group}
\index{group!alternating}
\nomenclature[A]{$A_n$}{alternating group}
The \emph{alternating group}~$A_{n+1}$ is the subgroup of~$S_{n+1}$ consisting
of symmetries of even signature (recall that the \emph{signature} of a symmetry
is the parity of the number of transpositions used to express it). It admits a
presentation with generators $P_1=\set{a_0,\ldots,a_{n-1}}$ and relations
\begin{align*}
  a_j^3&=1
  &
  (a_ia_j)^2&=1
\end{align*}
for $0\leq i<j<n-1$.
Here, a generator $a_i$ should be interpreted as the permutation
$(i(n-1)n)$, see~\cite[Section~6.3]{coxeter1972generators} for details.

\subsection{Braid groups and monoids}
\label{sec:braid-mon}
\index{braid!group}
\index{group!braid}
\index{monoid!braid}
\nomenclature[B+]{$B_n^+$}{braid monoid}
\nomenclature[B]{$B_n$}{braid group}
We introduce notations for the spaces $I=[0,1]$, $X=\R^2$ and $Y=I\times X$. Suppose given
$n\in\N$ and continuous functions
\[
  b_i
  :
  I
  \to
  Y
\]
with $0\leq i<n$, which are mutually disjoint (\ie for $t\in I$, $b_i(t)=b_j(t)$
implies $i=j$) and with fixed endpoints (say, $b_i(0)=(0,i,0)$ and
$b_i(1)=(1,i,0)$). This induces a subspace $\beta$ of $Y=I\times X$ defined as
\[
  \beta
  =
  \setof{(t,b_i(t))}{t\in I, 0\leq i<n}
  \pbox.
\]
Such a subspace is called a (\emph{geometric}) \emph{braid} with $n$ strands.
Graphically, a braid with three strands can be pictured as follows, where the first
coordinate is pictured vertically:
\[
  \fig[height=25mm]{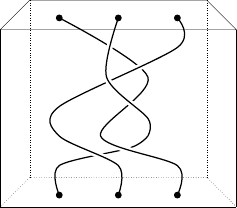}
\]
Braids are considered up to endpoint-preserving isotopy: we identify two
braids~$\beta$ and~$\beta'$ for which there exists a continuous function
\[
  h
  :
  I
  \to
  Y^Y
\]
such that $h(0)=\id_Y$, $h(1)(\beta)=\beta'$ and, for every $t\in I$, the
function $h(t)$ is a homeomorphism, whose restriction to
$(\set{0}\times X\sqcup\set{1}\times X)\subseteq Y$ is the identity, such that
$h(t)(\beta)$ is a braid. The set~$B_n$ of braids with $n$ strands (up to
isotopy) forms a group, where the composition is given by concatenation of
strands, and the identity is given by the braid induced by constant functions
$b_i(t)=(i,0)$.

An alternative description can be given as follows. Given a topological
space~$X$ (typically, $X=\R^2$ in the following), the $n$-element
\emph{configuration space}\index{configuration space} is
\[
  C_nX
  =
  \setof{(x_1,\ldots,x_n)}{\text{$x_i=x_j$ implies $i=j$}}
  \pbox.
\]
The $n$-element \emph{unlabeled configuration space}~$D_nX=C_nX/\Sigma_n$ is the
quotient of~$C_nX$ under the action of the symmetric group permuting the
coordinates. The $n$-strand \emph{braid group} can equivalently be defined as
the fundamental group of this space: $B_n=\pi_1(D_n\R^2)$.

Before presenting the braid group, we first present a submonoid of this group.
The (\emph{positive}) \emph{braid monoid}~$B_{n+1}^+$ with $n+1$ strands, admits
a presentation with $P_1=\set{a_0,\ldots,a_{n-1}}$, where $a_i$ can be pictured
as
\[
  \satex{br-ai}
\]
(with $i$ strands on the left), and the relations are
\begin{itemize}
\item $a_ia_{i+1}a_i=a_{i+1}a_ia_{i+1}$, for $0\leq i<n-1$,
\item $a_ia_j=a_ja_i$, for $0\leq i<i+1<j<n$.
\end{itemize}
For instance, with $n+1=4$, the generators $a_0$, $a_1$, $a_2$ respectively
correspond to the following braids with four strands:
\begin{align*}
  a_0&=\satex{br-a0}
  &
  a_1&=\satex{br-a1}
  &
  a_2&=\satex{br-a2}
\end{align*}
and the relations are
\begin{align*}
  \satex{br-a010}&=\satex{br-a101}
  &
  \satex{br-a121}&=\satex{br-a212}
  \\
  a_0a_1a_0&=a_1a_0a_1
  &
  a_1a_2a_1&=a_2a_1a_2
\end{align*}
and
\begin{align*}
  \satex{br-a02}&=\satex{br-a20}
  \\
  a_0a_2&=a_2a_0
\end{align*}
We can expect that the full braid group can be recovered by moreover adding the
other kind of crossings
\[
  \finv{a_i}
  =
  \satex{br-ai-inv}
\]
which are inverse to the crossings~$a_i$. Indeed, the braid group~$B_{n+1}$ is
the enveloping group of the monoid~$B^+_{n+1}$, from which a presentation can be
constructed, see \secr{env-group-pres}:
\begin{itemize}
\item generators are $a_i,\finv{a_i}$ for $0\leq i<n$,
\item relations are those of the braid monoid plus $a_i\finv{a_i}=1$ and
  $\finv{a_i}a_i=1$ for~$0\leq i<n$.
\end{itemize}
This presentation is due to Artin~\cite{artin1925theorie}.

From the above presentation, we see that the symmetric group can be obtained
from the braid group by identifying $a_i$ with $\finv{a_i}$. There is thus a
projection morphism $B_n\to S_n$ and we write $P_n$ for its kernel, which is
called the \emph{pure braid group}: the group~$P_n$ can be seen as the submonoid
of $B_n$ consisting of braids which become identities if we interpret them as
symmetries, and we have an exact sequence
\[
  1\to P_n\to B_n\to S_n\to 1
  \pbox.
\]
This group can also more directly be described as the fundamental group of a
configuration space: $P_n=\pi_1(C_n\R^2)$.
A presentation of $P_n$ can be given with generators $a_{ij}$ for
$0\leq i<j<n$, as well as their formal inverses~$\finv{a_{ij}}$. Considering
$P_n$ as a subgroup of~$B_n$, the generators $a_{ij}$ can be expressed in terms
of the generators of~$B_n$ as
\[
  a_{ij}
  =
  a_ia_{i+1}\ldots a_{j-2}a_{j-1}^2\finv{a_{j-2}}\ldots\finv{a_{i+1}}\finv{a_i}
  \pbox.
\]
Graphically,
\[
  \satex{br-pure}
\]
The relations satisfied by those are
\[
  a_{kl}a_{ij}\finv{a_{kl}}
  =
  \begin{cases}
    a_{ij}&\text{if $l<i$ or $j<k$,}\\
    \finv{a_{il}}a_{ij}a_{il}&\text{if $i<j=k<l$,}\\
    \finv{a_{ij}}\finv{a_{ik}}a_{ij}a_{ik}a_{ij}&\text{if $i<k<j=l$,}\\
    \finv{a_{il}}\finv{a_{ik}}a_{il}a_{ik}a_{ij}\finv{a_{ik}}\finv{a_{il}}a_{ik}a_{il}&\text{if $i<k<j<l$.}
  \end{cases}
\]

\subsection{Hyperoctahedral group}
\label{sec:hyperoctahedral}
\index{hyperoctahedral!group}
\index{group!hyperoctahedral}
The hyperoctahedral group~$B_n$ (or~$C_n$) is the group of symmetries of a cube
of dimension~$n$. It admits a presentation with generators $a_i$, with
$0\leq i<n$ and relations
\begin{itemize}
\item $a_0a_1a_0a_1=a_1a_0a_1a_0$
\item $a_ia_{i+1}a_i=a_{i+1}a_ia_{i+1}$ for $0<i<n-1$,
\item $a_ia_j=a_ja_i$ for $0\leq i<i+1<j<n$.
\end{itemize}
This group can also be described as the group of signed symmetries, \ie
$n{\times}n$ orthogonal matrices with integer entries, \ie $n{\times}n$ matrices
with coefficients in $\set{-1,0,1}$ containing exactly one non-null coefficient
on each row and each column, see \secr{hyperoctahedral-cat} for details.

\subsection{Progressive ribbon group}
\label{sec:progressive-ribbon-group}
\index{progressive ribbon!group}
\index{group!of progressive ribbons}
Given $n\in\N$, the monoid~$R^+_{n+1}$ of positive progressive ribbons with
$n+1$ strands admits a presentation with generators $a_i$, with $0\leq i<n$, and
$b_i$, with $0\leq i\leq n$, subject to the relations
\begin{itemize}
\item $a_ia_{i+1}a_i=a_{i+1}a_ia_{i+1}$ for $0\leq i<n-1$,
\item $a_ia_j=a_ja_i$ for $0\leq i<i+1<j<n$,
\item $b_ia_i=a_ib_{i+1}$ for $0\leq i<n$,
\item $b_{i+1}a_i=a_ib_i$ for $0\leq i<n$,
\item $b_ia_j=a_jb_i$ for $0\leq i<n$, $0\leq j\leq n$, and $i\not\in\set{i,i+1}$.
\end{itemize}
Graphically, an element of this monoid can be depicted as $n$ ribbons in the
space, the generators are respectively
\begin{align*}
  a_i&=\fig[height=1cm]{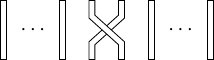}
  &
  b_i&=\fig[height=1cm]{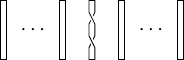}  
\end{align*}
and the relations have natural graphical interpretations, which the reader is
encouraged to draw. The group~$R_n$ of progressive ribbons is the enveloping
group of~$R_n^+$. Note that the relations satisfied by the generators $a_i$ are
precisely the relations of the braid group.  In fact, there is an obvious action
of the braid group $B_n$ on the set $\intset{n}=\set{0,\ldots,n-1}$ of braids
and $R_n$ can be obtained as the wreath product $R_n=\Z\wr B_n$ of the additive
group of integers with the braid group.
The groups $R_n$, as well as their presentation, are detailed
in~\cite{joyal1991geometry}, as well as in \secr{progressive-ribbons}.

\subsection{Dihedral groups}
\label{sec:dihedral-group}
\index{dihedral group}
\index{group!dihedral}
\nomenclature[D]{$D_n$}{dihedral group}
The dihedral group~$D_n$ is the group of symmetries of a regular polygon with
$n$ vertices. It admits a presentation by the 2-polygraph with 1-generators
$r_0,\ldots,r_{n-1}$ and $s_0,\ldots,s_{n-1}$, and relations
\begin{align*}
  r_0&=1&
  r_ir_j&=r_{i+j}&
  r_is_j&=s_{i+j}&
  s_ir_j&=s_{i-j}&
  s_is_j&=r_{i-j}
\end{align*}
for every indices $0\leq i,j<n$, with addition and subtraction taken
modulo~$n$. Alternatively, it can be presented by the 2-polygraph
\[
  \Pres{\star}{r,s}{r^n=1,s^2=1,rsrs=1}
\]
(the relationship with previous generators is given by $r=r_1$, $s=s_0$,
$r_i=r^i$ and $s_i=r^is$).
In terms of symmetries of a polygon with $n$ vertices, $s$ corresponds to a
symmetry and~$r$ a rotation of $2\pi/n$, and more generally $s_i$ and $r_i$
respectively correspond to a reflection along the $i$-th axis and a rotation
of~$2\pi i/n$:
\[
  \fig{dihedral-polygons}
\]
The case of $D_3$ is detailed in \cref{ex:S3-two-presentations}.

\subsection{Artin monoids and Coxeter groups}
\label{sec:artin-monoid}
\label{sec:coxeter-group}
\index{Artin!monoid}
\index{monoid!Artin}
\index{group!Artin}
\index{Coxeter group}
\index{group!Coxeter}
Given a finite set~$S$ of generators, a \emph{Coxeter matrix}~$M$ is a function
which to every pair of elements $(s,t)\in S\times S$ associates
$m_{st}\in\N\uplus\set\infty$ such that $m_{st}=m_{ts}$, $m_{ss}=1$, and
$m_{st}>1$ for every $s,t\in S$ with $s\neq t$. Such a matrix induces
\begin{itemize}
\item an \emph{Artin monoid} with presentation
  \newcommand{\altprod}[1]{\langle#1\rangle}
  \[
    \Pres{\star}{s}{\altprod{st}^{m_{st}}=\altprod{ts}^{m_{ts}}}_{s,t\in S}
  \]
  where $\altprod{st}^k$ denotes the alternating product of $s$ and $t$ of
  length $k$ starting with~$k$ (e.g., $\altprod{st}^5=ststs$) and, by convention,
  there is no relation when $m_{st}=\infty$ or $m_{ts}=\infty$,
\item an \emph{Artin group}, which is the group freely generated by the above
  monoid, and
\item a \emph{Coxeter group} with presentation as a monoid
  \[
    \Pres{\star}{s}{\altprod{st}^{m_{st}}=\altprod{ts}^{m_{ts}},ss=1}_{s,t\in S}
  \]
  or equivalently
  \[
    \Pres{\star}{s}{(st)^{m_{st}}=1}_{s,t\in S}
    \pbox.
  \]
\end{itemize}
Since the set~$S$ is finite, writing $n$ for its cardinal, we can assume that
it is of the form $S=\setof{i\in\N}{0\leq i<n}$, which is convenient in the following.
A Coxeter matrix~$M$ is often pictured as a Dynkin diagram consisting of a
labeled non-oriented simple graph with~$n$ vertices~$x_i$, with $0\leq i<n$,
with an edge between~$x_i$ and~$y_i$ labeled by $m_{ij}$. By convention, the looping
edges on vertices, as well as edges with label~$2$, are not drawn, and edges with
label~$3$ are not labeled. Artin monoids are further detailed in
\cref{Section:CoherentPresentationArtinTitsMonoids}.

For instance, the Coxeter matrix such that $m_{i(i+1)}=3$, and $m_{ij}=2$
for~$|i-j|>2$, can be represented by the diagram
\[
  \begin{tikzpicture}
    \foreach \i in {0,1,2,4,5} {
      \filldraw (\i,0) circle (.05);
    }
    \draw (3,0) node {$\ldots$};
    \draw (0,0) -- (2.7,0);
    \draw (3.3,0) -- (5,0);
  \end{tikzpicture}
\]
The Artin monoid is the braid monoid with $n+1$ strands (see \secr{braid-mon})
and the Coxeter group is the symmetric group on~$n+1$ elements (see
\secr{sym-group-pres}). Similarly, the hyperoctahedral group~$B_n$ is the
Coxeter group associated to
\[
  \begin{tikzpicture}
    \foreach \i in {0,1,2,4,5} {
      \filldraw (\i,0) circle (.05);
    }
    \draw (3,0) node {$\ldots$};
    \draw (0,0) -- (2.7,0);
    \draw (3.3,0) -- (5,0);
    \draw (.5,0) node[above] {$4$};
  \end{tikzpicture}
\]
(with $n$ vertices) and the dihedral group~$D_n$ is the Coxeter group associated to
\[
  \begin{tikzpicture}
    \foreach \i in {1,2,4,5} {
      \filldraw (\i,0) circle (.05);
    }
    \filldraw (0,.5) circle (.05) -- (1,0);
    \filldraw (0,-.5) circle (.05) -- (1,0);
    \draw (3,0) node {$\ldots$};
    \draw (1,0) -- (2.7,0);
    \draw (3.3,0) -- (5,0);
  \end{tikzpicture}
\]
(with $n$ vertices). A detailed presentation of these groups can be found in
several places such as~\cite{geck2000characters}.

\subsection{Plactic, Chinese, and Sylvester monoids}
\index{plactic monoid}
\index{monoid!plactic}
\index{Chinese monoid}
\index{monoid!Chinese}
\index{Sylvester monoid}
\index{monoid!Sylvester}
\label{sec:plactic-monoid}
\nomenclature[P]{$P_n$}{plactic monoid}
\nomenclature[C]{$C_n$}{Chinese monoid}
The \emph{plactic monoid}~$P_n$ of rank $n$ in type $A$ was introduced by Knuth
in~\cite{knuth1970permutations} and further developed by Lascoux and
Schützenberger~\cite{lascoux1981monoide}. It admits a presentation with
$\{1,\ldots,n\}$ as generators, and relations
\[
zxy=xzy, \; \text{for $1\leq x\leq y<z \leq n$}, 
\qquad
yzx=yxz, \; \text{for $1\leq x<y\leq z \leq n$}.
\]
Its elements are in bijection with semistandard Young tableaux~\cite{knuth1970permutations,lothaire2002algebraic}, from which stem applications in representation theory~\cite{fulton1997young, littelmann1996plactic}. Plactic monoids are further detailed in \cref{S:CoherentPresentationPlacticMonoids}.

Many variants of plactic monoids exist, such as the Chinese monoids.
These monoids were introduced in~\cite{duchamp1994plactic} through the following
presentation, called Chinese presentation. The \emph{Chinese monoid of rank
  $n>0$}, denoted by~$Ch_n$, is the quotient of the free monoid~$\{1,\ldots,n\}^{\ast}$ by
the congruence generated by the following relations:
\begin{equation}
  \label{ChineseRelations}
  zyx = zxy = yzx \quad \text{for all} \quad 1\leq x\leq y \leq z \leq n.
\end{equation}

The \emph{Sylvester monoids}~$Syl_n$~\cite{hivert2005algebra} are also a variant, with the same generators and relations
\[
  zxuy=xzuy,
\]
for $1\leq x\leq y<z \leq n$ and $u\in\freecat{\{1,\ldots,n\}}$.

\subsection{Thompson group~$F$}
\label{sec:Thompson-F}
\index{Thompson!group}
\index{group!Thompson}
\index{Thompson!monoid}
\index{monoid!Thompson}
We write $I\subseteq\R$ for the interval $I=[0,1]$. An element $x\in I$ is
called \emph{dyadic} when it is of the form $\frac m{2^{n}}$ for some
$m,n\in\N$.
The Thompson group~$F$~\cite{cannon1996introductory,mckenzie1973elementary} has
as elements the homeomorphisms $f:I\to I$ which
\begin{itemize}
\item are piecewise linear with slopes $2^n$ for some $n\in\Z$,
\item with a finite number of breakpoints,
\item such that the coordinates of the breakpoints are dyadic.
\end{itemize}
The multiplication~$gf$ of two elements $f:I\to I$ and $g:I\to I$ is given by
their composite~$g\circ f$. For instance, we have the three functions whose
graphs are depicted below are elements of~$F$, respectively called $x_0$, $x_1$,
and $x_2$:
\begin{align*}
  \fig{F-x0}
  &&
  \fig{F-x1}
  &&
  \fig{F-x2}  
\end{align*}
The explicit definition of~$x_0$ is
\[
  x_0(t)=
  \begin{cases}
    2t&\text{if $0\leq t\leq\frac14$}\\
    t+\frac14&\text{if $\frac14\leq t\leq\frac12$}\\
    \frac{t+1}2&\text{if $\frac12\leq t\leq 1$}
  \end{cases}
\]
and more generally, given $i\in\N$, one can define a function~$x_i:I\to I$
in~$F$ by
\[
  x_i(t)=
  \begin{cases}
    t&\text{if $0\leq t\leq t_i$}\\
    x_0(2^i(t-t_i))/2^i+t_i&\text{if $t_i\leq t\leq 1$}\\
  \end{cases}
\]
where $t_i=1-\frac1{2^i}$ and $i\neq 0$. The group~$F$
is generated by the elements~$x_i$ and can be described as the enveloping group
of the \emph{Thompson monoid~$F^+$}, which is presented by
\[
  F^+
  =
  \Pres{\star}{x_i}{x_{j+1}x_i=x_ix_j}_{i<j\in\N}
  \pbox.
\]
More explicitly, this means that, starting from the above presentation, we should add a
formal inverse~$\finv x_i$ to each generator $x_i$, as explained in
\secr{env-group-pres}. Insights on this presentation, as well as alternative
descriptions of~$F$, are given in~\secr{3-Thompson-F}. Note that the
generator~$x_2$ is superfluous, since the relation $x_2x_0=x_0x_1$ can be
replaced by $x_2=x_0x_1\finv x_0$. More generally, given~$i\geq 2$, one has
$x_i=x_0^{n-1}x_1(\finv x_0)^{n-1}$, which shows that the group is generated
by~$x_0$ and~$x_1$.
Following this trail, the group~$F$ also admits the following finite
presentation with only two generators and their inverses~\cite{cannon1996introductory}:
\[
  \Pres{\star}{a,b,\finv a,\finv b}{[a\finv b,\finv aba],[a\finv b,\finv a\finv abaa],a\finv a,\finv aa,b\finv b,\finv bb}
\]
where $a$ and $b$ respectively correspond to $x_0$ and $x_1$ in the previous
presentation, $[u,v]$ is a notation for the commutator $uv\finv u\finv v$, and
we simply write~$u$ for a relation $u=1$. The Thompson group~$F$ is
$\FP_\infty$~\cite{brown1984infinite}.

Two other variants of the group~$F$ were introduced by
Thompson~\cite{cannon1996introductory}. The group~$T$ is a ``cyclic variant''
whose elements are piecewise linear functions $f:S^1\to S^1$, where $S^1$
denotes the circle of unit perimeter. The group~$V$ is a ``symmetric variant''
consisting of functions~$f:S^1\to S^1$ which are right-continuous, bijective, and
piecewise-linear satisfying similar requirements as before. There are natural
inclusion morphisms $F\into T\into V$. Typical elements of~$T$ and~$V$ are
\begin{align*}
  \fig{T-ex}
  &&
  \fig{V-ex}
\end{align*}
Both~$T$ and~$V$ also admit finite presentations. Generalizations of those
groups were also introduced by using $n$-ary trees instead of binary
ones~\cite{higman1974finitely}, by using~$I^n$ instead of~$I$ for the domain of
the endomorphisms~\cite{brin2004higher}, and by using braidings instead of
symmetries~\cite{brin2007algebra}.

\subsection{Wirtinger monoids}
\label{sec:wirtinger}
\index{Wirtinger monoid}
\index{monoid!Wirtinger}
\index{presentation!Wirtinger}
A \emph{knot} is an embedding $\kappa:S^1\to\R^3$ of the circle into
space. Suppose given a diagram representing such a knot on the plane, e.g., the
trefoil knot:
\[
  \fig{trefoil-abc}
\]
We label the arcs (hereby $a$, $b$, $c$) and orient them following a
conventional orientation of~$S^1$. This induces a presentation of a monoid, with
arcs as generators, and a relation
\[
  ab=ca
\]
for each crossing
\[
  \fig{wirtinger-crossing}
\]
(the orientation of~$b$ or~$c$ is not relevant). For instance, the monoid
associated to the trefoil knot is
\[
  \Pres{\star}{a,b,c}{ab=ca,bc=ab,ca=bc}
  \pbox.
\]
The enveloping group
of the monoid is the fundamental group of the knot, \ie the fundamental group of
$\R^3\setminus\kappa(S^1)$, and the presentation is called the \emph{Wirtinger
  presentation}. This group only depends on the knot, up to isotopy (not on the
diagram used to construct it or its orientation for instance). As a side note,
this group is also the enveloping group of the monoid
\[
  \Pres\star{a,b}{bab=aba}
\]
which is not isomorphic to the above monoid.

Let us see an application of those presentations. Two knots are equivalent when
their complements are isomorphic. Moreover, the fundamental group of the complement is
an invariant of the equivalence classes. This shows that the trefoil is not
isomorphic to the unknot: for instance, one has an abelian fundamental group
whereas the other does not.

\subsection{Temperley-Lieb monoids}
\label{sec:temperley-lieb}
\index{Temperley-Lieb!monoid}
\index{Temperley-Lieb!algebra}
\index{monoid!Temperley-Lieb}
Given $n\in\N$, the $n$-th \emph{Temperley-Lieb monoid} is presented by the
2-polygraph with 1-generators~$d$ and $u_i$, indexed by $0\leq i<n$,
together with relations, for $0\leq i,j<n$:
\begin{itemize}
\item $du_i=u_id$,
\item $u_iu_i=du_i$,
\item $u_iu_ju_i=u_i$ whenever $|j-i|=1$,
\item $u_ju_i=u_iu_j$ whenever $|j-i|>1$.
\end{itemize}
Given a ring~$R$ and a parameter $\delta\in R$, the \emph{Temperley-Lieb
  algebra} $A_n(\delta)$ is the free $R$-algebra, generated by the $u_i$ as
above, satisfying the above relations with $d=\delta$.
More details can be found in~\cite{abramsky2009temperley, jones1985polynomial,
  kauffman1990invariant, temperley1971relations}, as well as \secr{temperley-lieb-cat}.

\subsection{Brauer monoids}
\label{sec:mon-brauer}
\index{Brauer!monoid}
\index{monoid!Brauer}
Given $n\in\N$, the $n$-th \emph{Brauer monoid}~\cite{brauer1937algebras,
  kudryavtseva2006presentations} is presented by the 2-polygraph with
$1$-generators~$a_i$ and $u_i$, indexed by~$0\leq i<n$, subject to the 12
families of relations, indexed by $0\leq i,j<n$,
\begin{itemize}
\item for $0\leq i<n$,
  \begin{align*}
  a_ia_i&=1
  &
  u_iu_i&=1
  &
  a_iu_i&=u_i
  &
  u_ia_i&=u_i    
  \end{align*}
\item for $|j-i|=1$,
  \begin{align*}
  a_ia_ja_i&=a_ja_ia_j
  &
  u_iu_ju_i&=u_i
  &
  a_iu_ja_i&=a_ju_ia_j
  &
  u_ia_ju_i&=u_i    
  \end{align*}
\item for $|j-i|>1$,
  \begin{align*}
  a_ja_i&=a_ia_j
  &
  u_ju_i&=u_iu_j
  &
  u_ja_i&=a_iu_j    
  \end{align*}
\item for $0\leq i<n-3$,
  \begin{align*}
    a_ia_{i+1}u_iu_{i+2}&=a_{i+2}a_{i+1}u_iu_{i+2}
  \end{align*}
\end{itemize}
In particular, the relations satisfied by the $a_i$ are precisely those of the
symmetric group, see \secr{sym-group}, and the relations satisfied by the $u_i$
are those of Temperley-Lieb monoids with $d=\unit{}$, see \secr{temperley-lieb}.
More details about Brauer algebras, as well as a graphical illustration of the
relations can be found in~\secr{chord}.

\subsection{Graphs}
\index{graph}

Consider the monoid presented by
\[
  \Pres\star{s,t}{ss=st=s,ts=tt=t}.
\]
A right action of this monoid on a set~$X$ is precisely a directed graph: the
vertices are elements $x\in X$ such that $x\cdot s=x\cdot t=x$, and other
elements $y\in X$ are edges with source $y\cdot s$ and target $y\cdot t$.

\subsection{Tseitin monoid}
\label{sec:tseitin}
\index{Tseitin monoid}
\index{monoid!Tseitin}
The \emph{Tseitin monoid} is the one whose set of 1-gene\-rators is
$\set{a,c,b,d,e}$, subject to the relations
\begin{align*}
  ac&=ca
  &
  bc&=cb
  &
  eca&=ce
  &
  ccae&=cca
  \\
  ad&=da
  &
  bd&=db
  &
  edb&=de
\end{align*}
The reason why it is interesting is that it admits a small presentation, and yet
it was shown to have undecidable word problem by
Tseitin~\cite{tseitin1958associative} (see \cref{sec:word-problem}), meaning
that, writing~$P$ for the above $2$-polygraph, there is no algorithm which,
given any two $1$-cells $u,v\in\freecat{P_1}$, answers whether $u\approx^Pv$
holds or not. In fact, a $1$\nbd-poly\-graph with only three relations
exhibiting this property can be found~\cite{matiyasevich2005decision}.

\section{Presentations of Categories}
An example of a presentation of a category which is not a monoid was already given
in \exr{walking-iso}: the walking isomorphism category. There are possible
variations on this construction of the ``walking something'', as we show below,
although most interesting categorical constructions (e.g., adjunctions) are
2\nbd-cate\-gorical and will require using 3-polygraphs.

\subsection{The walking retract}
\index{retract}
A \emph{retract} in a category~$C$ consists of a pair of morphisms $a:x\to y$
and $b:y\to x$ such that $b\circ a=\id_x$. In this case, $a$ is called a
\emph{section} and $b$ a \emph{retraction}. The \emph{walking retract} is the
category presented by
\[
  \Pres{x,y}{a:x\to y,b:y\to x}{ab=\id_x}
  \pbox.
\]
It can be described as the category
\[
  \xymatrix{
    x\ar@/^/[r]^a&\ar@/^/[l]^by\ar@(ur,dr)^c
  }
\]
with two objects $x,y$ and three non-trivial morphisms
\begin{align*}
  a&:x\to y
  &
  b&:y\to x
  &
  c&:y\to y  
\end{align*}
with non-trivial compositions
\begin{align*}
  ab&=\id_x
  &
  ac&=a
  &
  ba&=c
  &
  cc&=c
\end{align*}

\subsection{Interval objects}
\index{interval}
\index{cylinder}
Given an object~$x$ in a category, a \emph{cylinder} on~$x$ is an object~$y$
together with morphisms $s,t:x\to y$ and $p:y\to x$ such that
$p\circ s=\unit{x}=p\circ t$. A typical example of cylinder object on a
topological space~$X$ is the space $X\times I$, where $I=\clint 01$ is the
standard interval, with $s(x)=(x,0)$, $t(x)=(x,1)$ and $p(x,a)=x$,
see~\cite{kamps1997abstract} for details. The theory of an object together with
a cylinder is the category presented by
\[
  \Pres{x,y}{s:x\to y,t:x\to y,p:y\to x}{sp=\unit{x},tp=\unit{x}}
  \pbox.
\]
Note that it is a variation on previous case, since it consists of two morphisms
with a common retraction. Depending on the applications, the cylinder can be
equipped with various extra morphisms such as
\begin{itemize}
\item \emph{reversion}: a morphism $r:y\to y$ satisfying
  \begin{align*}
    r\circ s&=t&
    r\circ t&=s&
    r\circ r&=\unit{x}&
    s\circ r&=s
  \end{align*}
  typically, the endomorphism of~$X\times I$ defined by $r(x,a)=(x,1-a)$,
\item \emph{concatenation}: morphisms $c^-,c^+:y\to y$ satisfying the 
  relations
  \begin{align*}
    s\circ c^-&=s&t\circ c^+&=t&
    p\circ c^-&=p&p\circ c^+&=p
  \end{align*}
  and moreover, in the presence of a reversion
  \begin{align*}
    r\circ c^-&=c^+\circ r&r\circ c^+&=c^-\circ r
  \end{align*}
  typically, on a topological cylinder $X\times I$, we define $c^-(x,a)=(x,a/2)$
  and $c^+(x,a)=(x,(a+1)/2)$.
\end{itemize}

\subsection{The walking factorization}
\index{factorization}
A \emph{factorization} of a morphism~$c:x\to z$ consists of a pair of morphisms
$a:x\to y$ and $b:y\to z$ such that $c=b\circ a$. The \emph{walking
  factorization} is the category presented by
\[
  \Pres{x,y,z}{a:x\to y,b:y\to z,c:x\to z}{ab=c}
  \pbox.
\]
It can be described as the category with $\set{0,1,2}$ as objects, one
morphism~$i\to j$ whenever $i\leq j$ and no morphism $i\to j$ whenever $i>j$.

\subsection{The walking $n$-span}
\index{span!8-@$n$-}
Given a natural number $n$, consider the poset~$S_n$ whose elements are
\[
  S_n
  =
  \set{x_0^-,x_0^+,x_1^-,x_1^+,\ldots,x_{n-1}^-,x_{n-1}^+,x_n}
\]
ordered by $x_i^\alpha<x_j^\beta$ whenever $i<j$, for any
$\alpha,\beta\in\set{-,+}$, and $x_n$ is the maximal element. For~$n=0,1,2,3$,
the Hasse diagram of the poset is
\begin{align*}
  \xymatrix@C=1ex@R=2ex{
    x_0
  }
  &&
  \xymatrix@C=1ex@R=2ex{
    &x_1&\\
    x_0^-\ar[ur]&&\ar[ul]x_0^+
  }
  &&
  \xymatrix@C=1ex@R=2ex{
    &x_2&\\
    x_1^-\ar[ur]&&\ar[ul]x_1^+\\
    x_0^-\ar[u]\ar[urr]&&\ar[ull]\ar[u]x_0^+
  }
  &&
  \xymatrix@C=1ex@R=2ex{
    &x_3&\\
    x_2^-\ar[ur]&&\ar[ul]x_2^+\\
    x_1^-\ar[u]\ar[urr]&&\ar[ull]\ar[u]x_1^+\\
    x_0^-\ar[u]\ar[urr]&&\ar[ull]\ar[u]x_0^+\pbox.
  }  
\end{align*}
This poset can be seen as a category, still denoted $S_n$, with the above set of
objects and an arrow $x\to y$ when $x\leq y$. Given a category~$C$, a functor
$S_n\to C$ (\resp $S_n^\op\to C$) is called an \emph{$n$-cospan} (\resp
\emph{$n$-span})~\cite{batanin1998monoidal}. In particular, a $1$-span is a span
in the usual sense. A typical $n$-cospan in the category~$\Top$ is given by
sending $x_i^\alpha$ to the $i$-disk $D^i$, the inclusion
$x_{i-1}^-\to x_i^\alpha$ (\resp $x_{i-1}^+\to x_i^\alpha$) corresponding to the
inclusion of the $(i{-}1)$-disk into the lower (\resp upper) hemisphere of the
sphere $S^{i-1}$ bounding the disk $D^i$.
The category~$S_n$ admits a presentation with
\begin{itemize}
\item 0-generators: $x_i^-,x_i^+,x_n$ for $0\leq i<n$,
\item 1-generators: $\xymatrix{x_{i-1}^\alpha\to x_i^\beta}$ and
  $\xymatrix{x_{n-1}^\alpha\to x_n}$ for $0<i<n$ and $\alpha,\beta\in\set{-,+}$,
\item 2-generators:
  \begin{align*}
    \xymatrix@C=2ex@R=2ex{
      &x_i^\gamma&\\
      x_{i-1}^\beta\ar[ur]&\To&\ar[ul]x_{i-1}^{\beta'}\\
      &\ar[ul]x_{i-2}^\alpha\ar[ur]&\\
    }
    &&
    \xymatrix@C=2ex@R=2ex{
      &x_n&\\
      x_{n-1}^\beta\ar[ur]&\To&\ar[ul]x_{n-1}^{\beta'}\\
      &\ar[ul]x_{n-2}^\alpha\ar[ur]&\\
    }    
  \end{align*}
  for $1<i<n$ and $\alpha,\beta,\beta',\gamma\in\set{-,+}$.
\end{itemize}

\subsection{The globe category}
\index{globe}
\index{category!globe}
\nomenclature[On]{$\glob^{(n)}$}{category of globes of dimension $\le n$}
Given an integer $n \ge 0$, we write $\glob^{(n)}$ for the category with
$\set{0,\ldots,n}$ as objects and two non-trivial morphisms $\cossce ji,
\cottge ji : j\to i$
for every $j<i$, with composition given by $f\circ \cossce ij=\cossce ik$ and
$f \circ \cottge ij = \cottge ik$ for every morphism $f:j\to k$. Writing
$\cosce i$ and $\cotge i$ for $\cossce i{i+1}$ and $\cottge i{i+1}$,
it can be pictured as
\[
  \xymatrix{
    0\doubr{\cosce{0}}{\cotge{0}}& 1 \doubr{\cosce{1}}{\cotge{1}}&
    \cdots\doubr{\cosce{n-1}}{\cotge{n-1}} & n
    \pbox.
    }
  \]
This category is called the \emph{category of globes of dimension $\le n$}
and a presheaf on $\glob^{(n)}$ an \emph{$n$-globular set}, see
\cref{sec:gset}.
The category~$\glob^{(n)}$ admits the presentation
\[
  \Pres{0,\ldots,n}{\cosce i, \cotge i: i \to i+1}{\cosce{i+1}\cosce{i} =
\cotge{i+1}\cosce{i}, \cosce{i+1}\cotge{i} = \cotge{i+1}\cotge{i}}
  \pbox.
\]

\subsection{The augmented simplicial category}
\index{simplicial!category}
\index{category!simplicial}
\nomenclature[.D]{$\Simplaug$}{augmented simplicial category}
The \emph{augmented simplicial category} $\Simplaug$ is the category where an
object is a natural number~$n\in\N$ and a morphism $f:m\to n$ is a weakly
monotone function $f:\intset{m}\to\intset{n}$, where $\intset{n}$ denotes the
set $\set{0,\ldots,n-1}$. A presentation of this category is detailed
in \secr{simpl-cat}.

\subsection{Presentations of monoidal categories}
Many categories of interest are actually monoidal categories. It is generally
much easier to present them as monoidal categories (as opposed to as a category)
and then deduce a presentation of those as categories if one insists on having
that, see \secr{3-pres-2-pres}. For instance, one can construct a presentation
of the augmented simplicial category~$\Simplaug$, see \secr{ass-simpl-pres}, and
recover the presentation evoked in previous section, see~\exr{simpl-3-2-pres}.


\chapter{Examples of Coherent Presentations of Monoids}
\chaptermark{Coherent presentations of monoids}
\label{chap:SomeCoherentPresentations}

In this chapter, we present examples of coherent presentations of monoids in the sense developed in \cref{chap:2-Coherent}. In particular, we focus on families of monoids which occur in algebra and whose coherent presentations are computed using the rewriting method that extends Squier's and Knuth-Bendix's completion procedures into a homotopical completion-reduction procedure as presented in \cref{chap:2coh}.
\Cref{Theorem:CoherentAction} shows how coherent presentations of monoids can be used to make explicit the actions of monoids on small categories.
We apply this construction to the case of Artin monoids in \cref{Section:CoherentPresentationArtinTitsMonoids}.
In particular, we prove that the Tits-Zamolodchikov $3$-generators extend the Artin presentation into a coherent presentation (\cref{Theorem:ArtinCoherentPresentation}) and, as a byproduct, we give a constructive proof of a theorem of Deligne on the actions of an Artin monoid on a category~\cite[Theorem~1.5]{Deligne97}.
We also give coherent presentations of plactic and Chinese monoids in \cref{SS:CoherentPresentationPlacticChinese}.

\section{Artin Monoids}
\label{Section:CoherentPresentationArtinTitsMonoids}
We provide here coherent presentations of Artin monoids, already introduced
in \cref{sec:artin-monoid}. We also explain that these provide an explicit
description of actions of those monoids on categories. 
The construction is done in two stages.
Given an Artin monoid $B^+(W)$ on a Coxeter group $W$, first we consider the coherent presentation $\Gar_3(W)$ constructed on the Garside presentation of $B^+(W)$. Then we apply the homotopical completion-reduction on this coherent presentation by examining a family of generating triple confluences. We thus obtain a coherent presentation $\Art_3(W)$ of $B^+(W)$ that extends the Artin presentation.

\subsection{Action of monoids on categories}
\index{action!of a monoid!on a category}
Let~$M$ be a monoid seen as a $2$\nbd-cate\-gory with exactly one $0$-cell
$\star$, with the elements of $M$ as $1$-cells, $0$\nbd-compo\-sition given
by product in~$M$ and with identity $2$-cells only, see \cref{lem:mon-cat}.
An \emph{action}~$T$ of the monoid~$M$ on a category is a pseudofunctor
$T:M\to\Cat$. More explicitly, such an action is specified by
\begin{itemize}
\item a category $\C=T(\star)$,
\item an endofunctor $T(u):\C\fl\C$ for every element~$u$ of~$M$,
\item a natural isomorphism $T_{u,v}:T(u)T(v)\dfl T(uv)$ for every pair $(u,v)$
  of elements of~$M$ and a natural isomorphism $T_{\star}:1_{\C}\dfl T(1)$,
\end{itemize}
satisfying the following conditions:
\begin{itemize}
\item for every triple $(u,v,w)$ of elements of $M$, the following diagram
  commutes:
  \[
    \xymatrix @!C @R=1em @C=1em{
      & T(uv)T(w)
      \ar@2 @/^2ex/ [dr] ^-{T_{uv,w}}
      \\
      T(u)T(v)T(w)
      \ar@2 @/^2ex/ [ur] ^-{T_{u,v}T(w)}
      \ar@2 @/_2ex/ [dr] _-{T(u)T_{v,w}}
      && T(uvw)
      \\
      & T(u)T(vw)
      \ar@2 @/_2ex/ [ur] _-{T_{u,vw}}
    }
  \]
\item for every element $u$ of $M$, the following two diagrams commute:
  \begin{align*}
    \xymatrix @!C @R=1.5em @C=1em {
      & T(1)T(u)
      \ar@2 @/^2ex/ [dr] ^-{T_{1,u}}
      \\
      T(u)
      \ar@2 @/^2ex/ [ur] ^-{T_{\star} T(u)}
      \ar@{=} [rr] _-{}="tgt"
      && T(u)
    }
    &&
    \xymatrix @!C @R=1.5em @C=1em {
      & T(u)T(1)
      \ar@2 @/^2ex/ [dr] ^-{T_{u,1}}
      \\
      T(u) 
      \ar@2 @/^2ex/ [ur] ^-{T(u) T_{\star}}
      \ar@{=} [rr] _-{}="tgt"
      && T(u).
    }
  \end{align*}
\end{itemize}
Such an action corresponds to a \emph{2\nbd-repre\-sentation} of~$M$
in~$\Cat$, as defined by Elgueta in~\cite{Elgueta08}. We denote by
$\nRep2(M,\Cat)$ the category of actions of~$M$ on categories, equipped with
suitable morphisms, as detailed in~\cite[Section~5.1]{GaussentGuiraudMalbos15}.
The following result is proved in \cite[Theorem~5.1.6]{GaussentGuiraudMalbos15}:

\begin{theorem}
  \label{Theorem:CoherentAction}
  Suppose given a monoid~$M$ and a coherent presentation of~$M$ by a
  $(3,1)$-polygraph~$P$. We have an equivalence of categories
  \[
    \nRep2(M,\Cat)
    \equivto
    \nCat2(\pcat P,\Cat)
  \]
  between actions of~$M$ on categories and 2-functors from the
  $(2,1)$-category presented by~$P$ to~$\Cat$.
\end{theorem}

\noindent
Deligne already observed that this equivalence holds for Garside's presentation
of spherical Artin monoids~\cite[Theorem 1.5]{Deligne97}. 
In the rest of this section, we present an
application of \cref{Theorem:CoherentAction} to give an equational
description of an action of an Artin monoid on a category using
Zamolodchikov relations.

\subsection{Coxeter groups}
\index{Coxeter group}
\index{group!Coxeter}
Recall from \cref{sec:coxeter-group} that a \emph{Coxeter group} is a group $W$ that admits a presentation with a finite set $S$ of
generators and with one relation
\begin{equation}
  \label{cox2}
(st)^{m_{st}} = 1
\end{equation}
with $m_{st}\in\mathbb{N}\sqcup\set{\infty}$, for every $s$ and $t$ in $S$, with
the following requirements and conventions:
\begin{itemize}
\item $m_{st}=\infty$ means that there is, in fact, no relation between $s$ and $t$,
\item $m_{st}=1$ if and only if $s=t$.
\end{itemize}
The last requirement implies that $s^2=1$ holds in $W$ for every $s$ in $S$. As
a consequence, the group $W$ can also be seen as the monoid with the same
presentation. 
As in \cref{sec:coxeter-group}, we denote by $\langle st \rangle^n$ the element of length~$n$ in
the free monoid~$S^*$, obtained by multiplication of alternating copies of $s$
and~$t$:
\begin{align*}
  \langle st \rangle^0 &= 1,
  &
  \langle st \rangle^{n+1} &= s\langle ts \rangle^n
  \pbox.
\end{align*}
When $s\neq t$ and $m_{st}<\infty$, we use this notation and the relations
$s^2=t^2=1$ to write~\eqref{cox2} as a \emph{braid relation}:
\begin{equation}
    \label{cox2'}
  \langle st \rangle^{m_{st}} = \langle ts \rangle^{m_{st}}
  \pbox.
\end{equation}

A \emph{reduced expression} of an element $u$ of $W$ is a representative of
minimal length of $u$ in the free monoid~$S^*$, and we write $l(u)$ for the
length of any of the reduced expressions of~$u$, called the \emph{length}
of~$u$. The Coxeter group $W$ is finite if and only if it admits an element of
maximal length~\cite[Theorem~5.6]{BrieskornSaito72}. In that case, this element
is unique, called the \emph{longest element of $W$} and is denoted by
$w_0(S)$. For $I\subseteq S$, the subgroup of $W$ spanned by the elements of $I$
is denoted by $W_I$. It is a Coxeter group with generating set $I$. If $W_I$ is
finite, we denote by $w_0(I)$ its longest element.

\index{Artin!monoid}
\index{monoid!Artin}
\index{presentation!Artin}
\index{Artin!presentation}
We recall that the \emph{Artin monoid} associated to $W$ is the monoid denoted
by $B^+(W)$, generated by $S$ and subject to the braid
relations~\eqref{cox2'}. This presentation, seen as a $2$-polygraph, is denoted
by $\Art_2(W)$ and called \emph{Artin's presentation}.
This is the same as the one of $W$, except for the relations $s^2=1$.

\subsection{Length notation}
For every $u$ and $v$ in $W$, we have $l(uv)\leq l(u)+l(v)$ and we will use the
following notations
\begin{align*}
  \typedeux{1} 
  &\quad\Leftrightarrow\quad
  l(uv)=l(u)+l(v)
  \pbox,
  \\
  \typedeux{0}
  &\quad\Leftrightarrow\quad
  l(uv)<l(u)+l(v)
  \pbox.  
\end{align*}
We generalize the notation for a greater number of elements of $W$. For example,
in the case of three elements $u$, $v$, and $w$ of $W$, we write \typetrois{?}
when both equalities $l(uv)=l(u)+l(v)$ and $l(vw)=l(v)+l(w)$ hold. This case
splits in the following two mutually exclusive subcases:
\begin{align*}
  \typetrois{1} 
  &\quad\Leftrightarrow\quad
  \begin{cases}
    \typetrois{} \\
    l(uvw)=l(u)+l(v)+l(w)\pbox,
  \end{cases}
  \\
  \typetrois{0}
  &\quad\Leftrightarrow\quad
  \begin{cases}
    \typetrois{} \\
    l(uvw)<l(u)+l(v)+l(w)\pbox.
  \end{cases}  
\end{align*}

\subsection{Garside's coherent presentation}
\index{Garside presentation}
\index{presentation!Garside}
\label{SS:GarsideCoherentPresentation}
Let $W$ be a Coxeter group. We call \emph{Garside's presentation of $B^+(W)$}
the $2$-polygraph $\Gar_2(W)$ whose $1$-cells are the elements of
$W\setminus\{1\}$ and with one $2$-generator
\[
  \alpha_{u,v} : u|v \To uv
\]
whenever $l(uv)=l(u)+l(v)$ holds. Here, we write $uv$ for the product in $W$ and
$u|v$ for the product in the free monoid over $W$.  We denote by $\Gar_3(W)$ the
extended presentation of $B^+(W)$ obtained from $\Gar_2(W)$ by adjunction of one
$3$-generator
\[
  \xymatrix @!C @R=1em {
    & uv|w
    \ar@2 @/^/ [dr] ^-{\alpha_{uv,w}}
    \ar@3 []!<0pt,-15pt>;[dd]!<0pt,+15pt> ^-{A_{u,v,w}}
    \\
    u|v|w
    \ar@2 @/^/ [ur] ^-{\alpha_{u,v}|w}
    \ar@2 @/_/ [dr] _-{u|\alpha_{v,w}}
    && uvw
    \\
    & u|vw
    \ar@2 @/_/ [ur] _-{\alpha_{u,vw}}
  }
\]
for all $u$, $v$, and $w$ of $W\setminus\{1\}$ with \typetrois{1}.

\subsection{Homotopical completion-reduction of Garside's presentation}
\label{SS:GarsideCompletionReduction}
The coherent presentation $\Gar_3(W)$ can be computed by coherent
comple\-tion-reduction of the $2$-polygraph $\Gar_2(W)$, as we now explain,
see~\cite{GaussentGuiraudMalbos15}.

Let $<$ denote the strict order on the elements of the free monoid $W^*$ that
first compares their length as elements of $W^*$, then the length of their
components, starting from the right. The order relation $\leq$ generated by $<$
by adding reflexivity is a termination order on $\Gar_2(W)$: for every $2$-generator
$\alpha_{u,v}$ of $\Gar_2(W)$, we have $u|v>uv$. Hence the $2$-polygraph
$\Gar_2(W)$ terminates, so that its coherent completion is defined (see
\cref{Section:HomotopicalCompletion}). By applying the coherence
completion-reduction procedure (see
\cref{Subsection:HomotopicalCompletionReduction}), one can obtain a coherent
extension of the Garside presentation $\Gar_2(W)$, as detailed
in~\cite[Proposition~3.2.1]{GaussentGuiraudMalbos15}. The resulting
$(3,1)$-polygraph $\KBS(\Gar_2(W))$ has one $0$-generator, one $1$-generator for every
element of~$W\setminus\set{1}$, the $2$-generators
\begin{align*}
  \alpha_{u,v}:u|v&\To uv
  &
  \beta_{u,v,w}:u|vw&\To uv|w
\end{align*}
for every $u$ and $v$ of $W\setminus\{1\}$ with \typedeux{1} and every $u$, $v$,
and $w$ of $W\setminus\{1\}$ with \typetrois{0}, and the nine following families
of $3$-generators
\[
  \begin{array}{cc}
    \vcenter{\xymatrix @!C @C=1ex @R=1.5em {
      & uv|w
      \ar@2 @/^/ [dr] ^-{\alpha_{uv,w}}
      \ar@3 []!<0pt,-20pt>;[dd]!<0pt,20pt> ^-{A_{u,v,w}}
      \\
      u|v|w
      \ar@2 @/^/ [ur] ^-{\alpha_{u,v}|w}
      \ar@2 @/_/ [dr] _-{u|\alpha_{v,w}}
      && uvw
      \\
      & u|vw
      \ar@2 @/_/ [ur] _-{\alpha_{u,vw}}
    }}
  &
  \vcenter{\xymatrix @!C {
      u|v|w
      \ar@2 @/^3ex/ [rr] ^-{\alpha_{u,v}|w} ^{}="src"
      \ar@2 @/_/ [dr] _-{u|\alpha_{v,w}}
      && uv|w
      \\
      & u|vw
      \ar@2 @/_/ [ur] _-{\beta_{u,v,w}}
      \ar@3 "src"!<0pt,-15pt>;[]!<0pt,20pt> ^-{B_{u,v,w}}
    }}
  \\
    \vcenter{\xymatrix @!C @C=1ex @R=1.5em {
      & uv|wx
      \ar@2 @/^/ [dr] ^-{\beta_{uv,w,x}}
      \ar@3 []!<0pt,-20pt>;[dd]!<0pt,20pt> ^->>{C_{u,v,w,x}}
      \\
      u|v|wx
      \ar@2 @/^/ [ur] ^-{\alpha_{u,v}|wx}
      \ar@2 @/_/ [dr] _-{u|\beta_{v,w,x}}
      && uvw|x
      \\
      & u|vw|x
      \ar@2 @/_/ [ur] _-{\alpha_{u,vw}|x}
    }}
  &
  \vcenter{\xymatrix @!C @C=1em {
      u|v|wx
      \ar@2 @/^3ex/ [rrr] ^-{\alpha_{u,v}|wx} ^{}="s"
      \ar@2 @/_/ [dr] _-{u|\beta_{v,w,x}}
      &&& uv|wx
      \\
      & u|vw|x
      \ar@2 [r] _-{\beta_{u,v,w}|x} _{}="t"
      & uv|w|x
      \ar@2 @/_/ [ur] _-{uv|\alpha_{w,x}}
      \ar@3 "s"!<0pt,-15pt>;"t"!<0pt,15pt> ^-{D_{u,v,w,x}}
    }}
  \\
  \vcenter{\xymatrix @!C @C=1ex @R=1.5em {
      & uv|w|x
      \ar@2 @/^/ [dr] ^-{uv|\alpha_{w,x}}
      \ar@3 []!<0pt,-20pt>;[dd]!<0pt,20pt> ^->>{E_{u,v,w,x}}
      \\
      u|vw|x
      \ar@2 @/^/ [ur] ^-{\beta_{u,v,w}|x}
      \ar@2 @/_/ [dr] _-{u|\alpha_{vw,x}}
      && uv|wx
      \\
      & u|vwx
      \ar@2 @/_/ [ur] _-{\beta_{u,v,wx}}
    }}
  &
  \vcenter{\xymatrix @!C @C=-1.3em @R=1.5em {
      && uv|w|xy
      \ar@2 @/^/ [drr] ^-{uv|\alpha_{w,xy}}
      \\
      u|vw|xy
      \ar@2 @/^/ [urr] ^-{\beta_{u,v,w}|xy}
      \ar@2 @/_/ [dr] _-{u|\beta_{vw,x,y}}
      &&&& uv|wxy
      \\
      & u|vwx|y
      \ar@2 [rr] _-{\beta_{u,v,wx}|y} _-{}="t"
      && uv|wx|y
      \ar@2 @/_/ [ur] _-{uv|\alpha_{wx,y}}
      \ar@3 "1,3"!<0pt,-20pt>;"t"!<0pt,20pt> ^->>{F_{u,v,w,x,y}}
    }}
  \\
  \vcenter{\xymatrix @!C @C=1ex @R=1.5em {
      & uv|w|xy
      \ar@2 @/^/ [dr] ^-{uv|\beta_{w,x,y}}
      \ar@3 []!<0pt,-20pt>;[dd]!<0pt,20pt> ^->>{G_{u,v,w,x,y}}
      \\
      u|vw|xy
      \ar@2 @/^/ [ur] ^-{\beta_{u,v,w}|xy}
      \ar@2 @/_/ [dr] _-{u|\beta_{vw,x,y}}
      && uv|wx|y
      \\
      & u|vwx|y
      \ar@2 @/_/ [ur] _-{\beta_{u,v,wx}|y}
    }}
  &
  \vcenter{\xymatrix @!C {
      & uv | xy
      \ar@2 @/^/ [dr] ^-{\beta_{uv,x,y}}
      \\
      u | vxy 
      \ar@2 @/^/ [ur] ^-{\beta_{u,v,xy}}
      \ar@2 @/_3ex/ [rr] _-{\beta_{u,vx,y}} _{}="t"
      && uvx | y
      \ar@3 "1,2"!<0pt,-20pt>;"t"!<0pt,15pt> ^-{H_{u,v,x,y}}
    }}
\end{array}
\]
\[
  \vcenter{
    \xymatrix @R=1.5em {
      & uv | w = uv | xy
      \ar@2 @/^3ex/ [dr] ^-{\beta_{uv,x,y}}
      \ar@3 []!<0pt,-25pt>;[dd]!<0pt,25pt> ^-{I_{u,v,w,v',w'}}
      \\
      *\txt{$u | vw$ \\ $=$ \\ $u | v'w'$}
      \ar@2 @/^3ex/ [ur] ^-{\beta_{u,v,w}}
      \ar@2 @/_3ex/ [dr] _-{\beta_{u,v',w'}}
      && *\txt{$uvx | y$ \\ $=$ \\ $uv'x' | y$}
      \\
      & uv' | w' = uv' | x'y
      \ar@2 @/_3ex/ [ur] _-{\beta_{uv',x',y}}
    }
  }
\]
These $3$-generators are families indexed by all the possible elements of
$W\setminus\set{1}$, that can be deduced by the sources and targets of the
$2$-cells. For example, there is one $3$-generator $A_{u,v,w}$ for all elements
$u$, $v$, and $w$ in $W\setminus\set{1}$ with \typetrois{1}\!, and one $3$-generator
$F_{u,v,w,x,y}$ for all elements $u$, $v$, $w$, $x$, and $y$ in
$W\setminus\set{1}$ with \typecinq{0}{1}{1}{}{0}{}\!\!\!.

By considering a family of generating triple confluences associated to some of
the triple critical branchings of $\KBS(\Gar_2(W))$, one can reduce this family
of $3$-generators and obtain the following
result~\cite[Theorem~6.4.3]{GaussentGuiraudMalbos15}:

\begin{theorem}
  \label{Theorem:GarsideCoherentPresentation}
  For every Coxeter group $W$, the Artin monoid $B^+(W)$ admits $\Gar_3(W)$ as a
  coherent presentation.
\end{theorem}

\noindent
The $(3,1)$-polygraph $\Gar_3(W)$ is called \emph{Garside's coherent
  presentation} of the Artin monoid $B^+(W)$.

\subsection{Artin's coherent presentation of Artin monoids}
\label{SS:ArtinCoherentPresentation}
Let $W$ be a Coxeter group with a totally ordered set $S$ of generators. The
homotopical reduction method of \cref{Section:HomotopicalCompletionReduction} can be used
on Garside's coherent presentation $\Gar_3(W)$ to contract it into a smaller
coherent presentation associated to Artin's
presentation~\cite{GaussentGuiraudMalbos15}.

\index{Artin!presentation!coherent}
\index{presentation!Artin!coherent}

We consider the presentation of the Artin monoid $B^+(W)$ by
the $2$-polygraph $\Art_2(W)$ with one $0$-generator, the elements of $S$ as $1$-generator
and one $2$-generator
\[
  \gamma_{s,t} : \langle ts \rangle^{m_{st}} \To \langle st \rangle^{m_{st}},
\]
for all $t>s$ in $S$ such that $m_{st}$ is finite.
The following result extends the $2$-polygraph $\Art_2(W)$ into a coherent
presentation of the Artin monoid $B^+(W)$, called the \emph{Artin coherent
  presentation} of $B^+(W)$~\cite[Theorem 4.1.1]{GaussentGuiraudMalbos15}.
This coherent presentation is obtained using the homotopical reduction \cref{Subsection:GenericHomotopicalReduction} on Garside's coherent presentation
$\Gar_3(W)$.

\begin{theorem}
  \label{Theorem:ArtinCoherentPresentation}
  For every Coxeter group $W$, the Artin monoid $B^+(W)$ admits the coherent
  presentation $\Art_3(W)$ made of Artin's presentation $\Art_2(W)$ and one
  $3$-generator $Z_{r,s,t}$ for all elements $t>s>r$ of $S$ such that the subgroup
  $W_{\set{r,s,t}}$ is finite.
\end{theorem}

Artin's coherent presentation has exactly one $k$-cell,
$0\leq k\leq 3$, for every subset~$I$ of~$S$ of rank $k$ such that the subgroup
$W_I$ is finite.
The shapes of the $3$-generators $Z_{r,s,t}$ of the polygraph $\Art_3(W)$ are obtained by projection of the $3$-generators of the polygraph $\Gar_3(W)$ given in \cref{SS:GarsideCompletionReduction} and depend on the type of the Coxeter type of the parabolic subgroup $W_{\{r,s,t\}}$. According to the classification of finite Coxeter groups~\cite[Chapter~VI, \textsection{}4, Theorem~1]{BourbakiLie4-6}, there are five cases, shown below:
\begin{align*}
  \begin{tikzpicture}
    \foreach \i in {1, 2, 3} {
      \coordinate (\i) at (\i, 0) ;
    }
    \draw [fill = black] (1) circle (0.05) node [above] {$r$} ;
    \draw [fill = black] (2) circle (0.05) node [above] {$s$} ;
    \draw [fill = black] (3) circle (0.05) node [above] {$t$} ;
    \draw (1) -- (2) ;
    \draw (2) -- (3) ;
    \node at (2, -0.5) {$A_3$} ;
  \end{tikzpicture}
  &&
  \begin{tikzpicture}
    \foreach \i in {1, 2, 3} {
      \coordinate (\i) at (\i, 0) ;
    }
    \draw [fill = black] (1) circle (0.05) node [above] {$r$} ;
    \draw [fill = black] (2) circle (0.05) node [above] {$s$} ;
    \draw [fill = black] (3) circle (0.05) node [above] {$t$} ;
    \draw (1) -- node [above] {$\sm 4$} (2) ;
    \draw (2) -- (3) ;
    \node at (2, -0.5) {$B_3$} ;
  \end{tikzpicture}
  &&
  \begin{tikzpicture}
    \foreach \i in {1, 2, 3} {
      \coordinate (\i) at (\i, 0) ;
    }
    \draw [fill = black] (1) circle (0.05) node [above] {$r$} ;
    \draw [fill = black] (2) circle (0.05) node [above] {$s$} ;
    \draw [fill = black] (3) circle (0.05) node [above] {$t$} ;
    \draw (1) -- node [above] {$\sm 5$} (2) ;
    \draw (2) -- (3) ;
    \node at (2, -0.5) {$H_3$} ;
  \end{tikzpicture}
\end{align*}
\begin{align*}
  \begin{tikzpicture}
    \foreach \i in {1, 2, 3} {
      \coordinate (\i) at (\i, 0) ;
    }
    \draw [fill = black] (1) circle (0.05) node [above] {$r$} ;
    \draw [fill = black] (2) circle (0.05) node [above] {$s$} ;
    \draw [fill = black] (3) circle (0.05) node [above] {$t$} ;
    \node at (2, -0.5) {$A_1\times A_1\times A_1$} ;
  \end{tikzpicture}
  &&
  \begin{tikzpicture}
    \foreach \i in {1, 2, 3} {
      \coordinate (\i) at (\i, 0) ;
    }
    \draw [fill = black] (1) circle (0.05) node [above] {$r$} ;
    \draw [fill = black] (2) circle (0.05) node [above] {$s$} ;
    \draw [fill = black] (3) circle (0.05) node [above] {$t$} ;
    \draw (1) -- node [above] {$\sm p$} (2) ;
    \node at (2, -0.5) {$I_2(p)\times A_1 \:\: \sm{3\leq p<\infty}$} ;
  \end{tikzpicture}
\end{align*}
According to these types, the $3$-generators $Z_{r,s,t}$ have the following shapes
\begin{itemize}
\item type $A_3$:
  \[
    \xymatrix @!C @C=2em @R=1.5em{
      & strsrt
      \ar@2 [r] ^{s\gamma_{rt}s\gamma_{rt}^-}
      & srtstr
      \ar@2 [r] ^{sr\gamma_{st}r}
      \ar@3 []!<0pt,-40pt>;[dddd]!<0pt,40pt> ^{\:Z_{r,s,t}}
      & srstsr
      \ar@2@/^2ex/ [dr] ^-{\gamma_{rs}tsr}
      \\
      stsrst
      \ar@2@/^2ex/ [ur] ^-{st\gamma_{rs}t}
      &&&& rsrtsr
      \\
      tstrst
      \ar@2 [u] ^{\gamma_{st}rst}
      \ar@2 [d] _{ts\gamma_{rt}st}
      &&&& rstrsr
      \ar@2 [u] _{rs\gamma_{rt}sr}
      \\
      tsrtst
      \ar@2@/_2ex/ [dr] _-{tsr\gamma_{st}}
      &&&& rstsrs
      \ar@2 [u] _{rst\gamma_{rs}}
      \\
      & tsrsts
      \ar@2 [r] _{t\gamma_{rs}ts}
      & trsrts
      \ar@2 [r] _{\gamma_{rt}s\gamma_{rt}^-s}
      & rtstrs
      \ar@2@/_2ex/ [ur] _-{\:r\gamma_{st}rs}
    }
  \]
  \renewcommand{\objectstyle}{\scriptstyle} 
  \renewcommand{\labelstyle}{\scriptstyle} 
\item type $B_3$:
  \[
    \xymatrix @!C @R=1.5em @C=1.5em{
      & srtstrstr
      \ar@2 [r] ^{sr\gamma_{st}rs\gamma_{rt}}
      & srstsrsrt
      \ar@2 [r] ^{srst\gamma_{rs}t}
      \ar@3 []!<0pt,-65pt>;[dddddddd]!<0pt,65pt> ^{\displaystyle \: Z_{r,s,t}}
      & srstrsrst
      \ar@2@/^2ex/ [dr] ^-{\quad srs\gamma_{rt}srst}
      \\
      srtsrtstr
      \ar@2@/^2ex/ [ur] ^-{srts\gamma_{rt}^-str}
      &&&&
      srsrtsrst
      \ar@2 [d] ^-{\gamma_{rs}tsrst}
      \\
      strsrstsr
      \ar@2 [u] ^-{s\gamma_{rt}sr\gamma_{st}^-r}
      &&&&
      rsrstsrst
      \\
      stsrsrtsr
      \ar@2 [u] ^{st\gamma_{rs}tsr}
      &&&&
      rsrtstrst
      \ar@2 [u] _{rsr\gamma_{st}rst}
      \\
      tstrsrtsr
      \ar@2 [u] ^{\gamma_{st}rsrtsr}
      \ar@2 [d] _{ts\gamma_{rt}s\gamma_{rt}^-sr}
      &&&&
      rsrtsrtst
      \ar@2 [u] _{rsrts\gamma_{rt}^-st}
      \\
      tsrtstrsr
      \ar@2 [d] _{tsr\gamma_{st}rsr}
      &&&&
      rstrsrsts
      \ar@2 [u] _{rs\gamma_{rt}sr\gamma_{st}^-}
      \\
      tsrstsrsr
      \ar@2 [d] _-{tsrst\gamma_{rs}}
      &&&&
      rstsrsrts
      \ar@2 [u] _{rst\gamma_{rs}ts}
      \\
      tsrstrsrs
      \ar@2@/_2ex/ [dr] _-{tsrs\gamma_{rt}srs}
      &&&&
      \ar@2 [u] _-{r\gamma_{st}rs\gamma_{rt}s}
      rtstrstrs
      \\
      & tsrsrtsrs
      \ar@2 [r] _{t\gamma_{rs}tsrs}
      & trsrstsrs
      \ar@2 [r] _{\gamma_{rt}sr\gamma_{st}^-rs}
      & rtsrtstrs
      \ar@2@/_2ex/ [ur] _-{\quad rts\gamma_{rt}^-strs}
    }
  \]
  \renewcommand{\objectstyle}{\displaystyle}
  \renewcommand{\labelstyle}{\displaystyle}
\item type $A_1\times A_1\times A_1$:
  \[
    \xymatrix @!C @C=2.8em @R=1.5em{
      & str
      \ar@2 [r] ^{s\gamma_{rt}} _{}="src"
      & srt
      \ar@2@/^2ex/ [dr] ^{\gamma_{rs}t}
      \\
      tsr
      \ar@2@/^2ex/ [ur] ^{\gamma_{st}r}
      \ar@2@/_2ex/ [dr] _{t\gamma_{rs}}
      &&& rst
      \\
      & trs
      \ar@2 [r] _{\gamma_{rt}s} ^{}="tgt"
      & rts
      \ar@2@/_2ex/ [ur] _{r\gamma_{st}}
      \ar@3 "src"!<0pt,-15pt>;"tgt"!<0pt,15pt> ^{\:Z_{r,s,t}}
    }
  \]
\item type $H_3$:
  \renewcommand{\objectstyle}{\scriptstyle} 
  \renewcommand{\labelstyle}{\scriptstyle} 
  \[
    \xymatrix @!C @R=1em @C=1.25em{
      & srsrtsrstrsrsrt
      \ar@2 [r] _{}="src"
      & srsrtsrstsrsrst
      \ar@2@/^2ex/ [dr]
      \\
      srstrsrsrtsrsrt
      \ar@2@/^2ex/ [ur]
      &&&
      srsrtsrtstrsrst
      \ar@2[d]
      \\
      srstsrsrstsrsrt
      \ar@2 [u]
      &&&
      srsrtstrsrtsrst
      \ar@2 [d]
      \\
      srtstrsrtstrsrt
      \ar@2 [u]
      &&&
      srsrstsrsrtsrst
      \ar@2 [d]
      \\
      srtsrtstrsrtstr
      \ar@2 [u]
      &&&
      rsrsrtsrsrtsrst
      \\
      srtsrstsrsrstsr
      \ar@2 [u]
      &&&
      rsrstrsrsrtsrst
      \ar@2 [u]
      \\
      srtsrstrsrsrtsr
      \ar@2 [u]
      &&&
      rsrstsrsrstsrst
      \ar@2 [u]
      \\
      strsrsrtsrsrtsr
      \ar@2 [u]
      &&&
      rsrtstrsrtstrst
      \ar@2 [u]
      \\
      stsrsrstsrsrtsr
      \ar@2 [u]
      &&&
      rsrtsrtstrsrtst
      \ar@2 [u]
      \\
      tstrsrstsrsrtsr
      \ar@2 [u]
      \ar@2 [d]
      &&&
      rsrtsrstsrsrsts
      \ar@2 [u]
      \\
      tsrtsrstsrstrsr
      \ar@2 [d]
      &&&
      rsrtsrstrsrsrts
      \ar@2 [u]
      \\
      tsrtsrtstrstrsr
      \ar@2 [d]
      &&&
      rstrsrsrtsrsrts
      \ar@2 [u]
      \\
      tsrtstrsrtstrsr
      \ar@2 [d]
      &&&
      rstsrsrstsrsrts
      \ar@2 [u]
      \\
      tsrstsrsrstsrsr
      \ar@2 [d]
      &&&
      rtstrsrtstrsrts
      \ar@2 [u]
      \\
      tsrstrsrsrtsrsr
      \ar@2 [d]
      &&&
      rtsrtstrsrtstrs
      \ar@2 [u]
      \\
      tsrsrtsrstrsrsr
      \ar@2 [d]
      &&&
      rtsrstsrsrstsrs
      \ar@2 [u]
      \\
      tsrsrtsrstsrsrs
      \ar@2[d]
      &&&
      rtsrstrsrsrtsrs
      \ar@2 [u]
      \\
      tsrsrtsrtstrsrs
      \ar@2@/_2ex/ [dr]
      &&&
      trsrsrtsrsrtsrs
      \ar@2[u]
      \\
      & tsrsrtstrsrtsrs
      \ar@2 [r] ^{}="tgt"
      & tsrsrstsrsrtsrs
      \ar@2@/_2ex/ [ur]
      \ar@3 "src"!<0pt,-160pt>;"tgt"!<0pt,160pt> ^{\displaystyle \: Z_{r,s,t}}
    }
  \]
\item type $I_2(p)\times A_1$, $p\geq 3$:
  \renewcommand{\objectstyle}{\displaystyle}
  \renewcommand{\labelstyle}{\displaystyle}
  \[
    \xymatrix @!C @C=2.35em @R=1.25em{
      & st\langle rs \rangle^{p-1}
      \ar@2 [r] ^-*+<8pt>{s\gamma_{rt}\langle rs \rangle^{p-2}}
      & (\cdots)
      \ar@2 [r]
      \ar@3 []!<0pt,-20pt>;[dd]!<0pt,20pt> ^{\:Z_{r,s,t}}
      & \langle sr \rangle^p t
      \ar@2@/^3ex/ [dr]^*+<4pt>{\gamma_{rs}t}
      \\
      t\langle sr \rangle^p
      \ar@2@/^3ex/ [ur]^{\gamma_{st}\langle rs \rangle^{p-1}}
      \ar@2@/_3ex/ [dr]_{t\gamma_{rs}}
      &&&& \langle rs \rangle^p t
      \\
      & t\langle rs \rangle^p
      \ar@2 [r]_-*+<6pt>{\gamma_{rt}\langle sr \rangle^{p-1}}
      & rt\langle sr \rangle^{p-1}
      \ar@2 [r]_-*+<6pt>{r\gamma_{st}\langle sr \rangle^{p-2}}
      & (\cdots)
      \ar@2@/_3ex/ [ur] 
    }
  \]
\end{itemize}

\subsection{Action of braid monoids on categories}
Following \cref{Theorem:CoherentAction} that, up to equivalence, the
actions of a monoid~$M$ on categories are the same as the $2$-functors from
$\cl{P}$ to $\Cat$, where $P$ is any coherent presentation of~$M$.  As an
application of \cref{Theorem:ArtinCoherentPresentation}, we establish the
relationship between coherent presentations of Artin monoids and Deligne's
notion of an action on a category. In particular, Deligne's
Theorem~\cite[Theorem~1.5]{Deligne97} is equivalent to
\cref{Theorem:GarsideCoherentPresentation}.

The extended presentation $\Gar_3(W)$ is a coherent presentation of the Artin
monoid $B^+(W)$. We thus get Deligne's Theorem~\cite[Theorem~1.5]{Deligne97}
for any Artin monoid as a consequence of
\cref{Theorem:GarsideCoherentPresentation,Theorem:CoherentAction}. Moreover,
\cref{Theorem:ArtinCoherentPresentation} gives a similar result in terms of
Artin's coherent presentation $\Art_3(W)$,
formalizing~\cite[Paragraph~1.3]{Deligne97} on the actions of $B^+_4$ on a
category.

\section{Plactic and Chinese Monoids}
\index{plactic!monoid}
\index{monoid!plactic}
\label{S:CoherentPresentationPlacticMonoids}
\label{SS:CoherentPresentationPlacticChinese}

We provide here coherent presentations of plactic and Chinese monoids whose presentations are recalled in \cref{sec:plactic-monoid}.
First, we recall a few points concerning the combinatorics of plactic monoids.

\subsection{Rows, columns, and tableaux}
A \emph{row} is a non-decreasing word $x_1\ldots x_k$ in the free monoid $\{1,\ldots,n\}^\ast$, i.e., with $x_{1}\leq x_{2}\leq\ldots\leq x_k$.
A \emph{column} is a decreasing word $x_p\ldots x_1$ in~$\{1,\ldots,n\}^\ast$, i.e., with~$x_{p}>x_{p-1}>\ldots>x_1$.  We denote by $\col(n)$ the set of non-empty columns.
A row  $x_1\ldots x_k$ \emph{dominates} a row $y_1 \ldots y_l$, and we denote $x_1\ldots x_k\vartriangleright y_1 \ldots y_l$, if $k\leq l$ and $x_{i}>y_{i}$, for $1\leq i \leq k$.
Any word $w$ in $\{1,\ldots,n\}^\ast$ has a unique decomposition as a product of rows of maximal length~$u_{1}\ldots u_{k}$, and it is called a \emph{tableau}  if $u_{1}\vartriangleright u_{2} \vartriangleright \ldots \vartriangleright u_{k}$.

\subsection{Schensted's algorithm}
The \emph{Schensted algorithm} computes for each $w$ in~$\{1,\ldots,n\}^\ast$ a tableau, denoted by $P(w)$, called the \emph{Schensted tableau} of $w$, and constructed from the following steps \cite{Schensted61}.
Given $u$ a tableau written as a product of rows of maximal length $u=u_1\ldots u_k$ and $1\leq y \leq n$, it computes the tableau $P(uy)$ as follows.
If $u_{k}y$ is a row, the result is $u_{1}\ldots u_{k}y$.
If $u_{k}y$ is not a row, then suppose~\mbox{$u_{k}=x_{1}\ldots x_{l}$} with~$1\leq x_{i} \leq n$ and let $j$ minimal such that $x_{j}>y$, 
then the result is $P(u_{1}\ldots u_{k-1}x_{j})v_{k}$, where~\mbox{$v_{k}=x_{1}\ldots x_{j-1}yx_{j+1}\ldots x_{l}$}.
The tableau $P(w)$ is computed from the empty tableau and iteratively applying the Schensted algorithm. In this way, $P(w)$ is the row reading of the planar representation of the tableau computed by the Schensted algorithm.

\subsection{Knuth's 2-polygraph}
\index{plactic!congruence}
\label{Placticcongruence}
\nomenclature[P]{$P_n$}{plactic monoid}
Let us give an orientation to the Knuth relations set out in \cref{sec:plactic-monoid} with
respect to the lexicographic order, thus forming a $2$-polygraph, denoted by~$\Knuth_2(n)$, whose $1$-generators are $1,\ldots,n$ and the $2$-generators are
\begin{align*}
  \eta_{x,y,z} &: zxy \dfl xzy,&&\text{for $1\leq x\leq y < z \leq n$,}
  \\
  \varepsilon_{x,y,z} &: yzx \dfl yxz,&&\text{for $1\leq x<y\leq z\leq n$.}
\end{align*}
The congruence on the free monoid
$\{1,\ldots,n\}^\ast$ generated by this polygraph is called the
\emph{plactic congruence of rank $n$} and the $2$-polygraph~$\Knuth_2(n)$ is a
presentation of the plactic monoid~$P_n$,~\cite[Theorem 6]{knuth1970permutations}.

\subsection{Pre-column presentation}
One adds to the presentation~$\Knuth_{2}(n)$ one superfluous generator $c_{u}$
for any~$u$ in $\col(n)$. We denote by $\Colo_1(n)$ the set of \emph{column
  generators}~$c_{u}$ for any $u$ in~$\col(n)$ and by
\[
  \gamma_u : c_{x_{p}}\ldots c_{x_{1}} \dfl  c_u,
\]
the defining relation for the column generators $u=x_{p} \ldots x_{1}$ in
$\col(n)$ of length greater than $2$. In the free monoid $\Colo_1(n)^\ast$, the
Knuth relations can be written in the following form
\begin{align*}
  \eta_{x,y,z}^c : c_{z}c_{x}c_{y} &\dfl c_{x}c_{z}c_{y},
  &&\text{for $1 \leq x\leq y<z  \leq n$,}
  \\
  \varepsilon_{x,y,z}^c : c_{y}c_{z}c_{x} &\dfl c_{y}c_{x}c_{z},
  &&\text{for $1 \leq x< y\leq z \leq n$.}
\end{align*}
The $2$-polygraph $\Knuthcc_2(n)$ whose $1$-generators are columns and $2$-generators are
the defining relations $\gamma_u$ for columns generators and the Knuth relations
$\eta_{x,y,z}^c$ and $\varepsilon_{x,y,z}^c$ is a presentation of the monoid
$P_n$.  We define the $2$-polygraph~$\PreColo_2(n)$ with column generators and
the $2$-generators are
\begin{align*}
  \alpha'_{x,zy} : c_{x}c_{zy} &\dfl c_{zx}c_{y},
  &&\text{for $1 \leq x\leq y<z \leq n$,}
  \\
  \alpha'_{y,zx} : c_{y}c_{zx} &\dfl c_{yx}c_{z},
  &&\text{for $1 \leq x< y\leq z \leq n$,}
  \\
  \alpha'_{x,u} : c_{x}c_{u} & \dfl c_{xu},
  &&\text{for $xu\in \col(n)$ and $1\leq x \leq n$,}
\end{align*}
where $\alpha'_{x,zy}$ and $\alpha'_{y,zx}$ correspond
respectively to the Knuth relations $\eta_{x,y,z}^c$
and~$\varepsilon_{x,y,z}^c$.  The $2$-polygraph $\PreColo_2(n)$ is a
presentation of the monoid $P_n$, then called the \emph{pre-column
  presentation}.

\subsection{Coherent column presentation}
\index{presentation!column}
\index{column!presentation}
Given columns $u$ and $v$, if the planar representation of the Schensted tableau
$P(uv)$ is not the tableau obtained as the concatenation of the two columns $u$
and $v$, we write $\typedeux{0}$. In this case, the tableau $P(uv)$
contains at most two columns~\cite[Lemma 3.1]{cain2015finite}.
We write~$\typedeuxbase{$u$}{$v$}{01}$ (\resp ~$\typedeuxbase{$u$}{$v$}{02}$) if the tableau
$P(uv)$ has one column (\resp two columns). When $\typedeux{0}$, we define a
$2$-generator
\[
  \alpha_{u,v} :  c_uc_v \dfl c_{w}c_{w'}
\]
where $w=uv$ and $c_{w'}=1$, if $\typedeux{01}$, and $w$ and $w'$ are
respectively the left and right columns of the tableau $P(uv)$, if
$\typedeux{02}$.

\index{column presentation!coherent}
\index{presentation!column!coherent}
The $2$-polygraph~$\Colo_2(n)$ whose set of $1$-generators is $\Colo_1(n)$ and the
$2$-generators are the~$\alpha_{u,v}$ is a finite convergent presentation of the
monoid $P_n$, called the \emph{column presentation}~\cite{cain2015finite}.
Using the coherent completion procedure defined in \cref{Subsubsection:HomotopicalCompletion}, this poly\-graph is extended into the \emph{column coherent presentation}~$\Colo_3(n)$ of the monoid $P_n$~\cite[Theorem 3.2.2]{HageMalbos17}.
Its $3$-generators, given by the confluence diagrams of the critical
branchings of the $2$-polygraph~$\Colo_2(n)$, have the following hexagonal form
\[
  \xymatrix @!C @C=3em @R=1em {
    &
    c_{e} c_{e'}c_t
    \ar@2[r] ^{c_{e}\alpha_{e',t}}
    \ar@3 []!<30pt,-8pt>;[dd]!<30pt,8pt> ^{\mathcal{X}_{u,v,t}}
    &
    c_{e} c_{b}c_{b'}
    \ar@2@/^/[dr] ^{\alpha_{e,b}c_{b'}}
    \\
    c_uc_vc_t
    \ar@2@/^/[ur] ^{\alpha_{u,v}c_t}
    \ar@2@/_/[dr] _{c_u\alpha_{v,t}}
    &&&
    c_{a}c_{d} c_{b'}
    \\
    &
    c_uc_{w}c_{w'}
    \ar@2[r] _{\alpha_{u,w}c_{w'}}
    &
    c_{a}c_{a'}c_{w'}
    \ar@2@/_/[ur] _{c_{a}\alpha_{a',w'}}
  }
\]
for any columns $u$, $v$, and $t$ such that \typetroisbasep{$u$}{$v$}{$t$}{00}\!.

\subsection{Pre-column coherent presentation}
Using the homotopical reduction procedure of~\cref{Section:HomotopicalCompletionReduction},
the coherent presentation~$\Colo_3(n)$ can be reduced into a smaller coherent
presentation of~$P_n$ as follows. Firstly, we apply a homotopical reduction on the
$(3,1)$-polygraph~$\Colo_3(n)$ with a collapsible part defined by some of the
generating triple confluences of the $2$-polygraph $\Colo_2(n)$. In this way, we
reduce the coherent presentation~$\Colo_3(n)$ of the monoid~$P_n$ into the
coherent presentation~$\overline{\Colo}_3(n)$ of~$P_n$, whose underlying
$2$-polygraph is $\Colo_2(n)$ and the $3$-generators~$\mathcal{X}_{u,v,t}$ are those
of $\Colo_3(n)$, but with~$u$ is of length~$1$.  Then we reduce the coherent
presentation~$\overline{\Colo}_3(n)$ into a coherent presentation
$\PreColo_3(n)$ obtained from $\PreColo_2(n)$ by adjunction of the $3$-generator
$R_{\Gamma_3}(C'_{x,v,t})$ where
\[
  \xymatrix @C=3em @R=0.6em {
    & {c_{xv}c_{t}}
    \ar@3 []!<0pt,-10pt>;[dd]!<0pt,+10pt> ^{C'_{x,v,t}}
    \\
    {c_{x}c_{v}c_{t}}
    \ar@2@/^/ [ur] ^-{\alpha_{x,v}c_{t}}
    \ar@2@/_/ [dr] _-{c_{x}\alpha_{v,t}}
    \\
    & {c_{x}c_{w}c_{w'}}
    \ar@2 [r] _-{\alpha_{x,w}c_{w'}}
    & {c_{xv}c_{z_{l}\ldots z_{q+1}}c_{w'}}
    \ar@2 [luu] _-{c_{xv}\alpha_{z_{l}\ldots z_{q+1}, w'}}
  }
\]
with \typetroisbasep{$x$}{$v$}{$t$}{12}, and the $3$-generator
$R_{\Gamma_3}(D_{x,v,t})$ where
\[
  \xymatrix @!C @C=2.3em @R=0.6em {
    &
    c_{e} c_{e'}c_t
    \ar@2[r] ^{c_{e}\alpha_{e',t}}
    \ar@3 []!<40pt,-10pt>;[dd]!<40pt,10pt> ^{D_{x,v,t}} 
    &
    c_{e} c_{b}c_{b'}
    \ar@2@/^/[dr] ^{\alpha_{e,b}c_{b'}}
    \\
    c_xc_vc_t
    \ar@2@/^/[ur] ^{\alpha_{x,v}c_t}
    \ar@2@/_/[dr] _{c_x\alpha_{v,t}}
    &&&
    c_{a}c_{d} c_{b'}
    \\
    &
    c_xc_{w}c_{w'}
    \ar@2[r] _{\alpha_{x,w}c_{w'}}
    &
    c_{a}c_{a'}c_{w'}
    \ar@2@/_/[ur] _{c_{a}\alpha_{a',w'}}
  }
\]
with \typetroisbasep{$x$}{$v$}{$t$}{22} and where the homotopical reduction
$R_{\Gamma_3}$ eliminates a collapsible part $\Gamma_3$ of
$\overline{\Colo}_3(n)$. In this way, we prove that the $(3,1)$-polygraph $\PreColo_3(n)$ is a coherent presentation of the monoid~$P_n$~{\cite[Theorem 4.3.4]{HageMalbos17}}.
For instance, the coherent presentation $\Colo_3(2)$ has only one $3$-generator
\[
  \xymatrix @C=3em @R=0.6em {
    & {c_{21}c_{21}}
    \ar@3 []!<-10pt,-10pt>;[dd]!<-10pt,+10pt> ^{C'_{2,1,21}}
    \\
    {c_2c_1c_{21}}
    \ar@2@/^/ [ur] ^-{\alpha_{2,1}c_{21}}
    \ar@2@/_/ [dr] _-{c_{2}\alpha_{1,21}}
    \\
    & {c_2c_{21}c_1}
    \ar@2 [r] _-{\alpha_{2,21}c_1}
    & {c_{21}c_2c_1.}
    \ar@2 [luu] _-{c_{21}\alpha_{2,1}}
  }
\]
In this case, the $(3,1)$-polygraphs $\PreColo_3(2)$ and $\Colo_3(2)$ coincide.

\subsection{Knuth's coherent presentation}
The coherent presentation $\PreColo_3(n)$ can be reduced into a coherent
presentation of the monoid $P_n$ whose underlying $2$-polygraph is
$\Knuth_2(n)$. We define an extended presentation $\Knuth_{3}(n)$ of the monoid
$P_n$ obtained from $\Knuth_{2}(n)$ by adjunction of the following set of
$3$-generators
\[
  \setof{\mathcal{R}(C'_{x,v,t})}{\typetroisbasep{$x$}{$v$}{$t$}{12}}
  \cup
  \setof{\mathcal{R}(D_{x,v,t})}{\typetroisbasep{$x$}{$v$}{$t$}{22}}
  \pbox,
\]
where $\mathcal{R} :\freegpd{\overline{\Colo}_3(n)}\fl\freegpd{\Knuthcc_3(n)}$
is a Tietze transformation, see {\cite[Section 4.4]{HageMalbos17}}. This gives a coherent presentation of the plactic monoid on the Knuth generators.

\begin{theorem}[{\cite[Theorem 4.4.7]{HageMalbos17}}]
  \label{KnuthcoherentTheorem}
  For $n > 0$, the $(3,1)$-polygraph $\Knuth_3(n)$ is a coherent presentation of
  the monoid~$P_n$.
\end{theorem}

\subsection{Example}
For instance, the Knuth coherent presentation of the mo\-noid~$P_2$ has generators
$c_1$ and $c_2$ subject to the Knuth relations
\[
\eta_{1,1,2}^c : c_2c_1c_1 \dfl c_1c_2c_1
\qtand
\varepsilon_{1,2,2}^c : c_2c_2c_1 \dfl c_2c_1c_2
\]
and the following~$3$-generator
\[
  \xymatrix@!C@C=3em{
    c_2c_2c_1c_1
    \ar@2@/^3ex/ [r] ^{2\eta_{1,1,2}^c}="src"
    \ar@2@/_3ex/ [r] _{\varepsilon_{1,2,2}^{c}1}="tgt"
    & c_2c_1c_2c_1
    \ar@3 "src"!<-5pt,-15pt>;"tgt"!<-5pt,15pt>  ^-{\, C''}
    \pbox.
  }
\]
Note that the Knuth coherent presentation of the monoid $P_2$ corresponds to the
coherent presentation that one can compute directly using the fact that the
$2$-polygraph $\Knuth_2(2)$ is convergent.

\subsection{Chinese monoids}
\label{SS:ChineseMonoids}
\index{Chinese monoid}
\index{monoid!Chinese}
\nomenclature[C]{$Ch_n$}{Chinese monoid}

Using the completion-reduction method, as applied previously for the plactic monoid, we calculate a coherent presentation of the Chinese monoid $Ch_n$ of rank $n>0$ whose presentation is recalled in \cref{ChineseRelations}.
We do not give here the details of this construction, developed in \cite{HageMalbos22}, but only present the method to obtain the form of the $3$-generators of this coherent presentation.

Chinese relations \cref{ChineseRelations} generate the Chinese congruence, denoted by~$\sim_{\C_n}$,  and interpreted in~\cite{CassaigneEspieKrobNovelliHibert01} in terms of Chinese staircases.
A \emph{Chinese staircase} is a collection of boxes in right-justified rows,  filled with non-negative integers, whose rows and  columns are indexed with elements of~$\{1,\ldots,n\}$ from top to bottom and from right to left respectively, and where the $i$-th row contains $i$ boxes, for~$1 \leq i \leq n$.
We denote by $R(t)$ the reading of a Chinese staircase $t$, row by row from right to left and from top to bottom.  A Schensted-like insertion algorithm, denoted by~$\insr{r}$, is introduced in~\cite{CassaigneEspieKrobNovelliHibert01}, and consists in inserting an element of~$\{1,\ldots,n\}$ into  a Chinese staircase from the right.
From a word $w=x_1x_2\ldots x_k$, we associate a Chinese staircase~$\llbracket w\rrbracket$ obtained by insertion of~$w$ in the empty staircase~$\emptyset$ by application of~$\insr{r}$ step by step from left to right:
\[
\llbracket w\rrbracket:=
(\ldots((\emptyset \insr{r} x_1) \insr{r} x_2) \insr{r} \ldots ) \insr{r} x_k.
\]
Chinese staircases satisfy the \emph{cross-section property} for the congruence~$\sim_{\C_n}$, that is,  for all words $w$ and $w'$, $w\sim_{\C_n} w'$ if and only if the insertion algorithm yields the same Chinese staircase: $\llbracket w\rrbracket = \llbracket w'\rrbracket$,~\cite{CassaigneEspieKrobNovelliHibert01}.
The elements of the monoid $Ch_n$ can thus be identified with the Chinese staircases, which also form a monoid, whose product is defined by setting~$t \star_{r} t' := (t \insr{r} R(t'))$, for all Chinese staircases~$t$ and $t'$.

We construct a finite convergent presentation $\Chin_2(n)$ of the monoid~$Ch_n$,  whose $1$-generators are \emph{columns} on $\{1,\ldots,n\}$ of length at most~$2$ and \emph{square generators}, as defined in {\cite[Section 4.1]{HageMalbos22}}, and whose $2$-generators are 
\[
\gamma_{u,v}:c_uc_v\dfl c_{e}c_{e'},
\]
for all columns~$c_u$ and $c_v$ such that~$c_uc_v$  does not form a  Chinese staircase and~$c_u\star_{r} c_v$ is equal to the Chinese staircase composed by the columns~$c_e$ and~$c_{e'}$. Note that the polygraph $\Chin_2(n)$ is obtained from the relations \cref{ChineseRelations} by applying Tietze transformations, defined in \cref{sec:2-tietze}, which consist in adding column generators and associated relations.

By definition of the $2$-generators, the source of each critical
branching of~$\Chin_2(n)$ has the form~$c_uc_vc_t$, for columns~$c_u$,~$c_v$,~$c_t$
such that~$c_uc_v$ and~$c_vc_t$ are not Chinese staircases. Their confluence diagrams are then obtained by applying the $2$-generators~$\gamma$, see {\cite[Theorem 5.6]{HageMalbos22}} for a detailed proof.
Following Squier's homotopical theorem, \cref{Theorem:SquierCompletion2polygraphs}, the $2$-polygraph~$\Chin_2(n)$ extends into a coherent presentation $\Chin_3(n)$ of the  monoid $Ch_n$ by adjunction of $3$-generators with the following decagonal form
\[
  \xymatrix @!C @C=3ex @R=3ex {
    &
    c_{e} c_{e'}c_t
    \ar@2[r] ^{c_{e}\gamma_{e',t}} 
    \ar@3 []!<95pt,-10pt>;[dd]!<95pt,10pt> ^{\mathcal{X}_{u,v,t}} 
    &
    c_{e} c_{b}c_{b'}
    \ar@2[r] ^{\gamma_{e,b}c_{b'}} 
    &
    c_{s} c_{s'}c_{b'}
    \ar@2[r] ^{c_{s}\gamma_{s',b'}}
    &
    {c_s c_k c_{k'}}
    \ar@2@/^/[dr] ^{\gamma_{s,k}c_{k'}}
    \\
    c_uc_vc_t
    \ar@2@/^/[ur] ^{\gamma_{u,v}c_t}
    \ar@2@/_/[dr] _{c_u\gamma_{v,t}}
    &&&&&
    c_{l}c_{m} c_{k'}
    \\
    &
    c_uc_{w}c_{w'}
    \ar@2[r] _{\gamma_{u,w}c_{w'}}
    &
    c_{a}c_{a'}c_{w'}
    \ar@2[r] _{c_{a}\gamma_{a',w'}}
    &
    c_{a}c_{d}c_{d'}
    \ar@2[r] _{\gamma_{a,d}c_{d'}}
    & c_{l}c_{l'}c_{d'}
    \ar@2@/_/[ur] _{c_{l}\gamma_{l',d'}}
  }
\]
for all columns~$c_u,c_v,c_t$ such that~$c_u c_v$ and~$c_v c_t$ are not  normal forms with respect to~$\Chin_2(n)$, and where the $\gamma$ denote either a $2$-generator of $\Chin_2(n)$ or an identity. This proves the following result.

\begin{theorem}
  \label{Theorem:CoherenceChinois}
  For $n>0$, the $(3,1)$-polygraph $Col_3(n)$ is a finite coherent convergent
  presentation of the Chinese monoid~$Ch_n$.
\end{theorem}


\chapter{A Catalogue of 3-Polygraphs}
\label{chap:3ex}
In this chapter, we give some examples of presentations of
$2$-categories by $3$\nbd-poly\-graphs. In many examples, the presented
$2$-categories are in fact monoidal categories and, actually, PROs. For those,
we consider presentations by 3\nbd-poly\-graphs~$P$ with only one
0-generator~$\star$ and one
1\nbd-gene\-rator~$a$, so that we only need to provide the 2-generators (which
we simply call \emph{generators}) and 3-generators (which we call \emph{rules}
or \emph{relations}). Moreover, we simply write $n$ instead of
$a^n$ for a 1-cell in $\freecat{P_1}$.

\section{Braids and Symmetries}
\subsection{Positive braids}
\label{sec:braid-cat}
\index{braid!category}
\index{category!of positive braids}
\index{Yang-Baxter relation}
\nomenclature[B]{$\pBraids$}{category of positive braids}
The category~$\pBraids$ of \emph{positive braids} contains all positive braid
groups $B_n^+$, see \secr{braid-mon}. It has natural numbers as objects, every
positive braid $b\in B_n^+$ induces a morphism $n\to n$, for every $n\in\N$, and
composition and identities are induced by multiplication and units of the
monoids~$B_n^+$. Otherwise said, if we consider the monoids $B_n^+$ as
one-object categories (see \secr{monoid}), we have
\[
  \pBraids
  =
  \coprod_{n\in\N}B_n^+
  \pbox.
\]
The expected monoidal structure makes a PRO out of it. As such, it admits a
presentation by a 3-polygraph with one 2-generator $\gamma:2\to 2$, called
\emph{braiding} and pictured as
\[
  \satex{tau}
\]
one 3-generator corresponding to the Yang-Baxter relation
\begin{equation}
  \label{eq:br-yb}
  \satex{yb_l}
  \qTO
  \satex{yb_r}
  \,\pbox.
\end{equation}

\subsection{Braids}
\label{sec:pres-braids}
\index{braid!groupoid}
\index{category!of braids}
\index{polygraph!of braids}
\index{Burau representation}
\nomenclature[B]{$\Braids$}{category of braids}
The category~$\Braids$ of \emph{braids} is defined similarly as
\[
  \Braids
  =
  \coprod_{n\in\N}B_n
\]
where $B_n$ is the $n$-th braid group, see \secr{braid-mon}, so that this
category is a groupoid.
It admits a presentation with two 2-generators $\gamma,\finv{\gamma}:2\to 2$,
\begin{align*}
  \satex{tau}
  &&
  \satex{tau2}  
\end{align*}
and three relations: the Yang-Baxter relation~\eqref{eq:br-yb}, as well as
\begin{align*}
  \satex{tau-inv-l}
  &\TO
  \satex{tau-inv-c}
  &
  \satex{tau-inv-r}
  &\TO
  \satex{tau-inv-c}
  \,\pbox.
\end{align*}

As an application, consider the ring $R=\Z[t,t^{-1}]$ of Laurent polynomials in
one variable~$t$. The category~$\Vect R$ of finitely generated free $R$-modules
is monoidal with the usual tensor product. Writing~$P$ for the above presentation,
we interpret a 1-cell $n$ of $P$ as $R^n$, the 2-generators $\gamma$ and $\finv\gamma$ as the
morphisms $R^2\to R^2$ whose matrix representations are respectively
\begin{align*}
  [\gamma]&=
  \pa{
    \begin{matrix}
      1{-}t&t\\
      1&0
    \end{matrix}
  }
  &
  [\finv\gamma]&=
  \pa{
    \begin{matrix}
      0&1\\
      t^{-1}&1{-}t^{-1}
    \end{matrix}
  }  
\end{align*}
This interpretation can be checked to be compatible with the relations of the
presentation and thus induces a monoidal functor
$f:\Braids\to\Vect{R}^{\Vect{R}}$, which is known as the \emph{Burau
  representation}~\cite{burau1935zopfgruppen}.  As a side note, this
representation is not faithful~\cite{bigelow1999burau} (\ie there are distinct
braids with the same image), but others are, such as Lawrence-Krammer
representation~\cite{bigelow2001braid,krammer2002braid,lawrence1990homological}.

\subsection{Permutations}
\label{sec:pres-sym}
\label{sec:pres-perm}
\index{permutation}
\index{symmetric!category}
\index{category!of permutations}
\index{polygraph!of permutations}
\nomenclature[S]{$\Bij$}{category of permutations}
The category $\Sym$ of \emph{permutations} (or sometimes \emph{symmetries}) is the monoidal category
\[
  \Sym
  =
  \coprod_{n\in\N}S_n
\]
where $S_n$ is the $n$-th symmetric group considered as a one-object category, see \secr{sym-group}.
Alternatively, the morphisms $m\to n$ can be described as the bijections
$\intset{m}\to\intset{n}$ where $\intset{n}$ is the set $\set{0,\ldots,n-1}$. A
presentation for this category can be obtained from the presentation of
$\Braids$ by adding the relation
\[
  \satex{sym_l}
  \TO
  \satex{sym_r}\,\text.
\]
In this case, the relation $\satex{tau}\TO\satex{tau2}$ can be derived, which is
why we generally remove one of the two generators, and note the remaining one as
\begin{equation}
  \label{eq:gamma}
  \satex{gamma}
\end{equation}
which is often called \emph{transposition} in this context. The two relations
of the presentation are thus
\begin{align}
  \label{eq:sym-rel}
  \satex{yb-l}
  &\TO
  \satex{yb-r}
  &
  \satex{sym-l}
  &\TO
  \satex{sym-r}
  \,\pbox.
\end{align}

A notion of canonical form (which is in fact a normal form for the above
rewriting system) was presented in \secr{sym-pres}: it consists of morphisms of
the form
\[
  \satex{empty}
  \qquad\qquad\text{or}\qquad\qquad
  \satex{sym-nf}
\]
where the diagram on the left is the empty diagram, and $\psi$ is a canonical form.

\subsection{Free braided and symmetric monoidal categories}
\label{sec:free-sym}
\index{braided!monoidal category}
\index{monoidal category!braided}
\index{category!braided monoidal}
A \emph{braided strict mo\-no\-idal category} $(C,\otimes,\monunit,\gamma)$ is a
strict monoidal category $(C,\otimes,\monunit)$ equipped with an invertible natural
transformation of components
\[
  \gamma_{u,v}
  :
  v\otimes u
  \to
  u\otimes v
\]
called \emph{braiding}, which is compatible with the monoidal structure: for
every objects $u,v,w\in C$,
\begin{align*}
  \gamma_{u\otimes v,w}&=(u\otimes\gamma_{v,w})\circ(\gamma_{w,u}\otimes v)\pbox,&\gamma_{\monunit,u}&=\unit{u}\pbox,\\
  \gamma_{u,v\otimes w}&=(\gamma_{u,v}\otimes w)\circ(v\otimes\gamma_{u,w})\pbox,&\gamma_{u,\monunit}&=\unit{u}\pbox.
\end{align*}%
\index{symmetric!monoidal category}%
\index{monoidal category!symmetric}%
A \emph{symmetric monoidal category} is a braided monoidal category in which the
braiding moreover satisfies $\gamma_{v,u}\circ\gamma_{u,v}=\unit{u\otimes v}$ for every
objects $u,v\in C$, in which case it is called a \emph{symmetry}.
A \emph{braided monoidal functor} $f:C\to D$ between braided monoidal
categories~$C$ and~$D$, with $\gamma^C$ and $\gamma^D$ as respective braidings,
is a monoidal functor which preserves the braiding, \ie
$f(\gamma^C_{u,v})=\gamma^D_{fu,fv}$ for every objects $u,v\in C$. A
\emph{symmetric monoidal functor} is a braided monoidal functor between
symmetric monoidal categories.

Given a monoidal category~$C$, there always exists a free braided monoidal
category $\tilde C$: it is a braided monoidal category equipped with a monoidal
functor $C\to\tilde C$ such that, given a braided monoidal category~$D$, a
monoidal functor~$C\to D$ extends uniquely as a braided monoidal functor
$\tilde C\to D$. A similar statement holds for symmetric, instead of braided,
monoidal categories. When~$C$ is presented by a 3-polygraph, $\tilde C$ admits
the following presentation.

\begin{theorem}
  \index{free!braided monoidal category}
  \label{thm:free-bmc}
  \label{thm:free-smc}
  Suppose given a 3-polygraph~$P$ presenting a monoidal category $\pcat{P}$, \ie
  $P_0=\set{\star}$. The free braided monoidal category on~$\pcat{P}$ is
  presented by the 3-polygraph~$Q$ such that
  \begin{align*}
    Q_0&=\set{\star}\\
    Q_1&= P_1\\
    Q_2&= P_2\sqcup\setof{\gamma_{a,b}:ba\To ab,\gamma_{a,b}^-:ab\To ba}{a,b\in P_1}\\
    Q_3&= P_3\sqcup\setof{G_{a,b},G_{a,b}',L_{a,\alpha},R_{a,\alpha}}{a,b\in P_1,\alpha\in Q_2}
  \end{align*}
  with relations
  \begin{align*}
    G_{a,b}:\gamma_{a,b}^-\comp{}\gamma_{a,b}&\TO\unit{ab}
    &
    G_{a,b}':\gamma_{a,b}\comp{}\gamma_{a,b}^-&\TO\unit{ab}
    \\
    \intertext{for $a,b\in P_1$ and, for $\alpha:u'\to u$,}
    L_{a,\alpha}:\alpha a\comp{}\gamma_{a,u}&\TO\gamma_{u',a}\comp{}a\alpha
    &
    R_{a,\alpha}:a\alpha\comp{}\gamma_{u,a}&\TO\gamma_{a,u'}\comp{}\alpha a    
  \end{align*}
  Above, $\gamma_{a,u}:ua\To au$ is a notation for the morphism defined
  inductively on~$u$ by
  \begin{align*}
    \gamma_{a,\monunit}&=\unit{a}
    &
    \gamma_{a,b}&=\gamma_{a,b}
    &
    \gamma_{a,bu}&=b\gamma_{a,u}\comp{}\gamma_{a,b}u
  \end{align*}  
  and similarly, for $\gamma_{u,a}:au\To ua$.
  Graphically, $\gamma_{a,b}$, $\gamma_{a,u}$, and $\gamma_{u,a}$ are
  respectively depicted as
  \begin{align*}
    \satex{tau-ba}
    &&
    \satex{tau-au}
    &&
    \satex{tau-ua}    
  \end{align*}
  and the relations $L_{a,\alpha}$ and $R_{a,\alpha}$ are respectively
  \begin{align}
    \label{eq:sym-nat}
    \satex{tau-L-l}
    &\TO
    \satex{tau-L-r}
    &
    \satex{tau-R-l}
    &\TO
    \satex{tau-R-r}    
  \end{align}
  The braiding on~$\pcat{Q}$ is the one which is given on objects $a,b\in Q_1$
  by~$\gamma_{a,b}$.
  The free symmetric monoidal category is presented by the 3-polygraph obtained
  from~$Q$ by removing the $2$-generators $\gamma_{a,b}^-$ and the associated
  relations $G_{a,b}$ and $G_{a,b}'$, and adding the relations
  \[
    I_{a,b}
    :
    \gamma_{b,a}\comp{}\gamma_{a,b}
    \TO
    \unit{ab}
  \]
  indexed $a,b\in P_1$ (see also \cref{sec:free-PROP}).
\end{theorem}

\noindent
In particular, we see that $\Braids$ is the free braided monoidal category on
the terminal category: there is only one 1-generator~$a$, one invertible
2-generator $\gamma_{a,a}$, and the two relations $L_{a,\gamma_{a,a}}$ and
$R_{a,\gamma_{a,a}}$ both correspond to the Yang-Baxter
relation~\eqref{eq:br-yb}. Similarly, $\Sym$ is the free symmetric monoidal
category on the terminal category. It can also be shown that~$\Braids$ (\resp
$\Sym$) is the terminal braided (\resp symmetric) monoidal category.

In every presentation of a braided (\resp symmetric) monoidal category, the
generators $\gamma_{a,b}$ are definable and the relations $L_{a,\alpha}$ and
$R_{a,\alpha}$ (\resp and~$I_{a,b}$) are derivable. Up to Tietze equivalence, we
can thus suppose that every presentation of a braided (\resp symmetric) monoidal
category is of the form given in the above theorem. For this reason, in a
presentation~$P$ with only one $1$-generator, we often say that a $2$-generator
$\gamma:2\to 2$ is a \emph{symmetry}\index{symmetry} when it is involutive (\ie
satisfies the equation on the right of \cref{eq:sym-rel}) and satisfies
relations \eqref{eq:sym-nat} for every $2$-generator~$\alpha$ (in particular,
for $\alpha=\gamma$, the Yang-Baxter relation has to be satisfied).

\section{Monoids}
\label{sec:2pres-mon}
%
Consider the category $\Fun$ where an object is a natural number and a morphism
$f:m\to n$ is a function $f:\intset{m}\to\intset{n}$. Alternatively, this
category can be described as the skeleton of the full subcategory of~$\Set$ on
finite sets. It can be equipped with a tensor product similar to the one of the
simplicial category, see \secr{ass-simpl-pres}, thus making it a PRO. This
category contains interesting subcategories, with the same objects, closed under
tensor product, with the following morphisms:
\begin{itemize}
\item the category $\Funinj$ of injective functions,
\item the category $\Funsurj$ of surjective functions,
\item the category $\Bij$ of bijections,
\item the category $\Simplaug$ of non-decreasing functions,
\item the category $\Simplinj$ of injective non-decreasing functions,
\item the category $\Simplsurj$ of surjective non-decreasing functions.
\end{itemize}
Details for this section can be found in~\cite{burroni1993higher,lack2004composing,lafont2003towards}.

\subsection{Non-decreasing injections}
\label{sec:incr-inj}
\nomenclature[.D]{$\Simplinj$}{category of injective non-decreasing functions}
The PRO $\Simplinj$ of injective non-decreasing functions admits a presentation
with one generator $\eta:0\to 1$, called \emph{unit} and pictured as
\[
  \satex{eta}
\]
and no relation. The generator is interpreted as the terminal function $0\to 1$
and, for instance, the interpretation of the diagram on the left is the
injective non-decreasing function $f:\intset{4}\to\intset{6}$ whose graph is
depicted on the right:
\begin{align*}
  \satex{eta-ex}
  \qquad\qquad\qquad\quad
  \vcenter{
    \xymatrix@C=1ex@R=1.5ex{
      0\ar@{-}[d]&1\ar@{-}[d]&2\ar@{-}[dr]&3\ar@{-}[drr]\\
      0&1&2&3&4&5\pbox.
    }
  }  
\end{align*}
The opposite PRO $\Simplinj^\op$ can be described as the category whose morphisms are
non-decreasing partial surjective functions.

\subsection{Non-decreasing surjections}
\label{sec:incr-surj}
\nomenclature[.D]{$\Simplsurj$}{category of surjective non-decreasing functions}
The PRO $\Simplsurj$ of non-decreasing surjective functions (already encountered in \cref{sec:simpl-surj-coh}) admits a
presentation with one generator $\mu:2\to 1$, called \emph{multiplication} and
pictured as
\[
  \satex{mu}
\]
and one relation
\[
  \satex{mon-assoc-l}
  \TO
  \satex{mon-assoc-r}
\]
The generator can be interpreted as the terminal function $2\to 1$, whose
graph is
\[
  \xymatrix@C=1ex@R=1.5ex{
    0\ar@{-}[d]&\ar@{-}[dl]1\\
    0
  }
\]
and, for instance, the interpretation of the diagram on the left is depicted on
the right:
\[
  \satex{mu-ex}
  \qquad\qquad\quad
  \vcenter{
    \xymatrix@C=0.5ex@R=5ex{
      0\ar@{-}[d]&1\ar@{-}[dl]&2\ar@{-}[dll]&3\ar@{-}[dl]&4\ar@{-}[d]&5\ar@{-}[dl]&6\ar@{-}[dll]&7\ar@{-}[dlll]&8\ar@{-}[dllll]&9\pbox.\ar@{-}[dlllll]\\
      0&&1&&2
    }
  }
\]
This category is thus the theory for \emph{semigroups}: a monoidal functor from
$\Simplsurj$ to~$\Set$ (with the monoidal structure induced by cartesian
product) consists of a set equipped with an associative binary operation.
The rewriting system is convergent and normal forms can be described as the
canonical form given by
\[
  \satex{empty}
  \qquad\qquad\text{or}\qquad\qquad
  \satex{mu-cf-id}
  \qquad\qquad\text{or}\qquad\qquad
  \satex{mu-cf-mu}
\]
where $\psi$ is a canonical form (the diagram on the left is the empty diagram).

\subsection{Non-decreasing functions}
\index{simplicial!category}
\index{augmented simplicial category}
\index{polygraph!of monoids}
\label{sec:incr-fun}
\label{sec:pres-simpl}
\nomenclature[.D]{$\Simplaug$}{augmented simplicial category}
The PRO $\Simplaug$ (the augmented simplicial category) can be presented with two
generators $\eta:0\to 1$ and $\mu:2\to 1$, pictured as above, together with the
rules of \secr{incr-inj,sec:incr-surj}:
\[
  \satex{mon-assoc-l}
  \TO
  \satex{mon-assoc-r}
\]
as well as the additional rules
\begin{align*}
  \satex{mon-unit-l}
  &\TO
  \satex{mon-unit-c}
  &
  \satex{mon-unit-r}
  &\TO
  \satex{mon-unit-c}
  \,\pbox.
\end{align*}
These define a distributive law $\dlaw$ between the two previous categories so
that
\[
  \Simplaug
  =
  \Simplsurj\otimes_\dlaw\Simplinj
  \pbox.
\]
The rewriting system is convergent, see \secr{ass-simpl-pres}, and normal forms
can be described as the following canonical forms:
\begin{align*}
  \satex{empty}
  &&\text{or}&&
  \satex{mu-cf-id}
  &&\text{or}&&
  \satex{mu-cf-mu}
  &&\text{or}&&
  \satex{mu-cf-eta}  
\end{align*}
where $\psi$ is a canonical form; moreover, in the third case, we suppose that
$\psi$ is not of the form given in the fourth case.

An alternative notion of canonical form (which is not a normal form for a rewriting system, for similar reasons as in \cref{ex:cf-not-nf}) is
\begin{align*}
  \satex{empty}
  &&\text{or}&&
  \satex{mu-cf-mu}
  &&\text{or}&&
  \satex{mu-cf-eta}
  \,\pbox.
\end{align*}
This means that morphisms in $\Simplaug$ can be shown to be in bijection with
$2$-cells of the above form. Note that, here, the canonical form of an identity
does not only consist of wires.

From this presentation, one sees that $\Simplaug$ is the theory for
\emph{monoids}: the monoidal functors from $\Simplaug$ to a monoidal
category~$C$ correspond to monoids~$C$ (see \exr{ass-alg}). We can also deduce
that it is the free cocartesian category on the terminal category, see
\cref{sec:trs-3pol,sec:free-cart-cat}.

\index{Tamari lattice}
Writing $P$ for the above polygraph, the morphisms from $n$ to $1$ in
$\freecat P_2$, for some $n\in\N$ can be ordered by $\phi\leq\psi$ whenever
there exists a rewriting path $\phi\TO\psi$. The resulting poset is a
well-studied lattice called the $n$-th \emph{Tamari
  lattice}~\cite{friedman1967problemes}.

\subsection{Injections}
\nomenclature[F]{$\Funinj$}{category of injective functions}
The PRO $\Funinj$ admits a presentation with two generators
\begin{align*}
  \satex{gamma}
  &&
  \satex{eta}  
\end{align*}
satisfying the relations for symmetries (\secr{pres-sym}) and non-decreasing
injections (\secr{incr-inj})
\begin{align*}
  \satex{yb-l}
  &\TO
  \satex{yb-r}
  &
  \satex{sym-l}
  &\TO
  \satex{sym-r}
  \\
  \intertext{as well as the compatibility relations}
  \satex{eta-gamma-l}
  &\TO
  \satex{eta-gamma-r}
  &
  \satex{eta-gamma2-l}
  &\TO
  \satex{eta-gamma2-r}
  \,\pbox.
\end{align*}
Those define a distributive law $\dlaw$ between $\Sym$ and $\Simplinj$, so
that
\[
  \Funinj
  =
  \Sym\otimes_\dlaw\Simplinj
  \pbox.
\]
The resulting presentation is convergent with normal forms being given by
\begin{align*}
  \satex{f-cf-eta}
  &&\text{or}&&
  \satex{f-cf-gamma}
  \pbox.
\end{align*}
Note that the second (or the first) compatibility relation is redundant since we
have
\[
  \satex{eta-gamma2-l1}
  \OT
  \satex{eta-gamma2-c1}
  \TO
  \satex{eta-gamma2-r1}
  \,\pbox.
\]

\subsection{All functions}
\index{polygraph!of commutative monoids}
\label{sec:pres-com-mon}
\nomenclature[F]{$\Fun$}{category of functions}
The PRO~$\Fun$ corresponds to commutative monoids. It admits a presentation with
three generators
\begin{align*}
  \satex{mu}
  &&
  \satex{eta}
  &&
  \satex{gamma}
\end{align*}
and relations which consist of those for monoids (\secr{incr-fun})
\begin{align*}
  \satex{mon-assoc-l}
  &\TO
  \satex{mon-assoc-r}
  &
  \satex{mon-unit-l}
  &\TO
  \satex{mon-unit-c}
  &
  \satex{mon-unit-r}
  &\TO
  \satex{mon-unit-c}
\end{align*}
those for symmetries (\secr{pres-sym})
\begin{align*}
  \satex{yb-l}
  &\TO
  \satex{yb-r}
  &
  \satex{sym-l}
  &\TO
  \satex{sym-r}     
\end{align*}
compatibility relations
\begin{align}
  \label{eq:mon-sym-compat}
  \satex{mu-gamma-l}
  &\TO
  \satex{mu-gamma-r}
  &
  \satex{gamma-mu-l}
  &\TO
  \satex{gamma-mu-r}
  &
  \satex{eta-gamma-l}
  &\TO
  \satex{eta-gamma-r}
  &
  \satex{gamma-eta-l}
  &\TO
  \satex{gamma-eta-r}
\end{align}
and the commutativity relation
\begin{equation}
  \label{eq:mon-com}
  \satex{mon-com-l}
  \TO
  \satex{mon-com-r}
  \,\pbox.
\end{equation}
This rewriting system is convergent. A notion of canonical form (which does not
exactly coincide with normal forms for similar reasons as in \cref{ex:cf-not-nf}) can
be given by
\begin{align*}
  \satex{f-cf-eta}
  &&\text{or}&&
  \satex{f-cf-mu}
  \pbox.
\end{align*}

This monoidal category is the theory for \emph{commutative monoids}:
a symmetric monoidal functor to a symmetric monoidal
category~$C$ corresponds to commutative monoid in~$C$.
It can be obtained as a composite PROP $\Fun=\Funsurj\otimes_\dlaw\Funinj$ along
the expected distributive law.

\subsection{Partial functions}
\nomenclature[F]{$\Funp$}{category of partial functions}
Consider the PRO $\Funp$ where a morphism $f:m\to n$ is a partial function
$f:\intset{m}\to\intset{n}$. It admits a presentation with four generators
\begin{align*}
  \satex{mu}
  &&
  \satex{eta}
  &&
  \satex{gamma}
  &&
  \satex{eps}  
\end{align*}
satisfying the relations of \secr{pres-com-mon} together with
\begin{align*}
  \satex{eps-gamma-l}
  &\TO
  \satex{eps-gamma-r}
  &
  \satex{gamma-eps-l}
  &\TO
  \satex{gamma-eps-r}
  &
  \satex{mu-eps-l}
  &\TO
  \satex{mu-eps-r}
  &
  \satex{eta-eps-l}
  &\TO
  \satex{empty}
\end{align*}
The resulting rewriting system is convergent and canonical forms can be given by
\begin{align*}
  \satex{f-cf-eta}
  &&\text{or}&&
  \satex{f-cf-mu}
  &&\text{or}&&
  \satex{f-cf-eps}
  \,\pbox.  
\end{align*}
The category can be described as a composite PROP
$\Funp=\Funinj^\op\otimes_\dlaw\Fun$ along the expected distributive law, as
well as the composite PRO $\Fun=\Simplinj^\op\otimes_\dlaw\Fun$ along the
expected distributive law. In the last case, the associated factorization system
given by \propr{sfs-dlaw} is the usual factorization of a partial
function~$f:X\to Y$ as a partial non-decreasing injection (given by the
canonical partial function $X\to\operatorname{dom}(f)$ where
$\operatorname{dom}(f)\subseteq X$ is the domain of~$f$) followed by a total
function (the restriction of~$f$ to its domain).

Similar presentations exist for other variants of the PRO (e.g., partial
non-decreasing surjective functions). In particular, for the category of partial
injective non-decreasing functions~$\Simplpinj$, we obtain a decomposition as
\[
  \Simplpinj
  =
  \Simplinj^\op\otimes_\dlaw\Simplinj
  =
  \qCospan(\Simplinj)
  \pbox.
\]

\subsection{Symmetric monoids}
\label{sec:sym-mon}
\index{monoid!symmetric}
\index{symmetric!monoid}
The theory of \emph{symmetric monoids} is obtained from the theory of
commutative monoids of \secr{pres-com-mon} by removing the commutativity
relation~\eqref{eq:mon-com}: by \thmr{free-smc}, this is the free symmetric
monoidal category on~$\Simplaug$, and therefore it is the PROP of monoids. The
relations induce a distributive law $\dlaw$ so that this PRO can be obtained as
$\Bij\otimes_\dlaw\Simplaug$, see~\cite{lack2004composing}. A direct description of
this category can be given as the PRO where a morphism $f:m\to n$ is a function
$f:\intset{m}\to\intset{n}$ together with total order on each of the sets
$f^{-1}(i)$ for $0\leq i<n$, equipped with suitable
composition~\cite{pirashvili2002p}.

\subsection{Braided monoids}
\index{monoid!braided}
\index{braided!monoid}
As a variant of the theory for commutative monoids, one can consider
the one for
\emph{braided monoids}. It admits a presentation with four generators
\begin{align*}
  \satex{mu}
  &&
  \satex{eta}
  &&
  \satex{tau}
  &&
  \satex{tau2}  
\end{align*}
and the laws are similar to those of commutative monoids (\secr{pres-com-mon})
excepting that the transpositions satisfy the laws for
braidings (\secr{pres-braids}) instead of symmetries
(\secr{pres-sym}). As in \secr{sym-mon}, one could also considered the variant
without the commutativity relation, \ie the free braided monoidal category
on~$\Simplaug$, but in practice braided monoids are always understood commutative.

\subsection{Comonoids}
\label{sec:comonoid}
\index{comonoid}
Dual categories are also interesting. For instance, $\Simplaug^\op$ is the theory
representing comonoids and $\Fun^\op$ cocommutative comonoids.

\subsection{Free cartesian categories}
\index{free!cartesian category}
\label{sec:free-cart-cat}
The category~$\Fun^\op$ can be characterized as being the free cartesian
category on the terminal one (and of course the case of~$\Fun$ is dual, but the
opposite is more commonly used from this perspective). We briefly describe here
the situation and refer the reader to \cref{sec:trs-3pol}, where cartesian
categories are considered in detail.

Any cartesian category has an underlying symmetric monoidal category, where the
tensor product is induced by the cartesian product and unit is the terminal
object. Not every symmetric monoidal category can be obtained in this way: those
who can are characterized by the fact that every object is equipped with a
structure of commutative comonoid, in a natural way, see \cref{thm:mon-cart}.
From this observation, one can derive the characterization of free
cartesian categories given in \thmr{pres-free-cart}. In particular, the
category~$\Fun^\op$ is the free category on the terminal category.

\section{Distributive Laws}
\label{sec:dlaw}

\subsection{Distributive laws between monads}
\index{distributive law!of monads}
A monad\index{monad} is a particular case of a monoid, as already mentioned in
\exr{ass-alg}. Namely, a 2-functor from the 2-category~$\Simplaug$ (which
represents monoids, see \secr{pres-simpl}) to the 2\nbd-cate\-gory~$\Cat$ (of
categories, functors and natural transformations) amounts to the data of
\begin{itemize}
\item a category~$C$ (the image of $\star$),
\item a functor $T:C\to C$ (the image of~$a$),
\item natural transformations $\mu:T\circ T\To T$ and $\eta:\unit{C}\To T$ (the
  respective images of $\mu$ and $\eta$)
\end{itemize}
satisfying axioms so that $(T,\mu,\eta)$ is a monad on~$C$ (see \cref{sec:bicat-dlaw}).

As a variant, a theory corresponding to pairs of monads on a same
category can easily be constructed. The category $\Simplaug\sqcup\Simplaug$, the
coproduct of $\Simplaug$ with itself as monoidal categories
admit a presentation with two 1\nbd-gene\-rators $a$, $b$ and four
2\nbd-generators
\begin{align*}
  \mu_a:aa&\To a
  &
  \eta_a:\unit{}&\To a
  &
  \mu_b:bb&\To b
  &
  \eta_b:\unit{}&\To b
\end{align*}
respectively pictured as
\begin{align*}
  \satex{mu-a}
  &&
  \satex{eta-a2}
  &&
  \satex{mu-b}
  &&
  \satex{eta-b2}
\end{align*}
such that the pair $(\mu_a,\eta_a)$ satisfies the laws of monoids, see
\secr{pres-simpl}, as well as the pair $(\mu_b,\eta_b)$.

More interestingly, the previous theory can be modified in order to present a
pair of monads on a same category, together with a distributive law between
them. We recall that a \emph{distributive law} between two monads $T$ and $U$
consists of a natural transformation $\dlaw:T\circ U\To U\circ T$ satisfying
four suitable axioms, see \secr{bicat-dlaw} or~\cite{beck1969distributive}. This
theory~$\DLaw$ can be obtained from the one for $\Simplaug\sqcup\Simplaug$ by adding a
2-generator $\lambda:ba\To ab$, pictured as
\[
  \satex{dlaw}
\]
and the four relations
\begin{align*}
  \satex{dlaw-mu-a-r}
  &\TO
  \satex{dlaw-mu-a-l}
  &
  \satex{dlaw-mu-b-l}
  &\TO
  \satex{dlaw-mu-b-r}
  \\
  \satex{dlaw-eta-a-l}
  &\TO
  \satex{dlaw-eta-a-r}
  &
  \satex{dlaw-eta-b-l}
  &\TO
  \satex{dlaw-eta-b-r}
\end{align*}
This presentation is convergent.
It was shown by Beck~\cite{beck1969distributive} that a distributive law between
monads $T$ and $U$ induces a structure of monad on the composite endofunctor
$U\circ T$. This can be rephrased in the above setting as follows:

\begin{theorem}
  \label{thm:comp-monad}
  In $\DLaw$, the two following morphisms induce the structure of a monoid on $ab$:
  \begin{align*}
    \satex{dlaw-mu}
    &&
    \satex{dlaw-eta}
  \end{align*}
\end{theorem}
\begin{proof}
  Associativity is shown by the following derivation
  \[
    \satex{dlaw-assoc1}
    \TO
    \satex{dlaw-assoc2}
    \TO
    \satex{dlaw-assoc3}
    \TO
    \satex{dlaw-assoc4}
  \]
  and left and right neutrality by
  \[
    \satex{dlaw-unitl1}
    \TO
    \satex{dlaw-unitl2}
    \TO
    \satex{dlaw-unitl3}
    \qtand
    \satex{dlaw-unitr1}
    \TO
    \satex{dlaw-unitr2}
    \TO
    \satex{dlaw-unitr3}
  \]
  which concludes the proof.
\end{proof}

For instance, on~$\Set$, the monad of rings can be obtained as
$U\otimes_\dlaw T$ where~$T$ is the monad of free monoids, $U$ and the monad of
free abelian groups and $\dlaw:TU\To UT$ is the usual distributive law, which
sends the formal expression of the form $(a+b)(c+d)$ to $ac+ad+bc+bd$.

\subsection{Iterated distributive laws}
\label{sec:idlaw}
\index{distributive law!iterated}
It is useful to iterate this construction: given three monads $S$, $T$, and $U$ on a
category~$C$ and distributive laws between $S$ and $T$, $S$ and $U$, and $T$ and
$U$, we would like to have a monad structure on the composite $U\circ T\circ
S$. This is the case where the distributive laws satisfy a suitable compatibility
axiom~\cite{cheng2011iterated}. The theory~$\IDLaw$ for \emph{iterated
  distributive laws} axiomatizes such a situation. It has three 1-generators
$a$, $b$, and $c$, nine 2\nbd-generators
\begin{align*}
  \mu_a&:aa\To a&\mu_b&:bb\To b&\mu_c&:cc\To c&\dlaw_{ab}&:ba\To ab\\
  \eta_a&:\unit{}\To a&\eta_b&:\unit{}\To b&\eta_c&:\unit{}\To c&\dlaw_{ac}&:ca\To ac\\
  &&&&&&\dlaw_{bc}&:cb\To bc
\end{align*}
pictured as
\begin{align*}
  \satex{mu-a}
  &&
  \satex{eta-a2}
  &&
  \satex{mu-b}
  &&
  \satex{eta-b2}
  &&
  \satex{mu-c}
  &&
  \satex{eta-c2}
  &&
  \satex{dlaw-ab}
  &&
  \satex{dlaw-ac}
  &&
  \satex{dlaw-bc}  
\end{align*}
satisfying relations expressing that $(a,\mu_a,\eta_a)$, $(b,\mu_b,\eta_b)$, and
$(c,\mu_c,\eta_c)$ are mo\-noids, $\dlaw_{ab}$, $\dlaw_{bc}$, and $\dlaw_{ac}$ are
distributive laws, and the additional axiom
\[
  \satex{dlaw-yb-l}
  \TO
  \satex{dlaw-yb-r}
\]
reminiscent of the Yang-Baxter relation. This is illustrated in \secr{R-mat}.

In order to compose four or more monads, one might at first think that we should
need new axioms, but it is in fact enough to assume the above axioms for every
triple of monads in order to compose an arbitrary number of those,
see~\cite{cheng2011iterated} and \secr{iterated-dlaw}.

\subsection{Other distributive laws}
Variants of the notion of distributive law have been studied in the
literature. For instance, the expected notion of distributive law between a
monad and a comonad is studied in~\cite{power2002combining}. A weaker notion of
distributive law has also been studied by Street~\cite{street2009weak}: here the
two relations involving units have been replaced by the unique relation
\[
  \satex{dlaw-weak-l}
  \TO
  \satex{dlaw-weak-r}
\]
or, equivalently, the two relations
\begin{align*}
  \satex{dlaw-weak1-r}
  &\TO
  \satex{dlaw-weak1-l}
  &
  \satex{dlaw-weak2-r}
  &\TO
  \satex{dlaw-weak2-l}
  \,\pbox.
\end{align*}

\subsection{Linear non-linear terms}
\index{linear non-linear term}
A functor $T:\Fun\to\Set$ is an abstract way to encode a collection of
terms. Namely, for $n\in\N$, the set $Tn$ can be seen as the collection of terms
$t(x_0,\ldots,x_{n-1})$ with $n$ free variables, and given a morphism
$f:m\to n$, \ie a function $f:\intset{m}\to\intset{n}$, the function $Tf:Tm\to Tn$ is the
reindexing function, sending a term $t(x_0,\ldots,x_{n-1})$ to the term
$t(x_{f(0)},\ldots,x_{f(n-1)})$ obtained by replacing the variable $x_i$ by
$x_{f(i)}$, see for instance~\cite{fiore1999abstract}. Similarly, a functor
$T:\Bij\to\Set$ encodes a collection of ``linear'' terms: $Tn$ consists of terms
in which each $x_i$, for $0\leq i<n$, occurs exactly once and, for this reason,
we can only permute variables, and not merge two of them or forget one of them.

Now suppose that we are interested in a ``mixed'' situation where a term can
have both linear and non-linear variables. A linear variable can always be
considered as a non-linear one by forgetting about the fact that it should
occur exactly once. Those are naturally modeled by the following
theory~\cite{fiore2007towards,hyland2020lnl}, with two 1-generators $a$ and $b$,
six 2-generators
\begin{align*}
  \mu_a:aa&\To a
  &
  \eta_a:\unit{}&\To a
  &
  \sigma:b&\To a
  \\
  \gamma_{aa}:aa&\To aa
  &
  \gamma_{ab}:ba&\To ab
  &
  \gamma_{bb}:bb&\To bb
\end{align*}
pictured as
\begin{align*}
  \satex{mu-a}
  &&
  \satex{eta-a2}
  &&
  \satex{sigma-ba}
  &&
  \satex{gamma-aa}
  &&
  \satex{gamma-ba}
  &&
  \satex{gamma-bb}  
\end{align*}
such that the axioms of symmetries hold for whichever typing of the wires,
$(\mu_a,\eta_a,\gamma_{aa})$ satisfies the axioms of commutative monoids, and
the axioms
\begin{align*}
  \satex{lnl1-l}
  &\TO
  \satex{lnl1-r}
  &
  \satex{lnl2-l}
  &\TO
  \satex{lnl2-r}  
\end{align*}
hold for whichever typing of the wires. The object $a$ thus corresponds to a
non-linear variable (which can duplicated, erased, and exchanged with other
variables), the object $b$ to linear variable (which can only be exchanged with
other variables), and the morphism $\sigma$ to the fact that we can consider a
linear variable as a non-linear one.

\section{Bialgebras}
\subsection{Matrices}
\index{category!of matrices}
\index{matrix}
\nomenclature[M]{$\Mat R$}{category of matrices with coefficients in $R$}
Given a semiring~$R$, we write $\Mat R$ for the category whose objects are
natural numbers and a morphism $f:m\to n$ is an $n{\times}m$-matrix with coefficients
in~$R$, with usual composition and identities. This category is a PRO with the
usual direct sum of matrices.

When~$R$ is a ring, the category~$\Mat R$ is equivalent to the category of
finite-dimensional $R$-modules and $R$-linear maps.  When~$\kk$
is a field, the category~$\Mat\kk$ is equivalent to~$\Vect\kk$, the category of
vector spaces and $\kk$-linear maps. Writing~$K$ for the semiring of small cardinals,
with disjoint union as addition and cartesian product as product, the
category~$\Mat K$ is equivalent to the full subcategory of $\qSpan(\Set)$ (the
category of isomorphism classes of spans, see \secr{qspan}) on finite sets.

\subsection{Multirelations}
\index{multirelation}
Given two sets~$X$ and~$Y$, a \emph{multirelation} from~$X$ to~$Y$ is a function
$X\times Y\to\N$ such that the set
\[
  \setof{(x,y)\in X\times Y}{f(x,y)\neq 0}
\]
is finite. Any relation $R\subseteq X\times Y$ is canonically seen as the
multirelation~$f$ such that $f(x,y)=1$ if $(x,y)\in R$ and $f(x,y)=0$
otherwise. Conversely, a multirelation can be thought of as a relation with
multiplicities: given $(x,y)\in X\times Y$, $f(x,y)$ is called the
\emph{multiplicity} of the relation $(x,y)$. Given two multirelations $f:X\to Y$
and $g:Y\to Z$, their composite is given by
\[
  (g\circ f)(x,z)
  =
  \sum_{y\in Y}f(x,y)\times g(y,z)
\]
(note that the sum only involves a finite number of non-zero terms) and
the identity on a set~$X$ is such that $\unit{X}(x,x')=0$ if $x\neq x'$ and
$\unit{X}(x,x)=1$.
\nomenclature[MRel]{$\MRel$}{category of multirelations}
We write $\MRel$ for the category of sets and
multirelations. This is a monoidal category when equipped with the tensor
product induced by disjoint union.
A morphism $f:X\to Y$ in this category can be seen as a span
\begin{equation}
  \label{eq:mrel-span}
  \xymatrix@C=3ex@R=3ex{
    X&\ar[l]_-sR\ar[r]^-t&Y
  }
\end{equation}
where $R$ is a finite set, such that the image under~$f$ of $(x,y)\in X\times Y$
is the cardinal of the following set:
\[
  f(x,y)
  =
  \cardinal{\setof{r\in R}{\text{$x=s(r)$ and $y=t(r)$}}}
\]
making~$\MRel$ a subcategory of~$\qSpan(\Set)$, the category of sets and
isomorphism classes of spans of functions.

The full subcategory of~$\MRel$ on finite sets
is equivalent to $\Mat\N$, the category of matrices with
coefficients in~$\N$. Namely, a morphism $f:m\to n$ is an $n{\times}m$-matrix
with coefficients in~$\N$, which corresponds to the multirelation from
$\intset{m}$ to $\intset{n}$ such that the multiplicity of
$(i,j)\in\intset{m}\times\intset{n}$ is $f(i,j)$. Note that the category $\Mat\N$ is
isomorphic to $\qSpan(\Fun)$. This category can also be described as the PRO
where a morphism $f:m\to n$ is a morphism of free finitely generated monoids
$f:\N^m\to\N^n$, see~\cite{pirashvili2002p}.

\subsection{Bialgebras}
\label{sec:bialgebra}
\index{bialgebra}
\index{bimonoid}
\index{polygraph!of bialgebras}
The PRO~$\Mat\N$ admits a presentation with generators
\begin{align*}
  \mu&:2\To 1
  &
  \eta&:0\To 1
  &
  \delta&:1\To 2
  &
  \varepsilon&:1\To 0
  &
  \gamma&:2\To 2\\
  \intertext{pictured as}
  &\satex{mu}
  &&
  \satex{eta}
  &&
  \satex{delta}
  &&
  \satex{eps}
  &&
  \satex{gamma}  
\end{align*}
and relations expressing that $\gamma$ is a symmetry (\cref{sec:free-sym}),
$(\mu,\eta,\gamma)$ is a commutative monoid (\cref{sec:pres-com-mon}),
$(\delta,\varepsilon,\gamma)$ is a cocommutative comonoid (\cref{sec:comonoid})
and the compatibility relations
\begin{align*}
  \satex{bialg-l}
  &\TO
  \satex{bialg-r}
  &
  \satex{mu-eps-l}
  &\TO
  \satex{mu-eps-r}
  &
  \satex{eta-delta-l}
  &\TO
  \satex{eta-delta-r}
  &
  \satex{eta-eps-l}
  &\TO
  \satex{empty}
  \,\pbox.
\end{align*}
The interpretations of the generators are the following multirelations,
represented as matrices
\begin{align*}
  \mu&=
  \pa{
    \begin{matrix}
      1&1
    \end{matrix}
  }
  &
  \delta&=
  \pa{
    \begin{matrix}
      1\\
      1
    \end{matrix}
  }
  &
  \gamma&=
  \pa{
    \begin{matrix}
      0&1\\
      1&0
    \end{matrix}
  }  
\end{align*}
\ie graphically,
\begin{align*}
  \satex{mu-x}
  &&
  \satex{eta-x}
  &&
  \satex{delta-x}
  &&
  \satex{eps-x}
  &&
  \satex{gamma-x}
  \pbox.
\end{align*}
A notion of canonical form is given by
\begin{align*}
  \satex{empty}
  &&\text{or}&&
  \satex{mrel-cf-mu}
  &&\text{or}&&
  \satex{mrel-cf-eps}
  &&\text{or}&&
  \satex{mrel-cf-eta}
\end{align*}
up to some relations, for which normal forms can be given,
see~\cite{mimram2011structure}. For instance, the multirelation $f:3\to 4$ whose
graph is shown on the left (we link two points as many times as their
multiplicity in the multirelation) is represented by the canonical form on the
right:
\[
  \vcenter{
    \xymatrix@C=3ex@R=15ex{
      0\ar@{-}@/_/[d]\ar@{-}[d]\ar@{-}@/^/[d]\ar@{-}[drrr]&1&2\ar@{-}[d]\ar@{-}[dr]\\
      0&1&2&3
    }
  }
  \qquad\qquad\qquad\qquad
  \satex{mrel-cf-ex}
\]
This example should make more clear the correspondence between the category and its presentation,
and details can be found in various places~\cite{lack2004composing,lafont2003towards,mimram2011structure,pirashvili2002p}. This presentation can be obtained
from the one of~$\Fun$ given in \secr{pres-com-mon}, using the methods of
\secr{3-dlaws} for composing PROPs: we have
$\qSpan(\Fun)=\Fun^\op\otimes_\dlaw\Fun$ where the distributive law
$
\dlaw:\Fun\otimes_\Bij\Fun^\op\to\Fun\otimes_\Bij\Fun^\op
$
is given by pullback, see~\cite{lack2004composing}.

This is the theory for an algebraic structure called a \emph{bicommutative
  bialgebra}, or \emph{bimonoid}. More generally, if we drop the requirement
that the monoid and the comonoid structures should be commutative, we obtain the
theory of \emph{bialgebras} (which are said to be
commutative/cocommutative/bicommutative when the monoid/comonoid/both structures
are commutative). Bialgebras are generally considered in the category~$\Vect\kk$
of vector spaces over a fixed field~$\kk$.
For instance, given a monoid~$(M,\cdot,1)$, the vector space~$\kk M$ generated by
the set~$M$ is canonically a cocommutative bialgebra: writing $e_a$ with
$a\in M$ for a basis vector, the interpretations of the various morphisms are
given by
\begin{align*}
  \mu:\kk M\otimes\kk M&\to\kk M
  &
  \eta:1&\to\kk M
  \\
  e_a\otimes e_b&\mapsto e_{a\cdot b}
  &
  1&\mapsto e_1
  \\[2ex]
  \delta:\kk M&\to\kk M\otimes\kk M
  &
  \varepsilon:\kk M&\to 1
  &
  \gamma:\kk M\otimes\kk M&\to\kk M\otimes\kk M
  \\
  e_a&\mapsto e_a\otimes e_a
  &
  e_1&\mapsto 1
  &
  e_a\otimes e_b&\mapsto e_b\otimes e_a
\end{align*}

\subsection{Relations}
\index{relation}
\index{bialgebra!special}
\index{special bialgebra}
\nomenclature[Rel]{$\Rel$}{category of relations}
\label{sec:pres-rel}
\label{sec:special-bialgebra}
Consider the category~$\Rel$ whose objects are sets and where a morphism
from~$X$ to~$Y$ is a relation from~$X$ to~$Y$, \ie a subset
$R\subseteq X\times Y$. Explicitly, the composite of two relations
$R\subseteq X\times Y$ and $S\subseteq Y\times Z$ is the relation
$S\circ R\subseteq X\times Z$ defined by
\[
  S\circ R
  =
  \setof{(x,z)\in X\times Z}{\exists y\in Y, (x,y)\in R \land (y,z)\in S}
  \pbox.
\]
This category is monoidal with tensor product given on objects by disjoint union
of sets.
It may be described as a variant of category of~$\qSpan(\Set)$ where morphisms
are isomorphism classes of spans of the form~\eqref{eq:mrel-span}, such that $s$
and $t$ are jointly monic, \ie two distinct $x$ and $y$ give rise to distinct
pairs of images $(s(x),t(x))$ and $(s(y),t(y))$.

The full subcategory whose objects are finite sets is equivalent to the
category~$\Mat\Bool$ of matrices over the semiring of booleans (with $\lor$ as
addition and $\land$ as multiplication).
It admits a presentation obtained from the theory of bicommutative bialgebras by
further adding the relation
\begin{equation}
  \label{eq:rel-special}
  \satex{rel-l}
  \TO
  \satex{rel-r}
  \,\pbox.
\end{equation}
If we remember that the theory of bicommutative bialgebras is a presentation of
the category $\Mat\N$ of matrices with coefficients in~$\N$, adding the
relation \eqref{eq:rel-special} amounts to quotient coefficients in~$\N$ in
matrices by the relation $1+1=1$. Otherwise said, we describe~$\Mat\Bool$ as the
quotient of~$\Mat\N$ under the relation identifying two morphisms $f,g:m\to n$
which have the same non-zero coefficients, \ie when $f(i,j)=0$ iff $g(i,j)=0$ for
every $(i,j)\in\intset{m}\times\intset{n}$. This is the theory for bialgebras
which are called \emph{special}, \emph{relational}~\cite{hyland2000symmetric},
or \emph{qualitative}~\cite{mimram2011structure}.
Note that the category of relations can be described as the following pushout
in~$\MonCat$ (the category of monoidal categories an strict monoidal functors):
\nomenclature[MonCat]{$\MonCat$}{category of monoidal categories}
\[
  \xymatrix@C=3ex@R=3ex{
    \Funinj\sqcup\Funinj^\op\ar[d]\ar[r]&\ar@{.>}[d]\qSpan(\Funinj)\\
    \qCospan(\Funinj)\ar@{.>}[r]&\Mat\Bool
  }
\]
where the arrows are the canonical inclusions, from which the presentation
of~$\Mat\Bool$ can be deduced,
see \cite{fong2017universal} for details and generalizations of this situation.

\subsection{Matrices over $\Z/2\Z$}
The category~$\Mat{\Z/2\Z}$ of natural numbers and matrices with coefficients
over the ring $\Z/2\Z$ admits a presentation obtained from the theory of
bialgebras by adding the relation
\[
  \satex{asep-l}
  \TO
  \satex{asep-r}
\]
which amounts to quotient coefficients in~$\N$ by $1+1=0$. This is the theory
for bialgebras which are \emph{anti-separable}~\cite{lafont2003towards}. More
generally, matrices over $\Z/n\Z$ can be obtained by replacing the above
relation by a relation of the form
\[
  \satex{nsep-l}
  \TO
  \satex{nsep-r}
\]
where the diagram on the left contains $n-1$ morphisms $\satex{delta-small}$
and morphisms $\satex{mu-small}$.
An alternative confluent presentation is given
in~\cite[Section~3.2]{lafont2003towards} and shown to be terminating
in~\cite[Section~7]{guiraud2006termination}.

\subsection{Matrices over~$\Z$, Hopf algebras}
\label{sec:hopf}
\index{Hopf algebra}
\index{algebra!Hopf}
\index{polygraph!of Hopf algebras}
The category~$\Mat{\Z}$ of natural numbers and matrices over~$\Z$ admits a
presentation obtained from the theory of bialgebras by adding a generator
$\sigma:1\to 1$, called \emph{antipode}, pictured as
\[
  \satex{sigma}
\]
and interpreted as the $1{\times}1$-matrix $(\begin{matrix}-1\end{matrix})$,
together with the relations
\begin{align*}
  \satex{sigma-mu-l}
  &\TO
  \satex{sigma-mu-r}
  &
  \satex{sigma-eta-l}
  &\TO
  \satex{sigma-eta-r}
  &
  \satex{sigma-gamma-l}
  &\TO
  \satex{sigma-gamma-r}
  &
  \satex{sigma-mul-l}
  &\TO
  \satex{sigma-mul-r}
  \,\pbox.
  \\
  \satex{sigma-delta-l}
  &\TO
  \satex{sigma-delta-r}
  &
  \satex{sigma-eps-l}
  &\TO
  \satex{sigma-eps-r}
  &
  \satex{sigma-gamma2-l}
  &\TO
  \satex{sigma-gamma2-r}
\end{align*}
This is the theory for bicommutative \emph{Hopf algebras}. We have seen in
\secr{bialgebra} that every monoid induces a bialgebra; when this monoid is a
group~$G$, the bialgebra~$\kk G$ is canonically a Hopf algebra, the interpretation of
the antipode being given by inverses:
\begin{align*}
  \sigma:\kk G&\to\kk G\\
  e_a&\mapsto e_{-a}\,\text.
\end{align*}

\subsection{Matrices over an arbitrary semiring}
\label{sec:linear-bialg}
More generally, given a semiring~$R$, the PRO $\Mat R$ admits a presentation
containing the theory of matrices plus a generator
\[
  \satex{coef}
\]
for every $a\in R$, interpreted as the $1{\times}1$-matrix
$(\begin{matrix}a\end{matrix})$, and relations
\begin{align*}
  \satex{mat-mu-l}
  &\TO
  \satex{mat-mu-r}
  &
  \satex{mat-eta-l}
  &\TO
  \satex{mat-eta-r}
  &
  \satex{mat-mul-l}
  &\TO
  \satex{mat-mul-r}
  &
  \satex{mat-one-l}
  &\TO
  \satex{mat-one-r}
  \\
  \satex{mat-delta-l}
  &\TO
  \satex{mat-delta-r}
  &
  \satex{mat-eps-l}
  &\TO
  \satex{mat-eps-r}
  &
  \satex{mat-add-l}
  &\TO
  \satex{mat-add-r}
  &
  \satex{mat-zero-l}
  &\TO
  \satex{mat-zero-r}
  \\
  \satex{mat-gamma-l}
  &\TO
  \satex{mat-gamma-r}
  &
  \satex{mat-gamma2-l}
  &\TO
  \satex{mat-gamma2-r}
\end{align*}
see~\cite{lafont2003towards}. This is the theory for \emph{bicommutative
  $R$-linear bialgebras}. For instance the diagram on the left corresponds to
the linear transformation $\R^3\to\R^2$ associated to the matrix depicted on
the right:
\[
  \satex{mat-ex}
  \qquad\qquad\qquad
  \pa{
    \begin{matrix}
      4&3&8\\
      7&5&2
    \end{matrix}
  }
\]
Note that we recover the laws for Hopf algebras (\secr{hopf}) by setting
\[
  \satex{sigma}
  =
  \satex{none}
  \pbox.
\]

\subsection{Variants}
Other variants have been studied in the literature. We can mention presentations
of invertible, orthogonal, special orthogonal, unitary, and special unitary
matrices~\cite{lafont2003towards} and stochastic
matrices~\cite{fritz2009presentation}.

\section{Coefficients}
\label{sec:coefs}

\subsection{The free monoidal category on a category}
\index{free!monoidal category}
Given a category~$C$, the \emph{free monoidal category}~$\freemoncat{C}$ it
generates is the monoidal category whose monoid of objects is
$\freemoncat{C_0}$, the free monoid over the objects of~$C$, the monoid of
morphisms is $\freemoncat{C_1}$ the free monoid over the morphisms of~$C$, with
\[
  f_1\ldots f_n
  :
  a_1\ldots a_n
  \to
  b_1\ldots b_n
\]
whenever $f_i:a_i\to b_i$ for $1\leq i\leq n$, composition is given pointwise,
\ie
\[
  (g_1\ldots g_n)\circ(f_1\ldots f_n)
  =
  (g_1\circ f_1)\ldots(g_n\circ f_n)
\]
with identities $\unit{}\ldots\unit{}$ and tensor product is given by
composition in the free monoid, \ie concatenation.

There is an obvious functor $C\to\freemoncat{C}$, which is such that for every
functor~$C\to D$, where~$D$ a monoidal category, there is a unique strict
monoidal functor $\freemoncat{C}\to D$ making the following diagram commute:
\[
  \xymatrix@C=3ex@R=3ex{
    C\ar[d]\ar[r]&D\pbox.\\
    \freemoncat{C}\ar@{.>}[ur]
  }
\]

Given a presentation of a category~$C$ by a 2-polygraph~$P$, the monoidal
category~$\freemoncat{C}$ admits a presentation by the 3-polygraph~$Q$ with
\begin{align*}
  Q_0&=\set{\star}
  &
  Q_1&=P_0
  &
  Q_2&=P_1
  &
  Q_3&=P_2
\end{align*}
sometimes referred to as the \emph{suspension} of the $2$-polygraph~$P$.

\subsection{Monoid actions}
\label{sec:act-mon}
\index{monoid!action}
As a particular case, when~$M$ is a monoid considered as a category with one object, an algebra
for the theory~$\freemoncat{M}$ in a monoidal category~$C$, consists of an
\emph{action} of~$M$, \ie an object~$x$ of the monoidal category~$C$ together with a
morphism of monoids $M\to\Hom Cxx$.

For instance, given a monoid~$(M,\times,1)$, we can consider its standard
presentation, see \secr{2-std-pres}:
\[
  \Pres{\star}{a}{\alpha_{u,v}:uv\To(u\times v)}
\]
the resulting presentation of the free monoidal category generated by~$M$ has a
presentation with a generator
\[
  \satex{coef-mon}
\]
for every element $u\in M$, with relations
\[
  \satex{coef-mon-l}
  \TO
  \satex{coef-mon-r}
  \,\pbox.
\]
More generally, when a monoid~$M$ admits a presentation by a 2-polygraph~$P$,
the free monoidal category generated by~$M$ admits a presentation with
$2$-generators
\[
  \satex{coef}
\]
indexed by $a\in P_1$, together with a relation
\[
  \satex{act-mon-l}
  \TO
  \satex{act-mon-r}
\]
for every relation $a_1\ldots a_i\To b_1\ldots b_j$ in~$P_2$.

\subsection{Thee free linear category}
\label{sec:R-mat}
Suppose given a monoid~$M$ and write~$R$ for the semiring it freely
generates. As explained above, the monoid $M$ freely generates a monoidal
category, and we write here $\freemoncat{M}$ for the free symmetric monoidal
category generated by this monoidal category. We have seen in \secr{bialgebra}
that the PROP $\Mat\N$ can be seen as a composite PROP
$\Fun^\op\otimes_\dlaw\Fun$. More generally, the PROP~$\Mat R$ can be seen as
a composite PROP:
\[
  \Mat R
  =
  \Fun^\op\otimes_\Bij\freemoncat{M}\otimes_\Bij\Fun
  \pbox.
\]
The tensor product above is an iterated distributive law, in the sense of
\secr{idlaw}: it is induced by three distributive laws
\[
  \begin{array}{r@{\ :\ }r@{\ \to\ }l}
    \dlaw_1&\freemoncat{M}\otimes_\Bij\Fun^\op&\Fun^\op\otimes_\Bij\freemoncat{M}
    \\
    \dlaw_2&\Fun\otimes_\Bij\Fun^\op&\Fun^\op\otimes_\Bij\Fun
    \\
    \dlaw_3&\Fun\otimes_\Bij\freemoncat{M}&\freemoncat{M}\otimes_\Bij\Fun
  \end{array}
\]
such that the diagram
\[
  \xymatrix@!C=7ex@R=3ex{
    &\freemoncat{M}\otimes_\Bij\Fun\otimes_\Bij\Fun^\op\ar[rr]^*+{\scriptstyle\freemoncat{M}\otimes_\Bij\dlaw_2}&&\freemoncat{M}\otimes_\Bij\Fun^\op\otimes_\Bij\Fun\ar[dr]^{\dlaw_1\otimes_\Bij\Fun}&\\
    \Fun\otimes_\Bij\freemoncat{M}\otimes_\Bij\Fun^\op\ar[ur]^{\dlaw_3\otimes_\Bij\Fun^\op}\ar[dr]_{\Fun\otimes_\Bij\dlaw_1}&&&&\Fun^\op\otimes_\Bij\freemoncat{M}\otimes_\Bij\Fun\\
    &\Fun\otimes_\Bij\Fun^\op\otimes_\Bij\freemoncat{M}\ar[rr]_{\dlaw_2\otimes_\Bij\freemoncat{M}}&&\Fun^\op\otimes_\Bij\Fun\otimes_\Bij\freemoncat{M}\ar[ur]_{\Fun^\op\otimes_\Bij\dlaw_3}
  }
\]
commutes.
More generally, when the semiring~$R$ is not freely generated by a monoid, the
PROP $\Mat R$ admits a description of the form
\[
  \Mat R
  =
  (\Fun^\op\otimes_\Bij\freemoncat{M}\otimes_\Bij\Fun)/{\sim}
\]
meaning that it can be obtained from the above PROP by further quotienting by
the congruence~$\sim$ generated by the two relations
\begin{align*}
  \satex{mat-add-l}
  &\TO
  \satex{mat-add-r}
  &
  \satex{mat-zero-l}
  &\TO
  \satex{mat-zero-r}  
\end{align*}
see~\cite{bonchi2017interacting} for details.

\section{Frobenius Algebras}
\index{algebra!Frobenius}
\index{Frobenius algebra}
\index{polygraph!of Frobenius algebras}
In this section, we study Frobenius algebras, which are ``dual'' to bialgebras
in a sense that will be made precise. The axioms were first discovered by
Lawvere~\cite{lawvere1969ordinal} (and rediscovered
in~\cite{carboni1987cartesian}) and a nice introduction to the subject can be
found in Kock's book~\cite{kock2004frobenius}.

\subsection{Presentation}
\label{sec:frobenius-algebra}
The theory of \emph{Frobenius algebras} is the PRO presented by four generators
\begin{align*}
  \mu&:2\to 1
  &
  \eta&:0\to 1
  &
  \delta&:1\to 2
  &
  \varepsilon&:1\to 0
  \\
  \intertext{respectively pictured as}
  &\satex{mu}
  &&
  \satex{eta}
  &&
  \satex{delta}
  &&
  \satex{eps}
\end{align*}
such that $(\mu,\eta)$ is a monoid, $(\delta,\varepsilon)$ is a comonoid, and the
two following relations are satisfied
\begin{align}
  \label{eq:frob}
  \satex{frob-l}
  &\TO
  \satex{frob-c}
  &
  \satex{frob-r}
  &\TO
  \satex{frob-c}
\end{align}

The theory of \emph{commutative} Frobenius algebras is obtained by taking the
free symmetric algebra as described in \secr{free-sym} (\ie adding a generator
$\gamma:2\to 2$ pictured as usual, see \eqref{eq:gamma},
together with the axioms for symmetries~\eqref{eq:sym-rel}, the compatibility
relations between symmetry and the monoid structure~\eqref{eq:mon-sym-compat},
as well as the comonoid structure), and adding the commutativity relations
\begin{align*}
  \satex{gamma-mu}
  &\TO
  \satex{muu}
  &
  \satex{gamma-delta}
  &\TO
  \satex{deltaa}
  \,\pbox.
\end{align*}
The theory of \emph{special} or \emph{separable} Frobenius algebras (commutative
or not) can be obtained by further adding the relation
\begin{equation}
  \label{eq:frob-sep}
  \satex{rel-l}
  \TO
  \satex{rel-r}
\end{equation}
and the theory of \emph{extraspecial} Frobenius algebras by further adding to
the theory of special Frobenius algebras the relation
\[
  \satex{eta-eps-l}
  \TO
  \satex{empty}
  \,\pbox.
\]
Other variants can be found in~\secr{graph-cospan}.

Note that these axiomatizations are not claimed to be minimal. For instance,
associativity in presence of other axioms implies coassociativity:
\begin{align*}
  \satex{frob-ac1}
  \enskip&=\enskip
  \satex{frob-ac2}
  \enskip=\enskip
  \satex{frob-ac3}
  \enskip=\enskip
  \satex{frob-ac4}
  \\
  \enskip&=\enskip
  \satex{frob-ac5}
  \enskip=\enskip
  \satex{frob-ac6}  
\end{align*}
and the situation is similar for unitality and commutativity.

It is conjectured that the rewriting system for Frobenius algebras and its
variants is terminating and can be completed into a finite convergent rewriting
system~\cite{forest2018coherence}. Note that there are indexed critical pairs
such as
\[
  \satex{frob-icp-l}
  \OT
  \satex{frob-icp-c}
  \TO
  \satex{frob-icp-r}
\]
but we expect to be able to handle those using the techniques of
\cref{sec:3pol-indexed}.

Even in the absence of a notion of normal form, we can introduce the following
notion of canonical form for morphisms. We define a family of morphisms
$\phi_{m,g,n}$, indexed by $m,g,n\in\N$, defined by
\begin{equation}
  \label{eq:frob-nb}
  \begin{split}
  \phi_{m,g,n}
  &=
  \mu_m
  \comp{}
  (\delta\comp{}\mu)
  \comp{}
  (\delta\comp{}\mu)
  \comp{}
  \ldots
  \comp{}
  (\delta\comp{}\mu)
  \comp{}
  \delta_n
  \\
  &=
  \satex{frob-mgn-h}
  \end{split}
\end{equation}
(the diagram is drawn horizontally to save space)
where $\mu_m:m\to 1$ is a right comb of multiplications, defined inductively by
\begin{align*}
  \satex{comb0}&=\satex{eta}
  &
  \satex{combm1}&=\satex{comb-ind}
\end{align*}
and dually $\delta_n:1\to n$ is a right comb of comultiplications, and there are
$g$ occurrences of $\delta\comp{}\mu$ in the middle: $m$ and $n$ are
respectively the arity and coarity of the morphism and $g$ is called its
\emph{genus} for reasons explained below. It can then be shown that, in the
bicommutative case, any morphism rewrites to a tensor product of such morphisms
(up to composing with a symmetry).  In the case of special Frobenius algebras,
the normal forms are tensor products of morphism of the form $\phi_{m,0,n}$, and
in the case of extraspecial Frobenius algebras the normal forms are the same
excepting that~$\phi_{0,0,0}$ is not allowed to occur.

\subsection{Frobenius algebras}
A Frobenius algebra is typically considered in the category $\Vect\kk$: it
consists of a vector space~$A$ which is both an algebra (\ie equipped with
associative and unital morphisms $\mu:A\otimes A\to A$ and $\eta:\kk\to A$), a
coalgebra (\ie equipped with coassociative and counital morphisms
$\delta:A\to A\otimes A$ and $\varepsilon:A\to\kk$) satisfying the compatibility
axioms~\eqref{eq:frob}.

There are many alternative characterizations of those~\cite{kock2004frobenius,street2004frobenius}. For instance, it can also be defined as an
algebra~$(A,\mu,\eta)$, equipped with a non-degenerate bilinear form
$\sigma:A\otimes A\to A$ which is associative, in the sense that
$\sigma\circ(\mu\otimes A)=\sigma\circ(A\otimes\mu)$, see \secr{frob-pairing}
for details; or as an algebra~$(A,\mu,\eta)$ equipped with a linear form
$\varepsilon:A\to\kk$ such that $\varepsilon(ab)=0$ for every $a\in A$ implies
$b=0$. Namely, one can transform the first into the second definition, and
vice-versa, by defining $\varepsilon(a)=\sigma(1,a)$ and
$\sigma=\varepsilon\circ\mu$.

\subsection{Frobenius as cospans}
\label{sec:frob-cospan}
The category presented by the theory of special Frobenius algebras is the
category~$\qCoSpan(\Simplaug)$ of isomorphism classes of cospans of
non-decreasing functions. Namely, it is composed of the
presentation of~$\Simplaug^\op$ (comonoids), the presentation of~$\Simplaug$
(monoids), and compatibility axioms can be obtained as pushouts in~$\Simplaug$:
\begin{align*}
  \xymatrix@C=3ex@R=3ex{
    &\ar[dl]_{\mu 1}3\ar[dr]^{1\mu}&\\
    2\ar@{.>}[dr]_\mu&&2\ar@{.>}[dl]^\mu\\
    &1&
  }
  &&
  \xymatrix@C=3ex@R=3ex{
    &\ar[dl]_{\mu}2\ar[dr]^{\mu}&\\
    1\ar@{.>}[dr]_{\unit{}}&&1\ar@{.>}[dl]^{\unit{}}\\
    &1&
  }
  &&
  \xymatrix@C=3ex@R=3ex{
    &\ar[dl]_{1\mu}3\ar[dr]^{\mu 1}&\\
    2\ar@{.>}[dr]_\mu&&2\ar@{.>}[dl]^\mu\pbox.\\
    &1&
  }  
\end{align*}
The diagram in the middle corresponds to the separability
axiom~\eqref{eq:frob-sep}, whereas the two other to the Frobenius
axioms~\eqref{eq:frob}. The presentation can thus be deduced from the one of
$\Simplaug$, by using the fact that
$\qCospan(\Simplaug)=\Simplaug\otimes_\dlaw\Simplaug^\op$ where the distributive law
$\dlaw:\Simplaug^\op\otimes\Simplaug\to\Simplaug\otimes\Simplaug^\op$ is given by pushout.
Similarly, the theory of special commutative Frobenius algebras presents the
category~$\qCoSpan(\Fun)=\Fun\otimes_\dlaw\Fun^\op$ of isomorphism classes of cospans of functions. In
this sense it is dual to the theory~$\qSpan(\Fun)$ of commutative bialgebras described in \cref{sec:bialgebra} (by
analogy, one could be tempted to consider $\qSpan(\Simplaug)$ as a theory for
non-commutative bialgebras, but this is not well-defined because $\Simplaug$ does not
have all pullbacks).

\subsection{Corelations}
\nomenclature[Corel]{$\Corel$}{category of corelations}
\label{sec:corel}
\index{corelation}
The theory of extraspecial commutative Frobenius algebras corresponds to the
variant of the category~$\qCoSpan(\Fun)$ whose morphisms are jointly surjective
cospans (considered up to isomorphism), \ie cospans of the form
$\xymatrix@C=3ex@R=3ex{X\ar[r]|-f&Y&\ar[l]|-gZ}$ such that for every~$y\in Y$
there exists $x\in X$ such that $f(x)=y$ or there exists $z\in Z$ such
that~$g(z)=y$, see~\cite{coya2016corelations}. Dually to the case of relations,
see \secr{pres-rel}, these morphisms can alternatively be described as
follows. Given two sets~$X$ and~$Y$, a \emph{corelation}~$f:X\to Y$ is a
partition of~$X\sqcup Y$. Given two corelations $f:X\to Y$ and~$g:Y\to Z$, their
composite $g\circ f:X\to Z$ is defined as the restriction to $X\sqcup Z$ of the
finest partition of~$X\sqcup Y\sqcup Z$ which is coarser than $f$ (\resp $g$)
when restricted to $X\sqcup Y$ (\resp $Y\sqcup Z$). The identity corelation
$\unit{X}:X\to X$ is the diagonal corelation. We write $\Corel$ for the category
of natural numbers and corelations: this category is isomorphic to the category
of cospans described above.
For instance, we have the following composition of corelations
\[
  \fig{corel-comp}
\]
which corresponds to the following composition in the theory of extraspecial
commutative Frobenius algebras
\[
  \satex[angle=90,origin=c]{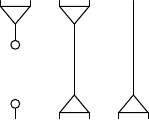}
  \qquad\quad
  \satex[angle=90,origin=c]{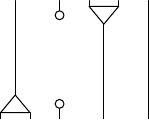}
  \qquad\quad
  \satex[angle=90,origin=c]{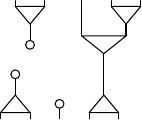}
\]

Any morphism $f:m\to n$ in $\qCoSpan(\Fun)$, which consists of a cospan
\[
  \xymatrix{m\ar[r]|{f_1}&p&\ar[l]|{f_2}n}
\]
up to isomorphism, can be uniquely be written as
\[
  f=f'\comp0\satex{eta-eps-small}\comp0\ldots\comp0\satex{eta-eps-small}
\]
where $f':m\to n$ is a corelation and the number of instances of
$\satex{eta-eps-small}$ indicates the number of elements of $\intset{p}$
which are neither in the image of~$f_1$ nor in the image of~$f_2$, \ie measures
the deficiency of surjectivity of~$f$.

\subsection{2-cobordisms}
\label{sec:ncob}
\label{sec:2cob}
\index{cobordism}
\index{topological quantum field theory}
We briefly recall the well-known description of the category
presented by the theory of commutative Frobenius algebras in geometrical
terms. A clear and detailed account of the situation can be found
in~\cite{kock2004frobenius}. Fix a natural number $n\in\N$. Given a smooth
oriented manifold with boundary~$\Sigma$, its boundary decomposes as
$\partial\Sigma=\partial^-\Sigma\sqcup\partial^+\Sigma$, where the two
components are determined according to the orientation, and we say that $\Sigma$
is an $n$-\emph{cobordism} from $\partial^-\Sigma$ to $\partial^+\Sigma$. One
can build a category $\Cob n$ whose objects are oriented smooth
$(n{-}1)$-manifolds, morphisms are $n$-cobordisms, considered up to
diffeomorphism, and composition is given by gluing (\ie taking pushouts) along
common boundaries. This category is symmetric monoidal, with tensor product
given on objects and morphisms by disjoint union, and people usually consider
\emph{topological quantum field theories}, which are symmetric monoidal functors
$\Cob n\to\Vect\kk$.

Here, we will be interested in $\Cob2$: an object is a
disjoint union of circles and a morphism consists of ``trousers'' between those,
such as on the left below
\begin{align}
  \label{eq:cob-ex}
  \fig{cob-ex}
  &&
  \satex{cob-ex}
  \,\pbox.
\end{align}
The category $\Cob 2$ admits the theory of commutative Frobenius algebras as
presentation. Namely, the generators are interpreted as
\begin{align*}
  \mu&=\fig{cob-mu}
  &
  \eta&=\fig{cob-eta}
  &
  \delta&=\fig{cob-delta}
  &
  \varepsilon&=\fig{cob-eps}
  &
  \gamma&=\fig{cob-gamma}
\end{align*}
so that, for instance, the morphism on the left of \eqref{eq:cob-ex} corresponds to the diagram on the
right. Verifying that the relations are satisfied in $\Cob2$ is
direct. Conversely, in order to show that they are sufficient, one can use the classical
result that the connected morphisms (\ie diffeomorphism classes of connected
compact oriented surfaces with boundaries) are characterized by their number of
inputs, outputs, and genus (roughly, the number of holes): these are thus in
bijection with normal forms~\eqref{eq:frob-nb} since the interpretation of
$\phi_{m,g,n}$ is a cobordism of genus~$g$ with~$m$ inputs and~$n$ outputs.

Non-commutative Frobenius algebras admit a similar description by
2\nbd-dimen\-sional thick tangles~\cite{lauda2005frobenius}.

\section{Linear Relations}
\nomenclature[LinRelk]{$\LinRel\kk$}{category of linear relations}
\label{sec:lin-rel}
\index{relation!linear}
Suppose fixed a field~$\kk$. The category $\LinRel\kk$ has vector spaces over $\kk$
as objects, a morphism $f:X\to Y$ is a subspace of the vector space $X\oplus Y$,
composition of $f:X\to Y$ and $g:Y\to Z$ is given by relational composition
\[
  g\circ f
  =
  \setof{x\oplus z\in X\oplus Z}{\exists y\in Y, x\oplus y\in f \land y\oplus z\in g}
\]
and identity on $X$ is the diagonal in $X\oplus X$. A morphism $f:X\to Y$ in
this category is a relation between the underlying sets~$X$ and~$Y$, which is
closed under addition and multiplication by scalars and thus called a
\emph{linear relation}. This category is typically used to provide semantics to
networks such as those found in electric circuits~\cite{baez2015compositional}
or control~\cite{baez2015categories, bonchi2017calculus}.

It is shown in \cite{baez2015categories, bonchi2017interacting} that this
category admits a presentation with generators
\begin{align*}
  \satex{mu}
  &&
  \satex{eta}
  &&
  \satex{delta}
  &&
  \satex{eps}
  &&
  \satex{mu3}
  &&
  \satex{eta3}
  &&
  \satex{delta3}
  &&
  \satex{eps3}
  &&
  \satex{gamma}
  &&
  \satex{coef2}    
\end{align*}
for $a\in\kk$ such that
\begin{itemize}
\item
  $(\satex{mu-small},\satex{eta-small},\satex{delta3-small},\satex{eps3-small},\satex{gamma-small},\satex{coef2-small})$
  is a bicommutative linear bialgebra (\cref{sec:linear-bialg}),
\item
  $(\satex{mu3-small},\satex{eta3-small},\satex{delta-small},\satex{eps-small},\satex{gamma-small},\satex{coef2-small})$
  is a bicommutative linear bialgebra (\cref{sec:linear-bialg}),
\item
  $(\satex{mu-small},\satex{eta-small},\satex{delta-small},\satex{eps-small},\satex{gamma-small})$
  is a bicommutative extraspecial Frobenius algebra (\cref{sec:frobenius-algebra}),
\item
  $(\satex{mu3-small},\satex{eta3-small},\satex{delta3-small},\satex{eps3-small},\satex{gamma-small})$
  is a bicommutative extraspecial Frobenius algebra (\cref{sec:frobenius-algebra}),
\item the following compatibility relations hold:
  \begin{align}
    \label{eq:linrel-compat}
    \satex{ih-compat1-l}
    &=
    \satex{ih-compat1-r}
    &
    \satex{ih-compat2-l}
    &=
    \satex{ih-compat2-r}    
  \end{align}
\end{itemize}
where $\satex{sigma2-small}$ is a notation for $\satex{coef2-1-small}$.
The generators should be interpreted as the following linear relations:
\begin{align*}
  \satex{mu-ih}
  &&
  \satex{eta-ih}
  &&
  \satex{delta-ih}
  &&
  \satex{eps-ih}
  &&
  \satex{gamma-ih}
  &&
  \satex{coef2-ih}
  \\
  \satex{mu3-ih}
  &&
  \satex{eta3-ih}
  &&
  \satex{delta3-ih}
  &&
  \satex{eps3-ih}
\end{align*}
meaning that the interpretation of $\satex{mu-vsmall}$ is the linear relation
\[
  \setof{(x\otimes y)\oplus(x+y)}{x,y\in\kk}
  \subseteq
  (\kk\otimes\kk)\oplus\kk
\]
and so on. For instance, the two compatibility relations \eqref{eq:linrel-compat} can be read as the fact
that the space of pairs $(x,y)$ such that $x+y=0$ coincides with the space of
pairs $(x,-x)$, with $x$ and $y$ ranging over $\kk$.
This theory is sometimes called the theory of \emph{interacting Hopf algebras}
over~$\kk$ and has applications to control theory.

In~\cite{bonchi2017interacting}, it is shown that this category can be obtained
as the following pushout in~$\MonCat$:
\[
  \xymatrix@C=3ex@R=3ex{
    \Vect\kk\sqcup\Vect\kk^\op\ar[d]\ar[r]&\qSpan(\Vect\kk)\ar@{.>}[d]\\
    \qCospan(\Vect\kk)\ar@{.>}[r]&\LinRel\kk
  }
\]
from which the above presentation can be obtained. Another possible description
is mentioned in \secr{lin-rel-dual}.

\section{Interchange}
\subsection{Interchange algebra}
\index{interchange!algebra}
\index{algebra!interchange}
The theory of \emph{interchange algebras}~\cite{lafont1997interaction,loday2006generalized}
models algebraic structures with two multiplications satisfying the
exchange law (which is sometimes also called the \emph{interchange law}) holds. It can
be presented by the 3-polygraph~$P$ with generators
\begin{align*}
  P_0&=\set\star
  &
  P_1&=\set{a}
  &
  P_2&=\set{\gamma:aa\To aa,\mu:aa\To a,\nu:aa\To a}
\end{align*}
and relations expressing that $\gamma$ is a symmetry (\cref{sec:free-sym})
together with the relation
\[
\satex{inter_l}
\TO
\satex{inter_r}
\]
where the white (\resp gray) triangle corresponds to $\mu$ (\resp $\nu$). An
interchange algebra is \emph{associative} when both $\mu$ and $\nu$ are. It is
\emph{unital} when equipped with two morphisms $\unit{}\To a$ acting as a left
and right unit for $\mu$ and $\nu$, respectively. For instance, an associative
unital interchange algebra in the cartesian category~$\Set$ is precisely a
$2$-category with only one $0$-cell and one $1$-cell. By the Eckmann-Hilton
argument~\cite{eckmann1962group}, this is the same as a commutative
monoid.

\subsection{Iterated monoidal categories}
\label{sec:it-mon-cat}
\index{iterated monoidal category}
\index{monoidal category!iterated}
A variant of interchange algebras with $n$ distinct monoid structures instead of
two is the following one.
An \emph{$n$-fold monoid} is a set equipped with $n$-distinct products which are
suitably compatible. More precisely, the corresponding theory can be presented
by the 3-polygraph~$P$ with generators $P_0=\set\star$, $P_1=\set{a}$ and
\[
  P_2=\set{\gamma:aa\To aa,\eta:\unit{}\To a,\mu_0:aa\To a,\ldots,\mu_{n-1}:aa\To a}
\]
such that for every $0\leq i<j<n$, $\eta$, $\mu_i$, and $\mu_j$ satisfy the
axioms of associative unital interchange algebras. Again, by the Eckmann-Hilton
argument, this theory is equivalent to the one of commutative monoids. However,
the theory of \emph{$n$-fold pseudomonoids} (where equalities are replaced by
coherent isomorphisms) is more interesting: its algebras in the cartesian
$2$-category~$\Cat$ are \emph{$n$-fold monoidal categories}. This structure is
the one one obtains by considering monoids in the category of monoids in the
category of monoids in the category~$\Cat$ with the monoidal structure
induced by cartesian product, see~\cite{balteanu2003iterated}.

\subsection{Interchange bialgebra}
\index{interchange!bialgebra}
An \emph{interchange bialgebra}~\cite{loday2006generalized} is a PROP which is
both an interchange algebra and an interchange coalgebra, with the following
compatibility relations
\begin{align*}
  \satex{ib_ww_l}
  &\TO
  \satex{ib_ww_r}
  &
  \satex{ib_bb_l}
  &\TO
  \satex{ib_bb_r}
  &
  \satex{ib_wb_l}
  &\TO
  \satex{ib_wb_r}
  &
  \satex{ib_bw_l}
  &\TO
  \satex{ib_bw_r}
  \,\pbox.
\end{align*}
The most studied variant of this notion is the one of interaction nets, which are
detailed in \secr{interaction-net}.

\section{Idempotent Objects}
\label{sec:3-Thompson-F}
\index{idempotent}
\index{Thompson!group}
\index{group!Thompson}
\index{Thompson!category}
\index{category!Thompson}
\newcommand{\TF}{\category{T}}
%
An \emph{idempotent object} in a monoidal category is an object $x$ equipped
with an isomorphism $x\otimes x\to x$. The theory for idempotent objects is the
PRO $\TF$, called the \emph{Thompson category}, presented by the 3-polygraph~$P$
with $P_0=\set\star$, $P_1=\set{a}$, with 2-generators $\mu:aa\To a$ and
$\delta:a\To aa$, depicted as
\begin{align*}
  \satex{mu}
  &&
  \satex{delta}  
\end{align*}
and subject to the relations
\begin{align*}
  \satex{mu-delta}
  &\TO
  \satex{sym-r}
  &
  \satex{bialg-qual-l}
  &\TO
  \satex{bialg-qual-r}
  \,\pbox.
\end{align*}
A morphism is a binary \emph{tree} (\resp a \emph{cotree}) when it has $a$ as
target (\resp source) and is a composite of generators $\mu$ (\resp $\delta$)
only. It is observed in~\cite{fiore2010abstract} that the monoid of
automorphisms $\TF(a,a)$ is isomorphic to the Thomson group~$F$, already
presented in \secr{Thompson-F}. Namely, we can recover the generators of the
usual presentation as defined by induction on~$i$ by
\[
  x_0
  =
  \satex{F-x0}
  \qqtand
  x_{i+1}
  =
  \satex{F-xi}
  \pbox.
\]
More generally, any morphism $a\to a$ of~$\TF$ decomposes as a cotree followed
by a tree respectively encoding the dyadic partitions of the input and of the
output of the corresponding morphism $I\to I$ in the Thompson group by
specifying when the interval should be split in two halves. For instance, the
morphism on the left corresponds to the function on the right:
\[
  \satex{F-ex}
  \qquad\qquad\qquad\qquad
  \begin{array}{rc}
    \satex[angle=-90,,scale=.27,trim=0 0 -32mm 0]{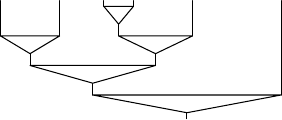}&\fig{F-ex}\\
    &\hspace{5mm}
    \satex[angle=180,scale=.5]{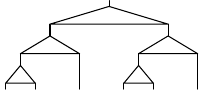}
  \end{array}
\]

By post-composition, the elements of~$F$ act on trees which are ``large enough''
(the action is partially defined). For instance, we have the following action
of~$x_0$:
\[
  \satex{F-x0-act-l}
  \overset\ast\TO
  \satex{F-x0-act-r}
\]
From this point of view, the generator~$x_0$ can be pictured as
\[
  \satex{mon-assoc-l}
  \overset{x_0}
  \to
  \satex{mon-assoc-r}
\]
since it will ``replace'' a prefix of a tree as on the left with a prefix as on
the right. Such a transformation is sometimes called an \emph{associative law}
and the group~$F$ can be described as the group of those with expected
composition. For instance, Mac Lane's pentagon is
\[
  \xymatrix@C=12ex@R=3ex{
    \satex{mon-cp1-0}\ar[dd]_{x_0}\ar[r]^{x_0x_0\finv{x_1}\finv{x_0}}&\satex{mon-cp1-1}\ar[dr]^{x_0}\\
    &&\satex{mon-cp1-2}\ar[d]^{x_1}\\
    \satex{mon-cp1-4}\ar[rr]_{x_0}&&\satex{mon-cp1-3}
  }
\]

Similarly, the Thompson group~$V$ can be recovered as the group of automorphisms
on~$a$ in the free symmetric monoidal category on~$\TF$.
This entails that~$V$ can be described as the group of automorphisms of the free
Cantor algebra on a singleton~\cite{brown1987finiteness}: we recall that a
\emph{Cantor algebra} is a set~$A$ equipped with a bijection
$\alpha:A\to A\times A$.

\section{Dualities}
\subsection{Adjunctions}
\label{sec:adjunction}
\index{adjunction}
An \emph{adjunction} consists of two functors~$f:C\to D$ and $g:D\to C$ together
with natural transformations $\eta:\unit{C}\To fg$ and
$\varepsilon:gf\To\unit{D}$ such that $(\eta f)\comp{}(f\varepsilon)=\unit{f}$
and $(g\eta)\comp{}(\varepsilon g)=\unit{g}$. In such a situation~$f$ is called
a \emph{left adjoint} to~$g$, and~$g$ a \emph{right adjoint} to~$f$.

The theory of \emph{adjunctions} is the 2-category~$\Adj$ presented by the
3\nbd-poly\-graph~$P$ with generators
\begin{align*}
  P_0&=\set{c,d}
  &
  P_1&=\set{f:c\to d,g:d\to c}
\end{align*}
and
\begin{align*}
  P_2&=\set{\eta:\unit{c}\To fg,\varepsilon:gf\To\unit{d}}  
\end{align*}
pictured as
\begin{align*}
  \eta&=\satex{adj-eta-lbl}
  &
  \varepsilon&=\satex{adj-eps-lbl}  
\end{align*}
and respectively called \emph{unit} and \emph{counit}, together with the two
relations
\begin{align}
  \label{eq:zigzag}
  \satex{adj1-l}
  &\TO
  \satex{adj1-r}
  &
  \satex{adj2-l}
  &\TO
  \satex{adj2-r}  
\end{align}
often called \emph{zigzag} or \emph{triangle identities}. This presentation is
convergent, the normal forms being the horizontal composites of $\eta$,
$\varepsilon$, and identities.
A 2-functor $F:\Adj\to\Cat$ corresponds precisely to an adjunction in the usual
sense.

The 2-category~$\Adj$ is studied in~\cite{schanuel1986free}. In particular, if
we consider the 2-category $\Adj(c,c)$, which is the full sub-2-category
of~$\Adj$ on the 0-cell~$c$, we have an isomorphism of 2-categories (or of
monoidal categories)
\[
  \Adj(c,c)
  \isoto
  \Simplaug
\]
Namely, the monoid of 1-cells is freely generated by~$fg:c\to c$ (thus
isomorphic to~$\N$), and one can define a structure of monoid on~$fg$ whose
multiplication and unit are respectively $f\varepsilon g$ and $\eta$, from which
the isomorphism can easily be deduced. For instance, the non-decreasing function
pictured on the left corresponds to the 2-cell on the right in the theory of
monoids
\begin{align*}
  \vcenter{
    \xymatrix@C=3ex@R=3ex{
      0\ar@{-}[d]&1\ar@{-}[dl]&2\ar@{-}[dll]&3\ar@{-}[dl]\\
      0&1&2
    }
  }
  &&
  \satex{adj-simpl-ex-l}  
\end{align*}
and to the following 2-cell in $\Adj(c,c)$:
\[
  \satex{adj-simpl-ex-r}
  \pbox.
\]
Similarly, we have $\Adj(d,d)\isoto\Simplaug^\op$, and $\Adj(c,d)$ and $\Adj(d,c)$
are the subcategories of~$\Simplaug$ whose objects are non-zero natural numbers and
morphisms are the last-element (\resp first-element) preserving
functions. Note that~$\Simplaug^\op$ is isomorphic to the subcategory of~$\Simplaug$
whose objects are non-zero natural numbers and morphisms are preserving both
first and last element.

\subsection{Duality}
\index{duality}
\index{dual}
\index{self-dual}
\index{autonomous category}
\index{category!autonomous}
\index{rigid category}
\index{category!rigid}
\index{compact closed category}
\index{category!compact closed}
\index{trace}
\index{pivotal category}
\index{category!pivotal}
A \emph{duality} in a monoidal category~$C$ is an adjunction in~$C$, considered
as a 2-category. It consists of two objects~$x$ and~$x^*$ together with
morphisms
\begin{align*}
  \eta_x&:\unit{}\to xx^*
  &
  \eps_x&:x^*x\to\unit{}  
\end{align*}
satisfying the zigzag relations~\eqref{eq:zigzag}. In this case $x$ is called
a \emph{left dual} of~$x^*$, and $x^*$ a \emph{right dual} of~$x$. Two left
(\resp right) duals of a given object are necessarily isomorphic. An object~$x$
is \emph{self-dual} when it admits a right dual~$x^*$ which is isomorphic
to~$x$.

A monoidal category is \emph{right-autonomous} (\resp \emph{left-autonomous})
when every object~$x$ admits a right dual $x^*$ (\resp a left dual
$\null^*x$). It is \emph{strictly} so when duals are chosen so that
$(xy)^*=y^*x^*$, $\unit{}^*=\unit{}$,
$\eta_{\unit{}}=\unit{\unit{}}=\eps_{\unit{}}$,
$\eta_{xy}=x\eta_yx^*\circ\eta_x$, and $\eps_{xy}=\eps_y\circ y^*\eps_xy$;
without loss of generality, we consider that this is always the case in the
following. It is \emph{autonomous} (or \emph{rigid}) when it is both left- and
right-autonomous.
A \emph{compact closed category} is a symmetric monoidal category which is
autonomous.
For instance, the category $\Vect\kk$ of $\kk$-vector spaces is compact closed,
the dual $x^*$ of a vector space~$x$ being its linear
dual.
A coherence theorem for compact closed categories was shown by Kelly and
Laplaza~\cite{kelly1980coherence}. In particular, in a compact closed category,
we have isomorphisms $\null^*x\isoto x^*$, $x^{**}\isoto x$, which we will be
considered as equalities in the following.
Given a morphism $f:x\to x$ of a compact closed category, the morphism
\[
  \trace(f)=\eps_{x^*}\circ fx^*\circ\eta_x
  \qquad\qquad\qquad
  \satex{cc-trace-f}
\]
is called its \emph{trace}: there is a general axiomatization of trace in
symmetric monoidal categories~\cite{joyal1996traced}, which every compact closed
category canonically possesses. In particular, the morphism $\trace(\id_x)$ is
often called the \emph{dimension} of~$x$. The category is \emph{loop-free} when
every object~$x$ is $1$-dimensional, \ie $\trace(x)=\unit{\unit{}}$.

A right-autonomous category is \emph{pivotal} when equipped with a monoidal
natural isomorphism $x\isoto x^{**}$, considered as an equality in the following.
A nice survey of the flavors of categories with duals can be found
in~\cite{selinger2010survey}.
  
\subsection{The free compact category}
\index{free!compact category}
Suppose given a symmetric monoidal category~$C$ presented by a
3-polygraph~$P$. The free compact category on~$C$ admits a presentation by the
3-polygraph~$Q$ with generators
\begin{align*}
  Q_0&=P_0=\set\star
  \\
  Q_1&=\setof{a^n}{\text{$a\in P_1$ and $n\in\Z$}}
  \\
  Q_2&=\setof{f:u^0\To v^0}{f:u\To v\in P_2}\sqcup\\
  &\quad\setof{\eta_{a^n}:\unit{}\To a^na^{n+1},\varepsilon_{a^n}:a^{n+1}a^n\To\unit{}}{a^n\in Q_1}
\end{align*}
where, for $u=a_1\ldots a_k$, we write $u^n$ for $a_1^n\ldots a_k^n$. The
generators $\eta_{a^n}$ and $\varepsilon_{a^n}$ are respectively pictured
\[
  \satex{eta-an}
  \qqtand
  \satex{eps-an}
  \pbox.
\]
The relations are those in~$P_3$ together with the zigzag
relations~\eqref{eq:zigzag} satisfied by $\eta_{a^n}$ and
$\varepsilon_{a^n}$ for every $a\in P_1$ and integer~$n$.

Here, the object $a^0$ corresponds to $a$ in the original category and, for
$n\in\N$, $a^n$ (\resp $a^{-n}$) corresponds to $a^{**\cdots*}$ (\resp
$\null^{*\cdots**}a$), with $*$ applied $n$ times.
Thanks to this presentation, one can for instance deduce that the canonical
monoidal functor from a monoidal category into its free compact category is
faithful~\cite{mimram2014towards}.

\subsection{The free compact closed category}
Suppose given a symmetric monoidal category~$C$ presented by a
3-polygraph~$P$. The free compact closed category on~$C$ admits a presentation
by the 3-polygraph~$Q$ with generators
\begin{align*}
  Q_0&=P_0=\set\star
  \\
  Q_1&=\setof{a,a^*}{a\in P_1}
  \\
  Q_2&=\setof{f:u\To v}{f:u\To v\in P_2}\sqcup\setof{\eta_a:\unit{}\To aa^*,\varepsilon_a:a^*a\To\unit{}}{a\in P_1}
\end{align*}
and relations being those in~$Q_3$ together with the zigzag relations for each
$\eta_a$ and $\varepsilon_a$.

The theory for a pair of adjoint endofunctors can be obtained from the
category~$\Adj$ of \secr{adjunction} by identifying the two objects, and its
presentation can be obtained from the one of~$\Adj$ by identifying the
0-generators~$c$ and~$d$. The above results show that this theory is the free
compact closed category on the terminal category.

\subsection{The free ``self-dual'' compact closed category}
\label{sec:free-sd-compact}
As a variant of the situation described in previous section, given a symmetric
monoidal category~$C$ presented by a polygraph~$P$, we call the \emph{free
  self-dual compact closed category} on~$C$, the category presented by the
$3$-polygraph~$Q$ with generators
\begin{align*}
  Q_0&=P_0=\set\star
  \\
  Q_1&=P_1
  \\
  Q_2&=\setof{f:u\To v}{f:u\To v\in P_2}\sqcup\setof{\eta_a:\unit{}\To aa,\varepsilon_a:aa\To\unit{}}{a\in P_1}
\end{align*}
and relations being those in~$Q_3$ together with the zigzag relations for each
$\eta_a$ and $\varepsilon_a$, and the relations
\begin{align*}
  \satex{sd-eta-l}
  &=
  \satex{sd-eta-r}
  &
  \satex{sd-eps-l}
  &=
  \satex{sd-eps-r}
  \\
  \gamma_{a,a}\circ\eta_a
  &=
  \eta_a
  &
  \eps_a\circ\gamma_{a,a}
  &=
  \eps_a
\end{align*}
for every 1-generator $a\in P_1$. In the resulting category, only the generators
are self-dual, for instance the dual of the $1$-cell~$ab$ is~$ba$. The
axiomatization of self-dual categories is quite subtle,
see~\cite{selinger2010autonomous}, and it is not clear that this construction is
an instance of those. However, it is quite useful in the following, so we use it
without claiming a universal property.

\subsection{1-cobordisms}
The category of $n$-cobordisms was presented in \secr{ncob}. In this section we
study the case $n=1$ and give a presentation of the monoidal
category~$\Cob1$. Its objects are disjoint unions of oriented $0$\nbd-dimen\-sional
manifolds, \ie points together with an orientation $-$ or $+$, and it admits a
presentation by the polygraph~$P$ with $P_0=\set{\star}$, $P_1=\set{-,+}$,
2-generators being
\begin{align*}
  \satex{gamma-m}
  &&
  \satex{gamma-mp}
  &&
  \satex{gamma-pm}
  &&
  \satex{gamma-p}
  &&
  \satex{cap-mp}
  &&
  \satex{cup-pm}  
\end{align*}
such that the four first generators equip the monoidal category with a symmetric
structure, see~\secr{free-sym}, and the two last generators satisfy the axioms
for dualities. In other words, $\Cob1$ is the free compact closed symmetric monoidal category on
one object.

\subsection{The Temperley-Lieb category}
\label{sec:temperley-lieb-cat}
\index{Temperley-Lieb!category}
\index{category!Temperley-Lieb}
\newcommand{\TL}{\category{TL}}
The \emph{Temperley-Lieb category} $\TL$, introduced and studied
in~\cite{abramsky2009temperley}, is the PRO generated by
\begin{align*}
  \satex{cap}
  &&
  \satex{cup}  
\end{align*}
subject to the relations
\begin{align*}
  \satex{adj1-l}
  &\TO
  \satex{adj1-r}
  &
  \satex{adj2-l}
  &\TO
  \satex{adj2-r}
  &
  \satex{tl-r}
  &\TO
  \satex{tl-l}
  \,\pbox.
\end{align*}
It is thus the free monoidal category containing a self-dual object, satisfying the
last relation above.
Note that this presentation is not terminating because of the loop
\[
  \satex{tl-loop1}
  \qTO
  \satex{tl-loop2}
  \qTO
  \satex{tl-loop1}
\]
it can however be shown to be quasi-terminating~\cite{alleaume2018rewriting}, in the sense of \cref{sec:quasi-termination}.
Given $n\in\N$, the monoid of endomorphisms $\TL(n,n)$ is generated by
\[
  d=\satex{tl-d}
  \qqtand
  u_i=\satex{tl-ui}
\]
with $0\leq i\leq n-2$, where $u_i$ is composed of the identity on $i$ on the
left and on $n-i-2$ on the right, subject to relations of \secr{temperley-lieb},
making it the $n$-th Temperley-Lieb monoid, see \secr{temperley-lieb}. For
instance, with $n=4$, we have the relations
\begin{align*}
  \satex{tl-ex1-l}&=\satex{tl-ex1-r}&
  \satex{tl-ex2-l}&=\satex{tl-ex2-r}&
  \satex{tl-ex3-l}&=\satex{tl-ex3-r}\,\pbox.
  \\
  u_0u_1u_0&=u_0
  &
  u_1u_1&=du_1
  &
  u_0u_2&=u_2u_0
\end{align*}

\subsection{Chord diagrams}
\label{sec:chord}
\index{chord diagram}
\index{Brauer!category}
\index{category!Brauer}
Consider the full subcategory~$C$ of $\Corel$ (see \secr{corel}) whose morphisms
are corelations $\xymatrix@C=3ex@R=3ex{X\ar[r]|-f&Y&\ar[l]|-gZ}$ which are
\emph{one to one}, meaning that for every $y\in Y$ the cardinal of
$f^{-1}(y)\cup g^{-1}(y)$ is precisely~$2$. A morphism can be represented by
drawing a line between two elements of the source~$X$ or the target~$Z$ which are sent
to the same element of~$Y$ by $f$ or $g$, as on the left below:
\begin{equation}
  \label{eq:chord-ex}
  \fig{chord-ex}
  \qquad\qquad\qquad\qquad
  \satex{chord-ex}
\end{equation}
such diagrams are sometimes called \emph{chord diagrams} or \emph{Brauer
  linkings}~\cite{brauer1937algebras}. The resulting PRO, which we call the
\emph{Brauer category}, has a presentation with generators
\begin{align*}
  \satex{cap}
  &&
  \satex{cup}
  &&
  \satex{gamma}  
\end{align*}
such that the last one is a symmetry (\cref{sec:free-sym}) and the relations
\begin{align}
  \label{eq:chord-rel}
  \begin{array}{r@{\ }l@{\qquad\qquad\quad}r@{\ }l@{\qquad\qquad\quad}r@{\ }l}
    \satex{adj1-l}
    &\TO
    \satex{adj1-r}
    &
    \satex{cap-com-l}
    &\TO
    \satex{cap-com-r}
    &
    \satex{cap-cup-l}
    &\TO
    \satex{empty}
    \\[4ex]
    \satex{adj2-l}
    &\TO
    \satex{adj2-r}
    &
    \satex{cup-com-l}
    &\TO
    \satex{cup-com-r}
  \end{array}
\end{align}
hold, see~\cite{hughes2008linking}. The morphism corresponding to the diagram on
the left of~\eqref{eq:chord-ex} is shown on the right. This category is the free
loop-free self-dual compact closed symmetric monoidal category on an object, in sense
of~\secr{free-sd-compact}. Given $n\in\N$, the monoid of endomorphisms on~$n$ is precisely
the $n$-th Brauer monoid, as described in~\secr{mon-brauer}. In terms of the
presentation of the monoid, the generators $a_i$ and $u_i$ are respectively
interpreted as
\begin{align*}
  a_i&=\satex{br-ai}
  &
  u_i&=\satex{br-ui}  
\end{align*}
and for instance, we have the derivation of the following relations:
\begin{align*}
  \satex{br-aua-l}
  &\TO
  \satex{br-aua-r}
  &
  \satex{br-uau-l}
  &\TO
  \satex{br-uau-r}
  &
  \satex{br-aauu-l}
  &\TO
  \satex{br-aauu-r}
  \,\pbox.
  \\
  a_0u_1a_0
  &=
  a_1u_0a_1
  &
  u_0a_1u_0
  &=
  u_0
  &
  a_0a_1u_0u_2
  &=
  a_2a_1u_0u_2
\end{align*}
A variant without loop-freeness can be obtained by dropping the last rule
of~\eqref{eq:chord-rel}. The resulting category corresponds to the subcategory
of $\qCospan(\Simplaug)$ whose morphisms satisfy a similar condition as above.

\subsection{Frobenius algebras and dualities}
\label{sec:frob-pairing}
A Frobenius algebra is always canonically equipped with a notion of duality and
can even be characterized in terms of this structure, see
\cite{kock2004frobenius} for details.

\begin{proposition}
  The following theories present the same category:
  \begin{enumerate}
  \item The theory of Frobenius algebras, generated by
    $\satex{mu-small}$, $\satex{eta-small}$,
    $\satex{delta-small}$, and $\satex{eps-small}$ such that
    \begin{itemize}
    \item $(\satex{mu-small},\satex{eta-small})$ is a monoid,
    \item $(\satex{delta-small},\satex{eps-small}$) is a
      comonoid,
    \item the Frobenius relations~\eqref{eq:frob} hold between
      $\satex{delta-small}$ and $\satex{mu-small}$.
    \end{itemize}
  \item The theory generated by $\satex{mu-small}$,
    $\satex{eta-small}$, $\satex{eps-small}$,
    $\satex{cap-small}$ such that
    \begin{itemize}
    \item $(\satex{mu-small},\satex{eta-small})$ is a monoid,
    \item the following relations hold:
      \begin{align*}
        \satex{cap-mu-l}
        &\TO
        \satex{cap-mu-r}
        &
        \satex{cap-eps-l}
        &\TO
        \satex{eta0}
        \OT
        \satex{cap-eps-r}
        \,\pbox.
      \end{align*}
    \end{itemize}
  \item The theory generated by $\satex{mu-small}$,
    $\satex{eta-small}$, $\satex{cap-small}$,
    $\satex{cup-small}$ such that
    \begin{itemize}
    \item $(\satex{mu-small},\satex{eta-small})$ is a monoid,
    \item $(\satex{cap-small},\satex{cup-small})$ is a duality,
    \item the following relation holds:
      \[
        \satex{cap-mu-l}\TO\satex{cap-mu-r}
      \]
      or equivalently the following relation holds:
      \[
        \satex{mu-cup-l}\TO\satex{mu-cup-r}
        \pbox.
      \]
    \end{itemize}
  \end{enumerate}
\end{proposition}
\begin{proof}
  The equivalence between theories can be derived using Tietze
  transformations. We sketch the equivalence between the first and last one. In
  the theory~(1), we can define
  \begin{align*}
    \satex{cap0}&=\satex{eta-delta-l}
    &
    \satex{cup0}&=\satex{mu-eps-l}    
  \end{align*}
  and conversely, in the theory (3), we can define
  \begin{align*}
    \satex{delta0}&=\satex{cap-mu-r}
    &
    \satex{eps0}&=\satex{eta-cup-l}    
  \end{align*}
  In both cases, the required relations are derivable.
\end{proof}

\noindent
Various extensions of this result are possible in order to take into account
variations on the notion of Frobenius algebra (commutative, special, etc.).

\label{sec:lin-rel-dual}
Following the same ideas as previously, in the theory of interacting Hopf
algebras (\secr{lin-rel}) a self-duality can be defined by
\begin{align*}
  \satex{cap00}&=\satex{linrel-cap}
  &
  \satex{cup0}&=\satex{linrel-cup}  
\end{align*}
and the structure can be axiomatized as follows taking these as
generators~\cite{baez2015categories}:

\begin{proposition}
  The following theories present the same category:
  \begin{enumerate}
  \item The theory of interacting Hopf algebras.
  \item The theory generated by $\satex{mu-small}$,
    $\satex{eta-small}$, $\satex{delta3-small}$,
    $\satex{eps3-small}$, $\satex{gamma-small}$,
    $\satex{coef2-small}$, $\satex{cap-small}$,
    $\satex{cup-small}$ such that
    \begin{itemize}
    \item
      $(\satex{mu-small},\satex{eta-small},\satex{delta3-small},\satex{eps3-small},\satex{gamma-small},\satex{coef2-small})$
      is a bicommutative linear bialgebra,
    \item $(\satex{cap-small},\satex{cup-small})$ is a duality,
    \item
      $(\satex{mu-small},\satex{eta-small},\satex{mu-t-small},\satex{eta-t-small},\satex{gamma-small})$
      is a bicommutative extraspecial Frobenius algebra,
    \item
      $(\satex{delta3-t-small},\satex{eps3-t-small},\satex{delta3-small},\satex{eps3-small},\satex{gamma-small})$
      is a bicommutative extraspecial Frobenius algebra.
    \end{itemize}
  \end{enumerate}
\end{proposition}

\subsection{Tangles}
\label{sec:tangles}
A category of tangles can be defined as a variation of the category of braids,
see \cref{sec:braid-mon,sec:braid-cat}, intuitively by allowing wires to
loop.
Given natural numbers $m,n\in\N$, a \emph{tangle}\index{tangle} from~$m$ to~$n$
is an embedding
\[
  t
  :
  T
  \to
  X
\]
where~$T$ is a 1-manifold with boundary and $X=\R^2\times[0,1]$, such that the
image of the boundary~$\partial T$ of~$T$ is of the form
\[
  t(\partial T)
  =
  \setof{(0,i,0)}{i\in\N, 0\leq i<m}
  \sqcup
  \setof{(0,i,1)}{i\in\N, 0\leq i<n}
\]
where the natural number $m$ (\resp $n$) is called the source (\resp target) of
the tangle. Graphically, a tangle from~$4$ to~$2$ can be pictured as
\[
  \fig[height=4cm]{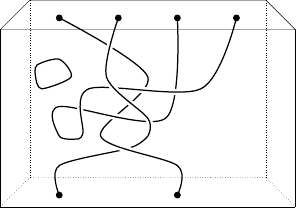}
\]
Tangles are considered up to endpoint-preserving isotopy: we identify two
tangles $t:T\to X$ and $t':T'\to X$
for which there exists a continuous map $h:[0,1]\to X^X$ such that
$h(0)=\id_X$, $h(1)\circ t=t'$ and, for every $t\in[0,1]$, $h(t)$ is a
homeomorphism such that $h(t)(\partial T)=h(0)(\partial T)$.
\nomenclature[Tang]{$\Tang$}{category of tangles}
We write~$\Tang$ for the PRO with tangles as morphisms with expected composition
(corresponding to linking wires) and tensor product (corresponding to
juxtaposition of diagrams).

The category~$\Tang$ is the free braided monoidal category on a self-dual
object, which moreover satisfies the first Reidemeister
move~\eqref{eq:reidemeister1} below, see~\cite{freyd1989braided,
  freyd1992coherence}. Explicitly, this means that it admits a presentation with
generators
\begin{align*}
  \satex{tau}
  &&
  \satex{tau2}
  &&
  \satex{cap}
  &&
  \satex{cup}  
\end{align*}
subject to the relations
\begin{itemize}
\item first Reidemeister move:
  \begin{equation}
    \label{eq:reidemeister1}
    \satex{knot1-l}\TO\satex{knot1-c}\OT\satex{knot1-r}
  \end{equation}
\item second Reidemeister move:
  \[
    \satex{knot2-l}\TO\satex{knot2-c}\OT\satex{knot2-r}
  \]
\item third Reidemeister move (aka Yang-Baxter rule):
  \[
    \satex{knot3-l}\TO\satex{knot3-r}
  \]
\item zigzag relations:
  \[
    \satex{adj1-l}
    \TO
    \satex{adj1-r}
    \OT
    \satex{adj2-l}
  \]
\item the naturality relations:
  \begin{align*}
    \satex{cap-tau-nat-l}
    &\TO
    \satex{cap-tau-nat-r}
    &
    \satex{cap-tau2-nat-l}
    &\TO
    \satex{cap-tau2-nat-r}
    \\
    \satex{tau2-cap-nat-l}
    &\TO
    \satex{tau2-cap-nat-r}
    &
    \satex{tau-cap-nat-l}
    &\TO
    \satex{tau-cap-nat-r}
    \\
    \satex{cup-tau2-nat-l}
    &\TO
    \satex{cup-tau2-nat-r}
    &
    \satex{cup-tau-nat-l}
    &\TO
    \satex{cup-tau-nat-r}
    \\
    \satex{tau-cup-nat-l}
    &\TO
    \satex{tau-cup-nat-r}
    &
    \satex{tau2-cup-nat-l}
    &\TO
    \satex{tau2-cup-nat-r}
  \end{align*}
  or equivalently, the sliding relations:
  \begin{align*}
    \satex{cap-tau-l}
    &\TO
    \satex{cap-tau-r}
    &
    \satex{cap-tau2-l}
    &\TO
    \satex{cap-tau2-r}
    \\
    \satex{cup-tau-l}
    &\TO
    \satex{cup-tau-r}
    &
    \satex{cup-tau2-l}
    &\TO
    \satex{cup-tau2-r}
  \end{align*}
\end{itemize}
If we replace relation~\eqref{eq:reidemeister1} by
\begin{equation}
  \label{eq:reidemeister1'}
  \satex{framed-l}
  \TO
  \satex{framed-r}  
\end{equation}
we present the category of \emph{framed tangles} or
\emph{ribbons}\index{ribbon}: in this category, morphisms correspond to pieces
of ribbon (instead of wire), \ie embeddings
\[
  r
  :
  T\times[0,1]
  \to
  X
\]
where $T$ is a $1$-manifold with boundary, e.g.,
\[
  \fig[height=2cm]{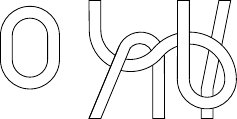}
\]
considered up to endpoint-preserving isotopy. Note that the relation \eqref{eq:reidemeister1} is not
satisfied since trying to strengthen the loop introduces a ``twist'' on the rope
\[
  \fig[height=2cm]{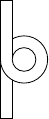}
  =
  \fig[height=2cm]{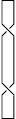}
\]
and is thus not the identity ribbon. If we remove the
relation~\eqref{eq:reidemeister1} (without adding \eqref{eq:reidemeister1'}), we
present the category of tangles up to regular isotopy, a variant of isotopy
which forces ribbons to always be flat against the plane, which prevents the
identity \eqref{eq:reidemeister1'} from holding:
\[
  \fig[height=2cm]{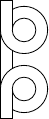}
  \qneq
  \fig[height=2cm]{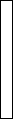}
  \,\pbox.
\]
A tangle from~$0$ to~$0$ is called a \emph{link}\index{link}, and a
\emph{knot}\index{knot} is a link $t:T\to X$ such that $T$ is the
$1$-sphere. The Reidemeister moves were originally introduced for
those~\cite{reidemeister1932knotentheorie}.

A variant of this category can be obtained by considering oriented tangles,
see~\cite{freyd1989braided, freyd1992coherence, kassel:quantum-groups,
  selinger2010survey}: the objects of this category are sequences of ``-'' and
``+'', and morphisms are oriented tangles up to endpoint-preserving planar
isotopy, e.g.,
\[
  \fig[height=3cm]{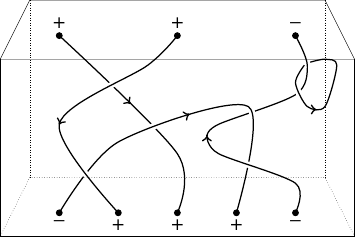}
\]
Of course, ribbon variants of this category can also be considered. The
resulting categories are pivotal (instead of being self-dual).

\subsection{First-order logic}
The dependencies between quantifiers in proofs for first-order logic are
characterized by the structure of free pivotal category on a bialgebra,
see~\cite{mimram2011structure}.

\section{Endomorphisms}
\subsection{The theory of endomorphisms}
\index{endomorphism}
\label{sec:end}
Consider the PRO whose morphisms are endomorphisms $f:n\to n$ consisting of
a list $(f_1,\ldots,f_n)$ of $n$~natural numbers. The composite of two morphisms
$f:n\to n$ and $g:n\to n$ is given by pointwise addition
$(f_1+g_1,\ldots,f_n+g_n)$, identities are lists $(0,\ldots,0)$, and tensor
product is given by concatenation of lists. This PRO admits a presentation with
one generator
\[
  \satex{sigma}
\]
and no relation, the generator being interpreted as the list $(1)$. An algebra
for this theory in a monoidal category consists of an object~$x$ together with
an \emph{endomorphism} $f:x\to x$ on~$x$.
This is a particular case of a presentation of the free monoidal category
generated by the monoid~$\N$, see \secr{act-mon}: since~$\N$ is free on one
generator, its action is entirely determined by an endomorphism (corresponding
to the action of~$1$).

\subsection{Actions of a set}
\label{sec:act-set}
Suppose fixed a set~$L$. The previous situation can be generalized by considering
the PRO where the morphisms are endomorphisms $f:n\to n$ consisting of a list
$(f_1,\ldots,f_n)$ of elements of $L^*$, the free monoid over~$L$: composition
is given by pointwise concatenation in~$L^*$ and tensor by concatenation of
lists. It admits a presentation with generators
\[
  \satex{coef}
\]
indexed by $a\in L$ and no relation. This is the theory for an \emph{action
  of~$L$}: it consists of an object~$x$ together with a family $f_a:x\to x$ of
endomorphisms of~$x$ indexed by $a\in L$.
Again, this is a particular case of a presentation of the free monoidal category
generated by the monoid~$\freecat{L}$ (the free monoid on the set~$L$). We
recover the theory of \secr{end} in the particular case where $L$ is a
singleton.

\subsection{Involutions}
The monoidal theory of involutions admits a presentation with one generator
\[
  \satex{sigma}
\]
and relation
\[
  \satex{sigma-inv-l}
  \TO
  \satex{sigma-inv-r}
\]
It is the free monoidal category on the monoid $\Z/2\Z$.

\subsection{The hyperoctahedral category}
\index{hyperoctahedral!category}
\label{sec:hyperoctahedral-cat}
The symmetric monoidal theory of involutions admits a presentation with two
generators
\begin{align*}
  \satex{sigma}
  &&
  \satex{gamma}  
\end{align*}
and relations
\begin{align*}
  \satex{sigma-inv-l}
  &\TO
  \satex{sigma-inv-r}
  &
  \satex{sigma-gamma-l}
  &\TO
  \satex{sigma-gamma-r}
  &
  \satex{sigma-gamma2-l}
  &\TO
  \satex{sigma-gamma2-r}  
\end{align*}
together with usual relations for symmetries, see \secr{pres-sym}.

Given $n\in\N$, the monoid of endomorphisms on~$n$ is in fact a group called the
\emph{hyperoctahedral group} or \emph{signed symmetric group}, and
noted~$C_n$, see \secr{hyperoctahedral}. It can be
described as the group of \emph{signed permutations} of a set with $n$~elements,
\ie $(n{\times}n)$-matrices where each column and each row contains one non-null
coefficient which is either $1$ or $-1$, with usual multiplication and
identities, the two generators being respectively interpreted as
\[
  \pa{
    \begin{matrix}
      -1
    \end{matrix}
  }
  \qqtand
  \pa{
    \begin{matrix}
      0&1\\
      1&0
    \end{matrix}
  }
  \pbox.
\]
It can also be described as the wreath product
\[
  C_n
  =
  S_2\wr S_n
\]
of symmetric groups.
The presentation of~$C_n$ given in \secr{hyperoctahedral} can be recovered as the
associated presentation of category (in the sense of \secr{3-pres-2-pres}), where the
generators~$a_0$ and $a_i$ (with $0<i<n$) respectively correspond to
\[
  \satex{Bn-a0}
  \qqtand
  \satex{Bn-ai}
  \,\pbox.
\]
For instance, the relation $a_0a_1a_0a_1=a_1a_0a_1a_0$ corresponds to
\[
  \satex{Bn-rel-ex-l}
  \qTO
  \satex{Bn-rel-ex-r}
  \,\pbox.
\]

In this category, every morphism factors uniquely as a symmetry (\ie a matrix
containing only $0$ and~$1$) followed by a diagonal matrix. This provides a
factorization system inducing a description of the category as
$\Bij\otimes_\dlaw\freemoncat{\N}$ where~$\freemoncat\N$ is the free monoidal
category over the monoid~$\N$, considered as a category, and $\dlaw$ is the
distributive law induced by the factorization system. Of course, this category
being self-dual, we also have a decomposition as
$\freemoncat\N\otimes_\dlaw\Bij$.

\subsection{Progressive ribbons}
\label{sec:progressive-ribbons}
\index{ribbon!progressive}
\index{progressive ribbon}
\index{balanced category}
\index{category!balanced}
\index{twist}
The category~$\category{R}$ is the subcategory of the category of ribbons up to
endpoint-preserving isotopy described in \cref{sec:tangles}, with the same objects, where
\begin{itemize}
\item we restrict to ribbons which are \emph{progressive}, \ie always ``go
  down'', so that the ribbon on the left is valid but not the one on the right:
  \begin{align*}
    \fig[height=12mm]{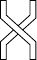}
    &&
    \fig[height=12mm]{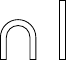}
    \,\pbox,
  \end{align*}
\item we restrict to ribbons which always show the ``same face'' at the
  boundary, so that the ribbon on the left is valid, but not the one on the
  right:
  \begin{align*}
    \fig[height=12mm]{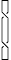}
    &&
    \fig[height=12mm]{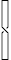}
    \,\pbox.
  \end{align*}
\end{itemize}
For instance, we have the following morphism $3\to 3$:
\[
  \fig[height=2cm]{ribbon-progressive-ex}
  \,\pbox.
\]
Note that the ribbons can be twisted. This category admits a presentation with
generators
\begin{align*}
  \satex{sigma}
  &&
  \satex{sigma2}
  &&
  \satex{tau}
  &&
  \satex{tau2}    
\end{align*}
such that the two first generators are mutually inverse and the two last
generators form a braiding (\cref{thm:free-bmc}).
For instance the above morphism corresponds to the diagram
\[
  \satex{ribbon-progressive-ex}
  \,\pbox.
\]
This category can be shown to be the free balanced category on an
object~\cite{joyal1991geometry}. We recall that a \emph{balanced category} is a
braided monoidal category equipped with a natural family of isomorphisms
$\theta_u:u\to u$, called \emph{twists}, such that $\theta_\monunit=\monunit$
and the following diagram commutes for every pair of objects $u$ and $v$:
\[
  \xymatrix@C=5ex@R=3ex{
    v\otimes u\ar[d]_{\theta_{v\otimes u}}\ar[r]^{\gamma_{u,v}}&u\otimes v\ar[d]^{\theta_u\otimes\theta_v{\displaystyle.}}\\
    v\otimes u&\ar[l]^{\gamma_{v,u}}u\otimes v
  }
\]
Graphically,
\[
  \fig[height=2cm]{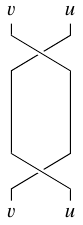}
  =
  \fig[height=2cm]{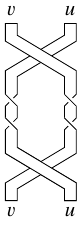}
  \,\pbox.
\]

Given $n\in\N$, the group of endomorphisms of~$n$ is the group $R_n$ of
progressive ribbons with $n$ strands described in
\secr{progressive-ribbon-group}:
\[
  \category{R}
  =
  \coprod_{n\in\N}R_n
\]
It can be obtained as the wreath product of the additive group of integers with
the $n$-th braid group:
\[
  R_n= \Z\wr C_n
  \pbox.
\]

\subsection{The pearl necklace}
\index{perl}
\index{polygraph!of pearls}
The \emph{pearl necklace} is the PRO, introduced in \cref{sec:pearls}, presented
by the 3-polygraph with three generators
\begin{align*}
  \satex{sigma}
  &&
  \satex{adj-eta}
  &&
  \satex{adj-eps}  
\end{align*}
subject to the four relations
\begin{align*}
  \satex{sigma-adj-eta-l}
  &\TO
  \satex{sigma-adj-eta-r}
  &
  \satex{sigma-adj-eps-l}
  &\TO
  \satex{sigma-adj-eps-r}
  &
  \satex{adj1-l}
  &\TO
  \satex{adj1-r}
  &
  \satex{adj2-l}
  &\TO
  \satex{adj2-r}
  \,\pbox.
\end{align*}
This is the free monoidal category on a self-dual object together with an
endomorphism, see also \cref{sec:free-sd-compact,sec:end}. This presentation was
introduced and studied in~\cite{guiraud2009higher} as a first example of a
convergent finite presentation which does not have a finite derivation type, see
\cref{sec:pearls}.

\subsection{Directed acyclic graphs}
\index{directed acyclic graph}
\label{sec:dag}
A \emph{directed acyclic graph}, or \emph{DAG}, is a graph
\[
  G
  =
  \xymatrix@C=3ex@R=3ex{G_0&\ar@<-.5ex>[l]_s\ar@<.5ex>[l]^tG_1}
\]
which is acyclic, \ie every path from a vertex to itself is empty.
Any~$n\in\N$ can canonically be seen as a DAG with the set
$\intset{n}=\set{0,\ldots,n-1}$ as set of vertices and no edge; we write
$\intset{n}$ for this graph.
A vertex~$x\in G_0$ of a DAG~$G$ is \emph{minimal} (\resp \emph{maximal}) when
there is no edge with~$x$ as target (\resp as source). Given a set $V\subseteq G_0$,
the \emph{restriction} of~$G$ to~$V$ is the subgraph of~$G$ with~$V$ as
vertices, whose edges are the edges~$G$ such that both their source and their target belong
to~$V$. Given a set $V\subseteq G_0$, the graph obtained by \emph{hiding}~$V$ in~$G$ is
the graph obtained from~$G$ by adding a new edge $x\to z$ for every vertex
$y\in V$ and pair of edges $x\to y$ and $y\to z$ in~$G$
and then restricting the resulting graph to~$G_0\setminus V$.

We can build a PROP $\category{DAG}$ of DAGs, where an object is a natural number
and a morphism from~$m$ to~$n$ is a cospan
\[
  \xymatrix@C=3ex@R=3ex{
    \intset{m}\ar[r]^f&G&\ar[l]_g\intset{n}
  }
\]
where~$G$ is a finite graph, $f$ and~$g$ are injective morphisms of graphs such
that for every $i\in\intset{m}$ (\resp $i\in\intset{n}$) the vertex $f(i)$
(\resp $g(i)$) is a minimal (\resp maximal) vertex in~$G$, and the images of~$f$
and~$g$ are disjoint. A vertex of~$G$ is a \emph{source} (\resp \emph{target},
\resp \emph{internal}) \emph{vertex} when it is in the image of~$f$ (\resp in
the image of~$g$, \resp neither in the image of~$f$ nor the image of~$g$). We
will picture by
\[
  \fig{dag-ex}
\]
a morphism with~$x$ and~$y$ as internal vertices, two source vertices, and three
target vertices. The composite of two morphisms $G:m\to n$ and $H:n\to o$ is
given by computing the pushout
\[
  \xymatrix@C=3ex@R=3ex{
    &&G\times_nH\\
    &G\ar@{.>}[ur]^{h'}&&\ar@{.>}[ul]_{g'}H\\
    \intset{m}\ar[ur]^f&&\ar[ul]_g\intset{n}\ar[ur]^h&&\ar[ul]_i\intset{o}
  }
\]
and then hiding $h'(g(\intset{n}))$ in the resulting graph, with the expected
cospan morphisms induced from $h'\circ f$ and $g'\circ i$. For instance, we have
the following composition of morphisms:
\[
  \fig{dag-comp-ex}
\]
The identity on~$n$ is the graph obtained from $\intset{n}\sqcup\intset{n}$ by
adding, for every $i\in\intset{n}$ a vertex from the first copy of~$i$ to the
second one, e.g., the identity on~$3$ is
\[
  \fig[height=15mm]{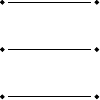}
\]
Tensor product is given on morphisms by disjoint union and symmetries are the
expected ones. It is shown in~\cite{fiore2013algebra} that this PROP admits a
presentation with generators
\begin{align*}
  \satex{mu}
  &&
  \satex{eta}
  &&
  \satex{delta}
  &&
  \satex{eps}
  &&
  \satex{gamma}
  &&
  \satex{sigma}  
\end{align*}
such that the five first generators satisfy the axioms for bicommutative
bialgebras (\cref{sec:bialgebra}) and the penultimate one is a symmetry
(\cref{sec:free-sym}). In this language, the composition of the two above
morphisms can be written
\[
  \satex{dag-comp-sd-ex}
\]
If we restrict the morphisms of the category $\category{DAG}$ to simple graphs,
\ie forbidding multiple edges with same source and same target, we can obtain a
presentation of the resulting category by further adding the
relation~\eqref{eq:rel-special} of special bialgebras.

\subsection{Posets}
\index{poset}
As a variant of previous situation, consider the PROP whose objects are integers
and morphisms $m\to n$ are cospans
\[
  \xymatrix@C=3ex@R=3ex{
    \intset{m}\ar[r]^f&E&\ar[l]_g\intset{n}
  }
\]
of finite posets, where $\intset{n}$ is the discrete poset with $n$ elements,
and the morphisms are injective non-decreasing functions with disjoint images,
such that the images of~$f$ (\resp $g$) are minimal (\resp maximal) elements of
the poset. Composition is obtained from the pushout by removing the elements
which are identified in the interface, as in \secr{dag}. This category admits
the same presentation as the category of simple graphs, with the extra
\emph{transitivity} relation
\[
  \satex{dag-trans-l}
  \TO
  \satex{dag-trans-r}
\]
see~\cite{mimram2015posets} for details.

\subsection{Frobenius}
\label{sec:graph-cospan}
We now consider a variant of the categories introduced in the previous section, as well
as a generalization of the theory of special commutative Frobenius algebras,
which is the category $\qCospan(\Fun)$ as explained \cref{sec:frob-cospan}.

\nomenclature[Circ]{$\Circ$}{category of circuits}
A formalization of ``circuits'' (such as electric circuits made of electronic
components) was introduced and studied in~\cite{rosebrugh2005generic} by
considering the category~$\Circ$, which is the full subcategory of the
category~$\qCospan(\FinGraph)$ whose objects are integers~$n\in\N$, seen as
graphs with $\intset{n}$ as vertices and no edge, morphisms being isomorphism
classes of cospans of finite graphs. A morphism $G:m\to n$ in this category is a
graph~$G$ (considered up to isomorphism) together with two functions
$f:\intset{m}\to G$ and $g:\intset{n}\to G$, composition is given by pushout and
identities are cospans of the form
$\xymatrix@C=3ex@R=3ex{\intset{n}\ar[r]&\intset{n}&\ar[l]\intset{n}}$ with both
morphisms being identities.
The resulting category is a PROP which admits a presentation with generators
\begin{align*}
  \satex{mu}
  &&
  \satex{eta}
  &&
  \satex{delta}
  &&
  \satex{eps}
  &&
  \satex{gamma}
  &&
  \satex{sigma}  
\end{align*}
such that the five first generators satisfy the axioms for special commutative
Frobenius algebras (\cref{sec:frobenius-algebra}) and the penultimate one is a
symmetry (\cref{sec:free-sym}).

Given a set~$L$ of \emph{labels} (thought of as the possible components of our
circuits), we can consider a variant of the previous category where the vertices
of graphs are labeled in~$L$. It is shown to have a similar presentation as
above, with the generator $\satex{sigma}$ being replaced by a family of
generators $\satex{coef}$, indexed by~$a\in L$. From \secr{act-set}, we deduce
that this is the theory for special Frobenius algebras equipped with an action
of~$L$~\cite[Proposition~3.2]{rosebrugh2005generic}. This category has found
many applications to modeling networks~\cite{baez2018props,baez2015compositional}.

\newcommand{\Net}[1]{\category{Net}_{#1}}

\section{Nets}
Some of the previous examples consist in categories whose morphisms are graphs of some sort (e.g.,
\cref{sec:dag} or \cref{sec:graph-cospan}). The kind of constructions performed there
can be adapted in order to build freely generated monoidal categories as
follows, by formalizing the networks occurring in string-diagrammatic
representations of morphisms. Those networks will comprise nodes,
which correspond to $2$-generators, ports which represent the inputs and outputs
of the generators and of the whole net, and wires which link ports together.

Throughout the section, we suppose fixed a \emph{signature} $2$-polygraph~$P$
with $P_0=\set{\star}$. We recall that, given an element $u=a_0\ldots a_{n-1}$
of~$\freecat{P_1}$, we write $\sizeof{u}=n$ for its length.

\subsection{Definition}
Following~\cite{lafont1997interaction}, a \emph{$P$-net} consists of
\begin{itemize}
\item a finite set~$N$ of \emph{nodes},
\item a labeling function $\ell:N\to P_2$,
\item finite totally ordered sets $X^-=\set{x_0^-,\ldots,x_m^-}$ and
  $X^+=\set{x_0^+,\ldots,x_n^+}$ of \emph{input} and \emph{output ports},
\item a labeling function $\ell:X\to P_1$, where
  \[
    X
    =
    X^-\sqcup X^+\sqcup X_\circ^-\sqcup X_\circ^+
  \]
  is the set of \emph{ports} and
  \begin{align*}
    X_\circ^-
    &=
    \setof{\nu_i^-}{\nu\in N, 0\leq i<\sizeof{\src1(\ell(\nu))}}
    \\
    X_\circ^+
    &=
    \setof{\nu_j^+}{\nu\in N, 0\leq j<\sizeof{\tgt1(\ell(\nu))}}
  \end{align*}
  are the sets of \emph{inner input} and \emph{output ports} respectively,
\item a finite set~$W$ of \emph{wires},
\item a labeling function~$\ell:W\to P_1$,
\item a \emph{boundary} $\partial w\subseteq X$ for every wire $w\in W$,
\end{itemize}
such that
\begin{itemize}
\item for every node $\nu\in N$ such that $\src1(\nu)=a_1\ldots a_p$ and
  $\tgt1(\nu)=b_1\ldots b_q$, we have $\ell(\nu_i^-)=a_i$ and
  $\ell(\nu_j^+)=b_j$ for $0\leq i<p$ and $0\leq j<q$,
\item for every $w\in W$, $\partial w$ contains $0$ or $2$ elements, and in the
  case $\partial w=\set{x,y}$ we have $\ell(w)=\ell(x)=\ell(y)$,
\item the sets $\partial w$ form a partition of $X$, \ie
  $X=\bigcup_{w\in W}\partial w$ and $w\neq w'$ implies
  $\partial w\cap\partial w'=\emptyset$.
\end{itemize}
The \emph{source} (\resp \emph{target}) of such a net is
$\ell(x_0^-)\ldots\ell(x_m^-)$ (\resp $\ell(x_0^+)\ldots\ell(x_n^+)$).

For instance, with $P_1=\set{a,b}$ and
$P_2=\set{\alpha:aaaa\To a,\beta:ab\To ab}$, the diagram on the left below
\[
  \satex{net-ex}
  \qquad\qquad
  \satex{net-ex-labels}
\]
can be encoded as the net (figured on the right) with
\begin{align*}
  N&=\set{\mu,\nu}
  &
  X^-&=\set{x,x'}
  &
  X^+&=\set{y}
  &
  W&=\set{w_1,\ldots,w_7}
\end{align*}
labels being
\begin{align*}
  \ell(\mu)&=\alpha
  &
  \ell(\nu)&=\beta
  &
  \ell(x)&=a
  &
  \ell(x')&=b
  &
  \ell(\mu_i^-)&=a
  &
  \ell(w_1)&=a
  &
  \ldots
\end{align*}
and boundaries being
\begin{align*}
  \partial w_1&=\set{\mu_0^-,\mu_1^-}
  &
  \partial w_2&=\set{\mu_2^-,x}
  &
  \partial w_3&=\set{\nu_1^-,x'}
  &
  \ldots&
  &
  \partial w_7&=\emptyset
\end{align*}
As it can be observed for $w_7$ above, wires with empty boundary encode
``loops''.

Nets are considered up to isomorphism, \ie renaming of ports, wires, and
nodes, preserving labels. Two nets can be composed by linking wires along inner boundary ports and
removing those ports; two nets can also be tensored by juxtaposition. We write
$\Net P$ for the resulting monoidal category, with $\freecat{P_1}$ as objects
and nets as morphisms.

\begin{proposition}
  The monoidal category~$\Net P$ is the free self-dual compact closed category
  on~$P$, in the sense of \secr{free-sd-compact}.
\end{proposition}

\noindent
Various other free constructions can be obtained by considering monoidal
subcategories obtained by restricting the notion of net.

\subsection{Traced categories}
\index{traced category}
\index{category!traced}
\label{sec:net-traced-cat}
The free traced category on~$P$, see~\cite{joyal1996traced}, can be obtained by
forbidding wires to link two input ports or two output ports. More formally, we
restrict the category~$\Net P$ to nets such that for every wire $w\in W$, its
boundary~$\partial w$ is either empty or of the form~$\partial w=\set{x,y}$ with
\[
  (x,y)
  \quad\in\quad
  \pa{X^-\times X_\circ^-}\cup\pa{X_\circ^+\times X_\circ^-}\cup\pa{X_\circ^+\times X^+}
\]
see~\cite{hasegawa2008finite} for details.

\subsection{Free symmetric monoidal categories}
\index{free!symmetric monoidal category}
Given a net, we define a relation~$\prec$ on nodes as the smallest transitive
relation such that $\mu\prec\nu$ whenever there exists a wire $w$ such that
$\partial w=\set{\mu_j^+,\nu_i^-}$ for some indices $i$ and $j$,
\ie some output port of~$\mu$ is linked to some input port of~$\nu$.
The free symmetric monoidal category on~$P$ can obtained by considering the
subcategory of~$\Net P$ whose morphisms are nets such that the relation $\prec$
is acyclic, the condition of \secr{net-traced-cat} is satisfied, and
$\partial w\neq\emptyset$ for every wire $w\in W$.

\subsection{Interaction nets}
\index{interaction!net}
\label{sec:interaction-net}
The morphisms in~$\Net P$ are called \emph{interaction nets} whenever all the
generators in~$P_2$ are of the form $\alpha:a_1\ldots a_n\To b$, \ie their
target consists of only one generator. An \emph{interaction rule}
$A:\phi\TO\psi$ consists of a pair of parallel nets $\phi:u\To v$ and
$\psi:u\To v$ which are of the form
\[
  \satex{ir-l}
  \TO
  \satex{ir-r}
  \vspace{-2ex}
\]
(in particular, $v$ is always an identity).
These were introduced by
Lafont~\cite{lafont1989interaction,lafont1997interaction} in order to provide a
distributed model of computation, where there is no global synchronization.
Because of the shape of the rules, the corresponding rewriting system
is confluent (as a rewriting system on nets) since there is no critical branching.
The signature~$P$ together with a set~$P_3$ of interaction rules can be encoded
as a 3\nbd-poly\-graph~$Q$ obtained from the free self-dual compact closed symmetric monoidal
category on~$P$, as described in \secr{free-sd-compact}, by adding the elements
of~$P_3$ as $3$-generators.

\index{interaction!combinator}
As a particular case, \emph{interaction
  combinators}~\cite{lafont1997interaction} are the following interaction
nets. The sets of $0$- and $1$-generators are respectively $P_0=\set{\star}$ and
$P_1=\set{a}$. The $2$-generators are
\begin{align*}
  \satex{ic-gamma}
  &&
  \satex{ic-delta}
  &&
  \satex{ic-eps}
\end{align*}
respectively called \emph{constructor} (noted~$\gamma$), \emph{duplicator}
(noted~$\delta$), and \emph{eraser} (noted~$\varepsilon$),
and the interaction rules are
\begin{align*}
  \satex{ic-gamma-gamma-l}
  &\TO
  \satex{ic-gamma-gamma-r}
  &
  \satex{ic-delta-delta-l}
  &\TO
  \satex{ic-delta-delta-r}
  \\
  \satex{ic-gamma-eps-l}
  &\TO
  \satex{ic-gamma-eps-r}
  &
  \satex{ic-delta-eps-l}
  &\TO
  \satex{ic-delta-eps-r}
  &
  \satex{ic-eps-eps-l}
  &\TO
  \satex{empty}
\end{align*}
\[
  \satex{ic-gamma-delta-l}
  \TO
  \satex{ic-gamma-delta-r}
\]
as well as symmetric ones.
Lafont's original notation was
\begin{align*}
  \satex{icl-gamma-gamma-l}
  &\TO
  \satex{icl-gamma-gamma-r}
  &
  \satex{icl-delta-delta-l}
  &\TO
  \satex{icl-delta-delta-r}
  &
  \satex{icl-gamma-delta-l}
  &\TO
  \satex{icl-gamma-delta-r}
  \\
  \satex{icl-gamma-eps-l}
  &\TO
  \satex{icl-gamma-eps-r}
  &
  \satex{icl-delta-eps-l}
  &\TO
  \satex{icl-delta-eps-r}
  &
  \satex{icl-eps-eps-l}
  &\TO
  \satex{empty}
  \,\pbox.
\end{align*}
which suggests that they can also be formalized without resorting to the compact
closure, at the cost of having to add two copies of each $2$-generator (one
going upward and one going downward).
Interaction nets are ``universal'' in the sense that they can simulate any
net~\cite[Theorem~1]{lafont1997interaction}. In particular, this implies that they are
Turing-complete.

\section{Simplicial and Cubical Categories}
Given a category~$\C$, we write $\hat\C$ for the category of \emph{presheaves}
on this category: the objects of $\hat\C$ are functors $\C^\op\to\Set$ and
morphisms are natural transformations.

\subsection{Simplicial categories}
\index{simplicial!category}
\index{augmented simplicial category}
\index{presimplicial!set}
\index{simplicial!set}
In \secr{2pres-mon}, we have already described presentations of monoidal
categories, whose associated presheaf categories are widely used in algebraic
topology:
\begin{itemize}
\item for $\Simplaug$ presented in \cref{sec:pres-simpl}, $\hatSimplaug$ is the
  category of \emph{augmented simplicial sets},
\item for $\Simplinj$ presented in \cref{sec:incr-inj}, $\hatSimplinj$ is the
  category of \emph{augmented presimplicial} or \emph{semisimplicial sets},
\item for $\Fun$ presented in \cref{sec:pres-com-mon}, $\hat\Fun$ is the
  category of \emph{augmented symmetric simplicial
    sets}~\cite{grandis2001finite}.
\end{itemize}
The non-augmented variants can be obtained by taking presheaves on the same
categories with the object~$0$ removed. For each of those, from their
presentation as a monoidal category, we can deduce a presentation as a category
by using the method of \secr{3-pres-2-pres} and thus obtain an algebraic
description of the associated presheaves.

For instance, from the presentation of the PRO~$\Simplinj$ given in
\cref{sec:incr-inj}, we can deduce the following presentation of it as a
category:
\[
  \Pres{n}{d^n_i:n\to n+1}{d_j^nd_i^{n+1}\To d_i^nd_{j+1}^{n+1}}_{0\leq i\leq j\leq n}
\]
(see \exr{simpl-3-2-pres} for details). An \emph{augmented presimplicial set}
thus consists of a family of sets~$(X_n)_{n\in\N}$ together with functions
\[
  \partial_i^n
  :
  X_{n+1}
  \to
  X_n
\]
with $0\leq i\leq n$, satisfying relations dual to those of the above presentation. The description of
non-augmented presimplicial sets is similar, excepting that we constrain $n$ to
be strictly positive (there is no $X_0$ nor $\partial^0_0$). For those, an element of~$X_{n+1}$ can
be interpreted geometrically as an $n$-simplex and the function~$\partial_i^n$
as the function associating to an $n$-simplex its $i$-th face, obtained by
removing its $i$-th vertex; more generally, a simplicial set can be thought of
as the result of gluing simplices. For instance, the space
\[
  \fig[height=3cm]{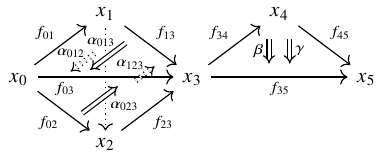}
\]
can be described by the simplicial set~$X$ with
\begin{align*}
  X_1&=\set{x_0,x_1,\ldots}
  &
  X_2&=\set{f_{01},f_{02},\ldots}
  &
  X_3&=\set{\alpha_{012},\ldots,\beta,\gamma}
  &
  X_4&=\set{A}  
\end{align*}
and $X_n=\emptyset$ for $n\geq 5$, face maps being
\begin{align*}
  \partial^2_0(\beta)&=f_{45}
  &
  \partial^2_1(\beta)&=f_{35}
  &
  \partial^2_2(\beta)&=f_{34}  
\end{align*}
and so on. A \emph{simplicial set}~$X$ is moreover equipped with
functions
\[
  \sigma_i^n
  :
  X_n
  \to
  X_{n+1}
\]
with $0\leq i<n$, sending an $n$-simplex to the corresponding $n{+}1$-simplex
with degenerated $i$-th face. For instance, given a $1$-simplex
$\xymatrix@C=3ex@R=3ex{x_0\ar[r]|{f}&x_1,}$ its images under $\sigma_0^3$ and
$\sigma_1^3$ are respectively
\begin{align*}
  \xymatrix@C=3ex@R=3ex{
    &x_0\ar@{}[d]|(.6){=}\ar[dr]^f&\\
    x_0\ar@{=}[ur]\ar[rr]_f&&x_1
  }
  &&
  \xymatrix@C=3ex@R=3
  ex{
    &x_1\ar@{}[d]|(.6){=}\ar@{=}[dr]&\\
    x_0\ar[ur]^f\ar[rr]_f&&x_1\pbox.
  }  
\end{align*}
Similarly, a \emph{symmetric} augmented simplicial set~$X$ is equipped with
functions
\[
  \gamma^n_f
  :
  X_n
  \to
  X_n
\]
indexed by bijections $f:\intset{n}\to\intset{n}$ sending an $n$-simplex to the
corresponding $n$-simplex with vertices renumbered according to~$f$.

\subsection{Cubical categories}
\index{cubical category}
\nomenclature[.C]{$\Cub$}{cubical category}
Developments similar to previous section can be performed using cubes instead of
simplices.
The \emph{precubical category} is the PRO~$\Cubinj$ generated by
\[
  \eta^-:0\to 1
  \qqtand
  \eta^+:1\to 0
\]
respectively depicted as
\begin{align*}
  \satex{etas}
  &&
  \satex{etat}  
\end{align*}
with no relation. As a category, $\Cubinj$ can be presented by
\[
  \Pres{n}{d^-_{n,i}:n\to n+1,d^+_{n,i}:n\to n+1}{d_{n,j}^\epsilon d_{n+1,i}^{\epsilon'}=d_{n,i}^{\epsilon'}d_{n+1,j+1}^\epsilon}_{0\leq i\leq j<n}
\]
with $\epsilon,\epsilon'\in\set{-,+}$.
A presheaf~$X$ on this category is called a \emph{precubical set}: it consists
of a family of sets~$(X_n)_{n\in\N}$, whose elements can be interpreted as
$n$-dimensional cubes, together with morphisms
\[
  \partial_i^\epsilon
  :
  X_{n+1}
  \to
  X_n
\]
with $0\leq i<n$ and $\epsilon\in\set{-,+}$, defined by
\[
  \partial_i^\epsilon
  =
  \unit{n-i}\otimes\eta^\epsilon\otimes\unit{n-i-1}
\]
and associating to each $(n{+}1)$-cube its source ($\epsilon=-$) or target
($\epsilon=+$) face in the $i$-th direction. For instance, the space
\[
  \fig[height=3.5cm]{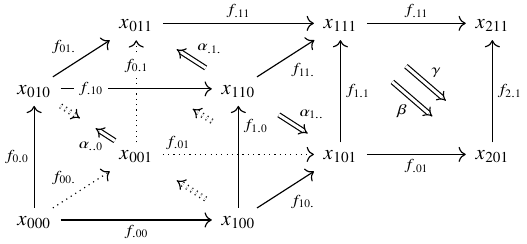}
\]
corresponds to the precubical set with
\begin{align*}
  X_0&=\set{x_{000},x_{001},\ldots}
  &
  \!
  X_1&=\set{f_{.00},f_{.01},\ldots}
  &
  \!\!
  X_2&=\set{\alpha_{..0},\alpha_{..1},\ldots}
  &
  \!\!\!
  X_3&=\set{A}  
\end{align*}
and $X_n=\emptyset$ for $n\geq 4$. Face maps are given by
\begin{align*}
  \partial^-_0(\beta)&=f_{1.1}
  &
  \partial^+_0(\beta)&=f_{2.1}
  &
  \partial^-_1(\beta)&=f_{\cdot 01}
  &
  \partial^+_1(\beta)&=f_{\cdot 11}  
\end{align*}
and so on.

Many variations on this category (and thus on associated presheaves) are
possible and considered in the literature, see~\cite{grandis2003cubical} for a
panorama which is briefly recalled here. The \emph{cubical category} $\Cub$ is
the PRO generated by
\begin{align*}
  \satex{etas}
  &&
  \satex{etat}
  &&
  \satex{eps}
\end{align*}
subject to the relations
\begin{align}
  \label{eq:cc-eta-eps}
  \satex{etas-eps}&\TO\satex{empty}
  &
  \satex{etat-eps}&\TO\satex{empty}
  \,\pbox.
\end{align}
A \emph{cubical set}~$X$ is a presheaf on this category and comes equipped with
morphisms
\[
  \sigma_i^n
  :
  X_n
  \to
  X_{n+1}
\]
with $0\leq i\leq n$, called \emph{degeneracies}, sending an $n$-cube to the
corresponding $(n{+}1)$-cube degenerated in dimension~$i$. For instance, the
images under $\sigma_0$ and~$\sigma_1$ of the $1$-cube
$\xymatrix@C=3ex@R=3ex{x\ar[r]|f&y}$ are respectively
\begin{align*}
  \xymatrix@C=3ex@R=3ex{
    y\ar@{=}[r]\ar@{}[dr]|{=}&y\\
    x\ar[u]^f\ar@{=}[r]&x\ar[u]_f
  }
  &&
  \xymatrix@C=3ex@R=3ex{
    x\ar[r]^f\ar@{}[dr]|{=}&y\\
    x\ar@{=}[u]\ar[r]_f&y\ar@{=}[u]\pbox.
  }  
\end{align*}
\index{symmetric!cubical category}
The \emph{symmetric cubical category}~$\Cubsym$ is the free symmetric monoidal
category on the PRO $\Cub$, it is generated by
\begin{align*}
  \satex{etas}
  &&
  \satex{etat}
  &&
  \satex{eps}
  &&
  \satex{gamma}  
\end{align*}
subject to the above relations~\eqref{eq:cc-eta-eps} and those of symmetries
(\cref{sec:pres-sym}). A presheaf~$X$ on it is called a \emph{symmetric cubical
  set} and is equipped with morphisms
\[
  \gamma^n_f
  :
  X_n
  \to
  X_n
\]
indexed by bijections $f:\intset{n}\to\intset{n}$ sending an $n$-cube to the
corresponding $n$-cube obtained by permuting the directions along $f$. For
instance the image of the $2$-cube on the left along the transposition
$\intset{2}\to\intset{2}$ is the $2$-cube on the right:
\begin{align}
  \label{eq:cc-2cube-trans}
  \vcenter{
    \xymatrix@C=3ex@R=3ex{
      x_{01}\ar[r]^{f_{.1}}&x_{11}\\
      x_{00}\ar[u]^{f_{0.}}\ar[r]_{f_{.0}}&\ar@{}[ul]|{\Toul}x_{10}\ar[u]_{f_{1.}}
    }
  }
  &&
  \vcenter{
    \xymatrix@C=3ex@R=3ex{
      x_{10}\ar[r]^{f_{1.}}&x_{11}\\
      x_{00}\ar[u]^{f_{.0}}\ar[r]_{f_{0.}}&\ar@{}[ul]|{\Todr}x_{01}\ar[u]_{f_{.1}}\pbox.
    }
  }
\end{align}
By definition, the \emph{cartesian cubical category}~$\Cubcart$ is the free cartesian category
on the PROP~$\Cubsym$: it is generated by
\begin{align*}
  \satex{etas}
  &&
  \satex{etat}
  &&
  \satex{eps}
  &&
  \satex{delta}
  &&
  \satex{gamma}  
\end{align*}
subject to the relations given in \cref{sec:free-cart-cat}, which
include~\eqref{eq:cc-eta-eps}. Note that this is the Lawvere theory for sets
which are bipointed (\ie equipped with two distinguished elements). A
presheaf~$X$ on this category is a symmetric cubical set equipped with morphisms
\[
  \delta^n_i
  :
  X_{n+1}
  \to
  X_n
\]
with $0\leq i<n$, which to an $(n{+}1)$-cube associate a diagonal $n$-cube (in
directions $i$ and $i+1$). For instance, the image under $\delta^1_0$ of the
$2$-cube on the left is a one cube as on the right
\begin{align*}
  \xymatrix@C=3ex@R=3ex{
    y_1\ar[r]&z\pbox.\\
    x\ar[u]\ar[r]&\ar@{}[ul]|{\Toul}\ar[u]y_0
  }
  &&
  \xymatrix@C=3ex@R=3ex{
    &z\\
    x\ar[ur]&
  }  
\end{align*}
The \emph{cubical category with connections}~$\Cubcon$\index{connection} is the PRO generated by
\begin{align*}
  \satex{etas}
  &&
  \satex{etat}
  &&
  \satex{eps}
  &&
  \satex{mus}
  &&
  \satex{mut}
\end{align*}
subject to relations~\eqref{eq:cc-eta-eps}, as well as
\begin{align*}
  \satex{mue-ass-l}
  &\TO
  \satex{mue-ass-r}
  &
  \satex{mue-unit-l}
  &\TO
  \satex{mue-unit-c}
  &
  \satex{mue-unit-r}
  &\TO
  \satex{mue-unit-c}
  \\
  \satex{mue-eps-l}
  &\TO
  \satex{mue-eps-r}
  &
  \satex{mue-abs-l}
  &\TO
  \satex{mue-abs-c}
  &
  \satex{mue-abs-r}
  &\TO
  \satex{mue-abs-c}
\end{align*}
for $\eta,\epsilon\in\set{-,+}$ with $\eta\neq\epsilon$. A presheaf~$X$ on this
category is a \emph{cubical category with connections} and is equipped with
morphisms
\[
  \kappa^\epsilon_i
  :
  X_n
  \to
  X_{n+1}
\]
called \emph{connections} which produce degenerated cubes, but in another way than
degeneracies. For instance, the images under~$\kappa^-_0$ and $\kappa^+_0$ of
the $1$-cube $\xymatrix@C=3ex@R=3ex{x\ar[r]|f&y}$ are shown on the left and
those of the $2$-cube (e.g., left of~\eqref{eq:cc-2cube-trans}) are shown on the
right (omitting 2-cells for clarity):
\begin{align*}
  \vcenter{
    \xymatrix@C=3ex@R=3ex{
      y\ar@{=}[r]^{\phantom{f}}&y\\
      x\ar[u]^f\ar[r]_f&y\ar@{=}[u]
    }
  }
  &&
  \vcenter{
    \xymatrix@C=3ex@R=3ex{
      y\ar[r]^f&y\\
      x\ar@{=}[u]\ar@{=}[r]_{\phantom{f}}&y\ar[u]_f
    }
  }
  &&
  \fig{2cube-degen-s}
  &&
  \fig{2cube-degen-t}  
\end{align*}
A symmetric variant can be obtained by taking the free symmetric monoidal
category (in practice, one moreover asks both monoid structures to be
commutative). The \emph{cubical category with reversions}\index{reversion}~$\Cubrev$ is the PRO
generated by
\begin{align*}
  \satex{etas}
  &&
  \satex{etat}
  &&
  \satex{eps}
  &&
  \satex{delta}
  &&
  \satex{gamma}
  &&
  \satex{sigma}  
\end{align*}
satisfying the relations of~$\Cubcon$ together with
\begin{align*}
  \satex{sigma-gamma-l}
  &\TO
  \satex{sigma-gamma-r}
  &
  \satex{sigma-eps-l}
  &\TO
  \satex{sigma-eps-r}
  &
  \satex{sigma-inv-l}
  &\TO
  \satex{sigma-inv-r}
  &
  \satex{etas-sigma-l}
  &\TO
  \satex{etas-sigma-r}
  &
  \satex{mus-sigma-l}
  &\TO
  \satex{mus-sigma-r}
  \pbox.
\end{align*}
A presheaf~$X$ on this category is equipped with morphisms
\[
  \rho^n_i
  :
  X_n
  \to
  X_n
\]
reversing cubes in the $i$-th direction. There is also a cartesian
variant~\cite{cohen2016cubical}, where one usually requires axioms corresponding
to the Lawvere theory of de Morgan algebras, \ie bounded distributive lattices
with an idempotent negation.

\section{Quantum Processes}
\index{ZX-calculus}
\index{ZW-calculus}
\newcommand{\bra}[1]{\langle #1|}
\newcommand{\ket}[1]{|#1\rangle}
Let us briefly mention that presentations of PROPs are intensively used nowadays
in the study of quantum processes, see~\cite{coecke2017picturing} for an
in-depth introduction. One of the most notable axiomatic approaches is the
\emph{ZX-calculus}~\cite{coecke2008interacting}, aiming at modeling operations
on qubits, which is a presented PROP together with a canonical interpretation in
the category~$\FdHilb$ of finite-dimensional Hilbert spaces and linear maps,
which is a PROP when equipped with the usual tensor product of vector
spaces. The ZX-calculus is generated by
\begin{align*}
  \satex{gamma}
  &&
  \satex{cap}
  &&
  \satex{cup}
  &&
  \satex{zx-alpha}
  &&
  \satex{zx-beta}
  &&
  \satex{zx-h}  
\end{align*}
where $\alpha^m_n$ (\resp $\beta^m_n$) have $m$ inputs and $n$ outputs, such that various axioms
are satisfied, among which the fact that the first generator induces a symmetry
and the second and the third a self-duality. The generating object of this PROP
is interpreted as~$\cplx^2$. Using the traditional notation, we write
$\ket 0,\ket 1$ for the standard basis (also called the $Z$-basis) of~$\cplx^2$
and $\ket -,\ket +$ for the Bell basis (also called the $X$-basis) defined by
$\ket -=\frac 1{\sqrt 2}(\ket 0-\ket 1)$ and
$\ket +=\frac 1{\sqrt 2}(\ket 0+\ket 1)$. We also write $\bra 0$ for the adjoint
of $\ket 0$, etc. The interpretation of the generators is then given by
\begin{align*}
  \intp{\satex{cap}}&=\ket 0^{\otimes 2}+\ket 1^{\otimes 2}\\
  \intp{\satex{cup}}&=\bra 0^{\otimes 2}+\bra 1^{\otimes 2}\\
  \intp{\satex{zx-alpha}}&=\ket 0^{\otimes n}\bra 0^{\otimes m}+\ce^{\ci\alpha}\ket 1^{\otimes n}\bra 1^{\otimes m}\\
  \intp{\satex{zx-beta}}&=\ket +^{\otimes n}\bra +^{\otimes m}+\ce^{\ci\alpha}\ket -^{\otimes n}\bra -^{\otimes m}\\
  \intp{\satex{zx-h}}&=\ket 0\bra ++\ket 1\bra -
\end{align*}
The original axiomatization is not \emph{complete} in
general~\cite{de2014zx}, meaning that the resulting functor to $\FdHilb$ is not
faithful. However, the rules can be completed so that it is the
case~\cite{ng2017universal}: the proof is based on an alternative complete
axiomatization called the
\emph{ZW-calculus}~\cite{hadzihasanovic2015diagrammatic}.


\chapter{A Syntactic Description of Free \pdfm{n}-Categories}
\chaptermark{Free $n$-categories}
\label{chap:syntactic_descr}
In this chapter, we provide an explicit description of the free $n$-category
$\freecat P$ generated by an $n$-polygraph~$P$. Variants of this construction
can be found in~\cite[Deuxième partie]{penon1999approche},
\cite[Section~7]{makkai2005wordcomp}, and
\cite[Section~4.1]{metayer2008cofibrant}, but this chapter is mostly inspired by
the work of Makkai~\cite{makkai2005wordcomp}, where the proofs of most
assertions can be found, formulated in a slightly different language, see
also~\cite[Section~2]{forest2021computational}.

We first provide, in \cref{sec:ncat-syntax}, a formal definition of the syntax
of $n$\nbd-cate\-gories, \ie a description of the morphisms in an $(n{+}1)$-category
freely generated by an $n$-polygraph, allowing reasoning by induction on its
terms to prove results on free categories.
It turns out that this syntax for $n$-categories, which corresponds to the one
introduced in \cref{chap:n-cat} and used throughout the book, is very
``redundant'', in the sense that there are many ways to express a composite of
cells which will give rise to the same result, and is sometimes not very practical for
this reason. In \cref{sec:ncat-alt-syntax}, we provide an alternative syntax,
which suffers less from these problems, by restricting compositions. Finally, in
\cref{sec:ncat-wp}, we briefly mention the word problem for free $n$-categories.

\section{A Syntax for \pdfm{n}-Categories}
\label{sec:ncat-syntax}

\subsection{A syntax for free \pdfm{n}-categories}
Suppose fixed an $n$-polygraph~$P$. We define two sets of terms, as the smallest
sets closed under certain operations.
\begin{itemize}
\item An \emph{expression} is
  \[
    x
    \qqqtor
    \unit{f}
    \qqqtor
    f\Tcomp{i}g
  \]
  where $x$ is an element of $\bigsqcup_{0\leq i\leq n}P_n$, $f$ and $g$ are
  expressions and $i$ is a natural number (which can be supposed to satisfy
  $0<i\leq n$ without loss of generality).
\item A \emph{type} is either
  \[
    \TObj
    \qqqtor
    \THom Tfg
  \]
  for some type~$T$ and expressions $f$ and $g$.
\end{itemize}

The expressions should be thought of as formal composites of generators; in
particular, $f\Tcomp{i}g$ corresponds to two formal $n$-cells $f$ and $g$
composed in dimension $n-i$ (thus the minus sign in the index) and types
represent either the set of objects ($\TObj$) or a particular hom-set.

\subsection{Type-theoretic syntax}
A \emph{judgment}\index{judgment} is an expression of one of the following
forms, with the following meanings:
\begin{itemize}
\item well-formed type:
  \[
    \vdash T
  \]
  for some type~$T$,
\item well-typed term:
  \[
    \vdash f:T
  \]
  for some term $f$ and type $T$,
\item equivalent types:
  \[
    \vdash T\Teq U
  \]
  for some types $T$ and $U$,
\item equivalent terms:
  \[
    \vdash f\Teq g:T
  \]
  for some terms $f$ and $g$ and type $T$.
\end{itemize}

An \emph{inference rule}\index{inference rule} is of the form
\[
  \inferrule{\vdash\Gamma_1\\\ldots\\\vdash\Gamma_n}{\vdash\Gamma}
\]
where $\vdash\Gamma_i$ and $\vdash\Gamma$ are judgments respectively, called the
\emph{premises} and the \emph{conclusion} of the inference rule. A judgment is
\emph{derivable} when it is the conclusion of an inference rule whose premises
are all derivable. An expression~$f$ is derivable, when the judgment
$\vdash f:T$ is derivable for some type~$T$, two terms $f$ and $g$ are
equivalent if $\vdash f\Teq g:T$ is derivable for some type $T$, and two
types~$T$ and~$U$ are equivalent when $\vdash T\Teq U$ is derivable.

We are going to describe a typing system~$\TS{P}$, \ie a set of inference rules,
on the above expressions and types. The notation $\TS P$ suggests that it
depends on~$P$: as we have seen, the expressions and types depend on~$P$, but
the rules will be uniform. This will allow us to define an
$n$\nbd-cate\-gory~$C^P$, where~$C^P_i$ is the set of derivable expressions of
dimension~$i$ (defined below), compositions are given by $\Tcomp{j}$ and
identities by $\unit{}$. The main property that this construction satisfies is the
following one.

\begin{theorem}
  The $n$-category~$C^P$ is isomorphic to $\freecat{P}$.
\end{theorem}

\noindent
The result will be proved by induction on~$n$: we thus suppose that the property
is satisfied for strictly smaller values of~$n$. In particular, given a cell
$f\in\freecat P_i$ for $0\leq i<n$, up to isomorphism, we can suppose that
$\freecat P_i=C^{\tpol iP}_i$, and thus that $f$ is the equivalence class of an
expression in $\TS{\tpol iP}$. Clearly, a term in $\TS{\tpol iP}$ is a term in
$\TS P$, a valid derivation of $\TS{\tpol iP}$ is a valid derivation of~$\TS P$
and finally the equivalence relation in~$\TS{\tpol iP}$ coincides with the one
of~$\TS P$ (by \lemr{gen-type}). In particular, given a cell $x\in P_{i+1}$, we
abusively write $\src ix$ and $\tgt ix$ for an expression representing the
source and target of~$x$ in~$\freecat P_i$. This convention allows us to
associate a type~$T_x$ to each generator $x\in P_i$ by
\[
  T_x=
  \begin{cases}
    \TObj&\text{if $i=0$,}\\
    \THom{T_{\src{i-1}x}}{\src{i-1}x}{\tgt{i-1}x}&\text{otherwise.}
  \end{cases}
\]

\subsection{Inference rules}
Our typing system consists of the following inference rules:
\begin{itemize}
\item rules for types:
  \[
    \inferrule{\null}{\vdash\TObj}
    \qquad\qquad
    \inferrule{\vdash f:T\\\vdash g:T}{\vdash\THom Tfg}
  \]
\item rules for terms:
  \[
    \inferrule{\vdash f:T\\\vdash T\Teq T'}{\vdash f:T'}
  \]
  \[
    \inferrule{\null}{\vdash x:T_x}
    \qquad
    \text{for $x\in\bigsqcup_iP_i$}
  \]
  \[
    \inferrule{\vdash f:T}{\vdash\unit{f}:\THom Tff}
    \qquad
    \inferrule{\vdash f:\THom T{g}{g'}\\\vdash f':\THom{T}{g'}{g''}}{\vdash f\Tcomp{1}f':\THom{T}{g}{g''}}
  \]
  \[
    \inferrule{\vdash f:\THom T{g}{g'}\\\vdash f':\THom{T'}{h}{h'}\\\vdash g\Tcomp{i}h:U\\\vdash g'\Tcomp{i}h':U}{\vdash f\Tcomp{(i+1)}f':\THom U{g\Tcomp{i}h}{g'\Tcomp{i}h'}}
  \]
\item equivalence on types:
  \[
    \inferrule{\vdash f\Teq f':T\\\vdash g\Teq g':T\\\vdash T\simeq T'}{\vdash(\THom Tfg)\Teq(\THom{T'}{f'}{g'})}
  \]
\item equivalence on terms:
  \[
    \inferrule{\vdash f:T}{\vdash f\Teq f:T}
    \qquad
    \inferrule{\vdash f\Teq g:T}{\vdash g\Teq f:T}
    \qquad
    \inferrule{\vdash f\Teq g:T\\\vdash g\Teq h:T}{\vdash f\Teq h:T}
  \]
  \[
    \inferrule{\vdash f\Teq f':T}{\vdash\unit{f}\Teq\unit{f'}:\THom Tff}
    \qquad
    \inferrule{\vdash f\Teq f':T\\\vdash g\Teq g':U\\\vdash f\Tcomp{i} g:V}{\vdash(f\Tcomp{i} g)\Teq(f'\Tcomp{i}g'):V}
  \]
  \[
    \inferrule{\vdash\unit{g}\Tcomp{i}f:T}{\vdash(\unit{g}\Tcomp{i}f)\Teq f:T}
    \qquad
    \inferrule{\vdash f\Tcomp{i}\unit{g}:T}{\vdash(f\Tcomp{i}\unit{g})\Teq f:T}
  \]
  \[
    \inferrule{\vdash(f\Tcomp{i}g)\Tcomp{i} h:T}{\vdash (f\Tcomp{i}g)\Tcomp{i} h\Teq f\Tcomp{i}(g\Tcomp{i}h):T}
  \]
  \[
    \inferrule{\vdash\unit{f\Tcomp{i}g}:T}{\vdash\unit{f\Tcomp{i}g}\Teq\unit{f}\Tcomp{(i+1)}\unit{g}:T}
  \]
  \[
    \inferrule{\vdash(f\Tcomp{j}f')\Tcomp{i}(g\Tcomp{j}g'):T}{\vdash(f\Tcomp{j}f')\Tcomp{i}(g\Tcomp{j}g')\Teq(f\Tcomp{i}g)\Tcomp{j}(f'\Tcomp{i}g'):T}{\text{\quad for $i<j$.}}
  \]
\end{itemize}

\begin{remark}
  Note that, in the last rule, the side condition $i<j$ is important. For
  instance, consider the $2$-polygraph corresponding to the following diagram:
  \[
    \xymatrix{
      x\ar@/^4ex/[rr]^a\ar@/^2ex/@{{}{ }{}}[rr]|{\alpha\Downarrow}\ar[r]|b&y\ar@/_4ex/[rr]_e\ar@/_2ex/@{{}{ }{}}[rr]|{\beta\Downarrow}\ar[r]|c&z\ar[r]|d&w
    }
  \]
  The expression on the left is derivable, but the one on the right is not:
  \begin{align*}
    (\alpha\Tcomp1\unit{d})&\Tcomp2(\unit{b}\Tcomp1\beta)
    &
    (\alpha\Tcomp2\unit{b})&\Tcomp1(\unit{d}\Tcomp2\beta)
    \pbox.
  \end{align*}
\end{remark}

\subsection{Admissible rules}
\index{admissible rule}
\index{inference rule!admissible}
A rule is \emph{admissible} in the system when, whenever the premises are
derivable, the conclusion is also derivable (with the above rules).
\begin{lemma}
  \label{lem:eq-type}
  The following rules are admissible
  \begin{align*}
    \inferrule{\vdash f\Teq g:T}{\vdash f:T}
    &&
    \inferrule{\vdash f:T\\\vdash f:U}{\vdash T\Teq U}
    &&
    \inferrule{\vdash f:T\\\vdash T\Teq U}{\vdash f:U}    
  \end{align*}
  In particular, the type of an expression is uniquely defined up to $\Teq$,
  and two equivalent expressions have equivalent types.
\end{lemma}

\noindent
This ensures that we can meaningfully consider terms and types up to
equivalence, and moreover, by definition of the equivalence, we have:

\begin{lemma}
  The relation $\Teq$ is a congruence on expressions.
\end{lemma}

\subsection{Dimension of cells}
The \emph{dimension}\index{dimension} $\dim(T)$ of a type~$T$ is the natural
number defined by
\begin{align*}
  \dim(\TObj)&=0
  &
  \dim(\THom Tfg)&=\dim(T)+1
\end{align*}
and the dimension of a derivable expression~$f$ is the dimension of $T$ for some
derivable judgment $\vdash f:T$. It is easily shown that two equivalent types
have the same dimension and, since by \lemr{eq-type} two equivalent expressions
have equivalent types, two equivalent expressions have the same dimension, \ie
the dimension is well-defined on equivalence classes. A derivable expression of
dimension $k$ is sometimes called a \emph{$k$-expression}.

\begin{lemma}
  \label{lem:gen-type}
  Any derivable expression involving a generator $x\in P_i$ has dimension at
  least~$i$.
\end{lemma}

\subsection{Construction of the free \pdfm{n}-category}
\index{free!8-category@$n$-category}
We define the category~$C^P$ as the category such that $C^P_i$, for
$0\leq i\leq n$, is the set of equivalence classes under $\simeq$ of derivable expressions of
dimension~$i$, composition of two $i$-cells $f$ and $g$ is defined by
$f\comp{i}g=f\Tcomp{(n-i)}g$, and the identity on $f$ is $\unit{f}$.

\begin{lemma}
  The $n$-category~$C^P$ is well-defined.
\end{lemma}

\begin{theorem}
  The $n$-category~$C^P$ is isomorphic to~$\freecat P$.
\end{theorem}

\section{Alternative Syntax for \pdfm{n}-Categories}
\label{sec:ncat-alt-syntax}

\subsection{Composition in maximal codimension}
\index{composable cells}
\label{sec:other-comp}
In an $n$-category~$C$, an $i$\nbd-cell $x$ is \emph{$k$-composable} with a
$j$-cell~$y$, with $0\leq k<i\wedge j$, whenever $\tgt kx=\src ky$. In this
case, we can extend the composition operation to cells which do not necessarily
have the same dimension, and define their \emph{$k$-composite} as
\[
  x\comp{k}y
  =
  \unit{i\vee j}(x)\comp{k}\unit{i\vee j}(y)
\]
which is a $(i\vee j)$-cell. It is moreover useful to adopt the following
convention. We write $x\comp{}y$ to mean that we compose $x$ and $y$ in the
maximal possible dimension:
\[
  x\comp{}y
  =
  x\comp{(i\wedge j)-1}y
  \pbox.
\]
We say that $x$ and $y$ are \emph{composable}, when this composite is defined,
\ie
\[
  \tgt{(i\wedge j)-1}x=\src{(i\wedge j)-1}y
  \pbox.
\]
Perhaps surprisingly, all the compositions can be recovered from those
compositions in maximal codimension:

\begin{lemma}
  Given cells $x$ and $y$ of respective dimensions~$i$ and~$j$, and
  $0\leq k<i\vee j-1$, we have
  \begin{align*}
    x\comp{k}y
    &=(x\comp{}\src{k+1}y)\comp{k+1}(\tgt{k+1}x\comp{}y)
    \\
    &=(\src{k+1}x\comp{}y)\comp{k+1}(x\comp{}\tgt{k+1}y)
  \end{align*}
  which allows to compute any composition by recurrence on~$k$.
\end{lemma}
\begin{proof}
  We have
  \begin{align*}
    x\comp{k}y
    &=(x\comp{k+1}\unit{\tgt{k+1}x})\comp{k}(\unit{\src{k+1}y}\comp{k+1}y)
    &&\text{identity is neutral}
    \\
    &=(x\comp{k}\src{k+1}y)\comp{k+1}(\tgt{k+1}x\comp{k}y)
    &&\text{exchange law}
    \\
    &=(x\comp{}\src{k+1}y)\comp{k+1}(\tgt{k+1}x\comp{}y)
    &&\text{definition of $\comp{}$}
  \end{align*}
  and similarly for the second equality.
\end{proof}

\noindent
In fact, the whole structure of $n$-category can be axiomatized using this
operation~\cite[Section~8]{makkai2005wordcomp}, as follows. Given an
$(i{+}1)$-cell~$x$, we write $\src{}x$ and $\tgt{}x$ instead of $\src{i}x$ and
$\tgt{i}x$, respectively.

\begin{proposition}
  Given $n\in\N\cup\set\omega$, an $n$-category~$C$ consists of an
  $n$-globular set equipped with composition and identity partial operations as
  follows. Two cells $x\in C_i$ and $y\in C_j$, with $0\leq i,j\leq n$ are said
  to be \emph{composable} when
  \[
    \tgt{(i\wedge j)-1}(x)=\src{(i\wedge j)-1}(y)
    \pbox.
  \]
  The operations are
  \begin{itemize}
  \item \emph{compositions}: for every composable cells $x\in C_i$ and
    $y\in C_j$, there is an $(i\vee j)$-cell
    \[
      x\comp{}y
    \]
  \item \emph{identities}: for every $x\in C_i$, with $0\leq i<n$, there is an
    $(i{+}1)$-cell
    \[
      \unit{x}
    \]
  \end{itemize}
  and should satisfy
  \begin{itemize}
  \item \emph{sources and targets of compositions}: for every composable cells
    $x\in C_i$ and $y\in C_j$, with $0\leq i,j\leq n$,
    \begin{align*}
      \src{}(x\comp{}y)&=
      \begin{cases}
        \src{}x\comp{}y&\text{if $i>j$,}\\
        \src{}x&\text{if $i=j$,}\\
        x\comp{}\src{}y&\text{if $i<j$,}
      \end{cases}
      &
      \tgt{}(x\comp{}y)&=
      \begin{cases}
        \tgt{}x\comp{}y&\text{if $i>j$,}\\
        \tgt{}y&\text{if $i=j$,}\\
        x\comp{}\tgt{}y&\text{if $i<j$,}
      \end{cases}
    \end{align*}
  \item \emph{sources and targets of identities}: for every $x\in C_i$, with
    $0\leq i<n$,
    \[
      \src{}(\unit{x})=x=\tgt{}(\unit{x})
    \]
  \item \emph{associativity}: for every $x\in C_i$, $y\in C_j$ and $z\in C_k$
    with $x$ and $y$ compatible, $y$ and $z$ compatible, and either $i=j\leq k$
    or $i=k\leq j$ or $j=k\leq i$,
    \[
      (x\comp{}y)\comp{}z=x\comp{}(y\comp{}z)
    \]
  \item \emph{distributivity}: for every $x\in C_i$, $y\in C_j$ and $z\in C_k$
    \[
      x\comp{}(y\comp{}z)=(x\comp{}y)\comp{}(x\comp{}z)
    \]
    if $i<j$ and $i<k$, and $(y,z)$, $(x,y)$, and $(x,z)$ are compatible, and
    \[
      (x\comp{}y)\comp{}z=(x\comp{}z)\comp{}(y\comp{}z)
    \]
    if $i>k$ and $j>k$, and $(x,y)$, $(x,z)$, and $(y,z)$ are compatible.
  \item \emph{unitality}: for every $x\in C_i$ and $y\in C_j$ with
    $0\leq i,j\leq n$,
    \begin{align*}
      \unit{x}\comp{}y
      &=
      \begin{cases}
        y&\text{if $i+1\leq j$,}\\
        \unit{x\comp{}y}&\text{if $i+1>j$,}
      \end{cases}
      &
      x\comp{}\unit{y}
      &=
      \begin{cases}
        x&\text{if $i\geq j+1$,}\\
        \unit{x\comp{}y}&\text{if $i<j+1$,}
      \end{cases}
    \end{align*}
    whenever $\unit{x}$ and $y$ (\resp $x$ and $\unit{y}$) are composable and
    $i<n$ (\resp $j<n$),
  \item \emph{commutativity}: for every $x\in C_i$ and $y\in C_j$, with
    $\tgt{k-1}x=\src{k-1}y$, where $k=(i\wedge j)-1$ is supposed to satisfy
    $k>0$,
    \[
      (x\comp{}\src{k}(y))\comp{}(\tgt{k}(x)\comp{}y)
      =
      (\src{k}(x)\comp{}y)\comp{}(x\comp{}\tgt{k}(y))
      \pbox.
    \]
  \end{itemize}
\end{proposition}

\section{The Word Problem for Free \pdfm{n}-Categories}
\label{sec:ncat-wp}

\subsection{The word problem}
\index{word problem!for 8-polygraphs@for $n$-polygraphs}
Given an $n$-polygraph~$P$, the \emph{word problem} for~$P$ consists in finding an algorithmic
answer to the following decision problem:
\begin{center}
  Given two derivable expressions $f$ and $g$, do we have $f\Teq g$?
\end{center}

It was shown by Makkai~\cite[Section~10]{makkai2005wordcomp} that this problem
is decidable. We provide here the main arguments of his construction, after
recalling the required notions and tools.

\subsection{Multisupport}
Given two multisets\index{multiset} $\mu$ and $\nu$ with common domain (see
\cref{sec:multisets}), we write
\begin{itemize}
\item $\mu+\nu$ for their pointwise sum (noted $\mu\mcup\nu$ in
  \cref{sec:multisets}),
\item $\mu-\nu$ for their pointwise difference (we only use this operation in
  situations where $\mu(x)\geq\nu(x)$ for every element $x$ of the domain),
\item $\mu\leq\nu$ whenever $\mu(x)\leq\nu(x)$ for every element~$x$ of the
  domain.
\end{itemize}

\index{multisupport}
\index{support!multiset}
The \emph{multisupport}~$\msupp f$ of a derivable
expression~$f$ is the multiset on~$\bigsqcup_{i\leq n}P_i$ consisting of all
occurrences of generators in~$f$. Formally, $\msupp f$ is defined by induction on~$f$ as
\[
  \begin{cases}
    \set{x}&\!\!\text{if $f=x$ is a generator of type $\TObj$,}\\
    \set{x}+\msupp g+\msupp h&\!\!\text{if $f=x$ is a generator of type $\THom Tgh$,}\\
    \msupp g&\!\!\text{if $f=\unit{g}$,}\\
    \msupp g-\msupp{\Ttgt i(g)}+\msupp h&\!\!\text{if $f=g\Tcomp i h$.}    
  \end{cases}
\]
This definition is invariant under equivalence of expressions:

\begin{lemma}
  \label{lem:msupp-eq}
  Given two derivable expressions $f$ and $g$ such that $f\Teq g$, we have
  $\msupp f=\msupp g$.
\end{lemma}

\subsection{Composition in maximal codimension}
We recall the convention introduced in \secr{other-comp}. Given two derivable
expressions $f$ and $g$, not necessarily of the same dimension, whose target and
source coincide in dimension $k=(\dim(f)\wedge\dim(g))-1$, we write~$f\comp{}g$
for their composite in dimension~$k$, \ie in their maximal
codimension. Formally, we suppose that
\[
  \Ttgt{(\dim(f)-k)}(f)=\Tsrc{(\dim(g)-k)}(g)
\]
and define
\[
  f\comp{}g
  =
  \unit f^{l-\dim(f)}\Tcomp{(l-k)}\unit g^{l-\dim(g)}
\]
which is a derivable expression of dimension $l=\dim(f)\vee\dim(g)$. By
convention, $\unit f^0=f$ and $\unit f^{i+1}=\unit{\unit f^i}$, and similarly
for $g$.

\subsection{Atoms and molecules}
We mutually define notions of atom and mole\-cule by recurrence on the
dimension~$k$ as follows. Given~$k\leq n$,
\begin{itemize}
\item an \emph{atom}\index{atom} of dimension~$k$ is a derivable expression of
  the form
  \begin{equation}
    \label{eq:atom}
    f_{k-1}\comp{}(f_{k-2}\comp{}(\ldots(f_1\comp{}x\comp{}g_1)\ldots)\comp{}g_{k-2})\comp{}g_{k-1}
  \end{equation}
  where $x$ is a generator of dimension~$k$ (the \emph{nucleus} of the atom) and
  each $f_i$ and $g_i$ is a molecule of dimension~$i$,
\item a \emph{molecule}\index{molecule} of dimension~$k$ is a derivable
  expression of the form
  \[
    f_1\comp{}\ldots\comp{}f_p
  \]
  where each $f_i$ is an atom of dimension~$k$ (by convention, the composition
  is bracketed on the right and, in the case $p=0$, it should be
  an expression of the form $\unit{f}$ with $f$ an expression of dimension
  $k-1$).
\end{itemize}
In the following, we sometimes say \emph{$k$-atom} (\resp \emph{$k$-molecule})
for an atom (\resp molecule) of dimension~$k$. Clearly, every $k$-atom is a
particular $k$\nbd-mole\-cule. Note that the notion of atom corresponds to the one of
rewriting step (\cref{sec:n-rewriting-step}) and the notion of molecule to the
one of rewriting path (\cref{sec:n-rewriting-path}).

Molecules provide canonical representatives of derivable expressions up to
equivalence, in the sense that each equivalence class contains at least a
molecule, see~\cite[Proposition~8.(12)]{makkai2005wordcomp} and
\cref{prop:rewr-path-cell}:

\begin{proposition}
  \label{prop:expr-molecule}
  Every $k$-expression is equivalent to a $k$-molecule.
\end{proposition}

\begin{example}
  In the polygraph~$P$ corresponding to the diagram
  \[
    \xymatrix{
      x\ar@/^/[r]^a\ar@/_/[r]_{a'}\ar@{}[r]|{\alpha\Downarrow}&y\ar@/^/[r]^b\ar@/_/[r]_{b'}\ar@{}[r]|{\beta\Downarrow}&z
    }
  \]
  the composite $\alpha\Tcomp2\beta$ is equivalent to the molecules
  \[
    (\unit{x}\comp{}\alpha\comp{}b)\comp{}(a'\comp{}\beta\comp{}\unit{z})
    \qqtand
    (a\comp{}\beta\comp{}\unit{z})\comp{}(\unit{x}\comp{}\alpha\comp{}b')
    \pbox.
  \]
\end{example}

\subsection{Equivalence of molecules}
Two molecules are equivalent when one can be obtained from the other by applying
the exchange relation on adjacent atoms.
Formally, we define an equivalence relation $\sim$ on $k$-molecules, as the
smallest equivalence relation such that,
\begin{enumerate}
\item for every molecules $f$ and $h$, and atoms $g$ and $g'$, we have
  \[
    g\Teq g'
    \qquad\text{implies}\qquad
    f\comp{}g\comp{}h\sim f\comp{}g'\comp{}h
  \]
  whenever all composites are defined,
\item for every $k$-atoms $f:\THom T{f'}{f''}$ and $g:\THom U{g'}{g''}$, for
  every $k$-atoms $f_1,g_1,f_2,g_2$ such that
  \begin{align}
    \label{eq:molecule-xch}
    f_1&= f\comp{}g'
    &
    g_1&= f''\comp{}g
    &
    f_2&= f\comp{}g''
    &
    g_2&= f'\comp{}g
  \end{align}
  and for every $k$-molecules $h$ and $h'$, we have
  \[
    h\comp{}f_1\comp{}g_1\comp{}h'
    \sim
    h\comp{}g_2\comp{}f_2\comp{}h'
    \pbox.
  \]
  Graphically,
  \[
    \vxym{
      \\
      \ar@/^8ex/[rr]\ar@/^4ex/@{{}{ }{}}[rr]|{h\Downarrow}\ar@/_8ex/[rr]\ar@/_4ex/@{{}{ }{}}[rr]|{h'\Downarrow}\ar@/^4ex/[r]|{f'}\ar@/^2ex/@{{}{ }{}}[r]|{f\Downarrow}\ar[r]|{f''}&\ar[r]|{g'}\ar@/_4ex/[r]|{g''}\ar@/_2ex/@{{}{ }{}}[r]|{g\Downarrow}&\\
      {}
    }
    \quad\sim\quad
    \vxym{
      \\
      \ar@/^8ex/[rr]\ar@/^4ex/@{{}{ }{}}[rr]|{h\Downarrow}\ar@/_8ex/[rr]\ar@/_4ex/@{{}{ }{}}[rr]|{h'\Downarrow}\ar[r]|{f'}\ar@/_2ex/@{{}{ }{}}[r]|{f\Downarrow}\ar@/_4ex/[r]|{f''}&\ar@/^4ex/[r]|{g'}\ar[r]|{g''}\ar@/^2ex/@{{}{ }{}}[r]|{g\Downarrow}&\text{~~~~.}\\
      {}
    }
  \]
\end{enumerate}
It is shown in \cite[Section~9]{makkai2005wordcomp}:

\begin{proposition}
  \label{prop:molecule-eq}
  Given two $k$-molecules $f$ and $g$, we have $f\Teq g$ if and only if
  $f\sim g$.
\end{proposition}

\subsection{Contexts}
\index{context}
Given $k>0$ and two molecules $h$ and $h'$ with common type~$T$,
with~$\dim(h)=\dim(h')=k-1$, a \emph{context} of type $\THom Th{h'}$ is a
$k$\nbd-expression in which a formal variable of type~$\THom Th{h'}$ occurs
exactly once. Formally, we consider the polygraph $P'$ obtained from~$P$, by
adding a new generator~$\ul x$, \ie $P'_k=P_k\sqcup\set{\ul x:h\to h'}$ and
$P_i'=P_i$ for $i\neq k$. A context~$c$ is then a derivable expression in
$\TS{P'}$ such that $\msupp c(\ul x)=1$. The notion of context corresponds to
the one already introduced in \cref{sec:contexts}.

To every $(k{+}1)$-atom~$f$ of the form \eqref{eq:atom}, with top dimensional
generator $x:\THom Th{h'}$ of dimension $k+1$, one can associate the $k$-context
\[
  f_{k-1}\comp{}(f_{k-2}\comp{}(\ldots(f_1\comp{}\ul x\comp{}g_1)\ldots)\comp{}g_{k-2})\comp{}g_{k-1}
\]
where~$\ul x$ is a formal variable of dimension~$k$ and of type~$T$. We then say
that $f$ can be obtained by \emph{substituting} the variable~$\ul x$ by the
generator~$x$ in~$c$, what we write $f=c[x]$.

\begin{example}
  Consider a $3$-atom $f_2\comp{}(f_1\comp{}x\comp{}g_1)\comp{}g_2$:
  \[
    \vcenter{
      \xymatrix@C=10ex{
        \\
        x\ar@/^6ex/[rrr]\ar@/^4ex/@{{}{ }{}}[rrr]|{f_2\Downarrow}\ar@/_6ex/[rrr]\ar@/_4ex/@{{}{ }{}}[rrr]|{g_2\Downarrow}\ar[r]|{f_1}&y\ar@/^2.5ex/[r]\ar@{}[r]|{h\Downarrow\overset x{\underset{\phantom x}\TO}\Downarrow h'}\ar@/_2.5ex/[r]&z\ar[r]|{g_1}&w\text.\\
        {}
      }
    }
  \]
  The associated context is $f_2\comp{}(f_1\comp{}\ul x\comp{}g_1)\comp{}g_2$:
  \[
    \vcenter{
      \xymatrix@C=10ex{
        \\
        x\ar@/^6ex/[rrr]\ar@/^4ex/@{{}{ }{}}[rrr]|{f_2\Downarrow}\ar@/_6ex/[rrr]\ar@/_4ex/@{{}{ }{}}[rrr]|{g_2\Downarrow}\ar[r]|{f_1}&y\ar@/^2.5ex/[r]\ar@{}[r]|{\ul x\Downarrow}\ar@/_2.5ex/[r]&z\ar[r]|{g_1}&w\text.\\
        {}
      }
    }
  \]
\end{example}

Equality on atoms can be characterized in the expected way,
see~\cite[Section~10.(6)]{makkai2005wordcomp}:

\begin{proposition}
  \label{prop:atom-eq}
  Given two $(k{+}1)$-atoms $f=c[x]$ and $g=d[y]$, we have $f\Teq g$ if and only
  if $x=y$ and $c\Teq d$.
\end{proposition}

\subsection{Finiteness of equivalence classes}
The main argument in order to show that the word problem is decidable is the
observation that there is only a finite number of molecules with given
generators of given multiplicities,
see~\cite[Lemma~10.(10)]{makkai2005wordcomp}:

\begin{proposition}
  \label{prop:eq-finite}
  Suppose given a multiset~$\mu$ on~$\bigsqcup_{i\leq n}P_i$. There is a finite
  number of molecules~$f$ such that $\msupp f\leq\mu$.
\end{proposition}

\subsection{Decidability results}
As expected, the proof of decidability of the word problem is performed by
recurrence on the dimension $k$ of the two cells $f$ and~$g$ we aim at
comparing. The base case is immediate, so that we focus on the inductive
case. By \propr{expr-molecule}, whose proof is constructive, we can suppose that
$f$ and~$g$ are $k$-molecules. In \propr{molecule-eq}, we have reduced the
equivalence of the molecules $f$ to the existence of a sequence of molecules
\begin{equation}
  \label{eq:molecule-xch-seq}
  f=f_1\sim f_2\sim\ldots\sim f_p=g
\end{equation}
such that $f_{i+1}$ can be obtained from~$f_i$ by exchanging adjacent atoms,
which can be tested using equivalence of $k$-atoms for some
decompositions. Moreover, in \propr{atom-eq}, we have reduced the equivalence of
two $k$-atoms to the equivalence of $(k{-}1)$\nbd-mole\-cules (the associated
contexts). In order to show our decidability result, it is therefore enough to
show that
\begin{itemize}
\item the decompositions~\eqref{eq:molecule-xch} in order to identify the
  possible exchange relations can be found in a finite search space,
\item the sequence of molecules~\eqref{eq:molecule-xch-seq} can be found in a
  finite search space.
\end{itemize}
Both properties follow from the facts that the multisupport is preserved under
equivalence (\lemr{msupp-eq}) and there is only a finite number of molecules
with given multisupport (\propr{eq-finite}).

\begin{theorem}
  The word problem for polygraphs is decidable.
\end{theorem}

\noindent
As a consequence, it can be shown that

\begin{theorem}
  Whether an expression is derivable or not is decidable.
\end{theorem}

\noindent
An actual implementation was performed by Forest~\cite{cateq,forest2021computational}.


\chapter{Complexes and Homology}
\label{chap:tools-homology}
\label{Chapter:ComplexesAndHomology}
In this chapter, we recall basic notions on modules, abelian resolutions, and
homology. We refer the reader to~\cite{lang2012algebra,MacLane95,Rotman09} for a
deeper presentation and the proofs of the given results. Throughout the chapter,
we fix a ring $R$, and we denote by $\unit{R}$ the multiplicative identity
in~$R$.

\section{Modules Over a Ring}

\subsection{Modules over a ring}
\index{module}
A \emph{left-$R$-module} consists of an abelian group
$(V,+)$ together with a \emph{left scalar multiplication}
\[
\cdot : R \times V \to V
\]
satisfying for all $r,s\in R$ and
$x,y\in V$ the following four relations:
\begin{align*}
  r\cdot (x+y) &= r \cdot x + r \cdot y,
  &
  (rs)\cdot x &= r\cdot(s\cdot x),
  \\
  (r+s)\cdot x &= r\cdot x + s\cdot x,
  &
  1_R\cdot x &= x.
\end{align*}
In the following, we often write $rx$ instead of $r\cdot x$.
We will say \emph{$R$-module}, or \emph{module}, for a left $R$-module. 
The notion of \emph{right $R$-module} can be defined similarly, based on a right
action $\cdot:G\times R\to G$, satisfying dual axioms.
All the notions presented in this appendix are defined in a similar fashion for right $R$-modules
  since every right $R$-module is a left $R^\op$-module, where $R^\op$ is the
  opposite ring. 

\subsection{}
If $V$ and $W$ are two left-$R$-modules, a \emph{morphism of $R$-modules} 
$f : V\to W$ is a morphism of abelian groups satisfying, for any $r\in R$
and $x\in V$,
\[
  f(r\cdot x) = r\cdot f(x)\pbox.
\]
The left-$R$-modules and their morphisms form a category denoted by $\Mod_R$.
\nomenclature[ModR]{$\Mod_R$}{category of $R$-modules}
We denote by~$\hom{\Mod_R}{V}{W}$ the abelian group of morphisms from $V$ to $W$.

For a morphism $f : V \fl W$, we will denote by 
\[
\ker f := \setof{x\in V}{f(x) =0},
\] 
the \emph{kernel} of $f$ and by 
\[
\im f := \setof{y \in W}{\text{$y = f(x)$ for some $x$ in $V$}},
\]
the \emph{image} of $f$. We will denote by
$\coker f := W / \im f$
the \emph{cokernel} of $f$.

\subsection{Exact sequences}
\index{exact sequence}
A pair of composable morphisms of modules
\[
  \svxym{
    V'\ar[r]^f&V\ar[r]^g&V''
  }
\]
is \emph{exact at $V$} if $\im f = \ker g$. A sequence of morphisms of modules
\[
  \xymatrix@C=4ex{
    \cdots\ar[r]& 
    V_{n+1}\ar[r]^-{d_{n+1}}&
    V_{n}\ar[r]^-{d_n}&
    V_{n-1}\ar[r]&
    \cdots
  }
\]
is \emph{exact} if each adjacent pair $(d_{i+1},d_i)$ of morphisms is exact at $V_i$.

\begin{example}
The sequences
\[
  \xymatrix@C=3.6ex{
    0\ar[r]&V\ar[r]^f&V',
  }
\quad
  \xymatrix@C=3.6ex{
    V\ar[r]^f&V'\ar[r]&0,
  }
\quad\text{and}\quad
  \xymatrix@C=3.6ex{
    0\ar[r]&V\ar[r]^f&V'\ar[r]&0
  }
\]
are exact if and only if the morphism $f$ is injective, surjective, and
bijective, respectively.
If the sequence $V' \ofl{f} V \ofl{g} V''$ is exact with $f$ surjective and
$g$ injective, then $V$ is the zero module.  
\end{example}

\subsection{Free modules}
\index{module!free}
\index{free!module}
A $R$-module $V$ is \emph{free} if it is a direct sum of copies of~$R$. If
$V=\coprod_{i\in I} Rx_i$, with $R\isoto Rx_i$, the set $\setof{x_i}{i\in I}$ is
called a \emph{basis} of~$V$. It follows that each element $x$ in $V$ has a
unique decomposition
\[
x = \sum_{i\in I} \lambda_i x_i,
\]
where $\lambda_i$ is in $R$ and almost all $\lambda_i$ are zero.
Given a set $X$, there exists a free $R$-module having $X$ as a basis, which is
usually denoted $\freemod RX$.

\begin{lemma}
  \label{lem:free-mod-free}
  Let $X=\setof{x_i}{i\in I}$ be a basis of a free module $V$. For every module
  $W$ and every map $f : X \to W$, there is a unique morphism of $R$-modules
  $\tilde{f} : V \to W$ extending $f$, \ie making the following triangle commute:
   \[
    \xymatrix{
      X 
      \ar@{->} [r]
      \ar[dr]_{f}
      &
      V
      \ar@{.>}[d]^{\tilde{f}}
      \\
      &
      W
    }
  \]
  where the morphism $X\to V$ is the canonical inclusion.
\end{lemma}

One shows that every $R$-module is a quotient of a free $R$-module. As a
consequence, any $R$-module $V$ may be described by \emph{generators} and
\emph{relations} in the following way. Given a free $R$-module $F$ with basis
$X$ and a surjective morphism of $R$-modules
$f : F \to V$, we say that $X$ is a set of \emph{generators of~$V$} and the kernel $\ker f$ is
called its submodule of \emph{relations}.

\subsection{Finitely generated modules}
\index{module!finitely generated}
\index{finitely!generated module}
A $R$-module $V$ is \emph{finitely generated} if there is a finite subset
$X=\set{x_1, x_2, \ldots, x_n}$ of $V$ such that for every element $x$ of $V$, there exist
$r_1, r_2, \ldots, r_n$ in $R$ with $x = r_1x_1 + r_2x_2 + \ldots +
r_nx_n$. Then the set $X$ is referred to as a
\emph{generating set} for $V$. The finite generators need not form a basis,
since they need not be linearly independent over $R$. An $R$-module $V$ is
finitely generated if and only if there is a surjective morphism:
\[
R^n \to V
\]
for some $n$. That is, $V$ is a quotient of a free module of finite rank.

\begin{proposition}
  \label{Proposition:UniversalPropertyFree}
  Let $F$, $V$, and $W$ be $R$-modules. If $F$ is free, $e : V \to W$ is a
  surjective morphism and $f :F \to W$ is any morphism, then there exists a
  morphism $\tilde{f} : F \to V$ making the following triangle commute
  \[
    \xymatrix@C=2.5em @R=2.5em{
      & 
      F
      \ar[d] ^-{f}
      \ar@{.>}[dl]  _-{\tilde{f}}
      \\
      V 
      \ar[r] _-{e}
      &
      W
      \ar[r]
      &
      0.
    }
  \]
\end{proposition} 

\noindent
As a consequence, for any free $R$-module, the functor
\[
\hom{\Mod_R}{F}{-} : \Mod_R \to \Ab
\]
is \emph{exact}, that is, for any exact sequence
\[
0\fl V' \fl V \fl V'' \fl 0
\]
of $R$-modules, the induced sequence 
\[
0\fl \hom{\Mod_R}{F}{V'} \fl \hom{\Mod_R}{F}{V} \fl \hom{\Mod_R}{F}{V''} \fl 0
\] 
is exact.

\subsection{Projective module}
\index{module!projective}
\index{projective!module}
\label{sec:projective-module}
A projective module is a module which behaves as the free module $F$ in
Proposition~\ref{Proposition:UniversalPropertyFree}. More explicitly, a
$R$-module $P$ is \emph{projective} if whenever $e : V \fl W$ is a surjective
morphism and $f : P \to W$ is any morphism, there exists a morphism
$\tilde{f} : P \to V$ making the following triangle commute:
\[
\xymatrix @C=2.5em @R=2.5em{
 & 
 P
 \ar[d] ^-{f} 
 \ar@{.>}[dl] _-{\tilde{f}}
 \\
 V 
 \ar[r] _-{e}
 &
 W
 \ar[r] 
 &
 0.
}
\] 
In particular, any free module is projective.
The following result gives several ways to characterize projective modules.

\begin{proposition}
  The following conditions are equivalent for a $R$-module~$P$:
  \begin{conditions}
  \item $P$ is projective,
  \item if $f : V \fl P$ is a surjective morphism, then there exists
    $h : P \fl V$ such that $fh=\id_P$,
  \item if $f : V \fl P$ is a surjective morphism, then
    $V\simeq P\oplus \ker f$,
  \item the functor $\hom{\Mod_R}{P}{-} :  \Mod_R \to \Ab$ is exact,
  \item\label{item:Projective} $P$ is a summand of a free module, that is, there
    exists a free $R$-module $F$ such that $F\isoto P\oplus Q$, for some
    $R$-module $Q$.
  \end{conditions}
\end{proposition}

\noindent
Note that the $R$-module $Q$ in \ref{item:Projective} is necessarily
projective.

\section{Chain Complexes}
\label{S:ChainComplexes}

We recall in this section the notion of chain complex which is the fundamental object of study in homological algebra.

\subsection{Chain complex}
\index{chain complex}
\label{paragr:def_chain_complex}
\label{paragr:def_ChR}
A \emph{chain complex} in the category $\Mod_R$ is a sequence
$(V_n)_{n\in\Z}$ of $R$-modules, together with a sequence
$(d_n)_{n\in\Z}$ of morphisms
\[
\xymatrix{
\cdots \ar[r]
& V_{n+1} \ar[r] ^-{d_{n+1}} 
& V_{n} \ar[r] ^-{d_{n}}
& V_{n-1} \ar[r] 
&
\cdots
}
\]
with the composite of adjacent morphisms being zero, that is $\bnd_n\bnd_{n+1}=0$, for every $n$ in $\Z$.
Such a complex is denoted $(V_\ast,d_\ast)$, or simply $V_\ast$ or $V$. We say that $V_n$ is the module of the complex $V_\ast$ in \emph{degree} $n$.

Given two chain complexes $(V_\ast,d_\ast)$ and $(V'_\ast,d'_\ast)$, a \emph{chain map}
\[
f:(V_\ast,d_\ast) \to (V'_\ast,d'_\ast)
\]
is a family of
morphisms of $R$-modules $(f_n : V_n \to V'_n)_{n\in\Z}$ making the following diagrams
commute
\[
\xymatrix@C=3em @R=3em{
\cdots \ar[r]
& V_{n+1} \ar[r] ^-{d_{n+1}} \ar[d]^{f_{n+1}}
& V_{n} \ar[r] ^-{d_n} \ar[d]^{f_{n}}
& V_{n-1} \ar[r] \ar[d]^-{f_{n-1}}
& \cdots
\\
\cdots \ar[r]
& V'_{n+1} \ar[r] ^-{d'_{n+1}}
& V'_{n} \ar[r] ^-{d'_n}
& V'_{n-1} \ar[r]
& \cdots\text.
}
\]
Chain complexes of $R$-modules and chain maps form a category that we will denote by $\Ch[R]$.
\nomenclature[ChR]{$\Ch[R]$}{category of chain complexes of $R$-modules}

A chain complex $V$ is \emph{positive} if $V_n=0$ for every $n<0$. A positive complex looks like a sequence 
\[
\xymatrix{
\cdots \ar[r]
& V_{n} \ar[r] ^-{d_n} 
& V_{n-1} \ar[r]
& \cdots \ar[r] 
& V_2 \ar[r] ^-{d_2} 
& V_1 \ar[r] ^-{d_1}
& V_0.
}
\]

\nomenclature[ChR0]{$\pCh[R]$}{category of positive chain complexes of $R$-modules}
Positive chain complexes of $R$-modules form a full subcategory of $\Ch[R]$ denoted by $\pCh[R]$.

\subsection{Homology}
\index{homology!of a chain complex}
\label{paragr:def_funct_Hn}
Let $V$ be a chain complex in $\Mod_R$, and $n\in \Z$. 
The morphisms $\bnd_n$ are called \emph{boundary maps} of $V$, and the elements of the module $V_n$ are called \emph{$n$-chains}. 
We denote by $Z_n(V)=\ker \bnd_n$ the module of \emph{$n$-cycles} of $V$ and by $B_n(V) = \im\bnd_{n+1}$ the module of \emph{$n$-boundaries} of $V$.

In the category $\Mod_R$, the equation  $\bnd_n\bnd_{n+1}=0$ is equivalent to the condition
$B_n(V) \subseteq Z_n(V)$, that is any $n$-boundary of $V$ is an $n$-cycle of $V$.
The \emph{$n$-th homology module of $V$} is the $R$-module defined as the quotient 
\[
  H_n(V) = Z_n(V) / B_n(V).
\]
An element of $H_n(V)$ is a coset $c+B_n(V)$, and called a \emph{$n$-th homology class} of~$V$. 

Given a chain map $f : V \to V'$, and $n\in \Z$, we define a morphism of modules
\[
H_n(f) : H_n(V) \to H_n(V')
\]
by setting $H_n(f)(c+B_n(V)) = f(c) + B_n(V')$. Then 
\[
H_n : \Ch[R] \fl \Mod_R
\]
is a functor called \emph{$n$-th homology functor}. We refer to {\cite[Proposition 6.8]{Rotman09}} for a detailed proof.

\subsection{Acyclic complex}
\index{chain complex!acyclic}
\index{acyclic!chain complex}
\index{exact sequence}
A complex $(V,d)$ is an \emph{exact sequence} if $H_n(V)=0$ for every $n\in\Z$. This means that all $n$-cycles are $n$-boundaries. In this way, homology of the complex $V$ measures the deviation of $V$ from being an exact sequence. An exact sequence is also called an \emph{acyclic complex}. 

\subsection{Chain homotopy}
\index{chain homotopy}
\label{paragr:chain_homot}
\index{chain complex!homotopy}
Let $f, g : V \to V'$ be two chain maps. We say that $f$ and $g$ are \emph{(chain) homotopic}, denoted by $f\simeq g$, if there exists a family of morphisms $(s_n : V_n \fl V'_{n+1})_{n\in \Z}$ such that the following relation holds for every $n\in \Z$
\[
f_n - g_n = \bnd'_{n+1}s_n + s_{n-1}\bnd_n.
\]
\[
\xymatrix@C=4em @R=4em{
\cdots \ar[r]
& V_{n+1} \ar[r] ^-{\bnd_{n+1}} \ar@<+.7ex>[d]^{f_{n+1}} \ar@<-.7ex>[d]_{g_{n+1}}
& V_n \ar[r] ^-{\bnd_n} \ar@<+.7ex>[d]^{f_{n}} \ar@<-.7ex>[d]_{g_{n}} \ar[dl] _(0.4){s_n}
& V_{n-1} \ar[r] \ar@<+.7ex>[d]^{f_{n-1}} \ar@<-.7ex>[d]_{g_{n-1}} \ar[dl] _(0.4){s_{n-1}}
& \cdots 
\\
\cdots \ar[r]
& V'_{n+1} \ar[r] ^-{\bnd'_{n+1}} 
& V'_n \ar[r] ^-{\bnd'_n} 
& V'_{n-1} \ar[r]  
& \cdots\text.
}
\]
A chain map $f : V \fl V'$ is \emph{null-homotopic} if $f \simeq 0$, where $0$ is the zero chain map.

It is easy to see that chain homotopy defines an equivalence relation between chain maps.
Moreover, two chain homotopic maps induce isomorphisms in homology as states by the following result.

Two complexes $V$ and $V'$ are of the \emph{same homotopy type}, or \emph{homotopic}, if there exist chain maps
\[
f : V \fl V',
\quad\text{and}\quad
g : V' \fl V,
\]
such that $gf \simeq 1_V$ and $fg\simeq 1_{V'}$. In that case, the chain maps $f$ and $g$ are called \emph{homotopy equivalences}.

\begin{proposition}\label{prop:homotopic_maps}
Let $f,g : V \fl V'$ be two chain maps such that $f\simeq g$, then for every $n\in \Z$,
\[
H_n(f) = H_n(g) : H_n(V) \fl H_n(V').
\]
\end{proposition}
\begin{proof}
Let $(s_n : V_n \fl V'_{n+1})_{n\in \mathbb{Z}}$ be a homotopy between $f$ and $g$. 
Given $x \in \ker \bnd_n$, we have 
\[
(f_n-g_n)(x) = \bnd'_{n+1}s_n(x) + s_{n-1}\bnd_n(x) = \bnd'_{n+1}s_n(x).
\]
Hence, $f_n(x)-g_n(x) \in \im d'_{n+1}$ and $f_n(x)=g_n(x)$ holds in $H_n(V')$.
\end{proof}

\subsection{Contracting homotopy}
\label{SS:ContractingHomotopy}
\index{contracting homotopy}
\index{homotopy!contracting}
A chain complex $(V,d)$ has a \emph{contracting homotopy} if its identity $1_V=(1_{V_n} : V_n \to V_n)_{n\in \Z}$ is null-homotopic. That is, there is a family of morphisms $(i_n : V_n \fl V_{n+1})_{n\in \Z}$ such that the following condition holds for every $n\in \Z$
\begin{equation}
\label{E:contractingHomotopy}
1_{V_n} = \bnd_{n+1}i_n + i_{n-1}\bnd_n.
\end{equation}

\begin{proposition}
\label{prop:contr-hom-exact}
A complex $V$ having a contracting homotopy is acyclic.
\end{proposition}
\begin{proof}
Prove that $H_n(V)=0$ for every $n\in\Z$. This amounts to proving the inclusion $Z_n(V) \subseteq B_n(V)$. This is a direct consequence of relation \eqref{E:contractingHomotopy}. Indeed, for any $x$ in $Z_n(V)$, we have $x = \bnd_{n+1}i_n(x)$. This proves the acyclicity of~$V$. 
\end{proof}

\section{Resolutions}
\label{S:Resolution}

\subsection{Resolutions}
\label{SS:ResolutionModules}
\index{resolution!of a module}

A \emph{resolution} of an $R$-module~$V$ is an exact sequence of $R$\nbd-modules
\begin{equation}
  \label{eq:resolution}
  \xymatrix@C=4ex{
    \cdots \ar[r]
    & V_{n} \ar[r] ^-{d_n} 
    & V_{n-1} \ar[r]
    & \cdots \ar[r] 
    & V_1 \ar[r] ^-{d_1}
    & V_0 \ar[r] ^-{\varepsilon}
    & V \ar[r]
    & 0.
  }
\end{equation}
From the definition, the morphism $\varepsilon$ is surjective and, for every
$n\in\N$, we have
\begin{equation}
\label{eq:resolution-laws}
\im d_1 = \ker \varepsilon,
\quad\text{and}\quad
\im d_{n+1} = \ker d_n.  
\end{equation}

\index{resolution!projective}
\index{projective!resolution}
Such a resolution is \emph{projective} (\resp \emph{free}) if all the
modules~$V_n$ are projective (\resp free). Given a natural number~$n$, a
\emph{partial resolution of length~$n$ of~$V$} is defined in a similar way but
with a bounded sequence $(V_k)_{0\leq k\leq n}$ of $R$-modules:
\[
\xymatrix{
V_{n} \ar[r] ^-{d_n} 
& V_{n-1} \ar[r]
& \cdots \ar[r] 
& V_1 \ar[r] ^-{d_1}
& V_0 \ar[r] ^-{\varepsilon}
& V \ar[r]
& 0.
}
\]

Note that by \cref{prop:contr-hom-exact} a way to prove that a complex \eqref{eq:resolution} is a resolution of~$V$ is to construct a
contracting homotopy
\[
\xymatrix{
\cdots
& V_{n+1} \ar[l]
& V_{n-1} \ar[l] _-{i_{n-1}} 
& \cdots \ar[l] 
& V_1 \ar[l]
& V_0 \ar[l] _-{i_0}
& V \ar[l] _-{i_{-1}}
}
\]
such that $d_0 i_{-1}=\id_{V}$.

\begin{proposition}
  \label{Proposition:ExistenceFreeResolution}
  Every $R$-module $V$ has a free resolution.
\end{proposition}
\begin{proof}
  We take $V_0=\freemod RV$ to be the free $R$-module generated by~$V$ and the
  morphism $\varepsilon:V_0\to V$ to be the morphism extending the identity
  on~$V$ (whose existence is asserted in \cref{lem:free-mod-free}). The morphism
  $\epsilon$ is clearly surjective, \ie satisfies $\im\varepsilon=V$. 
  
  We then
  define $V_1=\freemod R{\ker\varepsilon}$ to be the free module generated
  by~$\ker\varepsilon$ and $d_1$ to be the morphism extending the canonical
  inclusion $\ker\varepsilon\to V_0$: 
 
 \[
\xymatrix{
R[\ker\varepsilon] \ar[r] ^-{d_1} 
& V_{0} \ar[r] ^-{\varepsilon}
& V.\\
\ker\varepsilon
  \ar@{^(->}[u] \ar@{^(->}[ur]
}
\] 
 
  Similarly, by recurrence, if $V_n$ and
  $d_n:V_n\to V_{n-1}$ are defined for $n\geq 1$, we set $V_{n+1}=\freemod R{\ker d_n}$ and
  $d_{n+1}:V_{n+1}\to V_n$ to be the extension of the inclusion
  $\ker d_n\to V_n$. The relations \eqref{eq:resolution-laws} are easily seen to be
  satisfied.
\end{proof}

\begin{theorem}
\label{T:HomotopyResolutions}
Let
\[
\xymatrix{
\cdots \ar[r]&
P_{n} \ar[r] ^-{d_n} 
& P_{n-1} \ar[r]
& \cdots \ar[r] 
& P_1 \ar[r] ^-{d_1}
& P_0
}
\]
be a projective complex and
\[
\xymatrix{
V_{n} \ar[r] ^-{d_n} 
& V_{n-1} \ar[r]
& \cdots \ar[r] 
& V_1 \ar[r] ^-{d_1}
& V_0 
}
\]
be an acyclic complex.
Then, for every morphism $f : H_0(P_\ast) \fl H_0(V_\ast)$, there exists a chain map $\tilde{f} : P_\ast \fl V_\ast$ inducing $f$. 
Moreover, two chain maps inducing $f$ are homotopic.
\end{theorem}

We refer to~\cite[Theorem 4.1]{HILTONSTAMMBACH} for a detailed proof of this result. 
In the same way as the proof of \cref{Proposition:ExistenceFreeResolution}, we prove that every $R$-module $V$ has a projective resolution. Indeed, there exists a projective presentation of $V:\xymatrix@C=0.8em{
0 \ar[r] 
& V_1 \ar[r]
& P_0 \ar[r] 
& V \ar[r]
& 0.
}$ Then there exists a projective presentation of $V_1:
\xymatrix@C=0.8em{
0 \ar[r] 
& V_2 \ar[r]
& P_1 \ar[r] 
& V_1 \ar[r]
& 0, 
}$ and so on.
This gives us the following projective and acyclic complex $P_\ast$:
\[
\xymatrix{
\ar[r]
\cdots
&
P_{n} \ar[r] ^-{d_n} 
& P_{n-1} \ar[r]
& \cdots \ar[r] 
& P_1 \ar[r] ^-{d_1}
& P_0 
}
\]
where the morphism $d_n$ is defined as the composite
$P_n\to V_n\to P_{n-1}$.
As $H_0(P_\ast)=V$, this is a resolution a resolution of $V$.

From \cref{T:HomotopyResolutions} we deduce the following result, see {\cite[Theorem 4.3]{HILTONSTAMMBACH}} for a detailed proof.

\begin{proposition}
\label{Proposition:HomotopyResolutions}
Two projective resolutions of a module $V$ are of the same homotopy type.
\end{proposition}

\subsection{Quasi-isomorphisms and homotopy equivalences}
\label{paragr:def_qis}
\label{paragr:def_homot_equiv}
\index{homotopy!equivalence}
\index{quasi-!isomorphism}

Let $V$ and $V'$ be two chain complexes.
A chain map $f : V \to V'$ is a \ndef{quasi-isomorphism} if it induces an
isomorphism in homology, meaning that the induced map
\[
H(f) : H(V) \to H(V') 
\]
is an isomorphism of graded modules.

A sufficient condition for a chain map $f : V \to V'$ to be a
quasi-isomorphism is to be a \ndef{homotopy equivalence} in the sense that
it admits an inverse up to homotopy, that is, a chain map $g : V' \to V$
such that $gf$ and~$\unit{V}$ are homotopic, and $fg$ and $\unit{V'}$ are
homotopic. The fact that a homotopy equivalence is a quasi-isomorphism
follows from \cref{prop:homotopic_maps}.

\subsection{Resolutions as quasi-isomorphisms}
Note that a projective resolution \eqref{eq:resolution} induces a chain map
\[
\xymatrix@C=3em @R=3em{
  \cdots \ar[r]
  & V_n \ar[r] ^-{d_n} \ar[d]
  & V_{n-1} \ar[r] \ar[d]
  & \cdots \ar[r]
  & V_1 \ar[r] ^-{d_1} \ar[d]
  & V_0 \ar[d]^{\varepsilon}
  \\
  \cdots \ar[r]
  & 0 \ar[r]
  & 0 \ar[r]  
  & \cdots \ar[r]
  & 0 \ar[r]
  & V
}
\]
between the chain complex $(V_n,d_n)$ of projective $R$-modules and the chain
complex consisting of~$V$ in degree~$0$ and the $R$-module $0$ in other
degrees.

\subsection{Schanuel's lemma and finite homological type}
\label{sec:schanuel}
\label{SchanuelLemma}
\index{Schanuel's lemma}
\index{lemma!Schanuel}

Given two exact sequences of $R$-modules
  \[
    \begin{array}{c@{\ \to\ }c@{\ \to\ }c@{\ \to\ }c@{\ \to\ }c}
      0&K_1&P_1&V&0,\\[1ex]
      0&K_2&P_2&V&0,
    \end{array}
  \]
  where $P_1$ and $P_2$ are projective, Schanuel's lemma states that there is an isomorphism:
  \[
    K_1\oplus P_2
    \qisoto
    K_2\oplus P_1
    \pbox.
  \]
We refer the reader to \cite{Rotman09} for a detailed proof and applications of this result. The following proposition is an important generalization of this result which is useful for defining the finite homology type of modules in \cref{SS:FinitenessModulesCategory}.

\begin{proposition}[Generalized Schanuel's Lemma]
  \label{Proposition:GeneralisedSchanuelLemma}
  Given two exact sequences of $R$-modules
  \[
    \begin{array}{c@{\ \to\ }c@{\ \to\ }c@{\ \to\ }c@{\ \to\ }c@{\ \to\ }c@{\ \to\ }c@{\ \to\ }c@{\ \to\ }c}
      0&K&P_k&P_{k-1}&\cdots&P_1&P_0&V&0,\\[1ex]
      0&L&Q_k&Q_{k-1}&\cdots&Q_1&Q_0&V&0,
    \end{array}
  \]
  where all the $P_i$ and $Q_i$ are projective, we write
  \[
    P_{\textnormal{od}} :=\bigoplus_{\textnormal{$i$ odd}}P_i,
    \quad
    P_{\textnormal{ev}} := \bigoplus_{\textnormal{$i$ even}}P_i,
\quad
    Q_{\textnormal{od}} := \bigoplus_{\textnormal{$i$ odd}}Q_i,
\;\;\;\text{and}\;\;\;
    Q_{\textnormal{ev}} := \bigoplus_{\textnormal{$i$ even}}Q_i
    \pbox.
  \]
  The following properties hold:
  \begin{enumerate}
  \item If $k$ is even, then $K\oplus Q_{\textnormal{ev}} \oplus P_{\textnormal{od}} \simeq L\oplus Q_{\textnormal{od}} \oplus P_{\textnormal{even}}$.
  \item If $k$ is odd then $K\oplus Q_{\textnormal{od}} \oplus P_{\textnormal{ev}} \simeq L\oplus Q_{\textnormal{ev}} \oplus P_{\textnormal{od}}$.
  \end{enumerate}
\end{proposition}

Let us mention an important consequence of \cref{Proposition:GeneralisedSchanuelLemma} which will be used to define the finite homological type in \cref{S:CategoriesFiniteHomologicalType}.
\begin{corollary}
  \label{Proposition:ConsequenceGeneralisedSchanuelLemma}
  Given two exact sequences of $R$-modules
  \[
    \begin{array}{c@{\ \to\ }c@{\ \to\ }c@{\ \to\ }c@{\ \to\ }c@{\ \to\ }c@{\ \to\ }c@{\ \to\ }c@{\ \to\ }c}
      0&K&P_k&P_{k-1}&\cdots&P_1&P_0&V&0,\\[1ex]
      0&L&Q_k&Q_{k-1}&\cdots&Q_1&Q_0&V&0,
    \end{array}
  \]
  where all the $P_i$ and $Q_i$ are finitely generated and projective, then the
  $R$-module $K$ is finitely generated if and only if $L$ is finitely generated.
\end{corollary}

\section{Homology of Monoids}
\label{sec:mon-homology}

In this section we recall the homology of a monoid with integral coefficients. In particular, the homology of monoids is used in \cref{chap:2-homology} to define homological finiteness conditions for the existence of finite convergent presentations.

\subsection{Tensor product of modules}
\index{tensor product!of modules}
\index{module!tensor product}
Given a right $R$-module~$C$ and a left $R$-module~$D$, recall that their \emph{tensor
  product} $C\otimes_RD$ is the group obtained as the quotient of the free
abelian group over the set $C\times D$, whose elements are pairs noted
$x\otimes y$ with $x\in C$ and~$y\in D$, under the relations
\begin{align*}
  (x+x')\otimes y
  &=
  x\otimes y+x'\otimes y,
  &
  (x\cdot r)\otimes y
  &=
  x\otimes (r\cdot y),
  \\
  x\otimes(y+y')
  &=
  x\otimes y+x\otimes y',
\end{align*}
for $x,x'\in C$, $y,y'\in D$, and $r\in R$. This construction is functorial:
given a right $R$-module $C$ and a morphism $f:D\to D'$ of left $R$-modules,
there is an induced group morphism
\[
  C\otimes_Rf:C\otimes_R D\to C\otimes_R D'
\]
and this construction is compatible with composition and identities.

Given a monoid~$M$, the trivial $\Z M$-module $\Z$ is canonically a right module
with the right action given by $n\cdot u=n$ for $n\in\Z$ and $u\in M$. In the
following, we will mostly be interested in computing the tensor product of the
trivial module $\Z$ with a free $\Z M$-module $\freemod{\Z M}{X}$ on a set~$X$,
in which case the result is
\[
  \Z\otimes_{\Z M}\freemod{\Z M}{X}
  =
  \freemod\Z{X}
\]
where $\freemod\Z{X}$ is the free abelian group over~$X$ (which coincides with
the free $\Z$-module).

\subsection{Homology of monoids with integral coefficients}
\index{homology!of monoids}
\index{monoid!homology}

Let $M$ be a monoid. To a free resolution of the trivial $\Z M$-module  $\Z$ by left $\Z M$-modules
\[
\xymatrix@C=1.75em{
\cdots \ar[r]
& F_{n+1} \ar[r] ^-{d_{n+1}} 
& F_n \ar[r] 
& \cdots \ar[r]
& F_1 \ar[r] ^-{d_1}
& F_0 \ar[r] ^-{\varepsilon}
& \Z, 
}
\]
by tensoring with the trivial right $\Z M$-module~$\Z$, we associate the complex of $\Z$-modules:
\[
  \!
\xymatrix@C=1.70em{
\cdots \ar[r]
& \Z\otimes_{\Z M}F_{n+1} \ar[r] ^-{\widetilde{d}_{n+1}} 
& \Z\otimes_{\Z M}F_n \ar[r] 
& \cdots \ar[r]
& \Z\otimes_{\Z M}F_1 \ar[r] ^-{\widetilde{d}_1}
& \Z\otimes_{\Z M}F_0,
}
\]
where $\widetilde{d}_k$ denotes the map $\id_{\Z}\otimes_{\Z M} d_{k}$, for all $k\geq 1$.
Note that the $\Z$-module $\Z\otimes_{\Z M}F_{n}$ is obtained from $F_n$ by trivializing the action of $M$, that is $F_n$ quotiented by all relations $ux=x$, for $u$ in $M$ and $x$ in~$F_n$.
In particular, if $F_n=\freemod{\Z M}X$, then $\Z\otimes_{\Z M}F_{n}=\freemod{\Z}X$ is the free $\Z$-module on $X$.
We obtain a chain complex, because $d_nd_{n+1}=0$ induces that 
$\widetilde{\bnd}_{n}\widetilde{\bnd}_{n+1}=0$.

We define the \emph{$n$-th homology group} of $M$ with integral coefficient $\Z$ as the following $\Z$-module: 
\[
\mathrm{H}_n(M,\Z) = \ker \widetilde{\bnd}_n / \im \widetilde{\bnd}_{n+1},
\]
with the convention that $\widetilde{d}_0=0$. By definition, for any monoid $M$, we have $\mathrm{H}_0(M,\Z)\simeq \Z$. Following results of the previous section, the homology does not depend on the choice of the resolution used to compute it.

\begin{proposition}
\label{Proposition:HomologyInvariant}
The groups $\mathrm{H}_n(M,\Z)$ do not depend on a particular choice of a free resolution, but only on the monoid $M$. 
\end{proposition}


\chapter{Homology of Categories}
\label{chap:homology_categories}

In this chapter, we recall the classical notions of simplicial homology and
homology with coefficients for a category. We also recall the property of being
\emph{finite homological type}~$\FP_n$ for a category.

\section{Simplicial Homology and Nerve of a Category}
\label{sec:simpl-H}

\begin{paragraph}[Simplicial sets]
  \label{paragr:def_simpl_sets}
  \index{simplicial!category}
  \index{category!simplicial}
  \index{simplicial!set}
  \index{simplicial!set}
  \nomenclature[.D]{$\Simpl$}{simplicial category}
  We will denote by $\Simpl$ the simplicial category, that is, the full
  subcategory of the category of posets whose objects are the
  \[
    \Ordn{n} = \set{0 < 1 < \cdots < n}
  \]
  for $n \ge 0$. By definition, a \ndef{simplicial set} is a presheaf on
  $\Simpl$, \ie a functor $\Simpl^\op \to \Set$. We will denote by
  $\SSet$ the category of simplicial sets. If $X$ is a simplicial set, the
  set $X_n = X([n])$ is called the set of \ndef{$n$-simplices} of $X$.
\end{paragraph}

\begin{paragraph}[Homology of a simplicial set]
  \label{paragr:homology_SSet}
  \index{homology!of simplicial sets}
  \index{simplicial!set!homology}
  Let $n \ge 1$. For $i$ such that $0 \le i \le n$, denote by $\delta^i$ the
  unique injective order-preserving map
  \[ \delta^i_n : [n-1] \to [n] \]
  not reaching $i$. If $X$ is a simplicial set, we will denote by
  \[ d^n_i : X_n \to X_{n-1} \]
  the map $X(\delta^i_n)$.

  We define a functor
  \[
    \nc : \SSet \to \pCh[\Z]
  \]
  in the following way. If $X$ is a simplicial set, $\nc(X_n)$ is the free
  abelian group on~$X_n$ and, if $x$ is an $n$-simplex of $X$ with $n > 0$,
  we set
  \[
    d_n([x]) = \sum_{i=0}^n (-1)^n [d^i_n(x)].
  \]
  A classical calculation shows that $\nc(X)$ is indeed a chain
  complex. If $f : X \to Y$ is a morphism of simplicial sets, the morphism
  $\nc(f) : \nc(X) \to \nc(Y)$ is defined by sending $[x]$, for $x$ an
  $n$-simplex of $X$, to $[f_n(x)]$, where $f_n$ denotes~$f_{\Ordn{n}}$.

  The \ndef{homology of a simplicial set} $X$ is the homology of the chain
  complex~$\nc(X)$.
\end{paragraph}

\begin{paragraph}[The nerve functor]
  \label{paragr:nerve}
  \index{nerve}
  \nomenclature[N]{$N$}{nerve functor}
  As the category of posets embeds into the category~$\Cat$ of small
  categories, we get a functor $i : \Simpl \to \Cat$. This functor $i$ induces
  the so-called \ndef{nerve functor}
  \[ N : \Cat \to \SSet \]
  sending a category $C$ to the simplicial set $\Cat(i(\var), C)$.
  Explicitly, a $0$-simplex of $N(C)$ is an object of $C$ and an
  $n$-simplex, for $n > 0$, is a chain of $n$ composable arrows in $C$.
\end{paragraph}

\begin{paragraph}[Homology of categories]
  \label{paragr:homology_cat}
  \index{category!homology}
  \index{homology!of categories}
  If $C$ is a category, the \ndef{homology} of $C$ is the homology of its
  nerve $N(C)$.
  In particular, if $M$ is a monoid (or even more particularly a group), by
  considering $M$ as a category with one object, we get a notion of homology
  of $M$.
\end{paragraph}

\begin{theorem}
  If $M$ is a monoid, then the homology of $M$, seen as a category with one
  object, coincides with the homology $\mathrm{H}_n(M,\Z)$ of the monoid as
  defined in~\cref{sec:mon-homology}.
\end{theorem}

\begin{proof}
  If $M$ is a monoid, then the complex $\nc(N(M))$ is canonically isomorphic
  to the complex $\Z \otimes_{\Z M} B$, where $B$ is the bar resolution of $M$,
  see~\cite[Chapter~X, Section~5]{MacLane95}, whence the result.
\end{proof}

\section{Homology of Categories with Coefficients}

\subsection{Modules of a category}
\index{module!of a category}
\nomenclature[Mod(C)]{$\Mod(C)$}{category of modules over a category~$C$}
\label{SS:ModulesCategory}
A \emph{(left) module} of a category~$C$, or \emph{(left) $C$\nbd-module}, is a functor from~$C$ to the category~$\Ab$ of abelian groups~\cite{mitchell1972rings}. The $C$-modules and the natural transformations between them form an abelian category with enough projectives denoted by $\Mod(C)$. 

We denote by $\Z C$ the \emph{free $\Z$-category} on $C$ whose $0$-cells are the ones of~$C$ and for all $0$-cells $p$ and $q$ of $C$, the hom-set $\Z C(p,q)$ is the free $\Z$-module on the hom-set $C(p,q)$. It follows immediately from the definitions that the category~$\Mod(C)$ is isomorphic to the category of additive functors from $\Z C$ to $\Ab$.  

Given a $C$-module $M$, if $x\in M(p)$ for some $0$-cell $p$ of $C$, then we say that $x$ is an \emph{element} of $M$. The \emph{left action} on $M$ is defined for a $1$-cell $f:p\to q$ and $x$ in $M(p)$, by $f\cdot x = M(f)(x)$.
A family $X=(x_i)_{i\in I}$ of elements of $M$ is called a \emph{family of generators} for $M$ if every element $x$ of $M(p)$ can be written as
\[
x = \sum_{i\in I} \; f_i \cdot x_i,
\]
for $1$-cells $f_i: p_i \to p$ in $\Z C$, where all but a finite number of $f_i$ are zero. 
This amounts to say that the natural transformation 
\[
\varphi_X : \bigoplus_{i\in I} \Z C(p_i,-) \to M,
\]
with $x_i\in M(p_i)$ and which takes $1_{p_i}$ to $x_i$ is an epimorphism in $\Mod(C)$.
The family $X$ is a \emph{basis} for $M$ if the natural transformation $\varphi_X$ is an isomorphism. 
A $C$-module $M$ is \emph{free} if it has a basis. It is \emph{finitely generated} if it has a finite set of generators. A finitely generated module $M$ is thus the cokernel of a morphism of finitely generated  free modules:
\[
\xymatrix{
 \bigoplus_{i\in I} \Z C(q_i,-)\ar[r] &  \bigoplus_{j\in J} \Z C(p_j,-) \ar[r] & M \ar[r] & 0.
 }
\]

\subsection{Natural systems}
\label{SS:NaturalSystems}
\index{category!of factorizations}
\index{factorizations}
\index{natural system}

Given a small category $C$, the \emph{category of factorizations}~\cite{BauesWirsching85,Quillen72} is the category, denoted by $\fact{C}$, whose $0$-cells are the $1$-cells of~$C$ and whose $1$-cells from $w$ to $w'$ are pairs $(u,v)$ of $1$-cells of $C$ such that the
following diagram commutes in $C$:
\[
  \xymatrix @R=2em @C=2em{   
    \cdot
    \ar [r] ^-{w} _-{}="src"
    & \cdot
    \ar [d] ^-{v}
    \\
    \cdot
    \ar [u] ^-{u}
    \ar [r] _-{w'} ^-{}="tgt"
    & \cdot 
  }
\]
The triple $(u,w,v)$ is called a \emph{factorization of $w'$}. 
Composition in the category $\fact{C}$ is defined by pasting, \ie if $(u,v):w\fl w'$ and $(u',v'):w'\to w''$ are
$1$-cells of $\fact{C}$, then the composite $(u',v')(u,v)$ is defined by the pair $(u'u,vv')$:
\[
  \xymatrix @R=2em @C=2em{
    \cdot
    \ar [r] ^-{w} _-{}="src1"
    & \cdot
    \ar [d] ^-{v}
    \\
    \cdot
    \ar [u] ^-{u}
    \ar [r] ^-{w'} ^-{}="tgt1" _-{}="src2"
    & \cdot 
    \ar [d] ^-{v'}
    \\
    \cdot
    \ar [u] ^-{u'}
    \ar [r] _-{w''} ^-{}="tgt2"
    & \cdot
  }
\]
The identity of $w$ is the pair $(1_{\src{(w)}},1_{\tgt{(w)}})$:
\[
  \xymatrix @R=2em @C=2em{   
    \cdot
    \ar [r] ^-{w} _-{}="src"
    & \cdot
    \ar [d] ^-{1_{\tgt{}(w)}}
    \\
    \cdot
    \ar [u] ^-{1_{\src{}(w)}}
    \ar [r] _-{w} ^-{}="tgt"
    & \cdot 
  }
\]
An $\fact{C}$-module $D$ is called a \emph{natural system of abelian groups} on $C$. 
We will denote by $D_w$ the abelian group which is the
image of $w$ by the functor $D$. If there is no confusion, we denote by $uav$ the image of
an element $a$ of $D_w$ through the homomorphism of groups
$D(u,v) : D_w \to D_{w'}$. 
The category~$\Mod(\fact{C})$ is also denoted by~$\NatSys{C,\Ab}$.
\nomenclature[Nat(C,Ab)]{$\NatSys{C,\Ab}$}{category of natural systems}

\subsection{Homology of categories with coefficients}
\label{SS:HomologyWithCoefficients}
\index{homology!of categories!with coefficients}

The cohomology of categories with values in natural systems was defined by Baues and Wirsching in~\cite{BauesWirsching85}, and by Wells in~\cite{Wells79}. It generalizes Hochschild-Mitchell cohomology of categories with coefficients in bimodules \cite{mitchell1972rings}, the cohomology with coefficient in left modules \cite{Quillen72,Roos61}, and the cohomology of the classifying space, see~\cite{BauesWirsching85} for the correspondences.

Dually, we define the homology of a category $C$ with values in a \emph{contravariant natural system} $D$ on $C$, that is an $(FC)^\op$-module, as follows.
We consider the nerve $N(C)$ of $C$ defined in \cref{paragr:nerve}, with boundary maps denoted by $d_i:N_n(C)\fl N_{n-1}(C)$, for $0\leq i\leq n$.
For $s=(u_1,\dots,u_n)$ in $N_n(C)$, we denote by $\cl{s}$ the composite $1$-cell $u_1\cdots u_n$ of $C$. 
For $n$ in $\mathbb{N}$, the $n$-th chain group $C_n(C,D)$ is defined as the abelian group   
\[
C_n(C,D) \:=\: \bigoplus_{s\in N_n(C)} D_{\cl{s}}.
\]
We denote by $\iota_s$ the embedding of $D_{\cl{s}} \hookrightarrow C_n(C,D)$. 
We define a boundary map $d : C_n(C,D) \fl C_{n-1}(C,D)$ on the component $D_{\cl{s}}$ of $C_n(C,D)$ by setting
\[
d\iota_s 
	\:=\: \iota_{d_0(s)} u_{1*} 
		\:+\: \sum_{i=1}^{n-1} (-1)^i \iota_{d_i(s)} 
		\:+\: (-1)^n \iota_{d_n(s)} u_n^*\,, 
\] 
for $s=(u_1,\dots,u_n)$ in $N_n(C)$ and where $u_{1*}$ and $u_n^*$ denote the morphisms $D(u_1,1)$ and $D(1,u_n)$, respectively. The homology of the category $C$ with coefficients in $D$ is defined as the homology of the complex $(C_*(C,D),d_*)$:
\[
\mathrm{H}_*(C,D) \:=\: \mathrm{H}_*(C_*(C,D),d_*).
\]
Denoting $\mathrm{Tor}_*^{FC}(D,-)$ the left derived functor of the functor
$D\otimes_{FC}-$, we prove that there is an isomorphism
\[
\mathrm{H}_*(C,D) \:\simeq\: \mathrm{Tor}_*^{FC}(D,\mathbb{Z}).
\]
natural in $D$.

\section{Categories of Finite Homological Type}
\label{S:CategoriesFiniteHomologicalType}
\index{finite homological type}

\subsection{Modules of finite homological type}
\label{SS:FinitenessModulesCategory}

For~$n\in\Nbinfty$, a $C$-module~$M$ is of \emph{homological type~$\FP_n$} (where $\FP_n$ stands for finitely $n$-presented) if it admits a partial resolution of length~$n$ by finitely generated projective $C$-modules
\[
\xymatrix{
\:\cdots\:
	\ar [r]
& P_n
     \ar[r]
& P_{n-1}
     \ar [r]
& \:\cdots\:
	\ar [r]	
& P_0
	\ar [r]
& M
	\ar [r]
& 0.
}
\]
By general homological arguments, see~\cref{sec:schanuel}, we have the following characterization of the property~$\FP_n$ which is a consequence of 
\cref{Proposition:ConsequenceGeneralisedSchanuelLemma}.

\begin{proposition}
\label{L:FPn}
Let~$C$ be a category, $M$ be a $C$-module, and $n$ be a natural number. The following assertions are equivalent:  
\begin{enumerate}
\item $\!M$ is of homological type~$\FP_n$.
\item $\!M$ admits a partial resolution of length~$n$ by finitely generated free $C$\nbd-modules:
 \[
\xymatrix{
F_n \ar[r] 
& F_{n-1} \ar[r] 
& \;\ldots\; \ar[r] 
& F_0 \ar[r] 
& M.
}
\]
\item $\!M$ is finitely generated and, for every natural number~$k<n$ and every projective finitely generated partial resolution of~$M$ of length~$k$
\[
\xymatrix{
P_k
     \ar[r] ^-{d_k}
& P_{k-1}
     \ar [r]
& \:\cdots\:
	\ar [r]	
& P_0
	\ar [r]
& M
	\ar [r]
& 0,
}
\]
the $C$-module~$\ker d_k$ is finitely generated.  
\end{enumerate}
\end{proposition}

\subsection{Categories of finite homological type}
\label{SS:CategoriesFiniteHomologicalType}
The property for a category~$C$ to be of homological type~$\FP_n$ is defined according to a category of modules over one of the categories in the following diagram 
\[
\xymatrix @!C @R=1em { 
&& {\opp{C}} 
	\ar@{>->} @/^/ ^-{q_1} [dr]
\\ 
{\fact{C}}
	\ar@{->>} [r] ^-{\partial} 
& {\opp{C} \times C} 
	\ar@{->>} @/^/ [ur] ^-{p_1}
	\ar@{->>} @/_/ [dr] _-{p_2}
&& {\freegpd{C}}
\\
&& {C} 
	\ar@{>->} @/_/ [ur] _-{q_2}
}
\]
where~$\partial$ is the boundary map, $p_1$ and~$p_2$ are the projections of the cartesian product, $\freegpd{C}$ is the enveloping groupoid of~$C$, and~$q_1$ and~$q_2$ are the canonical inclusion morphisms. A category~$C$ is of \emph{homological type}
\begin{itemize}
\item \emph{$\FP_n$} if the constant natural system~$\Z$ is of type~$\FP_n$,
\item \emph{bi-$\FP_n$} if the $\opp{C}\times C$-module~$\Z C$ is of type~$\FP_n$,  
\item \emph{left-$\FP_n$} if the constant $C$-module~$\Z$ is of type~$\FP_n$,  
\item \emph{right-$\FP_n$} if the constant $\opp{C}$-module~$\Z$ is of type~$\FP_n$,  
\item \emph{top-$\FP_n$} if the constant $\freegpd{C}$-module~$\Z$ is of type~$\FP_n$. 
\end{itemize}

Using the fact that the property~$\FP_n$ is preserved by left Kan extensions~\cite[Lemma~5.1.4]{GuiraudMalbos12advances}, these finiteness homological properties of categories are related by the following implications~\cite[Proposition~5.2.4]{GuiraudMalbos12advances}:
\[
\xymatrix @!C @C=1em @R=1em { 
&& {\text{right-$\FP_n$}} 
	\ar@2 @/^/ [dr]
\\ 
{\FP_n}
	\ar@2 [r] 
& {\text{bi-$\FP_n$}} 
	\ar@2 @/^/ [ur] 
	\ar@2 @/_/ [dr] 
&& {\text{top-$\FP_n$.}}
\\
&& {\text{left-$\FP_n$}} 
	\ar@2 @/_/ [ur]
}
\]

If~$C$ is a groupoid, all of these implications are equivalences~\cite[Proposition~5.2.6]{GuiraudMalbos12advances}, but this is not the case in general. Indeed, Cohen constructed in~\cite{cohen1992monoid} a right-$\FP_{\infty}$ monoid which is not left-$\FP_{1}$: thus, top-$\FP_n$, left-$\FP_n$, and right-$\FP_n$ are not equivalent in general. Moreover, monoids with a finite convergent presentation are left-$\FP_{\infty}$ and right-$\FP_{\infty}$, see~\cite{anick1986homology,kobayashi1990complete,squier1987word}, but there exists a finitely presented monoid that is left-$\FP_{\infty}$ and right-$\FP_{\infty}$, but does not satisfy the homological finiteness condition~$\mathrm{FHT}$ introduced by Pride and Wang~\cite{kobayashi2001homotopical,pride00second}; since the properties~$\mathrm{FHT}$ and bi-$\FP_3$ are equivalent~\cite{KobayashiOtto03}, it follows that left-$\FP_n$ and right-$\FP_n$ do not imply bi-$\FP_n$ in general. 

Finally, we note the following consequence of the definition.

\begin{proposition}
If a category $C$ is of homological type $\FP_n$, for a natural number $n$, then the abelian group $\mathrm{H}_k(C,\mathbb{Z})$ is finitely generated for every $0\leq k\leq n$. 
\end{proposition}

\subsection{Monoids of finite homological type}
\index{finite homological type}

The notion of finite homological type for categories applies to monoids seen as categories with one object. In particular, the proofs to show that a monoid is of homological type left-$\FP_n$ in \cref{chap:HomologieSquierTheorem} are based on the following result, which is an immediate consequence of \cref{L:FPn}.

\begin{proposition}
  \label{lemme_pl_n}
  Let~$n$ be a natural number. A monoid~$ M$ is of homological type left-$\FP_n$ if and only if the trivial left $\Z M$-module~$\Z$ is of homological type left-$\FP_n$.
\end{proposition}

\chapter{Locally Presentable Categories}
\label{chap:loc_pres}

This appendix is a quick introduction to locally presentable
categories. We refer the reader to the classical book~\cite{adamek1994localp} for
a detailed presentation.

The notion of a locally presentable category is, in some sense, a
formalization of what is an algebraic structure. When category theory is
restricted to locally presentable categories, many things get simpler. In
particular, there are characterizations of adjoint functors purely in terms
of preservation of limits and colimits. Locally presentable categories also
play an important role in the theory of model categories through the
concept of combinatorial model categories.

There are many ways to define locally presentable categories. We start with
the presentation in terms of sketches, which are somehow categories encoding
the syntax of an algebraic structure. These sketches are used several times
in the body of the book. We then give the intrinsic categorical
characterization, defining on our way several notions that will be needed
for the theory of model categories. Finally, we give the syntactic
characterization.

\section{Sketches}
\label{sec:sketches}

\begin{paragraph}[Cones]
  \index{cone!inductive}
  \index{cone!projective}
  \index{cone!limit}
Let $\C$ be a category. By a \ndef{projective cone} in $\C$, we will mean a
triple $(F, X, \alpha)$, where $F : D \to \C$ is a functor from some
\emph{small} category $D$, $X$ is an object of $\C$, and $\alpha : X \tod F$
is a natural transformation, where $X$ is considered as a constant functor.
The diagram $F$ is then called the \ndef{base} of the cone and the object
$X$ the \ndef{tip} of the cone. We say that $\C$ is a \ndef{projective limit
cone} if the morphism $X \to \varprojlim F$ induced by $\alpha$ is an
isomorphism.

The notion of \ndef{inductive cone} is defined similarly by reversing the direction
of the natural transformation $\alpha$, that is, by considering a natural
transformation $\alpha : F \tod X$, and the notion of \ndef{inductive limit cone}
by asking for the induced morphism $\limind F \to X$ to be an isomorphism.
\end{paragraph}

\begin{paragraph}[Sketches]
  \index{sketch}
  \index{sketch!projective}
  \index{sketch!inductive}
  \index{model (of a sketch)}
A \ndef{sketch} is a triple $(S, P, I)$, where $S$ is a
\emph{small} category, $P$ a \emph{set} of projective cones in~$S$ and
$I$ a \emph{set} of inductive cones in~$S$. A sketch is
\ndef{projective} (\resp \ndef{injective}) if $I=\emptyset$ (\resp
$P=\emptyset$). By abuse of notation, we will often refer to the sketch $(S,
P, I)$ as $S$ only, and we will talk of the projective cones of $S$ (\resp
of the inductive cones of $S$) for the elements of $P$ (\resp of $I$).

A \ndef{model} of a sketch~$S$ in a category~$\C$ consists of a functor $F:
S\to\C$ sending every projective cone of $S$ to a projective limit cone in
$\C$ and every inductive cone of $S$ to an inductive limit cone in $\C$. A
\ndef{morphism of models of a sketch $S$} is just a natural transformation.
We write $\Models(S)$ for the category of models of a sketch $S$ in the
category of sets.
\nomenclature[Mod(S)]{$\Models(S)$}{category of models of a sketch~$S$}

A \ndef{morphism of sketches} from a sketch $S$ to a sketch
$S'$ is a functor $f : S \to S'$ sending every projective cone of $S$ to a
projective cone of $S'$ and every inductive cone of $S$ to an inductive cone
of $S'$. If $f$ is such a morphism, then precomposition by $f$ induces a
functor
  \[
    f^*
    :
    \Models(S')
    \to
    \Models(S)
    \pbox.
  \]
\end{paragraph}

\begin{paragraph}[Sketchable categories]
  \index{category!sketchable}
We say that a category $\C$ is \ndef{sketchable} (\resp \ndef{projectively
sketchable}) if there exists a sketch (\resp a projective sketch) $S$ such
that $\C$ is equivalent to $\Models(S)$.
\end{paragraph}

\interbreak

In this section, we will focus on projectively sketchable categories.

\begin{example}\label{exam:sketch_graph}
  The category of graphs is projectively sketchable. Indeed, this category
  is nothing but the category of functors from the category
\[
\xymatrix{
  [0] & \ar@<-.5ex>[l]_s \ar@<.5ex>[l]^t [1]
}
\]
to the category of sets, \ie the category of models of a sketch without any
projective or inductive cones. More generally, any diagram category is
projectively sketchable.
\end{example}

\begin{example}
  \newcommand\p[1]{\langle #1\rangle}
  \newcommand\s[1]{\vert #1\vert}
  The category of monoids is projectively sketchable. Indeed, consider the
  category $S$ with three objects
  \[ [0], \quad [1], \quad [2], \quad [3] \]
  and generated by the morphisms
  \[
    \xymatrix@C=1pc@R=1.5pc{
      & [2] \ar[dl]_{p_1} \ar[dr]^{p_2} \\
      [1] & & [1] \\
    }
    \qquad
    \qquad
    \xymatrix@C=1pc@R=1.5pc{
      & & [3] \ar[dll]_{q_1} \ar[d]^{q_2} \ar[drr]^{q_3} \\
      [1] & & [1] & & [1] \\
    }
  \]
  \[ \xymatrix{[3] \ar[r]^{q_{1,2}} & [2]} \qquad\qquad \xymatrix{[3]
  \ar[r]^{q_{2,3}} & [2]} \]
  \[ \xymatrix{[0] \ar[r]^u & [1]} \qquad\qquad \xymatrix{[2] \ar[r]^m & [1]} \]
  \[
  \xymatrix{[1] \ar[r]^{\p{u, \id}} & [2]}
  \qquad\qquad
  \xymatrix{[1] \ar[r]^{\p{\id, u}} & [2]}
  \]
  \[
  \xymatrix{[3] \ar[r]^{\p{m, \id}} & [2]}
  \qquad\qquad
  \xymatrix{[3] \ar[r]^{\p{\id, m}} & [2]}
  \]
  subject to the relations
  \[ 
    p_1q_{1,2} = q_1,
    \quad
    p_2q_{1,2} = q_2,
    \quad
    p_1q_{2,3} = q_2,
    \quad
    p_2q_{2,3} = q_3,
  \]
  \[
    p_1\p{u, \id} = u,
    \quad
    p_2\p{u, \id} = \id_{[1]},
    \quad
    p_1\p{\id, u} = \id_{[1]},
    \quad
    p_2\p{u, \id} = u,
  \]
  \[
    p_1\p{m, \id} = m q_{1,2}
    \quad
    p_2\p{m, \id} = q_3
    \quad
    p_1\p{\id, m} = q_1,
    \quad
    p_2\p{u, m} = m q_{2,3},
  \]
  and
  \[
  \xymatrix{
    [1] \ar@{=}[dr] \ar[r]^{\p{u, \id}} & [2] \ar[d]^m & [1] \ar[l]_{\p{\id,
    u}} \ar@{=}[dl] \\
    & [1] & \pbox{,}
  }
  \qquad
  \qquad
  \xymatrix{
    [3] \ar[r]^{\p{m, \id}} \ar[d]_{\p{\id,m}} & [2] \ar[d]^m \\
    [2] \ar[r]_m & [1] \pbox{.}
  }
  \]
  Endow $S$ with the three projective cones
  \[
    \xymatrix{ [0] }
    \qquad
    \qquad
    \xymatrix@C=1pc@R=1.5pc{
      & [2] \ar[dl]_{p_1} \ar[dr]^{p_2} \\
      [1] & & [1] \\
    }
    \qquad
    \qquad
    \xymatrix@C=1pc@R=1.5pc{
      & & [3] \ar[dll]_{q_1} \ar[d]^{q_2} \ar[drr]^{q_3} \\
      [1] & & [1] & & [1] \pbox{,} \\
    }
  \]
  the first one being indexed by the empty category.
  Then we claim that the data of a model $F : S \to \Set$ of the
  sketch $S$ in $\Set$ is equivalent to the data of a monoid of
  underlying set $M = F([1])$, multiplication $F(m)$ and unit given
  by~$F(u)$ (the value of $F([i])$, for $0 \le i \le 3$, being forced to be
  sent to $M^i$).

  More generally, any Lawvere theory defines a sketch whose models are the
  models of the starting Lawvere theory.
\end{example}

\begin{example}\label{exam:sketch_cat}
  The category $\Cat$ of small categories is projectively sketchable.
  Consider the full subcategory~$S\subseteq\Simpl^\op$ of the opposite
  category of the simplicial category (see~\cref{paragr:def_simpl_sets}) on
  the objects $[0]$, $[1]$, $[2]$, and $[3]$. If $[m] \hookto [n]$ is an
  inclusion of $\Simpl$ whose
  image avoids exactly $i_0, \dots, i_k$, the corresponding morphism in
  $\Simpl^\op$ will be denoted by $d_{i_0, \dots, i_k} : [n] \to [m]$. With
  this notation, the commutative diagrams
  \[
    \xymatrix@C=1pc@R=1.5pc{
      & [2] \ar[dl]_{d_2} \ar[dr]^{d_0} \\
      [1] \ar[dr]_{d_0} & & [1] \ar[dl]^{d_1} \\
      & [0] & \\
    }
    \qquad
    \qquad
    \xymatrix@C=1pc@R=1.5pc{
      & & [3] \ar[dll]_{d_{2,3}} \ar[d]^{d_{0,3}} \ar[drr]^{d_{0,1}} \\
      [1] \ar[dr]_{d_0} & & [1] \ar[dl]^{d_1} \ar[dr]_{d_0} & & [1]
      \ar[dl]^{d_1} \\
      & [0] & & [0] & \\
    }
  \]
  define two cones (which are actually limit cones) in $S$. One can show
  that the category of models of the sketch defined by $S$ and these two
  cones is equivalent to $\Cat$.
\end{example}

\begin{example}
  The category $\oCat$ of small \oo-categories is projectively sketchable
  (see Proposition~\ref{prop:ocatlimsketch}).
\end{example}

\begin{theorem}\label{thm:sketch-morphism}
  If $f : S \to S'$ is a morphism of projective sketches, then the
  functor
  \[
    f^*
    :
    \Models(S')
    \to
    \Models(S)
    \pbox.
  \]
  admits a left adjoint
  \[
    f_!
    :
    \Models(S)
    \to
    \Models(S')
    \pbox.
  \]
\end{theorem}

\begin{proof}
  See for instance \cite[Section~4, Theorem~4.1]{barr1985toposes}.
\end{proof}

\begin{example}
  There is an obvious inclusion from the sketch of graphs defined in
  Example~\ref{exam:sketch_graph} into the sketch of categories defined
  in Example~\ref{exam:sketch_cat}. This morphism of sketches induces the
  forgetful functor from small categories to graphs. The previous theorem
  shows the well-known fact that this forgetful functor admits a left
  adjoint, sending a graph to the free category on this graph.
\end{example}

\begin{proposition}\label{prop:mod_bicompl}
  A projectively sketchable category is complete and cocomplete.
\end{proposition}

\begin{proof}
  The completeness of such a category can be checked directly (limits are
  computed as in presheaves). Proving the cocompleteness is more involved,
  see~\cite[Example 3.11.8 and Theorem 1.38]{adamek1994localp}.
\end{proof}

\goodbreak

We end the section with a very powerful criterion due to Lair to prove that a
functor is monadic in terms of sketches.

\begin{theorem}\label{thm:Lair}
  Let $f : S \to S'$ be a morphism of projective sketches satisfying the two
  following conditions:
\begin{itemize}
  \item the base of any cone of $S'$ factors though $f$,
  \item every object of $S'$ not in the image of $f$ is the tip of a cone of
  $S'$.
\end{itemize}
  Then the induced functor $f^\ast : \Models(S') \to \Models(S)$ is monadic.
\end{theorem}

\begin{proof}
  By Theorem~\ref{thm:sketch-morphism}, the functor $f^\ast$ admits a
  right adjoint and the result is thus a particular case of \cite[Corollary
  1]{LairMonad}.
\end{proof}

\section{Locally Presentable Categories}
\label{sec:def_kappa-small}

\index{regular cardinal}
Fix $\kappa$ be a regular cardinal. Recall that the fact that $\kappa$ is
\ndef{regular} means that the category of sets of cardinal $< \kappa$ is closed
under colimits of cardinality $< \kappa$.

\begin{paragraph}[$\kappa$-filtered diagrams]
  \index{filtered diagram}
A small category $I$ is said to be \ndef{$\kappa$-small} if the cardinal of
its set of morphisms is strictly smaller than~$I$. An $\aleph_0$-small
category, $\aleph_0$ being the cardinal of countable sets, is called a
\ndef{finite category}. A diagram $F : I \to \C$ is said to be
$\kappa$-small if $I$ is $\kappa$-small. A \ndef{$\kappa$-filtered category}
is a category $I$ such that for every $\kappa$-small diagram $F : J \to I$,
there exists an object $i$ of $I$ and cone $F \tod i$. An
$\aleph_0$-filtered category is said to be \ndef{filtered}.
A \ndef{$\kappa$-filtered diagram} is a diagram
$F : I \to \C$ such that $I$ is $\kappa$-filtered.
\end{paragraph}

\begin{paragraph}[$\kappa$-presentable objects]
  \index{presentable object}
An object $X$ of a category~$\C$ is said to be
\ndef{$\kappa$\nbd-presentable} if the functor
\[ \C(X,-):\C\to\Set \]
preserves $\kappa$-filtered colimits, that is, if for every
$\kappa$-filtered diagram $F:I\to\C$ the canonical map
\[
  \colim_i\C(X,Fi)
  \to
  \C(X,\colim_iFi)
\]
is a bijection. An $\aleph_0$-presentable object is called a \ndef{finitely
presentable object}.
\end{paragraph}

\begin{paragraph}[Locally presentable categories]\label{paragr:def_lpc}
  \index{locally presentable category}
  \index{category!locally presentable}
A category~$\C$ is said to be \ndef{locally $\kappa$\nbd-presentable} if it is cocomplete
and if there exists a \emph{set} $S$ of $\kappa$-small objects of~$C$ such
that every object of~$\C$ can be obtained as a $\kappa$-filtered colimit of
objects in~$S$. A locally $\aleph_0$-presentable category is said to be \ndef{locally
finitely presentable}. Essentially by definition, if a category $\C$ is
locally $\kappa$-presentable, then it is locally $\kappa'$-presentable for
every regular cardinal~$\kappa' > \kappa$. A category~$\C$ is \ndef{locally
presentable} if it is locally $\kappa$-presentable for some regular cardinal
$\kappa$.
\end{paragraph}

\begin{theorem}
  A category is locally presentable if and only if it is equivalent to the
  category of models of some projective sketch. More precisely, if $\kappa$
  is a regular cardinal, a category is locally $\kappa$-presentable if and
  only it is the category of model of some projective sketch whose cones are
  over $\kappa$-small diagrams.
\end{theorem}

\begin{proof}
  See for instance~\cite[Corollary~1.52]{adamek1994localp}.
\end{proof}

\section{Essentially Algebraic Theories}
\label{sec:ess-alg-th}

\begin{paragraph}[Algebraic theory]
  \index{theory!algebraic}
  \index{algebraic theory}
  An \ndef{algebraic theory}~$P$ consists of
  \begin{itemize}
  \item a set~$P_0$ of \ndef{sorts},
  \item a set~$P_1$ of \ndef{operations} together with functions
    \begin{align*}
      \sce0:P_1\to\freecat{P_0}
      &&
      \tge0:P_1\to P_0      
    \end{align*}
    associating to an operation its \ndef{arity} and \ndef{coarity}, where
    $\freecat{P_0}$ denotes the free monoid on~$P_0$,
  \item a set~$P_2\subseteq\freecat{P_1}\times\freecat{P_1}$ of \ndef{relations}
    consisting of pairs of terms with the same arity and coarity.
  \end{itemize}
  An element $(u,v)\in\freecat{P_1}\times\freecat{P_1}$ is often written $u=v$.
  Note that this notion corresponds to the one of term rewriting system, as
  introduced in \cref{sec:trs}.
\end{paragraph}

\begin{paragraph}[Essentially algebraic theory]
  \index{theory!essentially algebraic}
  \index{essentially algebraic theory}
An \ndef{essentially algebraic theory} consists of an algebraic theory~$P$
together with a function
\[
  \dom
  :
  P_1
  \to
  \powerset(\freecat{P_1}\times\freecat{P_1})
\]
which to every operation in $a\in P_1$ associates a subset $\dom(a)$ of
$\freecat{P_1}\times\freecat{P_1}$, called the \ndef{domain} of~$a$, specifying
the relations under which the operation is defined. An operation $a\in P_1$ is
\ndef{total} when $\dom(a)=\emptyset$, and a term is total when it is composed of
total operations only. An essentially algebraic theory is required to satisfy the
following condition: for every operation $a\in P_1$ and relation
$(u,v)\in\dom(a)$, the terms~$u$ and~$v$ are total.
\end{paragraph}

\begin{paragraph}[Model of an essentially algebraic theory]
  \index{model!of an essentially algebraic theory}
A \ndef{model}~$M$ of an essentially algebraic theory consists of
\begin{itemize}
\item a set $M_s$ for every sort $s\in P_0$,
\item a partial function
  \[
    M_a
    :
    M_{s_1}\times\ldots\times M_{s_k}
    \to
    M_s
  \]
  for every operation
  \[
    a
    :
    s_1\ldots s_k
    \to
    s
  \]
  such that $M_a(m_1,\ldots,m_k)$ is defined if and only if
  \[
    M_u(m_1,\ldots,m_k)
    =
    M_v(m_1,\ldots,m_k)
  \]
  for every $(u,v)\in\dom(a)$.
\end{itemize}
In the definition above, the interpretation~$M_u$ for a term~$u$ is defined
by induction on~$u$ from the interpretation of operations in the expected
way.
\end{paragraph}

\begin{example}
  The essentially algebraic theory~$P$ of categories has two sorts $s_0$ and
  $s_1$ in $P_0$, and four operations in~$P_1$
  \begin{align*}
    s&:s_1\to s_0
    &
    t&:s_1\to s_0
    &
    c&:s_1,s_1\to s_1
    &
    e:s_0\to s_1    
  \end{align*}
  with domains
  \begin{align*}
    \dom(s)&=\emptyset
    &
    \dom(t)&=\emptyset
    &
    \dom(c)&=\set{t(x_1)=s(x_2)}
    &
    \dom(e)&=\emptyset    
  \end{align*}
  and relations
  \begin{align*}
    s\circ e(x_1)&=x_1&
    s\circ c(x_1,x_2)&=s(x_1)&
    c(e(x_1),x_2)&=x_2
    \\
    t\circ e(x_1)&=x_1&
    t\circ c(x_1,x_2)&=t(x_2)&
    c(x_1,e(x_2))&=x_1&
  \end{align*}
  and
  \[
    c(c(x_1,x_2),x_3)=c(x_1,c(x_2,x_3))
    \pbox.
  \]
  Note that the terms $t(x_1)$ and $s(x_2)$ occurring in the domain of~$c$ are
  total as required. The models of this theory are precisely the small
  categories.
\end{example}

\begin{theorem}
  A category is locally finitely presentable if and only if it is equivalent to
  the category of models of an essentially algebraic theory.
\end{theorem}

\begin{proof}
  See~\cite[Theorem~3.36]{adamek1994localp}.
\end{proof}

\begin{remark}
The above result can be generalized to locally $\kappa$-presentable categories
by considering theories~$P$ with operations~$a$ with an arity of the form
$\sce0(a)\in P_0^{\kappa_a}$ for some cardinal $\kappa_a<\kappa$ (we recover the
above case when all the $\kappa_a$ are finite).
\end{remark}


\chapter{Model Categories}
\label{chap:model-cat}

One of the goals of this book is to construct the so-called ``folk'' model
category structure on the category of strict $\omega$-categories. This is
achieved in \cref{chap:folk}. The notion of model category, introduced by
Quillen~\cite{quillen1967homotopical}, constitutes a very
general framework in which to study the homotopical properties of a
category endowed with a class of weak equivalences.

\Cref{chap:resolutions,chap:w-eq,chap:folk} introduce
along the way the main definitions of this theory. Nevertheless, for the
convenience of the reader, we gather in this appendix these main definitions
plus some complements. For more details, we refer the reader to
the classical books~\cite{HirMC, hovey2007model, quillen1967homotopical} or
to~\cite{riehl2014categorical} for a recent panorama on the subject.

\section{Definition}
 
We start with some preliminary definitions.

\begin{paragraph}[2-out-of-3 property]
  \index{2-out-of-3 property}
A class of maps $\clW$ in a category $\C$ is said to satisfy the
\ndef{2-out-of-3 property} if for any commutative triangle
\[
\xymatrix{
  X \ar[rr]^f \ar[dr]_h & & Y \ar[dl]^g \\
                    & Z
}
\]
in $\C$, if two morphisms among $f$, $g$, and $h$ are in $\clW$, then so is
the third one.
For instance, isomorphisms in a category satisfy the 2-out-of-3 property.
More generally, any reasonable notion of ``equivalence'' in a category
should satisfy this property.
\end{paragraph}

\begin{paragraph}[Retracts]
  \index{retract}
  Let $\C$ be a category. We say that a morphism $f : X \to Y$ of~$\C$ is a
  \ndef{retract} of a morphism $g : Z \to T$ of $\C$ if there exists a
  commutative diagram
  \[
    \xymatrix{
      X \ar[d]_j \ar[r] \ar@/^3ex/[rr]^{\unit{X}} & Z \ar[d]_g \ar[r] & X \ar[d]^f \\
      Y \ar[r] \ar@/_3ex/[rr]_{\unit{Y}} & T \ar[r] & Y
    }
  \]
  in $\C$. We say that a class of morphisms of $\C$ is \ndef{closed under
  retracts} if any retract of an element of the class belongs to the class.
\end{paragraph}

\begin{paragr}[Lifting properties]
  \index{lifting property}
  \index{lift}
  Let $f : X \to Y$ and $g : Z \to T$ be two morphisms of $\C$. One says
  that $f$ \ndef{has the left lifting property} with respect to $g$ or that
  $g$ \ndef{has the right lifting property} with respect to $f$ if for every
  commutative square
  \[
  \xymatrix{
    X \ar[d]_f \ar[r] & Z \ar[d]^g \\
    Y \ar[r] & T \\
  }
  \]
  there exists a \ndef{lift}, that is, a morphism $h : Y \to Z$ making
  the two triangles
  \[
  \xymatrix{
    X \ar[d]_f \ar[r] & Z \ar[d]^g \\
    Y \ar[r] \ar@{.>}[ur]^h & T \\
  }
  \]
  commute. More generally, one says that $f$ has the left lifting property
  with respect to a class of maps $\clI$ if it has the left lifting
  property with respect to every morphism in $\clI$, and similarly for the
  right lifting property.
  We will denote by $\lorth{\clI}$ and
  $\rorth{\clI}$ the class of maps having the left or right lifting property
  with respect to a class $\clI$.
\end{paragr}

\interskip

We can now give the definition of a model category:

\begin{paragraph}[Model category]
  \index{model category}
  \index{fibration}
  \index{cofibration}
  \index{weak equivalence}
  \index{trivial fibration}
  \index{fibration!trivial}
  \index{trivial cofibration}
  \index{cofibration!trivial}
  A \ndef{model category} is a category $\M$ endowed with three classes of
  maps: the \ndef{weak equivalences}, the \ndef{cofibrations}, and the
  \ndef{fibrations}; these data are required to satisfy the following
  axioms:
  \begin{enumerate}
    \item the category $\M$ is finitely complete and finitely cocomplete,
    \item the class of weak equivalences satisfies the 2-out-of-3 property,
    \item the class of weak equivalences, cofibrations, and fibrations are
      closed under retracts,
    \item cofibrations have the left lifting property with respect to
      \ndef{trivial fibrations} (that is, maps that are both a fibration and
      a weak equivalence); trivial cofibrations (that is, maps that are both
      a cofibration and a weak equivalence) have the left lifting property
      with respect to fibrations, and
    \item every map of $\M$ factors as a cofibration followed by a trivial
      fibration, and as a trivial cofibration followed by a fibration.
  \end{enumerate}
\end{paragraph}

\begin{example}
  \nomenclature[Top]{$\Top$}{category of topological spaces}
  \label{ex:model_Top}
  One of the motivating examples of Quillen is the following model category
  structure on the category $\Top$ of topological spaces:
  \begin{itemize}
  \item the weak equivalences are the \ndef{topological weak equivalences},
  that is, the maps $f:X\to Y$ which induce a bijection
    $\pi_0(f):\pi_0(X)\to\pi_0(Y)$ on path components and isomorphisms
    $\pi_n(f,x):\pi_n(X,x)\to\pi_n(Y,f(x))$ on homotopy groups for every
    $n \ge 1$ and every base point $x$ in $X$,
  \item the fibrations are the \ndef{Serre fibrations}, that is, the maps
  having the right lifting properties with respect to the inclusions
  \[
    \begin{split}
      D^n & \hookto D^n\times I \\
      x & \mapsto (x, 0)
    \end{split}
  \]
  of disks into cylinders, for $n \ge 1$,
  \item the cofibrations are the maps having the left lifting property with
  respect to maps that are both topological weak equivalences and Serre
  fibrations.
  \end{itemize}
\end{example}

\begin{example}\label{ex:model_Ch}
  Another motivating example of Quillen is the following model category
  structure on the category $\pCh[\Z]$ of chain complexes
  (see~\cref{paragr:def_ChR}):
\begin{itemize}
  \item the weak equivalences are the quasi-isomorphisms
  (see~\cref{paragr:def_qis}),
  \item the cofibrations are the monomorphisms $f$ such that, for every $n
  \ge 0$, the cokernel of $f_n$ is projective,
  \item the fibrations are the morphisms $f$ such that, for every $n > 0$,
  $f_n$ is an epimorphism.
  \end{itemize}
\end{example}

\begin{example}\label{ex:model_Cat}
  The category $\Cat$ of small categories can be endowed with the so-called
  ``folk'' model category structure:
  \begin{itemize}
    \item the weak equivalences are the equivalences of categories,
    \item the cofibrations are the functors injective on objects,
    \item the fibrations are the \ndef{iso-fibrations}, that is, the
    functors $f : C \to D$ such that for every object $x$ of $C$ and any
    isomorphism $v : f(x) \to y$ of $D$, there exists an isomorphism $u : x
    \to x'$ in $C$ such that $f(u) = v$.
  \end{itemize}
\end{example}

\begin{example}\label{ex:model_wCat}
  \cref{chap:resolutions,chap:w-eq,chap:folk} are devoted to the
  construction of the ``folk'' model category structure on the category
  $\ooCat$ of strict \oo-categories, see~\cref{thm:folk}.
\end{example}

\begin{paragraph}[Cofibrant and fibrant replacements]
\label{sec:cof-repl}
\index{cofibrant!replacement}
\index{fibrant!object}
Let $\M$ be a model category. An object $X$ of~$\M$ is said to be
\ndef{cofibrant} if the unique morphism $\varnothing \to X$ from the initial
object of~$\M$ to $X$ is a cofibration. Dually, one says that the object~$X$
is \ndef{fibrant} if the unique morphism  $X \to \ast$ from $X$ to the
terminal object of $\M$ is a fibration.

If $X$ is an object of $\M$, a \ndef{cofibrant replacement} of $X$ is a
cofibrant object~$QX$ of $\M$ endowed with a weak equivalence $QX \to X$. It
follows immediately from the axioms of model categories that cofibrant
replacements exist. Indeed, to produce one, it suffices to factor
the morphism $\varnothing \to X$ as a cofibration followed by a trivial
fibration. Dually, a fibrant replacement consists of a fibrant object~$RX$
together with a weak equivalence $X \to RX$.
\end{paragraph}

\begin{paragraph}[Combinatorial model categories]
  \index{combinatorial model category}
  \index{model category!combinatorial}
  \index{cofibration!generating}
  A model category $\M$ is said to be \ndef{combinatorial} if
  \begin{enumerate}
    \item the underlying category of $\M$ is locally presentable
    (see~\cref{paragr:def_lpc}),
    \item there exists \emph{sets} $I$ and $J$ of morphisms of $\M$ such that
    the class of cofibrations of $\M$ is $\lrorth{I}$ and the class of
    trivial cofibrations of $\M$ is $\lrorth{J}$.
  \end{enumerate}
  Sets $I$ and $J$, as in the definition, are called \ndef{generating
  cofibrations} and \ndef{generating trivial cofibrations}, respectively.
\end{paragraph}

\begin{remark}
  To build a model category, one often starts with a class of weak
  equivalences and sets $I$ and $J$ of candidates to be generators. This is
  what we did in~\cref{chap:resolutions,chap:w-eq,chap:folk} to build the
  folk model structure on $\ooCat$.
\end{remark}

\begin{example}
  The model category structure on $\pCh[\Z]$ defined in \cref{ex:model_Ch}
  can be proven to be combinatorial.
\end{example}

\begin{example}
  The model category structure on $\Cat$ defined in \cref{ex:model_Cat} is
  combinatorial. The set $I$ can be taken to consist of the three functors
  \[ \varnothing \into \cdot, \qquad
    \{\xymatrix@C=0.8pc{\cdot & \cdot}\}
    \into \{\xymatrix@C=0.8pc{\cdot \ar[r] & \cdot}\},
    \quad
    \{
    \xymatrix@C=2pc{
      \cdot
      \ar@/^2ex/[r]_{}="0"
      \ar@/_2ex/[r]_{}="1"
      &
      \cdot
    }
    \}
    \to
    \{
    \xymatrix@C=1.5pc{
      \cdot
      \ar[r]
      &
      \cdot
    }
    \},
  \]
  and the set $J$ to the functor
  \[
    \{x\} \into \{\xymatrix@C=1pc{x \ar[r]^{\sim} & y}\},
  \]
  where the symbol $\sim$ denotes an isomorphism.
\end{example}

\begin{example}
  The folk model structure on $\ooCat$ is combinatorial. A set of
  generating cofibrations is given by the inclusions of spheres into
  disks (this is the set $\setgencof$ of~\cref{paragr:cof_triv_fib}). It is
  harder to describe a set of generating trivial cofibrations, see the set
  $\setgentrivcof$ of~\cref{paragr:fib}.
\end{example}
\begin{remark}
  The model category structure on $\Top$ described in \cref{ex:model_Top} is
  not combinatorial because $\Top$ is not locally presentable. Nevertheless,
  there exist sets $I$ and $J$ as in the second point of the definition.
  For instance, one can take as $I$ the set of canonical inclusions of
  spheres into disks, and as $J$ the set of inclusions of disks into
  cylinders as in the definition of Serre fibrations.
  These sets $I$ and $J$ do not satisfy the assumptions of the statement of
  the ``small object argument'' we gave (see~\cref{prop:small_obj_arg_2}).
  Nevertheless, one can check that the argument still applies. One says that
  this model category structure is \ndef{cofibrantly generated}.
\end{remark}

\section{The Homotopy Category}

\begin{paragraph}[Localization]
\label{sec:localization}
\index{localizer}
\index{localization}
\index{homotopy!category}
\index{category!homotopy}
  A \ndef{localizer} (also called \ndef{relative category}) is a
  category~$\C$ endowed with a class $\clW$ of morphisms called \ndef{weak
  equivalences}. The \ndef{homotopy category} of a localizer $(\C, \clW)$ is
  the category $\loc{\C}{\clW}$ obtained from~$\C$ by formally inverting
  arrows in $\clW$. More precisely, the category $\loc{\C}{\clW}$ is the
  universal category endowed with a functor $p : \C\to\loc\C\clW$ sending
  the elements of $\clW$ to isomorphisms. Universal means that for any
  category~$\D$ equipped with a functor $F:\C\to\D$ sending the elements
  of~$\clW$ to isomorphisms, there is a unique functor $\widetilde
  F:\loc\C\clW\to\D$ making the triangle
  \[
    \xymatrix@C=2pc@R=2pc{
      \C\ar[d]_p\ar[r]^F&\D\\
      \loc\C\clW\ar@{.>}[ur]_{\widetilde F}
    }
  \]
  commute. We will often denote the category $\loc{\C}{\clW}$ by $\Ho(\C)$, making
  implicit the class $\clW$.
\end{paragraph}

\begin{paragraph}
  If $(\C, \clW)$ is a localizer with $\C$ a small category, then the category
  $\loc\C\clW$ can be described in the following way: the objects of
  $\loc\C\clW$ are the object of~$\C$ and its morphisms are sequences of
  zigzags of morphisms of~$\C$
  \[
    \xymatrix@C=2pc@R=2pc{
      \ar[r]^-{f_1}&&\ar[l]_-{f_2}&\ar[l]_-{f_3}\ar[r]^-{f_4}&\cdots&\ar[l]_-{f_n}
      \pbox,
    }
  \]
  with all the backward morphisms in~$\clW$, modded out by the smallest
  congruence such that
  \begin{align*}
    \xymatrix@C=1pc{
      \ar[r]^-{f_1}&\ar[r]^-{f_2}&
    }
    &=
    \xymatrix@C=1pc{
      \ar[rr]^-{f_2\circ f_1}&&
    }
    &
    \xymatrix@C=1pc{
      X\ar[r]^-{\unit{X}}&X
    }
    &=\unit{X}
    &
    \xymatrix@C=1pc{
      X\ar[r]^-{f}&Y&\ar[l]_-{f}X
    }
    &=\unit{X}
    \\
    \xymatrix@C=1pc{
      &\ar[l]_-{f_1}&\ar[l]_-{f_2}
    }
    &=
    \xymatrix@C=1pc{
      &&\ar[ll]_-{f_1\circ f_2}
    }
    &
    \xymatrix@C=1pc{
      X&\ar[l]_-{\unit{X}}X
    }
    &=\unit{X}
    &
    \xymatrix@C=1pc{
      Y&\ar[l]_-{f}X\ar[r]^-{f}&Y
    }
    &=\unit{Y}
    \pbox.
  \end{align*}
  With appropriate set-theoretic foundations, one can adapt this construction
  to locally small categories, although the localization of a locally small
  category is not locally small in general.
\end{paragraph}

\begin{paragraph}[Homotopy category]
  Any model category $\M$ has an underlying localizer $(\M, \clW)$ and thus
  a \ndef{homotopy category} $\Ho(\M)$.
\end{paragraph}

\interskip

One of the goals of the theory of model categories is to get a good
understanding of the homotopy category $\Ho(\M)$. In particular, it
admits a simpler description in terms of homotopies.

\begin{paragraph}[Homotopies]
\index{cylinder!object}
\index{path object}
\index{homotopy}

\newcommand{\cylobj}[1]{I#1}
\newcommand{\pathobj}[1]{#1^I}

Let $\M$ be a model category.
A \emph{cylinder object} for an object~$A$ of $\M$ is a factorization
\[
  \xymatrix@C=2pc@R=2pc{
    A +  A\ar[dr]_i\ar[rr]^-{(\id_A,\id_A)}&&A\\
    &\cylobj{A}\ar[ur]_s
  }
\]
of the codiagonal map $A + A\to A$ as a cofibration $i$ followed by a
weak equivalence~$s$. The components of the map $i$ are denoted
$i_0,i_1:A\to \cylobj{A}$.
Similarly, a \emph{path object} for $A$ is a factorization
\[
  \xymatrix@C=2pc@R=2pc{
    A\ar[rr]^-{(\id_A,\id_A)}\ar[dr]_r&&A\times A\\
    &\pathobj{A}\ar[ur]_p&
  }
\]
of the diagonal map $A\to A\times A$ as a weak equivalence~$r$ followed by
a fibration~$p$. The components of the map~$p$ are denoted
$p_0,p_1:\pathobj{A}\to A$.

Given morphisms $f,g:A\to B$ of $\M$, a \emph{left homotopy} from~$f$ to~$g$ is a
morphism $h:\cylobj{A}\to B$, for some cylinder object $\cylobj{A}$ of $A$, such
that $hi_0=f$ and $hi_1=g$. If such a homotopy exists, we say that $f$ and
$g$ are \emph{left homotopic}.
Dually, a \emph{right homotopy} from~$f$ to~$g$ is
a morphism $k:A\to\pathobj{B}$, for some path object~$\pathobj{B}$ of~$B$, such
that $p_0k=f$ and $p_1k=g$. If such a homotopy exists, we say that $f$ and
$g$ are \emph{right homotopic}.
\end{paragraph}

\newcommand{\cfrestr}[1]{#1_{\textnormal{cf}}}

\begin{proposition}
  Let $\M$ be a model category. In the full subcategory~$\cfrestr\M$ of~$\M$
  whose objects are both cofibrant and fibrant,
  \begin{enumerate}
  \item the relation ``being left homotopic'' coincides with the relation 
  ``being right homotopic'',
  \item the relation ``being (left or right) homotopic'' is a congruence.
  \end{enumerate}
\end{proposition}

\begin{proposition}
  \label{prop:ho-heq}
  The homotopy category $\Ho(\M)$ of a model category $\M$ is equivalent to
  the category $\cfrestr\M/{\sim}$ obtained from $\cfrestr\M$ by quotienting
  the morphisms by the relation ``being (left or right) homotopic''.
\end{proposition}

In particular, the homotopy category of a model category is locally small.

\section{Derived Functors}

{
\renewcommand\N{\mathcal{N}}

\begin{paragraph}[Homotopical functors]
  \index{homotopical!functor}
  \index{functor!homotopical}
\label{sec:homotopical-functor}
Let $(\C, \clW_\C)$ and $(\D, \clW_\D)$ be two localizers.
A functor $F:\C\to\D$ is said to be \ndef{homotopical} if it preserves weak
equivalences, that is, if it sends the weak equivalences of $\C$ to weak
equivalences of $\D$. In this case, by the universal property of the
localization, the functor $F$ induces a functor
\[
  \overline{F}
  :
  \Ho(\C)
  \to
  \Ho(\D)
\]
making the square
\[
  \xymatrix@C=2.5pc@R=2pc{
    \C\ar[d]_{p_\C}\ar[r]^F&\D\ar[d]^{p_\D}\\
    \Ho(\C)
    \ar@{.>}[r]_{\overline{F}}&
    \Ho(\D)
  }
\]
commute.
\end{paragraph}

\interskip

If $F : \C \to \D$ is not homotopical, there is in general no functor
$\overline{F}$ making the above square commute. Nevertheless, one can seek
for ``best approximations'' to this situation:

\begin{paragraph}[Derived functors]
\label{sec:derived-functor}
\index{derived functor}
\index{functor!derived}
  Let $(\C, \clW_\C)$ and $(\D, \clW_\D)$ be two localizers and let $F : \C
  \to \D$ be a functor. The \ndef{(total) left derived functor} of $F$, if
  it exists, is the universal pair consisting of a functor
  \[ \Lder F : \Ho(\C) \to \Ho(\D) \]
  and a natural transformation
  \[
  \xymatrix@C=2pc@R=2pc{
    \C\ar[d]_{p_\C}\ar[r]^F&\D\ar[d]^{p_\D}\\
    \Ho(\C)
    \ar@{.>}[r]_{\Lder{F}}
    \ar@{}[ur]_(0.35){}="a"_(0.65){}="b"
    \ar@2"a";"b"^{\lambda}
    &
    \Ho(\D)
    \pbox.
  }
  \]
  This means that for every other functor $G : \Ho(\C) \to \Ho(\D)$ and
  natural transformation
  \[
    \xymatrix@C=2pc@R=2pc{
      \C\ar[d]_{p_\C}\ar[r]^F&\D\ar[d]^{p_\D}\\
      \Ho(\C)
      \ar@{.>}[r]_G
      \ar@{}[ur]_(0.35){}="a"_(0.65){}="b"
      \ar@2"a";"b"^{\alpha}
      &
      \Ho(\D)
      \pbox,
    }
  \]
  there exists a unique natural transformation $\gamma : G \tod \Lder F$ such
  that $\alpha$ factors as $\alpha = \lambda \circ (\gamma \ast p_\C)$.
  By abuse of language, one often refers to $\Lder F$ as the left derived
  functor of~$F$.

  Similarly, the \ndef{(total) right derived functor} of $F$, if it exists, is the
  universal functor $\Rder F : \Ho(\C) \to \Ho(\D)$ endowed with a natural
  transformation
\[
  \xymatrix@C=2pc@R=2pc{
    \C\ar[d]_{p_\C}\ar[r]^F&\D\ar[d]^{p_\D}\\
    \Ho(\C)
    \ar@{.>}[r]_{\Rder{F}}
    \ar@{}[ur]_(0.35){}="a"_(0.65){}="b"
    \ar@2"b";"a"_{\rho}
    &
    \Ho(\D)
    \pbox.
  }
\]
\end{paragraph}

\begin{remark}
  If $F$ is homotopical, then $\overline{F}$ (endowed with the identity
  natural transformation) is both the left and the right derived functor of
  $F$.
\end{remark}

One important use of model categories is to provide tools to prove the
existence of derived functors and to compute them.

\begin{theorem}
  \label{thm:derived-replacement}
  Let $F : \M \to \N$ be a functor between model categories. Suppose that
  $F$ sends trivial cofibrations between cofibrant objects to weak
  equivalences. Then $F$ admits a left derived functor $\Lder F : \Ho(\M)
  \to \Ho(\N)$. Moreover, if $X$ is an object of $\M$, then $\Lder
  F(p_\M(X))$ is canonically isomorphic to $p_\N(F(QX))$, where $QX$ is a
  cofibrant replacement of $X$.

  Similarly, if $F$ sends trivial fibrations between fibrant objects to weak
  equivalences, then $F$ admits a right derived functor that can be
  computed using fibrant replacements.
\end{theorem}

\begin{paragraph}[Quillen pairs]
\index{Quillen!functor}
\index{Quillen!pair}
\index{Quillen!adjunction}
\index{functor!Quillen}
\index{adjunction!Quillen}
  Let $\M$ and $\N$ be two model categories and let
  \[
  F
  :
  \M
  \rightleftarrows
  \N
  :
  G
  \]
  be an adjunction. One says that $(F, G)$ is a \ndef{Quillen pair} or
  a \ndef{Quillen adjunction} if $F$ preserves cofibrations and trivial
  cofibrations. This is equivalent to asking that $G$ preserves fibrations
  and trivial fibrations. In this case, one also says that $F$ is a
  \ndef{left Quillen functor} and that $G$ is a \ndef{right Quillen
  functor}. Using the previous proposition, one can show that left Quillen
  functors admit left derived functors and right Quillen functors admit
  right derived functors.
\end{paragraph}

\begin{theorem}
  If
  \[
  F
  :
  \M
  \rightleftarrows
  \N
  :
  G
  \]
  is a Quillen pair, then
  \[
  \Lder F
  :
  \Ho(\M)
  \rightleftarrows
  \Ho(\N)
  :
  \Rder G
  \]
  is an adjunction.
\end{theorem}

\begin{paragraph}[Quillen equivalences]
  \index{Quillen!equivalence}
  \index{equivalence!Quillen}
Let
\[
  F
  :
  \M
  \rightleftarrows
  \N
  :
  G
\]
be a Quillen pair. One says that $(F, G)$ is a \ndef{Quillen equivalence} if
the adjunction
\[
  \Lder F
  :
  \Ho(\C)
  \rightleftarrows
  \Ho(\D)
  :
  \Rder G
\]
is an adjoint equivalence.
\end{paragraph}

\begin{proposition}
  Let
  \[
    F
    :
    \M
    \rightleftarrows
    \N
    :
    G
  \]
  be a Quillen pair. The following assertions are equivalent:
  \begin{enumerate}
    \item $(F, G)$ is a Quillen equivalence,
    \item for every cofibrant object $X$ of $\M$ and every fibrant object
    $Y$ of $\N$, a morphism $FX\to Y$ in~$\N$ is a weak equivalence if and
    only if the adjoint morphism $X\to GY$ in~$\M$ is a weak equivalence.
  \end{enumerate}
\end{proposition}
}

\bibliographystyle{plain}
\bibliography{papers}

\renewcommand{\nomname}{Index of notations}
\printnomenclature
\renewcommand{\indexname}{Index of terminology}
\printindex


\end{document}